\documentclass[11pt,oneside,letterpaper]{book}
\usepackage{amsfonts}
\usepackage{amssymb}
\usepackage{amsmath,amsthm,mathrsfs,bigints}
\usepackage{graphicx}
\usepackage{verbatim}
\usepackage{color}
\usepackage[pdfencoding=auto,psdextra,pdfpagelabels]{hyperref}
\usepackage{stmaryrd}
\usepackage{mathtools}
\usepackage{bbm}
\usepackage{enumerate}
\usepackage[mono=false]{libertine}
\usepackage[textsize=footnotesize]{todonotes}
\usepackage[russian,english]{babel}
\usepackage[OT2,T1]{fontenc}
\usepackage{placeins}

\DeclareSymbolFontAlphabet{\amsmathbb}{AMSb}
\DeclareSymbolFont{sfletters}{OML}{cmbrm}{m}{it}

\makeatletter
\renewenvironment{proof}[1][\proofname]{\par
 \pushQED{\qed}%
 \normalfont \topsep6\p@\@plus6\p@\relax
 \trivlist
 \item[\hskip\labelsep
    \scshape
  #1\@addpunct{.}]\ignorespaces
}{%
 \popQED\endtrivlist\@endpefalse
}
\makeatother

\usepackage[margin=1.05in]{geometry}

\newcommand{\N}{\mathcal N}

\newcommand{\err}{\mathbbm{e}}

\newcommand{\dd}{\mathrm{d}}
\newcommand{\E}{\amsmathbb E}
\newcommand{\Res}{\mathop{\,\rm Res\,}}
\newcommand{\W}{\amsmathbb W}

\newcommand{\bth}{\boldsymbol\theta}

\newcommand{\Qu}{\mathbf{B}}

\newcommand{\D}{\mathfrak{D}^2}

\newcommand{\Dr}{\mathfrak{D}}

\newcommand{\B}{\mathcal B}

\newcommand{\z}{\mathbf{z}}

\renewcommand{\P}{{\amsmathbb P_\N}}

\newcommand{\Z}{\mathscr{Z}}

\newcommand{\Gm}{\mathcal G}

\newcommand{\F}{\mathscr{F}}

\newcommand{\ii}{\mathbf{i}}

\newcommand{\x}{\boldsymbol x}
\newcommand{\y}{\boldsymbol y}

\newcommand{\eps}{\varepsilon}

\renewcommand{\L}{\mathfrak L}
\newcounter{lassumc}

\newcounter{lassumcprime}

\newtheorem{theorem}{Theorem}[chapter]
\newtheorem{proposition}[theorem]{Proposition}
\newtheorem{lemma}[theorem]{Lemma}
\newtheorem{corollary}[theorem]{Corollary}

\newtheorem{lassum}[lassumc]{Assumption}
\newtheorem{lassumprime}[lassumcprime]{Assumption}
\newtheorem{assumption}{Assumption}

\newcounter{thmintroc}

\newtheorem{theoremI}[thmintroc]{Theorem}

\newcounter{lassumtiling}

\newtheorem{assumptionT}[lassumtiling]{Assumption}

\newtheorem*{assumption*}{Assumption}

\newtheorem{conjecture}[theorem]{Conjecture}
\newtheorem{definition}[theorem]{Definition}

\theoremstyle{remark}
\newtheorem{remark}[theorem]{Remark}

\DeclareMathSymbol{\sbeta}{\mathord}{sfletters}{"0C}

\renewcommand{\i}{\mathbf{i}}

\newcommand{\I}{\mathcal{I}}

\newcommand{\Op}{\Upsilon}
\def\ra{{\rightarrow}}

\numberwithin{equation}{chapter} \numberwithin{section}{chapter}

  \title{{\LARGE \textbf{Macroscopic asymptotics in discrete beta-ensembles and random tilings}}}

\author{Ga{\"e}tan Borot\footnote{Institut f\"ur Mathematik und Institut f\"ur Physik, Humboldt-Universit\"at zu Berlin, Unter den Linden 6, 10099 Berlin, Germany | Email: \href{mailto:gaetan.borot@hu-berlin.de}{\texttt{gaetan.borot@hu-berlin.de}}} \and Vadim Gorin\footnote{Departments of Statistics and Mathematics, University of California, Berkeley, CA, United States | Email: \href{mailto:vadicgor@gmail.com}{\texttt{vadicgor@gmail.com}}} \and Alice Guionnet\footnote{Unité de Mathématiques Pures et Appliquées (UMPA), École Normale Supérieure de Lyon, 46 allée d’Italie, 69007 Lyon, France | Email: \href{mailto:alice.guionnet@ens-lyon.fr}{\texttt{alice.guionnet@ens-lyon.fr}}}}

\begin{document}

\frontmatter

\hypersetup{pageanchor=false}
\begin{titlepage}
\maketitle
\end{titlepage}

\chapter*{Abstract}

We carry out the asymptotic analysis of repulsive ensembles of $N$ particles which are discrete analogues of continuous 1d log-gases or $\sbeta$-ensembles of random matrix theory. The ensembles that we study have several groups of particles which can have different intensities of repulsion. They appear naturally in models of random domino and lozenge tilings, random partitions, $\mathscr{N} = 2$-supersymmetric gauge theory, asymptotic representation theory, discrete orthogonal polynomial ensembles, etc. We allow filling fractions (\textit{i.e.} particle counts in groups) to be either fixed, or free, or to vary while respecting affine constraints. We are
interested in the macroscopic behavior of the distribution of particles, captured by linear statistics, partition functions, and their finite-size corrections as $N \rightarrow \infty$.

We prove the law of large numbers and large deviations for the empirical measure around the equilibrium measure. To reach finite-size corrections we work under an off-criticality assumption. For fixed filling fractions, we prove an asymptotic expansion for the  partition function and for the cumulants of linear statistics, in particular establishing a central limit theorem. For varying filling fractions, we prove that the central limit theorem is perturbed by an additional discrete Gaussian component which oscillates with $N$. Our approach is based on developing the bootstrap strategy for Nekrasov equations, which replace Dyson--Schwinger equations for discrete ensembles. The results can be seen as a far-reaching common generalisation of \cite{BGG} (for discrete ensembles with a single group of particles in the one-band regime) and \cite{BG_multicut} (for continuous $\sbeta$-ensembles in the multi-band regime where the CLT + discrete Gaussian was identified), but there are many new ingredients and difficulties; many of them are caused by the saturation phenomenon, which was absent in the continuous ensembles.

We apply our general results to the study of uniformly random lozenge tilings on a large class of domains --- not necessarily planar, simply-connected, nor orientable. The domains we allow are obtained by gluing trapezoids along a common vertical line: the tiles along the vertical are then described by a discrete ensemble we can study. When the analogues of filling fractions are fixed and this domain is orientable, we show that the Gaussian fluctuations on the vertical extend to the whole liquid region and are governed there by the Gaussian free field, as predicted by the Kenyon--Okounkov conjecture. We also establish a modification of the Kenyon--Okounkov conjecture in the non-orientable case. Complementarily, we prove discrete Gaussian fluctuations for filling fractions, when they are not fixed.

A key ingredient of the asymptotic analysis is the continuous invertibility of an operator arising through the linearization of Nekrasov equations, which we establish by generalizing ideas from \cite{BGK}. We relate this operator to a vector-valued Riemann-Hilbert problem (RHP) and develop systematically the algebraic and geometric tools to solve it. We exploit the RHP to construct spectral curves and analyse the structure of its solutions even when explicit formulae are lacking. In particular, it gives ways to relate the leading covariance of (continuous or) discrete ensembles to Green functions and, in the case of random lozenge tilings, leads to the identification of the Gaussian free field behavior.

\clearpage

\mainmatter
\setcounter{tocdepth}{2}
\tableofcontents

\hypersetup{pageanchor=true}

\vfill
\subsubsection*{Acknowledgements}

In the period 2016-2020 G.B. benefited from the excellent work conditions at the Max-Planck-Institute for Mathematics, Bonn, and from the support of the Max-Planck-Society. He also thanks the MIT, ENS Lyon, Melbourne University and IH\'ES for their hospitality during visits where this work progressed, Bernhard M. M\"uhlherr for helpful remarks about reflection groups, and Tibor Szab\'o for the idea of the proof of Proposition~\ref{Prop:realiza}. A.G was partially supported by the ERC Project LDRAM : ERC-2019-ADG Project 884584. V.G. was partially supported by the NSF grants DMS-1407562, DMS-1664619, DMS-1855458 and DMS-2246449, and thanks Alexei Borodin for discussions. We also thank Leonid Petrov for his help with all the beautiful simulations of random tilings included in the text.

\chapter{Introduction}
\label{SIntro}


This book develops the tools necessary for the asymptotic analysis of a class of statistical physics ensembles of $N$ particles occupying discrete sites on the real line and subjected to a pairwise interaction resembling the two-dimensional Coulomb repulsion. We start with a gentle introduction to the problem in several layers: first with a general context (Section~\ref{Chap1Sec1}), second with some illustrations of results that could not be reached before (Section~\ref{Chap1Sec2}), third with a description of the overall strategy of proofs (Section~\ref{Chap1Sec3}). The organization of the book is summarised in Section~\ref{Chap1Sec4}.

\section{Background, motivations, and objectives}
\label{Chap1Sec1}

\subsection{Paradigm and continuous ensembles}

\label{Section_continuous}

In probability theory, the law of large numbers states that the average of the observations obtained from a large number of independent random samples converges to the true mean, if it exists. It was first derived by Bernoulli in the 17th century to study the average outcome of a large number of fair coin tosses. Understanding the law of large numbers or the typical behavior of a large number of random variables is since then an essential question in probability theory. The analysis of the error with respect to this limiting behavior is the next natural question. The study of the fluctuations of large systems goes back to the 18th century and the famous work of de Moivre \cite{demoivre} who was the first to establish a central limit theorem by showing that the sum of independent Bernoulli variables, once recentered with respect to their mean and renormalized, converges towards a Gaussian variable. In the early 19th century Laplace \cite{Laplace_book} extended this analysis and made the first steps towards universality by proving that the convergence of these fluctuations towards a Gaussian variable is much more general, in the sense that the random variables need not be Bernoulli variables. The next wave in the early 20th century created the modern central limit theorems with essentially optimal conditions named after Lindeberg and Lyapunov.

 A large part of the research in probability and statistics revolves around central limit theorems. The original theorem required the random variables to be independent, but it was later extended to other situations, for example mean-field ensembles where random variables interact via a weak potential. This includes ensembles with joint law of the form
\[
\dd \amsmathbb{P}_{N}(x_{1},\ldots,x_{N})=\frac{1}{\mathscr{Z}_{N}}\cdot e^{-NV(\frac{1}{N}\sum_{i=1}^{N }x_{i})}\cdot \prod_{i = 1}^{N} \dd \rho(x_{i})
\]
for some smooth potential $V$, some probability measure $\rho$, and a normalizing constant $\mathscr{Z}_{N}$ making $\amsmathbb{P}_{N}$ a probability measure. In this case, the law of large numbers can be established by proving a large deviation principle thanks to Cram\'er theorem and Varadhan lemma \cite{Dembo_Zeitouni}, and the fluctuations can then be derived by a precise Laplace method \cite{bolt1,bolt2}. In such ensembles, the interaction between variables remains somewhat weak, in the sense that their strength is comparable to the entropy.

 However, these classical approaches do not allow for the analysis of the fluctuations of strongly interacting variables, such as
 the eigenvalues of Gaussian matrices or particles encoding random tilings, which have required the development of new methods. Consider the simplest possible random matrix model, namely the Gaussian orthogonal ensemble ($\sbeta = 1$) or the Gaussian unitary ensemble $(\sbeta = 2)$. It consists of a random matrix $M_{N}$ of size $N\times N$, which is real symmetric ($\sbeta =1$) or complex Hermitian ($\sbeta = 2$) with centered Gaussian independent entries. If $\sbeta=1$ the entries have covariance $\frac{2}{N}$ on the diagonal and $\frac{1}{N}$ above the diagonal. If $\sbeta=2$ the entries are real on the diagonal with covariance $\frac{1}{N}$, and above the diagonal they have independent real and imaginary parts with covariance $\frac{1}{2N}$ each. The eigenvalues $(\lambda_{1},\ldots,\lambda_{N})$ of $M_{N}$ are $N$ real numbers that are non-linear functions of the entries, in particular they exhibit strong and complicated correlations. Yet, Wigner \cite{Wig58} showed that the empirical measure of these eigenvalues converges weakly towards the semi-circle law by estimating its moments $\sum_{i = 1}^N \lambda_{i}^{k}=\textnormal{Tr}\,M_{N}^k$. Bessis, Br\'ezin, Itzykson, Parisi and Zuber \cite{BIPZ,BIZ}  (see also \cite{Zvonkin}) related these moments with the enumeration of maps and derived the so-called topological expansion. There is an extensive literature devoted to the study of the fluctuations of the variables $(\textnormal{Tr}\,M_{N}^{k})_{k \geq 1}$, implying the convergence in the sense of moments towards a Gaussian vector even when the entries of $M_{N}$ are not Gaussian provided they admit sufficiently many finite moments, see \textit{e.g.} \cite{jonsson,AGZ,LP,BaiSil,PaSh}. A particularly striking feature of these results is that $\textnormal{Tr}\,M_{N}^{k}$ fluctuates very little: even though it is of order the dimension $N$, the difference with its mean is only of order $1$. These fluctuations are much smaller than the fluctuations for the sum of $N$ independent random variables, which are of order $\sqrt{N}$.

 The situation becomes more involved when one considers random matrices with joint law
\begin{equation}
\label{meanfield}\dd\amsmathbb{P}(M_{N})= \frac{1}{\mathscr{Z}_N} \cdot e^{-\frac{N\sbeta}{2} \textnormal{Tr}\,V(M_{N})}\cdot \dd M_{N},
\end{equation}
where $\dd M_{N}$ denotes the Lebesgue measure over the space of real symmetric ($\sbeta = 1$) or Hermitian ($\sbeta = 2$) matrices. If $V(x)= \frac{x^{2}}{2}$ we recover the law of the Gaussian ensembles described above and the entries are independent. For any potential $V$ which is not a quadratic polynomial, the entries of $M_{N}$ are not independent anymore and computing the moments of $\textnormal{Tr}\,M_{N}^{k}$ becomes challenging. An alternative approach to these matrix models is based on the joint law of the eigenvalues $(\lambda_1,\ldots,\lambda_N)$ of $M_N$ which can be computed by a classical change of variables, see \textit{e.g.} \cite[Chapter 3]{ME} or \cite[Section 2.5]{AGZ}
\begin{equation}
\label{beta-mod}
\dd \amsmathbb{P}_{N}(\lambda_{1},\ldots,\lambda_{N})=\frac{1}{\mathscr{Z}_{N}}\cdot\prod_{1 \leq i < j \leq N} |\lambda_{i}-\lambda_{j}|^{\sbeta} \cdot \prod_{i = 1}^{N} e^{-\frac{N\sbeta}{2} \sum_{i=1}^{N} V(\lambda_{i})} \dd \lambda_{i}\,.
\end{equation}
For arbitrary values of $\sbeta > 0$ (not necessarily equal to $1$ or $2$) this is the definition of the \emph{continuous $\sbeta$-ensemble}, that is considered as a statistical mechanics ensemble of $N$ particles at positions $\lambda_1,\ldots,\lambda_N$. Even if the interaction between the eigenvalues is mean-field in the sense that it depends only on the empirical measure of the eigenvalues of the $\lambda_{i}$s, a crucial difference with ensembles like \eqref{meanfield} is that the interaction is much stronger: it is a product of $N^{2}$ terms of order $1$. Such models are also examples of log-gases in dimension one \cite{Forresterbook} and various generalizations of them (with different repulsion function or with $\lambda_i$s living in $\mathbb{R}^d$ with $d > 1$) have been studied from the perspective of rigorous statistical mechanics, see \textit{e.g.} the recent overview \cite{Serfatyreview}.

Sticking to \eqref{beta-mod}, the convergence of the empirical measure $L_N =\frac{1}{N}\sum_{i=1}^{N}\delta_{\lambda_{i}}$ for the ensembles \eqref{beta-mod} as $N\rightarrow\infty$ is by now well known. It can be established via large deviations techniques, see \cite{Vo93,BoPaSh,arous1997large} or \cite[Section 2.6]{AGZ}, or concentration of measures, see \cite{MaMa}. The outcome is that $L_N$ converges weakly in probability towards a measure $\mu$ defined as the unique maximizer of
the functional
\begin{equation}\label{grfcontinuous}
\mathcal{I}[\nu]:=  \frac{\sbeta}{2}\iint \bigg(\log|x-y|- \frac{V(x)+V(y)}{2}\bigg)\dd\nu(x)\dd\nu(y).
\end{equation}
This maximizer $\mu$ is usually called \emph{equilibrium measure} \cite{AlCol,Deiftcours}. This result is similar
 to the law of large numbers. The question of the global fluctuations of the empirical measure\footnote{Asymptotic fluctuations of \emph{individual eigenvalues} have also been a focus of active research, but we do not address this direction in this book.} around this limit was first tackled mathematically in \cite{Johansson,AlPaSh}, although some heuristics can be found earlier in physics, see \textit{e.g.} \cite{Politzer}. More specifically, Johansson \cite{Johansson} showed that if $V$ is a strictly convex polynomial and $f$ a sufficiently smooth test function $f$, the linear statistics
\begin{equation}
\label{fluctlin}
\Delta \textsf{Lin}[f] := \sum_{i = 1}^{N} f(\lambda_i) - N\int f(x)\dd\mu(x)
\end{equation}
converges towards a Gaussian variable, which can be treated as (a version of) the central limit theorem for continuous $\sbeta$-ensembles. This Gaussian variable has mean proportional to $(1 - \frac{\sbeta}{2})$, in particular this mean vanishes for $\sbeta = 2$. The $\sbeta=2$ version of the statement can be directly related to the asymptotic formula for the determinant of a Toeplitz matrix given by the strong Szeg\"o limit theorem, see \cite{Johansson} for discussion. The asymptotic covariance of $\Delta \textsf{Lin}[f]$ with varying $f$ turns out to be universal: assuming the support of $\mu$ to be a single segment and apart from a proportionality factor $\frac{2}{\sbeta}$ it depends only on the endpoints of the segment, as was first observed in physics \cite{BeZe,Ben2}. In addition \cite{borodin2014clt,borodin2015general} emphasized that the covariance structure can be identified with the one of  a one-dimensional section of the two-dimensional \emph{Gaussian free field}.

The asymptotic analysis of $\Delta \textsf{Lin}[f]$ was extended in \cite{BG11,Shc2} by removing convexity and polynomiality assumptions on the potential $V$ and postulating instead that it is analytic, the equilibrium measure has a connected support and is off-critical in the sense that its density of $\mu$ vanishes like a square-root at the boundary of the support. Note that the fluctuations of critical ensembles are still not fully understood despite recent attempts, see \cite{BLS} and compare with \cite{ClKu,ClKrIt}. Similarly, while the results for less regular than analytic potentials are present in the literature, central limit theorems for continuous $\sbeta$-ensembles under minimal assumptions on the potential $V$ were not achieved so far, see \cite{BLS,lambert2019quantitative,dadoun2023asymptotics} for the best available results.

In related developments, \cite[Proposition 1.2]{BG11} showed that the partition function (\textit{i.e.} the normalization constant $\mathscr{Z}_{N}$ in \eqref{beta-mod}) and the moments of linear statistics admit an all-order asymptotic expansion as $N \rightarrow \infty$. This means that there exists $N$-independent sequences of real numbers $(\mathscr{F}^{[p]})_{p\geq 0}$ so that for every $P \geq 0$
 \begin{equation}\label{topexp}
 \frac{1}{N^{2}}\log\left(\frac{\mathscr{Z}_{N}}{\mathscr{Z}_{N}^{\textnormal{G}\beta\textnormal{E}}}\right)= \sum_{p=0}^{P} N^{-p}\mathscr{F}^{[p]}+ o(N^{-P}), \qquad N\rightarrow\infty,
 \end{equation}
where $\mathscr{Z}_{N}^{\textnormal{G}\sbeta\textnormal{E}}$ is the partition function for the Gaussian $\sbeta$-ensemble obtained by choosing $V(x)= \frac{x^2}{2}$ in \eqref{beta-mod}. The latter is known explicitly by rescaling the Selberg--Mehta formula \cite{Selberg,forrester2008importance}
\begin{equation}
\label{eq_Selberg_integral}
 \int_{\amsmathbb{R}^N} \prod_{1\leq i<j \leq N} |\lambda_i-\lambda_j|^{\sbeta} \cdot \prod_{i=1}^N e^{-N\frac{\sbeta \lambda_i^2}{4}}\dd\lambda_i =(2\pi)^{\frac{N}{2}} \cdot \prod_{j=1}^N \frac{\Gamma\big(1+\frac{j\sbeta}{2}\big)}{\Gamma\big(1+\frac{\sbeta}{2}\big)} \cdot \left(\frac{2}{N\sbeta}\right)^{\frac{N}{2} + \frac{\sbeta N(N - 1)}{4}}.
\end{equation}
In addition, for any $n$-tuple $f_1,\ldots,f_n$ of smooth test functions of polynomial growth, there is a $N$-independent sequence $(\mathscr{F}_{n}^{[p]}[f_1,\ldots,f_n])_{p \geq 0}$ such that
 \begin{equation}
 \label{topexp2}
 \amsmathbb{E}\left[\prod_{j = 1}^{n} \Delta \textsf{Lin}[f_j] \right] = \sum_{p = 0}^{P} N^{-p} \mathscr{F}_{n}^{[p]}[f_1,\ldots,f_n] + o(N^{-P}).
 \end{equation}
Earlier mathematical results of this kind for $\sbeta = 2$ and with different assumptions appeared in \cite{AlPaSh} for the moments and \cite{ErMcL} for the partition function. Importing ideas from the theory of integrable systems, theoretical physicists have progressively \cite{Krichever,ACM92,ACKM,Ake96,E1MM,CE06,KostovCFT,EORev,EBook} discovered that the all-order expansion of continuous $\sbeta$-ensembles can be understood solely in terms of the geometry of a Riemann surface called spectral curve and encoding the equilibrium measure of the ensemble. In particular, the universal character of the leading covariance $\mathscr{F}_2^{[0]}[f_1,f_2]$ first noticed in \cite{BeZe,Ben2} for the one-band case is much more general, it relates to the Gaussian free field (\textit{aka} the free boson) on the spectral curve, and we will comment in Section~\ref{laius} on the way it appears to us.

Many of the results about central limit theorems and asymptotic expansions were obtained by exploiting the so-called Dyson--Schwinger equations. This is a hierarchy of equations involving the moments of the variables $
\textsf{Lin}[f] = \sum_{i = 1}^{N} f(\lambda_i)$ with $f$ varying in a certain set of suitable test functions, which is derived by integrations by parts in the integrals defining expectation values. Assuming that the equilibrium measure is off-critical, these equations can be asymptotically solved and lead to expansions in $N$ for these moments. Carrying out this programme requires inverting a master operator which is the linearization of the Dyson--Schwinger equations near the law of large numbers. We recall this approach in Section~\ref{Theorem_correlators_intro} and refer to \cite{ST} or \cite[Chapter 1]{BGK} for an overview of some mathematical uses of Dyson--Schwinger equations and discussion of the master operator and its inversion, and to \cite{pain,EvKn,gorin2024dynamical} for some more recent developments.

The case where the support of the equilibrium measure is not connected, but consists of $K \geq 2$ pairwise disjoint segments $\amsmathbb{B}_1,\ldots,\amsmathbb{B}_K$, leads to new phenomenology. Indeed, if $\chi_k$ is a smooth function close to the indicator function of $\amsmathbb{B}_k$ the \emph{filling fractions} defined as
\[
N_k: = \sum_{i = 1}^N \chi_k(\lambda_i)
\]
should be the number of particles in a neighborhood of $\amsmathbb{B}_k$, and therefore an integer. Thus, it cannot fluctuate like a Gaussian variable, contradicting the asymptotic Gaussianity of \eqref{fluctlin}. Instead, it was shown in \cite{Shc2,BG_multicut} that the fluctuations of the filling fractions are asymptotically described by discrete Gaussians, which are random vectors with integer coordinates whose weights resemble the Gaussian density. As the centering of the discrete Gaussian component is $N$-dependent, it can only admit subsequential limits, with different subsequences leading to substantially different asymptotic random variables. When $f$ is different from a smooth indicator function, the linear statistic \eqref{fluctlin} is described by the sum of a Gaussian and a discrete Gaussian. The origin of this \emph{perturbed central limit theorem} (or "central limit theorem with interferences'') can be traced back to the numerical observations of \cite{Jurkiewicz} and the heuristic derivation of \cite{BDE}. In probabilistic terms it is described by conditioning the system to have fixed filling fractions $\boldsymbol{N}= (N_1,\ldots,N_K)$ and decomposing \eqref{fluctlin} as
\begin{equation}
\label{eq_decompolin}
\Delta\textsf{Lin}[f] = \bigg(\sum_{i = 1}^{N} f(\lambda_i) - N \int f(x)\dd\mu^{\hat{\boldsymbol{n}}}(x)\bigg) + N\int f(x)\dd(\mu - \mu^{\hat{\boldsymbol{n}}})(x),
\end{equation}
where $\mu^{\hat{\boldsymbol{n}}}$ is the equilibrium measure corresponding to the fixed filling fractions $\hat{\boldsymbol{n}} = \frac{\boldsymbol{N}}{\N}$, and obtained as a minimizer of $-\mathcal{I}$ in \eqref{grfcontinuous} constrained to have those fixed masses in $\amsmathbb{B}_1,\ldots,\amsmathbb{B}_K$. The first term in \eqref{eq_decompolin} corresponds to the system with fixed filling fractions and becomes asymptotically Gaussian. The second term gives the contribution of the discrete component. The two terms \eqref{eq_decompolin} are asymptotically independent, because the mean and variance of the first term eventually only depend on the macroscopic properties of the equilibrium measure, in which the precise dependence on $\hat{\boldsymbol{n}}$ is washed away.

\subsection{Random tilings and discrete ensembles}

\begin{figure}[t]
\begin{center}
\includegraphics[height=8cm]{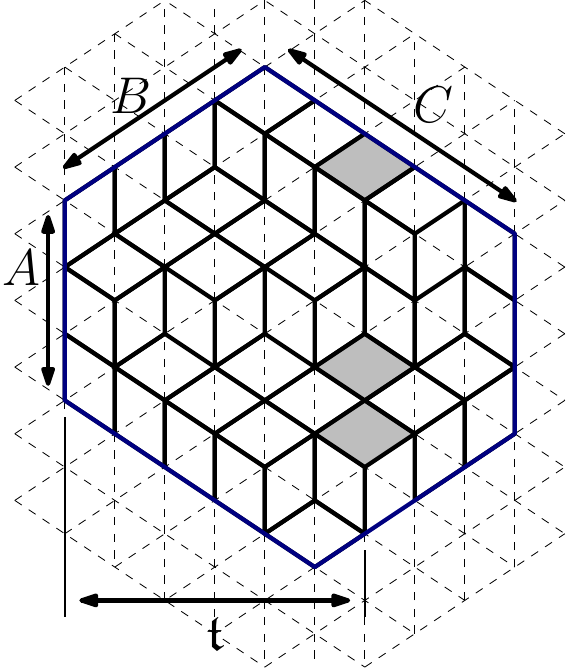}\hspace{2cm}
\includegraphics[height=8cm]{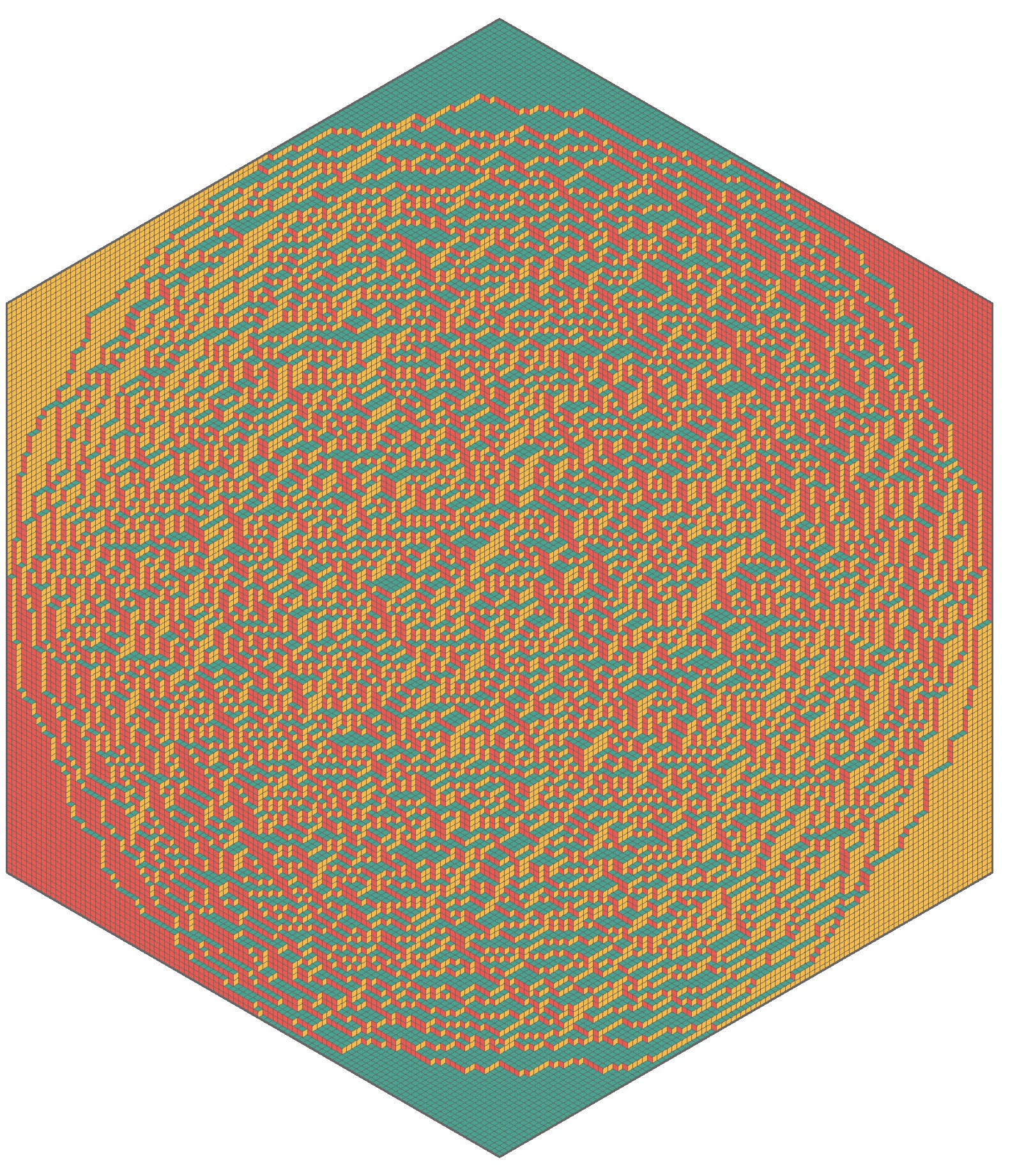}
\caption{\label{Fig_tiling_hex_222} Left: A lozenge tiling of $3\times 4 \times 5$ hexagon and its vertical section. Right: Simulation of a uniformly random lozenge tiling of $100\times 100 \times 100$ hexagon (three types of lozenges shown in colors)}
\end{center}
\end{figure}

Another well-known setting of random variables with strong interactions is given by random tilings --- a particular case of the dimer model, \textit{cf.} \cite{kenyon2009lectures}. The goal is to study the tiling, say by lozenges, of a given domain. When such a tiling is possible, it may be done in different ways and one considers a tiling picked uniformly at random within the set of possible ones. A lozenge tiling of a domain in the plane can be seen as the two-dimensional projection of a three-dimensional stepped surface, and therefore random tilings deliver one of the simplest models of random surfaces. One then would like to understand the typical shape of the stepped surface corresponding to a random tiling. This question has been most studied and best understood in the setting of a hexagonal domain. Consider an $A\times B\times C$ hexagon drawn on the regular triangular lattice. It can be tiled
using three types of elementary lozenges (three ways to glue together two adjacent lattice triangles), seen as two-dimensional projections of the three types of faces (green, yellow, red in the simulation of Figure~\ref{Fig_tiling_hex_222}) of the cube. We are interested in the asymptotic behavior of uniformly random tilings as $A,B,C \rightarrow \infty$ and would like to answer similar questions as before: Does the random tiling converge towards a deterministic limit shape? Can we describe the fluctuations around this limit shape? Is there an asymptotic expansion for the partition function? The same questions can be asked for lozenge tilings of more general domains.

For the hexagon, numerous results are already available and we refer to \cite{Vadimlecture} for a comprehensive overview. The computation of the partition function, \textit{i.e.} the number of tilings of a hexagon, dates back over more than a century and is given by the MacMahon formula \cite{macmahon1896,macmahon1915}, \textit{cf.} \eqref{eq_Selberg_integral}:
\begin{equation}
 \#\{\text{lozenge tilings of the }A\times B\times C\text{ hexagon}\}=\prod_{a=1}^A \prod_{b=1}^B \prod_{c=1}^C \frac{a+b+c-1}{a+b+c-2}.
\end{equation}
Much more recently \cite{CLP} proved that the stepped surface corresponding to random tiling converges (uniformly, in probability) towards a limit shape, which should be interpreted as the law of large numbers in this setting. Near the boundary of the domain all the tiles take the same orientation resulting in \emph{frozen regions} (these are flat regions for the stepped surface), whereas in the center of the domain the typical configuration shows a mixture of tiles in different directions resulting in \emph{liquid} regions (that are rough regions in the stepped surface). The boundary between these two types of regions is called the \emph{arctic curve} (Figure~\ref{Fig_tiling_hex_222}). The macroscopic fluctuations of the stepped surface inside the arctic curve are given by the two-dimensional Gaussian free field \cite{kenyon2008height,petrov2015asymptotics,BuGo2,duits2018global,berggren2024perfect}. The Gaussian free field is a two-dimensional Gaussian field whose covariance is given by the Green function of the Laplacian in the domain under consideration, with Dirichlet boundary conditions, and under a specific complex structure constructed from the limit shape, see \cite{sheffield2007gaussian,werner2020lecture}, and \cite[Lectures 11-12]{Vadimlecture}. We will review its complex-analytic aspects in Section~\ref{sec:complexstr}, see Definition~\ref{GFFdefinition}. The microscopic fluctuations of the tiles have also been studied and feature universal limit laws similar to those appearing in random matrices, but such questions are not in our focus.

\noindent

In this book we adopt the point of view that random tilings are discrete versions of random matrices. From the representation-theoretic perspective, the two are related by the semiclassical limit degenerating irreducible representations of a Lie group (tilings are in bijection with Gelfand--Tsetlin patterns enumerating bases of irreducible representations of $\textnormal{U}(N)$, see \cite[Section 2]{BP_lectures}) into orbital measures on the corresponding Lie algebra (the Lie algebra of $\textnormal{U}(N)$ is, up to multiplication by $\ii$, the space of Hermitian matrices), \textit{cf.} \cite[Section 19.3]{Vadimlecture}. From the perspective of integrable probability and special functions, tilings and matrices are described by very similar structures, such as determinantal point processes, Painlev\'e equations, Sine and Airy kernels, see \textit{e.g.} \cite{johansson2006random,baik2016combinatorics}.

To better understand the parallelism, first, observe that a lozenge tiling is uniquely determined by the positions of one type of lozenges, for instance the horizontal ones. It happens that the distribution of these horizontal tiles crossing a given vertical axis is very similar to the one of the $\sbeta$-ensemble \eqref{beta-mod} for $\sbeta = 2$, except for the fact that the positions of the lozenges take discrete values. More precisely, if we dissect the hexagon by a vertical line at distance $\mathfrak{t}$ from the left boundary, the total number $N$ of horizontal lozenges on this line is fixed by the geometric parameters $A,B,C,\mathfrak{t}$ (Figure~\ref{Fig_tiling_hex_222}). The distribution of these $N$ lozenges was first computed in \cite{CLP}, see also Section~\ref{Section_Hexagon}. For instance, in the case $\mathfrak{t}>\max(B,C)$ as in Figure~\ref{Fig_tiling_hex_2}, we have $N=B+C-\mathfrak{t}$. Introducing the coordinate system such that the lowest possible position for
horizontal lozenges on the $\mathfrak{t}$-th vertical line is $1$ and the highest one is $A+B+C-\mathfrak{t}$,
the formula for the distribution of the positions $\boldsymbol{\ell} = (\ell_1 < \cdots < \ell_N) \in \amsmathbb{Z}^N$ of the horizontal tiles reads
\begin{equation}
\label{eq_Hahn}
\begin{split}
 & \amsmathbb{P}_N(\ell_1,\ldots,\ell_N) =\frac{1}{\mathscr{Z}_N} \cdot \prod_{1 \leq i<j \leq N} (\ell_j - \ell_i)^2 \cdot \prod_{i=1}^N w(\ell_i), \\
 & w(\ell)=
 (A+B+C+1-\mathfrak{t}-\ell)_{\mathfrak{t}-B} \cdot (\ell)_{\mathfrak{t}-C},
\end{split}
\end{equation}
where $(a)_n =a(a+1)\cdots (a+n-1)$ is the Pochhammer symbol and $\mathscr{Z}_N$ is the normalizing
constant making $\amsmathbb{P}_N$ a probability measure. The distribution
\eqref{eq_Hahn} is known as the \emph{Hahn discrete orthogonal polynomial ensemble} and this formula was a starting point for many asymptotic results, see \cite{Johansson2,BKMM,Gor,BG,duits2018global}.

The key feature for our discussion is the similarity between \eqref{eq_Hahn} and \eqref{beta-mod} with $\sbeta=2$, which becomes apparent if we expand the weight using Stirling formula:
\[
w(\ell) = e^{\varpi N \log N - N V(\frac{\ell}{N}) + o(\log N)}.
\]
Here $\varpi$ is a constant and the potential $V$ is expressed via $\mathrm{Llog}(x) := x\log |x| - x$:
\begin{equation}
\label{eq_logpot}V(x) = \mathrm{Llog}(x) + \mathrm{Llog}(\hat{A} + \hat{B} + \hat{C} - \hat{\mathfrak{t}} - x) - \mathrm{Llog}(x + \hat{\mathfrak{t}} - \hat{C}) - \mathrm{Llog}(\hat{A} + \hat{C} - x).
\end{equation}
Let us nevertheless stress two important distinctions between \eqref{eq_Hahn} and \eqref{beta-mod}: the former deals with \emph{discrete} $\ell_i\in\amsmathbb{Z}$, while the latter deals with \emph{continuous} $\lambda_i\in\amsmathbb{R}$; the potential $V(x)$ is smooth in \eqref{beta-mod}, while only $\exp(-N V(x))$ is smooth in \eqref{eq_Hahn} and $V'(x)$ has logarithmic singularities.

\subsection{Discrete analogues of \texorpdfstring{$\sbeta$}{beta}-ensembles: our objects and goals}

The main goal of this book is to study a common generalization of \eqref{beta-mod} and \eqref{eq_Hahn} presented in details in Chapter~\ref{Chapter_Setup_and_Examples}. It is a discrete and multi-group version of the continuous $\sbeta$-ensembles, which is a probability measure on discrete particle configurations of weight
\begin{equation}
\label{eq_measure_general_intro}
\amsmathbb{P}_\N(\boldsymbol{\ell}) = \frac{1}{\mathscr{Z}_N} \cdot \prod_{1 \leq i < j \leq N} \frac{\Gamma(\ell_j - \ell_i + 1)\cdot \Gamma(\ell_j - \ell_i + \theta_{h(i),h(j)})}{\Gamma(\ell_j - \ell_i) \cdot \Gamma(\ell_j - \ell_i + 1 - \theta_{h(i),h(j)})} \cdot \prod_{i = 1}^{N} w_{h(i)}(\ell_i),
\end{equation}
where $N$ particles $\ell_1<\ell_2<\ldots<\ell_N$ are split in a monotonic way into $H$ groups with $h(i)\in\{1,2,\ldots,H\}$ denoting the group of the $i$-th particle. The differences $\ell_{i + 1} - \ell_i - \theta_{h,h}$ are assumed to be non-negative integers whenever the $i$-th and the $(i + 1)$-th particle are in the same group $h$. The repulsion intensity $\sbeta$ in \eqref{beta-mod} is replaced with the $H\times H$ matrix $2\boldsymbol{\Theta}$. We allow quite general symmetric positive semi-definite  matrices $\boldsymbol{\Theta}=(\theta_{g,h})_{1 \leq g,h \leq H}$ with positive diagonal. The weight $w(\ell)$ in \eqref{eq_Hahn} is replaced with group-dependent weights $w_h(\ell)$ on which we impose similar regularity conditions as for $w(\ell)$.

Note that for $\theta \in \{\frac{1}{2},1\}$ the interactions in \eqref{eq_measure_general_intro} are exactly the same as in continuous $\sbeta$-ensembles:
\[
 \frac{\Gamma(\ell_j - \ell_i + 1)\cdot \Gamma(\ell_j - \ell_i + \theta)}{\Gamma(\ell_j - \ell_i) \cdot \Gamma(\ell_j - \ell_i + 1 - \theta)}=(\ell_j-\ell_i)^{2\theta}.
\]
For other values of $\theta$ the match only occurs asymptotically:
\[
\frac{\Gamma\big(\ell_j - \ell_i + 1\big)\cdot \Gamma\big(\ell_j - \ell_i + \theta\big)}{\Gamma\big(\ell_j-\ell_i\big) \cdot \Gamma\big(\ell_j - \ell_i + 1 - \theta\big)} \sim (\ell_j - \ell_i)^{2\theta}, \qquad \text{as } (\ell_j-\ell_i)\rightarrow +\infty.
\]

The particular case of \eqref{eq_measure_general_intro} with all $\theta_{g,h}$ equal to a single $\theta > 0$ was first analyzed\footnote{However, even for this case, we go much further than \cite{BGG}.} in \cite{BGG}. The measures of the kind \eqref{eq_measure_general_intro} are widespread in probability and beyond.
\begin{itemize}
 \item For $H=1$, $\theta_{1,1}=1$, and various choices of $w_1(\ell)$, they appear in the study of random domino tilings of the Aztec diamond \cite{Johansson2}, random lozenge tilings of hexagons as in the previous section, last passage percolation \cite{Johansson3}, interacting particle systems in the TASEP family \cite{Johansson_shape}, stochastic systems of non-intersecting paths \cite{Johansson2,KoOCRo}, partition functions of Yang--Mills theory \cite{DK,LeMa}, and the representation theory of $\textnormal{U}(\infty)$ \cite{BoOl2007,Borodin2011}.
 \item For $H=1$ and general $\theta = \theta_{1,1}>0$, they appear in the analysis of Jack polynomials, of the $zw$-measures in the asymptotic representation theory, and of the Jack-deformation of the Plancherel measures on Young diagrams \cite{Fulman,Okounkov2005,BoOl2005, Matsumoto2008,Matsumoto2011,Olsh_hyper,DoFe,dimitrov2024global}.
 \item For general $H$ and matrices $\boldsymbol{\Theta}$ with entries in $\{0,\tfrac{1}{2},1\}$ they can describe sections of lozenge tilings of domains with very general topologies, as we explain in Chapter~\ref{Chap11}.
\end{itemize}

In the next section we will expand on some of these examples. We also remark that our approach to the analysis of \eqref{eq_measure_general_intro} involves interpolation in the matrix $\boldsymbol{\Theta}$ towards model cases, and subsequently, even if one is only interested in a special $\boldsymbol{\Theta}$, introducing and using more generic $\boldsymbol{\Theta}$s becomes inevitable.

\medskip

The objective of this book is to demonstrate how the phenomenology of global fluctuations for continuous $\sbeta$-ensembles outlined in Section~\ref{Section_continuous} extends (or not) to the discrete ensembles \eqref{eq_measure_general_intro}. The main achievements are:

\begin{enumerate}
 \item We develop the theory of equilibrium measures (including regularity properties), the law of large numbers, and large deviation principles for the discrete ensembles \eqref{eq_measure_general_intro} in Chapters~\ref{Chapterlarge} and \ref{Chapter_smoothness}.
 \item We find  the discrete substitute of the Dyson--Schwinger equations for \eqref{eq_measure_general_intro} --- in the spirit of the non-perturbative Dyson--Schwinger equations of Nekrasov --- in Chapter~\ref{ChapterNekra}.
 \item We obtain central limit theorems for the linear statistics \eqref{fluctlin} in discrete ensembles \eqref{eq_measure_general_intro} with fixed filling fractions, leading to Gaussian fluctuations in Chapter~\ref{Chapter_fff_expansions}.
 \item We establish the existence of several terms in the  $\N \rightarrow \infty$ asymptotic expansions of the partition function for discrete ensembles with fixed filling fractions in Chapter~\ref{Chapter_partition_functions}.
 \item We prove the appearance of discrete Gaussians in the asymptotics of filling fractions (as well as their necessary generalizations) and obtain the perturbed central limit theorems combining Gaussian and discrete Gaussian components in Chapter~\ref{Chapter_filling_fractions}.
 \item We use the ensembles \eqref{eq_measure_general_intro} to analyze random lozenge tilings of a wide class of domains, resulting in the description of their macroscopic fluctuations via the Gaussian free field in Chapter~\ref{Chap11}.
  \item We set up the functional and complex analysis tools to invert the operator appearing in Dyson--Schwinger / Nekrasov equations in Chapter~\ref{Chapter_SolvingN}, to analyze the spectral curves and compute the leading covariance in Chapter~\ref{Chapter_AG}. This is largely independent of the specific discrete ensembles and has an interest as well for the well-studied continuous ensembles.
\end{enumerate}

\section{Illustration of results}

\label{Chap1Sec2}

In this section we provide four different stochastic systems leading to the particular cases of distribution \eqref{eq_measure_general_intro}. Our goal is two-fold. First, we highlight how the general form of \eqref{eq_measure_general_intro} arises as a simplest common generalization of all these examples. Second, we indicate the main theorems about these systems which are corollaries of the results proven in this book.

\subsection{Lozenge tilings of the hexagon with a hole}
\label{Section_tiling_hole}

\begin{figure}[t]
\begin{center}
\includegraphics[width=0.5\textwidth]{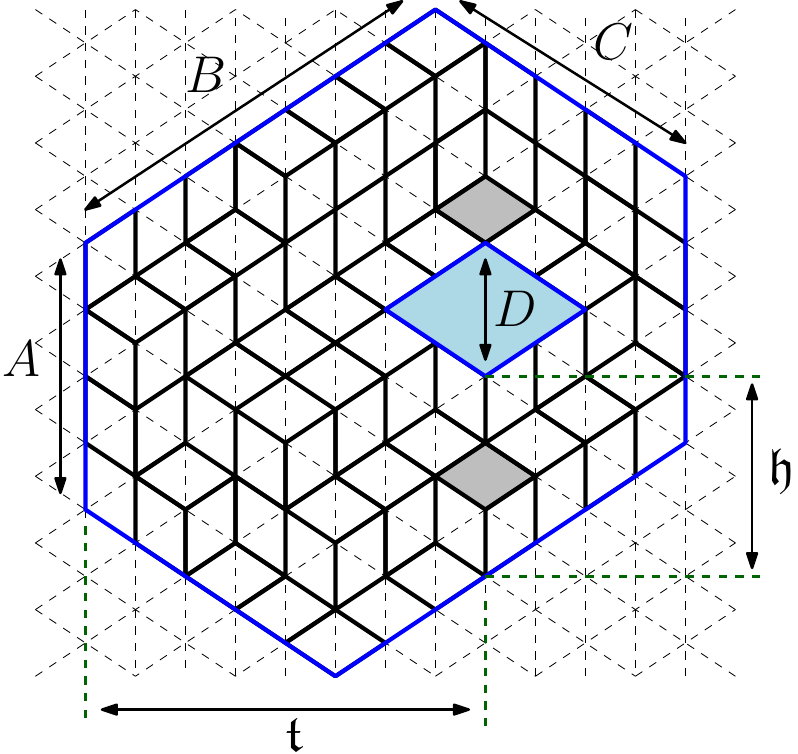} \hfill \includegraphics[width=0.43\linewidth]{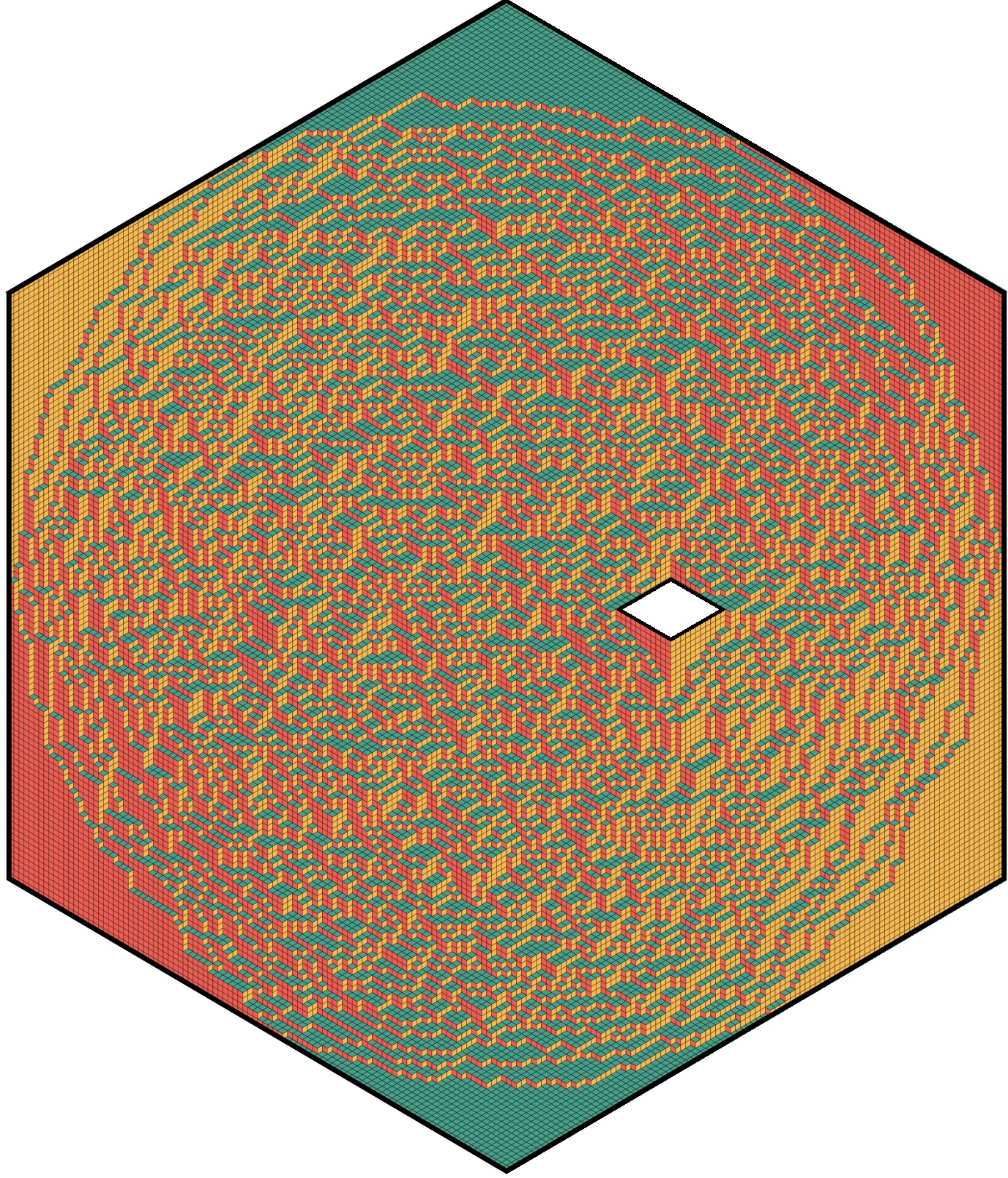}
\caption{\label{Fig_hex_hole_small} Left panel: Lozenge tiling of
the $4\times 7\times 5$ hexagon with a rhombic $2\times 2$ hole (shown in blue). The horizontal lozenges outside the hole and on the $\mathfrak{t}$-th vertical line are shown in
gray. Right panel: $100\times 100\times 100$ hexagon with a hole. We thank Leonid Petrov
for his help with this and other simulations.}
\end{center}
\end{figure}

Let us consider uniformly random lozenge tilings of the domain obtained by
cutting a
rhombic $D\times D$ hole in the hexagon, as shown in Figure~\ref{Fig_hex_hole_small}. This domain is not simply-connected and its random tilings were considered already in \cite{BGG,BuGo2}. We call $(\mathfrak{h},\mathfrak{t})$ the horizontal and vertical coordinates of the bottom of the hole compared to the left boundary and bottom point of the hexagon, and we are interested in the asymptotic behavior of tilings as $A,B,C,D,\mathfrak{t},\mathfrak{h}$ go to infinity at the same rate $\N$.

Let $\amsmathbb P_N$ be the probability distribution of the horizontal lozenges outside the
hole and on the $\mathfrak{t}$-th vertical line, induced by the uniform measure on all tilings of the
hexagon with the hole. Like in the hexagon, the number $N$ of horizontal lozenges on that line is fixed by the geometry. Focussing on the case $\mathfrak{t} > \max(B,C)$ as in Figure~\ref{Fig_hex_hole_small}, then the highest possible position is $A + B + C - \mathfrak{t}$, and the positions of the $N$ horizontal lozenges are specified by
\[
1 \leq \ell_1 < \cdots < \ell_N \leq A + B + C - \mathfrak{t}.
\]
It is explained in \cite[Section 9.2]{BGG}, see also Section~\ref{Section_hex_hole}, that the distribution of $\boldsymbol{\ell} = (\ell_i)_{i = 1}^{N}$ is
\begin{equation}
\label{eq_Hahn_cut}
 \amsmathbb P_N(\boldsymbol{\ell})=\frac{1}{\mathscr{Z}_N} \cdot \prod_{1\leq i<j\leq N} (\ell_j-\ell_i)^2 \cdot \prod_{i=1}^N w(\ell_i),
 \end{equation}
with the weight
\begin{equation}\label{weight:holey}w(\ell)=
 (A+B+C+1-\mathfrak{t}-\ell)_{\mathfrak{t}-B} \cdot (\ell)_{\mathfrak{t}-C}\cdot (\mathfrak{h} + 1 -\ell)_D^2.
 \end{equation}
One can split the particles in \eqref{eq_Hahn_cut} into $H = 2$ groups, depending on whether they are below or above the hole: the former means $\ell_i\in [a_1,b_1]= [1,\mathfrak{h}]$, the latter $\ell_i\in [a_2,b_2]=[ \mathfrak{h}+1 +D,A+B + C - \mathfrak{t}]$. Hence, we are in the setting of \eqref{eq_measure_general_intro} with $H=2$ and
\[
\boldsymbol{\Theta}= \left(\begin{array}{cc} 1& 1\\ 1 & 1 \end{array}\right).
\]

At this point, there are two possible situations to study. In \cite[Section 9.2]{BGG} and \cite{BuGo3}, we were \emph{fixing
the segment filling fractions}: we consider only the tilings with $N_1$ horizontal lozenges on the $\mathfrak{t}$-th vertical line below the hole and
$N_2=N-N_1$ of them above, where $N_1$ and $N_2$ are deterministic. In \cite{BGG} the fluctuations of linear statistics (as in \eqref{fluctlin}) encoding horizontal lozenges along the $\mathfrak{t}$-th vertical line were analyzed and shown to be Gaussian. In \cite{BuGo3}, this result was used to upgrade to the convergence of the full two-dimensional field of fluctuations of the stepped surface corresponding to the random tiling towards the Gaussian free field.

Another possibility is to leave $N_1$ and $N_2$ free: only their sum $N = N_1 + N_2$ is
deterministic while $N_1,N_2$ are random. In this situation we can ask a new question: what is the asymptotic
distribution of $N_1$ as the size of the hexagon and its hole become very large? In terms of the stepped surface, this is asking about the asymptotic distribution for the height of the hole. Chapter~\ref{Chap11} will answer this question with a discrete Gaussian distribution, see Corollary~\ref{Corollary_discrete_Gauss_in_hex} for which we present an informal version here.

\begin{theoremI} \label{Theorem_hex_hole_intro}
 Assume that $A$, $B$, $C$, $D$, $\mathfrak t$ grow linearly with a large parameter $\N$. Then, under a non-degeneracy condition (\textit{cf.} Definition~\ref{Definition_hex_hole_non_degenerate}), the random variable $N_1$ is asymptotically discrete Gaussian, in the sense that there exist three constants $m,r,\mathscr{Q}$ depending on $\frac{A}{\N},\frac{B}{\N},\frac{C}{\N},\frac{D}{\N},\frac{\mathfrak{t}}{\N},\frac{\mathfrak{h}}{\N}$ such that
\begin{equation}
\label{vareqna}
\lim_{\N\rightarrow\infty} \sup_{k\in \amsmathbb Z}\left|\amsmathbb{P}_{\N}[N_1=k_1] - \frac{1}{\mathscr{Q}} \exp\left( - \frac{(k_1- \N m)^2}{2r^2}\right)\right|=0.
\end{equation}
\end{theoremI}
In fact, we prove asymptotic Gaussianity and discrete Gaussianity in a much larger generality, for uniformly random lozenge tilings of many domains with various complicated topologies, see Chapter~\ref{Chap11} for further details.

\medskip

Theorem~\ref{Theorem_hex_hole_intro} is new. However, discrete Gaussian random variables were widely expected to appear in random tilings of non simply-connected domains, see \cite[Lecture 24]{Vadimlecture}. Some related rigorous results can be found in the literature. For dimers on Riemann surfaces (in the setting where no frozen regions appear), discrete Gaussians appear as increments of the height function along a loop of a non-trivial topology, see \cite{boutillier2009loop,Dubedator,dubedat2015asymptotics,kenyon2016asymptotics,berestycki2024dimers}. For domino tilings of domains with holes and special ``Temperleyan'' boundary conditions (which are not polygonal and again cannot lead to frozen regions), discrete Gaussian fluctuations were found in \cite{basok2023dimers,nicoletti2025temperleyan}. In the setting of periodically-weighted domino tilings of the Aztec diamond, the gaseous phases (which are significantly distinct from both liquid and frozen regions observed in uniformly random lozenge tilings) play the same role as holes and also lead to the appearance of discrete Gaussian components, see \cite{berggren2025gaussian}. Finally, \cite{basok2025kenyon} shows that if the limit of macroscopic fluctuations of tilings exists and satisfies certain very natural properties, then the discrete Gaussian component should necessarily be present. We believe that none of the methods of these articles could reach Theorem~\ref{Theorem_hex_hole_intro}.

\medskip

The non-degeneracy condition in Theorem~\ref{Theorem_hex_hole_intro} is explained in more detail in Section~\ref{sec:gen}. It is crucial for the applicability of our techniques, but it is less clear whether it is needed for the final conclusion. In particular, this condition includes an off-criticality assumption prohibiting the liquid region to touch the boundary of the domain along the vertical line passing through the center of the hole. The latter situation is non-generic, but for certain positions of the hole it does occur: we do not know whether asymptotic discrete Gaussianity survives in this case.

\subsection{Lozenge tilings of C-shaped domains}\label{Section_tiling_C}

\begin{figure}[t]
\begin{center}
\includegraphics[height=8cm]{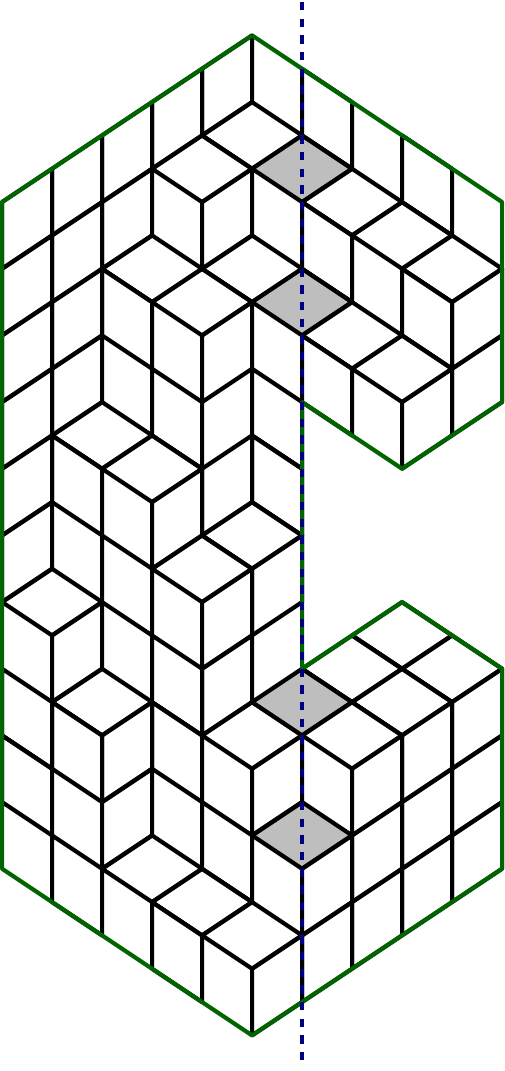}\qquad \qquad
\includegraphics[height=8cm]{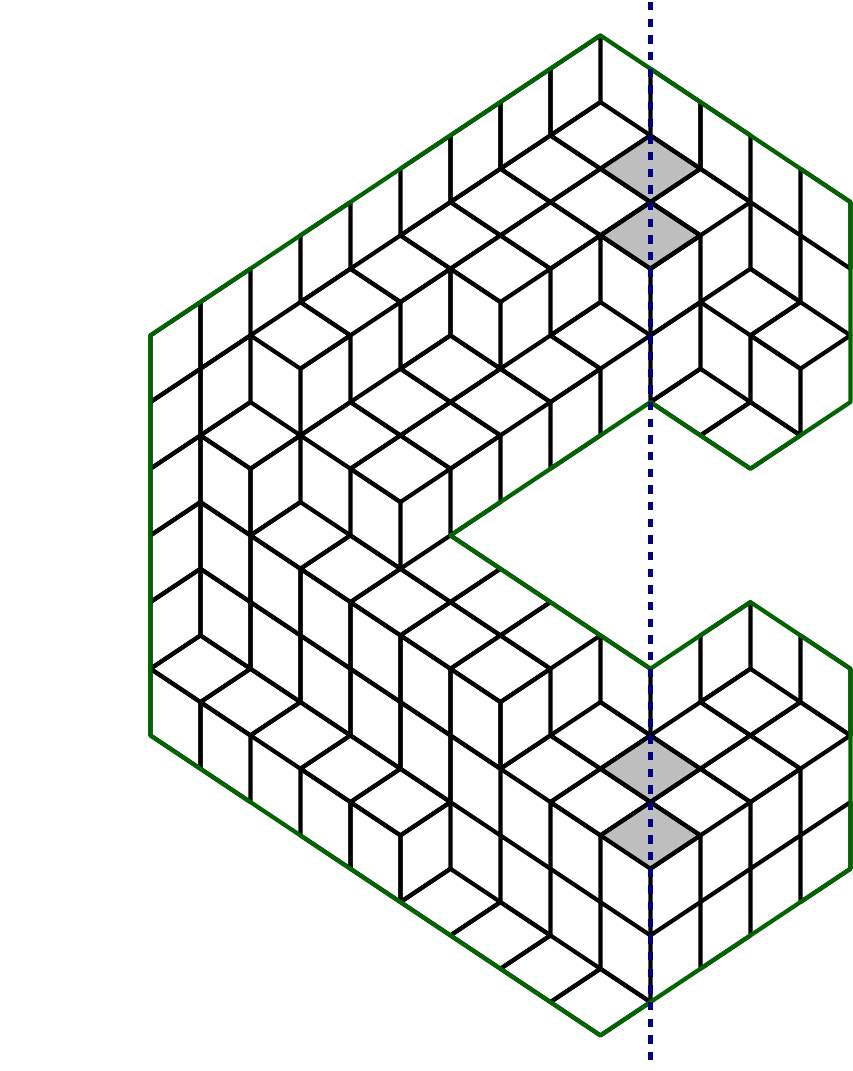}
\caption{\label{Fig:Cshapeint} Tilings of two different $C$-shaped regions. Gray horizontal lozenges are distributed as \eqref{eq_measure_general_intro}.}
\end{center}
\end{figure}

\begin{figure}[t]
\begin{center}
\includegraphics[height=8cm]{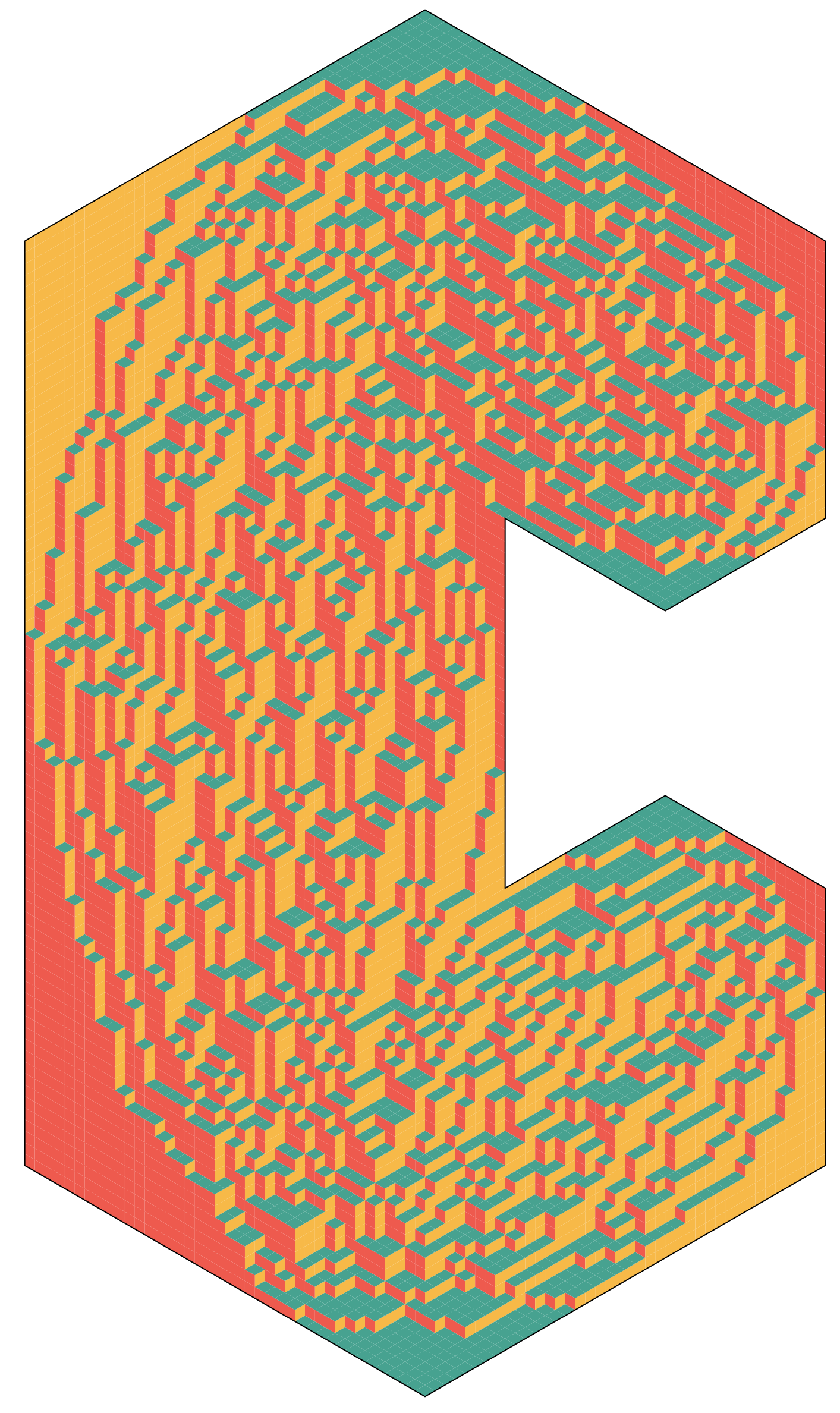}\qquad \qquad \qquad
\includegraphics[height=8cm]{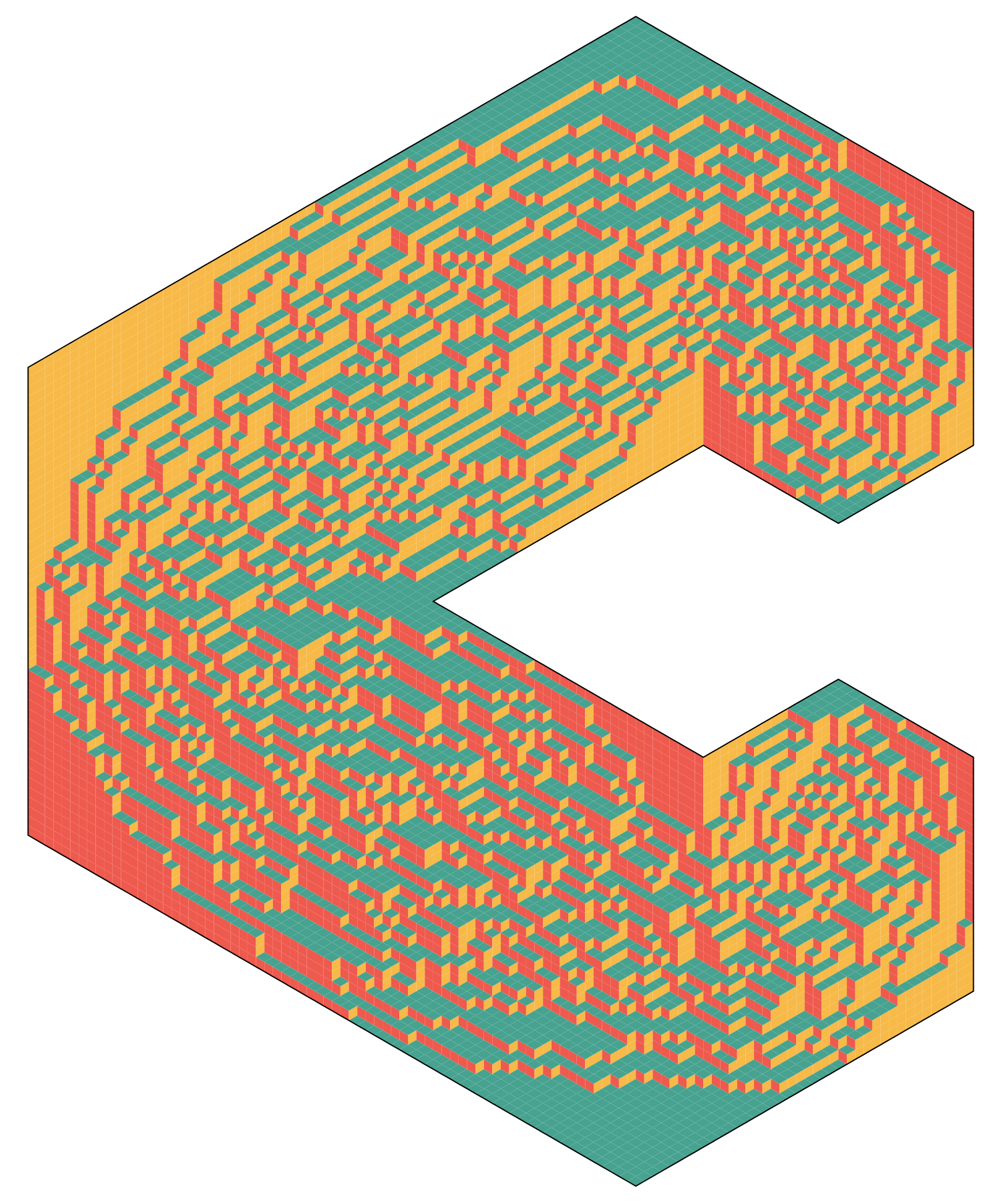}
\caption{\label{Fig:Cshape_simulation} Uniform lozenge tilings of domains in Figure~\ref{Fig:Cshapeint} resized $8x$. Simulation using \cite{Leo_ultimate}.}
\end{center}
\end{figure}

Our next example is still a uniformly random lozenge tiling model, but in a different domain: we consider either of the C-shaped domains of Figure~\ref{Fig:Cshapeint}. The figure shows finite domains and for asymptotic questions one should consider the larger similarly-shaped domains of the same proportions or, equivalently, fixed domains tiled with lozenges of smaller and smaller mesh sizes. A simulation is shown in Figure~\ref{Fig:Cshape_simulation}.

Let us first look at the horizontal lozenges along the vertical axis of the hole in the middle of the C-shape --- these are four gray lozenges on blue dashed lines in either of the panels of Figure~\ref{Fig:Cshapeint}. The lozenges are naturally split into two groups: $N_1$ of them are above the gap, and $N_2$ below the gap. In contrast to Section~\ref{Section_tiling_hole}, both $N_1$ and $N_2$ are deterministically fixed by the geometry of the C-shaped domain. We set $N=N_1+N_2$, denote $\boldsymbol{\ell} = (\ell_1<\ell_2<\cdots<\ell_N)$ the positions of the horizontal lozenges along the vertical axis and define the splitting of them into $H=2$ groups as
\[
 h(i)= \left\{\begin{array}{lll} 1 & & 1\leq i \leq N_1, \\ 2 & & N_1+1\leq i \leq N. \end{array}\right.
\]

A general counting argument explained in detail in Sections \ref{Section_simple_domains} and \ref{Section_gluing_def} leads to an explicit formula for the distribution of $\ell_1,\ldots,\ell_N$ of the form:

\begin{equation}
\label{eq_weight_C_shape}
 \amsmathbb P_N(\boldsymbol{\ell})=\frac{1}{\mathscr{Z}_N} \cdot \prod_{1\leq i<j\leq N_1} (\ell_j-\ell_i)^2 \cdot
 \prod_{N_1+1\leq i<j\leq N} (\ell_j-\ell_i)^2 \cdot
 \prod_{\substack{1 \leq i \leq N_1 \\ N_1 + 1 \leq j \leq N}}(\ell_j-\ell_i) \cdot
 \prod_{i=1}^N w_{h(i)}(\ell_i),
 \end{equation}
The weights $w_1(\ell)$ and $w_2(\ell)$ are products of Pochhammer symbols similar to \eqref{eq_Hahn} or \eqref{weight:holey}; their expressions are slightly different for the domains in the left and right panels of Figure~\ref{Fig:Cshapeint}. This distribution \eqref{eq_weight_C_shape} is of the form \eqref{eq_measure_general_intro} with interaction matrix
\begin{equation}
\label{Cintmatr}
\boldsymbol{\Theta} = \left(\begin{array}{cc} 1 & \tfrac{1}{2} \\[4pt] \tfrac{1}{2} & 1 \end{array}\right).
\end{equation}

Compared to the hexagon with a hole, the tilings of C-shaped domains are more difficult to approach with exact formulae. One manifestation is that the equilibrium measure for the former can be encoded through a solution of a single polynomial equation (whose coefficients depend in a non-trivial fashion on parameters of the hexagon and the hole), see Section~\ref{Section_hex_hole} and especially Lemma~\ref{Lemma_qpm_hex_hole}. In contrast, for the C-shaped domain one needs to solve a system of three coupled equations, and the solution is much less explicit, see Section~\ref{sec:Csha} for related discussions.

\begin{figure}[t]
\begin{center}
\includegraphics[height=7cm]{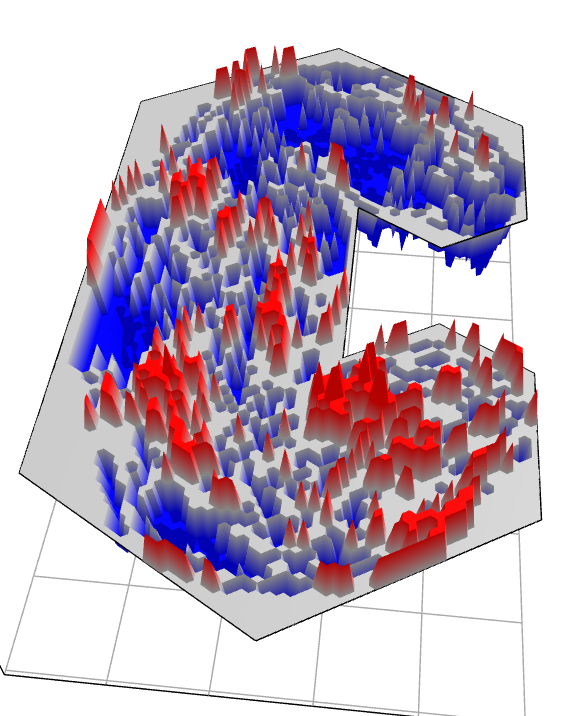}\qquad \qquad \qquad
\includegraphics[height=7cm]{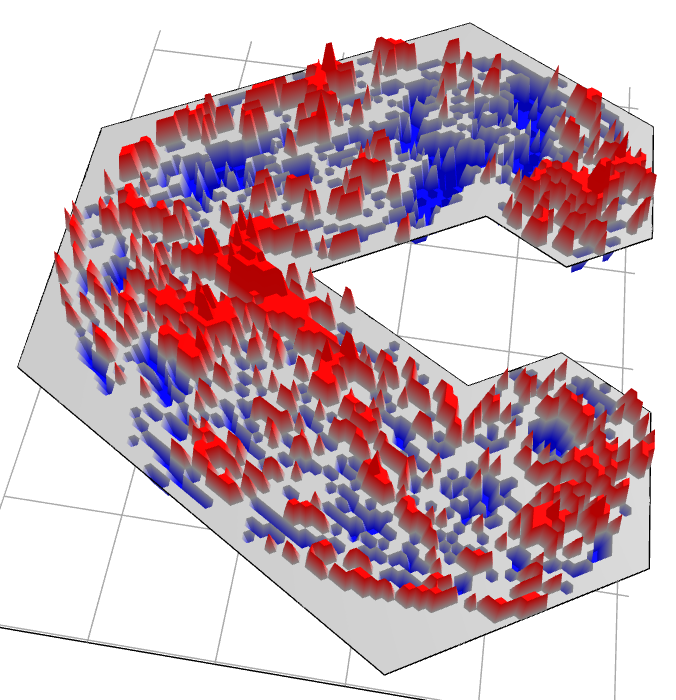}
\caption{\label{Fig:Cshape_simulation_fluct} Difference between heights of two independent samples as in Figure~\ref{Fig:Cshape_simulation}, divided by $\sqrt{2}$: for Gaussian fields, this has the same distribution as a single sample recentered around its mean.}
\end{center}
\end{figure}

Nevertheless, we can describe macroscopic fluctuations of tilings for the C-shaped domains. For that we use the language of height functions. The detailed definition and discussion is given in Section~\ref{sec:complexslopel}. For now, it is sufficient to define the height function of a tiling of a domain of linear scale $\N$, denoted $\mathsf{Ht}_\N(x,y)$, as the number of horizontal lozenges {\scalebox{0.16}{\includegraphics{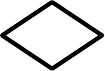}}} directly above the closest lattice point to $(\N x,\N y)$. In this language, the law of large numbers for tilings is the statement on convergence of the random variable $\frac{1}{\N}\mathsf{Ht}_\N(x,y)$ towards a deterministic limit shape $\mathfrak{Ht}(x,y)$. In the frozen regions (where asymptotically only one type of lozenges is present) the limit shape $\mathfrak{Ht}(x,y)$ is flat, while in the liquid region all three types of lozenges are present and $\mathfrak{Ht}(x,y)$ is curved. We are interested in the asymptotic behavior of $\mathfrak{Ht}(x,y)$ recentered with respect to its mean as $\N\rightarrow\infty$, see Figure~\ref{Fig:Cshape_simulation_fluct} for a simulation. Here is an informal version of our theorem on fluctuations, see Theorems~\ref{Theorem_GFF} and \ref{Theorem_GFF_general} for the detailed statements.
\begin{theoremI} \label{Theorem_C_shape_intro}
 Assume that the proportions of the C-shaped domain depend on a large parameter $\N$ and grow linearly with $\N$. Then, under a non-degeneracy condition, the random field $\sqrt{\pi}\big(\mathsf{Ht}_\N-\E[\mathsf{Ht}_\N]\big)$ converges as $\N \rightarrow \infty$ in
 the liquid region to the Gaussian free field in the complex structure defined by the Kenyon--Okounkov complex slope and with Dirichlet boundary
 conditions.
\end{theoremI}

The conjectural description of the fluctuations of heights for tilings via the Gaussian free field and the definition of the appropriate complex structure goes back to the work of Kenyon and Okounkov, see \cite[Section 2.3]{kenyon2007limit}. A detailed exposition and heuristics of the predictions can be found in \cite[Lectures 11-12]{Vadimlecture}. In Theorem~\ref{Theorem_GFF} we prove the Kenyon--Okounkov conjecture for lozenge tilings of a very large class of polygonal domains, which can be planar or embedded into orientable surfaces, might have complicated topologies, and always have non-trivial configurations of liquid and frozen regions. Our domains are obtained by gluing multiple elementary building blocks called ``trapezoids'' along a single vertical section; the distribution of horizontal lozenges along this section has the form \eqref{eq_measure_general_intro}, thus connecting this setting to the framework of discrete $\sbeta$-ensembles. We also go beyond that and in Theorem~\ref{Theorem_GFF_general} we deal with domains embedded into non-orientable surfaces (\textit{e.g.} a M\"obius strip). We show that in the non-orientable case the Kenyon--Okounkov conjecture needs to be adjusted: one should consider the oriented double covering of the domain and use the symmetrically conditioned Gaussian free field on the cover.

For the orientable case, there is  already an extensive literature related to the Kenyon--Okounkov conjecture, including developments based on determinantal processes and discrete harmonic functions in \cite{kenyon2000conformal,kenyon2001dominos,kenyon2008height} and more recently in \cite{li2013conformal}, \cite{berestycki2016note,berestycki2020dimers,berestycki2024dimers1,berestycki2024dimers}, \cite{russkikh2018dimers,russkikh2020dominos}, \cite{basok2023dimers,nicoletti2025temperleyan}; based on t-embeddings and origami maps in \cite{chelkak2023dimer,chelkak2021bipartite} and \cite{berggren2024perfect,berggren2024perfect2}; based on contour integrals and orthogonal polynomials in \cite{borodin2014anisotropic,duits2013gaussian,petrov2015asymptotics,duits2018global,berggren2025gaussian}; based dynamical loop equations in \cite{huang2020height,gorin2024dynamical,dimitrov2024global}; based on Schur generating functions in \cite{BuGo2, bufetov2018asymptotics, BuGo3, boutillier2021limit, ahn2020global, ahn2022lozenge}.

Despite numerous previous approaches and results, our Theorem~\ref{Theorem_GFF} is significantly different. The main novelty is our ability to handle the situations in which the tiled domain both has a non-trivial topology and leads to rich (that is, non-flat) limit shapes for the height functions. The particular case of the C-shaped domains in Theorem~\ref{Theorem_C_shape_intro} is also new\footnote{While the full details were not yet written down, it might be possible to adapt the method \cite{huang2020height} to the domain in the left panel of Figure~\ref{Fig:Cshapeint}, but not to the domain in the right panel.}. Furthermore, we have not seen detailed studies of the asymptotic heights for random tilings of non-orientable domains, as in our Theorem~\ref{Theorem_GFF_general}.

A minor shortcoming of our approach is the non-degeneracy condition in Theorem~\ref{Theorem_C_shape_intro}, which is a close relative of the off-criticality condition already present in the analysis of continuous $\sbeta$-ensembles and mentioned in Section~\ref{Section_continuous}. Fortunately, for tilings a large part of the off-criticality can be checked by convexity arguments (Lemma~\ref{Lemma_glued_tilings_assumptions}). Our analysis crucially relies on the off-criticality assumption, and this condition is somewhat generic, see \cite{Kuijlaarsgen} and \cite{colombo2024generic}, we do not know how to move forward when it fails. We expect that the statement remains true for certain test functions without this assumption, see \cite{BLS} for results in this direction for continuous $\sbeta$-ensembles

\subsection{The \texorpdfstring{$zw$}{zw}-measures}\label{zwintro}

For our next example, we look at the $zw$-measures, which are curious objects stemming from asymptotic representation theory and non-commutative harmonic analysis. Directly from \cite[Section 2]{Olsh_hyper}, the object of interest is a probability distribution on $N$ tuples of integers $\lambda_1\geq \lambda_2\geq \ldots\geq\lambda_N\in\amsmathbb{Z}$ depending on five parameters $\mathsf{z}_1,\mathsf{z}_2,\mathsf{w}_1,\mathsf{w}_2$ and $\theta$:
\begin{equation}
\label{eq_zw_measure_from_O}
 \frac{1}{\mathscr{Z}_N^{zw}} \cdot \prod_{1\leq i < j \leq N} \frac{\Gamma\big(\ell_j-\ell_i +1\big) \cdot \Gamma\big(\ell_j-\ell_i+\theta\big)}{\Gamma\big(\ell_j-\ell_i\big) \cdot \Gamma\big(\ell_j-\ell_i+1-\theta\big)} \cdot \prod_{\substack{1 \leq i \leq N \\ k = 1,2}} \frac{1}{\Gamma\big(\mathsf{z}_k-\ell_i+1\big)\cdot \Gamma\big(\mathsf{w}_k+\ell_i+(N-1)\theta+1\big)},
\end{equation}
where $\ell_i=\lambda_{N+1-i}+\theta(i-N)$. One needs to restrict the choices of the parameters to guarantee that \eqref{eq_zw_measure_from_O} is nonnegative and summable. For instance, assuming $\theta>0$, and $\mathsf{z}_1,\mathsf{z}_2,\mathsf{w}_1,\mathsf{w}_2\in\amsmathbb{C}$ such that $\mathsf{z}_2=\mathsf{z}_1^*$ and $\mathsf{w}_2=\mathsf{w}_1^*$ with $\mathrm{Re}(\mathsf{z}_1)>0$ and $\mathrm{Re}(\mathsf{w}_1)>0$ will work. Other choices are possible, for instance \eqref{eq_Hahn} can be matched with a special case of \eqref{eq_zw_measure_from_O} with all parameters being integers.

The first detailed investigation of $zw$-measures \eqref{eq_zw_measure_from_O} at $\theta=1$ goes back to the work of Borodin and Olshanski \cite{Olsh_zw,BoOl2005} where they appear as coefficients of decomposition of the generalized biregular representation\footnote{The precise formalism is that of spherical representations of the Gelfand pair $\big(\textnormal{U}(\infty)\times \textnormal{U}(\infty), \textnormal{U}(\infty)\big)$.} of the infinite-dimensional unitary group $\textnormal{U}(\infty)$. Here $\lambda_1\geq\dots\geq\lambda_N$ are identified with the highest weights parameterizing the irreducible representations of the unitary group $\textnormal{U}(N)$. The crucial property of \eqref{eq_zw_measure_from_O} underlying its relevance for the asymptotic representation theory is that it is coherent over $N$: the formula makes sense for all positive integers $N$ and the measure with $N=K-1$ can be linked to the one with $N=K$ through an explicit Markov transition kernel arising from the branching rules for the irreducible representations of $\textnormal{U}(K)$. The representation-theoretic constructions of \cite{Olsh_zw} can be extended to $\theta=\frac{1}{2}$ and $\theta=2$ in \eqref{eq_zw_measure_from_O}, corresponding to infinite-dimensional symmetric spaces $\textnormal{U}(\infty)/\textnormal{O}(\infty)$ and $\textnormal{U}(2\infty)/\textnormal{Sp}(\infty)$ related to orthogonal and symplectic groups, see \cite{okounkov1998asymptotics,okounkov2006limits} and references therein. The case of general $\theta>0$ and especially the form of the pair interaction is intimately related with the theory of Jack polynomials, see \cite{Olsh_hyper} and \cite{gorin2015multilevel}. It admits a $(q,t)$-deformation related to the Macdonald polynomials, see \cite{olshanski2021macdonaldb,olshanski2021macdonalda}. The continuous versions of the $zw$-measures \eqref{eq_zw_measure_from_O} were also studied, see \cite{Neretin_zw,neretin2003rayleigh} and \cite{borodin2001infinite}, and are related to the multivariate Bessel functions.

The $zw$-measures are described by \eqref{eq_measure_general_intro} with $H = 1$ group, intensity of repulsion $\theta_{1,1} = \theta > 0$ and weight $w_1$ given by the product of the four inverse Gamma functions. We remark that in order to see a non-trivial limit shape and fluctuations one should rescale the parameters $\mathsf{z}_1,\mathsf{z}_2,\mathsf{w}_1,\mathsf{w}_2$ linearly with $N$. This is unusual from the point of view of asymptotic representation theory of unitary groups, where these parameters are typically kept fixed with $N$. The equilibrium measure and the central limit theorem for linear statistics in this ensemble were previously studied in \cite[Section 9.4]{BGG}, see also Proposition~\ref{Proposition_ZW_equilibrium_measure}. Complimentary to them, we study in detail the partition function $\mathscr{Z}_N^{zw}$. It was evaluated in closed form by Olshanski \cite[(2.2)]{Olsh_hyper}
 \begin{equation}
 \label{eq_zw_partition_intro}
 \begin{split}
  \mathscr{Z}_N^{zw} & = \prod_{1\leq i < j \leq N} \frac{\Gamma\big( (j-i)\theta +1\big) \cdot \Gamma\big((j-i)\theta+\theta\big)}{\Gamma\big((j-i)\theta\big) \cdot \Gamma\big((j-i)\theta+1-\theta)} \\
  & \quad \times \prod_{i=1}^N \Bigg[ \frac{\Gamma\big(\mathsf{z}_1+\mathsf{z}_2+\mathsf{w}_1+\mathsf{w}_2+(i-1)\theta+1\big)}{\Gamma\big((i-1)\theta+1\big)}
 \cdot \prod_{k,l=1}^2 \frac{1}{\Gamma\big(\mathsf{z}_k+\mathsf{w}_l+(i-1)\theta+1\big)}\Bigg].
\end{split}
\end{equation}
This formula has the same importance for discrete ensembles as the Selberg formula \eqref{eq_Selberg_integral} has for the continuous $\sbeta$-ensembles. It can be related to the multivariate hypergeometric summation formula of \cite{gustafson1987multilateral}, or be derived by purely representation-theoretic means if $\theta \in \{\frac{1}{2},1,2\}$. Its $N=1$ case is the ${}_2 H_2$-summation formula going back by more than a century to the work of Dougall \cite{dougall1906vandermonde} and generalizing Gau\ss{} summation theorem for ${}_2F_1(a,b,c;z = 1)$ --- which corresponds furthermore to having one of the parameters $\mathsf{z}_1,\mathsf{z}_2,\mathsf{w}_1,\mathsf{w}_2$ being an integer. The asymptotic expansion of $\log(\mathscr{Z}_N^{zw})$ is then achieved in Proposition~\ref{Proposition_ZW_partition} by rewriting \eqref{eq_zw_partition_intro} in terms of Barnes double Gamma function and exploiting the known asymptotics expansion for the latter \cite{Spreafico}. The result takes the following form.
\begin{theoremI} \label{Theorem_zw_intro}
As $N \rightarrow \infty$ we have
\[
\log(\mathscr{Z}_N^{zw}) = N^2 \mathscr{F}^{[0],zw} + \theta N \ln N + N \mathscr{F}^{[1],zw}  + \bigg(\theta + \frac{1}{\theta} + 3\bigg)\ln N + \mathscr{F}^{[2],zw} + o(1),
\]
where $(\mathscr{F}^{[p],zw})_{p = 0}^{2}$ are explicit functions of the rescaled parameters $\frac{\mathsf{w}_1}{N},\frac{\mathsf{w}_2}{N},\frac{\mathsf{z}_1}{N},\frac{\mathsf{z}_2}{N}$. The asymptotic is valid and uniform when these parameters are in the range described in Definition~\ref{Definition_ZW_set} combined with \eqref{zwchoi}.
\end{theoremI}
In contrast, for general measures of the form \eqref{eq_measure_general_intro} there is no hope for a closed form evaluation of the partition function. Nevertheless, it does not prevent accessing their asymptotic expansion and we establish it in Chapter~\ref{Chapter_partition_functions}. Our key results in this direction are Theorems~\ref{Theorem_partition_one_band} and \ref{Theorem_partition_multicut}, which claim the existence of the large $N$ expansions of the logarithm of the partition function up to $o(1)$ order (in the fixed filling fractions situation and yet again assuming off-criticality). We provide formulae for the coefficients of the first two terms of orders $N^2$ and $N\log(N)$. Some formulae for the lower-order terms can be in principle collected from various pieces in our text, but we have not undertaken the task to assemble them into a nice and compact form: this remains an outstanding open question. Although we do not carry this out, it is possible to push our method to get the existence of asymptotic expansion to all-orders in $\frac{1}{N}$ like in \eqref{topexp}.

Olshanski formula for the partition function of the $zw$-discrete ensembles is an important ingredient in our proofs for general discrete ensembles \eqref{eq_measure_general_intro}. Indeed, our strategy is to interpolate between any given ensemble of interest and $zw$-measures, while controlling the change of the partition function in the process --- an overview will be given in Section~\ref{intro_interp}. In this way, we eventually get a formula for the logarithm of the ratio of the partition functions of two ensembles, similarly to \eqref{topexp}, and need to supplement it with the asymptotic formula for the $zw$-measures achieved in Proposition~\ref{Proposition_ZW_partition}. An analogous strategy was used for continuous ensembles in \cite{BG11,BG_multicut} with the Gaussian $\sbeta$-ensemble and its partition function \eqref{eq_Selberg_integral}.

\subsection{The discrete ensemble with Gaussian weight}\label{dis_Gauss_intro}

Our last example is the discrete ensemble with Gaussian weight\footnote{This should not be confused with discrete Gaussian random variables. These two classes intersect only for $N = 1$.}, which is a discretization of the Gaussian $\sbeta$-ensemble. In the notation of \eqref{eq_measure_general_intro}, it has a single group of particles with interaction intensity $\theta > 0$ and Gaussian weight depending on a parameter $\kappa>0$. The configurations in the ensemble are $N$-tuples of integers $\lambda_1\leq \lambda_2\le\dots\leq \lambda_N\in\amsmathbb{Z}$ encoded through
\begin{equation}
\label{le113}
\ell_i=\lambda_i+\theta\bigg(i - \frac{N+1}{2}\bigg)
\end{equation}
and with probability distribution
\begin{equation}
\label{eq_discrete_model_Gauss}
\amsmathbb{P}_\N(\boldsymbol{\ell}) = \frac{1}{\mathscr{Z}_N} \cdot \prod_{1 \leq i < j \leq N} \frac{\Gamma\big(\ell_j - \ell_i + 1\big)\cdot \Gamma\big(\ell_j - \ell_i + \theta\big)}{\Gamma\big(\ell_j - \ell_i\big) \cdot \Gamma\big(\ell_j - \ell_i + 1 - \theta\big)} \cdot \prod_{i = 1}^{N} e^{-\kappa \frac{\ell_i^2}{N}}.
\end{equation}
For $\theta = 1$ the distribution \eqref{eq_discrete_model_Gauss} was first studied in \cite{DK} in the context of pure Yang--Mills theory on the two-dimensional sphere, see also \cite{BoutetDeMonvelShcherbina1998,levy2015douglas,LeMa,lemoine2025two} for further developments. It was also studied in the context of non-intersecting random walks and discrete orthogonal polynomials in \cite{Liechty}. Still in the context of discrete orthogonal polynomials, a generalization of this model where $V(x) = \kappa x^2$ is replaced by an analytic potential is discussed in \cite[Chapter 3]{BleherL}, but the connection to discrete orthogonal polynomials would be lost for $\theta \neq 1$.

\begin{figure}[t]
\begin{center}
\includegraphics[height=4.8cm]{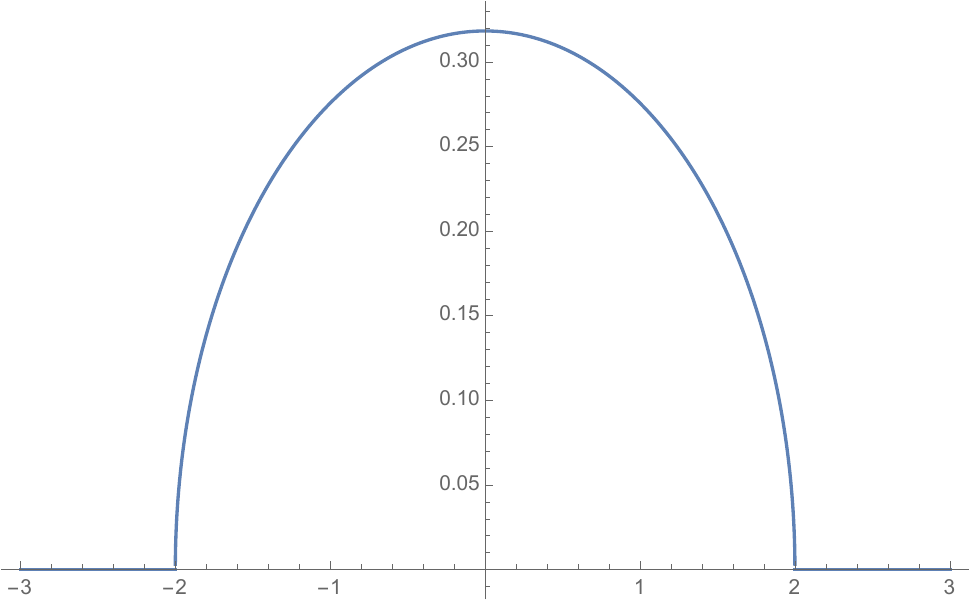} \hfill \includegraphics[height=4.8cm]{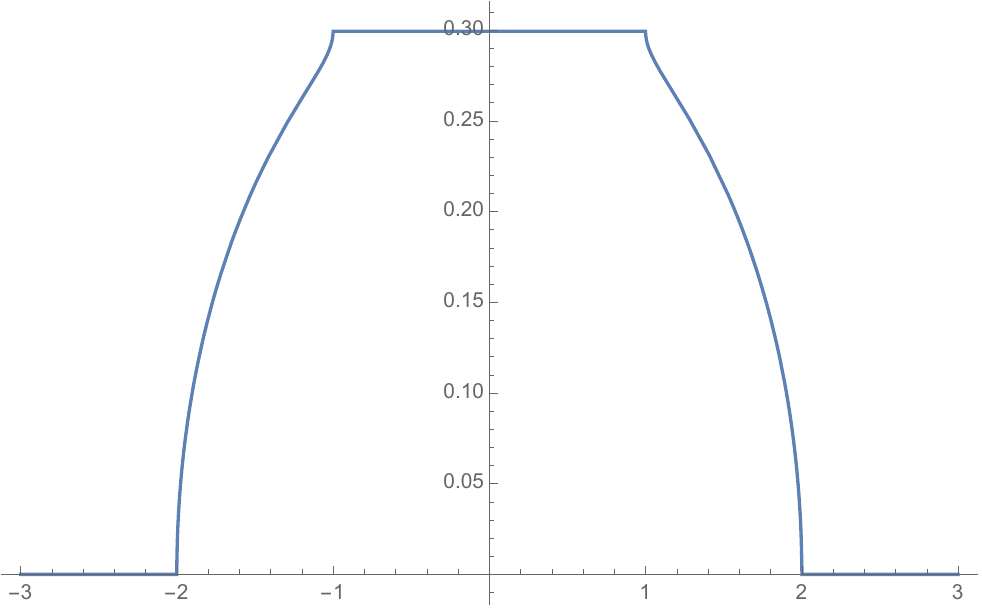}
 \caption{\label{Fig_Gaussian_intro} Equilibrium measure for the discrete ensemble with Gaussian weight. Left panel: $\theta=2$, $\kappa=1$ leading to the semicircle supported on $[-2,2]$. Right panel: $\theta \approx 3.34$, $\kappa \approx 1.69$ leading to $(\alpha,\beta) = (1,2)$.}
 \end{center}
\end{figure}

The remarkable feature discovered in \cite{DK} is a phase transition in the system, depending on the value of $\theta\kappa$. This transition is best described in the language of the equilibrium measure
\[
\mu = \lim_{N\rightarrow\infty} \frac{1}{N} \sum_{i=1}^N \delta_{\frac{\ell_i}{N}},
\]
whose detailed computation is reviewed in Section~\ref{mueqGaussian}. In the weak confinement phase $\sqrt{2\theta\kappa} < \pi$, the equilibrium measure is a semi-circle density (as for the continuous Gaussian $\sbeta$-ensemble) of support\footnote{In this book, $\beta$ will be the generic notation for right endpoints of bands. It has nothing to do with the intensity of repulsion $\sbeta$ in the continuous $\sbeta$-ensembles, for which we used a slightly different font. This unfortunate collision of notation only occurs in the introductory chapter.} $(-\beta,\beta)$, that is
\begin{equation}
 \mu(x)= \left\{\begin{array}{lll} \dfrac{\kappa}{\pi \theta}\sqrt{\beta^2 - x^2 } && |x| \leq \beta \\[10pt] 0 && |x| > \beta \end{array}\right. \qquad \textnormal{with}\quad \beta = \sqrt{\frac{2\theta}{\kappa}}.
\end{equation}
In contrast, in the strong confinement phase $\sqrt{2\theta\kappa} > \pi$ the density is more complicated:
\begin{equation}
\label{eq_insidesat_intro} \mu(x)= \left\{\begin{array}{lll}
 \dfrac{1}{\theta} & & |x| \leq \alpha \\[10pt]
  \dfrac{2 \sqrt{(\beta^2 - x^2)(x^2 - \alpha^2)}}{\pi \theta \beta |x|} \cdot \Pi_1\bigg[-\dfrac{\alpha^2}{x^2}\,;\,\dfrac{\alpha}{\beta}\bigg] & & |x| \in (\alpha,\beta) \\[10pt]
  0&& |x|> \beta \end{array}\right.
\end{equation}
Here $\Pi_1[\mathsf{c};\mathsf{k}]$ is an incomplete elliptic integral while the constants $\alpha$ and $\beta$ are found by solving two explicit equations involving $\theta$, $\kappa$, and complete elliptic integrals, see Proposition~\ref{Prop_Gaussian_LLN} for the full details and Figure~\ref{Fig_Gaussian_intro} for a plot of the densities.

Note that the minimal spacing between particles $\ell_i$ in \eqref{eq_discrete_model_Gauss} is $\theta$ and therefore the maximal possible density of the equilibrium measure is $\frac{1}{\theta}$. We refer to the segments in which the density takes the value $\frac{1}{\theta}$ as \emph{saturated regions}, while the maximal intervals in which the density is strictly inside $(0,\frac{1}{\theta})$ are called \emph{bands}. With this terminology, in the weak confinement phase the equilibrium measure of the discrete ensemble with Gaussian weight has a single band $(-\beta,\beta)$, while in the strong confinement phase there are two bands $(-\beta,-\alpha)$ and $(\alpha,\beta)$ separated by a saturated region $[-\alpha,\alpha]$. In this book we show that the phase transition also drastically affects the asymptotic behavior of linear statistics.

\begin{theoremI} \label{Theorem_Gauss_intro}
Suppose that $N\rightarrow\infty$ with $\theta>0$ and $\kappa>0$ fixed in the discrete ensemble with Gaussian weight \eqref{eq_discrete_model_Gauss}. Let $f(x)$ be an analytic function and consider the random variable
\begin{equation}
\label{eq_linear_Gauss_intro}
 \sum_{i = 1}^{N} f(\ell_i) - N\int f(x)\dd\mu(x).
\end{equation}
If $\sqrt{2\theta\kappa} < \pi$, then \eqref{eq_linear_Gauss_intro} converges in distribution to a Gaussian limit with explicit mean and covariance given in Proposition~\ref{CLTGAUSSIANWEIGHT}. If $\sqrt{2\theta\kappa} > \pi$, then \eqref{eq_linear_Gauss_intro} converges as $N\rightarrow\infty$ towards a sum of two independent random variables: the first one is a deterministic affine function of a discrete Gaussian random variable; the second one is Gaussian with explicit covariance given in Proposition~\ref{CLTGAUSSIANWEIGHT2}.
\end{theoremI}
The one-band case in Theorem~\ref{Theorem_Gauss_intro} is  the specialization to Gaussian weights  --- carried out in Section~\ref{Gausscorr} --- of Corollaries~\ref{Corollary_CLT} or \ref{Corollary_CLT_relaxed} which cover much more general ensembles of the form \eqref{eq_measure_general_intro}. It can also be extracted from \cite[Section 9.3]{BGG} but an explicit formula for the mean was not provided. The two-band case in Theorem~\ref{Theorem_Gauss_intro} is new. The covariance of the continuous Gaussian part of the fluctuations is the same than in a continuous $\sbeta$-ensembles with two bands $(-\beta,-\alpha)$ and $(\alpha,\beta)$. We do not compute the mean very explicitly, neither the parameters of the discrete Gaussian part of the fluctuations, even though all ingredients to do so are present in this book.  The far reaching generalization of Theorem~\ref{Theorem_Gauss_intro} covering all ensembles of the type \eqref{eq_measure_general_intro} under mild assumptions is achieved in Theorem~\ref{Theorem_linear_statistics_fluctuation_sat}.

The discrete Gaussian random variable in Theorem~\ref{Theorem_Gauss_intro} in the strong confinement phase is of similar nature as the discrete Gaussian random variable in Theorem~\ref{Theorem_hex_hole_intro}. Informally, it should be thought as the total number of particles in one of the bands, say, in the segment $[\N\alpha,\N\beta]$. If $\theta=1$ it is easy to make this identification precise by considering the random variable
\[
N_+=\#\{i\mid \ell_i>0\}.
\]
As $N\rightarrow\infty$, the discrete Gaussian in Theorem~\ref{Theorem_hex_hole_intro} becomes related to $N_+$ by a deterministic affine transformation. When $\theta\neq 1$, the situation is trickier, because the particles $\ell_i$ do not sit on fixed lattice sites. Instead each $\ell_i$ with fixed $i$ comes with its own $\theta$-dependent shift and $N_+$ is no longer appropriate. For instance, if $\theta=2$ and $N$ is odd, then near $0$ every second integer would be occupied by a particle and there would be two classes of configurations, depending whether $0$ or $1$ are occupied  by some $\ell_i$. The random variable $N_+$ mixes these two cases together, while distinguishing them is important for our construction of the discrete Gaussian random variable in Theorem~\ref{Theorem_Gauss_intro}. Instead of $N_+$, one can take the position of the particle $\ell_{\lfloor \frac{N}{2} \rfloor}$ as the random variable related by an affine transformation to the discrete Gaussian random variable of Theorem~\ref{Theorem_Gauss_intro}.

In our most general Theorem~\ref{Theorem_linear_statistics_fluctuation_sat} (and the closely related Theorem~\ref{Theorem_CLT_for_filling_fractions_saturation}) we need to rely on a mix of ordinary filling fractions, as in Theorem~\ref{Theorem_hex_hole_intro}, and positions of individual particles, as in Theorem~\ref{Theorem_Gauss_intro}, to define the limiting discrete Gaussian component of the fluctuations. This mixture is encoded in more elaborate random variables called \emph{extended filling fractions}. We demonstrate various ingredients of this approach in Section~\ref{Sec:condgaus} (in the case of Gaussian weight), Section~\ref{fluctuatnecmergitur}, and in Section~\ref{Section_Fillingfractionsatcasesec} (in the general case).

\section{Methodology} \label{method}
\label{Chap1Sec3}

We present the main ideas developed in this book, illustrated with the special case of the random tilings of the hexagon with a hole or the C-shaped domain introduced in Sections~\ref{Section_tiling_hole}-\ref{Section_tiling_C}. Some aspects of the hexagon with a hole were already treated in \cite{BGG} and \cite{BuGo3} in the case of fixed filling fractions, but many new ideas are necessary for the full treatment even in this basic case. We will also outline which additional difficulties arise in the generalization to more complicated settings. We only give a high-level overview of the arguments, delaying the technical details to the main text.

\subsection{Law of large numbers and concentration of measures}
We deal with tilings of Figure~\ref{Fig_hex_hole_small}  in the regime where the hexagon and the hole become large, \textit{i.e.} for some $\N$ going to infinity
\begin{equation}
\label{eq_rescaled_params_intro}
A=\hat A \N, \quad B=\hat B \N, \quad C= \hat C \N, \quad D= \hat D \N, \quad \mathfrak{t}=\hat{\mathfrak{t}} \N, \quad \mathfrak{h}=\hat{\mathfrak{h}}\N,
\end{equation}
with $\hat{\mathfrak{t}} >\max(\hat B,\hat C)$, so that the hole is in the right part of the hexagon. For simplicity of exposition, we assume that the asymptotic proportions $(\hat A, \hat B, \hat C, \hat D, \hat {\mathfrak t}, \hat{\mathfrak{h}})$ are fixed integers and $\N$ remains an integer as it tends to infinity. However, all the asymptotic bounds we develop are locally uniform in the proportions, and therefore, their extension to $\N$-dependent proportions (what we do in the main text) does not introduce additional difficulties. In particular, $\N$-dependent proportions allow us achieving irrational ratios, say, $\frac{A}{B} \rightarrow \sqrt{2}$, which would be impossible with fixed $\hat A$, $\hat B$ and integral $A$, $B$.

The number $N$ of horizontal lozenges on the vertical section along the center of the hole is
  \[N=B+C-D-t= (\hat B+\hat C - \hat D-\hat{\mathfrak{t}}) \N.\]
  The distribution $ \amsmathbb P_N$ of the $N$ horizontal lozenges $\boldsymbol{\ell} = (\ell_i)_{i = 1}^N$ is described in \eqref{eq_Hahn_cut}, copied here as
  \begin{equation}
\label{eq_Hahn_cut_2}
\amsmathbb{P}_{\N}(\boldsymbol{\ell}) = \frac{1}{\mathscr{Z}_N} \cdot \prod_{1\leq i<j\leq N} (\ell_j-\ell_i)^2 \cdot \prod_{i=1}^N (A+B+C+1-\mathfrak{t}-\ell_i)_{\mathfrak{t}-B} \cdot (\ell_i)_{\mathfrak{t}-C}\cdot (\mathfrak{h}+1-\ell_i)_D^2.
\end{equation}
We recall the splitting into two groups of $N_1$ and $N_2$ particles, belonging to the respective segments
\[
\ell_i\in [a_1,b_1]= [1,\mathfrak{h}] \qquad \textnormal{or}\qquad \ell_i\in [a_2,b_2]=[ \mathfrak{h}+D+1,A+B+C-\mathfrak{t}],
\]
and denote  $I_h=\{i\,\,|\,\,\ell_i\in [a_h,b_h]\}$ for $h=1,2$. We consider two different settings:
\begin{itemize}
 \item In the fixed filling fractions setting, $N_1$ and $N_2$ are deterministically fixed.
 \item In the fluctuating filling fractions setting, $N_1$ and $N_2$ are random variables and only $N = N_1+N_2$ is deterministic.
\end{itemize}

The law of large numbers describes the asymptotic positions $\boldsymbol{\ell}$ on average, and is formulated in terms of the empirical measure of the particles. For $h=1,2$ we set
\begin{equation}
\label{eq_empirical_measure_int} \mu_{\N,h} = \frac{1}{\N} \sum_{i \in I_h}
\delta_{\frac{\ell_i}{\N}}.
\end{equation}
We view $\mu_{\N,h}$ as a nonnegative measure on  $ [\frac{a_h}{\N},\frac{b_h}{\N}]=[\hat a_h,\hat b_h]$ with  mass  $\hat n_h=\tfrac{N_h}{\N}$ which is fixed in the fixed filling fraction case, and otherwise satisfies
$\hat n_1+\hat n_2= \hat{n} := \tfrac{N}{\N}$. Moreover,
since the distance between particles is at least one, the empirical measures satisfy for all $x<y$,
\[\mu_{\N,h}([x,y])\leq y-x +\frac{1}{\N}.\]
This implies that
any limit point of $\mu_{\N,h}$ as $\N \rightarrow \infty$ admits a density bounded by $1$ with respect to the Lebesgue measure. Let $(\mu_1,\mu_2)$ be the unique minimizer of the functional
\begin{equation} \label{eq_I_intro}
-{\mathcal I}[\mu_1,\mu_2] =\iint\left(-\log |x-y| +\frac{V(x)+V(y)}{2}\right) \dd(\mu_1+\mu_2)(x)\,\dd(\mu_1+\mu_2)(y)
\end{equation}
among  the set $\mathscr{P}_{\star}$ of pairs of nonnegative measures on $[\hat a_1,\hat b_1]$  and $[\hat a_2,\hat b_2]$ respectively, having densities bounded by $1$ and masses fixed to $(\hat{n}_1,\hat{n}_2)$ (in the fixed filling fractions case), or whose sum is fixed to be $\hat{n}$ (in the fluctuating filling fraction case). Here the potential is taken to be
 \begin{equation}
 \label{eq_potential_Hahn_intro}
\begin{split}
V(x) & =\mathrm{Llog}(\hat A+\hat C-x)-\mathrm{Llog}(\hat A+\hat B+\hat C-\hat{\mathfrak{t}}-x) + \mathrm{Llog}(x+\hat{\mathfrak{t}}-\hat C) -\mathrm{Llog}(x) \\
& \quad +2\big(\mathrm{Llog}(\hat D+\hat{\mathfrak{h}}-x)- \mathrm{Llog}(\hat{\mathfrak{h}}-x)\big),
\end{split}
\end{equation}
where $\textnormal{Llog}(x) := x\log |x| - x$.

\begin{theorem}[Law of large numbers]
\label{Theorem_LLN_intro}
For each $h \in \{1,2\}$, the recentered empirical measure $(\mu_{\N,h}-\mu_h)$ vanishes in probability for the weak topology as $\N \rightarrow \infty$. In other words, for any bounded continuous function $f$ on the real line
  \[\lim_{\N\ra\infty}\left| \int f(x)\dd\mu_{\N,h}(x)-\int f(x)\dd\mu_h(x)\right|=0\qquad \text{ in probability.}\]
   \end{theorem}
The potential $V(x)$ in \eqref{eq_potential_Hahn_intro} comes from the large $\N$ behavior of the weight $w(x)$ of \eqref{weight:holey}, in the sense that
\begin{equation}
\label{eq_Hahn_Stirlng_intro}
w(\ell)=  (A+B+C+1-\mathfrak{t}-\ell)_{\mathfrak{t}-B} \cdot  (\ell)_{\mathfrak{t}-C} \cdot  (\mathfrak{h}+1-\ell)_D^2 \approx e^{-{\N} V (\frac{\ell}{\N})}
\end{equation}
thanks to Stirling formula at least when $x$ remains away from the segment endpoints. But, at the closest integer sites outside segments the weight has a zero, while $V$ remains bounded. This error near endpoints is negligible for Theorem~\ref{Theorem_LLN_intro} as $\N \rightarrow \infty$. However, it does play a role when we proceed to the study of fluctuations, as we will use tools from complex analysis and the exact coincidence of zeros of various involved functions becomes crucial. For this reason our precise definition of the asymptotic data in Chapter~\ref{Chapter_Setup_and_Examples} are going to involve $\hat{a}_h,\hat{b}_h$ shifted by terms of order $\frac{1}{\N}$: the resulting parameters will be called $\hat{a}_h',\hat{b}'_h,\ldots$ \textit{etc.}, see Definition~\ref{def:eq_shifted_parameters}. To simplify the exposition in this introduction we ignore the shifts as long as we can.

By standard methods of potential theory \cite{DS,ST}, the existence and uniqueness of $(\mu_1,\mu_2)$ as the minimizer of $-{\mathcal I}$ follows from the upper semi-continuity of ${\mathcal I}$, the compactness of upper-level sets $\{{\mathcal I}\geq M\}$, and the concavity of $\mathcal I$, see Section~\ref{Section_Energy_functional}. The minimizer is characterized by the condition
\[
{\mathcal I}[\mu_1,\mu_2]- {\mathcal I}[\mu_1+ t \nu_1,\mu_2+ t\nu_2] \geq 0
\]
for any $t \in \amsmathbb{R}$ such that the density of the measure $\mu_h + t\nu_h$ remains between $0$ and $1$ for each $h \in \{1,2\}$. This translates into the existence of Lagrange multipliers $v_1,v_2 \in \amsmathbb{R}$ (with $v_1 = v_2$ in the case of fluctuating filling fractions) such that for $x \in [\hat{a}_h,\hat{b}_h]$ the effective potential
\begin{equation}
\label{effpotd}
V^{\textnormal{eff}}(x) = V(x) - 2\int \log |x -y|\dd(\mu_1 + \mu_2)(y)
\end{equation}
\begin{itemize}
\item is larger than $v_h$ on regions with density $0$ (voids);
\item is smaller than $v_h$ on regions with density $1$ (saturations);
\item is equal to $v_h$ on regions of intermediate density (bands).
\end{itemize}

The last property is best formulated in terms of the Stieltjes transform
\[
\mathcal{G}(z) = \int \frac{\dd(\mu_1 + \mu_2)(x)}{z -x}.
\]
By design, $\mathcal{G}(z)$ is a holomorphic function of $z$ in the complex plane away from the support of $\mu_1 + \mu_2$, such that $\Gm(z) \sim \frac{\hat{n}}{z}$ as $z \rightarrow \infty$. Let us first assume that $\mu_1 + \mu_2$ always has density smaller than $1$ --- this can be achieved for instance by choosing the parameters so that the hole is close to the right border of the hexagon. Then the support of $\mu_1 + \mu_2$ only consists of bands, and differentiating the equality $V^{\textnormal{eff}}(x) = v_h$  yields the equation
\begin{equation}
\label{eq_bandsrhp}
\forall x \in \textnormal{bands}\qquad \Gm(x^+) + \Gm(x^-) = V'(x).
\end{equation}
Here $f(x^{\pm})$ indicate the limit values of $f(z)$ as $z = x \pm \ii \epsilon$ and $\epsilon \rightarrow 0^+$. Multiplying by $\Gm(x^+) - \Gm(x^-)$ we learn that
\begin{equation}
\label{algeqn}
P(z) = \Gm(z)^2 - \Gm(z)V'(z)
\end{equation}
has no discontinuity when $z$ crosses the bands. If $V$ were a polynomial of degree $d$, then $P(z)$ is an entire function growing like $z^{d - 1}$ as $z \rightarrow \infty$, and by Liouville theorem it must be a polynomial of degree $d - 1$. This equation is well-known for continuous $\sbeta$-ensembles, it is the leading order as $N \rightarrow \infty$ of the Dyson--Schwinger equation stating that certain quadratic observables depending on $z \in \amsmathbb{C} \setminus \amsmathbb{R}$ have an expectation value which is a polynomial in $z$ of degree $d - 1$.

For discrete ensembles and in particular for \eqref{eq_Hahn_cut_2}, modifications of this classical theory are necessary to take into account the saturations  (which give rise to logarithmic singularities in $\Gm(z)$) and the non-polynomiality of the potential (which is manifest in \eqref{eq_potential_Hahn_intro}). Both obstacles can be overcome by exponentiation, \textit{i.e.} by looking at $\exp(\Gm(z))$ and $\exp(V'(x))$: the former (in contrast to $\Gm(z)$) has no jumps as $z$ crosses the real axis at a saturated point; the latter becomes a rational function of $z$. As a matter of fact, for the equilibrium measure of the hexagon with a hole the algebraic equation \eqref{algeqn} is replaced with the statement of Lemma~\ref{Lemma_qpm_hex_hole} that
\begin{equation} \label{eq_quadratic_equation_on_G_intro}
 \big(\hat{\mathfrak{t}}-\hat{C}+z\big) \cdot \big(\hat{b}_1-z\big)^2 \cdot \big(\hat{b}_2-z\big) \cdot \exp(\Gm(z)) + \big(\hat{A}+\hat{C} -z\big) \cdot \big(z-\hat{a}_1\big) \cdot \big(z-\hat{a}_2\big)^2 \cdot \exp(-\Gm(z))
\end{equation}
is a polynomial of degree $4$. Generalizations of this equation exist for discrete ensembles of the form \eqref{eq_measure_general_intro} provided the weight has some good analytic properties (\textit{cf.} Assumption~\ref{Assumptions_analyticity} in Chapter~\ref{Chapter_Setup_and_Examples}), and coincide with the leading order of the Nekrasov equations of Chapter~\ref{ChapterNekra}. We will use them in Chapter~\ref{Chapter_smoothness} to establish various regularity properties of the equilibrium measure such as H\"older continuity of the density or smooth dependence on the parameters.

Having discussed the equilibrium measure, Theorem~\ref{Theorem_LLN_intro} can be deduced from large deviation principles, \textit{cf.} \cite{Feral,arous1997large,Johansson_shape,Johansson2}. We will not derive such large deviations principles but rather a quantitative estimate expressing concentration of measures as it is more suitable to study the fluctuations afterwards. The idea is to prove concentration for the pseudo-distance $\Dr$ on $\mathscr{P}_{\star}$ defined through
\begin{equation*}
\D[\nu,\rho] :=-\iint \log  |x-y| \dd(\nu-\rho)(x)\dd(\nu-\rho)(y)
= \int_{0}^{+\infty} \left|
\int e^{\mathbf{i}px}\dd (\nu-\rho)(x)\right|^{2} \frac{\dd p}{p}.
\end{equation*}
The last equality is justified in Lemma~\ref{Lemma_I_quadratic_smoothing}, see also \cite{Deiftcours}. A technical difficulty is that $\Dr$ is infinite if applied to the empirical measures $\mu_{\N,h}$ since the latter are atomic. This is overcome by following the regularization procedure of  \cite{MaMa}, that is, replacing the empirical measure $\mu_{\N,h}$ with the convolution $u_\N *\mu_{\N ,h}$ with the uniform law $u_\N $ on $[0,\frac{1}{\N }]$. This gives an absolutely continuous measure with density bounded by $1$. However, it is important for our arguments for the pair of measures remains in $\mathscr{P}_\star$, which may not be the case after convolution since the support is general only included in $[\hat a_h,\hat b_h+\frac{1}{\\N }]$. So, we apply a last step consisting of moving the mass that spilled over $[\hat{b}_h,\hat{b}_h + \frac{1}{\N}]$ to
\[
\amsmathbb{A} = [\hat{a}_1,\hat{b}_1] \cup [\hat{a}_2,\hat{b}_2],
\]
so that the resulting pair of measures $(\tilde{\mu}_{\N,1},\tilde{\mu}_{\N,2})$ belongs to $\mathscr{P}_{\star}$. We then have the following estimate, appearing in greater details in Theorem~\ref{Prop_pseudodistance_bound}.

\begin{theorem}[Concentration of empirical measures]
\label{thm_concentration}
Let $\tilde\mu_\N=\tilde\mu_{\N,1}+\tilde\mu_{\N,2}$ and $\mu=\mu_1+\mu_2$. There exists $C>0$ such that for every $t\geq 0$,
\[
\amsmathbb{P}_\N\big[\Dr[\tilde\mu_{\N },\mu]\geq t\big] \leq e^{C\N \log \N
  - \N ^2t^2}.
\]
\end{theorem}
The idea for the proof of this theorem is that from \eqref{eq_Hahn_cut_2} and \eqref{eq_Hahn_Stirlng_intro} one can deduce
\[
\amsmathbb P_N(\ell) \leq  e^{C\N\log \N +\N^2({\mathcal I}[\tilde\mu_{\N,1},\tilde\mu_{\N,2}]-{\mathcal I}[\mu_1,\mu_2])}.
\]
Then, algebraic manipulations with the difference of $\mathcal{I}$s and the characterization of $\mu$ lead to
\[
\amsmathbb P_N(\ell) \leq e^{C\N\log \N -\N^2\D(\tilde\mu_{\N},\mu)},
\]
from which we deduce Theorem~\ref{thm_concentration}. From the compatibility of the pseudo-distance $\Dr$ with the weak topology, we infer in Lemma~\ref{Lemma_tail_bound_general} and Corollary~\ref{Corollary_a_priori_0}  the following bounds.
\begin{lemma}[Concentration of linear statistics] \label{Lemma_tail_bound_general_int}
 There exists
$C>0$ such that for any any positive integer $\N$, any $t >0$ and real numbers $\mathfrak{a} < \mathfrak{b}$, we have
\begin{equation*}
 \amsmathbb P_\N\Bigg[\exists f \in \mathscr{H}_{\textnormal{Lip},\frac{1}{2}}\quad  \min_{h=1,2} \bigg|\int_{\amsmathbb R}
 f(x) \dd(\mu_{\N,h} - \mu_h)(x)\bigg| \geq t |\!| f
 |\!|_{\frac{1}{2}}+ \tfrac{C}{\N} \big( |\!|f|\!|_\textnormal{Lip} +|\!|f|\!|_\infty \big) \Bigg]  \leq  e^{C \N\log \N-\N^2 \frac{t^2}{C}},
\end{equation*}
where $\mathscr{H}_{\textnormal{Lip},\frac{1}{2}}$ is the space of Lipschitz functions on $[\hat{a}_h,\hat{b}_h]$ with finite Sobolev $\frac{1}{2}$-norm (see \eqref{eq_Fourier}), and
\begin{equation}
\label{eq_ff_LLN_Intro}
\hspace{-0.3cm}\amsmathbb P_\N\bigg[\Big|\#\big\{i \,\,\big| \,\,\ell_i \in [\N \mathfrak{a}, \N \mathfrak{b}]\big\} - \N \mu([\mathfrak{a},\mathfrak{b}]) \Big|
\geq  C \cdot \N^{\frac{1}{2}} \log^2 \N\bigg] \leq e^{-\frac{\N}{C}\log^2
\N}.
\end{equation}
\end{lemma}

For more general discrete ensembles of the form \eqref{eq_measure_general_intro}, where the interactions, the potentials, and the restrictions on the filling fractions go beyond our running example of \eqref{eq_Hahn_cut_2}, Theorems~\ref{Theorem_LLN_intro}, \ref{thm_concentration} and Lemma~\ref{Lemma_tail_bound_general_int} continue to hold. The functional ${\mathcal I}$ is going to be replaced with a multi-group version of \eqref{eq_I_intro} given in \eqref{eq_functional_general} involving the coefficients of the interaction matrix $\boldsymbol{\Theta}$ while $\D$ is taken to be the quadratic part of the new $-\mathcal{I}$. Crucially, $\D$ remains nonnegative because we take the assumption that $\boldsymbol{\Theta}$ is positive semi-definite. The existence, uniqueness and characterization of the equilibrium measure in the general setting is developed in Theorem~\ref{Theorem_equi_charact_repeat_2}. Some of the ingredients for the proofs already appeared in \cite{DS} (addressing measures with upper bound on the density, but $H =1$) and  \cite{Hardy_Kuijlaars_min} (addressing $H > 1$ with a strictly convex matrix of interaction, but without density upper bound), but we need new ingredients as we want to allow filling fractions to vary while respecting a given set of affine constraints. In particular, positive semi-definiteness of $\boldsymbol{\Theta}$ is not sufficient and we require an additional assumption \eqref{eq_Theta_through_Theta_prime} fixing the interplay between positivity of $\boldsymbol{\Theta}$ and the chosen affine constraints. If filling fractions are uniquely fixed --- like in \cite{Hardy_Kuijlaars_min} --- this additional assumption is automatically verified.

 The generalization of the polynomial \eqref{eq_quadratic_equation_on_G_intro} is introduced in Chapter~\ref{Chapter_smoothness} and used to establish various regularity properties of the equilibrium measures which would be difficult to obtain only by potential-theoretic methods, see Theorem~\ref{Theorem_regularity_density}. The sum $\Gm(z)=\Gm_1(z)+\Gm_2(z)$ appearing in the exponentials in \eqref{eq_quadratic_equation_on_G_intro} is going to be replaced with a linear combination involving the coefficients of $\boldsymbol{\Theta}$, see Definition~\ref{GQdef}.

\subsection{Dyson--Schwinger equations and Nekrasov equations}
\label{sec_Nek_intro}

Theorem~\ref{Theorem_LLN_intro} explained that the empirical measures  $\mu_{\N,h}$ approximate the equilibrium measures $\mu_h$ as $\N\rightarrow\infty$ and the next step is to compute the fluctuations, \textit{i.e.} the leading asymptotics for the differences  $\mu_{\N,h}-\mu_h$. In contrast to the sequence of i.i.d. random variables where the difference would be of order $\N^{-\frac{1}{2}}$, for tilings and for discrete ensembles of the form \eqref{eq_measure_general_intro} the difference is of order $\N^{-1}$. In the continuous $\sbeta$-ensembles of Section~\ref{Section_continuous}, Johansson's study of the fluctuations of the empirical measure \cite{Johansson} was based on  Dyson--Schwinger equations (which are also named loop equations, master equations, Ward identities, \textit{etc.} depending on the context) for correlators. These equations come from integration by parts in expectation values, resulting in exact relations between products of linear statistics averaged over the law \eqref{beta-mod}. These equations are not closed \textit{per se} but one can show that they are \emph{asymptotically closed}: in the regime where the dimension $N$ is large and with enough concentration of measure properties, they have a triangular structure and can be uniquely solved under certain conditions. This structure was first revealed in the physics literature \cite{ACM92,ACKM} and used to compute all coefficients of asymptotic expansions like \eqref{topexp2} assuming the existence of such an expansion. The same structure was exploited in \cite{BG11} to set up a scheme of self-improving estimates allowing to establish this asymptotic expansion, and we will sketch it in Section~\ref{Theorem_correlators_intro}.

  In the setting of random tilings and more general discrete ensembles \eqref{eq_measure_general_intro}, integration by parts is not available and this remained a roadblock for the extensions of the methodology for many years. A breakthrough in a different field helped overcoming this problem: Nekrasov \cite{Nekrasovpaper,nekrasov2018quantum, nekrasov2023seiberg} found that in random partition models associated to supersymmetric quiver gauge theory, certain observables depending meromorphically on $z \in \amsmathbb{C}$ and called $qq$-characters have an expectation value which is polynomial in $z$. These equations come from studying the modification of partitions by addition of a box, and for this reason were named ``non-perturbative'' Dyson--Schwinger equations. The form of the equations for the specific models considered by Nekrasov gave \cite{BGG} enough hint to discover non-trivial meromorphic observables whose expectation value is holomorphic (\textit{i.e.} has no poles) in the discrete ensembles \eqref{eq_measure_general_intro} with $H = 1$ group of particles. We will call generically \emph{Nekrasov equations} such holomorphicity statements, which have been found for many discrete ensembles beyond the gauge-theoretic ones. They offer an appropriate replacement for Dyson--Schwinger equations for discrete ensembles. Albeit fully non-linear (instead of quadratic), they retain some asymptotically triangular structure and the methodology of \cite{BG11} can still be adapted to them in order to access to the asymptotic fluctuations of the empirical measures. This strategy has been first applied in \cite{BGG} and further developed in \cite{guionnet2019rigidity,dimitrov2019log,dimitrov2022asymptotics,gorin2024dynamical,dimitrov2024global,dimitrov2025multi}. Our book continues in this vein: we find a new set of Nekrasov equations (Chapter~\ref{ChapterNekra}) adapted to the general discrete ensembles \eqref{eq_measure_general_intro}, which we intensively use to carry out the asymptotic analysis. Compared to \cite{BGG}, we also need to import the amendments of the strategy developed in \cite{BG_multicut} to analyze situations with several bands and in \cite{BGK} to analyze ensembles with several groups of particles.

Let us describe the particular instance of the Nekrasov equation relevant for the hexagon with the hole. Using the weight $w(\ell)$ of \eqref{weight:holey}, we can decompose
\begin{equation}\label{nekcond}\frac{w(\ell+\frac{1}{2})}{w(\ell-\frac{1}{2})}=\frac{\Phi^+\left(\frac{\ell}{\N}\right)}{\Phi^-\left(\frac{\ell}{\N}\right)},\end{equation}
in terms of two polynomials
\begin{equation*}
\begin{split}
\Phi^+(z) & = \big(\hat{\mathfrak{t}}-\hat{C} + z- \tfrac{1}{2\N}\big)\cdot \big( \hat{A}+\hat{B}+\hat{C}- \hat{\mathfrak{t}}-z+ \tfrac{1}{2\N}\big) \cdot \big(\hat{\mathfrak{h}}-z+ \tfrac{1}{2\N}\big)^2 \\
\Phi^-(z)  & = \big(z-\tfrac{1}{2\N}\big)\cdot \big(\hat{A} + \hat{C} - z+\tfrac{1}{2\N}\big) \cdot \big(\hat{\mathfrak{h}}+\hat{D}- z+\tfrac{1}{2\N}\big)^2.
\end{split}
\end{equation*}
\begin{theorem}[Nekrasov equation]
\label{Theorem_Nekrasov_int}
For the discrete ensemble \eqref{eq_Hahn_cut_2}, the function
\begin{equation}
\label{Nekrasov_eqn_int}  R(z)= \Phi^-(z)
  \cdot \E\left[ \prod_{i=1}^N
  \bigg(1-\frac{1}{\N z-\ell_i+\frac{1}{2}} \bigg)
   \right]
  +
  \Phi^+(z)
    \cdot \E \left[ \prod_{i=1}^N
  \bigg(1+\frac{1}{\N z-\ell_i-\frac{1}{2}} \bigg)
   \right]
\end{equation}
is holomorphic on the complex plane, \textit{i.e.} it has no poles.
\end{theorem}
Theorem~\ref{Theorem_Nekrasov_int} can be traced back to \cite[Theorem 1.1]{BGG} and Theorem~\ref{Theorem_Nekrasov} will extend it to the general discrete ensembles \eqref{eq_measure_general_intro}. The hardest part is to figure out the correct form of the observable $R(z)$ --- this is where the work of Nekrasov gives some inspiration. Once the guess is made, checking cancellation of the poles between the terms in $R(z)$ is straightforward based on the explicit form of the measure \eqref{eq_Hahn_cut_2}.

\begin{corollary}[Holomorphicity of $R(z)$]
\label{Corollary_observable_polynomial}
The function $-\frac{1}{2}R(z)$ is a monic polynomial of degree $4$.
\end{corollary}
This a direct corollary of Theorem~\ref{Theorem_Nekrasov_int} and Liouville theorem: an entire function growing like $-2z^4$ as $z\rightarrow\infty$ must be a polynomial of degree $4$. The polynomiality of \eqref{eq_quadratic_equation_on_G_intro} from the previous section is readily deduced as $\N\rightarrow\infty$ limit using Theorem~\ref{Theorem_LLN_intro}. While versions of Theorem~\ref{Theorem_Nekrasov_int} are available for all discrete ensembles \eqref{eq_measure_general_intro},
 Corollary~\ref{Corollary_observable_polynomial} relies on the polynomiality of $\Phi^\pm(z)$ which is not true for general weights. Besides, for general weights the analog of the observable $R(z)$ may develop poles near rescaled endpoints of segments where particle live, that will bring some technical (but no conceptual) complications.

\subsection{Expansion of the correlators and central limit theorem for fixed filling fractions}
\label{Theorem_correlators_intro}

Keeping on with tilings of the hexagon with a hole, let us now describe the asymptotic theorems for fluctuations which can be derived from the Nekrasov equations when the filling fractions $(N_1,N_2)$ are deterministically fixed. We call $G(z)$  the Stieltjes transform of the empirical measure and $\Delta G(z)$ its recentering compared to the equilibrium one, that is for $z \in \amsmathbb{C} \setminus \amsmathbb{A}$
\begin{equation}\label{eq_stieldef_intro}
  G(z)=\sum_{i=1}^N \frac{1}{z-\frac{\ell_i}{\N}}, \qquad \mathcal{G}(z)=\int \frac{\dd\mu(x)}{z-x}, \qquad \Delta G(z)= G(z)-\N \mathcal{G}(z).
\end{equation}
As a corollary of Lemma~\ref{Lemma_tail_bound_general_int}, we have
\begin{equation}
   \label{eq_apriory_bound_intro}
    \amsmathbb{P}_\N\Bigg[\bigg|\frac{\Delta G(z)}{\N}\bigg|\geq \frac{t}{\sqrt{|\textnormal{Im}\,z|}}+\frac{C}{\sqrt{\N}}\bigg ]\leq  e^{C \N\log \N-\N^2 \frac{t^2}{C}}.
\end{equation}
The goal is to improve this result and get precise asymptotics for the moments of $\Delta G(z)$, implying that $\Delta G(z)$ becomes Gaussian as $\N\rightarrow\infty$. To this end it is crucial to assume that the ensemble is off-critical. This is made precise in Assumption~\ref{Assumptions_offcrit} in Section~\ref{Chapter_Setup_and_Examples} and roughly means that the density of the equilibrium measure has a square-root behavior at endpoints of bands, that these endpoints do not touch the endpoints of the segments, and that the inequalities for the effective potential are strict away from the bands. The latter condition is important when showing that off-criticality is an open condition, \textit{i.e.} that it is preserved under small perturbations of the parameters of the ensemble (Theorem~\ref{Theorem_off_critical_neighborhood}).  For the hexagon with a hole, Proposition~\ref{Proposition_hex_hole_offcrit} and Lemma~\ref{Lemma_symmetric_hex} describe a range of off-critical parameters.

\begin{theorem}[Central limit theorem for fixed filling fractions] \label{Theorem_CLT_G_Intro} For the ensemble \eqref{eq_Hahn_cut_2} with fixed filling fractions and under the off-criticality assumption, there exist holomorphic functions $W_{1}^{[1]}(z)$ and $W_1^{[2]}(z)$ of $z \in \amsmathbb{C} \setminus \amsmathbb{A}$ such that, as $\N\rightarrow\infty$ uniformly over $z$ away from $\amsmathbb{A}$
\begin{equation}
   \label{exp_1}\E\big[G(z)\big] = \N \mathcal{G}(z)  + W_1^{[1]}(z) + \frac{W_1^{[2]}(z)}{\N} + o\bigg(\frac{1}{\N}\bigg).
\end{equation}
Moreover, there exists a holomorphic function $\mathcal{F}(z_1,z_2)$ of $(z_1,z_2) \in (\amsmathbb{C} \setminus \amsmathbb{A})^2$ such that for any $n \geq 1$, as $\N \rightarrow \infty$ uniformly for $z_1,\ldots,z_n$ away from $\amsmathbb{A}$
\begin{equation}
\label{exp_2}\E\bigg[ \prod_{i=1}^n \big(G(z_i)-\E[G(z_i)]\big)\bigg]= \E\bigg[ \prod_{i = 1}^{n} \textnormal{\textsf{Gau\ss{}}}(z_i) \bigg] +o(1),
\end{equation}
 where $\big(\textnormal{\textsf{Gau\ss{}}}(z_i)\big)_{i = 1}^{n}$ is a centered random Gaussian vector with covariance matrix $(\mathcal{F}(z_i,z_j))_{i,j = 1}^{n}$.
\end{theorem}

For the hexagon with the hole \eqref{eq_Hahn_cut_2}, the asymptotic covariance $\mathcal{F}(z_1,z_2)$ is completely explicit \cite[Theorem 7.1]{BGG}. It depends only on the endpoints of two bands $(\alpha_1,\beta_1)$, $(\alpha_2,\beta_2)$ of the equilibrium measure, exactly matches the covariance of continuous $\sbeta$ ensembles having the two same bands, and can be related to a section of the Gaussian free field in an annulus, see \cite[Section 4.4]{BuGo3} and Section~\ref{sec:Oshaped}. Theorem~\ref{Theorem_CLT_G_Intro} admits a direct extension to general discrete ensembles \eqref{eq_measure_general_intro}  with fixed filling fractions under the off-criticality assumption: these are Theorems~\ref{Theorem_correlators_expansion} and \ref{Theorem_correlators_expansion_relaxed}. In general, the covariance depends on the endpoints of the bands and on the interaction matrix $\boldsymbol{\Theta}$, but the dependence can be quite complicated.

In continuous $(\sbeta = 2)$-ensembles the term $W_1^{[1]}(z)$ in the expansion \eqref{exp_1} --- responsible for the mean in the Gaussian --- actually vanishes. In the formulation of Theorem~\ref{Theorem_CLT_G_Intro} this term is non-zero even for $\theta = 1$ as it was observed in \cite[Remark 2.7]{BGG}.  But, we discover that it disappears if we use the shifted parameters $\hat{a}'_h,\hat{b}'_h$ of Definition~\ref{def:eq_shifted_parameters} instead of $\hat{a}_h,\hat{b}_h$ to define the equilibrium measure $\mu$ and the deterministic term $\mathcal{G}(z)$. This is shown in Corollary~\ref{Theorem_correlators_expansion_theta1} and Theorem~\ref{Theorem_correlators_expansion_relaxed_theta1} for general discrete ensembles when the intensity of repulsion within each group is $1$, that is when all diagonal entries of $\boldsymbol{\Theta}$ are equal to $1$. This condition always holds for tiling models.

\medskip

Let us outline the main ideas in the proof of Theorem~\ref{Theorem_CLT_G_Intro}. The first step is to expand the right-hand side of \eqref{Nekrasov_eqn_int}. We start with
\[
\Phi^{\pm}(z) = \phi^{\pm}(z) + \frac{\phi^{\pm,[1]}(z)}{\N} + \frac{\phi^{\pm,[2]}(z)}{\N^2} + \frac{\mathbbm{f}^{\pm}(z)}{\N^3},
\]
where $\phi^{\pm}(z)$ are asymptotic versions of $\Phi^\pm(z)$ given by
\begin{equation*}
\begin{split}
\phi^+(z) & = \big(\hat{\mathfrak{t}}-\hat{C} +z\big) \cdot \big( \hat A+\hat B+\hat C-
   t-z\big) \cdot \big(\hat{\mathfrak{h}}-z\big)^2, \\ \phi^-(z) & =  z  \cdot \big(\hat A +\hat C -z\big) \cdot \big(\hat{\mathfrak{h}}+\hat{D}-z\big)^2.
\end{split}
\end{equation*}
Again, we shall use in the text a slightly different definition for $\phi^\pm(z)$ (Definition~\ref{Definition_phi_functions}) containing shifts of order $\frac{1}{\N}$, which can be omitted for this introductory discussion. The functions $\phi^{\pm}(z)$ are linked with the potential \eqref{eq_potential_Hahn_intro} via
\begin{equation}
\label{eq_ratio_to_potential_intro} \forall x \in\amsmathbb{A} \qquad
\frac{\phi^{+}(x)}{\phi^{-}(x)} = \exp\big(-\partial_x V(x)\big).
\end{equation}
The asymptotic expansion of the Nekrasov equations is based on two functions:
\begin{equation}\label{eq_defq_intro}
 q^{\pm}(z)   :=  \phi^{-}(z) \cdot \exp\big(-\mathcal{G}(z)\big) \pm \phi^{+}(z) \cdot \exp\big( \mathcal{G}(z)\big).
\end{equation}
To get the series expansion of $R(z)$ from \eqref{Nekrasov_eqn_int}, we first observe that for any $\tau \in \{\pm 1\}$, we have (see the proof of Proposition~\ref{Proposition_Nek_1_asymptotic_form} for computational details)
\begin{equation}
\label{eq_x307}
 \prod_{i = 1}^N \bigg(1 + \frac{\tau}{\N z - \ell_i - \frac{\tau}{2}}\bigg)  = \exp\bigg[\frac{\tau G(z)}{\N} + \frac{\tau\partial_{z}^2 {G}(z)}{24 \N^3}  +
O\bigg(\frac{1}{\N^3}\bigg)\bigg],
\end{equation}
where $z$ is kept away from $\amsmathbb{A}$ and it should be remembered that $G(z)$ is of order $\N$. Taking $t = \N^{-\frac{1}{2}} \log\N$ in \eqref{eq_apriory_bound_intro}  we learn that  $\Delta G(z)=G(z)-\N \mathcal{G}(z)$ and $\partial_z \Delta G(z)$ are of order at most $\N^{\frac{1}{2}} \log \N$ with overwhelming probability. Therefore, using \eqref{eq_x307}, the right-hand side of \eqref{Nekrasov_eqn_int} expands into
\begin{equation}
 \label{eq_nek0} R(z)=q^{+}(z)- q^{-}(z) \cdot \frac{\E[\Delta G (z)]}{\N} +\frac{E(z)}{\N} + \frac{\mathbbm{c}(z)}{\N},
\end{equation}
where we set
\[
E(z) = \phi^{-,[1]}(z) \cdot e^{-\mathcal{G} (z)} + \phi^{+,[1]}(z) \cdot e^{\mathcal{G}(z)},
\]
and $\mathbbm{c}(z)$ is a reminder term upper-bounded by
\begin{equation}
\label{eq_error_bound_intro}
\big|\mathbbm{c}(z)\big| \leq  C \cdot \bigg(\frac{\E\big[|\Delta G(z)|^2\big]}{\N}+\frac{\E\big[|\partial_z \Delta G(z)|\big]}{\N} + o(1)\bigg)
\end{equation}
for some finite constant $C$, as long as $z$ is away from $\amsmathbb{A}$. The key observation is that the equation \eqref{eq_nek0} can be solved. Before explaining how, let us emphasize its structure:
\begin{itemize}
 \item The function  $\E\left[\Delta G (z)\right]$ is the unknown we would like to find;
 \item $q^\pm(z)$ and $E(z)$ are assumed to be known, the answer will be expressed in terms of them;
 \item The term $\mathbbm{c}(z)$ is expected to be negligible as $\N \rightarrow \infty$;
 \item $R(z)$ will not enter the answer, but the important information  is that $R(z)$ is a holomorphic function of $z$ in a complex neighborhood of $\amsmathbb{A}$.
\end{itemize}

The solvability of \eqref{eq_nek0} hinges on the structure of singularities of $q^{\pm}(z)$ and $R(z)$. We already know that $q^{+}(z)$ and $R(z)$ are entire functions. The off-criticality of the ensemble prescribes a decomposition
\begin{equation}\label{eq_off_crit_assumption_intro}
 q^{-}(z)= s(z) \cdot \sqrt{(z-\alpha_1)( z - \beta_1)(z-\alpha_2)(z - \beta_2)},
\end{equation}
where $(\alpha_h,\beta_h)$ is the band (here assumed unique) in $[\hat{a}_h,\hat{b}_h]$ of the equilibrium measure and $s(z)$ \emph{does not vanish} in the neighborhood of bands. For the hexagon with the hole, the presence of one band per segment and the decomposition \eqref{eq_off_crit_assumption_intro} can be verified for a wide range of parameters directly using
\[
(q^+(z))^2-(q^{-}(z))^2=4\,\phi^+(z) \cdot \phi^-(z)
\]
and the polynomiality of $q^+(z)$, see Section~\ref{Section_hex_hole} for the details. For the more general measures \eqref{eq_measure_general_intro}, the property of having one band per segment and the decomposition \eqref{eq_off_crit_assumption_intro} can be verified if the potentials are strictly convex (Proposition~\ref{prop:convexmueq}). Such convexity arguments do apply for tiling models (Lemma~\ref{Lemma_glued_tilings_assumptions}). Otherwise, it is a non-trivial assumption to make. We will see in Proposition~\ref{Proposition_density} that the density $\mu(x)$ of the equilibrium measure is closely related to $q^\pm(z)$ by
\begin{equation*}
 \forall x \in  \textnormal{bands} \qquad\tan\big(\pi \mu(x)\big)=\frac{ q^{-}(x^-) - q^{-}(x^+)}{2\ii \, q^+(x)},
\end{equation*}
where $q^{-}(x^\pm)$ are limits when approaching $x$ from the upper or lower half-planes. In other words, the above decomposition with non-vanishing $s$ means that  $\mu(x)=0 $ or $1$ outside bands, remains strictly between $0$ and $1$ in the interior of bands, and approaches those extreme values like a square-root near endpoints of the bands. For the analysis of general discrete ensembles the one-band condition can be relaxed as addressed in Section~\ref{Sec:nexdtflu}, but the remaining conditions are crucial and they form the off-criticality Assumption~\ref{Assumptions_offcrit} of Section~\ref{Chapter_Setup_and_Examples}.

We actually need an extension of \eqref{eq_off_crit_assumption_intro}: we would like to assume that $s(z)$ has no zeros in a neighborhood of the segments $\amsmathbb{A}$ rather than only in a neighborhood of the bands. This is not guaranteed by the off-criticality \textit{per se}, and in order to achieve this extension, we propose a \emph{localization} procedure in Chapter~\ref{Chapter_conditioning}. We observe that in the void regions there are no particles with probability exponentially close to $1$ (Theorem~\ref{Theorem_ldpsup}), while in saturated regions the particles are densely packed, \textit{i.e.} occupy every possible lattice site with probability exponentially close to $1$ (Theorem~\ref{Theorem_ldsaturated}).  This happens because of the strict inequalities for the effective potential $V_{\textnormal{eff}}$ that we require in the off-criticality assumption. This is further used to remove the majority of voids and saturations from the system and localize to new segments $\mathfrak{A}$ which are small neighborhoods of bands. After this procedure the functions $q^\pm(z)$ are asymptotically unchanged, and therefore the decomposition \eqref{eq_off_crit_assumption_intro} is preserved, but now we can safely assume that $\mathfrak{A}$ is inside the region where $s$ has no zeros, see Theorem~\ref{proposition_FFF_conditioning} and Theorem~\ref{proposition_fluct_conditioning} for further details.

The relevance of decomposition \eqref{eq_off_crit_assumption_intro} to solving equations of the kind \eqref{eq_nek0} is demonstrated by the following lemma.

\begin{lemma}[Solving Nekrasov equations with one band] \label{Lemma_invert_intro} Take four real parameters $\hat{a} <\alpha < \beta <\hat{b}$ and set $\sigma(z)=\sqrt{(z-\alpha)(z-\beta)}$. Assume that we have four functions $F(z)$, $A(z)$ (unknown), $E(z)$ and $s(z)$ (known) such that
\begin{equation}\label{eq_theeqa_intro}
s(z) \cdot \sigma(z) \cdot F(z)= E(z) +A(z),
\end{equation}
 where $F(z) =  O(\frac{1}{z^2})$ as $z \rightarrow \infty$, the functions $E(z),F(z)$ are holomorphic for $z \in \amsmathbb{C} \setminus [\hat{a},\hat{b}]$ while $s(z)$ and $A(z)$ are holomorphic in a neighborhood of $[\hat a,\hat b]$ and $s(z)$ has no zeros in this neighborhood. Then \eqref{eq_theeqa_intro} has a unique solution $F(z)$ given by
 \begin{equation}\label{eq_solution_intro}
F(z)=\Op\bigg[\frac{E}{s}\bigg](z) \qquad \textnormal{with}\quad \Op[f](z) := \frac{1}{\sigma(z)} \oint_{\gamma} \frac{\dd \zeta}{2\ii\pi} \frac{f(\zeta)}{z - \zeta},
\end{equation}
where the contour $\gamma$ surrounds $[\hat{a},\hat{b}]$ but leaves $z$ outside.
\end{lemma}
\begin{proof}
We rewrite \eqref{eq_theeqa_intro} as
 \[
 \sqrt{(z-\alpha)(z - \beta)} \cdot F(z)= \frac{E(z) + A(z)}{s(z)}.
 \]
 Then, by using Cauchy residue formula, taking a contour to be a small loop around $z$,  then deforming the contour into a new one surrounding $[\hat a,\hat b]$, and using that $\sqrt{(z-\alpha)(z - \beta)} \cdot F(z)$ is holomorphic outside $[\hat a,\hat b]$ and decays like $O(\frac{1}{z})$ as $z \rightarrow \infty$, we have
 \begin{equation*}
 \begin{split}
 \sqrt{(z-\alpha)(z - \beta)} \cdot F(z)&= \Res_{\zeta = z} \frac{\sigma(\zeta) F(\zeta) \dd \zeta}{\zeta - z} =-\oint_{\gamma} \frac{\dd \zeta}{2\ii\pi}\,\frac{E(\zeta) + A(\zeta)}{(\zeta - z)  \cdot s(\zeta)} \\
&=-\frac{1}{2 \pi \ii}\oint_{\gamma} \frac{\dd \zeta}{2\ii\pi} \frac{E(\zeta)}{(\zeta- z) \cdot s(\zeta)},
\end{split}
\end{equation*}
where in the last equality we used that $\frac{A(\zeta)}{s(\zeta)}$ is holomorphic in a neighborhood of $[\hat{a},\hat{b}]$ and therefore its contour integral vanishes.
\end{proof}
The equation \eqref{eq_nek0} we have to deal with resembles \eqref{eq_theeqa_intro} but has two important differences. First, \eqref{eq_theeqa_intro} is adapted to the one-band case, but we actually have two bands: $\sqrt{(z - \alpha)(z - \beta)}$ should be replaced with $\sqrt{(z-\alpha_1)(z-\beta_1)(z-\alpha_{2})(z-\beta_{2})}$. Second, there is the additional term $\frac{\mathbbm{c}(z)}{\N}$ in the right-hand side, which we expected to be negligible but not yet shown to be. In the two-band situation the following generalization of Lemma~\ref{Lemma_invert_intro} is proven by the same method.

\begin{lemma}[Solving Nekrasov equations, example with two bands] \label{Lemma_invert_two_cut_intro} Take eight real parameters satisfying
\[
\hat{a}_1 < \alpha_1 < \beta_1 < \hat{b}_1 < \hat{a}_2 < \alpha_2 < \beta_2 < \hat{a}_2,\qquad \amsmathbb{A} = [\hat{a}_1,\hat{b}_1] \cup [\hat{a}_2,\hat{b}_2]
\]
and introduce $\sigma(z)  =\sqrt{(z-\alpha_1)(z-\beta_1)(z-\alpha_2)(z-\beta_2)}$.  Assume that we have four functions $F(z)$, $A(z)$ (unknown), $E(z)$ and $s(z)$ (known) such that \eqref{eq_theeqa_intro} holds with $F(z) = O(\frac{1}{z^2})$ as $z \rightarrow \infty$, the functions $F(z),E(z)$ are holomorphic for $z \in \amsmathbb{C} \setminus \amsmathbb{A}$ while $s(z),A(z)$ are holomorphic in a neighborhood of $\amsmathbb{A}$ and $s(z)$ has no zeros in this neighborhood. Then, there is a unique solution $F(z) = \Op\big[\frac{E}{s}\big](z)$  satisfying the additional condition
\begin{equation}
\label{vanishintcond}
\forall h \in \{1,2\}\qquad \oint_{\gamma_h} F(z)\dd z=0,
\end{equation}
where $\gamma_h$ is a contour surrounding $[\hat{a}_h,\hat{b}_h]$.
\end{lemma}

The condition \eqref{vanishintcond} is new compared to Lemma~\ref{Lemma_invert_intro} and necessary for the uniqueness. Generally speaking, for the ensemble with $K$ bands we need $K$ conditions corresponding to integrals on contours surrounding each band to guarantee uniqueness. Lemma~\ref{Lemma_invert_intro} implicitly had one such condition, which was implied by $F(z)=O(\frac{1}{z^2})$.  When the unknown is $F(z) = \Delta G(z)$, this condition holds because empirical measures have the same filling fractions as the equilibrium measure $\mu$. If we were not fixing filling fractions deterministically, \eqref{vanishintcond} would no longer be true, and in fact the resulting indeterminacy hides a more complicated $\N$-dependence of $\Delta G(z)$ as exemplified by Theorem~\ref{Theorem_Gauss_intro}. The solution operator $\Op$ for Lemma~\ref{Lemma_invert_two_cut_intro} still has an explicit formula, \textit{cf.} Theorem~\ref{Theorem_Masterspecial}. However, for general discrete ensembles \eqref{eq_measure_general_intro} the relevant operator does not have any particularly nice formula, and we develop in Chapter~\ref{Chapter_SolvingN} an abstract theory of solutions based on Fredholm theory to establish its existence and its properties, in particular its continuity in suitable operator norms, building on ideas of \cite{BGK}.

Switching to the error term $\mathbbm{c}(z)$, the bound \eqref{eq_error_bound_intro} together with \eqref{eq_apriory_bound_intro} implies that $\mathbbm{c}(z)$ is of order at most $\log^2(\N)$, with the largest contribution coming from $\E[|\Delta G(z)|^2]$. This is not yet sufficient to get the asymptotics of $\E[\Delta G(z)]$ through \eqref{eq_nek0}: we want $\mathbbm{c}(z) = o(1)$ so that it does not contribute to the $\N\rightarrow\infty$ limit. However, Lemma~\ref{Lemma_invert_two_cut_intro} applied to \eqref{eq_nek0} together with continuity of the linear operator $\Op$ yields $\E[\Delta G(z)]=O(\log^2(\N))$, thus improving the previous estimate $\E[\Delta G(z)]=O\big(\N^{\frac{1}{2}}\log(\N)\big)$ we had in \eqref{eq_apriory_bound_intro}. In order to move forward, we also need to improve the bound on  $\E[\Delta G(z)^2]$. This requires studying the variance of $\Delta G(z)$, which will be anyway helpful as we eventually want to prove \eqref{exp_2}.

\medskip

For this purpose, we look for additional equations for the covariance by writing the Nekrasov equations for varying weights. This idea leads to the \emph{higher-order Nekrasov equations} of Corollary~\ref{Corollary_higherNek}. In more detail, we consider the measures $\amsmathbb{P}_\N^{(t)}$ and associated expectation value $\amsmathbb{E}^{(t)}$ constructed with the new weights
\[
w^{(t)}(\ell)
=w(\ell) \cdot \bigg(1 + \frac{t}{z_2 - \tfrac{\ell}{\N}}\bigg)	,
\]
and observe that for any random variable $\chi$
\[
\partial_t \E^{(t)}[\chi]|_{t=0}= \E\big[\chi \cdot (G(z_2)-\E[G(z_2)])\big].
\]
On the other hand we can write Nekrasov equations for these new measures and observe that this changes the functions $\Phi^\pm$ into
\[
\Phi^{\pm,(t)}(z) = \Phi^\pm(z) \cdot \bigg(1+\frac{t}{z_2- z \mp\frac{1}{2\N}}\bigg).
\]
Computing the $t$-derivative at $t = 0$  of the version of \eqref{eq_nek0} for the measures $\amsmathbb{P}_\N^{(t)}$, we get an equation similar to \eqref{eq_nek0}, but this time involving $\E[\Delta G(z) \cdot (G(z_2)-\E[G(z_2)])]$ instead of $\E[\Delta G(z)]$. The new equation is the $\N\rightarrow\infty$ expansion of the second-order Nekrasov equation. It will have an error term $\frac{\mathbbm{c}(z,z_2)}{\N}$ involving the third moments of $\Delta G$. We upper-bound the third moment using \eqref{eq_apriory_bound_intro} and use Lemma~\ref{Lemma_invert_two_cut_intro} to solve the equation and get the improved upper bound on the second moment of $\Delta G(z)$. This improved upper bound is sufficient to conclude $\mathbbm{c}(z)=o(1)$. Hence, the continuity of $\Op$ and Lemma~\ref{Lemma_invert_two_cut_intro} allows solving \eqref{eq_nek0} for $\amsmathbb{E}[\Delta G(z)]$ up to $o(1)$ and subsequently finish the computation of $W_1^{[1]}(z)$ in \eqref{exp_1}.

The next step in the proof of Theorem~\ref{Theorem_CLT_G_Intro} is to find that the covariance has a limit $\mathcal{F}(z_i,z_j)$. Given the technology which we have already developed, it would be sufficient to show that the remainder term $\frac{\mathbbm{c}_2(z,z_2)}{\N}$ can be neglected. This is done by iterating the same scheme of proof: we develop higher-order Nekrasov equations, feed into them \eqref{eq_apriory_bound_intro} as an initial bound and iterate until we get the optimal estimates. This is the scheme of self-improving estimates mentioned in Section~\ref{sec_Nek_intro} and originally introduced in \cite[Section 4.3 and Figure 2]{BG11}. The details will be carried out for general discrete ensembles in Sections~\ref{Section_asymptotic_Nekrasov}-\ref{secprosofggg}. Eventually, once all the estimates are done, the application of Lemma~\ref{Lemma_invert_two_cut_intro} to higher-order Nekrasov equations yields Theorem~\ref{Theorem_CLT_G_Intro}.

\medskip

The result of Theorem~\ref{Theorem_CLT_G_Intro} gives a central limit theorem for the empirical Stieltjes transforms. It can readily be used it to get the result for arbitrary linear statistics of the form \eqref{fluctlin} with analytic test-functions $f$, via the identity:
\begin{equation}
\label{eq_linear_statistic_integral_intro}
 \sum_{i=1}^N f(\lambda_i)=  \oint_{\gamma} \frac{\dd z}{2\ii\pi} G(z) f(z),
\end{equation}
where the contour $\gamma$ surrounds all connected components of $\amsmathbb{A}$. This is how we achieve the central limit theorem for linear statistics in Corollaries~\ref{Corollary_CLT} and \ref{Corollary_CLT_relaxed}. In particular, these central limit theorems can be used to prove Theorem~\ref{Theorem_C_shape_intro} in Section~\ref{Section_tiling_C} and the one-band case of Theorem~\ref{Theorem_Gauss_intro} in Section~\ref{dis_Gauss_intro}.

\subsection{Expansion of the partition function for fixed filling fractions}
\label{intro_interp}
The next step in our investigations is the analysis of the asymptotic expansion of the free energy, \textit{i.e.} the logarithm of the partition function, as $\N \rightarrow \infty$. We keep on illustrating the general methods on the example of the tiling of the hexagon with a hole with fixed filling fractions $N_1$ and $N_2$. Our goal is to expand $\log \mathscr{Z}_N$ up to $o(1)$, and show that the coefficients in this expansion are smooth functions of the normalized filling fractions $\hat{n}_1 = \frac{N_1}{\N}$ and $\hat{n}_2 = \frac{N_2}{\N}$. The result we prove takes the following form.

\begin{theorem}[Asymptotics of the free energy]
\label{Theorem_partition_hex_intro}
 Consider an off-critical ensemble \eqref{eq_Hahn_cut_2} with fixed filling fractions. For any $\eps > 0$ and $\N$ larger than some $\N_0(\eps)$, we have
 \begin{equation}
 \label{eq_partition_hex_intro}\Big|\log\Z_\N - \big(\I[\mu_1,\mu_2]\,\N^2 + (\hat{n}_1+\hat{n}_2)\N\log \N + \mathbbm{Rest}_1\,\N\big)\Big| < \eps,
 \end{equation}
 where $(\mu_1,\mu_2)$ is the equilibrium measure and  $\I[\mu_1,\mu_2]$ is the maximum of \eqref{eq_I_intro}, while $\mathbbm{Rest}_1$ is a twice differentiable function of the rescaled parameters $\hat{A},\hat{B},\hat{C},\hat{D},\hat{\mathfrak{t}},\hat{\mathfrak{h}},\hat{n}_1,\hat{n}_2$, which is bounded with $\N$ and whose first and second-order partial derivatives are bounded.
 \end{theorem}
The precise statement can be found in Theorems~\ref{Theorem_partition_one_band} and \ref{Theorem_partition_multicut} for general discrete ensembles. All these theorems are new. Following the proofs it is possible to collect terms and give a certain expression for the leading order of $\mathbbm{Rest}_1$, but it is rather involved and we have not tried to simplify it. Finding a conceptually simple expression for $\mathbbm{Rest}_1$ is left as an outstanding open problem. In contrast, for continuous $\sbeta$-ensemble there is a nice formula involving $(1 - \frac{\sbeta}{2})$ times the von Neumann entropy of the equilibrium measure, whose heuristic derivation goes back to Dyson \cite{DysonII} and proof can be found in \cite[Theorem 1]{Shc2} or \cite[Theorem 4]{BG_multicut}. Even without an explicit expression for  $\mathbbm{Rest}_1$, Theorem~\ref{Theorem_partition_hex_intro} has important probabilistic consequences explained in Section~\ref{Sec:nexdtflu}.

\medskip

The driving idea for the proof of Theorem~\ref{Theorem_partition_hex_intro} is to interpolate between the target ensemble (\eqref{eq_Hahn_cut_2} in our case) and a reference ensemble for which the partition function can be explicitly computed. The change in the partition function along the interpolation can be controlled by exploiting Theorem~\ref{Theorem_CLT_G_Intro} and its generalizations. The exactly solvable models usually have one band only and so we first restrict to target ensembles with one band. To estimate their partition function, we interpolate with a $zw$-ensemble discussed in Section~\ref{zwintro} and then use Theorem~\ref{Theorem_zw_intro}. In a second step, for target ensembles with having one band in each of their $H$ segments (we have $H=2$ in Theorem~\ref{Theorem_partition_hex_intro}), we interpolate with a collection of $H$ independent one-band discrete ensembles.

\medskip

For the first part, in Section~\ref{Section_Partition_onecut} we outline the interpolation between two general measures \eqref{eq_measure_general_intro}, which have the same interaction matrix $\boldsymbol \Theta$ but different weights $w(\ell)$. In order to illustrate the approach,  consider a toy partition function:
\begin{equation}
\label{Znzdeva}
 \Z_{N}^{(t)} = \sum_{\boldsymbol{\ell}} \prod_{1 \leq i<j \leq N} (\ell_j-\ell_i)^2 \cdot \prod_{i=1}^N w^{(t)}(\ell_i) \qquad t \in \{0,1\}.
\end{equation}
We would like to compute $\log(\Z_N^{(1)}) - \log(\Z_N^{(0)})$. It can be written
\[
\log(\Z_N^{(1)}) - \log(\Z_N^{(0)}) = \int_0^1 \big(\partial_t \log(\Z_{N}^{(t)})\big) \dd t,
\]
where the definition of $\Z_N^{(t)}$ of \eqref{Znzdeva} is extended to arbitrary $t \in [0,1]$ by
\[
w^{(t)}(\ell) = \big(w^{(0)}(\ell)\big)^{1 - t} \cdot (w^{(1)}(\ell)\big)^{t}.
\]
The key observation is that we can write the integrand as
\[
\partial_t \log(\Z_{N}^{(t)})= \dfrac{ \sum\limits_{\boldsymbol{\ell}}\prod\limits_{i<j} (\ell_j-\ell_i)^2 \cdot \prod\limits_{i=1}^N w^{(t)}(\ell_i)  \cdot \sum_{i=1}^N \log\left(\frac{w^{(1)}(\ell_i)}{w^{(0)}(\ell_i)}\right) }{\Z_N^{(t)}} =\E^{(t)}\Bigg[\sum_{i=1}^N  \log\bigg(\frac{w^{(1)}(\ell_i)}{w^{(0)}(\ell_i)}\bigg)\Bigg],
\]
where the expectation is taken with respect to the discrete ensemble whose partition function $\Z_N^{(t)}$ computes. The right-hand side has the form of a linear statistic, and therefore we can compute its asymptotic via Theorem~\ref{Theorem_CLT_G_Intro} and combining \eqref{exp_1} with \eqref{eq_linear_statistic_integral_intro}.

There are two technical details in this scheme. First, we need to make sure that $\log\big(\frac{w^{(1)}(\N z)}{w^{(0)}(\N z)}\big)$ is analytic for $z$ in the neighborhood of the segments $[\hat{a}_h,\hat{b}_h]$, as otherwise we will not be able to use \eqref{eq_linear_statistic_integral_intro}. In particular, this requires coincidence of the zeros of $w^{(0)}$ and $w^{(1)}$. Looking at \eqref{eq_Hahn}, \eqref{weight:holey}, \eqref{eq_zw_measure_from_O}, we observe that the weights of our interest do have zeros, and therefore we should be cautious about their locations. Second, in order to use (a generalization of) Theorem~\ref{Theorem_CLT_G_Intro}, we need to know off-criticality. In general, even if the equilibrium measures $\mu^{(0)}$ and $\mu^{(1)}$ are both off-critical, the equilibrium measure $\mu^{(t)}$ might fail to be so. However, if we additionally assume that $\mu^{(0)}$ and $\mu^{(1)}$ have exactly the same bands, voids, and saturations, then the characterization of the equilibrium measure can be used to show the identity $\mu^{(t)}=(1-t)\mu^{(0)}+ t\mu^{(1)}$ making off-criticality for all $t \in [0,1]$ obvious. For this reason, when we interpolate towards a $zw$-ensemble in Section~\ref{Section_Partition_onecut}, we use a fine tuning of their parameters allowing to match the bands, voids, and saturations of the target measure. We justify that such a tuning is always possible in Section~\ref{sec:Tuning}.

\medskip

For the second part, in Section~\ref{Section_Partition_multicut} we interpolate between a discrete ensemble \eqref{eq_measure_general_intro} with an arbitrary matrix $\boldsymbol{\Theta}$ and a discrete ensemble with a diagonal matrix $\boldsymbol{\Theta}$, simply dimming off-diagonal entries. In the latter, the particles in different segments $[a_h,b_h]$ are independent and the partition function factorizes into a product of $H$ one-segment partition functions. We illustrate the approach by dealing with a toy partition function having particles split in two groups $\ell_1,\dots,\ell_{N_1}$ and $\ell_{N_1+1},\dots,\ell_{N}$
\begin{equation*}
\begin{split}
  \Z_{N} & = \sum_{\boldsymbol{\ell}} \prod_{1\leq i<j\leq N} (\ell_j-\ell_i)^2 \cdot \prod_{i=1}^N w(\ell_i) \\
  & =\sum_{\boldsymbol{\ell}} \prod_{1\leq i<j\leq N_1} (\ell_j-\ell_i)^2 \cdot \prod_{N_1< i<j\leq N} (\ell_j-\ell_i)^2  \cdot \prod_{\substack{1 \leq i \leq N_1 \\ N_1 < j \leq N}}(\ell_j-\ell_i)^2 \cdot \prod_{i=1}^{N_1} w(\ell_i)  \cdot \!\!\prod_{i=N_1+1}^{N} w(\ell_i).
\end{split}
\end{equation*}
For the interpolation argument, we introduce $t$-dependent discrete ensembles with the corresponding partition function given by
\begin{equation*}
\begin{split}
\Z_{N}^{(t)} & =   \sum_{\boldsymbol{\ell}} \prod_{1\leq i<j\leq N_1} (\ell_j-\ell_i)^2 \cdot  \prod_{N_1< i<j\leq N} (\ell_j-\ell_i)^2 \cdot \prod_{\substack{1 \leq i \leq N_1 \\ N_1 < j \leq N}} \frac{\Gamma\big(\ell_j - \ell_i + 1\big)\cdot \Gamma\big(\ell_j - \ell_i + t\big)}{\Gamma\big(\ell_j - \ell_i\big) \cdot \Gamma\big(\ell_j  - \ell_i + 1 - t\big)} \\ & \qquad\qquad\times \prod_{i=1}^{N_1} w_1^{(t)}(\ell_i)  \prod_{i=N_1+1}^{N} w_2^{(t)}(\ell_i).
 \end{split}
\end{equation*}
This is a discrete ensemble with interaction matrix
\[
\boldsymbol{\Theta}^{(t)} = \left(\begin{array}{cc} 1 & t\\ t & 1\end{array}\right).
\]
For $t = 1$ we take $w_1^{(t = 1)} = w_2^{(t = 1)} = w$ as weights, so that $\Z_N^{(1)} = \Z_N$. For arbitrary $t \in [0,1)$ we choose weights $w^{(t)}_1$ and $w^{(t)}_2$ including the mean-field effect of the ratio of four Gamma functions --- see the proof of Theorem~\ref{Theorem_partition_multicut} for the exact formulae --- so that the equilibrium measure does not depend on $t$. This ensures off-criticality all the way through the interpolation. The rest of the argument is very similar to the first step. We write
\begin{equation*}
\begin{split}
 \log\bigg(\frac{\Z_{N}^{(1)}}{\Z_{N}^{(0)}}\bigg) & =\int_0^1 \big(\partial_t \log(\Z_{N}^{(t)})\big) \dd t \\
 & = \int_{0}^{1} \amsmathbb{E}^{(t)}\bigg[\sum_{1 \leq i \leq N_1}\sum_{N_1 < j \leq N} \big( (\log \Gamma)'(\ell_j - \ell_i + t) + (\log \Gamma)'(\ell_j - \ell_i + 1 -t)\big)\bigg]\dd t \\
 & = \int_{0}^{1} \amsmathbb{E}^{(t)}\bigg[ \sum_{1 \leq i \leq N_1}\sum_{N_1 < j \leq N} 2\,\ln(\ell_j - \ell_i) + \cdots\bigg]\dd t,
 \end{split}
\end{equation*}
where the $\cdots$ indicate the remainder in Stirling approximation, which only gives a $O(1)$ contribution to the result. Since the two segments to which $\ell_i$ (for $1 \leq i \leq N_1$) and $\ell_j$ (for $N_1 < j \leq N_2$) belong are separated, the function $(x,y) \mapsto \ln(y - x)$ extends to a holomorphic function for $x$ in a neighborhood of $[\hat{a}_1,\hat{b}_1]$ and $y$ in a neighborhood of $[\hat{a}_2,\hat{b}_2]$. Then, we can use the same trick as in \eqref{eq_linear_statistic_integral_intro} but with double contour integral: the expectation $\amsmathbb{E}\big[\sum_{i,j} \ln(\ell_j-\ell_i)\big]$ can be reduced to the expectation of a bilinear functional of $G(z_1)$ and $G(z_2)$, and eventually to the covariance of $\Delta G(z_1)$ and $\Delta G(z_2)$. The asymptotics of the latter can be understood by the general version of Theorem~\ref{Theorem_CLT_G_Intro}.

\medskip

The above sketch delivers the formula for $\log\Z_\N$ in \eqref{eq_partition_hex_intro} as a multiple integral of the asymptotic expansion in a general version of Theorem~\ref{Theorem_CLT_G_Intro}. The explicit terms $\I[\mu_1,\mu_2]\,\N^2$ and  $(\hat{n}_1+\hat{n}_2)\N\log\N$ appear in the direct simplifications in these integrals (for general discrete ensembles the coefficient of the $\N \ln \N$ term is rather $\sum_{h = 1}^{H} \theta_{h,h}\hat{n}_h$). For the lower-order terms corresponding to $ \mathbbm{Rest}_1\,\N$ such simplifications are yet to be found.

\subsection{Fluctuations of the filling fractions}
\label{Sec:nexdtflu}
In this section we outline the proof of Theorem~\ref{Theorem_hex_hole_intro} of Section~\ref{Section_tiling_hole} and its general versions presented in Chapter~\ref{Chapter_filling_fractions}. Our starting point is Theorem~\ref{Theorem_partition_hex_intro}, in which we would like to make the dependence on $N_1$ and $N_2$ more explicit. Since $N = N_1 + N_2$ is always fixed, we only keep the dependence in $N_1$. Recall the notations $\hat{n}_h = \frac{N_h}{\N}$ and $\hat{n} = \frac{N}{\N}$. Denoting $\Z_\N^{N_1}$ the partition function with filling fractions $(N_1,N - N_1)$, we rewrite \eqref{eq_partition_hex_intro} as
\begin{equation}
 \label{eq_partition_hex_intro_2}\log\Z_\N^{N_1} = \I[\mu_1^{\hat{n}_1},\mu_2^{\hat{n}_1}]\,\N^2 + \hat{n}\,\N\log \N + \mathbbm{Rest}_1(\hat{n}_1)\,\N + o(1).
\end{equation}
The first three terms depend smoothly on $\hat{n}_1$. If we switch from the setting of deterministically fixed filling fractions (as in Theorem~\ref{Theorem_partition_hex_intro}) to the setting of fluctuating filling fractions (as in Theorem~\ref{Theorem_hex_hole_intro}), then the distribution of random $N_1$  becomes
\begin{equation}
\label{eq_law_ff_intro}
 \amsmathbb{P}_{\N}[N_1=k_1]= \frac{  \Z_\N^{k_1} }{\sum_{k=0}^N \Z_\N^{k}},
\end{equation}
Unlike \cite{BG_multicut}, there are no additional combinatorial factors in this formula because the particles are ordered on the real line.

Let $(\mu^*_1,\mu^*_2)$ denote the equilibrium measure for the ensemble with fluctuating filling fractions: it is the maximizer of $\I[\mu_1^{\hat{n}_1},\mu_2^{\hat{n}_1}]$ over $\hat{n}_1 \in [0,\hat{n}]$. We would like to Taylor expand \eqref{eq_partition_hex_intro_2} near the optimal value $\hat{n}_1^*$ and then plug the result into \eqref{eq_law_ff_intro}. Since $\I[\mu_1^{\hat{n}_1},\mu_2^{\hat{n}_1}]$ reaches its maximum at $\hat{n}_1 = \hat{n}_1^*$, its first-order derivative vanish and we get
\begin{equation}
\begin{split}
 \label{eq_partition_hex_intro_3}
 \log\Z_\N^{N_1} & = \I[\mu_1^{\hat{n}_1^*},\mu_2^{\hat{n}_1^*}]\,\N^2 + \hat{n}\,\N\log \N  + \frac{1}{2} (N_1 - \N\hat{n}_1^*)^2 \cdot \big( \partial^2_{\hat{n}_1} \I[\mu_1^{\hat{n}_1},\mu_2^{\hat{n}_1}]\big)\big|_{\hat{n}_1 = \hat{n}_1^*} \\
 & \quad  + (N_1  - \N\hat{n}_1^*) \cdot \partial_{\hat{n}_1} \mathbbm{Rest}_1(\hat{n}_1)\big|_{\hat{n}_1 = \hat{n}_1^*}  + \cdots
\end{split}
\end{equation}
Thus, \eqref{eq_law_ff_intro} turns into
\begin{equation}
 \label{eq_partition_hex_intro_4} \amsmathbb{P}_{\N}[N_1=k_1] = \frac{1}{\mathscr{Q}} \exp\bigg( -\frac{B}{2} (k_1 - \N\hat{n}_1^* - u)^2 + \cdots \bigg),
 \end{equation}
 where
 \[
 B := -\partial^2_{\hat{n}_1} \I[\mu_1^{\hat{n}_1},\mu_2^{\hat{n}_1}]\big|_{\hat{n}_1 = \hat{n}_1^*}, \qquad u := B^{-1} \cdot \partial_{\hat{n}_1}\mathbbm{Rest}_1(\hat{n}_1)\big|_{\hat{n}_1 = \hat{n}_1^*},
\]
and $\mathscr{Q}$ is a normalization constant making \eqref{eq_partition_hex_intro_4} a probability distribution. An important technical ingredient is the law of large numbers for the filling fractions of \eqref{eq_ff_LLN_Intro}, which allows us  restricting the relevant values of $\frac{k_1}{\N}$ to a window of size $C \frac{\log^2\N}{\sqrt{\N}} = o(1)$ around $\hat{n}_1^*$ on an event of overwhelming probability, thus making the remainders $\cdots$ in the Taylor expansions negligible. In this regard, the asymptotic expansion \eqref{eq_partition_hex_intro_4} is enabled by the uniform smoothness of the expansion \eqref{eq_partition_hex_intro} with respect to $\hat{n}_1$ in a $\N$-independent neighborhood of $\hat{n}_1^*$. Finally, we need to know that $B > 0$, which is related to the minimizing property of $\hat{n}^*$. In general, $B$ is replaced with a quadratic form and Proposition~\ref{Proposition_Hessian_free_energy} shows that it is positive definite when restricted to the space of deformations of filling fractions obeying a  set of deterministic affine constraints prescribed in Equation $(\star)$ in Section~\ref{DataS} (for the previous example there is only one such constraint: $N_1 + N_2 = N$). Combining all these ingredients, we arrive at Theorem~\ref{Theorem_hex_hole_intro}.

\medskip

As we mentioned in Section~\ref{dis_Gauss_intro}, for the general discrete ensemble \eqref{eq_measure_general_intro} in the presence of saturations, the filling fractions should be replaced with more elaborate parameters --- extended filling fractions first introduced in Chapter~\ref{Chapter_conditioning} and analyzed in Section~\ref{Section_Fillingfractionsatcasesec}. The proof of their asymptotic discrete Gaussianity remains conceptually the same: we fix them deterministically, compute the asymptotics of the corresponding partition functions, and then plug into the general version of \eqref{eq_law_ff_intro}. The most challenging technical detail in this generalization is the careful bookkeeping for various lattices to which different particles belong, due to condition $\ell_{i+1}-\ell_i\in\{\theta_{h,h},\theta_{h,h}+1,\theta_{h,h}+2,\dots\}$ for consecutive particles in the $h$-th segment. Finally, the two-band case of Theorem~\ref{Theorem_Gauss_intro} and its analog for general discrete ensembles in Theorem~\ref{Theorem_linear_statistics_fluctuation_sat} are obtained by combining the asymptotics of the extended filling fractions with the asymptotic expansions for the linear statistics in the fixed filling fractions case explained in Section~\ref{Theorem_correlators_intro}.

\subsection{Two-dimensional fields of fluctuations for random lozenge tilings}

Let us outline the approach of Chapter~\ref{Chap11} for investigating the global fluctuations of the height function in random tilings and proving convergence to the two-dimensional Gaussian free field, as in Theorem~\ref{Theorem_C_shape_intro} of Section~\ref{Section_tiling_C}. In Section~\ref{Section_gluing_def} we present a class of polygonal domains on the triangular grid, which are obtained as gluings of trapezoids along a vertical line. A \emph{trapezoid} is a basic building block in our construction, see Section~\ref{Section_trapezoid} for detailed discussion. For instance, the hexagon in Figure~\ref{Fig_tiling_hex_222} can be obtained as a gluing of one left trapezoid and one right trapezoid along the line $x=\mathfrak{t}$. The same is true for the hexagon with a hole in Figure~\ref{Fig_hex_hole_small}, although the trapezoids become slightly different in order to create the hole. The domains of Figure~\ref{Fig:Cshapeint} are both obtained as gluings of three trapezoids along the dashed vertical line; adjusting the parameters of these trapezoids leads to different domains, as the comparison of two panels on this figure reveals.

The basic observation is that the exact enumeration for the number of lozenge tilings of a trapezoid (see Proposition~\ref{Proposition_number_tilings_trapezoid}) leads to the identification of the joint distribution of horizontal lozenges along the distinguished vertical line with an instance of the discrete ensemble \eqref{eq_measure_general_intro}. In this situation the matrix elements of $\boldsymbol{\Theta}$ are from $\{0,\tfrac{1}{2},1\}$ with $1$s on the diagonal. This is related to each trapezoid giving rise to a block submatrix filled with $\tfrac{1}{2}$, and the blocks are further added together, taking into account that each segment $[\hat a_h,\hat b_h]$ should eventually belong to exactly two trapezoids.

Hence, we can use all the machinery developed in this book and arrive at a generalization of Theorem~\ref{Theorem_CLT_G_Intro} to tilings of a wide class of polynomial domains. We still need to check off-criticality condition. But, due to some convexity properties of the potentials appearing in tilings, the possibility of criticality is even more restrained than in general discrete ensembles, see Lemma~\ref{Lemma_glued_tilings_assumptions}.

The analysis of discrete ensembles gives the central limit theorem for the fluctuations along the distinguished vertical section of the tiled domain for fixed filling fractions, and its perturbation by a discrete Gaussian component for fluctuating filling fractions. An additional step is required in order to extend the central limit theorem to the whole two-dimensional field of fluctuations of the height function in the domain. Here we use the approach first introduced in \cite{BuGo3} in the study of the hexagon with a hole. The idea is that for trapezoids, the theorems of \cite{BuGo3} can be used to show that once we know the macroscopic fluctuations along the boundary, we will know them everywhere inside as well. Therefore, using the fluctuations along the vertical lines as an input, we can extend them to the whole domain. A separate question is to identify the resulting field of the Gaussian fluctuations with the Gaussian free field to prove the Kenyon-Okounkov conjecture of \cite[Section 2.3]{kenyon2007limit} and \cite[Lectures 11-12]{Vadimlecture}. This is done by matching the structure of the covariance of trapezoid gluings with help of the well-known domain Markov property of the Gaussian free field. In \cite{BuGo3} this could be carried out on the basis of explicit formulae, which are not to be expected for the general domains analyzed in this book. Instead, we need to analyze more conceptually the interplay between the geometry of gluings of trapezoids and geometry of the spectral curve for the discrete ensembles \eqref{eq_measure_general_intro}, for which the tools of Chapter~\ref{Chapter_AG} become handy. This effort culminates with Theorems~\ref{Theorem_GFF} and \ref{Theorem_GFF_general} proving the Gaussian free field asymptotics.

Finally, we notice that the gluings of trapezoids do not have to be embedded into an orientable surface. In the non-orientable case, the main steps of our arguments remain the same, but the covariance identification at the end becomes different: we need to use orientable covers for the covariance of discrete ensembles in Corollary~\ref{cor:Greenfermion}. As a consequence, a symmetric conditioning of the Gaussian free field appears in Theorem~\ref{Theorem_GFF_general_nonor}. We are not aware of a previous occurrence of such symmetric Gaussian free fields in limits of ensembles of two-dimensional statistical mechanics. Dimer models on non-orientable surfaces have been studied in the past to access their partition function \cite{LuWuNonor,LuWuNonor2,Cimasoni,CimasoniPham}, but not to study fluctuations and their Gaussianity.

\subsection{Riemann--Hilbert problems, spectral curves and leading covariances}
\label{laius}

Unveiling the connection to Gaussian free fields for random tilings models requires geometric tools, which are at the center of Part~\ref{Part_Master_equation} and that we now outline. More generally, these tools are useful to carry out explicit computations of coefficients of the asymptotic expansions of correlators in continuous $\sbeta$-ensembles or in discrete ensembles. They lead to a characterization of the leading covariance (for linear statistics in the context of Theorems~\ref{Theorem_C_shape_intro}, \ref{Theorem_Gauss_intro}, \ref{Theorem_CLT_G_Intro} and the more general Theorems~\ref{Theorem_correlators_expansion} and \ref{Theorem_correlators_expansion_relaxed}) that is summarized in Theorem~\ref{thm_intro_FRHP} below. This characterization not only offers a way to get explicit formula for the leading covariance in some non-trivial cases, but also make it accessible for a structural analysis even when an explicit formula is lacking --- for general tiling models and beyond.

The construction of spectral curves plays the major role in this geometric approach. Before approaching the general situation for discrete models like \eqref{eq_measure_general_intro}, let us start by reviewing how it arises in the simplest situation, going back to the continuous $\sbeta$-ensembles mentioned in \eqref{beta-mod} of Section~\ref{Section_continuous} or a discrete ensemble with $H = 1$ group and intensity of interaction $\theta = \frac{\sbeta}{2}$. Here, the \emph{spectral curve} is a Riemann surface\footnote{Here 'curve' is meant in the sense of complex one-dimensional manifolds, which corresponds indeed to a Riemann surface.} $\Sigma$ on which the Stieltjes transform of the equilibrium measure
\[
\mathcal{G}(z) = \int \frac{\dd\mu(x)}{x - z}
\]
admits an analytic continuation. As we saw in \eqref{eq_bandsrhp}, in absence of saturations the minimizing property for the equilibrium measure implies a Riemann--Hilbert problem controlling the discontinuity of $\mathcal{G}(z)$ on the bands $\amsmathbb{B}$, in this case
\begin{equation}
\label{sdunasbs}\forall x \in \amsmathbb{B}\qquad \frac{\sbeta}{2}\big(\Gm(x^+) + \Gm(x^-)\big) = V'(x).
\end{equation}
The complex analysis methods to solve such equations date back to Sokhotski's 1873 doctoral dissertation \cite{Sokhotski} and to Plemelj \cite{Plemelj},  see \textit{e.g.} \cite{CoursEynard} for a recent implementation. The general solution in integral form can be found in \cite[Section 4.3]{Tricomi} or \cite{Musk}. For $V$ polynomial and in the regime where we have $K$ disjoint bands $(\alpha_1,\beta_1),\ldots,(\alpha_K,\beta_K)$, the solution of \eqref{sdunasbs} takes the form $\Gm(z) = \sbeta^{-1}V'(z) + s(z)\sigma(z)$ where $s(z)$ is a polynomial and
\begin{equation}
\label{spcurvehyperun} \sigma^2(z) = \prod_{k = 1}^{K} (z - \alpha_k)(z - \beta_k).
\end{equation}
Accordingly, $\Gm(z)$ is an algebraic function which can be analytically continued on a hyperelliptic Riemann surface of genus $K - 1$, and the latter is our spectral curve $\Sigma$. This Riemann surface is made by gluing two copies (also called \emph{sheets}) of the Riemann sphere along the bands, which are slit to form $K$ topological circles. Equation~\eqref{spcurvehyperun} is an algebraic presentation for this Riemann surface.  Moreover, the recursive solution of the Dyson--Schwinger equation of the continuous $\sbeta$-ensembles leading to \eqref{topexp2} can be recast solely in terms of the geometry of the spectral curve, giving rise to the theory of topological recursion by Chekhov, Eynard and Orantin \cite{E1MM,CE06,EORev,EBook}. For instance, the leading covariance of two Stieltjes transforms of the empirical measure, denoted $\mathcal{F}(z_1,z_2)$, satisfies
\begin{equation}
\label{FRHPS}
\forall x \in \amsmathbb{B} \qquad \frac{\sbeta}{2}\big(\mathcal{F}(x_1^+,z_2) + \mathcal{F}(x_1^-,z_2)\big) = - \frac{1}{(x_1 - z_2)^2}
\end{equation}
for any $z_2$ away from the bands. From there, one can show that $\frac{\sbeta}{2}\mathcal{F}(z_1,z_2)$ analytically continues to a meromorphic function on $\Sigma^2$ with a double pole on the diagonal having coefficient $1$. There are relatively explicit formulae for such objects on hyperelliptic Riemann surfaces, \textit{cf.} \cite{Ake96} or Section~\ref{Sec1241}, and it only depends on the endpoints of the bands. These properties for the leading covariance is a rather universal feature in many random matrix models, which is explained in physics by their relation with two-dimensional conformal field theory, see \textit{e.g.} \cite{KostovCFT}. In mathematical terms, it allows a neat description of the central limit theorem in the ensemble with fixed filling fractions. Indeed, after integration by parts, one can write
\begin{equation}
\label{covgree}
\textsf{Cov}[f,f] = \oint_{\gamma} \oint_{\gamma} \frac{\dd z_1 \dd z_2}{(2\ii\pi)^2} \mathcal{F}(z_1,z_2)\,f(z_1)f(z_2) = \frac{2}{\pi \sbeta} \int_{\amsmathbb{B}^2} \textnormal{Green}(x_1,x_2)\dd f(x_1) \dd f(x_2),
\end{equation}
where $\amsmathbb{B}$ are the bands and
\[
\textnormal{Green}(z_1,z_2) = - \frac{\sbeta/2}{4\pi} \int_{z_1^{*}}^{z_1}\int_{z_2^{*}}^{z_2} \mathcal{F}(\tilde{z}_1,\tilde{z}_2)\dd \tilde{z}_1 \dd \tilde{z}_2 - \frac{1}{2\pi} \log\left|\frac{z_1 - z_2}{z_1 - z_2^{*}}\right|.
\]
The analytic properties (in particular the double pole with coefficient $1$) of $\frac{\sbeta}{2}\mathcal{F}(z_1,z_2)$ lead to identify $\textnormal{Green}(z_1,z_2)$ with the \emph{Green function} on \emph{half of the spectral curve}. This half is a bordered Riemann surface which is best defined in the gluing presentation: we keep the upper half-plane in the first sheet, the lower half-plane in the second sheet, and on each sheet the real axis minus the bands is seen as the boundary. The Green function that we are talking about is the Green function for the Laplacian on this bordered surface with Dirichlet boundary conditions. In short, in the ensemble with fixed filling fractions, the field $\sqrt{\pi\sbeta/2} \cdot \Delta \textsf{Lin}[f]$ converges as $N \rightarrow \infty$ to the restriction to the bands of the derivative of the Gaussian free field on the half-spectral curve.

Now, let us return to discrete ensembles like \eqref{eq_measure_general_intro} (or continuous ensembles with several groups of particles). Then, we show that  \eqref{FRHPS} is replaced with a vector-valued Riemann--Hilbert problem for the $H$-tuple of Stieltjes transforms $(\mathcal{G}_{\mu_h}(z))_{h = 1}^{H}$ encoding the equilibrium measure of each segment. Namely for any $h \in [H]$ and $x$ in bands of the $h$-th segment, we have
\[
\theta_{h,h}\big(\Gm(x^+) + \Gm(x^-)\big) + \sum_{g \neq h} 2\theta_{h,g}\,\Gm_g(x) = V'_h(x).
\]
In the same spirit, the leading covariance appearing in the generalizations of Theorem~\ref{Theorem_CLT_G_Intro} satisfies the following generalization of \eqref{FRHPS}.
\begin{theorem}
\label{thm_intro_FRHP}
Consider an off-critical ensemble \eqref{eq_measure_general_intro} with fixed filling fractions. There exists functions $\mathcal{F}_{h_1,h_2}(z_1,z_2)$ indexed by $h_1,h_2 \in [H]$ depending only on $\boldsymbol{\Theta}$ and the endpoints of the bands, and having the following properties. First $\mathcal{F}_{h_1,h_2}(z_1,z_2) = \mathcal{F}_{h_2,h_1}(z_2,z_1)$. Second, for any $h_1,h_2 \in [H]$, any $x_1$ in bands of $[\hat{a}_{h_1},\hat{b}_{h_1}]$ and $z_2$ away from bands, we have
\[
\theta_{h_1,h_1}\big(\mathcal{F}_{h_1,h_2}(x_1^+,z_2) + \mathcal{F}_{h_1,h_2}(x_1^-,z_2)\big) + \sum_{g \neq h_1} 2\theta_{h_1,g}\,\mathcal{F}_{g,h_2}(x_1,z_2) = - \frac{\delta_{h_1,h_2}}{(x_1 - z_2)^2}.
\]
Third, for an analytic function $f$ near the bands, we have as $\N \rightarrow \infty$
\[
\amsmathbb{E}\big[(\Delta\textnormal{\textsf{Lin}}[f])^2\big] = \sum_{h_1,h_2 = 1}^{H} \oint_{\gamma_{h_1}}\oint_{\gamma_{h_2}} \frac{\dd z_1 \dd z_2}{(2\ii\pi)^2} \,\mathcal{F}_{h_1,h_2}(z_1,z_2)\,f(z_1)f(z_2) + o(1).
\]
where the contour $\gamma_h$ surrounds the $h$-th band.
\end{theorem}
The precise version of this statement is the combination of Theorem~\ref{Theorem_correlators_expansion_relaxed}, Theorem~\ref{thm:Bsym} and Proposition~\ref{thmBfund}. In Section~\ref{Section_Master_by_covariance} we explain that such Riemann--Hilbert problems closely relate to the solution operator $\Upsilon$ required to solve asymptotically Nekrasov equations (of which Lemma~\ref{Lemma_invert_intro} and \ref{Lemma_invert_two_cut_intro} gave the simplest instance). This is the origin of Theorem~\ref{thm_intro_FRHP}. Although these Riemann--Hilbert problems rarely afford explicit solutions, our point of view is that a lot can nevertheless be learned by analyzing their algebraic structure. To this end, Chapter~\ref{Chapter_AG} develops the construction of spectral curves on which the solutions of such Riemann--Hilbert problems admit an analytic continuation --- systematizing common arguments in the physics literature that were formalized in \cite{BESeifert}.

\medskip

\begin{figure}[t]
\begin{center}
\includegraphics[width=0.45\textwidth]{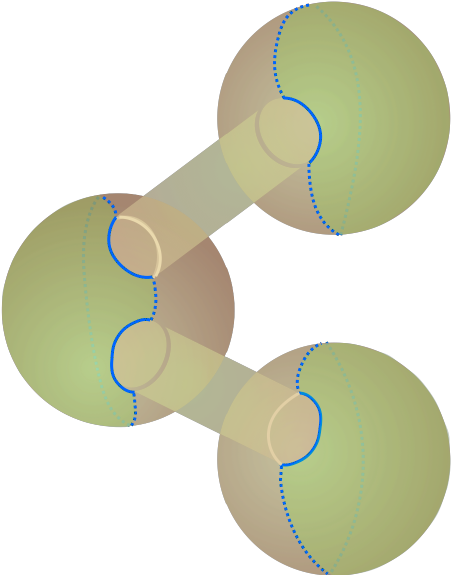}
\caption{\label{fig:Cgluefig} The spectral curve for the C-shaped domain is obtained by gluing together three Riemann spheres minus a disc (the sheets, conformally equivalent to Riemann spheres with slits). The half-spectral curve is obtained by keeping only the darker hemispheres, and is homemorphic to the C-shaped domain. More precise pictures describing the gluing will be provided in Chapter~\ref{Chapter_AG}.}
\end{center}
\end{figure}

Let us illustrate with the random lozenge tilings of the C-shaped domain of Section~\ref{Section_tiling_C} how to move forward with such results. The corresponding discrete ensemble has matrix of interaction \eqref{Cintmatr} of size $H = 2$. We have two empirical Stieltjes transforms
\[
\Delta G_1(z) = \sum_{i = 1}^{N_1} \frac{1}{z - \frac{\ell_i}{\N}} - \N\Gm_1(z) \qquad \textnormal{and}\qquad \Delta G_2(z) = \sum_{i = N_1 + 1}^{N} \frac{1}{z - \frac{\ell_i}{\N}} - \N\Gm_2(z),
\]
corresponding to each of the two groups of particles. Let $\amsmathbb{B}_h$ be the band in the $h$-th segment for $h \in \{1,2\}$. An easy consequence of Theorem~\ref{thm_intro_FRHP} is that for $\N \rightarrow \infty$ we have
\[
\forall h_1,h_2 \in \{1,2\}\qquad \amsmathbb{E}\big[\Delta G_{h_1}(z_1) \Delta G_{h_2}(z_2)\big] = \mathcal{F}_{h_1,h_2}(z_1,z_2) + o(1),
\]
and the leading covariance satisfies the Riemann--Hilbert problem
\begin{equation}
\label{RHPvecto2}
\begin{split}
\forall x_1 \in \amsmathbb{B}_1 & \quad \left\{\begin{array}{lll} \mathcal{F}_{1,1}(x_1^+,z_2) + \mathcal{F}_{1,1}(x_1^-,z_2) + \mathcal{F}_{1,2}(x_1,z_2) & = & - \dfrac{1}{(x_1 - z_2)^2}, \\[4pt] \mathcal{F}_{1,2}(x_1^+,z_2) + \mathcal{F}_{1,2}(x_1^-,z_2) + \mathcal{F}_{2,2}(x_1,z_2) & = & 0, \end{array}\right. \\[6pt]
\forall x_1 \in \amsmathbb{B}_2 & \quad \left\{\begin{array}{lll} \mathcal{F}_{2,1}(x_1^+,z_2) + \mathcal{F}_{2,1}(x_1^-,z_2) + \mathcal{F}_{1,1}(x_1,z_2) & = & 0, \\[4pt] \mathcal{F}_{2,2}(x_1^+,z_2) + \mathcal{F}_{2,2}(x_1^-,z_2) + \mathcal{F}_{1,2}(x_1,z_2) & = & - \dfrac{1}{(x_1 - z_2)^2}, \end{array}\right.
\end{split}
 \end{equation}
 and a symmetric one with respect to the variable $z_2$. To get our hands on these functions, we construct a spectral curve $\Sigma$ by gluing a central sheet to a second sheet along $\amsmathbb{B}_1$ and a third sheet along $\amsmathbb{B}_2$. The outcome is a Riemann sphere $\Sigma$ equipped with a meromorphic function $Z(\zeta)$ of degree $3$ giving the coordinate on each of the sheets (see Section~\ref{sec:Csha} for formulae). This construction mirrors the structure of the C-domain, which is the simply-connected domain obtained by gluing a central trapezoid to a second trapezoid along the bottom segment and to a third trapezoid along the top segment (see Figure~\ref{fig:Cgluefig} for an artistic view). The general procedure for the construction of spectral curves is explained in Section~\ref{sec:Algapproach}. Then, analyzing the Riemann--Hilbert problem \eqref{RHPvecto2} leads us to isolate the following property.

\begin{theorem}
\label{9table}
For the discrete ensemble of the C-shaped domain, the nine functions in the table
\begin{equation*}
\begin{array}{|c|c|c|}
\hline
(\mathcal{F}_{1,1} + \mathcal{F}_{2,1} + \mathcal{F}_{1,2} + \mathcal{F}_{2,2})(z_1,z_2) + \frac{1}{(z_1 - z_2)^2} & -(\mathcal{F}_{1,1} + \mathcal{F}_{1,2})(z_1,z_2) & - (\mathcal{F}_{2,1} + \mathcal{F}_{2,2})(z_1,z_2) \\[2pt]
\hline
-(\mathcal{F}_{1,1} + \mathcal{F}_{2,1})(z_1,z_2) & \mathcal{F}_{1,1}(z_1,z_2) + \frac{1}{(z_1 - z_2)^2} & \mathcal{F}_{2,1}(z_1,z_2) \\[2pt]
\hline
-(\mathcal{F}_{1,2} + \mathcal{F}_{2,2})(z_1,z_2) & \mathcal{F}_{1,2}(z_1,z_2) & \mathcal{F}_{2,2}(z_1,z_2) + \frac{1}{(z_1 - z_2)^2} \\[2pt]
\hline
\end{array}
\end{equation*}
are the nine branches of the meromorphic function on $\Sigma^2$
\begin{equation}
\label{0ZZ0}
\frac{1}{Z'(\zeta_1)Z'(\zeta_2)(\zeta_1 - \zeta_2)^2} = \partial_{z_1}\partial_{z_2}\log(\zeta_1 - \zeta_2)
\end{equation}
that we get when expressing it in terms of $z_1 = Z(\zeta_1)$ and $z_2 = Z(\zeta_2)$.
\end{theorem}
We explain in Section~\ref{sec:Csha} how to turn this compact statement into concrete ways of computing the leading covariance of $\textsf{Lin}_{h_1}[f_1]$ and $\textsf{Lin}_{h_2}[f_2]$ for analytic test functions $f_1,f_2$ and any choice of $h_1,h_2 \in \{1,2\}$. The logarithm in \eqref{0ZZ0} is reminiscent of the Green function, and indeed, we will combine the above description of $\mathcal{F}_{h_1,h_2}$ and the tools of \cite{BuGo2} to show that the fluctuations of the height field for the tilings of the C-shaped domain are asymptotically given by the Gaussian free field: this is a particular case of Theorem~\ref{Theorem_GFF_general}.

The key feature of Theorem~\ref{9table} responsible for the link with the Green function is the double pole with coefficient $1$ present only in the diagonal of the table. In Theorem~\ref{thm:Omegav} and Corollary~\ref{cor:Greennormal} we find the same phenomenon for all matrices $\boldsymbol{\Theta}$ coming from tiling models in orientable domains. Furthermore, we discover that the spectral curve is glued from sheets in the exact same way the tiled domain is glued from trapezoids. This is eventually responsible for the identification of the Green function defined from the spectral curve with the Green function appearing in Kenyon--Okounkov conjecture (Propositions~\ref{Proposition_uniformization_map} and \ref{Proposition_Green_function_explicit}). For tiling models of non-orientable domains, the situation is slightly different (Theorem~\ref{thm:spcurvenonbip}): the spectral curve is glued from sheets in the same way the orientation covering of the domain is glued from trapezoids, and double poles can also be found in certain off-diagonal entries in the analogue of Theorem~\ref{9table}. This time, the fine analysis of the pole structure reveals a link to the Green function symmetrized by the covering involution (Corollary~\ref{cor:Greenfermion}).

\section{Road map}
\label{Chap1Sec4}

\begin{figure}[t]
\begin{center}
\includegraphics[width=0.9\textwidth]{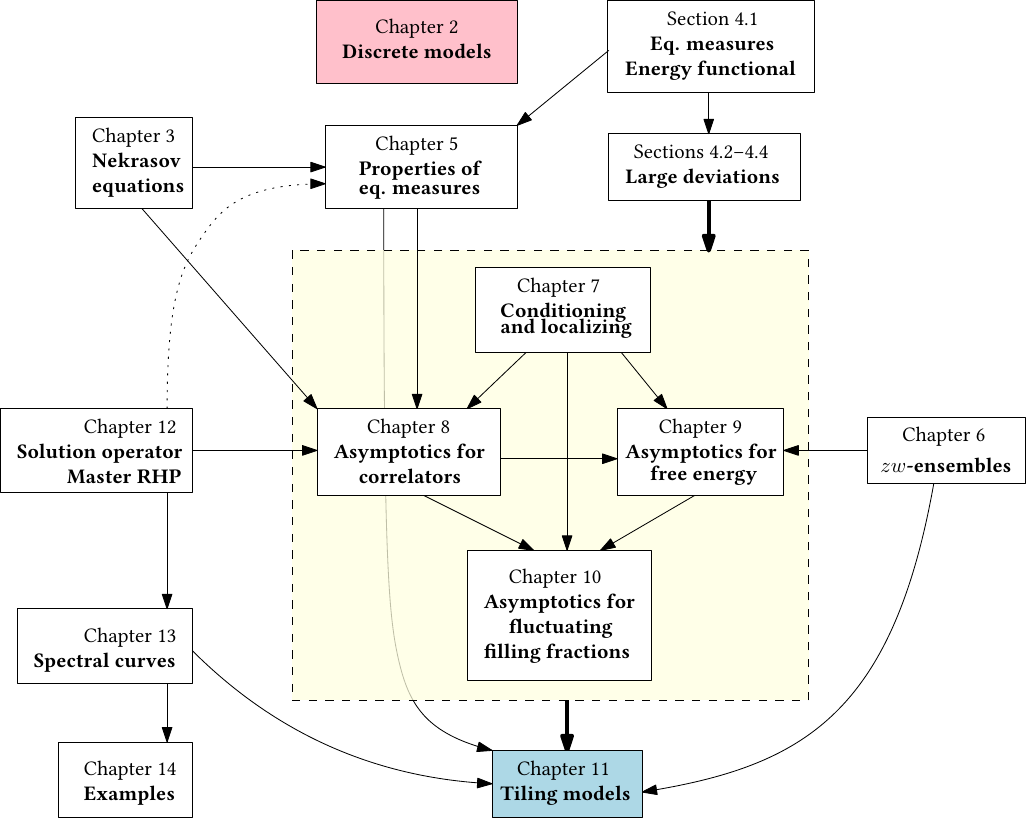}
\caption{\label{Fig:map} Dependencies of different chapters. The yellow box contains the core asymptotic results.}
\end{center}
\end{figure}

After this long introduction it is time to summarize the structure of this book (Figure \ref{Fig:map}). In Part~\ref{Part1}, we introduce our discrete ensembles and our assumptions in a precise way (Chapter~\ref{Chapter_Setup_and_Examples}), and state some of our main results. In Chapter~\ref{ChapterNekra} we present our first main tool: the Nekrasov equations. These equations will not only allow us to derive the asymptotics of the correlators  by solving them asymptotically (Chapter~\ref{Chapter_fff_expansions}), but it will be instrumental to study the regularity of the equilibrium measure (Chapter~\ref{Chapter_smoothness}).

 In Part~\ref{Part_Asymptotic}, we carry out the asymptotic analysis of the discrete ensembles. First, we discuss in Chapter~\ref{Chapterlarge} the weak convergence of the empirical distribution of the particles, that is the  law of large numbers, and characterize the limit distribution, namely the equilibrium measure, as the minimizer of the energy functional $-\mathcal{I}$. This characterization is essential in the analysis of the properties of this equilibrium measure, and also in our design of the various interpolation schemes that we invoke. We obtain quantitative estimates on this convergence in the form of concentration of measures estimates and derive large deviations estimates on the support giving conditions to ensure that no particles or no holes appear in void or saturated regions, respectively. This will be crucial to localize and condition our ensembles (Chapter~\ref{Chapter_conditioning}).

  In Chapter~\ref{Chapter_smoothness}, we study the regularity of the equilibrium measure, meaning the regularity of its density and the regularity when varying all the parameters of the ensembles. This is important in many respects. First it allows to make sure that if an ensemble is off-critical, it will remain off-critical if we change the parameters a little, a key point when we check the assumptions or proceed with interpolations. Second, it will be a key to show that the terms in the expansion of the partition functions depend smoothly on the filling fractions, an essential step to study the asymptotic law of the filling fractions and their extended version.

  In Chapter~\ref{Chapterzw} we study our reference discrete ensembles, namely the $zw$-ensembles, which are exactly solvable. We obtain precise estimates on their equilibrium measures and partition function.

  In Chapter~\ref{Chapter_conditioning}, we develop the methods to condition and localize the ensembles to a neighborhood of the bands of the equilibrium measures.

 We show in  Chapter~\ref{Chapter_fff_expansions} how to study the asymptotics of the linear statistics for discrete ensembles satisfying a set  of assumptions from Chapter~\ref{Chapter_Setup_and_Examples}, namely  Assumptions~\ref{Assumptions_Theta}, \ref{Assumptions_basic}, \ref{Assumptions_offcrit} and \ref{Assumptions_analyticity}, and the additional assumptions that the segment filling fractions are deterministically fixed and each segment $[\hat a_h,\hat b_h]$ has a unique band of the
equilibrium measure. This is based on Nekrasov equations that we solve asymptotically by linearization. This in turn hinges on the inversion of the so-called master operator, see Theorem~\ref{Theorem_Master_equation} whose proof is deferred to Chapter~\ref{Chapter_SolvingN}.

We deduce the expansion of the partition function in Chapter~\ref{Chapter_partition_functions} by interpolation  towards the $zw$-discrete ensemble. This allows us deriving the asymptotic behavior of the filling fractions and their extensions in Chapter~\ref{Chapter_filling_fractions}. This is the last step of the proofs of our main asymptotic results for the general discrete ensembles \eqref{eq_measure_general_intro}. All these results are specialized to Gaussian weights in Sections~\ref{mueqGaussian}, \ref{Sec:condgaus} and \ref{Gausscorr}.

We harvest the fruits of the previous chapters and apply our results to random tilings in Chapter~\ref{Chap11}. Here we reduce our assumptions to conditions that are easier to check. Furthermore, we prove that for fixed filling fractions, the global fluctuations in the whole two-dimensional tiling follow the Gaussian free field, thanks to the identification of the covariance given by solving asymptotically the Nekrasov equation with a Green function. When the filling fractions are not fixed, we also prove the appearance of the discrete Gaussian components in the global fluctuations of the heights of tilings.

Part~\ref{Part_Master_equation} contains the complex and algebraic geometry theory needed to invert the master operator linearizing the Nekrasov equations and to compute leading covariances. It can be read mostly independently of the two other parts. Chapter~\ref{Chapter_SolvingN} develops abstract continuity and smoothness properties of the master operator and connects the operator to a Riemann-Hilbert problem. Chapter~\ref{Chapter_AG} develops a general framework to study this Riemann--Hilbert problem via the construction of spectral curves and relates its fundamental solution (\textit{i.e.} the leading covariance for Part~\ref{Part_Asymptotic}) to Green functions and Gaussian free fields in favorable situations. This is the basis for the identification of the Gaussian free field in tiling models in Chapter~\ref{Chap11}. Chapter~\ref{Chap14} showcases examples where the constructions of the previous two sections can be made --- to some extent --- more explicit.

An index of all definitions, theorems, symbols, figures, \textit{etc.} is compiled in Part~\ref{apppart}.

\section{Some notations}
\label{Notasec}

If $k$ is a positive integer, we denote $[k] = \{1,\ldots,k\}$. The notation $\llbracket a,b\rrbracket$ stands for the set of integers $x$ such that $a \leq x \leq b$. If $x \in \amsmathbb{R}$, its integer part $\lfloor x \rfloor$ is the unique integer such that $\lfloor x \rfloor \leq x < \lfloor x \rfloor + 1$. We define likewise $\lceil x \rceil$ the unique integer such that $x < \lceil x \rceil \leq x + 1$. $\amsmathbb{R}_{> 0}$ (respectively $\amsmathbb{R}_{\geq 0}$, $\amsmathbb{R}_{< 0}$, $\amsmathbb{R}_{\leq 0}$) is the set of positive (respectively nonnegative, negative, nonpositive) real numbers. $\widehat{\amsmathbb{C}} = \amsmathbb{C} \cup \{\infty\}$ is the Riemann sphere. Its set of real points $\widehat{\amsmathbb{R}} = \amsmathbb{R} \cup \{\infty\}$ is topologically a circle. The conjugate of a complex number $z$ is $z^*$.
\label{index:indicfn}
If $X$ is a subset of a set $Y$ --- usually obvious from the context --- then $X^{\textnormal{c}} \subset Y$ denotes its complement subset. When $X$ is a subset of $\amsmathbb{C}$, $\mathring X$ denotes its topological interior, which is the maximal open set contained in $X$. If $X$ is a finite set, $\#X$ denotes its cardinality. The indicator function of a set $Y$ is denoted $\mathbbm{1}_Y$, it is equal to $1$ whenever its argument belongs to $Y$ and to $0$ otherwise. We sometimes write $\delta_{x,y}$ instead of $\mathbbm{1}_{\{y\}}(x)$, which is $1$ whenever $x=y$.

The Pochhammer symbol is defined as
\begin{equation}
\label{Pochdef}(x)_{n} := \frac{\Gamma(x + n)}{\Gamma(x)} = x(x + 1)\cdots (x + n - 1).
\end{equation}
The logarithm with base $e$ is denoted $\log$. The integrated logarithm functions will be frequently used:
\begin{equation}
\label{Llogdef}
\begin{split}
\mathrm{Llog}(x) & := \int_{0}^{x} \log |\xi|\,\dd \xi = x \log|x| - x, \\
\mathrm{LLlog}(x) & := \int_{0}^{x} \mathrm{Llog}(\xi)\,\dd \xi = \frac{x^2}{2}\log|x| -  \frac{3x^2}{4}.
\end{split}
\end{equation}
To handle multi-valued functions, we always make the following choice of branches. The square root $z \mapsto \sqrt{z}$ is the holomorphic function defined on $\amsmathbb{C}\setminus \amsmathbb{R}_{< 0}$ such that $\sqrt{1} = 1$. The complex logarithm $z \mapsto \log z$ is the holomorphic function defined on $\amsmathbb{C}\setminus \amsmathbb{R}_{< 0}$ such that $\log 1 = 0$. For the tangent reciprocal, we make the unusual choice $\arctan : \,\widehat{\amsmathbb{R}} \rightarrow [0,\pi]$, which is the unique continuous function satisfying $\tan(\arctan(t)) = t$. In particular, $\arctan(\infty) = \tfrac{\pi}{2}$.

We denote by $\boldsymbol{e}^{(1)},\ldots,\boldsymbol{e}^{(H)}$ the canonical basis of $\amsmathbb{R}^H$ or $\amsmathbb{C}^H$ and by $\langle \boldsymbol{v} \cdot \boldsymbol{w} \rangle$ the standard Euclidean scalar product of two vectors $\boldsymbol{v},\boldsymbol{w} \in \amsmathbb{R}^H$. If $\boldsymbol{v} \in \amsmathbb{C}^H$ and $p \in \amsmathbb{R}_{> 0}$, we denote
\[
|\!|\boldsymbol{v}|\!|_{\infty} = \max_{h \in [H]} |v_h|,\qquad |\!|\boldsymbol{v}|\!|_p = \sum_{h = 1}^{H} |v_h|^p.
\]
Likewise if $\boldsymbol{\Theta}$ is an $H \times H$ matrix we denote $|\!|\boldsymbol{\Theta}|\!|_{\infty} = \max\big\{|\theta_{g,h}|\,\,|\,\,g,h \in [H]\big\}$.

For a compactly-supported Lipschitz function $f: \,\amsmathbb{R} \rightarrow \amsmathbb{C}$, we define
\begin{equation}
\label{eq_Fourier}
 |\!| f|\!|_{\frac{1}{2}} = \bigg(\int_{\amsmathbb{R}} |s| |\hat{f}(s)|^2 \dd s\bigg)^{\frac{1}{2}},
 \quad \hat f(s) =\int_{\amsmathbb{R}} e^{\ii xs } f(x)\dd x,
\end{equation}
 \[
 |\!| f|\!|_{\infty} = \sup_{x \in \amsmathbb{R}} |f(x)|,
\quad
 \quad |\!| f|\!|_\textnormal{Lip}=\sup_{x\neq y}
 \left|\frac{f(x)-f(y)}{x-y}\right|.
 \]
 We define likewise the Fourier transform if $f$ is a measure on $\amsmathbb{R}$. If $X \subset \amsmathbb{R}$, the notation $\mathscr{L}^p(X)$ stands for the space of complex-valued measurable functions on $X$ whose $p$-th power is integrable with respect to Lebesgue.

If $f$ is a holomorphic function in a domain $U \setminus I$, where $U \subset \amsmathbb{C}$ is open and $I \subset \amsmathbb{R}$ is a finite union of segments, we denote for $x \in I$ by $f(x^{\pm})$ the limit $\lim_{\epsilon \rightarrow 0^+} f(x \pm {\ii}\epsilon)$ if it exists.

By absolutely continuous measure we mean a measure on $\amsmathbb{R}$ which is absolutely continuous with respect to the Lebesgue measure. For a signed absolutely continuous measure $\rho$, by $|\rho|$ we mean the positive measure whose density is the absolute value of the
density of $\rho$. If $\rho$ is a finite measure on $\amsmathbb{R}$, we denote
\[
\mathcal{G}_{\rho}(z) = \int_{\amsmathbb{R}} \frac{\dd\rho(x)}{z - x}
\]
its Stieltjes transform, which is a holomorphic function of $z \in \amsmathbb{C} \setminus \textnormal{supp}\,\rho$, where $\textnormal{supp}\,\rho$ is the support of $\rho$.

 If $A$ and $B$ depend on $\N$, we write: $A = O(B)$ if there exists a constant $C > 0$ independent of $\N$ such that $|A| \leq C |B|$ for $\N$ large enough; $A = o(B)$ if for any $\eps > 0$ there exists $\N_0(\eps) > 0$ such that $|A|\leq \eps |B|$ for $\N \geq \N_0(\eps)$.

In this book, $C$ will stand for various (large) universal finite constants, whose value can change from line to line. In particular they are always independent of $\N$ unless stated otherwise.

\part{THE DISCRETE ENSEMBLES}
 \label{Part1}

\chapter{Setup and main results}

\label{Chapter_Setup_and_Examples}

\section{A simple ensemble}

\label{Section_framework}

 In the simplest possible framework, we study random configurations of $N$
integers
 $\lambda_1\leq \lambda_2 \leq \ldots \leq \lambda_N$, whose probability
distribution depends on three pieces of information: a (possibly infinite) real interval $[a,b]$ called ``segment of definition'' or simply \emph{segment}, a real parameter $\theta>0$
 called \emph{intensity of interaction}, and a function $w: [a,b] \rightarrow \amsmathbb R_{\geq 0}$
called \emph{weight}. The distribution itself is defined via the positions of shifted \emph{particles}
\begin{equation}
\label{ilambdal}
\forall i \in [N]\qquad \ell_i=\lambda_i+ (i-1)\theta,
\end{equation}
where $[N]$ means the set $\{1,2,\ldots,N\}$ here and below.
The integrality of
$\lambda_i$s translates into
\begin{equation} \label{eq_lattice_restriction}
\forall i \in [N - 1] \qquad \ell_{i+1}-\ell_{i} \in \{\theta,\theta+1,\theta+2,\ldots\}.
\end{equation}
The probability of a configuration $\boldsymbol{\ell} = (\ell_i)_{i = 1}^{N}$ is non-zero only if all
$\ell_i$s belong to the interval $[a,b]$ and in the latter case is given by
\begin{equation} \label{eq_general_measure_one_cut}
 \amsmathbb{P}_{N}(\boldsymbol{\ell}) =\frac{1}{\Z_N} \cdot \prod_{1\leq i <j \leq N}
 \frac{1}{N^{2\theta}}\cdot \frac{\Gamma\big(\ell_j-\ell_i+1\big)\cdot \Gamma\big(\ell_j-\ell_i+\theta\big) }{\Gamma\big(\ell_j-\ell_i\big)\cdot
 \Gamma\big(\ell_j-\ell_i+1-\theta\big)}\cdot \prod_{i=1}^N w(\ell_i),
\end{equation}
where $\Z_N$ is the normalizing constant called \emph{partition function} which guarantees that
\[
\sum_{\boldsymbol{\ell}} \amsmathbb{P}_{N}(\boldsymbol{\ell})=1.
\]
If $\theta$ is an integer or a half-integer, then the ratios of Gamma functions in
\eqref{eq_general_measure_one_cut} can be simplified. In particular:
\begin{equation} \label{eq_general_measure_one_cut_beta_1}
\begin{split}
\textnormal{for}\,\,\theta = \frac{1}{2} & \qquad \amsmathbb{P}_{N}(\boldsymbol{\ell})=\frac{1}{\Z_N}\cdot \prod_{1\leq i <j \leq N}
\frac{\ell_j - \ell_i}{N} \cdot \prod_{i=1}^N w(\ell_i) \\
\textnormal{for}\,\,\theta = 1 & \qquad
 \amsmathbb{P}_{N}(\boldsymbol{\ell}) =\frac{1}{\Z_N} \cdot \prod_{1\leq i <j \leq N}\bigg(
\frac{\ell_j - \ell_i}{N}\bigg)^2 \cdot \prod_{i=1}^N w(\ell_i).
\end{split}
\end{equation}

The advantage of working with $\ell_i$ given by \eqref{ilambdal} instead of $\lambda_i$ is that the probability of a configuration has a mean-field nature: it is a function of the empirical measure $\sum_{i = 1}^{N} \delta_{\ell_i}$, in particular the contribution of $\ell_i$ does not depend on the index $i$. However, when $\theta \neq 1$ there is a trade-off which will often complicate our analysis in comparison to continuous ensembles: the sites available for $\ell_i$ do depend on $i$.

\section{General description of the discrete ensembles}
\label{Section_general_model}
In the general setting, we study a generalization of \eqref{eq_general_measure_one_cut} where the $N$ particles
are split into $H$ groups and the intensity of interaction can depend on the groups, which we motivated in Chapter~\ref{SIntro}.

\subsection{The data}
\label{DataS}
We introduce a real parameter $\N \geq 1$, called \emph{master parameter}. All the other data entering in the definition of the ensemble may depend on $\N$, and we are interested in the study of $\N \rightarrow \infty$ asymptotics. In particular, the number of particles in the system will be growing linearly in $\N$.

Let $H$ be a positive integer. For each $h \in [H]$ we are given a closed interval of the real line $[a_h,b_h]$
 and a weight $w_h: [a_h,b_h]\rightarrow \amsmathbb R_{\ge
0}$. We call these intervals \emph{segments}, they must be ordered, pairwise disjoint and have nonempty interior, \textit{i.e.}
\[
-\infty \leq a_1<b_1<a_2<b_2<\dots<a_H<b_H \leq + \infty.
\]
If $a_1 = -\infty$ one has to understand $[a_1,b_1]$ as $(-\infty,b_1]$, and similarly for $b_H$. We are given a symmetric matrix of intensities of interactions $\boldsymbol{\Theta} =
(\theta_{g,h})_{g,h = 1}^H$. We assume later that this matrix is positive semi-definite and has
positive diagonal. For instance, it could be a matrix filled with the same positive number
$\theta$, or a positive diagonal matrix.

\label{index:Nh}The configurations of the ensemble depend on an $H$-tuple of nonnegative integers $\boldsymbol{N} =(N_1,\ldots, N_H)$ which specify the number of particles in each segment, also called \emph{segment filling fractions}. We may impose extra constraints on them. The constraints we allow must form a
system of affine equations, that we denote \eqref{eq_equations_eqs}. In other words, we are given a family $(\mathfrak{r}_e)_{e = 1}^{\mathfrak{e}}$ of independent linear forms in $\amsmathbb{R}^H$, an $\mathfrak{e}$-tuple of real numbers $\boldsymbol{r} = (r_e)_{e = 1}^{\mathfrak{e}}$, and we impose
\begin{equation} \tag{$\star$} \label{eq_equations_eqs}
\forall e \in [\mathfrak{e}]\qquad \mathfrak{r}_e\bigg(\frac{\boldsymbol{N}}{\N}\bigg)= r_e.
\end{equation}
We denote $\mathfrak{r} : \amsmathbb{R}^H \rightarrow \amsmathbb{R}^{\mathfrak{e}}$ the linear map with components $\mathfrak{r}_1,\ldots,\mathfrak{r}_{\mathfrak{e}}$.

 \label{index:tot}The total number of particles $N := N_1+\cdots+N_H$ is always fixed as a parameter of the ensemble. In other words, \eqref{eq_equations_eqs}
 must imply the constraint
 \begin{equation}
 \label{eq_number_of_particles_contraint0}
\sum_{h = 1}^H \frac{N_h}{\N} = \frac{N}{\N}.
 \end{equation}
 Moreover, we assume that the linear forms $\mathfrak{r}_e$ have integral coefficients. In principle, we are not interested in the linear forms themselves, but rather only in the solutions to \eqref{eq_equations_eqs}. Hence, we have some freedom in choosing $\mathfrak{r}_1,\ldots,\mathfrak{r}_{\mathfrak{e}}$. In particular, \eqref{eq_number_of_particles_contraint0} might be the first constraint corresponding to $\mathfrak r_1$ and $r_1$, or it can be a corollary of other constraints. For instance, an allowed choice of constraints consists in setting $\mathfrak e=H$ and $\mathfrak{r}_h(X_1,\ldots,X_H)=X_h$ for $h \in [H]$, so that the $h$-th constraint fixes the $h$-th segment filling fraction to be $N_h$. In this situation $N_1+\cdots+N_H$ is certainly fixed and \eqref{eq_number_of_particles_contraint0} is a corollary of the constraints.

We could take $N = \N$ (\textit{i.e.} $r_1 = 1$), but we do not necessarily want to impose this. The reason to distinguish between the master parameter $\N$ and the total number of particles $N$ is that in certain situations some of the particles lie almost surely at deterministic locations. It then becomes convenient to remove these particles from the considerations. This changes $N$ but we can keep the same master parameter $\N$.

Integral solutions to \eqref{eq_equations_eqs}, their rescaled version, and real solutions of \eqref{eq_equations_eqs} all have a role to play:
\begin{equation} \label{eq_Lambda_def}
\begin{split}
\Lambda_{\star} & := \big\{\boldsymbol{N} \in \amsmathbb{Z}_{\geq 0}^H\quad |\quad \boldsymbol{N}\,\,\textnormal{satisfies}\,\,\eqref{eq_equations_eqs}\big\} \\
\hat{\Lambda}_{\star} & := \N^{-1}\Lambda_{\star} \\
\hat{\Lambda}_\star^{\amsmathbb{R}} & := \big\{\N^{-1}\,\boldsymbol{N} \quad \big|\quad \boldsymbol{N} \in \amsmathbb{R}_{\geq 0}^H\,\,\textnormal{satisfies}\,\,\eqref{eq_equations_eqs}\big\}.
\end{split}
\end{equation}
The difference between $\hat \Lambda_\star$ and $\hat{\Lambda}_\star^{\amsmathbb{R}}$ is in the integrality condition used for the former, but not for the latter. However in the latter we still impose nonnegativity. If each of the three sets $\Lambda_\star$, $\hat \Lambda_\star$, and $\hat{\Lambda}_\star^{\amsmathbb{R}}$ has precisely one element, we say that the ensemble has \emph{fixed filling fractions}.

\subsection{The state space}\label{Section_configuration_space}

Given all the above
information, the state space
$\W_\N$ is the set of $N$-tuples $\boldsymbol{\ell} = (\ell_1,\ldots,\ell_N)$ such that $\ell_1<\ell_2<\cdots<\ell_N$ and satisfying the following conditions
\begin{enumerate}
 \item For any $i \in [N]$, $\ell_i$ is a finite real number in $\bigcup_{i=1}^H [a_i,b_i]$. We denote $h(i) \in [H]$ the index of the unique segment such that $\ell_i\in[a_{h(i)},b_{h(i)}]$.
 \item The segment filling fractions $N_h =\#\big\{i \in [N]\,\,|\,\,\ell_i\in [a_h,b_h]\big\}$ indexed by $h \in [H]$ satisfy the constraints \eqref{eq_equations_eqs}.
 \item If for an index $i \in [N - 1]$, both $\ell_i$ and $\ell_{i+1}$ belong to the
 same interval $[a_h,b_h]$, then $\ell_{i+1}-\ell_i \in \theta_{h,h} + \amsmathbb{Z}_{\geq 0}$.
 \item \label{index:Iplusmoins}Let $I^-(h) = \min\{i \in [N]\,\,|\,\,\ell_i\in[a_h,b_h]\}$. For any $h \in [H]$, if $a_h$ is finite we should have $\ell_{I^-(h)} - a_h \in \amsmathbb Z_{\geq 0}$. Let $I^+(h) = \max\{i \in [N]\,\,|\,\, \ell_i\in[a_h,b_h]\}$. For any $h \in [H]$, if $b_h$ is finite we should have $b_h-\ell_{I^+(h)}\in \amsmathbb Z_{\geq 0}$. If $H=1$, $a_1=-\infty$, $b_1=+\infty$, none of these constraints make sense, and we instead assume that we are given $x_0 \in \amsmathbb{R}$ such that $\ell_1 \in \amsmathbb{L}:= x_0 + \amsmathbb{Z}$.
 \end{enumerate}

\noindent We observe that the third condition implies
\begin{equation}
\label{eq_segment_ff_relation}
 \forall h \in [H]\qquad b_h-a_h-\theta_{h,h} N_h\in \amsmathbb Z_{\geq 0}
\end{equation}
whenever $a_h$ and $b_h$ are finite. If the equations \eqref{eq_equations_eqs} uniquely fix the value of $N_h$, the condition \eqref{eq_segment_ff_relation} is not restrictive. Similarly, it is not restrictive if $\theta_{h,h}$ is a positive integer. But, more generally, \eqref{eq_segment_ff_relation} constrains the possible values of $N_h$; in particular, if $\theta_{h,h}$ is irrational, then $a_h$, $b_h$ and \eqref{eq_segment_ff_relation} uniquely determine the value of $N_h$\footnote{For irrational $\theta_{h,h}$, one can be tempted to achieve the possibility of varying $N_h$ with fixed $a_h,b_h$ by adjusting the definition and waiving the requirement $b_h-\ell_{I^+(h)}\in \amsmathbb{Z}_{\geq 0}$. This is, indeed, possible, yet it would lead to additional terms appearing in the expansions of the partition functions in Chapter~\ref{Chapter_partition_functions} and in the asymptotics of fluctuating filling fractions in Chapter~\ref{Chapter_filling_fractions}. We do not follow this route.}. Let us mention a slight difference in notations with \cite{BGG}: there the particles
belonged to the intervals $[a_h+1,b_h-1]$ while here they belong to $[a_h, b_h]$.

\subsection{The probability measure}

We assume that the data defining the ensemble is chosen such that the state space $\W_\N$ is nonempty; in particular, $\Lambda_\star$ has to be nonempty. The associated \emph{discrete ensemble} is the random ensemble specified by the probability measure $\P$ on $\W_\N$
\begin{equation} \label{eq_general_measure}
 \P(\boldsymbol{\ell})= \frac {1}{\Z_\N} \cdot \prod_{1\leq i<j \leq N}\bigg[
 \frac{1}{\N^{2\theta_{h(i),h(j)}}}\cdot \frac{\Gamma\big(\ell_j-\ell_i+1\big)\cdot
 \Gamma\big(\ell_j-\ell_i+\theta_{h(i),h(j)}\big)}{\Gamma\big(\ell_j-\ell_i\big)\cdot
 \Gamma\big(\ell_j-\ell_i+1-\theta_{h(i),h(j)}\big)}\bigg] \cdot
  \prod_{i=1}^N w_{h(i)}(\ell_i).
\end{equation}
In words, the interaction between pairs of particles given by the ratio of four Gamma functions now depends on the segments to which the particles belong, and similarly the expression
for the multiplicative weight may depend on the segment. We assumed implicitly that the \emph{partition function}
\begin{equation}
\label{eq_partition_function_definition}
 \Z_\N =\sum_{\ell \in \W_\N} \prod_{1\leq i<j \leq N}\bigg[
 \frac{1}{\N^{\,2\theta_{h(i),h(j)}}}\cdot \frac{\Gamma\big(\ell_j-\ell_i+1\big)\cdot
 \Gamma\big(\ell_j-\ell_i+\theta_{h(i),h(j)}\big)}{\Gamma\big(\ell_j-\ell_i\big)\cdot
 \Gamma\big(\ell_j-\ell_i+1-\theta_{h(i),h(j)}\big)}\bigg]\cdot
  \prod_{i=1}^N w_{h(i)}(\ell_i).
\end{equation}
is finite, so that $\P$ is a indeed a probability measure. Proposition~\ref{lem_finite_partition_function} will justify the finiteness of $\Z_\N$ under the natural Assumptions~\ref{Assumptions_Theta} and \ref{Assumptions_basic}. We will sometimes switch to the different notation $\ell_i^{h}$ for $h \in [H]$ and $i \in [N_h]$ in order to stress the various groups of particles. The correspondence between the two notations is
\begin{equation}
\label{eljah}\ell_{i}^{h} = \ell_{\sum_{g = 1}^{h - 1} N_{g} + i}.
\end{equation}

We can consider this ensemble equivalently as an ensemble of $H$ random Young diagrams $\boldsymbol{\lambda}^1,\ldots,\boldsymbol{\lambda}^{H}$, whose rows (labeled from the smallest to the largest) are related to the configuration $\boldsymbol{\ell}$ by $\ell_i^{h} = \lambda_i^h + i\theta_{h,h} + c_h$ for some constant $c_h$. Indeed, Condition 3. of Section~\ref{Section_configuration_space} turns into
\[
\forall i \in [N_h - 1]\qquad (\lambda_{i + 1}^h - \lambda_i^h) \in \amsmathbb{Z}_{\geq 0}.
\]
As in Section~\ref{Section_framework}, there is a clear advantage of working with $\boldsymbol{\ell}$ since the expression of the probability $\amsmathbb{P}_\N(\boldsymbol{\ell})$ only depends on the empirical measures of the configurations $\boldsymbol{\ell}^{1},\ldots,\boldsymbol{\ell}^{H}$. Yet, we must bear with complications in the description of allowed sites for a given particle $\ell_i^h$, as they now depend on the number of particles to its left (if $a_h$ is finite or $H = 1$ and $a_1,b_1$ are infinite) or to its right (if $h = 1$, $a_1$ is infinite and $b_1$ is finite) in this segment.

\section{List of the imposed assumptions}

\label{Section_list_of_assumptions}

Our ultimate goal is to study the asymptotic behavior of the partition function and various statistics of $\boldsymbol{\ell}$ under the law $\P$ of
\eqref{eq_general_measure} as $\N\rightarrow\infty$. For this we assume that all the data
specified in Section~\ref{Section_framework} depend on $\N$ but we impose
restrictions on this dependence. In general, if we talk about discrete ensembles, then we use variables $a_h$, $b_h$, $N_h$ and their rescaled version denoted with a hat
\begin{definition}
\label{def:eq_rescaled_parameters}
The rescaled parameters are
\begin{equation} \label{eq_rescaled_parameters}
 \hat a_h =\frac{a_h}{\N},\qquad \hat b_h =\frac{b_h}{\N},\qquad \hat n_h
=\frac{N_h}{\N}.
 \end{equation}
 \end{definition}
In our setting the variables with hat will usually be $O(1)$ as $\N \rightarrow \infty$ and for that reason are more adapted to the asymptotic analysis. If we discuss continuous limits --- such as the equilibrium measure --- we use shifted segment endpoints $\hat a'_h$ and $\hat b'_h$ instead.
\begin{definition}
\label{def:eq_shifted_parameters}
The left and right-shifted parameters are
\begin{equation}
\label{eq_shifted_parameters}
 \hat a'_h=\hat a_h - \frac{\theta_{h,h}-\frac{1}{2}}{\N}\qquad \textnormal{and} \qquad \hat b'_h=\hat b_h + \frac{\theta_{h,h}-\frac{1}{2}}{\N}.
\end{equation}
The shifted and rescaled domain of definition is denoted
\begin{equation}
\label{eq_measure_space}
 \amsmathbb{A} = \bigcup_{h = 1}^H
[\hat{a}'_{h},\hat{b}'_{h}]
\end{equation}
\end{definition}
The primary reason for these peculiar shifts is of technical nature: it is used in our arguments for a precise match of locations of certain zeros of holomorphic observables in discrete ensembles and in their continuous counterparts. They will be relevant when we discuss equilibrium measures. Since both rescaled segments and their shift will play a role, we will also need slightly enlarged segments.
\begin{definition}
\label{def:eq_enlarged_parameters}
The rescaled enlarged segment $\amsmathbb{A}_h^{\mathfrak{m}} = [\hat{a}_h^{\mathfrak{m}},\hat{b}_h^{\mathfrak{m}}]$ is defined by
\[
\hat{a}_h^{\mathfrak{m}} = \min(\hat{a}_h,\hat{a}_h') \qquad \textnormal{and}\qquad \hat{b}_h^{\mathfrak{m}} = \max(\hat{b}_h,\hat{b}_h').
\]
\end{definition}

\subsection{Interaction and segment filling fractions}

\begin{assumption} \label{Assumptions_Theta} There exists a constant $C > 0$ independent of $\N$ such that
\begin{enumerate}
\item $\boldsymbol{\Theta}$ is a real symmetric positive semi-definite $H \times H$ matrix.
\item $\forall h \in [H] \qquad \theta_{h,h}\geq \frac{1}{C}$.
\item $H \leq C$,\,\, $\mathfrak{e} \leq C$,\, and \,\, $|\!|\boldsymbol{\Theta}|\!|_{\infty} \leq C$.
\item There exists a real symmetric positive semi-definite $\mathfrak{e} \times\mathfrak{e}$ matrix $\boldsymbol{\Theta}'$ such that
\begin{equation}
\label{eq_Theta_through_Theta_prime}
\forall \boldsymbol{X} \in \amsmathbb{R}^H\qquad \boldsymbol{X}^{T} \cdot \boldsymbol{\Theta} \cdot \boldsymbol{X} \\ = \big( \mathfrak r_1(\boldsymbol{X}),
 \ldots, \mathfrak r_{\mathfrak{e}}(\boldsymbol{X}) \big)^{T} \cdot \boldsymbol{\Theta}' \cdot \big( \mathfrak r_1(\boldsymbol{X}),
 \ldots, \mathfrak r_{\mathfrak{e}}(\boldsymbol{X}) \big).
\end{equation}
\item We have $|\!|\boldsymbol{r}|\!|_{\infty} < C$, the coefficients of the linear forms $\mathfrak{r}_1,\ldots,\mathfrak{r}_{\mathfrak{e}}$ are bounded by $C$, and the image of the parallelepiped $\prod_{h = 1}^H \big[0,\,\frac{1}{\theta_{h,h}}(\hat b_h - \hat a_h)\big]$ under the map $\boldsymbol{X} \mapsto (\mathfrak{r}_e(\boldsymbol{X}))_{e = 1}^{\mathfrak{e}}$ contains the closed ball of radius $\frac{1}{C}$ centered at $\boldsymbol{r}$.
\item If for some $h\in [H]$ both $a_h$ and $b_h$ are finite and $\theta_{h,h} \notin \amsmathbb{Z}$, then we require that the equations \eqref{eq_equations_eqs} uniquely determine the value of the filling fraction $N_h$.
\end{enumerate}
\end{assumption}

Here $|\!|\cdot|\!|_{\infty}$ stands for the sup norm. The main role of the constant $C$ is to control the uniformity of our asymptotic results, and there will be several other assumptions of this kind. We will make statements like "this error or bound depends only on the constants in the assumptions", but it occurs so often that we sometimes will keep it implicit. In particular, $C$ constrains the dependence in $\N$ of the parameters of the ensemble. While it is natural that the rescaled segment endpoints are $\N$-dependent, it may seem artificial to allow a $\N$-dependence in the interaction matrix $\boldsymbol{\Theta}$. However, it would be even more artificial in the logic of our arguments to give a different status to the various parameters of the ensemble. Therefore, it is more appropriate to consider that all the parameters of the discrete ensemble depend on $\N$ unless specified otherwise, in the limits fixed by our assumptions.

Condition 1. and 4. guarantee that the energy functional driving the large $\N$ limit of the discrete ensemble is strictly convex --- \textit{cf.} Lemma~\ref{Lemma_I_quadratic_Fourier} for
details --- and this essential to start the asymptotic analysis. Conditions 1. and 4. will be automatic in tiling models of Chapter~\ref{Chap11} due to \eqref{ThetaRTR}: we have $\boldsymbol{\Theta}' = \frac{1}{2}\textnormal{\textbf{Id}}$ and part of the affine constraints are specifying the total number of free lozenges on the large base of each trapezoid. In general, if the equations \eqref{eq_equations_eqs} completely fix the segment filling fractions,
then Condition 4.\ becomes empty as the linear forms $\mathfrak{r}_e$ for $e \in [\mathfrak{e}]$ span the entire $\amsmathbb{R}^H$. It is plausible that Condition 4.\ can be weakened, the essential property being the aforementioned convexity. Condition 5.\ means that the constraints \eqref{eq_equations_eqs} on the segment filling
fractions are linearly independent and do not force the system to have intervals $[a_h,b_h]$ that are macroscopically empty or fully packed with particles (such configurations can exist but they should not be the only ones!), \textit{cf.} Lemma~\ref{Lemma_non_empty_interior}. This would leave no room for macroscopic fluctuations in this segment, and has little interest.
Condition 6.\ implies that any constraints on $N_h$ forced by \eqref{eq_segment_ff_relation} are already a part of \eqref{eq_equations_eqs}. Our results could be adapted if $\theta_{h,h}$ were rational and $N_h$ were allowed to vary; in this situation \eqref{eq_segment_ff_relation} implies that given $a_h$ and $b_h$ the possible values of $N_h$ form a lattice with spacing equal to the denominator of $\theta_{h,h}$. However, due to lack of natural examples we decided not to pursue this direction.

 A simple example of $\boldsymbol{\Theta}$ satisfying Assumption
\ref{Assumptions_Theta} is the matrix whose elements are all equal to the same parameter $\theta > 0$. Indeed, in this case we have $\mathfrak e=1$ and a single linear form $\mathfrak{r}_1(X_1,\ldots,X_H)=X_1+\dots+X_H$ corresponding to the constraint on the total number of particles \eqref{eq_number_of_particles_contraint0}. Thus, we can take $\boldsymbol{\Theta}'=\theta$, so that \eqref{eq_Theta_through_Theta_prime} takes the form
 $(X_1,\ldots,X_H) \cdot \boldsymbol{\Theta} \cdot
(X_1,\ldots,X_H)^T= \theta (X_1+\ldots+X_H)^2$. In this situation the interactions between all pairs of particles are the same. When the filling fractions are fixed, another example is $\boldsymbol{\Theta} = \textnormal{diag}(\theta,\ldots,\theta)$ for some $\theta > 0$. From the probabilistic point of view, the particles are split into $H$ independent groups. In this case
\[
\mathfrak e=H,\qquad {\mathfrak r}_e(X_1,\ldots,X_H)=X_e,\qquad \boldsymbol{\Theta}=\boldsymbol{\Theta}'.
\]

\subsection{Regularity of the weights at large \texorpdfstring{$\N$}{}}

\begin{assumption}\label{Assumptions_basic}
There exists a constant $C > 0$ independent of $\N$ making the following conditions hold.
\begin{enumerate}
\item $\frac{\N}{C} \leq N \leq \N C$.
\item $\forall g,h \in [H] \qquad |\hat b_g - \hat a_{h}| \geq \frac{1}{C}$.
\item For any $h \in [H]$, if $a_h$ is finite, then $|\hat{a}_h| < C$; if $b_h$ is finite, then $|\hat{b}_h| < C$.
\item For any $h \in [H]$, $w_h$ is positive, continuous in
$[a_h,b_h]$ and admits a decomposition
\begin{equation}
 \label{eq_weight_form}
  w_h(x)=\exp\Bigg(-\N V_h\bigg(\frac{x}{\N}\bigg)+\err_h(x)\Bigg).
 \end{equation}
There exists $\iota_h^{\pm} \in \amsmathbb{Z}_{\geq 0}$ and a twice-continuously differentiable function $U_{h}: \,\amsmathbb{A}_h^{\mathfrak{m}} \rightarrow \amsmathbb{R}$ such that
\begin{equation}
 V_h(x) = \iota_h^{-}\,\mathrm{Llog}(x - \hat{a}'_{h}) + \iota_{h}^+\,\mathrm{Llog}(\hat{b}'_h - x) + U_h(x).
 \label{eq_V_logs}
\end{equation}
where $\mathrm{Llog}(x) = x \log|x| - x$. We assume $\iota_h^{\pm} \leq C$. For each $D>0$, we require the existence of a constant $K(D) > 0$ independent of $h$ and $\N$ such that for any $h \in [H]$ and $x \in [-D,D] \cap \amsmathbb{A}_h^{\mathfrak{m}}$ we have
\begin{equation}
\label{V_bound} \max\big(|U_h(x)|, |\partial_xU_h(x)|, |\partial_x^2U_h(x)|\big) \leq K(D).
\end{equation}
The function $\err_h(x)$ satisfies
\begin{equation}
\label{error_bound_e}\forall x \in [a_{h}, b_{h}]\qquad |\err_h(x)|\leq C \log \N.
\end{equation}
\item \label{point5} If $\theta_{g,g}\neq 1$ for at least one $g \in [H]$, then we assume that for any $h \in [H]$, the function
 $\err_h(x)$ is piecewise monotonous, with at most $C$ intervals of monotonicity.
\item If $a_{1} = -\infty$, then we assume that
there exists $\eta_{1} > 0$ such that
\begin{equation}
\label{eq_confinement_left} \forall x \in (-\infty,-C)\qquad \Bigg(V_1(x) - \bigg[\eta_{1} + 2\sup_{\hat{\boldsymbol{n}} \in \hat{\Lambda}_\star}\bigg(\sum_{h=1}^H \hat{n}_{h}\theta_{1,h}\bigg)\bigg]\log|x|\Bigg)
> -C.
\end{equation}
In this situation we require $\iota_1^{-} = 0$, so that the corresponding term in \eqref{eq_V_logs} is absent.

If $b_H = +\infty$, then we assume that there exists $\eta_{H}>0$ such that
\begin{equation}
\label{eq_confinement_right}\forall x \in (C,+\infty)\qquad \Bigg(V_H(x) - \bigg[\eta_{H} + 2 \sup_{\hat{\boldsymbol{n}} \in \hat{\Lambda}_\star}\bigg(\sum_{h = 1}^H \hat{n}_{h}\theta_{H,h}\bigg)\bigg]\log|x|\Bigg) >
-C.
\end{equation}
In this situation we require $\iota_H^{+} = 0$, so that the corresponding term in \eqref{eq_V_logs} is absent.

If we simultaneously have $a_1=-\infty$, $b_H=+\infty$ and $\theta_{1,H}<0$, then we replace
\eqref{eq_confinement_left}-\eqref{eq_confinement_right} by the stronger assumptions
\begin{equation}
\label{eq_confinement_left_refined} \forall x \in (-\infty,-C)\qquad \Bigg(V_1(x) - \bigg[\eta_{1} + 2\sup_{\hat{\boldsymbol{n}} \in \hat{\Lambda}_\star}\bigg(\sum_{h=1}^{H-1} \hat{n}_{h}\theta_{1,h} + \frac{\hat{n}_H \theta_{1,H}}{2}\bigg)\bigg]\log|x|\Bigg)
> -C,
\end{equation}
\begin{equation}
\label{eq_confinement_right_refined}\forall x \in (C,+\infty)\qquad \Bigg(V_H(x) - \bigg[\eta_{H} + 2 \sup_{\hat{\boldsymbol{n}} \in \hat{\Lambda}_\star} \bigg(\frac{\hat{n}_1 \theta_{H,1}}{2} + \sum_{h = 2}^H \hat{n}_{h}\theta_{H,h}\bigg)\bigg]\log|x|\Bigg)
> -C.
\end{equation}
\end{enumerate}
\end{assumption}
Condition 1.\ means that the number of particles grows linearly with the master parameter $\N$. Condition 2.\ ensures that the intervals $[\hat{a}_{h},\hat{b}_{h}]$ keep a macroscopic size and remain at macroscopic distance from each other. Condition 3.\ prevents finite endpoints to escape to infinity as $\N$ grows, while simultaneously allowing $a_1$ or $b_H$ to be infinite; whether they are finite or infinite is in any case a property independent of $\N$.

Given $w_h(x)$, we have some freedom in Condition 4. to choose $V_h$, as we can move terms of order $O(\log \N)$ between $\N V_h$ and the $\err_h$ in \eqref{eq_weight_form}. The reason for the appearance of the error $\err_h$ in the definition is that we would like to be able to combine
the continuity of $V_h(x)$ with the condition $w_h(a_h-1)=w_h(b_h+1)=0$ that appears naturally in the random tiling models and which is helpful to carry out the asymptotic analysis through the last part of Theorem~\ref{Theorem_Nekrasov} and Remark~\ref{rem:vanishing_interpretation}. Away from the endpoints $a_h,b_h$, we can always absorb this error in $U_h$ and force $\err_h$ to vanish. It is important for us to allow $x\log |x|$ singularities in $V$, for two reasons. First, it often occurs in applications, in particular in random tiling models. Second, a macroscopic number of particles occupying consecutive sites will contribute to the potential felt by a nearby particle by such a singular term with integer $\iota = 2$. In our study we allow $\iota$ to be an arbitrary integer as the value $\iota = 2$ does not appear to be special. We can slightly move the positions of $x\log |x|$ singularities by modifying the error $\err_h$. In \eqref{eq_V_logs} a special choice is made by requiring these singularities to be at the shifted points $\hat a'_h$ and $\hat b'_h$ of Definition~\ref{def:eq_shifted_parameters}. This is a technical choice which reveals important for assembling together different parts of our arguments, in particular in Chapter~\ref{Chapter_fff_expansions} and Step 2 of the proof of Proposition~\ref{Proposition_Nek_1_asymptotic_form}. Besides, the use of shifted variables will allow us to recognize some simplifications in our asymptotic expansion results (Chapters~\ref{Chapter_fff_expansions} and \ref{Chapter_filling_fractions}) when $\theta_{g,g} = 1$ for every $g \in [H]$, in particular for random tiling models. A concrete example of a simplification can be observed in Lemma~\ref{lem:EulerMaclaurin}.

Condition 5.\ guarantees that the choice of splitting between $V_h$ and $\err_h$ in \eqref{eq_weight_form} is not too pathological. When combined with \eqref{V_bound}, the condition rules out some fast oscillations of $w_h$, yet at this point we do not know any interesting examples where we would not be able to choose $\err_h$ satisfying this condition. We use Condition 5.\ to simplify our large deviations argument in Chapter~\ref{Chapterlarge}.

Condition 6.\ is an expression of confinement. It guarantees that the energy functional (as in \eqref{eq_I_intro} and Definition~\ref{Definition_functional} below) is well-defined (Lemma~\ref{Lemma_I_welldef}), that the partition function of the ensemble is finite and the probability measure $\amsmathbb{P}_{\N}$ in \eqref{eq_general_measure} is well-defined (Proposition~\ref{lem_finite_partition_function}), and that the particles almost surely stay in a region of size $O(\N)$ (Remark~\ref{rkcompactsupport}). Due to the different intensities of interaction between the different groups of particles which can all have different populations, the confinement conditions are more involved than the ones usually met in Coulomb gases. Studying discrete ensembles without confinement is possible by importing ideas from \cite{Hardy_Kuijlaars_min,Hardynoconfin2} but introduces new technicalities which are beyond the scope of this book.
Later assumptions will in fact make $U_h$ real-analytic instead of twice-continuously differentiable, but we kept the set of assumptions separate for a better understanding of the logical necessities in the results.

\begin{remark}
\label{rembounder} The conditions \eqref{eq_V_logs} and \eqref{V_bound} imply for any $x \in [-D,D] \cap \amsmathbb{A}_h^{\mathfrak{m}}$ and $h \in [H]$
\begin{equation}
\label{eq_second derivative_bound}
\begin{split}
 |\partial_x V_h(x)| & \leq K'(D)\Big(1+ \big|\log|x-\hat a'_h|\big|+\big|\log|\hat b'_h - x|\big|\Big),\\
 |\partial_x^2 V_h(x)| & \leq K''(D)\bigg(1+ \frac{1}{|x-\hat a'_{h}|}+\frac{1}{|\hat{b}'_{h} - x|}\bigg).
 \end{split}
 \end{equation}
for some constants $K'(D),K''(D) > 0$.
\end{remark}

\subsection{Off-criticality of the equilibrium measure}\label{eqmeas}
\label{Offcrit_section}
Assumptions~\ref{Assumptions_Theta} and \ref{Assumptions_basic} allow us getting started on the asymptotic analysis. Chapter~\ref{Chapterlarge} develops the large deviation principles for the discrete ensembles of Section~\ref{Section_general_model} using only these assumptions. When we push the analysis further to reach an asymptotic expansion for the correlators (Chapter~\ref{Chapter_fff_expansions}), for the partition function (Chapter~\ref{Chapter_partition_functions}) and central limit theorems (Chapter~\ref{Chapter_filling_fractions}), we need additional assumptions on the weights $w_h(x)$, which are listed in this and the next section.

Before presenting the third assumption, we need to talk about \emph{equilibrium measures} and
related notations. Let $\mathscr{P}_h^{\hat n_h}$ be the set of absolutely
continuous nonnegative measures supported on $ [\hat a'_h,\hat b'_h]$ with mass
$\hat{n}_{h}$ and density bounded from above by $\frac{1}{\theta_{h,h}}$. We denote
\begin{equation}
\label{Pnhateq}
\mathscr{P}^{\hat{\boldsymbol{n}}} = \prod_{h = 1}^{H} \mathscr{P}_h^{\hat n_h},
\end{equation}
Elements of $\mathscr{P}^{\hat{\boldsymbol{n}}}$ can be considered in two equivalent ways: either as tuples of measures $\boldsymbol{\nu} = (\nu_h)_{h = 1}^{H}$, or as a single measure $\nu$ with support in $\amsmathbb{A}$. The relation between the two is $\nu = \sum_{h = 1}^{H} \nu_h$, and we will switch back and forth between the tuple $\boldsymbol{\nu}$ and the measure $\nu$ depending on the context. We set
\begin{equation}
\label{eq_ff_set_of_measures}
\mathscr{P}_{\star} = \bigcup_{\hat{\boldsymbol{n}} \in \hat{\Lambda}_\star^{\amsmathbb{R}}} \mathscr{P}^{\hat{\boldsymbol{n}}}.
\end{equation}
In other words, $\mathscr{P}_\star$ consists of those nonnegative measures with support in $\amsmathbb{A}$ and having density bounded from above by $\frac{1}{\theta_{h,h}}$, whose masses satisfy the constraints \eqref{eq_equations_eqs}.

\begin{definition} \label{Definition_functional} Given an $H$-tuple $\boldsymbol{\nu}=(\nu_h)_{h = 1}^H$ of nonnegative measures, with $\nu_h$ having support in $[\hat a'_h,\hat b'_h]$ and a density $\nu_h(x)$ with respect to the Lebesgue measure, we define
 \begin{equation}
 \label{eq_functional_general}
  \mathcal{I}[\boldsymbol{\nu}]=\sum_{g,h = 1}^{H} \int_{\hat a'_{g}}^{\hat b'_{g}} \int_{\hat a'_{h}}^{\hat b'_{h}} \theta_{g,h}\log|x-y|
 \, \nu_{g}(x)\nu_{h}(y)\dd x \dd y- \sum_{h=1}^{H} \int_{\hat a'_h}^{\hat b'_h} V_h(x) \nu_h(x)\dd x.
 \end{equation}
 whenever the integral in the right-hand side makes sense. We call $-\mathcal{I}$ the energy functional.
 \end{definition}

There are two possible sources of divergence for the integral in \eqref{eq_functional_general}: the logarithmic singularity and the presence of an infinite segment. As we work with measures of bounded density the logarithmic singularity is integrable, but infinite $\hat a'_1$ or $\hat b'_H$ can still cause divergences. Under Assumptions~\ref{Assumptions_Theta} and \ref{Assumptions_basic}, Lemma~\ref{Lemma_I_welldef} proves that the integral \eqref{eq_functional_general} is well-defined with values in $\amsmathbb{R} \cup \{-\infty\}$.
\begin{definition}
\label{Definition_equilibrium_measure} A minimizer of $-\mathcal I$ on
$\mathscr{P}_\star$ is called \emph{equilibrium measure} and denoted by
$\boldsymbol{\mu}=(\mu_h)_{h = 1}^H$. We write $\mu = \sum_{h = 1}^{H} \mu_h$, and for $x \in [\hat{a}_h',\hat{b}_h']$ we write $\mu_h(x)$ for the density of $\mu_h$.
\end{definition}

Since all data depend on $\N$, so does the equilibrium measure $\mu$. To access the leading asymptotics of our ensemble it would suffice to replace this $\N$-dependent equilibrium measure with its $\N \rightarrow \infty$ limit --- such a limit will typically exist up to extracting subsequences. But, as we aim at results going beyond leading asymptotics, it is preferable to keep the equilibrium measure $\N$-dependent. Otherwise, the way the parameters approach a limit in the $\N \rightarrow \infty$ would affect the description of the various subleading terms, and this would vary from ensemble to ensemble. In the context of continuous ensembles, working with $\N$-dependent equilibrium measures as a way to "resum many terms of an asymptotic expansion at once'' has also been used in \cite{BGKsinh,Kozlowski}.

 The definition of the equilibrium measure makes sense even if the linear forms
$(\mathfrak{r}_e)_{e = 1}^{\mathfrak{e}}$ appearing in \eqref{eq_equations_eqs} are such that
the space of integral solutions $\Lambda_\star$ is empty while the space of nonnegative real solutions $\hat{\Lambda}_\star^{\amsmathbb{R}}$ is nonempty. In particular, Theorem~\ref{Theorem_equi_charact} below still applies in such a case.

\begin{definition} \label{Definition_void_saturated} Given an equilibrium measure $\boldsymbol{\mu}$, we will use the following terminology for each $h \in [H]$.
\begin{itemize}
\item[$\bullet$] A \emph{void} is a maximal (for the inclusion) and closed interval in $[\hat{a}'_h,\hat{b}'_h]$ where
$\mu_h(x) =0$.
\item[$\bullet$] A \emph{saturation} is a maximal closed interval in
$[\hat{a}'_h,\hat{b}'_h]$ where $\mu_h(x) = \frac{1}{\theta_{h,h}}$.
\item[$\bullet$] A \emph{band} is a maximal open interval in $[\hat a'_h,\hat b'_h]$ where $0<\mu_h(x)<\frac{1}{\theta_{h,h}}$.
\end{itemize}
We denote by $\amsmathbb{V}_{h}$ (respectively $\amsmathbb{S}_{h}$, $\amsmathbb{B}_{h}$) the union of voids (respectively saturations, bands) in $[\hat{a}'_{h},\hat{b}'_{h}]$. If we speak of a \emph{saturated segment} (respectively, a \emph{void segment}), we mean a segment included in a saturation (respectively, in a void) but not necessarily maximal.
\end{definition}

If for instance $\mu_h$ has a single band inside $[\hat a'_h,\hat b'_h]$, then it can appear in four different configurations: \textit{saturation-band-saturation}, \textit{saturation-band-void}, \textit{void-band-saturation}, or \textit{void-band-void}.

\begin{definition}\label{def_eff_pot} Given an equilibrium measure $\boldsymbol{\mu}$, the $h$-th \emph{effective potential} is the function of $x \in [\hat{a}'_h,\hat{b}_h']$ defined through
 \begin{equation}
\label{eq_V_eff} V^{{\textnormal{eff}}}_h(x)= V_h(x) - \sum_{g =1}^H 2\theta_{h,g}
\int_{\hat a'_{g}}^{\hat b'_g} \log|x-y|\mu_{g}(y)\dd y.
\end{equation}
\end{definition}

The equilibrium measure can be characterized in the spirit of classical potential theory \cite{DS,ST} by the properties of its effective potential. We give a detailed proof in Section~\ref{Section_Energy_functional}.

\begin{theorem} \label{Theorem_equi_charact}
 Suppose that Assumptions~\ref{Assumptions_Theta} and \ref{Assumptions_basic} hold.
 Then the equilibrium measure $\mu$ has a compact
 support, \textit{i.e.} its support is inside $[-D,D]$ with $D>0$ depending only on the constants in the assumptions. Moreover, there exists an $H$-tuple of constants $\boldsymbol{v}$ such that for any $h \in [H]$ and $x \in [\hat{a}'_h,\hat{b}_h']$, we have
\begin{align}
V^{{\textnormal{eff}}}_h(x) \geq v_h& \quad \textnormal{if}\,\,x\,\,\textnormal{is in a void}, \label{eq_void_inequality_repeat}\\
 V^{{\textnormal{eff}}}_h(x) \leq v_h & \quad \textnormal{if}\,\,x\,\,\textnormal{is in a saturation},\label{eq_saturated_inequality_repeat}\\
 V^{{\textnormal{eff}}}_h(x) = v_h& \quad \textnormal{if}\,\,x\,\,\textnormal{is in a band}. \label{eq_band_equality_repeat}
\end{align}
If each interval $[\hat a'_h, \hat b'_h]$ has at
least one band, then $\boldsymbol{v}$ belongs to the $\amsmathbb{R}$-span of the $H$-dimensional vectors of coefficients of the linear forms $(\mathfrak{r}_{e})_{e = 1}^\mathfrak{e}$ in \eqref{eq_equations_eqs}, and the above conditions uniquely determine the equilibrium measure, \textit{i.e.} no other measure in $\mathscr{P}_\star$ satisfies them.
\end{theorem}
Our proof of the uniqueness of the measure satisfying \eqref{eq_void_inequality_repeat}-\eqref{eq_band_equality_repeat} relies on having at least one band in each interval $[\hat a'_h, \hat b'_h]$. If the entire segment $[\hat a'_h, \hat b'_h]$ is void or saturated, then the choice of the constant $v_h$ is not unique and the statement about the linear space to which $\boldsymbol{v}$ belongs would need further clarifications. We do not address such situations as they have little interest.

\medskip

We are now ready to formulate the \emph{off-criticality assumption}, \textit{cf.} \cite{BG11,BG_multicut}.

\begin{assumption} \label{Assumptions_offcrit}
 There exist $C > 0$, and a continuous function $\eps: \,\amsmathbb{R}_{\geq 0} \rightarrow \amsmathbb{R}_{\geq 0}$ such that $\eps(0) = 0$ and $\eps|_{\amsmathbb{R}_{> 0}} > 0$ making the following properties hold for any $\N \geq C$ and $h \in [H]$.
 \begin{enumerate}
 \item The segment $[\hat a'_h, \hat b'_h]$ contains at least one band. The bands in this segment are at distance at least $\frac{1}{C}$ of the endpoints. The total number of bands, voids and saturations of $\boldsymbol{\mu}$ is bounded by $C$, and the length of each of these regions is larger than $\frac{1}{C}$.
 \item $V^{{\textnormal{eff}}}_h(x) \geq v_h+\eps(\delta)$ for any $x \in \amsmathbb{V}_h$ at distance at least $\delta \geq 0$ from a boundary point of $\amsmathbb{V}_h$ which is neither $\hat a'_h$ nor $\hat b'_h$.
 \item $V^{{\textnormal{eff}}}_h(x) \leq v_h-\eps(\delta)$ for any $x \in \amsmathbb{S}_h$ at distance at least $\delta \geq 0$ from a boundary point of $\amsmathbb{S}_h$ which is neither $\hat a'_h$ nor $\hat b'_h$.
  \item For any $x \in \amsmathbb{B}_h$ at a distance at least $\delta$ from the boundary of $\amsmathbb{B}_h$, we have
  \[
  \eps(\delta) \leq \mu_h(x) \leq \frac{1}{\theta_{h,h}}- \varepsilon(\delta).
  \]
 \item For any point $x_0$ at the boundary of a void and a band in $[\hat a'_h, \hat b'_h]$, we have
 \[{ \frac{1}{C} \sqrt{|x-x_0|} \leq \mu_h(x) \leq C\sqrt{|x-x_0|}}\] for any $x$ in this band at distance at most $\frac{1}{C}$ from $x_0$.
 \item For any point $x_0$ at the boundary of a saturation and a band in $[\hat a'_h, \hat
 b'_h]$, we have \[{\frac{1}{C}\sqrt{|x-x_0|} \leq \frac{1}{\theta_{h,h}}-\mu_h(x) \leq C\sqrt{|x-x_0|}}\] for any $x$
in this band at distance at most $\frac{1}{C}$ from $x_0$.
 \end{enumerate}
\end{assumption}
In words, Conditions 1.\ and 2.\ mean that the inequalities of Theorem~\ref{Theorem_equi_charact} are strict inside voids and
saturations. Conditions 4.,5.,6.\ say that the density of the equilibrium measure does not vanish in the interior of the bands and has singularities of square-root type
 at the endpoints of bands. Such a behavior is believed to be generically true for many classes of potentials $V_h(x)$. In absence of saturations, this genericity property was proved for rescalings of real-analytic potentials $V_h(x)$ by \cite{Kuijlaarsgen} and the situation with saturations is discussed in \cite{colombo2024generic}. In practice, Assumption~\ref{Assumptions_offcrit} is easy to check with numerical simulations of the equilibrium measure. Part of it is automatic for strictly convex potentials, \textit{cf.} Proposition~\ref{prop:convexmueq}.

Similar off-criticality conditions are required for the continuous $\sbeta$-ensembles \cite{BG11,BG_multicut}. Conditions 2. and 4. guarantee that the different phases are pure in the sense that with overwhelming probability, there will be no particles in voids and particles are fully packed (no holes) in saturated segments. Conditions 5. and 6. conditions are technical conditions we need to prove the central limit theorem. We do not know whether a full central limit theorem holds without them; a partial result in this direction for continuous $\sbeta$-ensembles can be found in \cite{BLS}.

\subsection{Analyticity of the weights}\label{def_an}
\begin{assumption} \label{Assumptions_analyticity} There exist a constant $C>0$ independent of $\N$ and pairwise disjoint simply-connected open subsets $\amsmathbb M_h\subset \amsmathbb C$ independent of $\N$, indexed by $h \in [H]$ and containing $\big[\hat a_h-\frac{1}{C}, \hat b_h+\frac{1}{C}\big]$ with the following properties. For any $h \in [H]$, the weight in the $h$-th segment has a decomposition
 \begin{equation}
\label{eq_ansatzw}
w_h(x) = \exp\bigg(-\N U_h\Big(\frac{x}{\N}\Big)\bigg) \cdot \prod_{j=1}^{\iota_h^-} \frac{\N^{\,x - a_h + \rho^-_{h,j} - \frac{1}{2}}}{\Gamma\big(x - a_h + \rho^-_{h,j}\big)} \cdot \prod_{j=1}^{\iota_h^+} \frac{\N^{\,b_h-x + \rho^+_{h,j} - \frac{1}{2}}}{\Gamma\big(b_h - x + \rho^+_{h,j}\big)},
\end{equation}
with the following properties.
\begin{enumerate}
 \item $U_h(z)$ is a holomorphic function in $\left\{z\in \amsmathbb C \mid \textnormal{dist}(z,\amsmathbb M_h)\leq \frac{1}{2\N}\right\}$, and $|U_h(z)| \leq C$ in this domain.
 \item $\forall j \in [\iota_h^{\pm}]\qquad \rho^{\pm}_{h,j} \in \big(\frac{1}{C},C\big)$.
 \item The potential $V_h$ of Assumption~\ref{Assumptions_basic} is defined through
 \begin{equation}
 \label{eq_potential_weight_match}
 V_h(x)=U_h(x)+\iota_h^{-}\, \mathrm{Llog}(x - \hat{a}'_{h}) + \iota_{h}^+\, \mathrm{Llog}(\hat{b}'_h - x).
 \end{equation}
 \end{enumerate}
\end{assumption}
Let us compare \eqref{eq_ansatzw} with \eqref{eq_weight_form}. For this purpose we recall Stirling asymptotic formula
\begin{equation}
\label{eq_Stirling_basic}
 \log(\Gamma(z)) = \mathrm{Llog}(z) + \frac{1}{2}\log\bigg(\frac{2\pi}{z}\bigg)+o(1),
 \end{equation}
 where $o(1)$ tends to $0$ as $z \rightarrow +\infty$.  Applying this formula to \eqref{eq_ansatzw} and combining $O(\log x)$ remainders and finite shifts in $z$ into $\err_h(x)$, we arrive at \eqref{eq_weight_form} with $V_h(x)$ given by \eqref{eq_potential_weight_match}. Assumption~\ref{Assumptions_analyticity} will be intensively used in Chapters~\ref{ChapterNekra}, \ref{Chapter_smoothness} and \ref{Chapter_fff_expansions}--\ref{Chapter_filling_fractions}.
For many of our theorems, it would be sufficient to require analyticity of $U_h$ only in neighborhoods of bands, rather than neighborhoods of the entire segments $[\hat a_h,\hat b_h]$. Yet, we do not pursue this direction for lack of natural examples where this weaker assumption is useful.

\smallskip

 Our next task is to introduce several functions which will further enter into the exact functional relations between observables of the ensemble --- the Nekrasov equations written down in Chapter~\ref{ChapterNekra} --- that we exploit to establish asymptotic expansions as $\N \rightarrow \infty$. We also make explicit the asymptotic information about these functions which is implied by Assumption~\ref{Assumptions_analyticity}.

\begin{definition} \label{Definition_phi_functions} Given a weight $w_h(x)$ as in \eqref{eq_ansatzw} and $h\in[H]$ we define two holomorphic functions $\Phi^{\pm}_h(z)$ of $z \in \amsmathbb{M}_{h}$ satisfying
\begin{equation}
\label{eq_ratio_weights}
\forall x \in [\hat a_h,\hat b_h]\qquad \frac{w_{h}\big(\N x+\frac{1}{2}\big)}{w_{h}\big(\N x-\frac{1}{2}\big)} =
 \frac{\Phi^+_h(x)}{\Phi^-_h(x)}.
\end{equation}
They are specified by
\begin{equation}
\begin{split}
\label{eq_Phi_plus_def}
\Phi^+_h(z) & := e^{-\N[U_h(z+\frac{1}{2\N})-U_h(z-\frac{1}{2\N})]} \cdot \prod_{j=1}^{\iota_h^+} \bigg(\hat b_h - z-\frac{1}{2\N} + \frac{\rho^+_{h,j}}{\N}\bigg), \\
\Phi^-_h(z) & := \prod_{j=1}^{\iota_h^-} \bigg(z-\hat a_h -\frac{1}{2\N}+ \frac{\rho^-_{h,j}}{\N}\bigg).
\end{split}
\end{equation}
We also define holomorphic functions $\phi^\pm_h(z)$ of $z \in \amsmathbb{M}_{h}$ satisfying a $\N=\infty$-version of \eqref{eq_ratio_weights} (although they may still be $\N$-dependent as $V_h$ can depend on $\N$)
\begin{equation}
\label{eq_ratio_to_potential} \forall x \in [\hat{a}'_{h},\hat{b}'_{h}]\qquad
\frac{\phi^{+}_{h}(x)}{\phi^{-}_{h}(x)} = e^{-\partial_x V_h(x)}.
\end{equation}
They are specified by
\begin{equation}
\label{eq_phi_plus_minus_def}
\phi^+_h(z):= e^{-\partial_zU_h(z)} \cdot \big(\hat b'_h - z\big)^{\iota_h^+},\qquad
\phi^-_h(z):= \big(z-\hat a'_h\big)^{\iota_h^-}.
\end{equation}
\end{definition}

Let us emphasize that we use $\hat a_h$, $\hat b_h$ in the definition of $\Phi_h^\pm(z)$ and $\hat a'_h$, $\hat b'_h$ in the definition of $\phi^\pm_h(z)$. Given \eqref{eq_Phi_plus_def}, the property \eqref{eq_ratio_weights} is a corollary. Similarly, \eqref{eq_phi_plus_minus_def} implies \eqref{eq_ratio_to_potential}. While the product over $j$ parts of $\Phi^\pm_h(z)$ are uniquely fixed by \eqref{eq_ansatzw}, \eqref{eq_ratio_weights} and the request that $\Phi^\pm_h(z)$ has no poles, we have some freedom in distributing the part of the weight involving $U_h$ between $\Phi^+_h(z)$ and $\Phi^-_h(z)$ to obey \eqref{eq_ratio_weights}. With \eqref{eq_Phi_plus_def} we have chosen to make only $\Phi^+_h(z)$ dependent on $U_h(z)$, but other choices are also possible without incidence on the final results.

\begin{definition}
\label{Def_finer_Phi}
We introduce a further decomposition of holomorphic functions in $\amsmathbb{M}_h$
\begin{equation}
\label{Phifinerdec}
\Phi_h^{\pm}(z) = \phi_h^{\pm}(z) + \frac{\phi_{h}^{\pm,[1]}(z)}{\N} + \frac{\phi_h^{\pm,[2]}(z)}{\N^2} + \frac{\mathbbm{f}_h^{\pm}(z)}{\N^3}.
\end{equation}
where
\begin{equation}
\label{1storderphiexp}
\begin{split}
\phi_h^{+,[1]}(z) & =  \frac{\phi_h^+(z)}{\hat{b}'_h - z} \cdot \bigg(\sum_{j = 1}^{\iota_h^+} \big(\rho_{h,j}^+ - \theta_{h,h}\big)\bigg), \\
 \phi_h^{-,[1]}(z) & =  \frac{\phi_h^-(z)}{z - \hat{a}_h'} \cdot  \bigg(\sum_{j = 1}^{\iota_h^-} \big(\rho_{h,j}^- - \theta_{h,h}\big)\bigg), \\
\phi_h^{+,[2]}(z) & = \phi_h^{+}(z) \cdot \Bigg[-\frac{\partial_z^3 U_h(z)}{12} + \frac{1}{(\hat{b}'_h - z)^2} \cdot \bigg( \sum_{1 \leq j_1 < j_2 < \iota_h^+} (\rho_{h,j_1}^{+} - \theta_{h,h})(\rho_{h,j_2}^+ - \theta_{h,h})\bigg)\Bigg], \\
\phi_h^{-,[2]}(z) & = \frac{\phi_h^-(z)}{(z - \hat{a}_h')^2} \cdot \bigg(\sum_{1 \leq j_1 < j_2 \leq \iota_h^-} (\rho_{h,j_1}^{-} - \theta_{h,h})(\rho_{h,j_2}^- - \theta_{h,h})\bigg),
\end{split}
\end{equation}
and the error $\mathbbm{f}_h^{\pm}(z)$ is defined to make \eqref{Phifinerdec} holds.
\end{definition}
The expressions \eqref{1storderphiexp} arise by expanding $\Phi_h^{\pm}$ in the parameter $\frac{1}{\N}$ as $\N \rightarrow \infty$ beyond the leading order given by \eqref{eq_phi_plus_minus_def} always using the shifted parameters $\hat{a}_h',\hat{b}_h'$. There are no poles in $\phi_h^{\pm,[p]}$ for $p \in \{1,2\}$ because
 the denominator decreases the order of the zero (if any) already contained in  $\phi_h^{\pm}$. Let us collect helpful properties of the functions appearing in Assumption~\ref{Assumptions_analyticity} and Definition~\ref{Definition_phi_functions}.

\begin{lemma}\label{Lemma_phi_properties}
There exist a function $c\,:\,\amsmathbb{R}_{> 0} \rightarrow \amsmathbb{R}_{> 0}$ and for any compact $\amsmathbb{K}_h \subset \amsmathbb{M}_h$ a constant $C > 0$, both independent of $\N$, making the following properties hold for any $h \in [H]$ and $z \in \amsmathbb{K}_h$.
\begin{enumerate}
\item $\forall z \in \amsmathbb{K}_h\qquad \max\big(|\phi^{\pm}_{h}(z)|\,,\,|\phi^{\pm,[1]}_h(z)|\,,\,|\phi_h^{\pm,[2]}(z)|\,,\,|\mathbbm{f}_h^{\pm}(z)|\big) \leq C$.
\item For any $\eta > 0$ small enough independently of $\N$ and $x \in [\hat{a}_h' + \eta,\hat{b}_h' - \eta]$, we have $|\phi_h^{\pm}(x)| \geq c(\eta)$.
\item If $\iota_h^->0$, then $\phi_h^-(\hat{a}'_{h}) = 0$.
\item If $\iota_h^+>0$, then $\phi_h^+(\hat{b}'_{h}) = 0$.
\item If $\iota_h^{\pm} = 0$ or $\rho_{h,j}^{\pm} = \theta_{h,h}$ for every $j \in [\iota_h^\pm]$, then $\phi_h^{\pm,[1]}(z) = 0$.
\end{enumerate}
\end{lemma}

\begin{proof} The vanishing properties in 3.4.5. are obvious. The denominators in Definition~\ref{Def_finer_Phi} do not correspond to singularities because there are zeros of order $\iota_h^{\pm}$ in $\phi_h^{\pm}(z)$. Note that $\phi_h^{\pm,[1]}(z) = 0$ if $\iota_h^{\pm} = 0$, and the sums are absent in $\phi_h^{\pm,[2]}(z)$ if $\iota_h^{\pm} \leq 1$, that is
\[
\phi_h^{+,[2]}(z) = -\phi^+_h(z) \frac{\partial_z^3U_h(z)}{12}\,\, \textnormal{if}\,\,\iota_h^+ \leq 1 \qquad \textnormal{and} \qquad \phi_h^{-,[2]}(z) = 0\,\,\textnormal{if}\,\,\iota_h^- \leq 1.
\]
Since $U_h$ is holomorphic and uniformly bounded, so are the derivatives of $U_h$ thanks to the Cauchy residue formula. This implies the bounds in 1. and 2. except for the error term $\mathbbm{f}_h^{\pm}(z)$. For those we use a Taylor formula with remainder, which gives
\[
 \N \left[U_h\bigg(z+\frac{1}{2\N}\bigg)-U_h\bigg(z-\frac{1}{2\N}\bigg)\right] = \partial_zU_h(z) + \frac{\partial_z^3 U_h(z)}{12 \N^2} + O\bigg(\frac{1}{\N^3}\bigg),
\]
where the remainder is uniform in $z$ in any $\N$-independent compact of $\amsmathbb{M}_h$.
\end{proof}

\subsection{An auxiliary assumption}
\label{sec:auxiliary_assum}
An important step in our analysis is to study ensembles which satisfy certain extra assumptions that we now describe. Let us immediately emphasize that these extra assumptions do not restrict the generality of our results: we use them at some step of our arguments as they simplify the analysis, and at a later step we show how to get rid of them.

\begin{definition}
\label{bandlabel} We label the bands from $1$ to $K$ respecting their increasing order of appearance along the real line. The $k$-th band is denoted $(\alpha_k,\beta_k)$. It is contained in a segment $[\hat{a}_{h^k}',\hat{b}_{h^k}']$ for a unique $h^k \in [H]$. We call $\llbracket k^-(h),k^+(h)\rrbracket$ the set of indices of bands included in the $h$-th shifted segment $[\hat{a}_h',\hat{b}_h']$.
\end{definition}

\begin{assumption}
\label{Assumptions_extra} There exists a constant $C > 0$ independent of $\N$ making the following properties hold.
\begin{enumerate}
\item The constraints \eqref{eq_equations_eqs} fix deterministically the segment filling fractions $N_1,\ldots,N_H$.
\item For any $h \in [H]$, the equilibrium measure has only one band $(\alpha_{h},\beta_{h})$ in $[\hat{a}'_{h},\hat{b}'_{h}]$.
\item We have $|\hat{a}_{1}| < C$ and $|\hat{b}_{H}| < C$.
\item For any $h \in [H]$: if $\hat{a}'_h\in\amsmathbb{V}_{h}$, then we require $\iota_h^-=0$; if $\hat{b}'_h\in \amsmathbb{V}_{h}$, then we require $\iota_h^+=0$.
\item For any $h \in [H]$: if $\hat{a}'_{h} \in \amsmathbb{S}_{h}$, then we require $\iota_h^-=2$ and
 $\rho_{h,1}^-=1$, $\rho_{h,2}^-=\theta_{h,h}$. This implies in particular
\begin{equation}
\label{eq_Phi_minus_vanish}
\Phi_h^{-}\bigg(\hat{a}_{h} - \frac{1}{2\N}\bigg) = \Phi_h^{-}\big(\hat{a}'_{h}\big) =0, \qquad \phi_h^{-}\big(\hat{a}'_{h}\big)= \partial_z\phi_h^-(z)\big|_{z=\hat a'_h}= 0.
\end{equation}
If $\hat{b}'_{h} \in \amsmathbb{S}_{h}$, then we require $\iota_h^+=2$ and $
 \rho_{h,1}^+=1$, $\rho_{h,2}^+=\theta_{h,h}$. This implies in particular
  \begin{equation}
  \label{eq_Phi_plus_vanish}
  \Phi_h^{+}\bigg(\hat{b}_{h} + \frac{1}{2\N}\bigg) = \Phi_h^{+}\big(\hat{b}'_{h} \big) = 0, \qquad  \phi_h^{+}\big(\hat{b}'_{h}\big)= \partial_z \phi_h^+(z)\big|_{z=\hat b'_h} = 0.
  \end{equation}
If $\theta_{h,h}=1$, then in \eqref{eq_Phi_minus_vanish}-\eqref{eq_Phi_plus_vanish} the two zeros for $\Phi_h^\pm(z)$ merge into a double zero.
\item For any $h \in [H]$ and $z \in \amsmathbb{M}_{h}$, we have $|s_{h}(z)| \geq \frac{1}{C}$, where $s_h(z)$ is the auxiliary function defined in \eqref{eq_s_h_definition}, namely
\begin{equation} \label{eq_s_first_appearance}
\begin{split}
s_h(z) & = \frac{\phi_h^{+}(z)\cdot\exp\big(\sum_{g = 1}^H \theta_{h,g}\,\mathcal{G}_{\mu_{g}}(z)\big) - \phi_h^{-}(z)\cdot\exp\big(-\sum_{g = 1}^{H} \theta_{h,g}\,\mathcal{G}_{\mu_{g}}(z)\big)}{\big(z-\hat a'_h\big)^{\mathbbm{1}_{\amsmathbb{S}_h}(\hat{a}'_{h})} \cdot \big(z-\hat b'_h\big)^{\mathbbm{1}_{\amsmathbb{S}_h}(\hat{b}'_{h})} \cdot \prod_{k = k^{-}(h)}^{k^+(h)} \sqrt{(z - \alpha_k)(z - \beta_k)} },
\end{split}
\end{equation}
where $\mathcal{G}_{\mu_{g}}$ the Stieltjes transform of $\mu_{g}$, \textit{cf.} Section~\ref{Section_observables}.
\end{enumerate}
\end{assumption}
\begin{remark}
\label{remtheta11}We note that under Assumption~\ref{Assumptions_extra} and assuming $\theta_{h,h} = 1$, Condition 6. of Lemma~\ref{Lemma_phi_properties} is automatically satisfied, so the difference between $\Phi_h^{\pm}(z)$ and $\phi_h^{\pm}(z)$ is $O(\frac{1}{\N^2})$.
\end{remark}
In Chapter~\ref{Chapter_conditioning} we show that Assumption~\ref{Assumptions_extra} is not restrictive, namely we can always find an almost surely exhaustive family of conditionings of the ensemble which satisfy these extra assumptions. In other words, given an ensemble satisfying Assumptions~\ref{Assumptions_Theta}, \ref{Assumptions_basic}, \ref{Assumptions_offcrit} and \ref{Assumptions_analyticity}, we are always able to change the definition of the endpoints $(a_h,b_h)_{h = 1}^H$ so that the modified ensemble is undistinguishable from the original one with overwhelming probability but it satisfies Assumption~\ref{Assumptions_extra} as well. It turns out that formulae for asymptotic expansions will be simpler under this assumption. Nevertheless, our final results will not require Assumption~\ref{Assumptions_extra}.

\section{Observables of interest}
\label{Section_observables}

\subsection{Empirical measure and Stieltjes transforms}
\label{Def_Stiel_empi}
We summarize the definitions and notations for the observables that will be studied in this book, namely the empirical measure, its Stieltjes transform and its moments and cumulants. If $\boldsymbol{\nu}$ is an $H$-tuple of compactly-supported integrable measures on the real line, its Stieltjes transform is $\Gm_{\boldsymbol{\nu}}(z) = (\Gm_{\nu_h}(z))_{h = 1}^{H}$ where
\[
\forall h \in [H]\qquad \forall z \in \amsmathbb{C} \setminus \textnormal{supp}(\nu_h)\qquad \Gm_{\nu_h}(z) = \int_{\amsmathbb{R}} \frac{\dd \nu_h(x)}{z - x},
\]
\label{index:supp}where $\textnormal{supp}(\nu_h)$ is the support of $\nu_h$. The Stieltjes transform is a convenient tool to probe $\boldsymbol{\nu}$ against holomorphic test functions. Indeed, by Cauchy residue formula, it has the property that, for any holomorphic function $f_h$ in a complex neighborhood of $\textnormal{supp}\,\nu_h$ and any contour $\gamma_h$ which is included in this neighborhood, is oriented counterclockwise and surrounds $\textnormal{supp}\,\nu_h$, we have
\begin{equation}
\label{testest}
\int_{\amsmathbb{R}} f_h(x) \dd \nu_h(x) = \oint_{\gamma_h} \frac{\dd z}{2\ii\pi} f_h(z) \Gm_{\nu_h}(z).
\end{equation}
If we specialize $f$ to the constant function $1$, we get the mass of $\nu_h$, namely
\begin{equation}
\label{massnun}
\nu_h(\amsmathbb{R}) = \oint_{\gamma_h} \frac{\dd z}{2\ii\pi} \Gm_{\nu_h}(z) .
\end{equation}
This is equivalent to the statement that $\Gm_{\nu_h}(z) \sim \frac{\nu_h(\amsmathbb{R})}{z}$ as $z \rightarrow \infty$.

\begin{definition}
\label{def_empirical_mes}
Given a configuration $\boldsymbol{\ell} \in \W_\N$ in a discrete ensemble of Section~\ref{Section_general_model}, the \emph{empirical measure} is the $H$-tuple of nonnegative atomic measures $\boldsymbol{\mu}_{\N} = (\mu_{\N,h})_{h = 1}^{H}$ given by
\[
\forall h \in [H] \qquad \mu_{\N,h} = \frac{1}{\N} \sum_{i = 1}^{N_h} \delta_{\N^{-1}\ell_i^h}.
\]
This tuple $\boldsymbol{\mu}_{\N}$ is equivalently specified by the measure $\mu_{\N} = \sum_{h = 1}^{H} \mu_{\N,h} = \frac{1}{\N} \sum_{i = 1}^{N} \delta_{\N^{-1}\ell_i}$. The empirical resolvent is the tuple $\boldsymbol{G}(z) = (G_h(z))_{h = 1}^{H}$ defined for any $h \in [H]$ and $z \in\amsmathbb{C} \setminus [\hat{a}_h,\hat{b}_h]$ by
\begin{equation}
\label{Ghdefdef} G_h(z) = \sum_{i = 1}^{N_h} \frac{1}{z - \frac{\ell_i^h}{\N}} = \N \int_{\hat{a}_h}^{\hat{b}_h} \frac{\dd \mu_{\N,h}(x)}{z - x} = \N \Gm_{\mu_{\N,h}}(z).
\end{equation}
It is a holomorphic function of $z$ in this domain, such that for any $h \in [H]$, we have $G_h(z) \sim \frac{N_h}{z}$ as $z \rightarrow \infty$.

We sometimes recenter the empirical measure by subtracting the equilibrium measure $\boldsymbol{\mu}$ and introduce $\Delta \boldsymbol{G}(z) = \boldsymbol{G}(z) - \N \Gm_{\boldsymbol{\mu}}(z)$. In other words, for any $h \in [H]$ and $z \in \amsmathbb{C} \setminus \amsmathbb{A}_h^{\mathfrak{m}}]$ we have
\[
\Delta G_h(z) = G_h(z) - \N\Gm_{\mu_h}(z).
\]
\end{definition}
Under the law $\amsmathbb{P}_{\N}$ of the discrete ensemble --- \textit{cf.} \eqref{eq_general_measure} --- the empirical measure $\boldsymbol{\mu}_{\N}$ becomes an $H$-tuple of random nonnegative measures, and $\boldsymbol{G}(z)$ becomes a random $H$-tuple of holomorphic functions. We often test the --- recentered or not --- empirical measure against holomorphic test functions, using \eqref{testest} to express the outcome in terms of empirical Stieltjes transform.

\begin{definition}
\label{def_contourgammah}Consider a discrete ensemble satisfying Assumptions~\ref{Assumptions_analyticity} with complex domains $(\amsmathbb{M}_h)_{h = 1}^{H}$, and such that $a_1$ and $b_H$ are finite. For any $h \in [H]$, we denote $\gamma_h$ a contour included in $\amsmathbb{M}_h \setminus \amsmathbb{A}_h^{\mathfrak{m}}$ which is oriented counterclockwise and surrounds $\amsmathbb{A}_h^{\mathfrak{m}}$. Throughout this book, the contours $\gamma_h$ are always independent of $\N$ but may change from line to line --- in particular, whenever needed we can take contours closer and closer to $\gamma_h$ but still $\N$-independent.
\end{definition}

In Chapter~\ref{Chapter_smoothness} when we study abstractly the energy functional of Definition~\ref{Definition_functional}, we will consider the notion of variational datum, which is weaker than the notion of discrete ensemble, and for which only the segment $[\hat{a}_h',\hat{b}_h']$ is relevant. In this context, we only require that $\gamma_h$ of Definition~\ref{def_contourgammah} to surround $[\hat{a}_h',\hat{b}_h']$.

\subsection{Cumulants and correlators}

\label{Sec_cumu_corr}

If $(X_i)_{i = 1}^n$ is a $n$-tuple of complex random variables, we define their joint cumulant through
\begin{equation}
\label{eq_cumulant_def}
\E^{(\textnormal{c})}\big[X_1,\ldots,X_n\big] = \partial_{s_1} \cdots \partial_{s_n} \log \E\Bigg[\exp\bigg(\sum_{i= 1}^n s_i X_i\bigg)\Bigg]\Bigg|_{s_1 = \ldots = s_n = 0}.
\end{equation}
Joint cumulants are related to joint moments through
\[
\E \bigg[\prod_{i = 1}^n X_i \bigg] = \sum_{\amsmathbb{J} \vdash [n]} \prod_{J\, \in\, \amsmathbb{J}} \E^{(\textnormal{c})}\big[(X_j)_{j \in J}\big],
\]
where the sum ranges over all set partitions of $[n]$, \textit{i.e.} tuples of nonempty pairwise disjoint subsets of $[n]$. For instance
\begin{equation*}
\begin{split}
\E[X_1] & = \E^{(\textnormal{c})}[X_1], \\
\E\big[X_1X_2\big] & = \E^{(\textnormal{c})}[X_1,X_2] + \E^{(\textnormal{c})}[X_1] \cdot \E^{(\textnormal{c})}[X_2], \\
\E\big[X_1X_2X_3\big] & = \E^{(\textnormal{c})}[X_1,X_2,X_3] + \E^{(\textnormal{c})}[X_1,X_2]\cdot \E^{(\textnormal{c})}[X_3] + \E^{(\textnormal{c})}[X_1,X_3]\cdot \E^{(\textnormal{c})}[X_2] \\
& \quad + \E^{(\textnormal{c})}[X_2,X_3]\cdot \E^{(\textnormal{c})}[X_1] + \E^{(\textnormal{c})}[X_1]\cdot \E^{(\textnormal{c})}[X_2]\cdot \E^{(\textnormal{c})}[X_3].
\end{split}
\end{equation*}
For $n \geq 2$, adding a deterministic part to $X_1,\ldots,X_n$ does not change the value of the joint cumulant. In particular, if one of the $X_i$ is deterministic, we have $\E^{(\textnormal{c})}[X_1,\ldots,X_n] = 0$ for any $n \geq 2$. For more details on moments and cumulants see, \textit{e.g.}, \cite[Chapter 3]{PT}.

Most of our analysis of discrete ensembles beyond the leading order will be based on the analysis of the correlators, which we now define.

\begin{definition}
\label{def_correlators}
Let $n \geq 1$ and $h_1,\ldots,h_n \in [H]$. The \emph{$n$-point correlator} is defined as \emph{joint cumulant} of the $n$ random variables $\Delta
G_{h_1}(z_1),\ldots,\Delta G_{h_n}(z_n)$.
\begin{equation*}
\begin{split}
W_{n;h_1,\ldots,h_n}(z_1,\ldots,z_n) & = \E^{(\textnormal{c})}\big[\Delta
G_{h_1}(z_1),\ldots,\Delta G_{h_n}(z_n)\big] \\
& =\partial_{s_1}\cdots\partial_{s_n}\log \E\Bigg[ \exp\bigg(\sum_{i=1}^n s_i \Delta G_{h_i}(z_i)\bigg)\Bigg]\Bigg|_{s_1 = \cdots = s_n = 0}.
\end{split}
\end{equation*}
\end{definition}

\section{Main results}
\label{Section_results}
We state here the main results we prove for discrete ensembles in the subsequent chapters. We first recall the convergence of the $H$-tuple of empirical measures $\boldsymbol{\mu}_{\N} = (\mu_{\N,h})_{h = 1}^{H}$ from Definition~\ref{def_empirical_mes} towards the $H$-tuple of equilibrium measures $\boldsymbol{\mu} = (\mu_h)_{h = 1}^{H}$ from Definition~\ref{Definition_equilibrium_measure}. We recall that $\mu = \sum_{h = 1}^{H} \mu_h$.

\begin{theorem}[Law of large numbers] \label{Theorem_main_LLN} Suppose that Assumptions~\ref{Assumptions_Theta} and \ref{Assumptions_basic} hold.
Fix $h \in [H]$ and let $f$ be a bounded continuous function
defined on $\amsmathbb{A}_h^{\mathfrak{m}}$ for $\N$ large enough. Then
\[
 \lim_{\N\rightarrow\infty} \Bigg| \frac{1}{\N} \sum_{i=1}^{N_h} f\bigg(\frac{\ell_i^h}{\N}\bigg) -\int_{\hat
{a}'_h}^{\hat{b}_h'} f(x) \dd\mu_h(x)\Bigg|=0\qquad \text{in probability.}
\]
 In terms of the Stieltjes transform (Definition~\ref{def_empirical_mes}) for each $\N$-independent  $\varepsilon > 0$ and compact set $\amsmathbb{K}_h$ containing $[\hat{a}_h,\hat{b}_h]$  in its interior for $\N$ large enough, we have
\[
\lim_{\N\rightarrow\infty} \sup_{z\in \amsmathbb{C} \setminus \amsmathbb{K}_h}\,\P\Bigg[\bigg| \frac{G_h(z)}{\N} - \Gm_{\mu_h}(z)\bigg|>\eps\Bigg] =0.
\]
Both limits are uniform over the ensembles satisfying Assumptions~\ref{Assumptions_Theta} and \ref{Assumptions_basic} with fixed constant $C$.
\end{theorem}

This is an immediate consequence of  Lemma~\ref{Lemma_tail_bound_general} and Corollary~\ref{Corollary_a_priori_0}, which moreover contain quantitative estimates for the speed of convergence. We insist on the fact that $\mu_h$ may still depends on $\N$ through the potential $V_h$, the
intervals $[\hat{a}'_h,\hat{b}'_h]$ and the segment filling fractions $\hat{n}_h$ of Definitions~\ref{def:eq_rescaled_parameters}-\ref{def:eq_shifted_parameters}. Thus, if
one wants to turn Theorem~\ref{Theorem_main_LLN} into a statement of convergence
of linear statistics to a limit, then one needs to supply it with statements describing exactly the dependence of the parameters of the ensemble in $\N$. As we show in Chapter~\ref{Chapter_smoothness}, the dependence of $\mu_h$ on $V_h$, $\hat
a_h',\hat{b}_h',\hat{n}_h$ is continuous in an appropriate topology.

The next theorem describes global fluctuations of the empirical measure when the filling fractions are fixed and each segment has exactly one band. A particular case of this result was already derived in \cite{BGG}. It is established in Chapter~\ref{Chapter_fff_expansions} where the dependencies in the parameters of the ensembles are precised in view of the forthcoming analysis, see Theorem~\ref{Theorem_correlators_expansion}, Corollaries~\ref{Theorem_correlators_expansion_theta1}-\ref{Corollary_CLT}, Theorems~\ref{Theorem_correlators_expansion_relaxed}-\ref{Theorem_correlators_expansion_relaxed_theta1} and Corollary~\ref{Corollary_CLT_relaxed}. We only state here a summary of this last corollary.

\begin{theorem}[Central limit theorem for fixed filling fractions]
\label{Theorem_CLT_fff_intro}
Suppose that Assumptions~\ref{Assumptions_Theta}--\ref{Assumptions_analyticity} hold. In addition, assume that for each $h \in [H]$, $\mu_h$ has exactly one band and the equations \eqref{eq_equations_eqs} deterministically fix the filling fraction $\hat{n}_h$. Then, if $f(z)$ is a $\N$-independent holomorphic function of $z$ in a neighborhood of $\amsmathbb{A}$, the linear statistics
\begin{equation}
\label{thedsidun}
\sum_{i=1}^{N} f\bigg(\frac{\ell_i}{\N}\bigg) - \N \int f(x)\dd \mu(x)
\end{equation}
is asymptotically Gaussian as $\N \rightarrow \infty$, in the sense that its cumulants of order $\geq 3$ tend to $0$ as $\N \rightarrow \infty$. Furthermore, if $\theta_{g,g} = 1$ for every $g \in [H]$, then this Gaussian variable is asymptotically centered: the mean is $O(\frac{1}{\N})$ as $\N \rightarrow \infty$.
\end{theorem}
It is a feature of the discrete ensemble that the mean (present if some $\theta_{h,h} \neq 1$) and the covariance may depend on $\N$ but do so in a bounded way. Up to $o(1)$ we find the covariance of the linear statistics \eqref{thedsidun} in the form
\[
\textsf{Cov}[f,f] = \sum_{h_1,h_2 = 1}^{H} \oint_{\gamma_{h_1}} \oint_{\gamma_{h_2}} \frac{\dd z_1 \dd z_2}{(2\ii\pi)^2}\,\mathcal{F}_{h_1,h_2}(z_1,z_2)\,f(z_1)f(z_2),
\]
where $\gamma_{h}$ is a contour surrounding the band $(\alpha_h,\beta_h)$ in $[\hat{a}_h',\hat{b}'_h]$ and
\begin{equation}
\label{FFFFFFFFFFF}
\mathcal{F}_{h_1,h_2}(z_1,z_2) = \Upsilon_{h_1}\Bigg[\bigg(\frac{\sqrt{(* - \alpha_{h_2})(* - \beta_{h_2})}}{2(* - z_2)^2} \delta_{g,h_2}\bigg)_{g = 1}^{H}\,;\,\boldsymbol{\kappa} = \boldsymbol{0}\Bigg](z_1)
\end{equation}
is constructed via the linear operator $\Upsilon$ solving the so-called master problem and studied in Chapter~\ref{Chapter_SolvingN}. As a matter of fact, this is a symmetric function of $(h_1,z_1)$ and $(h_2,z_2)$ (Theorem~\ref{thm:Bsym}). It plays the role of a fundamental solution for the master problem, from which all other solutions can be reconstructed (\textit{cf.} Section~\ref{Section_Master_by_covariance}). This function only depends on $\boldsymbol{\Theta}$ and the endpoints of the bands. We give an explicit formula when  $\boldsymbol{\Theta}$ diagonal or has all entries equal (Theorem~\ref{Theorem_Masterspecial}). Such formulae are already known for continuous $\sbeta$-ensembles, but is scattered among many references, so we decided to include a unified derivation (Section~\ref{Section_Explicit}). Chapter~\ref{Chapter_AG} develops tools that lead to formulae for more general but still algebraic cases (this can occur only for a discrete set of $\boldsymbol{\Theta}$s). These algebraic formulae are illustrated with examples in Chapter~\ref{Chap14}, covering for instance the leading covariance for tilings in domains of shape E, C, S, G, O or in a M\"obius band. For truly general $\boldsymbol{\Theta}$ we do not have other formulae to propose than \eqref{FFFFFFFFFFF}.

Next, we describe our results when the assumptions of having one band per segment and deterministic filling fractions are relaxed. To this end, we first show that away from a small neighborhood of bands, the configurations have no particles or are densely packed with probability exponentially close to $1$ as $\N\rightarrow\infty$.
\begin{definition}
 Let $1 \leq p <q \leq N$. We say that the particles $\ell_p,\ell_{p + 1},\ldots,\ell_{q}$ are \emph{densely packed} if  there exists $h \in [H]$ such that they all belong to $[a_h,b_h]$ and $\ell_{i+1}-\ell_i=\theta_{h,h}$ for any $i \in \llbracket p,q-1 \rrbracket$.
\end{definition}

\begin{theorem} \label{Theorem_voids}
 Suppose that Assumptions~\ref{Assumptions_Theta}, \ref{Assumptions_basic} and \ref{Assumptions_offcrit}
 hold. If $[\mathfrak{a},\mathfrak{b}]$ is included in a void of the equilibrium measure, then for any $\eps>0$ there exists a real $C > 0$ independent of $\N$ such that the probability to have at least one particle $\ell_i$ inside $[\N(\mathfrak{a}+\eps),
 \N(\mathfrak{b}-\eps)]$ is at most $C^{-1}e^{-C\N}$.
 \end{theorem}

\begin{theorem} \label{Theorem_densely_packed}
 Suppose that Assumptions~\ref{Assumptions_Theta}, \ref{Assumptions_basic} and
 \ref{Assumptions_offcrit}
 hold. Let $[\mathfrak{a},\mathfrak{b}]$ be a saturation of the equilibrium measure. For each $\eps>0$ there
 exist $\N$-dependent integers $p < q$ and a real $C>0$ independent of $\N$ such that with
 probability greater than $1-C^{-1}e^{-C \N}$ the following two conditions hold:
 \begin{itemize}
  \item $\N \mathfrak{a} \leq \ell_{p}< \N(\mathfrak{a}+\eps)$ and $\N (\mathfrak{b}-\eps) <\ell_{q} \leq \N\mathfrak{b}$;
  \item The particles $\ell_p,\ell_{p+1},\ldots,\ell_{q}$ are densely packed.
 \end{itemize}
\end{theorem}
If $\mathfrak{a} = \hat{a}_h'$ (or $\mathfrak{b} = \hat{b}_h'$) then the densely packed configuration extends up to this endpoint, \textit{i.e.} the $\eps$-error is no longer needed. These two theorems will be proved as Theorem~\ref{Theorem_ldpsup} and Theorem~\ref{Theorem_ldsaturated} in the text. It should be possible to adapt our proofs to give them a more precise form, namely a large deviation principle for the presence of a particle with good rate function $V_h^{\textnormal{eff}} - v_h$, or for the presence of a hole with good rate function $v_h - V_h^{\textnormal{eff}}$. We have however not followed this route and the present form is sufficient for our purposes. In the case $H = 1$ such a large deviation principle with good rate functions has been established in \cite{DimitrovZhang} and we expect the same strategy to apply as well for $H > 1$.

If a discrete ensemble has several bands per segment but no saturation, we condition the ensemble to have particles only in new segments that are small neighborhoods of the bands. According to Theorem~\ref{Theorem_voids} we only pay the price of exponentially small errors as $\N \rightarrow \infty$. In particular, we are now back to having one band per segment, but the filling fractions of those new segments can fluctuate. Our next result states that these fluctuations are asymptotically described by a discrete Gaussian. A few comments are in order before making this precise.  By Assumption~\ref{Assumptions_offcrit}, the equilibrium measure $\mu$ has a finite number $K$ of bands. Let us denote them $(\alpha_k,\beta_k)$ for $k \in[K]$, in increasing order along the real line, and $h^k \in [H]$ the index of the segment in which the $k$-th band is included.
Let us choose $\eps > 0$ small enough independent of $\N$, depending only on the constants in the assumptions and such that the intervals $[\alpha_k-\eps,\beta_k+\eps]$ are pairwise disjoint for $k \in [K]$. We define the (random)\label{index:fluctn} \emph{fluctuating filling fractions} through
\begin{equation}
\label{Nbulleki} \forall k \in [K] \qquad N^{\circ}_k= \#\big\{i \in [N] \quad \big| \quad \ell_i \in (\N(\alpha_k-\eps),\N(\beta_k+\eps))\big\}, \qquad \hat{n}_k^{\circ} = \frac{N_k^{\circ}}{\N}.
\end{equation}
By Theorem~\ref{Theorem_main_LLN}, when $\N$ is large the rescaled filling fraction $\hat{n}_k^{\circ} $ is close to the deterministic value
\begin{equation}
\label{eq_def_equilibrium_ffi}
\hat{n}^{\boldsymbol{\mu}}_{k}=\frac{N^{\boldsymbol{\mu}}_{k}}{\N} := \mu([\alpha_k,\beta_k]).
\end{equation}
We will study the asymptotic distribution of the random integral vector $\boldsymbol{N}^{\circ} = (N^{\circ}_k)_{k = 1}^K$ as $\N\rightarrow\infty$. By Theorem~\ref{Theorem_voids}, with overwhelming probability as $\N \rightarrow \infty$ there are no particles outside the region ${\bigcup_{k=1}^K [\N(\alpha_k-\eps),\N(\beta_k+\eps)]}$. Thus the affine constraints \eqref{eq_equations_eqs} in Section~\ref{DataS} on $\boldsymbol{N} = (N_h)_{h = 1}^{H}$ can be rewritten as affine constraints on $\boldsymbol{N}^{\circ} = (N_{k}^{\circ})_{k = 1}^K$. They define an affine subspace $\L \subset \amsmathbb{R}^{K}$, which is parallel to a unique linear subspace $\L_0$ (the difference of two vectors in $\L$ is in $\L_0$). A \emph{discrete Gaussian} random variable $\textnormal{\textsf{\textbf{Gau\ss{}}}}_{\amsmathbb{Z}}[\Qu,
\L,\boldsymbol{u}]$ is a random variable taking values in the set
$\L \cap \amsmathbb{Z}^K $ and such that
\begin{equation} \label{eq_Discrete_Gaussian_maini}
\forall \boldsymbol{x} \in \L \cap \amsmathbb{Z}^{K}\qquad \amsmathbb{P}\big(\textnormal{\textsf{\textbf{Gau\ss{}}}}_{\amsmathbb{Z}}[\Qu,
\L, \boldsymbol u] = \x\big) = \frac{1}{\mathscr{Q}} \exp\left(-
\frac{1}{2}\Qu(\x - \boldsymbol u,\x - \boldsymbol u) \right),
\end{equation}
where $\mathscr{Q} > 0$ is the normalizing constant turning \eqref{eq_Discrete_Gaussian_maini} into a probability measure.

Next, we specify the parameters of the discrete Gaussian of relevance to the results. We first introduce
the equilibrium measure $\boldsymbol{\mu}^{\hat{\boldsymbol p}}$ which minimizes
 the energy functional $-\I$ defined in \eqref{eq_functional_general} over tuples of measures with filling fractions equal to a fixed $\hat{\boldsymbol{p}}$ in a small $\N$-independent neighborhood of $\hat{\boldsymbol{n}}^{\boldsymbol{\mu}}$. The quadratic form we need is $\boldsymbol{B}(\boldsymbol{x},\boldsymbol{x}) = \sum_{k,l = 1}^{K} B_{k,l} x_k x_l$ with
\begin{equation}
\label{eq_discrete_covariancei}
\forall k,l \in [K] \qquad B_{k,l} = - \partial_{\hat{p}_k}\partial_{\hat{p}_l} \I[\boldsymbol{\mu}^{\hat{\boldsymbol p}}]\big|_{\hat{\boldsymbol p}=\hat{\boldsymbol n}^{\boldsymbol{\mu}}}.
\end{equation}
Proposition~\ref{Proposition_Hessian_free_energy} shows that the restriction of $\Qu$ to $\L_0$ is positive definite. The definition of the vector $\boldsymbol u$ entering into \eqref{eq_Discrete_Gaussian_maini} is more delicate and rather implicit. Yet, it can be checked that
\begin{equation}
\label{eq_discrete_centeringi}
\boldsymbol{u}=\N \hat{\boldsymbol n}^{\boldsymbol{\mu}}+\log \N\,\boldsymbol{u}^{(1)} +\boldsymbol{u}^{(0)},
\end{equation}
where $\boldsymbol{u}^{(1)}$ is explicit and $\boldsymbol u^{(0)}$ is uniformly bounded, \textit{cf.} \eqref{eq_x234}. In particular, if $\theta_{h,h}$ does not depend on $h$, then $\boldsymbol{u}^{(1)} = 0$. We stress that $\Qu$ and $\boldsymbol{u}^{(1)},\boldsymbol{u}^{(0)}$ depend on the parameters of the ensemble, in particular they may depend on $\N$ but in a bounded way. As a consequence, the convergence in law towards such a random variable can only be described in terms of a notion of asymptotic equivalence that we now define.

\begin{definition}
\label{Definition_asymptotic_equivalencei}
Let $\boldsymbol{\chi}_{\N}^{(1)}$ and $\boldsymbol{\chi}_{\N}^{(2)}$ be two families (indexed by a large parameter $\N$) of discrete random $K$-dimensional real vectors. We say that the distribution of $\boldsymbol{\chi}_{\N}^{(1)}$ is asymptotically equal to the distribution of $\boldsymbol{\chi}_{\N}^{(2)}$ if
\begin{equation}
\label{eq_asymptotic_equivalencei}
 \lim_{\N\rightarrow\infty} \sup_{f \in \textnormal{Lip}_1} \big|\E[f(\boldsymbol{\chi}_{\N}^{(1)})]-\E[f(\boldsymbol{\chi}_{\N}^{(2)})]\big| =0,
\end{equation}
where $\textnormal{Lip}_1$ is the set $1$-Lipschitz functions real-valued functions on $\amsmathbb{R}^K$ that are upper-bounded by $1$. If \eqref{eq_asymptotic_equivalence} holds, we write $\chi_{\N}^{(1)} \stackrel{\textnormal{d}}{\sim} \chi_{\N}^{(2)}$.
\end{definition}

\begin{theorem}[Discrete Gaussian for fluctuating filling fractions] \label{Theorem_CLT_for_filling_fractionsi} If the discrete ensemble satisfies Assumptions~\ref{Assumptions_Theta}, \ref{Assumptions_basic}, \ref{Assumptions_offcrit}, \ref{Assumptions_analyticity} and its equilibrium measure has no saturations, then the distribution of $\boldsymbol{N}^{\circ}$ is asymptotically equal to the distribution of the discrete Gaussian $\textnormal{\textsf{\textbf{Gau\ss{}}}}_{\amsmathbb{Z}}[\Qu,
\L,\boldsymbol{u}]$ with the parameters discussed above.
\end{theorem}

This is proved in the text as Theorem~\ref{Theorem_CLT_for_filling_fractions}. Combining Theorem~\ref{Theorem_CLT_for_filling_fractionsi} with our previous Theorem~\ref{Theorem_CLT_fff_intro} on the fixed filling fractions case, we can obtain an asymptotic theorem for the fluctuations of linear statistics. Take a holomorphic function $f(z)$ defined in a complex neighborhood of the $k$-th band $(\alpha_k,\beta_k)$ and consider the random variable:
\begin{equation}
\label{eq_x235i}
\mathsf{Lin}_{k}[f] = \sum_{i = 1}^{N} \mathbbm{1}_{[\N(\alpha_k - \eps),\, \N(\beta_k + \eps)]}(\ell_i)\,f\bigg(\frac{\ell_i}{\N}\bigg).
\end{equation}
For $\eps$ small enough but fixed and $\N$ large enough, Theorem~\ref{Theorem_voids} shows that with overwhelming probability $\textnormal{\textsf{Lin}}_k[f]$ does not depend on $\eps$, because the particles remain in a small neighborhood of the bands. Its fluctuations have two origins: the fluctuation of filling fractions $\boldsymbol{N}^{\circ}$ leading to a varying number of non-zero terms in \eqref{eq_x235i}, and the fluctuation of the sum itself conditionally to the filling fractions. To describe the first, we introduce the deterministic quantity
\begin{equation}
\label{omegadefdef}
\forall k,m \in [K]\qquad \omega_{k,m}[f] = \int_{\alpha_k}^{\beta_k} f(x) \big(\partial_{\hat{p}_m} \mu^{\hat{\boldsymbol{p}}}(x)_{h^k}\big)\big|_{\hat{\boldsymbol{p}} = \hat{\boldsymbol{n}}^{\boldsymbol{\mu}}} \, \dd x.
\end{equation}
Equation~\eqref{formulaehhh} allows computing it in terms of the so-called first-kind functions, themselves defined via the fundamental solution $\boldsymbol{\mathcal{F}}$ of \eqref{FFFFFFFFFFF} or Definition~\ref{def:1stkind}).

\begin{theorem}
\label{Theorem_linear_statistics_fluctuation_ffi}
Suppose that Assumptions~\ref{Assumptions_Theta}, \ref{Assumptions_basic}, \ref{Assumptions_offcrit}, \ref{Assumptions_analyticity} hold and that the equilibrium measure has no saturations. Let $k_1,\ldots,k_L \in [H]$ be a $\N$-independent tuple and call $h_l = h^{k_l}$. Let $\boldsymbol{f}(z) = (f_l(z))_{l = 1}^{L}$ be functions such that $f_l(z)$ is holomorphic for $z$ in some $\N$-independent compact neighborhood of $[\alpha_{k_l} - \eps,\beta_{k_l} + \eps]$ where it is bounded independently of $\N$, for any $l \in [L]$. Let $\textnormal{\textsf{\textbf{Gau\ss{}}}}[\boldsymbol{f}]$ be a $L$-dimensional random Gaussian vector with covariance
\[
\textnormal{\textsf{Cov}}_{l_1,l_2}[\boldsymbol{f},\boldsymbol{f}] = \oint_{\gamma_{l_1}} \oint_{\gamma_{l_1}} \frac{\dd z_1\dd z_2}{(2\ii\pi)^2} \mathcal{F}_{l_1,l_2}(z_1,z_2) f_{l_1}(z_1) f_{l_2}(z_2)
\]
in terms of the fundamental solution \eqref{FFFFFFFFFFF} associated to the bands of $\mu$ and the matrix $(\theta_{h_{l_1},h_{l_2}})_{l_1,l_2 = 1}^{L}$. Then, as  $\N \rightarrow \infty$ we have
\[
\bigg( \textnormal{\textsf{Lin}}_{k_l}[f_l] - \N \int_{\alpha_{k_l}}^{\beta_{k_l}} f_l(x)\mu_{h_l}(x)\dd x \bigg)_{l = 1}^{L} \,\, \mathop{\sim}^{\textnormal{d}}\,\,\textnormal{\textbf{\textsf{Gau\ss{}}}}[\boldsymbol{f}] + \bigg( \sum_{m = 1}^{L} \big(\textnormal{\textsf{Gau\ss{}}}_{\amsmathbb{Z},m} - \N \hat{n}^{\boldsymbol{\mu}}_m\big) \cdot  \omega_{k_l,m}[f_l]\bigg)_{l = 1}^{L}.
\]
where the discrete Gaussian $\textnormal{\textsf{\textbf{Gau\ss{}}}}_{\amsmathbb{Z}} = \textnormal{\textsf{\textbf{Gau\ss{}}}}_{\amsmathbb{Z}}[\Qu,\L,\boldsymbol{u}]$ is the same as in Theorem~\ref{Theorem_CLT_for_filling_fractionsi}.
\end{theorem}
This is proved in Theorem~\ref{Theorem_linear_statistics_fluctuation_ff}. Last but not least, we extend this result in presence of saturations in Section~\ref{Section_Fillingfractionsatcasesec}. This is the general case where each of the $H$ segments can have a succession of voids, bands and saturations in any possible order. As part of the off-criticality Assumption~\ref{Assumptions_offcrit}, only voids and saturations are allowed to touch the endpoints of the segment. A new difficulty arises as \eqref{eq_x235i} will depend on $\eps$ if the $k$-th band is the neighbor of a saturation. There is even a second complication if $\theta_{h,h} \neq 1$, as the sites allowed for a particle in $[a_h,b_h]$ depend on the number of particles to its left in this segment. To overcome these difficulties, we construct a different discrete random variable $\boldsymbol{Q}^{\circ}$. Instead of looking at filling fractions in fixed segments as in \eqref{Nbulleki}, we let the endpoints of the segments slide randomly and count the number of particles in those random segments, in a way that is driven deterministically by $\boldsymbol{Q}^{\circ}$. We call $\overline{\boldsymbol{N}}^{\circ} = (\hat{\boldsymbol{a}},\hat{\boldsymbol{b}},\boldsymbol{N}^{\circ},\boldsymbol{Q}^{\circ})$ the tuple of \emph{extended filling fractions}. Its core properties are established in Lemma~\ref{Lemma_fluct_params_choice} and Theorem~\ref{Theorem_CLT_for_filling_fractions_saturation} which we can summarize as follows.

\begin{theorem}[Discrete Gaussian for extended filling fractions]
\label{discGextei}
Suppose that Assumptions~\ref{Assumptions_Theta}, \ref{Assumptions_basic}, \ref{Assumptions_offcrit}, \ref{Assumptions_analyticity} hold and that the equilibrium measure has at least one saturation. There exists $\eps > 0$ and random $K$-tuples of integers $\boldsymbol{Q}^{\circ},\boldsymbol{N}^{\circ}$ and $K$-tuple of real numbers $\hat{\boldsymbol{a}}^{\circ},\hat{\boldsymbol{b}}^{\circ}$ such that the following properties hold with probability exponentially close to $1$ as $\N \rightarrow \infty$.
\begin{enumerate}
\item $\boldsymbol{N}^{\circ},\hat{\boldsymbol{a}}^{\circ}$, $\hat{\boldsymbol{b}}^{\circ}$ are deterministic affine functions of $\boldsymbol{Q}^{\circ}$, with coefficients of the linear part depending only on the diagonal of $\boldsymbol{\Theta}$;
\item the segment $[\hat{a}_k^{\circ},\hat{b}_k^{\circ}]$ contains the $k$-th band and is contained in its $\eps$-neighborhood;
\item the segment $[\N\hat{a}_k^{\circ},\N\hat{b}_k^{\circ}]$ contains exactly $N_k^{\circ}$ particles, the leftmost and rightmost ones being at integer distance of its endpoints.
\end{enumerate}
Besides, the distribution of $\boldsymbol{Q}^{\circ}$ is asymptotically equal to a discrete Gaussian distribution $\textnormal{\textsf{\textbf{Gau\ss{}}}}_{\amsmathbb{Z}}$ as $\N \rightarrow \infty$, whose quadratic form $\boldsymbol{B}$ and centering vector $\boldsymbol{u}$ depend on the parameters of the ensemble.
\end{theorem}
The parameters $\boldsymbol{B}$ and $\boldsymbol{u}$ are obtained as in Theorem~\ref{Theorem_CLT_for_filling_fractionsi} (\textit{cf.} \eqref{eq_discrete_covariancei}-\eqref{eq_discrete_centeringi}) but where $\hat{\boldsymbol{p}}$ represents the value of $\frac{\boldsymbol{Q}^{\circ}}{\N}$ instead of the usual filling fractions. From Theorem~\ref{discGextei}, we derive the fluctuations of the linear statistics in Theorem~\ref{Theorem_linear_statistics_fluctuation_sat}, which take the following form.
 \begin{theorem}[Perturbed central limit theorem]
\label{Theorem_linear_statistics_fluctuation_sati}
In the setting of Theorem~\ref{Theorem_linear_statistics_fluctuation_ffi} but with at least one saturation, we have as $\N \rightarrow \infty$
\[
\bigg(\textnormal{\textsf{Lin}}_{k_l}[f_l] - \N\int_{\alpha_{k_l}}^{\beta_{k_l}} f_l(x)\mu_{h_l}(x)\dd x\bigg)_{l = 1}^{L} \,\,\mathop{\sim}^{\textnormal{d}} \,\,\textnormal{\textsf{\textbf{Gau\ss{}}}}[\boldsymbol{f}] + \big(\textnormal{\textsf{Shift}}_{h_l}[f_l]\big)_{l = 1}^{L}.
\]
The first part is the random Gaussian vector from Theorem~\ref{Theorem_linear_statistics_fluctuation_ffi}. The second part is a deterministic function of the discrete Gaussian vector $\textnormal{\textsf{\textbf{Gau\ss{}}}}_{\amsmathbb{Z}}$ taken from Theorem~\ref{discGextei}: it is specified in Definition~\ref{DhfL} as sum of a term like the second one in Theorem~\ref{Theorem_linear_statistics_fluctuation_ffi} --- with $\hat{n}^{\boldsymbol{\mu}}_m$ replaced by the equilibrium value of $\frac{\hat{Q}_m^{\circ}}{\N}$ --- and a term vanishing in case $\theta_{g,g} = 1$ for every $g \in [H]$.
\end{theorem}

An important intermediate result to reach Theorems~\ref{discGextei}-\ref{Theorem_linear_statistics_fluctuation_sati} is to estimate the partition functions up to order one. Such results are obtained by interpolating the discrete ensemble under consideration with discrete ensemble for which the partition functions can be estimated, see Theorem~\ref{Theorem_partition_multicut}.

\begin{theorem}[Asymptotics of free energy]
\label{Theorem_partition_multicuti}
Suppose that Assumptions~\ref{Assumptions_Theta}, \ref{Assumptions_basic}, \ref{Assumptions_offcrit}, \ref{Assumptions_analyticity} hold and that there is a single band per segment. Then the partition function \eqref{eq_partition_function_definition} admits an expansion
\begin{equation}
 \label{eq_partition_multicut_prime}
\lim_{\N \rightarrow \infty} \Bigg|\log\Z_\N - \bigg[\I[\boldsymbol{\mu}]\N^2 + \bigg( \sum_{h = 1}^{H} \theta_{h,h}\hat{n}_h\bigg) \N\log \N + \mathbbm{Rest}_1\N\bigg]\Bigg| = 0,
 \end{equation}
where $\mathbbm{Rest}_1$ is a twice-continuously differentiable function in the parameters of the ensemble, bounded and with first and second partial derivatives bounded independently of $\N$.
 \end{theorem}
We contented ourselves (and it is sufficient to get the perturbed central limit theorems) to obtain this result under the assumption of having one band per segment. It is always possible to reduce to this setting by localizing the ensembles near the bands, thus removing pieces of saturations. Up to exponentially small corrections, the partition function after removal is the partition function before removal times the interaction and weight factors for the fully packed configuration of particles in the saturated segments that were removed. The large $\N$ asymptotic expansion of the latter can be extracted from Lemma~\ref{Lemma_densely_packed_expansion} for the interaction part and \eqref{phiexp}-\eqref{gnexpmfgun} for the weight part. The final result is the same as in Theorem~\ref{Theorem_partition_multicuti}, only giving additional contributions to the  $\mathbbm{Rest}_1$ that we have anyway kept implicit.

\chapter{Nekrasov equations}
\label{ChapterNekra}
\section{First-order Nekrasov equation}
\label{SNEK1}

Our main tool for the analysis of the discrete ensembles presented in Chapter~\ref{Chapter_Setup_and_Examples} is a
system of $H$ equations which generalizes the Nekrasov equations of
\cite[Theorem~4.1]{BGG}. They find their origin in the context of supersymmetric gauge theories in \cite{Nekrasovpaper}, where they were called ``non-perturbative" Dyson--Schwinger equations. We stress that these equations are exact, in the sense that they are valid for each $\N$, before taking any asymptotics. In this chapter, we silently assume that the partition function $\Z_\N$ is finite --- by Proposition~\ref{lem_finite_partition_function} this is always true under Assumptions~\ref{Assumptions_Theta} and \ref{Assumptions_basic} which we are going to impose in the next chapters.

\begin{theorem}
\label{Theorem_Nekrasov}
Consider a discrete ensemble as in Section~\ref{Section_general_model}. Fix $h \in [H]$ and suppose that there exist two holomorphic functions $\Phi^{+}_h(z)$ and $\Phi^-_h(z)$ defined for $z$ in an open set $\amsmathbb{M}_h\subset \amsmathbb C$ which contains $\big[\hat{a}_h-\tfrac{1}{2\N},\hat{b}_h+\tfrac{1}{2\N}\big]$ but not $\big[\hat{a}_{g}-\tfrac{1}{2\N},\hat{b}_{g}+\tfrac{1}{2\N}\big]$ for $g \neq h$, and such that
 \begin{equation}
 \label{eq_ratio_def}
\forall x \in [\hat a_h, \hat b_h]\qquad \frac{w_h\big(\N x+\frac{1}{2}\big)}{w_h\big(\N
x-\frac{1}{2}\big)}= \frac{\Phi^+_{h}(x)}{\Phi^-_{h}(x)}.
 \end{equation}
Define
\begin{equation}
\label{Nekrasov_eqn} R_h(z)= \Phi^-_h(z)
 \cdot  \E \left[ \prod_{i=1}^N
 \bigg(1-\frac{\theta_{h,h(i)}}{\N z-\ell_i+\frac{1}{2}} \bigg)
  \right]
 +
 \Phi^+_h(z)
  \cdot \E \left[ \prod_{i=1}^N
 \bigg(1+\frac{\theta_{h,h(i)}}{\N z-\ell_i-\frac{1}{2}} \bigg)
  \right],
\end{equation}
Then $R_h(z)$ has no
singularities in $\amsmathbb{M}_h$ except, perhaps, for simple poles at
\,$\hat{a}_h-\tfrac{1}{2\N}$ or \,$\hat{b}_h+\tfrac{1}{2\N}$. Moreover
\begin{itemize}
\item if $\hat{a}_h$ is finite and $\Phi_h^-\big(\hat{a}_h-\tfrac{1}{2\N}\big)=0$, then $R_h(z)$ has no pole at
$\hat{a}_h-\frac{1}{2\N}$;
\item if $\hat{b}_{h}$ is finite and $\Phi_h^+\big(\hat{b}_h+\frac{1}{2\N}\big)=0$, then
$R_h(z)$ has no pole at $\hat{b}_h+\frac{1}{2\N}$.
\end{itemize}
\end{theorem}
\begin{remark}
\label{rem:vanishing_interpretation} If the weight $w_h$ is defined on a segment larger than $[a_h,b_h]$, the condition $\Phi_h^-\big(\hat{a}_h - \frac{1}{2\N}\big) = 0$ means $w_h(a_h - 1) = 0$ while the condition $\Phi_h^+(\hat{b}_h + \frac{1}{2\N}\big) = 0$ means $w_h(b_h + 1) = 0$.
\end{remark}
\begin{remark} \label{rem:allway} Like the Dyson--Schwinger equations for continuous ensembles, these equations are insensitive to filling fractions. In particular, they hold both for the fixed filling fraction or for fluctuating filling fractions. Note however that in discrete ensembles, due to the integrality condition \eqref{eq_segment_ff_relation} the vanishing conditions $\Phi_h^{-}(\hat{a}_h - \frac{1}{2\N}) = 0$ or $\Phi_h^{+}(\hat{b}_h + \frac{1}{2\N}) = 0$ cannot be achieved for generic $\theta_{h,h}$ without fixing filling fractions.\end{remark}
\begin{proof} We first assume that $a_1$ and $b_H$ are finite. In this situation $R_h(z)$ is a finite sum --- over all possible $\boldsymbol{\ell} = (\ell_i)_{i = 1}^{N}$ in the state space $\W_\N$ --- of meromorphic functions of $z\in \amsmathbb{M}_h$.
The possible singularities of $R_h(z)$ are simple poles at points $\frac{1}{\N}\big(\ell_i+ \frac{1}{2}\big)$ and
$\frac{1}{\N}\big(\ell_i-\frac{1}{2}\big)$ for some $\boldsymbol{\ell} \in\W_\N$ and $i \in [N]$. Let us compute the residue of $R_h(z)$ at such a point $z = \frac{m}{\N}$ in the segment
$[\hat{a}_h,\hat{b}_h]$. It arises from configurations $\boldsymbol{\ell} = (\ell_i)_{i = 1}^{N}$ such
that $m = \ell_i\pm \frac{1}{2}$ for some $i \in [N]$, therefore
\begin{equation}
\label{eq_x8}
\begin{split}
\Res_{z = \frac{m}{\N}} R_h(z)\dd z & = - \frac{\theta_{h,h}}{\N} \sum_{\boldsymbol{\ell} \in\W_\N} \sum_{i = 1}^{N} \delta_{\ell_i,m + \frac{1}{2}} \cdot \Phi^-_{h}\Big(\frac{m}{\N}\Big)\cdot \P(\boldsymbol{\ell}) \cdot \prod_{j\neq i}
 \bigg(1-\frac{\theta_{h,h(j)}}{m-\ell_j+\frac{1}{2}}\bigg) \\
 & \quad + \frac{\theta_{h,h}}{\N} \sum_{\boldsymbol{\ell}\in\W_\N} \sum_{i = 1}^{N} \delta_{\ell_i,m-\frac{1}{2}} \cdot \Phi^+_h\Big(\frac{m}{\N}\Big) \cdot
 \P(\boldsymbol{\ell}) \cdot \prod_{j\neq i}
 \bigg(1+\frac{\theta_{h,h(j)}}{m-\ell_j-\frac{1}{2}}\bigg).
\end{split}
\end{equation}
Note the difference between the two summation sets in \eqref{eq_x8}. Fix $i \in [N]$. If $\ell_i$ is
not the leftmost particle in $[a_h,b_h]$, then in the first sum appear configurations where
$\ell_{i-1}=m+\frac{1}{2} -\theta_{h,h}$. These are absent from the second sum because $(\ell_{i} - \ell_{i - 1} - \theta_{h,h}) \in \amsmathbb{Z}_{\geq 0}$.
However, each such term in the first sum is actually zero due to the numerator of the $j = (i - 1)$-th factor in the product. If $\ell_i$ is the
leftmost particle and $\ell_i = a_h = m + \frac{1}{2}$, no configuration with the same value of $\ell_i$ contribute to the second sum. But if
$\Phi^-_h\big(\hat{a}_{h}-\frac{1}{2\N}\big)=0$, the corresponding term in the first sum vanishes as well. Similar
considerations apply to the cases where $\ell_i$ is not the rightmost particle in $[a_h,b_h]$, or is the rightmost particle and is located at $b_h$. We conclude that it suffices to study the case when the terms of the two sums in \eqref{eq_x8} are in one-to-one
correspondence.

In order to analyze \eqref{eq_x8}, let us identify how the probability measure \eqref{eq_general_measure} changes when the position of one particle
$\ell_i=x$ is changed to $\ell_i=x-1$. Setting $\ell_j=y$,
 the $(i,j)$-th factor in the double product over pairs in \eqref{eq_general_measure} is multiplied by
\begin{equation}
\label{eq_x6}
\begin{split}
& \quad \frac{\Gamma\big(y-x+1\big)\cdot\Gamma\big(y-x+\theta_{h,h(j)}\big)}{\Gamma\big(y-x\big)\cdot\Gamma\big(y-x+1-\theta_{h,h(j)}\big)}
\cdot
\frac{\Gamma\big(y-x+1\big)\cdot\Gamma\big(y-x+2-\theta_{h,h(j)}\big)}{\Gamma\big(y-x+2\big)\cdot\Gamma\big(y-x+1+\theta_{h,h(j)}\big)} \\
& = \frac{(y-x)(y-x+1-\theta_{h,h(j)})}{(y-x+1)(y-x+\theta_{h,h(j)})},
\end{split}
\end{equation}
if $i<j$, and by
\begin{equation}
\label{eq_x7}
\begin{split}
& \quad \frac{\Gamma\big(x-y+1\big)\cdot\Gamma\big(x-y+\theta_{h,h(j)}\big)}{\Gamma\big(x-y\big)\cdot\Gamma\big(x-y+1-\theta_{h,h(j)}\big)}
\cdot
\frac{\Gamma\big(x-y-1\big)\cdot\Gamma\big(x-y-\theta_{h,h(j)}\big)}{\Gamma\big(x-y\big)\cdot\Gamma\big(x-y+\theta_{h,h(j)}-1\big)} \\
& = \frac{(x-y)(x-y+\theta_{h,h(j)}-1)}{(x-y-1)(y-r-\theta_{h,h(j)})},
\end{split}
\end{equation}
if $i> j$. Note that \eqref{eq_x6} and \eqref{eq_x7} are two forms of the same
rational expression. Using \eqref{eq_x6}, \eqref{eq_x7} and \eqref{eq_ratio_def}, we see that
\begin{equation*}
\begin{split}
& \quad \Phi^-_h(\tfrac{m}{\N}) \cdot
\P(\ell_1,\ldots\ell_{i-1},m+\tfrac{1}{2},\ell_{i+1},\ldots,\ell_N) \cdot \prod_{j\neq i}
 \bigg(1-\frac{\theta_{h,h(j)}}{m-\ell_j+\frac{1}{2}}\bigg)
\\
& =\Phi^+_h(\tfrac{m}{\N}) \cdot
 \P(\ell_1,\ldots\ell_{i-1},m-\tfrac{1}{2},\ell_{i+1},\ldots,\ell_N) \cdot \prod_{j\neq i}
 \bigg(1+\frac{\theta_{h,h(j)}}{m-\ell_j-\frac{1}{2}}\bigg).
\end{split}
\end{equation*}
Therefore, the contributions differing by a shift by $+1$ for the position of the $i$-th particle in the first and second sum in \eqref{eq_x8} cancel out
and the total residue is zero. Even if for given $\boldsymbol{\ell} \in \W_\N$, the configuration obtained by shifting the position of the $i$-th particle by $1$ fails to be in $\W_\N$, in such cases the formula for $\P(\ell)$ and the expressions in \eqref{eq_x8} still make sense but give zero.

\smallskip

If $a_1$ or $b_H$ is infinite, then the state space $\W_\N$ is infinite. Hence, we should additionally show that the sum over configurations involved in $R_h(z)$ and in \eqref{eq_x8} are convergent. We first remark that, due to Condition 4. in the definition of $\W_\N$ in Section~\ref{Section_configuration_space}, the set of possible locations of particles for configurations in $\W_\N$ is discrete. Take $z\in \amsmathbb M_h$ and assume that $z - \tfrac{1}{2\N}$ is at distance $\eta > 0$ from the aforementioned set of possible locations of particles. Then we have
\[
\Bigg|\prod_{i = 1}^N \bigg(1 \mp \frac{\theta_{h,h(i)}}{\N z - \ell_i \pm \frac{1}{2}}\bigg)\Bigg| \leq \bigg(1 + \frac{|\!|\boldsymbol{\Theta}|\!|_{\infty}}{\eta\N}\bigg)^{N}.
\]
This shows that the two random variables under the expectations in \eqref{Nekrasov_eqn} are uniformly bounded. Hence, both terms in \eqref{Nekrasov_eqn} are finite.
\end{proof}

\section{Higher-order Nekrasov equations}\label{higher}

Higher-order Nekrasov equations can be derived from the first-order Nekrasov equation by varying the weights. They will be used in Chapter~\ref{Chapter_fff_expansions} to gain information as $\N \rightarrow \infty$ about the observables introduced in Section~\ref{Section_observables}. Recall from Definition~\ref{def_empirical_mes} the notation for the Stieltjes transform of the empirical measure in the segments:
\[
\forall h \in [H]\qquad G_h(z) = \sum_{i = 1}^{N_h} \frac{1}{z - \frac{\ell_i^{h}}{\N}}.
\]
This is a random meromorphic function of $z \in \amsmathbb{C}$.

\begin{corollary}
\label{Corollary_higherNek} Consider a discrete ensemble as in Section~\ref{Section_general_model}. Fix $h \in [H]$ and suppose that there exist two holomorphic functions $\Phi^{+}_h(z)$ and $\Phi^{-}_h(z)$
 defined for $z$ in an open set $\amsmathbb{M}_h\subset \amsmathbb C$ which contains $\big[\hat{a}_h-\tfrac{1}{2\N},\hat{b}_h+\tfrac{1}{2\N}\big]$ but not $\big[\hat{a}_{g}-\tfrac{1}{2\N},\hat{b}_{g}+\tfrac{1}{2\N}\big]$ for $g \neq h$, and such that
 \begin{equation}
 \label{eq_ratio_def2}
\forall x \in [\hat a_h, \hat b_h]\qquad \frac{w_h\big(\N x+\frac{1}{2}\big)}{w_h\big(\N
x-\frac{1}{2}\big)}= \frac{\Phi^+_{h}(x)}{\Phi^-_{h}(x)}.
 \end{equation}
For any integer $n \geq 2$, segment indices $h,h_2,\ldots,h_n \in [H]$ and variables $z_j\in \amsmathbb{C}\setminus [\hat a_{h_j}, \hat b_{h_j}]$ indexed by $j \in \llbracket 2,n\rrbracket$, the following expression
\begin{equation}
\label{eq_Nekrasov_higher}
\begin{split}
& \quad R_{h;h_2,\ldots,h_n}(z;z_2,\ldots,z_n) \\
& = \sum_{J \subseteq \llbracket 2,n\rrbracket} \sum_{\tau \in \{\pm 1\}} \bigg(\prod_{j \in J} \frac{\delta_{h,h_j}}{z_j - z - \frac{\tau}{2\N}}\bigg) \cdot
 \Phi_h^{\tau}(z) \cdot \E^{(\textnormal{c})}\Bigg[\prod_{i = 1}^N\bigg(1 + \frac{\tau}{\N}\,\frac{\theta_{h,h(i)}}{z - \frac{\ell_i}{\N} - \frac{\tau}{2\N}}\bigg)\,,\,\big(G_{h_j}(z_j)\big)_{j \notin J}\Bigg]
\end{split}
\end{equation}
is a meromorphic function of the variable $z\in \amsmathbb M_h$. The only possible poles are located
\begin{itemize}
\item at $z=z_j \pm \frac{1}{2\N}$ for some $j \in \llbracket 2,n\rrbracket$ --- these poles are simple if all points $z_j \pm \frac{1}{2\N}$ are pairwise distinct, but multiplicities add up if some pairs coincide;
\item at $z = \hat{a}_h - \frac{1}{2\N}$. This pole is at most simple, and it is absent if $\Phi_h^-\big(\hat{a}_h-\tfrac{1}{2\N}\big) = 0$;
\item at $z = \hat{b}_h + \frac{1}{2\N}$. This pole is at most simple, and it is absent if $\Phi_h^+\big(\hat{b}_h+\tfrac{1}{2\N}\big) = 0$.
\end{itemize}
In addition, \eqref{eq_Nekrasov_higher} is a meromorphic function of the variables $z_j\in \amsmathbb{C}\setminus [\hat a_{h_j}, \hat b_{h_j}]$ for $j \in \llbracket 2,n\rrbracket$, with only possible poles $z_j=z \pm \frac{1}{2\N}$.
\end{corollary}
Remark~\ref{rem:allway} still applies: for $\theta_{h,h}$ generic the last two statements are only interesting in the fixed filling fractions setting. It is possible to give $R_{h;h_2,\ldots,h_n}(z;z_2,\ldots,z_n)$ an equivalent form making the structure of the signs clearer. It is obtained by using the identity
\[
\forall \tau \in \{\pm 1\}\qquad \frac{1}{z_j - z - \frac{\tau}{2\N}} = \frac{z_j - z +
\tfrac{\tau}{2\N}}{(z_j - z)^2 - \frac{1}{4\N^2}},
\]
and expanding the products in \eqref{eq_Nekrasov_higher}
\begin{equation}
\label{otherformR2}
\begin{split}
& \quad R_{h;h_2,\ldots,h_n}(z;z_2,\ldots,z_n)\\
 & = \sum_{J \subseteq \llbracket 2,n\rrbracket} \sum_{J' \subseteq J} \bigg(\prod_{j \in J} \frac{\delta_{h,h_j} }{(z - z_j)^2 - \frac{1}{4\N^2}}\bigg) \cdot \prod_{j \in J'} (z_{j} - z)  \\
& \qquad \qquad\qquad \times \sum_{\tau \in \{\pm 1\}} \bigg(\frac{\tau}{2\N}\bigg)^{\# J - \#J'} \cdot \Phi_h^{\tau}(z) \cdot \E^{(\textnormal{c})}\Bigg[\prod_{i = 1}^N\bigg(1 + \frac{\tau}{\N}\,\frac{\theta_{h,h(i)}}{z - \frac{\ell_i}{\N} - \frac{\tau}{2\N}}\bigg)\,,\,\big(G_{h_j}(z_j)\big)_{j \notin J}\Bigg].
\end{split}
\end{equation}

\begin{proof}[Proof of Corollary~\ref{Corollary_higherNek}]
The holomorphicity in $z_2,\ldots,z_n$ away from the announced poles is clear from the definition. In order to prove meromorphicity in $z$ we note that Nekrasov equations are valid for all weights satisfying the assumptions of
Theorem~\ref{Theorem_Nekrasov},
and it is clear from their derivation that they hold even if $w_h$
is complex-valued. Let $\boldsymbol{t} = (t_i)_{i = 2}^n$ be an auxiliary tuple of
variables. If the assumptions of Theorem~\ref{Theorem_Nekrasov} are satisfied for given weights $w_1,\ldots,w_H$, the ensemble with perturbed weights
\begin{equation}
\label{perturbed_weight} \forall g \in [H] \qquad w_{g}^{(\boldsymbol{t},\boldsymbol{z})}(\boldsymbol{\ell}) = w_{g}(\boldsymbol{\ell}) \cdot \prod_{j =
2}^{n} \bigg(1 + \frac{\delta_{g,h_j}\,t_{j}}{z_j - \tfrac{\ell}{\N}}\bigg),
\end{equation}
also satisfies Nekrasov equations, as we can take for each $g \in [H]$ and $z \in \amsmathbb{M}_{g}$
\[
\Phi_{g}^{\tau,(\boldsymbol{t},\boldsymbol{z})}(z) = \Phi_{g}^{\tau}(z) \cdot \prod_{j = 2}^{n}
\bigg(1 + \frac{\delta_{g,h_j} \,t_j}{z_j - z - \frac{\tau}{2\N}}\bigg),
\]
and note that $\Phi_{g}^{\tau,(\boldsymbol{t},\boldsymbol{z})}(z)|_{\boldsymbol{t} = 0} =
\Phi_{g}^{\tau}(z)$. We denote $\amsmathbb{P}_{\N}^{(\boldsymbol{t},\boldsymbol{z})}$ the complex
measure with total mass $1$ associated with this ensemble\footnote{In order to be able to normalize the measure to mass $1$, we should assume that the partition function does not vanish. Vanishing could occur for some values of $t_i$, since the weights are now complex. However, due to continuity in $t_i$, the partition function does not vanish for small enough values of $t_i$, and this is all we need.} and $\E^{(\boldsymbol{t},\boldsymbol{z})}$ the corresponding expectation value. Theorem~\ref{Theorem_Nekrasov} for this measure states that
\begin{equation}
\label{Nekreload}
 R_{h}^{(\boldsymbol{t},\boldsymbol{z})}(z) := \Phi_h^{\tau,(\boldsymbol{t},\boldsymbol{z})}(z) \cdot \E^{(\boldsymbol{t},\boldsymbol{z})}\Bigg[\prod_{i
= 1}^N \bigg(1 + \frac{\tau}{\N}\,\frac{\theta_{h,h(i)}}{z - \frac{\ell_i}{\N} -
\frac{\tau}{2\N}}\bigg)\Bigg]
\end{equation}
is a meromorphic function of $z \in \amsmathbb{M}_{h}$, with possible (simple) poles $z = \hat{a}_h - \frac{1}{2\N}$ and $z = \hat{b}_h + \frac{1}{2\N}$. We introduce
\begin{equation}
\label{eq_x89}
R_{h;h_2,\ldots,h_n}(z;z_2,\ldots,z_n):= \partial_{t_2} \cdots \partial_{t_n} R_{h}^{(\boldsymbol{t},\boldsymbol{z})}(z) \big|_{\boldsymbol{t} = 0}.
\end{equation}
We now proceed in two steps: we first show that this expression matches \eqref{eq_Nekrasov_higher}, and we then describe its analytic properties in the variable $z$.

Let us examine in detail \eqref{eq_x89} which is obtained by differentiating \eqref{Nekreload}. We have a sum over all subsets $J \subseteq \llbracket 2,n\rrbracket$ which index the derivatives that will hit $\Phi_{h}^{\tau,(\boldsymbol{t},\boldsymbol{z})}$. For any $\tau \in \{\pm 1\}$ we compute
\begin{equation}
\label{phider}
\bigg( \prod_{j \in J} \partial_{t_j}\bigg) \Phi_{h}^{\tau,(\boldsymbol{t},\boldsymbol{z})}(z)\Big|_{\boldsymbol{t} = 0} = \bigg(\prod_{j \in J} \frac{\delta_{h,h_j}}{z_j - z - \frac{\tau}{2\N}}\bigg) \cdot \Phi_{h}^{\tau}(z).
\end{equation}
The derivatives indexed by $j \notin J$ will then hit
$\E^{(\boldsymbol{t},\boldsymbol{z})}$, whose
$(\boldsymbol{t},\boldsymbol{z})$-dependence comes from the factors
\[
\prod_{g = 1}^H \prod_{i = 1}^{N_{g}} w_{g}^{(\boldsymbol{t},\boldsymbol{z})}(\ell_i^g).
\]
in the probability measure. We claim that
\begin{equation}
\label{notinJ}
\bigg(\prod_{j \notin J} \partial_{t_{j}}\bigg)
\E^{(\boldsymbol{t},\boldsymbol{z})}\bigg[\prod_{i = 1}^N \bigg(1 +
\frac{\tau}{\N}\,\frac{\theta_{h,h(i)}}{z - \frac{\ell_i}{\N} -
\frac{\tau}{2\N}}\bigg)
\bigg] \bigg|_{\boldsymbol{t} = 0} = \E^{(\textnormal{c})}\Bigg[\prod_{i = 1}^N \bigg(1 +
\frac{\tau}{\N}\,\frac{\theta_{h,h(i)}}{z - \frac{\ell_i}{\N} -
\frac{\tau}{2\N}}\bigg) ,\big(G_{h_j}(z_j))_{j \notin J}\Bigg].
\end{equation}
To justify this formula, we use the definition \eqref{eq_cumulant_def} with the random variable
\[
X_1=\prod\limits_{i = 1}^N \left(1 +
\frac{\tau}{\N}\,\frac{\theta_{h,h(i)}}{z - \frac{\ell_i}{\N} -
\frac{\tau}{2\N}}\right)
\]
and the remaining random variables being $G_{h_j}(z_j)$ for $j \notin J$. We differentiate explicitly in $s_1$ and get
\begin{equation}
\label{eq_x90}
\begin{split}
& \quad \E^{(\textnormal{c})}\Bigg[\prod_{i = 1}^N \bigg(1 +
\frac{\tau}{\N}\,\frac{\theta_{h,h(i)}}{z - \frac{\ell_i}{\N} -
\frac{\tau}{2\N}}\bigg) ,\big(G_{h_j}(z_j))_{j \notin J}\Bigg] \\
& = \bigg(\prod_{j \notin J} \partial_{t_j}\bigg) \left\{\frac{\E \left[ \prod\limits_{i = 1}^N \bigg(1 +
\frac{\tau}{\N}\,\frac{\theta_{h,h(i)}}{z - \frac{\ell_i}{\N} -
\frac{\tau}{2\N}}\bigg) \exp\left(\sum_{j\notin J} t_j G_{h_j}(z_j)\right) \right]} {\E \left[ \exp\left(\sum_{j\notin J} t_j G_{h_j}(z_j)\right) \right]}\right\}_{\boldsymbol{t} = 0},
\end{split}
\end{equation}
where $\E$ is expectation with respect to the original --- not $(\boldsymbol{t},\boldsymbol{z})$-deformed --- measure. We notice that
\[
1 +\delta_{h,h_j} \frac{t_{j}}{z_j - \frac{\ell}{\N}} \,\,\mathop{=}_{t_j \rightarrow 0}\,\, \exp\left( \frac{t_{j}\,\delta_{h,h_j}}{z_j - \tfrac{\ell}{\N}}\right)+ O(t_j^2),
\]
So, Equation \eqref{eq_x90} yields exactly \eqref{notinJ}. Putting together \eqref{phider} and \eqref{notinJ}, and summing over $\tau \in \{\pm 1\}$, we conclude that \eqref{eq_x89} is the same as \eqref{eq_Nekrasov_higher}.

\medskip

Next, we turn to the analytic properties of \eqref{eq_x89} in the variable $z$. We take a point $z_0\in \amsmathbb M_h$ and aim to show the holomorphicity in a small neighborhood of $z_0$ chosen distinct from $\hat{a}_h - \frac{1}{2\N}$ and $\hat{b}_h + \frac{1}{2\N}$. For this purpose we choose a simple contour $\gamma_h \subset \amsmathbb{M}_h$ oriented counterclockwise and such that
\begin{itemize}
 \item $z_0$ is inside $\gamma_h$;
 \item $\gamma_h$ does not intersect the discrete set $\big\{\N^{-1}(\ell^h_i \pm \frac{1}{2})\,\,|\,\,i \in [N_h]\,\,\textnormal{and}\,\,\boldsymbol{\ell} \in \W_\N\big\}$;
 \item $z_j \pm \frac{1}{2\N}$ remains outside of $\gamma_h$ for any $j \in \llbracket 2,n\rrbracket$;
 \item if $\Phi_h^-\big(\hat{a}_h-\tfrac{1}{2\N}\big)\neq 0$, then
$\hat{a}_h-\tfrac{1}{2\N}$ is outside $\gamma_h$;
\item if $\Phi_h^+\big(\hat{b}_h+\tfrac{1}{2\N}\big)\neq 0$, then $\hat{b}_h+\tfrac{1}{2\N}$ is outside $\gamma_h$.
\end{itemize}
By Cauchy integral formula, for any $z$ in a small neighborhood of $z_0$ we have
\[
 R_{h}^{(\boldsymbol{t},\boldsymbol{z})}(z) = \oint_{\gamma_h} \frac{\dd\zeta}{2\ii\pi}\, \frac{R_{h}^{(\boldsymbol{t},\boldsymbol{z})}(\zeta)}{\zeta-z}.
\]
Since the contour $\gamma_h$ avoids all singularities in the definition of $R_{h}^{(\boldsymbol{t},\boldsymbol{z})}(z)$, we can directly differentiate the last formula with respect to $t_2,t_3,\ldots,t_n$ under the integral sign and then set $t_2=\cdots=t_n=0$. The outcome is
\[
 R_{h;h_2,\ldots,h_n}(z;z_2,\ldots,z_n)=\ \oint_{\gamma_h} \frac{\dd \zeta}{2\ii\pi}\,\frac{R_{h;h_2,\ldots,h_n}(\zeta;z_2,\ldots,z_n)}{\zeta-z},
\]
and this shows that $R_{h;h_2,\ldots,h_n}(z;z_2,\ldots,z_n)$ is holomorphic a function of $z$ near $z_0$. Finally, the fact that there are at most simple poles at $z=\hat{a}_h-\tfrac{1}{2\N}$ or $z=\hat{b}_h+\tfrac{1}{2\N}$ is immediate from the definition \eqref{eq_Nekrasov_higher}.
\end{proof}

\part{ASYMPTOTIC ANALYSIS}
\label{Part_Asymptotic}

\chapter{Law of large numbers and large deviations principles}
\label{Chapterlarge}

In this chapter we perform the leading-order analysis of the discrete ensembles described in Chapter~\ref{Chapter_Setup_and_Examples}. We explain how a typical configuration looks like macroscopically, by establishing a form of the law of large numbers. We also produce large deviations bounds which show which configurations become very unlikely as $\N\rightarrow\infty$.

The following consequences of the Stirling asymptotic expansion of the Gamma function will be used repeatedly in this chapter. For any $\theta \in \amsmathbb{R}$ there exist constants $c_1,c_2 \in \amsmathbb{R}$ such that
\begin{equation}
\label{eq_Stirling_ratio}
\forall \xi \in (\max(\theta,1-\theta),+\infty)\qquad \exp\left(2\theta \log \xi + \frac{c_2}{\xi}\right) \leq
\frac{\Gamma(\xi + 1)\cdot \Gamma(\xi + \theta)}{\Gamma(\xi)\cdot \Gamma(\xi + 1 - \theta)} \leq \exp\left(2\theta \log \xi + \frac{c_1}{\xi}\right).
\end{equation}
Moreover, if we fix some $\theta_0>0$, then the constants $c_1,c_2$ can be chosen
uniformly for any $\theta\in [-\theta_0,\theta_0]$. We often use \eqref{eq_Stirling_ratio} with $\xi=|\N x-\N y|$, and it takes the form
\begin{equation}
\label{eq_Stirling_xy}
|x - y|^{2\theta} \cdot e^{\frac{c_2}{\N |x - y|}} \leq \frac{\Gamma\big(|\N x - \N y| + 1\big)\cdot\Gamma\big(|\N x - \N y| + \theta\big)}{\N^{2\theta}\cdot \Gamma\big(|\N x - \N y|\big)\cdot\Gamma\big(|\N x
- \N y| + 1 - \theta\big)} \leq |x - y|^{2\theta} \cdot e^{\frac{c_1}{\N |x -
y|}}
\end{equation}
for $|x-y|\geq \frac{1}{\N}\max(\theta,1-\theta)$.

This chapter is organized as follows. In Section~\ref{Section_Energy_functional} we
introduce the energy functional and its minimizer called equilibrium measure, and we establish the characterization of equilibrium measures announced in Theorem~\ref{Theorem_equi_charact}. In Section~\ref{section_rought_ann} we state the theorems of concentration of the empirical measure around the equilibrium measure as $\N\rightarrow\infty$. Their proofs appear in Section~\ref{Section_rough_proofs}. Finally, in Section~\ref{Section_ld_suport} we establish the large deviations principle for the particle counts in voids and saturations at large $\N$.

\section{Energy functional and equilibrium measure}

\label{Section_Energy_functional}
\subsection{Characterization of the equilibrium measure}

The energy functional $-\mathcal{I}$ was introduced in Section~\ref{Definition_functional} and will play a central role in this chapter:
\begin{equation}
\label{eq_functional_general_repeat}  \mathcal{I}[\boldsymbol{\nu}]=\sum_{g,h = 1}^{H} \int_{\hat a'_{g}}^{\hat b'_g} \int_{\hat a'_{h}}^{\hat b'_{h}} \theta_{g,h}\log|x-y|
 \, \nu_g(x)\nu_{h}(y)\dd x \dd y- \sum_{h=1}^{H} \int_{\hat a'_h}^{\hat b'_h} V_h(x) \nu_h(x)\dd x.
\end{equation}
In this section we study its basic properties, prove the existence and
uniqueness of its minimizer and Theorem~\ref{Theorem_equi_charact}. For the particular case where the
equations \eqref{eq_equations_eqs} in Section~\ref{DataS} deterministically fix the segment filling fractions, the results
are known, \textit{cf.} \cite{Hardy_Kuijlaars_min} and references therein.

\label{index:Pstar}Before starting, one should keep in mind the notations of Section~\ref{Offcrit_section}. In particular, we are interested in minimizing $-\mathcal{I}$ over the set $\mathscr{P}_\star$ of nonnegative measures $\nu$ on $\amsmathbb{A}=\bigcup_{h=1}^H [\hat a'_h,\hat b'_h]$ which are absolutely continuous, whose density on $[\hat a'_h,\hat b'_h]$ is bounded from above by $\frac{1}{\theta_{h,h}}$ for each $h \in [H]$, and whose masses satisfy the constraints \eqref{eq_equations_eqs}. Such measures $\nu$ are equivalently described as tuples of measures $\boldsymbol{\nu} = (\nu_h)_{h = 1}^{H}$, where $\nu_h$ stands for the restriction of $\nu$ to $[\hat{a}_h',\hat{b}_h']$. We equip $\mathscr{P}_\star$ with the weak topology: a sequence $(\boldsymbol{\nu}^{(n)})_{n \geq 1}$ converges to $\boldsymbol{\nu}$ if and only if for any bounded continuous function $f: \amsmathbb{A} \rightarrow \amsmathbb{R}$ we have
\[
 \lim_{n\rightarrow\infty}\sum_{h=1}^H \int_{\hat a'_h}^{\hat b'_h} f(x)\, \dd \nu^{(n)}_h(x) = \sum_{h=1}^H \int_{\hat a'_h}^{\hat b'_h} f(x)\, \dd \nu_h(x).
\]

\begin{lemma}
\label{Lemma_I_welldef} Under Assumptions~\ref{Assumptions_Theta} and
\ref{Assumptions_basic}, there exists $C>0$ depending only on the constants in the
assumptions and such that
\[
\forall \boldsymbol{\nu} \in \mathscr{P}_{\star}\qquad \mathcal{I}[\boldsymbol{\nu}] \leq C - \frac{1}{C} \int_{\amsmathbb{R}} \log\big(1 + |x|\big)\nu(x)\dd x.
\]
 In particular, if $\int_{\amsmathbb{R}} \log\big(1 + |x|\big)\nu(x)\dd x=+\infty$, then $\mathcal I [\boldsymbol{\nu}]=-\infty$.
\end{lemma}
\begin{proof} Take a measure $\boldsymbol{\nu}\in \mathscr{P}_\star$. Since $x \mapsto \log|x|$ is integrable near
$x=0$ and $\nu$ has a bounded density, the integral in the definition of $\mathcal
I[\boldsymbol{\nu}]$ is well-behaved near the diagonal $x=y$. In the rest of the proof we analyze
its tail behavior as $x$ and $y$ become large. We use the following inequality
\begin{equation}
\label{eq_ln_upperbound} \forall x,y \in \amsmathbb{R}\qquad \log|x - y| \leq \log\big(1 + |x|) + \log\big(1 + |y|\big).
\end{equation}

First, assume that $b_H$ is finite. We claim that there exists $C>0$ such that for any $
g,h\in [H]$ and any $(x,y) \in [\hat a'_g, \hat b'_g] \times [\hat a'_{h}, \hat b'_{h}]$ we have
\begin{equation}
\label{eq_ln_better_bound}
 \theta_{g,h} \log|x-y| \leq\theta_{g,h} \Big(\log\big(1 + |x|\big)+\log\big(1 + |y|\big)\Big) +C.
\end{equation}
Indeed, if $\theta_{g,h}>0$, then \eqref{eq_ln_better_bound} is an immediate corollary of
\eqref{eq_ln_upperbound}; if $\theta_{g,h}<0$, then $g\neq h$ and we can assume without loss of
generality that $h>g$, which implies that $y$ belongs to a bounded set, and $|x-y|$ is bounded
away from $0$. Therefore
\begin{equation}
\label{eq_ln_first}
 \log|x-y|>\log\big(1 + |x|\big)+\log\big(1 + |y|\big) +C
\end{equation}
for some constant $C$ independent of $x$ and $y$, which implies
\eqref{eq_ln_better_bound} --- perhaps with a different value of $C$. Integrating the bound \eqref{eq_ln_better_bound} we see that
\begin{equation}
\label{eq_I_upperbound}\mathcal{I}[\boldsymbol{\nu}] \leq \sum_{h=1}^H \int_{\hat a'_{h}}^{\hat
b'_h} \bigg(- V_{h}(x)+\sum_{g = 1}^H 2\,\theta_{h,g}\hat{n}_{g}\log\big(1 + |x|\big)
\bigg)\nu_{h}(x)\dd x + C.
\end{equation}
Since $b_H$ is finite, in the sum \eqref{eq_I_upperbound} all the terms for $h>1$ are integrals of a
bounded function over a compact set, and therefore are finite. This is also the case for $h = 1$ if $a_1$ is finite. If $a_1 = - \infty$, we use \eqref{eq_confinement_left} to get an upper bound of the $h = 1$ integrand. This yields
\[
\mathcal{I}[\boldsymbol{\nu}] \leq C - \frac{1}{C} \int_{\amsmathbb{R}}
 \log\big(1 + |x|\big)\, \nu(x)\dd x.
\]

\smallskip

Second, we treat the case $b_H = +\infty$. If $a_1$ is finite or $\theta_{1,H}\geq 0$, then
we can repeat the argument using \eqref{eq_confinement_right} instead of \eqref{eq_confinement_left}.

The last case to handle is $a_1 = -\infty$, $b_H=+\infty$ and $\theta_{1,H}<0$. As the diagonal elements of $\boldsymbol{\Theta}$ are positive, we must have $H>1$.
In this case the bound \eqref{eq_ln_better_bound} could fail for $x\in [\hat a'_1, \hat b'_1]$ and $y\in
[\hat a'_H, \hat b'_H]$. Instead we then use
\begin{equation}
\label{eq_ln_second}
\forall (x,y) \in [\hat a'_1, \hat b'_1] \times [\hat a'_H, \hat b'_H]\qquad
 \log|x-y|\geq \frac{1}{2}\log\big(1 + |x|\big)+\frac{1}{2}\log\big(1 + |y|\big) -C,
 \end{equation}
which implies --- for a different positive constant $C$ ---
\begin{equation}
\label{eq_ln_better_bound_infinite}
\forall (x,y) \in [\hat a'_1, \hat b'_1] \times [\hat a'_H, \hat b'_H]\qquad
 \theta_{1,H} \log|x-y| \leq\frac{\theta_{1,H}}{2} \Big(\log\big(1 + |x|)+\log\big(1 + |y|\big)\Big) +C,
  \end{equation}
and then the same argument works again by using \eqref{eq_confinement_left_refined} and
\eqref{eq_confinement_right_refined} instead of \eqref{eq_confinement_left} and
\eqref{eq_confinement_right}.
\end{proof}

\begin{definition}
\label{DEFI2} Let $\mathcal I^{(2)}$ denote the quadratic part of the energy functional, \textit{i.e.}
\begin{equation}
 \label{eq_functional_quadratic}
  \mathcal{I}^{(2)}[\boldsymbol{\nu}] := \sum_{g,h = 1}^{H} \int_{\hat a'_g}^{\hat b'_g} \int_{\hat a'_{h}}^{\hat b'_{h}} \theta_{g,h}\log|x-y|\,
 \nu_g(x)\nu_{h}(y)\,\dd x \dd y
 \end{equation}
 \end{definition}
We would like to extend the definition \eqref{eq_functional_quadratic} to integrable signed measures, \textit{i.e.} differences of two finite-mass nonnegative measures. Since $\log|x-y|$ is unbounded for both small and large $|x-y|$, this requires some care. If $\nu$ is a signed measure, we denote $\nu^{\pm}$ its positive and negative parts, and $|\nu| = \nu^+ - \nu^-$ its absolute value. We first define for $T>0$
\begin{equation}
 \label{eq_functional_quadratic_smoothed}
\begin{split}
 K_T(x,y) & := \int_{\frac{1}{T}}^{T} \frac{e^{-t}- e^{-(x-y)^2 t}}{2t}\,\dd t, \\
  \mathcal{I}^{(2)}_T[\boldsymbol{\nu}] & := \sum_{g,h = 1}^{H}  \int_{\hat a'_g}^{\hat b'_g} \int_{\hat a'_{h}}^{\hat b'_{h}} \theta_{g,h} K_T(x,y)\,
\dd \nu_g(x) \dd \nu_{h}(y).
 \end{split}
 \end{equation}
Note that for each fixed $T>0$ the kernel $K_T(x,y)$ is uniformly bounded. Therefore, \eqref{eq_functional_quadratic_smoothed} is always well-defined. We use this observation to give an alternative definition for $\mathcal{I}^{(2)}[\boldsymbol{\nu}]$, following the approach of \cite{arous1997large}:
\begin{equation}
 \label{eq_functional_quadratic_alternative}
 \mathcal{I}^{(2)}[\boldsymbol{\nu}]=\lim_{T\rightarrow\infty} \mathcal{I}^{(2)}_T[\boldsymbol{\nu}]
\end{equation}
whenever the limit exists.
\begin{lemma}\label{Lemma_I_quadratic_smoothing}
 Suppose that for any $h\in[H]$, $\nu_h$ is an integrable signed measure such that
\begin{equation}
 \label{eq_quadratic_absolute}
  \sum_{g,h = 1}^{H}  \int_{\hat a'_g}^{\hat b'_g} \int_{\hat a'_{h}}^{\hat b'_{h}} \Big|\theta_{g,h}\log|x-y|
  \Big|\, \dd|\nu_g|(x)\,\dd |\nu_{h}|(y) < +\infty,
 \end{equation}
 Then both \eqref{eq_functional_quadratic} and \eqref{eq_functional_quadratic_alternative} exist and coincide.
\end{lemma}
If the condition \eqref{eq_quadratic_absolute} fails, we can still use \eqref{eq_functional_quadratic_alternative} as a definition of $\I^{(2)}$, if the limit exists in $\amsmathbb{R} \cup \{\pm \infty\}$. It can however be problematic to use \eqref{eq_functional_quadratic} in such cases, since one needs to properly define the value of the absolutely divergent integral.
\begin{proof}[Proof of Lemma~\ref{Lemma_I_quadratic_smoothing}] Let us first show that
\begin{equation}
\label{eq_KT_to_log}
\lim_{T\rightarrow\infty} K_T(x,y)=\log|x-y|.
\end{equation}
Denoting $|x-y|=\xi$, we need to show
\begin{equation}
\label{eq_log_integral}
\forall \xi \in \amsmathbb{R}_{> 0}\qquad \int_{0}^{+\infty} \frac{e^{-t} - e^{-\xi^2t}}{2t}\,\dd t =\log \xi.
\end{equation}
This is justified by observing that \eqref{eq_log_integral} is true at $\xi=1$, and the $\xi$-derivative of both sides is $\frac{1}{\xi}$. The integrand in the definition of $K_T(x,y)$ is positive whenever $|x-y|>1$ and nonpositive whenever $|x-y|<1$. Therefore, \eqref{eq_KT_to_log} implies $|K_T(x,y)|\leq \big|\log|x-y|\big|$. Then, the existence of the limit in \eqref{eq_functional_quadratic_alternative} and its coincidence with \eqref{eq_functional_quadratic} follows from the dominated convergence theorem.
\end{proof}

\begin{lemma} \label{Lemma_I_quadratic_Fourier}
Suppose that Assumption~\ref{Assumptions_Theta} holds and let
$\nu$ be an integrable signed measure on
$\amsmathbb{A}$ such that
\begin{equation}
\label{eq_nu_mass}
\forall e \in [\mathfrak{e}]\qquad \mathfrak r_e\big((\nu([\hat a'_h,\hat b'_h]))_{h = 1}^H\big) = 0.
\end{equation}
Then the limit in \eqref{eq_functional_quadratic_alternative} is either nonpositive or $-\infty$, and it admits the Fourier representation \eqref{eq_Fourier}
\begin{equation}
\label{eq_I_quadratic_Fourier}
 \mathcal I^{(2)}[\boldsymbol{\nu}]=-\int_{0}^{+\infty} \frac{\dd s}{s}\bigg( \sum_{g,h = 1}^H \theta_{g,h}\,\widehat \nu_g(s) \widehat \nu_{h}^*(s) \bigg).
\end{equation}
 In particular, both sides of \eqref{eq_I_quadratic_Fourier} are nonpositive and they achieve the value $-\infty$ simultaneously.
\end{lemma}
\begin{proof}
 A direct computation of the Gaussian integral shows that
\begin{equation}
\label{eq_x189}
 \int_{\hat a'_g}^{\hat b'_g} \int_{\hat a'_{h}}^{\hat b'_{h}} e^{-(x-y)^2 t}\dd \nu_g(x)\dd \nu_{h}(y) = \frac{1}{\sqrt{4\pi t}} \int_{\amsmathbb{R}} \widehat \nu_g(s) \widehat \nu_{h}^*(s) e^{-\frac{s^2}{4t}}\dd s.
\end{equation}
 For $e \in [\mathfrak{e}]$, we set $\nu^{(e)}=\mathfrak r_e(\boldsymbol{\nu})$. This is an integrable signed measure on $\amsmathbb{A}$ of mass $0$. Hence, \eqref{eq_x189} implies for any $T > 0$ and $e_1,e_2 \in [\mathfrak{e}]$
\begin{equation}
 \label{eq_x190}
 \begin{split}
\int_{\amsmathbb R^2} K_T(x,y) \dd \nu^{(e_1)}(x)\dd \nu^{(e_2)}(y) & =-\int_{\frac{1}{T}}^{T} \frac{\dd t}{4\sqrt{\pi t^3}} \int_{\amsmathbb{R}} \widehat \nu^{(e_1)}(s) \widehat \nu^{(e_2)*}(s) e^{-\frac{s^2}{4t}} \dd s \\
& = -\int_{\amsmathbb{R}} \dd s\, \widehat \nu^{(e_1)}(s) \widehat \nu^{(e_2)*}(s)
 \int_{\frac{1}{T}}^{T} \frac{e^{-\frac{s^2}{4t}}\dd t}{4\sqrt{\pi t^3}} \\
 &= -\int_{\amsmathbb{R}} \frac{\dd s}{2|s|}\, \widehat \nu^{(e_1)}(s) \widehat \nu^{(e_2)*}(s)
 \int_{\frac{s^2}{4T}}^{ \frac{T s^2}{ 4}} \frac{e^{-t}\dd t}{\sqrt{\pi t}}.
\end{split}
\end{equation}
 Multiplying \eqref{eq_x190} by $(e_1,e_2)$-th matrix element of
 $\boldsymbol{\Theta}'$ introduced in Condition 4. of Assumption~\ref{Assumptions_Theta} and
 summing over all $e_1,e_2\in [\mathfrak e]$, we get
 \begin{equation}
 \label{eq_x191}
  \mathcal I^{(2)}_T[\boldsymbol{\nu}]= -\int_0^\infty \frac{\dd s}{s} \bigg( \sum_{g,h=1}^H
 \theta_{g,h}\, \widehat \nu_g(s) \widehat \nu_{h}^*(s) \bigg) \int_{\frac{s^2}{4T}}^{ \frac{Ts^2}{ 4}} \frac{\exp(-t)\dd t}{\sqrt{\pi t}}.
 \end{equation}
 Since $\boldsymbol{\Theta}$ is a positive semi-definite matrix, $\sum_{g,h = 1}^H
 \theta_{g,h}\, \widehat \nu_g(s)\widehat \nu_{h}^*(s)\geq 0$ and we conclude that $\mathcal I^{(2)}_T[\boldsymbol{\nu}]\leq 0$ and it is decreasing as $T\rightarrow\infty$. Sending $T \rightarrow \infty$ in \eqref{eq_x191} and using $\int_0^{+\infty} \frac{e^{-t}}{\sqrt{t}}\dd t=\sqrt{\pi}$ we get \eqref{eq_I_quadratic_Fourier}.
\end{proof}

In the next statements we need a square root for the matrix of interactions.
\begin{definition}
\label{defsqur}Since $\boldsymbol{\Theta}$ is a positive semi-definite symmetric matrix, we can find another positive semi-definite
 symmetric matrix $\boldsymbol{\Theta}^{1/2}$ such that $(\boldsymbol{\Theta}^{1/2})^2=\boldsymbol{\Theta}$.
\end{definition}

\begin{corollary} \label{Corollary_I_positive}
 Under assumptions and with the notations of Lemma~\ref{Lemma_I_quadratic_Fourier}, we always have
 $\mathcal I^{(2)}[\boldsymbol{\nu}]\leq 0$, with equality if and only if $\boldsymbol{\nu}$ is the zero $H$-tuple of measures.
\end{corollary}
\begin{proof}
Since the matrix $\boldsymbol{\Theta}$ is positive semi-definite, $\sum_{g,h=1}^H
 \theta_{g,h}\, \widehat \nu_g(s)\widehat \nu_{h}^*(s)\geq 0$ for any $s$
 and
 \eqref{eq_I_quadratic_Fourier} implies $\mathcal{I}^{2}[\boldsymbol{\nu}]\leq 0$.
 By the same argument and continuity of $\widehat \nu_h(s)$, the equality $\mathcal I^{(2)}[\boldsymbol{\nu}]=0$
 implies that
 \begin{equation}
 \label{eq_x9}
 \forall s \in \amsmathbb{R}_{\geq 0}\qquad \sum_{g,h=1}^H \theta_{g,h}\, \widehat \nu_g(s)\widehat \nu_{h}^*(s)=0.
 \end{equation}
Using the square root $\boldsymbol{\Theta}^{1/2}$, \eqref{eq_x9} becomes
 \[
 \forall s \in \amsmathbb{R}_{\geq 0}\qquad \big|\boldsymbol{\Theta}^{1/2}(\widehat{\boldsymbol{\nu}}(s))\big|^2_2=0,
 \]
 which means that for any $h \in [H]$ we have $\sum_{g=1}^H [\boldsymbol{\Theta}^{1/2}]_{h,g} \widehat \nu_{g}(s) =0$
 for any $s$.
 Since the Fourier transform is isometric, this implies $\boldsymbol{\Theta}^{1/2}(\boldsymbol{\nu}) =0$, which after multiplication by $\boldsymbol{\Theta}^{1/2}$ yields $\boldsymbol{\Theta}(\boldsymbol{\nu}) = 0$, \textit{i.e.}
 \begin{equation}
 \label{eq_x10}
\forall h \in [H]\qquad \sum_{g=1}^H \theta_{h,g}\,\nu_{g} =0.
\end{equation}
As the signed measures $\nu_1,\ldots,\nu_H$ have pairwise disjoint supports and $\theta_{h,h}>0$ for any $h \in [H]$, we deduce from \eqref{eq_x10} that $\nu_{h}=0$ for any $h \in [H]$.
\end{proof}

The next statement turns Corollary~\ref{Corollary_I_positive} into a quantitative bound on $\I^{(2)}$

\begin{corollary} \label{Corollary_I_bound}
 Under the assumptions of Lemma~\ref{Lemma_I_quadratic_Fourier}, for any $\varkappa > 0$ and $M >1$ we have
\begin{equation}
\label{eq_I_lower_bound}
\max_{h\in[H]} \frac{\pi}{2} \cdot [\boldsymbol{\Theta}^{1/2}]_{h,h} \cdot \big|\nu([\hat{a}'_{h},\hat{b}'_{h}])\big| \leq \Bigg(-\bigg(\log M + \max_{h\in[H]} \frac{(\hat b'_h-\hat a'_h)^2}{8}\bigg) \cdot \mathcal{I}^{(2)}[\boldsymbol{\nu}]\Bigg)^{\frac{1}{2}} + \bigg(\frac{\mathcal{J}_{\varkappa}[\boldsymbol{\nu}]}{\varkappa M^{\varkappa}}\bigg)^{\frac{1}{2}},
\end{equation}
where
\begin{equation}
\label{eq_J_alpha}
\mathcal J_{\varkappa}[\boldsymbol{\nu}] :=\int_{0}^{+\infty} \frac{\dd s}{s^{1-\varkappa}} \bigg( \sum_{g,h =1}^H
 \theta_{g,h} \widehat \nu_g(s) \widehat \nu_{h}^*(s) \bigg).
\end{equation}
\end{corollary}
We aim to apply this corollary for large values of $M$ in situations where we can bound $\mathcal J_{\varkappa}[\boldsymbol{\nu}]$ from above and $\nu([\hat{a}'_{h},\hat{b}'_{h}])$ from below. Eventually, this would lead to $\I^{(2)}[\boldsymbol{\nu}]$ being bounded away from $0$.

\begin{proof}[Proof of Corollary~\ref{Corollary_I_bound}]
 For fixed $h\in[H]$, using $\nu_{g}(x)=0$ for any $x\in[\hat a'_h,\hat b'_h]$ and $g \neq h$, and the Plancherel theorem for the scalar product of the indicator function $\mathbbm{1}_{[\hat a'_h,\hat b'_h]}$ with $\sum_{g=1}^H [\boldsymbol{\Theta}^{1/2}]_{h,g} \nu_{g}(x)$, we have
\begin{equation}
 \label{eq_x186}
 \begin{split}
 [\boldsymbol{\Theta}^{1/2}]_{h,h}\int_{\hat a'_h}^{\hat b'_h} \dd \nu_h(x) & = \int_{\hat a'_h}^{\hat b'_h} \bigg(\sum_{g=1}^H \,[\boldsymbol{\Theta}^{1/2}]_{h,g}\, \dd\nu_{g}(x)\bigg) \\
 & = \frac{1}{2\pi} \int_{\amsmathbb{R}} \frac{e^{\ii s \hat b'_h} - e^{\ii s \hat a'_h}}{\ii s} \cdot \bigg(\sum_{g=1}^H\, [\boldsymbol{\Theta}^{1/2}]_{h,g}\,\widehat \nu_{g}(s)\bigg)\dd s.
\end{split}
\end{equation}
We split the last integral in two parts: over $|s|<M$ and over $|s|>M$. For the first one we use the Cauchy--Schwarz inequality as follows:
\begin{equation}
\label{eq_x185}
\begin{split}
 & \quad \bigg| \frac{1}{2\pi} \int_{-M}^{M} \frac{e^{\ii s \hat b'_h}-e^{\ii s \hat a'_h}}{\ii s} \cdot \bigg( \sum_{g=1}^H \,[\boldsymbol{\Theta}^{1/2}]_{h,g}\,\widehat \nu_{g}(s)\bigg)\dd s\bigg| \\
 & \leq \frac{1}{\pi} \int_{0}^{M}  \bigg| \frac{e^{\ii s \hat b'_h}-e^{\ii s \hat a'_h}}{\sqrt{s}}\bigg| \cdot\Bigg| \frac{\sum_{g=1}^H \,[\boldsymbol{\Theta}^{1/2}]_{h,g}\,\widehat \nu_{g}(s)}{\sqrt{s}}\Bigg|\,\dd s \\
 & \leq \frac{1}{\pi} \bigg(\int_{0}^{M}  \frac{\big|e^{\ii s \hat b'_h}-e^{\ii s \hat a'_h}\big|^2}{s}\dd s\bigg)^{\frac{1}{2}}\cdot \bigg(\int_{0}^{M}  \frac{\big|\sum_{g=1}^H\, [\boldsymbol{\Theta}^{1/2}]_{h,g}\,\widehat \nu_{g}(s)|^2}{s}\dd s\bigg)^{\frac{1}{2}}.
\end{split}
\end{equation}
 For the integral in the first factor we bound the integrand by $\frac{4}{s}$ for $s>1$, and by
 \[
\frac{\big|\exp\big(\frac{\ii s}{2}(\hat b'_h-\hat a'_h)\big)-\exp\big(\frac{\ii s}{2}(\hat a'_h-\hat b'_h)\big)\big|^2}{s}= \frac{4}{s} \,\sin^2\bigg(\frac{s(\hat b'_h-\hat a'_h)}{2}\bigg) \leq s (\hat b'_h-\hat a'_h)^2
 \]
 for $s<1$, where we used the inequality $|\sin(x)|\leq |x|$. Altogether, we bound the first integral in the right-hand side of \eqref{eq_x185} from above by $4 \log M + \frac{1}{2}(\hat b'_h-\hat a'_h)^2$. For the second integral, we use $(\boldsymbol{\Theta}^{1/2})^2 = \boldsymbol{\Theta}$ to get
 \[
 \bigg|\sum_{g=1}^H \,[\boldsymbol{\Theta}^{1/2}]_{h,g}\,\widehat \nu_{g}(s)\bigg|^2 \leq \sum_{h=1}^{H} \bigg|\sum_{g=1}^H\, [\boldsymbol{\Theta}^{1/2}]_{h,g}\,\widehat \nu_{g}(s)\bigg|^2=\big| \boldsymbol{\Theta}^{1/2}(\widehat{\boldsymbol{\nu}}(s))\big|^2_2=\sum_{g,h=1}^H \theta_{g,h}\,\widehat\nu_g(s)\widehat\nu^*_{h}(s).
 \]
 Hence, the second integral is bounded from above by $-\I^{(2)}[\boldsymbol{\nu}]$. Combining the results for the integrals in two factors, we bound \eqref{eq_x185} from above by
\[
 \frac{2}{\pi}\Bigg(-\bigg(\log M + \frac{(\hat b'_h-\hat a'_h)^2}{8}\bigg)\cdot \I^{(2)}[\boldsymbol{\nu}]\Bigg)^{\frac{1}{2}}.
\]
 We switch to the part of \eqref{eq_x186} corresponding to $|s|>M$ and use the Cauchy--Schwarz inequality in a different way:
\begin{equation*}
\begin{split}
 &\quad \bigg|\frac{1}{2\pi} \int_{|s|>M}\frac{e^{\ii s \hat b'_h}-e^{\ii s \hat a'_h}}{\ii s} \cdot \bigg(\sum_{g=1}^H\, [\boldsymbol{\Theta}^{1/2}]_{h,g}\,\widehat \nu_{g}(s)\bigg)\dd s\bigg| \\
\leq & \quad \frac{1}{\pi} \int_{M}^{+\infty}  \bigg| \frac{e^{\ii s \hat b'_h}-e^{\ii s \hat a'_h}}{\sqrt{s^{1+\varkappa}}}\bigg| \cdot\bigg| \frac{\sum_{g=1}^H\, [\boldsymbol{\Theta}^{1/2}]_{h,g}\,\widehat \nu_{g}(s)}{\sqrt{s^{1-\varkappa}}}\bigg|\,\dd s \\
& \leq \frac{1}{\pi} \bigg(\int_{M}^{+\infty}  \frac{\big|e^{\ii s \hat b'_h}-e^{\ii s \hat a'_h}\big|^2}{s^{1+\varkappa}}\dd s\bigg)^{\frac{1}{2}} \cdot \bigg(\int_{M}^{+\infty}  \frac{\big|\sum_{g=1}^H\, [\boldsymbol{\Theta}^{1/2}]_{h,g}\,\widehat \nu_{g}(s)\big|^2}{s^{1-\varkappa}}\dd s\bigg)^{\frac{1}{2}}.
\end{split}
\end{equation*}
 We bound the first integral from above by $4\int_M^{+\infty} s^{-1-\varkappa}\dd s= \frac{4}{\varkappa M^{\varkappa}}$. For the second integral we argue like for the second factor in \eqref{eq_x185} and get $\mathcal J_{\varkappa}[\boldsymbol{\nu}]$.
 \end{proof}

The next lemma is an analogue of Lemma~\ref{Lemma_I_quadratic_Fourier} for $\mathcal J_{\varkappa}[\boldsymbol{\nu}]$.
\begin{lemma} \label{Lemma_J_alpha} Suppose that Assumption~\ref{Assumptions_Theta} holds and let
$\nu$ be an integrable signed measure on
$\amsmathbb{A}$. Whenever the following integrals are absolutely convergent, we have an alternative expression for $\mathcal J_{\varkappa}[\boldsymbol{\nu}]$ of \eqref{eq_J_alpha}:
 \[
 \mathcal J_{\varkappa}[\boldsymbol{\nu}]=\Gamma(\varkappa) \cos\bigg(\frac{\pi \varkappa}{2}\bigg) \sum_{g,h=1}^H \theta_{g,h} \int_{\hat a'_g}^{\hat b'_g} \int_{\hat a'_{h}}^{\hat b'_{h}} \frac{\dd\nu_g(x)\dd\nu_{h}(y)}{|x - y|^{\varkappa}}.
 \]
\end{lemma}
\begin{proof}
 Since the integrand in the definition of $\mathcal J_{\varkappa}[\boldsymbol{\nu}]$ is positive, we can write
 \begin{equation}
 \label{eq_x187}
\mathcal  J_{\varkappa}[\boldsymbol{\nu}]=\lim_{\eps \rightarrow 0^+} \int_{0}^{+\infty} \frac{e^{-s\eps}\,\dd s}{s^{1-\varkappa}} \left( \sum_{g,h=1}^H
 \theta_{g,h} \int_{\hat a'_g}^{\hat b'_g} e^{\ii x s}\dd\nu_g(x) \int_{\hat a'_{h}}^{\hat b'_{h}} e^{-\ii y s}\dd\nu_{h}(y)\right).
 \end{equation}
 We interchange the order of integration and integrate in $s$ first using the identity
\begin{equation}
 \label{eq_x188}
 \int_{0}^{+\infty} \frac{e^{-s\eps+\ii x s-\ii y s}\dd s}{s^{1-\varkappa}} = \frac{1}{\big(\eps+\ii (x-y)\big)^{\varkappa}} \int_{(\eps+\ii (x-y))\amsmathbb R_{>0}} \frac{e^{-r}\dd r}{r^{1-\varkappa}} = \frac{\Gamma(\varkappa)}{\big(\eps+\ii (x-y)\big)^{\varkappa}},
\end{equation}
where in the denominator we choose the principal branch when computing the $\varkappa$-th power, \textit{i.e.} such that at $x=y$ we get a positive real number. Plugging \eqref{eq_x188} into \eqref{eq_x187} we get
\[
 \mathcal{J}_{\varkappa}[\boldsymbol{\nu}]=\text{Re}(\mathcal{J}_{\varkappa}[\boldsymbol{\nu}]) =\Gamma(\varkappa) \left(\sum_{g,h = 1}^H \theta_{g,h} \int_{\hat a'_g}^{\hat b'_g}\int_{\hat a'_{h}}^{\hat b'_{h}} \lim_{\eps \rightarrow 0^+} \text{Re}\,\big(\eps+\ii(x - y)\big)^{-\varkappa} \,\dd\nu_g(x)\dd\nu_{h}(y)\right),
\]
 where we use the given absolute convergence of the integrals to justify the interchange of the $\eps \rightarrow 0^+$ limit and integration. It remains to notice
\[
 \lim_{\eps \rightarrow 0^+} \text{Re}\,\big(\eps+\ii (x-y)\big)^{-\varkappa} =\frac{\cos\big(\frac{\pi \varkappa}{2}\big)}{|x-y|^{\varkappa}}.\qedhere
\]
\end{proof}

An important property of $-\mathcal I$ for developing the large deviation principles is that its lower level sets are compact. In this situation one says that $-\mathcal I$ is a \emph{good rate function}, \textit{cf.} \cite{Dembo_Zeitouni} (note the minus sign due to our convention). Here is the formal statement.

\begin{lemma}\label{Lemma_good_rate}
Under Assumptions~\ref{Assumptions_Theta} and \ref{Assumptions_basic}, for any $L\in\amsmathbb R$ the set $\big\{\boldsymbol{\nu} \in \mathscr{P}_\star \mid \mathcal{I}[\boldsymbol{\nu}]\geq L\big\} $ is compact.
\end{lemma}
\begin{proof}
 For $M\in\amsmathbb R$ and $\boldsymbol{\nu} \in \mathscr{P}_\star$ we define
 \begin{equation}
\label{eq_I_M}
 \mathcal I_M[\boldsymbol{\nu}]=
 \sum_{g,h = 1}^{H} \int_{\hat a'_g}^{\hat b'_g} \int_{\hat a'_{h}}^{\hat b'_{h}} \max\bigg( \theta_{g,h} \log|x-y|
 - \frac{V_g(x) + V_{h}(y)}{2 \sum_{i=1}^H \hat n_i },-M\bigg)\, \nu_g(x) \nu_{h}(y) \dd x \dd y.
\end{equation}
Note that $\mathcal I_M[\boldsymbol{\nu}]$ is a decreasing function of $M$ and we have by the monotone convergence theorem
 \begin{equation}
 \label{eq_x210}
 \lim_{M\rightarrow+\infty} \mathcal I_M[\boldsymbol{\nu}]=\inf_{M > 0} \mathcal I_M[\boldsymbol{\nu}]=\mathcal I[\boldsymbol{\nu}],
\end{equation}
where in the last equality we used the identity
\[
\sum_{g,h = 1}^{H} \int_{\hat a'_g}^{\hat b'_g} \int_{\hat a'_{h}}^{\hat b'_{h}}
 \frac{V_g(x) + V_{h}(y)}{2 \sum_{i=1}^H \hat n_i}\, \nu_g(x) \nu_{h}(y) \dd x \dd y = \sum_{h=1}^{H} \int_{\hat a'_h}^{\hat b'_h} V_h(x) \nu_h(x) \dd x.
\]
We deduce from \eqref{eq_x210} that
\begin{equation}
\label{eq_x200}
\big\{\boldsymbol{\nu}\in \mathscr{P}_\star \,\, | \,\,\mathcal{I}[\boldsymbol{\nu}]\geq L\big\}=\bigcap_{M>0} \big\{\boldsymbol{\nu}\in \mathscr{P}_\star \,\, | \,\, \mathcal{I}_M[\boldsymbol{\nu}] \geq L\big\}.
\end{equation}
While the function $\boldsymbol{\nu} \mapsto \mathcal I[\boldsymbol{\nu}]$ on $\mathscr{P}_{\star}$ is not continuous due to the singularity of the logarithm, for any $M > 0$ the function $\boldsymbol{\nu} \mapsto \mathcal{I}_M[\boldsymbol{\nu}]$ is continuous with respect to the weak topology. Hence, all sets in \eqref{eq_x200} are closed. In addition, by Lemma~\ref{Lemma_I_welldef}, for any $L\in\amsmathbb R$ we can find a finite constant $C_L$ depending on $L$ and on the constant in Lemma~\ref{Lemma_I_welldef} such that
\begin{equation}\label{eq_x201}
\big\{\boldsymbol{\nu}\in \mathscr{P}_\star \,\, |\,\, \mathcal{I}[\boldsymbol{\nu}]\geq L\big\}\subseteq\bigg\{\boldsymbol{\nu} \in \mathscr{P}_\star \,\, \bigg| \,\,\int_{\amsmathbb{R}} \log\big(1 + |x|\big)\nu(x)\dd x\leq C_L\bigg\}.
\end{equation}
As a corollary of Prokhorov theorem, the right-hand side of \eqref{eq_x201} is compact. Hence, the level set $\big\{\boldsymbol{\nu}\in \mathscr{P}_\star \,\,|\,\,\mathcal{I}[\boldsymbol{\nu}] \geq L\big\}$ is a closed subset of a compact set, which implies its compactness.
\end{proof}

One consequence of Lemma~\ref{Lemma_good_rate} is the existence of minimizers for $-\mathcal I$.

\begin{proposition} \label{Lemma_maximizer}
Under Assumptions~\ref{Assumptions_Theta} and \ref{Assumptions_basic},  $-\mathcal I$ achieves its minimum value
at a unique measure in $\mathscr{P}_\star$.
\end{proposition}

\begin{definition}
\label{DefEqmestwoconv} We denote the minimizer in Proposition~\ref{Lemma_maximizer} by $\boldsymbol{\mu} = (\mu_h)_{h = 1}^H$ and call it the \emph{equilibrium measure}. Following our conventions $\mu = \sum_{h = 1}^{H} \mu_h$ refers to the corresponding measure on $\amsmathbb{A}$.
\end{definition}

\begin{proof}[Proof of Proposition~\ref{Lemma_maximizer}]
Define
\[
\mathfrak i=\sup_{\boldsymbol{\nu}\in\mathscr{P}_\star} \mathcal I[\boldsymbol{\nu}].
\]
Lemma~\ref{Lemma_I_welldef} implies that the values of $\mathcal I[\boldsymbol{\nu}]$ are bounded
from above and therefore $\mathfrak i$ is a finite number. There exists a sequence
$(\boldsymbol{\mu}^{(n)})_{n \geq 0}$ of measures in $\mathscr{P}_\star$ such that
\begin{equation}
\label{eq_approximation_of_sup}
 \lim_{n\rightarrow\infty} \mathcal I[\boldsymbol{\mu}^{(n)}]=\mathfrak i\,.
\end{equation}
By Lemma~\ref{Lemma_good_rate}, up to extracting a subsequence we can
assume that the $H$-tuple of measures $\boldsymbol{\mu}^{(n)}$ weakly converge to an $H$-tuple of measures $\boldsymbol{\mu}\in \mathscr{P}_\star$.
Using the functional $\mathcal I_M$ of \eqref{eq_I_M}, we can write for any $M>0$
\begin{equation}
\label{eq_x203}
\mathfrak i=\lim_{n\rightarrow\infty} \mathcal I[\boldsymbol{\mu}^{(n)}]\leq  \lim_{n\rightarrow\infty} \mathcal I_M [\boldsymbol{\mu}^{(n)}]=\mathcal I_M [\boldsymbol{\mu}],
\end{equation}
where the last equality follows from the continuity of $\mathcal I_M$. We send $M\rightarrow\infty$ in the right-hand side of \eqref{eq_x203} using the monotone convergence theorem and get
\[
 \mathfrak i\leq \mathcal I[\boldsymbol{\mu}].
\]
Since $\mathfrak i$ was defined as the supremum, we conclude that $\mathcal
I[\boldsymbol{\mu}]=\mathfrak i$ and, hence, $\boldsymbol{\mu}$ is the desired minimizer.

\medskip

We proceed to the uniqueness part. If $\boldsymbol{\mu}^{0}$ and $\boldsymbol{\mu}^1$ are two minimizers of $-\mathcal I$, we have $\mathcal
I[\boldsymbol{\mu}^0]=\mathcal I[\boldsymbol{\mu}^1]= \mathfrak i$. Then, denoting $\boldsymbol{\mu}^{t} = (1 - t)\boldsymbol{\mu}^0 + t\boldsymbol{\mu}^1$ and $\boldsymbol{\nu} =
\boldsymbol{\mu}^1 - \boldsymbol{\mu}^0$, we have for any $t \in [0,1]$
\begin{equation}
\label{eq_x11} 0 = \mathcal{I}[\boldsymbol{\mu}^{t}] - (1 - t)\mathcal{I}[\boldsymbol{\mu}^0] - t\mathcal{I}[\boldsymbol{\mu}^1] = - t(1 - t)
\sum_{g,h = 1}^{H}\theta_{g,h} \int_{\hat a'_g}^{\hat b'_g} \int_{\hat a'_h}^{\hat b'_h} \log|x - y|\, \nu_g(x)\nu_{h}(y) \dd x \dd y.
\end{equation}
Up to the prefactor $-t(1-t)$ we would recognize $\mathcal{I}^{(2)}[\boldsymbol{\nu}]$ after we make sure that the integrals are absolutely convergent as in Lemma~\ref{Lemma_I_quadratic_smoothing}. For the singularity of $\log|x-y|$ on the diagonal, note that $\mu^0$ and $\mu^1$ have bounded densities, and therefore, so does $|\nu|$. Hence, the part of the integral corresponding to $|x-y|\leq 1$ is absolutely convergent. For the singularity of $\log|x-y|$ as $|x-y| \rightarrow +\infty$, we use Lemma~\ref{Lemma_I_welldef} to deduce that
\[
\forall h \in [H]\qquad \int_{\hat a'_h}^{\hat b'_h} \log\big(1 + |x|\big)|\nu_h(x)|\dd x< +\infty,
\]
Then, \eqref{eq_ln_upperbound} implies the desired absolute convergence. We can now apply Corollary~\ref{Corollary_I_positive} and learn that the right-hand side of \eqref{eq_x11} is zero if and only if $\boldsymbol{\nu} = 0$. This shows $\boldsymbol{\mu}^0 = \boldsymbol{\mu}^1$, \textit{i.e.} there is a unique minimizer.
 \end{proof}

We are now in position to establish the characterization of the equilibrium measure announced in Theorem~\ref{Theorem_equi_charact}, which we restate here for convenience. The notions of voids, bands and saturations were introduced in Definition~\ref{Definition_void_saturated}, and the effective potential in Definition~\ref{def_eff_pot}:
\[
V^{\textnormal{eff}}_h(x) = V_h(x) - \sum_{g = 1}^{H} 2\theta_{h,g} \int_{\hat{a}_g'}^{\hat{b}_g'} \log|x - y|\mu_g(y)\dd y.
\]

\begin{theorem} \label{Theorem_equi_charact_repeat_2}
 Suppose that Assumptions~\ref{Assumptions_Theta} and \ref{Assumptions_basic} hold.
 Then the equilibrium measure $\mu$ has a compact
 support, \textit{i.e.} its support is inside $[-D,D]$ with $D>0$ depending only on the constants in the assumptions. Moreover, there exists an $H$-tuple of constants $\boldsymbol{v}$ such that for any $h \in [H]$ and $x \in [\hat{a}'_h,\hat{b}_h']$, we have
\begin{align}
V^{{\textnormal{eff}}}_h(x) \geq v_h,& \quad \textnormal{if} \,\,x\,\,\textnormal{is in a void},\label{eq_void_inequality_repeat_2}\\
 V^{{\textnormal{eff}}}_h(x) \leq v_h, & \quad \textnormal{if}\,\,x\,\,\textnormal{is in a saturation},\label{eq_saturated_inequality_repeat_2}\\
 V^{{\textnormal{eff}}}_h(x) = v_h,& \quad \textnormal{if}\,\,x \,\,\textnormal{is in a band}. \label{eq_band_equality_repeat_2}
\end{align}
If each interval $[\hat a'_h, \hat b'_h]$ has at
least one band of the equilibrium measure, then $\boldsymbol{v}$ belongs to the $\amsmathbb{R}$-span of the $H$-dimensional vectors of coefficients of the linear forms $(\mathfrak{r}_{e})_{e = 1}^\mathfrak{e}$ in \eqref{eq_equations_eqs}, and the above conditions determine a unique measure in $\mathscr{P}_{\star}$.
\end{theorem}

\begin{remark} \label{remark:continuous} Proposition~\ref{Lemma_maximizer} and Theorem~\ref{Theorem_equi_charact_repeat_2} remain true if one minimizes $-\mathcal{I}$ over all $H$-tuples of nonnegative measures without an upper bound on their density. Some little details of the proofs are slightly different in that case as we should not assume existence of a density compared to the Lebesgue measure, but handling this is standard in potential theory. The characterization is simpler because saturations are absent.
\end{remark}

\begin{proof}[Proof of Theorem~\ref{Theorem_equi_charact_repeat_2}]
\textsc{Step 1.} Let $\nu$ be an integrable signed measure on $\amsmathbb{A}$ supported on a compact subset of the union of
 bands of $\mu$ and such that
 \begin{equation}
 \label{eq_x12}
\forall e \in [\mathfrak{e}]\qquad  \mathfrak r_e\big((\nu([\hat a'_h, \hat b'_h]))_{h = 1}^H\big)=0.
\end{equation}
For any $\eps \in \amsmathbb{R}$ with $|\eps|$ small enough, $\mu +\eps\nu$ is an absolutely continuous nonnegative measure on $\amsmathbb{A}$. Due to \eqref{eq_x12} it also satisfies the conditions \eqref{eq_equations_eqs}. Therefore $\boldsymbol{\mu} + \eps\boldsymbol{\nu} \in \mathscr{P}_\star$ and $\mathcal{I}[\boldsymbol{\mu}+\eps\boldsymbol{\nu}]\leq \mathcal{I}[\boldsymbol{\mu}]$. On the other hand, using the definition of $\mathcal{I}$ and the effective potential we
 have
 \begin{equation}
 \label{eq_x125}
 \mathcal{I}[\boldsymbol{\mu}+\eps\boldsymbol{\nu}] = \mathcal{I}[\boldsymbol{\mu}]-\eps \sum_{h=1}^H \int_{\hat a'_h }^{\hat b'_h} V_h^{\textnormal{eff}}(x)\, \nu_h(x) \dd x+
 \eps^2 \sum_{g,h=1}^H \theta_{g,h} \int_{\hat a'_g}^{\hat b'_g} \int_{\hat a'_{h}}^{\hat b'_{h}}
 \log|x-y|\, \nu_g(x) \nu_{h}(y) \dd x \dd y.
 \end{equation}
Therefore, the validity of the inequality $\mathcal{I}[\boldsymbol{\mu}+\eps\boldsymbol{\nu}]\leq \mathcal{I}[\boldsymbol{\mu}]$ for arbitrarily small positive and negative values $\eps$ implies
 \begin{equation}
 \label{eq_x13}
 \sum_{h=1}^H \int_{\hat a'_h }^{\hat b'_h} V_h^{\textnormal{eff}}(x)\, \nu_h(x)\dd x=0.
 \end{equation}
We observe that $x \mapsto V_h^{\textnormal{eff}}(x)$ is a continuous function on $[\hat{a}_h',\hat{b}_h']$. If this function took two different values at points $x_1 \neq x_2$ belonging to bands --- possibly two different bands --- of the $h$-th segment, we could construct a zero-mass integrable signed measure $\nu$ supported on small neighborhoods of $x_1$
 and $x_2$ in $[\hat{a}'_h,\hat{b}'_h]$, and for which \eqref{eq_x13} would not hold --- \eqref{eq_x12} is guaranteed to hold because $\mathfrak{r}_e$ for $e \in [\mathfrak{e}]$ are linear forms and all the components of $\boldsymbol{\nu}$ have zero mass. We conclude that
 for each $h \in [H]$, the function $V_h^{\textnormal{eff}}$ assumes a constant value $v_h$ on the bands in $[\hat a'_h,
 \hat b'_h]$. Then, since $\nu$ is supported on bands, we can replace $V_h^{\textnormal{eff}}(x)$ by $v_h$ in \eqref{eq_x13} and get
 \begin{equation}
 \label{eq_x14}
 \sum_{h=1}^H v_h\,\nu([\hat a'_h, \hat b'_h])=0.
 \end{equation}
If $\mu$ has at least one band in each $[\hat a'_h, \hat b'_h]$, the mass $\nu([\hat a'_h, \hat b'_h])$
 can take arbitrary small values satisfying \eqref{eq_x12}. Then \eqref{eq_x14} will hold if
 and only if the vector $\boldsymbol{v} = (v_h)_{h = 1}^H$ belongs to the linear span of the coefficients of the linear
 forms $(\mathfrak r_e)_{e = 1}^{\mathfrak{e}}$.

Let us now prove the inequalities $V_h^{\textnormal{eff}}(x)\geq v_h$ in voids and $V_h^{\textnormal{eff}}(x)\leq v_h$ in
 saturations. They readily follow by allowing $\nu$ to be supported in voids/saturations
as in \cite{DS,ST}. We should still have $\boldsymbol{\mu}+\eps\boldsymbol{\nu} \in \mathscr{P}_\star$, \textit{i.e.} the perturbed measure should be nonnegative with density
bounded from above by $\frac{1}{\theta_{h,h}}$ in $[\hat a'_h, \hat b'_h]$. Hence, the restriction of $\nu$ to voids must be a nonnegative measure and the restriction of $\nu$ to saturations must be a nonpositive measure. Thus, we can no longer freely change the sign of $\nu$ and \eqref{eq_x13} is only nonnegative, rather than being equal to $0$. Then we only get inequalities, but not equalities to $v_h$.

In more details, take $x_1$ in a saturation or in a band and take $x_2$ in a void or in a band, with both $x_1,x_2 \in [\hat{a}_h',\hat{b}_h']$ for some $h \in [H]$. The properties \eqref{eq_void_inequality_repeat_2} and \eqref{eq_saturated_inequality_repeat_2} say that $V_h^{\textnormal{eff}}(x_1)\leq V_h^{\textnormal{eff}}(x_2)$. This inequality is proven by choosing $\nu_h$ to have a negative density in a small neighborhood of $x_1$ and to have a positive density in a small neighborhood of $x_2$ so that the total mass is zero, and setting $\nu(x)=0$ at points $x$ outside these neighborhoods. Combining \eqref{eq_x125} with the inequality $\mathcal{I}[\boldsymbol{\mu}+\eps \boldsymbol{\nu}]\leq \mathcal{I}[\boldsymbol{\mu}]$ we get
 \begin{equation}
 \label{eq_x209}
 \int_{\hat a'_h }^{\hat b'_h} V_h^{\textnormal{eff}}(x)\, \nu_h(x)\dd x\geq 0.
 \end{equation}
Since $V_h^{\textnormal{eff}}$ is continuous and $\nu_h([\hat a'_h,\hat b'_h])=0$, the inequality \eqref{eq_x209} cannot hold if $V_h^{\textnormal{eff}}(x_1)>V_h^{\textnormal{eff}}(x_2)$.

\medskip

\noindent \textsc{Step 2.} Let us prove the compactness of the support of the equilibrium measure. This needs to be checked only if $a_1$ or $b_H$ are infinite. We can assume without loss of generality that $b_H = +\infty$. In this situation we will prove two inequalities. First, for some constant $D > 0$, we have
 \begin{equation}
 \label{eq_potential_eff_lower bound}
\forall x \in (D,+\infty) \qquad V^{\textnormal{eff}}_H(x) \geq -D + \eta_H \log x,
 \end{equation}
 where $\eta_H>0$ is the constant in \eqref{eq_confinement_right} of Assumption~\ref{Assumptions_basic}. Second, for some constants $D'\in\amsmathbb R$ and $M > 0$ larger than the total mass of $\mu$, we have
 \begin{equation}
 \label{eq_potential_eff_upper bound}
\forall x\in[\hat a'_H+1, \hat a'_H+1+M \theta_{H,H}] \qquad V^{\textnormal{eff}}_H(x) \leq D'.
 \end{equation}
 The choice of the constants $D,D',M$ can be made depending only on the constants in Assumptions~\ref{Assumptions_Theta} and \ref{Assumptions_basic}.
 In \eqref{eq_potential_eff_upper bound} we silently assumed that $a_H$ is finite --- this is automatic if $H>1$, but might fail for $H=1$. If $a_H=-\infty$, then we can instead prove and use \eqref{eq_potential_eff_upper bound} for $x\in [0, 1+\theta_{H,H}]$.

 Before justifying \eqref{eq_potential_eff_lower bound} and \eqref{eq_potential_eff_upper bound}, let us show how they imply compactness. On the one hand, due to the choice of $M$ and the upper bound $\frac{1}{\theta_{H,H}}$ for the density of $\mu_H$, it is impossible for the entire segment $[\hat a'_H+1, \hat a'_H+1+M \theta_{H,H}]$ to be saturated. Hence, \eqref{eq_potential_eff_upper bound} together with \eqref{eq_void_inequality_repeat_2} and \eqref{eq_band_equality_repeat_2} imply $v_H\leq D'$. On the other hand, \eqref{eq_potential_eff_lower bound} implies that $ V^{\textnormal{eff}}_H(x)>D'$ for $x \geq x_0$ for some large enough choice of $x_0 > 0$. In particular, $ V^{\textnormal{eff}}_H(x)>v_H$ for $x\geq x_0$. Hence, by \eqref{eq_band_equality_repeat_2} and \eqref{eq_saturated_inequality_repeat_2}, saturations and bands cannot intersect $[x_0,+\infty)$; in other words, the entire segment $[x_0,+\infty)$ is included in a void. A similar argument holds for $V_1^{\textnormal{eff}}$ in the case $a_1 = -\infty$. Thus $\mu$ has support included in $[-x_0,x_0]$ for some $x_0> 0$ depending only on the constants in the assumptions.

 \smallskip

To justify the inequality \eqref{eq_potential_eff_lower bound}, we first assume that either $a_1$ is finite or $\theta_{H,1}\geq 0$. We proved \eqref{eq_ln_better_bound} under the assumption that either $\theta_{H,1}\geq 0$ or $b_H$ is finite, yet the proof for the case of finite $a_1$ is the same. Then, using the inequality \eqref{eq_ln_better_bound} there exist $C,C'>0$ such that for $x>\hat{a}'_H$
\begin{equation}
\label{eq_x15}
\begin{split}
 V^{\textnormal{eff}}_H (x) & = V_H(x)-2 \sum_{h=1}^H \theta_{H,h} \int_{\hat a'_h}^{\hat
 b'_h} \log|x-y|\, \mu_h(y)\dd y \\
 & \geq V_H(x)-2 \sum_{h=1}^H \theta_{H,h} \int_{\hat
 a'_h}^{\hat b'_h} \Big(\log\big(1 + |x|\big)+\log\big(1 + |y|\big)\Big)\mu_h(y)\dd y-C \\
 & \geq V_H(x)- 2 \log\big(1 + |x|)\bigg(\sum_{h=1}^H \theta_{H,h}\, \mu([\hat a'_h,\hat b'_h])\bigg)- C',
 \end{split}
 \end{equation}
In the last line we used the fact $\int_{\amsmathbb{R}} \log\big(1 + |y|\big)\mu(y)\dd y< + \infty$, which follows from Lemma
\ref{Lemma_I_welldef} by bounding $\I[\boldsymbol{\mu}]$ from below by the value of $\I$ on some prescribed measure. By Conditions 3. and \eqref{eq_confinement_right} in Assumption~\ref{Assumptions_basic}, we can find a constant $D > \hat{a}'_H$ such that the claimed inequality \eqref{eq_potential_eff_lower bound} holds.

If $a_1 =-\infty$ and $\theta_{H,1}<0$, then we instead have using \eqref{eq_ln_second} for $h=1$,
\begin{equation*}
\begin{split}
 V^{\textnormal{eff}}_H (x) & \geq V_H(x)-2 \sum_{h=2}^H \theta_{H,h} \int_{\hat
 a'_h}^{\hat b'_h} \Big(\log\big(1 + |x|\big)+\log\big(1 + |y|\big)\Big)\mu_h(y)\dd y \\
 & \quad\phantom{V_H(x)-2} - \theta_{H,1}
 \int_{\hat a'_1}^{\hat b'_1} \Big(\log\big(1 + |x|\big)+ \log\big(1 + |y|\big)\Big)\mu_h(y)\dd y - C \\ & \geq V_H(x)- 2 \log\big(1 + |x|\big)\bigg(\sum_{h=2}^H \theta_{H,h}\mu([\hat a'_h,\hat b'_h]) \bigg) - \log\big(1 + |x|\big) \theta_{H,1} \mu([\hat a'_1,\hat b'_1]) - C,
\end{split}
\end{equation*}
which again leads to the desired \eqref{eq_potential_eff_lower bound} with help of \eqref{eq_confinement_right_refined} in Assumption
\ref{Assumptions_basic}.

To justify the inequality \eqref{eq_potential_eff_upper bound}, we notice that $V_H$ is uniformly bounded on $[\hat a'_H+1, \hat a'_H+1+M \theta_{H,H}]$ by \eqref{V_bound} in Assumption~\ref{Assumptions_basic}. As for the integrals in the definition of $V^{\textnormal{eff}}_H (x)$, we have for any $h\in [H]$
\begin{equation}
\label{eq_x202}
 \bigg| \int_{\hat a'_h}^{\hat b'_h} \log|x-y|\, \mu_h(y) \dd y\bigg|\leq \frac{1}{\theta_{h,h}} \int_{-1}^1 \big|\log|\xi|\big| \dd \xi + \int_{\hat a'_h}^{\hat b'_h} \log\big(1 + |x|\big)\mu_h(y)\dd y + \int_{\hat a'_h}^{\hat b'_h} \log\big(1 + |y|\big)\mu_h(y)\dd y,
\end{equation}
where the first term comes from the integration over $y\in[x-1,x+1]$ and the next two terms correspond to the two terms in \eqref{eq_ln_upperbound} for $y$ outside $[x-1,x+1]$. All three terms in the right-hand side of \eqref{eq_x202} are uniformly bounded: the first one by integrability of the logarithm and boundedness of $\frac{1}{\theta_{h,h}}$ included in Condition 2. in Assumption~\ref{Assumptions_Theta}; the second one by the uniform bound on the mass of $\mu_h$ due to Condition 1. in Assumption~\ref{Assumptions_basic}; the third one by Lemma
\ref{Lemma_I_welldef} and by bounding $\I[\boldsymbol{\mu}]$ by the value of $\I$ on some prescribed measure. We conclude that $V^{\textnormal{eff}}_H (x)$ is uniformly bounded on $[\hat a'_H+1, \hat a'_H+1+M\theta_{H,H}]$, which implies the claimed inequality \eqref{eq_potential_eff_upper bound}.

\medskip

\noindent \textsc{Step 3.} We now turn to the uniqueness statement. Take $\boldsymbol{\mu} \in \mathscr{P}_{\star}$ satisfying all the conditions stated in Theorem~\ref{Theorem_equi_charact_repeat_2} and let $\boldsymbol{\mu}+\boldsymbol{\nu}$ be
another measure in $\mathscr{P}_\star$. We claim that in this case the second and
the third terms in the right-hand side of \eqref{eq_x125} for $\eps=1$ are both
negative. This would imply that $\boldsymbol{\mu}$ is indeed a minimizer, which is unique by
Proposition~\ref{Lemma_maximizer}. Indeed, the third term is negative by Corollary~\ref{Corollary_I_positive}. For the second term, due to \eqref{eq_x14} we can replace
$V_h^{\textnormal{eff}}(x)$ with $V_h^{\textnormal{eff}}(x)-v_h$ without changing the sum. Then, considering separately voids, bands, and saturations, we
see that the integrand is nonnegative on each of them, which finishes the proof.
\end{proof}
\begin{remark}
\label{remark:v_when_fully_void_or_sat}
 If a segment $ [\hat{a}'_h,\hat{b}_h']$, $h\in[H]$ is fully void, then \eqref{eq_void_inequality_repeat_2} in Theorem~\ref{Theorem_equi_charact_repeat_2} only says that the effective potential is lower-bounded on $ [\hat{a}'_h,\hat{b}_h']$. If a segment $ [\hat{a}'_h,\hat{b}_h']$ is fully saturated, then \eqref{eq_saturated_inequality_repeat_2} in Theorem~\ref{Theorem_equi_charact_repeat_2} only says that the effective potential is upper-bounded on $[\hat{a}'_h,\hat{b}_h']$. Both statements are easily seen to be true.
\end{remark}

\subsection{Finiteness of the partition function}

We eventually are in position to justify the finiteness of the partition function so that the discrete ensembles we wish to study exist. Recall from Definition~\ref{def_empirical_mes} the random empirical measure $\boldsymbol{\mu}_{\N} = (\mu_{\N,h})_{h = 1}^{H}$, where
\[
\mu_{\N,h} := \frac{1}{\N} \sum_{i = 1}^{N_h} \delta_{\N^{-1}\ell_i^h},\qquad \mu_{\N} := \sum_{h = 1}^{H} \mu_{\N,h}.
\]

\begin{proposition}
\label{lem_finite_partition_function} If Assumptions~\ref{Assumptions_Theta} and \ref{Assumptions_basic} hold, then for $\N$ large enough (compared to the constants in the assumptions), the partition function
\[
\Z_{\N} := \sum_{\boldsymbol{\ell} \in \amsmathbb{W}_{\N}} \prod_{1 \leq i < j \leq N} \frac{1}{\N^{2\theta_{h(i),h(j)}}}\cdot \frac{\Gamma\big(\ell_j - \ell_i + 1\big)\cdot\Gamma\big(\ell_j - \ell_i + \theta_{h(i),h(j)}\big)}{\Gamma\big(\ell_j - \ell_i\big)\cdot\Gamma\big(\ell_j - \ell_i + 1 - \theta_{h(i),h(j)}\big)} \cdot \prod_{i = 1}^{N} w_{h(i)}(\ell_i)
\]
is finite. Hence
\begin{equation}
\label{PNprobarepeat} \amsmathbb{P}_{\N}(\boldsymbol{\ell}) := \frac{1}{\Z_{\N}} \prod_{1 \leq i < j \leq N} \frac{1}{\N^{2\theta_{h(i),h(j)}}}\cdot\frac{\Gamma\big(\ell_j - \ell_i + 1\big)\cdot\Gamma\big(\ell_j - \ell_i + \theta_{h(i),h(j)}\big)}{\Gamma\big(\ell_j - \ell_i\big)\cdot\Gamma\big(\ell_j - \ell_i + 1 - \theta_{h(i),h(j)}\big)} \cdot \prod_{i = 1}^{N} w_{h(i)}(\ell_i)
\end{equation}
is a well-defined probability measure on the configuration set $\amsmathbb{W}_{\N}$.
\end{proposition}
\begin{proof}
We consider the regularized energy of the empirical measure
\begin{equation}
\label{eq_x207}
 \tilde {\mathcal I}[\boldsymbol{\mu}_\N]= 2 \sum_{1\leq i<j \leq N} \theta_{h(i),h(j)} \log\bigg( \frac{\ell_j - \ell_i}{\N}\bigg)- \sum_{i=1}^N V_{h(i)}\bigg(\frac{\ell_i}{\N}\bigg).
\end{equation}
The functional $\tilde {\mathcal I}[\boldsymbol{\mu}_\N]$ is almost the same as the value of the functional $\mathcal I$ on the empirical measure. The only two differences are in the integration only between $\hat a'_h$ and $\hat b'_h$ in the definition of $\mathcal I$ --- a few particle positions $\frac{\ell_i}{\N}$ might be outside these segments due to distinction between $[\hat a'_h,\hat b'_h]$ and $[\hat a_h, \hat b_h]$ --- and in the exclusion of the explosive terms $\log\big(\frac{\ell_i}{\N}-\frac{\ell_i}{\N}\big)$.

For some constant $C > 0$ depending only on the constants in the assumptions, the unnormalized measure\label{index:notnorm} $\widetilde{\amsmathbb{P}}_{\N} = \Z_\N \cdot \amsmathbb{P}_\N$ assigns to a configuration $\boldsymbol{\ell} \in \amsmathbb{W}_{\N}$ the weight
\begin{equation}
\label{prodijN}
\begin{split}
\widetilde{\amsmathbb{P}}_{\N}(\boldsymbol{\ell}) & := \prod_{1 \leq i < j \leq N} \frac{1}{\N^{2\theta_{h(i),h(j)}}}\cdot \frac{\Gamma\big(\ell_j - \ell_i + 1\big) \cdot \Gamma\big(\ell_j - \ell_i + \theta_{h(i),h(j)}\big)}{\Gamma\big(\ell_j - \ell_i\big) \cdot \Gamma\big(\ell_j - \ell_i + 1 - \theta_{h(i),h(j)}\big)}\cdot \prod_{i = 1}^{N} w_{h(i)}(\ell_i) \\
& \leq \exp\bigg(\mathcal{N}^2 \tilde{\mathcal I}[\boldsymbol{\mu}_{\N}] + C\N \log \N + \sum_{1 \leq i < j \leq N} \frac{c_1}{\ell_j - \ell_i}\bigg),
\end{split}
\end{equation}
where we have used the upper bound in \eqref{eq_Stirling_xy} for the pairwise interaction ($c_1$ comes from there), replaced according to Assumption~\ref{Assumptions_basic} the weight with the potential \eqref{eq_weight_form} up to the error term \eqref{error_bound_e} that contributes to the $\N\log \N$ term, and recognized the regularized energy of the empirical measure \eqref{eq_x207}. Besides, we have
\[
\sum_{1 \leq i < j \leq N} \frac{1}{\ell_j - \ell_i} \leq \frac{1}{\min_{h} |\theta_{h,h}|} \sum_{1 \leq i < j \leq N} \frac{1}{j - i} \leq C'\N \log \N.
\]
for some constant $C' > 0$.

Next, we claim that there exists another constant $C>0$ such that
\begin{equation}
\label{eq_x208}
\forall \boldsymbol{\ell} \in \amsmathbb{W}_{\N}\qquad \tilde{\mathcal I}[\boldsymbol{\mu}_N]\leq C - \frac{1}{C} \int_{\amsmathbb{R}} \log\big(1 + |x|\big)\dd\mu_\N(x).
\end{equation}
Indeed, we can repeat the proof of Lemma~\ref{Lemma_I_welldef}. The only part of this proof which used existence of the density of $\nu$ (obviously failing for $\mu_\N$) is the convergence of the double integral involving $\log|x-y|$ on the diagonal $x=y$, however, we simply excluded the diagonal when defining $\tilde {\mathcal I}$.

Using \eqref{prodijN} and \eqref{eq_x208} we deduce that for a larger constant $C > 0$
\[
\widetilde{\amsmathbb{P}}_{\N}(\boldsymbol{\ell}) \leq \exp(C\N\log \N) \cdot \prod_{i = 1}^{N} \bigg(1 + \bigg|\frac{\ell_i}{\N}\bigg|\bigg)^{-\frac{\N}{C}}.
\]
By design of the state space $\amsmathbb{W}_{\N}$ and by Assumption~\ref{Assumptions_basic}, there exists $\boldsymbol{d} \in \amsmathbb{R}^H$ such that $|\!|\boldsymbol{d}|\!|_{\infty} \leq \N C$ and for any $\boldsymbol{\ell} \in \amsmathbb{W}_{\N}$ and $i \in [N]$ we have
\[
\ell_i = \lambda_i + (i - 1)\theta_{h(i),h(i)} + d_{h(i)}\qquad \textnormal{with} \quad \lambda_i \in \amsmathbb{Z}.
\]
Therefore,
\[
\Z_{\N} := \sum_{\boldsymbol{\ell} \in \amsmathbb{W}_{\N}} \widetilde{\amsmathbb{P}}_{\N}(\boldsymbol{\ell}) \leq \exp(C\N\log \N) \cdot \prod_{i = 1}^{N} \Bigg( \sum_{\lambda_i \in \amsmathbb{Z}} \bigg(1 + \frac{1}{\N}\bigg|\lambda_i + (i - 1)\theta_{h(i),h(i)} + d_{h(i)}\bigg|\bigg)^{-\frac{\N}{C}}\Bigg),
\]
and each sum over $\lambda_i \in \amsmathbb{Z}$ is convergent for $\N > C$.
\end{proof}

\section{Example with Gaussian weights}
\label{mueqGaussian}

A recurring example in this book is the discrete ensemble with Gaussian weights. This is the ensemble depending on two real parameters $\theta, \kappa > 0$ and corresponding in the notations of Section~\ref{Chapter_Setup_and_Examples} to
\[
H = 1,\qquad [a_1,b_1] = \amsmathbb{R},\qquad \theta_{1,1}=\theta,\qquad w(x) = e^{-NV(\frac{x}{N})} =e^{- \kappa\frac{x^{2}}{N}}, \qquad \N = N.
\]
In this simple case, we describe more precisely the equilibrium measure introduced in the previous section and characterized in Theorem~\ref{Theorem_equi_charact_repeat_2}. A plot of the resulting formulae was given in Figure~\ref{Fig_Gaussian_intro}.
 To express it, we need the complete elliptic integrals of the first, second, and third kind, denoted respectively $K(\mathsf{k})$, $E(\mathsf{k})$ and $\Pi_1[\mathsf{c};\mathsf{k}]$, and depending on the elliptic modulus $\mathsf{k} \in [0,1)$ through:
\begin{equation}
\label{index:ellipticEK}\begin{split}
K(\mathsf{k}) & := \int_{0}^{1} \frac{\dd t}{\sqrt{(1 - \mathsf{k}^2t^2)(1 - t^2)}}, \\
E(\mathsf{k}) & := \int_{0}^{1} \dd t\,\sqrt{\frac{1 - \mathsf{k}^2t^2}{1 - t^2}}, \\
\Pi_1[\mathsf{c};\mathsf{k}] & := \int_{0}^{1} \frac{\dd t}{(1 + \mathsf{c}t^2)\sqrt{(1 - \mathsf{k}^2t^2)(1 - t^2)}}.
\end{split}
\end{equation}

\begin{proposition} \label{Prop_Gaussian_LLN}
The measure $\mu$ has a density which is given, if $\sqrt{2 \theta \kappa}\leq \pi$, by
\begin{equation}
\label{outsidesat} \mu(x)= \left\{\begin{array}{lll}
  \dfrac{\kappa}{\pi \theta}\sqrt{\dfrac{2\theta}{\kappa}- x^2 } & & |x| \leq \sqrt{\dfrac{2\theta}{\kappa}}, \\[10 pt] 0 & & \textnormal{otherwise,} \end{array}\right.
\end{equation}
 and if $\sqrt{2 \theta \kappa}>{\pi}$, by
\begin{equation}
\label{insidesat} \mu(x)= \left\{\begin{array}{lll}
 \dfrac{1}{\theta} & & |x| \leq \alpha \\[10pt]
  \dfrac{2 \sqrt{(\beta^2 - x^2)(x^2 - \alpha^2)}}{\pi \theta \beta |x|}\,\Pi_1\bigg[-\dfrac{\alpha^2}{x^2}\,;\,\mathsf{k}\bigg] & & |x| \in (\alpha,\beta) \\[10pt]
  0& & |x|\geq \beta \end{array}\right. ,
\end{equation}
 where $\beta$ and $\mathsf{k} = \frac{\alpha}{\beta}$ are determined by the system of two equations:
\begin{equation}
\label{theconstatint} \kappa \beta= 2K(\mathsf{k}),\qquad \kappa\theta = 2K(\mathsf{k})\big(2E(\mathsf{k}) - (1 - \mathsf{k}^2)K(\mathsf{k})\big).
\end{equation}
\end{proposition}
\begin{proof}
The existence and uniqueness of the equilibrium measure $\mu$ in this situation was proved in \cite[Theorem~5.4]{BGG}. The characterization of this minimizer is exposed later in Theorem~\ref{Theorem_equi_charact_repeat_2}. The computation of $\mu$ was done explicitly in \cite{DK}, where some formulae unfortunately have misprints which we correct below. Let us sketch the steps of this computation. The minimization problem implies that if $x$ is in a band, we must have the equality
\begin{equation}
\label{pvint} \textnormal{p.v.} \int \frac{2\theta\,\dd\mu(y)}{x - y} = V'(x) = 2 \kappa x,
\end{equation}
 where p.v. indicates the Cauchy principal value. To solve this problem, the usual strategy is to assume the location of bands and saturations are known, then obtain a measure solving the linear equation \eqref{pvint}. The remaining constraints posed by the positivity properties for this measure to satisfy the remaining constraints in the minimization problem is equivalent to certain constraints on the endpoints of bands and saturations. If for given $\theta$ and $\kappa$ one can find such endpoints matching all the constraints (among which the positivity of $\mu$), the solution of \eqref{pvint} we have found for this choice of endpoints must coincide with the unique and sought-for solution of the minimization problem.

If $\mu$ has no saturations, it must be equal to the solution of the minimization problem without the constraint $\mu \leq \frac{1}{\theta}$: this is the famous semi-circle law with density \eqref{outsidesat}. The constraint $\mu(x) \leq \frac{1}{\theta}$ for all $x$ requires $\sqrt{2\theta \kappa} \leq \pi$. If $\sqrt{2\theta \kappa} > \pi$, we rather assume that $\mu$ has two bands separated by a saturation --- and we will check at the end the validity of this ansatz. Since $V'$ is odd and the solution of the minimization problem is unique, $\mu$ must be invariant under $x \mapsto -x$. We can therefore assume that the saturation is of the form $[-\alpha,\alpha]$ while the bands are of the form $(-\beta,-\alpha)$ and $(\alpha,\beta)$ for some $\beta > \alpha > 0$. The equation \eqref{pvint} implies the following equation for the restriction $\tilde{\mu}$ of $\mu$ to the bands and $x \in (-\beta,-\alpha) \cup (\alpha,\beta)$
\[
\textnormal{p.v.} \int \frac{\theta\,\dd\tilde{\mu}(y)}{x - y} = \kappa x + \log\bigg(\frac{x - \alpha}{x + \alpha}\bigg).
\]
The solution of this equation in terms of the Stieltjes transform is given by \cite[Equation 31]{DK} where the logarithm appears with a wrong sign. The correct formula is
\begin{equation}
\label{Wmumf} \mathcal{G}_{\tilde{\mu}}(x) := \int \frac{\dd\tilde{\mu}(y)}{x - y} = \frac{1}{\theta}\Bigg(\kappa x + \log\bigg(\frac{x - \alpha}{x + \alpha}\bigg) + \int_{-\alpha}^{\alpha} \frac{\dd y}{x - y}\,\frac{\sigma(x)}{\sigma(y)}\Bigg),
\end{equation}
where
\[
\sigma(x) := \sqrt{(x^2 - \alpha^2)(x^2 - \beta^2)}
\]
 and the branch of the square root is fixed by requiring that $\sigma(x)$ is holomorphic for $x \in \amsmathbb{C}\setminus [-\beta,-\alpha] \cup [\alpha,\beta]$ and $\sigma(x) \sim x$ as $x \rightarrow \infty$. This is equivalent to
\[
\forall x \in \amsmathbb{R} \qquad \frac{\dd\tilde{\mu}(x)}{\dd x} = \mathbbm{1}_{(\alpha,\beta)}(|x|) \cdot\frac{\textnormal{sgn}(x)}{\theta \pi} \int_{-\alpha}^{\alpha} \frac{\dd y}{x - y}\,\frac{\sqrt{(\beta^2 - x^2)(x^2 - \alpha^2)}}{\sqrt{(\beta^2 - y^2)(y^2 - \alpha^2)}}.
\]
By setting $\mathsf{k} = \frac{\alpha}{\beta}$, changing the integration variable to $t = \frac{y}{\alpha}$ and setting $\tilde{x} = \frac{x}{\alpha}$, we find
\begin{equation}
\label{muposs}\frac{\dd\tilde{\mu}(x)}{\dd x} =\frac{ \textnormal{sgn}(x)\,\sqrt{(\beta^2 - x^2)(x^2 - \alpha^2)}}{\pi \theta \alpha \beta} \int_{0}^{1} \dd t\bigg(\frac{1}{\frac{x}{\alpha} - t} + \frac{1}{\frac{x}{\alpha} + t}\bigg)\,\frac{1}{\sqrt{(1 - t^2)(1 - \mathsf{k}^2t^2)}},
\end{equation}
which can be rewritten in terms of elliptic functions as \eqref{insidesat}. The constraint that $\mu$ has total mass $1$ reads
\begin{equation}
\label{teoungfdb}\int \frac{\dd \mu(y)}{z - y} \,\,\mathop{\sim}_{z \rightarrow \infty}\,\, \frac{1}{z}.
\end{equation}
From \eqref{Wmumf} we find that $\int \frac{\theta\,\dd \mu(y)}{z - y} =  \theta\mathcal{G}_{\tilde{\mu}}(z) - \log\big(\frac{z - \alpha}{z + \alpha}\big)$ behaves as $z \rightarrow \infty$ like
\[
\bigg(\kappa - \int_{-\alpha}^{\alpha} \frac{\dd y}{\sqrt{(\alpha^2 - y^2)(\beta^2 - y^2)}}\bigg)z + \bigg(\int_{-\alpha}^{\alpha} \dd y\,\frac{\frac{\alpha^2 + \beta^2}{2} - y^2}{\sqrt{(\alpha^2 - y^2)(\beta^2 - y^2)}}\bigg)\frac{1}{z} + O\bigg(\frac{1}{z^2}\bigg).
\]
The coefficients can be rewritten in terms of elliptic integrals, and comparing to \eqref{teoungfdb} gives the two constraints
\[
\frac{\beta}{2} = \frac{K(\mathsf{k})}{\kappa},\qquad \theta = \beta\big(2E(\mathsf{k}) - (1 - \mathsf{k}^2)K(\mathsf{k})\big).
\]
Eliminating $\alpha$ between the two constraints give the claimed \eqref{theconstatint}, which is also \cite[Equation 36]{DK}. The function
\[
\mathsf{k} \mapsto 4K(\mathsf{k})\big(E(\mathsf{k}) - (1 - \mathsf{k}^2)K(\mathsf{k})\big),\qquad
\]
is strictly increasing from the domain $\mathsf{k} \in (0,1)$ to its image $(\pi^2,+\infty)$. Therefore, for any $\theta,\kappa > 0$ such that $\sqrt{2\theta \kappa} > \pi$, there exists a unique $\mathsf{k} \in (0,1)$ satisfying \eqref{theconstatint} and therefore a unique choice of endpoints $\beta= \frac{2K(\mathsf{k})}{\kappa}$ and $\alpha = \mathsf{k}\beta$ making $\tilde{\mu}$ a nonnegative measure of mass $1$. The proof that this is the solution of the minimization problem if $\sqrt{2\theta \kappa} > \pi$ will be complete after we justify that the effective potential
\[
V^{\textnormal{eff}}(x) = V(x) - 2\theta \int \dd y\,\log|x - y|\dd\mu(y)
\]
is nonnegative in the voids and nonpositive in the saturation.

We already know that $V_{\textnormal{eff}}(x)$ is equal to a constant $v$ in the bands. For $x$ outside the bands we have
\[
\frac{(V^{\textnormal{eff}})'(x)}{2} = \kappa x + \log\bigg(\frac{x - \alpha}{x + \alpha}\bigg) - \theta \mathcal{G}_{\tilde{\mu}}(x).
\]
Starting from \eqref{Wmumf}, we see that
\[
(V^{\textnormal{eff}})'(x) = - \textnormal{p.v.} \int_{-\alpha}^{\alpha} \frac{\dd y\,f_{x}(y)}{x - y},\qquad f_{x}(y) = \frac{\sigma(x)}{\sigma(y)}.
\]
Let us examine the saturation first. For $x \in (-\alpha,\alpha)$, we have
\begin{equation}
\label{1110}
\begin{split}
\frac{(V^{\textnormal{eff}})'(x)}{2} & = \int_{-\alpha}^{\alpha} \dd y\,\frac{f_{x}(x) - f_{x}(y)}{x - y} - \textnormal{p.v.} \int_{-\alpha}^{\alpha} \frac{\dd y}{x - y} \\
& = \int_{-\alpha}^{\alpha} \dd y \bigg(\int_{0}^{1} \dd t\,f_{x}'(ty + (1 - t)x)\bigg) - \log\bigg(\frac{\alpha - x}{\alpha + x}\bigg) \\
& = \int_{0}^{1} \frac{\dd t}{t} \big(f_{x}(t\alpha + (1 - t)x) - f_{x}(-t\alpha + (1 - t)x)\big) - \log\bigg(\frac{\alpha - x}{\alpha + x}\bigg).
\end{split}
\end{equation}
Notice that $f_{x}(y)$ is an even function of $y$ which is increasing for $y \in [0,\alpha]$. If $x \in [0,\alpha]$, we have for any $t \in [0,1]$ that
\[
|t\alpha + (1 - t)x| \geq |-t\alpha + (1 - t)x|.
\]
Hence
\[
f_{x}(t\alpha + (1 - t)x) \geq f_{x}(-t\alpha + (1 - t)x),
\]
thus the integral term in \eqref{1110} is nonnegative. As the logarithmic term is also nonnegative, we obtain $(V^{\textnormal{eff}})'(x) \geq 0$. Besides, we also have
\begin{equation}
\label{ineein0} (V^{\textnormal{eff}})'(x) \geq -2\log\bigg(\frac{\alpha - x}{\alpha + x}\bigg) \rightarrow +\infty,
\end{equation}
therefore $\lim_{x \rightarrow \alpha^{-}} (V^{\textnormal{eff}})'(x) = +\infty$. Finally, notice that $(V^{\textnormal{eff}})'(-x) = - (V^{\textnormal{eff}})'(x)$. These properties imply that
\begin{equation}
\label{ineein1}\forall x \in (-\alpha,\alpha)\qquad V^{\textnormal{eff}}(x) < v.
\end{equation}
A similar argument for the voids, \textit{i.e.} $|x| > \beta$, shows that
\begin{equation}
\label{ineein12} \forall x \in \amsmathbb{R} \setminus [-\beta,\beta]\qquad V^{\textnormal{eff}}(x) > v.
\end{equation}
The two inequalities \eqref{ineein1} and \eqref{ineein12} conclude the proof.
\end{proof}

\begin{center}
\begin{figure}[t]
\scalebox{0.5}{\includegraphics{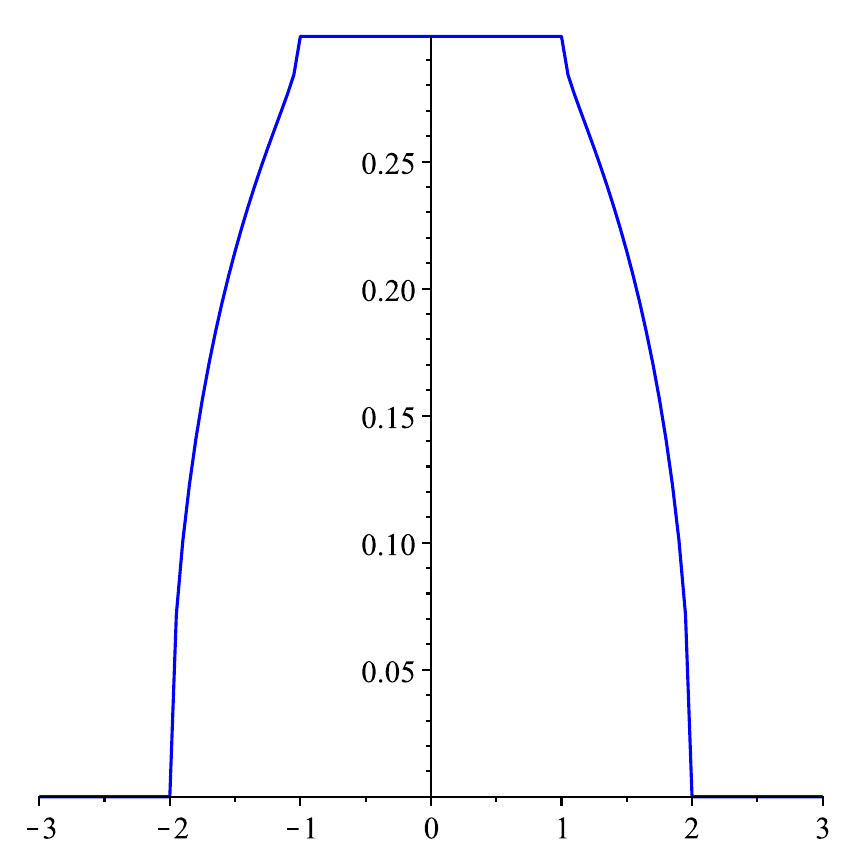}} \qquad
\scalebox{0.5}{\includegraphics{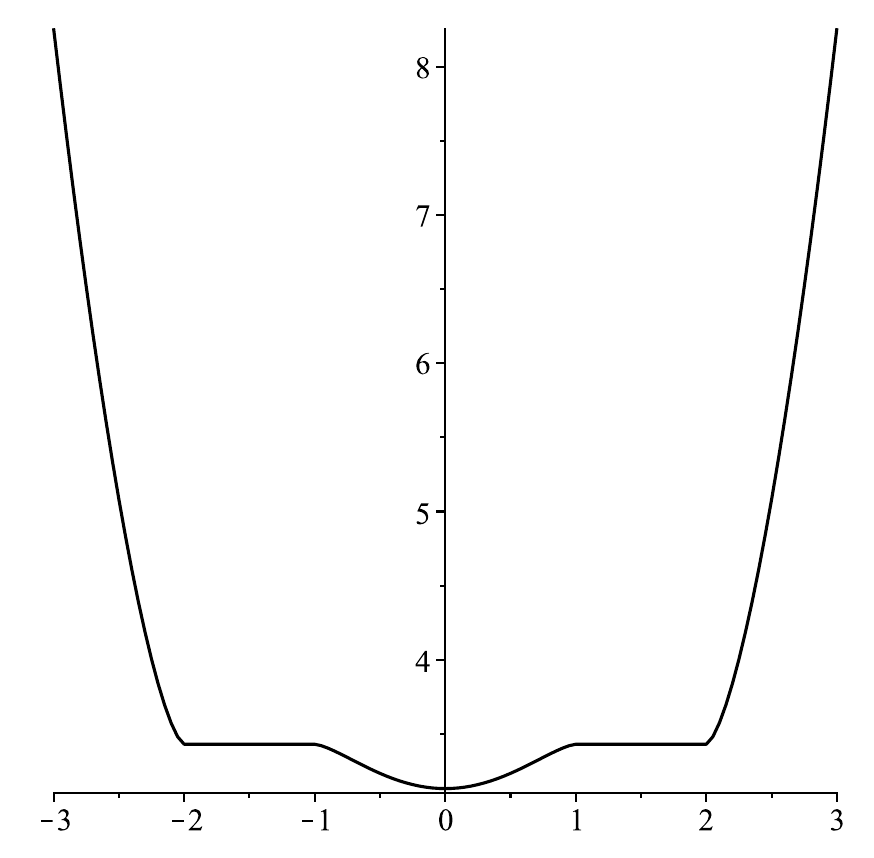}} \caption{\label{Fig_Gaussian} Density (left panel) plotted with formula \eqref{insidesat} and effective potential (right panel) for the equilibrium measure of the discrete ensemble with Gaussian weight with $\alpha = 1$ and $\beta = 2$, resulting in $\theta\approx 3.34$ and $c \approx 1.69$.}
\end{figure}
\end{center}
Proposition~\ref{Prop_Gaussian_LLN} indicates that
the equilibrium measure $\mu$ has two bands of density between $0$ and $\frac{1}{\theta}$ separated
by the saturation where $\mu(x)=\frac{1}{\theta}$. The examples in Sections \ref{Section_tiling_hole} and
\ref{Section_tiling_C} also have several bands, but there are caused by the geometry of
the models and the number of intervals of the support of the measure remains equal to the number of
bands. On the other hand, we see that for the discrete ensemble with Gaussian weight, the measure supported on a single interval develops two bands.

\section{Rough concentration estimates}
\subsection{Main results}
\label{section_rought_ann}
We now come back to the discrete ensembles of Section~\ref{Section_general_model}, which by Proposition~\ref{lem_finite_partition_function} are well-defined for $\N$ large enough compared to the constants in the assumptions. \emph{In all statements in the remaining of this book, we will constantly --- and often, silently --- assume that $\N$ is large enough in this sense.}

The parameters of a discrete ensemble determine an energy functional $-\mathcal{I}[\cdot]$ by \eqref{eq_functional_general_repeat} and thus a unique minimizer by Proposition~\ref{Lemma_maximizer}, which we have called \emph{equilibrium measure} and denoted $\boldsymbol{\mu}$. It is a deterministic measure \emph{which is a priori $\N$-dependent}, as all the parameters of the discrete ensemble can be. Besides, the empirical measure $\boldsymbol{\mu}_{\N} = (\mu_{\N,h})_{h = 1}^{H}$ is a random measure encoding configurations $\boldsymbol{\ell} \in \amsmathbb{W}_{\N}$ by
\[
\mu_{\N,h} := \frac{1}{\N} \sum_{i = 1}^{N_h} \delta_{\N^{-1}\ell_i^h},\qquad \mu = \sum_{h = 1}^{H} \mu_{\N,h}.
\]
The purpose of this section is to announce estimates of concentration of measures for the random empirical measure $\boldsymbol{\mu}_{\N}$ around the deterministic equilibrium measure $\boldsymbol{\mu}$ as $\N$ is large enough. For clarity of the exposition, most proofs will be deferred to the next Section~\ref{Section_rough_proofs}. Along the way, we derive the leading-order asymptotics for the partition function.

\begin{proposition}
\label{leadingZN}Under Assumptions~\ref{Assumptions_Theta} and \ref{Assumptions_basic}, there exists a constant $C > 0$ depending only on the constants in the assumptions, such that
\begin{equation}
\label{ZNestimateleading} |\log \Z_{\N} - \N^2 \mathcal{I}[\boldsymbol{\mu}]\big| \leq C \N \log \N.
\end{equation}
\end{proposition}

The key result is the following concentration of measures for $\boldsymbol{\mu}_\N - \boldsymbol{\mu}$.

\begin{definition}
\label{deffdeminorm} For a function $f : \amsmathbb{R} \rightarrow \amsmathbb{R}$, the Lipschitz norm is
\[
|\!| f |\!|_{\textnormal{Lip}} = \sup_{x \neq y} \bigg|\frac{f(x) - f(y)}{x - y}\bigg|.
\]
For $f \in \mathscr{L}^1(\amsmathbb{R})$ with Fourier transform $\widehat f$ --- \textit{cf.} our convention in Section~\ref{Notasec} --- we will make use of the norm
\[
|\!|f|\!|_{\frac{1}{2}}=\bigg(\int_{\amsmathbb{R}} |s| |\widehat f(s)|^2 \dd s\bigg)^{\frac{1}{2}}.
\]
We also introduce the space $\mathscr{H}_{\textnormal{Lip},\frac{1}{2}}$ of bounded real-valued functions on $\amsmathbb{R}$ with finite norms
 $|\!|f|\!|_{\textnormal{Lip}}$, and $|\!|f|\!|_{\frac{1}{2}}$.
\end{definition}

\begin{lemma} \label{Lemma_tail_bound_general}
 Suppose that Assumptions~\ref{Assumptions_Theta} and \ref{Assumptions_basic} hold. Then there exists
$C>0$ depending only on the constants in the assumptions, such that for any $t >0$ and $\N$ large enough
\begin{equation*}
\begin{split}
& \quad \P\Bigg[\exists f\in \mathscr{H}_{\textnormal{Lip},\frac{1}{2}}\quad \Big|\quad \min_{h \in [H]} \bigg|\int_{\amsmathbb R} f(x) \dd(\mu_{\N,h} - \mu_h)(x)\bigg| \geq t |\!| f |\!|_{\frac{1}{2}}+ \frac{C}{\N}\big( |\!|f|\!|_\textnormal{Lip} +|\!|f|\!|_\infty \big)\Bigg] \\
& \leq \exp\bigg(C \N\log \N- \frac{t^2\N^2}{C}\bigg).
\end{split}
\end{equation*}
\end{lemma}

Lemma~\ref{Lemma_tail_bound_general} readily implies the law of large numbers announced in Theorem~\ref{Theorem_main_LLN} and the following concentration of measures for the number of particles in macroscopic segments.

\begin{corollary}\label{Corollary_a_priori_0} Suppose that Assumptions~\ref{Assumptions_Theta} and \ref{Assumptions_basic}
hold. Then there exists
$C> 0$ depending only on the constants in the assumptions such that for any $\mathfrak{a},\mathfrak{b} \in \amsmathbb{R} \cup \{\pm \infty\}$ with $\mathfrak{a} < \mathfrak{b}$, we have
\[
\P\bigg[\Big|\#\big\{i \in [N] \,\,\big|\,\, \ell_i \in [\N \mathfrak{a},\N \mathfrak{b}]\big\} - \N \mu([\mathfrak{a},\mathfrak{b}])\Big| \geq C\N^{\frac{1}{2}} (\log\N)^2\bigg] \leq \exp\bigg(- \frac{\N(\log\N)^2}{C}\bigg).
\]
\end{corollary}
\begin{remark}
\label{rkcompactsupport}As the support of $\mu$ is contained in a compact that only depends on the constants in our assumptions (\textit{cf.} Theorem~\ref{Theorem_equi_charact_repeat_2}), Corollary~\ref{Corollary_a_priori_0} implies that particles stay in this compact with overwhelming probability.
\end{remark}

In the next chapters, Lemma~\ref{Lemma_tail_bound_general} will be mostly used to control moments and cumulants of the empirical measure expressed through the empirical resolvent, the resolvent of the equilibrium measure, and their difference, as in Definition~\ref{def_empirical_mes}
\[
G_{h}(z) = \sum_{i=1}^{N_h} \frac{1}{z-\frac{\ell_i^h}{\N}},\qquad
\Gm_{\mu_h}(z) =\int_{\hat{a}_h'}^{\hat{b}_h'} \frac{\dd\mu_h(x)}{z-x},\qquad \Delta G_{h}(z) =
G_{h}(z) - \N\mathcal{G}_{\mu_h}(z).
\]

\begin{corollary} \label{Corollary_a_priory_1} Suppose that Assumptions~\ref{Assumptions_Theta} and \ref{Assumptions_basic}
hold. Then for any $\delta>0$, there exists a constant $C >0$ depending only on $\delta$ and the constants in the assumptions, such that
for any $t > 0$ and any $z$ at distance greater than $\delta$ from $\amsmathbb{A}$ we have
\begin{equation}
\label{eq_concentration_Stieltjes}
\amsmathbb{P}_{\N}\Bigg[\bigg|\frac{\Delta G_h(z)}{\N}\bigg| \geq t\Bigg] \leq \exp\bigg(C \N\log \N -
\frac{t^2\N^2}{C}\bigg).
\end{equation}
In addition for any $\eps>0$, any positive integer $n$, any $h_1,\ldots,h_n \in [H]$ and $z_1,\ldots,z_n \in\amsmathbb C$ such that $\min_{i} \textnormal{dist}(z_i,\amsmathbb{A}) \geq \delta$, we have
 \begin{equation}
 \label{eq_weak_estimate} \E\Bigg[\bigg|
\prod_{i=1}^n \Delta G_{h_i}(z_i)\bigg|\Bigg] \leq \N^{n(\frac{1}{2}+\eps)} + n \bigg(\frac{\N}{\delta}\bigg)^n \exp\bigg(C \N\log \N -
\frac{\N^{1+2\eps}}{C}\bigg).
 \end{equation}
\end{corollary}
If $n$, $\delta$, and $\eps$ are fixed, then for large enough $\N$ the second term in \eqref{eq_weak_estimate} is uniformly smaller than the first term.
\begin{proof}[Proof of Corollary~\ref{Corollary_a_priory_1}] If $t<\frac{\log \N}{\N}$, then \eqref{eq_concentration_Stieltjes} is trivial, since its right-hand side is greater than $1$ for $C\geq 1$ and $\N\geq 2$. In order to deal with $t \geq \frac{\log \N}{\N}$, we take a smooth real-valued function $\tilde{f}$ which vanishes outside
$\bigcup_{h=1}^H [\hat a_h -\frac{\delta}{2}, \hat b_h + \frac{\delta}{2}]$ and is equal to $1$ on
$\amsmathbb{A}$. For
\[
f(x) := \frac{\tilde{f}(x)}{z - x},
\]
the norms $|\!|f|\!|_{\infty}$, $|\!|f|\!|_\textnormal{Lip}$, and $|\!|f|\!|_{\frac{1}{2}}$ are all bounded from above by a constant $C_1>0$. Hence, using $t \geq \frac{\log \N}{\N}$ we have for $\N$ large enough
\[
 2 C_1 t \geq t |\!| f |\!|_{\frac{1}{2}}+ \frac{C}{\N} \big( |\!|f|\!|_\textnormal{Lip} +|\!|f|\!|_\infty\big).
\]
Therefore, the event $\big\{|\frac{1}{\N}\Delta G_h(z)| \geq 2 C_1 t\big\}$ is a subset of the event in Corollary~\ref{Corollary_a_priori_0} and by a rescaling of $t$, \eqref{eq_concentration_Stieltjes} is a direct corollary of Lemma
\ref{Lemma_tail_bound_general}.

To justify the bound \eqref{eq_weak_estimate}, we write
\begin{equation}
\label{eq_x204}
\begin{split}
\E\Bigg[\bigg| \prod_{i=1}^n \Delta G_{h_i}(z_i)\bigg|\Bigg] & \leq \E\Bigg[\bigg|\prod_{i=1}^n \Delta G_{h_i}(z_i)\bigg| \cdot \prod_{j = 1}^{n} \mathbbm{1}_{|\Delta G_{h_j}(z_j)|< \N^{\frac{1}{2}+\eps}}\Bigg] \\
& \quad + \sum_{j=1}^n \E\Bigg[\bigg|\prod_{i=1}^n \Delta G_{h_i}(z_i)\bigg|\cdot \mathbbm{1}_{|\Delta G_{h_j}(z_j)|\geq \N^{\frac{1}{2}+\eps}}\Bigg].
\end{split}
\end{equation}
The first term in \eqref{eq_x204} is bounded from above by the first term in \eqref{eq_weak_estimate}. Using \eqref{eq_concentration_Stieltjes} with $t=\N^{-\frac{1}{2}+\eps}$ and the inequality
\[
\forall i \in [n]\qquad \big|\Delta G_{h_i}(z_i)\big|\leq \frac{\N}{\delta},
\]
which is implied by the assumptions on $z_i$, the second term in \eqref{eq_x204} is bounded from above by the second term in \eqref{eq_weak_estimate}.
\end{proof}

\begin{corollary}
\label{Corollary_Stieltjes_apriory} Suppose that Assumptions~\ref{Assumptions_Theta} and \ref{Assumptions_basic} hold. Then for any $\delta>0$, there exists a constant $C>0$, which depends only on $\delta$ and on the constants in the assumptions, such that any $\eps>0$, any $n \in \amsmathbb{Z}_{\geq 3}$, any $z_1,\ldots,z_n \in \amsmathbb{C}$ at distance at least $\delta$ from $\amsmathbb{A}$, and $\N$ large enough, we have
\[
\E\Bigg[\bigg|\prod_{i = 1}^{n} \Delta G_{h_i}(z_i) \bigg|\Bigg] \leq \N^{ (n - 2)(\frac{1}{2} + \eps)} \cdot \max_{i \in [n]} \E\big[ |\Delta G_{h_i}(z_i)|^2\big] + n \bigg(\frac{\N}{\delta}\bigg)^n\exp\bigg(C\N\log \N -
\frac{\N^{1+2\eps}}{C}\bigg).
\]
\end{corollary}
\begin{proof} We use \eqref{eq_x204} and bound from above the second term in the same way as in the previous proof. For the first term, we bound from above each factor except for two (say, the first two) with $\N^{\frac{1}{2} + \eps}$. For the remaining two factors, we get an upper bound using
\[
\big|\Delta G_{h_1}(z_1) \Delta G_{h_1}(z_1)\big| \leq \frac{\big|\Delta G_{h_1}(z_1)\big|^2+ \big|\Delta G_{h_2}(z_2)\big|^2}{2}.\qedhere
\]
\end{proof}

\subsection{Proofs}
 \label{Section_rough_proofs}
\label{rough_conc}

The goal of this section is to prove Proposition~\ref{leadingZN} and Lemma~\ref{Lemma_tail_bound_general}.

\begin{definition} \label{Definition_pseudodistance} Let $\nu^1,\nu^2$ be two nonnegative measures that are absolutely continuous of bounded density on $\amsmathbb{A}$, and as before, $\boldsymbol{\nu}^{1},\boldsymbol{\nu}^{2}$ the corresponding $H$-tuple of measures. We introduce
\begin{equation}
 \label{Ddefnu1nu2} \D\big[\boldsymbol{\nu}^1,\boldsymbol{\nu}^2\big] := -\mathcal{I}^{(2)}\big[\boldsymbol{\nu}^1-\boldsymbol{\nu}^2\big]= - \sum_{g,h = 1}^H \theta_{g,h} \iint_{\amsmathbb{R}^2} \log|x-y| (\nu_g^1-\nu_{g}^2)(x)(\nu_{h}^1 - \nu_{h}^2)(y)\,\dd x \dd y.
\end{equation}
\end{definition}

The basic properties of this quantity were described in Lemma~\ref{Lemma_I_quadratic_Fourier} and Corollary~\ref{Corollary_I_positive}. In particular, it is nonnegative\footnote{We will keep the notation $\mathfrak{D}^2$ even if the two measures do not satisfy the same constraints $(\star)$. In that case, contrarily to what the notation suggests we do not make statements about the sign of this quantity.} if $\nu^1,\nu^2$ satisfy both the same constraints $(\star)$, and will serve in this situation as a pseudo-distance to study the proximity for large $\N$ of the random empirical measure $\boldsymbol{\mu}_{\N}$ to the deterministic equilibrium measure $\boldsymbol{\mu}$ studied in Section~\ref{Section_Energy_functional}. However, $\mu_{\N}$ has atoms which can lead to infinite values of the right-hand side of \eqref{Ddefnu1nu2}. Therefore, we shall first approximate the empirical measures by absolutely continuous measures of bounded density. As another technical ingredient, we are going to construct a configuration whose empirical measure is close to the equilibrium measure $\boldsymbol{\mu}$. This construction requires some care: the empirical measures of configurations $\boldsymbol{\ell} \in \W_\N$ have integral --- after multiplication by $\N$ --- masses on the segments of the ensemble, while $\mu([\hat{a}_h',\hat{b}_h'])$ need not have any integrality property. Hence, we need to show that small perturbations of the filling fractions are always possible to match integrality conditions. This is the content of the following lemma.

We recall from Section~\ref{DataS} that filling fractions in the discrete ensembles satisfy affine conditions
\[
(\star) \qquad \forall e \in [\mathfrak{e}]\qquad \mathfrak{r}_e(\hat{\boldsymbol{n}}) = r_e,
\]
and the linear forms in the left-hand side have integral coefficients. Besides, we had introduced in \eqref{eq_Lambda_def} the set $\hat{\Lambda}_{\star}^{\amsmathbb{R}}$ of solutions $\hat{\boldsymbol{n}} \in \amsmathbb{R}_{\geq 0}^H$ of \eqref{eq_equations_eqs} and its subset $\hat{\Lambda}_{\star}$ of solutions that have integral coordinates after multiplication by $\N$.

\begin{lemma}
\label{Lemma_filling_modification} Under Assumption~\ref{Assumptions_Theta}, there exists $C>0$ depending only on constants in the assumption, such that for any $\boldsymbol{X} \in \amsmathbb{R}^H$ satisfying
\begin{equation*}
\begin{split}
\forall h \in [H] & \qquad 0 \leq X_h \leq \frac{\hat b'_h-\hat a'_h}{\theta_{h,h}}, \\
\forall e \in [\mathfrak{e}] & \qquad \big|\mathfrak{r}_e(\boldsymbol{X}) - r_e\big| \leq \eps,
\end{split}
\end{equation*}
there exists an $H$-tuple of real numbers $\tilde{\boldsymbol{X}} \in \hat{\Lambda}_{\star}^{\amsmathbb{R}}$ such that
\begin{equation}
\label{iugbnrugbr}\forall h \in [H] \qquad |\tilde{X}_h - X_h| < C \eps \quad \textnormal{and}\quad 0 \leq \tilde{X}_h \leq \frac{\hat b'_h-\hat a'_h}{\theta_{h,h}}.
\end{equation}
If moreover, $\hat{\Lambda}_\star \neq \emptyset$ and $\eps \geq \frac{1}{\N}$, then we can even choose $\tilde{\boldsymbol{X}} \in \hat{\Lambda}_{\star}$.
\end{lemma}
In the proof we need an auxiliary linear algebra statement.
\begin{lemma} \label{Lemma_non_empty_interior}
 Under Assumption~\ref{Assumptions_Theta}, there exists $C>0$, which depends only on the constants in the assumptions, such that $\hat{\Lambda}_{\star}^{\amsmathbb{R}}$ necessary has at least one element $\boldsymbol{X}^{\star}$ satisfying
\begin{equation}
 \forall h \in [H] \qquad \frac{1}{C} \leq X_h^{\star} \leq \frac{\hat b'_h-\hat a'_h}{\theta_{h,h}}-\frac{1}{C}.
\end{equation}
\end{lemma}
\begin{proof}[\textsc{Proof of Lemma~\ref{Lemma_non_empty_interior}}] We first assume that $\hat{a}_1$ and $\hat{b}_h$ are finite. Let $\Lambda^{\Box} \subset \hat{\Lambda}_{\star}^{\amsmathbb{R}}$ denote the set of points $\boldsymbol{X} = (X_h)_{h = 1}^H$ inside $\prod_{h = 1}^H \big[0,\frac{1}{\theta_{h,h}}(\hat{b}_h - \hat{a}_h)\big]$ and satisfy the equations \eqref{eq_equations_eqs}. This set is convex and by the fifth condition of Assumption~\ref{Assumptions_Theta} it is non-empty. By convexity and increasing $C$, if needed, it is sufficient for each $h\in[H]$ to find a two points $\boldsymbol{X},\boldsymbol{X}' \in\Lambda^\Box$ with
\[
X_h\geq \frac{1}{C}\qquad \text{and}\qquad X'_h\leq \frac{\hat b'_h-\hat a'_h}{\theta_{h,h}}-\frac{1}{C}.
\]
We show existence of such a $\boldsymbol{X}$ for $h=1$. The same argument will work for other $h$ up to index permutation, and for $\boldsymbol{X}'$ by reflection symmetry. By a compactness argument, it is sufficient to find $\boldsymbol{X}\in\Lambda^\Box$ with $X_1>0$. We argue by contradiction and suppose that for all $\boldsymbol{X}\in\Lambda^\Box$ we have $X_1=0$. Some other coordinates may be equal to the extreme allowed values for all points in $\Lambda^\Box$. Without loss of generality we can assume that for all $\boldsymbol{X} \in \Lambda^\Box$ we have
\begin{equation}
\label{eq_x230}
 X_1 =X_2=\cdots=X_K=0, \qquad X_{K+1}=\frac{\hat b'_{K+1}-\hat a'_{K+1}}{\theta_{K+1,K+1}}, \ldots, X_{K+M}=\frac{\hat b'_{K+M}-\hat a'_{K+M}}{\theta_{K+M,K+M}},
\end{equation}
for some given $K$ and $M$, and that there exists $\boldsymbol{X}^0 \in \Lambda^\Box$, such that
\[
 \forall g > K + M,\qquad 0<X^0_g<\frac{\hat b'_g-\hat a'_g}{\theta_{g,g}}
\]
For this latter fact we again used convexity of $\Lambda^\Box$.

Let us plug the values \eqref{eq_x230} into the linear forms $\mathfrak{r}_e$, $e\in [\mathfrak{e}]$ and consider the result as an affine function of the remaining variables $X_{K+M+1}, \ldots, X_H$. Assumption~\ref{Assumptions_Theta}.5 implies that the resulting $\mathfrak{e}$ functions are linearly independent. Simultaneously, $\boldsymbol{X}^{0}$ solve the equations \eqref{eq_equations_eqs}. Therefore, if we take small $\delta>0$ and replace $X_1=0$ with $X_1=\delta$, then there exists $(X^{\delta}_{g})_{g = K + M + 1}^{H}$ close to $(X^{0}_g)_{g = K + M + 1}^{H}$, and which still lead to a solution of \eqref{eq_equations_eqs}. This contradicts $X_1=0$ for all $\boldsymbol{X}\in \Lambda^\Box$, and finishes the proof.
\end{proof}

\begin{proof}[Proof of Lemma~\ref{Lemma_filling_modification}] It is enough to consider only the case when $\eps$ is small enough, as
otherwise the statement becomes trivial by taking a constant $C$ large enough. Let us decompose the linear space $\amsmathbb R^H$, equipped with its standard Euclidean structure, into an orthogonal direct sum $\mathscr{R} \oplus \mathscr{R}^{\perp}$, where
$\mathscr{R}= \textnormal{Ker}(\mathfrak{r})$ and
\[
\mathfrak{r} = (\mathfrak{r}_1,\ldots,\mathfrak{r}_\mathfrak{e}) : \amsmathbb{R}^H \longrightarrow \amsmathbb{R}^{\mathfrak{e}}.
\]
Assumption~\ref{Assumptions_Theta} guarantees that
$\dim(\mathscr{R}^\perp)=\mathfrak e$ and that
\begin{equation}
\label{eq_x75}
\forall \boldsymbol{X} \in \mathscr{R}^{\perp}\qquad C^{-1} |\!|\boldsymbol{X}|\!|_{\infty} \leq |\!|\mathfrak{r}(\boldsymbol{X})|\!|_{\infty} \leq C |\!|\boldsymbol{X}|\!|_{\infty}.
\end{equation}
 Using \eqref{eq_x75}, we can
find a vector $\Delta\boldsymbol{X} \in \mathscr{R}^\perp$ with $|\!|\boldsymbol{X}|\!|_{1} \leq C\eps$ such
that $\boldsymbol{X}+\Delta\boldsymbol{X}$ satisfies the constraints \eqref{eq_equations_eqs}. Note that the second condition in \eqref{iugbnrugbr} might fail for $\boldsymbol{X}+\Delta\boldsymbol{X}$ --- say, if $X_1=0$ and $\Delta X_1<0$. However, since $|\!|\boldsymbol{X}|\!|_{1} \leq C\eps$, we have
\[
-C\eps \leq X_h + \Delta X_h \leq \frac{\hat{b}'_h - \hat{a}_h}{\theta_{h,h}} + C\eps.
\]
Hence, using $\boldsymbol{X}^{\star}$ from Lemma~\ref{Lemma_non_empty_interior} and assuming $\eps$ to be small enough, the direction from $\boldsymbol{X}+\Delta\boldsymbol{X}$ to $\boldsymbol{X}^{\star}$ points towards the interior of parallelepiped $\prod_{h = 1}^{H} \big[0,\frac{1}{\theta_{h,h}}(\hat{b}_h' - \hat{a}_h')\big]$. Therefore, there exists $t \in (0,C' \eps)$ such that
\[
\tilde{\boldsymbol{X}}:=\boldsymbol{X}+\Delta\boldsymbol{X}+ t(\boldsymbol{X}^{\star}- \boldsymbol{X}-\Delta\boldsymbol{X})
\]
satisfies the second condition of \eqref{iugbnrugbr}. Since both $\boldsymbol{X}+\Delta\boldsymbol{X}$ and $\boldsymbol{X}^{\star}$ satisfy the constraints \eqref{eq_equations_eqs}, so does $\tilde{\boldsymbol{X}}$ and we are done with the first part of the statement.

For the second part, since the linear forms $\mathfrak{r}_1,\ldots,\mathfrak{r}_{\mathfrak{e}}$ have integral coefficients from the description in Section~\ref{DataS}, we can choose a $\amsmathbb{R}$-basis of $\mathscr{R}$ such that
the coordinates of all the basis vectors are integral after multiplication by $\N$. Assuming that $\hat{\Lambda}_\star\neq\emptyset$, this implies the existence of a $\frac{C}{\N}$-net of solutions of \eqref{eq_equations_eqs} in $\hat{\Lambda}_{\star}$, that is with coordinates in $\frac{1}{\N}\amsmathbb Z_{\geq 0}$. By $\frac{C}{\N}$-net here we mean that any point of $\Lambda_{\star}^{\amsmathbb{R}}$ is at distance at most $\frac{C}{\N}$ from a point in the net. Replacing $\tilde{\boldsymbol{X}}$ by (one of) the closest point in the net finishes the proof.
\end{proof}

Let us now introduce a procedure for the replacement of the empirical measure $\mu_{\N}$ with an absolutely continuous measure of bounded density $\tilde{\mu}_{\N}$. As the empirical measure, the measure $\tilde{\mu}_{\N}$ depends on a configuration $\boldsymbol{\ell} \in \amsmathbb{W}_{\N}$ and thus can be seen as a random measure under $\P$. The exact details of the procedure are not of central importance to us provided that it fulfills the properties listed in Lemma~\ref{mumutilde_lemma}.

\begin{definition}
\label{Definition_regularisation_measure_by_conv}
For any $h \in [H]$, we denote $\mu^*_{\N,h}$ the convolution of $\mu_{\N,h}$ with the uniform probability measure on $\big[0,\frac{1}{\N\theta_{h,h}}\big]$. We denote $\boldsymbol{\mu}^*_{\N} = (\mu_{\N,h})_{h = 1}^{H}$, and $\mu^*_{\N} = \sum_{h = 1}^{H} \mu_{\N,h}^*$ the corresponding measure on $\amsmathbb{R}$.
\end{definition}

\begin{lemma}\label{mumutilde_lemma}
Under Assumption~\ref{Assumptions_Theta}, there exists $C > 0$ depending only on the constants in the assumptions, such that for any $\boldsymbol{\ell} \in \amsmathbb{W}_\N$, there exists an absolutely continuous nonnegative measure $\tilde{\mu}_{\N}$ with the following properties:
\begin{itemize}
\item $\widetilde{\mu}_{\N} \in \mathscr{P}_{\star}$.
\item $|\widetilde{\mu}_{\N} - \mu_{\N}^*|(\amsmathbb{R}) \leq \frac{C}{\N}$.
\item $\widetilde{\mu}_{\N}$ and $\mu_{\N}^*$ coincide outside $[-C,C]$.
\end{itemize}
\end{lemma}
\begin{proof} Let $\boldsymbol{\ell} \in \amsmathbb{W}_\N$ and $h \in [H]$. The density of $\mu^*_{\N,h}$ only takes the two values $0$ or $\frac{1}{\theta_{h,h}}$ and its support is a union of $N_h$ segments of length $\theta_{h,h}$ inside $[\hat{a}_h,\hat{b}_h + \theta_{h,h}]$. Some of these segments may extend outside $[\hat{a}_h',\hat{b}_h']$. When this happens we redistribute the mass $[\hat{a}_h',\hat{b}_h']$ to the closest available regions so that the resulting measure has density either $0$ or $\frac{1}{\theta_{h,h}}$. This may not always be possible. Indeed, the maximal mass of a nonnegative measure with density bounded by $\frac{1}{\theta_{h,h}}$ and support included in $[\hat{a}_h',\hat{b}_h']$ is
\begin{equation}
\label{eq:maxx}
\frac{\hat{b}_h' - \hat{a}_h'}{\theta_{h,h}} = \frac{\hat{b}_h - \hat{a}_h}{\theta_{h,h}} + \frac{1}{\N} + \frac{1}{\N}\left(1 - \frac{1}{\theta_{h,h}}\right)
\end{equation}
while the total mass of $\mu_{N,h}^*$ (which is also the total mass of $\mu_{\N,h}$) is $\frac{\hat{b}_h - \hat{a}_h}{\theta_{h,h}} + \frac{1}{\N}$. If $\theta_{h,h} < 1$, the latter could exceed \eqref{eq:maxx}. In this case, we move the excess mass to another segment $[\hat{a}'_g,\hat{b}'_g]$, and Lemma~\ref{Lemma_filling_modification} guarantees that it can be done globally so that the resulting measure $\widetilde{\mu}_{\N}$ remains in $\mathscr{P}_{\star}$.
\end{proof}

We are now in position to prove the concentration of measure estimates for $\tilde{\boldsymbol{\mu}}_\N$ in the pseudo-distance $\D$.

\begin{theorem} \label{Prop_pseudodistance_bound} Under Assumptions~\ref{Assumptions_Theta} and \ref{Assumptions_basic}, there exists a constant $C$ depending only on the constants in the assumptions, such that for any $t > 0$ and $\N$ large enough we have
\begin{equation}
\label{eq_x205}
 \P\big[\Dr[\tilde{\boldsymbol{\mu}}_\N,\boldsymbol{\mu}]\geq t\big] \leq \exp(C\N\log\N - \N^2t^2).
\end{equation}
\end{theorem}
\begin{proof}
 The main point in the proof is to show that for any $\boldsymbol{\ell} \in \amsmathbb{W}_{\N}$ we have
\begin{equation}\label{bo1}
\P(\boldsymbol{\ell})\leq \exp\Big(\N^2\big(\mathcal{I}[\tilde{\boldsymbol{\mu}}_\N]-\mathcal{I}[\boldsymbol{\mu}]\big) + O(\N\log \N)\Big),
\end{equation}
with the measure $\tilde{\boldsymbol{\mu}}_{\N}$ depending on $\boldsymbol{\ell}$ constructed in Lemma~\ref{mumutilde_lemma} and an error independent of $\boldsymbol{\ell}$.

\medskip

\noindent \textsc{Step 1.} We first obtain a lower bound on the partition function $\Z_\N$ in
\eqref{eq_general_measure} by keeping only one particular configuration $\boldsymbol{p} \in \W_\N$.
For that take $\eps= \frac{1}{\N}$ and set $\boldsymbol{X}$ equal to the vector of filling fractions of the equilibrium measure
\[
\boldsymbol{X} = \big(\mu([\hat{a}_h',\hat{b}_h'])\big)_{h = 1}^{H}
\]
We can apply Lemma~\ref{Lemma_filling_modification} and --- given that a prerequisite to talk about a discrete ensemble is that $\hat{\Lambda}_{\star}$ is nonempty --- find a nearby $\tilde{\boldsymbol{X}}$ satisfying the integrality conditions. By moving small masses between the connected components of $\amsmathbb{A}$, we can construct an absolutely continuous nonnegative measure $\tilde \mu$ with support included in a $\N$-independent compact and density on $[\hat{a}_h',\hat{b}_h']$ bounded by $\frac{1}{\theta_{h,h}}$ for any $h \in [H]$, and satisfying
\[
\tilde{\boldsymbol{X}} = \big(\tilde{\mu}([\hat{a}_h',\hat{b}_h'])\big)_{h = 1}^{H}\qquad \textnormal{and}\qquad |\mu - \tilde{\mu}|(\amsmathbb{A}) = O\bigg(\frac{1}{\N}\bigg).
\]
The key difference between\footnote{These measures should be not confused with the empirical measure $\mu_{\N}$ and its regularization $\tilde{\mu}_{\N}$ constructed in Lemma~\ref{mumutilde_lemma}. The latter will play a role in the second part of the proof.} $\mu$ and $\tilde \mu$ is that the segment filling fractions for the latter are integral after multiplication by $\N$.

Before describing $\boldsymbol{p} \in \W_\N$, we first define the quantiles $\boldsymbol{q} = (q_1,\ldots,q_N)$
of $\tilde \mu$ through
\[
\forall i \in [N]\qquad \tilde \mu((-\infty, q_i])=\frac{i-\frac{1}{2}}{\N}.
\]
For any $h \in [H]$, we define
\[
I_{\boldsymbol{q}}^-(h) = \min\big\{i \in [N]\quad\big|\quad q_i \in [\hat a'_h,
\hat b'_h]\big\}\qquad I_{\boldsymbol{q}}^+(h) = \max\big\{i \in [N]\quad \big|\quad q_i \in [\hat a'_h,
\hat b'_h]\big\}.
\]
For each $h \in [H]$, the configuration $\boldsymbol{p}$ we are going to construct will exactly have
\[
N_h(\boldsymbol{p}) := \N\tilde{\mu}([\hat{a}_h',\hat{b}_h']) = I_{\boldsymbol{q}}^+(h) - I_{\boldsymbol{q}}^-(h) + 1
\]
particles in the $h$-th segment. If $a_h$ is finite, we set
\[
\forall i \in \llbracket I^-_{\boldsymbol{q}}(h),I^+_{\boldsymbol{q}}(h)\rrbracket \qquad p_i :=a_h+ (i-I^-_{\boldsymbol{q}}(h)) \theta_{h,h} + \big\lfloor \N (q_i -\hat a'_h) - (i-I^-_{\boldsymbol{q}}(h))\theta_{h,h} \big\rfloor.
\]
If $a_h=-\infty$, we must have $h = 1$. Then, if $b_1$ is finite we set
\[
\forall i \in \llbracket I^-_{\boldsymbol{q}}(1),I^+_{\boldsymbol{q}}(1)\rrbracket\qquad p_i := b_1+ (i-I^+_{\boldsymbol{q}}(1)) \theta_{1,1} + \big\lfloor \N(q_i - \hat b'_1) - (i-I^+_{\boldsymbol{q}}(1))\theta_{1,1}\big\rfloor,
\]
and if $b_1=+\infty$ we rather set
\[
\forall i \in \llbracket I^-_{\boldsymbol{q}}(1),I^+_{\boldsymbol{q}}(1)\rrbracket\qquad p_i :=(i-1)\theta_{1,1}+ \big\lfloor \N q_i-(i-1)\theta_{1,1}\big\rfloor.
\]
Since the density of $\tilde\mu_h$ is bounded from above by $\frac{1}{\theta_{h,h}}$, we have
\begin{equation}
\label{lowerqii}
q_{i + 1} - q_i \geq \frac{\theta_{h,h}}{\N}
\end{equation}
whenever both $i$ and $i+1$ belong to $\llbracket I^-_{\boldsymbol{q}}(h), I^+_{\boldsymbol{q}}(h)\rrbracket$, and therefore,
$p_{i+1}-p_i\in \theta_{h,h} + \amsmathbb{Z}_{\geq 0}$. We conclude that $\boldsymbol{p} \in
\W_\N$. Moreover
\begin{equation}
\label{eq_x17} \forall i \in [N]\qquad \bigg|\frac{p_i}{\N}- q_i\bigg| \leq \frac{1}{\N}+ \frac{\big|\theta_{h(i),h(i)}-\frac{1}{2}\big|}{\N}.
\end{equation}
In the right-hand side, the first term is an error coming from the integral parts in the definition of $\boldsymbol{p}$ and the second term comes from the difference between $\hat a_h$ or $\hat b_h$ and their shifted versions $\hat{a}'_h$ or $\hat{b}'_h$ --- see their definition in \eqref{eq_shifted_parameters}.

As we did in \eqref{eljah}, we record with the notation $p_i^h$ for $i \in [N_h(\boldsymbol{p})]$ and $h \in [H]$ the particles among $p_1,\ldots,p_N$ that belong to the $h$-th segment, in increasing order. We do the same with $q_i^h$ referring to quantiles belonging to the rescaled and shifted segment $[\hat{a}_h',\hat{b}_h']$. Then, \eqref{eq_x17} can be rewritten for any
\begin{equation}
\label{eq_x17bis}
\forall h \in [H]\qquad \forall i \in [N_h(\boldsymbol{p})]\qquad \bigg|\frac{p_i^h}{\N} - q_i^h\bigg| \leq \frac{1}{\N} + \frac{\big|\theta_{h,h}-\frac{1}{2}\big|}{\N}.
\end{equation}

Using the lower bound \eqref{eq_Stirling_xy} for the factors of the product over pairs in the unnormalized measure $\widetilde{\P}$, \textit{cf.} \eqref{eljah}, we get
\begin{equation}
\label{eq_ZLowerbound1}
\begin{split}
\Z_\N & \geq \widetilde{\P}(\boldsymbol{p}) \geq \exp\Bigg[-C \N\log\N + \!\! \sum_{1 \leq g < h \leq H}
\sum_{\substack{1 \leq i \leq N_g(\boldsymbol{p}) \\ 1 \leq j \leq N_{h}(\boldsymbol{p})}}
2\theta_{g,h}\log\bigg(\frac{p_{j}^{h} - p_i^g}{\N}\bigg) + \!\!
 \sum_{\substack{1 \leq h \leq H \\ 1 \leq i \leq N_h(\boldsymbol{p})}} \! \log w_h(p_i^h) \\
 & \qquad\qquad \qquad \quad\quad + \sum_{\substack{1 \leq h \leq H \\ 1 \leq i < j \leq N_h(\boldsymbol{p})}}
2\theta_{h,h}\log\bigg(\frac{p_{j}^h - p_i^h}{\N}\bigg) \Bigg].
\end{split}
\end{equation}
We now analyze each term in \eqref{eq_ZLowerbound1} separately. For $\log w_h$ we
use Assumption~\ref{Assumptions_basic}. The error $\err_h$ appearing in the weight $w_h$ contributes to the error $C \N
\log\N $, while for the potential part $V_h$ we write
\[
 \N \sum_{i=1}^{N_h(\boldsymbol{p})} V_h\bigg(\frac{p_i^h}{\N}\bigg)= \N^2 \sum_{i=2}^{N_h(\boldsymbol{p})-1}
 \int_{q_i^h}^{q_{i+1}^h} \Big(V_h(x)+ (x-\N^{-1}p_i^h)V_h'(\xi_x) \Big)\tilde\mu_{h}(x)\dd x + O(\N),
\]
where $\xi_x$ is a point between $x$ and $\frac{p_i^h}{\N}$. Here we could include the $i=1$ and $i=N_h(\boldsymbol{p})$ terms in the $O(\N)$-error because the potential \eqref{eq_V_logs} remains bounded near $\hat a'_h$ and $\hat b'_h$. The remaining $p_i^h$ for $i \neq 1,N_h(\boldsymbol{p})$ are in $[\N \hat a'_h + \frac{1}{2}, \N\hat b'_h - \frac{1}{2}]$. Then by Assumption~\ref{Assumptions_basic} we have $|V_h'(\xi_x)|=O(\log \N)$ uniformly for $x$ in the integration range, and using \eqref{eq_x17bis} and the fact that the mass of $\tilde{\mu}$ between two consecutive quantiles is $\frac{1}{\N}$, we get
\begin{equation}
\label{eq_x16}
\begin{split}
\N \sum_{i=1}^{N_h(\boldsymbol{p})} V_{h}\bigg(\frac{p_i^h}{\N}\bigg) & = \N^2 \sum_{i=2}^{N_h(\boldsymbol{p})-2}
\int_{q_i^h}^{q_{i+1}^h} \Big(V_h(x) +\big(q_{i+1}^h-q_i^h+O(\N^{-1})\big) \cdot O(\log \N)\Big)\tilde{\mu}_{h}(x)\dd x \\
& = \N^2 \int_{\hat a'_h}^{\hat b'_h} V_h(x)\tilde\mu_{h}(x)\dd x + \N
\Bigg(\sum_{i=2}^{N_h(\boldsymbol{p})-1} \big(q_{i+1}^h-q_i^h+O(\N^{-1})\big)\Bigg) \cdot O(\log \N) \\
 & = \N^2 \int_{\hat a'_h}^{\hat b'_h} V_h(x)\tilde\mu_{h}(x)\dd x +O(\N\log\N).
 \end{split}
 \end{equation}
A computation similar to \eqref{eq_x16} also works for the sum of
$\log\big(\frac{p_{j}^{h} - p_i^g}{\N}\big)$ featured in \eqref{eq_ZLowerbound1}, and we conclude
\begin{equation}
\label{eq_x18}
\begin{split}
& \quad \sum_{1 \leq g < h \leq H} \sum_{\substack{1 \leq i \leq N_g \\ 1 \leq j \leq N_{h}}} 2\theta_{g,h}
\log\bigg(\frac{p_{j}^{h} - p_i^g}{\N}\bigg) - \sum_{h = 1}^H \sum_{i = 1}^{N_h(\boldsymbol{p})} \log w_h(p_i^h) \\
& = \N^2\bigg(\iint_{\amsmathbb{R}^2} \sum_{1 \leq g \neq h \leq H} \theta_{g,h}\log|x -
y|\,\tilde\mu_{g}(x)\tilde\mu_{h}(y)\,\dd x \dd y -\sum_{h=1}^H \int_{\amsmathbb{R}}
V_h(x)\,\tilde\mu_{h}(x)\dd x \bigg) + O(\N\log \N).
\end{split}
\end{equation}
The terms corresponding to pairs of particles belonging to the same
$[\hat{a}'_{h},\hat{b}'_{h}]$ in the second line of \eqref{eq_ZLowerbound1}
should be treated separately because of the logarithmic singularity. Choose $d > 5+\frac{4}{\theta_{h,h}}$. Using the
monotonicity of the logarithm, \eqref{lowerqii} and \eqref{eq_x17bis}, we get
 \begin{equation}
 \label{eq_x19}
 \begin{split}
 & \quad \sum_{1\leq i<j \leq N_h(\boldsymbol{p})} \log\bigg(\frac{p_j^h - p_i^h}{\N}\bigg) =
 \sum_{\substack{1\leq i,j \leq N_h(\boldsymbol{p}) \\
 i+d<j}} \log\bigg(\frac{p_j^h - p_i^h}{\N}\bigg) +O(\N\log\N) \\
 &\geq \sum_{\substack{1\leq i,j \leq N_h(\boldsymbol{p}) \\
 i+d<j}} \log\bigg(q_j^h-q_i^h-\frac{3+2 \theta_{h,h}}{\N}\bigg)
 +O(\N\log\N)\\
 &\geq \N^2 \sum_{\substack{1\leq i,j \leq N_h(\boldsymbol{p}) \\
 i+d<j}} \int_{q_{j-1}^h}^{q_{j}^h} \int_{q_{i}^h}^{q_{i+1}^h} \log\bigg(x - y -\frac{3+2 \theta_{h,h}}{\N}\bigg)\tilde\mu_{h}(x)\tilde\mu_{h}(y)\,\dd x \dd y +O(\N\log\N) \\
&\geq \N^2 \iint_{y<x} \log(x-y)\tilde\mu_{h}(x)\tilde\mu_{h}(y)\dd x \dd y +O(\N\log\N).
\end{split}
\end{equation}
The constant $\tfrac{3+2 \theta_{h,h}}{\N}$ appears in the second line as an upper bound for $2\big(\tfrac{1}{\N}+\tfrac{|\theta_{h,h}-\frac{1}{2}|}{\N}\big)$ which comes from using \eqref{eq_x17bis} twice. The constant $d$ was chosen so that the arguments of logarithms are positive in \eqref{eq_x19}.

Inserting \eqref{eq_x18} and \eqref{eq_x19} into
\eqref{eq_ZLowerbound1}, we find the lower bound
\begin{equation}
\label{eq_x20_1} \Z_\N \geq \exp\big(\N^2\mathcal{I}[\tilde{\boldsymbol{\mu}}] + O(\N\log
\N)\big).
\end{equation}
By construction, $|\mu-\tilde\mu|$ is a measure of total mass $O(\frac{1}{\N})$, both $\mu$ and $\tilde\mu$ have a uniformly bounded density and their support is included in a compact independent of $\N$. Therefore, we can find a constant $C>0$ such that
\[
\big|{\mathcal{I}}(\boldsymbol{\mu})-{\mathcal{I}}(\tilde{\boldsymbol{\mu}})\big|\leq\frac{C}{\mathcal N}.
\]
and get
\begin{equation}
\label{eq_x20} \Z_\N \geq \exp\big(\N^2\mathcal{I}[\boldsymbol{\mu}] + O(\N\log \N)\big).
\end{equation}

\medskip

\noindent \textsc{Step 2.} We proceed to obtain an upper bound on $\P(\boldsymbol{\ell})$ by using \eqref{eq_Stirling_xy} to replace the
$\Gamma$-interaction by its Stirling approximation, and using the lower bound \eqref{eq_x20} for the partition function. We get
\begin{equation}
\label{eq_x22} \P(\boldsymbol{\ell}) \leq \exp\big(\N^2(\tilde {\mathcal I }[\boldsymbol{\mu}_{\N}] -
\mathcal{I}[\boldsymbol{\mu}]) + O(\N\log \N)\big)
\end{equation}
Here appears the regularized energy of the empirical measure introduced in \eqref{eq_x207}
\begin{equation}
\label{eq_x206}
 \tilde {\mathcal I}[\boldsymbol{\mu}_N]= 2 \sum_{1\leq i<j \leq N} \theta_{h(i),h(j)} \log\bigg( \frac{\ell_j - \ell_i}{\N}\bigg)- \sum_{i=1}^N V_{h(i)}\bigg(\frac{\ell_i}{\N}\bigg).
\end{equation}
We next estimate the cost of replacing the atomic measure $\mu_{\N}$ by its regularization constructed in Lemma~\ref{mumutilde_lemma}, \textit{i.e.} we aim to bound $\tilde{\mathcal I}[\boldsymbol{\mu}_\N]-\mathcal I[\tilde{\boldsymbol{\mu}}_\N]$. Take $h \in [H]$ and an integer $d>1+\frac{1}{\theta_{h,h}}$. Using Assumption
\ref{Assumptions_basic} and the fact that $|\tilde \mu_{\N,h}-\tilde \mu^*_{\N,h}|$ has total mass $O\big(\frac{1}{\N}\big)$, we have
\begin{equation}
\label{eq_x23}
\begin{split}
& \quad \int_{\hat a_h}^{\hat b_h} V_h(x) \dd\mu_{\N,h}(x) - \int_{\hat a'_h}^{\hat b'_h} V_h(x) \dd\tilde \mu_{\N,h}(x)  \\
& = \frac{1}{\N}\sum_{i=1+d}^{N_h-d} \Bigg[ V_h\bigg(\frac{\ell_i^h}{\N}\bigg)- \frac{\N}{\theta_{h,h}} \int_{\frac{\ell_i^h}{\N}}^{\frac{\ell_i^h+\theta_{h,h}}{\N}} \!\!\!V_h(x) \dd x \Bigg] + O\bigg(\frac{1}{\N}\bigg) \\
& \leq \frac{1}{\N^2}\sum_{i=1+d}^{N_h-d}\, \sup\big\{|V'_h(x)| \quad \big| \quad \N x \in [\ell_i^h,\ell_i^h+\theta_{h,h}]\big\} + O\bigg(\frac{1}{\N}\bigg) \\
& \leq \frac{C\,\log\N}{\N}.
\end{split}
\end{equation}
To get the second line we first extended $V_h$ by zero outside $[\hat{a}'_h,\hat{b}'_h]$, then replaced $\widetilde{\mu}_{\N,h}$ with $\mu_{\N,h}^*$ using $|\tilde \mu_{\N,h}-\tilde \mu^*_{\N,h}|$, and eventually moved the $d$ extreme terms on each side to the $O(\frac{1}{\N})$. The $N_h - 2d$ remaining terms involve values of $V_h$ away from the points $\hat{a}'_h$ and $\hat b'_h$ where $V'_h$ can explode logarithmically. Similarly, for $g\neq h$,
\begin{equation}
\label{eq_x24}
 \int_{\hat a_g}^{\hat b_g} \int_{\hat a_{h}}^{\hat b_{h}} \theta_{g,h}\log|x-y| \dd\mu_{\N,h}(x) \dd\mu_{\N,h}(y)
 - \int_{\hat a'_g}^{\hat b'_g} \int_{\hat a'_{h}}^{\hat b'_{h}} \theta_{g,h}\log|x-y|\tilde \mu_{\N,h}(x)\tilde{\mu}_{\N,h}(y)\,\dd x \dd y \leq \frac{C}{\N}.
\end{equation}
For the interaction of pairs in the segment $[\hat{a}_{h},\hat{b}_{h}]$, we have --- the first equality holds up to $O\big(\frac{\log \N}{\N}\big)$, which we omitted for lack of space but restored in the next lines ---
\begin{equation}
\label{eq_x25}
\begin{split}
& \quad  \iint_{x < y} \log(y - x)\dd\mu_{\N,h}(x)\dd\mu_{\N,h}(y) - \iint_{x < y} \log(y - x)\tilde{\mu}_{\N,h}(x)\tilde{\mu}_{\N,h}(y)\,\dd x \dd y \\
 = & \quad -\sum_{1\leq i < j \leq N_h} \frac{1}{ (\theta_{h,h})^2}
\int_{0}^{\frac{\theta_{h,h}}{\N}} \!\! \int_{0}^{\frac{\theta_{h,h}}{\N}} \log\bigg(1 + \frac{\N(u -
v)}{\ell_{j}^h - \ell_{i}^h}\bigg)\dd u\dd v -\sum_{i=1}^{N_h} \, \,
\iint\limits_{0<v<u<\frac{\theta_{h,h}}{\N}} \!\!\!\!\!\!
\log (u-v) \dd u\dd v \\
 \leq & \quad - \sum_{1\leq i < j \leq N_h}
\frac{1}{(\theta_{h,h})^2}\,\iint\limits_{0<v<u<\frac{\theta_{h,h}}{\N}} \log\Bigg[1-
\bigg(\frac{\N(v - u)}{\theta_{h,h}(j-i)}\bigg)^2 \Bigg]\dd u\dd v
-\sum_{i=1}^{N_h} \int_{0}^{\frac{\theta_{h,h}}{\N}}
\!\!\big(u\log (u)-u\big)\dd u  \\
 \leq & \quad 2 \sum_{1\leq i < j \leq N_h}
\frac{1}{(\theta_{h,h})^2}\,\iint\limits_{0<v<u< \frac{\theta_{h,h}}{\N}} \bigg(\frac{\N(v -
u)}{\theta_{h,h}(j-i)}\bigg)^2\dd u\dd v + O\bigg(\frac{\log\N}{\N}\bigg)=O\bigg(\frac{\log\N}{\N}\bigg).
\end{split}
\end{equation}
In the last inequality we have used $-\log(1-x^2)\leq 2 x^2$ for $|x|\leq \tfrac{1}{2}$ corresponding to $j-i\geq 2$ and the case $j-i=1$ was separately absorbed into $O\big(\frac{\log \N}{\N}\big)$. Therefore, combining \eqref{eq_x23}, \eqref{eq_x24}, \eqref{eq_x25} with \eqref{eq_x22} we get
\begin{equation}
\label{eq_x26} \P(\boldsymbol{\ell}) \leq \exp\big(\N^2(\mathcal{I}[\tilde{\boldsymbol{\mu}}_{\N}] -
\mathcal{I}[\boldsymbol{\mu}]) + O(\N\log \N)\big).
\end{equation}
We further write
\begin{equation}
\label{eq_x27}
\mathcal{I}[\tilde{\boldsymbol{\mu}}_{\N}] - \mathcal{I}[\boldsymbol{\mu}] =
-\mathfrak{D}^2[\tilde{\boldsymbol{\mu}}_{\N},\boldsymbol{\mu}] - \sum_{h=1}^H \int_{\hat a'_h}^{\hat b'_h} V^{{\textnormal{eff}}}_h(x)\dd(\tilde{\mu}_{\N,h} - \mu_h)(x).
\end{equation}
We claim that the second term in the right-hand side is nonpositive. Indeed, since both $\boldsymbol{\mu}$ and $\tilde{\boldsymbol{\mu}}_\N$ belong to $\mathscr{P}_\star$, so does $\boldsymbol{\mu}+\epsilon(\tilde{\boldsymbol{\mu}}_{\N} - \boldsymbol{\mu})$ for every $\epsilon \in (0,1)$. Since $\tilde\mu_\N -\mu$ has a uniformly bounded density and its support is contained in a $\N$-independent compact, and because $\boldsymbol{\mu}$ is a minimizer we have
\[
 0 \geq \mathcal I [\boldsymbol{\mu}+\epsilon(\tilde{\boldsymbol{\mu}}_{\N} - \boldsymbol{\mu})]- \mathcal I[\boldsymbol{\mu}]= -\epsilon
 \sum_{h=1}^H \int_{\hat a'_h}^{\hat b'_h} V^{{\textnormal{eff}}}_h(x)\dd(\tilde{\mu}_{\N,h} - \mu_h)(x) + o(\epsilon),
\]
which implies the claim by sending $\epsilon \rightarrow 0$.

 Hence, using \eqref{eq_x26}, \eqref{eq_x27} --- for clarity we restore explicitly the dependence of the measure $\tilde{\boldsymbol{\mu}}_\N$ on configurations $\boldsymbol{\ell}$ as $\tilde{\boldsymbol{\mu}}_\N^{(\boldsymbol{\ell})}$ in the notation --- we get for $t > 0$
\begin{equation}
\label{eq_x28}
\begin{split}
 &\P\big[ \mathfrak{D}^2[\tilde{\boldsymbol{\mu}}_\N,\boldsymbol{\mu}]\geq t^2\big] \\
 & \leq \sum_{\substack{\boldsymbol{\ell}\in \W_\N \\ \mathfrak{D}^2[\tilde{\boldsymbol{\mu}}_\N^{(\boldsymbol{\ell})},\boldsymbol{\mu}]\geq t^2}}
\exp\bigg(-\N^2 \mathfrak{D}^2[\tilde{\boldsymbol{\mu}}_\N^{(\boldsymbol{\ell})},\boldsymbol{\mu}] - \N^2\sum_{h=1}^H \int_{\hat a'_h}^{\hat b'_h} V^{{\textnormal{eff}}}_h(x)\dd(\tilde{\mu}_{\N,h}^{(\boldsymbol{\ell})} - \mu_h)(x)+ O(\N\log \N)\bigg)\\
& \leq \exp(C\N \log\N- t^2 \N^2) \cdot \sum_{\boldsymbol{\ell}\in\W_\N} \exp\bigg(-
\N^2\sum_{h=1}^H \int_{\hat a'_h}^{\hat b'_h} V^{{\textnormal{eff}}}_h(x)\dd(\tilde{\mu}_{\N,h}^{(\boldsymbol{\ell})} - \mu_h)(x)\bigg).
\end{split}
\end{equation}
We claim that the sum over $\boldsymbol{\ell}$ in this expression is bounded by $\exp(C'\N \log \N)$ for some $C' > 0$ depending only on the constants in the assumptions. Indeed, we have just explained that the expression under exponent is nonpositive. This implies that the sum over all possible configurations $\boldsymbol{\ell}$ of particles inside a compact $[-\N D, \N D]$ is bounded from above by the total number of such configurations, which for fixed $D$ is $\exp\big(O(\N\log \N)\big)$. For the configurations which contain particles at large positions $|\ell_i|\geq \N D$ for some $i \in [N]$, we rely on the inequality
\[
\forall x \in \amsmathbb{R}\qquad |x| > D \quad \Longrightarrow \quad V^{{\textnormal{eff}}}_h(x)\geq \frac{\eta_H}{2} \log |x|,
\]
which is implied by \eqref{eq_potential_eff_lower bound} for $D$ large enough. Besides, by construction in Lemma~\ref{mumutilde_lemma} the measures $\tilde{\mu}_{\N}$ and $\mu_\N^*$ coincide outside $[-D,D]$. The latter is the convolution of the empirical measure of $\boldsymbol{\ell}$ with the indicator of a small segment. Hence, a large $\ell_i$ yields a contribution to the integral of $V^{{\textnormal{eff}}}_h(x)$ which is positive and large, namely $O(\log \ell_i)$. As it is multiplied by $-\N^2$ and appears in an exponential, the contribution of these large $\ell_i$ is summable and the entire sum remains a $\exp\big(O(\N \log \N)\big)$.

We conclude that the desired inequality \eqref{eq_x205} is a consequence of \eqref{eq_x28}.
 \end{proof}

\begin{proof}[Proof of Proposition~\ref{leadingZN}] We want to derive the leading-order asymptotics of the partition function. Many ingredients already appeared in the previous proof. First, we have already obtained the lower bound $\Z_{\N} \geq \exp\big(\N^2 \mathcal{I}[\boldsymbol{\mu}] + C \N \log \N\big)$ in \eqref{eq_x20}. Second, we also have seen how the regularized energy of the empirical measure $-\tilde{\mathcal I}[\boldsymbol{\mu}_{\N}]$ of \eqref{eq_x206} can be replaced with the energy of the regularized empirical measure $-\mathcal I[\tilde{\boldsymbol{\mu}}_{\N}]$ up to an error $O(\N \log \N)$. Then, relying on the identity \eqref{eq_x27} and nonnegativity of $\mathfrak{D}^2$, we get
\begin{equation}
\begin{split}
\label{ZNfdddddgnfug}\Z_{\N} & = \sum_{\boldsymbol{\ell} \in \amsmathbb{W}_{\N}} \exp\big(\N^2 \tilde{\mathcal I}[\boldsymbol{\mu}_{\N}] + O(\N \log \N)\big) \\
& \leq \exp\big(\N^2 \mathcal{I}[\boldsymbol{\mu}] + O(\N \log \N)\big) \sum_{\boldsymbol{\ell} \in \amsmathbb{W}_{\N}} \exp\bigg(-\N^2\sum_{h = 1}^H \int_{\hat a'_{h}}^{\hat b'_{h}} V^{\textnormal{eff}}_{h}(x)\dd(\tilde{\mu}_{\N,h} - \mu_h)(x)\bigg)
\end{split}
\end{equation}
 Third, we have seen at the end of the previous proof that the last sum over $\W_\N$ is $\exp\big(O(\N\log \N)\big)$. This leads to the desired estimate \eqref{ZNestimateleading}.
\end{proof}

\begin{lemma}
\label{Lemma_linear_through_distance}
 Let $\nu$ be a compactly-supported signed measure on $\amsmathbb{A}$ with density $\mathscr{L}^2(\amsmathbb A)\cap \mathscr{L}^1(\amsmathbb A)$, and such that
\[
\forall e \in [\mathfrak{e}]\qquad \mathfrak r_e\big((\nu([\hat a'_h,\hat b'_h]))_{h = 1}^H\big) = 0.
\]
Then there exists $C>0$ depending only on the constants of Assumption
 \ref{Assumptions_Theta}, such that for any function $f$ with finite $|\!|f|\!|_{\frac{1}{2}}$ we have
 \[
 \biggl|\int_{\amsmathbb R} f(x) \nu(x)\dd x \biggr|^2 \leq C\,\D[\boldsymbol{\nu},\boldsymbol{0}]\, |\!|f|\!|^2_{\frac{1}{2}}.
 \]

\end{lemma}
\begin{proof}
\label{ehhindex}Choose a basis $\boldsymbol{u}^{(1)},\ldots,\boldsymbol{u}^{(H)}$ of eigenvectors of $\boldsymbol{\Theta}$ of Euclidean norm $1$, with respective eigenvalues $\lambda^{(1)},\ldots,\lambda^{(H)} \geq 0$. We fix $h \in [H]$ and assume without loss of generality that the function $f$ vanishes on $\bigcup_{g \neq h} [\hat a'_{g},\hat b'_{g}]$. We write
\begin{equation}
\label{eigenvalhhh}
\theta_{h,h}= \big\langle \boldsymbol{\Theta}(\boldsymbol{e}^{(h)}) \cdot \boldsymbol{e}^{(h)}\big\rangle=\sum_{g = 1}^H \lambda^{(g)} \big|\big\langle \boldsymbol{u}^{(g)} \cdot \boldsymbol{e}^{(h)}\big\rangle\big|^2 \leq H \max_{g \in [H]} \lambda^{(g)} |u^{(g)}_h|^2.
\end{equation}
where $\langle \cdot \rangle$ is the standard Euclidean scalar product and $\boldsymbol{e}^{(1)},\ldots,\boldsymbol{e}^{(H)}$ the canonical basis in $\amsmathbb{C}^H$. Let $\overline{\boldsymbol{u}}$ be an eigenvector realizing this maximum, and $\overline{\lambda}$ the corresponding eigenvalue. Due to Assumption~\ref{Assumptions_Theta}, we have
\begin{equation}
\label{eq_x74}
\overline{\lambda} |\overline{u}_h|^2 \geq \frac{1}{C} >0.
\end{equation}

Let $\nu_h$ be the restriction of $\nu$ to $[\hat{a}_h',\hat{b}'_h]$. We use the Fourier transform and the Plancherel theorem to compute
\begin{equation*}
\begin{split}
\int_{\amsmathbb R} f(x) \nu_h(x)\dd x & = \frac{1}{\overline{u}_h}\int_{\amsmathbb R} f(x) \left(\sum_{g=1}^H \overline{u}_g \nu_g(x)\right)\dd x \\
& = \frac{1}{2\pi \overline{u}_h}\int_{\amsmathbb R} \widehat{f}(s)\left( \sum_{g=1}^H \overline{u}_g
\widehat{\nu}_g(s)\right)\dd s \\
& = \frac{1}{2\pi \overline{u}_h \overline{\lambda}}\int_{\amsmathbb R}
|s|^{\frac{1}{2}}\,\widehat{f}(s) \cdot \frac{\big\langle
\widehat{\boldsymbol{\nu}}(s) \cdot \boldsymbol{\Theta}^{1/2}(\boldsymbol{\overline{u}})\big\rangle}{|s|^{\frac{1}{2}}}\dd s,
\end{split}
\end{equation*}
where $\boldsymbol{\Theta}^{1/2}$ is the square root of $\boldsymbol{\Theta}$ --- \textit{cf.} Definition~\ref{defsqur}. By the Cauchy--Schwarz inequality we deduce
\begin{equation}
\label{eq_x39} \bigg|\int_{\amsmathbb R} f(x) \nu_h(x)\dd x\bigg|^2\leq \frac{1}{4 \pi^2
\overline{u}_h^2 \overline{\lambda}} \cdot |\!|f|\!|_{\frac{1}{2}}^2 \cdot \int_{\amsmathbb R} \big|\big\langle
\widehat{\nu}(s) \cdot \boldsymbol{\Theta}^{1/2}(\overline{\boldsymbol{u}})\big\rangle \big|^2 \frac{\dd s}{|s|}
\end{equation}
Since $\boldsymbol{\Theta}^{1/2}$ is symmetric, we can further use the
Cauchy--Schwarz inequality in the form
\begin{equation}
\label{eq_x36}
\begin{split}
 \big|\big\langle
\widehat{\boldsymbol{\nu}}(s) \cdot \boldsymbol{\Theta}^{1/2}( \boldsymbol{\overline{u}})\big\rangle \big|^2 & = \big|\big\langle
\boldsymbol{\Theta}^{1/2}(\widehat{\boldsymbol{\nu}}(s)) \cdot \boldsymbol{\overline{u}}\big\rangle \big|^2 \leq \big\langle
\boldsymbol{\Theta}^{1/2}(\widehat{\boldsymbol{\nu}}(s)) \cdot \boldsymbol{\Theta}^{1/2}(\widehat{\boldsymbol{\nu}}(s))
\big\rangle^2 \big\langle \boldsymbol{\overline{u}} \cdot \boldsymbol{\overline{u}} \big\rangle^2
\\
& = \big\langle \boldsymbol{\Theta}(\widehat{\boldsymbol{\nu}}(s)) \cdot \widehat{\boldsymbol{\nu}}(s) \big\rangle^2.
\end{split}
\end{equation}
Combining \eqref{eq_x39} with \eqref{eq_x36} and \eqref{eq_x74} and then using Lemma~\ref{Lemma_I_quadratic_Fourier} we get
\begin{equation}
\label{eq_x35} \bigg|\int_{\amsmathbb R} f(x) \nu(x)\dd x\bigg|^2\leq \frac{1}{2 \pi^2\overline{u}_h^2 \overline{\lambda}} \cdot \mathfrak{D}^2[\boldsymbol{\nu},\boldsymbol{0}] \cdot |\!|f|\!|_{\frac{1}{2}}^2 \leq
 C\,\mathfrak{D}^2[\boldsymbol{\nu},\boldsymbol{0}]\, |\!|f|\!|_{\frac{1}{2}}^2.\qedhere
\end{equation}
\end{proof}

We can now complete the proofs of the two concentration of measures results announced in Section~\ref{section_rought_ann}, for the difference $\boldsymbol{\mu}_{\N} -\boldsymbol{\mu}$ between the empirical measure and the equilibrium measure in Lemma~\ref{Lemma_tail_bound_general}, and for the number of particles in a macroscopic segment in Corollary~\ref{Corollary_a_priori_0}.

\begin{proof}[Proof of Lemma~\ref{Lemma_tail_bound_general}] We fix $h \in [H]$ and assume without loss of generality that $f$ vanishes on $\bigcup_{g \neq h} [\hat a_g,\hat
b_g]$. The error done by replacing the empirical measure $\mu_{\N,h}$ with its regularization $\tilde{\mu}_{\N,h}$ from Lemma~\ref{mumutilde_lemma} can be estimated as
\begin{equation}
\label{eq_x37} \bigg|\int_{\amsmathbb R} f(x) \dd(\mu_{\N,h} - \tilde \mu_{\N,h})(x)\bigg| \leq \frac{C}{\N} \big(|\!|f|\!|_\textnormal{Lip}+|\!|f|\!|_{\infty}\big).
\end{equation}
 It remains to combine Lemma~\ref{Lemma_linear_through_distance} with $\boldsymbol{\nu}=\tilde{
\boldsymbol{\mu}}_\N-\boldsymbol{\mu}$, Equation~\eqref{eq_x37}, and Theorem~\ref{Prop_pseudodistance_bound}.
\end{proof}

\begin{proof}[Proof of Corollary~\ref{Corollary_a_priori_0}] First, assume that both $\mathfrak{a}$ and $\mathfrak{b}$ are finite, and denote
\[
P_{\N}[\mathfrak{a},\mathfrak{b}] = \#\big\{i \in [N]\quad \big| \quad\ell_i \in [\N \mathfrak{a},\N \mathfrak{b}]\big\}.
\]
Let $\delta > 0$ be small enough and take a smooth nonnegative function $f$, which is equal to $1$
on $[\mathfrak{a},\mathfrak{b}]$, bounded by $1$ everywhere, equal to $0$ outside $[\mathfrak{a}-\delta,\mathfrak{b}+\delta]$, and such that $|\!|f'|\!|_{\infty} \leq \frac{2}{\delta}$. In particular
\[
|\!| f |\!|_{\textnormal{Lip}} + |\!|f|\!|_{\infty} \leq \frac{2}{\delta} + 1.
\]
For the $\frac{1}{2}$-norm, we find
\begin{equation}
\label{eq_x38}
|\!| f |\!|_{\frac{1}{2}} = \bigg(\int_{\amsmathbb{R}} |\widehat{f'}(s)|^2 \frac{\dd s}{|s|} \bigg)^{\frac{1}{2}} = \bigg(-2\iint_{\amsmathbb{R}^2} \log|x - y|f'(x)f'(y)\dd x \dd y\bigg)^{\frac{1}{2}},
\end{equation}
where we used \eqref{eq_I_quadratic_Fourier} and $\int_{\amsmathbb{R}} f'(s) \dd s=0$ in the last equality. Since $f'(x)$ is non-zero only in $[\mathfrak{a}-\delta,\mathfrak{a}]\cup [\mathfrak{b},\mathfrak{b}+\delta]$, and $|f'(x)|\leq \frac{2}{\delta}$,
we further bound \eqref{eq_x38} by
\[
 \frac{\sqrt{8}}{\delta}\Bigg(\int_{\mathfrak{a}-\delta}^{\mathfrak{a}}\int_{\mathfrak{a}-\delta}^{\mathfrak{a}} \big|\log|x - y|\big|\dd x\dd y+
\int_{\mathfrak{b}}^{\mathfrak{b}+\delta}\int_{\mathfrak{b}}^{\mathfrak{b}+\delta} \big|\log|x - y|\big|\dd x\dd y+\int_{\mathfrak{a}-\delta}^{\mathfrak{a}}\int_{\mathfrak{b}}^{\mathfrak{b}+\delta} \big|\log|x - y|\big|\dd x\dd y
 \Bigg)^{\frac{1}{2}}.
\]
Using $\int_0^{\delta} \log x\,\dd x=\delta\big(1+\log(\frac{1}{\delta})\big)$, we obtain
\[
|\!|f|\!|_{\frac{1}{2}} \leq C \cdot \Bigg(\log(2+\mathfrak{b}-\mathfrak{a}) +\log\bigg(\frac{1}{\delta}\bigg)\Bigg)^{\frac{1}{2}}.
\]

Next, observe that
\begin{equation}
\label{nirzgbrgb}
\begin{split}
 \bigg|\int_{\amsmathbb{R}} \N f(x) \dd\mu_\N(x) - P_{\N}[\mathfrak{a},\mathfrak{b}]\bigg| & \leq \frac{2\delta\N}{\min_{h} \theta_{h,h}} + 2, \\
 \bigg|\int_{\amsmathbb{R}} \N f(x) \dd\mu(x) - \N \mu([\mathfrak{a},\mathfrak{b}])\bigg| & \leq \frac{2\delta\N}{\min_h \theta_{h,h}}.
\end{split}
\end{equation}
If we set $\frac{1}{\delta} = \N^{\frac{1}{2}}\log \N$ and $t = \N^{-\frac{1}{2}}\log \N$, we get
\[
t |\!| f |\!|_{\frac{1}{2}}+ \frac{C}{\N}\big( |\!|f|\!|_\textnormal{Lip} +|\!|f|\!|_\infty \big) \leq C'\Big(\N^{\frac{1}{2}}\log\N \cdot \big(\log(2+\mathfrak{b}-\mathfrak{a}) +\log( \N^{\frac{1}{2}} \log\N)\big)^{\frac{1}{2}} + 2\N^{\frac{1}{2}} \log \N + 1\Big)
\]
for a perhaps larger constant $C' > 0$. We then apply Lemma~\ref{Lemma_tail_bound_general} and use \eqref{nirzgbrgb} to get
\begin{equation}
\label{eq_x40}
 \P\bigg[\Big|P_{\N}[\mathfrak{a},\mathfrak{b}]- \N \mu([\mathfrak{a},\mathfrak{b}])\Big| \geq M_{\N}\bigg] \leq \exp\bigg(C \N \log\N- \frac{\N(\log\N)^2}{C}\bigg),
\end{equation}
where
\[
M_{\N} = \frac{4\N^{\frac{1}{2}}}{\min_h \theta_{h,h}\,\log \N} + 2 + C'\Big(\N^{\frac{1}{2}}\log\N \cdot \big(\log(2+\mathfrak{b}-\mathfrak{a}) +\log( \N^{\frac{1}{2}} \log\N)\big)^{\frac{1}{2}} + 2\N^{\frac{1}{2}} \log \N + 1\Big).
\]
Fix $D > 0$ large enough and assume that $\max(|\mathfrak{a}|,|\mathfrak{b}|)<D$. Then \eqref{eq_x40} implies
\begin{equation}
\label{eq_x41}
 \P\bigg[\Big|P_{\N}(\mathfrak{a},\mathfrak{b}) - \N\mu([\mathfrak{a},\mathfrak{b}])\Big|\geq C_D\, \N^{\frac{1}{2}} (\log \N)^2
\bigg]\leq \exp\bigg(-\frac{\N(\log\N)^2}{C}\bigg),
\end{equation}
where the constant $C_D$ depends only on $D$.

It remains to remove the dependence on $D$ in the bound. For that we fix a large $D$
for which the equilibrium measure $\mu$ is supported inside $[-D,D]$, and,
therefore, $\mu([-D,D]) = \frac{N}{\N}$. Then \eqref{eq_x41} with $\mathfrak{a}=\mathfrak{b}=D$ together with $P_\N[-\infty,\infty]=N$
 imply
\begin{equation}
\label{eq_x42}
 \P\Big[ P_{\N}[-\infty,-D]+ P_{\N}[D,+\infty] \geq C_D\,\N^{\frac{1}{2}}(\log\N)^2 \Big]\leq \exp\bigg(-\frac{\N(\log\N)^2}{C}\bigg).
\end{equation}
The statement of the corollary now follows by splitting an arbitrary interval $[\mathfrak{a},\mathfrak{b}]$ into $[\mathfrak{a},\mathfrak{b}] \cap [-D,D]$ and
$[\mathfrak{a},\mathfrak{b}]\setminus [-D,D]$.
\end{proof}

\section{Large deviations for the support}

\label{Section_ld_suport}
In the next chapters, it will become convenient to work under the auxiliary assumptions that each segment $[\hat
a_h,\hat b_h]$ for $h \in [H]$ is included in a $\N$-independent compact and contains only one band of the
equilibrium measure $\boldsymbol{\mu}$. This is part of the auxiliary Assumption~\ref{Assumptions_extra} on which the final results of this book do not depend. This section develops the basic mechanism for reducing a general case to this particular restrictive setting. For that we show that, with overwhelming probability as $\N$ becomes large, void segments have no particles and saturated segments are densely packed with particles. Hence, the configurations inside voids and saturations are essentially deterministic. This will allow us in Chapter~\ref{Chapter_conditioning} to condition by such configurations in voids and saturations and effectively continue the asymptotic analysis without them.

\subsection{Large deviations in voids}
\label{Sec : large_dev_voids}

\begin{theorem}\label{Theorem_ldpsup}
Suppose that Assumptions~\ref{Assumptions_Theta}, \ref{Assumptions_basic} and \ref{Assumptions_offcrit}
 hold. Choose $h \in [H]$, and a void interval $(\mathfrak{a},\mathfrak{b}) \subseteq [\hat{a}_h,\hat{b}_h]$ of the equilibrium measure. Assume that $\mathfrak{a}$ is a point separating a void from a band while $\mathfrak{b}$ is not (and is allowed to be $+\infty$). Then for any $\eps>0$, there exists $C>0$ depending only on $\eps$ and on the constants in the assumptions, such that
 \[
 \P\big[\exists i \in [N] \quad \big| \quad \ell_i \in [\N(\mathfrak{a}+\eps),\N\mathfrak{b}] \bigg] \leq C\exp\bigg(-\frac{\N}{C}\bigg).
\]
Similarly, if $\mathfrak{b}$ is a point separating a void from a band while $\mathfrak{a}$ is not (and is allowed to be $-\infty$), then
\[
\P\big[\exists i \in [N] \quad \big| \quad \ell_i \in [\N\mathfrak{a},\N(\mathfrak{b}-\eps)]\big]\leq C\exp\bigg(-\frac{\N}{C}\bigg).
\]
\end{theorem}
Note that Assumption~\ref{Assumptions_offcrit} (2) is necessary in general to make sure that no particle stay in the void regions, see \cite{FaGuSoWa} for sufficient conditions for a continuous $\sbeta$-ensemble to have particles in the void region.
Our assumptions guarantee that each void is adjacent to at least one band. Therefore, Theorem~\ref{Theorem_ldpsup} yields that asymptotically there are no particles in the
interior of voids. Since the equilibrium measure has compact support (Theorem~\ref{Theorem_equi_charact_repeat_2}), the result also shows that the rescaled location of particles $\frac{\ell_i}{\N}$ remain in a $\N$-independent compact with overwhelming probability. The rest of this subsection is devoted to the proof of Theorem~\ref{Theorem_ldpsup}, and all the lemmata follow the notations of the theorem, in particular we work with a fixed $h \in [H]$ such that $(\mathfrak{a},\mathfrak{b}) \subseteq [\hat{a}_h,\hat{b}_h]$. We only discuss the first part of the theorem, where $\mathfrak{a}$ separates a void from a band, as the second part is similar.

\begin{lemma} \label{Lemma_two_zetas} Suppose that Assumptions~\ref{Assumptions_Theta}, \ref{Assumptions_basic} and \ref{Assumptions_offcrit}
 hold. Assume that $\mathfrak{b}$ is finite. For any $\zeta_2 \in \big(0,\frac{\mathfrak{b} - \mathfrak{a}}{2}\big)$,
there exists $C > 0$ depending only on the constants in the assumptions and $\zeta_1 \in (0,\zeta_2)$ (which may depend on all the data), such that
\begin{equation}
\label{eq_x63}
 \inf_{x\in [\mathfrak{a}+\zeta_2, \mathfrak{b}]} \big(V^{\textnormal{eff}}_h(x) - V^{\textnormal{eff}}_h(\mathfrak{a}+\zeta_1)\big)\geq \frac{1}{C},
\end{equation}
and
\begin{equation}
\label{eq_x64}
 V^{\textnormal{eff}}_h(\mathfrak{a}+\zeta_1)=\inf_{x\in [\mathfrak{a}+\zeta_1,\mathfrak{b}]} V^{\textnormal{eff}}_h.
\end{equation}
Moreover, there exist $\zeta_1^\pm$ depending only on $\zeta_2$ and the constants in the
assumptions such that the choice of $\zeta_1$ can be made to obey $0<\zeta_1^-<\zeta_1<\zeta_1^+<\zeta_2$.
\end{lemma}
\begin{proof}
Note that $V^{\textnormal{eff}}_h$ is a continuous function:  for $V_h$ it is clear from \eqref{eq_V_logs} in Assumption~\ref{Assumptions_basic}; for the part involving the integral of logarithm against the equilibrium measure $\mu$ this follows from the integrability of the logarithm and the fact that $\mu$ has a bounded density. Moreover, if $V_h$ varies over the potentials meeting the requirements of Assumptions
\ref{Assumptions_Theta}, \ref{Assumptions_basic} and \ref{Assumptions_offcrit} with
fixed values of all the involved constants, then $V^{\textnormal{eff}}_h$ is even equicontinuous.

Assumption~\ref{Assumptions_offcrit} implies the existence of $\eta>0$ depending
only on the constants in the assumption and such that
\[
\inf_{x\in [\mathfrak{a}+\zeta_2, \mathfrak{b}]} \big(V^{\textnormal{eff}}_h(x) - V^{\textnormal{eff}}_h(\mathfrak{a})\big) >\eta.
\]
Choose $\frac{1}{C} = \frac{\eta}{2}$ and define $\zeta_1$ to be the minimal $\zeta$ such that
\[
\inf_{x\in [\mathfrak{a}+\zeta,\mathfrak{b}]} \big(V^{\textnormal{eff}}_h(x) -V^{\textnormal{eff}}_h(\mathfrak{a})\big)= \frac{\eta}{2}.
\]
Then \eqref{eq_x63} holds automatically, and the continuity of $V^{\textnormal{eff}}_h$ implies
that \eqref{eq_x64} also holds. Finally, the existence of the uniform $\zeta_1^\pm$ follows from the equicontinuity of $V^{\textnormal{eff}}_h(x)$.
\end{proof}

\begin{lemma} \label{Lemma_two_zetas_infinity} Suppose that Assumptions~\ref{Assumptions_Theta}, \ref{Assumptions_basic} and \ref{Assumptions_offcrit}
 hold. Assume that $h = H$ and $\mathfrak{b} = +\infty$. Take an arbitrary $\mathfrak{a}$ such that
$[\mathfrak{a},+\infty)$ is included in a void. Then, for any $c>0$ there exists an interval $(\zeta_1,\zeta_2) \subseteq \amsmathbb{R}_{> 0}$ and $\eps > 0$ such that
\[
\forall x \in ( \mathfrak{a}+\zeta_2,+\infty) \qquad V^{\textnormal{eff}}_H(x) - V^{\textnormal{eff}}_H(\mathfrak{a}+\zeta_1)> c + \eps \log\big(1 + |x|\big),
\]
and
\[
 V^{\textnormal{eff}}_H(\mathfrak{a}+\zeta_1)=\inf_{x\in [\mathfrak{a}+\zeta_1,+\infty)} V^{\textnormal{eff}}_H(x).
\]
There exists $C>0$ depending only on $\mathfrak{a},c$ and the constants in the assumptions, such that the choice of $(\zeta_1,\zeta_2)$ can be made to obey $\zeta_2<C$ and $\zeta_2-\zeta_1> \frac{1}{C}$.
\end{lemma}

An important difference between Lemmata~\ref{Lemma_two_zetas} and
\ref{Lemma_two_zetas_infinity} is that we fix $\zeta_2$ in the former and $c$
in the latter.

\begin{proof}[Proof of Lemma~\ref{Lemma_two_zetas_infinity}]
We define
\begin{equation*}
\begin{split}
\zeta_1 & = \min\Big\{\zeta \geq 0\quad \Big| \quad V^{\textnormal{eff}}_H(\mathfrak{a} + \zeta) = \inf_{x \in [\mathfrak{a} + \zeta,\mathfrak{b}]} V^{\textnormal{eff}}_H(x)\Big\}, \\
\zeta_2 & = \min\Big\{\zeta > \zeta_1 \quad \Big| \quad \forall x \in [\mathfrak{a} + \zeta,+\infty) \qquad V^{\textnormal{eff}}_H(x) - V^{\textnormal{eff}}_H(\mathfrak{a} + \zeta_1) > c + \eps \log\big(1 + |x|\big) \Big\}.
\end{split}
\end{equation*}
Such $\zeta_1,\zeta_2$ exist and are uniformly bounded from above due to the growth assumptions for the potential $V^{\textnormal{eff}}_H(x)$ as $x \rightarrow +\infty$ in Assumption~\ref{Assumptions_basic}. The equicontinuity of $V^{\textnormal{eff}}_H(x)$ on compact intervals implies that $\zeta_2-\zeta_1$ is uniformly bounded from below.
\end{proof}

Fix $\eps>0$ and choose $D>0$ independent of $\N$ and depending only on the constants in the assumptions, such that the
support of the equilibrium measure $\mu$ is inside $\big[-\frac{D}{2},\frac{D}{2}\big]$ --- \textit{cf.} Theorem~\ref{Theorem_equi_charact_repeat_2}. The exact value of $D$ will not be important.

\begin{definition}
\label{Definition_event_part_in_void} Let $\mathcal A_\eps\subseteq \W_\N$ denote the set of configurations such that
\begin{itemize}
\item there are at most $C \N^{\frac{1}{2}}(\log\N)^2$ particles $\ell_i$ in $[\N(\mathfrak{a}+\eps),\N\mathfrak{b}]\cup \big(-\infty,-\frac{\N D}{2}\big]\cup\big[\frac{\N D}{2},+\infty\big)$;
\item the event studied in Lemma~\ref{Lemma_tail_bound_general} with $t =\N^{-\frac{1}{2}}\log\N$ does not hold, \textit{i.e.}
\[
\forall f \in \mathscr{H}_{\textnormal{Lip},\frac{1}{2}} \quad \exists h \in [H] \qquad \bigg|\int_{\amsmathbb{R}} f(x)\dd(\mu_{\N,h} - \mu_h)(x)\bigg| \leq \frac{\log\N}{\N^{\frac{1}{2}}}\cdot |\!|f|\!|_{\frac{1}{2}} + \frac{C}{\N}\cdot \big(|\!|f|\!|_{\textnormal{Lip}} + |\!|f|\!|_{\infty}\big).
\]
\end{itemize}
For $0<\eps_1<\eps_2$, let $\mathcal{E}\big[\begin{smallmatrix} K_1 & K_2 \\ \eps_1 & \eps_2 \end{smallmatrix}\big] \subseteq \W_\N$ be the set of configurations with exactly $K_1$ particles in $[\mathfrak{a}+\eps_1, \mathfrak{b}]$ and $K_2 \leq K_1$ particles in $[\mathfrak{a}+\eps_2,\mathfrak{b}]$.
\end{definition}
Lemma~\ref{Lemma_tail_bound_general} and Corollary~\ref{Corollary_a_priori_0} imply the existence of $C>0$ such that
\[
\P[\mathcal{A}_\eps] \geq 1- C\exp\bigg(-\frac{\N(\log\N)^2}{C}\bigg).
\]
Therefore, it suffices to consider only the configurations from $\mathcal{A}_\eps$ in the proof of Theorem
\ref{Theorem_ldpsup}.

\begin{lemma} \label{Lemma_K_K_prime} Suppose that Assumptions~\ref{Assumptions_Theta}, \ref{Assumptions_basic} and \ref{Assumptions_offcrit} hold. Let $(\zeta_1,\zeta_2) \subseteq \amsmathbb{R}_{> 0}$ as in Lemma~\ref{Lemma_two_zetas}, or as in Lemma
\ref{Lemma_two_zetas_infinity} with $c=3$. Let $\delta>0$. There exists
$C>0$ depending only on the constants in the assumptions, such that if $K_1<\N^{\frac{1}{2}-\delta} K_2$, then
\[
 \P\Big[\mathcal{E}\big[\begin{smallmatrix} K_1 & K_2 \\ \zeta_1 & \zeta_2\end{smallmatrix}\big] \cap \mathcal A_{\zeta_1}\Big] \leq C\exp\bigg(-\frac{\N}{C}\bigg).
\]
\end{lemma}
\begin{proof}
 We assume $K_1\leq C \N^{\frac{1}{2}}(\log\N)^2$, as otherwise the probability vanishes since the set is empty. Take
$\boldsymbol{\ell}\in \mathcal{E}\big[\begin{smallmatrix} K_1 & K_2 \\ \zeta_1 & \zeta_2 \end{smallmatrix}\big] \cap \mathcal A_{\zeta_1}$. We write $I = \llbracket i_0,i_0 + K - 1\rrbracket$ for the indices of the particles in
$[\mathfrak{a}+\zeta_1,\mathfrak{b}]$
and $I^{\textnormal{c}} = [N]\setminus I$. We would like to bound $\P(\boldsymbol{\ell})$ by comparing it to the probability of the
configuration $\tilde{\boldsymbol{\ell}}$ deterministically obtained from $\boldsymbol{\ell}$ by the following rules.
\begin{itemize}
\item[$\bullet$] $\forall i \in I^{\textnormal{c}}\quad \tilde{\ell}_{i} = \ell_i$.
\item[$\bullet$] $\tilde{\ell}_{i_0}$ occupies the smallest available site in $[\N(\mathfrak{a}+\zeta_1),\N(\mathfrak{b}-\zeta_2)]$. Note that this position may be affected by the existence of a rescaled particle in the left
vicinity of $\mathfrak{a}+\zeta_1$, but only by a term of order $1$.
\item[$\bullet$] For any $i \in \llbracket i_0,i_0 + K_1 - 2\rrbracket$, we have $\tilde{\ell}_{i + 1} - \tilde{\ell}_{i} = \theta_{h,h}$, \textit{i.e.} the spacing between $\tilde{\ell}_{i + 1}$ and $\tilde{\ell}_i$ is minimal.
\end{itemize}
The configuration $\tilde{\boldsymbol{\ell}}$ has the same segment filling fractions as $\boldsymbol{\ell}$, and, therefore,
$\tilde{\boldsymbol{\ell}} \in \W_\N$. We decompose
\begin{equation}
\label{decompositionll} \frac{\P(\boldsymbol{\ell})}{\P(\tilde{\boldsymbol{\ell}})} = \prod_{k = 1}^3
P^{(k)}_{\N}(\boldsymbol{\ell}),
\end{equation}
with
\begin{equation*}
\begin{split}
P_\N^{(1)}(\boldsymbol{\ell})&= \prod_{\substack{i,j \in I \\ i < j}} \frac{\big(\ell_{j}-\ell_{i}\big)\cdot \Gamma\big(\ell_{j}-\ell_{i} + \theta_{h,h}\big)}{\Gamma\big(\ell_{j} - \ell_{i} + 1 - \theta_{h,h}\big)}\cdot \frac{\Gamma\big(\tilde{\ell}_{j} - \tilde{\ell}_{i} + 1 - \theta_{h,h}\big)}{(\tilde{\ell}_{j}-\tilde{\ell}_{i})\cdot \Gamma\big(\tilde{\ell}_{j}-\tilde{\ell}_{i} + \theta_{h,h}\big)}, \\
P_\N^{(2)}(\boldsymbol{\ell})&=\prod_{i \in I} \prod_{j \in I^\textnormal{c}} \bigg|\frac{\ell_i-\ell_j}{\tilde{\ell}_i-\ell_j}\bigg| \cdot \frac{\Gamma\big(|\ell_{i}-\ell_{j}|+\theta_{h,h(j)}\big)}{\Gamma\big(|\ell_{i}-\ell_{j}| + 1 - \theta_{h,h(j)}\big)} \cdot \frac{\Gamma\big(|\tilde{\ell}_{i}-\ell_{j}|+ 1 - \theta_{h,h(j)}\big)}{\Gamma\big(|\tilde{\ell}_i-\ell_j| + \theta_{h,h(j)}\big)}, \\
P_\N^{(3)}(\boldsymbol{\ell})&=\prod_{i\in I} \frac{w_h(\ell_i)}{w_h (\tilde \ell_i)}.
\end{split}
\end{equation*}
The absolute values in the second line avoids discussing the ordering of
particles appearing in the probability measure $\P$. We shall bound each of the three factors
separately.

\medskip

\noindent \textsc{Claim 1.} For $C>0$
\begin{equation}
\label{eq_Claim1} P_{\N}^{(1)}(\boldsymbol{\ell}) \leq \exp\Bigg(C K_1^2 \log\N +C K_1 \sum_{i\in
I} \log\bigg(\bigg|\frac{\ell_i}\N\bigg|+1\bigg)\Bigg).
\end{equation}
To justify this, we use \eqref{eq_Stirling_ratio} to replace the Gamma functions
involving the $\ell_i$s by $|\ell_i-\ell_{j}|^{2\theta_{h,h}}$, getting the error
\[
 \exp\bigg(C \sum_{\substack{i,j \in I \\ i< j}} \frac{1}{\ell_j-\ell_i}\bigg)\leq \exp\bigg(C
K_1 \sum_{k=1}^{K_1} \frac{1}{\theta_{h,h} \cdot k}\bigg)\leq \exp( C' K_1 \log K_1 ).
\]
We do the same for the factors involving the $\tilde \ell_i$s and get
\[
P_{\N}^{(1)}(\boldsymbol{\ell}) \leq \exp(C\, K_1 \log K_1) \cdot \prod_{\substack{i,j \in I \\ i < j}}
\frac{\Big(\frac{\ell_{j}}{\N}-\frac{\ell_{i}}\N\Big)^{2\theta_{h,h}}}{
\Big(\frac{\tilde \ell_{j}}{\N}-\frac{\tilde \ell_i}{\N}\Big)^{2\theta_{h,h}}}.
\]
For the numerators in the last formula we use the inequality $\log|x-y|\le
\log\big(1 + |x|)+\log\big(1 + |y|\big)$ and we bound the denominators of the last formula from
below by $\frac{\theta_{h,h}}{\N}$. As a result, we get \eqref{eq_Claim1}.

\medskip

\noindent \textsc{Claim 2.} The constant $D$ was chosen before the statement of Lemma~\ref{Lemma_K_K_prime}. With this choice, we claim that for any $\delta>0$ there exists $C>0$ such that
\begin{equation}
\label{eq_Claim2}
\begin{split}
 P_{\N}^{(2)}(\boldsymbol{\ell}) & \leq \exp\Bigg( C K_1 \N^{\frac{1}{2}+\frac{\delta}{2}} + 2 K_{\textnormal{far}} \N + C
\N^{\frac{1}{2}}
(\log\N)^2 \sum_{i\in I} \log\bigg(\bigg|\frac{\ell_i}\N\bigg|+1\bigg)\Bigg) \\
& \quad \times \exp\Bigg( 2 \N \sum_{i\in I} \sum_{g=1}^H \theta_{h,g}\int_{\hat
a'_{g}}^{\hat b'_{g}} \bigg[\log\bigg|\frac{\tilde
\ell_i}{\N}-x\bigg|-\log\bigg|\frac{\ell_i}{\N}-x\bigg|\bigg] \dd\mu_{g}(x) \Bigg),
\end{split}
\end{equation}
where $\boldsymbol{\mu}=(\mu_h)_{h = 1}^H$ is the equilibrium measure and $K_{\textnormal{far}}\leq K_1$
is the number of such $i \in I$ for which $|\ell_i|>D\N$. We
start by arguing in the same way as in Claim 1 to get a bound
\begin{equation}
\label{eq_x45} P_{\N}^{(2)}(\boldsymbol{\ell}) \leq \exp(CK_1\log \N) \cdot \prod_{\substack{i \in I \\
j \in I^\textnormal{c}}} \Bigg|\frac{\frac{\ell_{i}}{\N}-\frac{\ell_{j}}{\N}}{\frac{\tilde
\ell_{i}}{\N}-\frac{\ell_j}{\N}}\Bigg|^{2\theta_{h,h(j)}}.
\end{equation}

Let us split the particles in $I^\textnormal{c}$ into two groups: those for which
$|\ell_j|>D\N$ and all others. Denote the former $I^\textnormal{c}_{\textnormal{far}}$ and the latter
$I^\textnormal{c}_{\textnormal{close}}$. So, for any $j \in I^{\textnormal{c}}_{\textnormal{close}}$ we have
$|\ell_j|<D\N$. We observe that $\# I^\textnormal{c}_{\textnormal{far}} < \N^{\frac{1}{2}}(\log\N)^2$ by the definition of
$\mathcal A_{\zeta_1}$ and use for the particles from $I^\textnormal{c}_{\textnormal{far}}$ an elementary bound which holds
for some $C>0$
\[
\Bigg| \frac{\frac{\ell_i}{\N}-\frac{\ell_j}{\N} } {\frac{\tilde
\ell_i}{\N}-\frac{\ell_j}{\N}}\Bigg| \leq C\bigg(
\bigg|\frac{\ell_i}{\N}\bigg|+1\bigg)\quad \textnormal{whenever}\quad \bigg|\frac{\tilde\ell_i}{\N}\bigg|\leq
\frac{D}{2}\quad \textnormal{and} \quad \bigg|\frac{\ell_j}{\N}\bigg|>D.
\]
We conclude that
\begin{equation}
\label{eq_x53}
 \prod_{\substack{i \in I \\ j \in
I^\textnormal{c}_{\textnormal{far}}}} \Bigg|\frac{\frac{\ell_{i}}{\N}-\frac{\ell_{j}}{\N}}{\frac{\tilde
\ell_{i}}{\N}-\frac{\ell_j}{\N}}\Bigg|^{2\theta_{h,h(j)}} \leq \exp\Bigg(C
\N^{\frac{1}{2}}(\log\N)^2 \sum_{i\in I} \log\bigg(
\bigg|\frac{\ell_i}{\N}\bigg|+1\bigg)\Bigg).
\end{equation}
It remains to deal with $j\in I^\textnormal{c}_{\textnormal{close}}$. We split the
particles of $I$ as well into two groups: for $i\in I_{\textnormal{close}}$ we have
$|\ell_i|\leq D\N$ and for $i\in I_{\textnormal{far}}$ we have
$|\ell_i|>D\N$. We recall that $K_{\textnormal{far}}= \# I_{\textnormal{far}}$ and also set
$K_{\textnormal{close}}= \# I_{\textnormal{close}}$.

Let us start from $I_{\textnormal{close}}$, \textit{i.e.} deal with the case
$\max(\ell_i,\ell_j) \leq D \N$. We fix a small $\eta>0$ to
be specified later and introduce the function $L_{\eta}(x)$ through
\begin{equation*}
\begin{split}
 L_{\eta}'(x) & = \frac{x\,\mathbbm{1}_{[0,2D]}(|x|)}{x^2+\eta^2} + \frac{2D\,\mathbbm{1}_{[2D,2D + 1]}(|x|)}{4D^2+\eta^2}
\big((2D+1)\textnormal{sgn}(x) -x\big), \\
L_\eta(-2D-1) & = L_\eta(2D+1) = 0.
\end{split}
\end{equation*}
Observe that
\begin{equation}
\label{eq_x47}
\begin{split}
\forall x \in [-(2D+1),2D+1] &\qquad |L_{\eta}'(x)| \leq \frac{|x|}{x^2+\eta^2},\\
\forall x \in [-2D,2D] & \qquad L_\eta(x) = \frac{1}{2}\log(x^2+\eta^2) + C_D,
\end{split}
\end{equation}
where $C_{D} = L_{\eta}(2D) - \frac{1}{2}\log(4D^2 + \eta^2)$. We would like to replace the linear factors in \eqref{eq_x45} by exponentials of $L_\eta$. For that we use
\[
\forall x \in \amsmathbb{R}_{> 0}\qquad 0\leq \frac{1}{2}\log(x^2+\eta^2) -\log x =
\frac{1}{2}\log\bigg(1+\frac{\eta^2}{x^2}\bigg),
\]
which implies after replacing the spacings between particles by minimal possible ones
\begin{equation}
\label{eq_x46}\Bigg|\sum_{\substack{j \in I^{c}_{\textnormal{close}}}}
\log\bigg|\frac{\ell_{j} - \ell_i}{\N}\bigg| -
L_{\eta}\bigg(\bigg|\frac{\ell_{j} - \ell_i}{\N}\bigg|\bigg)\Bigg| \leq
 \sum_{k=1}^{N} \log\bigg(1 + \frac{\eta^2\N^2}{k^2\theta_{h,h}^2}\bigg) \leq
\N\int_{0}^{\hat{n}} \log\biggl(1+\frac{\eta^2}{x^2\theta_{h,h}^2}\bigg) \dd x \leq C'\eta\N,
\end{equation}
given that $\min_{h \in [H]} \theta_{h,h}$ is bounded away from $0$ by Condition 2. in Assumption~\ref{Assumptions_Theta}. Therefore,
\begin{equation}
\label{eq_x48}
\begin{split}
& \quad \prod_{\substack{i \in I_{\textnormal{close}} \\
j \in I^\textnormal{c}_{\textnormal{close}}}} \Bigg|
\frac{\frac{\ell_{i}}{\N}-\frac{\ell_{j}}{\N}}{\frac{\tilde \ell_{i}}{\N}-\frac{\ell_j}{\N}}\Bigg|^{2\theta_{h,h(j)}} \\
 & \le
 \exp\Bigg(CK_{\textnormal{close}}\,\eta\,\N+ \sum_{\substack{i \in I_{\textnormal{close}} \\
j \in I^\textnormal{c}_{\textnormal{close}}}}
2\theta_{h,h(j)}\bigg[L_{\eta}\bigg(\bigg|\frac{\ell_{i} - \ell_j}{\N}\bigg|\bigg)-
L_{\eta}\bigg(\bigg|\frac{\tilde
\ell_{i} - \ell_j}{\N}\bigg|\bigg)\bigg]\Bigg).
\end{split}
\end{equation}
Next we would like to use Lemma~\ref{Lemma_tail_bound_general} and replace the
sum over $j$ in \eqref{eq_x48} by the integral over the equilibrium measure. For that we need to bound
various norms of $L_{\eta}$. Equation \eqref{eq_x47} implies
\begin{equation}
|\!|L_{\eta}|\!|_{\textnormal{Lip}}\leq \frac{1}{2\eta}.
\end{equation}
Using \eqref{eq_x47} we get
 $|\!|L_{\eta}|\!|_{\infty}\leq C'_D |\log \eta|$ for small $\eta$, where $C'_D$ depends only
on $D$. We also estimate the $\frac{1}{2}$-norm of $L_\eta$ by using
\eqref{eq_x38} and the upper bound \eqref{eq_x47} for its derivative
\begin{equation*}
\begin{split}
|\!| L_{\eta} |\!|_{\frac{1}{2}}^2 & = \int_{\amsmathbb{R}} \frac{\big|\widehat{L_{\eta}'}(s)\big|^2}{|s|}\,\dd s = -2 \iint \limits_{|x|,|y| < 2D + 1} \log|x - y|\,L_{\eta}'(y)\,\dd x\dd y \\
& \leq 2 \iint\limits_{|x|,|y| < 2D + 1} \big|\log|x - y|\big| \cdot \frac{|x|\dd x}{x^2 + \eta^2} \cdot \frac{|y| \dd y}{y^2 + \eta^2} \\
& \leq 2 \iint\limits_{|x|,|y| < \frac{2D + 1}{\eta}} \big|\log \eta + \log|x - y|\big| \cdot \frac{|x|\dd x}{x^2 + 1} \cdot \frac{|y| \dd y}{y^2 + 1} \\
& \leq C\,\big|\log \eta\big| \bigg(\int_{-\frac{2D + 1}{\eta}}^{\frac{2D + 1}{\eta}} \frac{|x|\dd x}{x^2 + 1}\bigg)^2 + \iint\limits_{|x|,|y| < \frac{2D + 1}{\eta}} \frac{\big|\log|x - y| \big| \cdot |x||y|\dd x \dd y}{(x^2 + 1)(y^2 + 1)} \\
& \leq C\,\log^3\bigg(\frac{D}{\eta}\bigg).
\end{split}
\end{equation*}
We deduce from Lemma~\ref{Lemma_tail_bound_general} applied to $t= \N^{-\frac{1}{2}}\log \N$ and the definition of $\mathcal A_{\zeta_1}$ that for any $|u| < D$ and $g \in [H]$
\begin{equation}
\label{eq_x49}
\begin{split}
& \quad \Bigg|\frac{1}{\N} \sum_{\substack{j \in I^{c}_{\textnormal{close}}\\h(j)=g}}
L_{\eta}\bigg(\bigg|\frac{\ell_j}{\N} - u\bigg|\bigg) - \int_{\hat{a}'_g}^{\hat{b}'_g} L_{\eta}\big(|x -
u|\big)\dd\mu_{g}(x)\Bigg| \\
& \leq C \Bigg(\frac{1}{\N \eta}+\frac{1}{\N |\log\eta|} + \log^3\bigg(\frac{D}{\eta}\bigg)
\frac{\log
\N}{\sqrt{\N}}\Bigg)+\frac{1}{\N}\sum_{\substack{ j \notin I_{\textnormal{close}}^{c} \\ h(j) = g}} \Bigg|L_{\eta}\bigg(\bigg|\frac{\ell_j}{\N} - u\bigg|\bigg)\Bigg| \\
& \leq C \Bigg(\frac{1}{\N \eta} + \log^3\bigg(\frac{D}{\eta}\bigg) \frac{\log
\N}{\N^{\frac{1}{2}}}\Bigg)+C\,\frac{ \N^{\frac{1}{2}}(\log\N)^2}{\N |\log\eta|},
\end{split}
\end{equation}
where in the last inequality we used
\[
\#\big\{j \notin I_{\textnormal{close}} \,\, \big|\,\, h(j)=g\big\} \leq C
\N^{\frac{1}{2}} (\log\N)^2
\]
and the bound on $|\!|L_{\eta}|\!|_\infty$. We further specialize \eqref{eq_x49} to $u=\frac{\ell_i}{\N}$ or $u =
\frac{\tilde{\ell}_i}{\N}$.
 Choosing
$\eta=\N^{-\frac{1}{2}+\frac{\delta}{2}}$ with $\delta \in \big(0,\frac{1}{2}\big)$ and combining
with \eqref{eq_x48} --- and its analogue for the integral of logarithm rather than
sum --- we arrive at
\begin{equation}
 \label{eq_x54}
 \begin{split}
& \quad \prod_{\substack{i \in I_{\textnormal{close}} \\ j \in I^\textnormal{c}_{\textnormal{close}}}}
\Bigg|\frac{\frac{\ell_{i}}{\N}-\frac{\ell_{j}}{\N}}{\frac{\tilde \ell_{i}}{\N}-\frac{\ell_j}{\N}}\Bigg|^{2\theta_{h,h(j)}} \\
& \leq \exp\Bigg(CK_{\textnormal{close}} \N^{\frac{1}{2}+\frac{\delta}{2}} + 2 \N \sum_{\substack{i\in I_{\textnormal{close}} \\ g \in [H]}}
\theta_{h,g}\int_{\hat a'_{g}}^{\hat b'_{g}} \bigg[\log\bigg|\frac{
\ell_i}{\N}-x\bigg|-\log\bigg|\frac{\tilde \ell_i}{\N}-x\bigg|\bigg]\dd\mu_{g}(x)
\Bigg).
\end{split}
\end{equation}

It remains to study the case $(i,j) \in I_{\textnormal{far}} \times I^\textnormal{c}_{\textnormal{close}}$. For
the terms involving $\tilde \ell_i$, the same estimates as before work since
$|\tilde{\ell}_i|\leq D\N$. For the terms involving $\ell_i$, we use a bound
\[
\forall (i,j) \in I_{\textnormal{far}} \times I^{\textnormal{c}}_{\textnormal{close}} \qquad |\ell_i-\ell_j|<|2\ell_i|,
\]
and the fact that the segment filling fractions of $\boldsymbol{\ell} \in \mathcal{A}_{\zeta_1}$ need to be close to the segment filling
fractions of $\mu$, as follows from Lemma
\ref{Lemma_tail_bound_general} for an appropriate smooth approximation of the
indicator function of the interval $[\hat a_h, \hat b_h]$. We conclude that
\begin{equation}
\label{eq_x51}
 \prod_{\substack{i \in I_{\textnormal{far}} \\
j \in I^\textnormal{c}_{\textnormal{close}}}}
\Bigg|\frac{\frac{\ell_{i}}{\N}-\frac{\ell_{j}}\N}{\frac{\tilde \ell_{i}}{\N}-\frac{\ell_j}{\N}}\Bigg|^{2\theta_{h,h(j)}}
 \leq \exp\Bigg(CK_{\textnormal{far}} \N^{\frac{1}{2}+\frac{\delta}{2}} + 2 \N
\sum_{\substack{i\in I_{\textnormal{far}} \\ g \in [H]}} \theta_{h,g}\int_{\hat
a'_{g}}^{\hat b'_{g}} \bigg[\log\bigg|\frac{ 2\ell_i}{\N}\bigg|-\log\bigg|\frac{\tilde
\ell_i}{\N}-x\bigg|\bigg]\dd\mu_{g}(x)\Bigg).
\end{equation}
Note also that for $x$ in the support of $\mu$,
\[
\forall i \in I_{\textnormal{far}}\qquad 2 + \log\bigg|\frac{2\ell_i}{\N}\bigg|< \log \bigg|\frac{\ell_i}{\N}-x\bigg|.
\]
Therefore, \eqref{eq_x51} implies
\begin{equation}
\label{eq_x52}
\begin{split}
& \quad \prod_{\substack{i \in I_{\textnormal{far}} \\
k \in I^\textnormal{c}_{\textnormal{close}}}}
\Bigg|\frac{\frac{\ell_{i}}{\N}-\frac{\ell_{j}}{\N}}{\frac{\tilde \ell_{i}}{\N}-\frac{\ell_j}{\N}}\Bigg|^{2\theta_{h,h(j)}} \\
 \leq & \quad
 \exp\Bigg(C K_{\textnormal{far}} \N^{\frac{1}{2}+\frac{\delta}{2}}+2 K_{\textnormal{far}} \N + 2 \N\! \sum_{\substack{i\in I_{\textnormal{far}} \\ g \in [H]}} \theta_{h,g}\int_{\hat a_{g}}^{\hat b_{g}} \bigg[\log\bigg|\frac{
\ell_i}{\N}-x\bigg|-\log\bigg|\frac{\tilde \ell_i}{\N}-x\bigg|\bigg]
\dd\mu_{g}(x)\Bigg).
\end{split}
\end{equation}
Combining \eqref{eq_x45} with \eqref{eq_x53}, \eqref{eq_x54}, \eqref{eq_x52}, we get
the claimed inequality \eqref{eq_Claim2}.

\medskip

\noindent \textsc{Claim 3.} There exists $C>0$, such that
\begin{equation}
\label{eq_Claim3} P_{\N}^{(3)}(\boldsymbol{\ell}) \leq \exp\Bigg(C K_1 \log\N - \N \sum_{i\in
I}\bigg[V_h\bigg(\frac{\ell_i}{\N}\bigg) - V_h\bigg(\frac{\tilde\ell_i}{\N}\bigg)\bigg]\Bigg).
\end{equation}
As a matter of fact, this is a direct consequence of Assumption~\ref{Assumptions_basic}.

\medskip

\noindent \textsc{End of the proof of Lemma~\ref{Lemma_K_K_prime}.} We consider separately the cases of finite $\mathfrak{b}$ (using Lemma~\ref{Lemma_two_zetas}) and infinite $\mathfrak{b}$ (using Lemma~\ref{Lemma_two_zetas_infinity}).

For the first case note that $K_{\textnormal{far}}=0$ and $\log\ell_i$ for $i\in I$ is
uniformly bounded. Therefore, combining Claims 1, 2, and 3, we get
\[
 \P(\boldsymbol{\ell})\leq \P(\tilde{\boldsymbol{\ell}}) \cdot \exp\Bigg(CK_1^2\log\N+ CK_1 \N^{\frac{1}{2}+\frac{\delta}{2}}-\N\sum_{i\in
I}\bigg[V^{\textnormal{eff}}_h\bigg(\frac{\ell_i}{\N}\bigg)-V^{\textnormal{eff}}_h\bigg(\frac{\tilde
\ell_i}{\N}\bigg)\bigg]\Bigg).
\]
Further, using Assumption~\ref{Assumptions_basic} and the fact that $\mathfrak{a}+\zeta_1$ is bounded away both from
$\hat a'_h, \hat b'_h$ and from the support of $\mu$, we get
\[
 V^{\textnormal{eff}}_h\bigg(\frac{\tilde
\ell_i}{\N}\bigg)=V^{\textnormal{eff}}_h(\mathfrak{a}+\zeta_1)+\int_{\mathfrak{a}+\zeta_1}^{\frac{\tilde{\ell}_i}{\N}}
(V^{\textnormal{eff}}_h)'(x) \dd x \leq V^{\textnormal{eff}}_h(\mathfrak{a}+\zeta_1) +C \bigg(\frac{\tilde
\ell_i}{\N}- \mathfrak{a} - \zeta_1\bigg) \leq V^{\textnormal{eff}}_h(\mathfrak{a}+\zeta_1)+ \frac{C K}{\N},
\]
Then, using also $K_1<C \N^{\frac{1}{2}}(\log\N)^2$ we arrive to
\begin{equation}
\label{eq_x57}
 \P(\boldsymbol{\ell})\le\P(\tilde{\boldsymbol{\ell}}) \cdot \exp\Bigg(CK_1 \N^{\frac{1}{2}+\frac{\delta}{2}}-\N\sum_{i\in
I}\bigg[V^{\textnormal{eff}}_h\bigg(\frac{\ell_i}{\N}\bigg)-V^{\textnormal{eff}}_h(\mathfrak{a}+\zeta_1)\bigg]\Bigg).
\end{equation}
The definitions of $\zeta_1,\zeta_2$ imply that for each $i\in I$ we have
\[
V^{\textnormal{eff}}_h\bigg(\frac{\ell_i}{\N}\bigg)-V^{\textnormal{eff}}_h(\mathfrak{a}+\zeta_1)\geq 0,
\]
while for $K_2$ of these particles
\[
V^{\textnormal{eff}}_h\bigg(\frac{\ell_i}{\N}\bigg)-V^{\textnormal{eff}}_h(\mathfrak{a}+\zeta_1)\geq c >0.
\]
We conclude that for any $\boldsymbol{\ell}\in \mathcal{E}\big[\begin{smallmatrix} K_1 & K_2 \\ \zeta_1 & \zeta_2 \end{smallmatrix}\big]$
\begin{equation}
\label{eq_x55}
 \P(\boldsymbol{\ell})\le\P(\tilde{\boldsymbol{\ell}}) \cdot \exp\big(CK_1 \N^{\frac{1}{2}+\frac{\delta}{2}}- c K_2 \N \big).
\end{equation}
Observe that a given $\tilde{\boldsymbol{\ell}}$ might correspond at most to $(C\N)^{K_1}$ different choices
of $\boldsymbol{\ell}$. Thus, \eqref{eq_x55} implies
\begin{equation} \label{eq_x56}
 \P\Big[\mathcal{E}\big[\begin{smallmatrix} K_1 & K_2 \\ \zeta_1 & \zeta_2 \end{smallmatrix}\big]\Big] \leq \exp\big(CK_1\log\N+ CK_1 \N^{\frac{1}{2}+\frac{\delta}{2}}-c K_2\N \big).
\end{equation}
Using the inequality $K_1<\N^{\frac{1}{2}-\delta} K_2$, we see that \eqref{eq_x56} decays
exponentially fast as $\N\rightarrow\infty$, which finishes the proof of the finite $\mathfrak{b}$ case
of Lemma~\ref{Lemma_K_K_prime}.

\medskip

If $\mathfrak{b}$ is infinite we again combine Claims 1, 2, and 3. Instead of \eqref{eq_x57} we now
get, again using $K_1<C \N^{\frac{1}{2}}(\log\N)^2$, the inequality
\begin{equation}
\label{eq_x58}
\begin{split}
 \P(\boldsymbol{\ell}) & \leq \P(\tilde{\boldsymbol{\ell}}) \cdot \exp\Bigg(CK_1 \N^{\frac{1}{2}+\frac{\delta}{2}}-\N\sum_{i\in
I}\bigg[V^{\textnormal{eff}}_h\bigg(\frac{\ell_i}{\N}\bigg)-V^{\textnormal{eff}}_h(\mathfrak{a}+\zeta_1)\bigg] \\
& \qquad\qquad\qquad +C\N^{\frac{1}{2}}(\log\N)^2 \sum_{i\in I}
\log\bigg(\bigg|\frac{\ell_i}{\N}\bigg|+1\bigg)+2 K_{\textnormal{far}}\N\Bigg).
\end{split}
\end{equation}
 Recall that $\zeta_1$ and $\zeta_2$ are from Lemma~\ref{Lemma_two_zetas_infinity} with $c=3$. We have
\[
\forall i \in I\qquad V^{\textnormal{eff}}_h\bigg(\frac{\ell_i}{\N}\bigg)-V^{\textnormal{eff}}_h(\mathfrak{a}+\zeta_1)\geq 0.
\]
Exactly $K_2$ of these particles belong to $[\mathfrak{a}+\eps_2,\mathfrak{b}]$. Hence, denoting the set of indices of those $K_2$ particles by $I_2$,
\[
\forall i \in I_2\qquad V^{\textnormal{eff}}_h\bigg(\frac{\ell_i}{\N}\bigg)-V^{\textnormal{eff}}_h(\mathfrak{a}+\zeta_1)\geq 3+\eps
\log\bigg(\bigg|\frac{\ell_i}{\N}\bigg|+1\bigg).
\]
Therefore, using $K_{\textnormal{far}}\leq K_2$ and the fact that
$\log\big(\big|\frac{\ell_i}{\N}\big|+1\big)$ is uniformly bounded for $i\in I\setminus I_2$, we get
\begin{equation}
\label{eq_x59}
 \P(\boldsymbol{\ell})\le\P(\tilde{\boldsymbol{\ell}})\cdot \exp\Bigg(CK_1 \N^{\frac{1}{2}+\frac{\delta}{2}}-\N\sum_{i\in
I_2}\bigg[1+\bigg(\eps- C\frac{(\log\N)^2}{\N^{\frac{1}{2}}}\bigg)
\log\bigg(\bigg|\frac{\ell_i}{\N}\bigg|+1\bigg)\bigg]\Bigg).
\end{equation}
Assume that $\N$ is so large that $\eps- C \N^{-\frac{1}{2}}(\log\N)^2>
\frac{\eps}{2}$ in the last formula. We get
\begin{equation}
\label{eq_x60}
 \P(\boldsymbol{\ell})\le\P(\tilde{\boldsymbol{\ell}}) \cdot \exp\Bigg(CK_1 \N^{\frac{1}{2}+\frac{\delta}{2}}-\N\sum_{i\in
I_2}\bigg[1+\frac{\eps}{2} \log\bigg(\bigg|\frac{\ell_i}{\N}\bigg|+1\bigg)\bigg]\Bigg).
\end{equation}
Now we fix $\tilde{\boldsymbol{\ell}}$ and sum $\P(\boldsymbol{\ell})$ over all possible choices of $\boldsymbol{\ell}$
leading to this particular $\tilde{\boldsymbol{\ell}}$. We have $(C\N)^{K_1-K_2}\leq \exp(C'
K_1\log\N)$ possible choices for the positions of the particles from $I\setminus I_2$.
For each of the indices $i\in I_2$ we have infinitely many choices for the particle
position, but if we sum
\[
 \exp\Bigg(-\frac{\N\eps}{2} \log\bigg(\bigg|\frac{\ell_i}{\N}\bigg|+1\bigg)\Bigg)
\]
over such choices we get a bounded --- in fact, going to $0$ as $\N\rightarrow\infty$ ---
contribution. We conclude that for large $\N$
\begin{equation} \label{eq_x61}
 \P\Big[\mathcal{E}\big[\begin{smallmatrix} K_1 & K_2 \\ \zeta_1 & \zeta_2 \end{smallmatrix}\big]\Big] \le\exp\big(CK_1\log\N+ CK_1 \N^{\frac{1}{2}+\frac{\delta}{2}}-K_2\N \big).
\end{equation}
Because $K_1<\N^{\frac{1}{2}-\delta} K_2$, the last expression decays exponentially fast
as $\N\rightarrow\infty$, which finishes the proof of Lemma~\ref{Lemma_K_K_prime}.
\end{proof}

\begin{proof}[Proof of Theorem~\ref{Theorem_ldpsup}]
 Let us show that there exists $D>0$ such that
 \begin{equation}
 \label{eq_x62} \P\big[\exists i \in [N] \quad \big| \quad |\ell_i|>\N D\big]\leq C\exp\bigg(-\frac{\N}{C}\bigg).
\end{equation}
We only discuss the event $\ell_i>D\N$ as $\ell_i <-D\N$ can be handled similarly.

Take any $\mathfrak{a}_1$ such that $[\mathfrak{a}_1,+\infty)$ is included in a void. By Corollary
\ref{Corollary_a_priori_0} with overwhelming probability there are at most
$\N^{\frac{1}{2}}(\log\N)^2$ particles in $[\mathfrak{a}_1+\infty)$. Applying Lemma
\ref{Lemma_K_K_prime}
 with $\delta= \frac{1}{10}$, $\mathfrak{a}=\mathfrak{a}_1$, and summing over all possible choices of $K_1$ and $K_2$,
we find that for some $\mathfrak{a}_2>\mathfrak{a}_1$ with overwhelming probability there are at
most
\[
 \N^{\frac{1}{2}}(\log\N)^2 \cdot \N^{-\frac{1}{2}+\frac{1}{10}}=\N^{\frac{1}{10}}(\log\N)^2
\]
particles in $[\mathfrak{a}_2,+\infty)$. Applying Lemma~\ref{Lemma_K_K_prime} again with
$\delta= \frac{1}{10}$ and $\mathfrak{a}=\mathfrak{a}_2$, we find that for some $\mathfrak{a}_3>\mathfrak{a}_2$ with overwhelming
probability there at most
\[
\N^{\frac{1}{10}}(\log\N)^2 \cdot \N^{-\frac{1}{2}+\frac{1}{10}}<1
\]
particles in $[\mathfrak{a}_3,+\infty)$. Therefore we can choose $D=\mathfrak{a}_3$ and \eqref{eq_x62} is
proven.

\medskip

Now assume that the interval $[\mathfrak{a},\mathfrak{b}]$ in Theorem~\ref{Theorem_ldpsup} is finite and
--- without loss of generality --- that $\mathfrak{a}$ is a boundary point between a band and a void.
We then can again use the same iterative argument. Namely, by Corollary
\ref{Corollary_a_priori_0} for each $\eps>0$ with overwhelming probability
there are at most $\N^{\frac{1}{2}}(\log\N)^2$ particles in $[\mathfrak{a}+\eps,\mathfrak{b}]$. Then by
Lemma~\ref{Lemma_K_K_prime}
 with $\delta=\frac{1}{10}$, we deduce that for each $\eps>0$ there are at most $\N^{\frac{1}{10}}(\log\N)^2$ particles in
$[\mathfrak{a}+\eps,\mathfrak{b}]$ with probability exponentially close to $1$. Applying Lemma~\ref{Lemma_K_K_prime}
 with $\delta= \frac{1}{10}$ again, we deduce that for each $\eps>0$ with exponentially high
probability there are no particles in $[\mathfrak{a}+\eps,\mathfrak{b}]$.
\end{proof}

\subsection{Large deviations in saturations}
\label{Sec : large_dev_saturated}

Fix $h \in [H]$, take a configuration $\boldsymbol{\ell}\in \W_\N$ and an interval $[\mathfrak{a},\mathfrak{b}]\subset
[\hat a_h, \hat b_h]$, such that $\N(\mathfrak b-\mathfrak a)\geq \theta_{h,h}$, and define
\begin{equation}
\label{eq_particles_in_segment}
I=\big\{i \in [N] \quad \big|\quad \ell_i \in [\N\mathfrak{a},\N\mathfrak{b}]\big\}= \llbracket i_0+1,i_0+\#I\rrbracket.
\end{equation}

\begin{definition} \label{Definition_no_hole}We say that $\boldsymbol{\ell}$ has no holes in $[\mathfrak{a},\mathfrak{b}]$ if all spacings between particles indexed by $I$ and their immediate neighbors are minimal. In more detail, $I$ should be nonempty and
the following three conditions should be satisfied
\begin{enumerate}
\item Bulk: $\forall i \in \llbracket i_0+1,i_0 + \#I-1\rrbracket\quad \ell_{i + 1} = \ell_i + \theta_{h,h}$.
\item Left boundary: either\footnote{If $i_0=0$ and $a_1=-\infty$, then there is no way for this condition to be satisfied and there is no way for $\boldsymbol{\ell}$ to have no holes.} $\ell_{i_0+1}= a_h$ or $\ell_{i_0+1}=\ell_{i_0}+\theta_{h,h}$.

\item Right boundary: either\footnote{If $i_0+\#I=N$ and $b_H=+\infty$, then there is no way for this condition to be satisfied and no way for $\boldsymbol{\ell}$ to have no holes.}
$\ell_{i_0 + \#I}=b_h$ or $\ell_{i_0+\# I+1}= \ell_{i_0 + \#I}+\theta_{h,h}$.
\end{enumerate}
\end{definition}

\begin{theorem}\label{Theorem_ldsaturated}
Suppose that Assumptions~\ref{Assumptions_Theta}, \ref{Assumptions_basic} and \ref{Assumptions_offcrit} hold. Take $h \in [H]$ and a saturated interval $(\mathfrak{a},\mathfrak{b}) \subseteq [\hat{a}_h,\hat{b}_h]$ of the equilibrium measure $\boldsymbol{\mu}$. Assume that $\mathfrak{a}$ is a
point separating a saturation from a band, but $\mathfrak{b}$ is not.\footnote{As the equilibrium measure has finite mass, this situation requires $\mathfrak{b}$ to be finite.} Then for each $\eps>0$ there exists $C>0$ depending only on $\eps$ and the constants in the assumptions, such that
\[
\P\big[\textnormal{There are no holes in}\, [\mathfrak{a}+\eps,\mathfrak{b}]\big]\geq 1- C\exp\bigg(-\frac{\N}{C}\bigg).
\]
Similarly, if $\mathfrak{b}$ is a point separating a saturation from a band, but $\mathfrak{a}$
is not, then
\[
\P\big[\textnormal{There are no holes in}\, [\mathfrak{a},\mathfrak{b}-\eps]\big] \geq 1- C\exp\bigg(-\frac{\N}{C}\bigg).
\]
\end{theorem}

In the rest of this subsection we prove Theorem~\ref{Theorem_ldsaturated} in the case where $\mathfrak{a}$ separates a saturation from a band --- the second case being similar. All the lemmata follow the setting and notations of the theorem. The proof is in fact similar to that of Theorem
\ref{Theorem_ldpsup}, except that because we are looking at saturations, we have to deal with holes instead of particles. This is more delicate and gives rise to additional error terms. We should first give a precise definition of what we mean by a hole. This is easy in the case $\theta_{h,h}=1$: in this situation all particles in $[\hat a_h, \hat b_h]$ can only occupy pre-determined sites with integral spacing and we can define configuration of holes as those sites that are not occupied by particles. For general values of $\theta_{h,h}>0$, the ``allowed sites'' depend on the rest of the configuration so we need to be more careful.

\begin{figure}[t]
\centering
\includegraphics[width=0.35\linewidth]{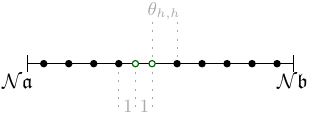}
\includegraphics[width=0.9\linewidth]{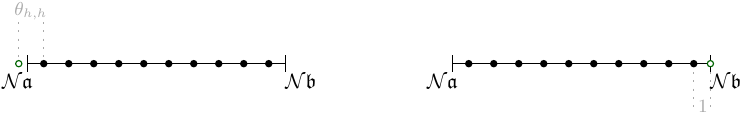}
 \caption{Top panel: two holes of the first type. Bottom panel: holes of the second and third types. \label{Fig_holes}}
\end{figure}

\begin{definition}\label{Definition_hole} Using the notation of \eqref{eq_particles_in_segment}, we say that $p\in\amsmathbb R$ is the position of a hole in $[\mathfrak{a},\mathfrak{b}]$ for the configuration $\boldsymbol{\ell} \in\W_\N$ if \emph{one of the following conditions} is satisfied.
\begin{itemize}
\item $\exists i \in \llbracket i_0 + 1,i_0 + \# I - 1\rrbracket\qquad \ell_i < p < p +\theta_{h,h} \leq \ell_{i + 1}\quad \textnormal{and} \quad p - \ell_i \in \amsmathbb{Z}_{> 0}$.
\item $\exists m \in \amsmathbb{Z}_{\geq 0} \qquad p=\ell_{i_0+1}-\theta_{h,h}-m\quad \textnormal{and} \quad p\geq \max\big(\N\mathfrak a-\theta_{h,h}, a_h, \ell_{i_0}+1\big)$.
\item $\exists m \in \amsmathbb{Z}_{\geq 0} \qquad p=\ell_{i_0+\#I}+m \quad \textnormal{and} \quad p\leq \max\big(\N \mathfrak b +1,b_h,\ell_{i_0+\#I+1}-\theta_{h,h}\big)$.
\end{itemize}
If $i_0=0$ or $i_0+\# I+1=N$, then the corresponding inequalities are ignored.
\end{definition}
The three types of holes are shown in Figure~\ref{Fig_holes}. We refer to holes of Definition~\ref{Definition_hole} as holes in $[\mathfrak a, \mathfrak b]$ although $\frac{p}{\N}$ may happen to be outside this segment; in particular, for each $\boldsymbol{\ell}\in\W_\N$ there is a well-defined total number of holes in $[\mathfrak a, \mathfrak b]$.

We remark that while the notion of ``having no holes'' (Definition~\ref{Definition_no_hole}) is somewhat canonical, to talk about holes themselves in Definition~\ref{Definition_hole} we made several \textit{ad hoc} choices. One choice is that particles at positions $0$ and $1+\theta_{h,h}$, according to Definition~\ref{Definition_hole}, lead to a hole at $1$, rather than at $\theta_{h,h}$. The exact inequalities we impose in the boundary cases also could have been different. Nevertheless, Definitions \ref{Definition_no_hole} and \ref{Definition_hole} are compatible with each other: if $I \neq \emptyset$, then a configuration $\boldsymbol{\ell} \in \W_\N$ has no holes according to Definition~\ref{Definition_no_hole} if and only if it has no holes according to Definition~\ref{Definition_hole}.

\medskip

The following lemma is proven in the same way as Lemma~\ref{Lemma_two_zetas}
\begin{lemma}
\label{Lemma_two_zetas_saturated} Suppose Assumptions~\ref{Assumptions_Theta}, \ref{Assumptions_basic} and \ref{Assumptions_offcrit} hold. For any $\zeta_2 \in \big(0,\frac{\mathfrak{b} - \mathfrak{a}}{2}\big)$, there exist $c>0$ depending only on the constants in the assumptions, and $\zeta_1 \in (0,\zeta_2)$ (which may depend on all the data) such that
\begin{equation}
\label{eq_x65} \sup_{x\in [\mathfrak{a}+\zeta_2, \mathfrak{b}]} \big(V^{\textnormal{eff}}_h(x) -
V^{\textnormal{eff}}_h(\mathfrak{a}+\zeta_1)\big) \leq c
\end{equation}
and
\begin{equation}
\label{eq_x66} V^{\textnormal{eff}}_h(\mathfrak{a}+\zeta_1)=\sup_{x\in [\mathfrak{a}+\zeta_1,\mathfrak{b}]}
V^{\textnormal{eff}}_h(x).
\end{equation}
Moreover, there exist $\zeta_1^\pm$ depending only on $\zeta_2$ and the constants in the
assumptions such that the choice of $\zeta_1$ can be made to obey $0<\zeta_1^-<\zeta_1<\zeta_1^+<\zeta_2$.
\end{lemma}

Fix $\eps>0$ and choose $D>0$ independent of $\N$ and depending only on the constants in the assumptions, such that the
support of the equilibrium measure $\mu$ is inside $\big[-\frac{D}{2},\frac{D}{2}\big]$ --- \textit{cf.} Theorem~\ref{Theorem_equi_charact_repeat_2}. The exact value of $D$ will not be important.

\begin{definition}
\label{Definition_event_hole_in_sat} Let $\mathcal A^{\textnormal{s}}_\eps\subset \W_\N$ denote the set of configurations such that
\begin{itemize}
\item there are at most $C \N^{\frac{1}{2}}(\log\N)^2$ holes in $[\mathfrak{a}+\eps,\mathfrak{b}]$;
\item there are no particles outside $(-D,D)$;
\item the event studied in Lemma~\ref{Lemma_tail_bound_general} with $t= \N^{-\frac{1}{2}}\log\N$ is not realized, \textit{i.e.}
\[
\forall f \in \mathscr{H}_{\textnormal{Lip},\frac{1}{2}} \quad \exists h \in [H] \qquad \bigg|\int_{\amsmathbb{R}} f(x)\dd(\mu_{\N,h} - \mu_h)(x)\bigg| \leq \frac{\log\N}{\N^{\frac{1}{2}}} \cdot |\!|f|\!|_{\frac{1}{2}} + \frac{C}{\N} \cdot \big(|\!|f|\!|_{\textnormal{Lip}} + |\!|f|\!|_{\infty}\big).
\]
\end{itemize}
For $0<\eps_1<\eps_2$, let $\mathcal{E}^{\textnormal{s}}\big[\begin{smallmatrix} K_1 & K_2 \\ \eps_1 & \eps_2 \end{smallmatrix}\big] \subset \W_\N$ be the set of configurations having $K_1$ holes in $[\mathfrak{a}+\eps_1, \mathfrak{b}]$ and $K_2 \leq K$ holes in $[\mathfrak{a}+\eps_2,\mathfrak{b}]$.
\end{definition}
Lemma~\ref{Lemma_tail_bound_general}, Corollary~\ref{Corollary_a_priori_0}, and Theorem~\ref{Theorem_ldpsup} imply the existence of a constant $C > 0$ such that
\[
\P\big[\mathcal A^{\textnormal{s}}_\eps\big] \geq 1- C\exp\bigg(-\frac{\N}{C}\bigg).
\]
Therefore, it suffices to consider only the configurations from $\mathcal A_\eps^{\textnormal{s}}$ in the proof of Theorem~\ref{Theorem_ldsaturated}.

\begin{lemma} \label{Lemma_K_K_prime_saturated} Suppose that Assumptions~\ref{Assumptions_Theta}, \ref{Assumptions_basic} and \ref{Assumptions_offcrit} hold. Let $(\zeta_1,\zeta_2) \subset \amsmathbb{R}_{> 0}$ be as in Lemma
\ref{Lemma_two_zetas_saturated}, and $\delta>0$. There exists $C>0$ depending only on the constants in the assumptions such that, if $K_1<\N^{\frac{1}{2}-\delta} K_2$, then
\[
 \P\Big[\mathcal{E}^{\textnormal{s}}\big[\begin{smallmatrix} K_1 & K_2 \\ \zeta_1 & \zeta_2 \end{smallmatrix}\big] \cap \mathcal A^{\textnormal{s}}_{\zeta_1}\Big] \leq C\exp\bigg(-\frac{\N}{C}\bigg).
\]
\end{lemma}
\begin{proof} Parts of the proof closely follow those of Lemma~\ref{Lemma_K_K_prime}, so we
omit some details. We assume $K_1 \leq C \N^{\frac{1}{2}}(\log\N)^2$, as otherwise the probability vanishes since the set is empty.
We take $\boldsymbol{\ell}\in \mathcal{E}^{\textnormal{s}}\big[\begin{smallmatrix} K_1 & K_2 \\ \zeta_1 & \zeta_2 \end{smallmatrix}\big] \cap \mathcal
A^{\textnormal{s}}_{\zeta_1}$, use the notation \eqref{eq_particles_in_segment} for the segment $[\mathfrak a+\zeta_1, \mathfrak b]$, and let $p_1 < \cdots < p_{K_1}$ be the position of the holes in $[\mathfrak{a} + \zeta_1,\mathfrak{b}]$. Denote $u=\min(p_1,\ell_{i_0+1})$ and $v=\max(p_{K_1},\ell_{i_0+\#I})$. Altogether we have $P = \# I$ particles and $K_1$ holes inside the segment $[u,v]$.

Our plan is to bound from above $\P(\boldsymbol{\ell})$ by comparing it to the probability of a modified
configuration $\tilde{\boldsymbol{\ell}}$. Outside $[u,v]$ the configurations $\boldsymbol{\ell}$ and
$\tilde{\boldsymbol{\ell}}$ coincide. Inside $[u,v]$ we obtain $\tilde{\boldsymbol{\ell}}$ from $\boldsymbol{\ell}$ by
pushing all the holes to the leftmost positions and pushing all
the particles to the rightmost positions. The number of particles and holes inside
$[u,v]$ remain unchanged. In more detail,
the new positions
of holes in $[\mathfrak{a} + \zeta_1,\mathfrak{b}]$ for $\tilde{\boldsymbol{\ell}}$ are denoted through $\tilde{p}_1, \ldots, \tilde{p}_{K_1}$ and we set
\begin{equation*}
\begin{split}
\forall k \in [K_1] \qquad & \tilde p_k=u+k-1 \\
\forall j \in [P] \qquad & \tilde \ell_{i_0+j}=v+\theta_{h,h}(j-P)
\end{split}
\end{equation*}
as in Figure~\ref{Fig_holes_pushed}. We aim at proving
\begin{equation}
\label{eq_x67}
 \P(\boldsymbol{\ell})\leq\P(\tilde{\boldsymbol{\ell}}) \cdot \exp\Bigg(CK \N^{\frac{1}{2}+\frac{\delta}{2}}+\frac{\N}{\theta_{h,h}}\sum_{i=1}^{K_1}
\bigg[V^{\textnormal{eff}}_h\bigg(\frac{p_i}{\N}\bigg)-V^{\textnormal{eff}}_h(\mathfrak{a}+\zeta_1)\bigg]\Bigg).
\end{equation}
Note that \eqref{eq_x67} would lead to the statement of Lemma
\ref{Lemma_K_K_prime_saturated} in the same way as \eqref{eq_x57} led to Lemma
\ref{Lemma_K_K_prime}, and therefore, in the rest of the proof we only justify
\eqref{eq_x67}.

\begin{figure}[t]
\centering
\includegraphics[width=0.45\linewidth]{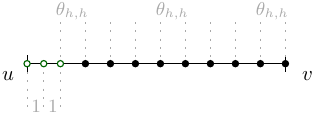}
 \caption{In $\tilde \ell$ particles $\tilde \ell_{i_0+j}$ are pushed to the right and holes $\tilde p_k$ are pushed to the left.\label{Fig_holes_pushed}}
\end{figure}

The configuration $\tilde{\boldsymbol{\ell}}$ can be obtained from $\boldsymbol{\ell}$ in $K_1$ steps, such that in the $k$-th step the
hole $p_k$ is moved to $\tilde p_k$. In this way we define a sequence of particle
configurations
\begin{equation}
\label{transfomrmum}
\boldsymbol{\ell}=\boldsymbol{\ell}^{(0)},\boldsymbol{\ell}^{(1)},\ldots,\boldsymbol{\ell}^{(K)}=\tilde{\boldsymbol{\ell}},
\end{equation}
such that $\boldsymbol{\ell}^{(k)}$ differs from $\boldsymbol{\ell}^{(k-1)}$ by the move of the hole $p_k$,
and study the ratio
\[
\frac{\P(\boldsymbol{\ell}^{(k-1)})}{\P(\boldsymbol{\ell}^{(k)})}.
\]
The argument is the same for
each $k$, so we only consider the $k=1$ case, \textit{i.e.} the move of the leftmost hole $p_1$.
Let $\ell_{i_0+1}, \ell_{i_0+2},\ldots \ell_{i_0+M}$ denote all the particles inside
$[\N(\mathfrak{a}+\zeta_1),p_1]$. These particles are densely
packed, \textit{i.e.}
\[
\forall m \in [M]\qquad \ell_{i_0+m}=\ell_{i_0+1}+(m-1)\theta_{h,h}.
\]
The move $p_1 \rightarrow \tilde{p}_1$ is equivalent to shifting all these particles to the right by $1$, \textit{i.e.}
$\ell_{i_0 + m} \rightarrow \ell_{i_0 + m} +1$ for $m \in [n]$. Hence, using Assumption~\ref{Assumptions_basic}
\begin{equation}
\label{eq_x68}
\begin{split}
 \frac{\P(\boldsymbol{\ell})}{\P(\boldsymbol{\ell}^{(1)})} & =\prod_{i=i_0 + 1}^{i_0 + M} \exp\Bigg(
 -\N V_h\bigg(\frac{\ell_i}{\N}\bigg)+ \N V_h\bigg(\frac{\ell_i+1}{\N}\bigg) + \err_h(\ell_i)-\err_h(\ell_i+1)\Bigg)
 \\
 & \quad \times \prod_{\substack{j\in
 [N] \\ j \notin \llbracket i_0+1,i_0+M\rrbracket}} \frac{|\ell_i-\ell_j|\cdot \Gamma\big(|\ell_i-\ell_j|+\theta_{h,h(j)}\big)}
 {\Gamma\big(|\ell_i-\ell_j|+1-\theta_{h,h(j)}\big)} \cdot
 \frac{\Gamma\big(|\ell_i+1-\ell_j|+1-\theta_{h,h(j)}\big)}{|\ell_i+1-\ell_j|\cdot\Gamma\big(|\ell_i+1-\ell_j|+\theta_{h,h(j)}\big)}.
\end{split}
\end{equation}
We now bound the various parts of \eqref{eq_x68}. First, note that since $\err_h(x)$ is
piecewise monotonous and uniformly bounded by $C\log\N$, we have a bound
\begin{equation}
\label{eq_x73}
 \prod_{i=i_0 + 1}^{i_0 + M} \exp\big(
 \err_h(\ell_i)-\err_h(\ell_i+1)\big) \leq \exp(C\log\N).
\end{equation}
Let us remark that if $\theta_{h,h}=1$, then the product in \eqref{eq_x73} is
telescoping because $\ell_{i + 1}=\ell_{i} +\theta_{h,h} = \ell_{i} + 1$ and we do not need the assumption on the monotonicity of $\err_h(x)$.

Next, we bound the factors containing the potential $V_h$. For this purpose, let us first assume that all relevant $\ell_i$s are at distance at least $\theta_{h,h}+1$ from the points $\N \hat a'_h$ and $\N\hat b'_h$. Then bounding $V''_h$ by Assumption~\ref{Assumptions_basic} (see Remark~\ref{rembounder}) and since $M \leq N$, we have
\begin{equation}
\label{eq_x69}
\begin{split}
& \quad \prod_{i=i_0 + 1}^{i_0 + M} \exp\Bigg(
 -\N V_h\bigg(\frac{\ell_i}{\N}\bigg)+
 \N V_h\bigg(\frac{\ell_i+1}{\N}\bigg)\Bigg) \\
 & \leq \prod_{i=i_0 + 1}^{i_0 + M} \exp\Bigg(
 \frac{\N}{\theta_{h,h}} \bigg[V_h\bigg(\frac{\ell_i+\theta_{h,h}}{\N}\bigg)-
 V_h\bigg(\frac{\ell_i}{\N}\bigg)+ \frac{C}{\N^2} \sup_{\N x\in [\ell_i,
 \ell_i+\max(\theta_{h,h},1)]} |V_h''(x)|
  \bigg] \Bigg) \\
  & \leq \exp\bigg(\frac{C}{\N} \sum_{k=1}^N \frac{\N}{k}\bigg)\cdot \prod_{m=1}^{M} \exp\Bigg(
 \frac{\N}{\theta_{h,h}}\bigg[V_h\bigg(\frac{\ell_{i_0+1}+m\theta_{h,h}}{\N}\bigg)-V_h\bigg(\frac{\ell_{i_0+1}+(m-1)\theta_{h,h}}{\N}\bigg)\bigg]\Bigg),
 \end{split}
 \end{equation}
 where $\frac{C}{\N} \sum_{k=1}^N \frac{\N}{k}$ is an upper bound for the terms involving $V''_h$.
 The product over $m$ in \eqref{eq_x69} is telescoping, and therefore, \eqref{eq_x69} simplifies
 to the following expression in terms of the locations of the holes
\begin{equation}
\label{eq_x71}
 \exp\Bigg(\frac{\N}{\theta_{h,h}}\bigg[V_h\bigg(\frac{p_1}{\N}\bigg)-V_h\bigg(\frac{\tilde p_1}{\N}\bigg)\bigg] +
 C\log\N\Bigg)
\end{equation}
For $\ell_i$ close to $\N \hat a'_h$ or $\N \hat b'_h$, the bound \eqref{eq_x69} needs an additional clarification because the $x\log|x|$ parts of $V_h(x)$ lead to exploding second derivatives $V''_h(x)$. In fact, $\ell_i$ in \eqref{eq_x69} cannot be close to $\N \hat a'_h$ because of the inequality $\frac{\ell_i}{\N} \geq \mathfrak a+\zeta_1$. It can be close to $\N \hat b'_h$, though, but the total number of $\ell_i$ at distance at most $\theta_{h,h}+1$ from $\N \hat b'_h$ is uniformly bounded. Hence, we can use an elementary bound
\[
\bigg|(\ell_i-\N \hat b'_h) \log\bigg(\frac{\ell_i}{\N}-\hat b'_h\bigg)\bigg|\leq C\log\N
\]
for the exploding part of $V_h(x)$. Therefore, \eqref{eq_x71} remains valid.

Jumping to the second line of \eqref{eq_x68}, for fixed $j$ the factor there is equal to
\[
 \frac{\ell_i-\ell_j}{\ell_i-\ell_j+1} \cdot
 \frac{\ell_i-\ell_j+1-\theta_{h,h(j)}}{\ell_i-\ell_j+\theta_{h,h(j)}}.
\]
both for $j>i_0+M$ and for $j\leq i_0$. Therefore, the second line of \eqref{eq_x68} is
\begin{equation}
\label{eq_x69s0}
\begin{split}
& \prod_{i=i_0 + 1}^{i_0 + M} \,\, \prod_{\substack{j\in
 [N] \\ j \notin \llbracket i_0+1,i_0+M\rrbracket}} \frac{\ell_i-\ell_j}{\ell_i-\ell_j+1} \cdot
 \frac{\ell_i-\ell_j+1-\theta_{h,h(j)}}{\ell_i-\ell_j+\theta_{h,h(j)}} \\
 & = \prod_{\substack{j\in [N] \\ j \notin \llbracket i_0+1,i_0+M\rrbracket}}
 \exp\Bigg(\sum_{m=0}^{M-1}\bigg[\log\bigg|\frac{\ell_{i_0+1}+m\theta_{h,h} - \ell_j}{\N}\bigg|
 -\log\bigg|\frac{\ell_{i_0+1}+m\theta_{h,h}+1 - \ell_j}{\N}\bigg| \\
 & \qquad\qquad\qquad\qquad +
 \log\bigg|\frac{\ell_{i_0+1}+m\theta_{h,h}+1-\theta_{h,h(j)} - \ell_j}{\N}\bigg|
 - \log\bigg|\frac{\ell_{i_0+1}+m\theta_{h,h}+\theta_{h,h(j)} - \ell_j}{\N}\bigg|\bigg]\Bigg).
\end{split}
\end{equation}
Let us explain how to bound the second line of \eqref{eq_x69s0} from above. Assume without loss of generality that $j\leq i_0$. Using the fact that the derivative of $x \mapsto \log\big(x-\frac{\ell_j}{\N}\big)$ is positive and monotonous for $x > \frac{\ell_j}{\N}$, we have for some $u \in (0,1)$
\[
\log\bigg(x-\frac{\ell_j}{\N}\bigg)-\log\bigg(x+\frac{1}{\N}-\frac{\ell_j}{\N}\bigg)=
-\frac{1}{\N}\log'\bigg(x+\frac{u}{\N}-\frac{\ell_j}{\N}\bigg)\leq - \frac{1}{\N}\log'\bigg(x-\frac{\ell_j}{\N}\bigg).
\]
Hence, using monotonicity of $x \mapsto \log'\big(x-\frac{\ell_j}{\N}\big)$ again, we have
\begin{equation*}
\begin{split}
& \quad \sum_{m=0}^{M-1}\Bigg(\log\bigg(\frac{\ell_{i_0+1}+m\theta_{h,h} - \ell_j}{\N}\bigg)
 -\log\bigg(\frac{\ell_{i_0+1}+m\theta_{h,h}+1 - \ell_j}{\N}\bigg)\Bigg) \\
 & \leq -\frac{1}{\N}\sum_{m=0}^{M-1}\log'\bigg(\frac{\ell_{i_0+1}+m\theta_{h,h} - \ell_j}{\N}\bigg) \\
 & \leq -\frac{1}{\theta_{h,h}}\int_{\frac{\ell_{i_0+1}}{\N}}^{\frac{\ell_{i_0+1}}{\N}+(m-1)\theta_{h,h}} \log'\bigg(x-\frac{\ell_j}{\N}\bigg)+C\log\N
 \end{split}
 \end{equation*}
and the integral in the last line is equal to
 \[
 \frac{1}{\theta_{h,h}}\Bigg(\log\bigg(\frac{\ell_{i_0+1 - \ell_j}}{\N}\bigg)-\log\bigg(\frac{\ell_{i_0+M} - \ell_j}{\N}\bigg)\Bigg) = \frac{1}{\theta_{h,h}}\Bigg(
 \log\bigg|\frac{\tilde p_1 - \ell_j}{\N}-\frac{\ell_j}{\N}\bigg|-\log\bigg|\frac{p_1 - \ell_j}{\N}\bigg|\Bigg).
\]
Doing similar computations for $j>i_0+M$ and for the third line of \eqref{eq_x69s0}, we upper bound \eqref{eq_x69s0} by
\begin{equation}
\label{eq_x69s}
\prod_{\substack{j\in
 [N] \\ j \notin \llbracket i_0+1,i_0+g\rrbracket}}
 \exp\Bigg(\frac{2 \theta_{h,h(j)}}{\theta_{h,h}}\bigg[
 \log\bigg|\frac{\tilde p_1 - \ell_j}{\N}\bigg| -\log\bigg|\frac{p_1 - \ell_j}{\N}\bigg|\bigg] +C\log\N \Bigg).
\end{equation}
It is somewhat inconvenient that the product in \eqref{eq_x69s} does not range over all
$j$, \textit{i.e.}, it avoids $M$ of the indices --- and $M$ might be large. However, we
can change that by observing
\begin{equation*}
\begin{split}
& \quad \prod_{\substack{i_0 + 1 \leq j \leq i_0 + M \\ \ell_j -\tilde p_1> 1}}
 \exp\bigg(2
 \log\bigg|\frac{\tilde p_1 - \ell_j}{\N}\bigg|-2\log\bigg|\frac{p_1 - \ell_j}{\N}\bigg|
 \bigg) \\
 & \geq \exp\Bigg(-C\log\N+\frac{2\N}{\theta_{h,h}}\int_{\frac{\tilde p_1}{\N}}^{\frac{ p_1}{\N}}\bigg(\log\bigg|\frac{\tilde p_1}{\N}-x\bigg|-\log\bigg|\frac{p_1}{\N}-x\bigg|\bigg)\dd x\Bigg) =\exp(-C\log\N),
\end{split}
\end{equation*}
where the inequality is obtained by treating the sum over $\ell_j$ --- which are all between $\tilde p_1$ and $p_1$ --- as a Riemann sum and replacing it by the corresponding integral. We conclude that the second line of \eqref{eq_x68} is bounded from above by
\begin{equation}
\label{eq_x70}
 \prod_{\substack{j \in [N] \\ |\tilde p_1 - \ell_j>1}}
 \exp\Bigg(\frac{2 \theta_{h,h(j)}}{\theta_{h,h}}\bigg[
 \log\bigg|\frac{\tilde p_1 - \ell_j}{\N}\bigg|-\log\bigg|\frac{p_1 - \ell_j}{\N}\bigg|\bigg] +C\log\N \Bigg).
\end{equation}
At this point we can repeat the argument in the proof of Lemma
\ref{Lemma_K_K_prime} (see Claim 2 there), replace the logarithm in \eqref{eq_x70} by $L_\eta$, and conclude
that the product over the particles $\ell_j$ can be approximated by the exponent of
the integral of the density of the equilibrium measure $\mu$. Combining with \eqref{eq_x71}, we reach the inequality
\begin{equation}
\label{eq_x72}
\begin{split}
 \frac{\P(\boldsymbol{\ell})}{\P(\boldsymbol{\ell}^{(1)})} & \leq
\exp\Bigg(\frac{\N}{\theta_{h,h}}\bigg[V_h\bigg(\frac{p_1}{\N}\bigg)-\sum_{g=1}^H
2\theta_{h,g} \int_{\hat a'_{g}}^{\hat b'_{g}}\log\bigg|x-\frac{p_1}{\N}\bigg|\mu_{g}(x)\dd x \\
& \quad\qquad -V_h\bigg(\frac{\tilde p_1}{\N}\bigg)+\sum_{g=1}^H 2\theta_{h,g} \int_{\hat a'_{g}}^{\hat
b'_{g}}\log\bigg|x-\frac{\tilde p_1}{\N}\bigg|\mu_{g}(x) \dd x\bigg] +
 C \N^{\frac{1}{2}+\frac{\delta}{2}}\Bigg),
 \end{split}
 \end{equation}
 where the $C \N^{\frac{1}{2}+\frac{\delta}{2}}$ remainder is like in \eqref{eq_x54}.

The argument for the $k$-th step for $k \geq 2$ of the moves \eqref{transfomrmum} is very similar and leads to the same bound as
\eqref{eq_x72}. Multiplying all these bounds over $k \in [K_1]$, we get
\begin{equation}
\label{eq_x106}
\begin{split}
 \frac{\P(\boldsymbol{\ell})}{\P(\tilde{\boldsymbol{\ell})}} & \leq
\exp\Bigg(\frac{\N}{\theta_{h,h}}\sum_{k=1}^{K_1} \bigg[V_h\bigg(\frac{p_k}{\N}\bigg)-\sum_{g=1}^H
2\theta_{h,g} \int_{\hat a'_{g}}^{\hat b'_{g}}\log\bigg|x-\frac{p_k}{\N}\bigg|\mu_{g}(x)\dd x \\
& \quad\qquad -V_h\bigg(\frac{\tilde p_k}{\N}\bigg)+\sum_{g=1}^H 2\theta_{h,g} \int_{\hat a'_{g}}^{\hat
b'_{g}}\log\bigg|x-\frac{\tilde p_k}{\N}\bigg|\mu_{g}(x) \dd x\bigg] + C K \N^{\frac{1}{2}+\frac{\delta}{2}}\Bigg).
 \end{split}
 \end{equation}
 It remains to replace $\frac{\tilde p_k}{\N}$ by $\mathfrak a+\zeta_1$ everywhere in the last inequality. Note that
 \[
 \forall k \in [K_1]\qquad \bigg|\frac{\tilde p_k}{\N}-(\mathfrak a+\zeta_1)\bigg|\leq C \N^{-\frac{1}{2}}(\log\N)^2
 \]
 by using the definition of $\tilde p_k$ and inequality $K_1\leq C \N^{\frac{1}{2}}(\log\N)^2$. In addition, observe that $V_h^\textnormal{eff}(x)$ is uniformly Lipschitz in a neighborhood of $a+\zeta_1$ because $\mu_h$ is saturated in such a neighborhood. We conclude that replacement of $\frac{\tilde p_k}{\N}$ by $\mathfrak a+\zeta_1$ leads to an error which can be absorbed into $C K_1 \N^{\frac{1}{2}+\frac{\delta}{2}}$ in \eqref{eq_x106}. Hence, we get \eqref{eq_x67} and finish the proof of Lemma~\ref{Lemma_K_K_prime_saturated}.
\end{proof}

\begin{proof}[Proof of Theorem~\ref{Theorem_ldsaturated}]
We repeat word by word the proof of Theorem~\ref{Theorem_ldpsup} using Lemma
\ref{Lemma_K_K_prime_saturated} instead of Lemma~\ref{Lemma_K_K_prime}. \end{proof}

\chapter{Regularity of the equilibrium measure}
\label{Chapter_smoothness}

In this chapter we investigate the regularity properties of the equilibrium measure, such as the smoothness of its density and its smooth dependence on the parameters of the ensemble. Beyond its intrinsic importance and interest, this study also allows us to derive the asymptotic expansion of the free energy around the equilibrium values of various parameters --- a crucial step for the estimation carried out in Chapter~\ref{Chapter_filling_fractions} of the partition function of the discrete ensembles. In addition, we prove the stability of the off-criticality assumption under small perturbations of the parameters, which is crucial for conditioning procedure in Chapter~\ref{Chapter_conditioning} and for the interpolation procedure in Chapter~\ref{Chapter_partition_functions}.

Although in the most of this book we only use the equilibrium measures that come from the discrete ensembles of Section~\ref{Section_general_model}, it is worth considering \textit{per se} the notion of an equilibrium measure --- \textit{i.e.} a minimizer of the energy functional $-\I$, which can be defined from a smaller amount of data. We call it \emph{variational datum} in order to make a clear distinction with the setting of discrete ensembles, and it will have its own set of imposed assumptions which we describe in the next paragraph. All of them are particular cases or subsets of the assumptions from Section~\ref{Section_list_of_assumptions} that concern directly the variational data. As this list was rather long, we recollect in Section~\ref{allassuml} all the properties needed specifically in our study of variational data and their equilibrium measure.

\section{Setup for equilibrium measures}
\label{Section_parameters}
\subsection{Variational datum}
\label{Section_variational_data}
A variational datum consists of a positive integer $H$ and four groups of parameters.
\begin{enumerate}
\item An $H\times H$ matrix $\boldsymbol{\Theta}$, called matrix of interactions.
\item $2H$ points on the real line $-\infty < \hat{a}'_1< \hat{b}'_1< \hat{a}'_2<\cdots<\hat{a}'_{H}<\hat{b}'_H < +\infty$, defining the segments.
\item $H$ functions $V_h:\,[\hat{a}'_{h},\hat{b}'_{h}]\rightarrow\amsmathbb R$, called potentials.
\item An $H$-tuple of positive real numbers $\hat{\boldsymbol{n}} = (\hat{n}_h)_{h = 1}^H$ called segment filling fractions.
\end{enumerate}
We continue to denote
\[
\amsmathbb{A} = \bigcup_{h = 1}^H [\hat{a}'_{h},\hat{b}'_{h}].
\]

The data defining discrete ensembles in Section~\ref{Section_general_model} naturally gives a variational datum which has a dependence in the parameter $\N$ and must satisfy integrality conditions. The latter are obtained from the fact that a segment $[a_h,b_h]$ contain $\N \hat n_h$ of particles, and by using the formulae of Section~\ref{Section_list_of_assumptions} to express $a_h$ and $b_h$ through $\hat a'_h$ and $\hat b'_h$:
\begin{equation}
\label{eq_integ111} \forall h \in [H]\qquad \N\hat{n}_{h} \in \amsmathbb{Z}_{\geq 0} \quad \textnormal{and} \quad (\N \hat{b}'_{h} - \N \hat{a}'_{h} +\theta_{h,h}- \N \theta_{h,h}\hat{n}_{h}) \in \amsmathbb{Z}_{>0}.
\end{equation}
These integrality conditions are not needed when we define the equilibrium measure in an abstract way as a solution to a minimization problem, hence, they do not play any role for the statements of theorems of this chapter. Yet, they will reappear in the proofs, because we will approximate arbitrary equilibrium measures by equilibrium measures coming from auxiliary discrete ensembles and use the control given by Nekrasov equations (Chapter~\ref{ChapterNekra}) on the latter to study the regularity of the former. From the potential-theoretic point of view, this detour via the discrete ensembles to study variational data may seem artificial. But, we would not know how to prove in full generality the regularity properties established in this chapter without approximations by discrete ensembles.

Of course, we can apply the regularity results of this chapter to the $\N$-dependent equilibrium measures associated to discrete ensembles. In this case, the properties established in this chapter hold for each $\N$ large enough individually --- before studying any $\N \rightarrow \infty$ asymptotics. Besides, many instances the uniformity constants that appear in our results can be chosen independently of $\N$ --- this is usually clear from the context.

\subsection{List of imposed assumptions}
\label{allassuml}

The assumptions that we will consider for variational data will be named by letters, to distinguish them from the assumptions that we considered for discrete ensembles.

\begin{lassum}
\label{Assumption_A}There exists a constant $C > 0$ such that
\begin{enumerate}
\item $\boldsymbol{\Theta}$ is a real symmetric positive semi-definite $H \times H$ matrix;
\item $\forall h \in [H] \qquad \theta_{h,h} \geq \frac{1}{C}$;
\item $H \leq C$ and $|\!|\boldsymbol{\Theta}|\!|_{\infty} \leq C$;
\item $\forall h \in [H] \quad \frac{1}{C} \leq \hat{n}_h \leq \frac{\hat{b}'_{h} - \hat{a}'_{h}}{\theta_{h,h}} - \frac{1}{C}$.
\end{enumerate}
\end{lassum}

Assumption~\ref{Assumption_A} is equivalent to Assumption~\ref{Assumptions_Theta} of Chapter~\ref{Chapter_Setup_and_Examples} in the particular case when $a_1$ and $b_H$ are finite, the $\hat{a}'_h,\hat{b}'_h$ are obtained from $\hat{a}_h,\hat{b}_h$ by the shifts \eqref{eq_shifted_parameters}, and the filling fractions are deterministically fixed to the values $\hat{n}_{h}$. In other words, Assumption~\ref{Assumption_A} corresponds to the case when equations \eqref{eq_equations_eqs} from Section~\ref{DataS} fix all the filling fractions.

\begin{lassum}
\label{Assumption_B} There exists a constant $C > 0$ such that
\begin{enumerate}
\item $\hat{a}'_{1}\geq -C$ and $\hat{b}'_{H} \leq C$;
\item $\forall g,h \in [H] \quad |\hat{b}'_{g} - \hat{a}'_{h}| \geq \frac{1}{C}$;
\item for any $h \in [H]$, there exist nonnegative integers $\iota_{h}^{\pm} \leq C$ and holomorphic functions $U_{h}(z)$ defined for $z$ in an open subset $\amsmathbb{M}_{h} \subset \amsmathbb{C}$ containing the complex $\frac{1}{3C}$-neighborhood of $[\hat{a}'_{h},\hat{b}'_{h}]$, such that
\[
\forall x \in [\hat{a}'_{h},\hat{b}'_{h}] \qquad V_h(x)=\iota_h^{-}\,\mathrm{Llog}(x-\hat a'_h) + \iota_h^+\,\mathrm{Llog}(\hat b'_h-x) +U_h(x).
\]
where we recall the notation $\textnormal{Llog}(x) = x\log |x| - x$.
\item $\forall h \in [H]\quad \forall z \in \amsmathbb{M}_h\quad |U_{h}(z)| \leq C$.
\end{enumerate}
\end{lassum}
Assumption~\ref{Assumption_B} is a subset of Assumptions~\ref{Assumptions_basic} and \ref{Assumptions_analyticity} relevant for the equilibrium measure. The holomorphicity of $U_h$ is taken from Assumption~\ref{Assumptions_analyticity} and the remainder from Assumption~\ref{Assumptions_basic}. Note that the growth conditions of Assumption~\ref{Assumptions_basic} are not needed, since $\hat a'_1$ and $\hat b'_H$ are taken to be finite in variational data.

\begin{definition} \label{Definition_phi_functions_2} For each $h \in [H]$, we define two holomorphic functions $\phi^+_h$ and $\phi_h^-$ in $\amsmathbb{M}_{h}$, which satisfy
\begin{equation}
\label{eq_ratio_to_potential_2} \forall x \in [\hat{a}'_{h},\hat{b}'_{h}]\qquad
\frac{\phi^{+}_{h}(x)}{\phi^{-}_{h}(x)} = e^{-\partial_x V_h(x)}.
\end{equation}
More specifically, we set:
\begin{equation}
\label{eq_phi_plus_minus_def_2}
\forall z \in \amsmathbb{M}_h\qquad \phi^+_h(z):= \big(\hat b'_h - z\big)^{\iota_h^+}\cdot e^{-\partial_zU_h(z)},\qquad
\phi^-_h(z):= \big(z-\hat a'_h\big)^{\iota_h^-}.
\end{equation}
\end{definition}

This definition is a copy of Definition~\ref{Definition_phi_functions}. As mentioned below it, there is some freedom in choosing $\phi^{\pm}_h(z)$ to obey \eqref{eq_ratio_to_potential_2}, as for each $h\in [H]$ both functions can be multiplied by an arbitrary function $f(z)$, which only needs to be holomorphic without zeros $\amsmathbb M_h$. We have removed this freedom by the asymmetric choice \eqref{eq_phi_plus_minus_def_2}. The following useful properties of $\phi^\pm_h(z)$ are implied by the definitions of these functions and Assumption~\ref{Assumption_B}. These are the same properties as those in Lemma~\ref{Lemma_phi_properties}.

\begin{lemma}\label{Lemma_phi_properties_2}
There exist a constant $C > 0$ and a function $c\,:\,\amsmathbb{R}_{> 0} \rightarrow \amsmathbb{R}_{> 0}$ such that, for any $h \in [H]$.
\begin{enumerate}
\item $\forall z \in \amsmathbb{M}_h \quad \big|\phi^{\pm}_h(z)\big| \leq C$;
\item for $\eta > 0$ small enough and $x \in [\hat{a}'_{h} + \eta,\hat{b}'_{h} - \eta]$ we have $|\phi_h^{\pm}(x)| \geq c(\eta)$;
\item if $\iota_h^->0$, then $\phi_h^-(\hat{a}'_{h}) = 0$; if $\iota_h^+>0$, then $\phi_h^+(\hat{b}'_{h}) = 0$.
\end{enumerate}
\end{lemma}

If Assumptions~\ref{Assumption_A} and \ref{Assumption_B} are satisfied, then the conclusion of Proposition~\ref{Lemma_maximizer} and Theorem~\ref{Theorem_equi_charact_repeat_2}, all the results of Section~\ref{Section_Energy_functional}, as well as Lemma~\ref{Lemma_linear_through_distance} are valid --- because Assumptions~\ref{Assumption_A} and \ref{Assumption_B} include all relevant parts of Assumptions~\ref{Assumptions_Theta} and \ref{Assumptions_basic}. This tells us in particular that minus the energy functional
\begin{equation}
\label{eq_functional_general_copy}
\mathcal{I}[\boldsymbol{\mu}] = \sum_{g,h = 1}^{H} \int_{\hat{a}'_{g}}^{\hat{b}'_{g}}\int_{\hat{a}'_{h}}^{\hat{b}'_{h}} \theta_{g,h}\log|x - y| \mu_g(x)\mu_{h}(y)\dd x \dd y-\sum_{h=1}^H \int_{\hat a'_h}^{\hat b'_h} V_h(x)
\mu_h(x)\dd x
\end{equation}
admits a unique minimizer $\boldsymbol{\mu} = (\mu_h)_{h = 1}^{H}$ among $H$-tuples of absolutely continuous nonnegative measures $\mu_h$ with support in $[\hat{a}'_{h},\hat{b}'_{h}]$ and density bounded by $\frac{1}{\theta_{h,h}}$ there, and such that $\mu_h([\hat{a}_{h},\hat{b}_{h}]) = \hat{n}_{h}$, for any $h \in [H]$. The $H$-tuple of constants $\boldsymbol{v}$ appearing in the characterization do not have to satisfy any constraint since all the segment filling fractions were fixed. An equivalent description of $\boldsymbol{\mu}$ is via the measure $\mu = \sum_{h = 1}^H \mu_h$ on $\amsmathbb{A}$. We call $\boldsymbol{\mu}$ or $\mu$ the equilibrium measure of the variational datum. The following functions will play an important role in the study of $\boldsymbol{\mu}$.

\begin{definition}
\label{GQdef} We introduce the Stieltjes transform of the components of the equilibrium measure
\[
\forall h \in [H]\qquad \Gm_{\mu_{h}}(z) = \int_{\hat{a}'_h}^{\hat{b}'_h} \frac{\mu_h(x)\dd x}{z - x}.
\]
and $\Gm_{\mu}(z) = \sum_{h = 1}^{H} \Gm_{\mu_h}(z)$. We also introduce a $2H$-tuple auxiliary functions $(q_h^{\pm}(z))_{h = 1}^{H}$ defined for $\tau \in \{\pm 1\}$ and $z \in \amsmathbb{M}_h \setminus \amsmathbb{A}$ by
\begin{equation}
\label{eq_q_pm}
q_h^{\tau}(z) := \phi_h^{+}(z)\cdot \exp\bigg(\sum_{g = 1}^H \theta_{h,g}\,\Gm_{\mu_{g}}(z)\bigg) + \tau\,\phi_h^{-}(z) \cdot \exp\bigg(- \sum_{g = 1}^H \theta_{h,g}\,\Gm_{\mu_{g}}(z)\bigg).
\end{equation}
\end{definition}

Under Assumptions~\ref{Assumption_A} and \ref{Assumption_B}, we will see in Remark~\ref{Remark_finitely_many_bands} that the equilibrium measure has a finite number of bands, labeled from $1$ to $K$ respecting their increasing order of appearance along the real line. As in Definition~\ref{bandlabel}, we denote $(\alpha_k,\beta_k)$ the $k$-th band, $h^k \in [H]$ the index of the segment in which it is included, and $\llbracket k^-(h),k^+(h)\rrbracket \subseteq [K]$ the set of indices of the bands contained in the $h$-th segment. We also recall that $\amsmathbb{V}_h$, $\amsmathbb{B}_h$ and $\amsmathbb{S}_h$ refer respectively to the voids, the bands and the saturations of $\mu$ in $[\hat{a}_h',\hat{b}_h']$, \textit{cf.} Definition~\ref{Definition_void_saturated}.

\begin{definition}
\label{GQdef2} For each $h \in [H]$, we introduce the function
\begin{equation}
\label{eq_sigma_h_definition}
\sigma_h(z)= \prod_{k = k^-(h)}^{k^+(h)} \sqrt{(z - \alpha_k)(z - \beta_k)}
 \end{equation}
 for $z$ in $\amsmathbb{C}\setminus \bigcup_{k = k^-(h)}^{k^+(h)} [\alpha_k,\beta_k]$, where we use the branch of the square root such that $\sigma_h(z)$ is a holomorphic function in this domain and such that $\sigma_h(z) \sim z^{k^+(h) - k^-(h) + 1}$ as $z \rightarrow \infty$. Next, we define
\[
\psi_h(z) :=\psi_h^+(z) \cdot \psi_h^-(z),
\]
where
\begin{equation}
\label{psihdefdef}
 \psi_h^-(z) =\begin{cases} z-\hat a'_h & \textnormal{if }\iota_h^- \geq 2 \textnormal{ and } \hat a'_h\in \amsmathbb{S}_h,\\ 1 & \textnormal{otherwise}. \end{cases}\qquad
 \psi_h^+(z) =\begin{cases} z-\hat b'_h & \textnormal{if }\iota_h^+ \geq 2 \textnormal{ and } \hat b'_h\in \amsmathbb{S}_h,\\ 1 & \textnormal{otherwise}. \end{cases}
\end{equation}
We finally set
\begin{equation}
\label{eq_s_h_definition}
 s_h(z) := \frac{q_h^{-}(z)}{\sigma_h(z) \cdot \psi_h(z)}.
\end{equation}
\end{definition}
\begin{remark} \label{remextrasssh}
Suppose that we have  either $\iota_h^{-} = 2$ if $\hat{a}_h'$ is saturated, or $\iota_h^- = 0$, and either $\iota_h^+ = 2$ if $\hat{b}_h'$ is saturated, or $\iota_h^+ = 0$ (in the context of discrete ensembles, such conditions appear in Assumption~\ref{Assumptions_extra}). These properties lead to a simplified definition
\[
\psi_h(z) = (z - \hat{a}_h')^{\mathbbm{1}_{\amsmathbb{S}_h}(\hat{a}'_h)} \cdot (z - \hat{b}'_h)^{\mathbbm{1}_{\amsmathbb{S}_h}(\hat{b}'_h)}.
\]
Then, $s_h(z)$ from \eqref{eq_s_h_definition} matches the definition of $s_h(z)$ given in \eqref{eq_s_first_appearance}.
\end{remark}

From the definitions, $\Gm_{h}(z)$ is \textit{a priori} holomorphic for $z$ in $\amsmathbb{C} \setminus (\amsmathbb{B}_h \cap \amsmathbb{S}_h)$, while $q_h^{\pm}(z)$ and $s_h(z)$ are holomorphic functions of $z$ in $\amsmathbb{C} \setminus (\amsmathbb{B}_h \cap \amsmathbb{S}_h)$. As soon as we know both $q_h^+(z)$ and $q_h^-(z)$, we can reconstruct $\Gm_{\mu_{h}}(z)$ and, hence, the equilibrium measure itself. The properties of these two functions are very different: as we will see, the function $q_h^+(z)$ is holomorphic in $\amsmathbb M_h$, while $q_h^-(z)$ is not --- it has discontinuities. On the other hand, after we divide by $\sigma_h(z)$, this discontinuity disappears and $s_h(z)$ is again holomorphic in $\amsmathbb{M}_h$. The additional division by $\psi_h(z)$ is introduced in the definition of $s_h(z)$ in order to get rid of possible zeros of $q_h^-(z)$ at the endpoints $\hat a'_h$, $\hat b'_h$. In fact $q_h^+(z)$ can also have zeros at the same points. This seemingly technical modification is crucial for the arguments developed in Chapter~\ref{Chapter_fff_expansions}, which require the regularity of $\frac{1}{s_h(z)}$ in $\amsmathbb M_h$.

\medskip

We say that the equilibrium measure is off-critical if it satisfies the following assumption, copied from Assumption~\ref{Assumptions_offcrit}. Recall the definition of the effective potential:
\begin{equation}
\label{eq_V_eff_copy}
V^{\textnormal{eff}}_{h}(x) = V_h(x) - \sum_{g= 1}^H 2\theta_{h,g}\int_{\hat{a}'_{g}}^{\hat{b}'_{g}} \log|x - y|\,\mu_{g}(y)\dd y.
\end{equation}

\begin{lassum}\label{Assumption_C}
 There exist $C > 0$ and a continuous function $\eps: \,\amsmathbb{R}_{\geq 0} \rightarrow \amsmathbb{R}_{\geq 0}$ such that $\eps(0) = 0$ and $\eps|_{\amsmathbb{R}_{> 0}} > 0$ making the following properties hold for any $h \in [H]$.
 \begin{enumerate}
 \item The segment $[\hat a'_h, \hat b'_h]$ contains at least one band. The bands in this segment are at distance at least $\frac{1}{C}$ of the endpoints. The total number of bands, voids and saturations of $\boldsymbol{\mu}$ is bounded by $C$, while the length of each of these regions is larger than $\frac{1}{C}$.
 \item $V^{{\textnormal{eff}}}_h(x) \geq v_h+\eps(\delta)$ for any $x \in \amsmathbb{V}_h$ at distance at least $\delta$ from a boundary point of $\amsmathbb{V}_h$ which is neither $\hat a'_h$ nor $\hat b'_h$.
 \item $V^{{\textnormal{eff}}}_h(x) \leq v_h-\eps(\delta)$ for any $x \in \amsmathbb{S}_h$ at distance at least $\delta$ from a boundary point of $\amsmathbb{S}_h$ which is neither $\hat a'_h$ nor $\hat b'_h$.
  \item For any $x \in \amsmathbb{B}_h$ at a distance at least $\delta$ from the boundary of $\amsmathbb{B}_h$, we have
  \[
  \eps(\delta) \leq \mu_h(x) \leq \frac{1}{\theta_{h,h}} - \varepsilon(\delta).
  \]
 \item For any point $x_0$ separating a void from a band in $[\hat a'_h, \hat b'_h]$, we have
 \[{\frac{1}{C}\sqrt{|x-x_0|} \leq \mu_h(x) \leq C\sqrt{|x-x_0|}}\] for any $x$ in this band at distance at most $\frac{1}{C}$ from $x_0$.
 \item For any point $x_0$ separating a saturation from a band in $[\hat a'_h, \hat
 b'_h]$ we have \[{\frac{1}{C}\sqrt{|x-x_0|} \leq \frac{1}{\theta_{h,h}} -\mu_h(x) \leq C\sqrt{|x-x_0|}}\] for any $x$
in this band at distance at most $\frac{1}{C}$ from $x_0$.
 \end{enumerate}
\end{lassum}

\section{Spatial regularity: results}
\label{RegSecReg}
The main result of this section is the smoothness of the density of the equilibrium measure and analytic properties of the auxiliary functions $q^\pm_h(z)$. The proofs are postponed to Section~\ref{ProofRegSecReg}.

\begin{theorem} \label{Theorem_regularity_density}
Let $\boldsymbol{\mu}$ be the equilibrium measure of a variational datum obeying Assumptions~\ref{Assumption_A} and \ref{Assumption_B}. Then the following properties hold for any $h \in [H]$.
 \begin{enumerate}
 \item[(i)] The density $\mu_h(x)$ is an $\frac{1}{2}$-H\"older continuous function of $x$ everywhere except, possibly, at $x \in \{\hat{a}_h',\hat{b}_h'\}$.
 \item[(ii)] For any $x\in\amsmathbb R\setminus \{\hat{a}_h',\hat{b}_h'\}$, the limits
 \[
 \Gm_{\mu_h}(x^{\pm}) = \lim_{\epsilon \rightarrow 0^+} \Gm_{\mu_h}(x \pm {\ii}\epsilon)
 \]
 exist. Furthermore, we have $\textnormal{Im}(\Gm_{\mu_h}(x^{\pm}))=\mp \pi \mu_h(x)$.
 \item[(iii)] The function $q_h^+(z)$ from Definition~\ref{GQdef} is meromorphic for $z \in \amsmathbb M_h$. Its only possible singularities are simple poles at $\hat a'_h$ (respectively $\hat b'_h$), which can appear only if $\iota^-_h=0$ (respectively $\iota^+_h=0$).
 \item[(iv)] The function $ (q_h^-(z))^2$ is meromorphic in $\amsmathbb M_h$. Its only possible singularities in $\amsmathbb{M}_h$ are simple or double poles at $\hat a'_h$ (respectively $\hat b'_h$), which can appear only if $\iota^-_h=0$ (respectively $\iota^+_h=0$). Besides, the endpoints of a band of $\mu_h$ in $(\hat a'_h,\hat b'_h)$ is necessarily a zero of $(q_h^-(z))^2$.
 \item[(v)] The function $s_h(z)$ from Definition~\ref{GQdef2} is meromorphic for $z \in \amsmathbb M_h$. Its only possible singularities in $\amsmathbb{M}_h$ are simple poles at $\hat a'_h$ (respectively $\hat b'_h$), which can appear only if $\iota^-_h=0$ (respectively $\iota^+_h=0$).
 \end{enumerate}
\end{theorem}
\begin{remark} \label{Remark_finitely_many_bands}
 Property (iv) implies that the equilibrium measure has finitely many bands, since a meromorphic function (which is not identically equal to $0$) must have finitely many zeros on a finite interval. Since $\mu_h$ has continuous density, it implies that there are also finitely many voids and saturations. However, we do not claim that the total number of bands is bounded by a constant depending only on the constants in the assumptions. Later, we will make use of such a uniform boundedness for the number of bands, voids and saturations, explaining why we include it in Assumption~\ref{Assumption_C} 1. --- see also Assumption~\ref{Assumptions_offcrit} 1. in the context of discrete ensembles.
\end{remark}

Assumptions~\ref{Assumption_A} and \ref{Assumption_B} include finiteness of $\hat a'_1$, $\hat b'_H$ and deterministically fixed filling fractions, while the setting of Section~\ref{Section_list_of_assumptions} did not require it. Yet, the conclusion of Theorem~\ref{Theorem_regularity_density} remains true in the more general setting of Section~\ref{Section_list_of_assumptions}. Indeed, first, the equilibrium measure necessarily has a compact support by Theorem~\ref{Theorem_equi_charact_repeat_2}, hence, we can always restrict our attention to finite segments. Second, each equilibrium measure which is a minimizer of the energy functional without fixing filling fractions is simultaneously a minimizer for the same energy functional with filling fractions fixed to their optimal values.

We complement Theorem~\ref{Theorem_regularity_density} with two auxiliary statements. The first one expresses the density of the equilibrium measure through the functions $q^{\pm}_h(z)$ while the second one provides uniform upper and lower bounds on $q^{\pm}_h(z)$. We recall the notation
\[
f(x^{\pm})=\lim\limits_{\epsilon \rightarrow 0^+}f(x\pm \ii \epsilon),
\]
whenever the right-hand side is well-defined.
\begin{proposition}
\label{Proposition_density} Under Assumptions~\ref{Assumption_A} and \ref{Assumption_B}, the density of the equilibrium measure can be computed for $h \in [H]$ through
\begin{equation}\label{eq_density_tan}
 \forall x \in (\hat{a}'_{h},\hat{b}'_{h})\qquad\tan\big(\pi\theta_{h,h}\,\mu_h(x)\big)=\frac{ q^{-}_h(x^-) - q^{-}_h(x^+)}{2\ii \, q^+_h(x)}.
\end{equation}
For $x$ inside a band, the density is also given by
\begin{equation}
 \label{eq_density_exponent}
1-\exp\big(2\ii\pi \theta_{h,h}\, \mu_h(x)\big)=\frac{2\, q^{-}_h(x^+)}{q^{+}_h(x)+q^{-}_h(x^+)}.
\end{equation}
\end{proposition}

\begin{proposition}
\label{Proposition_q_bounds}
 Consider a variational datum satisfying Assumptions~\ref{Assumption_A}, \ref{Assumption_B} and \ref{Assumption_C}. There exist \mbox{$\eps,C>0$} which depend only on the constants in the assumptions such that for any $h\in[H]$ and with the notation $\amsmathbb{B}_h^\eps$ for the complex $\eps$-neighborhood of the bands of $\mu_h$,
\begin{equation}
\label{eq_q_bounds}
 \sup_{\textnormal{dist}(z,\amsmathbb{B}_h) \leq \varepsilon} |q^{\pm}_h(z)| \leq C, \qquad \sup_{\textnormal{dist}(z,\amsmathbb{B}_h) \leq \varepsilon} |s_h(z)| \leq C,
\qquad
\inf_{\textnormal{dist}(z,\amsmathbb{B}_h) \leq \varepsilon} |s_h(z)| \geq  \frac{1}{C}.
\end{equation}
\end{proposition}
In fact, the first two bounds of \eqref{eq_q_bounds} require only Assumptions~\ref{Assumption_A} and \ref{Assumption_B}, while for the validity of the third bound we have to require Assumption~\ref{Assumption_C}. Indeed, \eqref{eq_density_exponent} shows that if the density $\mu_h(x)$ is close to $0$ in an interior point of a band, then so is $s_h(x)$.

Theorem~\ref{Theorem_regularity_density}, Proposition~\ref{Proposition_density} and the finiteness of the number of bands, voids and saturations (the so-called finite-gap ansatz) are essentially known if the constraints $\mu_{h} \leq \frac{1}{\theta_{h,h}}$ are absent and the potential has no singularity (\textit{i.e.} $\iota_h^{-} = \iota_h^+ = 0$) \cite{DeiftKMCL}. They are also known in presence of the upper-constraint on the density if $H = 1$ and the potential has no singularity \cite{Kuij1}.

In absence of upper constraint on the density, such results were also obtained in \cite{BGK}.  In particular, Theorem~\ref{Theorem_regularity_density} (i) and (ii) were derived in the course of studying a random matrix ensemble for which $\mu$ represents the spectral measure in the limit of large sizes. Large deviation theory in the random matrix ensemble shows that the average of the spectral measure approximates the equilibrium measure, while the Dyson--Schwinger equation which is exact for finite sizes allows some control in the large size limit, sufficient to derive the desired properties for $\mu$. The limiting form of this equation is an analogue of property (iii) of Theorem~\ref{Theorem_regularity_density}.

We use in Section~\ref{ProofRegSecReg} a similar approach to prove (iii) and approximate $\mu$ by the equilibrium measure of an auxiliary discrete ensemble which we control well thanks to the Nekrasov equations from Chapter~\ref{ChapterNekra}. Yet, in presence of the constraint $\mu_h(x) \leq \frac{1}{\theta_{h,h}}$, the situation is more complicated. The main reason is that an analog of the random matrix ensemble featuring this constraint is a discrete ensemble whose definition imposes integrality conditions like \eqref{eq_integ111}. They force all the parameters of Section~\ref{Section_parameters} to depend on the total number of particles, here driven by $\N$. Hence, it becomes impossible to make approximations with discrete ensembles, while keeping the parameters such as endpoints $\hat a'_h$, $\hat b'_h$ fixed. Therefore, we need to proceed in two steps. First, in Lemma~\ref{Lemma_regularity_one_band} we provide a discrete approximation argument in the special case $H=1$, $\theta_{1,1}=1$, and rational $\hat b'_1-\hat a'_1$ and $\hat n$. Next, we add a second layer of approximations by showing that the properties of Theorem~\ref{Theorem_regularity_density} are stable under limits in families of variational data --- we do this in Lemma~\ref{Lemma_continuity_in_parameters}. Eventually, both lemmata are combined to finish the proof of Theorem~\ref{Theorem_regularity_density}.

Note that if the properties (i) and (ii) of Theorem~\ref{Theorem_regularity_density} were known, then (iii), (iv), and (v) could be deduced from the characterization of the equilibrium measure of Theorem~\ref{Theorem_equi_charact_repeat_2} in a straightforward way. However, we do not have in the desired generality\footnote{It is conceivable that such an argument can be reached by generalizing ideas of \cite{DeiftKMCL,Kuij1}.} a direct argument proving (i) and (ii), and therefore, we proceed in an opposite direction by first showing that (iii) holds and then deducing (i) and (ii) as a corollary. This is the reason why we need to rely on the Nekrasov equations of some well-chosen families of discrete ensembles and take limits within these families.

We conclude this section with a useful result: for strictly convex potentials, parts of the off-criticality Assumption~\ref{Assumption_C} hold automatically. For polynomial potentials with a single group of particles and unconstrained equilibrium measures (\textit{i.e.} minimizers of the energy functional $-\mathcal{I}$ without constraints on the maximal density) this is known by an argument of  \cite{BoPaSh} or \cite{Johansson}, and we extend the validity of this argument to our context.

\begin{proposition}
\label{prop:convexmueq}
Consider a variational datum satisfying Assumptions \ref{Assumption_A} and \ref{Assumption_B}. Suppose that for any $g,h \in [H]$, we have $\theta_{g,h} \geq 0$ and there exists a constant $c > 0$ such that $V_h''(x) > c > 0$ for any $x \in (\hat{a}_h',\hat{b}_h')$. Then, for each $h \in [H]$ either the $h$-th segment is void or saturated, or the following properties simultaneously hold:
\begin{itemize}
\item[(i)] $\mu_h$ has at least one band;
\item[(ii)] the support of $\mu_h$ is connected;
\item[(iii)] if there is a void, it cannot be reduced to an interior point of the $h$-th segment;
\item[(iv)] the property 2., the lower bound in 4., and the property 5. of Assumption~\ref{Assumption_C} hold.
\end{itemize}
The constants in (iv) depend on the constants in the assumptions and on a positive lower bound for the length of the bands.
\end{proposition}

\section{Spatial regularity: proofs}

\label{ProofRegSecReg}

We turn to the proofs of Theorem~\ref{Theorem_regularity_density} and Propositions~\ref{Proposition_density}, \ref{Proposition_q_bounds} and \ref{prop:convexmueq}.

\begin{lemma} \label{Lemma_regularity_one_band} Theorem~\ref{Theorem_regularity_density} without (vi) holds when $H=1$, $\theta_{1,1}=1$, $(\hat b'_1-\hat a'_1)\in\amsmathbb Q$ and $\hat n_1 \in\amsmathbb Q$.
\end{lemma}
\begin{proof}
Consider a variational datum with $H=1$, specified by intensity of interaction $\theta > 0$, endpoints $\hat{a}',\hat{b}'$, potential $V$ and filling fraction $\hat{n}$. We omit the subscripts throughout the proof, \textit{i.e.} $\hat{a}'_1 = \hat a'$, $\hat b_1'=\hat b'$, \textit{etc}. From Assumption~\ref{Assumption_B} we also have nonnegative integers $\iota^{\pm}$ and a complex domain $\amsmathbb{M}$. The rationality assumption \eqref{eq_integ111} and $\theta =1$ imply the existence of three sequences of positive integers $N^{(m)}$, $M^{(m)}$, $\N^{(m)}$ indexed by $m \in \amsmathbb{Z}_{> 0}$ and tending to $+\infty$ as $m \rightarrow \infty$, such that
\[
\forall m \in \amsmathbb{Z}_{> 0}\qquad \hat{n} = \frac{N^{(m)}}{\N^{(m)}} \qquad \textnormal{and} \qquad \hat{b} - \hat{a} = \frac{N^{(m)} + M^{(m)}}{\N^{(m)}}.
\]
Eventually $m$ will go to infinity. Let us construct an auxiliary discrete ensemble, with $H=1$, $N^{(m)}$ particles and master parameter $\N^{(m)}$. The expectation value corresponding to its law will be denoted $\E^{(m)}$. The defining segment for this discrete ensemble is $[a^{(m)},b^{(m)}]$ with
\begin{equation*}
\begin{split}
& a^{(m)} = \N^{(m)}\hat{a}^{(m)},\qquad \quad \! b^{(m)} = \N^{(m)}\hat{b}^{(m)}, \\
& \hat a^{(m)} = \hat a' + \frac{1}{2\N^{(m)}},\qquad \hat{b}^{(m)} = \hat b' - \frac{1}{2\N^{(m)}},
\end{split}
\end{equation*}
which takes into account the shift \eqref{eq_shifted_parameters}, and the weight is
\[
w^{(m)}(\ell) = \frac{\exp\Big[-\N^{(m)}\, U\big(\frac{\ell}{\N^{(m)}}\big)\Big]}{\bigg[\Gamma\big(\ell - \N^{(m)}\hat{a} + 1\big) \left(\dfrac{e}{\N^{(m)}}\right)^{\ell - \N^{(m)}\hat{a} + 1}\bigg]^{\iota^{-}}\!\!\! \cdot \bigg[\Gamma\big(\N^{(m)}\hat{b} - \ell + 1\big) \left(\dfrac{e}{\N^{(m)}}\right)^{\N^{(m)}\hat{b} - \ell + 1}\bigg]^{\iota^+}}.
\]
where $e = \exp(1)$. Using the Stirling approximation in the form
\[
\log \Gamma(x+c) \mathop{=}_{x \rightarrow +\infty} x\log x -x + O(\log x)
\]
for a fixed $c$, we deduce that the auxiliary ensemble satisfies Assumption~\ref{Assumptions_basic} with $\N=\N^{(m)}$ and same potential $V$ as in the variational datum. Hence, the equilibrium measure of the auxiliary discrete ensemble is the same as the equilibrium measure of Theorem~\ref{Theorem_regularity_density} and does not depend on $m$. We also find that the auxiliary ensemble satisfies Assumption~\ref{Assumptions_analyticity} with the functions of Definition~\ref{Definition_phi_functions} given by
\begin{equation*}
\begin{split}
\Phi^{+,(m)}(z) = & \big(\hat{b}' - z\big)^{\iota^+}\cdot\exp\Bigg(-\N^{(m)}\bigg[U\bigg(z + \frac{1}{2\N^{(m)}}\bigg) - U\bigg(z - \frac{1}{2\N^{(m)}}\bigg)\bigg]\Bigg), \\
\Phi^{-,(m)}(z) = & \big(z - \hat{a}'\big)^{\iota^-}.
\end{split}
\end{equation*}

Then, the Nekrasov equation of Theorem~\ref{Theorem_Nekrasov} states that for random $(\ell_1,\ldots, \ell_{N^{(m)}})$ under the law of the auxiliary discrete ensemble, the function
\begin{equation}
\label{eq_NekAux} R^{(m)}(z):= \sum_{\tau \in \{\pm 1\}} \Phi^{\tau,(m)}(z)\cdot \E^{(m)}\Bigg[\prod_{i = 1}^{N^{(m)}} \bigg(1 + \frac{\tau}{\N^{(m)}z - \ell_i - \frac{\tau}{2}}\bigg)\Bigg]
\end{equation}
is meromorphic for $z \in \amsmathbb M$, and its only possible singularities are simple poles at $z = \hat{a}'$ if $\iota^{-} = 0$, and at $z = \hat{b}'$ if $\iota^{+} = 0$. Let $\gamma \subset \amsmathbb{M}$ be a contour surrounding $[\hat{a}',\hat{b}']$ as described in Definition~\ref{def_contourgammah}. Applying Cauchy integral formula to the holomorphic function $(z - \hat{a}')(z - \hat{b}') \cdot R^{(m)}(z)$, we write for $z$ inside $\gamma$:
\begin{equation}
\label{eq_RCauchy} R^{(m)}(z) = \oint_{\gamma} \frac{\dd\zeta}{2\ii\pi}\, \frac{(\zeta - \hat{a}')(\zeta - \hat{b}')}{(z - \hat{a}')(z - \hat{b}')}\cdot \frac{R^{(m)}(\zeta)}{\zeta - z}.
\end{equation}
Besides, it follows from Lemma~\ref{Lemma_tail_bound_general} (\textit{cf.} also Corollary~\ref{Corollary_a_priory_1}) that we have for $\tau \in \{\pm 1\}$
\[
\lim_{m \rightarrow \infty} \E^{(m)}\bigg[\prod_{i = 1}^{N^{(m)}} \bigg(1 + \frac{\tau}{\N^{(m)}\zeta - \ell_i - \frac{\tau}{2}}\bigg)\bigg] = \exp\big(\tau\,\Gm_{\mu}(\zeta)\big).
\]
and the convergence is uniform over $\zeta\in\gamma$. On the other hand we have the uniform convergence
\[
\lim_{m \rightarrow \infty} \Phi^{\tau,(m)}(z) = \phi^{\tau}(z)
\]
for $z\in\amsmathbb M$, with $\phi^\tau$ as in Definition~\ref{Definition_phi_functions_2}:
\begin{equation}
\label{eq_phifdsf}\phi^{+}(z) = (\hat{b}' - z)^{\iota^{+}}\cdot e^{-\partial_z U(z)},\qquad \phi^{-}(z) = (z - \hat{a}')^{\iota^{-}}.
\end{equation}
We then deduce from \eqref{eq_NekAux} that $R^{(m)}(z)$ has a uniform limit as $m \rightarrow \infty$ on the contour $\gamma$, and, therefore, by \eqref{eq_RCauchy} also everywhere inside $\gamma$. Let $S(z)$ denote the limiting function: $S(z):=\lim_{m\rightarrow\infty} R^{(m)}(z)$. Since $R^{(m)}(z)$ is meromorphic in $\amsmathbb M$ with only possible singularity simple poles at $z = \hat{a}'$ if $\iota^-=0$ and at $z = \hat{b}'$ if $\iota^+=0$,
so is $S(z)$ by taking the limit of the contour integral \eqref{eq_RCauchy} with an appropriate choice of $\gamma$. Simultaneously, for $z\in\amsmathbb M\setminus[\hat a',\hat b']$, we have by taking the limit of \eqref{eq_NekAux}
\begin{equation}
\label{eq_limitpsiXi}\phi^+(z) \cdot \exp\big(\Gm_{\mu}(z)\big) + \phi^{-}(z) \cdot \exp\big(-\Gm_{\mu}(z)\big) = S(z).
\end{equation}
We conclude that $S(z)$ coincides with $q^{+}(z)$, thus, proving Theorem~\ref{Theorem_regularity_density}, (iii).

Let us treat \eqref{eq_limitpsiXi} as an equation on an unknown function $\Gm_{\mu}(z)$:
\begin{equation}
\label{eq_Stiel_equation}\phi^+(z)\cdot \Xi(z) + \frac{\phi^{-}(z)}{\Xi(z)} = S(z),\qquad \Xi(z) := \exp\big(\Gm_{\mu}(z)\big).
\end{equation}
The solutions of \eqref{eq_Stiel_equation} are, for all $z \in \amsmathbb M \setminus [\hat{a}',\hat{b}']$,
\begin{equation}
\label{eq_Xieq} \Xi^{\pm}(z) = \frac{S(z) \pm \sqrt{\Delta(z)}}{2\phi^{+}(z)}\qquad \textnormal{with}\quad \Delta(z) = (S(z))^2- 4\phi^{-}(z)\phi^{+}(z).
\end{equation}
Since $\mathcal G_\mu(z)$ is a holomorphic function of $z$ outside the support of $\mu$, the sign $\pm$ cannot depend on $z$. Since $\phi^+(z)$ defined in \eqref{eq_phifdsf} has no zeros for $z \in (\hat{a}',\hat{b}')$, we deduce that the limits
\[
\Xi(x^{\pm}) = \lim_{\epsilon \rightarrow 0^+} \Xi(x \pm {\ii}\epsilon)
\]
exist for $x \in (\hat{a}',\hat{b}')$. Note that the two solutions, $\Xi^{\pm}(z)$ in \eqref{eq_Xieq} satisfy the relation
\begin{equation}
\label{eq_Xiinvrel}\Xi^{+}(z)\cdot \Xi^{-}(z) = \frac{\phi^{-}(z)}{\phi^{+}(z)} = \frac{ (\hat{b}' - z)^{\iota^{+}}}{ (z - \hat{a}')^{\iota^{-}}} e^{-\partial_z U(z)},
,
\end{equation}
and, thus, $\Xi(x^{\pm})\neq 0$ for $x\in (\hat a',\hat b')$. This proves the first part of Theorem~\ref{Theorem_regularity_density}, (ii), \textit{i.e.}\ that $\Gm_\mu(z)=\log \Xi(z)$ has finite limits on the real line. Note that the equality $\textnormal{Im}\big(\Gm_{\mu}(x^\pm)\big) = \mp \pi \mu(x)$ (second part of Theorem~\ref{Theorem_regularity_density}, (ii)) is a general property of the Stieltjes transform of any measure.

Since $S(z)$, $\phi^{\pm}(z)$, $\Delta(z)$ are holomorphic in $\amsmathbb M$, \eqref{eq_Xieq} implies that $\Xi(x^{\pm})$ is $\frac{1}{2}$-H\"older continuous inside $(\hat a',\hat b')$, and therefore, so is $\pi \mu(x)= \textnormal{Im}(\log \Xi(x^-))$, proving Theorem~\ref{Theorem_regularity_density}, (i).

Next, we prove part (iv) in Theorem~\ref{Theorem_regularity_density}. From \eqref{eq_limitpsiXi} and \eqref{eq_Xieq} we have
\[
 \Delta(z)=\Big(\phi^+(z)\cdot\exp\big(\Gm_{\mu}(z)\big) - \phi^{-}(z)\cdot\exp\big(-\Gm_{\mu}(z)\big)\Big)^2.
\]
We conclude that $ (q^{-}(z))^2= \Delta(z)$. Since $\Delta(z)=(S(z))^2- 4\phi^{-}(z)\phi^{+}(z)$ is holomorphic in $\amsmathbb M$ except for the possible double poles at $\hat a'$ and $\hat b'$, we conclude that so is $(q^-(z))^2$. It remains to show that all the endpoints of the bands of $\mu$ are necessary zeros of $\Delta(z)$.

Note that $\Xi(z)=\exp\big(\Gm_\mu(z)\big)$ defined for $z$ in the complex domain $\amsmathbb{M}$ is continuous near a point $z = x \in\amsmathbb R$ whenever $x$ belongs to a void or a saturation of $\mu$ --- indeed, $\Gm_\mu(z)$ makes an additive jump by $2\ii\pi \mu(x)=2\ii\pi$ when crossing the real axis through the saturation, which gives a factor $1$ upon exponentiation. On the other hand, $\Xi(z)$ makes a jump whenever $z$ crosses the real axis through a band.

Let us now compare with \eqref{eq_Xieq}. It identifies $\Xi(z)$ with a branch of a bi-valued holomorphic function. The above argument shows that any endpoint of a band of $\mu$ is a ramification point for this function, as the value of $\Xi(z)$ is changed when we loop around and make a jump crossing the real axis on one side, but not on the other side. On the other hand, since $S(z)$ and $\phi^+(z)$ are holomorphic in $\amsmathbb M$, the ramification points of \eqref{eq_Xieq} are those of $\sqrt{\Delta(z)}$, which are zeros of $\Delta(z)$.

\medskip

The final step is to prove Theorem~\ref{Theorem_regularity_density}, (v). For this purpose, we study the jump of the function $s(z)$ from Definition~\ref{GQdef2} as $z$ crosses the real axis at a point $x$. There are several cases to consider.
\begin{enumerate}
 \item If $x$ is outside $[\hat a',\hat b']$ or if $x$ is inside a void of $\mu$ in $(\hat a',\hat b')$, then $\Gm_{\mu}$ is holomorphic near the point $x$. Therefore, so is $q^{-}(z)$. Since $\sigma(z)$ is holomorphic outside bands and $\psi(z)$ is holomorphic everywhere, we conclude that $s(z)$ is holomorphic near $z = x$.
 \item If $x$ is inside a saturation of $\mu$, then the limits of $\exp\big(\Gm_{\mu}(z)\big)$ as $z$ approaches $x$ from the upper and from the lower half-planes coincide since the imaginary parts of $\Gm_{\mu}$ differ by a multiple of $2\pi$.
 Hence, $\exp\big(\Gm_{\mu}(z)\big)$ is holomorphic near $x$, and we again conclude that $s(z)$ is holomorphic near $z = x$.
 \item If $x$ is inside a band, then we compute:
\begin{equation}
\label{eq_x171}
 s(x^+)= \frac{\phi^+(x) \cdot \exp\big(\Gm_{\mu}(x)-\ii \pi \mu(x)\big)-\phi^-(x)\cdot \exp\big(-\Gm_\mu(x)+\ii\pi\mu(x)\big)}{\sigma(x^+) \cdot \psi(x)},
\end{equation}
\begin{equation}
\label{eq_x172}
 s(x^-)= \frac{\phi^+(x) \cdot \exp\big(\Gm_{\mu}(x)+\ii \pi \mu(x)\big)-\phi^-(x) \cdot \exp\big(-\Gm_\mu(x)-\ii\pi\mu(x)\big)}{\sigma(x^-) \cdot \psi(x)},
\end{equation}
where $\Gm_\mu(x)$ is given by the principal value of the integral, which is well-defined since $\mu(x)$ is $\frac{1}{2}$-H\"older continuous.
We claim that \eqref{eq_x171} and \eqref{eq_x172} coincide. To see this, we use the properties $\sigma(x^+)=-\sigma(x^-)$ of the square root, $\frac{\phi^+(x)}{\phi^-(x)}= e^{-\partial_xV(x)}$ given in \eqref{eq_ratio_to_potential_2}, and
\[
\partial_x V(x)= \Gm_\mu(x^+) + \Gm_{\mu}(x^-)
\]
obtained by differentiating with respect to $x$ the characterization of $\mu$ of Theorem~\ref{Theorem_equi_charact_repeat_2}. This differentiation is a legitimate operation, since we already proved the $\frac{1}{2}$-H\"older continuity of the density. We conclude that $s(z)$ has the same limits as $z$ approaches $x$ from the upper and lower half-planes and, therefore, $s(z)$ is holomorphic near $z = x$.
\end{enumerate}

We have shown that $s(z)$ is holomorphic everywhere except, perhaps, at the endpoints of the bands or at the points $\hat a'$ and $\hat b'$. In fact, the function $s(z)$ is bounded near each of the endpoints of bands, since each of them is a simple zero of $(\sigma(z))^2$ and simultaneously a zero of $(q^-(z))^2$ by (iv). We conclude that the singularities are removable and $s(z)$ is necessary holomorphic near all endpoints of bands. Eventually, for the segment endpoints $\hat a'$ and $\hat b'$, the previous arguments show that $s(z) \cdot \psi(z)$ is either holomorphic or has a simple pole at $\hat a'$ (the same for $\hat b'$). It remains to observe that whenever $\psi^-(z)=z-\hat a'$, \textit{i.e.} if $\iota^- > 1$ and $\hat{a}'$ belongs to a saturation, both terms in the definition of $q^-(z)$ vanish at $z=\hat a'$ and, hence, the ratio $\frac{q^-(z)}{z-\hat a'}$ is regular at $z = \hat a'$. Similarly, whenever $\psi^+(z)=z-\hat b'$, \textit{i.e.} if $\iota^+ > 1$ and $\hat{b}'$ belongs to a saturation, there is no singularity of $s(z)$ at $z = \hat b'$.
\end{proof}

The general case of Theorem~\ref{Theorem_regularity_density} is deduced from Lemma~\ref{Lemma_regularity_one_band} with the help of an approximation lemma.

\begin{lemma} \label{Lemma_continuity_in_parameters}
Assume that we have a sequence of variational data, indexed by $m \in \amsmathbb{Z}_{> 0}$ in such a way that Assumptions~\ref{Assumption_A} and \ref{Assumption_B} hold for each $m$ with $m$-independent number of segments $H$, complex domains $\amsmathbb{M}_1,\ldots,\amsmathbb{M}_H$ and integers $\iota_1^{\pm},\ldots,\iota_H^{\pm}$. Let $\boldsymbol{\mu}^{(m)}$ denote the corresponding equilibrium measure. In addition, assume that the variational data converge as $m\rightarrow\infty$ to limiting values which still satisfying Assumptions~\ref{Assumption_A} and \ref{Assumption_B} --- for the potential $U_h(x)$ we use the topology of the uniform convergence in the $m$-independent sets $\amsmathbb{M}_h$, for $h \in [H]$. If the five conclusions of Theorem~\ref{Theorem_regularity_density} hold for $\boldsymbol{\mu}^{(m)}$ for any $m \in \amsmathbb{Z}_{> 0}$, then the density of $\mu^{(m)}$ converges to the density of a limit measure $\mu$ uniformly on any compact and the five conclusions of Theorem~\ref{Theorem_regularity_density} hold for $\boldsymbol{\mu}$ as well.
\end{lemma}

In Section~\ref{Section-ParReg} we are going to show that under the off-criticality Assumption~\ref{Assumption_C} the dependence of the equilibrium measure on the parameters is in fact differentiable. However, some of the arguments in that section use Theorem~\ref{Theorem_regularity_density} as an ingredient. We give here an independent proof of Lemma~\ref{Lemma_continuity_in_parameters} to avoid a circular reasoning

\begin{proof}[Proof of Lemma~\ref{Lemma_continuity_in_parameters}]

All parameters of the $m$-th variational datum are written with an exponent $(m)$, while the limiting values --- once shown to exist --- are written without exponent.

Since the extreme endpoints $a^{\prime\,(m)}_1$, $b^{\prime\, (m)}_H$ and the segment filling fractions $\hat{\boldsymbol{n}}^{(m)}$ converge to finite values as $m \rightarrow \infty$, there exists a constant $D > 0$ such that the nonnegative measures $\mu_1^{(m)},\ldots,\mu_H^{(m)}$ are supported inside n $m$-independent large segment $[-D,D]$ and have mass bounded by $D$. Hence, $(\boldsymbol{\mu}^{(m)})_{m \geq 1}$ admits converging subsequences. Let $\boldsymbol{\mu}$ be one of the limiting points, realized as the limit along the subsequence $(m_l)_{l \geq 1}$ . Eventually, we are going to show that there is a unique limiting point, and this will imply existence of the limit $\lim_{m\rightarrow\infty}\boldsymbol{\mu}^{(m)}=\boldsymbol{\mu}$.

Fix $h\in [H]$ and choose a contour $\gamma_h \subset \amsmathbb M_h$ surrounding $[\hat{a}^{\prime\,(m)}_h,\hat{b}^{\prime\,(m)}_h]$ for $m$ large enough, as described in Definition~\ref{def_contourgammah}. The weak convergence of $\mu^{(m_l)}_h$ implies that the corresponding Stieltjes transforms $\Gm_{\mu^{(m_l)}_h}(z)$ converge to $\Gm_{\mu_h}(z)$ uniformly over $z\in\gamma_h$ as $l \rightarrow \infty$. Hence, the functions $q^{\pm,(m_l)}_h(z)$ converge uniformly over $z\in\gamma_h$ as $l \rightarrow \infty$ to the functions $q^{\pm}_h(z)$ defined for $\tau \in \{\pm 1\}$ by
\begin{equation}
\label{eq_q_pm_limit}
q_h^{\tau}(z) = \phi_h^{+}(z) \cdot \exp\bigg(\sum_{g = 1}^H \theta_{h,g}\,\Gm_{\mu_{g}}(z)\bigg) + \tau\,\phi_h^{-}(z)\cdot \exp\bigg(- \sum_{g = 1}^H \theta_{h,g}\,\Gm_{\mu_{g}}(z)\bigg).
\end{equation}
At this point we can repeat the arguments in the proof of Lemma~\ref{Lemma_regularity_one_band}. Since $q^{+,(m_l)}_h(z)$ and $(q^{-,(m_l)}_h(z))^2$ are meromorphic functions of $z\in \amsmathbb M_h$ with only possible poles at $\hat a^{\prime\,(m_l)}_h$ and $\hat b^{\prime\,(m_l)}_h$, Cauchy integral formula implies that $q^{+}_h(z)$ and $(q^{-}_h(z))^2$ are also meromorphic in $\amsmathbb M_h$ with only possible poles at $\hat a'_h$ and $\hat b'_h$. Then we can treat \eqref{eq_q_pm_limit} with $\tau=+1$ as a quadratic equation with holomorphic in $\amsmathbb M_h$ coefficients for the unknown function $\Xi(z)=\exp\big(\theta_{h,h}\,\Gm_{\mu_h}(z)\big)$. We further conclude that $\Xi(z)$ has finite limits as $z$ approaches real $x\in (\hat{a}'_h,\hat{b}'_h)$, which depend on $x$ in n $\frac{1}{2}$-H\"older continuous way. Identifying the imaginary parts of these limits with the density of $\mu_h$, we deduce the five conclusions of Theorem~\ref{Theorem_regularity_density} for the measure $\mu_h$. We also find that the density of $\mu^{(m_{l})}_h$ converges uniformly on any compact to the density of $\mu_h$ --- this is again because we can express this density through integrals over $\gamma_h$ by using Cauchy integral formula.

It remains to show that there is at most one limiting point. For that we first notice that, for any $h \in [H]$, by Theorem~\ref{Theorem_regularity_density}, (iv) the number of bands of $\hat \mu_h$ is bounded from above by the number of zeros of $q^{-}_h(z)$ in $[\hat a'_h, \hat b'_h]$. The latter number is finite, as a holomorphic function that does not identically vanish in a compact set has finitely many zeros in this compact. Thus, we do not need to worry about pathological situations. Further, since the endpoints of the bands are found as simple zeros of holomorphic functions, Rouch\'e theorem implies that these endpoints for the bands of $\mu^{(m_l)}$ converge as $l \rightarrow \infty$ to those of the bands of $\mu$.

Let us now show that any limiting point $\boldsymbol{\mu}$ satisfies the characterization of Theorem~\ref{Theorem_equi_charact_repeat_2}. Since there is a unique measure with fixed filling fractions satisfying this characterization, this implies that there is only one limiting point, and that the previous stated convergences hold for the whole sequence $\mu^{(m)}$ as $m \rightarrow \infty$.

The weak convergence of $\boldsymbol{\mu}^{(m_l)}$ towards $\boldsymbol{\mu}$ as $l \rightarrow \infty$, the fact that all the involved measures have uniformly bounded density and the local integrability near $0$ of the function $x \mapsto \log|x|$ implies the pointwise convergence of the effective potentials \eqref{eq_V_eff} of $\mu^{(m_l)}$ to the effective potentials of $\boldsymbol{\mu}$:
\[
\forall x \in \amsmathbb{R}\qquad \lim_{l\rightarrow\infty} V^{\textnormal{eff},(m_l)}_h(x)= V^{\textnormal{eff}}_h(x).
\]
Hence the characterization of Theorem~\ref{Theorem_equi_charact_repeat_2} for $\boldsymbol{\mu}$ is obtained as $l\rightarrow\infty$ limit of the similar characterizations of $\boldsymbol{\mu}^{(m_l)}$ as equilibrium measures.
\end{proof}

We now have all ingredients to finish the proof of the most general case in Theorem~\ref{Theorem_regularity_density}.

\begin{proof}[Proof of Theorem~\ref{Theorem_regularity_density}]
 We first note that by Lemma~\ref{Lemma_continuity_in_parameters} we can waive the rationality restrictions in Lemma~\ref{Lemma_regularity_one_band} by approximating irrational parameters by rational ones. Further, we would like to allow the values of the interaction parameter $\theta_{1,1}=\theta$ different from $\theta=1$, while still keeping $H=1$. For this, let $\mathcal{I}$ be the functional for intensity of interaction $\theta$ and mass $\hat{n}$, and let $\mathcal{I}^1$ be the functional \eqref{eq_functional_general} for intensity of interaction $1$ and mass $\hat{n}' = \theta\hat{n}$, with all other parameters the same as for $\mathcal I$. If $\nu$ is a nonnegative measure with mass $\hat{n}$ and density bounded by $\frac{1}{\theta}$, then $\nu^1 = \theta \nu$ is a nonnegative measure with mass $\hat{n}'$ and density bounded by $1$ and we have
 \[
 \mathcal{I}[\nu] = \frac{\mathcal{I}^1[\nu^1]}{\theta}.
 \]
 Therefore, the minimizer $\mu$ of $-\mathcal{I}$ is related to the minimizer $\mu^1$ of $-\mathcal{I}^1$ by $\mu^1 = \theta \mu$. Note that the functions $\phi^\pm(z)$ and $q^\pm(z)$ are the same for $\mu$ and for $\mu^1$. So, the five properties listed in Theorem~\ref{Theorem_regularity_density} for $H=1$ and general value of $\theta>0$ follow by this procedure from the $\theta=1$ case.

Having fully established Theorem~\ref{Theorem_regularity_density} for $H=1$, we now proceed to the $H>1$ case. We fix $h\in[H]$. The idea of the extension to $H>1$ is that for any $h \in [H]$ the restriction $\mu_h$ of $\mu$ to $[\hat a'_h,\hat b'_h]$ itself is an equilibrium measure with $H=1$.
In order to make this precise we define a modified potential
\[
V^{\textnormal{mod}}_h(x) = V_h(x) - \sum_{g \neq h} 2\theta_{h,g} \int_{\hat{a}'_{g}}^{\hat{b}'_{g}} \log\big((x - y)\cdot \textnormal{sgn}(h - g)\big)\,\mu_{g}(y) \dd y
\]
and decompose it as
\[
 V^{\textnormal{mod}}_h(x)= \iota_h^{-}\,\mathrm{Llog}(x - \hat a'_h) + \iota_h^{+}\,\mathrm{Llog}(\hat b'_h - x) + U^{\textnormal{mod}}_h(x),
\]
where
\[
U^{\textnormal{mod}}_h(x) = U_h(x) - \sum_{g \neq h} 2\theta_{h,g} \int_{\hat{a}'_{g}}^{\hat{b}'_{g}} \log\big((x - y)\cdot\textnormal{sgn}(h - g)\big)\,\mu_{g}(y) \dd y.
\]
The reason to have replaced the absolute value with a sign function is that, since $U_h(x)$ is holomorphic in a simply-connected complex domain $\amsmathbb M_h$ including $[\hat a'_h,\hat b_h']$, but not $[\hat a'_{g},\hat b'_{g}]$, $g\neq h$, we can choose appropriate branches of the logarithm, so that $U^{\textnormal{mod}}_h(x)$ is also holomorphic in $\amsmathbb M_h$. Let $\mu^{\textnormal{mod}}$ denote the equilibrium measure for the parameters $H=1$, intensity of interaction $\theta:=\theta_{h,h}$, segment $[\hat{a}',\hat{b}'] := [\hat a'_h, \hat b'_h]$, potential $V := V^{\textnormal{mod}}_h$, and filling fraction $\hat{n} = \mu_h([\hat a'_h,\hat b'_h])$. The energy functional $-\I$ for this variational datum is --- up to addition of a constant --- the same as the original $H>1$ energy functional restricted to the set of measures coinciding with $\mu_{g}$ on segments $[\hat a'_{g}, \hat b'_{g}]$ for every $g \neq h$, so that its variable is effectively a single measure on $[\hat{a}',\hat{b}'] = [\hat a'_h,\hat b'_h]$. Hence $\mu^{\textnormal{mod}}$ coincides with $\mu_h$.

The properties (i) and (ii) of Theorem~\ref{Theorem_regularity_density} hold for the measure $\mu_h$ since we have already proven them for $\mu^{\textnormal{mod}}$. For the last two properties notice that
\[
 \partial_x V^{\textnormal{mod}}_h(x) =\partial_x V_h(x) - \sum_{g \neq h} 2\theta_{h,g}\, \Gm_{\mu_{g}}(x).
\]
Hence, by Definition~\ref{Definition_phi_functions_2} and formula \eqref{eq_q_pm}, the functions $q^{\pm}_h(z)$ for the original $H>1$ ensemble and for the modified $H=1$ ensemble differ by multiplication with the function
\[
 \exp\left(\sum_{g \neq h} \theta_{h,g} \,\Gm_{\mu_{g}}(x)\right),
\]
which is holomorphic and has no zeros in $\amsmathbb M_h$. Therefore, the properties (iii), (iv), and (v) of Theorem~\ref{Theorem_regularity_density} for $\mu_h$ and $\mu^{\textnormal{mod}}$ are equivalent statements.
\end{proof}

\begin{proof}[Proof of Proposition~\ref{Proposition_density}]

 The only part of the definition of $q^{\pm}_h$, which is not holomorphic near the point $x\in(\hat a_h,\hat b_h)$, is the function $\Gm_{\mu_h}(z)$; it has the upper- and lower-boundary values
 \[
\forall x \in \amsmathbb{R}\qquad \Gm_{\mu_h}(x^{\pm}) = \lim_{\eps \rightarrow 0^+} \Gm_{\mu_h}(x\pm \ii \eps)= \textnormal{p.v.}\, \Gm_{\mu_h}(x)\mp \ii\pi\mu_h(x).
 \]
 Here $\textnormal{p.v.}\, \Gm_{\mu_h}(x)$ is the Cauchy principal value of the defining integral for the Stieltjes transform, which exists because the density of $\mu_h$ is $\frac{1}{2}$-H\"older continuous. Hence, using continuity of $q^+_h$ near $x$, we write
\begin{equation*}
\begin{split}
& \quad \lim_{\eps \rightarrow 0^+} \frac{q^{-}_h(x-\ii \eps) - q^{-}_h(x+\ii\eps)}{2 \, q^+_h(x)} \\
 & = \lim_{\eps \rightarrow 0^+} \frac{q^{-}_h(x-\ii \eps) - q^{-}_h(x+\ii\eps)}
{q^{+}_h(x-\ii \eps) + q^{+}_h(x+\ii\eps)} \\
& = \frac{ 2\ii \sin\big(\pi \theta_{h,h} \mu_h(x)\big) \Big[\phi_h^{+}(x)\cdot\exp\big(\sum_{g = 1}^H \theta_{h,g}\,\Gm_{\mu_{g}}(x)\big) + \phi_h^{-}(x)\cdot\exp\big(- \sum_{g = 1}^H \theta_{h,g}\,\Gm_{\mu_{g}}(x)\big)\Big] } {2\cos\big(\pi \theta_{h,h} \mu_h(x)\big)\Big[\phi_h^{+}(x)\cdot\exp\big(\sum_{g = 1}^{H} \theta_{h,g}\,\Gm_{\mu_{g}}(x)\big) + \phi_h^{-}(x)\cdot\exp\big(- \sum_{g = 1}^H \theta_{h,g}\,\Gm_{\mu_{g}}(x)\big)\Big]} \\
& = \ii \tan\big(\pi \theta_{h,h} \mu(x)\big),
\end{split}
\end{equation*}
which proves \eqref{eq_density_tan}. For \eqref{eq_density_exponent}, we use the identity
\begin{equation}
\label{eq_characterization_derivative}
 \phi_h^{+}(x)\cdot \exp\bigg(\sum_{g = 1}^H \theta_{h,g}\,\Gm_{\mu_{g}}(x)\bigg) = \phi_h^{-}(x)\cdot\exp\bigg(- \sum_{g = 1}^{H} \theta_{h,g}\,\Gm_{\mu_{g}}(x)\bigg),
\end{equation}
which is valid for $x$ inside the band, where it is equivalent to $q^+_h(x^+)-q^{+}_h(x^-)=0$ proven in Theorem~\ref{Theorem_regularity_density}. This gives
\begin{equation*}
\begin{split}
& \quad \frac{2\,q^-_h(x^+)}{q^+_h(x^+) + q^-_h(x^+)} \\
& = \tfrac{\phi_h^+(x) \cdot \exp\big(\sum_{g = 1}^{H} \theta_{h,g} \,\Gm_{\mu_g}(x) - \ii \pi \theta_{h,h}\,\mu_h(x)\big) - \phi_h^+(x) \cdot \exp\big(-\sum_{g = 1}^{H} \theta_{h,g}\,\Gm_{\mu_g}(x) + \ii \pi \theta_{h,h}\,\mu_h(x)\big)}{\phi_h^+(x) \cdot \exp\big(\sum_{g = 1}^{H} \theta_{h,g}\,\Gm_{\mu_g}(x) - \ii \pi \theta_{h,h}\,\mu_h(x)\big)} \\
& = 1 - \exp\big(2\ii\pi \theta_{h,h}\,\mu_h(x)\big). \qedhere
\end{split}
\end{equation*}
\end{proof}

\begin{proof}[Proof of Proposition~\ref{Proposition_q_bounds}]
 By Theorem~\ref{Theorem_regularity_density} the following three functions are holomorphic for $z$ in $\amsmathbb{M}_{h}$:
 \[
 (z-\hat a'_h)(z-\hat b'_h) \cdot q^{+}_h(z),\qquad (z-\hat a'_h)^2(z-\hat b'_h)^2 \cdot (q^-_h(z))^2,\qquad (z-\hat a'_h)(z-\hat b'_h) \cdot s_h(z).
 \]
Choose a contour $\gamma_h \in \amsmathbb M_h$ surrounding $[\hat a'_h, \hat b'_h]$ as described in Definition~\ref{def_contourgammah}. By the maximum principle for harmonic functions the maximum over the compact domain bounded by $\gamma_h$ for each of the three functions
\[
\big|(z-\hat a'_h)(z-\hat b'_h) \cdot q^{+}_h(z)\big|,\qquad \big|(z-\hat a'_h)^2(z-\hat b'_h)^2 \cdot (q^{-}_h(z))^2\big|,\qquad \big|(z-\hat a'_h)(z-\hat b'_h) \cdot s_h(z)\big|
\]
 is attained on $\gamma_h$. On the other hand, on $\gamma_h$ these three functions are uniformly bounded from above directly from their definitions --- up to changing slightly $\gamma_h$. For $\partial_z U_h(z)$ entering into the definitions, we should additionally use a corollary of Cauchy integral formula, namely that the derivative of a uniformly bounded holomorphic function is uniformly bounded on a slightly smaller domain. Since for small enough $\eps$ the points $\hat a'_h$ and $\hat b'_h$ are bounded away from $\amsmathbb B_h^\eps$, this proves the upper bounds in \eqref{eq_q_bounds}.

For the lower bound on $s_h(z)$, note that by Cauchy integral formula and the upper bound we have just proven for the derivative $\partial_z s_h(z)$ is uniformly bounded in a complex neighborhood of $[\hat a'_h, \hat b'_h]$. Hence, it suffices to prove that $s_h(z)$ is bounded away from zero inside each band of $\mu_h$. Let $(\alpha,\beta)$ be such a band.
We rewrite \eqref{eq_density_exponent} as
\begin{equation}
\label{eq_x174}
1-\exp\big(2\ii\pi \theta_{h,h}\, \mu_h(x)\big)=\frac{\sigma_h(x^+) \cdot \psi(x) \cdot s_h(x)}{ \phi_h^{+}(x) \cdot \exp\big(-\ii\pi \theta_{h,h} \mu_{h}(x) + \sum_{g = 1}^H \theta_{h,g}\,\Gm_{\mu_{g}}(x) \big)}.
\end{equation}
Note that \eqref{eq_characterization_derivative} implies
\[
 \exp\bigg(\sum_{g = 1}^{H} 2\theta_{h,g}\,\Gm_{\mu_{g}}(x)\bigg)=\frac{\phi^-_h(x)}{\phi^+_h(x)}.
\]
Combining with \eqref{eq_x174} we conclude that
\begin{equation}
\label{eq_x175}
|s_h(x)|=\bigg|\frac{1-\exp\big(2\ii\pi \theta_{h,h}\,\mu_h(x)\big)}{\sigma_h(x^+) \cdot \psi(x)}\bigg|\cdot \sqrt{\big|\phi^-_h(x) \phi^+_h(x)\big|}.
\end{equation}
The first factor in the right-hand side of \eqref{eq_x175} is bounded away from zero by Conditions 4, 5, and 6 of Assumption~\ref{Assumption_C}. The second factor in the right-hand side of \eqref{eq_x175} is bounded away from $0$ by Definition~\ref{Definition_phi_functions_2} and Assumption~\ref{Assumption_B}.4, implying uniform bounds for $\partial_x U_h(x)$.
\end{proof}

\begin{proof}[Proof of Proposition~\ref{prop:convexmueq}]
Fix $h \in [H]$. Since $\mu_h$ has continuous density by Theorem~\ref{Theorem_regularity_density}, having the $h$-th segment neither fully void nor fully saturated implies the existence of at least one band, justifying (i). In this situation, the assumptions imply that the $h$-th effective potential
\[
V_h^{\textnormal{eff}}(x) = V_h(x) - \sum_{g = 1}^{H} 2\theta_{g,h} \int_{\hat{a}_g'}^{\hat{b}_g'} \log|x - y|\dd \mu(y)
\]
is a smooth strictly convex function of $x \in (\hat{a}_h',\hat{b}'_h) \setminus \textnormal{supp}(\mu_h)$, uniformly so depending on the constants in the assumption (it is a sum of convex functions, the first one being uniformly strictly convex). By extending this function to have constant value $v_h$ on the support, we get a continuous convex function. The locus where this function reaches its minimum $v_h$ must be a segment. Since the support coincides with this locus up to a zero-measure set, it must actually be a segment, proving (ii). Then, the uniform strict convexity outside of the support implies that we must have $\partial_x V_{h}^{\textnormal{eff}}(x) < c_L < 0$ to the left of the support and $\partial_x V_h^{\textnormal{eff}} > c_R > 0$ to the right of the support. This implies Assumption~\ref{Assumption_C}.2.

The derivative of the effective potential in a band $(\alpha,\beta)$ vanishes. This can be rewritten as a functional equation for the Stieltjes transform of the restriction of $\mu_h$ to the band
\begin{equation}
\label{tricoeq}
\forall x \in (\alpha,\beta)\qquad \int_{\alpha}^{\beta} \frac{\mu_h(y)\dd y}{x - y} = \tilde{V}'(x).
\end{equation}
The right-hand side involves the modified potential
\begin{equation}
\label{eq_x292}
\tilde{V}(x) = \frac{1}{\theta_{h,h}}\bigg(V_h(x) - 2\theta_{h,h} \int_{S} \log|x - y|\dd\mu_h(y) - \sum_{g \neq h} 2\theta_{h,g} \int_{\hat{a}_g'}^{\hat{b}'_g} \log|x - y|\dd\mu_g(y)\bigg),
\end{equation}
where $S = \textnormal{supp}(\mu_h) \setminus (\alpha,\beta)$. Note that $\tilde{V}(x)$ is again uniformly strictly convex. The function $\tilde{V}'(x)$ diverges at most logarithmically near the endpoints (due to the possible presence of saturations in $S$ adjacent to the band, or to the possible existence of such singularities in $V_h$ if $\alpha = \hat{a}'_h$ or $\beta = \hat{b}'_h$), so $V$ is of class $L^p$ for any $p > 0$. The solution of \eqref{tricoeq} is described by Tricomi \cite[Eqn. (18), p179]{Tricomi}. Keeping in mind that the density of $\mu_h$ is bounded, this yields
\begin{equation}
\label{Tricomun}
\mu_h(x) = \frac{\sqrt{(\beta - x)(x - \alpha)}}{2\pi^2}\int_{\alpha}^{\beta} \frac{\tilde{V}'(x) - \tilde{V}'(y)}{x - y} \frac{\dd y}{\sqrt{(\beta - y)(y - \alpha)}}.
\end{equation}
The uniform strict convexity of $\tilde{V}$ implies that the integrand is larger than a uniform positive constant. Hence, the condition 5.\ and the lower bound of 4.\ from Assumption~\ref{Assumption_C} are satisfied with constants depending only on $c$, the constants in the assumptions and $\beta-\alpha$. If all bands have length larger than $c' > 0$, then the constants in this statement can be made to depend only on $c$, $c'$, and the constants in the assumptions. This completes the justification of (iv).

Recall from Definition~\ref{Definition_void_saturated} that bands are open intervals, while voids and saturations are closed intervals (possibly reduced to a point). In this terminology, it is possible to have two consecutive bands $(\alpha,\beta)$ and $(\beta,\alpha')$ with $\{\beta\}$ being a void, \textit{i.e.} $\mu_h(\beta) = 0$. In this situation, the equation \eqref{tricoeq} is valid by continuity on the whole interval $(\alpha,\beta')$ and we conclude again from \eqref{Tricomun} by strict convexity of $\tilde{V}$ that the density of $\mu_h$ cannot vanish at $\beta$. This shows (iii).
\end{proof}

\section{Parametric Lipschitz regularity}
\label{Section-ParReg}

In this section we prove that equilibrium measures depend in a Lipschitz way on the parameters in the variational datum, and in particular is continuous with respect to those parameters. In each of the following lemmata, we have two variational data, denoted $\textnormal{I}$ and $\textnormal{II}$, two corresponding equilibrium measures $\boldsymbol{\mu}^{\textnormal{I}}$ and $\boldsymbol{\mu}^{\textnormal{II}}$, \textit{etc}. Throughout this section the analyticity of $U_h(x)$ of Assumption~\ref{Assumption_B} is not necessary and all the results continue to hold (with the same proofs) for $U_h(x)$ which is only twice-continuously differentiable, as in the Assumption~\ref{Assumptions_basic} for discrete ensembles.

\begin{lemma}
\label{Lemma_continuity_interactions}
Take two variational data satisfying Assumptions~\ref{Assumption_A} and \ref{Assumption_B}, differing only by their matrices of interactions $\boldsymbol{\Theta}^{\textnormal{I}}$ and $\boldsymbol{\Theta}^{\textnormal{II}}$. There exists $M > 0$ depending only on the constants in the assumptions, such that for any measurable function $f$ with finite $|\!|f|\!|_{\frac{1}{2}}$, we have
\begin{equation}
 \label{eq_continuity_interactions}
 \bigg| \int_{\amsmathbb R} f(x)\dd(\mu^{\textnormal{I}}-\mu^{\textnormal{II}})(x)\bigg|\leq M \cdot |\!|f|\!|_{\frac{1}{2}} \cdot \sqrt{|\!|\boldsymbol{\Theta}^{\textnormal{I}}- \boldsymbol{\Theta}^{\textnormal{II}}|\!|_{\infty}}.
 \end{equation}
\end{lemma}

In the following study of the dependence on the parameters from the second, third and fourth group of Section~\ref{Section_variational_data}, we assume that $\boldsymbol{\mu}^{\textnormal{I}}$ and $\boldsymbol{\mu}^{\textnormal{II}}$ satisfy the off-criticality Assumption~\ref{Assumption_C}, and deduce Lipschitz regularity. Yet, the first parts of our proofs would still give $\frac{1}{2}$-H\"older continuity statements without this assumption. We prefer to state the stronger results under the stronger Assumption~\ref{Assumption_C}, as they are used directly in Section~\ref{Section_Smoothnessparam}.

\begin{lemma}
\label{Lemma_continuity_potential}
 Take two variational data satisfying Assumptions~\ref{Assumption_A} and \ref{Assumption_B}, differing only by their potentials $V^{\textnormal{I}}_h$ and $V^{\textnormal{II}}_h$ on the $h$-th segment for a single $h \in [H]$. Assume furthermore that $\boldsymbol{\mu}^{\textnormal{I}}$ satisfies Assumption~\ref{Assumption_C}. Then, there exists $M > 0$ depending only on the constants in the assumptions, such that for any measurable function $f$ with finite $|\!|f|\!|_{\frac{1}{2}}$, we have \begin{equation}
 \label{eq_potential_inequality}
\bigg| \int_{\amsmathbb{R}} f(x)\dd(\mu^{\textnormal{I}} - \mu^{\textnormal{II}})(x) \bigg| \leq M \cdot |\!|f|\!|_{\frac{1}{2}} \cdot \Big(|\!|V^{\textnormal{I}}_{h} - V^{\textnormal{II}}_{h}|\!|_{\infty} + \sup_{[\hat{a}_h' + \frac{1}{M},\hat{b}_h' - \frac{1}{M}]} \big|(V_h^{\textnormal{I}} - V_{h}^{\textnormal{II}})'\big|\Big).
\end{equation}
\end{lemma}

 \begin{lemma}
\label{Lemma_continuity_filling_fractions}
Take two variational data satisfying Assumptions~\ref{Assumption_A} and \ref{Assumption_B}, differing only by their segment filling fractions $\hat{n}_h^{\textnormal{I}}$ and $\hat{n}_{h}^{\textnormal{II}}$ for a single $h \in [H]$. Assume furthermore that $\boldsymbol{\mu}^{\textnormal{I}}$ satisfies Assumption~\ref{Assumption_C}. Then, there exists $M > 0$ depending only on the constants in the assumptions, such that for $|\hat{n}_{h}^{\textnormal{I}} - \hat{n}_{h}^{\textnormal{II}}| < \frac{1}{M}$ and for any measurable function $f$ with finite $|\!|f|\!|_{\frac{1}{2}}$ and $|\!|f|\!|_{\infty}$ norms we have
 \begin{equation}
 \label{eq_continuity_filling_fractions}
 \bigg| \int_{\amsmathbb R} f(x)\dd(\mu^{\textnormal{I}}-\mu^{\textnormal{II}})(x)\bigg|\leq M \cdot \big(|\!|f|\!|_{\frac{1}{2}} + |\!|f|\!|_{\infty}\big)\cdot |\hat{n}_{h}^{\textnormal{I}} - \hat{n}_{h}^{\textnormal{II}}|. \end{equation}
\end{lemma}

\begin{lemma}
\label{Lemma_continuity_end_points_void}
Take two variational data satisfying Assumptions~\ref{Assumption_A} and \ref{Assumption_B}, differing only by their segment endpoints $\hat{a}_{h}^{\prime\,\textnormal{I}}$ and $\hat{a}_{h}^{\prime\,\textnormal{II}}$ for a single $h \in [H]$. Assume furthermore that $\boldsymbol{\mu}^{\textnormal{I}}$ satisfies Assumption~\ref{Assumption_C} and that $\hat a^{\prime\,\textnormal{I}}_h$ belongs to a void of $\boldsymbol{\mu}^\textnormal{I}$. Then there exists $M > 0$ depending only on the constants in the assumptions, such that if $|\hat{a}_{h}^{\prime\,\textnormal{I}} - \hat{a}_{h}^{\prime\,\textnormal{II}}| < \frac{1}{M}$, then $\boldsymbol{\mu}^{\textnormal{I}} = \boldsymbol{\mu}^{\textnormal{II}}$.
\end{lemma}
Note that in the above lemma, we must have $\iota_h^-=0$ for both variational data, as otherwise the potentials $V_h^\textnormal{I}$ and $V_h^\textnormal{II}$ specified by Assumption~\ref{Assumption_B}.3 cannot be equal on their common domain of definition.

\begin{lemma}
\label{Lemma_continuity_end_points_saturated} Take two variational data satisfying Assumptions~\ref{Assumption_A} and \ref{Assumption_B} with $\iota_h^{-,\textnormal{I}}=\iota_h^{-,\textnormal{II}}=2$, differing by their segment endpoints $\hat{a}_{h}^{\prime\,\textnormal{I}}$ and $\hat{a}_{h}^{\prime\,\textnormal{II}}$ for a single $h \in [H]$, and by the potentials
\begin{equation*}
\begin{split}
V_h^{\textnormal{I}}(x) & = -2\, \mathrm{Llog}(x - \hat{a}_{h}^{\prime\,\textnormal{I}}) + V_h^0(x), \\
V_h^{\textnormal{II}}(x) & = -2 \,\mathrm{Llog}(x - \hat{a}_{h}^{\prime\,\textnormal{II}}) + V_{h}^{0}(x).
\end{split}
\end{equation*}
Assume furthermore that $\boldsymbol{\mu}^{\textnormal{I}}$ satisfies Assumption~\ref{Assumption_C} and $\hat{a}_{h}^{\prime\,\textnormal{I}}$ belongs to a saturation of $\boldsymbol{\mu}^{\textnormal{I}}$. Then, there exists $M > 0$ depending only on the constants in the assumptions, such that for $|\hat{a}_{h}^{\prime\,\textnormal{I}} - \hat{a}_{h}^{\prime\,\textnormal{II}}| < \frac{1}{M}$ and for any measurable function $f$ with finite $|\!|f|\!|_{\frac{1}{2}}$ we have
\begin{equation}
\label{eq_continuity_end_points_saturated}
 \bigg| \int_{\amsmathbb R} f(x)\dd(\mu^{\textnormal{I}}-\mu^{\textnormal{II}})(x)\bigg| \leq M \cdot \big(|\!|f|\!|_{\frac{1}{2}} + |\!|f|\!|_{\infty}\big) \cdot |\hat{a}_{h}^{\prime\,\textnormal{I}} - \hat{a}_{h}^{\prime\,\textnormal{II}}|.
\end{equation}
\end{lemma}
Let us explain why $\iota_h^-=2$ in Lemma~\ref{Lemma_continuity_end_points_saturated}. When we restrict an equilibrium measure to a smaller interval by removing a saturated segment, we need to absorb interactions with this segment into the potential; direct integration over the saturated segment leads to the appearance of $-2\,\mathrm{Llog}(x-\hat{a}')$ term in the potential, if $\hat{a}'$ was the left endpoint of the removed saturated segment. Comparing to the form of the potential specified in Assumption~\ref{Assumption_B}.3, this automatically creates a $x\log |x|$ singularity with prefactor $\iota^- = 2$. In Lemma~\ref{Lemma_continuity_end_points_saturated} we want to compare potentials having same type of singularities, and therefore we are forced to assume $\iota_h^{-,\,\textnormal{I}} = \iota_h^{-,\,\textnormal{II}} = 2$ from the start.

Results similar to Lemmata~\ref{Lemma_continuity_end_points_void} and \ref{Lemma_continuity_end_points_saturated} also hold for the right endpoints $\hat{b}_{h}'^{,\textnormal{I}}$ and $\hat{b}_{h}'^{,\textnormal{II}}$. The reason for the complicated setup in Lemma~\ref{Lemma_continuity_end_points_saturated} is that the weight $w_h(x)$ in the discrete ensemble
usually vanishes at the $a_h$ which are endpoints of saturations. In particular, the vanishing is useful when we apply the Nekrasov equations (\textit{cf.} Theorem~\ref{Theorem_Nekrasov}), and it also shows up when we condition on positions of densely packed particles in saturated segments in Chapter~\ref{Chapter_conditioning}. The vanishing leads to the appearance of the $(x-\hat a'_h)\log |x-\hat a'_h|$ term in the formula for $V_h(x)$ of Assumption~\ref{Assumption_B}.3. Therefore, we do not want to change $\hat a'_h$ without simultaneously changing the potential $V_h(x)$.

\begin{proof}[Proof of Lemma~\ref{Lemma_continuity_interactions}] Hereafter we denote, for $\textnormal{J} \in \{\textnormal{I},\textnormal{II}\}$, by $\I^\textnormal{J}, \mathfrak D^\textnormal{J}, V_h^{\textnormal{eff},J}$ the functional, distance, $h$-th effective potential associated with the variational datum $\textnormal{J}$, the potential $V_h^\textnormal{J}$, and the equilibrium measure $\boldsymbol{\mu}^\textnormal{J}$. First, we consider the case of matrices $\boldsymbol{\Theta}^\textnormal{I}$ and $\boldsymbol{\Theta}^\textnormal{II}$ having
 the same diagonal elements. Set $\boldsymbol{\nu}=\boldsymbol{\mu}^\textnormal{I}-\boldsymbol{\mu}^\textnormal{II}$. Using the definition of $\boldsymbol{\mu}^\textnormal{I}$, we have
\begin{equation}
\label{eq_x116}
\begin{split}
0 & \leq \I^\textnormal{I}[\boldsymbol{\mu}^\textnormal{I}]-\I^\textnormal{I}[\boldsymbol{\mu}^\textnormal{II}]\\
& = \sum_{g,h=1}^H
\int_{\hat a'_{g}}^{\hat b'_{g}} \int_{\hat a'_{h}}^{\hat b'_{h}}
\log|x-y|\,
\big(\theta_{g,h}^\textnormal{I}\dd(\mu_g^\textnormal{II}+\nu_g)(x)\dd(\mu_h^\textnormal{II}+\nu_h)(y)-\theta_{g,h}^\textnormal{II}\dd
\mu^\textnormal{II}_g(x)\dd\mu^\textnormal{II}_h(y)\big)\\
& \quad -\sum_{h=1}^H \int_{\hat a'_{h}}^{\hat b'_{h}} V_h(x)
\nu_h(x) \dd x\\
& = -\mathfrak{D}^{2\,\textnormal{{II}}}[\boldsymbol{\mu}^\textnormal{I},\boldsymbol{\mu}^\textnormal{II}]-\sum_{h=1}^{H} \int_{\hat a'_{h}}^{\hat b'_{h}}
V^{\textnormal{eff},\rm II}_{h}(x) \nu_h(x)\dd x \\
 & \quad + \sum_{g,h= 1}^{H} \int_{\hat
a'_{g}}^{\hat b'_{g}} \int_{\hat a'_{h}}^{\hat b'_{h}}
(\theta_{g,h}^\textnormal{I}-\theta_{g,h}^\textnormal{II}) \log|x-y|\,
\dd(\mu^\textnormal{II}_g+\nu_g)(x)\dd(\mu^\textnormal{II}_h+\nu_h)(y).
\end{split}
\end{equation}
We claim that
\begin{equation}
\label{eq_x76}
\forall h \in [H]\qquad \int_{\hat a'_{h}}^{\hat b'_{h}} V^{\textnormal{eff},\textnormal{II}}_{h}(x) \nu_h(x) \dd x\geq 0.
\end{equation}
Indeed, if we let $\boldsymbol{v}^{\textnormal{II}}$ be the constants of Theorem
\ref{Theorem_equi_charact_repeat_2} for $\boldsymbol{\mu}^{\textnormal{II}}$, we can replace
$V^{\textnormal{eff},\textnormal{II}}_{h}(x)$ with $V^{\textnormal{eff},\textnormal{II}}_{h}(x)-v_{h}^{\textnormal{II}}$ in
\eqref{eq_x76}. The integral is unchanged since $\boldsymbol{\mu}^\textnormal{I}$ and $\boldsymbol{\mu}^{\textnormal{II}}$ have the
same segment filling fractions. Studying separately the cases of bands, voids, and saturations
of $\boldsymbol{\mu}^{\textnormal{II}}$ using Theorem~\ref{Theorem_equi_charact_repeat_2}, we conclude that
\begin{equation}
\label{eq_x77}
\forall h \in [H] \quad \forall x\in [\hat a'_{h},\hat b'_{h}]\qquad \big(V^{\textnormal{eff},\textnormal{II}}_{h}(x)-v_{h}^\textnormal{II}\big) \nu_h(x) \geq 0,
\end{equation}
which implies \eqref{eq_x76}. On the other hand, the
absolute value of the third term in \eqref{eq_x116} is bounded from above by $M \cdot |\!|\boldsymbol{\Theta}^{\textnormal{I}} - \boldsymbol{\Theta}^{\textnormal{II}}|\!|_{\infty}$ for some $M > 0$ because of the uniform bound on the density of the equilibrium measure. Therefore
\begin{equation}
\label{eq_x117}
\mathfrak{D}^{2\,\textnormal{II}}[\boldsymbol{\mu}^\textnormal{I},\boldsymbol{\mu}^\textnormal{II}] \leq M \cdot\,|\!|\boldsymbol{\Theta}^{\textnormal{I}} - \boldsymbol{\Theta}^{\textnormal{II}}|\!|_{\infty}.
\end{equation}
Using Lemma~\ref{Lemma_linear_through_distance}, the inequality
\eqref{eq_continuity_interactions} follows from \eqref{eq_x117}.

\smallskip

We proceed to the case of matrices $\boldsymbol{\Theta}^\textnormal{I}$, $\boldsymbol{\Theta}^\textnormal{II}$ with different diagonal elements. The difference is that now $\boldsymbol{\mu}^\textnormal{I}$
is constrained by $\mu_h^\textnormal{I}(x)\leq (\theta_{h,h}^\textnormal{I})^{-1}$ while $\boldsymbol{\mu}^\textnormal{II}$ is constrained by
$\mu_h^\textnormal{II}(x)\leq (\theta_{h,h}^\textnormal{II})^{-1}$. The inequality $0\le
\I^\textnormal{I}[\boldsymbol{\mu}^\textnormal{I}]-\I^\textnormal{I}[\boldsymbol{\mu}^\textnormal{II}]$ may then fail to be valid, because $\boldsymbol{\mu}^\textnormal{II}$ does not belong to the class of measures over which we minimize $-\I^\textnormal{I}$ to get $\boldsymbol{\mu}^\textnormal{I}$. We fix this by changing the matrix elements from $\boldsymbol{\Theta}^\textnormal{I}$ to $\boldsymbol{\Theta}^\textnormal{II}$
one by one, and proving the bound \eqref{eq_continuity_interactions} on each step
separately. Without loss of generality we can assume that in one step all
matrix elements of $\boldsymbol{\Theta}^\textnormal{I}$ and $\boldsymbol{\Theta}^\textnormal{II}$ are the same except for
$\theta_{1,1}$, and we have $\theta^\textnormal{I}_{1,1}<\theta^\textnormal{II}_{1,1}$. Then the inequality
\[
\I^I[\boldsymbol{\mu}^\textnormal{I}]-\I^I[\boldsymbol{\mu}^\textnormal{II}] \geq 0
\]
does hold. We also should be more careful with the second term in the third line of \eqref{eq_x116}. Indeed, it might be no longer true
that $
V^{\textnormal{eff},\textnormal{II}}_{1}(x) \nu(x)\geq 0$ for $x$ in a saturation of $\boldsymbol{\mu}^\textnormal{II}$ inside $[\hat a'_1,\hat b'_1]$, as $\nu(x)$ might be positive there. However, we still have for $x$ in a saturation of $\boldsymbol{\mu}^\textnormal{II}$,
\[\nu(x)\le
\frac{1}{\theta_{1,1}^\textnormal{I}} - \frac{1}{\theta_{1,1}^\textnormal{II}} \leq C\,
|\theta^\textnormal{I}_{1,1}-\theta^\textnormal{II}_{1,1}| \quad \textnormal{and}\quad V^{\textnormal{eff},\textnormal{II}}_{1}(x)-v_1^{\textnormal{II}}\leq 0,
\]
and therefore, the second term in \eqref{eq_x116} is bounded from above by $M \cdot
|\theta^\textnormal{I}_{1,1}-\theta^\textnormal{II}_{1,1}|$ for some $M > 0$. We conclude that \eqref{eq_x117} still holds and
the argument goes through.
\end{proof}

\begin{proof}[Proof of Lemma~\ref{Lemma_continuity_potential}] Since it does not involve potentials, $\mathfrak{D}^2$ is common to the two variational data. Set $\boldsymbol{\nu}=\boldsymbol{\mu}^\textnormal{I}-\boldsymbol{\mu}^{\textnormal{II}}$. Using the definition of $\boldsymbol{\mu}^\textnormal{I}$ we have
\begin{equation}
\label{eq_x79}
\begin{split}
0 & \leq \I^\textnormal{I}[\boldsymbol{\mu}^\textnormal{I}]-\I^\textnormal{I}[\boldsymbol{\mu}^{\textnormal{II}}] \\
& = \sum_{g_1,g_2 = 1}^H \int_{\hat
a'_{g_1}}^{\hat b'_{g_1}} \int_{\hat a'_{g_2}}^{\hat b'_{g_2}} \theta_{g_1,g_2}
\log|x-y|\, \big[\dd(\mu^{\textnormal{II}}_{g_1}+\nu_{g_1})(x)\dd(\mu^{\textnormal{II}}_{g_2}+\nu_{g_2})(y) -\dd
\mu^{\textnormal{II}}_{g_1}(x)\dd\mu^{\textnormal{II}}_{g_2}(y)\big] \\
& \quad -\sum_{g=1}^H \int_{\hat a'_{g}}^{\hat b'_{g}} V_{g}^{\textnormal{I}}(x)
\nu_g(x) \dd x \\
& = - \D[\boldsymbol{\mu}^\textnormal{I},\boldsymbol{\mu}^{\textnormal{II}}]-\int_{\hat a'_h}^{\hat b'_h}
\big(V^{\textnormal{I}}_h(x)-V^{\textnormal{II}}_h(x)\big)\nu_h(x) \dd x -\sum_{g=1}^{H} \int_{\hat a'_{g}}^{\hat
b'_{g}} V^{\textnormal{eff},\textnormal{II}}_{g}(x) \nu_g(x)\dd x,
\end{split}
\end{equation}
where $V^{\textnormal{eff},\textnormal{II}}_g$ is the effective potential of \eqref{eq_V_eff} on the $g$-th segment
corresponding to the variational datum $\textnormal{II}$. As in the proof of Lemma~\ref{Lemma_continuity_interactions}, one checks that for every $g \in [H]$
\[
 \int_{\hat a'_{g}}^{\hat b'_{g}} V^{\textnormal{eff},\textnormal{II}}_{g}(x) \nu_g(x) \dd x \geq 0.
\]
Combining the last inequality with \eqref{eq_x79}, we deduce
\begin{equation}
\label{eq_x80} \D[\boldsymbol{\mu}^\textnormal{I},\boldsymbol{\mu}^{\textnormal{II}}]\leq \int_{\hat a'_h}^{\hat b'_h}
\big(V_h^{\textnormal{II}}(x)-V_h^{\textnormal{I}}(x)\big)\nu_h(x)\dd x
\end{equation}
Recall that Assumption~\ref{Assumption_C} guarantees the existence of a small $\delta > 0$ depending only on the constants in the assumptions such that $[\hat{a}'_{h},\hat{a}'_{h} + \delta] \cup [\hat{b}'_{h} - \delta,\hat{b}'_{h}]$ is void or saturated in $\boldsymbol{\mu}^{\textnormal{I}}$. We next claim that the same is true for $\boldsymbol{\mu}^{\textnormal{II}}$ if $V_h^{\textnormal{I}}-V_h^{\textnormal{II}}$ is small enough:

\vspace{0.2cm}
\noindent \textsc{Claim.} There exists $\eta > 0$ depending only on the constants in the assumptions such that $\big[\hat{a}'_{h},\hat{a}'_{h} + \frac{\delta}{2}\big] \cup \big[\hat{b}'_{h} - \frac{\delta}{2},\hat{b}'_{h}\big]$ is void or saturated in $\boldsymbol{\mu}^{\textnormal{I}}$ and $\boldsymbol{\mu}^{\textnormal{II}}$ as soon as $|\!|V^{\textnormal{I}}_h - V^{\textnormal{II}}_h|\!|_{\infty} \leq \eta$.
\vspace{0.2cm}

We temporarily admit the claim and continue with the proof of Lemma~\ref{Lemma_continuity_potential}. The claim implies that we have $\nu_h(x) = 0$ for any $x \in \big[\hat{a}'_{h},\hat{a}'_{h} + \frac{\delta}{2}\big] \cup \big[\hat{b}'_{h} - \frac{\delta}{2},\hat{b}'_{h}\big]$.
Let us choose a smooth function $\chi$ on the real line taking values in $[0,1]$, vanishing outside $\big[\hat{a}'_{h} + \frac{\delta}{4},\hat{b}'_{h} - \frac{\delta}{4}\big]$, and equal to $1$ on $\big[\hat{a}'_{h} + \frac{\delta}{2},\hat{b}'_{h} - \frac{\delta}{2}\big]$. Then \eqref{eq_x80} together with Lemma~\ref{Lemma_linear_through_distance} imply
\begin{equation}
\label{633zwei}\D[\boldsymbol{\mu}^{\textnormal{I}},\boldsymbol{\mu}^{\textnormal{II}}] \leq \int_{\hat{a}'_{h}}^{\hat{b}'_{h}} \chi(x) \cdot \big(V^{\textnormal{I}}_h(x) - V^{\textnormal{II}}_h(x)\big)\dd\nu_h(x) \leq M \cdot |\!|\chi \cdot (V^{\textnormal{I}}_h - V^{\textnormal{II}}_h)|\!|_{\frac{1}{2}}\cdot\mathfrak{D}[\boldsymbol{\mu}^{\textnormal{I}},\boldsymbol{\mu}^{\textnormal{II}}]
\end{equation}
for some constant $M > 0$. By the Plancherel theorem for the Fourier transform
\begin{equation*}
\begin{split}
|\!|F|\!|_{\frac{1}{2}}^2 & = \int_{\amsmathbb{R}} |s|\,|\widehat{F}(s)|^2\,\dd s \\
& \leq \frac{1}{2}\bigg(\int_{\amsmathbb{R}} |\widehat{F}(s)|^2\dd s + \int_{\amsmathbb{R}} |s\,\widehat{F}(s)|^2\,\dd s\bigg) = \pi\bigg(\int_{\amsmathbb{R}} |F(x)|^2\,\dd x + \int_{\amsmathbb{R}} |F'(x)|^2\,\dd x\bigg)
\end{split}
\end{equation*}
Using this inequality for $F = \chi \cdot (V^{\textnormal{I}}_h - V^{\textnormal{II}}_h)$ in \eqref{633zwei}, we conclude that
\[
\mathfrak{D}[\boldsymbol{\mu}^{\textnormal{I}},\boldsymbol{\mu}^{\textnormal{II}}] \leq M \cdot \Big( |\!|V^{\textnormal{I}}_h - V^{\textnormal{II}}_h|\!|_{\infty} + \sup_{x \in [\hat{a}'_{h} + \frac{1}{M},\hat{b}'_{h} - \frac{1}{M}]} \big|(V^{\textnormal{I}}_h - V^{\textnormal{II}}_h)'(x)|\Big)
\]
for a perhaps different constant $M > 0$. Using Lemma~\ref{Lemma_linear_through_distance} we arrive to the desired bound \eqref{eq_potential_inequality}.

It remains to justify the claim. Let $\eta > 0$ and assume that $ |\!|V^{\textnormal{I}}_h - V^{\textnormal{II}}_h|\!|_{\infty} \leq \eta$. Since $\boldsymbol{\nu}=\boldsymbol{\mu}^\textnormal{I}-\boldsymbol{\mu}^\textnormal{II}$, \eqref{eq_x80} implies that $\D[\boldsymbol{\mu}^{\textnormal{I}},\boldsymbol{\mu}^{\textnormal{II}}] \leq 2\eta$. Using Lemma~\ref{Lemma_linear_through_distance}, we deduce that for any $\tilde{\eta} > 0$ and $g \in [H]$ we have
\[
\forall x,y \in [\hat{a}'_{h},\hat{b}'_{h}]\qquad \int_{\hat{a}'_{g}}^{\hat{b}'_{g}} \frac{1}{2}\log\big((x - y)^2 + \tilde{\eta}^2\big) \dd\nu_g(x) \leq M_0(\tilde{\eta})\,\eta^{\frac{1}{2}}
\]
for some constant $M_0(\tilde{\eta}) > 0$. On the other hand, for any $g \in [H]$
\begin{equation}
\label{636zwei} \int_{\hat{a}'_{g}}^{\hat{b}'_{g}} \bigg(\log|x - y| - \frac{1}{2}\log\big((x - y)^2 + \tilde{\eta}^2\big)\bigg)\dd\nu_g(x) = \int_{\hat{a}'_{g}}^{\hat{b}'_{g}} \frac{1}{2}\,\log\bigg(\frac{(x - y)^2}{(x - y)^2 + \tilde{\eta}^2}\bigg)\,\dd\nu_g(x).
\end{equation}
The integrand in the last expression converges to $0$ as $\tilde{\eta} \rightarrow 0$ while
for $\tilde{\eta} \leq 1$
\[
\bigg|\log\bigg(\frac{(x - y)^2}{(x - y)^2 + \tilde{\eta}^2}\bigg)\bigg| \leq \log\bigg(1 + \frac{1}{(x - y)^2}\bigg)
\]
is integrable. We deduce by the dominated convergence theorem that \eqref{636zwei} converges to $0$ as $\tilde{\eta} \rightarrow 0$. We conclude that
\[
\sup_{y \in [\hat{a}'_{h},\hat{b}'_{h}]} \bigg|\int_{\hat{a}'_{g}}^{\hat{b}'_{g}} \log|x - y|\,\dd\nu_g(x)\bigg|
\]
can be made arbitrarily small by choosing $\tilde{\eta}$ and $\eta$ small enough depending only on the constants in the assumptions. Consequently, $|\!| V^{\textnormal{eff},\textnormal{I}}_{h} - V^{\textnormal{eff},\textnormal{II}}_h|\!|_{\infty}$ also becomes arbitrarily small. In combination with Conditions 2. and 3. of Assumption~\ref{Assumption_C} for $\boldsymbol{\mu}^{\textnormal{I}}$ and the characterization of voids, bands and saturations in Theorem~\ref{Theorem_equi_charact_repeat_2}, we obtain the claim.
\end{proof}

\begin{proof}[Proof of Lemma~\ref{Lemma_continuity_filling_fractions}] The value of the constants $M$ that will appear in this proof can vary from line to line, and it is understood as usual that they can be chosen to depend only on the constants in the assumptions. The energy functional $-\I$ is common to the two variational data, but the equilibrium measures $\boldsymbol{\mu}^{\textnormal{I}}$ and $\boldsymbol{\mu}^{\textnormal{II}}$ minimize it over two different sets of measures since the segment filling fractions differ. Likewise, the pseudo-distance $\mathfrak{D}$ is common to the two variational data. Let $\boldsymbol{v}^{\textnormal{I}}$ be the $H$-tuple of constants of Theorem~\ref{Theorem_equi_charact_repeat_2} for $\boldsymbol{\mu}^{\textnormal{I}}$, and set
\[
\forall g \in [H] \quad \forall x \in [\hat{a}'_g,\hat{b}'_g]\qquad U^{\textnormal{eff},\textnormal{I}}_{g}(x)=V^{\textnormal{eff},\textnormal{I}}_{g}(x)-v_{g}^{\textnormal{I}}.
\]
For any measure $\nu$ on $ \amsmathbb{A}$
 \begin{equation}
 \label{eq_x84}
 \I[\boldsymbol{\nu}]=\I[\boldsymbol{\mu}^{\textnormal{I}}]- \sum_{g = 1}^{H} \int_{\hat{a}'_g}^{\hat{b}_g'}
 U^{\textnormal{eff},\textnormal{I}}_g(x)\, \dd(\nu_g-\mu_g^\textnormal{I})(x)-\D[\boldsymbol{\nu},\boldsymbol{\mu}^\textnormal{I}]-\sum_{g=1}^H
 v_{g}^{\textnormal{I}}\Big(\nu([\hat a'_{g},\hat b'_{g}])- \hat{n}_{g}^{\textnormal{I}}\Big).
 \end{equation}
 Take a measure $\boldsymbol{\nu}$ whose voids, bands, and saturations are the same as for
 the equilibrium measure $\boldsymbol{\mu}^\textnormal{I}$, and whose segment filling fractions are $\hat{\boldsymbol{n}}^{\textnormal{II}}$. Since $\boldsymbol{\mu}^{\textnormal{II}}$ is a minimizer of $-\I$ over the set of measures with the same segment filling fractions --- among other conditions ---, using
 \eqref{eq_x84} we get
 \begin{equation}
 \label{eq_x85}
 0\le
 \I[\boldsymbol{\mu}^\textnormal{II}]-\I[\boldsymbol{\nu}]=\D[\boldsymbol{\nu},\boldsymbol{\mu}^\textnormal{I}]-\D[\boldsymbol{\mu}^\textnormal{II},\boldsymbol{\mu}^\textnormal{I}]+
 \sum_{g = 1}^{H} \int_{\hat{a}'_g}^{\hat{b}_g'} U^{\textnormal{eff},\textnormal{I}}_g(x)\, \dd(\nu_g-\mu_g^\textnormal{II})(x).
 \end{equation}
 We claim that the last term in \eqref{eq_x85} is nonpositive. Indeed, considering
 separately $x$ in voids, bands, and saturations of $\boldsymbol{\nu}$ --- equivalently,
 of $\boldsymbol{\mu}^\textnormal{I}$ --- and using Theorem~\ref{Theorem_equi_charact_repeat_2}, we see that the
 product of the integrand and the density of the measure is nonpositive everywhere. Therefore,
 \begin{equation}
 \label{eq_x86}
 \D[\boldsymbol{\mu}^\textnormal{II},\boldsymbol{\mu}^\textnormal{I}]\leq
 \D[\boldsymbol{\nu},\boldsymbol{\mu}^\textnormal{I}].
 \end{equation}

Next, we make a particular choice of $\boldsymbol{\nu}$. By Assumption~\ref{Assumption_C}, there exists $\eps > 0$ and a segment $[\mathfrak{a},\mathfrak{b}] \subset [\hat a'_h, \hat b'_h]$ of positive length, both depending only on the constants in the assumptions, such that
\[
\forall x \in [\mathfrak{a},\mathfrak{b}]\qquad \eps < \mu^{\textnormal{I}}(x) < \frac{1}{\theta_{h,h}} - \eps.
\]
We set
\[
\nu = \mu^\textnormal{I}+ \frac{3(\hat n^\textnormal{II}_{h}-\hat n^\textnormal{I}_{h})}{\mathfrak{b} - \mathfrak{a}} \textnormal{Leb}_{\big[\frac{2\mathfrak{a}+ \mathfrak{b}}{3},\frac{\mathfrak{a} + 2\mathfrak{b}}{3}\big]},
\]
where $\textnormal{Leb}_{[r,s]}$ is the Lebesgue measure of density $1$ on $[r,s]$ and $0$ outside. Thanks to the assumption on this segment, $\boldsymbol{\nu}$ has the same voids, bands and saturations as $\boldsymbol{\mu}^{\textnormal{I}}$, and has density between $0$ and $\frac{1}{\theta_{h,h}}$ on $[\hat{a}'_{h},\hat{b}'_{h}]$ provided that
\[
|\hat{n}^{\textnormal{I}}_h - \hat{n}^{\textnormal{II}}_{h}| \leq \frac{\eps(\mathfrak{b} - \mathfrak{a})}{3}.
\]
By construction, the filling fractions are the same for $\boldsymbol{\nu}$ and $\boldsymbol{\mu}^\textnormal{II}$, and we compute
\[
\D[\boldsymbol{\nu},\boldsymbol{\mu}^{\textnormal{I}}] = \theta_{h,h}\cdot |\hat{n}_{h}^{\textnormal{I}} - \hat{n}_{h}^{\textnormal{II}}|^2 \cdot \bigg(\frac{3}{2} + \log 3 - \log(\mathfrak{b} - \mathfrak{a})\bigg).
\]
We conclude by \eqref{eq_x86} that
\begin{equation}
\label{640ein} \D[\boldsymbol{\mu}^{\textnormal{II}},\boldsymbol{\mu}^{\textnormal{I}}] \leq M \cdot (\hat{n}_{h}^{\textnormal{I}} - \hat{n}_{h}^{\textnormal{II}})^2.
\end{equation}
for some constant $M > 0$. Note that as $\boldsymbol{\mu}^{\textnormal{II}}$ and $\boldsymbol{\mu}^{\textnormal{I}}$ do not have the same segment filling fractions, we cannot use directly Lemma~\ref{Lemma_linear_through_distance}.
Let us introduce an auxiliary $H$-tuple of nonnegative measures $\boldsymbol{\varrho}$, with compact support, such that $\varrho_g = 0$ for $g \neq h$ and $\varrho_h$ has bounded density and total mass $1$. We have
\begin{equation}
\label{ppmpmp}
\begin{split}
\D\big[\boldsymbol{\mu}^{\textnormal{I}} + (\hat{n}_{h}^{\textnormal{II}} - \hat{n}_{h}^{\textnormal{I}})\boldsymbol{\varrho},\mu^{\textnormal{II}}\big]
 & = \D[\boldsymbol{\mu}^{\textnormal{I}},\boldsymbol{\mu}^{\textnormal{II}}] + (\hat{n}_{h}^{\textnormal{I}} - \hat{n}_{h}^{\textnormal{II}})^2 \cdot\D[\boldsymbol{\varrho},\boldsymbol{0}]
 \\ & \quad - 2(\hat{n}_{h}^{\textnormal{II}} - \hat{n}_{h}^{\textnormal{I}}) \bigg(\sum_{g = 1}^{H} \theta_{h,g} \int_{\hat{a}'_{g}}^{\hat{b}'_{g}}\int_{\hat{a}'_{h}}^{\hat{b}'_{h}} \log|x - y|\,\dd\varrho_g(x)\dd(\mu^{\textnormal{I}}_h - \mu^{\textnormal{II}}_h)(y) \bigg).
\end{split}
\end{equation}
We use this equation in two steps. Recall that Assumption~\ref{Assumption_C} guarantees the existence of $\delta > 0$ depending only on the constants in the assumptions, such that $[\hat{a}'_{h},\hat{a}'_{h} + \delta] \cup [\hat{b}'_{h} - \delta,\hat{b}'_{h}]$ is void or saturated in $\boldsymbol{\mu}^{\textnormal{I}}$. We can always take $\delta$ small enough.
\vspace{0.2cm}

\noindent \textsc{Claim.} There exists $\eta > 0$ depending only on $\delta$ and the constants in the assumptions, such that $\big[\hat{a}'_{h},\hat{a}'_{h} + \frac{\delta}{2}\big] \cup \big[\hat{b}'_{h} - \frac{\delta}{2},\hat{b}'_{h}\big]$ is void or saturated in $\boldsymbol{\mu}^{\textnormal{I}}$ and in $\boldsymbol{\mu}^{\textnormal{II}}$ as soon as $|\hat{n}_{h}^{\textnormal{I}} - \hat{n}_{h}^{\textnormal{II}}| \leq \eta$.

\vspace{0.2cm}

We temporarily admit the claim and finish the proof of Lemma~\ref{Lemma_continuity_filling_fractions}. Choose $\varrho$ in \eqref{ppmpmp} so that it is supported in $\big[\hat{a}'_{h},\hat{a}'_{h} +\frac{\delta}{8}\big]$. Choose a large enough constant $D > 0$, so that all the measures of interest are supported inside $[-D,D]$. Choose a smooth function $\chi$ on the real line taking values in $[0,1]$, which vanishes on $\big[\hat{a}'_{h} - \frac{\delta}{4},\hat{a}'_{h} + \frac{\delta}{4}\big]$ and outside $[-D - \delta,D + \delta]$, and is equal to $1$ on $\big[-D,\hat{a}'_{h} - \frac{\delta}{2}\big] \cup \big[\hat{a}'_{h} +\frac{\delta}{2},
\hat{b}'_{h}\big]$. Since $\mu^{\textnormal{I}}_h - \mu^{\textnormal{II}}_h$ vanishes on $\big[\hat{a}'_{h},\hat{a}'_{h} + \frac{\delta}{2}\big]$, we have for any $g \in [H]$
\begin{equation}
\label{642ein}\int_{\hat{a}'_{g}}^{\hat{b}'_{g}}\!\!\int_{\hat{a}'_{h}}^{\hat{b}'_{h}} \log|x - y|\,\dd\varrho_g(x)\dd(\mu^{\textnormal{I}}_h - \mu^{\textnormal{II}}_h)(y) = \int_{\hat{a}'_{g}}^{\hat{b}'_{g}}\!\!\int_{\hat{a}'_{h}}^{\hat{b}'_{h}} \log|x - y|\,\chi(y)\,\dd\varrho_g(x)\dd\big(\mu^{\textnormal{I}}_h - \mu^{\textnormal{II}}_h + (\hat{n}^{\textnormal{II}}_{h} - \hat{n}^{\textnormal{I}}_{h})\varrho_h\big)(y).
\end{equation}
Note that the function
\[
y \mapsto \chi(y) \int_{\hat{a}'_{g}}^{\hat{b}'_{g}} \log|x - y|\,\dd\varrho_g(x)
\]
is continuously differentiable on the real line and has a compact support. Using that $\mu^{\textnormal{I}} - \mu^{\textnormal{II}} + (\hat{n}_{h}^{\textnormal{II}} - \hat{n}_{h}^{\textnormal{I}})\varrho$ has total mass $0$, we apply Lemma~\ref{Lemma_linear_through_distance} to the right-hand side of \eqref{642ein} and bound its absolute value from above by
\begin{equation}
\label{643ein} M \cdot\mathfrak{D}\big[\boldsymbol{\mu}^{\textnormal{I}} + (\hat{n}_{h}^{\textnormal{II}} - \hat{n}_{h}^{\textnormal{I}})\boldsymbol{\varrho},\boldsymbol{\mu}^{\textnormal{II}}\big]
\end{equation}
for some $M > 0$. Combining \eqref{643ein} with \eqref{ppmpmp} and \eqref{640ein}, we get
\begin{equation}
\label{644ein} \D\big[\boldsymbol{\mu}^{\textnormal{I}} + (\hat{n}^{\textnormal{II}}_{h} - \hat{n}_{h}^{\textnormal{I}})\boldsymbol{\varrho},\boldsymbol{\mu}^{\textnormal{II}}\big] \leq M\cdot \mathfrak{D}\big[\boldsymbol{\mu}^{\textnormal{I}} + (\hat{n}_{h}^{\textnormal{II}} - \hat{n}_{h}^\textnormal{I})\boldsymbol{\varrho},\boldsymbol{\mu}^{\textnormal{II}}\big] + M \cdot (\hat{n}_{h}^{\textnormal{I}} - \hat{n}_{h}^{\textnormal{II}})^2.
\end{equation}
We treat \eqref{644ein} as a quadratic inequality for the unknown nonnegative quantity $\mathfrak{D}\big[\boldsymbol{\mu}^{\textnormal{I}} + (\hat{n}_{h}^{\textnormal{II}} - \hat{n}_{h}^{\textnormal{I}})\boldsymbol{\varrho},\boldsymbol{\mu}^{\textnormal{II}}\big]$. Then \eqref{644ein} implies
\[
\mathfrak{D}\big[\boldsymbol{\mu}^{\textnormal{I}} + (\hat{n}_{h}^{\textnormal{II}} - \hat{n}_{h}^{\textnormal{I}})\boldsymbol{\varrho},\boldsymbol{\mu}^{\textnormal{II}}\big] \leq M \cdot |\hat{n}_{h}^{\textnormal{I}} - \hat{n}_{h}^{\textnormal{II}}|
\]
for some $M > 0$. Combining this inequality with Lemma~\ref{Lemma_linear_through_distance}, we get
\[
\bigg|\int_{\hat{a}'_{h}}^{\hat{b}'_{h}} f(x)\,\dd\big(\mu_h^{\textnormal{I}} - \mu_h^{\textnormal{II}} + (\hat{n}_{h}^{\textnormal{II}} - \hat{n}_{h}^{\textnormal{I}})\varrho_h\big)(x)\bigg| \leq M\cdot |\!|f|\!|_{\frac{1}{2}}\,|\hat{n}_{h}^{\textnormal{I}} - \hat{n}_{h}^{\textnormal{II}}|.
\]
Clearly, we also have
\[
\bigg|\int_{\hat{a}'_{h}}^{\hat{b}'_{h}} f(x)\,\dd\big((\hat{n}_{h}^{\textnormal{II}} - \hat{n}_{h}^{\textnormal{I}})\varrho_h\big)(x)\bigg| \leq M\cdot |\!|f|\!|_{\infty}\,|\hat{n}_{h}^{\textnormal{I}} - \hat{n}_{h}^{\textnormal{II}}|
\]
for some $M > 0$. Combining the two previous equations leads to the desired \eqref{eq_continuity_filling_fractions}.

The proof will be complete after we justify the claim. As the argument is very similar to the one used for the claim in the proof of Lemma~\ref{Lemma_continuity_potential}, we omit some details. Using \eqref{eq_x86} and \eqref{ppmpmp} we get
\[
\D\big[\boldsymbol{\mu}^{\textnormal{I}} + (\hat{n}_{h}^{\textnormal{II}} - \hat{n}_{h}^{\textnormal{I}})\boldsymbol{\varrho},\boldsymbol{\mu}^{\textnormal{II}}\big] \leq M \cdot |\hat{n}_h^{\textnormal{I}} - \hat{n}_{h}^{\textnormal{II}}|
\]
for some $M > 0$, which implies that the effective potentials constructed from the measures $\boldsymbol{\mu}^{\textnormal{I}} + (\hat{n}_{h}^{\textnormal{II}} - \hat{n}_{h}^{\textnormal{I}})\boldsymbol{\varrho}$ and $\boldsymbol{\mu}^{\textnormal{II}}$ are close to each other. It follows that the effective potentials of $\boldsymbol{\mu}^{\textnormal{I}}$ and $\boldsymbol{\mu}^{\textnormal{II}}$ are also close to each other. Then, using Conditions 2.\ and 3.\ in Assumption~\ref{Assumption_C} and the characterization of voids, bands and saturations in Theorem~\ref{Theorem_equi_charact_repeat_2}, we obtain the claim. \end{proof}

\begin{proof}[Proof of Lemma~\ref{Lemma_continuity_end_points_void}]
 Both $\boldsymbol{\mu}^\textnormal{I}$ and $\boldsymbol{\mu}^\textnormal{II}$ are minimizers of energy functionals; the only difference between these minimization problems is that the first one allows measures to have mass
 inside $[\hat a^{\prime\,\textnormal{I}}_{h}, \hat b^{\prime\,\textnormal{I}}_{h}]$, while the second one allows measures to have mass
 inside $[\hat a^{\prime\,\textnormal{II}}_{h}, \hat b^{\prime\,\textnormal{II}}_{h}]$, where $\hat b^{\prime\,\textnormal{I}}_h=\hat b^{\prime\,\textnormal{II}}_h$, but
 $\hat a^{\prime\,\textnormal{I}}_h\neq \hat a^{\prime\,\textnormal{II}}_h$.

 We claim that $\boldsymbol{\mu}^\textnormal{I}$ satisfies the characterization of Theorem~\ref{Theorem_equi_charact_repeat_2} for the second variational datum which would imply $\boldsymbol{\mu}^\textnormal{I}=\boldsymbol{\mu}^\textnormal{II}$ by uniqueness. If $\hat a^{\prime\,\textnormal{I}}_h \in (\hat a^{\prime\,\textnormal{II}}_h - \frac{1}{M},\hat a^{\prime\,\textnormal{II}}_h)$ and $\hat a^{\prime\,\textnormal{II}}_h$ belongs to a void of $\boldsymbol{\mu}^\textnormal{I}$ --- the latter is automatic for $M>0$ large enough due to Assumption~\ref{Assumption_C} --- then there is nothing to check, as each characterizing condition of Theorem~\ref{Theorem_equi_charact_repeat_2} for the second variational datum is implied by the same condition for the first variational datum. On the other hand, if $\hat a^{\prime\,\textnormal{II}}_h < \hat a^{\prime\,\textnormal{I}}_h$, then we need to check that $V^{\textnormal{eff},\textnormal{I}}_h(x)\geq v_h^{\textnormal{I}}$ for $x\in [\hat a^{\prime\,\textnormal{II}}_h,\hat a^{\prime\,\textnormal{I}}_h]$. This follows from the uniform continuity of $V^{\textnormal{eff},\textnormal{I}}_h(x)$ with constants depending only on those in the assumptions, and the fact that $V^{\textnormal{eff},\textnormal{I}}_h(x)-v_h^{\textnormal{I}}$ is positive and bounded away from $0$ at $x=\hat a^{\prime\,\textnormal{I}}_h$ by Assumption~\ref{Assumption_C}.2.
 \end{proof}

\begin{proof}[Proof of Lemma~\ref{Lemma_continuity_end_points_saturated}] First, consider the case $\hat a_h^{\prime\,\textnormal{I}} < \hat a_{h}^{\prime\,\textnormal{II}}$. This guarantees that for $\hat{a}^{\prime\,\textnormal{II}}_{h} - \hat{a}^{\prime\,\textnormal{I}}_h < \eps$, the whole segment $[\hat{a}_{h}^{\prime\,\textnormal{I}},\hat{a}_{h}^{\prime\,\textnormal{II}}]$ is saturated for $\boldsymbol{\mu}^{\textnormal{I}}$. Let $\lambda_h$ be the Lebesgue measure of density $\frac{1}{\theta_{h,h}}$ supported on $[\hat{a}^{\prime\,\textnormal{I}}_{h},\hat{a}_{h}^{\prime\,\textnormal{II}}]$, let $\lambda_g = 0$ for $g \neq h$, and set $\tilde{\boldsymbol{\mu}} = \boldsymbol{\mu}^{\textnormal{I}} - \boldsymbol{\lambda}$. This measure is supported on the segments corresponding to the set of parameters $\textnormal{II}$, and has new segment filling fractions $\tilde{n}_{g} = \hat{n}_{g}$ for $g \neq h$ and $\tilde{n}_{h} = \hat{n}_{h} - \frac{1}{\theta_{h,h}}(\hat{a}_{h}^{\prime\,\textnormal{II}} - \hat{a}_{h}^{\prime\,\textnormal{I}})$. Let us formulate a minimization problem solved by $\tilde{\boldsymbol{\mu}}$. From the characterization of $\boldsymbol{\mu}^{\textnormal{I}}$ in Theorem~\ref{Theorem_equi_charact_repeat_2} which we restrict to the segments corresponding to the parameters $\textnormal{II}$, there exists an $H$-tuple of constants $\boldsymbol{v}$ such that for any $g_0 \in [H]$ and $x \in [\hat{a}_{g_0}^{\prime\,\textnormal{I}},\hat{b}_{g_0}^{\prime\,\textnormal{I}}]$, we have
\[
\tilde{V}_{g_0}(x) - \sum_{g = 1}^H 2\theta_{g_0,g} \int_{\hat{a}_{g}^{\prime\,\textnormal{II}}}^{\hat{b}_{g}^{\prime\,\textnormal{II}}} \log|x - y|\dd\tilde{\mu}_{g}(y)\,\,\begin{cases} \leq v_{g_0}^{\textnormal{I}} & \textnormal{if} \,\,\tilde{\mu}_{g_0}(x) = \theta_{g_0,g_0}^{-1},\\ = v_{g_0}^{\textnormal{I}} & \textnormal{if}\,\,0 < \tilde{\mu}_{g_0}(x) < \theta_{g_0,g_0}^{-1}, \\ \geq v_{g_0}^{\textnormal{I}} & \textnormal{if}\,\,\tilde{\mu}_{g_0}(x) = 0, \end{cases}
\]
where we have introduced a new potential
\begin{equation*}
\begin{split}
\tilde{V}_{g_0}(x) & = V_{g_0}^{\textnormal{I}}(x) - 2\theta_{h,g_0} \int_{\hat{a}_{h}^{\prime\,\textnormal{I}}}^{\hat{a}_{h}^{\prime\,\textnormal{II}}} \log|x - y| \lambda(y) \dd y= V_{g_0}^{\textnormal{I}}(x) + 2\Delta V(x) = V_{g_0}^{\textnormal{II}}(x) + 2\delta_{h \neq g_0} \Delta V(x) ,\\
\Delta V(x) & = (\hat{a}_{h}^{\prime\,\textnormal{II}} - x)\log|\hat{a}_{h}^{\prime\,\textnormal{II}} - x| + (x - \hat{a}_{h}^{\prime\,\textnormal{I}})\log|x - \hat{a}_{h}^{\prime\,\textnormal{I}}| + (\hat{a}^{\prime\,\textnormal{II}}_{h} - \hat{a}^{\prime\,\textnormal{I}}_{h}).
\end{split}
\end{equation*}
Therefore, $\tilde{\boldsymbol{\mu}}$ is the equilibrium measure for parameters that differ from those of $\boldsymbol{\mu}^{\textnormal{II}}$ only by the potential and the filling fractions of the $h$-th segment, and we have
\begin{equation*}
\begin{split}
\sup_{g \in [H]} \sup_{x \in [\hat{a}_{g}^{\prime\,\textnormal{II}},\hat{b}_{g}^{\prime\,\textnormal{II}}]} \big|(\tilde{V}_{g} - V^{\textnormal{II}}_{g})(x)\big| & = 2\,\sup_{g \neq h} \sup_{x \in [\hat{a}_{g}^{\prime\,\textnormal{II}},\hat{b}_{g}^{\prime\,\textnormal{II}}]} |\Delta V(x)| \leq M \cdot (\hat{a}_{h}^{\prime\,\textnormal{II}} - \hat{a}_{h}^{\prime\,\textnormal{I}}), \\
\sup_{g \in [H]} \sup_{x \in [\hat{a}_{g}^{\prime\,\textnormal{II}},\hat{b}_{g}^{\prime\,\textnormal{II}}]} \big|(\tilde{V}_{g} - V_{g}^{\textnormal{II}})'(x)\big| & = 2\,\sup_{g \neq h} \sup_{x \in [\hat{a}_{g}^{\prime\,\textnormal{II}},\hat{b}_{g}^{\prime\,\textnormal{II}}]} \big|(\Delta V)'(x)\big| \leq M \cdot (\hat{a}_{h}^{\prime\,\textnormal{II}} - \hat{a}_{h}^{\prime\,\textnormal{I}}),\\
|\tilde{n}_{h} - \hat{n}_{h}^{\textnormal{II}}| = 2\,\frac{\hat{a}_{h}^{\prime\,\textnormal{II}} - \hat{a}_{h}^{\prime\,\textnormal{I}}}{\theta_{h,h}} & \leq M \cdot (\hat{a}_{h}^{\prime\,\textnormal{II}} - \hat{a}_{h}^{\prime\,\textnormal{I}}),
\end{split}
\end{equation*}
for some constant $M > 0$. Recall that $\boldsymbol{\mu}^{\textnormal{I}}$ satisfies Assumption~\ref{Assumption_C}. Its restriction $\tilde{\mu} = \mu|_{\amsmathbb{A}^{\textnormal{II}}}$ also satisfies Assumption~\ref{Assumption_C}, because by construction of $\tilde{V}$ the effective potential of $\tilde{\boldsymbol{\mu}}$ is equal to the restriction to $\amsmathbb{A}^{\textnormal{II}}$ of the effective potential of $\boldsymbol{\mu}^{\textnormal{I}}$. Thus, we can apply Lemma~\ref{Lemma_continuity_potential} and \ref{Lemma_continuity_filling_fractions}. We deduce the existence of a perhaps larger constant $M > 0$ such that, if $|\hat{a}_{h}^{\prime\,\textnormal{I}} - \hat{a}_{h}^{\prime\,\textnormal{II}}| < \frac{1}{M}$, then the following bound holds
\begin{equation}
\label{fguim}\bigg|\int_{\amsmathbb{R}} f(x)\dd(\tilde{\mu} - \mu^{\textnormal{II}})(x)\bigg| \leq M \cdot \big(|\!|f|\!|_{\infty} + |\!|f|\!|_{\frac{1}{2}}\big) \cdot \big|\hat{a}_{h}^{\prime\,\textnormal{I}} - \hat{a}_h^{\prime\,\textnormal{II}}\big|.
\end{equation}
 The proof is completed by inserting \eqref{fguim} in the inequality
\begin{eqnarray*}
\bigg|\int_{\amsmathbb{R}} f(x)\dd(\mu^{\textnormal{I}} - \mu^{\textnormal{II}})(x)\bigg| & \leq & \bigg|\int_{\amsmathbb{R}} f(x)\dd\lambda(x)\bigg| + \bigg|\int_{\amsmathbb{R}} f(x)\dd(\tilde{\mu} - \mu^{\textnormal{II}})(x)\bigg| \\
& \leq & M \cdot |\!|f|\!|_{\infty} \cdot \big|\hat{a}_{h}^{\prime\,\textnormal{I}} - \hat{a}_{h}^{\prime\,\textnormal{II}}\big| + \bigg|\int_{\amsmathbb{R}} f(x)\dd(\tilde{\mu} - \mu^{\textnormal{II}})(x)\bigg|,
\end{eqnarray*}
which holds for some constant $M > 0$.

\smallskip

For the second case $\hat a_h^{\prime\,\textnormal{I}} < \hat a_{h}^{\prime\,\textnormal{II}}$ the argument is the same with $\tilde{\mu}$ now defined as a sum of $\mu^{\textnormal{I}}$ and the Lebesgue measure of density $\frac{1}{\theta_{h,h}}$ supported on $[\hat a_h^{\prime\,\textnormal{I}}, \hat a_{h}^{\prime\,\textnormal{II}}]$. The only difference is that to show the off-criticality of $\tilde \mu$ we also need to appeal to the uniform continuity of the effective potential on the segment $[\hat a_h^{\prime\,\textnormal{I}}, \hat a_{h}^{\prime\,\textnormal{II}}]$, as we did in the proof of Lemma~\ref{Lemma_continuity_end_points_void}.
\end{proof}

\section{Preservation of the off-criticality}

\label{Section_off_criticality}

In this section we show that the off-criticality condition of Assumption~\ref{Assumption_C}
is preserved under small changes of variational data. In contrast to Section~\ref{Section-ParReg}, the analyticity of $U_h$ is significantly used here.

\begin{theorem} \label{Theorem_off_critical_neighborhood}
Take two variational data satisfying Assumptions~\ref{Assumption_A} and \ref{Assumption_B}. Assume that for any $h \in [H]$ the complex domain $\amsmathbb{M}_h$ in Assumption~\ref{Assumption_B} is common to the two variational data. Then, there exists $\eps > 0$ depending only on the constants in the assumptions such that, if the variational datum $\textnormal{I}$ satisfies the off-criticality Assumption~\ref{Assumption_C} and for any $g,h \in [H]$ the inequalities
\begin{equation*}
\begin{split}
\max\big(|\hat a_h^{\prime\,\textnormal{I}}-\hat a_h^{\prime\,\textnormal{II}}|\,,\,|\hat b_h^{\prime\,\textnormal{I}}-\hat b_h^{\prime\,\textnormal{II}}|\,,\,|\hat n_h^{\textnormal{I}}-\hat n_h^\textnormal{II}|\,,\,|\theta_{g,h}^{\textnormal{I}} - \theta_{g,h}^{\textnormal{II}}|\big) & \leq \eps, \\
 \int_{\max (\hat a_h^{\prime\,\textnormal{I}},\hat a_h^{\prime\,\textnormal{II}})}^{\min(\hat b_h^{\prime\,\textnormal{I}},\hat b_h^{\prime\,\textnormal{II}})} \big|V_h^{\textnormal{I}}(x)-V_h^\textnormal{II}(x)\big| \dd x & \leq \eps, \\
\sup_{z\in \amsmathbb{M}_h} \big|\phi^{\pm,\textnormal{I}}_h(z)-\phi^{\pm,\textnormal{II}}_h(z)\big| & \leq \eps,
\end{split}
\end{equation*}
then the variational datum $\textnormal{II}$ also satisfies Assumption~\ref{Assumption_C}.
\end{theorem}
\begin{proof}
 We first show that $\boldsymbol{\mu}^\textnormal{I}$ and $\boldsymbol{\mu}^\textnormal{II}$ are close to each other. For the changes in
 the filling fractions $\hat n_h$ and endpoints $\hat a'_h$, $\hat b'_h$, this is the content of Lemmata
 \ref{Lemma_continuity_filling_fractions}--\ref{Lemma_continuity_end_points_saturated}.
 For the change in the potential we cannot directly use Lemma~\ref{Lemma_continuity_potential},
 as it used a different topology on potentials. However, \eqref{eq_x80} in its proof is enough for us, as it shows that
 $\D[\boldsymbol{\mu}^\textnormal{I},\boldsymbol{\mu}^{\textnormal{II}}]< C\eps$. Combining all these propositions together, we obtain the following statement:
 For each $\delta>0$ there exists $\eps_0>0$, such that if $\eps<\eps_0$, then for
 every measurable function $f$ on the real line with finite $|\!|f|\!|_{\frac{1}{2}}$ and $|\!|f|\!|_{\infty}$ norms we have
 \begin{equation}
 \label{eq_x110}
 \bigg|\int_{\amsmathbb R} f(x)\dd(\mu^{\textnormal{I}}-\mu^\textnormal{II})(x)\bigg|\leq \delta\cdot
 \big(|\!|f|\!|_{\frac{1}{2}}+|\!|f|\!|_{\infty}\big).
 \end{equation}
 In particular \eqref{eq_x110} implies --- by following the same argument as the one
 used in the proof of the Claim inside Lemma
 \ref{Lemma_continuity_potential} --- that the effective potentials of $\boldsymbol{\mu}^{\textnormal{I}}$
 and $\boldsymbol{\mu}^{\textnormal{II}}$ become arbitrary close as $\eps\rightarrow 0$. Together with the characterization
 of the equilibrium measure of Theorem~\ref{Theorem_equi_charact_repeat_2}, this implies that all
 off-criticality properties of Assumption~\ref{Assumption_C} hold for
 $\boldsymbol{\mu}^{\textnormal{II}}$ away from the endpoints of the bands of $\boldsymbol{\mu}^{\textnormal{I}}$. In particular, if the segment
 $[u,v]$ is void (or saturated) in $\boldsymbol{\mu}^{\textnormal{I}}$ and $c>0$, then for small
 enough $\eps$ depending only $c$ and the constants in the assumptions, the segment $[u+c,v-c]$ is void (or saturated) in $\boldsymbol{\mu}^{\textnormal{II}}$.

 For the off-criticality properties near the endpoints of
 the bands we need a separate argument that we now present. Choose one band of $\boldsymbol{\mu}^{\textnormal{I}}$,
 that we denote $[\alpha^{\textnormal{I}},\beta^{\textnormal{I}}]$, and which is included in $[\hat{a}_h',\hat{b}_h']$ for some $h \in [H]$.
 Fix a contour $\gamma_h \subset \amsmathbb{M}_h$ surrounding $[\alpha^{\textnormal{I}},\beta^{\textnormal{I}}]$ as described in Definition~\ref{def_contourgammah}, and such that all other bands of $\boldsymbol{\mu}^\textnormal{I}$ are outside $\gamma_h$. In particular, intersections of
 $\gamma_h$ with the real axis
 are inside voids or saturations of $\boldsymbol{\mu}^\textnormal{I}$.
 Since we have already shown that the voids (respectively, saturations) of $\boldsymbol{\mu}^\textnormal{I}$ and $\boldsymbol{\mu}^{\textnormal{II}}$ are close to each other, if $\eps$ is small enough, the intersections of $\gamma_h$ with the real axis are
 also inside voids or saturations of $\boldsymbol{\mu}^{\textnormal{II}}$.

 Let $q_h^{\pm,\textnormal{I}}$ and $q_h^{\pm,\textnormal{II}}$ denote the functions $q_h^\pm(z)$ of Definition~\ref{GQdef}
 corresponding to $\boldsymbol{\mu}^\textnormal{I}$ and $\boldsymbol{\mu}^\textnormal{II}$, respectively. We claim that for each
 $\delta>0$ there exists $\eps_0>0$ such that for all $\eps \in (0,\eps_0)$
 \begin{equation}
 \label{eq_x111}
 \sup_{z \in \gamma_h} \big|q_h^{+,\textnormal{I}}(z)-q_h^{+,\textnormal{II}}(z)\big|<\delta.
 \end{equation}
 Indeed, given the definition of $q_h^+(z)$, for $z$ outside the real axis \eqref{eq_x111} is a direct corollary of
 \eqref{eq_x110}. For $z$ on the real axis we combine \eqref{eq_x111} with the fact that
neighbors of $z$ on the real axis belong to a void (or a saturation) simultaneously for $\boldsymbol{\mu}^{\textnormal{I}}$ and $\boldsymbol{\mu}^{\textnormal{II}}$, and therefore they contribute in the same way to the Stieltjes
 transforms $\Gm_{\mu^{\textnormal{I}}(z)}$, $\Gm_{\mu^\textnormal{II}(z)}$ entering the
 definition of $q_h^{+,\textnormal{I}}(z)$, $q_h^{+,\textnormal{II}}(z)$.
 Recall that $q_h^{+,\textnormal{I}}(z)$ and $q_h^{+,\textnormal{II}}(z)$ are holomorphic by Theorem~\ref{Theorem_regularity_density}. Hence, by the maximum principle for harmonic functions, \eqref{eq_x111}
 implies the same uniform bound for $z$ in the compact domain bounded by $\gamma_h$.

 Next, adjusting $\gamma_h$ if necessary, and using the lower bound on $s_h(z)$ of Proposition~\ref{Proposition_q_bounds}, we
 observe that the only zeros of $(q_h^{-,\,\textnormal{I}}(z))^2$ are at $\alpha^{\textnormal{I}}$ and $\beta^{\textnormal{I}}$. Recall that
 \begin{equation}
 \label{eq_x198}
\forall \textnormal{J} \in \{\textnormal{I},\textnormal{II}\}\qquad (q_h^{-,\,\textnormal{J}}(z))^2=(q_{h}^{+,\,\textnormal{J}}(z))^2- 4\phi^{+,\,\textnormal{J}}_h(z)
 \phi^{-,\textnormal{J}}_h(z).
 \end{equation}
 Therefore, combining \eqref{eq_x111} and its extension inside $\gamma_h$ with the Rouch\'{e} theorem, we conclude that for
 small enough $\eps$, the function $(q_h^{-,\rm II}(z))^2$ has only two zeros inside $\gamma_h$, which are simple
 and close to $\alpha^{\textnormal{I}}$ and $\beta^{\textnormal{II}}$, respectively. Let us denote these zeros
 $\alpha^{\textnormal{II}}$ and $\beta^{\textnormal{II}}$, respectively. By Theorem~\ref{Theorem_regularity_density}, the zeros are endpoints of the band of $\boldsymbol{\mu}^\textnormal{II}$.

 By Theorem~\ref{Theorem_regularity_density}, $s^{\textnormal{II}}_h(z)$ is holomorphic inside the contour $\gamma_h$.
 By \eqref{eq_x111}, \eqref{eq_x198} and the (just proven) closeness of $\alpha^{\textnormal{II}}$ and $\beta^{\textnormal{II}}$ to $\alpha^{\textnormal{I}}$ and $\beta^{\textnormal{I}}$, respectively, $s^{\textnormal{II}}_h(z)$ is uniformly close to $s^{\textnormal{I}}_h(z)$ for $z \in \gamma_h$. In particular, $s^{\textnormal{II}}_h(z)$ is bounded away from $0$ on $\gamma_h$.
 Since $(q_h^{-,\,\textnormal{II}}(z))^2$ has only two zeros inside $\gamma_h$, the expression $s^{\textnormal{II}}_h(z)$
 has none. Applying the maximum principle to $-|s^\textnormal{II}_h(z)|$, this implies that $s^\textnormal{II}_h(z)$ is bounded away from $0$ by a uniform
 constant inside $\gamma_h$.

 At this point we can use the formula \eqref{eq_density_exponent}. More precisely, we use its corollary given in \eqref{eq_x175} and valid for $x\in (\alpha^{\textnormal{II}},\beta^{\textnormal{II}})$:
 \begin{equation}
 \label{eq_x176}
\big|s_h^\textnormal{II}(x)\big|= \big|1-\exp\big(2\ii\pi\theta_{h,h}^{\textnormal{II}}\,\mu^\textnormal{II}_h(x)\big)\big| \cdot
\frac{
\sqrt{\big|\phi^{-,\,\textnormal{II}}_h(x) \cdot \phi^{+,\,\textnormal{II}}_h(x)\big|}}{|\psi^\textnormal{II}(x)|} \cdot \frac{1}{\big|\sigma_h^\textnormal{II}(x^+)\big|}.
\end{equation}
Since the three functions
\[
\big|s_h^{\textnormal{II}}(z)\big|,\qquad \frac{\sqrt{\big|\phi^{-,\,\textnormal{II}}_h(x) \cdot \phi^{+,\,\textnormal{II}}_h(x)\big|}}{\big|\psi^\textnormal{II}(x)\big|},\qquad \left|\frac{\sigma^\textnormal{II}_h(x)}{\sqrt{(x-\alpha^\textnormal{II})(\beta^\textnormal{II}-x)}}\right|
\]
are bounded away from $0$ and infinity --- for the last one we use the already proven fact that endpoints of the bands for the parameters $\rm II$ are close to those for the parameters $\rm I$ --- the identity \eqref{eq_x176} implies that for
$C>0$
\[
C^{-1}\,\sqrt{\big|x-\alpha^{\textnormal{II}}\big|} \leq \big|1-\exp\big(2\ii\pi\theta^{\textnormal{II}}_{h,h}\,
 \mu_h^{\textnormal{II}}(x)\big)\big| \leq C\,\sqrt{\big|x-\alpha^{\textnormal{II}}\big|}
\]
in a neighborhood of $\alpha^{\textnormal{II}}$, and
\[
C^{-1}\,\sqrt{\big|x-\beta^{\textnormal{II}}\big|} \leq \big|1-\exp\big(2\ii\pi\theta_{h,h}^{\textnormal{II}}\,
 \mu_h^{\textnormal{II}}(x)\big)\big| \leq C\,\sqrt{\big|x-\beta^\textnormal{II}\big|}
\]
in a neighborhood of $\beta^{\textnormal{II}}$. This implies the square-root behavior of
the density $\mu_h^{\textnormal{II}}(x)$ near the endpoints of the bands, as claimed in Assumption~
\ref{Assumption_C}.5 and 6.

\smallskip

It now remains to establish the bound of Assumption~\ref{Assumption_C}.2 and~\ref{Assumption_C}.3  on the
effective potential in voids and saturations --- still in a neighborhood of
$\alpha^{\textnormal{II}}$ or $\beta^{\textnormal{II}}$, because away from these points the bound is implied by \eqref{eq_x110} and the bound for $\boldsymbol{\mu}^\textnormal{I}$. We will only consider a neighborhood of $\alpha^{\textnormal{II}}$, the argument for $\beta^{\textnormal{II}}$ is similar. First, suppose that $\alpha^{\textnormal{II}}$ is an endpoint of a void rather than a saturation. The $\frac{1}{2}$-H\"older continuity of the density $\mu^{\textnormal{II}}(x)$ implies that the derivative of the effective potential $\partial_{x}V^{\textnormal{eff},\textnormal{II}}_h(x)$ is a continuous function of $x$ in a neighborhood of $\alpha^{\textnormal{II}}$. Inside the band,
\textit{i.e.} to the right of $\alpha^{\textnormal{II}}$, this derivative vanishes by Theorem
\ref{Theorem_equi_charact_repeat_2}. Therefore, the behavior of the effective potential for
$x$ in a left neighborhood of $\alpha^{\textnormal{II}}$ is determined by the sign of its second
derivative.

Fix a small $c>0$ and write
\begin{equation}
\label{eq_x113}
\forall x \leq \alpha^{\textnormal{II}}\qquad V^{\textnormal{eff},\textnormal{II}}_h(x)=-\theta_{h,h}^{\textnormal{II}} \int_{\alpha^\textnormal{II}}^{\alpha^\textnormal{II}+c} \mu^\textnormal{II}_h(y)\log
 (y-x) \dd y + \Delta V(x),
\end{equation}
where $\Delta V(x)$ is twice-continuously differentiable in a neighborhood of
$\alpha^{\textnormal{II}}$ with uniformly bounded second derivative. Introduce the function
\[
\forall x < 0 \qquad f(x) :=-\partial_{x}^2\left(\int_{0}^{c} \sqrt{y}\log(y-x) \dd y\right)
=\int_{0}^{c} \frac{\sqrt{y}\dd y}{(y-x)^2}.
\]
This function is positive integrable and goes to $+\infty$
as $x \rightarrow 0^-$. Therefore, for $x$ in a small enough left neighborhood of $\alpha^{\textnormal{II}}$,
\eqref{eq_x113} implies the existence of $C>0$ such that
\begin{equation}
\label{eq_x114}\partial_{x}^2V^{\textnormal{eff},\textnormal{II}}_h(x)>-C + \frac{f(x-\alpha^{\textnormal{II}})}{C}.
\end{equation}
The key feature of \eqref{eq_x114} is that the right-hand side is strictly positive for $x$ near $\alpha^{\textnormal{II}}$.
Integrating \eqref{eq_x114} twice, we get the desired lower bound for the effective
potential near $\alpha^{\textnormal{II}}$.

\smallskip

We proceed to the second case and $\alpha^{\textnormal{II}}$ is an endpoint of a saturated
region. The effective potential is still continuously differentiable in a
neighborhood of $\alpha^{\textnormal{II}}$ and we need to analyze its second derivative. We
replace \eqref{eq_x113} with
\begin{equation}
\label{eq_x115}
\forall x \leq \alpha^{\textnormal{II}}\qquad V^{\textnormal{eff},\textnormal{II}}_h(x)=-\theta_{h,h}^{\textnormal{II}} \int_{\alpha^{\textnormal{II}}}^{\alpha^{\textnormal{II}}+c} \bigg(\mu^{\textnormal{II}}_h(y)-\frac{1}{\theta_{h,h}}\bigg)\log
 (y-x) \dd y + \Delta V(x),
\end{equation}
where $\Delta V(x)$ is twice-continuously differentiable in a neighborhood of $\alpha^{\textnormal{II}}$. At this point we repeat the argument from the first case.
\end{proof}

\section{Parametric smooth regularity}
\label{Section_Smoothnessparam}
In this section we would like to upgrade the estimates of Propositions
\ref{Lemma_continuity_potential},
\ref{Lemma_continuity_filling_fractions},
\ref{Lemma_continuity_end_points_saturated} and prove
smoothness of the equilibrium measure with respect to the involved
parameters. This upgrade will rely on the results obtained in Part~\ref{Part_Master_equation} concerning the continuous invertibility --- via an operator $\boldsymbol{\Upsilon}$ --- of certain systems of equations. The construction of the operator $\boldsymbol{\Upsilon}$ and the study of its properties is completely independent of the asymptotic analysis that we are carrying in Part~\ref{Part_Asymptotic} of the book, justifying that it is discussed separately. The results from Part~\ref{Part_Master_equation} that we use in the present section will be mentioned when they appear.

We prove the differentiability under some additional requirements. The first one is that each segment $[\hat a'_h, \hat b'_h]$ has exactly one band of the equilibrium measure. This is not really restrictive, as the variational data can be localized near bands (\textit{cf.} Chapter~\ref{Chapter_conditioning}). The second requirement is the off-criticality Assumption~\ref{Assumption_C}, which was already imposed in Section~\ref{Section-ParReg}. It is a crucial ingredient of our arguments, yet there might be situations in which the differentiability is still true in the critical cases. However, it is impossible to completely drop out both one band and off-criticality assumptions, as we expect non-differentiable behavior of the density of the equilibrium measure at phase transitions, \textit{e.g.} at values of parameters where two bands merge together.

\medskip

For the clarity of the exposition, we proceed in two steps: Proposition
\ref{Proposition_differentiability_filling_fraction} deals in details with
derivatives with respect to the segment filling fractions only and has slightly more restrictive assumptions, while Theorem
\ref{Theorem_differentiability_full} discusses the full differentiability with respect to all
parameters and without this restrictive assumption.

\begin{proposition}
\label{Proposition_differentiability_filling_fraction}
Consider a variational datum satisfying Assumptions~\ref{Assumption_A}, \ref{Assumption_B} and \ref{Assumption_C}, with filling fractions $\hat{\boldsymbol{n}}^{(0)}$. Additionally, assume that the following properties hold for any $h \in [H]$.
\begin{enumerate}
\item The equilibrium measure has only one band $[\alpha_h,\beta_h]$ in $[\hat{a}'_{h},\hat{b}'_{h}]$.
\item The function $s_h$ of Definition~\ref{GQdef2} satisfies $|s_h(z)| \geq \frac{1}{C}$ for any $z \in \amsmathbb{M}_{h}$.
\item If $\hat{a}'_{h} \in \amsmathbb{S}_{h}$, then we require a double zero: $\phi_h^{-}(\hat{a}'_{h})=[\phi_h^{-}]'(\hat{a}'_{h}) = 0$. Likewise, if $\hat{b}'_{h} \in \amsmathbb{S}_{h}$, then we require a double zero: $\phi_h^{+}(\hat{b}'_{h})=[\phi_h^{+}]'(\hat{b}'_{h}) = 0$.
\end{enumerate}
Then, there exists $\delta >0$ depending only on the constants in the assumptions such, for any variational datum differing from the original one only by the vector of filling fraction $\hat{\boldsymbol{n}}$ which satisfies $|\!|\hat{\boldsymbol{n}} - \hat{\boldsymbol{n}}^{(0)}|\!|_{\infty} \leq \delta$, the properties in \ref{Assumption_A}, \ref{Assumption_B} and \ref{Assumption_C} are satisfied. Moreover, if we let $\mu^{\hat{\boldsymbol{n}}}$ be the corresponding equilibrium measure, then the $\amsmathbb{C}^H$-valued function
 \[
 \hat{\boldsymbol{n}} \mapsto \big(\Gm_{\mu_h^{\hat{\boldsymbol{n}}}}(z)\big)_{h = 1}^{H}
\]
 is smooth, with limits in the definition of derivatives that are uniform over $z$ in compact subsets of
 $\amsmathbb C\setminus \amsmathbb{A}$. We have the formulae,
 \begin{equation}
 \label{formulaehhh}
\forall g,h \in [H]\quad \forall z \in \amsmathbb{C} \setminus [\alpha_h,\beta_h] \qquad \partial_{\hat{n}_{g}} \Gm_{\mu^{\hat{\boldsymbol{n}}}_h}(z) = \Op_h\big[0\,;\,\boldsymbol{e}^{(g)}\big] = \mathfrak{c}^{\textnormal{1st}}_{h;g}(z)
\end{equation}
in terms of the operator $\boldsymbol{\Upsilon}$ constructed in Section~\ref{sec:Masterpb} ($\boldsymbol{e}^{(1)},\ldots,\boldsymbol{e}^{(H)}$ is the canonical basis of $\amsmathbb{C}^H$), or equivalently the functions $\mathfrak{c}^{\textnormal{1st}}$ introduced in Definition~\ref{def:1stkind}. Simultaneously, the endpoints of the bands are also smooth functions of $\hat{\boldsymbol{n}}$. All the partial derivatives with respect to $\hat{\boldsymbol{n}}$ at any order depend on the intensity of interactions $\boldsymbol{\Theta}$ and the $H$-tuple of potentials $\boldsymbol{V}$ in a continuous way --- for the latter the topology induced by the norm appearing in the right-hand side of \eqref{eq_potential_inequality} is used.
\end{proposition}

\begin{proof}[Proof of Proposition~\ref{Proposition_differentiability_filling_fraction}] It is immediate from the definitions that Assumptions~\ref{Assumption_A} and \ref{Assumption_B} are preserved under small perturbations of the segment filling fractions. For Assumption~\ref{Assumption_C} we proved this in
 Theorem~\ref{Theorem_off_critical_neighborhood}. Preservation of the additional three assumptions in the statement of the proposition follows from the continuous dependence of endpoints of the bands and functions $s_h(z)$ on the segment filling fractions, which we have established inside the proof of Theorem~\ref{Theorem_off_critical_neighborhood}. In particular, the one-band condition is preserved, as was argued previously around \eqref{eq_x198}.

We proceed to the computation of derivatives in the filling fractions. Let $\hat{\boldsymbol{n}}$ close enough to $\hat{\boldsymbol{n}}^{(0)}$. Choose $\Delta\hat{\boldsymbol{n}} \in \amsmathbb{R}^H$. For $\epsilon \in \amsmathbb{R}$ with $|\epsilon|$ small enough, we set $\hat{\boldsymbol{n}}^{(\epsilon)} = \hat{\boldsymbol{n}} + \epsilon \Delta\boldsymbol{n}$ and denote $\boldsymbol{\mu}^{(\epsilon)}$ the equilibrium measure corresponding to variational datum obtained from the original one by just changing the segment filling fractions to $\hat{\boldsymbol{n}}^{(\epsilon)}$. The equilibrium measure of the original variational datum and all the corresponding auxiliary functions is obtained for $\epsilon = 0$, and denoted without an exponent $(0)$. For each $h \in [H]$, choose a sequence $(\amsmathbb{M}_h^{(m)})_{m \geq 1}$ of compact subsets of $\amsmathbb{C} \setminus [\hat{a}_h',\hat{b}_h']$ whose union is $\amsmathbb{C} \setminus [\hat{a}'_h,\hat{b}_h']$. Proposition
\ref{Lemma_continuity_filling_fractions} implies that for $|\epsilon| \neq 0$ small enough
\[
\forall m \geq 1 \quad \forall h \in [H]\qquad \sup_{z\in \amsmathbb{M}_h^{(m)}} \frac{1}{\epsilon}\big|\Gm_{\mu^{(\epsilon)}_h}(z)-\Gm_{\mu_h}(z)\big| < + \infty.
\]
Hence, for each $m \geq 1$, Montel theorem tells us that the family of functions $\frac{1}{\epsilon}\big(\Gm_{\mu^{(\epsilon)}_h} -\Gm_{\mu_{h}})$ on $\amsmathbb{M}_h^{(m)}$ indexed by $\epsilon \neq 0$ small enough is precompact in the uniform convergence topology as $\epsilon \rightarrow 0$. Therefore, we can extract a sequence $(\epsilon_l)_l$ tending to $0$ as $l \rightarrow 0$ and functions $\Gm_{\mu_h}'(z)$ which are holomorphic for $z \in \amsmathbb{C} \setminus [\hat{a}'_h,\hat{b}_h']$ for each $h \in [H]$, such that for each $m \geq 1$
\begin{equation}
\label{Gmprimeseq}\forall h \in [H]\qquad \lim_{l\rightarrow\infty} \frac{1}{\epsilon_l}\big(\Gm^{\mu_h^{(\epsilon_l)}}(z) - \Gm_{\mu_h}(z)\big) = \Gm'_{\mu_h}(z),
\end{equation}
uniformly over $z$ in $\amsmathbb{M}_h^{(m)}$, and hence over $z$ in any compact subset of $\amsmathbb C\setminus [\hat a'_{h},\hat b'_{h}]$.

By writing $\Gm_{\mu_h}'(z)$ we did not presume that $\epsilon \mapsto \Gm_{\mu_h^{(\epsilon)}}(z)$ is differentiable at $\epsilon = 0$; this is only a notation for the limit \eqref{Gmprimeseq}, and in particular this limit may depend on the chosen sequence $(\epsilon_l)_l$. The next step is to show that this is not the case, implying differentiability at $\epsilon = 0$. We will do so by identifying the limit point $\Gm'_{\mu_h}(z)$ in \eqref{Gmprimeseq}, and more precisely showing that it is uniquely characterized as the solution of a certain set of equations.

We introduce the functions
$q_{h}^{+,(\epsilon)}(z)$ as in Definition~\ref{GQdef}, so that for any $h \in [H]$
\begin{equation}
\label{eq_x101}
 q_{h}^{+,(\epsilon)}(z)=\phi^+_h(z)\cdot \exp\left( \sum_{g=1}^H \theta_{h,g} \, \mathcal{G}_{\mu^{(\epsilon)}_g}(z)\right) + \phi^-_h(z) \cdot \exp\left(-\sum_{g=1}^H \theta_{h,g}\, \mathcal{G}_{\mu_g^{(\epsilon)}}(z)\right).
\end{equation}
The limit
\begin{equation}
\label{limqhplus}
\lim_{l \rightarrow \infty} \frac{1}{\epsilon_l}\big(q_h^{+,(\epsilon_l)}(z) - q_h^{+}(z)\big)
\end{equation}
exists and is uniform for $z$ in any compact of $\amsmathbb M_h\setminus [\hat a'_h,\hat b'_h]$ by its definition and the corresponding statement for $\Gm_{\mu^{(\epsilon)}_g}$ for any $g \in [H]$ --- \textit{cf.} \eqref{eq_x101}. On the other hand, $z\mapsto q_h^{+,(\epsilon)}(z)$ is meromorphic for $z \in \amsmathbb{M}_h$ by Theorem~\ref{Theorem_regularity_density}, with fixed locations for possible poles, thus its values for $z$ in a neighborhood of $[\hat a'_h,\hat b'_h]$ can be expressed via Cauchy integral formula. Studying the increase rate of this expression between $\epsilon = \epsilon_{l}$ and $\epsilon = 0$, we conclude that the limit \eqref{limqhplus} exists uniformly for $z$ in any compact of $\amsmathbb{M}_h$, and we denote it $(q_h^{+})'(z)$. Again, this is only a notation for this sequential limit and does not presume of the differentiability of $\epsilon \mapsto q_h^{+,(\epsilon)}(z)$ at $\epsilon = 0$. The Leibniz rule to compute limits of sequential increase rates is nevertheless valid and we get from \eqref{eq_x101} the expression
\begin{equation}
\label{eq_x102} (q_h^{+})'(z) = \left(\sum_{g=1}^H \theta_{h,g}\,\Gm'_{\mu_{g}}(z)\right)\,q_h^{-}(z).
\end{equation}
Since $(\alpha_h,\beta_h)$ is the unique band of $\mu$ in $[\hat{a}_h,\hat{b}_h]$, the auxiliary tuple of functions $(s_h)_{h = 1}^{H}$ of Definition~\ref{GQdef2} associated to the equilibrium measure $\boldsymbol{\mu}$ is
\[
\forall h \in [H]\qquad s_h(z) = \frac{q_h^{-}(z)}{\sqrt{(z - \alpha_h)(z - \beta_h)} \cdot \psi_h(z)}.
\]
Dividing \eqref{eq_x102} by $s_h(z) \cdot \psi_h(z)$, we get
\begin{equation}
\label{eq_x177} \forall h \in [H]\qquad \sqrt{(z-\alpha_h)(z-\beta_h)} \left( \sum_{g=1}^H \theta_{h,g}\,\Gm'_{\mu_g}(z)\right) = \frac{(q_{h}^{+})'(z)}{s_h(z) \cdot \psi_h(z)}.
\end{equation}
We claim that the right-hand side of \eqref{eq_x177} is holomorphic inside $\amsmathbb{M}_h$. Indeed, we already know that $q_{h}^{+,(\epsilon)}(z)$ and $s_h(z)$ are holomorphic and also $s_h(z)$ is bounded away from $0$ by our additional assumption. We may have a pole either at $z = \hat a'_h$ or $z = \hat b'_h$ because of $\psi_h(z)$ in the denominator. However, if $\hat a'_h$ is an endpoint of a saturation, then at this point $\phi^-_h(z)$ has a double zero by our additional assumption. So, the functions
\[
\exp\left(-\sum_{g=1}^H \theta_{h,g}\, \Gm_{\mu^{(\epsilon)}_{h}}(z)\right)\qquad \textnormal{respectively} \qquad \exp\left(\sum_{g=1}^H \theta_{h,g}\, \Gm_{\mu^{(\epsilon)}_{g}}(z)\right)
\]
get a simple zero, respectively a simple pole, at $z = \hat{a}_h'$. We conclude that $q_{h}^{+,(\epsilon)}(z)$ has a zero at $z = \hat a'_h$. The location of this zero is independent of $\epsilon$, and, hence, by Rouch\'e theorem, the sequential derivative $(q_h^{+})'(z)$ also admits a zero at $z = \hat{a}'_h$. Therefore, the ratio in the right-hand side of \eqref{eq_x177} has no singularity at $z = \hat a'_h$. The argument for $\hat b'_h$ is the same. In addition, studying the sequential increase rate of the condition that the segment filling fractions of $\boldsymbol{\mu}^{(\epsilon)}$ are $\hat{\boldsymbol{n}} + \epsilon \Delta\hat{\boldsymbol{n}}$, we get
\begin{equation}
\label{eq_x103}
\forall h \in [H]\qquad \oint_{\gamma_h} \frac{\dd z}{2\ii\pi}\,(\Gm_{\mu_h})'(z) = \Delta\hat{n}_h.
\end{equation}

We treat \eqref{eq_x177}, \eqref{eq_x103} as a system of $2H$ equations on the unknown $H$-tuple of functions $(\Gm'_{\mu_h})_{h = 1}^{H}$. In these equations the exact form of the right-hand side of \eqref{eq_x177} is irrelevant; only its holomorphicity inside $\amsmathbb M_h$ is important. Systems of equations of this form are studied independently of the rest of this book in Part~\ref{Part_Master_equation}. In particular, it will be established in Chapter~\ref{Chapter_SolvingN} (\textit{cf.} Theorem~\ref{Theorem_Master_equation_12}) that they admit a unique solution. For the present use we only need the specialization of Theorem~\ref{Theorem_Master_equation_12} to $\boldsymbol{E}(z) = 0$, which provides the unique solution of \eqref{eq_x177},\eqref{eq_x103} in the form
\begin{equation}
\label{eq_x104}
\forall h \in [H] \qquad \Gm'_{\mu_h}(z)=\Op_h[\boldsymbol{0}\,;\,\Delta\hat{\boldsymbol{n}}],
\end{equation}
where $\boldsymbol{\Op} = (\Op_h)_{h = 1}^H$ is a certain $H$-tuple of continuous linear operators depending smoothly on the parameters involved. The formulae \eqref{eq_x104} imply that the sequential limit in \eqref{Gmprimeseq} is unique, and hence the function $\epsilon \mapsto \Gm_{\mu_h^{(\epsilon)}}(z)$ is differentiable at $\epsilon = 0$ and its derivative is \eqref{eq_x104}.

Although it may seem that the right-hand side of \eqref{eq_x104} does not depend
on the segment filling fractions, it is not quite true. Indeed, the equations \eqref{eq_x177} and, hence, the operator $\boldsymbol{\Op}$
depend on the $H$-tuple of band endpoints $\boldsymbol{\alpha},\boldsymbol{\beta}$ which, in turn,
depend on the segment filling fractions. Let us show that the endpoints of the bands depend on the segment filling fractions in a
differentiable way --- it will imply that $\boldsymbol{\Op}$ also depends on the segment filling
fractions in a differentiable way. Let $\mathfrak{a}^{(\epsilon)}$ be an endpoint
of the $h$-th band of the equilibrium measure $\mu^{(\epsilon)}$. In
Theorem~\ref{Theorem_regularity_density} we have
shown that $\mathfrak{a}^{(\epsilon)}$ is a zero of $(q_{h}^{-,(\epsilon)}(z))^2$; this zero is necessary simple, since $s_h^{(\epsilon)}(z)$ is bounded away from $0$. The function
$(q_{h}^{-,(\epsilon)}(z))^2$ is holomorphic for $z$ in a neighborhood of $[\hat a'_h,\hat b'_h]$, and it
depends on $\epsilon$ in a differentiable way (in the uniform topology), as follows from the
differentiability of $\epsilon \mapsto \Gm_{\mu^{(\epsilon)}_h}(z)$. As the zero is simple, the
derivative in $z$ of $(q_{h}^{-,(\epsilon)}(z))^2$ is non-zero at
$\mathfrak{a}^{(\epsilon)}$. Therefore, the implicit function theorem guarantees that
$\epsilon \mapsto\mathfrak{a}^{(\epsilon)}$ is differentiable.

At this point, \eqref{eq_x104} readily gives the smoothness of
$\epsilon \mapsto \Gm_{\mu_h^{(\epsilon)}}(z)$. Indeed, making the dependence on partial fractions in \eqref{eq_x104}
explicit, we can rewrite it as
\begin{equation}
\label{eq_x105}
\forall g,h \in [H]\qquad \partial_{\hat{n}_{g}} \Gm_{\mu_h^{\hat{\boldsymbol{n}}}}(z)=\Op_{h}^{\hat{\boldsymbol{n}}}
 [\boldsymbol{0}\,;\,\boldsymbol{e}^{(g)}](z).
 \end{equation}
 where the notation $\Op_h^{\hat{\boldsymbol{n}}}$ stresses that we use the operator $\Upsilon_h$ for the endpoints associated to $\mu^{\hat{\boldsymbol{n}}}$. The alternative expression via first-kind functions $\mathfrak{c}^{\textnormal{1st}}_{h;g}(z)$ is proved in Corollary~\ref{cor:1stkindapp}. Theorem~\ref{Theorem_Master_equation_12} explains that the operator $\Op_h$ depends smoothly on the $H$-tuples $\boldsymbol{\alpha},\boldsymbol{\beta}$ of endpoints of the bands. Therefore, the right-hand side of \eqref{eq_x105} is differentiable
 with respect to $\hat{\boldsymbol{n}}$, hence, so is the left-hand side, and we get a similar
expression for the second-order partial derivatives of $\Gm_{\mu^{\hat{\boldsymbol{n}}}_h}(z)$ through first-order partial derivatives of $\Op^{\hat{\boldsymbol{n}}}$. Continuing, we get the expressions for the third-order partial derivatives of $\Gm_{\mu^{\hat{\boldsymbol{n}}}_h}(z)$, \textit{etc}. This eventually shows smoothness of $\Gm_{\mu^{\hat{\boldsymbol{n}}}_h}(z)$ with respect to $\hat{\boldsymbol{n}}$ in a small enough neighborhood of $\boldsymbol{\hat{n}}^{(0)}$.
\end{proof}

We proceed to the most general differentiability statement. In order to allow variation of the potential and intensity of interaction, we consider families of variational data depending on an auxiliary real parameter denoted $u$.

\begin{theorem}
\label{Theorem_differentiability_full}
Let $p \in \amsmathbb{Z}_{\geq 0}$. Consider a variational datum satisfying Assumptions~\ref{Assumption_A}, \ref{Assumption_B} and \ref{Assumption_C}, with segment filling fractions $\hat{\boldsymbol{n}}^0$, and such that for each $h \in [H]$ the segment $[\hat{a}^{\prime\, 0}_h,\hat{b}^{\prime\,0}_h]$ contains a single band of the equilibrium measure. There exists a constant $\delta > 0$ depending only on the constants in the assumptions, with the following properties.

Consider any family of variational data parameterized by $\boldsymbol{t} = (\hat{\boldsymbol{a}}',\hat{\boldsymbol{b}}',\hat{\boldsymbol{n}},u)$ varying in the $\delta$-neighborhood of $\boldsymbol{t}^0 = (\hat{\boldsymbol{a}}^{\prime\,0},\hat{\boldsymbol{b}}^{\prime\,0},\hat{\boldsymbol{n}}^{0},0)$ in $\amsmathbb{R}^{3H + 1}$, such that the first three entries of $\boldsymbol{t}$ represent the segment endpoints and filling fractions and the variational datum for $\boldsymbol{t} = \boldsymbol{t}^0$ coincides with the original one. Suppose that these variational data satisfy Assumption~\ref{Assumption_B} for a $\boldsymbol{t}$-independent tuple of integers $\boldsymbol{\iota}^{\pm}$ and $\boldsymbol{t}$-independent complex domains $(\amsmathbb{M}_h)_{h = 1}^{H}$, and that for any $h \in [H]$ the regular part of the potential $U_h^{\boldsymbol{t}}(z)$ is a $p$-times continuously differentiable function of $\boldsymbol{t}$ --- with any partial derivative up to order $p$ being a holomorphic function of $z \in \amsmathbb{M}_h$.

Then, denoting $\boldsymbol{\mu}^{\boldsymbol{t}}$ the equilibrium measure corresponding to the $\boldsymbol{t}$-dependent variational datum, for any $h \in [H]$ the function
\[
\boldsymbol{t} \longmapsto \Gm_{\mu^{\boldsymbol{t}}_h}(z)
\]
is $p$ times continuously differentiable, with any partial derivative up to order $p$ being a holomorphic function of $z \in \amsmathbb{C} \setminus [\hat{a}_h^{\prime\,\boldsymbol{t}},\hat{b}_h^{\prime\,\boldsymbol{t}}]$. Simultaneously, for any $\boldsymbol{t}$ as above there remains exactly one band per segment and the endpoints of the bands are $p$ times continuously differentiable functions of $\boldsymbol{t}$. All the derivatives depend on the matrix of interactions $\boldsymbol{\Theta}$ in a continuous way.
\end{theorem}
\begin{proof} The proof splits into two parts: first, we reduce the most general statement of Theorem~\ref{Theorem_differentiability_full} to a minor extension of Proposition~\ref{Proposition_differentiability_filling_fraction} by restricting the equilibrium measure to smaller segments; the endpoints of these segments remain fixed and do not vary with $\boldsymbol{t}$. Second, we explain how this extension is proved by the same method as the one we used in Proposition~\ref{Proposition_differentiability_filling_fraction}.

We write $\boldsymbol{\mu}^{0} : = \boldsymbol{\mu}^{(\boldsymbol{t}^0)}$ and more generally indicate with an exponent ${}^0$ on all the quantities associated with the original variational data. Due to the results of Sections~\ref{Section-ParReg} and \ref{Section_off_criticality}, the $\boldsymbol{t}$-dependent variational datum in the statement of Theorem~\ref{Theorem_differentiability_full} satisfies Assumptions~\ref{Assumption_A}, \ref{Assumption_B} and \ref{Assumption_C} for $\boldsymbol{t}$ close enough to $\boldsymbol{t}^0$ and the nature --- void or saturated --- of the endpoints of the segments do not change with $\boldsymbol{t}$. The latter implies $\mathbbm{1}_{\amsmathbb{S}_h}(\hat{a}_h') = \mathbbm{1}_{\amsmathbb{S}_h^0}(\hat{a}_h^{\prime\,0})$, \textit{etc}. which we will be used at several places in the proof to avoid writing an exponent ${}^0$.

We start by choosing two $H$-tuples $\Delta\hat{\boldsymbol{a}},\Delta\hat{\boldsymbol{b}} \in \amsmathbb{R}_{\geq 0}^{H}$ for which there exists $\eta > 0$ satisfying the following conditions for any $h\in[H]$.
\begin{itemize}
 \item $\hat a_h^{\prime\,0} +\Delta \hat a_h < \hat b_h^{\prime\,0} - \Delta \hat b_h$,
 \item The entire segment $[\hat a_h^{\prime\,0},\hat a_h^{\prime\,0}+\Delta \hat{a}_h+\eta]$ is void, or is saturated. The entire segment $[\hat{b}_h^{\prime\,0} - \Delta \hat b_h - \eta,\hat{b}_h^{\prime\,0}]$ is void, or is saturated.
 \item $\forall x \in [\hat{a}_h^{\prime\,0} + \Delta \hat{a}_h,\hat{b}_h^{\prime\,0} - \Delta \hat{b}_h]\quad s_h(x) \neq 0$.
 \end{itemize}
It is possible to make such a choice owing to the one-band-per-segment condition and the lower bound on $|s_h|$ from Proposition~\ref{Proposition_q_bounds}.

As mentioned before, take $\delta > 0$ small enough. For each $\boldsymbol{t}$ in a $\delta$-neighborhood of $\boldsymbol{t}^0$, we are going to define a modified variational datum, which has $H$ segments $[\hat{a}_h^{\prime\,0} + \Delta \hat{a}_h,\hat{b}_h^{\prime\,0} - \Delta \hat{b}_h]$ and unchanged intensities of interactions $\boldsymbol{\Theta}^{\boldsymbol{t}}$. The other parameters are chosen so that the modified equilibrium measure on the $h$-th segment is the restriction of $\mu_h^{\boldsymbol{t}}$ to the modified segments. This can be achieved by compensating the filling fractions and the potential for the saturated parts that have been removed. In greater details, for $h \in [H]$ the $h$-th modified filling fractions is
\begin{equation}
\label{nhtmodd}
\hat{n}_{h}^{\boldsymbol{t},\textnormal{mod}} = \hat{n}^h - \frac{\hat{a}_h' - \hat{a}_h^{\prime\,0} - \Delta\hat{a}_h}{\theta_{h,h}^{\boldsymbol{t}}}\mathbbm{1}_{\amsmathbb{S}_h}(\hat{a}_h') - \frac{\hat{b}_h^{\prime\,0} - \Delta\hat{b}_h - \hat{b}_h'}{\theta_{h,h}^{\boldsymbol{t}}} \mathbbm{1}_{\amsmathbb{S}_h}(\hat{b}_h'),
\end{equation}
and the modified potential is
\[
V_h^{\boldsymbol{t},\textnormal{mod}}(x) = V_h^{\boldsymbol{t}}(x) - \sum_{g = 1}^{H}\ \frac{2\,\theta_{h,g}^{\boldsymbol{t}}}{\theta_{g,g}^{\boldsymbol{t}}}\left( \mathbbm{1}_{\amsmathbb{S}_g}(\hat{a}_g') \int_{\hat{a}'_g}^{\hat{a}_g^0 + \Delta \hat{a}_g} \log|x - y| + \mathbbm{1}_{\amsmathbb{S}_g}(\hat{b}_g') \int_{\hat{b}_g^{\prime\,0} - \Delta \hat{b}_g}^{\hat{b}'_{g}} \log|x - y|\dd y\right).
\]
This also leads to modifications in the function appearing in Assumption~\ref{Assumption_C}. For any $h \in [H]$, we have
\begin{equation}
\label{phihmodmd}
\phi_h^{-,\boldsymbol{t},\textnormal{mod}}(z) = \phi_h^{-,\boldsymbol{t}}(z) \cdot \left\{ \begin{array}{lll} \dfrac{1}{(z - \hat{a}_h)^{\iota_h^-}} & & \textnormal{if} \quad \hat{a}_h' \in \amsmathbb{V}_h, \\[15pt] \dfrac{(z - \hat{a}_h^{\prime\,0} - \Delta \hat{a}_h)^2}{(z - \hat{a}_{h}')^{\iota_h^-}} && \textnormal{if}\quad \hat{a}_h' \in \amsmathbb{S}_h, \end{array} \right.
\end{equation}
while $\phi_h^{+,\boldsymbol{t},\textnormal{mod}}(z)$ is equal to the product of $\phi_h^{+,\boldsymbol{t}}(z)$ with all the following factors depending on the nature of the $h$-th segment endpoints
\begin{equation}
\label{phihmodmd2}
\begin{split}
 \quad (z - \hat{a}_h')^{-\iota_h^-} & \qquad \textnormal{if}\,\,\hat{a}_h' \in \amsmathbb{V}_h, \\
 (z - \hat{a}_h')^{2 - \iota_h^-} & \qquad \textnormal{if}\,\, \hat{a}_h' \in \amsmathbb{S}_h, \\
 \left(\frac{z - \hat{b}_h^{\prime 0} + \Delta \hat{b}_h}{z - \hat{b}_h'}\right)^2 & \qquad \textnormal{if}\,\, \hat{b}_h' \in \amsmathbb{S}_h,
\end{split}
\end{equation}
and with the product over $g \neq h$ of all the following factors appearing if the $g$-th segment endpoints are saturated
\begin{equation}
\label{phihmodmd3}
\begin{split}
\exp\left(\frac{2\,\theta_{h,g}^{\boldsymbol{t}}}{\theta_{g,g}^{\boldsymbol{t}}} \int_{\hat{a}'_g}^{\hat{a}_g^{\prime\,0} + \Delta \hat{a}_g} \frac{\dd y}{z - y}\right) & \qquad \textnormal{if}\,\, \hat{a}_g' \in \amsmathbb{S}_g, \\
\exp\left(\frac{2\,\theta_{h,g}^{\boldsymbol{t}}}{\theta_{g,g}^{\boldsymbol{t}}} \int_{\hat{b}_g^{\prime\,0} - \Delta \hat{b}_g}^{\hat{b}_g'} \frac{\dd y}{z - y}\right) & \qquad \textnormal{if}\,\, \hat{b}_g' \in \amsmathbb{S}_g'.
\end{split}
\end{equation}

We observe that the modifications in the functions $\phi^{\pm}_h(z)$ are being offset with similar modifications in $\sum_{g = 1}^H \theta_{h,g}\,\Gm_{\mu_{h}}(z)$, so that eventually both terms in the definition of $q^-_h(z)$ are multiplied by the same factor. As a result, for any $h \in [H]$ the function $s_h^{\boldsymbol{t},\textnormal{mod}}(z)$ is obtained by multiplying $s_h^{\boldsymbol{t}}(z)$ with is an explicit holomorphic function of $z$ in $\amsmathbb{M}_h$ bounded away from $0$ on the interval $[\hat a_h^{\prime\,0}+\Delta \hat a_h, \hat b_h^{\prime 0} - \Delta \hat b_h]$. Indeed, the only potentially vanishing factors are $ (z-\hat a^{\prime\,0}_h- \Delta \hat a_h)^2$ in $\phi^{-,\boldsymbol{t},\textnormal{mod}}_h(z)$ and $(x-\hat b^{\prime\,0}_h+\Delta \hat b_h)^2$ in $\phi^{+,\boldsymbol{t},\textnormal{mod}}_h(z)$ if the corresponding endpoints belong to a saturation. The former is offset by the pole in $\exp\big(-\theta_{h,h}\,\Gm_{\mu_{h}^{\boldsymbol{t},\textnormal{mod}}}(z)\big)$ and by the division with $\psi_h(z)$ in the definition of $s_h(z)$; the latter is offset by the pole in $\exp\big(\theta_{h,h}\,\Gm_{\mu_{h}^{\boldsymbol{t},\textnormal{mod}}}(z)\big)$ and by the division with $\psi_h(z)$ in the definition of $s_h(z)$. The conclusion from this discussion is that for any $h \in [H]$, the auxiliary functions $s_h^{\boldsymbol{t},\textnormal{mod}}(z)$ for the modified variational are bounded away from $0$, like in the second additional assumption of Proposition~\ref{Proposition_differentiability_filling_fraction}.

\medskip

Let us summarize the properties of the modified variational data, for $\delta$ chosen small enough. First, they satisfy by design the Assumptions~\ref{Assumption_A}, \ref{Assumption_B} and \ref{Assumption_C}. Second, they satisfy the three additional assumptions of Proposition~\ref{Proposition_differentiability_filling_fraction} --- something which may have failed before the modification. Third, the endpoints of the segments are independent of $\boldsymbol{t}$: the variation of $\hat a'_h$ and $\hat b'_h$ got re-expressed into changes of the filling fractions and potentials $V_h$. Fourth, we retrieve $\boldsymbol{\mu}^{\boldsymbol{t}}$ from its modification by adding back saturated segments following a procedure independent of the segment filling fractions $\hat{\boldsymbol{n}}$: the existing saturations in the modification are simply enlarged as the $h$-th modified segment is enlarged to $[\hat{a}_h,\hat{b}_h]$ for each $h \in [H]$. Hence, $\Gm_{\mu_h^{\boldsymbol{t}}}(z)$ is a differentiable function of $\boldsymbol{t}$ if and only if $\Gm_{\mu_h^{\boldsymbol{t},\textnormal{mod}}}(z)$ is a differentiable function of $\boldsymbol{t}$. In the rest of the proof we establish the latter. This in fact closely follow the proof of Proposition~\ref{Proposition_differentiability_filling_fraction}, so we omit some details.

\medskip

 Choose $\boldsymbol{t}$ close enough to $\boldsymbol{t}^0$, a real vector $\Delta\boldsymbol{t} = (\Delta\hat{\boldsymbol{a}},\Delta\hat{\boldsymbol{b}},\Delta\boldsymbol{n},\Delta u) \in \amsmathbb{R}^{3H + 1}$ and look at the subfamily of \emph{modified} --- as explained above --- variational data corresponding to $\boldsymbol{t}^{(\epsilon)} = \boldsymbol{t} + \epsilon \Delta\boldsymbol{t}$ parameterized by a real $\epsilon$ in a small enough neighborhood of $0$. To work with a shorter notation we denote with an exponent ${}^{(\epsilon)}$ all the associated quantities, and without any exponent all the quantities associated with $\epsilon = 0$. The Lemmata~\ref{Lemma_continuity_potential}, \ref{Lemma_continuity_filling_fractions},
 \ref{Lemma_continuity_end_points_void},
 \ref{Lemma_continuity_end_points_saturated} guarantee that $\Gm_{\mu_h^{(\epsilon)}}(z)$
 depends on $\epsilon$ in a Lipschitz way, and therefore, there exist a sequence $(\epsilon_l)_l$ tending to $0$ such that the following limit exist
 \begin{equation}
\label{eq_x157}
\forall h \in [H]\qquad \lim_{l\rightarrow\infty} \frac{1}{\epsilon_l} \big(\Gm_{\mu_h^{(\epsilon_l)}}(z) - \Gm_{\mu_h}(z)\big) = \Gm'_{\mu_h}(z).
\end{equation}
Again, the $'$ is a notation for the sequential limit and does not presume of differentiability. In order to identify the in limit in \eqref{eq_x157}, we use the auxiliary functions of the modified and $\epsilon$-dependent variational datum
\begin{equation}
\label{eq_x158}
 q_{h}^{+,(\epsilon)}(z)=\phi^{+,(\epsilon)}_h(z) \cdot \exp\left(\sum_{g=1}^H \theta_{h,g}\,\Gm_{\mu_g^{(\epsilon)}}(z)\right) + \phi^{-,(\epsilon)}_h(z) \cdot \exp\left(-\sum_{g=1}^H \theta_{h,g}\,\Gm_{\mu_g^{(\epsilon)}}(z)\right).
\end{equation}
The difference with \eqref{eq_x101} in the proof of Proposition~\ref{Proposition_differentiability_filling_fraction} is that $\phi^{\pm,(\epsilon)}_h$ now depends on $\epsilon$ through the $\epsilon$-dependent position of the non-modified endpoints that may occur in numerators and denominators --- \textit{cf.} \eqref{phihmodmd}, \eqref{phihmodmd2}, \eqref{phihmodmd3} with $\boldsymbol{t}$ replaced with $\boldsymbol{t} + \epsilon \Delta \boldsymbol{t}$. This dependence is clearly differentiable in $\epsilon$. Taking the $l \rightarrow \infty$ limit in the increase rate of \eqref{eq_x158} between $\epsilon = \epsilon_l$ and $\epsilon=0$ we get by Leibniz rules
\begin{equation}
\label{eq_x158bis}
\sqrt{(z-\alpha_h)(z-\beta_h)} \left(\sum_{g=1}^{H}
\theta_{h,g}\,\Gm'_{\mu_g}(z)\right) = \frac{\partial_{\epsilon}q_{h}^{+,(\epsilon)}(z)\big|_{\epsilon=0}}{s_h(z) \cdot \psi_h(z)} + E_h(z),
\end{equation}
where
\[
E_h(z) = -\frac{\partial_{\epsilon} \phi^{+,(\epsilon)}_h(z)\big|_{\epsilon=0}\exp\big(\sum_{g=1}^H \theta_{h,g}\,\Gm_{\mu_g}(z)\big) +
 \partial_{\epsilon}\phi^{-,(\epsilon)}_h(z)\big|_{\epsilon=0} \cdot \exp\big(-\sum_{g=1}^H \theta_{h,g}\,\Gm_{\mu_g}(z)\big)}{s_h(z) \cdot \psi_h(z)},
\]
where $(\alpha_h,\beta_h)$ is the band of $\boldsymbol{\mu}$ in the $h$-th segment. In addition, since the segment filling fractions have an affine dependence in $\epsilon$ as per \eqref{nhtmodd}, we have
\begin{equation}
\label{eq_x159}
\forall h \in [H]\qquad \oint_{\gamma_h} \frac{\dd z}{2\ii\pi}\,\Gm'_h(z) = \Delta \hat{n}_h - \frac{\Delta \hat{a}_h}{\theta_{h,h}} \cdot \mathbbm{1}_{\amsmathbb{S}_h}(\hat{a}_h') - \frac{\Delta \hat{b}_h}{\theta_{h,h}} \cdot \mathbbm{1}_{\amsmathbb{S}_h}(\hat{b}_h') := \kappa_h.
\end{equation}
The last two terms come from the change in the filling fraction caused by restricting from $[\hat{a}_h',\hat{b}_h']$ to $[\hat a_h^{\prime\,0}+\Delta \hat a_h, \hat b_h^{\prime\,0} +\Delta \hat b_h]$. Observe that the second line in \eqref{eq_x158bis} is holomorphic inside $\amsmathbb{M}_h$ --- compared to the proof of Proposition~\ref{Proposition_differentiability_filling_fraction}, it is important here that the endpoints of the segments $[\hat a_h^0+\Delta \hat a_h, \hat b_h^0 +\Delta \hat b_h]$ do not change with $\epsilon$, as otherwise a pole might have appeared because of the division by $\psi_h(z)$. The tuples $\boldsymbol{\kappa} = (\kappa_h)_{h = 1}^{H}$ and $\boldsymbol{E}(z) = (E_h(z))_{h =1}^{H}$ are explicit and do not depend on the choice of the sequence
$(\epsilon_l)_l$. Hence, the equations \eqref{eq_x158bis}, \eqref{eq_x159} fit into the
framework of Chapter~\ref{Chapter_SolvingN} --- up to taking for each $h \in [H]$ smaller complex domains $\amsmathbb{M}_{h}$ so that the possible poles of the third line of \eqref{eq_x158bis} at $\hat a'_h$ or $\hat b'_h$ are kept outside this set. Therefore, these equations admit a
unique solution given in the form
\begin{equation}
\label{fomur}
\forall h \in [H] \qquad \Gm'_{\mu_h}(z) =\Op_h\big[\boldsymbol{E}\,;\,\boldsymbol{\kappa}\big](z).
\end{equation}
This formula shows that all the sequential limits are the same, and therefore, \eqref{eq_x157} computes the $\epsilon \rightarrow 0$ limit, \textit{i.e.} $\Gm^{\mu_h^{(\epsilon)}}(z)$ is differentiable at $\epsilon = 0$. Since the operators $\Upsilon_1,\ldots,\Upsilon_H$ depend smoothly on all the parameters, we also have $p$-times continuous differentiability with respect to $\boldsymbol{t}$ by differentiating further \eqref{fomur} in the same way as we did with \eqref{eq_x104}.
\end{proof}

\section{Derivatives of the free energy}
\label{Section_derivatives_free_energy}

We now explain two consequences of Theorem~\ref{Theorem_differentiability_full} pertaining to the derivatives of the minimum of the energy functional $-\I$ from Definition~\ref{Definition_functional}. The first one, Proposition~\ref{Proposition_Hessian_free_energy} below, will be used in Chapter~\ref{Chapter_filling_fractions} when we deal with ensembles having fluctuating filling fractions and study the asymptotic behavior of filling fractions as random variables. The second one, first stated in a simpler form in Proposition~\ref{Proposition_quadratic_potential_differentiability} and in greater generality and abstraction in Proposition~\ref{proposition_Energy_series} below, will be used in Chapter~\ref{Chapter_partition_functions} when we interpolate between discrete ensembles with fixed filling fractions and a reference discrete ensemble.

\begin{proposition}
\label{Proposition_Hessian_free_energy}
Consider a variational datum satisfying Assumptions~\ref{Assumption_A}, \ref{Assumption_B} and \ref{Assumption_C} with segment filling fractions $\boldsymbol{n}^{0}$, such that each segment contains a single band of the equilibrium measure. Consider also affine equations \eqref{eq_equations_eqs} as in Section~\ref{DataS} --- with variable $\frac{N_h}{\N}$ replaced with $\hat{\boldsymbol{n}}$ --- which we assume to be satisfied by $\hat{\boldsymbol{n}}^0$, to fulfill the conditions 4.5.6. of Assumption~\ref{Assumptions_Theta}, and to admit other solutions than $\hat{\boldsymbol{n}}^0$.

Then, there exist constants $\delta > 0$ and $C > 0$ depending only on the constants in the assumptions with the following property. For $\hat{\boldsymbol{n}}$ varying in a $\delta$-neighborhood of $\hat{\boldsymbol{n}}^0$ in the space of solutions of ($\star$), denoting $\boldsymbol{\mu}^{\hat{\boldsymbol{n}}}$ the equilibrium measure of the variational datum differing from the previous one only by the segment filling fractions which are equal to $\hat{\boldsymbol{n}}$, the function
\[
\hat{\boldsymbol{n}} \longmapsto -\I[\hat{\boldsymbol{n}}]
\]
is twice-continuously differentiable, its Hessian is positive definite, and we have
\begin{equation}
\label{soungnzn} \frac{1}{C} \cdot |\!|\hat{\boldsymbol{n}} - \hat{\boldsymbol{n}}^0 |\!|_{\infty}^2 \leq -\I[\boldsymbol{\mu}^{\hat{\boldsymbol{n}}}] \leq C \cdot |\!|\hat{\boldsymbol{n}} - \hat{\boldsymbol{n}}^0 |\!|_{\infty}^2.
\end{equation}
\end{proposition}
\begin{proof}
 We first claim that the density of the equilibrium measure as a function of the segment filling fractions
 \begin{equation}
\label{eq_x118}
\hat{\boldsymbol{n}} \longmapsto \mu^{\hat{\boldsymbol{n}}}(x)
\end{equation}
 is twice-continuously differentiable for all $x \in \amsmathbb{A}$ except at the endpoints of the
 bands --- where the derivatives explode. We already proved in Proposition~\ref{Proposition_differentiability_filling_fraction} and Theorem~\ref{Theorem_differentiability_full} that the endpoints of the bands depend on $\hat{\boldsymbol{n}}$ close enough to $\hat{\boldsymbol{n}}^0$. For $x$ outside the bands, the density is constant so it is clearly differentiable with respect to $\hat{\boldsymbol{n}}$. For $x$ in a band, for any $h \in [H]$ we expressed in Proposition~\ref{Proposition_density} the density $\mu_h(x)$ through the auxiliary functions $q^{-,\hat{\boldsymbol{n}}}_h(x^+)$ and $q^{+,\hat{\boldsymbol{n}}}_h(x)$ of Definition~\ref{GQdef}. We now use Theorem~\ref{Theorem_differentiability_full}. Because both $(q^{-,\hat{\boldsymbol{n}}}_h(z))^2$ and $q^{+,\hat{\boldsymbol{n}}}_h(z)$ are holomorphic functions of $z$ in a complex neighborhood of the $h$-th segment, the twice-continuous differentiability of $\hat{\boldsymbol{n}} \mapsto (\Gm_{\mu^{\hat{\boldsymbol{n}}}_h}(z))_{h = 1}^{H}$ implies the same regularity for $\boldsymbol{n} \mapsto (q^{-,\hat{\boldsymbol{n}}}_h(z))^2$ and $\hat{\boldsymbol{n}} \mapsto q^{+,\hat{\boldsymbol{n}}}_h(z)$ for $z \in \amsmathbb{C}$ close to but outside the $h$-th band, hence by Cauchy integral formula for all $z$ in a complex neighborhood of the $h$-th band. Since the only zeros of $q^{-,\hat{\boldsymbol{n}}}_h(z)$ near the $h$-th band are endpoints of this band due to the lower bound for $s_h(z)$ in Proposition~\ref{Proposition_q_bounds}, we conclude that both $q^{-,\hat{\boldsymbol{n}}}_h(x^+)$ and $q^{+,\hat{\boldsymbol{n}}}_h(x)$ are twice-continuously differentiable with respect to $\hat{\boldsymbol{n}}$ for any $x$ in the $h$-th band. By \eqref{eq_density_exponent}, we deduce that $\hat{\boldsymbol{n}} \mapsto \mu^{\hat{\boldsymbol{n}}}_h(x)$ is also twice continuously differentiable.

 The same argument implies that the first-order partial derivatives of the densities in \eqref{eq_x118} have
inverse square-root singularities at the endpoints of the bands. In particular, they are still integrable in $x$ near these endpoints. We conclude that the definition \eqref{eq_functional_general} of $\I$ can be differentiated once and for each $i \in[H]$ we have
 \begin{eqnarray*}
\partial_{\hat{n}_i} \I[\boldsymbol{\mu}^{\hat{\boldsymbol{n}}}] & = & \partial_{\hat{n}_i}\left(
 \sum_{g,h = 1}^{H} \int_{\hat a'_{g}}^{\hat b'_{g}} \int_{\hat a'_{h}}^{\hat
 b'_{h}} \theta_{g,h}\,\log|x-y| \mu^{\hat{\boldsymbol{n}}}_{g}(x)\mu^{\hat{\boldsymbol{n}}}_{h}(y)\dd x\dd y
-\sum_{h=1}^H \int_{\hat a'_h}^{\hat b'_h} V_h(x)\mu^{\hat{\boldsymbol{n}}}_h(x)\dd x \right) \\
& = & -\sum_{h=1}^H \int_{\hat a'_h}^{\hat b'_h} V^{\textnormal{eff},\hat{\boldsymbol{n}}}_h(x) \big(\partial_{\hat n_i}\mu^{\hat{\boldsymbol{n}}}_{h}(x) \big) \dd x = -\sum_{h=1}^H \int_{\hat a'_h}^{\hat b'_h} v^{\hat{\boldsymbol{n}}}_h\,\big(\partial_{\hat{n}_i}\mu^{\hat{\boldsymbol{n}}}_{h}(x)\big)\dd x \\
& = & -\sum_{h=1}^H v^{\hat{\boldsymbol{n}}}_h\,\partial_{\hat n_i} \left(\int_{\hat a'_h}^{\hat b'_h} \mu^{\hat{\boldsymbol{n}}}_{h}(x) \dd x\right)= - v_i^{\hat{\boldsymbol{n}}},
 \end{eqnarray*}
 where we used the effective potential and the constants from the characterization Theorem~\ref{Theorem_equi_charact_repeat_2} for the variational datum with segment filling fractions $\hat{\boldsymbol{n}}$ and the fact that $\partial_{\hat{n}_i}\mu^{\hat{\boldsymbol{n}}}_{h}(x)$ is zero outside of the $h$-th band for any $h \in [H]$.

 Taking the second derivative, we get for any $i,j \in [H]$
\[
 \partial_{\hat n_i}\partial_{\hat n_j} \I[\boldsymbol{\mu}^{\hat{\boldsymbol{n}}}]=-\partial_{\hat n_j} v_i^{\hat{\boldsymbol{n}}}= -\partial_{\hat n_j} V^{\textnormal{eff},\hat{\boldsymbol{n}}}_i(x),
\]
where $x$ is an arbitrary point in the $i$-th band of $\mu^{\hat{\boldsymbol{n}}}$, that can be taken independent of $\hat{\boldsymbol{n}}$. Using the definition, we can further write for any $i,j \in [H]$
\begin{equation}
\label{eq_x119}
\begin{split}
\partial_{\hat{n}_{j}} V^{\textnormal{eff},\hat{\boldsymbol{n}}}_i(x) & = -\partial_{\hat{n}_{j}}\bigg( \sum_{h=1}^H 2\theta_{h,i}\int_{\hat a'_h}^{\hat
 b'_h} \log|x-y|\mu^{\hat{\boldsymbol{n}}}_h(y)\dd y\bigg) \\
 & = -\sum_{h=1}^H 2\theta_{h,i}\int_{\hat a'_h}^{\hat
 b'_h} \log|x-y| \big( \partial_{\hat{n}_{j}}\mu^{\hat{\boldsymbol{n}}}_h(y)\big)\,\dd y.
 \end{split}
 \end{equation}
Since this expression does not depend on the choice of $x$ in the $i$-th band and since we have
\[
1=\int_{\hat a'_i}^{\hat b'_i} \big(\partial_{\hat{n}_i} \mu^{\hat{\boldsymbol{n}}}_i(x)\big)\dd x,
\]
we can multiply \eqref{eq_x119} by the function $\partial_{\hat{n}_i} \mu^{\hat{\boldsymbol{n}}}_{i}(x)$ --- whose support is inside the $i$-th band --- and integrate to get
\begin{equation}
\label{eq_x120}
 \partial_{\hat n_i}\partial_{\hat n_j} \I[\boldsymbol{\mu}^{\hat{\boldsymbol{n}}}]= \sum_{h=1}^H 2\theta_{h,i}\int_{\hat a'_i}^{\hat b'_{i}}\int_{\hat a'_h}^{\hat
 b'_h} \log|x-y| \big(\partial_{\hat{n}_i} \mu^{\hat{\boldsymbol{n}}}_i(x)\big) \big( \partial_{\hat{n}_j} \mu^{\hat{\boldsymbol{n}}}_h(y)\big)\dd x \dd y.
\end{equation}
On the other hand, if $g \neq i$, then
\[
0=\int_{\hat a'_{g}}^{\hat b'_{g}}\,
\partial_{\hat{n}_{i}}\mu^{\hat{\boldsymbol{n}}}_{g}(x) \dd x.
\]
Hence, multiplying \eqref{eq_x119} by
 $\partial_{\hat{n}_{i}} \mu^{\hat{\boldsymbol{n}}}_{g}(x)$ and integrating, we get
\begin{equation}
\label{eq_x121} 0= \sum_{h=1}^H 2\theta_{h,g}\int_{\hat a'_{g}}^{\hat b'_{g}}\int_{\hat
a'_h}^{\hat
 b'_h} \log|x-y| \big( \partial_{\hat{n}_{i}} \mu^{\hat{\boldsymbol{n}}}_{g}(x)\big) \cdot \big( \partial_{\hat{n}_{j}} \mu^{\hat{\boldsymbol{n}}}_h(y)\big)\dd x \dd y, \quad g\neq i.
\end{equation}
Summing \eqref{eq_x120} and \eqref{eq_x121} over all $g \neq i$, we finally get for any $i,j \in [H]$
\begin{equation}
\label{eq_x182}
 \partial_{\hat{n}_{i}}\partial_{\hat{n}_{j}} \I[\boldsymbol{\mu}^{\hat{\boldsymbol{n}}}]= \sum_{g,h=1}^{H} 2\theta_{h,g}\int_{\hat a'_{g}}^{\hat b'_{g}}\int_{\hat a'_{h}}^{\hat
 b'_{h}} \log|x-y|  \big(\partial_{\hat{n}_{i}} \mu^{\hat{\boldsymbol{n}}}_{g}(x)\big) \cdot \big(\partial_{\hat{n}_{j}} \mu^{\hat{\boldsymbol{n}}}_{h}(y)\big)\dd x \dd y.
\end{equation}

At this point we can use Corollaries~\ref{Corollary_I_positive} and \ref{Corollary_I_bound} to show that the Hessian of $\hat{\boldsymbol{n}} \mapsto -\I[\boldsymbol{\mu}^{\hat{\boldsymbol{n}}}]$ is positive definite. As we assume that ($\star$) does not have a unique solution, the tangent space to the space of solutions of ($\star$) in $\amsmathbb{R}^H$ is non-zero. Take any non-zero vector $\boldsymbol{\zeta} \in \amsmathbb{R}^H$ of unit Euclidean length in this tangent space. We want to show that
\begin{equation}
\label{eq_x183}
 \sum_{i,j=1}^H \zeta_i\zeta_j \partial_{\hat{n}_{i}}\partial_{\hat{n}_{j}} \I[\boldsymbol{\mu}^{\hat{\boldsymbol{n}}}]\leq - C \sum_{h = 1}^{H} \zeta_h^2 =-C.
\end{equation}
for some $C > 0$. By \eqref{eq_x182}, the left-hand side of \eqref{eq_x183} can be written in terms of Definition~\ref{DEFI2} as
\begin{equation}
\label{eq_x194}
2\, \I^{(2)}[\boldsymbol{\nu}]\qquad \textnormal{where}\quad \nu(x) = \sum_{i=1}^H \zeta_i \partial_{\hat{n}_{i}} \mu^{\hat{\boldsymbol{n}}}(x).
\end{equation}
By Corollary~\ref{Corollary_I_positive} the last expression is nonpositive. In order to bound it away from $0$,
we apply Corollary~\ref{Corollary_I_bound} to the $H$-tuple of integrable signed measures $\boldsymbol{\nu}$. The required condition \eqref{eq_nu_mass} holds precisely because we chose $\boldsymbol{\zeta}$ in the tangent space of the solutions to ($\star$) and we have
\[
\forall i,h \in [H]\qquad \int_{\hat a'_h}^{\hat b'_h} \big(\partial_{\hat{n}_{i}} \mu^{\hat{\boldsymbol{n}}}(x)\big) \dd x=\delta_{i,h}.
\]
The same identities and the fact that $\boldsymbol{\zeta}$ has unit norm guarantee that the left-hand side of \eqref{eq_I_lower_bound} is larger than some $C_1 > 0$ depending only on the constants in the assumptions.

We now turn to the last term of the right-hand side of \eqref{eq_I_lower_bound}, which involves the functional $\mathcal{J}_{\varkappa}$ of \eqref{eq_J_alpha} and Lemma~\ref{Lemma_J_alpha}. We observe that by the arguments at the beginning of the proof, the function $\partial_{\hat{n}_{i}} \mu^{\hat{\boldsymbol{n}}}(x)$ is uniformly bounded away from the endpoints of the bands, while near the endpoint of a band it diverges like an inverse square root. Choosing $\varkappa = \frac{1}{2}$ in Corollary~\ref{Corollary_I_bound} and using the integrability of the function $(x,y) \mapsto \frac{1}{\sqrt{|xy(x - y)|}}$ near $0$ --- here $\frac{1}{\sqrt{|x-y|}}$ arises in Lemma~\ref{Lemma_J_alpha}, while $\frac{1}{\sqrt{|x|}}$ and $\frac{1}{\sqrt{|y|}}$ come from the upper bounds for the functions $\partial_{\hat{n}_{i}} \mu^{\hat{\boldsymbol{n}}}(x)$ --- we conclude that there exists $C>0$, depending on the constants in the assumptions and such $\mathcal{J}_{\frac{1}{2}}[\boldsymbol{\nu}] < C_2$. Hence, we can choose $t$ large enough in Corollary~\ref{Corollary_I_bound} so that the second term in the right-hand side of \eqref{eq_I_lower_bound} is at most $\frac{C_1}{2}$, where $C_1$ is the constant from the previous paragraph. Then, \eqref{eq_I_lower_bound} implies
\[
 \frac{C_1}{2}\leq \frac{2}{\pi} \left(-\bigg(\frac{|\!|\hat{\boldsymbol{b}}' - \hat{\boldsymbol{a}}'|\!|_{\infty}^2}{8} + \log t\bigg)\I^{(2)}[\boldsymbol{\nu}]\right)^{\frac{1}{2}},
\]
which is the desired \eqref{eq_x183} and the upper bound in \eqref{soungnzn}. The lower bound is immediate from \eqref{eq_x194} and the regularity of the functions $\partial_{\hat{n}_{i}} \mu^{\boldsymbol{\hat{n}}}_h(x)$ for $i,h \in [H]$.
\end{proof}

\begin{proposition}
\label{Proposition_quadratic_potential_differentiability}
Consider a variational datum satisfying Assumptions~\ref{Assumption_A}, \ref{Assumption_B} and \ref{Assumption_C} with $H=1$ and admitting a single band --- we drop the indices $h$ from the notations. Let $\Delta V(z)$ be a holomorphic function in a complex neighborhood of $[\hat a', \hat b']$. Consider the family of variational data parameterized by a real $\epsilon$ in a small enough neighborhood of $0$ and differing from the original variational datum only by the potential $V^t(x) = V(x) + \epsilon \Delta V(x)$. Then, denoting with exponent ${}^{\epsilon}$ all the associated quantities, $\epsilon \mapsto \I^\epsilon[\boldsymbol{\mu}^{\epsilon}]$ is differentiable at $\epsilon = 0$ and its derivative can be computed as
 \begin{equation}
 \label{eq_t_derivative}
\partial_{\epsilon}\I^\epsilon[\mu^\epsilon]\big|_{\epsilon =0}=-\int_{\hat a'}^{\hat b'} \Delta V(x)\,\mu^0(x)\dd x.
 \end{equation}
\end{proposition}
\begin{proof}[Proof of Proposition~\ref{Proposition_quadratic_potential_differentiability}]
Note that because $H=1$, the single filling fraction is fixed. Repeating the proof of Theorem~\ref{Theorem_differentiability_full}, but differentiating with respect to $\epsilon$ in $V^\epsilon$ rather than with respect to segment filling fractions or endpoints, we conclude that $\epsilon \mapsto \Gm_{\mu^{\epsilon}}(z)$ is smooth for any $z \in \amsmathbb{C} \setminus [\hat{a}',\hat{b}']$. Arguing as in Proposition~\ref{Proposition_Hessian_free_energy} we further conclude that the density $\mu^\epsilon(x)$ is differentiable in $\epsilon$ at $\epsilon=0$ for any $x$ except the endpoints of the band of $\mu^0$, and that its derivative at $\epsilon = 0$ diverges at worst like an inverse square root near these points. For $x$ outside the band, we have $\partial_{\epsilon}\mu^\epsilon(x)|_{\epsilon=0}=0$. We can then differentiate the definition of $\I^\epsilon[\mu^\epsilon]$ as
\begin{equation*}
\begin{split}
\partial_{\epsilon}\I^\epsilon[\mu^{\epsilon}]\big|_{\epsilon=0} & = \partial_{\epsilon}\bigg(
 \int_{\hat a'}^{\hat b'} \int_{\hat a'}^{\hat
 b'} \theta \log|x-y|\mu^{\epsilon}(x) \mu^{\epsilon}(y)\dd x \dd y
-\int_{\hat a'}^{\hat b'} \big(V(x)+\epsilon\Delta V(x)\big)\mu^{\epsilon}(x)\dd x\bigg)\bigg|_{\epsilon=0} \\
& = -\int_{\hat a'}^{\hat b'} V^{\textnormal{eff},0}(x) \big(\partial_{\epsilon}\mu^{\epsilon}(x)\big)_{\epsilon = 0}\dd x -\int_{\hat a'}^{\hat b'} \Delta V(x) \mu^0(x) \dd x \\
& = -v\int_{\hat a'}^{\hat b'} \big(\partial_{\epsilon}\mu^{\epsilon}(x)\big)|_{\epsilon = 0} \dd x -\int_{\hat a'}^{\hat b'} \Delta V(x) \mu^0(x)\,\dd x \\
& =  -\int_{\hat a'}^{\hat b'} \Delta V(x) \mu^0(x)\dd x,
\end{split}
\end{equation*}
where $V^\textnormal{eff,0}(x)$ and $v$ are respectively the effective potential and the constant from the characterization of $\mu^0$ in Theorem~\ref{Theorem_equi_charact_repeat_2}. In the last equality we used the fact the total mass of $\mu^\epsilon$ is independent of $\epsilon$.
\end{proof}

\begin{proposition}
\label{proposition_Energy_series}
Consider a variational datum satisfying Assumptions~\ref{Assumption_A}, \ref{Assumption_B} and \ref{Assumption_C} and having a single band per segment. For each $h \in [H]$, take a function $\Delta V_h$ which is holomorphic in a complex neighborhood of $[\hat{a}'_h,\hat{b}_h']$. Then consider the family of variational data parameterized by a real $\epsilon$ in a small neighborhood of $0$, differing from the original variational datum only by the potentials $V_h^\epsilon = V_h + \epsilon \Delta V_h$ for any $h \in [H]$. Then, denoting with an exponent ${}^{\epsilon}$ all the associated quantities, we have the second-order Taylor expansion
\begin{equation}
 \label{eq_t_derivative_2}
\I^{\epsilon}[\mu^\epsilon]\mathop{=}_{\epsilon \rightarrow 0} \I^0[\mu^0]+ \I'[\mu^0] \epsilon + \I''[\mu^0] \epsilon^2 +o(\epsilon^2),
 \end{equation}
The error $o(\epsilon^2)$ is uniform as long as $|\Delta V_h(z)|$ is uniformly bounded for any $h \in [H]$ and the constants in Assumptions~\ref{Assumption_A}, \ref{Assumption_B}, and \ref{Assumption_C} remain the same. The factors $\I'[\mu^0]$ and $\I''[\mu^0]$ do not depend on $\epsilon$ and are uniformly bounded under the same assumptions.

Furthermore, if we vary the original variational datum as in Proposition~\ref{Proposition_differentiability_filling_fraction} and Theorem~\ref{Theorem_differentiability_full}, then the factors $\I'[\mu^0]$ and $\I''[\mu^0]$ depend smoothly on the variation parameters and the remainder $o(\epsilon^2)$ remains uniform.
\end{proposition}

It will be clear from the proof that a similar Taylor expansion holds to any fixed order under the same assumptions. We limited ourselves to state it up to order $2$ because this is enough to obtain in Chapters~\ref{Chapter_partition_functions} and \ref{Chapter_filling_fractions} the asymptotic expansion of $\log \Z_\N$ up to $o(1)$ as $\N \rightarrow \infty$.

\begin{proof}
 The first derivative in $\epsilon$ is computed for $H=1$ in \eqref{eq_t_derivative} and for $H>1$ the proof remains exactly the same. Next, note that we can represent
\[
 \forall h \in [H]\qquad \int_{\hat a'_h}^{\hat b'_h} \Delta V_h(x)\,\mu^0_h(x)\dd x =\frac{1}{2\pi \ii} \oint_{\gamma_h} \Delta V_h(z)\,\Gm_{\mu^0_h}(z) \dd z,
\]
This allows us to use the differentiability of $\Gm_{\mu_h^{\epsilon}}(z)$ proven in Proposition~\ref{Proposition_differentiability_filling_fraction} and Theorem~\ref{Theorem_differentiability_full}. We conclude that the $\epsilon$-derivative of $\I^\epsilon[\mu^\epsilon]$ depends smoothly on $\epsilon$, on filling fractions, on the endpoints of the segments, and on a possible auxiliary parameter $t$ that the variational datum may depend on, as in Theorem~\ref{Theorem_differentiability_full}. This implies the desired second-order Taylor expansion \eqref{eq_t_derivative_2}.
\end{proof}

\section{Parametric regularity with varying filling fractions}
In the setup of Section~\ref{Section_parameters} the filling fractions $\hat n_h$ were fixed and taken to be a part of the variational data. The equilibrium measure was explicitly depending on the filling fractions and we studied this dependence in Lemma~\ref{Lemma_continuity_filling_fractions}, Proposition~\ref{Proposition_differentiability_filling_fraction} and Proposition~\ref{Proposition_Hessian_free_energy}. This setup differs from the setting of discrete ensembles in Section~\ref{Section_general_model}, where the filling fractions were allowed to fluctuate subject to linear constraints \eqref{eq_equations_eqs}, which may not fix them deterministically. Similarly, in regard to the equilibrium measure, in Theorem~\ref{Theorem_equi_charact_repeat_2} the filling fractions were not prescribed but instead were optimizing themselves in the range allowed by \eqref{eq_equations_eqs} to minimize the energy functional $-\mathcal{I}$. In this section, we relax the setting of Section~\ref{Section_parameters} and allow filling fractions to vary. The spatial regularity results of Section~\ref{RegSecReg} remain unchanged, because any equilibrium measure with varying filling fractions is also an equilibrium measure with fixed filling fractions (chosen to be the optimal ones). Further, in the off-critical situation of Assumption~\ref{Assumption_C}, one can use Proposition~\ref{Proposition_Hessian_free_energy} to conclude that the optimal filling fractions smoothly depend on the remaining variational data. Therefore, the parametric smooth regularity of Sections~\ref{Section-ParReg} and \ref{Section_Smoothnessparam} continues to hold in the setting of varying filling fractions.

But, in the present section we want to work without the off-criticality Assumption~\ref{Assumption_C} and redo the arguments of Section~\ref{Section-ParReg} in the setting of varying filling fractions. Without off-criticality the results are weaker: rather than Lipschitz regularity we only get $\frac{1}{2}$-H\"older regularity. In Chapter~\ref{Chap11} where these results are applied, we will only use the implied continuity.

\subsection{Setup with varying filling fractions}

In comparison to Section~\ref{Section_variational_data}, we no longer fix the values of $\hat{\boldsymbol{n}} = (\hat{n}_h)_{h = 1}^H$. Instead, as in Section~\ref{DataS}, we are given a family $(\mathfrak{r}_e)_{e = 1}^{\mathfrak{e}}$ of linear forms in $\amsmathbb{R}^H$ with integral coefficients, and an $\mathfrak{e}$-tuple of real numbers $\boldsymbol{r} = (r_e)_{e = 1}^{\mathfrak{e}}$. The filling fractions $\hat{\boldsymbol{n}}$ are subject to the linear constraints \eqref{eq_equations_eqs}, which read
\begin{equation} \tag{$\star$}
\forall e \in [\mathfrak{e}]\qquad \mathfrak{r}_e(\hat{\boldsymbol{n}})= r_e.
\end{equation}
As before, $\mathfrak{r} : \amsmathbb{R}^H \rightarrow \amsmathbb{R}^{\mathfrak{e}}$ denotes the linear map with components $\mathfrak{r}_1,\ldots,\mathfrak{r}_{\mathfrak{e}}$. Assumption~\ref{Assumption_A} is replaced with the following.

\begin{lassumprime}
\label{Assumption_A_prime} There exists a constant $C > 0$ such that
\begin{enumerate}
\item $\boldsymbol{\Theta}$ is a real symmetric positive semi-definite $H \times H$ matrix;
\item $\forall h \in [H] \qquad \theta_{h,h} \geq \frac{1}{C}$;
\item $H \leq C$, $\mathfrak{e} \leq C$ , and $|\!|\boldsymbol{\Theta}|\!|_{\infty} \leq C$;
\item There exists a real symmetric positive semi-definite $\mathfrak{e} \times\mathfrak{e}$ matrix $\boldsymbol{\Theta}'$ such that
\begin{equation}
\label{eq_Theta_through_Theta_prime_repeat}
\forall \boldsymbol{X} \in \amsmathbb{R}^H\qquad \boldsymbol{X}^{T} \cdot \boldsymbol{\Theta} \cdot \boldsymbol{X} \\ = \big( \mathfrak r_1(\boldsymbol{X}),
 \ldots, \mathfrak r_{\mathfrak{e}}(\boldsymbol{X}) \big)^{T} \cdot \boldsymbol{\Theta}' \cdot \big( \mathfrak r_1(\boldsymbol{X}),
 \ldots, \mathfrak r_{\mathfrak{e}}(\boldsymbol{X}) \big).
\end{equation}
\item We have $|\!|\boldsymbol{r}|\!|_{\infty} < C$, the coefficients of the linear forms $\mathfrak{r}_1,\ldots,\mathfrak{r}_{\mathfrak{e}}$ are bounded by $C$, and the image of the parallelepiped $\prod_{h = 1}^H \big[0,\,\frac{1}{\theta_{h,h}}(\hat b'_h - \hat a'_h)\big]$ under the map $\boldsymbol{X} \mapsto (\mathfrak{r}_e(\boldsymbol{X}))_{e = 1}^{\mathfrak{e}}$ contains the closed ball of radius $\frac{1}{C}$ centered at $\boldsymbol{r}$.
\end{enumerate}
\end{lassumprime}
The rest of the setup is identical to Section~\ref{Section_variational_data}.

\subsection{H\"older parametric regularity and preservation of off-criticality}

In this section we have two variational data, denoted $\textnormal{I}$ and $\textnormal{II}$ and two corresponding equilibrium measures $\boldsymbol{\mu}^{\textnormal{I}}$ and $\boldsymbol{\mu}^{\textnormal{II}}$ defined as the unique minimizers of the energy functional $-\mathcal{I}$ among measures whose segment filling fractions $\boldsymbol{n}$ satisfy the affine constraints
\[
\forall e \in [\mathfrak{e}]\qquad \mathfrak{r}_e(\boldsymbol{n}) = r_e^{\textnormal{I}/\textnormal{II}}
\]
The integral parameters (like $H$, $\mathfrak{e}$ and the linear forms $\mathfrak{r}_e$ that have integral coefficients) are implicitly taken to be the same for the two variational data. As in Section~\ref{Section-ParReg}, the analyticity of $U_h(x)$ of Assumption~\ref{Assumption_B} is not necessary and all the results continue to hold for $U_h(x)$ which is only twice-continuously differentiable.

\begin{lemma}
\label{Lemma_continuity_varyinf_ff}
Take two variational data satisfying Assumptions~\ref{Assumption_A_prime} and \ref{Assumption_B}, differing by their matrices of interactions $\boldsymbol{\Theta}^{\textnormal{I}}$ and $\boldsymbol{\Theta}^{\textnormal{II}}$, by their segment endpoints $(\hat{a}_{h}^{\prime\,\textnormal{I}},\hat{b}_{h}^{\prime\,\textnormal{I}})_{h=1}^{H}$ and $(\hat{a}_{h}^{\prime\,\textnormal{II}},\hat{b}_{h}^{\prime\,\textnormal{II}})_{h=1}^H$, by their potentials $(V^{\textnormal{I}}_h)_{h=1}^H$ and $(V^{\textnormal{II}}_h)_{h=1}^H$, by the right-hand sides $(r_e^{\textnormal{I}})_{e=1}^{\mathfrak{e}}$ and $(r_e^{\textnormal{II}})_{e=1}^{\mathfrak{e}}$ of the affine constraints \eqref{eq_equations_eqs}. There exists $M > 0$ depending only on the constants in the assumptions, such that for any measurable function $f$ with finite $|\!|f|\!|_{\frac{1}{2}}$ and $|\!|f|\!|_{\infty}$ we have
\begin{equation}
 \label{eq_continuity_varying_ff}
 \begin{split}
&  \bigg| \int_{\amsmathbb R} f(x)\dd(\mu^{\textnormal{I}}-\mu^{\textnormal{II}})(x)\bigg| \leq M \cdot \big(|\!|f|\!|_{\frac{1}{2}}+|\!|f|\!|_{\infty}\big)  \\
& \quad \times \Biggl[ \sum_{e=1}^{\mathfrak e} \bigl|r_e^{\textnormal{I}}-r_e^{\textnormal{II}}\bigr| + \sum_{g,h=1}^H \bigl|\theta^{\textnormal{I}}_{g,h}- \theta^{\textnormal{II}}_{g,h}\bigr| + \sum_{h=1}^H \Big(\bigl|\hat{a}_{h}^{\prime\,\textnormal{I}}-\hat{a}_{h}^{\prime\,\textnormal{II}}\bigr|+ \bigl|\hat{b}_{h}^{\prime\,\textnormal{I}}-\hat{b}_{h}^{\prime\,\textnormal{II}}\bigr| + \sup_{x} \bigl|V_h^{\textnormal{I}}(x) - V_{h}^{\textnormal{II}}(x)\bigr|\Big)\Biggr]^{\frac{1}{2}},
 \end{split}
 \end{equation}
 where the supremum ranges over $x \in [\max(\hat{a}_{h}^{\prime\,\textnormal{I}},\hat{a}_{h}^{\prime\,\textnormal{II}}),\min(\hat{b}_{h}^{\prime\,\textnormal{I}},\hat{b}_{h}^{\prime\,\textnormal{II}})]$.
\end{lemma}
Let us emphasize the square root in the last formula, which makes the conclusion weaker than those of Lemmata~\ref{Lemma_continuity_interactions}, \ref{Lemma_continuity_potential}, \ref{Lemma_continuity_filling_fractions}, \ref{Lemma_continuity_end_points_void}, \ref{Lemma_continuity_end_points_saturated}. On the other hand, we also have less assumptions than in those lemmata. An important consequence of Lemma~\ref{Lemma_continuity_varyinf_ff} is that the property of preservation of off-criticality of Theorem~\ref{Theorem_off_critical_neighborhood} still holds in the situation where filling fractions are not completely fixed \textit{a priori}.
\begin{corollary}
\label{co:Lemregwithout}
In the setting of Lemma~\ref{Lemma_continuity_varyinf_ff}, there exists $\varepsilon > 0$ depending only on the constants in the assumptions with the following property. If $\boldsymbol{\mu}^{\textnormal{I}}$ satisfies the off-criticality Assumption~\ref{Assumptions_offcrit} and for any $g,h \in [H]$ and $e \in [\mathfrak{e}]$
\begin{equation}
\label{ineqesps}
\begin{split}
\max\big(|\hat a_h^{\prime\,\textnormal{I}}-\hat a_h^{\prime\,\textnormal{II}}|\,,\,|\hat b_h^{\prime\,\textnormal{I}}-\hat b_h^{\prime\,\textnormal{II}}|\,,\,|\theta_{g,h}^{\textnormal{I}} - \theta_{g,h}^{\textnormal{II}}|\,,\,|r_e^{\textnormal{I}}-r_e^\textnormal{II}|\big) & \leq \eps, \\
 \int_{\max (\hat a_h^{\prime\,\textnormal{I}},\hat a_h^{\prime\,\textnormal{II}})}^{\min(\hat b_h^{\prime\,\textnormal{I}},\hat b_h^{\prime\,\textnormal{II}})} \big|V_h^{\textnormal{I}}(x)-V_h^\textnormal{II}(x)\big| \dd x & \leq \eps, \\
\sup_{z\in \amsmathbb{M}_h} \big|\phi^{\pm,\textnormal{I}}_h(z)-\phi^{\pm,\textnormal{II}}_h(z)\big| & \leq \eps,
\end{split}
\end{equation}
then $\boldsymbol{\mu}^{\textnormal{II}}$ also satisfies the off-criticality Assumption~\ref{Assumptions_offcrit}.
\end{corollary}

\begin{proof}
The vector of filling fractions $\boldsymbol{n}^{\textnormal{I}/\textnormal{II}}$ for the measure $\boldsymbol{\mu}^{\textnormal{I}/\textnormal{II}}$ satisfy the affine constraints $(\star)$ with right-hand side $(r_e^{\textnormal{I}/\textnormal{II}})_{e = 1}^{\mathfrak{e}}$. As it is the unique minimizer of the energy functional $-\mathcal{I}$ among measures whose filling fractions satisfy these affine constraints, it is also uniquely minimizes $-\mathcal{I}$ among measures whose filling fractions are equal to $\boldsymbol{n}^{\textnormal{I}/\textnormal{II}}$. Therefore, $\boldsymbol{\mu}^{\textnormal{I}/\textnormal{II}}$ are also measured concerned with the previous results of this chapter modulo the additional input of the values $\boldsymbol{n}^{\textnormal{I}/\textnormal{II}}$, in particular Theorem~\ref{Theorem_off_critical_neighborhood}. To apply this theorem, we just need to check that $\boldsymbol{n}^{\textnormal{II}}$ stays arbitrarily close to $\boldsymbol{n}^{\textnormal{I}}$ if we assume the inequalities \eqref{ineqesps} with $\varepsilon > 0$ small enough. This is guaranteed by Lemma~\ref{Lemma_continuity_varyinf_ff}, which tells us that $\boldsymbol{n}^{\textnormal{II}}$ has $\frac{1}{2}$-H\"older dependence in the parameters of the variational datum $\textnormal{II}$.
\end{proof}

\subsection{Proof of Lemma~\ref{Lemma_continuity_varyinf_ff}}

\noindent \textsc{Step 1.} We start from the case when only $(V_h)_{h=1}^H$ differs between the variational data $\textnormal{I}$ and $\textnormal{II}$, while all other parameters are the same. The argument is similar to the first part of the proof of Lemma~\ref{Lemma_continuity_potential}. We denote, for $\textnormal{J} \in \{\textnormal{I},\textnormal{II}\}$, by $-\I^\textnormal{J}$ and $V_h^{\textnormal{eff},J}$ the energy functional and the $h$-th effective potential associated with the variational datum $\textnormal{J}$, by $V_h^\textnormal{J}$ the potential, by $\boldsymbol{\mu}^\textnormal{J}$ the equilibrium measure. Denote $\boldsymbol{\nu}:=\boldsymbol{\mu}^\textnormal{I}-\boldsymbol{\mu}^\textnormal{II}$. We claim that
\begin{equation}
\label{eq_x286}
 \sum_{h=1}^H \int_{\hat a'_{h}}^{\hat b'_{h}} V^{\textnormal{eff},\textnormal{II}}_{h}(x) \nu_h(x) \dd x \geq 0.
\end{equation}
The proof of this inequality is slightly different from similar parts of the arguments in Section~\ref{Section-ParReg} (\textit{cf.} \eqref{eq_x76}) as we no longer can directly rely on Theorem~\ref{Theorem_equi_charact_repeat_2} --- the step where we subtracted $v_h$ from the potential does not work if $\boldsymbol{\mu}^\textnormal{II}$ does not have bands and filling fractions are allowed to vary. Instead, we note that by convexity of the set of measures $\mathscr{P}_\star$ of \eqref{eq_ff_set_of_measures}, $\boldsymbol{\mu}^\textnormal{I}\in \mathscr{P}_\star$ and $\boldsymbol{\mu}^\textnormal{II} \in \mathscr{P}_\star$ implies $\boldsymbol{\mu}^\textnormal{II}+\eps \boldsymbol{\nu}\in \mathscr{P}_\star$ for all $\eps \in (0,1)$. As in \eqref{eq_x125}, we have
\[
 \I[\boldsymbol{\mu}^\textnormal{II}+\eps\boldsymbol{\nu}] = \I^{\textnormal{II}}[\boldsymbol{\mu}^{\textnormal{II}}]-\eps \sum_{h=1}^H \int_{\hat a'_h }^{\hat b'_h} V_h^{\textnormal{eff},\textnormal{II}}(x)\, \nu_h(x) \dd x+
 \eps^2 \sum_{g,h=1}^H \theta_{g,h} \int_{\hat a'_g}^{\hat b'_g} \int_{\hat a'_{h}}^{\hat b'_{h}}
 \log|x-y|\, \nu_g(x) \nu_{h}(y) \dd x \dd y.
\]
This expression should be smaller than $\I^\textnormal{II}[\boldsymbol{\mu}^\textnormal{II}]$ because $\boldsymbol{\mu}^{\textnormal{II}}$ is the minimizer of $-\I^{\textnormal{II}}$. Sending $\eps \rightarrow 0^+$, we arrive at \eqref{eq_x286}.

Next, using the definition of $\boldsymbol{\mu}^\textnormal{I}$ as a minimizer of $-\I^\textnormal{I}$ over the measures in $\mathscr{P}_\star$ and recalling Definition~\ref{Definition_pseudodistance} for $\mathfrak{D}^2$ which is common to the two ensembles, we have:
\begin{equation}
\label{eq_x285}
\begin{split}
0 & \leq \I^\textnormal{I}[\boldsymbol{\mu}^\textnormal{I}]-\I^\textnormal{I}[\boldsymbol{\mu}^{\textnormal{II}}] \\
& = \sum_{g,h = 1}^H \int_{\hat
a'_{g}}^{\hat b'_{g}} \int_{\hat a'_{h}}^{\hat b'_{h}} \theta_{g,h}
\log|x-y|\, \big([(\mu^{\textnormal{II}}_{g}+\nu_{g})(x)\mu^{\textnormal{II}}_{h}+\nu_{h})(y) -
\mu^{\textnormal{II}}_{g}(x)\dd\mu^{\textnormal{II}}_{h}(y)\bigr]\dd x\dd y \\
& \quad -\sum_{h=1}^H \int_{\hat a'_{h}}^{\hat b'_{h}} V_{h}^{\textnormal{I}}(x)
\nu_h(x) \dd x \\
& = - \D[\boldsymbol{\mu}^\textnormal{I},\boldsymbol{\mu}^{\textnormal{II}}]-\sum_{h=1}^H\int_{\hat a'_h}^{\hat b'_h}
\big((V^{\textnormal{I}}_h(x)-V^{\textnormal{II}}_h(x)\big)\nu_h(x) \dd x -\sum_{h=1}^{H} \int_{\hat a'_{h}}^{\hat
b'_{h}} V^{\textnormal{eff},\textnormal{II}}_{h}(x) \nu_h(x)\dd x.
\end{split}
\end{equation}
Thanks to \eqref{eq_x286}, we deduce
\begin{equation}
 \D[\boldsymbol{\mu}^\textnormal{I},\boldsymbol{\mu}^{\textnormal{II}}]\leq \sum_{h=1}^H\int_{\hat a'_h}^{\hat b'_h}
\big((V_h^{\textnormal{II}}(x)-V_h^{\textnormal{I}}(x)\big)\nu_h(x)\dd x.
\end{equation}
Since $\nu_h$ is a difference of two positive measures of uniformly bounded total masses, we conclude
\begin{equation}
 \D[\boldsymbol{\mu}^\textnormal{I},\boldsymbol{\mu}^{\textnormal{II}}]\leq C \sum_{h=1}^H\sup_{x\in [\hat a'_h,\hat b'_h]}
\big(|V_h^{\textnormal{II}}(x)-V_h^{\textnormal{I}}(x)\bigr|,
\end{equation}
for a similarly uniform constant $C>0$. Using Lemma~\ref{Lemma_linear_through_distance}, we deduce \eqref{eq_continuity_varying_ff}.

\bigskip

\noindent \textsc{Step 2.} In the general case, the new difficulty is that the set of measures $\mathscr{P}_\star$ differs between the variational data $\textnormal{I}$ and $\textnormal{II}$, for three reasons: the segments where the measures are supported may differ, the maximal density $\theta_{h,h}^{-1}$ may differ, and the right-hand side $\boldsymbol{r}$ of the constraints \eqref{eq_equations_eqs} may differ. We call $\mathscr{P}_\star^\textnormal{I}$ and $\mathscr{P}_\star^\textnormal{II}$ the two corresponding sets. It only makes sense to evaluate $\I^\textnormal{I}$ on the measures of $\mathscr{P}_\star^\textnormal{I}$ and to evaluate $\I^\textnormal{II}$ on the measures from $\mathscr{P}_\star^\textnormal{II}$ and the new ingredient is to show that we can move from between these two sets of measures up to a small error.

\begin{lemma} \label{Lemma_classes_close} Denote
\[
 \mathfrak{d} = \sum_{e=1}^{\mathfrak e}\big(|r_e^{\textnormal{I}}-r_e^{\textnormal{II}}\bigr| + \sum_{h=1}^H \Big( \big(|\theta^{\textnormal{I}}_{h,h}- \theta^{\textnormal{II}}_{h,h}\bigr| + \big(|\hat{a}_{h}^{\prime\,\textnormal{I}}-\hat{a}_{h}^{\prime\,\textnormal{II}}\bigr| + \big(|\hat{b}_{h}^{\prime\,\textnormal{I}}-\hat{b}_{h}^{\prime\,\textnormal{II}}\bigr|\Big).
\]
In the setting of Lemma~\ref{Lemma_continuity_varyinf_ff} there exists $C>0$ (depending only on the constants in Assumptions \ref{Assumption_A_prime} and \ref{Assumption_B}) and two measures $\boldsymbol{\mu}^{\textnormal{I} \rightarrow \textnormal{II}}$ and $\boldsymbol{\mu}^{\textnormal{II} \rightarrow \textnormal{I}}$, such that:
\begin{align}
\label{eq_x287}&\boldsymbol{\mu}^{\textnormal{I} \rightarrow \textnormal{II}}\in \mathscr{P}_\star^\textnormal{II}, \qquad \sum_{h=1}^H \int_{\hat a'_h}^{\hat b'_h}\left|\mu^{\textnormal{I} \rightarrow \textnormal{II}}_h(x)- \mu^{\textnormal{I}}_h(x)\right|\dd x < C \cdot \mathfrak{d},\\
&\boldsymbol{\mu}^{\textnormal{II} \rightarrow \textnormal{I}}\in \mathscr{P}_\star^\textnormal{I}, \qquad \sum_{h=1}^H \int_{\hat a'_h}^{\hat b'_h} \left|\mu^{\textnormal{II} \rightarrow \textnormal{I}}_h(x)- \mu^{\textnormal{II}}_h(x)\right|\dd x < C \cdot \mathfrak{d}, \label{eq_x288}
\end{align}
where $\hat a'_h= \min(\hat{a}_{h}^{\prime\,\textnormal{I}},\hat{a}_{h}^{\prime\,\textnormal{II}})$ and $\hat b'_h=\max(\hat{b}_{h}^{\prime\,\textnormal{I}},\hat{b}_{h}^{\prime\,\textnormal{II}})$.
\end{lemma}
\begin{proof}
 Starting from $\boldsymbol{\mu}^{\textnormal{I}}$, we construct $\boldsymbol{\mu}^{\textnormal{I} \rightarrow \textnormal{II}}$ by the following modifications:
 \begin{enumerate}
  \item Multiply the density $\mu^{\textnormal{I}}_h(x)$ by $\frac{\theta^{\textnormal{I}}_{h,h}}{\theta^{\textnormal{II}}_{h,h}}$ for each $h\in[H]$.
  \item Replace the density by $0$ everywhere outside $\bigcup_{h=1}^H [ \max(\hat{a}_{h}^{\prime\,\textnormal{I}},\hat{a}_{h}^{\prime\,\textnormal{II}}),\min(\hat{b}_{h}^{\prime\,\textnormal{I}},\hat{b}_{h}^{\prime\,\textnormal{II}})]$.
 \end{enumerate}
 Clearly, the bound \eqref{eq_x287} is satisfied after these transformations. Note that the right-hand sides of \eqref{eq_equations_eqs} also got deformed as a result, but increments are at most $C\cdot \mathfrak{d}$. We next use Lemma~\ref{Lemma_filling_modification}, which says that it is possible to slightly deform the filling fractions, as to satisfy the constraints \eqref{eq_equations_eqs} corresponding to the variational data $\textnormal{II}$. The second condition in \eqref{iugbnrugbr} guarantees that this deformation can be achieved by slightly increasing/decreasing densities $\mu_h(x)$, and we eventually arrive at a measure from $\mathscr{P}_\star^\textnormal{II}$. The first condition in \eqref{iugbnrugbr} implies the inequality \eqref{eq_x287} and we are done constructing $\boldsymbol{\mu}^{\textnormal{I} \rightarrow \textnormal{II}}$. The argument for $\boldsymbol{\mu}^{\textnormal{II} \rightarrow \textnormal{I}}$ is the same.
\end{proof}

\noindent \textsc{Step 3.} We now use Lemma~\ref{Lemma_classes_close} to finish the proof of Lemma~\ref{Lemma_continuity_varyinf_ff} by upgrading the arguments of Step 1 for the general case. The inequality \eqref{eq_x285} is replaced with
\begin{equation}
\label{eq_x290}
\begin{split}
0 & \leq \I^\textnormal{I}[\boldsymbol{\mu}^\textnormal{I}]-\I^\textnormal{I}[\boldsymbol{\mu}^{\textnormal{II} \rightarrow \textnormal{I}}] \\
& = - \D[\boldsymbol{\mu}^\textnormal{I},\boldsymbol{\mu}^{\textnormal{II} \rightarrow \textnormal{I}}]-\sum_{h=1}^H\int_{\hat a^{\prime\,\textnormal{I}}_h}^{\hat b^{\prime\,\textnormal{I}}_h}
\big((V^{\textnormal{I}}_h(x)-V^{\textnormal{II}}_h(x)\big)\nu_h(x) \dd x -\sum_{h=1}^{H} \int_{\hat a^{\prime\,\textnormal{I}}_{h}}^{\hat
b^{\prime\,\textnormal{I}}_{h}} V^{\textnormal{eff},\textnormal{II} \rightarrow \textnormal{I}}_{h}(x) \nu_h(x)\dd x,
\end{split}
\end{equation}
where $\boldsymbol{\nu}=\boldsymbol{\mu}^\textnormal{I}-\boldsymbol{\mu}^{\textnormal{II} \rightarrow \textnormal{I}}$ and
\[
V^{{\textnormal{eff},\textnormal{II} \rightarrow \textnormal{I}}}_h(x)= V^{\textnormal{II}}_h(x) - \sum_{g =1}^H 2\theta^{\textnormal{I}}_{h,g}
\int_{\hat a^{\prime\,\textnormal{I}}_{g}}^{\hat b^{\prime\,\textnormal{I}}_g} \log|x-y|\mu^{\textnormal{II} \rightarrow \textnormal{I}}_{g}(y)\dd y.
\]
Let us analyze the last term in \eqref{eq_x290}. For that we introduce $\widetilde{\boldsymbol{\nu}}=\boldsymbol{\mu}^{\textnormal{I} \rightarrow \textnormal{II}}-\boldsymbol{\mu}^{\textnormal{II}}$ and write for each $h \in [H]$
\begin{equation}
\label{eq_x289}
\begin{split}
\sum_{h = 1}^H \int_{\hat a^{\prime\,\textnormal{I}}_{h}}^{\hat
b^{\prime\,\textnormal{I}}_{h}} V^{\textnormal{eff},\textnormal{II} \rightarrow \textnormal{I}}_{h}(x) \nu_h(x)\dd x & = \sum_{h = 1}^H
 \int_{\hat a^{\prime\,\textnormal{II}}_{h}}^{\hat
b^{\prime\,\textnormal{II}}_{h}} V^{\textnormal{eff},\textnormal{II}}_{h}(x) \widetilde{\nu}_h(x)\dd x \\
& \quad + \sum_{h = 1}^H \int_{\hat a^{\prime\,\textnormal{II}}_{h}}^{\hat
b^{\prime\,\textnormal{II}}_{h}} \big(V^{\textnormal{eff},\textnormal{II} \rightarrow \textnormal{I}}_{h}(x)-V^{\textnormal{eff},\textnormal{II}}_{h}(x)\big) \widetilde{\nu}_h(x)\dd x \\
&\quad + \sum_{h = 1}^H
\int V^{\textnormal{eff},\textnormal{II} \rightarrow \textnormal{I}}_{h}(x) (\nu_h(x)-\widetilde {\nu}_h(x))\dd x,
\end{split}
\end{equation}
where the last integration is over the support of $\nu_h - \widetilde {\nu}_h$. The first sum in the right-hand side is nonnegative by the same argument as for \eqref{eq_x286}. For the second sum, note that
\begin{equation*}
\begin{split}
 & \quad  \int_{\hat a^{\prime\,\textnormal{II}}_{h}}^{\hat
b^{\prime\,\textnormal{II}}_{h}} \left(V^{\textnormal{eff},\textnormal{II} \rightarrow \textnormal{I}}_{h}(x)-V^{\textnormal{eff},\textnormal{II}}_{h}(x)\right) \widetilde{\nu}_h(x)\dd x \\
&= 2 \sum_{g =1}^H \left[\theta^{\textnormal{II}}_{h,g}
\int_{\hat a^{\prime\,\textnormal{II}}_{g}}^{\hat b^{\prime\,\textnormal{II}}_g}\mu^{\textnormal{II}}_{g}(y)-\theta^{\textnormal{I}}_{h,g}
\int_{\hat a^{\prime\,\textnormal{I}}_{g}}^{\hat b^{\prime\,\textnormal{I}}_g}\mu^{\textnormal{II} \rightarrow \textnormal{I}}_{g}(y)\right] \left[ \int_{\hat a^{\prime\,\textnormal{II}}_{h}}^{\hat
b^{\prime\,\textnormal{II}}_{h}} \log|x-y| \widetilde{\nu}_h(x)\dd x\right]\dd y.
\end{split}
\end{equation*}
Since the second $[\cdot]$ factor is uniformly bounded, \eqref{eq_x288} implies that the absolute value of the second sum in \eqref{eq_x289} is upper-bounded by $C\cdot \mathfrak{d}'$, where
\begin{equation*}
\begin{split}
 \mathfrak{d}' & = \sum_{e=1}^{\mathfrak e}\big|r_e^{\textnormal{I}}-r_e^{\textnormal{II}}\big| + \sum_{g,h = 1}^H \big|\theta^{\textnormal{I}}_{g,h}- \theta^{\textnormal{II}}_{g,h}\big| \\
 & \quad + \sum_{h=1}^H \Big( |\hat{a}_{h}^{\prime\,\textnormal{I}}-\hat{a}_{h}^{\prime\,\textnormal{II}}\bigr|+ |\hat{b}_{h}^{\prime\,\textnormal{I}}-\hat{b}_{h}^{\prime\,\textnormal{II}}\bigr| + \sup_{x\in [\max(\hat{a}_{h}^{\prime\,\textnormal{I}},\hat{a}_{h}^{\prime\,\textnormal{II}}),\min(\hat{b}_{h}^{\prime\,\textnormal{I}},\hat{b}_{h}^{\prime\,\textnormal{II}})]} \big|V_h^{\textnormal{I}}(x) - V_{h}^{\textnormal{II}}(x)\big| \Big).
\end{split}
\end{equation*}
The third sum in \eqref{eq_x289} is also upper-bounded by $C\cdot \mathfrak{d}'$, because $V^{\textnormal{eff},\textnormal{II} \rightarrow \textnormal{I}}_{h}$ is uniformly bounded and $\nu_h-\widetilde {\nu}_h$ is small by Lemma~\ref{Lemma_classes_close}.

Further bounding the second sum in the right-hand side of \eqref{eq_x290} as in Step 1, we get
$\D[\boldsymbol{\mu}^\textnormal{I},\boldsymbol{\mu}^{\textnormal{II} \rightarrow \textnormal{I}}]\leq C\cdot \mathfrak{d}'$. Hence, by Lemma~\ref{Lemma_linear_through_distance}
\[
 \bigg| \int_{\amsmathbb R} f(x)\dd(\mu^{\textnormal{I}}-\mu^{\textnormal{II} \rightarrow \textnormal{I}})(x)\bigg|\leq C\dot |\!|f|\!|_{\frac{1}{2}} \sqrt{\mathfrak{d}'}.
\]
Besides, \eqref{eq_x288} implies
\[
 \bigg| \int_{\amsmathbb R} f(x)\dd(\mu^{\textnormal{II}}-\mu^{\textnormal{II} \rightarrow \textnormal{I}})(x)\bigg|\leq C\cdot |\!|f|\!|_{\infty} \cdot \mathfrak{d}'.
\]
Using the triangular inequality, and noticing that $\mathfrak{d}'$ is uniformly bounded so that $\mathfrak{d}' \leq C'\sqrt{\mathfrak{d}'}$ for some uniform constant $C' > 0$, we arrive at \eqref{eq_continuity_varying_ff}.

\medskip

\begin{remark}
\label{sec:unrestricted}
As in Remark~\ref{remark:continuous}, we can also associate to variational data the ''unrestricted equilibrium measure'', that is the minimizer of $-\mathcal{I}$ over all $H$-tuples $\boldsymbol{\mu} = (\mu_h)_{h = 1}^{H}$ of nonnegative measures with support in $[\hat{a}'_{h},\hat{b}'_{h}]$ and such that $\mu_h([\hat{a}_{h},\hat{b}_{h}]) = \hat{n}_{h}$, for any $h \in [H]$, \textit{i.e.} not restricting to measures with density bounded by $\theta_{h,h}^{-1}$. The unrestricted equilibrium measure will be used in the proof of Theorem~\ref{thm:Bsym} (we will only need its existence, uniqueness and characterization in terms of the effective potential). In all the rest of the book we only work with the equilibrium measure studied including the upper bound on the density.

The unrestricted equilibrium is relevant for the continuous analogue of our discrete ensembles, \textit{i.e.} statistical mechanics ensembles defined like in Section~\ref{Section_general_model} but dropping all the discrete aspects, in particular the $\ell_i$ can vary in the real segments $[a_h,b_h]$. These continuous ensembles are much simpler to define than the discrete ensembles and have already been studied, see \textit{e.g.} \cite{BGK,BEO}. In particular, in \cite{BGK} the analogue of Lemma~\ref{Lemma_continuity_filling_fractions}, of Theorem~\ref{Theorem_off_critical_neighborhood} but allowing only filling fractions to change, and of Proposition~\ref{Proposition_Hessian_free_energy} have been proved. The interested reader can follow the proofs of the present chapter and observe many simplifications when the equilibrium measure has no saturations. This gives more complete results about the unrestricted equilibrium measure that are worth mentioning, namely the analogue of Lemmata~\ref{Lemma_continuity_interactions}, \ref{Lemma_continuity_potential} (continuity of the unrestricted equilibrium measure with respect to all parameters), Theorem~\ref{Theorem_off_critical_neighborhood} (preservation of off-criticality under variations of all parameters), Theorem~\ref{Theorem_differentiability_full} (smooth dependence of the Stieltjes transform under variation of all parameters) hold, \textit{etc}.
\end{remark}

\chapter{Reference ensembles and \texorpdfstring{$zw$}{zw}-measures}
\label{Chapterzw}

In Chapter~\ref{Chapter_partition_functions}, we will determine the asymptotic expansion of the partition function of the general discrete ensemble with fixed segment filling fractions by interpolation with reference ensembles having a single band. Therefore, we need to analyze the asymptotics for such a reference ensemble, where the position of the band and the nature (void or saturated) of its left and right neighboring regions can be tuned at will.
The $zw$-measures provide such a reference. They fit within the framework of discrete ensembles of Chapter~\ref{Chapter_Setup_and_Examples}, and generalize the projection of the uniform measure on lozenge tilings of a hexagon onto a vertical section (corresponding to $\theta = 1$). This chapter is devoted to the study of the $zw$-measures, more precisely a specialization of them that depends on four parameters (in addition to $\theta$ and $N$) necessary for the aforementioned tuning. We use an explicit product formula for their partition function, which allows us to derive easily their asymptotic expansion. We also compute explicitly their equilibrium measure and establish the corresponding phase diagram in terms of the four parameters.

\section{The \texorpdfstring{$zw$}{zw}-measure and its partition function}

The $zw$-measure depends on $\theta>0$ and four complex parameters $\mathsf{z}_1,\mathsf{z}_2,\mathsf{w}_1,\mathsf{w}_2$. It is a --- \textit{a priori}, complex --- measure on $N$-tuples
 \begin{equation}
 \label{eq_zw_state_space}
 \ell_1<\ell_2<\cdots<\ell_N,\qquad
 \ell_{i+1}-\ell_i \in \theta + \amsmathbb Z_{\geq 0},\qquad \ell_N\in\amsmathbb Z,
 \end{equation}
 which is given by the formula
\begin{equation}
 \label{eq_zw_1}
\frac{1}{\Z^{zw}_{N}}\cdot \prod_{1\leq i<j\leq N} \frac{\Gamma\big(\ell_j-\ell_i+1\big) \cdot\Gamma\big(\ell_j-\ell_i+\theta\big) }{\Gamma\big(\ell_j-\ell_i\big)\cdot
 \Gamma\big(\ell_j-\ell_i+1-\theta\big)} \cdot \prod_{i=1}^N \prod_{k = 1}^2 \frac{1}{\Gamma\big(\mathsf{z}_k+1-\ell_i\big)\cdot
\Gamma\big(\mathsf{w}_k+(N-1)\theta+1 +\ell_i\big)}.
\end{equation}

We will choose the parameters $\mathsf{z}_1,\mathsf{z}_2,\mathsf{w}_1,\mathsf{w}_2$ in such a way that \eqref{eq_zw_1} is
nonnegative for all relevant values of $\ell_i$ and has finite total mass. There are several classes of
possible choices here, and we refer to \cite{Olsh_hyper} and references therein for a discussion.

At $\theta=1$, the measure \eqref{eq_zw_1} for a specific integer choices of
$(\mathsf{z}_1,\mathsf{z}_2,\mathsf{w}_1,\mathsf{w}_2)$ appears in random lozenge tilings of the hexagon, as first observed in
\cite{CLP}. Such measures were introduced in another context by Kerov \cite{Kerov_zw}. For
$\theta \in \big\{\frac{1}{2},1,2\big\}$ the $zw$-measures are of great interest in the asymptotic
representation theory, \textit{cf.} \textit{e.g.}\ \cite{Olsh_zw,BO_zw,Neretin_zw}.
The case of a general $\theta$ is treated in details in\footnote{The $\ell_i$ here corresponds to the weakly decreasing sequence of integers $\nu_i = \ell_{N - i + 1} + (i - 1)\theta$ in \cite{Olsh_hyper}.} \cite{Olsh_hyper}; in particular, the
partition function was computed there.

\begin{proposition} \cite[(2.2)]{Olsh_hyper} \label{Proposition_Olsh_ZW} Suppose $\,\text{Re}(\mathsf{z}_1 + \mathsf{z}_2 + \mathsf{w}_1 + \mathsf{w}_2 + 1) > 0$ and set
\begin{equation}
 \label{eq_zw_normalization}
 \begin{split}
\Z^{zw}_{N} & = \prod_{1\leq i<j\leq N} \frac{\Gamma\big(\theta(j-i)+1\big) \cdot \Gamma\big(\theta(j-i)+\theta\big) }{\Gamma\big(\theta(j-i)\big)\cdot
 \Gamma\big(\theta(j-i)+1-\theta\big)} \\
 & \quad \times \prod_{i=1}^N \frac{\Gamma\big(\mathsf{z}_1+\mathsf{z}_2+\mathsf{w}_1+\mathsf{w}_2+(i-1)\theta+1\big)}{\Gamma\big((i-1)\theta+1\big) \prod_{k,\ell \in \{1,2\}} \Gamma\big(\mathsf{z}_k + \mathsf{w}_{\ell} +(i-1)\theta+1\big)}.
 \end{split}
 \end{equation}
 If the parameters are chosen so that this quantity does not vanish, then the measure \eqref{eq_zw_1} has mass $1$.
\end{proposition}

\section{Definition of the \texorpdfstring{$zw$}{zw}-discrete ensemble}
\label{SChoice}

We will now introduce a discrete ensemble in the sense of Chapter~\ref{Chapter_Setup_and_Examples} whose measure is proportional to the $zw$-measure with a certain choice of parameters. In the remaining of the text, we will call it generically the \emph{$zw$-discrete ensemble}.
It is specified by five real numbers $\theta>0$, $A_1,A_2,B_1,B_2$ and a positive integer $N$ satisfying the integrality condition
\begin{equation}
\label{integrazwwww}
B_1-(A_1+(N-1)\theta) \in \amsmathbb{Z}_{> 0}.
\end{equation}
We then set
\begin{equation}
 \label{zwchoi} \mathsf{z}_1 :=B_1,\quad \mathsf{z}_2 :=B_2,\quad \mathsf{w}_1 :=-A_1-\theta(N-1),\quad \mathsf{w}_2 :=-A_2-\theta(N-1).
\end{equation}
and consider configurations of $N$ particles satisfying
\begin{equation} \label{eq_ZW_model_state_space}
 \ell_i\in[A_1, B_1],\quad \ell_1\in A_1 + \amsmathbb Z_{\geq 0}, \quad \ell_{N}\in B_1 - \amsmathbb Z_{\geq 0},
 \quad \ell_{i+1}-\ell_i\in \theta + \amsmathbb Z_{\geq 0},
\end{equation}
Note that we do not require $B_1$ to be an integer, and so $\ell_N$ does not
have to be an integer unlike in \eqref{eq_zw_state_space}. These conditions differ from \eqref{eq_zw_state_space} by a translation common to all particles, therefore this change does not spoil the applicability of Proposition~\ref{Proposition_Olsh_ZW}.

When we take a large parameter $\N$ and consider a family of $zw$-discrete ensembles in which the data $N, A_1,A_2,B_1,B_2$ depend on $\N$, we use the notations
\begin{equation}
\label{eq_x134} \hat{n} = \frac{N}{\N},\qquad \hat{A}_1 = \frac{A_1}{\N},\qquad \hat{A}_2 = \frac{A_2}{\N}, \qquad \hat{B}_1 = \frac{B_1}{\N},\qquad \hat{B}_2 = \frac{B_2}{\N}.
\end{equation}
With the notations of Chapter~\ref{Chapter_Setup_and_Examples}, the defining segment is $[\hat{a},\hat{b}] = [\hat{A}_1,\hat{B}_1]$. We also need shifted parameters as in Definition~\ref{def:eq_shifted_parameters}
\[
\hat A'_k=\hat A_k - \frac{\theta-\frac{1}{2}}{\N},\quad \hat B'_k=\hat B_k + \frac{\theta-\frac{1}{2}}{\N} \qquad k \in \{1,2\}.
\]
We consider two different kinds of assumptions on the parameters and the precise definition of the ensemble depends on the kind we consider. The two kinds will have in common that
\[
[A_1,B_1]\subseteq[\min(A_2,B_2),\max(A_2,B_2)]
\]
but they will differ by the inequality $A_2\leq B_2$ or $A_2\geq B_2$ that is imposed. Both kinds will be useful for the tuning in Section~\ref{sec:Tuning}.

\subsection{First kind}
\label{S:sec_first_type}
If $\hat{A}_2 \leq \hat{A}_1$ and $\hat{B}_1 \leq \hat{B}_2$, we consider the probability measure
\begin{equation}
\label{zwdiscretemd}
\amsmathbb{P}_{\N}(\boldsymbol{\ell}) = \frac{1}{\Z_{\N}^{\textnormal{ref}}} \cdot \prod_{1 \leq i < j \leq N} \frac{1}{\N^{2\theta}}\cdot \frac{\Gamma\big(\ell_j - \ell_i + 1\big)\cdot\Gamma\big(\ell_j - \ell_i + \theta\big)}{\Gamma\big(\ell_j - \ell_i\big)\cdot\Gamma\big(\ell_j - \ell_i + 1 - \theta\big)} \cdot\prod_{i = 1}^N w(\ell_i)
\end{equation}
 on the state space \eqref{eq_ZW_model_state_space} with weight
\begin{equation}
 \label{weightzw} w(x) = \frac{ 2\pi\cdot\N^{B_1+B_2 - A_1-A_2 + 2}}{\Gamma\big(B_1 + 1 - x\big)\cdot\Gamma\big(B_2 + 1 - x\big)\cdot\Gamma\big(x + 1 - A_1\big)\cdot \Gamma\big(x + 1 - A_2\big)}.
\end{equation}
It differs from \eqref{eq_zw_1} only by a convenient normalizing factor of $\N$ and $\pi$. By comparison with Proposition~\ref{Proposition_Olsh_ZW}, the following value of the partition function makes \eqref{zwdiscretemd} a probability measure.
\begin{equation}
\label{ZNfirstty}\begin{split}
\Z_{\N}^{\textnormal{ref}} & = (2\pi)^{N}\cdot \N^{-\theta N^2 + (B_1 + B_2 - A_1 - A_2 + 2 + \theta)N}\cdot \Z_{N}^{zw} \\
& = (2\pi)^{N}\cdot \N^{-\theta N^2 + (B_1 + B_2 - A_1 - A_2 + 2 + \theta)N} \cdot \prod_{1 \leq i < j \leq N} \frac{\Gamma\big(\theta(j - i) + 1\big)\cdot\Gamma\big(\theta(j -i) + \theta\big)}{\Gamma\big(\theta(j - i)\big)\cdot\Gamma\big(\theta(j - i) + 1 - \theta\big)} \\
& \quad \times \prod_{i = 1}^N \frac{\Gamma\big(B_1 + B_2 - A_1 - A_2 - (2N - i - 1)\theta + 1\big)}{\Gamma\big((i - 1)\theta + 1\big) \cdot \prod_{k,\ell \in \{1,2\}} \Gamma\big(B_k - A_\ell - (i - 1)\theta + 1\big)},
\end{split}
\end{equation}
where we have used the change of index $i \rightarrow N + 1 - i$ in the factors corresponding to $k,\ell \in \{1,2\}$.

\subsection{Second kind}
\label{S:sec_second_type}
We sometimes need to consider the measure \eqref{zwdiscretemd} for another range of parameters, namely
\[
\hat{B}_2 \leq \hat{A}_1 \qquad \textnormal{and} \qquad \hat{B}_1 \leq \hat{A}_2.
\]
If we directly plug these numbers into \eqref{zwdiscretemd} and \eqref{ZNfirstty}, then the formulae might become singular, due to the poles of Gamma functions at negative integers. Yet, in the end the singularities in numerators and denominators compensate each other and the $zw$-discrete ensemble still makes sense. In order to see that, we need to rearrange the factors in the formulae.

\begin{proposition}
\label{Prop_62}
Suppose that $B_2 \leq {A}_1\leq {B}_1 \leq {A}_2$ and consider the measure \eqref{zwdiscretemd} on the state space \eqref{eq_ZW_model_state_space} with weight
\begin{equation}
\label{weightzwminus} w(x) =\frac{\N^{B_1 + B_2 - A_1 - A_2+2}}{2\pi} \cdot \frac{\Gamma\big(x - B_2\big)\cdot \Gamma\big(A_2 - x\big)}{\Gamma\big(B_1 + 1 - x\big)\cdot \Gamma\big(x + 1 - A_1\big)}.
\end{equation}
Then the following value of the partition function makes \eqref{zwdiscretemd} with weight \eqref{weightzwminus} a probability measure.
\begin{equation}
\label{ZNsecondty}\begin{split}
\Z_{\N}^{\textnormal{ref}}
& = \frac{\N^{-\theta N^2 + (B_1 + B_2 - A_1 - A_2 +2 + \theta)N}}{(2\pi)^N} \cdot \prod_{1 \leq i < j \leq N} \frac{\Gamma\big(\theta(j - i) + 1\big) \cdot\Gamma\big(\theta(j -i) + \theta\big)}{\Gamma\big(\theta(j - i)\big)\cdot\Gamma\big(\theta(j - i) + 1 - \theta\big)} \\
& \times \prod_{i = 1}^N \frac{\Gamma\big(A_2 - B_2 + (i - 1)\theta\big)\cdot\Gamma\big(A_2 - B_1 + (i - 1)\theta\big)\cdot\Gamma\big(A_1 - B_2 + (i - 1)\theta\big)}{\Gamma\big((i-1)\theta +1\big)\cdot\Gamma\big(A_1 + A_2 - B_1 - B_2 + (2N - i - 1)\theta\big)\cdot\Gamma\big(B_1 - A_1 - (i - 1)\theta + 1\big)}.
\end{split}
\end{equation}
\end{proposition}
\begin{proof}
Due to ${B}_2 \leq {A}_1\leq {B}_1 \leq {A}_2$, the weight \eqref{weightzwminus} is nonnegative at any $\ell_i$ as in \eqref{eq_zw_state_space}. Let us show that the weight only differs from \eqref{zwdiscretemd} by a prefactor independent of the configuration. We recall that \mbox{$\ell_i \in A_1 + (i - 1)\theta + \amsmathbb{Z}_{\geq 0}$} to rewrite
\begin{equation}
\label{thesimplm}\begin{split}
& \quad \frac{\Gamma\big(B_2 - A_1 - (i - 1)\theta + 1\big)}{\Gamma\big(B_2 + 1 - \ell_i\big)}\cdot\frac{\Gamma\big(A_1 - A_2 + (i - 1)\theta + 1\big)}{\Gamma\big(\ell_i + 1 - A_2\big)} \\
& = \frac{(B_2 + 1 - \ell_i)(B_2 + 2 - \ell_i)\cdots (B_2 - A_1 - (i - 1)\theta)}{(\ell_i - A_2)(\ell_i - 1 + A_2)\cdots (A_1 - A_2 + (i - 1)\theta + 1)} \\
& = \frac{(\ell_i - B_2 - 1)(\ell_i - B_2 - 2)\cdots (A_1 - B_2 + (i - 1)\theta)}{(A_2 - \ell_i)(A_2 - \ell_i + 1)\cdots (A_2 - A_1 - (i - 1)\theta - 1)} \\
& = \frac{\Gamma\big(\ell_i - B_2\big)}{\Gamma\big(A_1 - B_2 + (i - 1)\theta\big)}\cdot\frac{\Gamma\big(A_2 - \ell_i\big)}{\Gamma\big(A_2 - A_1 - (i - 1)\theta\big)}
\end{split}
\end{equation}
and deduce that the value of \eqref{weightzwminus} for $x = \ell_i$ as above is equal to the value of \eqref{weightzw} for $x = \ell_i$, multiplied by $\frac{\Z^{(i)}}{4\pi^2}$ with
\begin{equation*}
\begin{split}
\Z^{(i)} & = \Gamma\big(B_2 - A_1 - (i - 1)\theta + 1\big)\cdot \Gamma\big(A_1 - A_2 + (i - 1)\theta + 1\big) \\
& \quad \times \Gamma\big(A_1 - B_2 + (i - 1)\theta\big)\cdot \Gamma\big(A_2 - A_1 - (i - 1)\theta\big).
\end{split}
\end{equation*}

Next, we apply Proposition~\ref{Proposition_Olsh_ZW}. Although the condition $\,\text{Re}(\mathsf{z}_1 + \mathsf{z}_2 + \mathsf{w}_1 + \mathsf{w}_2 + 1) > 0$ of Proposition~\ref{Proposition_Olsh_ZW} is no longer valid for the measures of this section, the conclusion of the proposition remains valid, as can be seen by an analytic continuation argument for the values of $\mathsf{z}_1$, $\mathsf{z}_2$, $\mathsf{w}_1$, $\mathsf{w}_2$ which are allowed to be complex numbers but should be chosen so that all the arguments of Gamma functions avoid being negative integers. Note that we deal with finite sums. In contrast, the most complicated case in \cite{Olsh_hyper} is when the state space is countable and the condition $\,\text{Re}(\mathsf{z}_1 + \mathsf{z}_2 + \mathsf{w}_1 + \mathsf{w}_2 + 1) > 0$ is necessary for convergence of the infinite sums; for us, analytic continuation does not need an additional justification. We conclude that
\begin{equation}
\label{eq_x3}
 \Z_{\N}^{\textnormal{ref}} = \frac{\N^{-\theta N^2 + (B_1 + B_2 - A_1 - A_2 +2 + \theta)N}}{(2\pi)^N} \cdot \Z_{N}^{zw} \cdot \prod_{i = 1}^{N} \Z^{(i)}
\end{equation}
and it remains to transform the last formula into \eqref{ZNsecondty}. We will assume in the following computation that $A_2$, $B_2$ are generic, so that no arguments of Gamma functions are negative integers at the intermediate steps; by an analytic continuation argument, the final formula \eqref{ZNsecondty} is valid without this restriction. The formula \eqref{eq_x3} is
\begin{equation*}
\begin{split}
& \frac{\N^{-\theta N^2 + (B_1 + B_2 - A_1 - A_2 +2 + \theta)N}}{(2\pi)^N}\cdot \prod_{1\leq i<j\leq N} \frac{\Gamma\big(\theta(j-i)+1\big)\cdot \Gamma\big(\theta(j-i)+\theta\big) }{\Gamma\big(\theta(j-i)\big)\cdot
 \Gamma\big(\theta(j-i)+1-\theta\big)} \\
& \times \prod_{i=1}^N \tfrac{\Gamma(B_1+B_2-A_1-A_2-(2N-i-1)\theta+1)\cdot \Gamma(A_1 - A_2 + (i - 1)\theta + 1)\cdot \Gamma(A_1 - B_2 + (i - 1)\theta)\cdot \Gamma(A_2 - A_1 - (i - 1)\theta) }{\Gamma((i-1)\theta+1)\cdot \Gamma(B_1 -A_1 -(i-1)\theta+1)\cdot \Gamma(B_1 -A_2 -(i-1)\theta+1) \cdot \Gamma(B_2 -A_2 -(i-1)\theta+1)}.
\end{split}
\end{equation*}
We further transform the last formula using the reflection equation for the Gamma function:
\[
\Gamma(x)\cdot\Gamma(1-x)=\frac{\pi}{\sin(\pi x)}.
\]
As a result, we get \eqref{ZNsecondty} multiplied by
\begin{equation}
\label{eq_x4}
 \prod_{i=1}^N \frac{\sin\big(\pi(B_1 -A_2 -(i-1)\theta)\big)\cdot\sin\big(\pi(B_2 -A_2 -(i-1)\theta)\big)}{\sin\big(\pi(B_1+B_2-A_1-A_2-(2N-i-1)\theta)\big)\cdot\sin\big(\pi(A_1 - A_2 + (i - 1)\theta)\big)}.
\end{equation}
We claim that \eqref{eq_x4} equals $1$. In order to see that, we recall that $B_1 - A_1 - (N - 1)\theta \in \amsmathbb{Z}_{\geq 0}$ and use $\sin(x+k\pi)=(-1)^k\sin(x)$ for $k\in\amsmathbb Z$. Hence, \eqref{eq_x4} can be transformed into
\begin{equation*}
 \prod_{i=1}^N \frac{\sin\big(\pi(A_1 -A_2 -(i-N)\theta)\big)\cdot \sin\big(\pi(B_2 -A_2 -(i-1)\theta)\big) \,  }{\sin\big(\pi(B_2-A_2-(N-i)\theta)\big)\cdot\sin\big(\pi(A_1 - A_2 + (i - 1)\theta)\big) }=1. \qedhere
\end{equation*}
\end{proof}

\subsection{Remaining data of the discrete ensemble}

We now check that the $zw$-discrete ensembles fit in the framework of Chapter~\ref{Chapter_Setup_and_Examples}; for that we need to introduce all the remaining data appearing in various assumptions and then check the validity of the assumptions.

\begin{definition}\label{Definition_ZW_set}
We denote $\mathcal{D}^+$ the set of $(\theta,\hat{n},\hat{A}_1,\hat{A}_2,\hat{B}_1,\hat{B}_2) \in \amsmathbb{R}_{> 0}^2 \times \amsmathbb{R}^4$ cut out by the inequalities
\begin{equation}
\label{fnunf1}
\hat{B}_1 - \hat{A}_1 > \theta\hat{n} \quad \textnormal{and} \quad \hat{A}_1 > \hat{A}_2 \quad \textnormal{and}\quad \hat{B}_2 > \hat{B}_1.
\end{equation}
Likewise, we denote $\mathcal{D}^-$ the subset of $\amsmathbb{R}_{> 0}^2 \times \amsmathbb{R}^4$ cut out by the inequalities
\begin{equation}
\label{fnunf2}
\hat{B}_1 - \hat{A}_1 > \theta\hat{n} \quad \textnormal{and}\quad \hat{A}_2 > \hat{B}_1 \quad \textnormal{and}\quad \hat{A}_1>\hat{B}_2.
\end{equation}
We set $\mathcal{D} = \mathcal{D}^+ \cup \mathcal{D}^-$.
\end{definition}
The first inequality in \eqref{fnunf1} and \eqref{fnunf2} means that the number of particles $N$ we allow is macroscopically smaller than the maximal number of sites in $[A_1,B_1]$ that can be occupied by a particle. We introduce the following functions on $\mathcal{D}$
\begin{equation}
\label{Ddefzw}
\hat{D}_{k,\ell} = \hat{B}_k - \hat{A}_{\ell} - \theta\hat{n},\qquad \hat{D} = \hat{B}_1 + \hat{B}_2 - \hat{A}_1 - \hat{A}_2 - 2\theta\hat{n}.
\end{equation}
On $\mathcal{D}^+$ (first kind), they are all positive, while on $\mathcal{D}^{-}$ (second kind) $\hat{D}_{1,1}$ is positive and $\hat{D}_{1,2}$, $\hat{D}_{2,1}$, $\hat{D}_{2,2}$, $\hat{D}$ are negative.

\begin{lemma}
\label{Checkzw} Fix $(\theta,\hat{n},\hat{A}_1,\hat{A}_2,\hat{B}_1,\hat{B}_2) \in \mathcal{D}$. Assume that absolute values of the six parameters are smaller than a constant $C>0$ and that the $\frac{1}{C}$-neighborhood of $(\theta,\hat{n},\hat{A}_1,\hat{A}_2,\hat{B}_1,\hat{B}_2)$ is contained in $\mathcal{D}$. Then Assumptions~\ref{Assumptions_Theta}, \ref{Assumptions_basic} and \ref{Assumptions_analyticity} are satisfied for the corresponding $zw$-discrete ensemble with constants depending only on this $C$.
\end{lemma}
As we will see in Corollary~\ref{Corollary_ZW_off-critical_1}, the $zw$-ensemble also satisfies the off-criticality Assumption~\ref{Assumptions_offcrit} under an extra condition, which is generically true.

\begin{proof} We only detail the case of $\mathcal{D}^+$ as the treatment of $\mathcal{D}^-$ is similar.
 Checking all conditions of Assumption~\ref{Assumptions_Theta} is straightforward, because $H=1$ and the single matrix element of $\boldsymbol{\Theta}$ is $\theta$. The first three conditions of Assumption~\ref{Assumptions_basic} are also clear. The last condition of Assumption~\ref{Assumptions_basic} is empty, since the segment $[\hat A_1,\hat B_1]$ is uniformly bounded. Let us specify the parameters in Assumption~\ref{Assumptions_analyticity}. We omit the lower indices, since $H=1$; for instance, we write $\iota^- = \iota_1^-$. We have:
\[
 \iota^+=\iota^-=1 \qquad \textnormal{and} \qquad \rho^-=\rho^+=1,
\]
and the regular part of the potential is
\begin{equation}
\label{eq_ZW_U}
 U(z)=\frac{1}{\N} \Big(\log\Gamma\big(\N \hat B_2- \N \hat z+1\big)+\log\Gamma\big(\N \hat z-\N \hat A_2+1\big)-(\N \hat B_2- \N\hat A_2+1)\log \N-\log(2\pi)\Big).
\end{equation}
It is immediate from the definition that all three conditions of Assumption~\ref{Assumptions_analyticity} are satisfied. We also compute
\[
\frac{w\big(\N x + \frac{1}{2}\big)}{w\big(\N x - \frac{1}{2}\big)} = \frac{\big(\hat{B}_1 - x + \frac{1}{2\N}\big)\big(\hat{B}_2 - x + \tfrac{1}{2\N}\big)}{\big(x - \hat{A}_1 + \tfrac{1}{2\N}\big)\big(x - \hat{A}_2 + \tfrac{1}{2\N}\big)}.
\]
In fact, this computation remains valid for both the first and the second kind of $zw$-discrete ensembles. Hence, according to Definition~\ref{Definition_phi_functions}, the functions $\Phi^{\pm}(z)$ are specified by the formulae
\begin{equation}
\label{phipluszw}
\begin{split}
\Phi^{+}(z) & = \frac{\big(\hat{B}_1 - z + \tfrac{1}{2\N}\big)\big(\hat{B}_2 - z + \tfrac{1}{2\N}\big)}{z - \hat{A}_2 + \tfrac{1}{2\N}}, \\
\Phi^{-}(z) & = \big(z - \hat{A}_1 + \tfrac{1}{2\N}\big).
\end{split}
\end{equation}
The only difference between the two kinds is that $\Phi^+(z)$ and $\Phi^{-}(z)$ are positive for $x \in [\hat{A}_1,\hat{B}_1]$ in the first kind, and negative for the second kind. This sign will be irrelevant in our analysis. It remains to check the conditions 4.\ and 5.\ on $w(x)$ in Assumption~\ref{Assumptions_basic}; they are essentially implied by the validity of Assumption~\ref{Assumptions_analyticity}. In greater details, using \eqref{eq_potential_weight_match} we have:
\begin{equation}
\label{eq_ZW_V}
 V(x)=U(x) + \mathrm{Llog}(x-\hat A'_1) + \mathrm{Llog}(\hat{B}'_1 - x).
\end{equation}
Hence, the function $\err(x)$ of Assumption~\ref{Assumptions_basic} is given by
\begin{equation*}
\begin{split}
 \err(x) & =- \log\Gamma\big(x-\N \hat A_1+1\big)- \log\Gamma\big(\N \hat B_1-x+1\big)+(\N \hat B_1-\N \hat A_1+1)\log \N \\
 & \quad + \N\,\mathrm{Llog}\bigg(\frac{x}{\N} - \hat{A}_1'\bigg) + \N\,\mathrm{Llog}\bigg(\hat{B}_1' - \frac{x}{\N}\bigg)
 \end{split}
\end{equation*}
We apply the Binet integral formula --- see, \textit{e.g.}, \cite[Section 1.9]{Erd}:
\begin{equation}
\label{eq_Binet_formula}
\forall z > 0\qquad \log \Gamma(z)=  \mathrm{Llog}(z) + \frac{1}{2}\log\bigg( \frac{2\pi}{z}\bigg) + \int_0^{+\infty} \frac{2\arctan\big(\frac{t}{z}\big)}{e^{2\pi t}-1} \dd t.
\end{equation}
This formula is first used to justify
\[
 \log\Gamma(z) \mathop{=}_{z \rightarrow +\infty} \mathrm{Llog}(z) + O(\log z),
\]
and this results for $x\in [\N \max(\hat A_1,\hat A_1')+1,\N \min(\hat B_1,\hat B_1')-1]$ in the estimate
\begin{equation*}
\begin{split}
 \err(x) & = - \mathrm{Llog}(x - \N\hat{A}_1 + 1) - \mathrm{Llog}(\N \hat{B}_1 -x  + 1) + (\N \hat B_1-\N \hat A_1+1)\log \N \\ & \quad + \N\,\mathrm{Llog}\bigg(\frac{x}{\N} - \hat{A}_1'\bigg) + \N\,\mathrm{Llog}\bigg(\hat{B}_1' - \frac{x}{\N}\bigg) + O(\log \N)
\end{split}
\end{equation*}
with uniform error as $\N \rightarrow \infty$. Simplifying with help of $\mathrm{Llog}(\N \xi) = \xi \N \log \N + \N\,\mathrm{Llog}(\xi)$ we conclude that $\err(x)=O(\log \N)$, as desired for Assumption~\ref{Assumptions_basic}.
It remains to show the piecewise monotonicity of $\err(x)$.
For that we differentiate \eqref{eq_Binet_formula} twice and get
\begin{equation}
\label{eq_Binet_formula_derivative}
 \partial_z^2 \log \Gamma(z)= \frac{1}{z} + \frac{1}{2z^2} +\frac{1}{z^3}\int_0^{+\infty} \frac{4t\,\dd t }{\big(e^{2\pi t}-1\big)\big(1+\frac{t^2}{z^2}\big)} \dd t
 -\frac{1}{z^5}\int_0^{+\infty} \frac{4t^3\,\dd t}{\big(e^{2\pi t}-1\big)\big(1+\frac{t^2}{z^2}\big)^2}.
\end{equation}
Let us denote the last two terms in \eqref{eq_Binet_formula_derivative} $f(z)$.
For $z\rightarrow+\infty$, we have
\begin{equation}
\label{eq_Bernoulli}
 f(z) \mathop{=}_{z \rightarrow +\infty} \frac{1}{z^3} \int_0^{+\infty} \frac{4t\,\dd t}{e^{2\pi t}-1}+ O\bigg(\frac{1}{z^5}\bigg)=\frac{1}{6z^3} + O\bigg(\frac{1}{z^5}\bigg),
\end{equation}
In particular,
\[
 \partial_z^2 \log\Gamma(z) \,\,\mathop{=}_{z \rightarrow +\infty} \,\,\frac{1}{z} + \frac{1}{2z^2} + O\bigg(\frac{1}{z^3}\bigg).
\]
Then, for $x\in [\N \max(\hat A_1,\hat A_1')+1,\N \min(\hat B_1,\hat B_1')-1]$ we have
\begin{equation*}
\begin{split}
 \partial^2_x\err(x)&= - \frac{1}{x-\N \hat A_1+1} - \frac{1}{2(x-\N \hat A_1+1)^2} - \frac{1}{\N \hat B_1-x+1} - \frac{1}{2(\N \hat B_1-x+1)^2} \\
 & \quad + \frac{1}{x-\N\hat A'_1}+\frac{1}{\N \hat B'_1-x} - f\big(x-\N \hat A_1+1\big)-f\big(\N \hat B_1-x+1\big) \\=
 & \quad (1- \theta)\bigg(\frac{1}{(x-\N \hat A_1)^2} + \frac{1}{(\N \hat B_1-x)^2}\bigg) + O\bigg(\frac{1}{(x-\N \hat A_1)^3}+\frac{1}{(\N \hat B_1-x)^3}\bigg).
\end{split}
\end{equation*}
If $\theta\neq 1$, then outside finite neighborhoods of $\N\hat A_1$ and $\N\hat B_1$, we have $(1-\theta) \partial^2_x\err(x)>0$. Therefore, $\partial_x\err(x)$ can have at most one zero, implying that $\err(x)$ is piecewise monotone. In a finite neighborhood of $\N \hat A_1$ and $\N \hat B_1$, the function $\err(x)$ is also piecewise monotone by general smoothness arguments.

For $\theta=1$ we need to be more careful. In this case we further transform the expression for $\partial^2_x\err_x$:
\begin{equation*}
\begin{split}
 \partial^2_x\err(x) & =\frac{1}{2(x-\N \hat A_1+1)(x-\N\hat A'_1)} - \frac{1}{2(x-\N \hat A_1+1)^2}- f\big(x-\N \hat A_1+1\big)\\
  &\quad +\frac{1}{2(\N \hat B_1-x+1)(\N \hat B'_1-x)}- \frac{1}{2(\N \hat B_1-x+1)^2}-f\big(\N \hat B_1-x+1\big) \\
 & = \frac{1}{4(x-\N \hat A_1+1)^2(x-\N\hat A'_1)} - f\big(x-\N \hat A_1+1\big)\\
 &\quad + \frac{1}{4(\N \hat B_1-x+1)^2(\N \hat B'_1-x)}-f\big(\N \hat B_1-x+1\big).
 \end{split}
\end{equation*}
Using \eqref{eq_Bernoulli} we conclude that the last expression is positive outside finite neighborhoods of $\N \hat A_1$ and $\N \hat B_1$ and, therefore, $\partial_x\err(x)$ can have at most one zero there, again.
\end{proof}

The off-criticality of the equilibrium measure of Assumption~\ref{Assumptions_offcrit} will be addressed in Section~\ref{Eqmeszw}.

\section{Asymptotic expansion of the partition function}
\label{AsymZWpart} We recall the definition of the \emph{Barnes double zeta function} and refer to \cite{Spreafico} for its properties. For $\text{Re}(s) > 2$ we define
\[
\chi_2(s|b_1,b_2) = \sum_{\substack{m_1,m_2 \geq 0 \\ (m_1,m_2) \neq (0,0)}} \frac{1}{(b_1m_1 + b_2m_2)^{s}}.
\]
It is extended to a meromorphic function of $s\in\amsmathbb C$ by an analytic continuation. For instance,
\[
\chi'_2(0|1,1) = \zeta'(-1)\qquad \textnormal{and}\qquad \chi'_2(0|2,1) = \frac{\zeta'(-1)}{2} + \frac{5 \log 2}{24}
\]
in terms of the Riemann zeta function. The result for the asymptotic expansion of the partition function is best expressed using the shifted variables
\begin{equation}
\label{Dpar}
\begin{split}
\hat{D}' & = \hat{B}_1' + \hat{B}_2' - \hat{A}_1' - \hat{A}_2' - 2\theta\hat{n} = \hat{D} + \frac{4\theta - 2}{\N}  \\
\hat{D}_{k,\ell}' & = \hat{B}'_k - \hat{A}_\ell' - \theta\hat{n} = \hat{D}_{k,\ell} + \frac{2\theta -1}{\N},
\end{split}
\end{equation}
because they will appear in the equilibrium measure.

\begin{proposition}
\label{Proposition_ZW_partition} In the setting of Lemma~\ref{Checkzw}, uniformly in any compact subset of $\mathcal{D}$, we have the asymptotic expansion as $\N \rightarrow \infty$
\[
\log \Z^{\textnormal{ref}}_{\N} = \N^2 \F^{[0],\textnormal{ref}} + \theta \hat{n}\,\N \log \N + \N\F^{[1],\textnormal{ref}} + \frac{1}{12}\bigg(\theta + \frac{1}{\theta} - 3\bigg)\log \N + \F^{[2],\textnormal{ref}} + O\bigg(\frac{1}{\N}\bigg).
\]
The coefficients $(\F^{[p],\textnormal{ref}})_{p = 0}^{2}$ are given using the notation \eqref{Ddefzw} by the following expressions
\begin{equation*}
\begin{split}
\F^{[0],\textnormal{ref}} & = \frac{1}{\theta}\bigg[\mathrm{LLlog}(\theta\hat{n}) + \mathrm{LLlog}(\hat{D}' + \theta\hat{n}) - \mathrm{LLlog}(\hat{D}')  - \!\! \sum_{k,\ell = 1}^{2} \big(\mathrm{LLlog}(\hat{D}'_{k,\ell} + \theta\hat{n}) - \mathrm{LLlog}(\hat{D}_{k,\ell}')\big)\bigg] \\
\F^{[1],\textnormal{ref}} & = -  \hat{n}\log\big(2\pi \Gamma(\theta)\big) + \frac{\theta -1}{2\theta}\bigg[\mathrm{Llog}(\theta\hat{n}) - 5\big(\mathrm{Llog}(\hat{D}' + \theta\hat{n}) - \mathrm{Llog}(\hat{D}')\big) \\
& \qquad + \sum_{k,\ell = 1}^{2}  3 \big(\mathrm{Llog}(\hat{D}'_{k,\ell} + \theta\hat{n}) - \mathrm{Llog}(\hat{D}'_{k,\ell})\big)\bigg] \\
\F^{[2],\textnormal{ref}} & = \frac{37(\theta - 1)^2 - \theta}{12\theta} \log \bigg|\frac{\hat{D}' + \theta\hat{n}}{\hat{D}'}\bigg| - \sum_{k,\ell = 1}^{2} \frac{13(\theta - 1)^2 - \theta}{12\theta} \log\bigg|\frac{\hat{D}'_{k,\ell} + \theta \hat{n}}{\hat{D}'_{k,\ell}}\bigg| \\
& \quad + \frac{1}{12}\bigg(\theta + \frac{1}{\theta} - 3\bigg)\log\hat{n} + \chi_2'(0|\theta^{-1},1).
\end{split}
\end{equation*}
\end{proposition}
\begin{remark}
\label{Remark_Dplus_Dmoins}
These expressions are valid both in $\mathcal{D}^+$ and $\mathcal{D}^-$; this is the reason for having chosen differently the normalizing factors of $2\pi$ in \eqref{weightzw} and \eqref{weightzwminus}. If the parameters are in $\mathcal{D}^+$ we can omit the absolute values. For our following developments the precise expressions for $(\F^{[p],\textnormal{ref}})_{p = 0}^{2}$ are not that important, but it is crucial that they are smooth functions of all the involved parameters.
\end{remark}
\begin{proof}[Proof of Proposition~\ref{Proposition_ZW_partition}] We are going to rewrite the partition function in terms of Barnes double Gamma functions, whose asymptotics are well-known\label{index:Gamma2fn}. Following the notations of \cite{Spreafico}, this function $\Gamma_2(x|b_1,b_2)$ satisfies
\begin{equation}
\label{eq_Double_Gamma_functional}
\frac{\Gamma_2(x|b_1,b_2)}{\Gamma_2(x + b_2|b_1,b_2)}= \Gamma\bigg(\frac{x}{b_1}\bigg)\cdot \frac{b_1^{\frac{x}{b_1}-\frac{1}{2}}}{\sqrt{2\pi}}, \qquad \Gamma_2(1|b_1,b_2) = 1.
\end{equation}
It admits the asymptotic expansion as $x \rightarrow +\infty$ \cite[Proposition~8.11]{Spreafico}
\begin{equation}
\label{eq_expansion_Gamma2}
\begin{split}
\log \Gamma_2(x|b_1,b_2) & = -\frac{\mathrm{LLlog}(x)}{2b_1b_2} + \frac{1}{2}\bigg(\frac{1}{b_1} + \frac{1}{b_2}\bigg) \mathrm{Llog}(x) \\
&  \quad - \frac{1}{12}\bigg(3 + \frac{b_1}{b_2} + \frac{b_2}{b_1}\bigg)\log(x) - \chi_2'(0|b_1,b_2) + O\bigg(\frac{1}{x}\bigg),
\end{split}
\end{equation}
with uniformly small remainder for $(b_1,b_2)$ in any compact subset of $\amsmathbb{R}_{> 0}^2$. We have the function
\[
\mathrm{LLlog}(x) = \frac{x^2}{2} \log |x| - \frac{3x^2}{4},
\]
which is such that $\mathrm{LLlog}' = \mathrm{Llog}$.

If the parameters of the ensemble belong to $\mathcal{D}^+$, we can decompose the partition function \eqref{ZNfirstty} as
\begin{equation}
\label{firstXN}
\begin{split}
\Z_{\N}^{\textnormal{ref}} & = \N^{-\theta N^2 + (B_1 + B_2 - A_1 - A_2 + 2 + \theta)N}\cdot \tilde{\gimel}_N\cdot \frac{\gimel_{N}\big(B_1 + B_2 - A_1 - A_2 - 2\theta(N - 1)\big)}{\prod_{k,\ell \in \{1,2\}} \gimel_{N}\big(B_k - A_{\ell} - (N - 1)\theta\big)} \\
& = \N^{\theta \hat{n}^2 \N^2  + \hat{D}'\hat{n} \N^2  + (4 - 3\theta)\hat{n}\N} \cdot \tilde{\gimel}_{N} \cdot \frac{\gimel\big(\N \hat{D}' - 2(\theta - 1)\big)}{\prod_{k,\ell \in \{1,2\}} \gimel_N\big(\N \hat{D}_{k,\ell}' - (\theta - 1)\big)}.
\end{split}
\end{equation}
in terms of the auxiliary quantities
\begin{equation}
\label{phiexp}
\begin{split}
\gimel_{N}(x) & := \prod_{i = 0}^{N - 1} \Gamma\big(x + i\theta + 1\big) = \frac{(2\pi)^{\frac{N}{2}} \cdot \theta^{\frac{\theta N^2}{2} + (\frac{1 - \theta}{2} + x)N} \cdot \Gamma_2\big(\frac{x + 1}{\theta}\big|\frac{1}{\theta},1\big)}{\Gamma_2\big(\frac{x + 1}{\theta} + N\big|\frac{1}{\theta},1\big)}\,, \\
\tilde{\gimel}_N & := \prod_{1 \leq i < j \leq N} \frac{\Gamma\big(\theta(j - i) + 1\big)\cdot \Gamma\big(\theta(j - i) + \theta\big)}{\Gamma\big(\theta(j - i)\big)\cdot\Gamma\big(\theta(j - i) + 1 - \theta\big)}\cdot \prod_{i = 1}^{N} \frac{1}{\Gamma\big((i - 1)\theta + 1\big)}.
\end{split}
\end{equation}
We compute the telescoping product for $N \geq 1$
\begin{equation}
 \label{eq_x197}
 \begin{split}
\frac{\tilde{\gimel}_{N + 1}}{\tilde{\gimel}_{N}} & = \frac{1}{\Gamma\big(\theta N + 1\big)}\cdot \prod_{i = 1}^N\bigg[ \frac{\Gamma\big(\theta(N + 1 - i)+1\big)}{\Gamma\big(\theta(N + 1 - i)\big)}\cdot \frac{\Gamma\big(\theta(N + 2 - i)\big)}{\Gamma\big(\theta(N - i) + 1\big)}\bigg] \\
& = \frac{1}{\Gamma\big(\theta N+1\big)} \cdot \frac{\Gamma\big(\theta N+1\big)}{\Gamma(1)} \cdot\frac{\Gamma\big(\theta(N+1)\big)}{\Gamma(\theta)}= \frac{\Gamma\big(\theta(N+1)\big)}{\Gamma(\theta)}.
\end{split}
\end{equation}
Using $\tilde{\gimel}_1 = 1$ and \eqref{eq_Double_Gamma_functional} we obtain
\begin{equation}
\label{eq_x199}
\tilde{\gimel}_N = \prod_{i = 1}^{N} \frac{\Gamma(\theta i)}{\Gamma(\theta)} = \frac{(2\pi)^{\frac{N}{2}} \cdot \theta^{\frac{\theta N^2}{2} + \frac{(\theta - 1)N}{2}}}{\Gamma(\theta)^N \cdot \Gamma_2\big(N + 1\big|\frac{1}{\theta},1\big)}.
\end{equation}
Recall that $N = \N \hat{n}$ and we want to send $\N \rightarrow \infty$ while keeping $\hat{n}$ positive bounded away from $0$. We deduce from \eqref{eq_expansion_Gamma2} and \eqref{eq_x199} the asymptotic expansion
\begin{equation}
\label{eq_expansion_gamma_tilde}
\begin{split}
\log \tilde{\gimel}_N & = \frac{\theta \hat{n}^2}{2}\,\N^2\log\N + \frac{\mathrm{LLlog}(\theta\hat{n})}{\theta}\,\N^2 + \frac{\theta - 1}{2}\,\hat{n}\N\log\N \\
& \quad +\left( \frac{\theta - 1}{2\theta}\,\mathrm{Llog}(\theta\hat{n}) -  \hat{n} \log \Gamma(\theta) + \frac{\hat{n}}{2}\log(2\pi)\right)\N \\
& \quad + \frac{1}{12}\bigg(\theta + \frac{1}{\theta} - 3\bigg)(\log \N + \log\hat{n}) + \chi_2'(0|\theta^{-1},1) + O\bigg(\frac{1}{\N}\bigg).
\end{split}
\end{equation}
The remainder $O\big(\frac{1}{\N}\big)$ in the last formula is uniform over $(\theta,\hat{n})$ in any compact of $\amsmathbb{R}_{> 0}^2$. We also deduce from \eqref{eq_expansion_Gamma2} and \eqref{phiexp} the asymptotic expansion
\begin{equation}
\label{gnexpmfgun}\begin{split}
&\quad \log \gimel_N(\N\xi +\xi') \\
& = \bigg(\xi + \frac{\theta\hat{n}}{2}\bigg)\hat{n}\,\N^2\log\N + \frac{\mathrm{LLlog}(\xi + \theta\hat{n}) - \mathrm{LLlog}(\xi)}{\theta}\,\N^2 \\
& \quad + \bigg(\xi' + \frac{1 - \theta}{2}\bigg)\hat{n}\N\log \N + \left[\bigg(\xi' + \frac{1 - \theta}{2}\bigg)\frac{\mathrm{Llog}(\xi + \hat{n}\theta) - \mathrm{Llog}(\xi)}{\theta} + \frac{\hat{n}}{2}\log(2\pi) \right]\N \\
& \quad + \frac{1}{2\theta}\left[\bigg(\xi' + \frac{1 - \theta}{2}\bigg)^2 - \frac{\theta^2 + 1}{12}\right]\log\bigg(\frac{\xi + \theta\hat{n}}{\xi}\bigg)+ O\bigg(\frac{1}{\N}\bigg),
\end{split}
\end{equation}
which is uniform for $(\xi,\xi',\theta,\hat{n}) \in \amsmathbb{R}^4$ for any compact in the region $\xi,\hat{n},\theta > 0$. We should use \eqref{gnexpmfgun} five times with the $\N$-dependent variables
\[
(\xi,\xi') = \big(\hat{D}',-2(\theta -1)\big) \qquad \textnormal{and} \qquad (\xi,\xi') = \big(\hat{D}'_{k,\ell},-(\theta - 1)\big)\quad \textnormal{with} \,\, k,\ell \in \{1,2\}.
\]
We then have to collect the terms by order of magnitude in $\N$ and this leads to the claimed formulae in Proposition~\ref{Proposition_ZW_partition} without the absolute values in the logarithms (they are irrelevant in the $\mathcal{D}^+$ case). In particular, due to the power of $\N$ in prefactor of the partition function, we observe that the $\N^2\log \N$ cancel each other and the $\N\log\N$ terms appear with coefficient $\theta\hat{n}$.

If the parameters are in $\mathcal{D}^-$, the partition function \eqref{ZNsecondty} decomposes as
\begin{equation*}
\begin{split}
Z_{\N}^{\textnormal{ref}} & = \N^{-\theta N^2 + (B_1 + B_2 - A_1 - A_2 + 2 + \theta)N}\cdot \tilde{\gimel}_{N}
\\ & \quad \times \frac{\gimel_{N}\big(A_2 - B_2 - 1\big) \cdot \gimel_{N}\big(A_2 - B_1 - 1\big) \cdot \gimel_{N}\big(A_1 - B_2 - 1\big)}{\gimel_N\big(A_1 + A_2 - B_1 - B_2 + (N - 1)\theta - 1\big) \cdot \gimel_N\big(B_1 - A_1 - (N - 1)\theta\big)} \\
& = \N^{\theta\hat{n}^2\N^2 + \hat{D}' \hat{n} \N^2 + (4 - 3\theta)\hat{n}\N} \cdot \tilde{\gimel}_{N} \cdot \gimel_{N}\big(-\N(\hat{D}_{2,2}' + \theta\hat{n}) + 2(\theta - 1) + 1\big) \\
& \quad \times \frac{\gimel_N\big(-\N(\hat{D}_{1,2}' + \theta\hat{n}) + 2(\theta -1)+1\big) \cdot \gimel_N\big(-\N(\hat{D}_{2,1}' + \theta\hat{n}) + 2(\theta - 1) + 1\big)}{\gimel_{N}\big(-\N(\hat{D}' + \theta \hat{n}) + 3(\theta - 1)\big) \cdot \gimel_N\big(\N \hat{D}_{1,1}' - (\theta -1)\big)}.
\end{split}
\end{equation*}
The asymptotic expansion is obtained by a similar algebra and the result takes the same form, this time with the absolute values taking into account the fact that $\hat{D}'$ and $\hat{D}'_{k,\ell}$ for $(k,\ell) \neq (1,1)$ are negative (at least for $\N$ large enough).
\end{proof}

In the course of the proof we have observed while computing the quantity $\tilde{\gimel}_{N}$ the telescoping of the pairwise interaction for consecutive particles. As it will be useful later, we record the following result.

\begin{lemma} \label{Lemma_densely_packed_expansion} With the notation $\hat m= \frac{M}{\N}$, we have the asymptotic expansion as $\N \rightarrow \infty$
\begin{equation*}
\begin{split}
& \log\Bigg(\prod_{1 \leq i < j \leq M} \frac{1}{\N^{2\theta}}\cdot \frac{\Gamma\big(\theta(j - i) + 1\big)\cdot\Gamma\big(\theta(j - i) + \theta\big)}{\Gamma\big(\theta(j - i)\big)\cdot\Gamma\big(\theta(j - i) + 1 - \theta\big)}\Bigg) \\
& = \bigg(\frac{1}{\theta} \int_{0}^{\theta\hat{m}} \int_{0}^{\theta\hat{m}} \log|x - y|\,\dd x\,\dd y\bigg)\N^2 +\theta \hat m \N \log \N \\
& \quad + \hat{m}\log\bigg(\frac{2\pi}{\Gamma(\theta)}\bigg) \N + \frac{1}{6}\bigg(\theta + \frac{1}{\theta} - 3\bigg)\log \N \\
& \quad + \frac{1}{6}\bigg(\theta + \frac{1}{\theta} - 3\bigg)\log \hat{m} + \frac{\log(2\pi)}{2} + \chi_2'(0|\tfrac{1}{\theta},1) + O\bigg(\frac{1}{\N}\bigg),
\end{split}
\end{equation*}
which is uniform for $\hat{m}$ in any compact of $\amsmathbb{R}_{> 0}$.
\end{lemma}
\begin{proof}
By comparison with \eqref{phiexp} and \eqref{eq_x199}, the quantity inside the logarithm is equal to
\[
 \N^{-\theta M (M-1)}\cdot \tilde \gimel_M \cdot \gimel_M(0) = \frac{\N^{-\theta M(M - 1)} \cdot (2\pi)^{M} \cdot \theta^{\,\theta M^2} \cdot \Gamma_2\big(\frac{1}{\theta}\big|\frac{1}{\theta},1)}{\Gamma(\theta)^{M}\cdot \Gamma_2\big(M + 1\big|\frac{1}{\theta},1) \cdot \Gamma_2\big(M + \frac{1}{\theta}\big|\frac{1}{\theta},1\big)} .
\]
Using the asymptotics \eqref{eq_expansion_Gamma2} for $\Gamma_2$ and the special value $\Gamma_2\big(\frac{1}{\theta}\big|\frac{1}{\theta},1\big) = \sqrt{2\pi}$ from \cite[Proposition~8.12]{Spreafico} we obtain the desired result except that the leading term appears as $\frac{2}{\theta}\,\mathrm{LLlog}(\hat{m}\theta)$. A direct computation shows that this coincides with the double integral announced above.
\end{proof}

\section{Simplified equilibrium measure}

\label{Eqmeszw}

We would like to study the equilibrium measure $\mu$ (\textit{cf.} Proposition~\ref{Lemma_maximizer}) corresponding to the $zw$-discrete ensemble of Section~\ref{SChoice}. The measure itself is quite complicated because of the delicate definition of $V(x)$ in \eqref{eq_ZW_U}, \eqref{eq_ZW_V} involving logarithms of the Gamma function. Therefore, we would like to first introduce and investigate the equilibrium measure $\widetilde \mu$ for a simplified potential obtained by replacing Gamma functions with their Stirling approximations; the latter measure is explicit and its computation can be traced back to \cite{CLP} and \cite[Section 9.4]{BGG}. Ultimately, we can deduce all the needed properties of $\mu$ from those of $\widetilde \mu$ by appealing to the parametric smoothness results of Chapter~\ref{Chapter_smoothness}.

\medskip

Let us specify the variational datum (\textit{cf.} Section~\ref{Section_variational_data}) for the simplified equilibrium measure $\widetilde \mu$. We set $H=1$. All indices $h$ can be then dropped. Next,
\begin{equation}
\label{eq_ZW_data_1}
 \theta_{11}=\theta,\qquad \hat n=\frac{N}{\N},\qquad \hat a'= \hat A'_1= \frac{A_1}{\N}-\frac{\theta-\frac{1}{2}}{\N}, \qquad \hat b'=\hat B'_1=\frac{B_1}{\N}+\frac{\theta-\frac{1}{2}}{\N}.
\end{equation}
For the sake of the symmetry we also define
\begin{equation}
\label{eq_ZW_data_11}
 \hat A'_2= \frac{A_2}{\N}-\frac{\theta-\frac{1}{2}}{\N}, \qquad\hat B'_2=\frac{B_2}{\N}+\frac{\theta-\frac{1}{2}}{\N}.
\end{equation}
With these notations, the potential is defined by
\begin{equation}
\label{eq_ZW_data_2}
 \widetilde V(x)= \mathrm{Llog}(x-\hat A'_1)+ \mathrm{Llog}(\hat B'_1-x) + \mathrm{Llog}(x-\hat A'_2) + \mathrm{Llog}(\hat B'_2-x).
\end{equation}
Let $\widetilde \mu$ denote the equilibrium measure corresponding to the above variational datum and let $\Gm_{\widetilde \mu}(z)$ be the corresponding Stieltjes transform
\[
\Gm_{\widetilde \mu}(z) = \int_{\hat{a}'}^{\hat{b}'} \frac{\widetilde \mu(x) \dd x}{z - x}.
\]
The result will be formulated with the shifted parameters $\hat{D}_{k,\ell}',\hat{D}'$ of \eqref{Dpar}.
\begin{theorem} \label{Proposition_ZW_equilibrium_measure}
 Consider the variational datum \eqref{eq_ZW_data_1}-\eqref{eq_ZW_data_2} and suppose that the sextuple
 \[
 (\theta, \hat n, \hat A_1', \hat A_2', \hat B_1', \hat B_2')
 \]
 belongs to the domain $\mathcal D$ of Definition~\ref{Definition_ZW_set}. Then the simplified equilibrium measure $\widetilde \mu$ has a single band $(\alpha,\beta)$ and its density inside the band is given by
 \begin{equation}
 \label{eq_ZW_density_final}
 \forall x \in (\alpha,\beta)\qquad \widetilde \mu(x) =\frac{1}{\pi\theta}\,\arctan\Bigg(\hat D' \frac{\sqrt{(x - \alpha)(\beta - x)}}{q^+(x)}\Bigg),
 \end{equation}
 where $\arctan(z)$ takes values between $0$ and $\pi$ for $z\in [-\infty,+\infty]$ with $\arctan(\pm \infty)=\tfrac{\pi}{2}$,
 \begin{equation}
\label{eq_x130}
\begin{split} q^+(z) & = 2z^2 - (\hat{A}'_1 + \hat{A}'_2 + \hat{B}'_1 + \hat{B}'_2)z + \hat{B}'_1\hat{B}'_2 + \hat{A}'_1\hat{A}'_2 - \theta\hat{n}(\theta\hat{n} + \hat{D}') \\
\alpha & = \frac{1}{(\hat{D}')^2}\Big(-(\hat{A}'_1 + \hat{A}'_2 + \hat{B}'_1 + \hat{B}'_2)(\hat{D}' + \theta\hat{n})\theta\hat{n} + (\hat{B}'_1\hat{B}'_2 - \hat{A}'_1\hat{A}'_2)(\hat{D}' + 2\theta\hat{n})- \sqrt{\Delta} \Big), \\
 \beta & = \frac{1}{(\hat{D}')^2}\Big(-(\hat{A}'_1 + \hat{A}'_2 + \hat{B}'_1 + \hat{B}'_2)(\hat{D}' + \theta\hat{n})\theta\hat{n} + (\hat{B}'_1\hat{B}'_2 - \hat{A}'_1\hat{A}'_2)(\hat{D}' + 2\theta\hat{n}) + \sqrt{\Delta}\Big), \\
 \Delta & = 4\theta\hat{n}(\hat{D}' + \theta\hat{n})\hat{D}'_{1,1}\hat{D}'_{1,2}\hat{D}'_{2,1}\hat{D}'_{2,2} > 0.
\end{split}
\end{equation}
Besides, if the parameters are in $\mathcal{D}^{\tau}$ with $\tau \in \{\pm 1\}$, then
\begin{itemize}
\item[$\bullet$] if $\tau\,q^+(\hat{A}_1') > 0$, then $[\hat{A}'_1,\alpha]$ is a void;
\item[$\bullet$] if $\tau\,q^+(\hat{A}_1') < 0$, then $[\hat{A}'_1,\alpha]$ is a saturation;
\item[$\bullet$] if $\tau\,q^+(\hat{B}_1') > 0$, then $[\beta,\hat{B}'_1]$ is a void;
\item[$\bullet$] if $\tau\,q^+(\hat{B}_1') < 0$, then $[\beta,\hat{B}'_1]$ is a saturation.
\end{itemize}
\end{theorem}
We do not impose any integrality conditions on $(\theta, \hat n, \hat A_1', \hat A_2', \hat B_1', \hat B_2')$, they are allowed to be arbitrary sextuples of reals in the set $\mathcal D$. The formulae \eqref{eq_x130} imply that the length of the band is
 \begin{equation}
\label{width}
\beta - \alpha = \frac{4\sqrt{\theta\hat{n}(\hat{D'} + \theta\hat{n})\hat{D}'_{1,1}\hat{D}'_{1,2}\hat{D}'_{2,1}\hat{D}'_{2,2}}}{(\hat{D}')^2}.
\end{equation}

\begin{proof}[Proof of Theorem~\ref{Proposition_ZW_equilibrium_measure}] The auxiliary functions of Definition~\ref{GQdef} associated with $\widetilde{\mu}$ are denoted $\widetilde{q}^{\pm}(z)$.

\medskip

\noindent \textsc{Step 1.} We first compute $\Gm_{\widetilde \mu}(z)$. Let us define
\begin{equation}
\label{eq_qpluszw} q^+(z) = (\hat{B}'_1 - z)(\hat{B}'_2 - z)\cdot \exp\big(\theta \,\Gm_{\widetilde \mu}(z)\big) + (z - \hat{A}'_1)(z - \hat{A}'_2)\cdot \exp\big(-\theta\,\Gm_{\widetilde \mu}(z)\big).
\end{equation}
This is not exactly the function $\widetilde{q}^+(z)$ that Definition~\ref{GQdef} associates to $\widetilde{\mu}$. Indeed, we have
\begin{equation}
\label{eq_x2}
e^{- \partial_z \widetilde V'(z)} =\frac{(\hat B'_1-z)(\hat B'_2-z)}{(z-\hat A'_1)(z-\hat A'_2)}= \frac{\left(\frac{(\hat B'_1-z)(\hat B'_2-z)}{(z-\hat A'_2)}\right)} {(z-\hat A'_1)}.
\end{equation}
and comparing with Definition~\ref{Definition_phi_functions_2}, $\widetilde{\phi}^+(z)$ should be the numerator and $\widetilde{\phi}^-(z) = (z - \hat{A}'_1)$ the denominator of this last expression. We conclude that $q^+(z)$ from \eqref{eq_qpluszw} differs from $\widetilde{q}^+(z)$ by a multiplication with the linear factor $(z - \hat{A}'_2)$, whose only zero $\hat A'_2$ is outside $[\hat A'_1,\hat B'_1]$. Applying Theorem~\ref{Theorem_regularity_density}, we conclude that $q^+(z)$ is holomorphic for $z$ in a complex neighborhood of $[\hat A'_1,\hat B'_1]$. It is manifest in \eqref{eq_qpluszw} that $q^+(z)$ is also holomorphic for $z \in \amsmathbb{C} \setminus [\hat A'_1,\hat B'_1]$. We conclude that $q^+(z)$ is an entire function.

On the other hand, we have an asymptotic expansion
\begin{equation}
\label{eq_Wmuzwasym}\Gm_{\widetilde \mu}(z) \mathop{=}_{z \rightarrow +\infty} \frac{\hat{n}}{z} + \frac{\int_{\amsmathbb{R}} x\mu(x)\dd x}{z^2} + O\bigg(\frac{1}{z^3}\bigg),
\end{equation}
Plugging into \eqref{eq_qpluszw} we conclude that $q^+(z)=O(z^2)$ as $z\rightarrow\infty$. By Liouville theorem, this implies that $q^+(z)$ is a polynomial of degree $2$. This polynomial is identified by inserting \eqref{eq_Wmuzwasym} into \eqref{eq_qpluszw} to reach the $z \rightarrow \infty$ asymptotics up to $o(1)$ in the right-hand side of \eqref{eq_qpluszw}. We find
\begin{equation}
\label{eq_q_ZW_explicit}
q^+(z) = 2z^2 - (\hat{A}'_1 + \hat{A}'_2 + \hat{B}'_1 + \hat{B}'_2)z + \hat{B}'_1\hat{B}'_2 + \hat{A}'_1\hat{A}'_2 + \theta\hat{n}\big(\hat{A}'_1 + \hat{A}'_2 - \hat{B}'_1 - \hat{B}'_2 + \theta\hat{n}\big),
\end{equation}
which coincides with \eqref{eq_x130}. We further view the combination of \eqref{eq_qpluszw} and \eqref{eq_q_ZW_explicit} as a quadratic equation for $\exp\big(\theta\,\Gm_{\widetilde \mu}(z)\big)$. Its solution is
\begin{equation}
\label{eq_Wzw} \exp\big(\theta\, \Gm_{\widetilde \mu}(z)\big) = \frac{q^+(z) \pm \sqrt{(q^+(z))^2 - 4(\hat{B}'_1 - z)(\hat{B}'_2 - z)(z - \hat{A}'_1)(z - \hat{A}'_2)}}{2(\hat{B}'_1 - z)(\hat{B}'_2 - z)},
\end{equation}
where $q^+(z)$ is given by \eqref{eq_q_ZW_explicit}. To find the correct branch of the square root in this formula (equivalently, the correct $\pm$ sign), we analyze again the behavior as $z \rightarrow \infty$. We have
\begin{equation}
\label{eq_x132}
\begin{split}
L(z) &:= (q^+(z))^2 - 4(\hat{B}'_1 - z)(\hat{B}'_2 - z)(z - \hat{A}'_1)(z - \hat{A}'_2) \\
& = (\hat{D}')^2z^2 +2\big((\hat{A}'_1 + \hat{A}'_2 + \hat{B}'_1 + \hat{B}'_2)(\theta\hat{n} + \hat{D}')\theta\hat{n} + (\hat{A}'_1\hat{A}'_2 - \hat{B}'_1\hat{B}'_2)(2\theta\hat{n} + \hat{D}')\big)z \\
& \quad + \theta^4\hat{n}^4 + 2\hat{D}'\theta^3\hat{n}^3 + \big( (\hat{D}')^2 - 2\hat{A}'_1\hat{A}'_2 - 2\hat{B}'_1\hat{B}'_2\big)\theta^2\hat{n}^2 \\
& \quad - 2\hat{D}'(\hat{A}'_1\hat{A}'_2 + \hat{B}'_1\hat{B}'_2)\hat{n}\theta + (\hat{A}'_1\hat{A}'_2 - \hat{B}'_1\hat{B}'_2)^2 \\
& = (\hat{D}')^2z^2 + O(z).
\end{split}
\end{equation}
Let us choose the branch of the square root in \eqref{eq_q_ZW_explicit} which maps positive real numbers to positive real numbers. Then \eqref{eq_x132} implies
\begin{equation}
\label{eq_x135}
\exp\big(\theta\, \Gm_{\widetilde \mu}(z)\big) \mathop{=}_{z \rightarrow +\infty} \frac{q^+(z) \pm \hat{D}'z + O(1) }{2(\hat{B}'_1 - z)(\hat{B}'_2 - z)}.
\end{equation}
On the other hand, directly from the definitions we have
\begin{equation}
\label{eq_x136}
\exp\big(\theta\,\Gm_{\widetilde \mu}(z)\big) \mathop{=}_{z \rightarrow +\infty} 1 + \frac{\theta\hat{n}}{z} + O\bigg(\frac{1}{z^2}\bigg).
\end{equation}
Comparing \eqref{eq_x135} with \eqref{eq_x136} and the definition of $\hat D'=\hat{B}'_1 + \hat{B}'_2 - \hat{A}'_1 - \hat{A}'_2 - 2\theta\hat{n}$, we conclude that the correct sign is $-$. The final result is that
\begin{equation}
\label{eq_Wzw_final} \exp\big(\theta\, \Gm_{\widetilde \mu}(z)\big) = \frac{q^+(z) - \sqrt{(q^+(z))^2 - 4(\hat{B}'_1 - z)(\hat{B}'_2 - z)(z - \hat{A}'_1)(z - \hat{A}'_2)}}{2(\hat{B}'_1 - z)(\hat{B}'_2 - z)},
\end{equation}
with the branch of the square root mapping positive real numbers to positive real numbers.

\medskip

\noindent \textsc{Step 2.} We use \eqref{eq_Wzw_final} to identify the band of the equilibrium measure $\widetilde \mu$. We first make the preliminary remark that for $x \in \amsmathbb{A} = [\hat{A}'_1,\hat{B}'_1]$, the Stieltjes transform is related to the density of the equilibrium measure via
\begin{equation}
\label{eq_x137}
\widetilde \mu(x) = -\frac{1}{\pi}\,\textnormal{Im}\,\Gm_{\widetilde \mu}(x^+)=\frac{1}{\pi}\, \textnormal{Im}\,\Gm_{\widetilde \mu}(x^-)
\end{equation}
and we have $0 \leq \widetilde \mu(x) \leq \frac{1}{\theta}$. In particular, $x \in \amsmathbb{A}$ belongs to a void of $\mu$ if and only if $\exp\big(\theta\,\Gm_{\widetilde{\mu}}(x^+)\big)$ is real positive, and to a saturation if and only if $e^{\theta\,\Gm_{\widetilde{\mu}}(x^+)}$ is real negative. Using \eqref{eq_x137} we conclude that a point $x \in \amsmathbb A$ belongs to a band if and only if the expression in \eqref{eq_Wzw_final} is non-real, which means that
\begin{equation}
\label{eq_discrizw} L(x) := q^+(x)^2 - 4(\hat{B}'_1 - x)(\hat{B}'_2 - x)(x - \hat{A}'_1)(x - \hat{A}'_2) < 0.
\end{equation}
and $x$ is the endpoint of a band in the equality case. The expression \eqref{eq_discrizw}, as simplified in \eqref{eq_x132}, is a polynomial of degree $2$ in $x$ with discriminant
\[
 16\theta\hat{n}(\hat{D'} + \theta\hat{n})\hat{D}'_{1,1}\hat{D}'_{1,2}\hat{D}'_{2,1}\hat{D}'_{2,2} > 0.
\]
Therefore, $L(x)$ is negative if $x \in [\alpha,\beta]$ given by the formulae of \eqref{eq_x130}. Consequently, $\widetilde \mu$ has at most one band. Because $\widetilde \mu(x)$ has to be zero if $x \notin \amsmathbb A$, we also deduce that $\hat A'_1 \leq \alpha < \beta \leq \hat B'_1$. Therefore, $\widetilde \mu$ has a unique band $(\alpha,\beta)\subset \amsmathbb A$.

\medskip

\noindent \textsc{Step 3.} We now develop the formula for the density $\widetilde \mu(x)$ by using \eqref{eq_density_tan} of Proposition~\ref{Proposition_density}. We define
\begin{equation}
\label{eq_qminuszw} q^-(z) = (\hat{B}'_1 - z)(\hat{B}'_2 - z)\cdot \exp\big(\theta \,\Gm_{\widetilde \mu}(z)\big) - (z - \hat{A}'_1)(z - \hat{A}'_2)\cdot\exp\big(-\theta\,\Gm_{\widetilde \mu}(z)\big).
\end{equation}
We conclude that (similarly to $q^+(z)$) the function $q^-(z)$ differs from $q^-_1(z)$ of Definition~\ref{Definition_phi_functions_2} by the linear factor $(z - \hat{A}'_2)$. The definitions imply that $q^-(z)^2-q^{+}(z)^2$ is a polynomial, and we have already proven the polynomiality of $q^+(z)$. We conclude that $(q^{-}(z))^2$ is a polynomial. Analyzing the $z$ behavior as $z \rightarrow \infty$ using \eqref{eq_x136}, we also have that
\[
q^-(z) = -\hat D' z + O(1).
\]
Hence, $(q^{-}(z))^2$ is a degree two polynomial in $z$ with leading coefficient $(\hat D')^2$. By Theorem~\ref{Theorem_regularity_density} both endpoints of the band, $\alpha$ and $\beta$, are roots of $q^{-}(z)$. Therefore, we reach the formula:
\begin{equation}
 q^{-}(z)=-\hat D' \sqrt{(z-\alpha)(z-\beta)},
\end{equation}
where the branch of the square root is chosen again to map $\amsmathbb R_{>0}$ to itself. We can now apply Proposition~\ref{Proposition_density} leading to
\[
\widetilde \mu(x) = \frac{1}{\pi\theta}\,\arctan\bigg(\frac{q^{-}(x^-)-q^{-}(x^+)}{2\ii q_+(x)}\bigg)=\frac{1}{\pi\theta}\,\arctan\bigg(\hat D'\frac{\sqrt{(\beta - x)(x - \alpha)}}{q_+(x)}\bigg).
\]

\medskip

\noindent \textsc{Step 4.} It remains to examine $\widetilde \mu$ outside the band. There are two segments $[\hat A'_1,\alpha]$ and $[\beta, \hat B'_1]$; each of them can be either void or saturated. We first suppose that the parameters of the ensemble are in $\mathcal{D}^+$; in this case the denominator in \eqref{eq_Wzw} is positive. We deduce from \eqref{eq_Wzw}, \eqref{eq_x137}, and discussion after it that if $q^+(\hat{A}'_1) <0$, then $\exp\big(\theta\, \Gm_{\widetilde \mu}(\hat A'_1)\big)$ is a negative real number and, therefore, $[\hat{A}'_1,\alpha]$ is a saturation. On the other hand, if $q^+(\hat A'_1)>0$, then for small $\eps>0$, we see that
$\exp\big(\theta\, \Gm_{\widetilde \mu}(\hat A'_1+\eps)\big)$ is a positive real number and, therefore, $[\hat{A}'_1,\alpha]$ is a void. The same discussion holds for $\hat{B}'_1$: if $q^+(\hat{B}'_1)>0$, then $[\beta, \hat B'_1]$ is a void, and if $q^+(\hat{B}'_1)<0$, then $[\beta, \hat B'_1]$ is a saturation. The cases $q^+(\hat{A}'_1) = 0$ or $q^+(\hat{B}'_1) = 0$ are critical. They correspond to $\hat A'_1 = \alpha$ and $\hat B'_1 = \beta$, respectively, because $L(\hat A'_1)=0$ or $L(\hat B'_1)=0$ in these cases, as seen directly from the definition \eqref{eq_x132}.

If the parameters of the ensemble are in $\mathcal{D}^-$, the denominator of \eqref{eq_Wzw} becomes negative for $x\in \amsmathbb A$ because $\hat{B}'_2 < \hat{A}'_1$. Thus, the above discussion is still valid, except that the inequalities determining whether the endpoints of the segments are in a void or a saturation are reversed.
\end{proof}

\begin{definition}
\label{Definition_off_critical_set_ZW}
 $\mathcal{D}^{\textnormal{off}}$ is the open subset of $\mathcal D$ of Definition~\ref{Definition_ZW_set}, consisting of points such that
\begin{equation}
\label{eq_D_off_def}
\big(\hat{D}_{1,1}\hat{D}_{2,1}+\theta\hat n (\hat A_2-\hat A_1)\big) \cdot \big(\hat{D}_{1,1}\hat{D}_{1,2}-\theta\hat{n}(\hat B_2-\hat B_1)\big) \neq 0.
\end{equation}
\end{definition}
\begin{corollary} \label{Corollary_ZW_off-critical_1}
 Consider the variational datum from \eqref{eq_ZW_data_1}, \eqref{eq_ZW_data_2} and suppose that $(\theta, \hat n, \hat A_1', \hat A_2', \hat B_1', \hat B_2')$ belongs to a compact of $\mathcal{D}^{\textnormal{off}}\subset \mathcal D$. Then the corresponding equilibrium measure $\widetilde \mu$ is off-critical, \textit{i.e.} it satisfies Assumption~\ref{Assumption_C} in which the constants can be chosen to depend only on the said compact. In particular, for parameters in a given compact of $\mathcal{D}^{\textnormal{off}}$, the $zw$-ensemble satisfies Assumptions~\ref{Assumptions_Theta}, \ref{Assumptions_basic}, \ref{Assumptions_offcrit} and \ref{Assumptions_analyticity} for $\N$ large enough.
 \end{corollary}
\begin{proof} We first note that the product \eqref{eq_D_off_def} can be rewritten as
\[
 \big(\hat{D}'_{1,1}\hat{D}'_{2,1}+\theta\hat n (\hat A'_2-\hat A'_1)\big) \cdot \big(\hat{D}'_{1,1}\hat{D}'_{1,2}-\theta\hat{n}(\hat B'_2-\hat B'_1)\big) = q^+(\hat A_1') \cdot q^+(\hat B_1').
\]
Hence, for the parameters varying inside a given compact of $\mathcal{D}^{\textnormal{off}}$ both $q^+(\hat A_1')$ and $q^+(\hat B_1')$ are bounded away from $0$. Since the endpoints of the band $[\alpha,\beta]$ are the roots of $L(z)$ given by \eqref{eq_x132}, we conclude that $\alpha$ is bounded away from $\hat A_1'$ and $\beta$ is bounded away from $\hat B_1'$. This gives a lower bound on the lengths of voids and saturations and verifies Condition 1. in Assumption~\ref{Assumption_C}.

In order to verify the next conditions in Assumption~\ref{Assumption_C}, we compute the derivative of the effective potential:
\[
\forall x \in [\hat A'_1,\hat B'_1] \qquad \partial_x \widetilde V^\textnormal{eff}(x)= \partial_x \widetilde V(x)- 2\theta\, \Gm_{\widetilde \mu}(x)=\log\bigg( \frac{(x-\hat A'_1)(x-\hat A'_2)}{(B_1'-x)(B_2'-x)}\bigg)- 2\theta\, \Gm_{\widetilde \mu}(x),
\]
where the integral in the definition of $\Gm_{\widetilde \mu}(x)$ should be understood with a Cauchy principal value, \textit{i.e.}
\[
 \Gm_{\widetilde \mu}(x)=\frac{1}{2}\big( \Gm_{\widetilde \mu}(x^+)+ \Gm_{\widetilde \mu}(x^-)\big).
\]
We would like to understand the sign of $\partial_x \widetilde V^\textnormal{eff}(x)$ in voids/saturations $[\hat A_1',\alpha]$ and $[\beta, \hat B_1']$. For this purpose, we exponentiate it and use \eqref{eq_Wzw_final}. We get
\[
e^{-\partial_x \widetilde V^\textnormal{eff}(x)} = \frac{\Big(q^+(x) - \sqrt{(q^+(x))^2 - 4(\hat{B}'_1 - x)(\hat{B}'_2 - x)(x - \hat{A}'_1)(x - \hat{A}'_2)}\Big)^2}{4(\hat{B}'_1 - x)(\hat{B}'_2 - x)(x-\hat A'_1)(x-\hat A'_2)}.
\]
 The last identity can be rewritten as
\begin{equation}
\label{eq_x138}
e^{-\partial_x \widetilde V^\textnormal{eff}(x)} =\big(f(x) - \sqrt{(f(x))^2 - 1}\big)^2,
\end{equation}
\[
 f(x)=\frac{q^+(x)}{2\sqrt{(\hat{B}'_1 - x)(\hat{B}'_2 - x)(x-\hat A'_1)(x-\hat A'_2)}}.
\]
At this point we have four cases depending on whether we deal with parameters in $\mathcal{D}^+$ or $\mathcal{D}^-$, and whether we deal with a neighborhood of $\hat A'_1$ or $\hat B'_1$. All cases are similar and we provide details only for the parameters in $\mathcal{D}^+$, analyzing the effective potential in $[\hat A_1',\alpha]$. In this situation, by Theorem~\ref{Proposition_ZW_equilibrium_measure}, if $q^+(\hat A_1')>0$, then $[\hat A_1',\alpha]$ is a void; otherwise, it is a saturation.

In the former case, $q^+(x)$ is positive, we have $f(x)>1$ for all $x\in (\hat A_1',\alpha)$, and $f(\alpha)=1$. In fact, the comparison of $f(x)$ with $1$ is equivalent to comparison of $L(x)$ with $0$ in the proof of Theorem~\ref{Proposition_ZW_equilibrium_measure}. Using \eqref{eq_x138} we conclude that $0< e^{-\partial_x \widetilde V^\textnormal{eff}(x)} <1$, which implies $\partial_x \widetilde V^\textnormal{eff}(x)>0$ for any $x\in (\hat A_1',\alpha)$. In other words, the effective potential is strictly monotone, thus, satisfying Condition 2. of Assumption~\ref{Assumption_C}.

In the latter case, $q^{-}(x)$ is negative, we have $f(x)<-1$ for any $x\in (\hat A_1',\alpha)$, and $f(\alpha)=-1$. Using \eqref{eq_x138} we conclude that $e^{-\partial_x \widetilde V^\textnormal{eff}(x)}>1$, which implies $\partial_x \widetilde V^\textnormal{eff}(x)<0$ for any $x\in (\hat A_1',\alpha)$. Thus, the effective potential is again strictly monotone, thus, satisfying Condition 3. in Assumption~\ref{Assumption_C}.

The remaining conditions 4.5.6. in Assumption~\ref{Assumption_C} all immediately follow from the formula \eqref{eq_ZW_density_final} for the density of $\widetilde \mu$.

Assumption~\ref{Assumptions_offcrit} corresponds to Assumption~\ref{Assumption_C} for the original (non-simplified) equilibrium measure of the ensemble. Since the latter come from variational data differing from the simplified one by $O(\frac{1}{\N})$-corrections, we can use Theorem~\ref{Theorem_off_critical_neighborhood} to conclude that off-criticality for the simplified equilibrium measure implies off-criticality for the non-simplified one. The other Assumptions~\ref{Assumptions_Theta}, \ref{Assumptions_basic} and \ref{Assumptions_analyticity} were already checked in Lemma~\ref{Checkzw}.
\end{proof}

We now switch from the simplified equilibrium measure $\widetilde \mu$, corresponding to the potential $\widetilde V(x)$ of \eqref{eq_ZW_data_2} to the true equilibrium measure $\mu$ for $zw$-measures, corresponding to the potential specified by \eqref{eq_ZW_U} and \eqref{eq_ZW_V}.

\begin{corollary} \label{Corollary_ZW_off-critical_2}
 Consider the variational datum for $zw$-measures given by the segments \eqref{eq_ZW_data_1} and potential from \eqref{eq_ZW_U} and \eqref{eq_ZW_V}. Suppose that $(\theta, \hat n, \hat A_1, \hat A_2, \hat B_1, \hat B_2)$ vary in a given compact subset of $\mathcal{D}^{\textnormal{off}}\subset \mathcal D$. For $\N$ large enough depending only on the compact, the corresponding equilibrium measure $\mu$ is off-critical, \textit{i.e.} it satisfies Assumption~\ref{Assumption_C} with constants in the assumption depending only on the said compact.
\end{corollary}
\begin{proof} Using the Stirling formula, we see that the difference between the potential $V(x)$ of \eqref{eq_ZW_U}, \eqref{eq_ZW_V} and the potential $\widetilde V(x)$ of \eqref{eq_ZW_data_2} becomes arbitrary close to the constant $\hat A'_1-\hat B'_1+\hat A'_2-\hat B'_2$
 as $\N\rightarrow\infty$. Corollary~\ref{Corollary_ZW_off-critical_1} says that the simplified equilibrium measure $\widetilde \mu$ is off-critical. Hence, using Theorem~\ref{Theorem_off_critical_neighborhood} and the observation that shifting of the potential by a constant does not change the equilibrium measure, we conclude that the true equilibrium measure $\mu$ is also off-critical.
\end{proof}

We end this section by rewriting the leading term of the expansion of Proposition~\ref{Proposition_ZW_partition} in terms of the equilibrium measure.

\begin{proposition}
\label{Proposition_ZW_partition_leading} In the setting of Proposition~\ref{Proposition_ZW_partition}, the leading term of the partition function satisfies
\begin{equation} \label{eq_x224}
\F^{[0],\textnormal{ref}} = \I[\mu]+ \frac{\mathbbm{Rest}_0^{\textnormal{ref}}}{\N}+ o\bigg(\frac{1}{\N^{2}}\bigg),
\end{equation}
as $\N \rightarrow \infty$. Here $\mu$ is the equilibrium measure, $-\I$ is the energy functional that it minimizes, $\mathbbm{Rest}_0^{\textnormal{ref}}$ depends on $\N$ and are twice-continuously differentiable with respect to $(\hat{n},\hat{A}_1,\hat{A}_2,\hat{B}_1,\hat{B}_2)$ for $(\theta, \hat{n},\hat{A}_1,\hat{A}_2,\hat{B}_1,\hat{B}_2)\in\mathcal D$ with derivatives of orders $0$, $1$, $2$, uniformly bounded in $\N$. Finally, $o(\frac{1}{\N^{2}})$ is a term which after multiplication by $\N^{2}$ tends to $0$ as $\N\rightarrow\infty$ uniformly over the parameters varying in compact subsets of $\mathcal D$.
\end{proposition}
 Since $\mathbbm{Rest}_0^{\textnormal{ref}}$ depends on $\N$, one could have thought of absorbing the $o\big(\frac{1}{\N^2}\big)$ term inside it. However, the difference between these terms is that we claim bounds on partial derivatives for $\mathbbm{Rest}_0^{\textnormal{ref}}$, but not for $o\big(\frac{1}{\N^2}\big)$ remainder.
\begin{proof}
Let  $-\widetilde \I$ be the energy functional corresponding to the simplified potential \eqref{eq_ZW_data_2} and minimized by $\widetilde \mu$. We claim that
 \begin{equation}
 \label{eq_ZW_leading_partition}
 \begin{split}
& \quad \widetilde \I [\widetilde \mu] =  \F^{[0],\textnormal{ref}} \\
& = \frac{1}{\theta}\bigg[\mathrm{LLlog}(\theta\hat{n}) + \mathrm{LLlog}(\hat{D}' + \theta\hat{n}') - \mathrm{LLlog}(\hat{D}')  - \sum_{k,\ell =1}^2 \big(\mathrm{LLlog}(\hat{D}'_{k,\ell} + \theta\hat{n}) - \mathrm{LLlog}(\hat{D}_{k,\ell}')\big)\bigg]
\end{split}
\end{equation}
The identity \eqref{eq_ZW_leading_partition} readily implies \eqref{eq_x224} by comparing $\widetilde \I [\widetilde \mu]$ to $\I [\mu]$ through an application of Proposition~\ref{proposition_Energy_series}, in which we put $t=\tfrac{1}{\N}$ and combine the $t$ and $t^2$ terms together in a single error term of the same kind as $\frac{\mathbbm{Rest}_0^{\textnormal{ref}}}{\N}$ in the statement of Proposition~\ref{Proposition_ZW_partition_leading}.

It should be possible to directly compute $\widetilde \I [\widetilde \mu]$ to prove \eqref{eq_ZW_leading_partition}. However, rather than executing this computation we proceed in an indirect way. We view \eqref{eq_ZW_leading_partition} as an identity of two explicit functions of six parameters $(\theta, \hat{n},\hat{A}'_1,\hat{A}'_2,\hat{B}'_1,\hat{B}'_2)\in\mathcal D$. These functions do not depend on $\N$. Take now a $zw$-discrete ensemble with parameters $\big(\hat{n}^{(\N)},\hat{A}^{\prime\,(\N)}_1,\hat{A}^{\prime\,(\N)}_2,\hat{B}^{\prime\,(\N)}_1,\hat{B}^{\prime\,(\N)}_2\big)$ depending on $\N$ and tending to $(\hat{n},\hat{A}'_1,\hat{A}'_2,\hat{B}'_1,\hat{B}'_2)$ as $\N \rightarrow \infty$. We denote $\Z_\N^{(\N)}$ its partition function, and compute from Proposition~\ref{Proposition_ZW_partition}
\begin{equation}
\label{eq_x225}
\lim_{\N \rightarrow \infty} \frac{\log \Z_\N^{(\N)}}{\N^2} = \F^{[0],\textnormal{ref}}.
\end{equation}
On the other hand, Proposition~\ref{leadingZN} yields that
 \begin{equation}
\label{eq_x226}
\lim_{\N\rightarrow\infty}\bigg( \frac{\log \Z_\N^{(\N)}}{\N^2} - \I^{(\N)}[\mu^{(\N)}]\bigg)=0,
\end{equation}
where $\mu^{(\N)}$ is the equilibrium measure of the $zw$-ensemble for the particular value of $\N$ and $-\I^{(\N)}$ is the corresponding energy functional that it minimizes. Let $\widetilde \mu^{(\N)}$ denote the simplified equilibrium measure obtained by replacing the potential of the $\N$-dependent discrete $zw$-ensemble with its simplified version \eqref{eq_ZW_data_2} still using the $\N$-dependent parameters. The corresponding energy functional is denoted $-\widetilde \I^{(\N)}$. The relation between $\mu^{(\N)}$ and $\widetilde \mu^{(\N)}$ is the same as between $\mu$ and $\widetilde \mu$. Applying Proposition~\ref{proposition_Energy_series} with $\epsilon=\frac{1}{\N}$, we conclude that
\begin{equation}
\label{eq_x227}
\lim_{\N\rightarrow\infty} \big(\I^{(\N)}[\mu^{(\N)}] -\widetilde \I^{(\N)}[\widetilde \mu^{(\N)}]\big) =0.
\end{equation}
As $\N\rightarrow\infty$, the variational datum associated to the equilibrium measure $\mu^{(\N)}$ converges to the variational datum associated with the equilibrium measure $\widetilde \mu$: the approximation errors are $O\big(\frac{1}{\N}\big)$. Since these measures themselves were explicitly computed in Theorem~\ref{Proposition_ZW_equilibrium_measure}, we conclude that
\[
 \lim_{\N\rightarrow\infty} \widetilde \I^{(\N)}[\widetilde \mu^{(\N)}] =\widetilde \I[\widetilde \mu].
\]
Combining with \eqref{eq_x226}, \eqref{eq_x227},
 we conclude that
\begin{equation}
\label{eq_x228}
 \lim_{\N\rightarrow\infty} \frac{\log \Z_\N^{(\N)}}{\N^2} = \lim_{\N\rightarrow\infty} \I^{(\N)}[\mu^{(\N)}] = \widetilde \I[\widetilde \mu].
\end{equation}
Comparing \eqref{eq_x225} with \eqref{eq_x228} proves the claimed \eqref{eq_ZW_leading_partition}.
\end{proof}

\section{Adjusting the parameters}
\label{sec:Tuning}
In this section we continue analyzing the simplified equilibrium measures $\widetilde \mu$. Our task is to show that we can reach arbitrary types of the equilibrium measure and arbitrary lengths of the band by properly choosing the parameters of the variational datum. Note that eventually in the following chapters we are going to deal with the true equilibrium measure $\mu$, rather than its approximation $\widetilde \mu$. However, as $\mu$ and $\widetilde \mu$ become arbitrary close as $\N\rightarrow\infty$, it is sufficient to develop all the nice properties only for $\widetilde \mu$.

Our strategy for this section is to assume that we are given the parameters $\theta$, $\hat n$, the band $(\alpha,\beta)$, and one of the four types of the equilibrium measure \emph{void-band-void}, \emph{void-band-saturated}, \emph{saturated-band-void}, \emph{saturated-band-saturated}. We would like to choose the parameters $(\hat A_1', \hat A_2', \hat B_1', \hat B_2')$ of the $zw$-variational datum, so that $\widetilde \mu$ has the prescribed band and type. An additional requirement is that we want to have control over $\alpha-\hat A'_1$ and $\hat B'_1-\beta$, \textit{i.e.}, we want to guarantee that these numbers are neither too large nor too small.

There is more than one way to match the parameters of the $zw$-variational datum with the desired $(\alpha,\beta)$ and type. We present three different ways in the following three propositions. Each of them has some restrictions on the location and the length of the band and the type for which it works, but together these three propositions cover all the situations that we need. Let us emphasize again that in all three propositions we deal only with variational data and the equilibrium measure $\widetilde \mu$ of Theorem~\ref{Proposition_ZW_equilibrium_measure}. In particular, we are not imposing any integrality conditions on any of the involved parameters.

\begin{proposition}
\label{Proposition_zw_band_and_type_1}
 Fix $\eps>0$. There exists $\delta>0$ such that for each of the four types \emph{void-band-void}, \emph{void-band-saturated}, \emph{saturated-band-void}, \emph{saturated-band-saturated}, each $\hat{n} \in \big(\eps,\frac{1}{\eps}\big)$ and $\theta \in \big(\eps,\frac{1}{\eps}\big)$, each $\hat A_1'$, $\hat B_1'$
 satisfying
 \[
 \max\big(\hat B'_1- \hat A'_1- \theta \hat n,| \hat B'_1-\hat A'_1 -2\theta \hat n|\big) > \eps \qquad \textnormal{and} \qquad \max\big(|\hat A'_1|,|\hat B'_1|\big) < \frac{1}{\eps}
 \]
 and each interval $(\alpha,\beta)$ such that
 \[
 0 < \alpha - \hat A'_1 < \delta \qquad \textnormal{and} \qquad 0 < \hat B'_1 - \beta < \delta,
 \]
 we can find $\hat A'_2$ and $\hat B'_2$ such that the equilibrium measure $\widetilde \mu$ of Theorem~\ref{Proposition_ZW_equilibrium_measure} has the prescribed type and the band $(\alpha,\beta)$. Moreover, $\hat A'_2$ and $\hat B'_2$ can be chosen to depend smoothly on $\alpha,\beta, \hat A'_1, \hat B'_1, \hat n, \theta$.

 For the types \emph{void-band-void} and \emph{saturated-band-saturated}, if we additionally set $\alpha - \hat A'_1 = \hat B'_1 - \beta$, it is sufficient to require $|\hat B'_1-\hat A'_1 - 2\theta \hat n|\neq 0$ instead of $>\eps$.
\end{proposition}
The most important restriction in Proposition~\ref{Proposition_zw_band_and_type_1} is $|\hat B'_1-\hat A'_1 -2\theta \hat n|>\eps$ for the types \emph{void-band-saturated} and \emph{saturated-band-void}. Since $\alpha$ and $\beta$ are close to $\hat A'_1$ and $\hat B'_1$, this proposition is useful only in the situation when $| \beta-\alpha -2\theta \hat n|$ is bounded away from $0$.

\begin{proposition}
\label{Proposition_zw_band_and_type_2}
 Fix $\eps>0$. There exist $\eta>0$ and $\delta>0$, such that for each of the two types
 \emph{void-band-void} or \emph{saturated-band-void}, each $\hat{n} \in \big(\eps,\frac{1}{\eps})$ and $\theta \in \big(\eps,\frac{1}{\eps}\big)$, each $\hat A'_1 \in \big(-\frac{1}{\eps},\frac{1}{\eps}\big)$, and each interval $(\alpha,\beta)$
 such that
 \[
 0<\alpha-\hat A'_1<\delta\qquad \textnormal{and}\qquad \theta \hat \hat{n} +\eps< \beta-\alpha< 4 \theta \hat n,
 \]
we can find $\hat A'_2<\hat A'_1$ and $\hat B'_1=\hat B'_2=\hat B$ such that the simplified equilibrium measure $\widetilde \mu$ of Theorem~\ref{Proposition_ZW_equilibrium_measure} has the prescribed type and the band
 $(\alpha,\beta)$. In addition, we can enforce $\hat{B} - \beta > \eta$. Moreover, $\hat A'_2$ and $\hat B$,
 can be chosen to depend smoothly on $\alpha,\beta, \hat A_1', \hat n, \theta$.
\end{proposition}

\begin{proposition}
\label{Proposition_zw_band_and_type_3}
 Fix $\eps>0$. There exists $\eta>0$ and $\delta>0$, such that for each of the two types
 \emph{void-band-void}, \emph{void-band-saturated}, each $\hat{n} \in \big(\eps,\frac{1}{\eps}\big)$ and $\theta \in \big(\eps,\frac{1}{\eps}\big)$, each $\hat B'_1 \in \big(-\frac{1}{\eps},\frac{1}{\eps}\big)$, and each interval $(\alpha,\beta)$ such that
 \[
 0<\hat B'_1-\beta<\delta \qquad \textnormal{and}\qquad \theta \hat{n} +\eps< \beta-\alpha< 4 \theta \hat n,
 \]
 we can find $\hat B'_2>\hat B'_1$ and $\hat A'_1=\hat A'_2=\hat A$
 such that the simplified equilibrium measure $\widetilde \mu$ of Theorem~\ref{Proposition_ZW_equilibrium_measure} has the prescribed type and the band
 $(\alpha,\beta)$. In addition, we can enforce $\alpha-\hat{A}> \eta$. Moreover, $\hat A$ and $\hat B'_2$,
 can be chosen to depend smoothly on $\alpha,\beta, \hat B'_1, \hat n, \theta$.
\end{proposition}

 The most important restriction in the two last propositions is $\theta \hat \hat{n} +\eps< \beta-\alpha< 4\theta \hat n-\eps$, which constrains the possible length of the band.

\begin{proof}[Proof of Proposition~\ref{Proposition_zw_band_and_type_1}]
 We study the equilibrium measure $\widetilde \mu$ given in Theorem~\ref{Proposition_ZW_equilibrium_measure} for real parameters $\theta$, $\hat n$, $\hat A'_1$, $\hat
 A'_2$, $\hat B'_1$, $\hat B'_2$. Note that if we shift the four parameters $\hat A'_1$, $\hat A'_2$, $\hat B'_1$, $\hat
B'_2$ by the same number $c$, then the equilibrium measure $\widetilde \mu$ and its band
$(\alpha,\beta)$ are also translated by $c$. Therefore, it suffices to consider
only the symmetric case $\hat A'_1=-\hat B'_1$. We use the shortcut notation $T := \theta\hat{n}$. Recall that the function $q^+(z)$ is not exactly the function $\widetilde{q}^+(z)$ that Definition~\ref{GQdef} associates to $\widetilde{\mu}$, but rather given by \eqref{eq_x130}.

 Take $S> \frac{T}{2}$ and $R\in\amsmathbb R$ and set
 \begin{equation}
 \label{eq_x145}
 \hat A'_1=-S, \qquad \hat B'_1=S,\qquad \hat A'_2=-R,\qquad \hat B'_2=R
 \end{equation}
 to start with. We tune these parameters in such a way that $q^+(-S) = q^+(S)=0$, which means by \eqref{eq_x130} that
 \begin{equation}
 \label{eq_x146}
 R=-\frac{1}{2} \frac{(S- T)^2+S^2}{S- T}.
 \end{equation}
 The requirement $|\hat B'_1-\hat A'_1 - 2\theta \hat n|\neq 0$ in Proposition~\ref{Proposition_zw_band_and_type_1} guarantees that the
 denominator in \eqref{eq_x146} does not vanish.
 Note that $|R|>S>0$, and therefore, $(\theta, \hat n, -S, -R, S,R)\in \mathcal D$. Moreover, by definition of the quantity $L(x)$ in \eqref{eq_discrizw}, it also vanishes at $x = \pm S$ in this case and, therefore, the band $(\alpha,\beta)$ is exactly $(\alpha,\beta) = (-S,S)$.

Next, we verify
 that by slightly changing $\hat A_2$ and $\hat B_2$ from the values given by \eqref{eq_x145}, we can get
 arbitrary endpoints of the band $(\alpha,\beta)$ satisfying
 $0<\alpha-\hat A'_1$ and $0 < \hat B'_1-\beta<\delta$. For this purpose we take $\mathfrak{a},\mathfrak{b} \in \amsmathbb{R}$ small enough in absolute value and set
\[
\hat A'_1=-S,\qquad \hat B'_1=S,\qquad \hat A'_2=-R + R\mathfrak{a},\qquad \hat B'_2=R+R\mathfrak{b}.
\]
Although we could directly analyze the formulae \eqref{eq_x130} for $\alpha$ and $\beta$, we find simpler to look at the equation $L(x)=0$ which they solve. We recall the definition of the function $L(x)$ given in \eqref{eq_x132} and study its dependence on $\mathfrak a$ and $\mathfrak b$ for $x$ near $\pm S$. Near $x = S$ we have
\begin{equation}
\begin{split}
\label{eq_x148}
 4 (S- T)^2\cdot L(S+\xi) & = (T^2 - 2TS + 2S^2)^2 \cdot \big((2S -T)\mathfrak{a} + T\mathfrak{b}\big)^2 + 4T^2(2S - T)^2(2S\xi + \xi^2) \\
& \quad + Q_1(\mathfrak{a},\mathfrak{b})\xi + Q_2(\mathfrak{a},\mathfrak{b})\xi^2,
\end{split}
\end{equation}
where $Q_1$ and $Q_2$ are linear and quadratic polynomials in $(\mathfrak{a},\mathfrak{b})$ whose coefficients are uniformly bounded as soon as $\eps$ is fixed. For small $\mathfrak{a},\mathfrak{b}$, the right endpoint $\beta$ of the band is found as $S$ plus the closest to $0$ root of \eqref{eq_x148} in the variable $\xi$. Therefore,
\begin{equation}
\label{eq_x149}
 \beta=S - \frac{(T^2 - 2ST + 2S^2)^2 \cdot \big((2S - T)\mathfrak{a} + T\mathfrak{b}\big)^2}{8T^2(2S - T)^2S} + Q_3(\mathfrak{a},\mathfrak{b}),
 \end{equation}
where $Q_3(t\mathfrak{a},t\mathfrak{b}) = O(t^3)$ as $t \rightarrow 0$. Similarly, near $x = -S$ we have
\begin{equation*}
\begin{split}
4(S - T)^2\cdot L(-S + \xi) & = (T^2 - 2ST + 2S^2)^2\cdot \big(T\mathfrak{a} + (2S - T)\mathfrak{b}\big) + 4T^2(2S - T)^2(\xi^2 - 2S\xi) \\
& \quad + \widetilde Q_1(\mathfrak{a},\mathfrak{b})\zeta + \widetilde{Q}_2(\mathfrak{a},\mathfrak{b})\zeta^2,
\end{split}
\end{equation*}
where $\widetilde Q_1$ and $\widetilde{Q}_2$ are different expressions with the same property as $Q_1$ and $Q_2$ before. Adding to $-S$ the closest to $0$ root of this expression in the variable $\xi$ yields
\begin{equation}
\label{eq_x150} \alpha = -S + \frac{(T^2 - 2ST + 2S^2)^2\cdot \big(T\mathfrak{a} + (2S - T)\mathfrak{b}\big)^2}{8T^2(2S - T)^2S} + \widetilde Q_3(\mathfrak{a},\mathfrak{b}),
\end{equation}
with $\widetilde Q_3(t\mathfrak{a},t\mathfrak{b}) = O(t^3)$ as $t \rightarrow 0$.

\medskip

Our assumption $\hat B'_1- \hat A'_1- \theta \hat n>\eps$ translates into $2S - T>\eps$. Hence, the denominators in \eqref{eq_x149} and \eqref{eq_x150} are bounded away from $0$. Further, the assumption $|\hat B_1-\hat A_1-2\theta\hat n|> \eps$ translates into $|S- T|> \frac{\eps}{2}$, which guarantees that the expressions
\begin{equation}
\label{eq_x151} (2S - T)\mathfrak{a} + T\mathfrak{b} \qquad \textnormal{ and }\qquad T\mathfrak{a} + (2S - T)\mathfrak{b}
\end{equation}
are two independent linear forms in $(\mathfrak{a},\mathfrak{b})$. Therefore,
 using the formulae \eqref{eq_x149}, \eqref{eq_x150}, which are
approximations for the explicit formulae \eqref{eq_x130} for $\alpha$ and $\beta$,
for any small enough $\delta>0$ and any fixed $(\alpha^0,\beta^0)$ for which $\alpha^0 \in (-S,-S + \delta)$ and $\beta^0 \in (S - \delta,S)$ we can find \emph{four} ordered pairs
$(\mathfrak{a},\mathfrak{b})$ such that $(\alpha,\beta)$ of \eqref{eq_x149},\eqref{eq_x150} coincides with $(\alpha^0,\beta^0)$. These four possible choices arise because of squaring in these
formulae: if one choice $(\mathfrak{a},\mathfrak{b})$ is valid, so is $(-\mathfrak{a},-\mathfrak{b})$.
As $S- T$ approaches $0$, we no longer can reach arbitrary values of $\alpha$ and $\beta$, because the linear forms \eqref{eq_x151} become collinear. However, as they approach the same linear form, we can still reach the symmetric case
\[
\alpha-\hat A'_1=\hat B'_1-\beta.
\]
This time we have only \emph{two} different choices for $(\mathfrak{a},\mathfrak{b})$, which differ by a global sign.

Let us determine the types of the equilibrium measure that correspond to each of these four (or two) choices. This is done according to Theorem~\ref{Proposition_ZW_equilibrium_measure} by looking at the signs of
\begin{equation*}
\begin{split}
q^+(\hat{B}'_1) & = \frac{(S - T)^2 + S^2}{2(S - T)}\big((2S - T)\mathfrak{a} + T\mathfrak{b}\big), \\
q^+(\hat{A}'_1) & = -\frac{(S - T)^2 + S^2}{2(S - T)}\big(T\mathfrak{a} + (2S - T)\mathfrak{b}\big).
\end{split}
\end{equation*}
We see that if $|S - T| > \eps$, then the four aforementioned choices realize the four types \emph{void-band-void}, \emph{void-band-saturated}, \emph{saturated-band-void}, \emph{saturated-band-saturated}; otherwise, the two aforementioned choices realize the \emph{void-band-void} and \emph{saturated-band-saturated} types.
\end{proof}

\begin{proof}[Proof of Proposition~\ref{Proposition_zw_band_and_type_2}]
As in the previous proof, shifting $(\hat{A}'_1,\hat{B}'_1,\hat{A}'_2,\hat{B}'_2)$ by a constant merely translates the equilibrium measure $\widetilde \mu$ and its band, so we can set $\hat{A}'_1 = 0$ without loss of generality. We still use the notation $T = \theta\hat{n}$, and set
\[
 \hat B'_1=\hat B'_2=\hat B,\qquad \hat A'_2=-\frac{(\hat B-T)^2}{T} + \mathfrak{a},
\]
where we arranged the parameters so that $q^+(0)$ vanishes at $\mathfrak{a} = 0$. Note that $q^+(0)$ is a linear function of the variable $\mathfrak a$. By choosing small enough $\mathfrak{a}$ we can guarantee either
$q^{+}(0)>0$ or $q^+(0)<0$, thus (\textit{cf.} Theorem~\ref{Proposition_ZW_equilibrium_measure})
forcing $\widetilde \mu$ to have either a void or a saturation near $0$. This also implies,
that the left endpoint $\alpha$ of the band equals $0$ at $\mathfrak{a}=0$. This can be also
checked directly as the equation \eqref{eq_x132} transforms into
\begin{equation}
\label{eq_x152}
L(x) = T^2\mathfrak{a}^2 + \frac{2\big(2\hat{B}^2(\hat{B} - T)^2 - T(\hat{B} - T)(3\hat{B} - T)\mathfrak{a} + T^2\mathfrak{a}^2 \big)}{T}x + \frac{(\hat{B}^2 - T^2 - T\mathfrak{a})^2}{T^2} x^2,
\end{equation}
which has a root at $x = 0$ if $\mathfrak{a}=0$.
On the other hand,
\[
q^+(\hat B)= \frac{\hat{B} - T}{T}\big(\hat{B}(\hat{B} - T) - T\mathfrak{a}\big)
\]
is positive for $\mathfrak{a}$ small enough, forcing $\widetilde \mu$ to have a void near $\hat{B}$. Therefore, by Theorem~\ref{Proposition_ZW_equilibrium_measure}, we realize one of the two
types \emph{void-band-void} or \emph{saturated-band-void}, as desired.

The left endpoint of the band is found as the closest to $0$ root of
\eqref{eq_x152} in the variable $x$, which gives
\begin{equation}
\label{eq_x154} \alpha \mathop{=}_{\mathfrak{a} \rightarrow 0} \,\,\frac{T^3}{ 4 \hat B^2 (\hat B- T)^2}\,\mathfrak{a}^2 + O(\mathfrak{a}^3).
\end{equation}
The length of the band \eqref{width} becomes
\begin{equation}
\label{eq_x153}
\begin{split}
\beta - \alpha & = \frac{4T(\hat{B} - T)\big(\hat{B}(\hat{B} - T) - T\mathfrak{a}\big)\sqrt{\hat{B}^2 - T\mathfrak{a}}}{(\hat{B}^2 - T^2 - T\mathfrak{a})^2} \\
& \mathop{=}_{\mathfrak{a} \rightarrow 0} \,\,\frac{4T\hat{B}^2}{(\hat{B} + T)^2}\bigg(1 + \frac{(\hat{B} - T)T}{2\hat{B}^2(\hat{B} + T)}\,\mathfrak{a} + O(\mathfrak{a}^2)\bigg).
\end{split}
\end{equation}
Combining \eqref{eq_x154} with \eqref{eq_x153}, we can obtain arbitrary (small)
values for $\alpha$ and arbitrary values for $(\beta-\alpha)$ in $(T,4T)$, where the last restriction arises because of the inequality
\[
T < \frac{4T \hat B^2 }
 {(\hat B+T)^2}< 4T,\quad \textnormal{for}\,\,\hat B> T.
\]
Note that since \eqref{eq_x154} contains $\mathfrak{a}^2$, we can have both signs for $\mathfrak{a}$,
and thus realize the \textit{void-band-void} and \textit{saturated-band-void} types of the equilibrium measure.

It remains to prove the lower bound asserted in the proposition holds for $\hat{B}'_1 - \beta$. From \eqref{eq_x154} and \eqref{eq_x153} we deduce that
\[
 \beta \mathop{=}_{\mathfrak{a} \rightarrow 0} \,\,\frac{4T\hat{B}^2}{(\hat{B} + T)^2}+O(\mathfrak{a}).
\]
Hence,
\[
 \hat B'_1-\beta=\hat B-\beta \mathop{=}_{\mathfrak{a} \rightarrow 0}\,\, \frac{\hat B (\hat{B} + T)^2- 4T\hat{B}^2}{(\hat{B} + T)^2}+O(\mathfrak{a})= \frac{\hat B (\hat{B} -T)^2}{(\hat{B} + T)^2}+O(\mathfrak{a}).
\]
Because $\hat B$ and $\hat T$ are bounded away from $0$ and $\infty$ and $\hat{B} -T>\beta-\alpha-T>\eps$ by our assumption, the last expression is bounded away from $0$ for small $\mathfrak a$.
\end{proof}

\begin{proof}[Proof of Proposition~\ref{Proposition_zw_band_and_type_3}]
We repeat the arguments of Proposition~\ref{Proposition_zw_band_and_type_2}, applying the symmetry $(\hat
A'_1,\hat A'_2,x) \longleftrightarrow (-\hat B'_1,-\hat B'_2,-x)$.
\end{proof}

\chapter{Conditioning and localizing discrete ensembles}
\label{Chapter_conditioning}

In this chapter we condition the discrete ensemble on an event of overwhelming probability as $\N$ becomes large, so that we still get a discrete ensemble that fits the general framework of Chapter~\ref{Chapter_Setup_and_Examples}, but with new parameters. Roughly speaking, we condition the ensemble so that particles are frozen in saturations and that there are no particles in void regions. The conditioned ensemble satisfies all the assumptions of the original ensemble, but in addition it satisfies Assumption~\ref{Assumptions_extra}, and in particular it has exactly one band per segment; these properties will become important in Chapter~\ref{Chapter_fff_expansions} to perform the asymptotic analysis of the correlators of the empirical measure. The asymptotics of the partition function are then recovered in Chapter~\ref{Chapter_partition_functions}  by interpolation with an ensemble closely related to the $zw$-measures studied in Chapter~\ref{Chapterzw}. Finally, the asymptotics of the original ensemble are recovered by de-conditioning.

Our general strategy is to start with an ensemble satisfying the Assumptions~\ref{Assumptions_Theta}--\ref{Assumptions_analyticity} and construct a new one that satisfies these assumptions as well as Assumption~\ref{Assumptions_extra}, in two steps: conditioning and localizing. We begin by selecting a \emph{localization segment} around each band of the equilibrium measure of the initial ensemble. The endpoints of these segments are chosen to be close to the endpoints of the bands, which eventually leads to all the conditions of Assumption~\ref{Assumptions_extra} being automatically satisfied. The \emph{conditioning} procedure consists in freezing, outside the localization segments, the particles of the saturations at deterministic locations and prohibiting particles from being in voids. The \emph{localization} procedure consists in forgetting the positions of the frozen particles. The new ensemble describes the distribution of configurations of particles within the localization segments, which extend (by adding particles at deterministic locations) to configurations in the conditioned initial ensemble.

Compared to continuous one-dimensional Coulomb gases, the conditioning procedure is complicated by the fact that the locations allowed for the particles in our setting are discrete. Moreover, if some of the diagonal entries of $\boldsymbol{\Theta}$ are not equal to $1$, then the allowed locations depend on the number of particles to the left, so we are forced to choose localization segments that depend on segment filling fractions. In the case where the filling fractions are not fixed in the initial ensemble --- \textit{i.e.}, if the equations \eqref{eq_equations_eqs} in Section~\ref{DataS} have more than one solution --- the above conditioning is performed simultaneously with fixing the values of the filling fractions. Due to the fluctuating nature of the lattices in which the particles live for general matrices $\boldsymbol{\Theta}$, it would be impossible to freeze the positions of the particles in saturations without fixing the filling fractions.

\section{Fixed filling fractions case}
\label{fixedsecfill}
We recall that for the segment endpoints of the discrete ensemble we use $a_h,b_h$ and their rescaled version $\hat a_h$, $\hat b_h$ (Definition~\ref{def:eq_rescaled_parameters}), while for the quantities involving the continuous space and the equilibrium measure, we use the shifted rescaled parameters $\hat a'_h$, $\hat b'_h$ (Definition~\ref{def:eq_shifted_parameters}).
\[
\amsmathbb{A} = \bigcup_{h = 1}^{H} [\hat{a}_h',\hat{b}_h']
\]
involves exclusively the shifted parameters. This shift causes a small hiatus when talking about particles and holes in voids and saturations, because particles live in the segments of the ensemble, while voids and saturations are notions pertaining to the equilibrium measure, which is supported on $\amsmathbb{A}$. We will make the distinction clear when it is necessary.

\begin{theorem}
\label{proposition_FFF_conditioning}
 Suppose that a discrete ensemble $X$ satisfies Assumptions~\ref{Assumptions_Theta}--\ref{Assumptions_analyticity}, that the equations \eqref{eq_equations_eqs} deterministically fix the filling fractions and that each segment $[\hat a_h^X, \hat b_h^X]$ contains exactly one band $(\alpha_h,\beta_h)$ of the equilibrium measure.

 There exists $\eps_0>0$ depending only on the constants in the assumptions and such that for any $\eps \in (0,\eps_0)$, any choice of pairwise disjoint segments $[\hat a_h^Y, \hat b_h^Y]$ indexed by $h \in [H]$ and satisfying
\[
\forall h \in [H]\qquad [\alpha_h-\eps ,\beta_h+\eps]\subset [\hat a_h^Y, \hat b_h^Y]\subset [\alpha_h-2\eps ,\beta_h+2\eps]\subset [\hat a_h^X, \hat b_h^X],
\]
and an \emph{integrality condition} spelled out below, there exists another discrete ensemble $Y$ with same number of segments $H$, same intensities of interactions $\boldsymbol{\Theta}$, and defined on segments $[\hat a_h^Y, \hat b_h^Y]$ indexed by $h \in [H]$, and there exists a constant $C > 0$ depending only on $\eps$ and on the constants in the assumptions for the ensemble $X$, having the following properties for $\N$ large enough depending only on the constants in the assumptions.
 \begin{enumerate}
 \item The ensemble $Y$ satisfies Assumptions~\ref{Assumptions_Theta}--\ref{Assumptions_analyticity} as well as Assumption~\ref{Assumptions_extra}, with constants that only depend on the constants in the assumptions for the ensemble $X$ and on $\eps$.
 \item There is an event $\mathcal{A}$ of probability at least $1-C\exp\big(-\frac{\N}{C}\big)$ in the ensemble $X$, such that the ensemble $Y$ is obtained by localizing to $\bigcup_{h = 1}^{H} [\hat{a}_h^Y,\hat{b}_h^Y]$ the ensemble $X$ conditioned on $\mathcal{A}$. In particular, the equilibrium measure for $Y$ coincides with the restriction to $\amsmathbb{A}^{Y}$ of the equilibrium measure for $X$.
 \item On $\mathcal{A}$, for any $h \in [H]$, each of the segments $[a_h^X,a_h^Y]$ either has no particles (if $[\hat{a}_h^{\prime\,X},\hat{a}_h^{\prime\,Y}]$ is void for $X$) or has no holes as in Definition~\ref{Definition_no_hole} (if $[\hat{a}_h^{\prime\,X},\hat{a}_h^{\prime\,Y}]$ is saturated for $X$). A similar property holds for the segments $[b_h^Y,b_h^X]$.
 \item For any $h \in [H]$, the filling fractions of the discrete ensemble $Y$ --- which are deterministic conditionally on $\mathcal{A}$ due to property 3. --- are obtained by subtracting from $\hat{n}_h^X$ the count of all particles in $[a_h^X,a_h^Y) \cup (b_h^Y,b_h^X]$ divided by $\N$ (there is something to subtract only if the corresponding rescaled and shifted segments are saturated for the ensemble $X$).
 \end{enumerate}
The aforementioned \emph{integrality condition} reads, for any $h \in [H]$
\begin{enumerate}
  \item If $a_h^X \neq -\infty$ and $[\hat a_h^{\prime\,X},\hat a_h^{\prime\,Y}]$ is void for $X$, we require $\N(\hat a_h^Y-\hat a_h^X)\in\amsmathbb Z_{\geq 0}$. If $a_1^X=-\infty$, we require that $\N \hat{a}_1^Y$ belongs to the lattice $\amsmathbb{L}$ in which $\ell_1$ varies.
  \item If $[\hat a_h^{\prime\,X},\hat a_h^{\prime\,Y}]$ is saturated for $X$, we require $\N(\hat a_h^Y-\hat a_h^X)\in \theta_{h,h} \amsmathbb Z_{\geq 0}$.
  \item Similar conditions are imposed on $\hat{b}_{h}^X,\hat{b}_h^Y$.
\end{enumerate}
We say that $Y$ is the localization of $X$ to $\amsmathbb{A}^{Y}$.
\end{theorem}

\begin{remark} \label{rem:misreg}In general the equilibrium measure of the ensemble $Y$ is not exactly the restriction to the equilibrium measure of the ensemble $X$ to $\amsmathbb{A}^{Y}$, but is $O(\frac{1}{\N})$ close. The reason is that the regular part of the potential --- the function $U_h(x)$ ot \eqref{eq_V_logs}, \eqref{eq_potential_weight_match} of the assumptions in Section~\ref{Section_list_of_assumptions} or in Assumption~\ref{Assumption_B} of Section~\ref{allassuml} and used to define the equilibrium measure  --- depends on the chosen segments; this is related to the assignment of $\err_h(x)$ in \eqref{eq_weight_form} and the choices we made to resolve it. The fact that there is a small mismatch is important and will show up for instance in the proof of Theorem~\ref{Theorem_correlators_expansion_relaxed}. We will revisit this mismatch in Section~\ref{Section_alternative_localization} and Corollary~\ref{mueqYYtilde}.
\end{remark}

\begin{proof}[Proof of Theorem~\ref{proposition_FFF_conditioning}]
We first consider any $H$-tuples of segment endpoints $\hat{\boldsymbol{a}}^{Y},\hat{\boldsymbol{b}}^{Y}$ satisfying the above integrality conditions and such that
\[
\forall h \in [H]\qquad [\alpha_h-\eps ,\beta_h+\eps]\subset [\hat a_h^Y, \hat b_h^Y]\subset [\alpha_h-2\eps ,\beta_h+2\eps]\subset [\hat a_h^X, \hat b_h^X].
\]
Thanks to Assumption~\ref{Assumptions_offcrit}, by choosing small $\eps>0$ we can guarantee that for each $h \in [H]$, the segment $[\hat a_h^{\prime\,X},\hat a_h^{\prime\,Y}]$ (respectively $[\hat b_h^{\prime\,Y},\hat b_h^{\prime\,X}]$) is either void or saturated for the equilibrium measure of the ensemble $X$. Then, by Assumption~\ref{Assumptions_offcrit} and Theorem~\ref{Theorem_ldpsup}, whenever $[\hat a_h^{\prime\,X},\hat a_h^{\prime\,Y}]$ (respectively $[\hat b_h^{\prime\,Y},\hat b_h^{\prime\,X}]$) is void, the corresponding segment $[a_h^X,a_h^Y]$ (respectively $[b_h^Y,b_h^X]$) has no particles of the ensemble $X$ with overwhelming probability. And by Theorem~\ref{Theorem_ldsaturated}, if $[\hat a_h^{\prime\,X},\hat a_h^{\prime\,Y}]$ (respectively $[\hat b_h^{\prime\,Y},\hat b_h^{\prime\,X}]$) is saturated, the corresponding segment $[a_h^X,a_h^Y]$ or $[b_h^Y, b_h^X]$, respectively, has no holes of the ensemble $X$ with overwhelming probability. We claim that conditioning on the intersection of these events of overwhelming probability and forgetting about the position of particles in $\bigcup_{h = 1}^{H} [a_h^X,a_h^Y) \cup (b_h^Y,b_h^X]$ gives an ensemble $Y$ with segments $[\hat a_h^Y,\hat b_h^Y]$ indexed by $h \in [H]$ meeting all the requirements of Theorem~\ref{proposition_FFF_conditioning}, as long as $\eps>0$ is small enough.

\smallskip

Let us first see that the ensemble $Y$ is indeed a discrete ensemble in the sense of Section~\ref{Section_general_model}. The idea is rather simple: for each $h \in [H]$, the potential $V_h$, the functions $\Phi^{\pm}_h$, and their leading order $\phi_h^{\pm}$ for the ensemble $Y$ will be obtained from those of $X$ by absorbing --- via additional terms or factors --- the interactions with particles of $\bigcup_{g = 1}^{H} [a_g^{X}, a_g^{Y}) \cup (b_g^{Y},b_g^{X}]$. Given the amount of data that specifies a discrete ensemble, the detailed procedure is a bit long to describe. We offer two descriptions for this. The first one is sequential: we explain a sequence of steps to perform in order to get the data of the ensemble $Y$ from the data of the ensemble $X$. The second one, postponed to Section~\ref{Section_alternative_localization}, is global: we give the final formulae for the data of the ensemble $Y$. The interested reader can test their understanding of the procedure by showing that the equivalence of the two descriptions.

The sequential description of the procedure is as follows. For each $h \in [H]$ such that $[\hat a_h^{\prime\,X},\hat a_h^{\prime\,Y}]$ is saturated, we apply a first list of updates on the functions associated with $g \in [H]$ in the ensemble $X$ in order to obtain the corresponding functions in the ensemble $Y$.

\medskip

\noindent \textsc{Weights.} The weight $w_{g}(x)$ is multiplied by the term coming from the interaction with the particles occupying all the sites available in $[a_h^X,a_h^Y)$, \textit{i.e.} by
 \begin{equation}
 \label{eq_x178}
  \prod_{\begin{smallmatrix}\ell= a_h^X +i \theta_{h,h}: \\ i\in \amsmathbb Z_{\geq 0}, \,\, \ell \leq a_h^Y-\theta_{h,h} \end{smallmatrix}} \frac{1}{\N^{2\theta_{h,g}}}\cdot \frac{\Gamma\big(|x-\ell|+1\big)\cdot\Gamma\big(|x-\ell|+\theta_{h,g}\big)}{\Gamma\big(|x-\ell|\big)\cdot \Gamma\big(|x-\ell|+1-\theta_{h,g}\big)}.
 \end{equation}
 Note that for $g=h$ we always have $x>\ell$ and the product \eqref{eq_x178} telescopes to
  \begin{equation}
 \label{eq_x179}
  \frac{1}{\N^{2 \N (\hat a_h^Y-\hat a_h^X)}} \cdot \frac{\Gamma\big(x- \N\hat a_h^X+1\big)\cdot\Gamma\big(x-\N a_h^X+\theta_{h,h}\big)}{\Gamma\big(x- \N\hat a_h^Y+1\big)\cdot\Gamma\big(x-\N \hat a_h^Y+\theta_{h,h}\big)}.
 \end{equation}

 \medskip

\noindent \textsc{Potentials.} The potential $V_{g}(x)$ gets an additional leading term
\begin{equation}
\label{eq_Vhshift_leading} \frac{2\theta_{h,g}}{\theta_{h,h}}\big(\textnormal{Llog}(x - \hat{a}_h^{\prime\,Y}) - \mathrm{Llog}(x - \hat{a}_h^{\prime\,X})\big).
\end{equation}
For $x$ outside the segment $[\hat a_h^{\prime\,X},\hat a_h^{\prime\,Y}]$, the expression \eqref{eq_Vhshift_leading} can be transformed into
\begin{equation}
\label{Vhsfhit} -\frac{2\theta_{h,g}}{\theta_{h,h}} \int_{\hat a_h^{\prime\,X}}^{\hat a_h^{\prime\,Y}} \log|x-y|\,\dd y.
\end{equation}
By \eqref{eq_Stirling_xy} it matches the leading asymptotic term of the logarithm of \eqref{eq_x178}. For $x$ inside $[\hat a_h^{\prime\,X},\hat a_h^{\prime\,Y}]$, the equality does not apply. But since we only need to deal with $x\in [\min (\hat a_{g}^Y, \hat a_{g}^{\prime\,Y}), \max (\hat b_{h}^Y, \hat b_{h}^{\prime\,Y})]$, this mismatch can be eventually absorbed into $\err_h(x)$ terms in \eqref{eq_weight_form}.

Beyond that, the potential $V_{g}(x)$ gets another subleading term. For $h\neq g$ this term is
\begin{equation}
\label{eq_weight_subleading_1}
\begin{split}
& \quad \frac{2\theta_{h,g}}{\theta_{h,h}} \int_{\hat a_h^{\prime\,X}}^{\hat a_h^{\prime\,Y}} \log|x-y|\,\dd y
 \\ & -\frac{1}{\N}\log\Bigg( \prod_{\begin{smallmatrix}\ell = a_h^X +i \theta_{h,h}: \\ i\in \amsmathbb Z_{\geq 0}, \,\, \ell \leq a_h^Y-\theta_{h,h} \end{smallmatrix}} \frac{1}{\N^{2\theta_{h,g}}} \cdot \frac{\Gamma\big(|\N x-\ell|+1\big)\cdot\Gamma\big(|\N x-\ell|+\theta_{h,g}\big)}{\Gamma\big(|\N x-\ell|\big)\cdot\Gamma\big(|\N x-\ell|+1-\theta_{h,g}\big)}\Bigg).
 \end{split}
\end{equation}
Using \eqref{eq_Stirling_xy}, we see that the last term is decaying as $O\big(\frac{1}{\N}\big)$ as $\N\rightarrow\infty$. Note that the first term in \eqref{eq_weight_subleading_1} cancels with \eqref{Vhsfhit}: we think of \eqref{Vhsfhit} as the leading contribution and \eqref{eq_weight_subleading_1} is a small correction to it. For $h=g$ the subleading term that we add is slightly different, it is
\begin{equation}
\label{eq_weight_subleading_2}
 2\,\mathrm{Llog}(x-\hat a_h^{\prime\,X})-\frac{1}{\N}\log\bigg(\frac{\Gamma\big(\N x- \N\hat a_h^X+1\big)\cdot\Gamma\big(\N x-\N a_h^X+\theta_{h,h}\big)}{\N^{2\N x- 2\N\hat a_h^X+\theta_{h,h}}}\bigg).
\end{equation}
In other words, it accounts only for a half of the factors in \eqref{eq_x179}, because another half is treated separately in our form of the weight \eqref{eq_ansatzw}. The Stirling formula \eqref{eq_Stirling_basic} implies that \eqref{eq_weight_subleading_2} decays as $O\big(\frac{1}{\N}\big)$ as $\N\rightarrow\infty$.

\medskip

\noindent \textsc{Regular part of the potentials.} The updates in $U_g(x)$ can be inferred directly from the ones of $V_g(x)$ according to \eqref{eq_potential_weight_match}.

\medskip

\noindent \textsc{Functions $\Phi^{\pm}_g$.} The ratio of \eqref{eq_x178} at
 $x=\N z + \frac{1}{2}$ and $x=\N z - \frac{1}{2}$ is a telescoping product, which is equal for $g < h$ to
 \begin{equation}
\label{ratmoins} \frac{\Gamma\Big(\frac{\N\hat{a}_h^{Y} - \frac{1}{2} - \N z}{\theta_{h,h}}\Big)}{\Gamma\Big(\frac{\N\hat{a}_h^X - \frac{1}{2} - \N z}{\theta_{h,h}}\Big)}\cdot\frac{\Gamma\Big(\frac{\N\hat{a}_{h}^{X} + \frac{1}{2} - \N z}{\theta_{h,h}}\Big)}{\Gamma\Big(\frac{\N\hat{a}_{h}^{Y} + \frac{1}{2} - \N z}{\theta_{h,h}}\Big)}\cdot\frac{\Gamma\Big(\frac{\N\hat{a}_{h}^{X} + \theta_{h,g} - \frac{1}{2}- \N z}{\theta_{h,h}}\Big)}{\Gamma\Big(\frac{\N\hat{a}_{h}^{Y} + \theta_{h,g} - \frac{1}{2} - \N z }{\theta_{h,h}}\Big)}\cdot\frac{\Gamma\Big(\frac{\N \hat{a}_{h}^{Y} - \theta_{h,g} + \frac{1}{2} - \N z}{\theta_{h,h}}\Big)}{\Gamma\Big(\frac{\N\hat{a}_{h}^{X} - \theta_{h,g} + \frac{1}{2} - \N z}{\theta_{h,h}}\Big)},
\end{equation}
and for $g \geq h$ to
\begin{equation}
\label{ratplus} \frac{\tilde{\Gamma}\Big(\frac{\N z - \frac{1}{2} - \N\hat{a}_{h}^{Y}}{\theta_{h,h}}\Big)}{\tilde{\Gamma}\Big(\frac{\N z - \frac{1}{2} - \N\hat{a}_{h}^{X}}{\theta_{h,h}}\Big)}\cdot \frac{\tilde{\Gamma}\Big(\frac{\N z + \frac{1}{2} - \N\hat{a}_{h}^{X}}{\theta_{h,h}}\Big)}{\tilde{\Gamma}\Big(\frac{\N z + \frac{1}{2} - \N\hat{a}_{h}^{Y}}{\theta_{h,h}}\Big)}\cdot\frac{\tilde{\Gamma}\Big(\frac{\N z - \frac{1}{2} + \theta_{h,g} - \N \hat{a}_{h}^{X}}{\theta_{h,h}}\Big)}{\tilde{\Gamma}\Big(\frac{\N z - \frac{1}{2} + \theta_{h,g} - \N \hat{a}_{h}^{Y}}{\theta_{h,h}}\Big)}\cdot\frac{\tilde{\Gamma}\Big(\frac{\N z + \frac{1}{2} - \theta_{h,g} - \N \hat{a}_{h}^{Y}}{\theta_{h,h}}\Big)}{\tilde{\Gamma}\Big(\frac{\N z + \frac{1}{2} - \theta_{h,g} - \N \hat{a}_{h}^{X}}{\theta_{h,h}}\Big)},
\end{equation}
where $\tilde{\Gamma}(x) = \Gamma(x + 1)$. For $g = h$ this simplifies to
  \begin{equation}
 \label{rateq}  \frac{\big(z- \hat a_h^X+\frac{1}{2\N}\big)\big(z-a_h^X+\frac{\theta_{h,h}}{\N}-\frac{1}{2\N}\big)}{\big(z-\hat a_h^Y+\frac{1}{2\N}\big)\big(z- \hat a_h^Y+\frac{\theta_{h,h}}{\N}-\frac{1}{2\N}\big)}.
  \end{equation}
If $g < h$, then we multiply $\Phi_{g}^{+}(z)$ by \eqref{ratmoins}. If $g > h$, then we multiply $\Phi_{g}^{+}(z)$ by \eqref{ratplus}. If $g = h$, then we multiply $\Phi_{g}^{-}(z)$ by the denominator of \eqref{rateq} and we multiply $\Phi_{g}^{+}(z)$ by the numerator of \eqref{rateq}.

\medskip

\noindent \textsc{Functions $\phi_g^{\pm}$.} Differentiating \eqref{eq_Vhshift_leading}, we see that $e^{-\partial_xV_{g}(x)}$ is multiplied by
\begin{equation}
\label{eq_phi_addon}
 \exp\bigg(-\frac{2\theta_{h,g}}{\theta_{h,h}}\big(-\log|x-\hat a_h^{\prime\,X}|+\log|x-\hat a_h^{\prime\,Y}|\big)\bigg).
\end{equation}
Removing absolute values, for each $g\neq h$ we multiply $\phi^+_{g}(z)$ by
\begin{equation}
\label{eq_x211}
 \exp\Bigg(-\frac{2\theta_{h,g}}{\theta_{h,h}}\log\bigg(\frac{z-\hat a_h^{\prime\,Y}}{z-\hat a_h^{\prime\,X}}\bigg)\Bigg).
\end{equation}
If $g=h$, then \eqref{eq_phi_addon} simplifies to
\begin{equation}
\label{eq_x139}
 \bigg(\frac{x-\hat a_h^{\prime\,X}}{x-\hat a_h^{\prime\,Y}}\bigg)^2.
\end{equation}
In this situation we multiply $\phi^-_{g}(z)$ by $(z-\hat a_{g}^{\prime\,Y})^2$ and we multiply $\phi^+_{g}(z)$ by $(z-\hat a_{g}^{\prime\,X})^2$.

Additionally, we need to multiply $\phi^+_{g}(z)$ by the exponent of the derivative of \eqref{eq_weight_subleading_1} for each $g\neq h$ and by the exponent of \eqref{eq_weight_subleading_2} for $g=h$. These factors do not have zeros or singularities in a small complex neighborhood of the segments of the ensemble $Y$ and behave as $1+O\big(\frac{1}{\N}\big)$.

\medskip

A similar list of updates is carried out to account for the effect of each saturated segment of the form $(\hat b_h^{\prime\,Y},\hat b_h^{\prime\,X}]$. One notable difference is such segments affects the updates of $\Phi_{g}^+(z)$ and $\phi_{g}^+(z)$ for all $g$, rather but $\Phi_{g}^-(z)$ and $\phi_{g}^-(z)$ at this step. Indeed, even for $h=g$, the functions $\Phi^-_{g}(z)$ and $\phi^-_{g}(z)$ are uniquely fixed by \eqref{eq_Phi_plus_def}, \eqref{eq_phi_plus_minus_def} in a way which does not change when we account for the segment $(\hat b_h^{\prime\,Y},\hat b_h^{\prime\,X}]$ as we obtain the ensemble $Y$ from the ensemble $X$. As a result, the updated $\Phi_{g}^{+}(z)$ gets zeros at $z=\hat{b}_{g}^Y + \frac{1}{2\N}$ and at $z = \hat{b}^{\prime\,Y}_{g}=\hat b_{g}^Y + \frac{1}{\N}\big(\theta_{g,g}-\frac{1}{2}\big)$, while the updated $\phi_{g}^{+}(z)$ gets a double zero at $z=\hat{b}^{\prime\,Y}_{g}$.

\smallskip

In addition we need to apply a second list of updates for each $g \in [H]$, no matter whether the segments $[\hat a_{g}^{\prime\,X},\hat a_h^{\prime\,Y})$ and $(\hat b_h^{\prime\,Y},\hat b_h^{\prime\,X}]$ are saturated or void.

\medskip

\noindent \textsc{Potentials.} Our treatment of Gamma functions in the weight depends on whether the argument of these Gamma functions vanishes near the endpoints of the segments or not. According to Assumption~\ref{Assumptions_analyticity} we absorb them into $U_h(x)$ in the latter case and deal with them separately in the former case. Due to this feature, for each $g \in [H]$ we need to additionally add to $V_{g}(x)$ the term:
\begin{equation}
\begin{split}
\label{eq_weight_subleading_3}
& - \frac{1}{\N}\log\Bigg( \prod_{j=1}^{\iota_h^-} \frac{\N^{\N x - \N \hat a_h^X + \rho^-_{h,j} - \frac{1}{2}}}{\Gamma\big(\N x - \N \hat a_h^X + \rho^-_{h,j}\big)} \cdot \prod_{j=1}^{\iota_h^+} \frac{\N^{\N \hat b_h^X-\N x + \rho^+_{h,j} - \frac{1}{2}}}{\Gamma\big(\N \hat b_h^X - \N x + \rho^+_{h,j}\big)}\Bigg)
\\
& -\iota_h^{-}\, \mathrm{Llog}(x - \hat{a}^{\prime\,X}_{h}) - \iota_{h}^+\, \mathrm{Llog}(\hat{b}^{\prime\,X}_h - x).
\end{split}
\end{equation}
Note that the Stirling formula \eqref{eq_Stirling_basic} implies that \eqref{eq_weight_subleading_3} decays as $O\big(\frac{1}{\N}\big)$ as $\N\rightarrow\infty$.

\medskip

\noindent \textsc{Functions $\Phi^{\pm}$.} For each $g\in[H]$, we need to compensate the linear factors in \eqref{eq_Phi_plus_def} for $\Phi^-_{g}(z)$ and in \eqref{eq_phi_plus_minus_def} for $\phi^-_{g}(z)$ of the ensemble $X$. For that we multiply \emph{both} $\Phi_{g}^{+}(z)$ and $\Phi_{g}^{-}(z)$ by
\begin{equation}
 \label{eq_x5}
\prod_{j=1}^{\iota_{g}^{-X}} \bigg(x-\hat a_{g}^X -\frac{1}{2\N}+ \frac{\rho^{-,X}_{g,j}}{\N}\bigg)^{-1},
\end{equation}
and we multiply \emph{both} $\phi^+_h(z)$ and $\phi^-_h(z)$ by
\begin{equation}
\label{eq_x21}
\big(z-\hat a_{g}^{\prime\,X}\big)^{-\iota_{g}^{-,X}}.
\end{equation}
There is a certain asymmetry caused by the choices made in Definition~\ref{Definition_phi_functions}. This additional step does not involve the right endpoints $\hat b_{g}^{X}$ or their shifted version $\hat b_{g}^{\prime\,X}$.

After all these updates have been carried out, we obtain the data defining the discrete ensemble $Y$.

Let us now check that all the desired assumptions are satisfied for the ensemble $Y$. Checking that it satisfies Assumptions~\ref{Assumptions_Theta}, \ref{Assumptions_basic} and \ref{Assumptions_analyticity} is straightforward, as all the conditions are implied by similar conditions for the ensemble $X$, and we omit this part. Several parts of Assumption~\ref{Assumptions_extra} are also straightforward to check. The Conditions 1.2.3. are immediate from the construction. For Condition 4., note that we multiply the weight by new $\Gamma$-factors only in the saturated situation and, therefore, for the endpoints in the void, the weight in the ensemble $Y$ does not have the corresponding $\Gamma$-part in \eqref{eq_ansatzw}, \textit{i.e.} we absorb everything into $U_h$. For Condition 5., the required vanishing is clear from \eqref{rateq} and \eqref{eq_x139}: $\Phi^-_h(z)$ and $\phi^-_h(z)$ are multiplied by the respective inverses of these two factors; similarly for the right endpoints and the corresponding vanishing of $\Phi^+_h(z)$ and $\phi^+_h(z)$.

It remains to check Assumption~\ref{Assumptions_offcrit} and Condition 6. of Assumption~\ref{Assumptions_extra}. Both of them are properties of the equilibrium measure. In order to establish them, we introduce an auxiliary variational datum $\widetilde Y$ by taking the variational datum of the ensemble $Y$ and removing the $O(\frac{1}{\N})$ terms in the potential given by \eqref{eq_weight_subleading_1}, \eqref{eq_weight_subleading_2}, and \eqref{eq_weight_subleading_3}. Note that the variational problem for the equilibrium measure for $\widetilde Y$ matches the one for the ensemble $X$ when restricted to the shifted and rescaled segments of $Y$ and absorbing the interaction with the removed saturated segments of the ensemble $X$ into the potential of $\widetilde Y$. Hence, the equilibrium measure of $\widetilde Y$ differs from the equilibrium measure of the ensemble $X$ by removal of saturated segments of the form $[\hat a_h^{\prime\,X},\hat a_h^{\prime\,Y})$ or $(\hat b_h^{\prime\,Y},\hat b_h^{\prime\,X}]$. The equilibrium measure of $\widetilde Y$ then satisfies Assumption~\ref{Assumptions_offcrit} because this assumption holds for the ensemble $X$. For Condition 6. in Assumption~\ref{Assumptions_extra}, we claim that for any $g \in [H]$, the function $s_{g}^{\widetilde Y}(z)$ differs from $s_{g}^X(z)$ by factors bounded and bounded away from $0$ in a complex neighborhood of $[\hat a_{g}^Y, \hat b_{g}^Y]$. Recall from Definition~\ref{GQdef2}
\begin{equation} \label{eq_s_copy}
s_{g}(z) = \frac{\phi_{g}^{+}(z)\cdot \exp\big(\sum_{h = 1}^H \theta_{g,h}\,\mathcal{G}_{\mu_{h}}(z)\big) - \phi_{g}^{-}(z)\cdot \exp\big(-\sum_{h = 1}^{H} \theta_{g,h}\,\mathcal{G}_{\mu_{h}}(z)\big)}{\sqrt{(z - \alpha_g)(z - \beta_g)} \cdot (z-\hat a'_{g})^{\mathbbm{1}_{\amsmathbb{S}_g}(\hat{a}'_{g})} \cdot (z-\hat b'_{g})^{\mathbbm{1}_{\amsmathbb{S}_g}(\hat{b}'_{g})}},
\end{equation}
where we used the fact that in the present situation, each segment contains a single band.
 When we pass from $X$ to $\widetilde Y$, for each $h$ so that $[\hat a_{h}^{\prime\,X},\hat a_{h}^{\prime\,Y}]$ is saturated, the first term of the second line of \eqref{eq_s_copy} is multiplied by
\begin{equation}
\label{eq_x212}
 \exp\bigg(-\theta_{g,h} \int_{a_{h}^{\prime\,X}}^{a_{h}^{\prime\,Y}} \frac{\dd x}{\theta_{h,h}(z-x)}\bigg)=
 \exp\Bigg(\frac{\theta_{g,h}}{\theta_{h,h}} \log\bigg( \frac{z-a_{h}^{\prime\,Y}}{z-a_{h}^{\prime\,X}}\bigg)\Bigg),
\end{equation}
while the second term is multiplied by the inverse of the same factor. If $h\neq g$, then $\phi^+_{g}(z)$ is simultaneously multiplied by \eqref{eq_x211}, while $\phi^-_{g}(z)$ is not multiplied by anything; hence, both terms in the second line of \eqref{eq_s_copy} are multiplied by the inverse of \eqref{eq_x212}, which is clearly bounded and bounded away from $0$. If $h=g$, then we instead multiply $\phi^-_{g}(z)$ by $(z-\hat a_{g}^{\prime\,Y})^2$ and we multiply $\phi^+_{g}(z)$ by $(z-\hat a_{g}^{\prime\,X})^2$; in this situation both terms in the second line of \eqref{eq_s_copy} are multiplied by $(z-a_{h}^{\prime\,Y})(z-a_{h}^{\prime\,X})$. The first factor in the last product cancels with the linear factor in the first line of \eqref{eq_s_copy} for $\widetilde Y$, while the second factor is bounded and bounded away from $0$. Similar arguments apply for the contribution of the saturated segments $[\hat b_{h}^{\prime\,X},\hat b_{h}^{\prime\,Y}]$. Finally, the multiplication of both $\phi^{\pm}_{g}(z)$ by \eqref{eq_x21} again gives a factor bounded and bounded away from $0$ and the claim is proven.

Once the claim is established, we can use Proposition~\ref{Proposition_q_bounds}: it says that for the ensemble $X$ the $H$-tuple of functions $(s_{h}^X(z))_{h = 1}^H$ is bounded away from $0$ in a small neighborhood of the bands. Hence, by choosing $\eps$ to be small enough, we can guarantee that for the auxiliary variational datum $\widetilde Y$ the functions $s_h^{\widetilde Y}(z)$ are bounded away from $0$ in a neighborhood of $[\hat a_{h}^{Y},\hat b_{h}^{Y}]$ for any $h \in [H]$.

\smallskip

We have checked Assumption~\ref{Assumptions_offcrit} and Condition 6. in Assumption~\ref{Assumptions_extra} for $\widetilde Y$. The difference between the potential for $\widetilde Y$ and for $Y$ vanishes as $\N\rightarrow\infty$ (we will give a more precise estimate in Corollary~\ref{mueqYYtilde}, but it is not needed to conclude with the present argument). Using Lemma~\ref{Lemma_continuity_potential} and Theorem~\ref{Theorem_off_critical_neighborhood} we can deduce that Assumption~\ref{Assumptions_offcrit} and Condition 6. of Assumption~\ref{Assumptions_extra} are also valid for the ensemble $Y$ as long as $\N$ is large enough.
\end{proof}

It is useful to record the change of the partition function between the ensemble $X$ and its conditioning/localization $Y$. We have for a $C>0$
\begin{equation}
\label{eq_changeofZ}
\begin{split}
\frac{\Z_{\N}^{X}}{\Z_\N^Y} & = \prod_{i \in \amsmathbb{F}} w^{X}(x_i)\, \cdot \prod_{\substack{i < j \\ i,j \in \amsmathbb{F}}} \frac{1}{\N^{2\theta_{h(i),h(j)}}}\cdot \frac{\Gamma\big(x_j - x_i + 1\big)\cdot\Gamma\big(x_j - x_i + \theta_{h(i),h(j)}\big)}{\Gamma\big(x_j - x_i\big)\cdot\Gamma\big(x_j - x_i + 1 - \theta_{h(i),h(j)}\big)} \\
& \quad \times \Bigg[1 + O\bigg(\exp\bigg(-\frac{\N}{C}\bigg)\bigg)\Bigg],
\end{split}
\end{equation}
where $\amsmathbb{F} \subseteq [N]$ is the set of indices corresponding to frozen particles which are present in the ensemble $X$ but not in the ensemble $Y$, and $x_i$ is the position of the $i$-th particle for $i \in \amsmathbb{F}$. We can also compare the energy functionals of the two ensembles.
\begin{lemma}
\label{Lemma_energyconditioner}Under the conditions of Theorem~\ref{proposition_FFF_conditioning}, let $\boldsymbol{\mu}^{X}$ and $-\mathcal{I}^{X}$ (respectively $\boldsymbol{\mu}^{Y}$ and $-\mathcal{I}^{Y}$) be the equilibrium measure and energy functional of the ensemble $X$ (respectively $Y$). We have as $\N \rightarrow \infty$
\begin{equation*}
\begin{split}
\mathcal{I}^{X}[\boldsymbol{\mu}^{X}] & \,\,\mathop{=}_{\N \rightarrow \infty}\,\, \mathcal{I}^{Y}[\boldsymbol{\mu}^{Y}] + \sum_{g,h = 1}^{H} \theta_{g,h} \iint_{(\amsmathbb{A}^{X}_g \setminus \amsmathbb{A}^{Y}_g)\times (\amsmathbb{A}^{X}_{h} \setminus \amsmathbb{A}^{Y}_{h})} \log|x - y|\mu_g^{X}(x)\mu_{h}^{X}(y)\dd x\dd y \\
& \quad \,\,\,\,- \sum_{h = 1}^{H} \int_{\amsmathbb{A}^{X}_h\setminus \amsmathbb{A}_h^{Y}} V^{X}_{h}(x)\mu_{h}^{X}(x)\dd x + \frac{\mathbbm{Rest}_1^{\mathcal{I}}}{\N} + o\bigg(\frac{1}{\N^2}\bigg),
\end{split}
\end{equation*}
where $\amsmathbb{A}^{X}_h := [\hat a_h^{\prime\,X}, \hat b_h^{\prime\,X}]$ and $\amsmathbb{A}^{Y}_h := [\hat a_h^{\prime\,Y}, \hat b_h^{\prime\,Y}]$, respectively. The term $\mathbbm{Rest}_1^{\mathcal{I}}$ depends on $\N$, but is uniformly bounded. The remainder $o\big(\frac{1}{\N^2}\big)$ is uniformly small for the ensembles satisfying Assumptions~\ref{Assumptions_Theta}--\ref{Assumptions_analyticity} with the same choice of constants. Moreover, if we allow the filling fractions and endpoints of the segments to vary, as in Theorem~\ref{Theorem_differentiability_full}, then $\mathbbm{Rest}_1^{\mathcal{I}}$ is twice-differentiable in these parameters with bounded second derivatives. \label{RestI:index}
\end{lemma}
\begin{proof} Recall the auxiliary variational datum $\widetilde Y$ from the proof of Theorem~\ref{proposition_FFF_conditioning}. The desired statement would follow from two identities
\begin{equation*}
\begin{split}
 \mathcal{I}^{X}[\boldsymbol{\mu}^{X}] & = \mathcal{I}^{\widetilde Y}[\boldsymbol{\mu}^{\widetilde Y}] + \sum_{g,h = 1}^{H} \theta_{g,h} \iint_{(\amsmathbb{A}^{X}_g \setminus \amsmathbb{A}^{Y}_g)\times (\amsmathbb{A}^{X}_{h} \setminus \amsmathbb{A}^{Y}_{h})} \log|x - y|\mu_g^{X}(x)\mu_{h}^{X}(y)\dd x\dd y \\
 & \quad - \sum_{h = 1}^{H} \int_{\amsmathbb{A}^{X}_h\setminus \amsmathbb{A}_h^{Y}} V^{X}_{h}(x)\mu_{h}^{X}(x)\dd x \\
 \mathcal{I}^{\widetilde Y}[\boldsymbol{\mu}^{\widetilde Y}] & \mathop{=}_{\N \rightarrow \infty} \mathcal{I}^{Y}[\boldsymbol{\mu}^{Y}] + \frac{\mathbbm{Rest}_1^{\mathcal{I}}}{\N} +o\bigg(\frac{1}{\N^{2}}\bigg).
 \end{split}
\end{equation*}

The second identity comes from the construction of the ensemble $\widetilde Y$ and Proposition~\ref{proposition_Energy_series} in which we put $\epsilon=\tfrac{1}{\N}$ and combine $\epsilon$ and $\epsilon^2$ terms into a single one. For the first identity we observe on \eqref{Vhsfhit} that the potentials of the ensemble $\widetilde Y$ and of the ensemble $X$ are related by
\begin{equation}
\label{VysandVx}
\forall h \in [H]\qquad V^{X}_{h}(x) = V^{\widetilde Y}_{h}(x) - \sum_{g = 1}^H 2\theta_{h,g} \int_{\amsmathbb{A}^{X}_{g} \setminus \amsmathbb{A}^Y_{g}} \log|x - y|\mu_{g}^{X}(y) \dd y,
\end{equation}
while the equilibrium measure $\mu^{\widetilde Y}$ is the restriction of the equilibrium measure $\mu^{X}$ to $\amsmathbb{A}^Y$. The claimed identity is therefore obtained by inserting the decomposition
\[
\forall h \in [H]\qquad \mu^{X}_h = \mu^{X}_h|_{\amsmathbb{A}^{X}_h \setminus \amsmathbb{A}^{Y}_h} + \mu^{\widetilde Y}_h
\]
into $\mathcal{I}^{X}[\boldsymbol{\mu}^{X}]$ which is quadratic in the equilibrium measure, and recombining the terms to isolate $\mathcal{I}^{\widetilde Y}[\boldsymbol{\mu}^{\widetilde Y}]$.
\end{proof}

\section{Example with saturation: Gaussian weights}
\label{Sec:condgaus}
Our next aim is to proceed with the conditioning procedure for ensembles which have more than one band per
segment $[\hat a_h, \hat b_h]$ or fluctuating filling fractions. Before presenting the most general case in
 the next section, we describe the procedure in the case of the discrete ensemble with Gaussian weight introduce in Section~\ref{dis_Gauss_intro}. Recall it has $H = 1$ segment $(a_1,b_1) = (-\infty,+\infty)$ and weight $w(x) = e^{-\kappa \frac{x^2}{\N}}$. Recall also that $\ell_1\in \amsmathbb Z -\theta\frac{N-1}{2}$ and $\ell_N\in\amsmathbb Z +\theta\frac{N-1}{2}$. We set the large parameter $\N$ equal to the number $N$ of particles but we keep using the letters $\N$ and $N$ for their natural meaning to facilitate the comparison with the forthcoming results. Proposition~\ref{Prop_Gaussian_LLN} computed the equilibrium measure $\mu$ of the ensemble and shows that for $\sqrt{2\theta}\kappa > \pi$, the ensemble has two bands $(-\alpha,-\beta) \cup (\beta,\alpha)$ and a saturation $[-\beta,\beta]$ separating them.

Our aim is two-fold. First, we want to define the filling fractions for each of the two bands, which is
complicated by the presence of the saturation. Second, we want to restrict the ensemble to a small neighborhood of the bands. The latter is not quite necessary for the Gaussian case due to the simplicity of the weight $w(x)$. However, in the general situation, this is the only way we have found to guarantee that the auxiliary Assumption~\ref{Assumptions_extra} holds.

\begin{figure}[t]
\center \scalebox{1.3}{\includegraphics{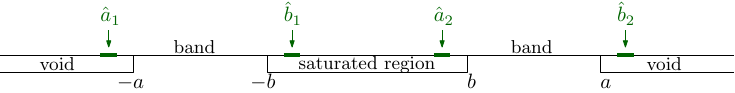}} \caption{Configuration of points appearing in Proposition~\ref{proposition_Gaussian_conditioning} for the discrete ensemble with Gaussian weight. \label{Fig_Gauss_bands}}
\end{figure}

\begin{proposition}\label{proposition_Gaussian_conditioning}
 Let $\kappa,c > 0$ and suppose $\sqrt{2\theta \kappa} - \pi > c$ and $\min(\theta,\kappa) > c$. There exists $\eps_1,\eps_2>0$ small enough and $C>0$ large enough depending only on $c$ and satisfying the following conditions for all $\N$ large enough. Take any six parameters $-\infty < \mathfrak{a}_1<\mathfrak{b}_1< \mathfrak{a}_2< \mathfrak{b}_2 < +\infty$ and $\mathfrak{n}_1,\mathfrak{n}_2 \in \frac{1}{\N}\amsmathbb{Z}_{\geq 0}$, such that:
 \begin{enumerate}
  \item $\N \mathfrak{b}_1 -\N\mathfrak{a}_1 - \theta (\N\mathfrak{n}_1-1)\in \amsmathbb Z_{\geq 0}$ and $\,\N \mathfrak{b}_2 -\N\mathfrak{a}_2 - \theta (\N\mathfrak{n}_2-1)\in \amsmathbb Z_{\geq 0}$,
  \item $\N\mathfrak{a}_1\in\amsmathbb Z-\theta\frac{N-1}{2}$ and $\,\N\mathfrak{b}_2 \in\amsmathbb Z+\theta\frac{N-1}{2}$,
  \item $\mathfrak{a}_1 \in (-\alpha-2\eps_1,-\alpha-\eps_1)$,\, $\mathfrak{b}_1 \in (-\beta+\eps_1,-\beta +2\eps_1)$,\, $\mathfrak{a}_2 \in (\beta-2\eps_1,\beta-\eps_1)$,\, $\mathfrak{b}_2 \in (\alpha+\eps_1,\alpha+2\eps_1)$,
  \item $\big|\mathfrak{n}_1-\mu([\mathfrak{a}_1,\mathfrak{b}_1])\big|<\eps_2$ and $\big|\mathfrak{n}_2-\mu([\mathfrak{a}_2,\mathfrak{b}_2])\big|<\eps_2$,
  \item $\N\mathfrak{n}_1+ \frac{\N}{\theta}(\mathfrak{a}_2-\mathfrak{b}_1)-1+ \N\mathfrak{n}_2=N$.
 \end{enumerate}
Let us define the event $\mathcal{A}$ by the following two properties:
 \begin{enumerate}
 \item[A.] There are no particles outside the interval $[\N \mathfrak{a}_1,\N\mathfrak{b}_2]$.
 \item[B.] There are no holes in the interval $[\N\mathfrak{b}_1,\N\mathfrak{a}_2]$, \textit{i.e.} it is occupied by $\frac{\N}{\theta}(\mathfrak{a}_2-\mathfrak{b}_1)+1$ particles densely packed at distances $\theta$ from each other.
 \end{enumerate}
Then, $\mathcal{A}$ has probability at least $1 - C\exp(-\frac{\N}{C})$, and the discrete ensemble with Gaussian weight and conditioned on $\mathcal{A}$ and on the event that the $\N\mathfrak{n}_1$-th particle sits at position $\N \mathfrak{b}_1$ --- given $\mathcal A$, this implies as well that the $(\N \mathfrak{n}_1+ \frac{\N}{\theta} (\mathfrak{a}_2-\mathfrak{b}_1)\big)$-th particle sits at position $\N\mathfrak{a}_2$ --- and localized to $\mathfrak{A} := [\mathfrak{a}_1,\mathfrak{b}_1] \cup [\mathfrak{a}_2,\mathfrak{b}_2]$ yields an ensemble satisfying Assumptions~\ref{Assumptions_Theta}, \ref{Assumptions_basic}, \ref{Assumptions_offcrit}, \ref{Assumptions_analyticity} and the additional Assumption~\ref{Assumptions_extra}.
\end{proposition}
The configurations of the points of the real line is shown in Figure~\ref{Fig_Gauss_bands}. In words, Conditions 1.\ and 2.\ on the five parameters are integrality conditions dictated by the lattices on which the particles live; Condition 3.\ says that we localize the ensemble in small enough neighborhoods of the bands; Condition 4.\ means that the filling fraction $\mathfrak{n}_1$ for the segment $[\mathfrak{a}_1,\mathfrak{b}_1]$ is close to the optimal one predicted by the equilibrium measure $\mu$; Condition 5.\ guarantees that the filling fractions $(\mathfrak{n}_1,\mathfrak{n}_2)$ are compatible with the total number $N$ of particles in the initial ensemble. The properties A.\ and B.\ mean that the particles outside $\mathfrak{A}$ sit at deterministic locations and, therefore, we can absorb their interactions with the particles in $\mathfrak{A}$ into a modification of the weight, potential, \textit{etc}. of the initial ensemble, as we did in the proof of Theorem~\ref{proposition_FFF_conditioning}.

The first step of the proof of Proposition~\ref{proposition_Gaussian_conditioning} is to check the off-criticality Assumption~\ref{Assumptions_offcrit} so that we can use the large deviations estimates of Chapter~\ref{Chapterlarge}.

\begin{lemma}
\label{Lemma_Discr_Gauss_Off_critical}
Let $\kappa,c > 0$ and assume $|\sqrt{2\theta \kappa} - \pi| \geq c$ and $\min(\theta,\kappa) \geq c$. The equilibrium measure $\mu$ of Proposition~\ref{Prop_Gaussian_LLN} for the discrete ensemble with Gaussian weight is off-critical, \textit{i.e.} $\mu$ satisfies Assumption~\ref{Assumptions_offcrit} for constants that only depend on $c$.
\end{lemma}
\begin{proof} If $\sqrt{2\theta \kappa} < \pi$, then $\mu$ is a semi-circle law and it is clearly off-critical. If $\sqrt{2\theta \kappa} - \pi > c$, then $\mu$ has two bands separated by a saturation, and the elliptic modulus $\mathsf{k}$ appearing in the formula for $\mu$ is larger than a positive constant depending only on $c$. Hence, the length $2\beta$ of the saturation is bounded from below by a positive constant depending only on $c$. This proves Condition 1.\ in Assumption~\ref{Assumptions_offcrit}. Condition 3.\ on the behavior of the effective potential in the saturation is implied by integrating the inequality \eqref{ineein0} appearing in the proof of Proposition~\ref{Prop_Gaussian_LLN} from $\beta$ to $x \in [0,\beta)$. Checking the analogous Condition 2.\ in voids is left to the reader. The remaining conditions 4., 5. and 6. can be checked by explicit computations with the formula \eqref{insidesat} for $\mu$, which we omit.
\end{proof}

\begin{proof}[Proof of Proposition~\ref{proposition_Gaussian_conditioning}]
Take $\eps_{1}>0$ and \emph{arbitrary} finite $\mathfrak a_1^0<\mathfrak b_1^0<\mathfrak a_2^0<\mathfrak b_2^0$ such that
\begin{equation*}
\begin{split}
 [-\alpha-\eps_{1},-\beta +\eps_{1}] & \subset(\mathfrak a_1^0,\mathfrak b_1^0)\subset [-\alpha-2\eps_{1},-\beta +2\eps_{1}], \\
[\beta-\eps_{1},\alpha+\eps_{1}] & \subset(\mathfrak a_2^0, \mathfrak b_2^0) \subset [\beta-2\eps_{1},\alpha+2\eps_{1}].
\end{split}
\end{equation*}
Therefore, $(-\infty,\mathfrak a_1^0]$ and $[\mathfrak b_2^0,+\infty)$ are void for the equilibrium measure $\mu$, while $[\mathfrak b_1^0,\mathfrak a_2^0]$ is saturated. Note that we are not imposing any integrality conditions on $\mathfrak a_1^0$, $\mathfrak b_1^0$, $\mathfrak a_2^0$, $\mathfrak b_2^0$.

We now consider an auxiliary variational datum --- following the terminology of Chapter~\ref{Chapter_smoothness} --- corresponding to the conditioned discrete ensemble with Gaussian weight. We deal with $H = 2$ segments $[\mathfrak a_1^0,\mathfrak b_1^0]$ and $[\mathfrak a_2^0,\mathfrak b_2^0]$ in this order on the real line, the $2\times 2$ matrix $\boldsymbol{\Theta}$ with all matrix elements equal to $\theta$ and the potentials
\begin{equation}
\label{eq_x87}
 V_1(x)=V_2(x)= \kappa x^2 -2 \int_{\mathfrak b_1^0}^{\mathfrak a_2^0} \log|x-y| \dd y.
\end{equation}
We set the filling fractions to be $\mathfrak n_1^0=\mu([\mathfrak a_1^0,\mathfrak b_1^0])$ and $\mathfrak n_2^0=\mu([\mathfrak a_2^0,\mathfrak b_2^0])$. Then by Theorem~\ref{Theorem_equi_charact_repeat_2} the equilibrium measure for this auxiliary variational datum is the same as the restriction of $\mu$ to $[\mathfrak a_1^0,\mathfrak b_1^0] \cup [\mathfrak a_2^0,\mathfrak b_2^0]$. Hence, it satisfies Assumptions~\ref{Assumption_A}, \ref{Assumption_B} and \ref{Assumption_C} of Chapter~\ref{Chapter_smoothness}. In particular, the measure is off-critical. By Proposition~\ref{Proposition_differentiability_filling_fraction}, the auxiliary equilibrium measure remains off-critical for the filling fractions in a neighborhood of $\mathfrak n_1^0$, $\mathfrak n_2^0$.

The auxiliary equilibrium measure also satisfies the lower bounds on $s_1(z)$, $s_2(z)$ of Proposition~\ref{Proposition_q_bounds}. Decreasing $\eps_1$, and observing as in the proof of Theorem~\ref{proposition_FFF_conditioning} that $s_1(z)$ and $s_2(z)$ are multiplied by bounded factors in this procedure, we can assume that the lower bounds hold in small complex neighborhoods of the intervals $[-\alpha-2\eps_1, -\beta+2\eps_1]$ and $[\beta-2\eps_1, \alpha+2\eps_1]$. At the same time, according to Theorem~\ref{Theorem_differentiability_full}, the Stieltjes transform $\Gm_{\mu}(z)$ for $z \in \amsmathbb{C} \setminus \big([\mathfrak{a}_1^0,\mathfrak{b}_1^0] \cup [\mathfrak{a}_2^0,\mathfrak{b}_2^0]\big)$ and the endpoints of the bands depend smoothly on the above segment filling fractions. This is also true for the functions $s_1(z)$ and $s_2(z)$, so that the lower bounds on $s_1(z)$ and $s_2(z)$ remain valid for filling fractions in a neighborhood of $\mathfrak n_1^0$, $\mathfrak n_2^0$. As a consequence, we conclude that there exist small $\eps_1,\eps_2>0$ such that for the auxiliary variational datum on the segments $[\mathfrak a_1^0,\mathfrak b_1^0]$ and $[\mathfrak a_2^0,\mathfrak b_2^0]$, potential given by \eqref{eq_x87} and segment filling fractions $\mathfrak{n}_1,\mathfrak{n}_2$ satisfying
\[
\big|\mathfrak n_1-\mu([\mathfrak{a}_1^0,\mathfrak{b}_1^0])\big|<\eps_2 \qquad \textnormal{and}\qquad \big|\mathfrak n_2-\mu([\mathfrak{a}_2^0,\mathfrak{b}_2^0])\big|<\eps_2,
\]
the equilibrium measure have two bands included in $(-\alpha-\eps_1, -\beta+\eps_1)$ and $(\beta-\eps_1, \alpha+\eps_1)$, respectively and the corresponding functions $s_1(z)$, $s_2(z)$ are uniformly bounded away from $0$ in $[-\alpha-2\eps_1, -\beta+2\eps_1]$ and $[\beta-2\eps_1, \alpha+2\eps_1]$, respectively.

We claim that such $\eps_1$, $\eps_2$ satisfy the conditions of Proposition~\ref{proposition_Gaussian_conditioning}.
Indeed, let us fix $\mathfrak n_1,\mathfrak n_2$ as above and pick points $\mathfrak a_1,\mathfrak b_1,\mathfrak a_2, \mathfrak b_2$ satisfying Conditions 1.\ to 5.\ and so that
\[
\mathfrak{a}_1 \in (\mathfrak{a}_1^0,-\alpha - \eps_1),\quad \mathfrak{b}_1 \in (-\beta + \eps_1,\mathfrak{b}_1^0),\quad \mathfrak{a}_2 \in (\mathfrak{a}_2^0,\beta - \eps_1),\quad \mathfrak{b}_2 \in (\alpha + \eps_1,\mathfrak{b}_2^0).
\]
We define shifted points
\begin{equation*}
\forall h \in \{1,2\}\qquad {\mathfrak a}'_h= {\mathfrak a}_h - \frac{\theta - \frac{1}{2}}{\N} \quad \textnormal{and} \quad {\mathfrak b}'_h= {\mathfrak b}_h + \frac{\theta-\frac{1}{2}}{\N}.
\end{equation*}

The segment filling fractions of the conditioned ensemble with Gaussian weight --- with conditioning on the event $\mathcal A$ and on the position of the $\N \mathfrak{b}_1$-th particle --- localized to the intervals $[\mathfrak a_1,\mathfrak b_1]$ and $[\mathfrak a_2,\mathfrak b_2]$ are readily translated into the segment filling fractions of the auxiliary equilibrium measure by an affine transformation, since on each of the four intervals $[\mathfrak a_1^0, \mathfrak a'_1]$, $[\mathfrak b'_1,\mathfrak b_1^0]$, $[\mathfrak a_2^0,\mathfrak a'_2]$, $[\mathfrak b'_2,\mathfrak b_2^0]$ the equilibrium measure is either void or saturated. Simultaneously, outside $[\N \mathfrak a_1,\N \mathfrak b_2] $ as $\N$ becomes large the discrete ensemble with Gaussian weight has no particles with overwhelming probability by Theorem~\ref{Theorem_ldpsup} and in $[\N \mathfrak b_1, \N\mathfrak a_2]$ there are no holes with overwhelming probability by Theorem~\ref{Theorem_ldsaturated}. This shows that the event $\mathcal A$ has overwhelming probability as $\N \rightarrow \infty$.

Next, we should check that the conditioned ensemble satisfies Assumptions~\ref{Assumptions_Theta}--\ref{Assumptions_analyticity} and \ref{Assumptions_extra}. This is done using exactly the same approach as in the proof of Theorem~\ref{proposition_FFF_conditioning}: we perform the same transformations of the weights and potentials and introduce an additional variational datum $\widetilde Y$ on the intervals $[\mathfrak a'_1,\mathfrak b'_1]$, $[\mathfrak a'_2,\mathfrak b'_2]$ by ignoring the subleading corrections as in \eqref{eq_weight_subleading_1}, \eqref{eq_weight_subleading_2}, and \eqref{eq_weight_subleading_3}.
The equilibrium measure of $\widetilde Y$ is the restriction to $[\mathfrak a'_1,\mathfrak b'_1]\cup[\mathfrak a'_2,\mathfrak b'_2]$ of the auxiliary equilibrium measure on $[\mathfrak a_1^0,\mathfrak b_1^0] \cup [\mathfrak a_2^0,\mathfrak b_2^0]$ discussed at the beginning of the present proof. This implies the off-criticality of the equilibrium measure of $\widetilde Y$ and the lower bounds on $s_1(z)$ and $s_2(z)$, which follow from similar lower bounds for the auxiliary equilibrium measure (as in the verification of Condition 6. of Assumption~\ref{Assumptions_extra} in Theorem~\ref{proposition_FFF_conditioning}). Using Lemma~\ref{Lemma_continuity_potential} and
Theorem~\ref{Theorem_off_critical_neighborhood} we conclude that the same properties hold for the conditioned and localized ensemble with Gaussian weight. Checking the remaining parts of Assumptions~\ref{Assumptions_Theta}--\ref{Assumptions_analyticity} and \ref{Assumptions_extra} is routine and we omit it.
\end{proof}

\section{General case with fluctuating filling fractions}
\label{fluctuatnecmergitur}
In this section we present a generalization of the conditioning and localization procedure that we used for Gaussian weights, based on the same principles, for discrete ensembles where the constraints ($\star$) do not uniquely fix the segment filling fractions. Our strategy is to restrict the ensemble to a small neighborhood of bands, because in such neighborhoods Assumption~\ref{Assumptions_extra} is satisfied. One technical difficulty is to carefully deal with the various integrality conditions relative to the lattices on which the particles are allowed to vary.

We recall from Definition~\ref{bandlabel} that bands of an equilibrium measure are denoted $(\alpha_k,\beta_k)$ for $k \in [K]$, that $h^k \in [H]$ is the index of the segment in which this band appears, and that the set of indices of bands in the $h$-th segment is $\llbracket k^-(h),k^+(h)\rrbracket$.

\begin{theorem}\label{proposition_fluct_conditioning}
Consider a discrete ensemble $X$ satisfying Assumptions~\ref{Assumptions_Theta}, \ref{Assumptions_basic}, \ref{Assumptions_offcrit} and \ref{Assumptions_analyticity}, with segments $[\hat{a}_h,\hat{b}_h]$ indexed by $h \in [H]$ and equilibrium measure $\boldsymbol{\mu}$. There exists two small numbers $\eps_1,\eps_2>0$ and a large number $C>0$, all depending only on the constants in the assumptions and satisfying the following conditions for $\N$ large enough.

Consider a $3K$-tuple of real parameters $(\mathfrak{a}_k$, $\mathfrak{b}_k$, $\mathfrak{n}_k)_{k = 1}^{K}$ such that $\mathfrak a_1 < \mathfrak{b}_1<\mathfrak a_2<\dots< \mathfrak a_k<\mathfrak b_k$ and $[\mathfrak{a}_{k},\mathfrak{b}_{k}] \subset [\hat{a}_{h^k}',\hat{b}_{h^k}']$ for any $k \in [K]$. Define the shifted parameters
 \[
\forall k \in [K]\qquad {\mathfrak a}'_k= {\mathfrak a}_k - \frac{\theta_{h^k,h^k}-\frac{1}{2}}{\N}\quad \textnormal{and}\quad {\mathfrak b}'_k= {\mathfrak b}_k + \frac{\theta_{h^k,h^k}-\frac{1}{2}}{\N}.
 \]
Assume the following properties.
 \begin{enumerate}
  \item $\forall k\in[K] \quad \N \mathfrak{n}_{k}\in\amsmathbb Z_{\geq 0} \quad \textnormal{and}\quad \N \mathfrak{b}_k -\N \mathfrak{a}_k - \theta_{h^k,h^k} (\N\mathfrak n_k-1)\in \amsmathbb Z_{\geq 0}$.

  \item For any $h\in[H]$
    \begin{itemize}
    \item if $[\hat a'_h,\mathfrak a'_{k^-(h)}]$ is void, then $\N(\mathfrak a_{k^-(h)}-\hat a_h)\in\amsmathbb Z_{\geq 0}$, otherwise $[\hat a'_h,\mathfrak a'_{k^-(h)}]$ is saturated and $\N(\mathfrak a_{k^-(h)}-\hat a_h)\in\theta_{h,h}\amsmathbb Z_{\geq 0}$;
    \item if $[\mathfrak b'_{k^+(h)},\hat b'_h]$ is void, then $\N(\hat b_h-\mathfrak b_{k^+(h)})\in\amsmathbb Z_{\geq 0}$, otherwise $[\mathfrak b'_{k^+(h)},\hat b'_h]$ is a saturated and then $\N(\hat b_h-\mathfrak b_{k^+(h)})\in\theta_{h,h}\amsmathbb Z_{\geq 0}$;
    \item for any $k \in \llbracket k^-(h),k^+(h)-1\rrbracket$, if $[\mathfrak b'_{k},\mathfrak a'_{k+1}]$ is void, then $\N(\mathfrak a_{k+1}-\mathfrak b_k)\in\amsmathbb Z_{\geq 0}$, otherwise $[\mathfrak b'_{k},\mathfrak a'_{k+1}]$ is a saturated and then $\N(\mathfrak a_{k+1}-\mathfrak b_k)\in\theta_{h,h}\amsmathbb Z_{\geq 0}$.
    \end{itemize}
If $a_1 = -\infty$ the first condition is omitted for $h = 1$. If $b_H = +\infty$ the second condition is omitted for $h = H$. If $H = 1$, $a_1 = -\infty$ and $b_1 = +\infty$, then we require that $\N\mathfrak a_1$ belongs to the lattice $\amsmathbb{L}$ where the leftmost particle of the system varies.
  \item $\forall k \in [K] \quad [\alpha_{k} - \eps_1,\beta_k + \eps_1] \subset [\mathfrak a_k,\mathfrak b_k] \subset [\alpha_k - 2\eps_1,\beta_k + 2\eps_1]$.
   \item $\forall k\in[K] \quad \big|\mathfrak{n}_{k}-\mu([\mathfrak a'_k,\mathfrak b'_k])\big|<\eps_2$.
   \item The $H$-tuple $(\widetilde{\mathfrak{n}}_h)_{h = 1}^{H}$ satisfies the constraints ($\star$), where for any $h \in [H]$ the number $\widetilde{\mathfrak{n}}_h$ is defined as the sum of the following terms: $\mathfrak{n}_k$ for $k \in \llbracket k^-(h),k^+(h)\rrbracket$; the mass $\frac{v - u}{\theta}$ of all the saturations $[u,v] \subseteq [\hat{a}'_h,\hat{b}'_h]$ appearing in Condition 2.; an extra $\big(-\frac{1}{\N}\big)$ for each $k \in \llbracket k^-(h),k^+(h) - 1\rrbracket$ such that $[\mathfrak{b}'_k,\mathfrak{a}'_{k + 1}]$ is saturated.
 \end{enumerate}
Define an event $\mathcal{A}$ by the two properties
 \begin{enumerate}
 \item[A.] The void segments appearing in Condition 2. do not have particles.
 \item[B.] The saturated segments appearing in Condition 2. do not have holes.
 \end{enumerate}
Then, the event $\mathcal{A}$ has probability at least $1 - C\exp\big(-\frac{\N}{C}\big)$, and the discrete ensemble $X$ conditioned on the intersection of $\mathcal{A}$ with the event that for any $k \in [K]$ there are $\N \mathfrak{n}_k$ particles in $[\N \mathfrak a_k, \N \mathfrak b_k]$, and localized to $\mathfrak{A} = \bigcup_{k = 1}^K [\mathfrak{a}_{k},\mathfrak{b}_{k}]$, yields a discrete ensemble obeying Assumptions~\ref{Assumptions_Theta}, \ref{Assumptions_basic}, \ref{Assumptions_offcrit}, \ref{Assumptions_analyticity}, and the additional Assumption~\ref{Assumptions_extra}.
\end{theorem}
\begin{proof}
The proof is a combination of the arguments of Theorem~\ref{proposition_FFF_conditioning} and \ref{proposition_Gaussian_conditioning}. We first choose auxiliary non-overlapping segments $[\mathfrak{a}^0_k,\mathfrak{b}^0_k]$ satisfying the integrality Condition 2. and such that for any $k \in [K]$, a small neighborhood of the $k$-th band of $\boldsymbol{\mu}$ is included in $[\mathfrak{a}^0_k, \mathfrak{b}^0_k]$ and a small neighborhood of $[\mathfrak{a}^0_k, \mathfrak{b}^0_k]$ is included in $[\hat a_{h^k},\hat b_{h^k}]$. In addition, we choose auxiliary filling fractions $(\mathfrak{n}_k^0)_{k = 1}^K$ obeying Conditions 1. and 5. in the statement. We introduce a new ensemble $Y^0$ by first conditioning $X$ on the intersection of $\mathcal{A}$ with the event $[\mathfrak{a}_{k}^0,\mathfrak{b}_{k}^0]$ has $\mathfrak{n}_{k}^0$ particles, and then localizing it to $\mathfrak{A}^0 = \bigcup_{k = 1}^K [\mathfrak{a}_{k}^0,\mathfrak{b}_{k}^0]$, as in Theorem~\ref{proposition_FFF_conditioning}. The new ensemble has $K$ rather than $H$ segments and its matrix of interactions $\boldsymbol{\Theta}$ is obtained by pre-composing the matrix of intensities of interactions of the ensemble $X$ with the map $k\mapsto h^k$ in rows and column indices. For the weights $(w_k)_{k = 1}^{K}$, potentials $(V_k)_{k = 1}^{K}$, and holomorphic functions $(\Phi^\pm_k)_{k = 1}^{K}$ and $(\phi^\pm_k)_{k = 1}^{K}$ we first define them as pre-compositions of the ones of the ensemble $X$ with $k\mapsto h^k$ in row indices, and then, exactly in the same way as in Theorem~\ref{proposition_FFF_conditioning}, we absorb into them the interaction with each saturation outside $\mathfrak{A}^0$. By design, the ensemble $Y^0$ satisfies Assumptions~\ref{Assumptions_Theta} and \ref{Assumptions_basic} and, therefore, it admits an equilibrium measure $\mu^{\mathfrak{n}^0}$ from Theorem~\ref{Theorem_equi_charact_repeat_2}. As usual, we set
\[
\forall k \in [K]\qquad {\mathfrak a}^{\prime\,0}_k := {\mathfrak a}^0_k - \frac{\theta_{h^k,h^k}-\frac{1}{2}}{\N} \quad \textnormal{and} \quad {\mathfrak b}^{\prime\,0}_k= {\mathfrak b}^0_k + \frac{\theta_{h^k,h^k}-\frac{1}{2}}{\N}.
\]
In addition to $Y^0$, we consider an auxiliary variational datum $\widetilde Y^0$, which is obtained by ignoring all subleading corrections to the potentials, similarly to $\widetilde Y$ obtained in the proof of Theorem~\ref{proposition_FFF_conditioning} by ignoring the lower-order terms \eqref{eq_weight_subleading_1}, \eqref{eq_weight_subleading_2}, and \eqref{eq_weight_subleading_3}. Note two small differences with the proof of Proposition~\ref{proposition_Gaussian_conditioning}: we did not impose integrality conditions on $\mathfrak{a}^0_k$ ,$\mathfrak{b}^0_k$ there, and there were no subleading corrections to the potentials on this step --- they did not arise because of the simplicity of the Gaussian weight.

Let us momentarily assume that our choice of segment filling fractions is
\[
\forall k \in [K] \qquad \mathfrak{n}_k^0 = \mu([\mathfrak{a}_k^{\prime\,0},\mathfrak{b}_k^{\prime\,0}]).
\]
Then the comparison with the characterization of the equilibrium measure in Theorem~\ref{Theorem_equi_charact_repeat_2} shows that for the auxiliary $\widetilde Y^0$, its equilibrium measure $\boldsymbol{\mu}^{\mathfrak{n}^0}$ is the restriction of the equilibrium measure $\boldsymbol{\mu}$ of $X$ to $\bigcup_{k = 1}^K [\mathfrak{a}_{k}^{\prime\,0},\mathfrak{b}_{k}^{\prime\,0}]$. Therefore, in this situation $\widetilde Y^0$ is off-critical, \textit{i.e.} it satisfies Assumption~\ref{Assumptions_offcrit}. Arguing as in Theorem~\ref{proposition_FFF_conditioning}, we further deduce that $Y^0$ --- whose potential differs from the one of $\widetilde Y^0$ by terms tending to $0$ as $\N\rightarrow\infty$ --- satisfies Assumptions~\ref{Assumptions_Theta}, \ref{Assumptions_basic}, \ref{Assumptions_offcrit}, \ref{Assumptions_analyticity} \emph{and} the additional Assumption~\ref{Assumptions_extra} for $\N$ large enough.

Thanks to the smooth dependence with respect to the segment filling fractions established in Proposition~\ref{Proposition_differentiability_filling_fraction}, there exist $\eps_1,\eps_2 > 0$ such that, for any other choice of segment filling fractions $(\mathfrak{n}^0_k)_{k = 1}^{K}$ satisfying
\begin{equation}
\label{inugngggg} \max_{k \in [K]} \big|\mathfrak{n}^0_k-\mu([\mathfrak{a}^{'0}_k,\mathfrak{b}^{'0}_k])\big|<\eps_2,
\end{equation}
the ensemble $Y$ defined like $Y^0$ but with arbitrary segment filling fractions satisfying \eqref{inugngggg}, still satisfies Assumptions~\ref{Assumptions_Theta}, \ref{Assumptions_basic}, \ref{Assumptions_offcrit}, \ref{Assumptions_analyticity} and the additional Assumption~\ref{Assumptions_extra}. Besides, its equilibrium measure $\boldsymbol{\mu}^{\mathfrak{n}^0}$ has precisely $K$ bands, each of them inside the $\eps_1$-neighborhood of the corresponding band of $\boldsymbol{\mu}$, and the list of saturations and voids for $\boldsymbol{\mu}^{\mathfrak{n}^0}$ --- in increasing order along the real line --- is in order-preserving bijection with that of the restriction of $\boldsymbol{\mu}$ to $\bigcup_{k = 1}^K [\mathfrak{a}_{k}^{\prime\,0},\mathfrak{b}_{k}^{\prime\,0}]$. Besides, we can always take a smaller $\eps_1$ if necessary, so that for any $k \in [K]$ the points $\mathfrak{a}_k^0$ and $\mathfrak{b}_k^0$ lie at distance greater than $2\eps_1$ from the endpoints of bands.

Such a choice of $(\eps_1,\eps_2)$ is suitable for what we want to prove. Indeed, let us now take parameters $(\mathfrak a_k,\mathfrak b_k,\mathfrak{n}_k)_{k = 1}^{K}$ satisfying all the assumptions of Theorem~\ref{proposition_fluct_conditioning}. In particular, Condition 4.\ for $\mathfrak{n}$ implies the condition \eqref{inugngggg}, because
\[
\forall k \in [K]\qquad \mathfrak{n}^0_k = \mathfrak{n}_k + \mu([\mathfrak{a}_{k}^{\prime\,0},\mathfrak{a}'_k]) + \mu([\mathfrak{b}'_{k},\mathfrak{b}_{k}^{\prime\,0}]),
\]
Our choices guarantee the inequalities $\mathfrak{a}_{k}^{\prime\,0}<\mathfrak{a}'_k<\mathfrak{b}'_{k}<\mathfrak{b}_{k}^{\prime\,0}$. Thanks to our previous discussions, we can apply Theorems~\ref{Theorem_ldpsup} and \ref{Theorem_ldsaturated} which say that with overwhelming probability, the ensemble $Y$ has no particles (respectively, no holes) in the segments outside $\N\mathfrak{A}$ which are included after usual rescaling and shifting in voids (respectively, saturations) of the equilibrium measure. Therefore, the initial ensemble $X$ conditioned on $\mathcal{A}$ and on the event that $[\N\mathfrak{a}_{k},\N\mathfrak{b}_{k}]$ has $\N\mathfrak{n}_{k}$ particles and localized to $\mathfrak{A}$ coincides with this ensemble $Y$ conditioned on the similar event $\mathcal{A}$ and localized to $\mathfrak{A}$. Due to the arguments of Theorem~\ref{proposition_FFF_conditioning}, since $Y$ satisfies Assumptions~\ref{Assumptions_Theta}, \ref{Assumptions_basic}, \ref{Assumptions_offcrit}, \ref{Assumptions_analyticity} and the additional Assumption~\ref{Assumptions_extra}, so does its localization to $\mathfrak{A}$, implying all the desired properties.
\end{proof}

\section{Revisiting the localization of ensembles}

\label{Section_alternative_localization}

In Theorem~\ref{proposition_FFF_conditioning} we defined the localization of an ensemble $X$, initially defined on segments $[a_h^X,b_h^X]$, to smaller segments $[a_h^Y,b_h^Y]$, for $h \in [H]$. We have described the data specifying the localized ensemble $Y$ by indicating the list of updates of the data specifying $X$ that eventually leads to the data specifying $Y$. This list is long but the logic behind it is apparent. We give here the final formulae for the data specifying $Y$, which offer a direct overview on $Y$ although the logic is now hidden. Besides, this will allow us in Section~\ref{Sec_Mismatch} to be more precise in the asymptotic comparison of the ensembles $Y$ and $X$, in particular, regarding their equilibrium measure.

We recall that the ensemble $Y$ has $H$ segments, each containing a single band of the equilibrium measure, and the intensities of interactions in the ensembles $X$ and $Y$ are the same and simply denoted $(\theta_{g,h})_{g,h = 1}^{H}$. If $\hat{a}_g^{\prime\,X}$ is in a saturation, we call
\[
\amsmathbb{F}_g^- = \big\{a_g^X,a_g^X + \theta_{g,g},a_g^X + 2\theta_{g,g},\ldots,a_g^Y - \theta_{g,g}\big\}
\]
the set of positions of the particles filling the segment $[a_g^X,a_g^Y)$. If $\hat{a}_g^{\prime\,X}$ is in a void, we set $\amsmathbb{F}_g^{-} = \emptyset$. We define likewise the set $\amsmathbb{F}_g^+$ of positions of the particles filling the segment $(b_g^Y,b_g^X]$ (if any). We will encounter indicator functions $\mathbbm{1}_{\amsmathbb{S}_h}(\hat{a}_h')$: since $a_h^{\prime\,X}$ is in a void (respectively, saturation) if and only if $\hat{a}_h^{\prime\,Y}$ is in a void (respectively, saturation), we have written in the argument of the indicator $\hat{a}_h'$ to refer indifferently to $\hat{a}_h^{\prime\,X}$ or $\hat{a}_h^{\prime\,Y}$.

\medskip

\noindent \textsc{Preliminary ingredients.} For each $h \in [H]$, we will need the functions
\begin{equation}
\label{Inthpmell}
\textnormal{Int}_{h}^{\pm}(\ell) := \prod_{g \neq h} \prod_{\ell' \in \amsmathbb{F}_g^{\pm}} \frac{1}{\N^{2\theta_{g,h}}}\cdot\frac{\Gamma\big(|\ell - \ell'| + 1\big) \cdot\Gamma\big(|\ell - \ell'| + \theta_{g,h}\big)}{\Gamma\big(|\ell - \ell'|\big)\cdot \Gamma\big(|\ell - \ell'| + 1 - \theta_{g,h}\big)},
\end{equation}
describing the interaction factors between the frozen particles away from the $h$-th segment with varying particles from the $h$-th segment, and their ratios
\[
\textnormal{RInt}_{h}^{\pm}(z) := \frac{\textnormal{Int}_{h}^{\pm}\big(\N z + \frac{1}{2}\big)}{\textnormal{Int}_{h}^{\pm}\big(\N z - \frac{1}{2}\big)}.
\]
More explicitly:
\begin{equation*}
\begin{split}
& \quad \textnormal{RInt}_{h}^{-}(z) \\
& = \prod_{g < h} \frac{\Gamma\Big(\frac{\N z - a_g^X + \frac{1}{2}}{\theta_{g,g}} + 1\Big)}{\Gamma\Big(\frac{\N z - a_g^Y + \frac{1}{2}}{\theta_{g,g}} + 1\Big)}\cdot\frac{\Gamma\Big(\frac{\N z - a_g^X - \frac{1}{2} + \theta_{g,h}}{\theta_{g,g}} + 1\Big)}{\Gamma\Big(\frac{\N z - a_g^Y - \frac{1}{2} + \theta_{g,h}}{\theta_{g,g}} + 1\Big)}\cdot \frac{\Gamma\Big(\frac{\N z - a_m^Y - \frac{1}{2}}{\theta_{g,g}} +1\Big)}{\Gamma\Big(\frac{\N z - a_g^X - \frac{1}{2}}{\theta_{g,g}} + 1\Big)}\cdot\frac{\Gamma\Big(\frac{\N z - a_g^Y + \frac{1}{2} - \theta_{g,h}}{\theta_{g,g}} + 1\Big)}{\Gamma\Big(\frac{\N z - a_g^X + \frac{1}{2} - \theta_{g,h}}{\theta_{g,g}} + 1\Big)} \\
&\quad \times \prod_{g > h} \frac{\Gamma\Big(\frac{a_g^Y - \N z + \frac{1}{2}}{\theta_{g,g}}\Big)}{\Gamma\Big(\frac{a_g^X - \N z + \frac{1}{2}}{\theta_{g,g}}\Big)}\cdot\frac{\Gamma\Big(\frac{a_g^Y - \N z + \frac{1}{2} - \theta_{g,h}}{\theta_{g,g}}\Big)}{\Gamma\Big(\frac{a_g^X - \N z + \frac{1}{2} - \theta_{g,h}}{\theta_{g,g}}\Big)}\cdot \frac{\Gamma\Big(\frac{a_g^X - \N z + \frac{1}{2}}{\theta_{g,g}}\Big)}{\Gamma\Big(\frac{a_g^Y - \N z + \frac{1}{2}}{\theta_{g,g}}\Big)}\cdot\frac{\Gamma\Big(\frac{a_g^X - \N z - \frac{1}{2} + \theta_{g,h}}{\theta_{g,g}}\Big)}{\Gamma\Big(\frac{a_g^Y - \N z - \frac{1}{2} + \theta_{g,h}}{\theta_{g,g}}\Big)},
\end{split}
\end{equation*}
The expression for $\textnormal{Rint}_{h}^{+}(z)$ is obtained by replacing in the above expression the points $(a_g^X,a_g^Y)$ with the points $(b_g^Y + \theta_{g,h},b_g^X + \theta_{g,h})$ for any $g \neq h$.

\medskip

\noindent \textsc{The parameters $\iota$ and $\rho$.} For any $h \in [H]$, the ensemble $X$ has its set of parameters $\iota_h^{\pm,X}$ and $\rho_{h,j}^{\pm,X}$ with $j \in [\iota_h^{\pm,X}]$, while the ensemble $Y$ has parameters
\[
\iota_h^{-,Y} = 2\mathbbm{1}_{\amsmathbb{S}_h}(\hat{a}_h') \qquad \textnormal{and}\qquad \iota_h^{+,Y} = 2 \mathbbm{1}_{\amsmathbb{S}_h}(\hat{b}_h').
\]
Besides, if $\iota_h^{\pm,Y} = 2$ we have $\rho_{h,1}^{\pm,Y} = 1$ and $\rho_{h,2}^{\pm ,Y} = \theta_{h,h}$.

\medskip

\noindent \textsc{The weights.} For the ensemble $Y$, the weight on the $h$-th segment for $h \in [K]$ reads
\begin{equation*}
\begin{split}
w_h^Y(\ell) & = \exp\Bigg(-\N U_h^X\bigg(\frac{\ell}{\N}\bigg)\Bigg) \cdot \prod_{j = 1}^{\iota_h^{-,X}} \frac{\N^{\ell - a_h^X + \rho_{h,j}^{-,X} - \frac{1}{2}}}{\Gamma\big(\ell - a_h^X + \rho_{h,j}^{-,X}\big)} \cdot \prod_{j = 1}^{\iota_h^{+,X}} \frac{\N^{b_h^X - \ell + \rho_{h,j}^{+,X} - \frac{1}{2}}}{\Gamma\big(b_h^X - \ell + \rho_{h,j}^{+,X}\big)} \\
&\quad \times \textnormal{Int}^-_{h}(\ell) \cdot \bigg[\frac{1}{\N^{2(a_h^Y - a_h^X)}}\cdot\frac{\Gamma\big(\ell - a_h^X + 1\big)\cdot\Gamma\big(\ell - a_h^X + \theta_{h,h}\big)}{\Gamma\big(\ell - a_h^Y + 1\big)\cdot\Gamma\big(\ell - a_h^Y + \theta_{h,h}\big)}\bigg]^{\mathbbm{1}_{\amsmathbb{S}_h}(\hat{a}_h')} \\
& \quad \times \textnormal{Int}^+_{h}(\ell) \cdot \bigg[\frac{1}{\N^{2(b_h^X - b_h^Y)}}\cdot \frac{\Gamma\big(b_h^X - \ell + 1\big)\cdot\Gamma\big(b_h^X - \ell + \theta_{h,h}\big)}{\Gamma\big(b_h^Y - \ell + 1\big)\cdot\Gamma\big(b_h^Y - x + \theta_{h,h}\big)}\bigg]^{\mathbbm{1}_{\amsmathbb{S}_h}(\hat{b}_h')}.
\end{split}
\end{equation*}
The first line is simply $w_h^X(\ell)$ (see Assumption~\ref{Assumptions_analyticity}).

\medskip

\noindent \textsc{The potentials.} We recall that for any ensemble satisfying Assumption~\ref{Assumptions_basic}, the potential is decomposed in regular and singular part as follows
\[
\forall h \in [H]\qquad V_h(x) = U_h(x) + \iota_h^-\, \mathrm{Llog}(x - \hat{a}_h^{\prime}) + \iota_h^+\,\mathrm{Llog}\big(\hat{b}_h^{\prime} - x).
\]
For the ensemble $Y$ the full potential on the $h$-th segment is
\begin{equation}
\label{TheVVX}
\begin{split}
V_h^Y(x) & = V_h^X(x) + 2\, \mathrm{Llog}(x - \hat{a}_h^{\prime\,Y}) \mathbbm{1}_{\amsmathbb{S}_h}(\hat{a}_h') + 2\,\mathrm{Llog}(\hat{b}_h^{\prime\,Y} - x)\mathbbm{1}_{\amsmathbb{S}_h}(\hat{b}_h') \\
& \quad - \iota_h^{-,X}\, \mathrm{Llog}(x - \hat{a}_h^{\prime\,X}) - \iota_h^{+,X}\,\mathrm{Llog}(\hat{b}_h^{\prime\,X} - x) - \frac{1}{\N}\log(\textnormal{Int}_h^{-}(\N x) \cdot \textnormal{Int}_h^+(\N x)\big)  \\
& \quad - \frac{1}{\N} \log\Bigg(\prod_{j = 1}^{\iota_h^{-,X}} \frac{\N^{\N x - a_h^X + \rho_{h,j}^{-,X} - \frac{1}{2}}}{\Gamma\big(\N x - a_h^X + \rho_{h,j}^{-,X}\big)} \cdot \prod_{j = 1}^{\iota_h^{+,X}} \frac{\N^{b_h^X - \N x + \rho_{h,j}^{+,X} - \frac{1}{2}}}{\Gamma\big(b_h^X - \N x + \rho_{h,j}^{+,X}\big)}\Bigg) \\
& \quad - \frac{1}{\N}\Bigg[ \log\bigg(\frac{\Gamma\big(\N x - a_h^X + 1\big) \cdot \Gamma\big(\N x - a_h^X + \theta_{h,h}\big)}{\N^{2(\N x - a_h^X) + \theta_{h,h}}}\bigg)\Bigg] \mathbbm{1}_{\amsmathbb{S}_h}(\hat{a}_h') \\
& \quad - \frac{1}{\N}\Bigg[\log\bigg(\frac{\Gamma\big(b_h^X - \N x + 1\big)\cdot \Gamma(b_h^X - \N x + \theta_{h,h})}{\N^{2(b_h^X - \N x) + \theta_{h,h}}}\bigg)\Bigg]\mathbbm{1}_{\amsmathbb{S}_h}(\hat{b}'_h),
\end{split}
\end{equation}
and its regular part is
\begin{equation}
\label{theUYX}
\begin{split}
U_h^Y(x) & = U_h^X(x) - \frac{1}{\N}\log\Bigg(\prod_{j = 1}^{\iota_h^{-,X}} \frac{\N^{\N x - a_h^X + \rho_{h,j}^{-,X} - \frac{1}{2}}}{\Gamma\big(\N x - a_h^X + \rho_{h,j}^{-,X}\big)} \cdot \prod_{j = 1}^{\iota_h^{+,X}} \frac{\N^{b_h^X - \N x + \rho_{h,j}^{+,X} - \frac{1}{2}}}{\Gamma\big(b_h^X - \N x + \rho_{h,j}^{+,X}\big)}\Bigg) \\
& - \frac{1}{\N}\log\big(\textnormal{Int}_h^-(\N x)\big) - \frac{\mathbbm{1}_{\amsmathbb{S}_h}(\hat{a}_h')}{\N}\log\bigg(\frac{\Gamma(\N x - a_h^X + 1)\cdot \Gamma(\N x - a_h^X + \theta_{h,h})}{\N^{2(\N x - a_h^X) + \theta_{h,h}}}\bigg) \\
& - \frac{1}{\N} \log\big(\textnormal{Int}_h^+(\N x)\big) - \frac{\mathbbm{1}_{\amsmathbb{S}_h}(\hat{b}_h')}{\N} \log\bigg(\frac{\Gamma(b_h^X - \N x + 1)\cdot \Gamma(b_h^X - \N x + \theta_{h,h})}{\N^{2(\N x - a_h^X) + \theta_{h,h}}}\bigg).
\end{split}
\end{equation}

\medskip

\noindent \textsc{Auxiliary functions.} We recall the decomposition
\[
\forall h \in [H]\qquad \frac{w_h\big(\N z + \frac{1}{2}\big)}{w_h\big(\N z - \frac{1}{2}\big)} = \frac{\Phi_h^+(z)}{\Phi_h^-(z)},
\]
where $\Phi_h^{\pm}(z)$ are holomorphic in a neighborhood of the $h$-th segment and are specified by the formulae in Assumption~\ref{Assumptions_analyticity}. These factors for the ensemble $Y$ are:

\begin{equation*}
\begin{split}
\Phi_h^{+,Y}(z) & = \Phi_h^{+,X}(z) \cdot \Bigg[\prod_{j = 1}^{\iota_{h}^{-,X}} \bigg(z - \hat{a}_h^X + \frac{\rho_{h,j}^{-,X} - \frac{1}{2}}{\N}\bigg)^{-1}\Bigg] \cdot \textnormal{RInt}_{h}^-(z) \cdot \textnormal{RInt}_h^+(z) \\
& \quad  \times \Bigg[\bigg(z - \hat{a}_h^X + \frac{1}{2\N}\bigg)\bigg(z - \hat{a}_h^X + \frac{\theta_{h,h} - \frac{1}{2}}{\N}\bigg)\Bigg]^{\mathbbm{1}_{\amsmathbb{S}_h}(\hat{a}_h')} \\
& \quad  \times \Bigg[\frac{\big(\hat{b}_h^Y - z + \frac{1}{2\N}\big)\big(\hat{b}_h^Y - z + \frac{\theta_{h,h} - \frac{1}{2}}{\N}\big)}{\big(\hat{b}_h^X - z + \frac{1}{2\N}\big)\big(\hat{b}_h^X - z + \frac{\theta_{h,h} - \frac{1}{2}}{\N}\big)}\Bigg]^{\mathbbm{1}_{\amsmathbb{S}_h}(\hat{b}_h')} \\
\Phi_h^{-,Y}(z) & = \Bigg[\bigg(z - \hat{a}_h^Y + \frac{1}{2\N}\bigg)\bigg(z - \hat{a}_h^Y + \frac{\theta_{h,h} - \frac{1}{2}}{\N}\bigg)\Bigg]^{\mathbbm{1}_{\amsmathbb{S}_h}(\hat{a}_h')}.
\end{split}
\end{equation*}

\medskip

\noindent \textsc{Leading order of auxiliary functions.} Recall that in Definition~\ref{Definition_phi_functions} we specified leading factors in $\Phi_{h}^{\pm}$, which we denoted $\phi_h^{\pm}$ and which satisfy
\[
\frac{\phi_h^+(z)}{\phi_h^{-}(z)} = e^{-\partial_zV_h(z)}.
\]
For the ensemble $Y$ they are given by
\begin{equation*}
\begin{split}
& \quad \phi_h^{+,Y}(z) = \phi_h^{+,X}(z) \cdot \exp\Big[ \big(\log(\textnormal{Int}_h^+)\big)'(\N z) + \big(\log(\textnormal{Int}_h^-)\big)'(\N z)\big)\Big] \\
&\quad\quad \times \exp\Bigg(\sum_{j = 1}^{\iota_h^{-,X}} \N\big(\log \N - \mathsf{\Psi}(\N z - a_h^X + \rho_{h,j}^{-,X})\big) + \sum_{j = 1}^{\iota_h^{+,X}} \N\big(-\log \N + \mathsf{\Psi}(b_h^X - \N z + \rho_{h,j}^{+,X})\big)\Bigg) \\
&\quad \qquad \times \Bigg[\frac{\exp\big[\N\big(2\log \N - \mathsf{\Psi}(b_h^X - \N z + 1) - \mathsf{\Psi}(b_h^X - \N z + \theta_{h,h})\big]}{(\hat{b}_h^{\prime\,X} - z)^{\iota_{h}^{+,X}}}\Bigg]^{\mathbbm{1}_{\amsmathbb{S}_h}(\hat{b}_h')} \\
&\quad \quad \times \Big[\exp\big[\N\big(-2\log \N + \mathsf{\Psi}(\N z - a_h^X + 1) + \mathsf{\Psi}(\N z - a_h^X + \theta_{h,h})\big]\Big]^{\mathbbm{1}_{\amsmathbb{S}_h}(\hat{a}_h')},
\end{split}
\end{equation*}
where\label{index:digam} $\mathsf{\Psi} = (\log \Gamma)'$ is the digamma function, and
\[
\phi_h^{-,Y}(z) = \big(z - \hat{a}_h^{\prime\,Y}\big)^{2\mathbbm{1}_{\amsmathbb{S}_h}(\hat{a}_h')}.
\]

\section{Mismatch of equilibrium measures after localization}
\label{Sec_Mismatch}

We mentioned in Section~\ref{fixedsecfill} that the conditioned and localized ensemble delivered by Theorem~\ref{proposition_FFF_conditioning} have an equilibrium measure which is not exactly the restriction of the equilibrium measure of the original ensemble, due to our definition of singular and regular parts of the potential which depends on segment endpoints. In Theorem~\ref{proposition_fluct_conditioning}, the original ensemble $X$ has fluctuating filling fractions. But, after we have chosen segments and filling fractions as described in Theorem~\ref{proposition_fluct_conditioning}, we can also speak of the equilibrium measure of $X$ conditioned to have those fixed filling fractions: this is the minimizer of the functional $-\mathcal{I}$ over the space of measures obeying these filling fraction constraints. And for the same reason, the restriction of this minimizer to $\amsmathbb{A}^Y$ does not exactly coincide with the equilibrium measure of the localized ensemble $Y$ delivered by Theorem~\ref{proposition_fluct_conditioning}. If one attempts obtaining concrete formulae for specific ensembles where the localization procedure must be used (\textit{e.g.} in Section~\ref{Gausscorr} and Chapter~\ref{Chap11}), it is important to understand better this mismatch. In this direction the following result will be useful.

\begin{proposition}
\label{mueqYYtilde}
In the setting of Theorem~\ref{proposition_FFF_conditioning}, let $\mu^X$ be the equilibrium measure of the ensemble $X$, $\mu^Y$ the equilibrium measure of $Y$. Let $h \in [H]$ and set $\amsmathbb{A}_h^{Y} = \amsmathbb{A}^Y \cap [\hat{a}_h',\hat{b}_h']$. For any $\delta > 0$, there exists constants $C_1,C_2 > 0$ depending only on $\delta$ and the constants in the assumptions such that, for any $z \in \amsmathbb{C}$ at distance at least $\delta$ from $\amsmathbb{A}_h^Y$, we have
\[
\bigg|\int_{\amsmathbb{A}^{Y}_h} \frac{\mu_h^X(x) - \mu^{Y}_h(x)}{z - x} \dd x \bigg| \leq \frac{C_1}{\N}.
\]
The upper bound is actually
\[
\bigg|\int_{\amsmathbb{A}^{Y}_h} \frac{\mu_h^X(x) - \mu^{Y}_h(x)}{z - x} \dd x \bigg| \leq \frac{C_2}{\N^2}.
\]
in case the following three conditions hold simultaneously
\begin{itemize}
\item the endpoints $\hat{a}_h'$ and $\hat{b}'_h$ are in a void, or $\theta_{h,h} = 1$;
\item $\iota_h^{-,X} = 0$, or $\rho_{h,j}^{-,X} = \theta_{h,h}$ for every $j \in [\iota_h^{-,X}]$;
\item $\iota_h^{+,X} = 0$, or $\rho_{h,j}^{+,X} = \theta_{h,h}$ for every $j \in [\iota_h^{+,X}]$.
\end{itemize}
Similar results hold in the setting of Theorem~\ref{proposition_fluct_conditioning}, provided we take $\mu^X$ to be the measure obeying fixed filling fractions constraints and compare it with $\mu^Y$ corresponding to the same constraints.
\end{proposition}
Note that the second condition is the same as the one met at the end of Lemma~\ref{Lemma_phi_properties}. The proof of Proposition~\ref{mueqYYtilde} is based on the asymptotic analysis of the formulae of Section~\ref{Section_alternative_localization}. At first sight they look complicated, but many simplifications occur if we only keep the leading and the first subleading contribution as $\N \rightarrow \infty$ as shown in the following lemma.

\begin{lemma}
\label{goodsimz}
We have as $\N \rightarrow \infty$
\begin{equation}
\label{VyavecVx}
\begin{split}
V_h^Y(x) & = V_h^{X}(x)  - \sum_{g = 1}^{H} 2\theta_{g,h} \int_{\amsmathbb{A}^{X} \setminus \amsmathbb{A}^Y} \log|x - y|\dd \mu^{X}(y) + \varpi \\
& \quad + \frac{1}{\N}\bigg((\theta_{h,h} - 1)\mathbbm{1}_{\amsmathbb{S}_h}(\hat{a}'_h) + \sum_{j = 1}^{\iota_h^{-,X}} (\rho_{h,j}^{-,X} - \theta_{h,h}) \bigg)\log(x - \hat{a}^{\prime\,X}_h)  \\
& \quad + \frac{1}{\N}\bigg((\theta_{h,h} - 1)\mathbbm{1}_{\amsmathbb{S}_h}(\hat{b}'_h) + \sum_{j = 1}^{\iota_h^{+,X}} (\rho_{h,j}^{X,+} - \theta_{h,h})\bigg)\log(\hat{b}^{\prime\,X}_h - x)  \\
& \quad + \frac{1}{\N} \sum_{g \neq h} \frac{\theta_{g,g}}{\theta_{g,h}}(\theta_{g,g} - 1) \Bigg[- \mathbbm{1}_{\amsmathbb{S}_g}(\hat{a}'_g) \, \log\bigg|\frac{x - \hat{a}_g^{\prime\,X}}{x - \hat{a}_g^{\prime\,Y}}\bigg| + \mathbbm{1}_{\amsmathbb{S}_g}(\hat{b}'_g) \,\log\bigg|\frac{x - \hat{b}_g^{\prime\,X}}{x - \hat{b}_g^{\prime\,Y}}\bigg|\Bigg] +  O\bigg(\frac{1}{\N^2}\bigg),
\end{split}
\end{equation}
where $\varpi$ is independent of $x$ and the error is uniform for $x \in [\hat{a}_h^{\prime\,Y},\hat{b}_h^{\prime\,Y}]$.
\end{lemma}

\begin{proof}[Proof of Lemma~\ref{goodsimz}]
We analyze the large $\N$ behavior of each term in \eqref{TheVVX}. The symbol $\varpi$ will be used for various quantities independent of $x$ whose value is irrelevant. We start with the terms involving the various forms of $x - \hat{a}_h$. Let $\hat{a} \in \amsmathbb{R}$ and $\rho,\theta > 0$, and set $\hat{a} = \hat{a}' + \frac{\theta - 1/2}{\N}$. The Stirling expansion for the Gamma function yields\footnote{Even though the right-hand side is a function of $\hat{a}$ and $\rho$ independent of $\theta$, the parameter $\theta$ appears in the right-hand side because we use the shifted endpoint $\hat{a}'$ which depends on $\hat{a}$ and $\theta$.} for $x > \hat{a}$
\begin{equation}
\label{Stirlingsimpli}
\log\left(\frac{\Gamma\big(\N(x - \hat{a}) + \rho\big)}{\N^{\N(x - \hat{a}) + \rho - \frac{1}{2}}} \right) = \N\,\mathrm{Llog}(x - \hat{a}') + (\rho -  \theta)\log(x - \hat{a}') + \frac{\log(2\pi)}{2} + O\bigg(\frac{1}{\N}\bigg)
\end{equation}
If we are given $\varepsilon > 0$ independent of $\N$ and such that $\min(\rho,\theta,x - \hat{a}) > \varepsilon$, the absolute value of the error $O(\frac{1}{\N})$ is bounded by $\frac{C}{\N}$ for a constant $C > 0$ that can be chosen to depend only on $\varepsilon$. Applying this to $(\hat{a},\theta)  = (\hat{a}_h^{X},\theta_{h,h})$ and $\rho = \rho_{h,j}^{-,X}$ for $j \in [\iota_h^-]$ or $\rho = 1$ or $\rho = \theta_{h,h}$ treats the corresponding terms in the third and fourth line of \eqref{TheVVX}. Combining them with the terms involving $\hat{a}^{\prime\,X}_h$ or $\hat{a}^{\prime\,Y}_h$ from the first two lines, we get a contribution
\begin{equation}
\label{2ugnssqs1}
\begin{split}
& \quad 2\,\mathbbm{1}_{\amsmathbb{S}_h}(\hat{a}'_h)\big(\mathrm{Llog}(x - \hat{a}_h^{\prime\,Y}) - \mathrm{Llog}(x - \hat{a}_h^{\prime\,X})\big) \\
& + \frac{1}{\N}\bigg( \mathbbm{1}_{\amsmathbb{S}}(\hat{a}'_h) \cdot (\theta_{h,h} - 1) + \sum_{j = 1}^{\iota_h^{-,X}} (\rho_{h,j}^{-,X} - \theta_{h,h})\bigg)\log(x - \hat{a}_{h}^{\prime\,X}) + \varpi_h + O\bigg(\frac{1}{\N^2}\bigg),
\end{split}
\end{equation}
The identity
\[
\mathrm{Llog}(x - \hat{a}_h^{\prime\,Y}) - \mathrm{Llog}(x - \hat{a}_h^{\prime\,X}) = - \int_{\hat{a}_h^{\prime\,Y}}^{\hat{a}_h^{\prime\,X}} \log|x - y| \dd y.
\]
helps us recognizing the corresponding integral term in the first line of the claimed formula \eqref{VyavecVx}, and therefore \eqref{2ugnssqs1} reproduces the terms involving $(x - \hat{a}_h^{\prime\,X})$ and $(x - \hat{a}_h^{\prime\,Y})$ in \eqref{VyavecVx}. The claimed terms involving $(\hat{b}_h^{\prime\,X} - x)$ or $(\hat{b}_h^{\prime\,Y} - x)$ are obtained in a similar way.

The only term that we are so far missing in \eqref{VyavecVx} is the sum over $g \neq h$ in the first line. It will come from the asymptotics of the $\textnormal{Int}_h^{\pm}(\N x)$-terms in \eqref{TheVVX}. We first look at the factors with $g > h$ in  $\textnormal{Int}_h^{-}(\N x)$. Such a factor is present is $\hat{a}_g^{\prime}$ is saturated, and it is then equal to
\[
\prod_{i = 0}^{M - 1} \frac{1}{\N^{2\theta_{g,h}}} \cdot \frac{\Gamma\big(\N(\hat{a}_g^X - x) + i\theta_{g,g} + 1\big)\cdot \Gamma\big(\N(\hat{a}_g^X - x) + i\theta_{g,g} + \theta_{g,h}\big)}{\Gamma\big(\N(\hat{a}_g^X - x) + i \theta_{g,g}\big) \cdot \Gamma\big(\N(\hat{a}_g^X - x) + i \theta_{g,g} + 1 - \theta_{g,h}\big)},
\]
where $\N(\hat{a}_g^Y - \hat{a}_g^{X}) = \hat{a}_g^{\prime\,Y} - \hat{a}_{g}^{\prime\,X} = M\theta_{g,g}$. By comparison with \eqref{phiexp} with the substitution $\theta = \theta_{g,g}$, this product is equal to
\begin{equation}
\label{1N2tM}
\frac{1}{\N^{2\theta_{g,h} M}} \cdot  \frac{\gimel_{M}\big(\N(\hat{a}_g^{\prime\, X} - x) + \theta_{g,g} - \frac{1}{2}\big) \cdot \gimel_M\big(\N(\hat{a}_g^{\prime\,X} - x) + \theta_{g,h} + \theta_{g,g} - \frac{3}{2}\big)}{\gimel_{M}\big(\N(\hat{a}_g^{\prime\,X} - x) + \theta_{g,g} - \frac{3}{2}\big)  \cdot \gimel_{M}\big(\N(\hat{a}_g^{\prime\,X} - x) - \theta_{g,h} +  \theta_{g,g} - \frac{1}{2}\big)}.
\end{equation}
We recall from \eqref{gnexpmfgun} the asymptotic expansion as $\N \rightarrow \infty$ of $\gimel_{M}(\N \xi + \xi')$  with $M = \N \hat{m}$ and $\xi, \hat{m}$ positive away from $0$, $\xi'$ bounded. In particular if we have another bounded parameter  $\xi''$ we obtain
\begin{equation*}
\begin{split}
\log\left(\frac{\gimel_{M}(\N \xi + \xi')}{\gimel_M(\N \xi + \xi'')}\right) & = (\xi' - \xi'')\hat{m}\N \log \N + \frac{(\xi' - \xi'')}{\theta_{g,g}}\big( \mathrm{Llog}(\xi + \theta_{g,g}\hat{m}) - \mathrm{Llog}(\xi)\big) \N \\
& \quad + \frac{(\xi' - \xi'')(\xi' + \xi'' + 1 - \theta_{g,g})}{2\theta_{g,g}}\log\left(\frac{\xi + \theta_{g,g}\hat{m}}{\xi}\right) + \varpi_0 + O\bigg(\frac{1}{\N}\bigg),
\end{split}
\end{equation*}
where $\varpi_0$ is independent of $\xi$. We apply this for $\xi = (\hat{a}_g^{\prime\,X} - x)$, which leads to $\xi + \theta_{g,g}\hat{m} = \hat{a}_g^{\prime\, Y}$, and to
\[
(\xi',\xi'') = \bigg(\theta_{g,g} - \frac{1}{2},\theta_{g,g} - \frac{3}{2}\bigg) \qquad \textnormal{or} \qquad (\xi',\xi'') = \bigg(\theta_{g,g} + \theta_{g,h} - \frac{3}{2},\theta_{g,g} - \theta_{g,h} - \frac{1}{2}\bigg).
\]
This implies that $-\frac{1}{\N}$ times the logarithm of \eqref{1N2tM} is
\begin{equation}
\label{offcterm}
\frac{2\theta_{g,h}}{\theta_{g,g}}\big(\mathrm{Llog}(\hat{a}_g^{\prime\,X} - x) - \mathrm{Llog}(\hat{a}_g^{\prime\,Y} - x)\big) + \frac{(\theta_{g,g} - 1)\theta_{g,h}}{\theta_{g,g}\,\N} \log\bigg(\frac{\hat{a}_g^{\prime\,X} - x}{\hat{a}_g^{\prime\,Y} - x}\bigg) + \varpi_h^{-} + O\bigg(\frac{1}{\N^2}\bigg).
\end{equation}
After using again
\[
\mathrm{Llog}(\hat{a}_g^{\prime\,X} - x) - \mathrm{Llog}(\hat{a}_g^{\prime\,Y} - x) = -\int_{\hat{a}_g^{\prime\,X}}^{\hat{a}_g^{\prime\,Y}} \log|x - y|\,\dd y,
\]
to identify the corresponding integral term in the first line of \eqref{VyavecVx}, we recognize all terms of \eqref{offcterm} as contributions in \eqref{VyavecVx}. The remaining terms in \eqref{VyavecVx} can be obtained by a similar analysis of the $g < h$ factors of $\textnormal{Int}^{-}(\N x)$, and of all factors of $\textnormal{Int}^+(\N x)$.
\end{proof}

\begin{proof}[Proof of Proposition~\ref{mueqYYtilde}]
We start with the setting of Theorem~\ref{proposition_FFF_conditioning} for which Lemma~\ref{goodsimz} was designed. In the proof of Theorem~\ref{proposition_FFF_conditioning} we had introduced a variational datum $\widetilde{Y}$ whose equilibrium measure was the restriction of the equilibrium measure of $X$ to $\amsmathbb{A}^{Y}$. Comparing with \eqref{VysandVx}, for any $h \in [H]$ the $h$-th potential for $\widetilde{Y}$ matches the first line in \eqref{VyavecVx}. Lemma~\ref{goodsimz} tells us it differs from the potential $V_h^Y(x)$ by an irrelevant constant and a $O(\frac{1}{\N})$ term which is a smooth function of $x$, or even a $O(\frac{1}{\N^2})$ term if the extra condition in the corollary is satisfied. Then, Theorem~\ref{Theorem_differentiability_full} guarantees that for any $h \in [H]$, the Stieltjes transform of $\mu^X|_{\amsmathbb{A}_h^Y} - \mu^{Y}$ is $O(\frac{1}{\N})$ (or $O(\frac{1}{\N^2})$ under the extra condition) for $z$ away from $\amsmathbb{A}_h^Y$ as well. In the setting of Theorem~\ref{proposition_fluct_conditioning} where the original ensemble $X$ may have fluctuating filling fractions, apart from more complicated choices of segments $\amsmathbb{A}^Y$ that also depend on filling fractions, the parameters of $X$ are updated exactly as we described above to get the ensemble $Y$, so the same conclusion applies.
\end{proof}

\chapter{Asymptotics with fixed filling fractions: correlators}
\label{Chapter_fff_expansions}

In this chapter we study the $n$-point correlators as $\N\rightarrow\infty$ under Assumptions~\ref{Assumptions_Theta}, \ref{Assumptions_basic}, \ref{Assumptions_offcrit} and \ref{Assumptions_analyticity}, and the additional assumptions that the segment filling fractions are deterministically fixed and each segment $[\hat a_h,\hat b_h]$ has a unique band of the
equilibrium measure. The main results are summarized in Section
\ref{Section_cumulant_expansion} after we introduce various notations. These results are required for Chapter~\ref{Chapter_filling_fractions} where we handle the more general case where we do not impose the two additional assumptions.

To recap the definition of the observables (\textit{cf.} Section~\ref{section_rought_ann}) we deal with the Stieltjes transform of the empirical measure recentered with the equilibrium measure, which by convention we multiply by $\N$:
\begin{equation}
\label{decomposstil}\forall h \in [H]\qquad \Delta G_h(z) := G_h(z) - \N\Gm_{\mu_h}(z),
\end{equation}
where
\[
G_{h}(z) = \sum_{i=1}^{N_h} \frac{1}{z-\frac{\ell_i^h}{\N}},\qquad \Gm_{\mu_h}(z) =\int_{\hat a'_h}^{\hat b'_h} \frac{\mu_h(x)\dd x}{z-x}.
\]
The $(\ell_i^h)_{i = 1}^{N_h}$ are the particles of the discrete ensemble that belong to the $h$-th segment, and $\mu_h$ is the equilibrium measure. The computations in this chapter will more directly concern the cumulants of these Stieltjes transforms. Namely, for $n \geq 1$ and $h_1,\ldots,h_n \in [H]$, the $n$-point correlator is defined as the cumulant of the $n$ random variables $\Delta
G_{h_1}(z_1),\ldots,\Delta G_{h_n}(z_n)$.
\begin{equation} \label{eq_correlators_def}
\begin{split}
W_{n;h_1,\ldots,h_n}(z_1,\ldots,z_n) & = \E^{(\textnormal{c})}\big[\Delta
G_{h_1}(z_1),\ldots,\Delta G_{h_n}(z_n)\big] \\
& =\partial_{s_1}\cdots\partial_{s_n}\log \E\Bigg[ \exp\bigg(\sum_{i=1}^n s_i \Delta G_{h_i}(z_i)\bigg)\Bigg]\Bigg|_{s_1 = \cdots = s_n = 0}.
\end{split}
\end{equation}
For $n \geq 2$ one can use $(G_{h_i}(z_i))_{i = 1}^n$ in this formula instead of $(\Delta G_{h_i}(z_i))_{i = 1}^n$ since the recenterings $(\Gm_{\mu_{h_i}}(z_i))_{i = 1}^n$ are deterministic. The $n$-point correlators are homogeneous polynomials of degree $n$ in the moments of $\Delta\boldsymbol{G}(z)$. For a discrete ensemble satisfying the basic Assumptions~\ref{Assumptions_Theta}, the large deviation estimates of Corollary~\ref{Corollary_a_priory_1} show that, for any $\delta,\varepsilon > 0$ independent of $\N$, any $h \in [N]$ and $z \in \amsmathbb{C}$ at distance at least $\delta$ away from $[\hat{a}_h,\hat{b}_h]$, the recentered empirical Stieltjes transform $\Delta G_h(z)$ is in absolute value smaller that $C\N^{\frac{1}{2}
+ \varepsilon}$ for some constant $C > 0$ independent of $\N$ and $h$, with
overwhelming probability. So, the order of magnitude of the first correlator $W_{1;h}(z) = \E[\Delta G_h(z)]$ is at
most $\N^{\frac{1}{2} + \varepsilon}$ and for $n \geq 2$, the order of magnitude of
$W_{n;h_1,\ldots,h_n}(z_1,\ldots,z_n)$ is at most $\N^{\frac{n}{2} + \varepsilon}$, for $z,z_1,\ldots,z_n$ staying uniformly away from the segments defining the ensemble.

In this chapter we improve these \textit{a priori} estimates by proving that the correlators are in fact much smaller as $\N \rightarrow \infty$, and in fact provide their asymptotic expansions with the help of the Nekrasov equations. As a corollary, we show that the random variables $\Delta G_h(z)$ are asymptotically Gaussian as $\N\rightarrow\infty$ and compute their asymptotic covariance.

\section{Preliminaries}

\label{Section_notations}

In this section we introduce various notations used in the formulation of the asymptotic expansions of $W_{n;h_1,\ldots,h_n}(z_1,\ldots,z_n)$.

\subsection{The \texorpdfstring{$\langle \bth_h \cdot \boldsymbol{F}\rangle $}{} notation}
\label{Theopr}

We are already familiar with the notation of $H$-tuples by bold characters. If $\boldsymbol{F} = (F_h)_{h = 1}^H$ is an $H$-dimensional vector, and
$P$ is a polynomial, we will denote for any $h \in [H]$
\[
\big\langle P(\bth_h)\cdot \boldsymbol{F}\big\rangle = \sum_{g = 1}^H P(\theta_{h,g})\, F_{g}.
\]
For instance, we have
\[
\big\langle \bth_h(\bth_h-1)\cdot \boldsymbol{F}\big\rangle = \sum_{g = 1}^H \theta_{h,g}(\theta_{h,g}-1)\,F_{g}.
\]

In this chapter we will deal with $n$-point correlators, which are actually $H^n$-tuples of functions. We will also denote them in bold characters:
\[
\boldsymbol{W}_n(z_1,\ldots,z_n) = \big(W_{n;h_1,\ldots,h_n}(z_1,\ldots,z_n)\big)_{h_1,\ldots,h_n = 1}^{H}
\]
We sometimes want to form partial tuples, by fixing the values of $h_i$ for some $i$ only, while the other indices are left to vary. In that case, we still use the bold notation, but indicate with $\bullet$ the indices that vary to form a tuple. If $h,h_1 \in [H]$ are fixed values and $z,z_1,z_2,z_3,z_4,z$ complex variables, here are two examples
\begin{equation*}
\begin{split}
\boldsymbol{W}_{2;\bullet,h_2}(z_1,z_2) & = \big(W_{2;h_1,h_2}(z_1,z_2)\big)_{h_1 = 1}^{H}, \\
\boldsymbol{W}_{4;h,h,\bullet,\bullet}(z,z,z_3,z_4) & = \big(W_{4;h,h,h_3,h_4}(z,z,z_3,z_4)\big)_{h_3,h_4 = 1}^{H}.
\end{split}
\end{equation*}
The first one is an $H$-tuple, the second one an $H^2$-tuple. This can be used in combination with the $\langle \cdot \rangle$ notation. For instance, we will meet expressions like
\[
\big\langle \bth_{h}\cdot \boldsymbol{W}_{2;\bullet,h_2}(z,z_2)\big\rangle=\sum_{g=1}^H \theta_{h,g}\, W_{2;g,h_2}(z,z_2).
\]
We see that the $\bullet$ is replaced with the dummy index $g$ which is summed over $[H]$, as prescribed by $\langle \cdot \rangle$.

It also happens that we want to sum over several dummy indices, so we need a new layer of notations. If $n$ is a positive integer, $P_1,\ldots,P_n$ are polynomials, $h_1,\ldots,h_n \in [H]$ are fixed and $\boldsymbol{F}$ is an $H^n$-tuple, we set
\begin{equation}
\begin{split}
\label{bipP1}
\big\langle P_1(\bth_{h_1})\otimes P_2(\bth_{h_2}) \otimes \cdots \otimes P_n(\bth_{h_n}) \cdot \boldsymbol{F}\big\rangle =
\sum_{g_1,\ldots,g_n =1}^H \bigg(\prod_{i = 1}^n
P_i(\theta_{h_i,g_i})\bigg)\,F_{g_1,\ldots,g_n}.
\end{split}
\end{equation}
For instance, if $\boldsymbol{F} = (F_{g,h})_{g,h = 1}^{H}$ is an $H^2$-tuple we have
\begin{equation*}
\begin{split}
\big\langle \bth_{h_1} \otimes \bth_{h_2} \cdot \boldsymbol{F}\big\rangle & = \sum_{g_1,g_2 = 1}^{H} \theta_{h_1,g_1} \theta_{h_2,g_2}\, F_{g_1,g_2}, \\
\big\langle \bth_h^{\otimes 2} \cdot \boldsymbol{F} \big\rangle & = \sum_{g_1,g_2=1}^{H} \theta_{h,g_1} \theta_{h,g_2} \,F_{g_1,g_2}.
\end{split}
\end{equation*}
If we want to insist on the indices that have to be summed up, we can also write these expressions as $\langle \bth_{h_1} \otimes \bth_{h_2} \cdot \boldsymbol{F}_{\bullet,\bullet} \rangle$ and $\langle \bth_{h}^{\otimes 2} \cdot \boldsymbol{F}_{\bullet,\bullet} \rangle$. When doing so, it is implicit that the index that is summed over and appearing in the polynomial $P_i$ in \eqref{bipP1} is the one replacing the $i$-th $\bullet$ starting from the left. Another example is
\[
\big\langle \bth_h^{\otimes 2} \cdot \boldsymbol{F}_{h,\bullet,\bullet} \big \rangle = \sum_{g_1,g_2 = 1}^{H} \theta_{h,g_1}\theta_{h,g_2}\,F_{h,g_1,g_2}.
\]

\subsection{Master problems and the solution operator \texorpdfstring{$\boldsymbol{\Op}$}{Op}}

\label{TheopUpsec2}

We write the asymptotic expansions of $W_{n;h_1,\ldots,h_n}(z_1,\ldots,z_n)$ in terms of the solutions of a system of linear equations that we call \emph{master problem}. The solution is encoded in an operator $\boldsymbol{\Op}$, which already appeared in the proofs of Section~\ref{Section_Smoothnessparam}. This operator is constructed and studied in details in Chapter~\ref{Chapter_SolvingN}, where a close relation between the master problem and a Riemann--Hilbert problem is also discussed. The master problem that we need to study is specified by the following data.
\begin{itemize}
\item[(i)] The positive integer $H$ and the real symmetric matrix $\boldsymbol{\Theta}$ of size $H$.
\item[(ii)] Pairwise disjoint real segments $\amsmathbb{A}^{\mathfrak{m}}_h := [\hat{a}_h^{\mathfrak{m}},\hat{b}_h^{\mathfrak{m}}]$ and sub-segments $[\alpha_h,\beta_h] \subset (\hat{a}_h^{\mathfrak{m}},\hat{b}_h^{\mathfrak{m}})$ indexed by $h \in [H]$.
\item[(iii)] Pairwise disjoint open subsets $\amsmathbb{M}_h \subset \amsmathbb{C}$ indexed by $h \in [H]$ and containing $\amsmathbb{A}^{\mathfrak{m}}_h$ but not $\bigcup_{g \neq h} \amsmathbb{A}_g^{\mathfrak{m}}$.
\item[(iv)] An $H$-tuple of functions $\boldsymbol{E}(z)=(E_h(z))_{h=1}^H$ such that $E_h$ is holomorphic in
$\amsmathbb{M}_h\setminus \amsmathbb{A}_h^{\mathfrak{m}}$ for any $h \in [H]$.
\item[(v)] An $H$-tuple of real numbers $\boldsymbol{\kappa} = (\kappa_h)_{h = 1}^H$.
\end{itemize}
\begin{definition}
\label{Definition_master_prob}
Given these data, the master problem asks to find an $H$-tuple of functions $\boldsymbol{F}(z)=(F_h(z))_{h=1}^H$ having the following properties for any $h \in [H]$.
\begin{enumerate}
\item $F_h(z)$ is holomorphic in $\amsmathbb C\setminus \amsmathbb{A}_h^{\mathfrak{m}}$.
\item There exists a holomorphic function $A_h(z)$ of $z \in \amsmathbb{M}_h$ such that
\begin{equation}
\label{eq_Master_equationpre}
 A_h(z) = \sqrt{(z-\alpha_h) (z-\beta_h)} \,\big\langle \bth_h \cdot \boldsymbol{F}(z)\big\rangle + E_h(z).
 \end{equation}
 \item $F_h(z) = \frac{\kappa_h}{z}+O\big(\frac{1}{z^2}\big)$ as $z \rightarrow \infty$. This can be equivalently restated as
\[
F_h(z)\,\, \mathop{=}_{z \rightarrow \infty} \,\, O\bigg(\frac{1}{z}\bigg) \qquad \textnormal{and} \qquad \oint_{\gamma_h} \frac{\dd z}{2\ii\pi}\, F_h(z) = \kappa_h,
\]
where $\gamma_h \subset \amsmathbb{M}_h$ is a contour surrounding $\amsmathbb{A}_h^{\mathfrak{m}}$ and having counterclockwise orientation.
\end{enumerate}
\end{definition}
The only information we have about $A_h(z)$ is that it is holomorphic on $\amsmathbb{M}_h$. This is a non-trivial property because $E_h$ is only holomorphic on $\amsmathbb{M}_h \setminus \amsmathbb{A}_h^{\mathfrak{m}}$, and the key result is that this property is sufficient to reconstruct uniquely the function $\boldsymbol{F}(z)$.

\begin{theorem}
\label{Theorem_Master_equation}
Consider data specifying a master problem. Assume Conditions 1.,2.,3. of Assumption~\ref{Assumptions_Theta} for these data, and additionally assume the existence of a constant $C > 0$ for which the distance between any pair of points in $\bigcup_{h = 1}^{H} \big\{\hat{a}_h^{\mathfrak{m}},\alpha_h,\beta_h,\hat{b}^{\mathfrak{m}}_h\big\}$ is larger than $\frac{1}{C}$ and the absolute value of all these points is smaller than $C$. Then, the master problem \eqref{eq_Master_equationpre} has a unique solution
 \begin{equation}
 \label{eq_solution_Master_equationpre}
 \forall h \in [H]\quad \forall z \in \amsmathbb{C} \setminus \amsmathbb{A}_h^{\mathfrak{m}}  \qquad F_h(z) = \Op_h\big[\boldsymbol{E}\,;\,\boldsymbol{\kappa}\big](z).
 \end{equation}
 For any $h \in [H]$, the solution operator $\Op_h$ is linear in $\boldsymbol{E}$ and $\boldsymbol{\kappa}$, depends smoothly on the parameters $\boldsymbol{\Theta}$ and $\boldsymbol{\alpha},\boldsymbol{\beta}$, and depends on the segments $(\amsmathbb{A}_g^{\mathfrak{m}})_{g = 1}^{H}$ and the domains $(\amsmathbb{M}_g)_{g = 1}^{H}$ only through the domain of definition (the precise meaning is explained in Chapter~\ref{Chapter_SolvingN}).

 Besides, it is continuous in the following sense. For any $2H$-tuple of \label{index:Kgset}compact subsets $(\amsmathbb{K}_g,\amsmathbb{K}_g^-)_{g = 1}^{H}$ such that
 \[
 \forall g \in [H]\qquad \amsmathbb{K}_g^- \subset \mathring{\amsmathbb{K}}_g \subset \amsmathbb{K}_{g} \subset \amsmathbb{C} \setminus \amsmathbb{A}_g^{\mathfrak{m}},
 \]
 there exists a constant $C' > 0$ depending only on these compacts and the constants in the assumptions such that, for any $\boldsymbol{E},\boldsymbol{\kappa}$ as in (iv)-(v), we have
 \begin{equation}
\label{ContinuityUp_eqpre}\max_{h \in [H]} \sup_{z \in \amsmathbb{K}_h} \big|\Op_h\big[\boldsymbol{E}\,;\,\boldsymbol{\kappa}\big](z)\big| \leq C' \cdot \Big( \max_{h \in [H]} \sup_{z \in \amsmathbb{K}_h^-} |E_h(z)| + |\!|\boldsymbol{\kappa}|\!|_{\infty}\Big).
\end{equation}
Besides, if $E_h(z)$ is holomorphic in $\amsmathbb{M}_h$ for each $h \in [H]$, then
\[
\forall h \in [H]\qquad \Upsilon_h[\boldsymbol{E}\,;\,\boldsymbol{0}] = 0.
\]
\end{theorem}
\begin{proof}
The theorem is restated in Chapter~\ref{Chapter_SolvingN} as Theorem~\ref{Theorem_Master_equation_12} and proved there.
\end{proof}

We can use the symbol * as a placeholder for the variable of the function the operator acts on:
\[
\forall h \in [H]\qquad \Upsilon_h\big[\boldsymbol{E}\,;\,\boldsymbol{\kappa}\big](z) = \Upsilon_h\big[\boldsymbol{E}(*)\,;\,\boldsymbol{\kappa}\big](z).
\]
This will be helpful when we apply $\Upsilon_h$ to functions that depend on several variables. For instance, if $\boldsymbol{E}(z_1,z_2)$ is a function of two variables, we write $\Upsilon_h\big[\boldsymbol{E}(*,z_2)\,;\,\boldsymbol{\kappa}\big]$ to indicate the application of the operator $\Upsilon_h$ to the function $z \mapsto \boldsymbol{E}(z,z_2)$ for fixed $z_2$, and $\Upsilon_h\big[\boldsymbol{E}(*,*)\,;\,\boldsymbol{\kappa}\big]$ to indicate the application of the operator $\Upsilon_h$ to the function $z \mapsto \boldsymbol{E}(z,z)$.

For general matrices $\boldsymbol{\Theta}$, there is no explicit formula for the solution operator $\boldsymbol{\Op}$. In Section~\ref{Section_Master_by_covariance}, we show that for certain classes of $\boldsymbol{E}(z)$ it admits an integral
representation in terms of the \emph{fundamental} solution $\boldsymbol{\mathcal{F}}(z_1,z_2) = \big(\mathcal{F}_{h_1,h_2}(z_1,z_2)\big)_{h_1,h_2 = 1}^{H}$ defined by
\begin{equation}
\label{eq_covariancepre}
\forall h_1,h_2 \in [H] \qquad \mathcal{F}_{h_1,h_2}(z_1,z_2) =
\Op_{h_1}\Bigg[\frac{\sqrt{(*-\alpha_{h_2})(*-\beta_{h_2})}}{2(*- z_2)^2}\boldsymbol{e}^{(h_2)}\,;\,\boldsymbol{\kappa} = \boldsymbol{0}\Bigg](z_1),
\end{equation}
An important property of the fundamental solution proved in Theorem~\ref{thm:Bsym} is its symmetry
\[
\mathcal{F}_{h_1,h_2}(z_1,z_2) = \mathcal{F}_{h_2,h_1}(z_2,z_1)
\]
For a general $\boldsymbol{\Theta}$ there is no explicit formula for $\boldsymbol{\mathcal{F}}(z_1,z_2)$, we only know that it depends (smoothly) on $\boldsymbol{\Theta},\boldsymbol{\alpha},\boldsymbol{\beta}$ and does not depend on the segments $(\amsmathbb{A}_h^{\mathfrak{m}})_{h = 1}^{H}$ and the complex domains $(\amsmathbb{M}_{g})_{g = 1}^{H}$. For simple matrices $\boldsymbol{\Theta}$ (\textit{cf.} Section~\ref{Section_Explicit}) and matrices $\boldsymbol{\Theta}$ coming from the tiling models (\textit{cf.} Theorems~\ref{thm:Omegav} and \ref{thm:spcurvenonbip}), the fundamental solution can admit a more explicit form or there is at least an algorithm to construct it in an algebraic way. This will be further illustrated in Chapter~\ref{Chap14}.

The data of a discrete ensemble satisfying Assumptions~\ref{Assumptions_Theta}, \ref{Assumptions_basic}, \ref{Assumptions_offcrit}, \ref{Assumptions_analyticity} and having a single band per segment specifies the \textit{a priori} $\N$-dependent data (i)-(ii)-(iii) for a master problem. The segments are given in Definition~\ref{def:eq_shifted_parameters}-\ref{def:eq_enlarged_parameters}. For each $h \in [H]$ we take $(\alpha_h,\beta_h)$ in (ii) to be the band in the $h$-th segment $[\hat{a}_h',\hat{b}_h']$. The $H$-tuples $\boldsymbol{E}(z)$ and $\boldsymbol{\kappa}$ play the role of source terms in the master problem, and they will be specified in each instance where we need to solve master problems. The assumptions of Theorem~\ref{Theorem_Master_equation} are valid due to Assumptions~\ref{Assumptions_basic} and \ref{Assumptions_offcrit}.

\subsection{Auxiliary functions}
\label{Subsection_Auxiliary_functions}

We collect here the definitions of several
useful functions depending on the data of a discrete ensemble. They are defined in terms of
the Stieltjes transform of the equilibrium measure $\boldsymbol{\mu}$ and the functions appearing in the expansion of $\Phi_h^{\pm}$ in Definitions~\ref{Definition_phi_functions} and \ref{Def_finer_Phi} and involved in Assumption~\ref{Assumptions_analyticity}.
\[
\Phi^{\pm}_{h}(z) = \phi_h^{\pm}(z) + \frac{\phi_h^{\pm,[1]}(z)}{\N} + \frac{\phi_h^{\pm,[2]}(z)}{\N^2} + \frac{\mathbbm{f}_h^{\pm}(z)}{\N^3}.
\]
Uniform bounds for these functions were provided in Lemma~\ref{Lemma_phi_properties}. The following functions $Q^\pm_h(z)$ appear in the large $\N$ asymptotic expansion of the Nekrasov equation of Theorems~\ref{Theorem_Nekrasov} and Corollary~\ref{Corollary_higherNek}.

\begin{equation*}
\begin{split}
 & \quad Q_h^{\pm}(z) \\
 & :=\Phi_{h}^{+}(z)\cdot\exp\bigg(\big\langle \bth_{h}\cdot \Gm_{\boldsymbol{\mu}}(z)\big\rangle + \Big\langle \frac{\bth_h(\bth_h - 1)}{2\N}\cdot \partial_{z}\Gm_{\boldsymbol{\mu}}(z)\Big\rangle + \Big\langle \frac{\bth_h(4\bth_h^2 - 6\bth_h + 3)}{24\N^2}\cdot \partial_{z}^2\Gm_{\boldsymbol{\mu}}(z)\Big\rangle \bigg) \\
 & \quad \pm \Phi_{h}^{-}(z)\cdot\exp\bigg(-\big\langle \bth_{h}\cdot \Gm_{\boldsymbol{\mu}}(z)\big\rangle + \Big\langle \frac{\bth_h(\bth_h - 1)}{2\N}\cdot \partial_{z}\Gm_{\boldsymbol{\mu}}(z)\Big\rangle - \Big\langle\frac{\bth_h(4\bth_h^2 - 6\bth_h + 3)}{24\N^2}\cdot \partial_{z}^2\Gm_{\boldsymbol{\mu}}(z)\Big\rangle\bigg).
\end{split}
\end{equation*}
Further expanding $Q_h^\pm(z)$ in the powers of $\tfrac{1}{\N}$ we decompose
\begin{equation}
\label{eq_Q_plus_determ_expansion}
Q_h^{\pm}(z)= q_h^{\pm}(z) + \frac{q_h^{\pm,[1]}(z)}{\N} + \frac{q_h^{\pm,[2]}(z)}{\N^2} + \frac{\mathbbm{q}^{\tau}_h(z)}{\N^3}.
\end{equation}
The first term has previously appeared in \eqref{eq_q_pm}
\begin{equation}
\label{eq_x95_1}
q_h^{\pm}(z) = \phi_{h}^{+}(z)\cdot e^{\langle \bth_{h}\cdot \Gm_{\boldsymbol{\mu}}(z)\rangle} \pm \phi_{h}^{-}(z)\cdot e^{-\langle\bth_{h}\cdot \Gm_{\boldsymbol{\mu}}(z)\rangle}.
\end{equation}
The next two terms are first defined here:
\begin{equation}
\label{eq_x215}
\begin{split}
q_h^{\pm,[1]}(z) & = \bigg(\frac{\phi_h^{+}(z)}{2} \big\langle\bth_h(\bth_h - 1) \cdot \partial_{z}\Gm_{\boldsymbol{\mu}}(z)\big\rangle + \phi_h^{+,[1]}(z)\bigg)e^{\langle \bth_{h}\cdot \Gm_{\boldsymbol{\mu}}(z)\rangle} \\
&\quad \pm \bigg(\frac{\phi_h^{-}(z)}{2} \big\langle\bth_h(\bth_h - 1) \cdot \partial_{z}\Gm_{\boldsymbol{\mu}}(z)\big\rangle + \phi_h^{-,[1]}(z)\bigg) e^{-\langle\bth_{h}\cdot \Gm_{\boldsymbol{\mu}}(z)\rangle}, \\
q_h^{\pm,[2]}(z) & = \bigg(\frac{\phi_h^{+}(z)}{24}\Big[\big\langle \bth_h(4\bth_h^2 - 6\bth_h + 3) \cdot \partial_{z}^2\Gm_{\boldsymbol{\mu}}(z)\big\rangle + 3\big\langle\bth_h(\bth_h - 1)\cdot \partial_{z}\Gm_{\boldsymbol{\mu}}\big\rangle^2\Big] \\& \qquad + \frac{\phi_h^{+,[1]}(z)}{2}\big\langle\bth_h(\bth_h - 1)\cdot \partial_{z}\Gm_{\boldsymbol{\mu}}(z)\big\rangle + \phi_h^{+,[2]}(z)\bigg)e^{\langle\bth_h\cdot \Gm_{\boldsymbol{\mu}}(z)\rangle} \\& \quad \pm \bigg(\frac{\phi_h^{-}(z)}{24}\Big[-\big\langle\bth_h(4\bth_h^2 - 6\bth_h + 3) \cdot \partial_{z}^2\Gm_{\boldsymbol{\mu}}(z)\big\rangle + 3\big\langle\bth_h(\bth_h - 1)\cdot \partial_{z}\Gm_{\boldsymbol{\mu}}\big\rangle^2\Big] \\
& \qquad\quad + \frac{\phi_h^{-,[1]}(z)}{2}\,\big\langle\bth_h(\bth_h - 1)\cdot \partial_{z}\Gm_{\boldsymbol{\mu}}(z)\big\rangle + \phi_h^{-,[2]}(z)\bigg)e^{-\langle\bth_h\cdot \Gm_{\boldsymbol{\mu}}(z)\rangle}.
\end{split}
\end{equation}
Alternative expressions are
\begin{equation}
\label{qhsmallexp}
\begin{split}
q_h^{\pm,[1]}(z) & = \frac{q_h^{\pm}(z)}{2}\big\langle\bth_h(\bth_h - 1) \cdot \partial_z \Gm_{\boldsymbol{\mu}}(z) \big\rangle + \phi_h^{+,[1]}(z)\cdot e^{\langle\bth_h \cdot \Gm_{\boldsymbol{\mu}}(z)\rangle} \pm \phi_h^{-,[1]}(z) \cdot e^{-\langle \bth_h \cdot \Gm_{\boldsymbol{\mu}}(z)\rangle}, \\
q_h^{\pm,[2]}(z) & = \frac{q_h^{\mp}(z)}{24}\big\langle \bth_h(4\bth_h^2 - 6\bth_h + 3) \cdot \partial_z^2 \Gm_{\boldsymbol{\mu}}(z)\big\rangle + \frac{q_h^{\pm}(z)}{8} \big\langle \bth_h(\bth_h - 1)\cdot \partial_z \Gm_{\boldsymbol{\mu}}(z)\big\rangle^2 \\
& \quad + \frac{1}{2} \big\langle \bth_h(\bth_h - 1) \cdot \partial_z \Gm_{\boldsymbol{\mu}}(z)\big\rangle \left( \phi_h^{+,[1]}(z) e^{\langle \bth_h \cdot \Gm_{\boldsymbol{\mu}}(z)\rangle} \pm \phi_h^{-,[1]}(z) e^{-\langle \bth_h \cdot \Gm_{\boldsymbol{\mu}}(z)\rangle}\right) \\
& \quad + \phi_h^{+,[2]}(z) \cdot e^{\langle \bth_h \cdot \Gm_{\boldsymbol{\mu}}(z)\rangle}  \pm \phi_h^{-,[2]}(z) \cdot e^{-\langle \bth_h \cdot \Gm_{\boldsymbol{\mu}}(z)\rangle}.
\end{split}
\end{equation}
The exact form of the remainder $\mathbbm{q}^{\pm}_h(z)$ is irrelevant for our purposes, it is only important to know that it is uniformly bounded as $\N\rightarrow\infty$ for $z$ belonging to compact subsets of $\amsmathbb{M}_{h}\setminus \amsmathbb{A}_{h}^{\mathfrak{m}}$ for any $h \in [H]$. The relevance of the functions $Q_h^\pm(z)$ and their asymptotic expansions comes from the fact that, for any $h \in [H]$, the expressions under the expectations in the Nekrasov equation of Theorem~\ref{Theorem_Nekrasov} can be transformed into

\begin{equation}
\label{expobserv0}
\begin{split}
 & \quad \prod_{i = 1}^N \bigg(1 \pm \frac{1}{\N}\,\frac{\theta_{h,h(i)}}{z - \frac{\ell_i}{\N} \mp \frac{1}{2\N}}\bigg) \\
& = \exp\bigg(\pm \big\langle \bth_h\cdot \Gm_{\boldsymbol{\mu}}\big\rangle + \frac{1}{2\N}\big\langle\bth_h(\bth_h - 1) \cdot \partial_{z}\Gm_{\boldsymbol{\mu}}\big\rangle \pm \frac{1}{24\N^2}\big\langle\bth_h(4\bth_h^2 - 6\bth_h + 3)\cdot \partial_{z}^2\Gm_{\boldsymbol{\mu}}\big\rangle\bigg) \\
 & \quad \times \big(\textnormal{polynomial in }\Delta G(z)\big) +O(\N^{-\frac{5}{2}+\eps}).
\end{split}
\end{equation}
This is explained in more detail in Section~\ref{Section_asymptotic_Nekrasov} below.

\medskip

Last, we recall the expression of the auxiliary function $s_h(z)$ from Definition~\ref{GQdef2}. In the specific context of Assumption~\ref{Assumptions_extra} (see \eqref{eq_s_first_appearance} or Remark~\ref{remextrasssh}) and with only one band $(\alpha_h,\beta_h)$ in the $h$-th segment, this is
\begin{equation}
\label{eq_s_again}
s_h(z) = \frac{q_h^-(z)}{ \sqrt{(z-\alpha_{h})(z-\beta_{h})} \cdot (z-\hat a'_h)^{\mathbbm{1}_{\amsmathbb{S}_h}(\hat{a}'_{h})} \cdot (z-\hat b'_h)^{\mathbbm{1}_{\amsmathbb{S}_h}(\hat{b}'_{h})}}.
\end{equation}

We have established in Theorem~\ref{Theorem_regularity_density} the key properties of the functions $q_h^{\pm}(z)$, namely that $q_h^{+}(z)$ and $s_h(z)$ are holomorphic up to possible singularities at the endpoints $\hat{a}_h',\hat{b}'_h$. But, we will use those functions in a specific context where such singularities are absent. To state these properties it is helpful to introduce a modified version of $s_h(z)$, without the linear factors in the denominator.
\begin{definition}
\label{GQdef3pre} We introduce
\begin{equation}
\label{eq_s_circ}
s^{\circ}_h(z):=  \frac{q_h^-(z)}{\sqrt{(z-\alpha_{h})(z-\beta_{h})}}= (z-\hat a'_h)^{\mathbbm{1}_{\amsmathbb{S}_h}(\hat{a}'_{h})} \cdot (z-\hat b'_h)^{\mathbbm{1}_{\amsmathbb{S}_h}(\hat{b}'_{h})}\cdot s_h(z).
\end{equation}
\end{definition}

\begin{lemma}
\label{Lemm84}
Under Assumptions~\ref{Assumptions_Theta}, \ref{Assumptions_basic}, \ref{Assumptions_offcrit}, \ref{Assumptions_analyticity} and \ref{Assumptions_extra},  for any $h \in [H]$ the functions $q_h^+(z)$, $s_h^{\circ}(z)$ and $q_h^+(z)/s_h^{\circ}(z)$ are holomorphic functions of $z$ in a complex neighborhood of $\amsmathbb{A}_h^{\mathfrak{m}}$.
\end{lemma}
\begin{proof}
We know by Theorem~\ref{Theorem_regularity_density} that $s_h(z)$ is holomorphic and has no zeros in a small neighborhood of $\amsmathbb{A}_h^{\mathfrak{m}}$. Accordingly, the only possible zeros of $s_h^{\circ}(z)$ in this neighborhood are at $z = \hat{c}_h'$ with $\hat{c}_h'  \in\{\hat{a}_h',\hat{b}_h'\}$, coming from the factors we multiplied $s_h(z)$ with in \eqref{eq_s_circ}. More precisely, there is a simple zero at $\hat{c}_h'$ if this point is saturated, otherwise $s_h^{\circ}(\hat{c}_h') \neq 0$.

We also know by Theorem~\ref{Theorem_regularity_density} that $q_h^+(z)$ is holomorphic in a small neighborhood of $\amsmathbb{A}_h^{\mathfrak{m}}$, except perhaps for simple poles at $\hat{a}_h'$ when $\iota_h^- = 0$, or at $\hat{b}_h'$ when $\iota_h^+ = 0$. Actually, if $\iota_h^- = 0$ Assumption~\ref{Assumptions_extra} forces $\hat{a}_h'$ to be in a void. Then $\mathcal{G}_{\mu_h}(z)$ is holomorphic near $z = \hat{a}_h'$, and so must $q_h^{\pm}(z)$ be. The other situation tolerated by Assumption~\ref{Assumptions_extra} is $\iota_h^+ = 2$, in which case $\hat{a}_h'$ is the endpoint of a saturation, creating a simple zero for $\exp\big(\theta_{h,h}\mathcal{G}_{\mu_h}(z)\big)$  at $z = \hat{a}_h'$. Then, in the expression
\[
q_h^+(z) = \phi_h^+(z) \cdot e^{\langle \bth_h \cdot \mathcal{G}_{\boldsymbol{\mu}}(z) \rangle} + \phi_h^-(z) \cdot e^{-\langle\bth_h \cdot \mathcal{G}_{\boldsymbol{\mu}}(z) \rangle}
\]
the first term has a simple zero from the exponential and the second term has a simple zero coming from the double zero of $\phi_h^-(z)$ compensated by the simple pole of the exponential. This creates a simple zero for $q_h^+(z)$. We conclude that in any case $q_h^+(z)$ is holomorphic near $z = \hat{a}_h'$, and that in the situation where $s_h^{\circ}(z)$  has a zero at $z = \hat{a}_h'$ it cancels with the zero of $q_h^+(z)$ to make the ratio $q_h^+(z)/s_h^{\circ}(z)$ holomorphic near $z = \hat{a}_h'$. A similar argument shows that $q_h^+(z)/s_h^{\circ}(z)$ is holomorphic near $z = \hat{b}_h'$.
\end{proof}

\section{Main results for the correlators}

\label{Section_cumulant_expansion}

We now present our asymptotic results for the cumulants in the fixed filling fraction ensembles with one band per segment. The presentation is split into two parts: we first work under additional Assumption~\ref{Assumptions_extra}, and then remove it using the conditioning procedure of Chapter~\ref{Chapter_conditioning}, at the expense of some of the formulae becoming less explicit.

\subsection{Results under Assumption~\ref{Assumptions_extra}}
\label{Section_under_5}
All the notations appearing in the following theorems have been introduced in Section~\ref{Section_notations}.

\begin{theorem}
\label{Theorem_correlators_expansion} Consider a discrete ensemble satisfying Assumptions~\ref{Assumptions_Theta}, \ref{Assumptions_basic}, \ref{Assumptions_offcrit}, \ref{Assumptions_analyticity} and the additional Assumption~\ref{Assumptions_extra}. For $h,h_1,h_2 \in [H]$, we set
\begin{equation}
\label{eq:8585}
\begin{split}
W_{2;h_1,h_2}^{[0]}(z_1,z_2) & = \mathcal{F}_{h_1,h_2}(z_1,z_2), \\
W_{1; h}^{[1]}(z) & = \Op_{h}\Bigg[\bigg(\frac{q^{+,[1]}_{g}}{s_g^{\circ}}\bigg)_{g = 1}^H\,;\,\boldsymbol{0}\Bigg](z_1), \\
W_{1;h}^{[2]}(z) & = \Op_{h}\Bigg[\bigg(\frac{1}{s_g^{\circ}}\Big\{q^{+,[2]}_{g} + q_g^{-,[1]} \big\langle\bth_g \cdot \boldsymbol{W}_{1}^{[1]}\big\rangle + \frac{q_g^+}{2}\big\langle\bth_g\cdot \boldsymbol{W}_{1}^{[1]}\big\rangle^{2} \\
& \qquad\qquad + \frac{q_g^+}{2} \big\langle \bth_g(\bth_g - 1)\cdot \partial \boldsymbol{W}_{1}^{[1]}\big\rangle
+ \frac{q_g^+}{2} \big\langle \bth_g^{\otimes 2}\cdot \boldsymbol{\mathcal{F}}_{\bullet,\bullet}(*,*) \big\rangle \Big\}\bigg)_{g =
1}^H;\, \boldsymbol{0}\Bigg](z).
\end{split}
\end{equation}
where $z \in \amsmathbb{M}_{h}$ and $z_i \in \amsmathbb{M}_{h_i}$ for $i \in \{1,2\}$.

Then, for any $\eps >0$ and compact sets $\amsmathbb{K}_h \subset \amsmathbb C\setminus [\hat a^{\mathfrak m}_h,\hat b^{\mathfrak m}_h]$ indexed by $h \in [H]$, there exist a sequence of positive constants $(C_n)_{n \geq 1}$ depending only $\eps$, the compacts, and the constants in the assumptions, such that for any $n \in \amsmathbb{Z}_{> 0}$, $h_1,\ldots,h_n \in [H]$ and $z_i \in \amsmathbb{K}_{h_i}$ for $i \in [n]$, we have
\begin{align}
\label{eq_x30} \big|W_{1;h}(z) - W_{1;h}^{[1]}(z) - \N^{-1}W_{1;h}^{[2]}(z) \big| & \leq C_1\,\N^{-\frac{3}{2} + \varepsilon}, \\
\label{eq_x31} \max_{ l = 1,2} \big|W_{1;h}^{[l]}(z)\big| & \leq C_1, \\
\label{eq_x32} \big|W_{2;h_1,h_2}(z_1,z_2) - W_{2;h_1,h_2}^{[0]}(z_1,z_2)\big| & \leq
C_2\,\N^{-\frac{1}{2} + \varepsilon}, \\
\label{eq_x33} \big|W_{2;h_1,h_2}^{[0]}(z_1,z_2)\big| & \leq C_2
\end{align}
and for $n \geq 3$
\begin{equation}
\label{eq_x34} \big|W_{n;h_1,\ldots,h_n}(z_1,\ldots,z_n)\big| \leq C_n\,\N^{-\frac{1}{2} + \varepsilon}.
\end{equation}
\end{theorem}

In the case $\theta_{g,g} = 1$ for every $g \in [H]$ (as will be the case for tilings in Chapter~\ref{Chap11}) there are many simplifications. Among others, the first subleading correction is smaller by order of magnitude $\frac{1}{\N^2}$ (instead of $\frac{1}{\N}$) compared to the leading order.

\begin{corollary}
\label{Theorem_correlators_expansion_theta1} Assume in the setting of Theorem~\ref{Theorem_correlators_expansion} that we have $\theta_{g,g} = 1$ for every $g \in [H]$. Let $h \in [H]$ and set
\begin{equation}
\label{W12W12}
\hat{W}_{1;h}^{[2]}(z) = \Op_{h}\Bigg[\bigg(\frac{1}{s_g^{\circ}}\Big\{\hat{q}^{+,[2]}_{g} + \frac{q_g^+}{2} \big\langle \bth_g^{\otimes 2}\cdot \boldsymbol{\mathcal{F}}_{\bullet,\bullet}(*,*) \big\rangle \Big\}\bigg)_{g = 1}^{H}\,;\,\boldsymbol{0}\bigg],
\end{equation}
where
\begin{equation}
\label{W12W120}
\hat{q}_{h}^{+,[2]}(z) = \frac{q_h^-(z)}{24} \big\langle \bth_h(4\bth_h^2 - 6\bth_h + 3) \cdot \partial_z^2 \mathcal{G}_{\boldsymbol{\mu}}(z)\big\rangle -   \frac{\partial_z^3 U_h(z)}{12} \cdot \phi_h^{+}(z) \cdot e^{\langle \bth_h \cdot \mathcal{G}_{\boldsymbol{\mu}}(z)\rangle}.
\end{equation}
Then, \eqref{eq_x30} and \eqref{eq_x31} can be replaced with
\[
\big|W_{1;h}(z) - \N^{-1}\hat{W}_{1;h}^{[2]}(z)\big| \leq C_1 \N^{-\frac{3}{2} + \varepsilon}, \qquad \big|\hat{W}_{1;h}^{[2]}(z)\big| \leq C_1.
\]
\end{corollary}

To state our central limit theorem we will use the following mode of convergence.
\begin{definition}
\label{Definition_convergence_moments} Let $L$ be a positive integer. If $\boldsymbol{\chi}^{(1)}$ and $\boldsymbol{\chi}^{(2)}$ are two random $L$-dimensional vectors depending on a large parameter $\N$ and having moments of all orders, we say that $\boldsymbol{\chi}^{(1)}$ approximates $\boldsymbol{\chi}^{(2)}$ in the sense of moments if
\[
\forall k \in \amsmathbb{Z}_{> 0} \quad \forall l_1,\ldots,l_k \in [M]\qquad \lim_{\N \rightarrow \infty} \bigg|\amsmathbb{E}\bigg[\prod_{j = 1}^{k} \chi_{l_j}^{(1)}\bigg] - \amsmathbb{E}\bigg[\prod_{j = 1}^{k} \chi_{l_j}^{(2)}\bigg]\bigg| = 0.
\]
In this case we write $\boldsymbol{\chi}^{(1)} \mathop{\sim}^{\text{mom}} \boldsymbol{\chi}^{(2)}$. In this definition, we do not require the convergence to be uniform in $k,l_1,\ldots,l_k$.
\end{definition}

\begin{corollary}[Central limit theorem]
\label{Corollary_CLT} Consider a discrete ensemble satisfying Assumptions~\ref{Assumptions_Theta}, \ref{Assumptions_basic}, \ref{Assumptions_offcrit}, \ref{Assumptions_analyticity} and the additional Assumption~ \ref{Assumptions_extra}. Let $L \in \amsmathbb{Z}_{> 0}$ and a $L$-tuple of integers $\boldsymbol{h} \in [H]^L$ both independent of $\N$, let a (possibly $\N$-dependent) $L$-tuple of functions $\boldsymbol{f}(z)$ such that $f_l(z)$ is a holomorphic function of $z \in \amsmathbb{M}_{h_l}$ for any $l \in [L]$. Assume there exists a constant $C > 0$ such that $\max_{l} \sup_{z} |f_l(z)| \leq C$.

We define the random $L$-dimensional Gaussian vector $\textnormal{\textbf{\textsf{Gau\ss{}}}}[\boldsymbol{f}]$ with mean and covariance
\begin{equation}
\label{eq_x144}
\begin{split}
\textnormal{\textsf{Mean}}_l[\boldsymbol{f}] & = \oint_{\gamma_{h_l}} \frac{ \dd z}{2\ii\pi}\,
f_l(z)\,\,W_{1;h_l}^{[1]}(z) , \\
\textnormal{\textsf{Cov}}_{l_1,l_2}[\boldsymbol{f},\boldsymbol{f}] & = \oint_{\gamma_{h_{l_1}}} \oint_{\gamma_{h_{l_2}}} \frac{\dd z_1\,\dd z_2}{(2\ii\pi)^2}\,
 f_{l_1}(z_1) f_{l_2}(z_2)
\,\mathcal{F}_{h_{l_1},h_{l_2}}(z_1,z_2)
\end{split}
\end{equation}
If $\boldsymbol{\ell}$ is $\P$-distributed with filling fractions deterministically as per Assumption~\ref{Assumptions_extra}, then as $\N \rightarrow \infty$
\begin{equation}
\label{eq_x142}
\Bigg(\sum_{i = 1}^{N_{h_l}} f_{l}\bigg(\frac{\ell^{h_l}_{i}}{\N}\bigg) - \N
\int_{\hat{a}_{h_l}'}^{\hat{b}_{h_l}'} f_{l}(x)\mu_{h_l}(x)\dd x\Bigg)_{l=1}^L \mathop{\sim}^{\textnormal{mom}} \quad \textnormal{\textbf{\textsf{Gau\ss{}}}}[\boldsymbol{f}],
\end{equation}
and the convergence in the sense of Definition~\ref{Definition_convergence_moments} is uniform for fixed choice of constants in the assumptions. If furthermore $\theta_{g,g} = 1$ for every $g \in [H]$, the result holds with $\textnormal{\textsf{Mean}}_l[\boldsymbol{f}] = 0$ for every $l \in [L]$.
\end{corollary}
On the one hand, we stress that the covariance expression through the matrix of functions $\boldsymbol{\mathcal{F}}$ given by \eqref{eq_covariancepre} depends only on the matrix $\boldsymbol{\Theta}$ and the endpoints of the bands $(\alpha_h,\beta_h)_{h = 1}^{H}$. On the other hand, the mean depends on the parameters of the ensemble in a more delicate way. It is zero if the diagonal of $\boldsymbol{\Theta}$ is $1,\ldots,1$ only because we have defined the equilibrium measure using the potential $V$ of \eqref{eq_potential_weight_match} and the \emph{shifted} segment $[\hat{a}'_h,\hat{b}'_h]$ of Definition~\ref{def:eq_shifted_parameters}. This contrasts with \cite[Remark 2.7]{BGG}, where a different definition of the equilibrium measure was used, leading to a non-zero mean even in the $H = 1$ and $\theta_{1,1} =1$ case treated there.

Recall that a random $L$-dimensional vector $\boldsymbol{\xi}$ is Gaussian if and only if for each deterministic $\boldsymbol{t} \in \amsmathbb{R}^L$ the scalar random variable $\sum_{l = 1}^{L} t_l \xi_l$ is Gaussian. Consequently, it is sufficient to prove Corollary~\ref{Corollary_CLT} for $L = H$ and $h_l = l$ for any $l$. The more sophisticated presentation of Corollary~\ref{Corollary_CLT} serves the direct comparison with the corresponding statements of central limit theorems in Chapter~\ref{Chapter_filling_fractions}, where the asymptotic distribution is no longer Gaussian and such a simplification is not possible.

\subsection{Comments and proofs of Corollaries~\ref{Theorem_correlators_expansion_theta1} and \ref{Corollary_CLT}}

Theorem~\ref{Theorem_correlators_expansion} is proved in Sections~\ref{Section_asymptotic_Nekrasov}-\ref{secprosofggg}. The strategy is to write down the $n$-th order Nekrasov equation of Corollary~\ref{Corollary_higherNek} and consider it as a functional relation between $\boldsymbol{W}_{n}$ and all other cumulants. Analyzing recursively in $n$ this tower of equations, we will justify that other cumulants actually give negligible contributions compared to the leading order of $\boldsymbol{W}_{n}$, and moving the error term to the right-hand side we find a linear equation for the leading asymptotics of $\boldsymbol{W}_{n}$. The equation is of the form discussed in Section~\ref{TheopUpsec} and it can be solved by an application of Theorem~\ref{Theorem_Master_equation}, using the operator $\boldsymbol{\Op}$. This is what we will call repeatedly ``solving Nekrasov equation''. For instance, we will find that $\boldsymbol{W}_{2}(z,z_2)$ satisfies for any $h,h_2 \in [H]$ and $(z,z_2) \in \amsmathbb{M}_h \times \amsmathbb{M}_{h_2}$
\begin{equation}
\label{eq_x140}
\frac{q_{h}^-(z)\delta_{h,h_2} }{2(z_2 - z)^2} + q_{h}^{-}(z) \,\big\langle \bth_{h}\cdot \boldsymbol{W}_{2;\bullet,h_2}(z,z_2)\big\rangle \,\,\mathop{=}_{\N \rightarrow \infty}\,\, R_{h;h_2}(z;z_2) + o(1),
\end{equation}
where $R_{h;h_2}(z;z_2)$ is such that
\[
\frac{R_{h;h_2}(z;z_2)}{(z-\hat a'_h)^{\mathbbm{1}_{\amsmathbb{S}_h}(\hat{a}'_{h})} \cdot (z-\hat b'_h)^{\mathbbm{1}_{\amsmathbb{S}_h}(\hat{b}_h')}}
\]
is a holomorphic function of $z \in \amsmathbb{M}_h$, \textit{cf.} later Proposition~\ref{Proposition_Nek_2_asymptotic_form} combined with the bounds of Proposition~\ref{Proposition_uniform_moments_bound}. Recalling \eqref{eq_s_again} and \eqref{eq_s_circ}, we divide \eqref{eq_x140} by the function $s_h^\circ(z)$ of Definition~\ref{GQdef3pre} and get
\begin{equation}
\label{eq_x141}
\frac{\sqrt{(z-\alpha_{h})(z-\beta_{h})}}{2(z_2 - z)^2}\delta_{h,h_2} + \sqrt{(z-\alpha_{h})(z-\beta_{h})}\,\big\langle \bth_{h}\cdot \boldsymbol{W}_{2;\bullet,h_2}(z,z_2)\big\rangle \,\, \mathop{=}_{\N \rightarrow \infty}\,\, \widetilde{R}_{h;h_2}(z;z_2) + o(1),
\end{equation}
with $\widetilde R_{h;h_2}(z;z_2)$ holomorphic. It is important here that $s_h(z)$ has no zeros for $z \in \amsmathbb M_h$, which is the reason why we work with Assumption~\ref{Assumptions_extra}. Thus, \eqref{eq_x141} takes the form of a master problem, and the application of Theorem~\ref{Theorem_Master_equation} yields the identity $\boldsymbol{W}_{2}(z_1,z_2) = \boldsymbol{\mathcal{F}}(z_1,z_2) + o(1)$ in terms of the fundamental solution defined in Section~\ref{TheopUpsec2}.

An important step in our arguments is to show that the remainder in \eqref{eq_x140} is indeed $o(1)$ as $\N\rightarrow\infty$. The key feature that makes our analysis possible is that, even at the stage where not all negligible contributions have been recognized as negligible, solving the Nekrasov equation by Theorem~\ref{Theorem_Master_equation} gives improved error bounds, and doing so recursively allows reaching optimal bounds. For that it is essential to use the continuity of $\boldsymbol{\Op}$ provided by Theorem~\ref{Theorem_Master_equation}, \textit{i.e.} applying it to a negligible term still gives a negligible term. This systematic strategy was first described in \cite{BG11}, in the context of (continuous) $\sbeta$-ensembles. The order of expansion put forward in Theorem~\ref{Theorem_correlators_expansion} is
just what is needed to obtain a central limit theorem for fluctuations of linear
statistics and filling fractions. It will be clear from the proof that the method we employ can be pushed first to show $\boldsymbol{W}_n = O(\N^{2 - n})$ as $\N \rightarrow \infty$ for any $n \geq 2$, and then to obtain an all-order asymptotic expansion of the form
\[
\boldsymbol{W}_{n}(z_1,\ldots,z_n) \,\,\mathop{=}_{\N \rightarrow \infty}\,\, \sum_{k = \delta_{n,1}}^{K} \N^{2 - n - k} \boldsymbol{W}_{n}^{[k]}(z_1,\ldots,z_n) + o(\N^{2 - n - K})
\]
for any fixed $K > 0$ as $\N \rightarrow \infty$, where the $k$-th term depends on $\N$ but remains bounded. The expressions for the coefficients of expansion become increasingly cumbersome as $(n,k)$ increase, but they can be computed by induction on $n + k$. This fits in a general scheme of asymptotic analysis of ``arbitrary'' Dyson--Schwinger equations, whose presentation in full generality is beyond the scope of this book.

\medskip

Accepting momentarily Theorem~\ref{Theorem_correlators_expansion_theta1}, we now prove the two corollaries.

\begin{proof}[Proof of Corollary~\ref{Theorem_correlators_expansion_theta1}]
We take a closer look at the expressions \eqref{eq:8585} for $\boldsymbol{W}_1^{[1]}(z)$ and $\boldsymbol{W}_1^{[2]}(z)$.  Since we are working with Assumption~\ref{Assumptions_extra} and  $\theta_{g,g} = 1$ for every $g \in [H]$, Remark~\ref{remtheta11} tells us that $\phi_h^{\pm,[1]}(z) = 0$.  Then
\[
q_h^{+,[1]}(z) = \frac{q_h^+(z)}{2} \cdot \bigg( \sum_{g \neq h} \theta_{h,g}(\theta_{h,g} - 1) \partial_z \mathcal{G}_{\mu_g}(z)\bigg),
\]
where the $g = h$ term could be excluded since $\theta_{h,h} = 1$. This is a holomorphic function of $z$ in a neighborhood of $\amsmathbb{A}_h^{\mathfrak{m}}$, and by Lemma~\ref{Lemm84} it remains so after division by $s_h^{\circ}(z)$. Then, by the last property stated in Theorem~\ref{Theorem_Master_equation} we get $0$ after applying the operator $\boldsymbol{\Upsilon}$, resulting in $\boldsymbol{W}_1^{[1]}(z) = 0$. Then, in the expression for $\boldsymbol{W}_1^{[2]}(z)$, only the contribution of $q_g^{+,[2]}$ and $\boldsymbol{\mathcal{F}}(*,*)$ remain. Yet, in the function $q_h^{+,[2]}(z)$ given by \eqref{qhsmallexp} for the reason already invoked, the terms which are products of $q_h^+(z)$ with a holomorphic function of $z$ will not contribute to $\boldsymbol{W}_1^{[2]}(z)$. And, since $\theta_{h,h} = 1$ the scalar products involving $\bth_h - 1$ are sums over $g \neq h$ involving $\mathcal{G}_{\mu_g}(z)$ and thus give holomorphic functions of $z$ near $\amsmathbb{A}_h^{\mathfrak{m}}$. So, only the first term --- involving $\partial_z^2 \mathcal{G}_{\boldsymbol{\mu}}(z)$ --- and the last line of \eqref{qhsmallexp} effectively contribute to $\boldsymbol{W}_1^{[2]}(z)$. In this last line, $\phi_h^{-,[2]}(z)$ is in fact zero if we compare to Definition~\ref{Def_finer_Phi} and use the fact that under Assumption~\ref{Assumptions_extra} and $\theta_{h,h} = 1$, either $\iota_h^+ = 0$ or $\iota_h^+ = 2$ with $\rho_{h,1} = 1 = \theta_{h,h} = \rho_{h,2}$. For the same reason, Definition~\ref{Def_finer_Phi} gives $\phi_h^{+,[2]}(z) = - \frac{1}{12}\partial_z^3 U_h(z) \cdot \phi_h^+(z)$. This yields the formulae \eqref{W12W12}-\eqref{W12W120}.
\end{proof}

\begin{proof}[Proof of Corollary~\ref{Corollary_CLT} conditionally to Theorem~\ref{Theorem_correlators_expansion}] First, for each $l \in [L]$ we take a point $z_l \in \gamma_{h_l}$ and set $f_l(z)=\frac{1}{z_l-z}$. Then the vector \eqref{eq_x142} coincides by definition in \eqref{decomposstil} with $\big(\Delta G_{h_l}(z_l)\big)_{l=1}^L$. The asymptotics of its mean is computed in \eqref{eq_x30} as $\big(W_{1;h_l}^{[1]}(z_l)\big)_{l=1}^L$, the asymptotics of covariance is computed in \eqref{eq_x32} as $\big(\mathcal{F}_{h_{l_1},h_{l_2}}(z_{l_1},z_{l_2})\big)_{l_1,l_2 = 1}^{L}$. Besides, \eqref{eq_x33} implies that all cumulants of order greater than $2$ vanish at $\N \rightarrow \infty$ and this characterizes Gaussian random vectors --- see, \textit{e.g.}, \cite[Chapter 3]{PT}. The more general case is obtained by representing for any $l \in [L]$
\[
\sum_{i = 1}^{N_{h_{l}}} f_{l}\bigg(\frac{\ell^{h_l}_{i}}{\N}\bigg) - \N
\int_{\hat{a}_{h_l}'}^{\hat{b}_{h_l}'} f_{l}(x)\mu_{h_l}(x)\dd x= \oint_{\gamma_{h_l}} \frac{\dd z}{2\ii\pi}\,f_{l}(z)\,\Delta G_{h_l}(z).
\]
Hence, using Fubini theorem, all joint moments of \eqref{eq_x142} are integrals of the joint moments of $\Delta G_{h}(z)$. For instance, the mean of the $l$-th component of \eqref{eq_x142} is represented as
\begin{equation}
\label{eq_x213}
 \oint_{\gamma_{h_l}} \frac{\dd z}{2\ii\pi}\,f_{l}(z) \,\E\big[\Delta G_{h_l}(z)\big]
\end{equation}
and the covariance of the $l_1$-th and $l_2$-th components is
\begin{equation}
\label{eq_x214}
\oint\limits_{\gamma_{h_{l_1}}} \oint\limits_{\gamma_{h_{l_2}}}  \frac{\dd z_{1}\,\dd z_{2}}{(2\ii\pi)^2} \,f_{l_1}(z_{1}) f_{l_2}(z_{2}) \E \Big[\big(\Delta G_{h_{l_1}}(z_{1})-\E[\Delta G_{h_{l_1}}(z_{1})]\big) \cdot \big(\Delta G_{h_{l_2}}(z_{2})-\E[\Delta G_{h_{l_2}}(z_{2})]\big)\Big].
\end{equation}
Sending $\N\rightarrow\infty$ in \eqref{eq_x213}, \eqref{eq_x214} using the first part of the proof leads to the expressions of \eqref{eq_x144} for the asymptotic mean and variance. For higher moments, we use similar integral expressions and the fact that the integrals of jointly Gaussian random variables are also Gaussian.
\end{proof}

\subsection{Results without Assumption~\ref{Assumptions_extra}}

We now produce analogues of Theorem~\ref{Theorem_correlators_expansion} and Corollary~\ref{Corollary_CLT} which do not rely on Assumption~\ref{Assumptions_extra} but still in the case of fixed segment filling fractions. The idea is to use the conditioning procedure of Chapter~\ref{Chapter_conditioning} to reduce it to Theorem~\ref{Theorem_correlators_expansion}.

\begin{theorem}
\label{Theorem_correlators_expansion_relaxed} Consider a discrete ensemble satisfying Assumptions~\ref{Assumptions_Theta}, \ref{Assumptions_basic}, \ref{Assumptions_offcrit}, and \ref{Assumptions_analyticity} and such that the equations \eqref{eq_equations_eqs} deterministically fix the segment filling fractions and each segment contains a single band.

For any $\varepsilon>0$ and compact sets $\amsmathbb{K}_h\subset \amsmathbb C\setminus \amsmathbb{A}_h^{\mathfrak{m}}$ indexed by $h\in[H]$, there exists a sequence of positive constants $(C_n)_{n \geq 1}$ depending only on $\eps$, the compacts and the constants in the assumptions, and there exists for each $h \in [H]$ two holomorphic functions $\widetilde{W}_{1;h}^{[1]}(z)$ and $\widetilde{W}_{2;h}^{[2]}(z)$ of $z \in \amsmathbb{C} \setminus (\alpha_h,\beta_h)$ which are $O\big(\frac{1}{z^2}\big)$ as $z \rightarrow +\infty$, such that for any $n \in \amsmathbb{Z}_{> 0}$, $h_1,\ldots,h_n \in [H]$ and any $z_i\in \amsmathbb{K}_{h_i}$ for $i \in [n]$, we have
\begin{align}
\label{eq_x80_2} \big|W_{1;h_1}(z_1) - \widetilde{W}_{1;h_1}^{[1]}(z_1) - \N^{-1} \widetilde{W}_{1;h_1}^{[2]}(z_1)\big| & \leq C_1\,\N^{-\frac{3}{2} + \varepsilon}, \\
\label{eq_x81_2}\big|W_{2;h_1,h_2}(z_1,z_2) - \mathcal{F}_{h_1,h_2}(z_1,z_2)\big| & \leq
C_2\,\N^{-\frac{1}{2} + \varepsilon}, \\
\label{eq_x82_2}\forall n \geq 3,\qquad \big|W_{n;h_1,\ldots,h_n}(z_1,\ldots,z_n)\big| & \leq
C_n\,\N^{-\frac{1}{2} + \varepsilon},
\end{align}
where $\boldsymbol{\mathcal{F}}(z_1,z_2)$ is as in Section~\ref{TheopUpsec2}. The functions $\widetilde{W}_{1;h}^{[p]}(z)$  for $p \in \{1,2\}$ depend smoothly on parameters, in the following way.

Let $p \in \amsmathbb{Z}_{\geq 0}$. Suppose that we allow segment filling fractions and endpoints $(\hat{\boldsymbol{a}},\hat{\boldsymbol{b}},\hat{\boldsymbol{n}})$ to vary in some open of $\amsmathbb{R}^3$. Suppose that the weights $(w_h(x))_{h = 1}^{H}$, the parameters $\rho_{h,j}$, the regular part of the potentials $(U_h(x))_{h = 1}^{H}$ are $p$ times continuously differentiable with respect to $(\hat{\boldsymbol{a}},\hat{\boldsymbol{b}},\hat{\boldsymbol{n}})$ while the integers $(\iota_h^{\pm})_{h = 1}^{H}$ remain constant, and that Assumptions~\ref{Assumptions_Theta}, \ref{Assumptions_basic}, \ref{Assumptions_offcrit}, \ref{Assumptions_analyticity} hold for uniform constants and the assumption of one band per segment remain satisfied. Then $\widetilde{W}_{1;h}^{[p]}(z)$ for $p \in \{1,2\}$ are $p$ times continuously differentiable with respect to $(\hat{\boldsymbol{a}},\hat{\boldsymbol{b}},\hat{\boldsymbol{n}})$ for $z$ away from the segments, and all the partial derivatives up to order $p$ are uniformly bounded for $z$ away from the segments, with constants determined by the constants in the assumptions.
\end{theorem}

\begin{theorem}
\label{Theorem_correlators_expansion_relaxed_theta1}
Consider the setting of Theorem~\ref{Theorem_correlators_expansion_relaxed} and assume $\theta_{g,g} = 1$ for every $g \in [H]$. Then, the same conclusion holds with $\boldsymbol{W}_1^{[1]}(z) = 0$.
\end{theorem}

At the end of this section we give the proof of Theorems~\ref{Theorem_correlators_expansion_relaxed} and \ref{Theorem_correlators_expansion_relaxed_theta1}, again conditionally on Theorem~\ref{Theorem_correlators_expansion} that is proved later. While filling fractions and endpoints are not arbitrary real numbers in discrete ensembles due to the integrality constraints described in Section~\ref{Section_general_model}, in the first place due to $\N \hat n_h\in \amsmathbb{Z}_{> 0}$, yet Theorem~\ref{Theorem_correlators_expansion_relaxed} claims that $\widetilde{W}_{1;h}^{[p]}(z)$ for $p \in \{1,2\}$ exist as smooth functions of continuously varying parameters $(\hat{\boldsymbol{a}},\hat{\boldsymbol{b}},\hat{\boldsymbol{n}})$, provided the weights $w_g$, the parameters $\rho_{g,j}$ and the regular part of the potential $U_g$ do for every $g \in [H]$. In the proof, we give a formula  \eqref{expsfguzbguwg} for $\widetilde{W}_{1;h}^{[p]}$ in terms of $W_{1;h}^{[p],Y}$  coming from Theorem~\ref{Theorem_correlators_expansion} but for a conditioned and localized ensemble $Y$ delivered from the original ensemble by Theorem~\ref{proposition_fluct_conditioning}. In case $\boldsymbol{\Theta}$ has diagonal $1,\ldots,1$, this formula yields $\widetilde{W}_{1;h}^{[1]}(z) = O(\frac{1}{\N})$. For this reason, in Theorem~\ref{Theorem_correlators_expansion_relaxed_theta1} we redefine $\widetilde{W}_{1;h}^{[1]}(z) = 0$ and rather absorb the aforementioned $O(\frac{1}{\N})$ as an additional $O(1)$ term in $\widetilde{W}_{1;h}^{[2]}$: this is formula \eqref{expsfguzbguwg22}.

\begin{corollary}[Central limit theorem]
\label{Corollary_CLT_relaxed}
Consider a discrete ensemble satisfying Assumptions~\ref{Assumptions_Theta}, \ref{Assumptions_basic}, \ref{Assumptions_offcrit}, \ref{Assumptions_analyticity}. Assume additionally that the equations \eqref{eq_equations_eqs} deterministically fix the segment filling fractions and that the equilibrium measure contains exactly one band per segment. Let $L \in \amsmathbb{Z}_{> 0}$ and a $L$-tuple of integers $\boldsymbol{h} \in [H]^L$ both independent of $\N$, let a (possibly $\N$-dependent) $L$-tuple of functions $\boldsymbol{f}(z)$ such that $f_l(z)$ is a holomorphic function of $z \in \amsmathbb{M}_{h_l}$ for any $l \in [L]$. Assume there exists a constant $C > 0$ such that $\max_{l} \sup_{z} |f_l(z)| \leq C$.

We define the random $L$-dimensional centered Gaussian vector $\textnormal{\textsf{\textbf{Gau\ss{}}}}_{\textbf{0}}[\boldsymbol{f}]$ with covariance
\begin{equation}
\label{eq_x162}
\textnormal{\textsf{Cov}}_{l_1,l_2}[\boldsymbol{f}] = \oint_{\gamma_{h_{l_1}}} \oint_{\gamma_{h_{l_2}}} \frac{\dd z_1\,\dd z_2}{(2\ii\pi)^2}
 f_{l_1}(z_1) f_{l_2}(z_2)
\,\mathcal{F}_{h_{l_1},h_{l_2}}(z_{l_1},z_{l_2}).
\end{equation}
If $\boldsymbol{\ell}$ is $\P$-distributed, then as $\N \rightarrow \infty$
\begin{equation}
\label{eq_x163}
\left(\sum_{i = 1}^{N_{h_l}} f_{l}\bigg(\frac{\ell^{h_l}_{i}}{\N}\bigg)-\E\left[\sum_{i = 1}^{N_{h_l}} f_{l}\bigg(\frac{\ell^{h_l}_{i}}{\N}\bigg)\right]\right)_{l=1}^L\,\,\mathop{\sim}^{\textnormal{mom}} \,\,\textnormal{\textsf{\textbf{Gau\ss{}}}}_{\textbf{0}}[\boldsymbol{f}]
\end{equation}
and the convergence in the sense of Definition~\ref{Definition_convergence_moments} is uniform for fixed choice of constants in the assumptions. Besides, for any $m \in [M]$, we have the asymptotic expansion
\begin{equation}
\label{eq_x164}
\E\left[\sum_{i = 1}^{N_{h_l}} f_{l}\bigg(\frac{\ell^{h_l}_{i}}{\N}\bigg)\right] = \N
\int_{\hat{a}_{h_l}'}^{\hat{b}_{h_l}'} f_{l}(x)\mu_{h_l}(x)\dd x + \mathbbm{Rest}[f_l] + O(\N^{-\frac{3}{2}+\eps}),
\end{equation}
where for fixed constants in the assumptions, $\mathbbm{Rest}[f_l]$ is uniformly bounded and the remainder $O(\N^{-\frac{3}{2} + \eps})$ is uniform as $\N \rightarrow \infty$. The term $\mathbbm{Rest}[f_l]$ depends smoothly on parameters in the following way.

Let $p \in \amsmathbb{Z}_{\geq 0}$. Suppose that we allow segment filling fractions and endpoints $(\hat{\boldsymbol{a}},\hat{\boldsymbol{b}},\hat{\boldsymbol{n}})$ to vary in some open of $\amsmathbb{R}^3$. Suppose that the weights $(w_h(x))_{h = 1}^{H}$, the parameters $(\rho_{h,j}^{\pm})_{h,j}$, the regular part of the potentials $(U_h(x))_{h = 1}^{H}$ are $p$ times continuously differentiable with respect to $(\hat{\boldsymbol{a}},\hat{\boldsymbol{b}},\hat{\boldsymbol{n}})$ while the integers $(\iota_h^{\pm})_{h = 1}^{H}$ remain constant, and that Assumptions~\ref{Assumptions_Theta}, \ref{Assumptions_basic}, \ref{Assumptions_offcrit}, \ref{Assumptions_analyticity} hold for uniform constants and the assumption of one band per segment remain satisfied. Then $\mathbbm{Rest}[f_l]$ is $p$ times continuously differentiable with respect to $(\hat{\boldsymbol{a}},\hat{\boldsymbol{b}},\hat{\boldsymbol{n}})$ for $z$ away from the segments, and all the partial derivatives up to order $p$ are uniformly bounded for $z$ away from the segments, with constants determined by the constants in the assumptions.

If furthermore $\theta_{g,g} = 1$ for every $g \in [H]$, the term $\mathbbm{Rest}[f_l]$ and all its partial derivatives are actually $O(\frac{1}{\N})$. In particular
\begin{equation}
\label{neatCLT}
\left(\sum_{i = 1}^{N_{h_l}} f_{l}\bigg(\frac{\ell^{h_l}_{i}}{\N}\bigg) - \N \int_{\hat{a}_{h_l}'}^{\hat{b}_{h_l}'} f_l(x)\mu_{h_l}(x)\dd x \right)_{l = 1}^{L}  \,\,\mathop{\sim}^{\textnormal{mom}} \,\,\textnormal{\textsf{\textbf{Gau\ss{}}}}_{\textbf{0}}[\boldsymbol{f}].
\end{equation}
\end{corollary}

The reduction of Corollary~\ref{Corollary_CLT_relaxed} to Theorem~\ref{Theorem_correlators_expansion_relaxed} is exactly the same as the reduction of Corollary~\ref{Corollary_CLT} to Theorem~\ref{Theorem_correlators_expansion} and we omit the proof. The covariance in Corollary~\ref{Corollary_CLT_relaxed} is exactly the same as in Corollary~\ref{Corollary_CLT}. However, a relatively simple formula for the mean is not available without Assumption~\ref{Assumptions_extra} and that is the reason why we recentered with respect to the mean in \eqref{eq_x163}, but not in \eqref{eq_x142}. We have collected
\[
\mathbbm{Rest}[f_l] = \oint_{\gamma_{h_l}} \frac{\dd z}{2\ii\pi} \widetilde{W}_{1;h_l}^{[1]}(z) f_l(z) + \frac{1}{\N} \oint_{\gamma_{h_l}} \frac{\dd z}{2\ii\pi} \widetilde{W}_{1;h_l}^{[2]}(z) f_l(z)
\]
in terms of the functions of Theorem~\ref{Theorem_correlators_expansion_relaxed} or \ref{Theorem_correlators_expansion_relaxed_theta1}.  In case $\theta_{h,h} = 1$, the order $1$ term is absent, resulting in the simpler statement \eqref{neatCLT} of the central limit theorem.

\subsection{Proof of Theorems~\ref{Theorem_correlators_expansion_relaxed} and \ref{Theorem_correlators_expansion_relaxed_theta1} conditional to Theorem~\ref{Theorem_correlators_expansion}}

We rely on Theorem~\ref{proposition_fluct_conditioning} (with $K=H$ in the present situation) and use $\eps_1,\eps_2>0$ from there, which can be chosen to be independent of $(\hat{\boldsymbol{a}}, \hat{\boldsymbol{b}}, \hat{\boldsymbol{n}})$ as long as the latter vary in a small fixed open set. We further choose the $2H$-tuple $(\mathfrak{a}_h,\mathfrak{b}_h)_{h = 1}^{H}$ as in Theorem~\ref{proposition_fluct_conditioning}, in such a way that the differences $\mathfrak a_h-\hat a_h$ and $\hat b_h - \mathfrak b_h$ do not depend on $(\hat{\boldsymbol{a}}, \hat{\boldsymbol{b}}, \hat{\boldsymbol{n}})$. If as in Theorem~\ref{proposition_fluct_conditioning} we use an exponent $X$ to refer to the original ensemble and $Y$ to refer to the ensemble localized to $\mathfrak{A} = \bigcup_{h = 1}^{H} [\mathfrak{a}_h,\mathfrak{b}_h]$, then we have:
\begin{equation}
\label{eq_x147}
\begin{split}
 G_h^X(z) & = G_h^Y(z) + \mathbbm{1}_{\amsmathbb{S}_h}(\hat a'_h)\left(\frac{1}{z - \hat{a}_h} + \frac{1}{x - \big(\hat{a}_h + \frac{\theta_{h,h}}{\N}\big)} + \cdots + \frac{1}{z - \big(\mathfrak{a}_h - \frac{\theta_{h,h}}{\N}\big)}\right) \\
& \quad + \mathbbm{1}_{\amsmathbb{S}_h}(\hat b'_h)\left( \frac{1}{z - \big(\mathfrak{b}_h + \frac{\theta_{h,h}}{\N}\big)} + \frac{1}{z - \big(\mathfrak{b}_h + 2\frac{\theta_{h,h}}{\N}\big)} + \cdots + \frac{1}{z - \hat{b}_h}\right) + \textnormal{error},
\end{split}
\end{equation}
where the error term is non-zero (but bounded) with exponentially small probability. Note that cumulants of order at least two are unchanged if we add deterministic parts to random variables. Both sums in the right-hand side of \eqref{eq_x147} are deterministic, hence \eqref{eq_x81_2} and \eqref{eq_x82_2} follow from \eqref{eq_x32} and \eqref{eq_x34} established in Theorem~\ref{Theorem_correlators_expansion} for the localized ensemble $Y$ which satisfies the extra Assumption~\ref{Assumptions_extra}. These deterministic sums can be compared to an integral by the Euler--Maclaurin formula, which we will use in the following form.

\begin{lemma}
\label{lem:EulerMaclaurin}
Choose $\theta > 0$, real parameters $\hat{a},\mathfrak{a},\hat{b},\mathfrak{b}$, and a twice continuously differentiable function $f$ in the range required by the formulae below. We introduce shifted parameters
\[
\hat{a}' = \hat{a} + \frac{\frac{1}{2} - \theta}{\N},\qquad \mathfrak{a}' = \mathfrak{a} + \frac{\frac{1}{2} - \theta}{\N},\qquad \hat{b}' = \hat{b} + \frac{\theta - \frac{1}{2}}{\N},\qquad \mathfrak{b}' = \mathfrak{b} + \frac{\theta - \frac{1}{2}}{\N}.
\]
If $\N(\mathfrak{a} - \hat{a})$ is $\theta$ times a positive integer $M$, we have
\begin{equation}
\label{EulerMaclaurinleft}
\begin{split}
& \quad \sum_{i = 0}^{M - 1} f\bigg(\hat{a} + \frac{i\theta}{\N}\bigg) - \frac{\N}{\theta} \int_{\hat{a}'}^{\mathfrak{a}'} f(x)\dd x \\
& = \frac{\theta -1}{2\theta}\big(f(\mathfrak{a}') - f(\hat{a}')\big) + \frac{2\theta^2 - 6\theta + 3}{24\theta \N}\big(f'(\mathfrak{a}') - f'(\hat{a}')\big) + O\bigg(\frac{1}{\N^2}\bigg).
\end{split}
\end{equation}
If $\N(\hat{b} - \mathfrak{b})$ is $\theta$ times a positive integer $M$, we have
\begin{equation}
\label{EulerMaclaurinright}
\begin{split}
& \quad \sum_{i = 1}^{M} f\bigg(\mathfrak{b} + \frac{i\theta}{\N}\bigg) - \frac{\N}{\theta} \int_{\mathfrak{b}'}^{\hat{b}'} f(x)\dd x \\
& = \frac{\theta -1}{2\theta}\big(f(\mathfrak{b}') - f(\hat{b}')\big) + \frac{2\theta^2 - 6\theta + 3}{24\theta \N}\big(f'(\hat{b}') - f'(\mathfrak{b}')\big) + O\bigg(\frac{1}{\N^2}\bigg).
\end{split}
\end{equation}
If $\N(\mathfrak{a} - \mathfrak{b})$ is $\theta$ times a positive integer $M$, we have
\begin{equation}
\label{EulerMaclaurinmiddle}
\begin{split}
& \quad \sum_{i = 1}^{M - 1} f\bigg(\mathfrak{b} + \frac{i\theta}{\N}\bigg) - \frac{\N}{\theta} \int_{\mathfrak{b}'}^{\mathfrak{a}'} f(x)\dd x \\
& = \frac{\theta -1}{2\theta}\big(f(\mathfrak{b}') + f(\mathfrak{a}')\big) + \frac{2\theta^2 - 6\theta + 3}{24\theta\N}\big(f'(\mathfrak{a}') - f'(\mathfrak{b}')\big) + O\bigg(\frac{1}{\N^2}\bigg).
\end{split}
\end{equation}
\end{lemma}
\begin{proof}
We focus on the first formula. The classical Euler--Maclaurin formula says that
\[
\sum_{i = 0}^{M - 1} f\bigg(\hat{a} + \frac{i\theta}{\N}\bigg) = \frac{\N}{\theta} \int_{\hat{a}}^{\mathfrak{a}} f(x)\dd x + \frac{f(\hat{a}) - f(\mathfrak{a})}{2} - \frac{\theta}{12\N}\big(f'(\hat{a}) - f'(\mathfrak{a})\big) + O\bigg(\frac{1}{\N^2}\bigg).
\]
We rewrite it in terms of left-shifted parameters
\begin{equation}
\label{EulerMaclaurinleftp}
\begin{split}
\sum_{i = 0}^{M - 1} f\bigg(\hat{a} + \frac{i\theta}{\N}\bigg) & = \frac{\N}{\theta} \int_{\hat{a}'}^{\mathfrak{a}'} f\bigg(x + \frac{\theta - \frac{1}{2}}{\N}\bigg)\dd x +  \frac{1}{2}\Bigg[f\bigg(\hat{a}' + \frac{\theta -\frac{1}{2}}{\N}\bigg) - f\bigg(\mathfrak{a}' + \frac{\theta - \frac{1}{2}}{\N}\bigg)\Bigg] \\
& \quad - \frac{\theta}{12\N} \Bigg[f'\bigg(\hat{a}' + \frac{\theta - \frac{1}{2}}{\N}\bigg) - f'\bigg(\mathfrak{a}' + \frac{\theta - \frac{1}{2}}{\N}\bigg)\Bigg] + O\bigg(\frac{1}{\N^2}\bigg).
\end{split}
\end{equation}
Let $F$ be a primitive function of $f$. We obtain after a Taylor expansion
\begin{equation*}
\begin{split}
& \quad \frac{\N}{\theta} \int_{\hat{a}'}^{\mathfrak{a}'} f\bigg(x + \frac{\theta - \frac{1}{2}}{\N}\bigg) \dd x = \frac{\N}{\theta}\Bigg[ F\bigg(\mathfrak{a}' + \frac{\theta - \frac{1}{2}}{\N}\bigg) - F\bigg(\hat{a}' + \frac{\theta - \frac{1}{2}}{\N}\bigg)\Bigg] \\
& = \frac{\N}{\theta}\big(F(\mathfrak{a}') -  F(\hat{a}')\big) + \frac{2\theta - 1}{2\theta} \big( f(\mathfrak{a}') - f(\hat{a}')\big) + \frac{(2\theta - 1)^2}{8\theta \N} \big(f'(\mathfrak{a}') - f'(\hat{a}')\big) + O\bigg(\frac{1}{\N^2}\bigg).
\end{split}
\end{equation*}
Likewise
\[
\frac{1}{2}\Bigg[f\bigg(\hat{a}' + \frac{\theta -\frac{1}{2}}{\N}\bigg) - f\bigg(\mathfrak{a}' + \frac{\theta - \frac{1}{2}}{\N}\bigg)\Bigg] = -\frac{1}{2}\big(f(\mathfrak{a}') - f(\hat{a}') \big) + \frac{1 - 2\theta}{2\N}\big(f'(\mathfrak{a}') - f'(\hat{a}')\big) + O\bigg(\frac{1}{\N^2}\bigg), \\
\]
and
\[
- \frac{\theta}{12\N} \Bigg[f'\bigg(\hat{a}' + \frac{\theta - \frac{1}{2}}{\N}\bigg) - f'\bigg(\mathfrak{a}' + \frac{\theta - \frac{1}{2}}{\N}\bigg)\Bigg] = \frac{\theta}{12\N}\big(f'(\mathfrak{a}') - f'(\hat{a}')\big) + O\bigg(\frac{1}{\N^2}\bigg).
\]
Collecting all the terms in \eqref{EulerMaclaurinleftp} we arrive to the claimed \eqref{EulerMaclaurinleft}. The two other formulae \eqref{EulerMaclaurinright}-\eqref{EulerMaclaurinmiddle} are obtained in a similar way.
\end{proof}

We can apply \eqref{EulerMaclaurinleft}-\eqref{EulerMaclaurinright} to $f(x) = \frac{1}{z - x}$, step $\theta = \theta_{h,h}$ and endpoints $\mathfrak{a}_h,\mathfrak{a}_h,\mathfrak{b}_h,\mathfrak{b}$ to rewrite \eqref{eq_x147}. In order to establish \eqref{eq_x80_2}, we also write with shifted parameters $\mathfrak a'_h$ and $\mathfrak b'_h$ from Theorem~\ref{proposition_fluct_conditioning}
\begin{equation}
\label{eq_x160}
\begin{split}
 \mathcal G_{\mu_h^X}(z) & =\mathcal G_{\mu_h^Y}(z) + \mathbbm{1}_{\amsmathbb{S}_h}(\hat{a}'_h) \int_{\hat a'_h}^{\mathfrak a'_h} \frac{\dd x}{\theta_{h,h}(z-x)} + \mathbbm{1}_{\amsmathbb{S}_h}(\hat b'_h) \int_{\mathfrak b'_h}^{\hat b'_h} \frac{\dd x}{\theta_{h,h}(z-x)} \\
 & \quad + \int_{\mathfrak a'_h}^{\mathfrak b'_h} \frac{\mu_h^X(x)-\mu_h^Y(x)}{z-x}\dd x.
 \end{split}
\end{equation}
The existence of the last term was explained in Section~\ref{Sec_Mismatch} and it was shown in Proposition~\ref{mueqYYtilde} to be $O(\frac{1}{\N})$ as $\N \rightarrow \infty$. We subtract \eqref{eq_x160} multiplied by $\N$ from the expectation value of \eqref{eq_x147} to get an expression of $\E\big[\Delta G_h^X(z)\big]$ in terms of $\E\big[\Delta G_h^Y(z)\big]$, and use \eqref{eq_x30} for the asymptotic expansion of the latter. We also use the Euler--Maclaurin approximations \eqref{EulerMaclaurinleft}-\eqref{EulerMaclaurinright} with $\theta = \theta_{h,h}$ and $f(x) = \frac{1}{z - x}$ to handle the deterministic sums in \eqref{eq_x147}.  As a result, we get \eqref{eq_x80_2} with
\begin{equation}
\label{expsfguzbguwg}
\begin{split}
\widetilde{W}_{1;h}^{[1]}(z) &= W_{1;h}^{[1],Y}(z) + \N \int_{\mathfrak{a}'_h}^{\mathfrak{b}'_h} \frac{\mu_h(x) - \mu_h^{Y}(x)}{z - x} \dd x \\
& \quad + \frac{\theta_{h,h} - 1}{2\theta_{h,h}}  \left(\frac{\mathbbm{1}_{\amsmathbb{S}_h}(\hat{a}'_h)}{z - \mathfrak{a}'_h} - \frac{\mathbbm{1}_{\amsmathbb{S}_h}(\hat{a}'_h)}{z - \hat{a}_h'} + \frac{\mathbbm{1}_{\amsmathbb{S}_h}(\hat{b}'_h)}{z - \mathfrak{b}'_h} - \frac{ \mathbbm{1}_{\amsmathbb{S}_h}(\hat{b}'_h)}{z - \hat{b}_h'}\right), \\
\widetilde{W}_{1;h}^{[2]}(z) & = W_{1;h}^{[2],Y}(z)  \\
& \quad + \frac{2\theta_{h,h}^2 - 6\theta_{h,h} + 3}{24\theta_{h,h}}\left(\frac{\mathbbm{1}_{\amsmathbb{S}_h}(\hat{a}'_h)}{(z - \mathfrak{a}'_h)^2} - \frac{\mathbbm{1}_{\amsmathbb{S}_h}(\hat{a}'_h)}{(z - \hat{a}_h')^2} - \frac{\mathbbm{1}_{\amsmathbb{S}_h}(\hat{b}'_h)}{(z - \mathfrak{b}'_h)^2} + \frac{\mathbbm{1}_{\amsmathbb{S}_h}(\hat{b}'_h)}{(z - \hat{b}_h')^2}\right).
\end{split}
\end{equation}
It remains to show that  for $p \in \{1,2\}$, the functions $\widetilde{W}_{1;h}^{[p]}(z)$ and their partial derivatives with respect to filling fractions and endpoints are uniformly bounded as $\N\rightarrow\infty$. We recall that $z$ is kept away from the segments defining the domain.

The two terms $W_{1;h}^{[1],Y}(z)$ and $W_{1;h}^{[2],Y}(z)$ come from Theorem~\ref{Theorem_correlators_expansion}, where they are expressed through the solution operator $\boldsymbol{\Op}$. This operator depends on the endpoints of the bands $(\boldsymbol{\alpha}^Y,\boldsymbol{\beta}^Y)$ for the ensemble $Y$ in a smooth way by Theorem~\ref{Theorem_Master_equation}.  Because all the data for the ensemble $Y$ depend smoothly on the data for the ensemble $X$, so do $(\boldsymbol{\alpha}^Y,\boldsymbol{\beta}^Y)$ by Theorem~\ref{Theorem_differentiability_full}. The functions on which we apply the operator $\boldsymbol{\Upsilon}$ involve $\phi^\pm_h(z)$, $\phi^{\pm,[1]}_h(z)$, $\phi_h^{\pm,[2]}(z)$ and the Stieltjes transforms of the equilibrium measure of $Y$. Using Theorem~\ref{Theorem_differentiability_full}, we conclude again that it depends smoothly on $(\hat{\boldsymbol{a}},\hat{\boldsymbol{b}},\hat{\boldsymbol{n}})$. The same is true for the terms present if there are saturations: by design in Theorem~\ref{proposition_fluct_conditioning} the differences $(\mathfrak a_h-\hat a_h)$ and $(\hat b_h - \mathfrak b_h)$ were chosen constant over the small open set in which $(\hat{\boldsymbol{a}},\hat{\boldsymbol{b}},\hat{\boldsymbol{n}})$ varies, and $\mathfrak{a}_h',\mathfrak{b}_h'$ are related to $\mathfrak{a}_h,\mathfrak{b}_h$ by a constant shift.

Finally, we have a term in \eqref{expsfguzbguwg} which is precisely $\N$ times the difference between Stieltjes transforms of the equilibrium measures of the ensemble $Y$ and of the ensemble $X$ restricted to $[\mathfrak{a}_h',\mathfrak{b}_h']$. This difference comes from the subleading corrections arising in the comparison between the potentials for $X$ and $Y$. Lemma~\ref{goodsimz} tells us that this difference is $O(\frac{1}{\N})$. Using Binet's integral formula  \eqref{eq_Binet_formula}  instead of asymptotic expansions in the proof of Lemma~\ref{goodsimz} allows checking that the errors remain of the same order of magnitude after taking partial derivatives. Controlling the resulting variation of the Stieltjes transform of equilibrium measures with help of Proposition~\ref{Theorem_differentiability_full}, we deduce its partial derivatives exist and remain $O(\frac{1}{\N})$ as well. This concludes the proof of Theorem~\ref{Theorem_correlators_expansion_relaxed}.

Finally, if $\theta_{g,g} = 1$ for every $g \in [H]$, many simplifications occur. First, Corollary~\ref{Theorem_correlators_expansion_theta1} where $W_{1;h}^{[1],Y}(z)$ comes from tells us that it vanishes for any $h \in [H]$. Second, the saturation contributions in \eqref{expsfguzbguwg} vanish. Third, Proposition~\ref{mueqYYtilde} tells us that the Stieltjes transform of the restriction of $\mu - \mu^Y$ to $[\mathfrak{a}_h',\mathfrak{b}_h']$ is of order $O(\frac{1}{\N^2})$ for any $h \in [H]$. Then, \eqref{expsfguzbguwg} leaves us with $\widetilde{W}_{1;h}^{[1]}(z) = O(\frac{1}{\N})$. Thus, instead of \eqref{expsfguzbguwg} we can take as new definition $\widetilde{W}_{1;h}^{[1]}(z) = 0$ and
\begin{equation}
\label{expsfguzbguwg22}
\begin{split}
\widetilde{W}_{1;h}^{[2]}(z)  & = W_{1;h}^{[2],Y} + \N \int_{\mathfrak{a}_h'}^{\mathfrak{b}_h'} \frac{\mu_h(x) - \mu_h^{Y}(x)}{z - x}\,\dd x  \\
& \quad - \frac{1}{24}\left(\frac{\mathbbm{1}_{\amsmathbb{S}_h}(\hat{a}'_h)}{(z - \mathfrak{a}'_h)^2} - \frac{\mathbbm{1}_{\amsmathbb{S}_h}(\hat{a}'_h)}{(z - \hat{a}_h')^2} - \frac{\mathbbm{1}_{\amsmathbb{S}_h}(\hat{b}'_h)}{(z - \mathfrak{b}'_h)^2} + \frac{\mathbbm{1}_{\amsmathbb{S}_h}(\hat{b}'_h)}{(z - \hat{b}_h')^2}\right),
\end{split}
\end{equation}
and these functions have the required properties for Theorem~\ref{Theorem_correlators_expansion_relaxed_theta1}.

\section{Asymptotic form of the Nekrasov equations}

\label{Section_asymptotic_Nekrasov}

The aim of this section is to rewrite the Nekrasov equations of Theorem~\ref{Theorem_Nekrasov} and Corollary~\ref{Corollary_higherNek} as polynomial equations on the joint cumulants $W_{n;h_1,\ldots,h_n}(z_1,\ldots,z_n)$ of the Stieltjes transforms, introduced in Definition~\ref{def_correlators}. We heavily rely on the $\bth_h$ and $\langle \cdot \rangle$ notations of Section~\ref{Theopr}. All the necessary auxiliary functions which are not specified in the statements have been introduced in Section~\ref{Subsection_Auxiliary_functions}.

\subsection{Rewriting the first-order Nekrasov equation}

\begin{proposition} \label{Proposition_Nek_1_asymptotic_form} Consider a discrete ensemble satisfying Assumptions~\ref{Assumptions_Theta}, \ref{Assumptions_basic}, \ref{Assumptions_offcrit}, \ref{Assumptions_analyticity} and the additional Assumption~\ref{Assumptions_extra}. For any $h\in[H]$ and $z\in \amsmathbb M_h\setminus \amsmathbb{A}^{\mathfrak{m}}_h$, we have
\begin{equation}
 \label{eq_First_Nekr_asymptotic}
 \begin{split}
& \quad Q_h^-(z)\,\big\langle \bth_h\cdot \boldsymbol{W}_1(z)\big\rangle + \frac{Q_h^+(z)}{2\N}\big\langle\bth_h\cdot \boldsymbol{W}_1(z)\big\rangle^2+
 \frac{Q_h^+(z)}{2\N}\big\langle \bth_h^{\otimes 2} \cdot \boldsymbol{W}_2(z,z)\big\rangle \\
 & \quad + \frac{Q_h^+(z)}{2\N}\big\langle \bth_h(\bth_h - 1)\cdot \partial_{z} \boldsymbol{W}_1(z)\big\rangle + q_h^{+,[1]}(z) + \frac{q_h^{+,[2]}(z)}{\N} \\
 & = \widetilde R_h(z)+ \mathbbm{Err}_h^{(1)}(z)+ \mathbbm{Err}_h^{(2)}(z),
\end{split}
\end{equation}
Here, $\widetilde R_h(z)$ is a holomorphic function of $z\in \amsmathbb M_h$. If $\hat a'_h \in \amsmathbb{S}_h$, then $\widetilde R_h(\hat a'_h)=0$; if $\hat b'_h\in \amsmathbb{S}_h$, then $\widetilde R_h(\hat b'_h)=0$.

Besides, for any $\eps > 0$ and compact $\amsmathbb{K}_h \subset \amsmathbb{M}_h \setminus \amsmathbb{A}_{h}^{\mathfrak{m}}$, there exists a constant $C > 0$ depending only on $\eps$, these compacts and the constants in the assumptions, such that for $\N$ large enough
\begin{equation}
\label{eq_Err1_bound}
\sup_{z \in \amsmathbb{K}_h} \big|\mathbbm{Err}_h^{(1)}(z)\big| \leq C \N^{-\frac{3}{2} + \eps}.
\end{equation}
The function $\mathbbm{Err}_h^{(2)}(z)$ is a finite sum --- with a number of terms which does not depend on anything --- of expressions of the form
\begin{equation}
\label{FERR2}
\frac{F(z)}{\N^{q}} \cdot \prod_{i=1}^k \E\Bigg[\prod_{j=1}^{l_i} \partial_z^{d_{i.j}} \Delta G_{g_{i,j}}(z) \Bigg] ,
\end{equation}
for some $q \in \amsmathbb{Z}_{\geq 2}$ and $k,l_1,\ldots,l_k \in \amsmathbb{Z}_{> 0}$ satisfying
  \begin{equation}
  \label{eq_Err2_bound2}
  q \leq \sum_{i=1}^k l_i \leq q+1,
  \end{equation}
 some function $F(z)$ which is holomorphic for $z \in \amsmathbb{K}_h$ and bounded from above in absolute value by $C$, and some $d_{i,j} \in \{0,1\}$ and $g_{i,j} \in [H]$.
\end{proposition}

The remainder term $\mathbbm{Err}_h^{(2)}(z;z_2)$ has an explicit expression, whose only virtue is to allow checking the inequality \eqref{eq_Err2_bound2}:
 \begin{equation}
\label{Err2qqq} \begin{split}
& \quad -\mathbbm{Err}_h^{(2)}(z) \\
& = \frac{Q_h^-(z)}{6\N^2}\big\langle \bth_h\cdot \boldsymbol{W}_1(z)\big\rangle^3 + \frac{Q_h^-(z)}{2\N^2}\big\langle \bth_h\cdot \boldsymbol{W}_1(z)\big\rangle\big\langle \bth_h^{\otimes 2}\cdot \boldsymbol{W}_{2}(z,z)\big\rangle
 + \frac{Q_h^-(z)}{6\N^2}\big\langle \bth_h^{\otimes 3}\cdot\boldsymbol{W}_3(z,z,z)\big\rangle \\
 & \quad + \frac{Q_h^+(z)}{24\N^3}\big\langle \bth_h\cdot \boldsymbol{W}_1(z)\big\rangle^4 + \frac{Q_h^+(z)}{8\N^3}\big\langle \bth_h^{\otimes 2} \cdot \boldsymbol{W}_2(z,z)\big\rangle^2 + \frac{Q_h^+(z)}{24\N^3}\big\langle\bth_h^{\otimes 4}\cdot \boldsymbol{W}_4(z,z,z,z)\big\rangle \\
& \quad + \frac{Q_h^+(z)}{6\N^3}\big\langle\bth_h\cdot \boldsymbol{W}_1(z)\big\rangle\big\langle \bth_h^{\otimes 3}\cdot \boldsymbol{W}_3(z,z,z)\big\rangle
 + \frac{Q_h^+(z)}{4\N^3}\big\langle \bth_h^{\otimes 2} \cdot \boldsymbol{W}_2(z,z)\big\rangle\big\langle \bth_h \cdot \boldsymbol{W}_1(z)\big\rangle^2
\\
& \quad + \frac{Q_h^-(z)}{2\N^2}\big\langle \bth_h\cdot \boldsymbol{W}_1(z)\big\rangle\big\langle \bth_h(\bth_h - 1)\cdot \partial_{z} \boldsymbol{W}_1(z)\big\rangle + \frac{Q_h^-(z)}{4 \N^2}\big\langle \bth_h \otimes \bth_h(\bth_h - 1) \cdot \partial_{z} \boldsymbol{W}_{2}(z,z)\big\rangle.
\end{split}
\end{equation}

Before proving Proposition~\ref{Proposition_Nek_1_asymptotic_form}, let us explain how we think about \eqref{eq_First_Nekr_asymptotic}. We treat \eqref{eq_First_Nekr_asymptotic} as a linear equation on $\boldsymbol{W}_1(z)$, $(\boldsymbol{W}_1(z))^{\otimes 2}$, and $\boldsymbol{W}_2(z,z)$ which all appear in the left-hand side, with the first term responsible for the leading contribution and the other terms giving subleading corrections. In the right-hand side, about the term $\widetilde R_h(z)$ we only know its analytic properties, but it turns out to be enough to guarantee that this term does not matter for the asymptotic expansion. The error terms $\mathbbm{Err}_h^{(1)}(z)$ and $\mathbbm{Err}_h^{(2)}(z)$ will be shown to be so small that they do not contribute to the asymptotic expansion at the precision that we are interested in. For $\mathbbm{Err}_h^{(1)}(z)$ this is clear from \eqref{eq_Err1_bound}; for $\mathbbm{Err}_h^{(2)}(z)$ this would eventually follow from $q\geq 2$ in \eqref{eq_Err2_bound2} after we prove that $\Delta \boldsymbol{G}(z)$ remains finite as $\N\rightarrow\infty$. The latter statement however needs an additional argument, because Corollary~\ref{Corollary_Stieltjes_apriory} is not enough to guarantee this.

\begin{proof}[Proof of Proposition~\ref{Proposition_Nek_1_asymptotic_form}] \phantom{s}

\medskip

\noindent \textsc{Step 1.} We start by perform the Taylor expansion of the right-hand side of \eqref{Nekrasov_eqn}.
We use an approximation based on the Taylor series expansion of $\log(1+x)$ and valid for $|x| \leq \frac{1}{2}$:
\[
1 + x = \exp\bigg(x - \frac{x^2}{2} + \frac{x^3}{3} + \mathbbm{c}(x)\bigg)\qquad\textnormal{with}\quad |\mathbbm{c}(x)| \leq 4x^4.
\]
Let us fix $\delta > 0$ small enough and independent of $\N$, and take $z$ at distance at least $\delta$ from $\bigcup_{h = 1}^{H} [\hat{a}_h,\hat{b}_h]$.
We can write for any $\tau \in \{\pm 1\}$ as $\N \rightarrow \infty$
\begin{equation}
\label{eq_x50}
\begin{split}
& \quad \prod_{i = 1}^N \bigg(1 + \frac{\tau}{\N}\,\frac{\theta_{h,h(i)}}{z - \frac{\ell_i}{\N} - \frac{\tau}{2\N}}\bigg) \\
& = \exp\Bigg(\sum_{i = 1}^N \bigg(\frac{\tau}{\N}\,\frac{\theta_{h,h(i)}}{z -
\frac{\ell_i}{\N} - \frac{\tau}{2\N}} -
\frac{1}{2\N^2}\,\frac{(\theta_{h,h(i)})^2}{(z - \frac{\ell_i}{\N} -
\frac{\tau}{2\N})^{2}} + \frac{\tau}{3\N^3}\,\frac{(\theta_{h,h(i)})^3}{(z -
\frac{\ell_i}{\N} - \frac{\tau}{2\N})^3} \bigg)+
O\bigg(\frac{1}{\N^3}\bigg)\Bigg).
\end{split}
\end{equation}
We also use an approximation
\[
\frac{1}{z - \frac{\ell_i}{\N} -
\frac{\tau}{2\N}} = \frac{1}{z -
\frac{\ell_i}{\N}}+\frac{\tau}{2\N} \frac{1}{(z -
\frac{\ell_i}{\N})^{2}}+\frac{1}{4\N^2} \frac{1}{(z -
\frac{\ell_i}{\N})^{3}}+ O\bigg(\frac{1}{\N^{3}}\bigg),
\]
valid for the same range of $z$, so that
\[
 \frac{1}{\N}\sum_{i=1}^N \frac{\theta_{h,h(i)}}{z -
\frac{\ell_i}{\N} - \frac{\tau}{2\N}} = \bigg\langle \bth_h \cdot \bigg(\frac{\boldsymbol{G}(z)}{\N} -\frac{\tau \partial_{z}\boldsymbol{G}(z)}{2\N^2}+\frac{\partial_z^2 \boldsymbol{G}(z)}{8\N^3}\bigg)\bigg\rangle+ O\bigg(\frac{1}{\N^{3}}\bigg).
\]
Similarly, we have
\begin{equation*}
\begin{split}
 \frac{1}{\N^2}\sum_{i=1}^N \frac{(\theta_{h,h(i)})^2}{(z -
\frac{\ell_i}{\N} - \frac{\tau}{2\N})^2} & = \bigg\langle (\bth_h)^2 \cdot \bigg(-\frac{\partial_z \boldsymbol{G}(z)}{\N^2} +\frac{\tau \partial_{z}^2 \boldsymbol{G}(z)}{2\N^3}\bigg)\bigg\rangle+ O\bigg(\frac{1}{\N^{3}}\bigg), \\
 \frac{1}{\N^3}\sum_{i=1}^N \frac{(\theta_{h,h(i)})^3}{(z -
\frac{\ell_i}{\N} - \frac{\tau}{2\N})^3} & =\frac{1}{2 \N^3} \big\langle (\bth_h)^3 \cdot \partial_z^2 \boldsymbol{G}(z)\big\rangle+ O\bigg(\frac{1}{\N^{3}}\bigg).
\end{split}
\end{equation*}
Plugging into \eqref{eq_x50}, we obtain
\begin{equation}
\label{eq_x91}
\begin{split}
 \prod_{i = 1}^N \bigg(1 + \frac{\tau}{\N}\,\frac{\theta_{h,h(i)}}{z - \frac{\ell_i}{\N} - \frac{\tau}{2\N}}\bigg) & = \exp\bigg(\frac{\tau}{\N}\big\langle \bth_h \cdot \boldsymbol{G}(z)\big\rangle + \frac{1}{2\N^2}\big\langle \bth_h(\bth_h - 1) \cdot \partial_{z}\boldsymbol{G}(z)\big\rangle \\
& \quad \phantom{\exp}\quad + \frac{\tau}{24 \N^3} \big\langle \bth_h(4\bth_h^2 - 6\bth_h +
3)\cdot \partial_{z}^2 \boldsymbol{G}(z)\big\rangle +
\frac{\mathbbm{c}_{1}}{\N^3}(z)\bigg),
\end{split}
\end{equation}
and $|\mathbbm{c}_{1}(z)|$ is bounded by some constant depending on $\delta$ but not on $\N$. We further transform \eqref{eq_x91} by decomposing
\[
\boldsymbol{G}_h(z) = \N \Gm_{\boldsymbol{\mu}}(z) + \Delta\boldsymbol{G}(z)
\]
as in Definition~\ref{def_correlators}, and expanding the exponential through a Taylor series
approximation
\[
\forall x \in [-1,1]\qquad \exp(x) = 1 + x + \frac{x^2}{2} + \frac{x^3}{6} + \frac{x^4}{24} + \mathbbm{c}_2(x) \quad \textnormal{with}\,\, |\mathbbm{c}_2(x)| \leq \frac{e |x|^5}{120}.
\]
We obtain
\begin{equation}
\label{eq_x92}
\begin{split}
 & \quad \prod_{i = 1}^N \bigg(1 + \frac{\tau}{\N}\,\frac{\theta_{h,h(i)}}{z - \frac{\ell_i}{\N} - \frac{\tau}{2\N}}\bigg)
 \\ & = \exp\bigg(\tau \big\langle \bth_h\cdot \Gm_{\boldsymbol{\mu}}\big\rangle + \frac{1}{2\N} \big\langle \bth_h(\bth_h - 1) \cdot \partial_{z}\Gm_{\boldsymbol{\mu}}\big\rangle + \frac{\tau}{24\N^2} \big\langle \bth_h(4\bth_h^2 - 6\bth_h + 3)\cdot \partial_{z}^2\Gm_{\boldsymbol{\mu}}\big\rangle \bigg) \\
 & \quad \times \bigg(1 + \frac{\tau}{\N}\big\langle \bth_h \cdot \Delta \boldsymbol{G}(z)\big\rangle + \frac{1}{2 \N^2} \big\langle \bth_h\cdot \Delta \boldsymbol{G}(z)\big\rangle^2 + \frac{\tau}{6 \N^3}\big\langle \bth_h\cdot \Delta \boldsymbol{G}(z)\big\rangle^3
 + \frac{1}{24\N^4} \big\langle \bth_h\cdot \Delta \boldsymbol{G}(z)\big\rangle^4 \\ & \qquad + \frac{1}{2\N^2} \big\langle \bth_h(\bth_h - 1)\cdot \partial_{z}\Delta \boldsymbol{G}(z)\big\rangle + \frac{\tau}{2\N^3} \big\langle \bth_h\cdot \Delta \boldsymbol{G}(z)\big\rangle \, \langle \bth_h(\bth_h - 1)\cdot \partial_{z}\Delta \boldsymbol{G}(z)\rangle + \mathbbm{c}_3(z)\bigg),
\end{split}
\end{equation}
where we claim that the remainder $\mathbbm{c}_3(z)$ is $O(\N^{-\frac{5}{2}+\eps})$ uniformly for fixed $\eps$ and $\delta$. Indeed, Corollary~\ref{Corollary_a_priory_1} yields that $|\Delta G_h (z)|=O(\N^{\frac{1}{2}+\eps})$ for any $h \in [H]$, in the sense that for any $k \in \amsmathbb{Z}_{> 0}$, the $k$-th moment is $O(\N^{\frac{k}{2} + k \eps})$. Moreover, using Cauchy integral formula in the form
\begin{equation}
\label{derCauchy}
\partial_{z}^d \Delta G_h(z) = \frac{(-1)^{d}}{d!} \oint \frac{\dd \zeta}{2\ii \pi}\,\frac{\Delta G_h(\zeta)}{(z - \zeta)^{d + 1}},
\end{equation}
where the contour of integration surrounds $z$ and no other singularity of the integrand, we conclude that $\partial_z^d \Delta G_h(z) =O(\N^{\frac{1}{2}+\eps})$ for each $d \geq 1$. Hence, the right-hand side of \eqref{eq_x92} includes explicitly all the terms which could potentially be larger than $O(\N^{-\frac{5}{2}+\eps})$: the next terms would be of the form
\[
\frac{(\Delta \boldsymbol{G}(z))^{\otimes 5}}{\N^5},\qquad \frac{(\partial_z \Delta \boldsymbol{G}_h(z))^{\otimes 2}}{\N^4},\qquad \frac{(\Delta \boldsymbol{G}(z))^{\otimes 2} \otimes \partial_z \Delta \boldsymbol{G}(z)}{\N^4},
\]
and others of even smaller order of magnitude in $\N$.

Returning to Definition~\ref{def_correlators} for the correlators and using the expansion \eqref{eq_Q_plus_determ_expansion} of the functions $Q_h^{\pm}$,
the expansion \eqref{eq_x92} implies for any $h \in [H]$
\begin{equation} \label{eq_x93}
\begin{split}
 & \quad \Phi^-_h(z) \cdot \E \Bigg[ \prod_{i=1}^N
 \bigg(1-\frac{\theta_{h,h(i)}}{\N z-\ell_i+\frac{1}{2}} \bigg)
  \Bigg]
 +
 \Phi^+_h(z) \cdot \E \Bigg[ \prod_{i=1}^N
 \bigg(1+\frac{\theta_{h,h(i)}}{\N z-\ell_i-\frac{1}{2}} \bigg)\Bigg] \\
 &=
  \frac{Q_h^-(z)}{\N} \big\langle \bth_h\cdot \boldsymbol{W}_1(z)\big\rangle + \frac{Q_h^+(z)}{2\N^2}\big\langle\bth_h\cdot W_1(z)\big\rangle^2+
 \frac{Q_h^+(z)}{2\N^2}\big\langle \bth_h^{\otimes 2} \cdot \boldsymbol{W}_2(z,z)\big\rangle \\
 & \quad + \frac{Q_h^+(z)}{2\N^2}\big\langle \bth_h(\bth_h - 1) \cdot \partial_{z} \boldsymbol{W}_1(z)\big\rangle + q_h^+(z) + \frac{q_h^{+,[1]}(z)}{\N} + \frac{q_h^{+,[2]}(z)}{\N^2} - \frac{\mathbbm{Err}_h^{(2)}(z)}{\N} + O\big(\N^{-\frac{5}{2}+\eps}\big),
 \end{split}
\end{equation}
where $\mathbbm{Err}_h^{(2)}(z)$ is defined by \eqref{Err2qqq}.

\medskip

\noindent \textsc{Step 2.} We are now ready to define $\widetilde R_h(z)$ using $R_h(z)$ of Theorem~\ref{Theorem_Nekrasov} recalled in the first line of \eqref{eq_x93}:
\begin{equation}
 \label{eq_x94}
\begin{split}
\widetilde R_h(z) & := \N R_h(z) - \N q_h^+(z) \\
& \quad - \mathbbm{1}_{\amsmathbb{V}_h}(\hat{a}_h') \cdot \left(\Res_{\zeta=\hat a_h - \frac{1}{2\N}} R_h(\zeta)\dd \zeta\right) \cdot \frac{\N}{z - \hat a_h + \frac{1}{2\N}} - \mathbbm{1}_{\amsmathbb{S}_h}(\hat{a}_h') \cdot \N R_h(\hat a'_h) \cdot \frac{z-\hat b'_h}{\hat a'_h - \hat b'_h} \\ & -  \mathbbm{1}_{\amsmathbb{V}_h}(\hat{b}_h')\cdot \left(\Res_{\zeta =\hat b_h + \frac{1}{2\N}} R_h(\zeta)\dd \zeta\right) \cdot \frac{\N}{z-\hat b_h - \frac{1}{2\N}}
 - \mathbbm{1}_{\amsmathbb{S}_h}(\hat{b}_h') \cdot \N R_h(\hat b'_h) \cdot \frac{z-\hat a'_h}{\hat b'_h - \hat a'_h}.
\end{split}
\end{equation}
We further set $\mathbbm{Err}_h^{(1)}(z)+\mathbbm{Err}_h^{(2)}(z)$ to be a function which makes the equation \eqref{eq_First_Nekr_asymptotic} true and check that all properties claimed in Proposition~\ref{Proposition_Nek_1_asymptotic_form} are satisfied.

We first check holomorphicity of $\widetilde{R}_h(z)$. By Theorem~\ref{Theorem_Nekrasov}, $R_h(z)$ is holomorphic for $z \in \amsmathbb{M}_h$ except perhaps for simple poles at $z = \hat a_h - \frac{1}{2\N}$ or $z = \hat b_h + \frac{1}{2\N}$. If $\hat{a}_h'$ is in a saturation, Assumption~\ref{Assumptions_extra} guarantees that $R_h(z)$ has no pole at $z = \hat{a}_h - \frac{1}{2\N}$. If $\hat{a}_h'$ is in a void, the simple pole of $R_h(z)$ is cancelled by the second line in \eqref{eq_x94}. By Assumption~\ref{Assumptions_offcrit} the point $\hat{a}_h'$ cannot be in a band, so we have covered all cases. A similar discussion can be made at $z = \hat{b}_h + \frac{1}{2\N}$. Furthermore, according to Theorem~\ref{Theorem_regularity_density}-(iii) the only possible singularities of $q_h^+(z)$ in $\amsmathbb{M}_h$ are located at $z = \hat{a}_h'$ if $\iota_h^- = 0$ and at $\hat{b}_h'$ if $\iota_h^+ = 0$. If $\hat{a}'_h$ is in a saturation, Assumption~\ref{Assumptions_extra} requires $\iota_h^- = 2$ so $q_h^+(z)$ is regular at $z = \hat{a}_h'$, while if $\hat{a}_h'$ is in a void, $q_h^+(z)$ is manifestly regular near $z = \hat{a}_h$ due to its definition in \eqref{eq_x95_1}. A similar argument holds near $\hat{b}_h'$. We conclude that $\widetilde{R}_h(z)$ from \eqref{eq_x94} is holomorphic in the whole domain $\amsmathbb{M}_h$.

Let us examine the behavior of $q_h^-$ at the endpoints. If $\hat{a}_h'$ is in a saturation, $e^{\langle \bth_h \Gm_{\boldsymbol{\mu}}(z)\rangle}$ has a simple zero at $z = \hat{a}'_h$. Then, in the definition \eqref{eq_x95_1} of $q_h^-(z)$: the first term vanishes at $z = \hat{a}'_h$; the exponential in the second term has a simple pole, but it is multiplied by the double zero of $\phi^-_h(z)$ required in Assumption~\ref{Assumptions_extra}. As a result, $q_h^-(\hat{a}_h') = 0$. We show $q_h^-(\hat{b}_h') = 0$ in a similar way except that it is now $e^{-\langle \bth_h \cdot \Gm_{\boldsymbol{\mu}}(z)\rangle}$ which has a simple zero and $\phi_h^+(z)$ which has a double zero by Assumption~\ref{Assumptions_extra}.

\medskip

\noindent \textsc{Step 3.} We now turn to the properties of the error terms. The properties of $\mathbbm{Err}_h^{(2)}(z)$ that were announced in the statement of Proposition~\ref{Proposition_Nek_1_asymptotic_form} are clear from the expression given in \eqref{eq_x93}. To justify the upper bound \eqref{eq_Err1_bound} for $\mathbbm{Err}_h^{(1)}(z)$, comparing \eqref{eq_x93} and \eqref{eq_x94} with the definition of $R_h(z)$ in Theorem~\ref{Theorem_Nekrasov} and $\widetilde R_h(z)$ in \eqref{eq_First_Nekr_asymptotic}, we conclude that $\mathbbm{Err}_h^{(1)}(z)$ is a sum of two terms: the first one is $\N$ times the $O(\N^{-\frac{5}{2}+\eps})$ remainder in \eqref{eq_x93} and satisfies \eqref{eq_Err1_bound} right away; the second one is the contribution of the last four terms in \eqref{eq_x94} and we need to show that these terms are small.

The term in \eqref{eq_x94} involving the residue at $\hat a_h -\frac{1}{2\N}$ is only present if $\hat{a}'_h$ is in a void. We recall that $R_h(z)$ coincides --- by definition in Theorem~\ref{Theorem_Nekrasov} --- with the first line of \eqref{eq_x93}, that is
\begin{equation}
\label{RhRhRh}
R_h(z) = \Phi^-_h(z) \cdot \E \Bigg[ \prod_{i=1}^N
 \bigg(1-\frac{\theta_{h,h(i)}}{\N z-\ell_i+\frac{1}{2}} \bigg)
  \Bigg]
 +
 \Phi^+_h(z) \cdot \E \Bigg[ \prod_{i=1}^N
 \bigg(1+\frac{\theta_{h,h(i)}}{\N z-\ell_i-\frac{1}{2}} \bigg)\Bigg].
\end{equation}
The second term in \eqref{RhRhRh} does not have any residue at $z = \hat{a}_h - \frac{1}{2\N}$. For the first term, the magnitude of its residue is bounded from above by
\begin{equation}
\label{eq_x95}
\begin{split}
& \P\big[\exists i \in [N] \quad|\quad \ell_i = \N \hat a_h\big] \cdot \Phi^-_h\bigg(\hat a_h -\frac{1}{2\N}\bigg) \\
& \times \bigg(1+\frac{1 - \delta_{h,1}}{\N}\cdot \frac{|\!|\boldsymbol{\Theta}|\!|}{\hat a_h -\hat b_{h-1}} +\frac{1 - \delta_{h,H}}{\N} \cdot \frac{|\!|\boldsymbol{\Theta}|\!|}{\hat a_{h+1} -\hat a_h} \bigg)^N \cdot \prod_{j=1}^{N-1} \bigg(1+\frac{1}{j}\bigg),
\end{split}
\end{equation}
where the first factor in the second line of \eqref{eq_x95} comes from bounding from above the factors with $h(i)\neq h$ under the expectation and the second factor from those with $h(i)=h$. Simplifying, we deduce the existence of $D >0$, such that the magnitude of the residue is bounded from above by
\[
 D \N \cdot \P\big[\exists i \in [N] \quad|\quad \ell_i = \N \hat a_h\big].
\]
In \eqref{eq_x94} the residue term only contributes if $\hat{a}'_h$ is in a void. In this situation, we can use Theorem~\ref{Theorem_ldpsup} which says the probability that there is a particle at $\N \hat{a}_h$ is exponentially small as $\N$ becomes large. Repeating the same argument for the right endpoint $\hat b_h$, we get the desired bound for the terms involving residues in \eqref{eq_x94}.

We proceed to the term in \eqref{eq_x94} which is present if $\hat{a}'_h$ is in a saturation and involves $R_h(\hat{a}_h')$. Fix $M \in \amsmathbb{Z}_{> 0}$ and denote $\mathcal{A}_M \subseteq \W_\N$ the event of having particles at the sites
\begin{equation}
\label{filledunu}\N \hat a_h, \N\hat a_h + \theta_{h,h}, \ldots, \N \hat a_h + M\theta_{h,h}.
\end{equation}
If $\hat{a}_h'$ is in a saturation, choosing $\eps' > 0$ small enough independent of $\N$, the event $\mathcal{A}_{\lfloor \eps' \N \rfloor}$ has probability exponentially close to $1$ as $\N$ becomes large, thanks to Theorem~\ref{Theorem_ldsaturated}. The two products under the expectation values in \eqref{RhRhRh} are uniformly bounded by a polynomial in $\N$ as in \eqref{eq_x95}, so up to exponentially small corrections as $\N \rightarrow \infty$ we can condition by the event $\mathcal{A}_M$ in the expectation values. On this event, we have a telescoping product
\begin{equation}
\label{telepou}
\begin{split}
 \prod_{m=0}^{M}
 \bigg(1-\frac{\theta_{h,h}}{\N z-\N \hat a_h-m\theta_{h,h}+\frac{1}{2}}\bigg) & =
 \prod_{m=0}^{M}
 \frac{\N z-\N \hat a_h-(m+1) \theta_{h,h}+\frac{1}{2}}{\N z-\N \hat a_h-m \theta_{h,h}+\frac{1}{2}}
 \\ & = \frac{\N z-\N \hat a_h- (M+1)\theta_{h,h}+\frac{1}{2}}{\N z-\N \hat a_h+\frac{1}{2}},
\end{split}
\end{equation}
representing the contribution of the particles \eqref{filledunu} to the product under the expectation value in the first term of $R_h(z)$ --- \textit{cf.} \eqref{RhRhRh}. There, it is multiplied with $\Phi^-_h(z)$. If $\theta_{h,h} = 1$, then in the right-hand side of \eqref{telepou} the denominator $\N \hat{a}'_h - \N \hat a_h + \frac{1}{2} = 1 - \theta_{h,h}$ vanishes, but in this situation we also have $\hat{a}_h - \frac{1}{2\N} = \hat{a}_h'$ and Assumption~\ref{Assumptions_extra} required that $\Phi^-_h(z)$ has a double zero at $z = \hat{a}'_h$, so the second term in $R_h(z)$ vanishes. If $\theta_{h,h} \neq 1$, the said denominator does not vanish and Assumption~\ref{Assumptions_extra} still requires a simple zero $\Phi_h^-(\hat{a}_h') = 0$. So the contribution of $\mathcal{A}_M$ in the first term of $R_h(\hat{a}_h')$ always vanishes, and thus $R_h(\hat{a}'_h)$ is always exponentially small as $\N \rightarrow \infty$.

Simultaneously, on the event $\mathcal{A}_{M}$ we have the telescoping product
\begin{equation}
\label{telepou2}
\begin{split}
 \prod_{m=0}^M
 \bigg(1+\frac{\theta_{h,h}}{\N z-\N \hat a_h-m\theta_{h,h}-\frac{1}{2}} \bigg) & = \prod_{m=0}^M
\frac{\N z-\N \hat a_h-(m-1)\theta_{h,h}-\frac{1}{2}}{\N z-\N \hat a_h-m\theta_{h,h}-\frac{1}{2}} \\
& = \frac{\N z-\N \hat a_h+\theta_{h,h}-\frac{1}{2}}{\N z-\N \hat a_h-M\theta_{h,h}-\frac{1}{2}},
\end{split}
\end{equation}
representing the contribution of the particles \eqref{filledunu} to the product under the expectation value in the second term of $R_h(z)$. As $z=\hat a'_h$ the numerator in the right-hand side of \eqref{telepou2} vanishes. Therefore, the contribution of $\mathcal{A}_M$ in the second term of $R_h(\hat{a}_h')$ vanishes. We eventually find that $R_h(\hat{a}_h')$ is exponentially small as $\N \rightarrow \infty$ and so is the term involving $\mathbbm{1}_{\amsmathbb{S}_h}(\hat{a}'_h)$ in \eqref{eq_x94}, once we remember that $z$ is kept at distance larger than $\delta$ from the segments of the ensemble. By a similar argument the same is true as well for the term involving $\mathbbm{1}_{\amsmathbb{S}_h}(\hat{b}'_h)$ in \eqref{eq_x94}. All in all, this justifies the estimate $\mathbbm{Err}_h^{(1)}(z) = O(\N^{-\frac{3}{2} + \eps})$.
\end{proof}

\subsection{Rewriting the second-order Nekrasov equation}

\begin{proposition} \label{Proposition_Nek_2_asymptotic_form} Consider a discrete ensemble satisfying Assumptions~\ref{Assumptions_Theta}, \ref{Assumptions_basic}, \ref{Assumptions_offcrit}, \ref{Assumptions_analyticity} and the additional Assumption~\ref{Assumptions_extra}. For any $h,h_2\in[H]$, any $z\in \amsmathbb M_h\setminus \amsmathbb{A}_h^{\mathfrak{m}}$ and $z_2\in\amsmathbb C\setminus \amsmathbb{A}_{h_2}^{\mathfrak{m}}$, we have
\begin{equation}
\label{eq_Second_Nekr_asymptotic}
\begin{split}
& \quad Q_{h}^{-}(z)\big\langle \bth_{h}\cdot \boldsymbol{W}_{2;\bullet,h_2}(z,z_2)\big\rangle + \frac{\delta_{h,h_2}\,Q_{h}^{-}(z)}{2(z_2 - z)^2} \\
& = \widetilde R_{h;h_2}(z;z_2)+ \mathbbm{Err}_{h;h_2}^{(1)}(z;z_2)+ \mathbbm{Err}_{h;h_2}^{(2)}(z;z_2).
\end{split}
\end{equation}
Here, $\widetilde R_{h;h_2}(z;z_2)$ is a meromorphic function of $z\in \amsmathbb M_h$, with only possible singularities being simple poles at $z=z_2 \pm \frac{1}{2\N}$. It is also holomorphic in $z_2\in\amsmathbb C\setminus \amsmathbb{A}_{h_2}^{\mathfrak{m}}$ except for the singularity we just mentioned. If $\hat a'_h \in \amsmathbb{S}_h$, then $\widetilde R_{h;h_2}(\hat a'_h;z_2)=0$; if $\hat b'_h\in \amsmathbb{S}_h$, then $\widetilde R_{h;h_2}(\hat b'_h;z_2)=0$.

Besides, for any $\N$-independent $\eps>0$ and disjoint pair of $\N$-independent compacts $\amsmathbb{K}_h \subset \amsmathbb{M}_h \setminus \amsmathbb{A}_h^{\mathfrak{m}}$ and $\amsmathbb{K}^{(2)}_{h_2} \subset \amsmathbb{C} \setminus \amsmathbb{A}_{h_2}^{\mathfrak{m}}$, there exists a constant $C>0$ depending only on $\eps$, these compacts and the constants in the assumptions, such that
   \begin{equation}
    \sup_{\substack{(z,z_2) \in \amsmathbb{K}_h \times \amsmathbb{K}_{h_2}^{(2)}}} \big|\mathbbm{Err}_{h;h_2}^{(1)}(z;z_2)\big|\leq C \N^{-\frac{1}{2}+\eps}.
   \end{equation}
The function $\mathbbm{Err}_{h;h_2}^{(2)}(z;z_2)$ is a finite sum --- with a number of terms which does not depend on anything --- of expressions of the form
  \begin{equation}
  \label{eq_Err2_bound_two}
  \frac{F(z)}{\N^{q}} \cdot \prod_{i=1}^k \E \left[ \big(\Delta G_{h_2}(z_2)\big)^{\tilde{l}_{i}} \cdot \prod_{j=1}^{l_i} \partial_z^{d_{i,j}} \Delta G_{g_{i,j}}(z) \right],
  \end{equation}
  for some $q,k,l_1,\ldots,l_k \in \amsmathbb{Z}_{> 0}$ and $\tilde{l}_1,\ldots,\tilde{l}_k \in \{0,1\}$ satisfying\footnote{In fact, there is exactly one $\tilde{l}_{i}$ equal to $1$ and the other are equal to $0$. So, \eqref{eq_Err2_bound_two_2} is the same constraint as \eqref{eq_Err2_bound2}.}
  \begin{equation}
  \label{eq_Err2_bound_two_2}
  q+1\leq \sum_{i=1}^k (l_i+\tilde{l}_i)\leq q+2,
  \end{equation}
  some function $F(z)$ holomorphic for $z \in \amsmathbb{K}_h$ and bounded from above in absolute value by $C$, and some $d_{i,j} \in \{0,1\}$ and $g_{i,j} \in [H]$.
\end{proposition}
The remainder term $\mathbbm{Err}_{h;h_2}^{(2)}(z;z_2)$ has an explicit expression on which the claimed properties can be checked:
\begin{equation}
\label{ERR2FORM}
\begin{split}
& \quad -\mathbbm{Err}_{h;h_2}^{(2)}(z;z_2)\\
&=\frac{Q_{h}^+(z)}{2\N}\big\langle \bth_{h}(\bth_{h} - 1)\cdot \partial_{z} \boldsymbol{W}_{2;\bullet,h_2}(z,z_2)\big\rangle
+\frac{Q_{h}^+(z)}{2\N}\big\langle \bth_{h}^{\otimes 2}\cdot \boldsymbol{W}_{3;\bullet,\bullet,h_2}(z,z,z_2)\big\rangle \\
& \quad + \frac{Q_h^+(z)}{\N}\big\langle \bth_h \cdot \boldsymbol{W}_1(z)\big\rangle\big\langle \bth_h \cdot \boldsymbol{W}_{2;\bullet,h_2}(z,z_2)\big\rangle
\\
&\quad +\frac{Q_{h}^-(z)}{6\N^2}\big\langle \bth_{h}^{\otimes 3}\cdot \boldsymbol{W}_{4;\bullet,\bullet,\bullet,h_2}(z,z,z,z_2)\big\rangle
+ \frac{Q_{h}^-(z)}{2\N^2}\big\langle\bth_{h}^{\otimes 2}\cdot \boldsymbol{W}_{2;\bullet,\bullet}(z,z)\big\rangle\big\langle \bth_{h}\cdot \boldsymbol{W}_{2;\bullet,h_2}(z,z_2)\big\rangle
 \\& \quad + \frac{Q_h^-(z)}{2\N^2} \big\langle \bth_h \cdot \boldsymbol{W}_1(z) \big\rangle^2 \big\langle \bth_h \cdot \boldsymbol{W}_{2;\bullet,h_2}(z,z_2)\big\rangle
 + \frac{Q_h^-(z)}{2\N^2}\big\langle \bth_h \cdot \boldsymbol{W}_1(z)\big\rangle\big\langle \bth_h^{\otimes 2} \cdot \boldsymbol{W}_{3;\bullet,\bullet,h_2}(z,z,z_2)\big\rangle.
\end{split}
\end{equation}
Similarly to Proposition~\ref{Proposition_Nek_1_asymptotic_form}, we think of \eqref{eq_Second_Nekr_asymptotic} as a linear equation for the unknown $\boldsymbol{W}_2$, and the terms in the right-hand side of \eqref{eq_Second_Nekr_asymptotic} will not contribute to the leading asymptotic term of the solution to this equation as $\N\rightarrow\infty$.

\begin{proof}[Proof of Proposition~\ref{Proposition_Nek_2_asymptotic_form}]
 We start from the second-order Nekrasov equation, that is Corollary~\ref{Corollary_higherNek} for $n=2$. It gives information on the analytic properties of the function $R_{h;h_2}(z;z_2)$, given in the alternative form \eqref{otherformR2}:
\begin{equation}
\label{R2R2R2}
\begin{split}
R_{h;h_2}(z;z_2) & = \sum_{\tau \in \{\pm 1\}} \Phi_h^{\tau}(z)\cdot \E^{(\textnormal{c})}\Bigg[\prod_{i = 1}^N\bigg(1 + \frac{\tau}{\N}\,\frac{\theta_{h,h(i)}}{z - \frac{\ell_i}{\N} - \frac{\tau}{2\N}}\bigg)\,,\, G_{h_2}(z_2)\Bigg] \\
 \\ & \quad + \delta_{h,h_2} \sum_{\tau \in \{\pm 1\}} \frac{z_2 - z + \frac{\tau}{2\N}}{ (z_2 - z)^2 - \frac{1}{4\N^2}} \cdot
 \Phi_h^{\tau}(z) \cdot \E\Bigg[\prod_{i = 1}^N\bigg(1 + \frac{\tau}{\N}\,\frac{\theta_{h,h(i)}}{z - \frac{\ell_i}{\N} - \frac{\tau}{2\N}}\bigg)\Bigg].
\end{split}
\end{equation}
Recalling the function $R_h(z)$ from \eqref{RhRhRh}, let us introduce
\begin{equation}
\label{eq_x97}
\begin{split}
R^{(1)}_{h;h_2}(z;z_2)&:= R_{h;h_2}(z;z_2) - R_h(z) \cdot \delta_{h,h_2} \frac{z_2 - z}{ (z_2 - z)^2 - \frac{1}{4\N^2}} \\
& = \sum_{\tau \in \{\pm 1\}}\,
 \Phi_h^{\tau}(z)\cdot \E^{(\textnormal{c})}\Bigg[\prod_{i = 1}^N\bigg(1 + \frac{\tau}{\N}\,\frac{\theta_{h,h(i)}}{z - \frac{\ell_i}{\N} - \frac{\tau}{2\N}}\bigg)\,,\, G_{h_2}(z_2)\Bigg]
 \\& \quad + \frac{\delta_{h,h_2}}{2\N} \frac{1}{ (z_2 - z)^2 - \frac{1}{4\N^2}} \sum_{\tau \in \{\pm 1\}}\tau
 \Phi_h^{\tau}(z)\cdot\E\Bigg[\prod_{i = 1}^N\bigg(1 + \frac{\tau}{\N}\,\frac{\theta_{h,h(i)}}{z - \frac{\ell_i}{\N} - \frac{\tau}{2\N}}\bigg)\Bigg].
\end{split}
\end{equation}
We use the Taylor expansion \eqref{eq_x92} under the expectation values. Since the covariance of two random variables vanishes if one of them is deterministic, we can replace $G_{h_2}(z_2)$ with $\Delta G_{h_2}(z_2)$ in \eqref{eq_x97}, and then use $|\Delta G_g (z)|=O(\N^{\frac{1}{2}+\eps})$ for any $g \in [H]$, $z$ uniformly away from $[\hat{a}_g,\hat{b}_g]$ and $\eps > 0$ arbitrarily small, as in the proof of Proposition~\ref{Proposition_Nek_1_asymptotic_form}. This leads to
\begin{equation}
\label{eq_x96}
\N \cdot R^{(1)}_{h,h_2}(z,z_2) = Q_{h}^{-}(z)\big\langle \bth_{h}\cdot \boldsymbol{W}_{2;\bullet,h_2}(z,z_2)\big\rangle + \delta_{h,h_2}\frac{Q_{h}^{-}(z)}{2(z_2 - z)^2} -\mathbbm{Err}_{h;h_2}^{(2)}(z) + O\big( \N^{-\frac{1}{2}+\eps}\big),
\end{equation}
where
\begin{equation}
 \label{eq_x127}
\begin{split}
& \quad -\mathbbm{Err}_{h;h_2}^{(2)}(z) \\
&= \frac{Q_h^+(z)}{2 \N} \cdot \E^{(\textnormal{c})} \Big[\big\langle \bth_h\cdot \Delta \boldsymbol{G}(z)\big\rangle^2\,,\, \Delta G_{h_2}(z_2)\Big]+ \frac{Q_h^-(z)}{6 \N^2}\cdot \E^{(\textnormal{c})}\Big[\big\langle \bth_h\cdot \Delta \boldsymbol{G}(z)\big\rangle^3\,,\, \Delta G_{h_2}(z_2)\Big]
\\ &+\frac{Q_h^+(z)}{2\N}\cdot\E^{(\textnormal{c})} \Big[\big\langle \bth_h(\bth_h - 1)\cdot \partial_{z}\Delta \boldsymbol{G}(z)\big\rangle \,,\, \Delta G_{h_2}(z_2)\Big]
\\ &+\frac{Q_h^-(z)}{2\N^2}\cdot\E^{(\textnormal{c})} \Big[\big\langle \bth_h\cdot \Delta \boldsymbol{G}(z)\big\rangle\big\langle \bth_h(\bth_h - 1)\cdot \partial_{z}\Delta \boldsymbol{G}(z)\big\rangle\, ,\, \Delta G_{h_2}(z_2) \Big].
\end{split}
\end{equation}
After expressing it in terms of correlators, we retrieve the expression \eqref{ERR2FORM} for $\mathbbm{Err}_{h;h_2}^{(2)}(z;z_2)$, while its claimed properties can be checked immediately by re-expressing cumulants in terms of moments --- \textit{cf.} Section~\ref{Sec_cumu_corr}. We further set
\begin{equation}
\begin{split}
& \quad \widetilde R_{h;h_2}(z;z_2) \\
& := \N \cdot R^{(1)}_{h;h_2}(z;z_2) \\
 & \quad - \mathbbm{1}_{\amsmathbb{V}_h}(\hat{a}'_h) \cdot \left(\Res_{\zeta =\hat a_h - \tfrac{1}{2\N}} R^{(1)}_{h;h_2}(\zeta;z_2)\dd \zeta\right) \cdot \frac{\N}{z-\hat a_h + \frac{1}{2\N}} - \mathbbm{1}_{\amsmathbb{S}_h}(\hat{a}'_h) \cdot \N R^{(1)}_{h;h_2}(\hat a'_h;z_2) \cdot \frac{z-\hat b'_h}{\hat a'_h - \hat b'_h} \\
 & \quad - \mathbbm{1}_{\amsmathbb{V}_h}(\hat{b}'_h) \cdot \left( \Res_{z'=\hat b_h + \tfrac{1}{2\N}} R^{(1)}_{h;h_2}(\zeta;z_2) \dd \zeta\right) \cdot \frac{\N}{z-\hat b_h - \frac{1}{2\N}} -\mathbbm{1}_{\amsmathbb{S}_h}(\hat{b}'_h) \cdot \N R^{(1)}_{h;h_2}(\hat b'_h;z_2) \cdot \frac{z-\hat a'_h}{\hat b'_h - \hat a'_h}.
\end{split}
\end{equation}
We can then argue as in Step 2. of the proof of Proposition~\ref{Proposition_Nek_1_asymptotic_form} that $\widetilde{R}_{h;h_2}(z;z_2)$ has all the desired analytic properties, using the key information on the analytic properties of $R_{h;h_2}(z;z_2)$ from Corollary~\ref{Corollary_higherNek} as well as Assumptions~\ref{Assumptions_offcrit} and \ref{Assumptions_extra}. We define $\mathbbm{Err}_1(z)$ to be the sum of $O\big( \N^{-\frac{1}{2}+\eps}\big)$ in \eqref{eq_x96} and the four terms of the difference $ \widetilde R_{h;h_2}(z;z_2)-\N R^{(1)}_{h;h_2}(z;z_2)$. It can be then shown that these four terms are exponentially small as $\N \rightarrow \infty$ following the same method as in Step 3. of Proposition~\ref{Proposition_Nek_1_asymptotic_form}, relying on Theorems~\ref{Theorem_ldpsup} and \ref{Theorem_ldsaturated} saying that with overwhelming probability there are no particles in voids and no holes in saturations.
\end{proof}

\subsection{Rewriting the higher-order Nekrasov equations}

\begin{proposition} \label{Proposition_Nek_higher_asymptotic_form} Consider a discrete ensemble satisfying Assumptions~\ref{Assumptions_Theta}, \ref{Assumptions_basic}, \ref{Assumptions_offcrit}, \ref{Assumptions_analyticity} and the additional Assumption~\ref{Assumptions_extra}. For any $n \geq 3$, $h,h_2,\ldots,h_n \in [H]$, any $z\in \amsmathbb M_h\setminus \amsmathbb{A}_h^{\mathfrak{m}}$ and $z_i\in\amsmathbb C\setminus \amsmathbb{A}_{h_i}^{\mathfrak{m}}$ for $i \in \llbracket 2,n\rrbracket$, we have
\begin{equation}
\label{eq_Higher_Nekr_asymptotic}
\begin{split}
& \quad Q_{h}^{-}(z)\big\langle \bth_{h}\cdot \boldsymbol{W}_{n;\bullet,h_2,h_3,\ldots,h_n}(z,z_2,z_3,\ldots,z_n)\big\rangle \\
& = \widetilde R_{h;h_2,\ldots,h_n}(z;z_2,\ldots,z_n)+ \mathbbm{Err}_{h;h_2,\ldots,h_n}^{(1)}(z;z_2,\ldots,z_n)+ \mathbbm{Err}^{(2)}_{h;h_2,\ldots,h_n}(z;z_2,\ldots,z_n).
\end{split}
 \end{equation}
Here, $R_{h;h_2,\ldots,h_n}(z;z_2,\ldots,z_n)$ is a meromorphic function of $z\in \amsmathbb M_h$ with only possible singularities being simple poles at $z=z_i \pm \frac{1}{2\N}$ for $i \in \llbracket 2,n\rrbracket$. It is also holomorphic in $z_i\in\amsmathbb C\setminus \amsmathbb{A}^{\mathfrak{m}}_{h_i}$ for $i \in \llbracket 2,n\rrbracket$ except for the singularities we just mentioned. If $\hat a'_h \in \amsmathbb{S}_h$, then $\widetilde R_{h;h_2,\ldots,h_n}(\hat a'_h; z_2, \ldots z_n)=0$; if $\hat b'_h\in \amsmathbb{S}_h$, then $\widetilde R_{h;h_2,\ldots,h_n}(\hat b'_h; z_2,\ldots,z_n)=0$.

Besides, for any $\N$-independent $\eps > 0$, any $\N$-independent compact sets $\amsmathbb{K}_h \subset \amsmathbb{M}_h \setminus \amsmathbb{A}_{h}^{\mathfrak{m}}$ and $\amsmathbb{K}^{(i)}_{h_i} \subset \amsmathbb{C} \setminus (\amsmathbb{A}_{h_i}^{\mathfrak{m}} \cup \amsmathbb{K}_h)$ indexed by $i \in \llbracket 2,n\rrbracket$, there exists a constant $C > 0$ depending only on $\eps$, these compacts and the constants in the assumptions, such that
\[
\sup_{\substack{z\in \amsmathbb{K}_h \\ z_i \in \amsmathbb{K}_{h_i}^{(i)}}} \big|\mathbbm{Err}^{(1)}_{h;h_2,\ldots,h_n}(z;z_2,\ldots,z_n)\big|\leq C \N^{-\frac{1}{2}+\eps}.
\]
The function $\mathbbm{Err}^{(2)}_{h;h_2,\ldots,h_n}(z;z_2,\ldots,z_n)$ is a finite sum --- with a number of terms that only depends on $n$ --- of expressions of the form
  \begin{equation}
  \label{eq_Err2_bound_multi}
  \frac{F(z)}{\N^{q}} \cdot \prod_{i=1}^k \E\Bigg[ \prod_{m=2}^n \big(\Delta G_{h_m}(z_m)\big)^{\tilde{l}_{i,m}} \cdot \prod_{j=1}^{l_i} \big(\partial_z^{d_{i,j}} \Delta G_{g_{i,j}}(z)\big) \Bigg],
  \end{equation}
  for some $q,k \in \amsmathbb{Z}_{> 0}$, $l_1,\ldots,l_k \in \amsmathbb{Z}_{\geq 0}$, $\tilde{l}_{i,m} \in \{0,1\}$ satisfying\footnote{In fact, for each $m \in \llbracket 2,n\rrbracket$, there is exactly one $\tilde{l}_{i,m}$ equal to $1$ and the other are equal to $0$.}
  \begin{equation}
  \label{eq_Err2_bound_multi_2}
  q+1\leq \sum_{i=1}^k \bigg(l_i+\sum_{m=2}^{n} \tilde{l}_{i,m}\bigg)\leq q+n,
  \end{equation}
some function $F(z)$ holomorphic for $z \in \amsmathbb{K}_h$ and bounded from above in absolute value by $C$, and some $d_{i,j} \in \{0,1\}$ and $g_{i,j} \in [H]$.
\end{proposition}

\begin{proof} We continue with the strategy used for the proofs of Propositions~\ref{Proposition_Nek_1_asymptotic_form} and \ref{Proposition_Nek_2_asymptotic_form}, using this time the $n$-th order Nekrasov equation. Take $n \geq 2$ and consider the function $R_{h;h_2,\ldots,h_n}(z;z_2,\ldots,z_n)$ of Corollary~\ref{Corollary_higherNek}. Writing separately the terms corresponding to $\#J = 0$, $\# J = n - 1$ and $\# J = n - 2$, we obtain
\begin{equation}
 \label{eq_x98}
\begin{split}
 & \quad R_{h;h_2,\ldots,h_n}(z;z_2,\ldots,z_n) \\
 & = \sum_{\tau \in \{\pm 1\}} \Phi_h^{\tau}(z)\cdot \E^{(\textnormal{c})}\Bigg[\prod_{i = 1}^N\bigg(1 + \frac{\tau}{\N}\,\frac{\theta_{h,h(i)}}{z - \frac{\ell_i}{\N} - \frac{\tau}{2\N}}\bigg)\,,\,\big(G_{h_j}(z_j)\big)_{j = 2}^{n}\Bigg] \\
 & \quad + \sum_{\tau \in \{\pm 1\}} \bigg(\prod_{j=2}^n \frac{\delta_{h,h_j}}{z_j - z - \frac{\tau}{2\N}}\bigg) \cdot
 \Phi_h^{\tau}(z)\cdot \E\Bigg[\prod_{i = 1}^N\bigg(1 + \frac{\tau}{\N}\,\frac{\theta_{h,h(i)}}{z - \frac{\ell_i}{\N} - \frac{\tau}{2\N}}\bigg)\Bigg] \\
 &\quad + \sum_{m=2}^n \sum_{\tau \in \{\pm 1\}}\bigg(\prod_{j \neq m} \frac{\delta_{h,h_j}}{z_j - z - \frac{\tau}{2\N}}\bigg) \cdot
 \Phi_h^{\tau}(z)\cdot \E^{(\textnormal{c})}\Bigg[\prod_{i = 1}^N\bigg(1 + \frac{\tau}{\N}\,\frac{\theta_{h,h(i)}}{z - \frac{\ell_i}{\N} - \frac{\tau}{2\N}}\bigg)\,,\,G_{h_m}(z_m)\Bigg] \\
 & \quad + \sum_{\substack{J \subset \llbracket 2,n\rrbracket\\ 0< \# J < n-2 \\ \tau \in \{\pm 1\}}} \bigg(\prod_{j \in J} \frac{\delta_{h,h_j}}{z_j - z - \frac{\tau}{2\N}}\bigg) \cdot \Phi_h^{\tau}(z)\cdot \E^{(\textnormal{c})}\Bigg[\prod_{i = 1}^N\bigg(1 + \frac{\tau}{\N}\,\frac{\theta_{h,h(i)}}{z - \frac{\ell_i}{\N} - \frac{\tau}{2\N}}\bigg)\,,\,\big(G_{h_j}(z_j)\big)_{j \notin J}\Bigg].
\end{split}
\end{equation}
Recalling the functions $R_h(z)$ from \eqref{RhRhRh} and $R^{(1)}_{h;h_m}(z;z_m)$ from \eqref{R2R2R2}-\eqref{eq_x97}, we introduce
\begin{equation*}
\begin{split} R^{(1)}_{h;h_2,\ldots,h_n}(z;z_2,\ldots,z_n) & = R_{h;h_2,\ldots,h_n}(z;z_2,\ldots,z_n) - R_h(z) \cdot \prod_{j=2}^n \delta_{h,h_j} \frac{z_j-z}{(z_j - z)^2 - \tfrac{1}{4\N^2}} \\
& - \sum_{m = 2}^{n} R^{(1)}_{h;h_m}(z;z_m) \cdot \prod_{j \neq m} \delta_{h,h_j} \frac{z_j-z}{(z_j - z)^2 - \frac{1}{4\N^2}}.
\end{split}
\end{equation*}
The result of the subtraction is to cancel the second and third lines in \eqref{eq_x98} up to small remainders. Hence,
\begin{equation}
\label{eq_x99}
\begin{split}
& \quad R^{(1)}_{h;h_2,\ldots,h_n}(z;z_2,\ldots,z_n) \\
& = \sum_{\tau \in \{\pm 1\}}
 \Phi_h^{\tau}(z)\cdot\E^{(\textnormal{c})}\Bigg[\prod_{i = 1}^N\bigg(1 + \frac{\tau}{\N}\,\frac{\theta_{h,h(i)}}{z - \frac{\ell_i}{\N} - \frac{\tau}{2\N}}\bigg)\,,\,\big(G_{h_j}(z_j)\big)_{j = 2}^n\Bigg] \\ & \quad + \sum_{\substack{J \subset \llbracket 2,n\rrbracket \\ 0< \# J <n-2 \\ \tau \in \{\pm 1\}}}\bigg(\prod_{j \in J} \frac{\delta_{h,h_j}}{z_j - z - \frac{\tau}{2\N}}\bigg) \cdot
 \Phi_h^{\tau}(z)\cdot \E^{(\textnormal{c})}\Bigg[\prod_{i = 1}^N\bigg(1 + \frac{\tau}{\N}\,\frac{\theta_{h,h(i)}}{z - \frac{\ell_i}{\N} - \frac{\tau}{2\N}}\bigg)\,,\,\big(G_{h_j}(z_j)\big)_{j \notin J}\Bigg]
\\ &\quad + O\bigg(\frac{1}{\N^{2}}\bigg) + O\bigg(\frac{1}{\N}\bigg) \cdot \Bigg\{\sum_{m=2}^n \Bigg| \,\E^{(\textnormal{c})}\Bigg[\prod_{i = 1}^N\bigg(1 + \frac{\tau}{\N}\,\frac{\theta_{h,h(i)}}{z - \frac{\ell_i}{\N} - \frac{\tau}{2\N}}\bigg)\,,\,G_{h_m}(z_m)\Bigg]\,\Bigg|\Bigg\}.
\end{split}
\end{equation}
Here and throughout this proof, the $O(\cdot)$ remainders, as function of $z,z_2,\ldots,z_n$, might have singularities at $z=z_j \pm \tfrac{1}{2\N}$ for $j \in \llbracket 2,n\rrbracket$, but they are uniform over $(z,z_2,\ldots,z_n)$ in $\N$-independent compact subsets of the complement of this singular locus in $(\amsmathbb{M}_h \setminus \amsmathbb{A}_h^{\mathfrak{m}}) \times \prod_{j = 2}^{n} \big(\amsmathbb{C} \setminus \amsmathbb{A}_{h_j}^{\mathfrak{m}}\big)$.

We continue this process of simplification in \eqref{eq_x99} using the functions from the $m$-th order Nekrasov equation with $m < n$. For this we split the second line of \eqref{eq_x99} as
\begin{equation}
\label{eq_x107}
\begin{split}
& \sum_{\substack{J \subset \llbracket 2,n\rrbracket \\ 0< \# J <n-2 \\ \tau \in \{\pm 1\}}} \bigg(\prod_{j \in J} \delta_{h,h_j} \frac{z_j - z}{(z_j - z) - \frac{1}{4\N^2}}\bigg) \cdot
 \Phi_h^{\tau}(z)\cdot \E^{(\textnormal{c})}\Bigg[\prod_{i = 1}^N\bigg(1 + \frac{\tau}{\N}\,\frac{\theta_{h,h(i)}}{z - \frac{\ell_i}{\N} - \frac{\tau}{2\N}}\bigg)\,,\,\big(G_{h_j}(z_j)\big)_{j \notin J}\Bigg] \\
 & + \sum_{\substack{J \subseteq \llbracket 2,n\rrbracket\\ 0< \# J <n-2 \\ \tau \in \{\pm 1\}}} \bigg( \prod_{j \in J} \frac{\delta_{h,h_j} }{(z_j - z) - \frac{1}{4\N^2}} \bigg) \cdot O\bigg(\frac{1}{\N}\bigg) \cdot \E^{(\textnormal{c})}\Bigg[\prod_{i = 1}^N\bigg(1 + \frac{\tau}{\N}\,\frac{\theta_{h,h(i)}}{z - \frac{\ell_i}{\N} - \frac{\tau}{2\N}}\bigg)\,,\,\big(G_{h_j}(z_j)\big)_{j \notin J}\Bigg]
\end{split}
\end{equation}
For each fixed set $J$, we recognize in the first line of \eqref{eq_x107} the first line of \eqref{eq_x99} with the $(n-\#J)$ variables $z,(z_j)_{j\notin J}$, up to multiplication the rational function in prefactor. Hence, by further subtractions, we can reduce the second line of \eqref{eq_x99} to lower-order terms. Let us denote the resulting function $R^{(2)}_{h;h_2,\ldots,h_n}(z;z_2,\ldots,z_n)$, so that
\begin{equation}
\label{eq_x108}
\begin{split}
 & \quad R^{(2)}_{h;h_2,\ldots,h_n}(z;z_2,\ldots,z_n) \\
 & = \sum_{\tau \in \{\pm 1\}} \Phi_h^{\tau}(z)\cdot\E^{(\textnormal{c})}\Bigg[\prod_{i = 1}^N\bigg(1 + \frac{\tau}{\N}\,\frac{\theta_{h,h(i)}}{z - \frac{\ell_i}{\N} - \frac{\tau}{2\N}}\bigg)\,,\,\big(G_{h_j}(z_j)\big)_{j = 2}^{n}\Bigg]
\\ & \quad + \sum_{\substack{J \subset \llbracket 2,n\rrbracket\\ 0< \#J <n-1 \\ \tau \in \{\pm 1\}}} O\bigg(\frac{1}{\N}\bigg) \cdot \E^{(\textnormal{c})}\Bigg[\prod_{i = 1}^N\bigg(1 + \frac{\tau}{\N}\,\frac{\theta_{h,h(i)}}{z - \frac{\ell_i}{\N} - \frac{\tau}{2\N}}\bigg)\,,\,\big(G_{h_j}(z_j)\big)_{j \notin J}\Bigg] + O\bigg(\frac{1}{\N^{2}}\bigg).
\end{split}
\end{equation}
 Since we constructed the function $R^{(2)}_{h;h_2,\ldots,h_n}(z;z_2,\ldots,z_n)$ by iteratively summing and sub\-trac\-ting the ho\-lo\-mor\-phic terms of Corollary~\ref{Corollary_higherNek}, it satisfies the same analytic properties as $R_{h;h_2,\ldots,h_n}(z;z_2,\ldots,z_n)$ in this corollary.

The next step is to perform the Taylor expansion of the product under the expectation values. We already did so in \eqref{eq_x92}, but this time we might need to go further. More precisely, we would like to keep expanding, until the bound $\Delta G_h(z)=O(\N^{\frac{1}{2}+\eps})$ of Corollary~\ref{Corollary_a_priory_1} guarantees that the remainder in transformed \eqref{eq_x108} is $O(\N^{-\frac{3}{2}+\eps})$. We also replace $G_{h_j}(z_j)$ by $\Delta G_{h_j}(z_j)$ inside the cumulants of order $2$ and higher, which does not change the result since these random variables differ by a deterministic part. The outcome is:
\begin{equation}
\begin{split}
\label{eq_x109} & \quad R^{(2)}_{h;h_2,\ldots,h_n}(z;z_2,\ldots,z_n)
\\& = \frac{Q_h^-(z)}{\N} \big\langle \bth_{h}\cdot \boldsymbol{W}_{n;\bullet,h_2,h_3,\ldots,h_n}(z,z_2,z_3,\ldots,z_n)\big\rangle -\frac{1}{\N} \mathbbm{Err}_{h;h_2,\ldots,h_n}^{(2)}(z; z_2,\ldots,z_n)+ O\big(\N^{-\frac{3}{2}+\eps}\big),
\end{split}
\end{equation}
where the function $\mathbbm{Err}^{(2)}_{h;h_2,\ldots,h_n}(z;z_2,\ldots,z_n)$ combines most of the terms of the Taylor expansion. We claim that it is a sum of terms of the form \eqref{eq_Err2_bound_multi} satisfying the inequalities \eqref{eq_Err2_bound_multi_2}. Indeed, the general form \eqref{eq_Err2_bound_multi} is obtained as in \eqref{eq_x92}. For the bound on degrees of \eqref{eq_Err2_bound_multi_2}, we observe that the first term of the expansion of the first line in \eqref{eq_x108} involves the $n$-th cumulant and the $\frac{1}{\N}$ prefactor which cancels with $\frac{1}{\N}$ in front of $\mathbbm{Err}^{(2)}_{h;h_2,\ldots,h_n}(z;z_2,\ldots,z_n)$ in \eqref{eq_x109}. Thus, for this term the second inequality of \eqref{eq_Err2_bound_multi_2} is actually an equality. Each next term in the expansion involves increasing powers of the components of $\Delta \boldsymbol{G}(z)$ --- or its derivatives --- and increasing powers of $\frac{1}{\N}$ by the same amount, hence, the equality is preserved. If we re-express the cumulants in terms of centered moments, the degrees are preserved. If we Taylor expand the second line in \eqref{eq_x108} the argument is similar: for $\#J=n-k$ with $k \in \llbracket 2,n-1\rrbracket$, the first term in the expansion involves the $k$-th order cumulant and a prefactor $\frac{1}{\N^{2}}$; for the next terms the order of the moment and the power of $\frac{1}{\N}$ grow in accord. We finally set
\begin{equation*}
\begin{split}
\widetilde R_{h;h_2,\ldots,h_n}(z;z_2,\ldots,z_n) & := \N \cdot R^{(2)}_{h;h_2,\ldots,h_n}(z;z_2,\ldots,z_n) \\
& \quad - \mathbbm{1}_{\amsmathbb{V}_h}(\hat{a}'_h) \cdot \left(\Res_{\zeta=\hat a_h - \tfrac{1}{2\N}} R^{(2)}_{h;h_2,\ldots,h_n}(\zeta;z_2,\ldots,z_n)\dd \zeta\right) \cdot \frac{\N}{z-\hat a_h + \frac{1}{2\N}} \\
& \quad
 -\mathbbm{1}_{\amsmathbb{S}_h}(\hat{a}'_h) \cdot \N R^{(2)}_{h;h_2,\ldots,h_n}(\hat a'_h;z_2,\ldots,z_n) \cdot \frac{z-\hat b'_h}{\hat a'_h - \hat b'_h}
 \\& \quad - \mathbbm{1}_{\amsmathbb{V}_h}(\hat{b}'_h) \cdot \left( \Res_{z'=\hat b_h + \tfrac{1}{2\N}} R^{(2)}_{h;h_2,\ldots,h_n}(\zeta;z_2,\ldots,z_n)\dd \zeta \right) \cdot \frac{\N}{z-\hat b_h - \frac{1}{2\N}} \\
 & \quad -\mathbbm{1}_{\amsmathbb{S}_h}(\hat{b}'_h) \cdot \N R^{(2)}_{h;h_2,\ldots,h_n}(\hat b'_h;z_2,\ldots,z_n) \cdot\frac{z-\hat a'_h}{\hat b'_h - \hat a'_h}.
\end{split}
\end{equation*}

By construction, the function $\widetilde R_{h;h_2,\ldots,h_n}(z;z_2,\ldots,z_n)$ satisfies all the analytic properties claimed in Proposition~\ref{Proposition_Nek_higher_asymptotic_form}. The function $\mathbbm{Err}^{(1)}_{h;h_2,\ldots,h_n}(z;z_2,\ldots,z_n)$ is the sum of $\N \cdot O\big( \N^{-\frac{3}{2}+\eps}\big)$ from \eqref{eq_x109} and the four terms of the difference $ \widetilde R_{h;h_2}(z;z_2,\ldots,z_n)-\N R^{(2)}_{h;h_2,\ldots,z_n}(z;z_2,\ldots,z_n)$. These terms are exponentially small as $N\rightarrow\infty$ by the same argument as in the proof of Propositions~\ref{Proposition_Nek_1_asymptotic_form} and \ref{Proposition_Nek_2_asymptotic_form} based on the large deviation estimates of Theorems~\ref{Theorem_ldpsup} and \ref{Theorem_ldsaturated} for the position of particles and holes.
\end{proof}

\section{Upper bound on the moments of \texorpdfstring{$\Delta G_h(z)$}{DeltaGh(z)}.}
\label{Secgg83gn5g}
In the next section we are going to show that Propositions~\ref{Proposition_Nek_1_asymptotic_form}, \ref{Proposition_Nek_2_asymptotic_form}, and \ref{Proposition_Nek_higher_asymptotic_form} imply Theorem~\ref{Theorem_correlators_expansion}. An important step in the argument is to establish a bound showing that the $\mathbbm{Err}^{(2)}$ terms in these propositions are negligible as $\N\rightarrow\infty$. In this section we prove such a bound, by improving the much weaker --- and not sufficient for a direct conclusion --- bound of Corollary~\ref{Corollary_a_priory_1}. Essentially, we show that the only way for all the asymptotic expansions of Propositions~\ref{Proposition_Nek_1_asymptotic_form}, \ref{Proposition_Nek_2_asymptotic_form}, and \ref{Proposition_Nek_higher_asymptotic_form} to hold simultaneously is that the $\mathbbm{Err}^{(2)}$ terms are indeed small.

\begin{proposition} \label{Proposition_uniform_moments_bound}
 Consider a discrete ensemble satisfying Assumptions~\ref{Assumptions_Theta}, \ref{Assumptions_basic}, \ref{Assumptions_offcrit}, \ref{Assumptions_analyticity} and the additional Assumption~\ref{Assumptions_extra}. For any tuple of compact sets $\amsmathbb{K}_h\subset \amsmathbb C\setminus \amsmathbb{A}_h^{\mathfrak{m}}$ indexed by $h \in [H]$, there exists a $\N$-independent sequence of positive constants $(C_n)_{n \geq 1}$ depending on these compacts and the constants in the assumptions, such that, for any $n \in \amsmathbb{Z}_{> 0}$, $h_1,\ldots,h_n \in [H]$ and $(z_1,\ldots,z_n) \in \prod_{i = 1}^{n} \amsmathbb{K}_{h_i}$ we have
\begin{equation}
\label{eq_x168}
 \E\Bigg[\bigg| \prod_{i=1}^n \Delta G_{h_i}(z_i) \bigg|\Bigg] \leq C_n.
\end{equation}
\end{proposition}
\begin{proof} \phantom{s}

\noindent \textsc{Step 1.} As a warm-up and to demonstrate how the method works, we prove
\begin{equation}
\label{uzgbfugbg}
\forall h \in [H] \qquad \big|\E\big[\Delta G_h(z)\big]\big|=\big|W_{1;h}(z)\big| \leq C_1.
\end{equation}
Contrarily to the final inequality \eqref{eq_x168} we aim at, the absolute value in \eqref{uzgbfugbg} is outside the expectation value, not inside. We look at \eqref{eq_First_Nekr_asymptotic} in Proposition~\ref{Proposition_Nek_1_asymptotic_form} and make the following observations, which hold for each $\eps > 0$.
 \begin{itemize}
 \item By \eqref{eq_Q_plus_determ_expansion}, we have $Q_h^-(z)=q_h^-(z)+O\big(\frac{1}{\N}\big)$.
 \item By Corollary~\ref{Corollary_a_priory_1}, we have $\frac{Q_h^+(z)}{2\N}\langle \bth_h(\bth_h - 1) \cdot \partial_{z} \boldsymbol{W}_1(z)\rangle=O(\N^{-\frac{1}{2}+\eps})$.
 \item By definition in \eqref{eq_x215}, we have $\frac{1}{\N}q_h^{+,[2]}(z)=O\big(\frac{1}{\N}\big) =O(\N^{-\frac{1}{2}+\eps})$.
 \item Applying Corollary~\ref{Corollary_a_priory_1}, we bound all terms \eqref{FERR2} in $\mathbbm{Err}^{(2)}_h(z)$ with help of the inequality \eqref{eq_Err2_bound2}. We see that all terms are $O(\N^{\frac{1-q}{2} + \eps})$, which is $O(\N^{-\frac{1}{2} + \eps})$ since $q \geq 2$.
 \end{itemize}
 All the $O(\cdot)$ errors mentioned above are uniform over $z$ belonging to compact subsets of $\amsmathbb M_h\setminus \amsmathbb{A}_h^{\mathfrak{m}}$. Hence, \eqref{eq_First_Nekr_asymptotic} can be transformed into
\begin{equation}
 \label{eq_x112}
 \begin{split}
 & \quad q_h^-(z)\big\langle \bth_h\cdot \boldsymbol{W}_1(z)\big\rangle + q_h^{+,[1]}(z) + \frac{Q_h^+(z)}{2\N}\big\langle\bth_h\cdot \boldsymbol{W}_1(z)\rangle^2 +
 \frac{Q_h^+(z)}{2\N}\big\langle \bth_h^{\otimes 2} \cdot \boldsymbol{W}_2(z,z)\big\rangle
 \\ &=\,\, \widetilde R_h(z)+ O(\N^{-\frac{1}{2}+\eps}).
\end{split}
\end{equation}
Let us emphasize one particular term which got absorbed into $O(\N^{-\frac{1}{2}+\eps})$ in the last formula: this is the $O(\frac{1}{\N})$ remainder from the expansion of $Q_h^-(z)$ multiplied by $\langle \bth_h\cdot \boldsymbol{W}_1(z)\rangle $. To see that this indeed contributes as $O(\N^{-\frac{1}{2}+\eps})$ we use again Corollary~\ref{Corollary_a_priory_1}.

At this stage we cannot yet ignore the third and fourth terms in the first line of \eqref{eq_x112}: Corollary~\ref{Corollary_a_priory_1} only shows that they are $O(\N^{\eps})$. We divide \eqref{eq_x112} by the function $s_h^{\circ}$ introduced in Definition~\ref{GQdef3pre}. The result is
\begin{equation}
 \label{eq_x124}
\begin{split}
& \quad \sqrt{(z-\alpha_{h})(z-\beta_{h})}\,\big\langle \bth_h\cdot \boldsymbol{W}_1(z)\big\rangle + \frac{q_h^{+,[1]}(z)}{s^{\circ}_h(z)} \\
& \quad + \frac{Q_h^+(z)}{2\N s^{\circ}_h(z)}\big\langle\bth_h\cdot \boldsymbol{W}_1(z)\big\rangle^2 +
 \frac{Q_h^+(z)}{2\N s^\circ_h(z)}\big\langle \bth_h^{\otimes 2} \cdot \boldsymbol{W}_2(z,z)\big\rangle
  \\ & = \frac{\widetilde R_h(z)}{s^{\circ}_h(z)}+ O(\N^{-\frac{1}{2}+\eps}).
\end{split}
\end{equation}
We implicitly used Condition 6. of Assumption~\ref{Assumptions_extra}, guaranteeing that the remainder $O(\N^{-\frac{1}{2}+\eps})$ remains of the same order after division by $s_h^{\circ}$. Proposition~\ref{Proposition_Nek_1_asymptotic_form} guarantees that $\frac{\widetilde R_h(z)}{s^{\circ}_h(z)}$ is a holomorphic function of $z\in\amsmathbb M_h$. Then we are in position to apply Theorem~\ref{Theorem_Master_equation} and find
\begin{equation}
\label{eq_x126}
\begin{split}
& \quad \boldsymbol{W}_{1}(z) \\
& = \boldsymbol{\Op}\Bigg[\bigg(\frac{q_h^{+,[1]}}{s^{\circ}_h} + \frac{Q_h^+}{2\N s^{\circ}_h}\big\langle\bth_h\cdot \boldsymbol{W}_1\rangle^2+
 \frac{Q_h^+}{2\N s^{\circ}_h}\big\langle \bth_h^{\otimes 2} \cdot \boldsymbol{W}_2(*,*)\big\rangle+ O(\N^{-\frac{1}{2}+\eps})\bigg)_{h=1}^H\,;\, \boldsymbol{0}\Bigg](z).
 \end{split}
\end{equation}
The second argument of the operator $\boldsymbol{\Op}$ is zero, because we are in the fixed filling fractions case, so that we have for any $h \in [H]$ as $z\rightarrow\infty$
\begin{equation*}
\begin{split}
\quad & G_h(z) = \frac{\N \hat n_h}{z}+O\bigg(\frac{1}{z^2}\bigg) \qquad \textnormal{and} \qquad \Gm_{\mu_h}(z) = \frac{\hat n_h}{z}+O\bigg(\frac{1}{z^2}\bigg) \\
\Longrightarrow  \quad & W_{1;h}(z)=\E\big[\Delta G_h(z)\big]= O\bigg(\frac{1}{z^2}\bigg).
\end{split}
\end{equation*}
The first argument of the operator $\boldsymbol{\Op}$ in \eqref{eq_x126} is $O(\N^{\eps})$ by Corollary~\ref{Corollary_a_priory_1}, hence, so is $\boldsymbol{W}_1(z)$ by the continuity properties of $\boldsymbol{\Op}$. This immediately gives a better bound for the second term in the first argument of $\boldsymbol{\Op}$, implying that it is $o(1)$ as $\N\rightarrow\infty$. The third term involves $\boldsymbol{W}_2$ and we are going to use the second-order Nekrasov equation \eqref{eq_Second_Nekr_asymptotic} of Proposition~\ref{Proposition_Nek_2_asymptotic_form} in order to bound it.

Applying Corollary~\ref{Corollary_a_priory_1}, we see that $\mathbbm{Err}^{(2)}_{h;h_2}(z,z_2)$ in Proposition~\ref{Proposition_Nek_2_asymptotic_form} is $O(\N^{1 - \frac{q}{2} + \eps})$ from the inequalities \eqref{eq_Err2_bound_two_2}, which is $O(\N^{\frac{1}{2} + \varepsilon})$ as $q \geq 1$. Hence, \eqref{eq_Second_Nekr_asymptotic} can be rewritten as
\begin{equation}
 q_{h}^{-}(z)\big\langle \bth_{h}\cdot \boldsymbol{W}_{2;\bullet,h_2}(z,z_2)\big\rangle + \delta_{h,h_2}\frac{q_{h}^{-}(z)}{2(z_2 - z)^2}
 \,\, = \, \, \widetilde R_{h;h_2}(z;z_2)+ O(\N^{\frac{1}{2}+\eps}).
\end{equation}
We have replaced $Q_{h}^{-}(z)$ by its leading term $q_{h}^{-}(z)$ to obtain the last formula: by \eqref{eq_Q_plus_determ_expansion} the error made in doing so is $O\big(\frac{1}{\N}\big)$ multiplied by $\langle \bth_{h}\cdot \boldsymbol{W}_{2;\bullet,h_2}(z,z_2)\rangle$ and, using Corollary~\ref{Corollary_a_priory_1}, we can absorb it into $O(\N^{\frac{1}{2}+\eps})$.

Dividing by $s^{\circ}_h(z)$ and using the operator $\Op_h$, this implies
\begin{equation}
\label{W2nund} \boldsymbol{W}_{2;h,h_2}(z,z_2)=\Op_h\Bigg[\frac{\sqrt{(*-\alpha_{h_2})(*-\beta_{h_2})}}{2(* - z_2)^2}\boldsymbol{e}^{(h_2)} + O(\N^{\frac{1}{2}+\eps})\,;\, \boldsymbol{0}\Bigg](z).
\end{equation}
In the statement of Theorem~\ref{Theorem_Master_equation} we required the function in the $h$-th component of the first argument of the operator $\boldsymbol{\Op}$ to be holomorphic in $\amsmathbb M_h\setminus \amsmathbb{A}_h^{\mathfrak{m}}$. Hence, the formula \eqref{W2nund} is \textit{a priori} valid only if we take $z_2$ to be outside $\bigcup_{h=1}^H \amsmathbb M_h$. This is not restrictive, since one could replace $\amsmathbb{M}_h$ with small neighborhoods of $[\hat{a}_h,\hat{b}_h]$, so that the functions in $O(\cdot)$ are still controlled uniformly on any $\N$-independent compacts of the complement in the complex plane of the defining segments of the ensemble. By the continuity properties of $\boldsymbol{\Op}$, it follows that $\boldsymbol{W}_2=O(\N^{\frac{1}{2}+\eps})$, justifying
\begin{equation}
\label{eq_x165}
\forall h,h_2 \in [H] \qquad \sup_{(z_1,z_2) \in \amsmathbb{K}_{h} \times \amsmathbb{K}_{h_2}} \big|W_{2;h_1,h_2}(z_1,z_2)\big| \leq C \N^{\frac{1}{2}+\eps}.
\end{equation}
with notations from Proposition~\ref{Proposition_uniform_moments_bound}.

Therefore, the third term in the first argument of $\boldsymbol{\Op}$ in \eqref{eq_x126} is $\frac{1}{\N} \cdot O(\N^{\frac{1}{2}+\eps})=O(\N^{-\frac{1}{2}+\eps})$ and becomes $o(1)$ as $\N\rightarrow\infty$. Hence,
\begin{equation}
\label{eq_x166}
\boldsymbol{W}_{1}(z)= \boldsymbol{\Op}\Bigg[\bigg(\frac{q_h^{+,[1]}}{s^{\circ}_h} + O(\N^{-\frac{1}{2}+\eps})\bigg)_{h=1}^H\,;\, \boldsymbol{0}\Bigg](z),
\end{equation}
In particular, $\E\big[\Delta \boldsymbol{G}(z)\big] = \boldsymbol{W}_{1}(z)$ is bounded as $\N\rightarrow\infty$.

\medskip

\noindent \textsc{Step 2.} We now repeat the same approach for the higher moments: we obtain some bounds through the higher-order Nekrasov equations, and then plug into the lower-order Nekrasov equations to get even stronger bounds. We will no longer need the first-order Nekrasov equation of Proposition~\ref{Proposition_Nek_1_asymptotic_form}. We rewrite the second and higher-order Nekrasov equations of Propositions~\ref{Proposition_Nek_2_asymptotic_form} and \ref{Proposition_Nek_higher_asymptotic_form} in a uniform way as $\N\rightarrow\infty$:
\begin{equation}
 \label{eq_Higher_Nekr_asymptotic_uniform}
 \begin{split}
 & \quad q_{h}^{-}(z)\big\langle \bth_{h}\cdot \boldsymbol{W}_{n;\bullet,h_2,h_3,\ldots,h_n}(z,z_2,z_3,\ldots,z_n)\big\rangle \\
 & =  \widetilde R_{h;h_2,\ldots,h_n}(z;z_2,\ldots,z_n)+ \mathbbm{Err}^{(3)}_{h;h_2,\ldots,h_n}(z;z_2,\ldots,z_n) +O(1).
 \end{split}
\end{equation}
The function $\mathbbm{Err}^{(3)}_{h;h_2,\ldots,h_n}(z;z_2,\ldots,z_n)$ contains the contributions of the error term in $Q_h^-(z)=q_h^-(z)+O(\frac{1}{\N})$, as well as $\mathbbm{Err}^{(2)}_{h;h_2,\ldots,h_n}(z;z_2,\ldots,z_n)$. The latter is a finite sum of terms in the form \eqref{eq_Err2_bound_multi} satisfying the inequality \eqref{eq_Err2_bound_multi_2}. Equation~\ref{eq_Higher_Nekr_asymptotic_uniform} is then a corollary of Proposition~\ref{Proposition_Nek_2_asymptotic_form} for $n=2$ and Proposition~\ref{Proposition_Nek_higher_asymptotic_form} for $n\geq 3$. Note that if we assume that the points $z_j, z_j\pm \tfrac{1}{2\N}$ are kept away from $\amsmathbb M_h$ for $j \in \llbracket 2,n\rrbracket$, then the $O(1)$ term in \eqref{eq_Higher_Nekr_asymptotic_uniform} is uniform in $z$ belonging to compact subset of $\amsmathbb M_h\setminus\amsmathbb{A}_{h}^{\mathfrak{m}}$.

Using the operator $\boldsymbol{\Op}$ acting in variable $z$, \eqref{eq_Higher_Nekr_asymptotic_uniform} implies that for $n\geq 2$.
\begin{equation}
\label{eq_x167}
 \boldsymbol{W}_{n;\bullet,h_2,h_3,\ldots,h_n}(z,z_2,z_3,\ldots,z_n)=\boldsymbol{\Op}\Bigg[\bigg(O(1)- \frac{\mathbbm{Err}^{(3)}_{h;h_2,\ldots,h_n}(*;z_2,\ldots,z_n)}{s^\circ_h}\bigg)_{h=1}^H\,;\,\boldsymbol{0}\Bigg](z).
\end{equation}
As we demonstrate in the next step, the desired bound \eqref{eq_x168} follows from \eqref{eq_x167}.

\medskip

\noindent \textsc{Step 3.} Our next goal is to prove the following statement by induction in $t$.

\medskip

\noindent \textsc{Claim.} For any $n \in \amsmathbb{Z}_{\geq 2}$, $t \in \amsmathbb{Z}_{\geq 0}$ and $\eps>0$, any compact sets $\amsmathbb K_h\subset \amsmathbb C\setminus\amsmathbb{A}_h^{\mathfrak{m}}$ indexed by $h \in [H]$, and for any $h_1,\ldots,h_n \in [H]$ we have as $\N \rightarrow \infty$
\begin{equation}
\label{eq_inductive_bound}
 \sup_{\substack{z_i\in \amsmathbb{K}_{h_i} \\ i \in [n]}} \E \left[ \prod_{i=1}^n \Delta G_{h_i}(z_i)\right] = O\big(\N^{\frac{n - t}{2} +\eps}\big)+O(1).
\end{equation}
Note that in the right-hand side, the first term is larger if $t\leq n$ while the second one is larger if $t>n$.
\smallskip

The base case $t=0$ is given in Corollary~\ref{Corollary_a_priory_1}. For the induction step, we assume that \eqref{eq_inductive_bound} is valid for some particular value of $t$ and aim to prove it for $t+1$.

By using Cauchy integral formula for the derivatives (\textit{cf.} \eqref{derCauchy}) the bound \eqref{eq_inductive_bound} implies a similar bound which some of $\Delta G_{h_i}(z_i)$ factors are replaced with their derivatives of arbitrary fixed order. Hence, under \eqref{eq_inductive_bound}, the magnitude of each term in $\mathbbm{Err}^{(2)}_{h;h_2,\ldots,h_n}$ in \eqref{eq_x167} of the form
 \begin{equation}
   \frac{F(z)}{\N^{q}} \, \prod_{i=1}^k \E \left[ \prod_{m=2}^n \big(\Delta G_{h_m}(z_m)\big)^{\tilde{l}_{i,m}} \prod_{j =1}^{l_i} \partial_z^{d_{i,j}} \Delta G_{g_{i,j}}(z) \right]
 \end{equation}
 is bounded as
\[
 O\big( \N^{\frac{p - t}{2}-q+\eps} \big)+O(1), \qquad \textnormal{for} \quad p = \sum_{i=1}^k \bigg(l_i+\sum_{m=2}^{n} \tilde{l}_{i,m}\bigg).
\]
Since $q+1\leq p\leq q+n$ by \eqref{eq_Err2_bound_multi_2} and $q\geq 1$, the last expression is
\[
 O\big( \N^{\frac{ n - t - q}{2}+\eps} \big)+O(1)= O\big( \N^{\frac{n - t - 1}{2}+\eps} \big)+O(1).
\]
Hence, using the continuity properties of the operator $\Op_h$ and \eqref{eq_x167}, we conclude
\begin{equation}
\label{eq_x169}
 \boldsymbol{W}_{n}(z_1,z_2,z_3,\ldots,z_n)= O\big( \N^{\frac{n - t - 1}{2}+\eps} \big)+O(1).
\end{equation}
Because centered moments are (homogeneous) polynomials in cumulants, \eqref{eq_x169} then implies
\begin{equation}
 \E \left[\prod_{i=1}^n \Delta G_{h_i }(z_i)\right] =O\big( \N^{\frac{n - t - 1}{2}+\eps}\big)+O(1),
\end{equation}
which is the induction step.

\medskip

\noindent \textsc{Step 4.} We use \eqref{eq_inductive_bound} for an even positive integer $n$, $t=n+1$, $h_1=h_2=\cdots=h_n=h$ and
\[
z_1=z_2=\cdots=z_{\frac{n}{2}}=z \qquad \textnormal{and}\qquad z_{\frac{n}{2}+1}=z_{\frac{n}{2}+2}=\cdots=z_n=z^*,
\]
 and use $\Delta G_h(z^*)= (\Delta G_h(z))^*$ to conclude
\begin{equation}
\label{eq_x173}
 \E\big[|\Delta G_{h }(z)|^n\big] = O(1).
\end{equation}
For odd $n$ we can use $|\Delta G_h(z)|^n\leq |\Delta G_h(z)|^{n+1}+1$ to reduce to even $n$ case and conclude that \eqref{eq_x170} holds for any $n \in \amsmathbb{Z}_{> 0}$. Finally, we observe
\[
 \E\Bigg[ \bigg| \prod_{i=1}^n \Delta G_{h_i}(z_i) \bigg|\Bigg] \leq \sum_{i=1}^n \E \big[ |\Delta G_{h_i}(z_i)|^n\big] .
\]
Together with \eqref{eq_x173} this implies \eqref{eq_x168}.
\end{proof}

\section{Proof of the asymptotic expansion (Theorem~\ref{Theorem_correlators_expansion})}
\label{secprosofggg}

We are ready to prove the asymptotic expansions \eqref{eq_x30}--\eqref{eq_x34}. The bounds of \eqref{eq_x31} and \eqref{eq_x33} follow from the definitions and properties of the solution operator $\boldsymbol{\Op}$ and we omit details here.

\smallskip

\noindent \textsc{Step 1.} Let us prove the bound for $\boldsymbol{W}_n$ with $n\geq 3$, given in \eqref{eq_x34}. Starting from \eqref{eq_Higher_Nekr_asymptotic}, using Proposition~\ref{Proposition_uniform_moments_bound}, $q\geq 1$ in \eqref{eq_Err2_bound_multi}, and recalling the asymptotic expansion for $Q_h^-(z)$ of \eqref{eq_Q_plus_determ_expansion}, we have
\begin{equation}
 q_{h}^{-}(z)\big\langle \bth_{h}\cdot \boldsymbol{W}_{n;\bullet,h_2,h_3,\ldots,h_n}(z,z_2,z_3,\ldots,z_n)\big\rangle+O(\N^{-\frac{1}{2}+\eps}) =  \tilde R_{h;h_2,\ldots,h_n}(z;z_2,\ldots,z_n).
\end{equation}
Dividing by $s_h^\circ(z)$ and using the solution operator $\boldsymbol{\Op}$ on the resulting function of the variable $z$, this implies
\[
\boldsymbol{W}_{n;\bullet,h_2,h_3,\ldots,h_n}(z,z_2,z_3,\ldots,z_n)=\boldsymbol{\Op}\big[O(\N^{-\frac{1}{2}+\eps})\,;\,\boldsymbol{0}\big](z).
\]
Using the continuity properties of $\boldsymbol{\Op}$ proven in Theorem~\ref{Theorem_Master_equation}, we obtain \eqref{eq_x34}.

\medskip

\noindent \textsc{Step 2.} We proceed to asymptotic expansion for $\boldsymbol{W}_2$. Starting from \eqref{eq_Second_Nekr_asymptotic}, using Proposition~\ref{Proposition_uniform_moments_bound}, observing as in Step 1 that $\mathbbm{Err}^{(2)}_{h;h_2}(z;z_2)$ is $O\big(\frac{1}{\N}\big)$, and replacing $Q_h^-$ by its leading order, we obtain
\[
 q_{h}^{-}(z)\big\langle \bth_{h}\cdot \boldsymbol{W}_{2;\bullet,h_2}(z,z_2)\big\rangle + \delta_{h,h_2}\frac{q_{h}^{-}(z)}{2(z_2 - z)^2}+ O(\N^{-\frac{1}{2}+\eps})= \widetilde R_{h;h_2}(z;z_2).
\]
Dividing by $s_h^\circ(z)$ and using the solution operator $\boldsymbol{\Op}$ operator in the resulting function of the variable $z$, we conclude
\[
 \boldsymbol{W}_{2;\bullet,h_2}(z,z_2)=\boldsymbol{\Op}\Bigg[\frac{\sqrt{(*-\alpha_{h_2})(*-\beta_{h_2})}}{2(* - z_2)^2}\boldsymbol{e}^{(h_2)}\,;\, \boldsymbol{0}\Bigg](z) + O(\N^{-\frac{1}{2}+\eps}),
\]
which matches \eqref{eq_x32}, \eqref{eq:8585} and identifies the leading order of $\boldsymbol{W}_{2}(z_1,z_2)$ with the fundamental solution $\boldsymbol{\mathcal{F}}(z_1,z_2)$ of \eqref{eq_covariancepre}.

\medskip

\noindent \textsc{Step 3.} Finally, let us prove the asymptotic expansion for $\boldsymbol{W}_1$ given in \eqref{eq_x30}. For this purpose we produce a higher-order version of \eqref{eq_x126}. Dividing \eqref{eq_First_Nekr_asymptotic} by $s^\circ_h(z)$ and recalling the asymptotic expansions \eqref{eq_Q_plus_determ_expansion}, we get
\begin{equation}
\label{eq_First_Nekr_asymptotic_final}
\begin{split}
& \quad \sqrt{(z-\alpha_h)(z-\beta_h)}\, \big\langle \bth_h\cdot \boldsymbol{W}_1(z)\big\rangle + \frac{q_h^{+,[1]}(z)}{s^\circ_h(z)} + \frac{q_h^{-,[1]}(z)}{\N s^\circ_h(z)}\big\langle \bth_h\cdot \boldsymbol{W}_1(z)\rangle + \frac{q_h^+(z)}{2\N s^\circ_h(z)}\big\langle\bth_h\cdot \boldsymbol{W}_1(z)\rangle^2 \\
 & \quad +
 \frac{q_h^+(z)}{2\N s^\circ_h(z)}\big\langle \bth_h^{\otimes 2} \cdot \boldsymbol{W}_2(z,z)\big\rangle + \frac{q_h^+(z)}{2\N s^\circ_h(z)}\big\langle \bth_h(\bth_h - 1)\cdot \partial_{z} \boldsymbol{W}_1(z)\big\rangle
 + \frac{q_h^{+,[2]}(z)}{\N s^\circ_h(z)} \\
 & = \frac{\widetilde R_h(z)}{s^\circ_h(z)}+ \frac{\mathbbm{Err}^{(1)}_h(z)}{s^\circ_h(z)}+ \frac{\mathbbm{Err}^{(3)}_h(z)}{s^\circ_h(z)},
\end{split}
\end{equation}
where $\mathbbm{Err}^{(2)}$ and further terms of the asymptotic expansions of $Q^{\pm}_h(z)$ have been absorbed into $\mathbbm{Err}^{(3)}(z)$. Since $q\geq 2$ in \eqref{FERR2}, Proposition~\ref{Proposition_uniform_moments_bound} implies that
\[
\frac{\mathbbm{Err}^{(3)}_h(z)}{s^\circ_h(z)}=O\bigg(\frac{1}{\N^{2}}\bigg).
\]
Hence, using Proposition~\ref{Proposition_uniform_moments_bound} again to bound all other terms, \eqref{eq_First_Nekr_asymptotic_final} implies
\begin{equation}
 \sqrt{(z-\alpha_h)(z-\beta_h)}\,\big\langle \bth_h\cdot \boldsymbol{W}_1(z)\big\rangle  + \frac{q_h^{+,[1]}(z)}{s^\circ_h(z)} + O\bigg(\frac{1}{\N}\bigg) = \frac{\widetilde R_h(z)}{s^\circ_h(z)}.
\end{equation}
Applying Theorem~\ref{Theorem_Master_equation}, it follows that
\begin{equation}
 \boldsymbol{W}_{1}(z)=\boldsymbol{\Op}\Bigg[ \bigg(\frac{q_h^{+,[1]}}{s^\circ_h}\bigg)_{h=1}^H\,;\,\boldsymbol{0} \Bigg](z)+ O\bigg(\frac{1}{\N}\bigg).
\end{equation}
This gives the first term in the asymptotic expansion of $\boldsymbol{W}_{1}(z)$ and matches the definition of $\boldsymbol{W}_{1}^{[1]}$ given in \eqref{eq:8585}. Plugging back into \eqref{eq_First_Nekr_asymptotic_final} and using the expansion for $\boldsymbol{W}_2(z,z)$ established on Step 2, we improve \eqref{eq_First_Nekr_asymptotic_final} to
\begin{equation}
\label{eq_First_Nekr_asymptotic_final_solution_step}
\begin{split}
& \sqrt{(z-\alpha_h)(z-\beta_h)}\,\big\langle \bth_h\cdot \boldsymbol{W}_1(z)\big\rangle + \frac{q_h^{+,[1]}(z)}{s^\circ_h(z)} +\frac{q_h^{-,[1]}(z)}{\N s^\circ_h(z)}\big\langle \bth_h\cdot \boldsymbol{W}_1^{[1]}(z)\big\rangle + \frac{q_h^+(z)}{2\N s^\circ_h(z)}\big\langle\bth_h\cdot \boldsymbol{W}_1^{[1]}(z)\big\rangle^2 \\ & +
 \frac{q_h^+(z)}{2\N s^\circ_h(z)}\big\langle \bth_h^{\otimes 2} \cdot \boldsymbol{\mathcal{F}}(z,z)\big\rangle + \frac{q_h^+(z)}{2\N s^\circ_h(z)}\big\langle \bth_h(\bth_h - 1)\cdot \partial_{z} \boldsymbol{W}_1^{[1]}(z)\big\rangle + \frac{q_h^{+,[2]}(z)}{\N s^\circ_h(z)} + O(\N^{-\frac{1}{2}+\eps}) \\
 & =  \frac{\widetilde R_h(z)}{s^\circ_h(z)}.
\end{split}
\end{equation}
Applying the solution operator $\boldsymbol{\Op}$ for this function of the variable $z$ variable again, we arrive at \eqref{eq_x30} with $\boldsymbol{W}^{[1]}_{1}(z)$ and $\boldsymbol{W}^{[2]}_{1}(z)$ given by \eqref{eq:8585}. This concludes the proof of Theorem~\ref{Theorem_correlators_expansion}.

\section{Application to the discrete ensemble with Gaussian weight}
\label{Gausscorr}

\subsection{Main results}

Let us illustrate the main results of this chapter by obtaining explicitly the expansion of correlators for the discrete ensemble with Gaussian weight. This is the discrete ensemble with $H = 1$, $a_1 = - \infty$, $b_1 = +\infty$, intensity of repulsion $\theta > 0$ and weight
\begin{equation}
\label{Gweig}
w(x) = e^{-\kappa \frac{x^2}{N}},
\end{equation}
where $\kappa > 0$. We set $\N = N$ and consider the $N \rightarrow \infty$ regime with $\theta$ and $\kappa$ independent of $N$. The equilibrium measure for this ensemble was computed in Proposition~\ref{Prop_Gaussian_LLN}. It is a semi-circle law supported on $[-\beta,\beta]$ in the weak confinement phase $\sqrt{2\theta\kappa} < \pi$, or a measure supported on a segment $[-\beta,\beta]$ with two bands $(-\beta,-\alpha)$ and $(\alpha,\beta)$ separated by a saturation $[-\alpha,\alpha]$ in the strong confinement phase $\sqrt{2\theta\kappa} > \pi$. The ensemble always satisfies Assumptions~\ref{Assumptions_Theta}, \ref{Assumptions_basic} and \ref{Assumptions_analyticity}, and Lemma~\ref{Lemma_Discr_Gauss_Off_critical} shows that it satisfies as well the off-criticality Assumption~\ref{Assumptions_offcrit} for $\sqrt{2\theta \kappa} \neq \pi$.

We first show how to derive from Theorem~\ref{Theorem_correlators_expansion} and Corollary~\ref{Corollary_CLT} the following form of the central limit theorem in the weak confinement phase.
\begin{proposition}
\label{CLTGAUSSIANWEIGHT}
Assume $\sqrt{2\theta \kappa} < \pi$.  Let $f$ be an analytic test function and denote
\begin{equation}
\label{sigmasqbeta}
\sigma(z) = \sqrt{z^2 - \beta^2}\qquad \textnormal{with}\quad \beta = \sqrt{\frac{2\theta}{\kappa}}.
\end{equation}
The random variable
\[
\sum_{i = 1}^{N} f\bigg(\frac{\ell_i}{N}\bigg) - N\int_{-\beta}^{\beta} f(x)\dd\mu(x)
\]
is approximated in the sense of moments as $N \rightarrow \infty$ by a Gaussian random variable, with mean
\begin{equation}
\label{Meanfmean}
\textnormal{\textsf{Mean}}[f] = \frac{\theta -1}{2\theta}\left(\frac{f(-\beta) + f(\beta)}{2} - \frac{1}{\pi} \int_{-\beta}^{\beta} \frac{f(x) \dd x}{|\sigma(x)|}  - \frac{\kappa}{\pi^2}\cdot \textnormal{p.v.}\int_{-\beta}^{\beta} \int_{-\beta}^{\beta} \frac{f(x) \cdot y \cdot \textnormal{T}\big[\kappa|\sigma(y)|\big]}{(x - y)|\sigma(x)|} \dd x \dd y \right)\!,
\end{equation}
where $\textnormal{T}(u) = u\big(u\frac{\cos(u)}{\sin(u)} - 1\big)$ and p.v. is the Cauchy principal value, and with covariance
\begin{equation}
\begin{split}
\textnormal{\textsf{Cov}}[f,f] & = \oint_{\gamma} \oint_{\gamma} \frac{\dd z_1\dd z_2}{(2\ii\pi)^2} \cdot \frac{f(z_1)f(z_2)}{2\theta (z_1 - z_2)^2}\left(-1 + \frac{z_1z_2 - \beta^2}{\sigma(z_1)\sigma(z_2)}\right) \\
& = \frac{1}{4\pi^2} \int_{-\beta}^{\beta}\int_{-\beta}^{\beta} f'(x_1)f'(x_2) \log\left(\frac{\beta^2 - x_1x_2 - \sqrt{(\beta^2 - x_1^2)(\beta^2 - x_2^2)}}{\beta^2 - x_1x_2  + \sqrt{(\beta^2 - x_1^2)(\beta^2 - x_2^2)}}\right) \dd x_1 \dd x_2.
\end{split}
\end{equation}
\end{proposition}
The first formula for the covariance already appears in \cite{BGG}. The insightful reader may recognize in the integrand of the last formula the Green function of the upper half-plane $\amsmathbb{H}$ with Dirichlet boundary condition after applying the uniformization map
\[
z \mapsto (z + \sigma(z)) \in \amsmathbb{H}.
\]
This is a more general phenomenon that was mentioned in Chapter~\ref{SIntro}, will be elucidated in Chapter~\ref{Chapter_AG} and applied to random tilings in Chapter~\ref{Chap11}. We will see in Section~\ref{sec:comparcd} that the formulae of Proposition~\ref{CLTGAUSSIANWEIGHT} are actually the same as in the continuous $\sbeta$-ensembles, except for the term in the mean involving $\textnormal{T}$ and which is a novelty of the discrete ensemble (only present for $\theta \neq 1$).

In the strong confinement phase, we have two bands separated by a saturation, so we cannot yet apply Theorem~\ref{Theorem_correlators_expansion} and \ref{Corollary_CLT} and will need to wait for Chapter~\ref{Chapter_filling_fractions} to get (the form of) the correct answer. However, we can handle with those theorems the ensemble delivered by Proposition~\ref{proposition_Gaussian_conditioning}: this is an ensemble where we impose fully packed particles in a fixed bit of the saturation and condition on fixed value of  filling fractions (at least not too far from their equilibrium value) on the two segments on both sides of the removed one. The equilibrium measure continues to have two bands $(\alpha_1,\beta_1)$ and $(\alpha_2,\beta_2)$, but now we can have $\alpha_1 \neq -\beta_2$ and  $\beta_1 \neq -\alpha_2$ because the choice of filling fractions can break the symmetry with respect to the origin.
\begin{proposition}
\label{CLTGAUSSIANWEIGHT2} Assume $\sqrt{2\theta \kappa} > \pi$ and consider the ensemble delivered in Proposition~\ref{proposition_Gaussian_conditioning} by localizing and conditioning on filling fractions not far from their equilibrium value. Introduce
\[
\sigma(z) = \sqrt{(z - \alpha_1)(z -\beta_1)(z - \alpha_2)(z - \beta_2)}
\]
depending on the endpoints of its two bands. The random variable
\[
\sum_{i = 1}^{N} f\bigg(\frac{\ell_i}{N}\bigg) - N\int f(x)\dd\mu(x)
\]
is approximated in the sense of moments by a Gaussian random variable. The mean is zero if $\theta = 1$. For general $\theta > 0$ the covariance is
\[
\textnormal{\textsf{Cov}}[f,f] = \oint_{\gamma} \oint_{\gamma} \frac{\dd z_1 \dd z_2}{(2\ii\pi)^2}\,\frac{f(z_1)f(z_2)}{4\theta \sigma(z_1)\sigma(z_2)}\left[\left(\frac{\sigma(z_1) - \sigma(z_2)}{z_1 - z_2}\right)^2 - (z_1 + z_2)(z_1 + z_2 - \mathsf{e}_1) + \mathsf{c}\right],
\]
where
\begin{equation*}
\begin{split}
\mathsf{e}_1 & = \alpha_1 + \beta_1 + \alpha_2 + \beta_2, \\
\mathsf{c} & = -(\alpha_1 + \beta_2)(\beta_1 + \alpha_2) + (\beta_2 - \beta_1)(\alpha_2 - \alpha_1)\frac{E(\mathsf{k})}{K(\mathsf{k})}, \\
\mathsf{k} & = \sqrt{\frac{(\beta_2 - \alpha_1)(\alpha_2 - \beta_1)}{(\beta_2 - \beta_1)(\alpha_2 - \alpha_1)}}.
\end{split}
\end{equation*}
\end{proposition}
A formula for the mean for $\theta \neq 1$ can in principle be extracted from \eqref{expsfguzbguwg} and the details of the localization procedure to compute $\mu^{Y}$ taking into account Sections~\ref{Section_alternative_localization}-\ref{Sec_Mismatch}, but we will not address this. The next Section is devoted to the proof of these formulae and some more, and we compare the outcome to the continuous $\sbeta$-ensembles in Section~\ref{sec:comparcd}.

\subsection{Expansion of correlators}
\label{SexpcorGAUSS}

We start by specifying all the auxiliary functions introduced in Section~\ref{Subsection_Auxiliary_functions}. For the Gaussian weight \eqref{Gweig} the definitions in Section~\ref{Section_list_of_assumptions} lead to a smooth potential
\[
U(x) = V(x) = \kappa x^2,
\]
and there is no error function when we compare potential and weights: $\mathbbm{e}(x) = 0$. The remaining functions appearing in the factorization of the weights in Definition~\ref{Definition_phi_functions} are
\begin{equation*}
\begin{split}
\Phi^+(z) & = e^{-\N[ U(z + \frac{1}{2\N}) - U(z - \frac{1}{2\N})]} = e^{-2\kappa z} = e^{-\partial_z V(z)} = \phi^+(z), \\
\Phi^-(z) & = \phi^-(z) = 1,
\end{split}
\end{equation*}
and we have $\iota^- = \iota^+ = 0$.

The ensemble never satisfies the extra Assumption~\ref{Assumptions_extra}. The segment is infinite but the large deviation bounds of Theorem~\ref{Theorem_ldpsup} imply that conditioning and localizing the ensemble on some finite segment $[-\hat{D}\N,\hat{D}\N]$ such that $\hat{D} > \frac{2\theta}{\kappa}$ only changes the partition function and the correlators (after integration against polynomial test functions) by exponentially small corrections as $\N \rightarrow \infty$. One does not need much of the sophistication of Chapter~\ref{Chapter_conditioning} for this, except to choose segments fulfilling integrality conditions. On this finite segment the ensemble satisfies the additional Assumption~\ref{Assumptions_extra} if $\sqrt{2\theta \kappa} < \pi$. If $\sqrt{2\theta \kappa} > \pi$ it still does not because we have two bands in a single segment. But, localizing and conditioning to fixed filling fractions with Proposition~\ref{proposition_Gaussian_conditioning}, we get an ensemble which does satisfy Assumptions~\ref{Assumptions_Theta}, \ref{Assumptions_basic}, \ref{Assumptions_offcrit}, \ref{Assumptions_analyticity} and \ref{Assumptions_extra}.

\medskip

Let us continue for $\sqrt{2\theta\kappa} < \pi$. Since the equilibrium measure is a semi-circle supported on $(-\beta,\beta)$, we have
\begin{equation}
\label{GGZG}\mathcal{G}(z) = \frac{\kappa(z - \sigma(z))}{\theta},\qquad \sigma(z) = \sqrt{z^2 - \beta^2}.
\end{equation}
The auxiliary functions of Section~\ref{Subsection_Auxiliary_functions} are $q^{\pm}(z) = e^{-2\kappa z + \theta \mathcal{G}(z)} \pm e^{-\theta \mathcal{G}(z)}$. Together with \eqref{GGZG} this yields
\[
q^+(z) = 2e^{-\kappa z}\,\mathrm{cosh}[\kappa \sigma(z)],\qquad q^-(z) = -2e^{-\kappa z} \mathrm{sinh}[\kappa \sigma(z)],\qquad s^{\circ}(z) = \frac{q^-(z)}{\sigma(z)}.
\]
Then, \eqref{qhsmallexp} yields
\begin{equation*}
\begin{split}
q^{\pm,[1]}(z) & = \frac{\theta(\theta - 1)}{2} q^{\pm}(z)  \cdot \partial_z \mathcal{G}(z), \\
q^{\pm,[2]}(z) & = \frac{\theta(4\theta^2 - 6\theta + 3)}{24} q^{\mp}(z) \cdot \partial_z^2 \mathcal{G}(z) + \frac{\theta^2(\theta - 1)^2}{8} q^{\pm}(z) \cdot \big(\partial_z \mathcal{G}(z)\big)^2.
\end{split}
\end{equation*}
Theorem~\ref{Theorem_correlators_expansion} gives us the expansion for the correlators
\begin{equation}
\label{firstexpanug}
\begin{split}
W_1(z) & = \N \mathcal{G}(z) + W_1^{[1]}(z) + \N^{-1} W_1^{[2]}(z) + O(\N^{-2}), \\
W_2(z_1,z_2) & = \mathcal{F}(z_1,z_2) + O(\N^{-1}).
\end{split}
\end{equation}
Since $H = 1$, the leading order for $W_2$ is given by Proposition~\ref{W20diagonal}:
\begin{equation}
\label{Fz1z2}
\mathcal{F}(z_1,z_2) = \frac{1}{2\theta(z_1 - z_2)^2}\left(-1 + \frac{z_1z_2 - \frac{2\theta}{\kappa}}{\sigma(z_1)\sigma(z_2)}\right).
\end{equation}
For the coefficients in $W_1$ as given by Theorem~\ref{Theorem_correlators_expansion}, we need the solution operation $\Upsilon$. It is particularly simple as $H = 1$ (\textit{cf.} Theorem~\ref{Theorem_Masterspecial}), namely
\[
\Upsilon[E;0](z) = \frac{1}{\theta \sigma(z)} \oint_{\gamma} \frac{\dd \zeta}{2\ii\pi}\,\frac{E(\zeta)}{\zeta - z},
\]
where $\gamma$ is a counterclockwise contour around the segment and $z$ is outside the contour. Then
\begin{equation*}
\begin{split}
W_1^{[1]}(z) & = \frac{1}{\theta \sigma(z)} \oint_{\gamma} \frac{\dd \zeta}{2\ii\pi} \frac{q^{+,[1]}(\zeta)}{(\zeta - z)s^{\circ}(\zeta)}, \\
W_1^{[2]}(z) & = \frac{1}{\sigma(z)} \oint_{\gamma} \frac{\big(q^{+,[2]}(\zeta) + \frac{q^+(\zeta)}{2} \mathcal{F}(\zeta,\zeta)\big)\dd\zeta}{2\ii\pi(\zeta - z)s^{\circ}(\zeta)} + O(\theta - 1),
\end{split}
\end{equation*}
where $O(\theta - 1)$ is a term vanishing for $\theta = 1$.  In particular, the first non-zero subleading order in $W_1$ is $\N^0 W_1^{[1]}$ if $\theta \neq 1$ and $\N^{-1}W_1^{[2]}$ if $\theta = 1$. Let us compute them more explicitly. We have
\begin{equation}
\label{W11zgauss}
\begin{split}
W_1^{[1]}(z) = \frac{\theta - 1}{2\theta\sigma(z)} \oint_{\gamma} \frac{\dd \zeta}{2\ii\pi}\,\frac{1 - \sigma'(\zeta)}{\zeta - z}\,\kappa\sigma(\zeta)\,\mathrm{coth}[\kappa \sigma(\zeta)].
\end{split}
\end{equation}
Observing that $\sigma'(z) = \frac{z}{\sigma(z)}$, we decompose
\begin{equation}
\label{ratdecompostri}
(1 - \sigma'(\zeta))\kappa\sigma(\zeta)\,\mathrm{coth}[\kappa \sigma(\zeta)] = (1 - \sigma'(\zeta)) + (\kappa\sigma(\zeta))^2 \textnormal{S}[\kappa \sigma(\zeta)] - \kappa^2 \zeta \sigma(\zeta) \textnormal{S}[\kappa \sigma(\zeta)],
\end{equation}
where $\textnormal{S}(t) = t\,\mathrm{coth}(t) - 1$ is even holomorphic function of $t$ in the domain $|t| < \pi$. An alternative formula for this function is the series expansion
\begin{equation}
\label{Lexpseries}
\textnormal{S}(t) = \sum_{n = 1}^{\infty} \frac{2^{2n}\textnormal{B}_{2n}}{(2n)!}\,t^{2n},
\end{equation}
\label{index:Bernou}where $\textnormal{B}_{n}$ are the Bernoulli numbers; it is adapted for computations when $\kappa$ is small. The maximum of $x \mapsto \kappa|\sigma(x)|$ for $|x| \leq \frac{2\theta}{\kappa}$ is equal to $\sqrt{2\theta \kappa}$, which is strictly smaller than $\pi$ because we study the one-band regime. Therefore, the second term in \eqref{ratdecompostri} is holomorphic inside the contour $\gamma$ and thus does not contribute to the contour integral in \eqref{W11zgauss}. The first term in \eqref{ratdecompostri} has no pole outside the contour and behaves as $O(\frac{1}{\zeta^2})$ as $\zeta \rightarrow \infty$. Then, in the contour integral \eqref{W11zgauss} we can move the contour to $\infty$ and this only picks the simple pole at $\zeta = z$, and we find
\begin{equation}
\label{W11zzz1}
W_1^{[1]}(z) = \frac{\theta - 1}{2\theta}\left(\frac{\sigma'(z) - 1}{\sigma(z)} + \oint_{\gamma} \frac{\dd \zeta}{2\ii\pi}\,\frac{\sigma(\zeta)}{\sigma(z)}\,\frac{\kappa^2 \zeta \,\textnormal{S}[\kappa \sigma(\zeta)]}{z - \zeta}\right).
\end{equation}
The last integral could be evaluated as series using the series expansion \eqref{Lexpseries} for $\textnormal{S}$ and the series expansion for the functions $(\sigma(\zeta))^{2n + 1}$, but we will not follow this route.

If $f$ is an analytic function, Corollary~\ref{Corollary_CLT} gives us Proposition~\ref{CLTGAUSSIANWEIGHT} with mean
\[
\textnormal{\textsf{Mean}}[f] = \oint_{\gamma} \frac{\dd z}{2\ii\pi}W_1^{[1]}(z) f(z),
\]
where we should insert \eqref{W11zzz1}. The expression $\sigma'(z)/\sigma(z)$ has simple pole at $z = \pm \beta$ with residue $\frac{1}{2}$, hence contributes to a term $\frac{f(-\beta) + f(\beta)}{2}$. We can rewrite the contribution of the last two terms of $W_1^{[1]}$ by squeezing the contour integral to the band. For the contribution of $-1/\sigma(z)$ this is easy and we get the second term of \eqref{Meanfmean} noticing that $\sigma(x^{\pm}) = \pm i |\sigma(x)|$. For the contribution of the third term in \eqref{W11zzz1}, we first squeeze the $\zeta$-contour integral to the band remembering that $z$ is outside the contour and that $\textnormal{S}$ is an even function:
\[
\oint_{\gamma} \frac{\dd \zeta}{2\ii\pi} \frac{\sigma(\zeta)}{\sigma(z)}\,\frac{\kappa^2 \zeta\, \textnormal{S}[\kappa \sigma(\zeta)]}{z - \zeta} = - \int_{-\beta}^{\beta} \frac{\dd y}{\pi}\,\frac{|\sigma(y)|}{\sigma(z)}\,\frac{\kappa^2 y \,\textnormal{S}[\ii \kappa |\sigma(y)|]}{z - y}.
\]
We then write $\textnormal{S}[\textnormal{i}t] = t^{-1}\textnormal{T}[t]$ with $\textnormal{T}(t) = t\big(t\cdot \mathrm{cot}(t) - 1)$. Integrating this against $\frac{f(z)}{2\ii\pi}$ and squeezing the contour to the band leads to the last term in \eqref{Meanfmean}. This finishes the proof of Proposition~\ref{CLTGAUSSIANWEIGHT}.

\medskip

For $\theta = 1$ we have $W_1^{[1]}(z) = 0$ and we should push the computation to $W_1^{[2]}(z)$. First, we need the specialization of $\mathcal{F}(z_1,z_2)$ at $z_1 = z_2$ and $\theta = 1$, which follows from \eqref{Fz1z2}
\begin{equation}
\label{Fzz}
\mathcal{F}(z,z) = \frac{1}{2\kappa (\sigma(z))^4}.
\end{equation}
Second, we observe that $\sigma(z)\sigma''(z) = 1 - \frac{z^2}{(\sigma(z))^2}$. Then
\begin{equation}
\label{W12zzz}
\begin{split}
W_1^{[2]}(z) & = \frac{1}{\sigma(z)} \oint_{\gamma} \frac{\dd\zeta}{2\ii\pi(\zeta - z)}\left(\frac{\kappa}{24}\bigg(\frac{\zeta^2}{(\sigma(\zeta))^2} - 1\bigg) - \frac{\kappa \sigma(\zeta)\,\mathrm{coth}[\kappa \sigma(\zeta)]\,\mathcal{F}(\zeta,\zeta)}{2\kappa}\right) \\
& = \frac{\kappa}{24\sigma(z)}\bigg(1 - \frac{z^2}{(\sigma(z))^2}\bigg) + \frac{\mathcal{F}(z,z)}{2\kappa \sigma(z)} + \frac{1}{\sigma(z)} \oint_{\gamma} \frac{\dd \zeta}{2\ii\pi(z - \zeta)}\,\frac{\textnormal{S}[\kappa \sigma(\zeta)]}{4\kappa^2(\sigma(\zeta))^4}.
\end{split}
\end{equation}
To get this last formula we decomposed the hyperbolic cotangent using the function $\textnormal{S}$, and for all other terms we could move the contour to infinity picking the residue at $\zeta = z$ and getting no contribution from infinity because the integrands were $O(\frac{1}{\zeta^2})$ as $\zeta \rightarrow \infty$. In the second term we have refrained from replacing $\mathcal{F}(z,z)$ by its expression \eqref{Fzz} for easier comparison in the next paragraph, but we did so in the last term.

\medskip

Let us turn briefly to $\sqrt{2\theta \kappa} > \pi$. We can also apply Theorem~\ref{Theorem_correlators_expansion} for the conditioned ensemble. It tells us that $W_1^{[1]}(z) = 0$ for $\theta = 1$, and that for general $\theta > 0$ the covariance is the fundamental solution $\mathcal{F}(z_1,z_2)$ computed for the matrix of interactions
\[
\boldsymbol{\Theta} = \left(\begin{array}{cc} \theta & \theta \\ \theta & \theta \end{array}\right)
\]
and the bands $(\alpha_1,\beta_1)$ and $(\alpha_2,\beta_2)$ of the conditioned ensemble (they depend on all the parameters, including the filling fractions $\mathfrak{n}_1,\mathfrak{n}_2$, chosen in Proposition~\ref{proposition_Gaussian_conditioning}). This is evaluated in Corollary~\ref{CoH2fund}, concluding the proof of Proposition~\ref{CLTGAUSSIANWEIGHT2}.

\subsection{Comparison with the Gaussian \texorpdfstring{$\sbeta$}{beta}-ensembles}
\label{sec:comparcd}

For a meaningful comparison we restrict to the case without saturation, that is $\sqrt{2\theta\kappa} < \pi$. The expression for the finite size corrections \eqref{W11zzz1} and \eqref{W12zzz} to $W_1$ seems to be new. Interestingly, they differ from the finite size corrections to $W_1$ in the continuous $\sbeta$-ensembles with repulsion $\theta = \frac{\sbeta}{2}$. The continuous ensemble in question was mentioned in Chapter~\ref{Section_continuous}, it is defined by the probability measure on $N$ particles $\boldsymbol{\lambda} = (\lambda_1,\ldots,\lambda_N) \in \amsmathbb{R}^N$
\begin{equation}
\label{lawgnuzn}
\dd \amsmathbb{P}_{N}(\boldsymbol{\lambda}) = \frac{1}{\mathscr{Z}_{N}} \cdot \prod_{1 \leq i < j \leq N} |\lambda_i - \lambda_j|^{2\theta} \cdot \prod_{i = 1}^{N} e^{-N \kappa \lambda_i^2} \dd\lambda_i.
\end{equation}
Its equilibrium measure is the semi-circle law supported on $[-\beta,\beta]$ with $\kappa \beta^2 = 2\theta$ (note this $\beta$ is different from the intensity of repulsion $\sbeta = 2\theta$) like in \eqref{sigmasqbeta}. Here is the analogue of Proposition~\ref{CLTGAUSSIANWEIGHT}.

\begin{proposition} \cite[Theorem 2.4]{Johansson}
\label{CLTGAUSSIANWEIGHTCC}
Assume $\theta$ and $\kappa$ positive. Let $f$ be an analytic test function and keep the notations $\sigma(z)$ and $\beta$ of \eqref{sigmasqbeta}. The random variable
\[
\sum_{i = 1}^{N} f(\lambda_i) - N \int_{-\beta}^{\beta} f(x)\dd \mu(x)
\]
is approximated in the sense of moments as $N \rightarrow \infty$ by a Gaussian random variable, with mean
\[
\textnormal{\textsf{Mean}}[f] = \frac{\theta -1}{2\theta}\left(\frac{f(-\beta) + f(\beta)}{2} - \frac{1}{\pi} \int_{-\beta}^{\beta} \frac{f(x) \dd x}{|\sigma(x)|} \right),
\]
and covariance like in Proposition~\ref{CLTGAUSSIANWEIGHT}.
\end{proposition}

We can rederive Proposition~\ref{CLTGAUSSIANWEIGHTCC} and discuss more precisely the first subleading order. It will be sufficient for this discussion to work with the first three correlators, defined by\footnote{Unlike in \eqref{eq_correlators_def}, we did not recenter $W_1(z)$ with respect to $N\Gm_{\mu}(z)$. This is more convenient to write the exact Dyson--Schwinger equations.}
\begin{equation*}
\begin{split}
W_1(z) & = \sum_{i = 1}^{N} \amsmathbb{E}\bigg[\frac{1}{z - \lambda_i}\bigg], \\
W_2(z_1,z_2) & =  \sum_{i_1,i_2 = 1}^{N} \left(\amsmathbb{E}\bigg[ \frac{1}{(z_1 - \lambda_{i_1})(z_2 - \lambda_{i_2})}\bigg] - \amsmathbb{E}\bigg[ \frac{1}{z_1 - \lambda_{i_1}}\bigg] \cdot \amsmathbb{E}\bigg[ \frac{1}{z_2 - \lambda_{i_2}}\bigg]\right), \\
W_3(z_1,z_2,z_3) & = \sum_{i_1,i_2,i_3 = 1}^{N} \left(\amsmathbb{E}\bigg[\prod_{p = 1}^3 \frac{1}{z_p - \lambda_{i_p}}\bigg] - \sum_{p = 1}^{3} \amsmathbb{E}\bigg[\frac{1}{z_p - \lambda_{i_p}} \bigg] \cdot \amsmathbb{E}\bigg[\frac{1}{(z_q - \lambda_{i_q})(z_r - \lambda_{i_r})}\bigg]\right. \\
& \qquad\qquad  \qquad \left.+ 2\prod_{p = 1}^{3} \amsmathbb{E}\bigg[\frac{1}{z_p - \lambda_{i_p}}\bigg]\right),
\end{split}
\end{equation*}
with the law \eqref{lawgnuzn} and where $\{p,q,r\} = \{1,2,3\}$. The existence of an asymptotic expansion
\begin{equation}
\label{thebetaexp}
\begin{split}
W_1(z) & = N\Gm_{\mu}(z) + W_1^{[1]}(z) + \frac{W_1^{[2]}(z)}{N} + O\bigg(\frac{1}{N^2}\bigg), \\
W_2(z_1,z_2) & = \mathcal{F}(z_1,z_2) + O\bigg(\frac{1}{N}\bigg), \\
W_3(z_1,z_2,z_3) & = O\bigg(\frac{1}{N}\bigg),
\end{split}
\end{equation}
as $N \rightarrow \infty$ is justified in \cite{BG11}. Formulae for the expansion of correlators for this ensemble have been discussed by many authors, see \textit{e.g.} \cite{Shakirov,Marchalbeta,ForresterWitte}. They can easily be derived by writing down the first two Dyson--Schwinger equations of the ensemble:
\begin{equation}
\label{SDeqnbeta}
\begin{split}
0 & = W_2(z,z) + \big(W_1(z)\big)^2 + \frac{\theta - 1}{\theta}\,\partial_z W_1(z) - \frac{2N\kappa}{\theta}\,z(W_1(z) - N), \\
0 & = W_3(z,z,z_2) + 2W_1(z)W_2(z,z_2) + \frac{\theta - 1}{\theta} \partial_z W_2(z,z_2) - \frac{2N\kappa}{\theta} z W_2(z,z_2) \\
& \quad + \frac{1}{\theta}\,\partial_{z_2}\bigg(\frac{W_1(z) - W_1(z_2)}{z - z_2}\bigg).
\end{split}
\end{equation}
Inserting the expansions \eqref{thebetaexp} in the first Dyson--Schwinger equation of \eqref{SDeqnbeta} gives after a short calculation
\begin{equation}
\label{W11zzz1c}
\begin{split}
\Gm_{\mu}(z) & = \frac{\kappa(z - \sigma(z))}{\theta}, \\
W_1^{[1]}(z) & = -\frac{(\theta - 1) \partial_z \mathcal{G}_{\mu}(z)}{\theta(2\mathcal{G}_{\mu}(z) - \frac{2\kappa}{\theta} z)} = \frac{\theta - 1}{2\theta}\cdot\frac{\sigma'(z) - 1}{\sigma(z)}, \\
W_1^{[2]}(z) & = \frac{\mathcal{F}(z,z)}{2\kappa \sigma(z)} + O(\theta - 1),
\end{split}
\end{equation}
where the $O(\theta - 1)$ is a term vanishing if $\theta = 1$. In turn, inserting \eqref{thebetaexp} in the second Dyson--Schwinger equations yields
\[
\mathcal{F}(z,z_2) = \frac{1}{2\theta\Gm_{\mu}(z) - 2\kappa z} \cdot \partial_{z_2}\left(\frac{\Gm_{\mu}(z) - \Gm_{\mu}(z_2)}{z - z_2}\right) = \frac{1}{2\theta(z - z_2)^2}\left(- 1 + \frac{z_1z_2 - \frac{2\theta}{\kappa}}{\sigma(z_1)\sigma(z_2)}\right).
\]
The leading terms $\Gm_{\mu}(z)$ and $\mathcal{F}(z_1,z_2)$ coincide with those in the discrete ensemble, \textit{cf.} \eqref{GGZG} and \eqref{Fz1z2}, and the latter gives the formula for the covariance in Proposition~\ref{CLTGAUSSIANWEIGHTCC}. The function $W_1^{[1]}(z)$ in the discrete ensemble \eqref{W11zzz1} has two terms. The first one is exactly $W_1^{[1]}(z)$ of the Gaussian $\sbeta$-ensemble given by \eqref{W11zzz1c}, and the same manipulations that led to the formula for the mean in Proposition~\ref{CLTGAUSSIANWEIGHT} give the formula for the mean in Proposition~\ref{CLTGAUSSIANWEIGHTCC}. The second term in \eqref{W11zzz1} is due to the presence of trigonometric functions instead of rational functions. This can be traced back to the fact that in the discrete ensemble, there is some symmetry between $\mu(x)$ (density for particles) and $\frac{1}{\theta} - \mu(x)$ (density for holes), which is manifest in $\sin(\pi \mu(x))$ but not in its rational analogue $\pi\mu(x)$, and here $\pi \mu(x) = \kappa|\sigma(x)|$.

For $\theta = 1$, Proposition~\ref{CLTGAUSSIANWEIGHT} for the discrete case and Proposition~\ref{CLTGAUSSIANWEIGHTCC} for the continuous case give the same result: we have $W_1^{[1]}(z) = 0$ in both cases. The leading term in $W_1(z)$ is rather $N^{-1}W_1^{[2]}(z)$. The $W_1^{[2]}(z)$ in the discrete ensemble \eqref{W12zzz} is a sum of three terms. The second one is exactly $W_1^{[2]}(z)$ of the Gaussian $\sbeta$-ensemble given by \eqref{W11zzz1c}, the third one again comes from the presence of trigonometric functions instead of rational functions, but we also have an extra term (the first one in \eqref{W12zzz}) involving a rational function of $z,\sigma(z)$. We recognize from the coefficient $\frac{1}{24}$ in \eqref{expobserv0} that it stems from derivatives of correlators appearing in the expansions of the observables involved in the Nekrasov equations (that replace the Dyson--Schwinger equations).

For general discrete ensembles even when $\theta = 1$ and in absence of saturations, this indicates that the all-order asymptotic expansion of the correlators is \emph{not governed} by Eynard--Orantin topological recursion applied to the spectral curve (in our case, the curve of equation $y^2 = x^2 - \frac{2}{\kappa}$) \emph{beyond the leading order}: we have seen explicitly that the first subleading term in $W_1$ differs from the one in the Gaussian unitary (\textit{i.e.} $\sbeta = 2$) ensemble. This invalidates the results of the heuristic derivation of \cite{Eynardpart,Eynardpart2} beyond leading order. This discrepancy between discrete and continuous ensembles has two origins: the presence of trigonometric functions (instead of rational functions) in the computations of the expansion for the discrete ensemble, but also the presence of more complicated higher-order derivative terms. This mismatch is not totally surprising:  for continuous $\sbeta$-ensembles Dyson--Schwinger equations encode Virasoro constraints (see \textit{e.g.} \cite{Shakirov}), while for discrete ensembles, the Nekrasov equations rather encode constraints related to a deformation of the Virasoro algebra \cite{HouYang,NieriZen} (in the latter reference one should consider $\theta = -\frac{\epsilon_1}{\epsilon_2}$ and $\frac{1}{N} = \sqrt{\epsilon_1\epsilon_2}$). The question of identifying the precise algorithm (even for $\theta = 1$) involving only the geometry of the spectral curve which should replace the Eynard--Orantin recursion in this context remains open.

\chapter{Asymptotics with fixed filling fractions: partition function}
\label{Chapter_partition_functions}

The aim of this chapter is to develop the large $\N$ asymptotic expansion of the logarithm of the partition
functions for discrete ensembles with segment filling fractions $\hat{\boldsymbol{n}}$ deterministically fixed
by the equations ($\star$) and having only one band per segment. Like in Chapter~\ref{Chapter_fff_expansions} for the correlators, we could in principle derive this asymptotic expansion up to any order $o(\N^{-m})$ as $\N \rightarrow \infty$, but we restrict ourselves to do so up to $o(1)$.

In a first stroke, we reduce the understanding of the partition function of the one-band case to a single $zw$-discrete ensemble already analyzed in Proposition~\ref{AsymZWpart}, by interpolating between their weights. In a second stroke, in the multi-band case, we decouple the particles of the various segments by interpolating the interactions between $\theta_{g,h}$ and $0$ for $g \neq h$. For both steps, this procedure allows computing the ratio of the two partition functions related by the interpolation in terms of correlators of the interpolating family of discrete ensembles. The asymptotic expansion of these correlators was studied in Chapter~\ref{Chapter_fff_expansions}, on which we rely.

When the number of particles near each band is not deterministically fixed --- this is the case either when there are several bands per segment, or when the equations ($\star$) do not deterministically fix the segment filling fractions --- the analysis is more complicated as the correlators do not a priori enjoy the $\frac{1}{\N}$-expansion described in Chapter~\ref{Chapter_fff_expansions}. This more general situation will be treated in the next Chapter~\ref{Chapter_filling_fractions}.

\section{One-band case: interpolation to \texorpdfstring{$zw$}{zw}-discrete ensembles}\label{Section_Partition_onecut}
\label{Sec91}
In this section we work with $H=1$ segment and omit the lower indices $1$ in the parameters: for instance, the segment $[\hat a_1,\hat b_1]$ is written as $[\hat a,\hat b]$, \textit{etc.}

In this section and in the next one, we are going to deal with families of discrete ensembles depending on segment filling fractions and endpoints, which in the $H=1$ case become the triplet $\boldsymbol{t} = (\hat{a},\hat{b},\hat{n})$. The rest of the data will depend on this triplet, and we will assume that they are $p$ times continuously differentiable for some $p \in \amsmathbb{Z}_{\geq 0}\cup \{\infty\}$. There is a tricky point here: $(\hat a, \hat b, \hat n)$ must satisfy integrality conditions, namely
\begin{equation}
\label{integra}
\N\hat n \in\amsmathbb Z_{>0} \qquad \textnormal{and}\qquad \N(\hat b - \hat a) -\theta (\N \hat n -1)\in \amsmathbb Z_{\geq 0},
\end{equation}
as otherwise the state space $\W_\N$ of Section~\ref{Section_configuration_space} would be empty. So, what we really mean is that the rest of the data (weights, potentials, auxiliary functions, \textit{etc.}) is encoded in formulae having a $p$-times continuously differentiable dependence in $\boldsymbol{t} = (\hat n, \hat a, \hat b)$ varying in some open of $\amsmathbb{R}^3$, while an actual discrete ensemble is obtained only if this triplet additionally satisfies the integrality conditions \eqref{integra}. For instance having a formula for the weight $w(x)$ involving a factorial of the form $(x+\N\hat n)!$ is not allowed, because the factorial is not defined for non-integral $x+\N\hat n$. On the other hand, if we first rewrite the same expression as $\Gamma(x+\N\hat n +1)$, then plugging non-integral $x+\N\hat n$ becomes possible and meaningful. Such shifted Gamma functions often appear in our weights, \textit{cf.} Assumption~\ref{Assumptions_analyticity}.

\begin{theorem}
\label{Theorem_partition_one_band}
Consider a family of discrete ensembles with $H=1$ parameterized by $\boldsymbol{t} = (\hat{a},\hat{b},\hat{n})$ in some open set of $\amsmathbb{R}^3$, where $[\hat{a},\hat{b}]$ represents the defining segment and $\hat{n}$ the segment filling fraction. We suppose that the weights satisfy Assumption~\ref{Assumptions_analyticity} with integers $\iota^{\pm}$ and complex domain $\amsmathbb{M}$ independent of $\boldsymbol{t}$ while the regular part of the potential $U(z)$ and the parameters $\rho_j^{\pm}$ for $j \in [\iota^{\pm}]$ are twice-continuously differentiable with respect to $\boldsymbol{t}$ for any $z \in \amsmathbb{M}$ --- with partial derivatives up to order $2$ which are holomorphic functions of $z \in \amsmathbb{M}$. Assume that the discrete ensemble for $\boldsymbol{t}^0$ satisfies Assumptions~\ref{Assumptions_Theta}, \ref{Assumptions_basic}, \ref{Assumptions_offcrit}, \ref{Assumptions_analyticity} and that the equilibrium measure has a single band.

Then, there exists a constant $C > 0$ depending only on the constants in the assumptions, such that for any $\varepsilon > 0$ and for any $\boldsymbol{t}$ in the $\frac{1}{C}$-neighborhood of $\boldsymbol{t}^0$ and satisfying the integrality constraints \eqref{integra}, for $\N$ large enough in a way that depends only on $\varepsilon$ and $C$ we have
\begin{equation}
\label{eq_partition_one_band}
\Big|\log \Z_\N -\big( \I[\boldsymbol{\mu}]\N^2  + \theta \hat n\, \N \log \N + \mathbbm{Rest}_1 \N \big) \Big| \leq  \eps,
\end{equation}
and the absolute values of $\mathbbm{Rest}_1$ and its first and second partial derivatives with respect to $\boldsymbol{t}$ are bounded by $C$.
\end{theorem}
The term $\mathbbm{Rest}_1$ depends on $\N$ and, therefore, there is some freedom in moving the parts of the asymptotic expansion of $\log\Z_\N$ between $\mathbbm{Rest}_1$ and the remainder in \eqref{eq_partition_one_band}, which is bounded from above by $\eps$. The difference between $\mathbbm{Rest}_1$ and the remainder is that we prove bounds on partial derivatives with respect to $\boldsymbol{t}$ for $\mathbbm{Rest}_1$, but not for the remainder.

For the proof of Theorem~\ref{Theorem_partition_one_band}, the strategy is to interpolate
between our ensemble and the $zw$-discrete ensemble, controlling the change of the
partition function in the process. The expression \eqref{eq_partition_one_band} is therefore obtained as an integral of the asymptotic expansions
of Chapter~\ref{Chapter_fff_expansions}. It would be interesting --- although cumbersome and not attempted here --- to find more explicit formulae for the term $\mathbbm{Rest}_1$.

For the interpolation argument we need the following lemma. We recall that there are four possible types of the equilibrium measure in the one-band case: \emph{void-band-void}, \emph{void-band-saturation}, \emph{saturation-band-void}, and \emph{saturation-band-saturation}.

\begin{lemma}
\label{Lemma_off_critical_interpolation}
Consider two discrete ensembles satisfying Assumptions~\ref{Assumptions_Theta}, \ref{Assumptions_basic}, \ref{Assumptions_offcrit}, \ref{Assumptions_analyticity} and the Conditions 1.,2.,3.,4.,5. of the additional Assumption~\ref{Assumptions_extra} with $H = 1$. Assume the two ensembles differ only by the weights $w^{(0)}$ and $w^{(1)}$ and the regular parts of the potential $U^{(0)}$ and $U^{(1)}$ appearing in Assumption~\ref{Assumptions_analyticity}. Their respective unique band are denoted $(\alpha^{(0)},\beta^{(0)})$ and $(\alpha^{(1)},\beta^{(1)})$.
 For any $\eta>0$ there exists $\eps>0$ depending only on the constants in the assumptions and $\eta$, such that if
 \begin{equation}
 \label{eq_smallbandpert} |\alpha^{(0)}-\alpha^{(1)}|+|\beta^{(0)}-\beta^{(1)}|<\eps,
 \end{equation}
 then for all $u \in [0,1]$, the discrete ensemble with weight $w^{(u)}(x)= \big(w^{(0)}(x)\big)^{1-u} \cdot \big(w^{(1)}(x)\big)^{u}$ and regular part of the potential $U^{(u)} = (1-u)U^{(0)} + uU^{(1)}$ also satisfies Assumptions \ref{Assumptions_Theta}, \ref{Assumptions_basic}, \ref{Assumptions_offcrit}, \ref{Assumptions_analyticity} and the Conditions 1.,2.,3.,4.,5. of the additional Assumption~\ref{Assumptions_extra}, with constants that can be chosen independently of $u$. Denoting $(\alpha^{(u)},\beta^{(u)})$ the unique band of the corresponding ensemble, we have
\begin{equation}
\label{eq_x180}
\sup_{u \in [0,1]} \big( |\alpha^{(u)}-\alpha^{(0)}|+|\beta^{(u)}-\beta^{(0)}|\big)<\eta.
\end{equation}
\end{lemma}
It is unclear whether the ensemble with weight $w^{(u)}(x)$ for $u \in (0,1)$ satisfies Condition 6. in Assumption~\ref{Assumptions_extra}, even if we assume that it does at $u=0$ and $u=1$. For this reason, we did not talk about this Condition 6. in the lemma.

\begin{proof}[Proof of Lemma~\ref{Lemma_off_critical_interpolation}]
 Observe that the ensemble with weight $w^{(u)}(x)$ satisfies Assumption~\ref{Assumptions_Theta}, as the intensity of interaction $\theta > 0$ and filling fraction $\hat{n}$ are unchanged. It also satisfies Assumptions~\ref{Assumptions_basic} and \ref{Assumptions_analyticity} with potential $V^{(u)}(x) :=(1-u) V^{(0)}(x)+u V^{(1)}(x)$. On the other hand, checking Assumption~\ref{Assumptions_offcrit} requires some efforts. For that, let $\mu^{(u)}$
be the equilibrium measure for the $u$-dependent ensemble. We are going to construct
 \emph{another} equilibrium measure $\widetilde{\mu}^{(u)}$, associated to a variational datum --- rather than a discrete ensemble, as discussed in Section~\ref{Section_parameters} --- which will be off-critical by construction, and
 which is close enough to $\mu^{(u)}$ so that we can use Theorem~\ref{Theorem_off_critical_neighborhood}.

Let us perform a dilation by $p > 0$ and translation by $q \in \amsmathbb{R}$ on the variational datum associated to $\mu^{(1)}$. The latter involves the shifted endpoints $\hat a'=\hat a - \frac{1}{\N}\big(\theta-\frac{1}{2}\big)$ and $\hat b'=\hat b + \frac{1}{\N}\big(\theta-\frac{1}{2}\big)$, and after translation and dilation we want to take
\begin{equation}
\label{diltra} \hat a'_{(p,q)} = \frac{\hat{a}' - q}{p}, \qquad \hat b'_{(p,q)}=\frac{\hat{b} - q}{p},\qquad V_{(p,q)}(x) = \frac{V^{(1)}(px + q)}{p},\qquad \hat{n}_{(p,q)} = \frac{\hat{n}}{p},
\end{equation}
while the intensity of interactions $\theta$ is unchanged. The density of the associated equilibrium measure is $\mu_{(p,q)}(x) = \mu^{(1)}(px + q)$. Because both $\mu^{(0)}$ and $\mu^{(1)}$ have a single band, there is a unique choice of $(p,q)$ such that the bands of $\mu^{(0)}$ and $\mu_{(p,q)}$ coincide. By \eqref{eq_smallbandpert} we have \begin{equation}
\label{boundpqd}
\max\big(|1-p|,|q|\big)<C_0\,\eps
\end{equation}
for a constant $C_0>0$ depending only on the constants in Assumption~\ref{Assumptions_offcrit} for $\mu^{(0)}$. We fix $(p,q)$ to these values.

It is inconvenient for the interpolation argument that the dilated and translated data live on a segment different from $[\hat a', \hat b']$. Hence, there is an extra step to make the segments coincide. We choose a small enough $\delta>0$ such that $[\hat a'+\delta, \hat b'-\delta] \subset [\hat a'_{(p,q)}, \hat b'_{(p,q)}]$ and bands of both $\mu^{(0)}$ and $\mu_{(p,q)}$ are included in $[\hat a'+2\delta, \hat b'-2\delta]$. Such a choice is possible because the distance between the band and the endpoints of the segment must be bounded away from $0$ due to Assumption~\ref{Assumptions_offcrit}.
For any $u \in [0,1]$, we define a variational datum by keeping $\theta$ as intensity of interactions, taking $\widetilde{\amsmathbb{A}} = [\hat{a}' + \delta,\hat{b}' - \delta]$ as the segment, and choosing the segment filling fraction and the potential as
\begin{equation}
\label{potence}
\begin{split}
\widetilde{n}^{(u)} & = (1 - u)\big(\hat{n} - \mu^{(0)}(\widetilde{\amsmathbb{A}}^{\textnormal{c}})\big) + u\big(\hat{n}_{(p,q)} - \mu_{(p,q)}(\widetilde{\amsmathbb{A}}^{\textnormal{c}})\big), \\
\widetilde{V}^{(u)}(x) & = (1 - u)\bigg(V^{(0)}(x) - 2\theta \int_{\widetilde{\amsmathbb{A}}^{\textnormal{c}}} \log|x - y|\mu^{(0)}(y) \dd y\bigg) \\ & \quad + u\bigg(V_{(p,q)}(x) - 2\theta \int_{\widetilde{\amsmathbb{A}}^{\textnormal{c}}} \log|x - y|\mu_{(p,q)}(y) \dd y\bigg), \\
\end{split}
\end{equation}
where $\widetilde{\amsmathbb{A}}^{\textnormal{c}}$ is the complement of $\widetilde{\amsmathbb{A}}$. We denote $\widetilde{\mu}^{(u)}$ the corresponding equilibrium measure. By construction, $\widetilde{\mu}^{(0)}$ is the restriction of $\mu^{(0)}$ to $\widetilde{\amsmathbb{A}}$ and $\widetilde{\mu}^{(1)}$ is the restriction of $\mu_{(p,q)}$ to $\widetilde{\amsmathbb{A}}$, and both equilibrium measures have the same band. The characterization of Theorem~\ref{Theorem_equi_charact_repeat_2} then implies that
\[
\widetilde{\mu}^{(u)} = (1 - u)\widetilde{\mu}^{(0)} + u\widetilde{\mu}^{(1)}.
\]

We claim that the $u$-dependent variational data satisfies Assumptions~\ref{Assumption_A}, \ref{Assumption_B}, and \ref{Assumption_C} of Section~\ref{allassuml} for constants that can be chosen independently of $u \in [0,1]$. The only point worth clarifying in the claim is the definition of the functions $\widetilde{U}^{(u)}(x)$ in the decomposition
\[
 \widetilde{V}^{(u)}(x)=\iota^{-}\,\mathrm{Llog}(x-\hat{a}' - \delta) + \iota^+\,\mathrm{Llog}(\hat{b}' - \delta -x) + \widetilde{U}^{(u)}(x)
\]
of Condition 3. in Assumption~\ref{Assumption_B}. Indeed, the logarithmic terms in the definition of $\widetilde V^{(u)}(x)$ could have created singularities in the holomorphic part $\widetilde U^{(u)}(x)$. Let us show that this, in fact, does not happen. The argument depends on the type of the equilibrium measure, which is common to all the involved measures $\mu^{(0)}$, $\mu^{(1)}$, $\mu_{(p,q)}$, and $\widetilde \mu^{(u)}$ for $u \in [0,1]$.

If the type is \emph{void-band-void}, then the integrals over $\widetilde{\amsmathbb{A}}^{\textnormal{c}}$ in the definition of $\widetilde V^{(u)}(x)$ vanish and $\widetilde U^{(u)}(x)=U^{(u)}(x)$ has no singularities. If the type is \emph{void-band-saturation}, then for $x\in \widetilde{\amsmathbb{A}}$ we have
\begin{equation*}
\begin{split}
 \widetilde{V}^{(0)}(x)&=V^{(0)}(x) - 2\theta\int_{\widetilde{\amsmathbb{A}}^{\textnormal{c}}} \log|x - y|\mu^{(0)}(y)\dd y \\
 & = V^{(0)}(x) - 2\int_{\hat{b}' - \delta}^{\hat b'} \log(y - x)\dd y \\
 &= V^{(0)}(x) - 2\,\mathrm{Llog}(\hat b'-x) + 2\,\mathrm{Llog}(\hat{b}' - \delta - x).
\end{split}
\end{equation*}
Because $\iota^+=2$, we can further transform $V^{(0)}(x)- 2\,\mathrm{Llog}(\hat b'-x) =\widetilde{U}^{(0)}(x)$. Similarly, we can decompose
\[
 \widetilde{V}^{(1)}(x)=V_{(p,q)}(x) - 2\theta \int_{\widetilde{\amsmathbb{A}}^{\textnormal{c}}} \log|x - y|\mu_{(p,q)}(y)\dd y = \widetilde{U}_{(p,q)}(x)+2\,\mathrm{Llog}(\hat{b}' - \delta -x).
\]
Plugging into the definition of $\widetilde{V}^{(u)}(x)$, we get
\[
 \widetilde{V}^{(u)}(x)=(1-u)\widetilde{V}^{(0)}(x)+u \widetilde{V}^{(1)}(x)= 2\,\mathrm{Llog}(\hat{b}' - \delta -x) + (1-u)\widetilde{U}^{(0)}(x)+u \widetilde{U}_{(p,q)}(x).
\]
Hence, $\widetilde U^{(u)}(z)=(1-u)\widetilde{U}^{(0)}(z)+u\widetilde{U}_{(p,q)}(z)$ is a holomorphic of $z \in \amsmathbb{M}$, as desired. The \emph{saturation-band-void} and \emph{saturation-band-saturation} types are treated similarly.

We have now established that the variational datum for $\widetilde \mu^{(t)}$ satisfies Assumptions~\ref{Assumption_A}, \ref{Assumption_B}, and \ref{Assumption_C} of Section~\ref{allassuml}. Since the variational datum for $\mu^{(u)}$ and for $\widetilde{\mu}^{(u)}$ are close, for $\varepsilon$ small in the sense specified by the assumptions of Theorem~\ref{Theorem_off_critical_neighborhood}, we can then apply this theorem and conclude that $\mu^{(u)}$ satisfies the off-criticality Assumption~\ref{Assumption_C} --- equivalently, Assumption~\ref{Assumptions_offcrit} --- with constants independent of $u \in [0,1]$.

In addition, this argument shows that $\mu^{(u)}$ has the same type as $\mu^{(0)}$ and $\mu^{(1)}$. This readily implies the validity of Conditions 1.,2.,3.,4.,5. of Assumption~\ref{Assumptions_extra} for $\mu^{(u)}$.

\medskip

It remains to justify \eqref{eq_x180}, \textit{i.e.} that the band $(\alpha^{(u)},\beta^{(u)})$ remains arbitrary close to $(\alpha^0,\beta^0)$ uniformly for $u \in [0,1]$ as $\eps \rightarrow 0$. We start by observing that as $\eps \rightarrow 0$, the dilation and translation factors $p$ and $q$ --- which were chosen to make the bands of $\mu^{(0)}$ and $\mu_{(p,q)}$ coincide --- converge respectively to $1$ and $0$. Hence, the endpoints of the band of $\widetilde{\mu}^{(u)}$ converge to $\alpha^{(0)}$ and $\beta^{(0)}$. Besides, the results of Section~\ref{Section-ParReg} imply that $\Gm_{\widetilde{\mu}^{(u)}}(z) - \Gm_{\mu^{(u)}}(z)$ tends to $0$ for $z$ away from $[\hat{a}',\hat{b}']$. Theorem~\ref{Theorem_regularity_density}-(iv) shows that the endpoints of the bands can be reconstructed as zeros of the holomorphic function $(q^-(z))^2$. Hence, we can use Rouch\'e theorem to conclude that the endpoints of the bands for $\widetilde{\mu}^{(u)}$ are close to the endpoints of the bands for $\mu^{(u)}$.
\end{proof}

\begin{proof}[Proof of Theorem~\ref{Theorem_partition_one_band}]
 The proof is slightly different for the various types of the equilibrium measure and we
only present a detailed discussion of the \textit{saturation-band-void} type.
 We call $X$ the discrete ensemble with parameter $\boldsymbol{t} = (\hat{a},\hat{b},\hat{n})$ given in the assumptions, and $X^0$ the one with parameter $\boldsymbol{t}^0$. We will add exponents ${}^0$ to refer to the quantities associated to the ensemble $X^0$ but to make notations lighter we do not write an exponent ${}^{\boldsymbol{t}}$ for those referring to the ensemble $X$. Whenever we mention twice-continuous differentiability of a quantity with respect to the variables $\boldsymbol{t} = (\hat{a},\hat{b},\hat{n})$, the word of caution at the beginning of Section~\ref{Sec91} apply, and the statement will implicitly include the statement that all its partial derivatives up to order two with respect to $\boldsymbol{t}$ are bounded uniformly by a constant depending only on the constants in the assumptions.

\medskip

\noindent \textsc{Step 1.} Our plan is to interpolate between $X$ and a $zw$-discrete ensemble. We start by specifying the target for interpolation. From Theorems~\ref{Theorem_off_critical_neighborhood} and \ref{Theorem_differentiability_full}, we deduce the existence of $c > 0$ depending only on the constants in the assumptions such that if $\eps' > 0$ is small enough and if
\begin{equation}
\label{eq_espnigh} |\!|\boldsymbol{t} - \boldsymbol{t}^0|\!|_{\infty} < \eps',
\end{equation}
then the ensemble $X$ satisfies Assumptions~\ref{Assumptions_Theta}, \ref{Assumptions_basic}, \ref{Assumptions_offcrit}, \ref{Assumptions_analyticity} and has a unique band $(\alpha,\beta)$, obeying
\begin{equation}
\label{eq_Lipschitzband} |\alpha - \alpha^0| + |\beta - \beta^0| < c \varepsilon'.
\end{equation}
We consider two cases in which we use Propositions~\ref{Proposition_zw_band_and_type_1} and \ref{Proposition_zw_band_and_type_2}, respectively:
\begin{enumerate}
 \item \emph{Long band:} $\beta^0-\alpha^0> \frac{5}{2}\theta \mu^0\big([\alpha^0,\beta^0]\big)$;
 \item \emph{Short band:} $\beta^0-\alpha^0< \frac{7}{2}\theta \mu^0\big([\alpha^0,\beta^0]\big)$.
\end{enumerate}
If both inequalities hold, then we can use either of the propositions. Fix $\delta>0$ small enough.

For the \emph{long band case}, we choose two nonnegative integers $K_\pm$ independent of $\N$, such that
\begin{equation}
\label{eq_x216}
 \alpha^0-\frac{2\delta}{3} < \hat a^0 + \frac{\theta K_-}{\N}< \alpha^0 - \frac{\delta}{3} \qquad \textnormal{and} \qquad \beta^0+\frac{\delta}{3} <\hat b^0- \frac{K_+}{\N} < \beta^0+\frac{2\delta}{3}.
\end{equation}
If $\N$ is large enough, then such $K_\pm$ exist. Note that $K_-$ is multiplied by $\theta$, but $K_+$ is not: indeed, in the \emph{saturated-band-void} type, with overwhelming probability the particles are densely packed with spacings $\theta$ in the left part of the interval while in the right part of the interval there are no particles and thus allowed sites have spacing $1$. Let us choose parameters of the $zw$-discrete ensemble as
\begin{equation}
\label{ABchoice} \hat A_1^0 =\hat a^0 + \frac{\theta K_-}{\N},\qquad \hat B_1^0 =\hat b^0 - \frac{K_+}{\N}, \qquad \hat n^{zw,0}=\hat n^0 - \frac{K_- - 1}{\N},
\end{equation}
The subtracted term in the filling fraction originates from the saturation to the left of the band. Then, using \eqref{eq_Lipschitzband} and choosing $\delta$ and $\eps'$ small enough, we can use Proposition~\ref{Proposition_zw_band_and_type_1} to define parameter $\hat{A}_2^0,\hat{B}_2^0$ so that the band of the auxiliary measure $\widetilde{\mu}$ of that proposition is exactly $(\alpha^0,\beta^0)$. For arbitrary $\boldsymbol{t}$ satisfying \eqref{eq_espnigh}, we call $X^{zw}$ the discrete $zw$-ensemble with parameters
\[
\hat{A}_1 = \hat{a} + \frac{\theta K_-}{\N},\qquad \hat{B}_1 = \hat{b}_1 - \frac{K_+}{\N},\qquad \hat{n}^{zw} = \hat{n} - \frac{K_- - 1}{\N},\qquad \hat{A}_2 = \hat{A}_2^0,\qquad \hat{B}_2 = \hat{B}_2^0.
\]
For this ensemble, let us compare the auxiliary equilibrium measure $\widetilde \mu$ to the true equilibrium measure $\mu^{zw}$ as in Corollary~\ref{Corollary_ZW_off-critical_2}. This corollary and Theorem~\ref{Theorem_differentiability_full} imply that the discrete ensemble $X^{zw}$ satisfies Assumptions~\ref{Assumptions_Theta}, \ref{Assumptions_basic}, \ref{Assumptions_offcrit}, \ref{Assumptions_analyticity} and the band $(\alpha^{zw},\beta^{zw})$ of the true equilibrium measure $\mu^{zw}$ at $\boldsymbol{t} = \boldsymbol{t}^0$ becomes arbitrary close to $(\alpha^0,\beta^0)$ as $\N$ becomes large. Moreover, if we further decrease $\eps'$ in \eqref{eq_espnigh}, then using Theorem~\ref{Theorem_differentiability_full} again, we can guarantee that the band changes very little as $\boldsymbol{t}$ vary and, therefore, the band $(\alpha^{zw},\beta^{zw})$ remains close to $(\alpha,\beta)$ for all $\boldsymbol{t}$ satisfying \eqref{eq_espnigh}.

For the \emph{short band case} we choose only one nonnegative integer $K_-$ independent of $\N$ such that
\[
\alpha^0 - \frac{2\delta}{3} < \hat{a}^0 + \frac{\theta K_-}{\N} < \alpha^0 - \frac{\delta}{3}.
\]
We specify the parameters of the $zw$-discrete ensemble as
\begin{equation}
\label{ABchoiceshort}\hat A_1^0=\hat a^0 + \frac{\theta K_-}{\N}, \qquad \hat n^{zw,0}=\hat n^0 - \frac{K_- - 1}{\N},
\end{equation}
and compared to the long band case, there is new detail to take into account in order to choose the remaining parameters. Proposition~\ref{Proposition_zw_band_and_type_2} guarantees the possibility to match arbitrary band positions with the band of the measure $\widetilde{\mu}$ associated to a $zw$-variational datum as in Theorem~\ref{Proposition_ZW_equilibrium_measure} by tuning $\hat{A}_2^0$ and $\hat B_1^0=\hat B_2^0$ to arbitrary real values. However, $zw$-discrete ensembles require the integrality condition \eqref{integrazwwww} for $\hat{B}_1^0$, namely
\begin{equation}
\label{zwinteg}
\N\big(\hat{B}_1^0 - \hat{A}_1^0 - \theta\hat{n}^{zw,0}\big) \in \theta + \amsmathbb{Z}_{> 0}.
\end{equation}
Given the choice we made for $\hat{A}_1^0$ and the fact that the segment $[\hat{a}^0,\hat{b}^0]$ of the discrete ensemble $X$ already satisfies the integrality condition \eqref{eq_segment_ff_relation}
\begin{equation}
\label{Xinteg}
\N\big(\hat{b}^0 - \hat{a}^0 - \theta \hat{n}^0\big) \in \amsmathbb{Z}_{\geq 0},
\end{equation}
this yields
\[
\N\big(\hat{b}^0 - \hat{B}_1^0\big) \in \amsmathbb{Z}.
\]
So, in a first stroke we use Proposition~\ref{Proposition_zw_band_and_type_2} with the values of $\hat{A}_1^0$,$\hat{n}^{zw,0}$ chosen by \eqref{ABchoiceshort} to find real values of $\hat{A}_2^0$ and $\hat{B}_1^0 = \hat{B}_2^0$ such that the band of $\widetilde{\mu}$ from this proposition matches the band $(\alpha^0,\beta^0)$. In a second stroke, we change the value of $\hat{B}_1^0 = \hat{B}_2^0$ by a $O\big(\frac{1}{\N}\big)$ so that it becomes of the form
\[
\hat{B}_1^0 = \hat{B}_2^0 = \hat{b}^0 - \frac{K_+}{\N}
\]
for some nonnegative integer $K_+$ independent of $\N$. Then, for any $\boldsymbol{t}$ satisfying \eqref{eq_espnigh}, we take as discrete ensemble $X^{zw}$ for the short band case the $zw$-ensemble with parameters
\[
\hat{A}_1 = \hat{a} + \frac{\theta K_-}{\N},\qquad \hat{n}^{zw} = \hat{n} - \frac{K_- - 1}{\N},\qquad \hat{A}_2 = \hat{A}_2^0,\qquad \hat{B}_1 = \hat{B}_2 = \hat{B}_1^0 = \hat{B}_2^0.
\]
Note that in the long band case, the choice of $\hat{A}_1,\hat{B}_1,\hat{n}^{zw}$ satisfied the integrality conditions \eqref{zwinteg} allowing the definition of the discrete ensemble $X^{zw}$ whenever $\hat{a},\hat{b},\hat{n}$ satisfied the integrality conditions \eqref{Xinteg} --- not only for $\boldsymbol{t} = \boldsymbol{t}^0$. The discussion then becomes similar to the long band case: choosing $\eps'$ small enough and using Corollary~\ref{Corollary_ZW_off-critical_2} and Theorem~\ref{Theorem_differentiability_full}, we can again guarantee that the ensemble $X^{zw}$ satisfies Assumptions~\ref{Assumptions_Theta}, \ref{Assumptions_basic}, \ref{Assumptions_offcrit}, \ref{Assumptions_analyticity} and that the band $(\alpha^{zw},\beta^{zw})$ of its (true) equilibrium measure $\mu^{zw}$ remains close to $(\alpha,\beta)$ for all $\boldsymbol{t}$ satisfying \eqref{eq_espnigh}.

\medskip

\noindent \textsc{Step 2.} We would like to interpolate between the ensembles $X$ and $X^{zw}$ as in Lemma~\ref{Lemma_off_critical_interpolation}. This lemma requires the ensembles to be defined on the same segments and to satisfy Conditions 1.,2.,3.,4.,5. of Assumption~\ref{Assumptions_extra}. \textit{A priori}, the ensembles $X$ and $X^{zw}$ do not meet these requirements. We are going to fix this by conditioning and localizing these ensembles on the same segment near the band as in Theorem~\ref{proposition_FFF_conditioning}. We momentarily set $\boldsymbol{t} = \boldsymbol{t}^0$. We choose $\delta'>0$ and four nonnegative integers $M_\pm,M^{zw}_\pm$ that may depend on the choice of $\delta$ in Step 1. but are independent of $\N$ and have the following properties.
\begin{enumerate}
 \item $\displaystyle \hat a^0+\theta\frac{M_-}{\N} =\hat A_1^0 + \theta \frac{M_-^{zw}}{\N}$\quad and\quad $\displaystyle\hat b - \frac{M_+}{\N} =\hat B_1 - \frac{M_+^{zw}}{\N}$.\\[4pt]  They represent the left and right endpoints of the segments to which we are going to localize the ensembles at $\boldsymbol{t} = \boldsymbol{t}^0$. The equalities mean that the ensembles $X$ and $X^{zw}$ are going to be localized to the same segments.
 \item $\displaystyle \alpha^0-\frac{\delta'}{3} < \hat a^0 + \theta \frac{M_-}{\N} < \alpha^0 - \frac{\delta'}{3}$\quad and \quad $\displaystyle \beta^0+\frac{\delta'}{3} <\hat b^0-\frac{M_+}{\N}< \beta^0+\frac{2\delta'}{3}.$
   \\[4pt] This condition means that we localize to a small neighborhood of bands, but still keep a little bit of distance to the endpoints of the bands.
 \item $\displaystyle \frac{M_-}{\N} > \delta' +\frac{K_-}{\N}$ and $\displaystyle \frac{M_+}{\N}> \delta' + \frac{K_+}{\N}$. \\[4pt] These conditions guarantee that the endpoints of the segments move at least by $\delta'$ in the localization procedure. Localization might have led to singularities in the potential at \emph{old} endpoints, and this condition guarantees that these singularities are at least some uniform distance away from the new endpoints.
 \end{enumerate}
These choices are always possible due to the way we have constructed $X^{zw}$ at $\boldsymbol{t} = \boldsymbol{t}^0$ in the previous step. We now set for arbitrary $\boldsymbol{t}$
\[
\hat{a}^{\textnormal{loc}} := \hat{a} + \frac{\theta M_-}{\N} = \hat{A}_1 + \frac{\theta M_-^{zw}}{\N} \qquad \hat{b}^{\textnormal{loc}} := \hat{b} - \frac{M_+}{\N} = \hat{B}_1 - \frac{M_+^{zw}}{\N},\qquad \hat{n}^{\textnormal{loc}} = \hat{n} - \frac{M_- - 1}{\N}
\]
We pick a perhaps smaller $\eps'$ for \eqref{eq_espnigh}, such that \eqref{eq_Lipschitzband} guarantees the inequalities
\[
\alpha-\frac{3\delta'}{4} < \hat a + \frac{\theta M_-}{\N} < \alpha - \frac{\delta'}{4}\qquad \textnormal{and} \qquad \beta+\frac{\delta'}{4} <\hat b- \frac{M_+}{\N} < \beta+\frac{3\delta'}{4}.
\]
For any $\boldsymbol{t}$ in the $\eps'$-neighborhood of $\boldsymbol{t}^0$, we call $Y$ (respectively $Y^{zw}$) the localization of the ensemble $X$ (respectively $X^{zw}$) to the segment $[\hat a^{\textnormal{loc}},\hat b^{\textnormal{loc}}]$ following the procedure of Theorem~\ref{proposition_FFF_conditioning}. By construction, the new ensembles $Y$ and $Y^{zw}$ satisfy Assumptions~\ref{Assumptions_Theta}, \ref{Assumptions_basic}, \ref{Assumptions_offcrit}, \ref{Assumptions_analyticity} and the additional Assumption~\ref{Assumptions_extra} --- we are not going to need Condition 6. in Assumption~\ref{Assumptions_extra}, but having it makes no hurt. We can even choose a smaller $\eps' > 0$ independent of $\N$ such that for any $\boldsymbol{t}$ in the $\eps'$-neighborhood of $\boldsymbol{t}^0$, we can apply Lemma~\ref{Lemma_off_critical_interpolation} to the discrete ensembles $Y$ and $Y^{zw}$. This will be useful in the next step.

\medskip

\noindent \textsc{Step 3.} Before being able to compare the partition functions of $X$ and $X^{zw}$, we first need to compare the partition functions of $Y$ and $Y^{zw}$. We can do this through an interpolation argument. We let $w^{(0)}(x)$ and $w^{(1)}(x)$ denote the weights of the ensembles $Y$ and $Y^{zw}$, respectively. We can apply Lemma~\ref{Lemma_off_critical_interpolation} to $Y$ and $Y^{zw}$ because their bands are close and they have the same intensity of interaction, segment, segment filling fraction, integers $\iota^{\pm}$ and parameters $\rho_{j}^{\pm}$ --- because they both satisfy the additional Assumption~\ref{Assumptions_extra}. Thus, as in Lemma~\ref{Lemma_off_critical_interpolation} we consider the weight
\[
w^{(u)}(x)=(w^{(0)}(x))^{1-u}\cdot (w^{(1)}(x))^u,
\]
with interpolation parameter $u \in [0,1]$ and denote $Y^{(u)}$ the corresponding discrete ensemble. All the quantities associated with this ensemble will be denoted with an exponent ${}^{(u)}$.

Notice that the number of particles $N^{\textnormal{loc}} = \N \hat{n}^{\textnormal{loc}}$ is independent of $u \in [0,1]$. We have
\begin{equation}
\label{eq_x217}
\log\Bigg(\frac{\Z_{\N}^{(1)}}{\Z_{\N}^{(0)}}\Bigg) = \int_{0}^{1} \dd u\,\partial_{u}\log \Z_{\N}^{(u)} = \int_{0}^{1}\dd u\,\E^{(u)}\bigg(\sum_{i = 1}^{N^{\textnormal{loc}}} \partial_{u} \log w^{(u)}(\ell_i)\bigg).
\end{equation}
In the setting of the interpolation of Lemma~\ref{Lemma_off_critical_interpolation}, we see that
\begin{equation}
\label{ftwt} f(z) := \frac{\partial_{u} \log w^{(u)}(\N z)}{\N} = \frac{1}{\N}\,\log\bigg(\frac{w^{(1)}(\N z)}{w^{(0)}(\N z)}\bigg)=U^{(0)}(z)-U^{(1)}(z)=V^{(0)}(z)-V^{(1)}(z)
\end{equation}
is independent of $u$, where $V^{(0)}$, $U^{(0)}$ and $V^{(1)}$, $U^{(1)}$ are the potential and its holomorphic part for the ensembles $Y$ and $Y^{zw}$, respectively. In the last two identities of \eqref{ftwt} we use the fact that, due to Assumptions~\ref{Assumptions_analyticity} and \ref{Assumptions_extra} satisfied by the two ensembles $Y$ and $Y^{zw}$, the products involving Gamma functions in \eqref{eq_ansatzw} are exactly the same in $w^{(1)}$ and $w^{(0)}$ by Assumption~\ref{Assumptions_extra} and cancel in the ratio; similarly, the logarithmic terms in $V^{(0)}$ and $V^{(1)}$ cancel. Therefore, $f(z)$ is holomorphic in a complex neighborhood of $[\hat a^{\textnormal{loc}}, \hat b^{\textnormal{loc}}]$. Hence, choosing a complex contour $\gamma\subset \amsmathbb M$ surrounding $[\hat a^{\textnormal{loc}},\hat b^{\textnormal{loc}}]$, we can rewrite \eqref{eq_x217} as
\begin{equation} \label{eq_x218}
\log\Bigg(\frac{\Z_{\N}^{(1)}}{\Z_{\N}^{(0)}}\Bigg) =\frac{\N}{2\pi \ii} \oint_{\gamma} \bigg( \int_{0}^{1}\E^{(u)}\big[ G(z)\big] \dd u\bigg)\, f(z) \dd z,
\end{equation}
where
\[
G(z) = \sum_{i = 1}^{N^{\textnormal{loc}}} \frac{1}{z - \frac{\ell_i}{\N}}
\]
is the Stieltjes transform of the unnormalized empirical measure for $Y^{(u)}$. By Lemma~\ref{Lemma_off_critical_interpolation}, all the intermediate ensembles $Y^{(u)}$ satisfy Assumptions~\ref{Assumptions_Theta}, \ref{Assumptions_basic}, \ref{Assumptions_offcrit} and \ref{Assumptions_analyticity}, with constants that can be chosen independent of $u \in [0,1]$. Therefore, we can apply Theorem~\ref{Theorem_correlators_expansion_relaxed} to the integrand in \eqref{eq_x218}. Letting $\boldsymbol{t}^{\textnormal{loc}} = (\hat{a}^{\textnormal{loc}},\hat{b}^{\textnormal{loc}},\hat{n}^{\textnormal{loc}})$, for any $\eps > 0$ independent of $\N$, we have as $\N \rightarrow \infty$
 \begin{equation} \label{eq_x219}
\log\Bigg(\frac{\Z_{\N}^{(1)}}{\Z_{\N}^{(0)}}\Bigg) =\frac{\N^2}{2\pi \ii} \oint_{\gamma} \Bigg( \int_{0}^{1} \bigg(\int_{\hat a^{\textnormal{loc}}}^{\hat b^{\textnormal{loc}}} \frac{\mu^{(u)}(x)\,f(z)}{z-x} \dd x\bigg) \dd u\Bigg)\dd z + \mathbbm{Rest}_{1}^{[1]}\N+ O(\N^{-\frac{1}{2}+\eps}).
\end{equation}
The term $\mathbbm{Rest}_1^{[1]}$ is obtained by integrating in $u$ and $z$ the term $\widetilde W_{1}^{[1],(u)}(z)$ of \eqref{eq_x80_2} multiplied by $f(z)$; we do not need the exact formula, it is only important that this term depends smoothly on $\boldsymbol{t}^{\textnormal{loc}}$. Exchanging the order of integration and computing the integral over $z$, we simplify \eqref{eq_x219} to
 \begin{equation} \label{eq_x219_2}
\log\Bigg(\frac{\Z_{\N}^{(1)}}{\Z_{\N}^{(0)}}\Bigg) =\N^2 \int_{0}^{1} \bigg(\int_{\hat a^{\textnormal{loc}}}^{\hat b^{\textnormal{loc}}} f(x)\mu^{(u)}(x)\dd x\bigg)\dd u +  \mathbbm{Rest}_{1}^{[1]} \N + O(\N^{-\frac{1}{2}+\eps}).
\end{equation}
We use Proposition~\ref{Proposition_quadratic_potential_differentiability} with $V(x)$ being the potential of $Y^{(u)}$ and with $\Delta V(x)$ being the difference of the potentials of $Y^{zw}$ and $Y$, that is $\Delta V(x)=V^{(1)}(x)-V^{(0)}(x)=-f(x)$. Hence, \eqref{eq_x219_2} becomes
 \begin{equation} \label{eq_x221}
\log\Bigg(\frac{\Z_{\N}^{Y^{zw}}}{\Z_{\N}^{Y}}\Bigg) = \log\Bigg(\frac{\Z_{\N}^{(1)}}{\Z_{\N}^{(0)}}\Bigg) = \big(\I^{(1)}[\mu^{(1)}]-\I^0[\mu^{(0)}]\big) \N^2 +  \mathbbm{Rest}_{1}^{[1]}\N+ O(\N^{-\frac{1}{2}+\eps}).
\end{equation}
Since $\boldsymbol{t}^{\textnormal{loc}}$ is an affine function of $\boldsymbol{t}$, the term $\mathbbm{Rest}_1^{[1]}$ depends smoothly on $\boldsymbol{t}$ as well.

\medskip

\noindent \textsc{Step 4.} We also need to compare the partition functions of $X$ and $X^{zw}$, and of $Y$ and $Y^{zw}$. Recall the procedure of Theorem~\ref{proposition_FFF_conditioning}: $Y$ is obtained from $X$ by conditioning on having $M_-$ densely packed particles in the interval $[\N\hat a, \N\hat a^{\textnormal{loc}})$ and on having no particles in $(\N\hat b^{\textnormal{loc}}, \N\hat b]$. The event we condition on has probability exponentially close to $1$ and, therefore, its effect on the change of the partition function is negligible. What does change the partition function is that we remove the contribution of the $M_-$ particles in $[\N\hat a, \N\hat a^{\textnormal{loc}})$ from the formula of the weight as well as their pairwise interaction, as expressed in \eqref{eq_changeofZ}. In this particular case, the partition function of $X$ is equal to the partition function of $Y$ multiplied by
\begin{equation}
\label{eq_x220}
\frac{\Z_\N^{X}}{\Z_\N^{Y}} = \prod_{1\leq i < j \leq M_-} \frac{1}{\N^{2\theta}}\cdot \frac{\Gamma\big(\theta(j-i)+1\big) \cdot \Gamma(\theta(j-i)+\theta)}{\Gamma\big(\theta(j-i)\big) \cdot \Gamma\big(\theta(j-i)+1-\theta\big)} \cdot \prod_{i=1}^{M_-} w\big(\N \hat a +(i-1)\theta\big),
\end{equation}
The asymptotic expansion of \eqref{eq_x220} as $\N\rightarrow\infty$ can be obtained by combining Lemma~\ref{Lemma_densely_packed_expansion} with Assumption~\ref{Assumptions_analyticity}, Stirling formula \eqref{eq_Stirling_basic}, and the approximation of integrals by Riemann sums. As a result, we transform \eqref{eq_x220} as $\N\rightarrow\infty$ into
\begin{equation}
\label{eq_x223}
\begin{split}
\frac{\Z_\N^{X}}{\Z_\N^Y} & =\exp\Bigg(\bigg(\frac{1}{\theta} \int_{\hat a}^{\hat a^{\textnormal{loc}}} \int_{\hat a}^{\hat a^{loc}} \log|x - y|\,\dd x\,\dd y - \frac{1}{\theta} \int_{\hat a}^{\hat a^{\textnormal{loc}}} V(x) \dd x \bigg)\N^2 +\theta (\hat n-\hat n^{\textnormal{loc}}) \N \log \N \\
&\qquad \qquad + \mathbbm{Rest}_1^{[2]} \N + o(1) \Bigg),
\end{split}
\end{equation}
where $V$ is the potential of the ensemble $X$ and $\mathbbm{Rest}_1^{[2]}$ depends smoothly on $\boldsymbol{t}$. We further use Lemma~\ref{Lemma_energyconditioner} to rewrite the leading $O(\N^2)$ terms in terms of the energy functionals of the equilibrium measure. We conclude that
\begin{equation}
\label{eq_x268}
 \log\Bigg(\frac{\Z_\N^{X}}{\Z_\N^Y}\Bigg) \,\, \mathop{=}_{\N \rightarrow \infty}\,\, \big(\I^{X}[\mu^{X}]-\I^Y[\mu^Y]\big) \N^2  +\theta (\hat n-\hat n^{\textnormal{loc}}) \N \log \N
+ \mathbbm{Rest}_1^{[3]} \N + o(1).
\end{equation}
The comparison of the partition functions of $X^{zw}$ and $Y^{zw}$ is similar, and we obtain
\[
 \log\Bigg(\frac{\Z_\N^{X^{zw}}}{\Z_\N^{Y^{zw}}}\Bigg)\,\, \mathop{=}_{\N \rightarrow \infty}\,\, \big(\I^{X^{zw}}[\mu^{X^{zw}}]-\I^{Y^{zw}}[\mu^{Y^{zw}}]\big) \N^2 +\theta (\hat n^{zw}-\hat n^{\textnormal{loc}}) \N \log \N + \mathbbm{Rest}_1^{[4]}\N + o(1).
\]
where the term $\mathbbm{Rest}_1^{[4]}$ depends smoothly on the parameters $\hat{a},\hat{b},\hat{n}$ of the $zw$-ensemble.

\smallskip

\noindent \textsc{Step 5.} We combine the results of the previous two steps to get as $\N \rightarrow \infty$
\begin{equation*}
\begin{split}
 \log\Bigg(\frac{\Z_\N^X}{\Z_\N^{X^{zw}}}\Bigg) & = \log\Bigg(\frac{\Z_\N^X}{\Z_\N^Y} \cdot \frac{\Z_\N^Y}{\Z_\N^{Y^{zw}}} \cdot \frac{\Z_\N^{Y^{zw}}}{\Z_\N^{X^{zw}}}\Bigg) \\
 & = \big(\I^{X}[\mu^{X}]-\I^{X^{zw}}[\mu^{X^{zw}}]\big)\N^2 +\theta (\hat n -\hat n^{zw})\N \log \N \\
 & \quad + \big(\mathbbm{Rest}_1^{[3]} - \mathbbm{Rest}_1^{[1]} -\mathbbm{Rest}_1^{[4]}\big)\N + o(1).
\end{split}
\end{equation*}
Propositions~\ref{Proposition_ZW_partition} and \ref{Proposition_ZW_partition_leading} together yield as $\N \rightarrow \infty$
\[
 \log \Z_\N^{X^{zw}} = \I^{X^{zw}}[\mu^{X^{zw}}]\,\N^2+ \theta \hat n^{zw}\N \log \N + \mathbbm{Rest}_1^{[6]}\N + o(1).
\]
Combining the last two asymptotic expansions, we get the desired \eqref{eq_partition_one_band}.
\end{proof}

\section{Multi-band case: interpolation to independence}
\label{Section_Partition_multicut}

We now extend the results of Section~\ref{Section_Partition_onecut} to the case of several segments, \textit{i.e.} $H>1$.

\begin{theorem}
\label{Theorem_partition_multicut}
 Consider a family of discrete ensembles with $H\geq 1$ segments, parameterized by $\boldsymbol{t} = (\hat{a}_h,\hat{b}_h,\hat{n}_h)_{h = 1}^{H}$ in some open set of $\amsmathbb{R}^{3H}$, where for each $h \in [H]$, $[\hat{a}_h,\hat{b}_h]$ and $\hat{n}_h$ represent the $h$-th segment and its segment filling fraction. We suppose that the tuples of integers $(\iota_h^{\pm})_{h = 1}^{H}$ and complex domains $(\amsmathbb{M}_h)_{h = 1}^{H}$ independent of $\boldsymbol{t}$, while for any $h \in [H]$ the regular part of the potential $U_h(z)$ for $z \in \amsmathbb{M}_h$ and the reals $(\rho_{h,j}^{\pm})_{j = 1}^{\iota_h^{\pm}}$ are twice-continuously differentiable with respect to $\boldsymbol{t}$ --- with partial derivatives of $U_h(z)$ derivatives up to order $2$ which are holomorphic functions of $z \in \amsmathbb{M}_h$. Assume that the discrete ensemble for $\boldsymbol{t}^0$ satisfies Assumptions~\ref{Assumptions_Theta}, \ref{Assumptions_basic}, \ref{Assumptions_offcrit}, Assumptions~\ref{Assumptions_analyticity} and has one band per segment.

 Then, there exists a constant $C > 0$ depending only on the constants in the assumptions, such that for any $\varepsilon > 0$ and any $\boldsymbol{t}$ in a $\frac{1}{C}$-neighborhood of $\boldsymbol{t}^0$ and satisfying the integrality constraints \eqref{eq_segment_ff_relation}, for $\N$ large enough in a way that depends only on $\varepsilon$ and $C$ we have
  \begin{equation}
 \label{eq_partition_multicut} \Bigg|\log\Z_\N - \bigg[\I[\boldsymbol{\mu}]\N^2 + \bigg( \sum_{h = 1}^{H} \theta_{h,h}\hat{n}_h\bigg) \N\log \N + \mathbbm{Rest}_1\N\bigg]\Bigg| \leq \varepsilon.
 \end{equation}
Besides, the absolute values of $\mathbbm{Rest}_1$ and its partial derivatives with respect to $\boldsymbol{t}$ up to order $2$ are bounded by $C$.
 \end{theorem}
 Similarly to Theorem~\ref{Theorem_partition_one_band}, there is some freedom in moving the parts of the asymptotic expansion between $\mathbbm{Rest}_1\N$ and the remainder in \eqref{eq_partition_multicut}. The difference between these terms is that we prove bounds on partial derivatives with respect to $\boldsymbol{t}$ for $\mathbbm{Rest}_1$ but not for the remainder bounded from above by $\eps$. The term $\mathbbm{Rest}_1$ will be obtained in the proof by combining the similar term in Theorem~\ref{Theorem_partition_one_band} with an integral of asymptotic expansions of Theorem~\ref{Theorem_correlators_expansion_relaxed}. As in Theorem~\ref{Theorem_partition_one_band}, it is in principle possible to obtain an explicit expression for this term, but we do not attempt to do so.

\begin{proof}[Proof of Theorem~\ref{Theorem_partition_multicut}]
The idea of the argument is to interpolate between the ensemble with arbitrary matrix of interactions $\boldsymbol{\Theta}$ and an
 ensemble with a diagonal matrix of interactions, controlling the change in the partition function. The latter ensemble splits into $H$ independent ensembles, and therefore, we can
 use Theorem~\ref{Theorem_partition_one_band} to reach its partition function. If we only changed the matrix elements
 of $\boldsymbol{\Theta}$ without adjusting other parameters of the ensemble, then the equilibrium measure would be changing in
 a non-trivial way, and it would be hard to check if the off-criticality Assumption~\ref{Assumptions_offcrit} is preserved in the interpolation.
  In order to avoid this, we will be changing the potentials together with $\boldsymbol{\Theta}$, so that the equilibrium measure never
  changes.

We call $X$ the ensemble which implicitly depends on parameters $\boldsymbol{t}$ close enough to $\boldsymbol{t}^0$. To keep short notations, all the quantities associated with the ensemble $X$ are denoted without exponent ${}^{\boldsymbol{t}}$. We are going to define a family of discrete ensembles $X^{(u)}$ parameterized by $u \in [0,1]$, such that $X^{(1)} = X$. We set $\boldsymbol{\Theta}^{(1)}=\boldsymbol{\Theta}$ and $\boldsymbol{\Theta}^{(0)}=\textnormal{diag}(\theta_{1,1},\ldots,\theta_{H,H})$ is the diagonal part of $\boldsymbol{\Theta}$. For arbitrary $u \in [0,1]$ we take as matrix of intensities of interactions for $X^{(u)}$
\[
\boldsymbol{\Theta}^{(u)}=u\boldsymbol{\Theta}^{(1)}+(1-u)\boldsymbol{\Theta}^{(0)}.
\]
The diagonal of $\boldsymbol{\Theta}^{(u)}$ does not depend on $u$, while the off-diagonal matrix elements of $\boldsymbol{\Theta}^{(u)}$ have the form $u\theta_{g,h}$. We take the segments and segment filling fractions of $X^{(u)}$ to be the ones of $X$, independently of $u \in [0,1]$. In particular, the total number of particles is the same for all $u$. This implies in particular that the state space $\W_\N$ for $X^{(u)}$ is actually independent of $u$, only the probability measure varies.

Let us define the holomorphic function
\begin{equation}
 \label{eq_x231}
\forall h \in [H] \quad \forall z \in \amsmathbb{M}_h\qquad \Delta V_h(z)=- \sum_{g \neq h} 2\theta_{h,g} \int_{\hat a'_{g}}^{\hat b'_{g}}
 \log\big(((x-y)\textnormal{sgn}(h - g)\big)\mu_{g}(y) \dd y.
\end{equation}
where $\boldsymbol{\mu} = (\mu_h)_{h = 1}^{H}$ is the equilibrium measure of $X$. We complete the definition of the ensemble $X^{(u)}$ by specifying the potential and its regular part in the ensemble $X^{(u)}$ for $h \in [H]$
\begin{equation}
\label{eq_x232}
\begin{split}
 w_h^{(u)}(x) & =w_h(x) \cdot \exp\Bigg((u-1)\N \Delta V_h\bigg(\frac{x}{\N}\bigg)\Bigg), \\
 V_h^{(u)}(x) & = V_h(x) + (1-u)\Delta V_h(x), \\
 U_h^{(u)}(x) & = U_h(x) + (1-u) \Delta V_h(x).
\end{split}
\end{equation}
The auxiliary functions for Assumptions~\ref{Assumptions_analyticity} follow from these expressions by Definition~\ref{Definition_phi_functions}. Explicitly:
\begin{equation*}
\begin{split}
\Phi_h^{+,(u)}(z) & = \Phi_h^{+}(z) \cdot e^{(1-u)\N[\Delta V_h(z + \frac{1}{2\N}) - \Delta V_h(z - \frac{1}{2\N})]}, \\
\phi_h^{+,(u)}(z) & = \phi_h^{+}(z) \cdot e^{(1-u)\partial_z \Delta V_h(z)},
\end{split}
\end{equation*}
while $\Phi_h^{-,(u)}(z) = \Phi_h^{-}(z)$ and $\phi_h^{-,(u)}(z) = \phi_h^{-}(z)$ are independent of $u \in [0,1]$.

We claim that the equilibrium measure of $X^{(u)}$ does not depend on $u \in [0,1]$ and
coincides with $\boldsymbol{\mu}$. This can be seen by applying Theorem~\ref{Theorem_equi_charact_repeat_2} and
noting that the effective potentials of Definition~\ref{def_eff_pot} computed with the measure $\boldsymbol{\mu}$ have the same characterizing properties as the effective potentials for the ensemble $X^{(u)}$. Indeed, the change in $\boldsymbol{\Theta}^{(u)}$ is exactly compensated by the change in $V_h^{(u)}$. In particular, the ensemble $X^{(u)}$ remains off-critical for any $u \in [0,1]$. It follows that the ensemble $X^{(u)}$ satisfies Assumptions~\ref{Assumptions_Theta}, \ref{Assumptions_basic}, \ref{Assumptions_offcrit} and \ref{Assumptions_analyticity} with constants that can be chosen independently of $u \in [0,1]$.

We then have
\begin{equation}
\begin{split}
 \label{eq_x170}
  \log\Bigg(\frac{\Z_\N^{(1)}}{\Z_\N^{(0)}}\Bigg) & =\int_0^1 \partial_{u}\log\Z_\N^{(u)} \dd u \\
 & =\int_0^1 \E^{(u)} \Bigg[\sum_{i=1}^{N} \partial_{u}\log\Bigg(w^{(u)}_{h(i)} \cdot \prod_{\substack{i < j \\ h(i) \neq h(j)}} \frac{\N^{1 - 2u\theta_{h(i),h(j)}} \cdot \Gamma\big(\ell_j - \ell_i + u \theta_{h(i),h(j)}\big)}{\Gamma\big(\ell_j - \ell_i + 1 - u\theta_{h(i),h(j)}\big)}\Bigg)\dd u \Bigg] \\
  & =\int_0^1 \E^{(u)}\bigg[ \N \sum_{i=1}^{N} \Delta V_{h(i)}\bigg(\frac{\ell_i}{\N}\bigg)+ \sum_{\substack{1\leq i<j \leq N \\ h(i) \neq h(j)}} S\bigg(\frac{\ell_j-\ell_i}{\N};\, u,\theta_{h(i),h(j)}\bigg)\bigg] \dd u,
\end{split}
\end{equation}
where
\[
S(z;u,\theta) :=\partial_u\log\bigg( \frac{\N^{1-2u\theta}\cdot \Gamma\big(\N z+u\theta\big)}{\Gamma\big(\N z+1-u\theta)}\bigg) = \theta\big(-2\log \N + \mathsf{\Psi}(\N z + u \theta) + \mathsf{\Psi}(\N z +1- u \theta)\big).
\]
and $\mathsf{\Psi} = (\log \Gamma)'$ is the digamma function. For $\N$ large enough, $S(z;u,\theta)$ is a holomorphic function of $z$ in the region where $\text{Re}(z)$ is positive and bounded away from $0$ uniformly with $\N$. We can rewrite \eqref{eq_x170} in terms of the Stieltjes transform $G_h(z)$ of the empirical measure:
\begin{equation}
 \label{eq_x169_2}
\begin{split}
   \log\Bigg(\frac{\Z_\N^{(1)}}{\Z_\N^{(0)}}\Bigg) & = \int_0^1 \bigg( \sum_{h=1}^H \N \oint_{\gamma_h} \frac{\dd z}{2\ii\pi}\,
  \Delta V_h(z) \E^{(u)} \big[G_h(z)\big] \\
  & \qquad + \sum_{1\leq g< h \leq H} \oint_{\gamma_g}\oint_{\gamma_h} \frac{\dd z\,\dd\zeta}{(2\ii\pi)^2}\, S\big(\zeta - z;u,\theta_{g,h}\big)
   \E^{(u)}\big[ G_g(z) G_{h}(\zeta)\big] \bigg) \dd u.
   \end{split}
   \end{equation}
   At this stage, we can use Theorem~\ref{Theorem_correlators_expansion_relaxed} to obtain the asymptotic expansion of \eqref{eq_x169_2}. For the first contour integral we use \eqref{eq_x80_2}. For the second double contour integral, we rewrite
\[
  \E^{(u)} \big[ G_g(z) G_{h}(\zeta) \big]= W_{1;g}^{(u)}(z) W_{1;h}^{(u)}(\zeta) +W_{2;g,h}^{(u)}(z,\zeta),
\]
where $W_{n;h_1,\ldots,h_n}^{(u)}$ are the correlators in the ensemble $X^{(u)}$, and \eqref{eq_x80_2} for the first term and \eqref{eq_x81_2} for the second term. Since the ensemble $X^{(u)}$ satisfy Assumptions~\ref{Assumptions_Theta}, \ref{Assumptions_basic}, \ref{Assumptions_offcrit}, \ref{Assumptions_analyticity} with constants independent of $u \in [0,1]$, the remainders in the asymptotic expansions of the correlators in Theorem~\ref{Theorem_correlators_expansion_relaxed} are uniform and can be integrated over $u$, as well as $z$ and $\zeta$. We find as $\N \rightarrow \infty$
\begin{equation}
\label{eq_x229}
\begin{split}
\log\Bigg(\frac{\Z_\N^{(1)}}{\Z_\N^{(0)}}\Bigg) & =  \sum_{h=1}^H \int_0^1 \dd u \left( \int_{\hat{a}'_h}^{\hat{b}_h'} \Delta V_h(x) \mu_h(x) \dd x\right) \N^2 \\
& \quad + \sum_{1\leq g < h \leq H} \int_{0}^{1} \dd u \left( \int_{\hat{a}'_g}^{\hat{b}'_g} \int_{\hat{a}_h'}^{\hat{b}_h'} S(x - y;u,\theta_{g,h}) \mu_{g}(x)\mu_h(y) \dd x \dd y\right) \\
& \quad + \mathbbm{Rest}_1^{[7]} \N + O(\N^{-\frac{1}{2}+\eps}).
\end{split}
\end{equation}
We recall that the equilibrium measure $\boldsymbol{\mu}$ of the ensemble $X^{(u)}$ is independent of $u \in [0,1]$, so the integral over $u$ can be removed in the first line. Here $\mathbbm{Rest}_1^{[7]}$ is a uniformly bounded remainder, which is twice-continuously differentiable with respect to $\boldsymbol{t}$ and whose partial derivatives up to order $2$ are also uniformly bounded. Recalling the definition of $\Delta V_h(x)$ from \eqref{eq_x231}, we see that the $\N^2$ terms in \eqref{eq_x30} cancel out and we finally get as $\N \rightarrow \infty$:
\begin{equation}
\label{eq_x239}
 \log\Bigg(\frac{\Z_\N^{(1)}}{\Z_\N^{(0)}}\Bigg) = \mathbbm{Rest}_1^{[8]}\N + o(1).
\end{equation}
Eventually, we are interested in $\log \Z_\N^{(1)}$. By design, in the ensemble $X^{(0)}$ particles in different segments do not interact. Therefore, its partition function factorizes as a product of $H$ partition function, each one being in the framework of Theorem~\ref{Theorem_partition_one_band}. We get as $\N \rightarrow \infty$
\[
\log \Z_\N^{(0)} = \sum_{h = 1}^{H} \I^{(0),h}[\mu_h] \N^2 + \mathbbm{Rest}_1^{[9]}\N + o(1),
\]
where $-\I^{(0),h}$ is the energy functional for the ensemble $X^{(0)}$ restricted to the $h$-th segment $[\hat{a}_h,\hat{b}_h]$. Since $\boldsymbol{\Theta}^{(0)}$ is diagonal, we clearly have
\[
\sum_{h = 1}^{H} \I^{(0),h}[\mu_h] = \I^{(0)}[\boldsymbol{\mu}].
\]
As $\boldsymbol{\mu}$ is the equilibrium measure for the ensemble $X^{(1)}$, and the potentials of $X^{(0)}$ and $X^{(1)} = X$ are related by \eqref{eq_x232}, we conclude that $\I^{(0)}[\boldsymbol{\mu}] = \I^{(1)}[\boldsymbol{\mu}] = \I[\boldsymbol{\mu}]$. Combining with Combining with \eqref{eq_x239}, we get
\begin{equation*}
 \log\Z_\N^{(1)} = \I[\boldsymbol{\mu}]\,\N^2 + \mathbbm{Rest}_1\,\N + o(1).\qedhere
\end{equation*}
\end{proof}

\chapter{Fluctuations of filling fractions and saturated parameters}

\label{Chapter_filling_fractions}

In Chapters~\ref{Chapter_fff_expansions} and \ref{Chapter_partition_functions} we studied the asymptotic expansions of correlators, linear statistics, and partition functions for the discrete ensembles in which the filling fractions are deterministically fixed and the ensemble has one band $(\alpha_h,\beta_h)$ in the $h$-th segment $[\hat a_h',\hat b_h']$ for each $h\in [H]$. In the present chapter we waive these restrictions and proceed to the most general case of arbitrarily fluctuating filling fractions and of arbitrary many bands per segment.

In our approach, we assign new filling fractions to each band of the equilibrium measure. We then proceed in two independent steps: we study the fluctuations of these filling fractions and the asymptotic expansions for the conditioned ensemble when the new filling fractions are fixed. The latter has already been computed in the previous chapters. For the former, we show in this chapter that the asymptotics of the filling fractions is governed by \emph{discrete} Gaussian distributions. The main technical work was done in Chapter~\ref{Chapter_partition_functions}: the asymptotic expansion for partition functions will give the weights of the asymptotic distributions. In the present chapter, in addition to collecting ingredients from the previous chapters, we also establish proper notations. The notations are more complicated in the cases when the equilibrium measure has saturations, as it even becomes more tricky to define filling fractions in such cases, due to the fluctuating nature of the allowed sites for particles. Hence, we first present the results for the case without saturations and only then proceed to the most general setting.

\smallskip

\section{Discrete Gaussian distributions and convergence mode}
In this section we explain what we mean by a discrete Gaussian random variable and convergence to such a random variable.
\label{index:Qufor}
Fix $d \in \amsmathbb{Z}_{> 0}$, a symmetric real bilinear form $\Qu$ on $\amsmathbb R^{d}$, an affine subspace
 $\L\subseteq \amsmathbb{R}^{d}$ and a vector $\boldsymbol u \in\amsmathbb R^d$. We denote $\L_0 \subseteq \amsmathbb{R}^d$ the unique linear subspace parallel to $\L$. We assume
that the restriction of $\Qu$ to $\L_0$ is
positive definite and that $\L \cap \amsmathbb{Z}^d \neq \emptyset$.
\begin{definition}
\label{Definition_disc_Gauss}
A \emph{discrete Gaussian} random variable $\textnormal{\textsf{\textbf{Gau\ss{}}}}_{\amsmathbb{Z}}[\Qu,
\L,\boldsymbol{u}]$ is a random variable taking values in the set
$\L \cap \amsmathbb{Z}^d $ and such that
\begin{equation} \label{eq_Discrete_Gaussian_main}
\forall \boldsymbol{x} \in \L \cap \amsmathbb{Z}^{d}\qquad \amsmathbb{P}\big(\textnormal{\textsf{\textbf{Gau\ss{}}}}_{\amsmathbb{Z}}[\Qu,
\L, \boldsymbol u] = \x\big) = \frac{1}{\mathscr{Q}} \exp\left(-
\frac{1}{2}\Qu(\x - \boldsymbol u,\x - \boldsymbol u) \right),
\end{equation}
where $\mathscr{Q} > 0$ is the normalizing constant turning \eqref{eq_Discrete_Gaussian_main} into a probability measure.
\end{definition}

Note that translating $\mathfrak{L}$ by an integer vector does not change the distribution of $\textnormal{\textsf{\textbf{Gau\ss{}}}}_{\amsmathbb{Z}}[\Qu, \L,\boldsymbol{u}]$, and that $\boldsymbol{u}$ itself does not have to belong to $\mathfrak{L}$. Besides, the shifted random variable $\textnormal{\textsf{\textbf{Gau\ss{}}}}_{\amsmathbb{Z}}[\Qu,\mathfrak{L},\boldsymbol{u}] - \boldsymbol{u}$ has the same distribution as $\textnormal{\textsf{\textbf{Gau\ss{}}}}_{\amsmathbb{Z}}[\Qu,\mathfrak{L} - \boldsymbol{u},0]$. We introduce the following definition in order to study convergence of distributions involving discrete Gaussians.
\begin{definition}
\label{Definition_asymptotic_equivalence}
Let $\boldsymbol{\chi}_{\N}^{(1)}$ and $\boldsymbol{\chi}_{\N}^{(2)}$ be two families (indexed by a large parameter $\N$) of discrete random $d$-dimensional real vectors. We say that the distribution of $\boldsymbol{\chi}_{\N}^{(1)}$ is asymptotically equal to the distribution of $\boldsymbol{\chi}_{\N}^{(2)}$ if
\begin{equation}
\label{eq_asymptotic_equivalence}
 \lim_{\N\rightarrow\infty} \sup_{f \in \textnormal{Lip}_1} \big|\E\big[f(\boldsymbol{\chi}_{\N}^{(1)})\big]-\E\big[f(\boldsymbol{\chi}_{\N}^{(2)})\big]\big| =0,
\end{equation}
where $\textnormal{Lip}_1$ is the set of $1$-bounded $1$-Lipschitz functions
\[
 \textnormal{Lip}_1=\big\{f:\amsmathbb{R}^d \rightarrow \amsmathbb{R} \quad \big|\quad \forall \boldsymbol{x},\boldsymbol{y} \in \amsmathbb{R}^d \quad |f(\x)-f(\y)| \leq |\!|\x-\y|\!|_{\infty}\,\,\textnormal{ and }\, |f(\x)|\leq 1\big\}.
\]
If \eqref{eq_asymptotic_equivalence} holds, we write $\chi_{\N}^{(1)} \stackrel{\textnormal{d}}{\sim} \chi_{\N}^{(2)}$.
\end{definition}

\smallskip

\noindent If $\boldsymbol{\chi}_{\N}^{(2)}=\boldsymbol{\chi}^{(2)}$ does not depend on $\N$, then $\boldsymbol{\chi}_{\N}^{(1)} \stackrel{\textnormal{d}}{\sim} \boldsymbol{\chi}_{(2)}$ is equivalent to the standard definition of convergence in distribution. On the other hand, the definition we use allows random variables to escape to infinity as $\N\rightarrow\infty$, but still remain close to each other. Our definition is closely related to convergence to $0$ of the Vaserstein metric $\mathcal W_1$ between distributions of $\boldsymbol{\chi}_{\N}^{(1)}$ and $\boldsymbol{\chi}_{\N}^{(2)}$: the difference is that the Vaserstein metric is only defined for random variables with finite expectations, while our definition avoids dealing with the first moments by imposing the boundedness condition on the test functions $f$. In this chapter, we do not establish tail bounds for finite but large deviations of filling fractions, and therefore we do not have a sufficient control to have convergence

\section{In absence of saturations}

\label{Section_FFF_no_sat}

\subsection{Definition of fluctuating filling fractions}

In this section we deal with a discrete ensemble from Section~\ref{Section_general_model} satisfying Assumptions~\ref{Assumptions_Theta}, \ref{Assumptions_basic}, \ref{Assumptions_offcrit}, \ref{Assumptions_analyticity} and such that its equilibrium measure $\boldsymbol{\mu}$ has no saturation. Following previous notations, the particles belong to the segments $[\N \hat a_h,\N\hat b_h]$ for $h\in[H]$, the (possibly random) number of particles in the $h$-th segment is denoted $N_h = \N \hat n_h$ and these numbers obey the affine constraints \eqref{eq_equations_eqs} in Section~\ref{DataS}.

By Assumption~\ref{Assumptions_offcrit}, $\boldsymbol{\mu}$ has a finite number $K$ of bands. Let us further denote $(\alpha_k,\beta_k)$ the bands of $\boldsymbol{\mu}$ indexed by $k\in[K]$ so that $\beta_{k}<\alpha_{k+1}$ for any $k \in [K - 1]$. They contain the mass
\begin{equation}
\label{eq_def_equilibrium_ff}
\hat{n}^{\boldsymbol{\mu}}_{k}=\frac{N^{\boldsymbol{\mu}}_{k}}{\N} := \mu_{h^k}([\alpha_k,\beta_k]).
\end{equation}
We recall that $h^k \in [H]$ is the index of the segment containing the $k$-th band, and that $\llbracket k^-(h),k^+(h)\rrbracket := \{k \in [K]\,\,|\,\,h^k = h\}$.

Let us choose $\eps > 0$ small enough independent of $\N$ and depending only on the constants in the assumptions, such that the intervals $[\alpha_k-\eps,\beta_k+\eps]$ are pairwise disjoint for $k \in [K]$. We define the (random) \emph{fluctuating filling fractions} through
\begin{equation}
\label{Nbullek} \forall k \in [K] \qquad N^{\circ}_k= \#\big\{i \in [N] \quad \big| \quad \ell_i \in (\N(\alpha_k-\eps),\N(\beta_k+\eps))\big\}, \qquad \hat{n}_k^{\circ} = \frac{N_k^{\circ}}{\N}.
\end{equation}
We aim at describing the asymptotic distribution of the random integral vector $\boldsymbol{N}^{\circ} = (N^{\circ}_k)_{k = 1}^K$ as $\N\rightarrow\infty$. Theorem~\ref{Theorem_CLT_for_filling_fractions} below shows that this vector is asymptotically equal to a discrete Gaussian distribution in the sense of Definition~\ref{Definition_asymptotic_equivalence}. Let us first describe the parameters of this distribution.

\subsection{Discrete Gaussian part}

\label{Sec:gausds}
By Theorem~\ref{Theorem_ldpsup}, with overwhelming probability as $\N \rightarrow \infty$ there are no particles outside the region ${\bigcup_{k=1}^K [\N(\alpha_k-\eps),\N(\beta_k+\eps)]}$, hence
\begin{equation}
\label{eq_x181}
\forall h \in [H]\qquad \N \hat{n}_h = \sum_{k = k^-(h)}^{k^+(h)} N^{\circ}_{k}.
\end{equation}
Therefore, on an event of overwhelming probability, the affine constraints \eqref{eq_equations_eqs} in Section~\ref{DataS} on $\boldsymbol{N} = (N_h)_{h = 1}^{H}$ can be rewritten as affine constraints on the numbers $\boldsymbol{N}^{\circ} = (N_{k}^{\circ})_{k = 1}^K$. They define an affine subspace $\L \subset \amsmathbb{R}^{K}$ and we call $\L_0$ the underlying linear subspace.

Next, we define the quadratic bilinear form $\Qu$ on $\amsmathbb R^K$ that will enter into \eqref{eq_Discrete_Gaussian_main}. For this purpose, we let an auxiliary vector of filling fractions $\hat{\boldsymbol p} = (\hat{p}_{k})_{k = 1}^{K}$ vary in a small neighborhood of $\hat{\boldsymbol n}^{\boldsymbol{\mu}}$ defined in \eqref{eq_def_equilibrium_ff}, and consider a variational datum consisting of
\begin{itemize}
\item the segments $[\alpha_k - \varepsilon,\beta_k + \varepsilon]$ indexed by $k\in[K]$;
\item the matrix of interactions $(\theta_{h^k,h^{l}})_{k,l = 1}^{K}$;
\item for each $k \in [K]$, the potential given by the restriction of the potential $V_{h^k}$ of the discrete ensemble to $[\alpha_k - \varepsilon,\beta_k + \varepsilon]$;
\item the filling fractions $\hat{\boldsymbol{p}}$.
\end{itemize}
For $\hat{\boldsymbol p} =\hat{\boldsymbol n}^{\boldsymbol{\mu}}$, this variational datum satisfies Assumptions~\ref{Assumption_A}, \ref{Assumption_B} and \ref{Assumption_C} with all the parameters $\iota^\pm$ set to $0$, because the discrete ensemble satisfies the corresponding Assumptions~\ref{Assumptions_Theta}, \ref{Assumptions_basic}, \ref{Assumptions_offcrit}, \ref{Assumptions_analyticity} and has no saturations. The associated equilibrium measure is the restriction of $\mu$ to $\bigcup_{k = 1}^{K} [\alpha_k - \varepsilon,\beta_k + \varepsilon]$. By the stability property established in Theorem~\ref{Theorem_off_critical_neighborhood}, for $\hat{\boldsymbol p}$ in a small enough and $\N$-independent neighborhood of $\hat{\boldsymbol n}^{\boldsymbol{\mu}}$, this variational datum still satisfies Assumptions~\ref{Assumption_A}, \ref{Assumption_B}, \ref{Assumption_C} and there is an associated equilibrium measure $\boldsymbol{\mu}^{\hat{\boldsymbol p}}$, which minimizes the same energy functional $-\I$ as the one for $\hat{\boldsymbol{n}}^{\boldsymbol{\mu}}$, but over tuples of measures with filling fractions equal to $\hat{\boldsymbol p}$. We further define the quadratic form $\Qu$ as minus the Hessian of $\hat{\boldsymbol{p}} \mapsto \I[\boldsymbol{\mu}^{\hat{\boldsymbol p}}]$ evaluated at $\hat{\boldsymbol{p}} = \hat{\boldsymbol{n}}^{\boldsymbol{\mu}}$:
\begin{equation}
\label{eq_discrete_covariance}
\forall k,l \in [K] \qquad \Qu_{k,l} = - \partial_{\hat{p}_k}\partial_{\hat{p}_l} \I[\boldsymbol{\mu}^{\hat{\boldsymbol p}}]\big|_{\hat{\boldsymbol p}=\hat{\boldsymbol n}^{\boldsymbol{\mu}}}.
\end{equation}
Proposition~\ref{Proposition_Hessian_free_energy} shows that the restriction of $\Qu$ to $\L_0$ is positive definite and its values on the intersection of the unit sphere in $\amsmathbb R^K$ with $\L_0$ are uniformly bounded and bounded away from $0$ as $\N\rightarrow\infty$.

\medskip

The definition of the vector $\boldsymbol{u}$ entering into \eqref{eq_Discrete_Gaussian_main} is more delicate. We start by defining a polynomial function $P(\x)$ of $\x\in \L\subseteq \amsmathbb R^K$, which has degree $2$ in $x_1,\ldots,x_K$, through the formula
\begin{equation}
\label{eq_x233}
P(\x)=-\frac{1}{2}\Qu(\x-\N \hat{\boldsymbol n}^{\boldsymbol{\mu}}, \x-\N \hat{\boldsymbol n}^{\boldsymbol{\mu}}) + \sum_{k=1}^K x_k \big(\theta_{h^k,h^k} \log \N + \mathbbm{Shift}_k\big).
\end{equation}
The constants $\mathbbm{Shift}_k$ in this formula are defined as partial derivatives in filling fractions of the sum of two terms: first, the remainder $\mathbbm{Rest}_1$ in \eqref{eq_partition_multicut} for the expansion of the partition function of the discrete ensemble localized to intervals $[\alpha_k-\eps,\beta_k+\eps]$ and with fixed $K$ filling fractions on these segments;\footnote{Two remarks are in order. First, the endpoints of the $k$-th segment are not precisely $[\alpha_k-\eps,\beta_k+\eps]$, but small modifications of those, so that the integrality conditions of Section~\ref{Section_configuration_space} (or of Theorem~\ref{proposition_fluct_conditioning}) are valid. Second, we did not present any simple expressions for $\mathbbm{Rest}_1$ in \eqref{eq_partition_multicut}, and therefore, we do not have such expressions for the constants $\mathbbm{Shift}_k$, their definition remains implicit.} second, the sum of the terms $\mathcal{I}'$ and $\mathcal{I}''$ of Proposition~\ref{proposition_Energy_series} applied with $\epsilon = \frac{1}{\N}$ and used to compare the energy functionals of the original ensemble and the localized one\footnote{The difference between the functionals arises from the (subleading as $\N \rightarrow \infty$) change in the potential in the localization procedure, itself arising from our treatment of Gamma functions in the weights, \textit{cf.} \eqref{eq_weight_subleading_3}.}. In consequence, $\mathbbm{Shift}_k$ is bounded as $\N$ becomes large and the upper bound only depends on the constants in the assumptions.

By completing the square, the same function $P(\x)$ can be re-expressed as
\begin{equation}
\label{eq_x234}
 P(\x)=-\frac{1}{2}\Qu(\x-\boldsymbol u, \x-\boldsymbol u) + c,
\end{equation}
where $c$ does not depend on $\x$. The determination of $\boldsymbol u$ through \eqref{eq_x233}-\eqref{eq_x234} may not be unique, because in this identity $\x$ ranges only over $\L$ and not over $\amsmathbb{R}^K$. In particular, identifying $\Qu$ with a $K\times K$ matrix, which is invertible when restricted to $\L_0$, and letting $\boldsymbol{\Pi}$ denote the orthogonal (with respect to the standard scalar product in $\amsmathbb R^K$) projection operator onto $\L_0$, we can have $\boldsymbol u$ explicitly given by
\begin{equation}
\label{eq_discrete_shift}
 \boldsymbol u=\N \hat{\boldsymbol n}^{\boldsymbol{\mu}} + \Qu^{-1}\boldsymbol{\Pi}\left(\theta_{h^k,h^k} \log \N + \mathbbm{Shift}_k\right)_{k=1}^K
\end{equation}
Both $\Qu$ and $\boldsymbol u$ may depend on $\N$; in particular they clearly depend on the equilibrium measure. The center $\boldsymbol u$ of the discrete Gaussian distribution always contains a term of order $\N$ coming from the filling fractions of the equilibrium measure and, in general, it also contains a correction of order $\log \N$. The latter becomes absent if all $\theta_{h,h}$ are equal for $h \in [H]$, because $(1,\ldots,1)$ is always orthogonal to $\L_0$, corresponding to the fact that the total number of particles $N_1^{\circ}+\cdots+N_K^{\circ}$ is always deterministically fixed by the constraints \eqref{eq_equations_eqs}.

\subsection{Results on the fluctuations of filling fractions and linear statistics}

Denoting $\textnormal{\textsf{\textbf{Gau\ss{}}}}_{\amsmathbb{Z}} := \textnormal{\textsf{\textbf{Gau\ss{}}}}_{\amsmathbb{Z}}[\Qu,\L,\boldsymbol u]$ the discrete Gaussian random variable with the parameters described in Section~\ref{Sec:gausds}, we are now in position to describe the asymptotic distribution of the filling fractions.

\begin{theorem} \label{Theorem_CLT_for_filling_fractions} If the discrete ensemble satisfies Assumptions~\ref{Assumptions_Theta}, \ref{Assumptions_basic}, \ref{Assumptions_offcrit}, \ref{Assumptions_analyticity} and its equilibrium measure has no saturations, then we have an asymptotic equality of distributions as $\N \rightarrow \infty$
\begin{equation}
\label{eq_discrete_filling_fractions_limit} \boldsymbol{N}^{\circ} \stackrel{\textnormal{d}}{\sim} \textnormal{\textsf{\textbf{Gau\ss{}}}}_{\amsmathbb{Z}},
\end{equation}
where the convergence in the sense of Definition~\ref{Definition_asymptotic_equivalence} is uniform for fixed constants in the assumptions.
\end{theorem}

Combining Theorem~\ref{Theorem_CLT_for_filling_fractions} with our previous results on the fixed filling fractions case, we also obtain an asymptotic theorem on the fluctuations of linear statistics. Take a holomorphic function $f(z)$ defined in a complex neighborhood of the $k$-th band $(\alpha_k,\beta_k)$ and consider the random variable:
\begin{equation}
\label{eq_x235}
\mathsf{Lin}_{k}[f] = \sum_{i = 1}^{N} \mathbbm{1}_{[\N(\alpha_k - \eps),\N(\beta_k + \eps)]}(\ell_i)\,f\bigg(\frac{\ell_i}{\N}\bigg).
\end{equation}
 If $\eps$ is small and $\N$ is large, then Theorem~\ref{Theorem_ldpsup} shows that with overwhelming probability, the random variable of \eqref{eq_x235} does not depend on the choice of $\eps$, because all the particles $\ell_i$ remain in small neighborhoods of bands. We are interested in the fluctuations of \eqref{eq_x235} as $\N\rightarrow\infty$. They naturally split into two components: fluctuations of the segment filling fractions $\boldsymbol N^\circ$ leading to the varying number of non-zero terms in \eqref{eq_x235}, and fluctuations of the sum itself conditionally to the fixed filling fractions. More precisely, recalling that $\boldsymbol{\mu}^{\hat{\boldsymbol{n}}^{\circ}}$ is the equilibrium measure for the variational datum as in the original ensemble, but with segment filling fractions fixed to be $\hat{\boldsymbol{n}}^{\circ} = \frac{\boldsymbol{N}^{\circ}}{\N}$, we decompose \eqref{eq_x235} as a sum of two terms:
 \begin{equation*}
 \begin{split}
 \mathsf{Lin}_k^{(1)}[f] & = \sum_{i = 1}^{N} \mathbbm{1}_{[\N(\alpha_k - \eps),\N(\beta_k + \eps)]}(\ell_i)\,f\bigg(\frac{\ell_i}{\N}\bigg) -\N \int_{\alpha_k}^{\beta_k} f(x) \mu_{h^k}^{\hat{\boldsymbol{n}}^{\circ}}(x) \dd x, \\
 \mathsf{Lin}_k^{(2)}[f] & =  \N \int_{\alpha_k}^{\beta_k} f(x) \mu^{\hat{\boldsymbol{n}}^{\circ}}_{h^k}(x) \dd x \\
 \end{split}
 \end{equation*}
 The random variable $(\mathsf{Lin}_k^{(1)}[f])_{k = 1}^{K}$ is asymptotically Gaussian by Corollary~\ref{Corollary_CLT_relaxed} and we will see that it is asymptotically independent of $(\mathsf{Lin}_k^{(2)}[f])_{k = 1}^{K}$. The variable $\mathsf{Lin}_k^{(2)}[f]$ is random because $\boldsymbol{N}^{\circ}$ is a random variable whose asymptotic distribution is discrete Gaussian random vector as already described in Theorem~\ref{Theorem_CLT_for_filling_fractions}. Accordingly, the asymptotic distribution of $\mathsf{Lin}_k^{(2)}[f]$ will be described by probing this discrete Gaussian random vector in a certain direction, which is the subject of the following definition.

\begin{definition}
\label{Definition_omega}
For each $k \in [K]$ and function $f$ defined a neighborhood of the $k$-th band, we introduce the deterministic $K$-dimensional vector
\[
\boldsymbol{\omega}_{k}[f] = \left(\int_{\alpha_k}^{\beta_k} f(x) \big(\partial_{\hat{p}_l} \mu_{h^k}^{\hat{\boldsymbol{p}}}(x)\big)|_{\hat{\boldsymbol{p}} = \hat{\boldsymbol{n}}^{\boldsymbol{\mu}}}\, \dd x\right)_{l = 1}^{K}.
\]
In line with the notations of Chapter~\ref{Chapter_fff_expansions}, if $\boldsymbol{v}$ is another $K$-dimensional vector we will use the notation
\[
\langle \boldsymbol{v} \cdot \boldsymbol{\omega}_{k} \rangle = \sum_{l = 1}^{K} v_l\, \omega_{k,l}.
\]
\end{definition}
In this formula appear the partial derivatives of the density of the equilibrium measure with respect to filling fractions, evaluated at the deterministic value $\hat{\boldsymbol{n}}^{\boldsymbol{\mu}}$. We have already been working such derivatives in the proof of Proposition~\ref{Proposition_Hessian_free_energy}.

\begin{theorem}
\label{Theorem_linear_statistics_fluctuation_ff}
 Consider a discrete ensemble satisfying Assumptions~\ref{Assumptions_Theta}, \ref{Assumptions_basic}, \ref{Assumptions_offcrit}, \ref{Assumptions_analyticity} and suppose that the equilibrium measure has no saturations. Let $\eps>0$ be small enough, so that the $\eps$-neighborhoods of the bands are pairwise disjoint. Let $L \in \amsmathbb{Z}_{> 0}$ and a $L$-tuple of integers $\boldsymbol{k} \in [H]^{L}$ both independent of $\N$, let a (possibly $\N$-dependent) $L$-tuple of functions $\boldsymbol{f}(z)$ such that $f_l(z)$ is a holomorphic function of $z$ in a $\N$-independent complex neighborhood of $[\alpha_{k_l}-\eps,\beta_{k_l}+\eps]$ for any $l \in [L]$. Assume there exists a constant $C > 0$ such that $\max_{l} \sup_{z} |f_l(z)| \leq C$.

 Let $\textnormal{\textsf{\textbf{Gau\ss{}}}}[\boldsymbol{f}]$ be a $L$-dimensional random Gaussian vector with covariance
 \begin{equation}
\label{eq_x240}
\textnormal{\textsf{Cov}}_{l_1,l_2}[\boldsymbol{f},\boldsymbol{f}] = \oint_{\gamma_{k_{l_1}}} \oint_{\gamma_{k_{l_2}}}
 f_{l_1}(z_1) f_{l_2}(z_2)\,\mathcal{F}_{k_{l_1},k_{l_2}}(z_1,z_2)\,\frac{\dd z_1 \dd z_2}{(2\ii\pi)^2},
\end{equation}
 and mean equal to the sum of $\mathbbm{Rest}[f_l]$ from \eqref{eq_x164} and $\mathbbm{Rest}'[f_l]$ from \eqref{eq_x257}\footnote{The appearance of this term is caused by the centering with respect to localized equilibrium measure in Corollary~\ref{Corollary_CLT_relaxed} versus the equilibrium measure of the original ensemble in our treatment through \eqref{eq_x256}. The term is just the difference of two centerings.}, and independent of the $K$-dimensional discrete Gaussian vector $\textnormal{\textsf{\textbf{Gau\ss{}}}}_{\amsmathbb{Z}}$ appearing in Theorem~\ref{Theorem_CLT_for_filling_fractions}.

 Then, we have as $\N \rightarrow \infty$
 \begin{equation}
 \label{Lin1km}
\left( \textnormal{\textsf{Lin}}_{k_l}[f_l] - \N \int_{\alpha_{k_l}}^{\beta_{k_l}} f_l(x)\mu_{h^{k_l}}(x)\dd x \right)_{l = 1}^{L} \,\,\stackrel{\textnormal{d}}{\sim} \,\,\textnormal{\textsf{\textbf{Gau\ss{}}}}[\boldsymbol{f}] + \Big( \big\langle (\textnormal{\textsf{\textbf{Gau\ss{}}}}_{\amsmathbb{Z}} - \N \hat{\boldsymbol{n}}^{\boldsymbol{\mu}}) \cdot \boldsymbol{\omega}_{k_l}[f_l] \big\rangle \Big)_{l = 1}^{L},
 \end{equation}
where the convergence in the sense of Definition~\ref{Definition_asymptotic_equivalence} is uniform for fixed constants in the assumptions.
 \end{theorem}

The right-hand side of \eqref{Lin1km} contains an implicit shift in both terms: the mean for the Gaussian part, and the bounded shift $(\mathbbm{Shift}_k)_{k=1}^K$ entering into the definition \eqref{eq_discrete_shift} of the center $\boldsymbol u$ for the discrete Gaussian part which affects $\boldsymbol{\omega}_{k}$. It would be interesting to find explicit formulae for these shifts.

 Taking $f_l(x)=1$ for any $l \in [L]$, we recover the fluctuations of the filling fractions of Theorem~\ref{Theorem_CLT_for_filling_fractions}. Indeed, the Gaussian component in that case $\textnormal{\textsf{\textbf{Gau\ss{}}}}[\boldsymbol{f}]$ is deterministically zero: it has zero covariance because $\mathcal{F}_{l_1,l_2}(z_1,z_2) = O(\frac{1}{z_1^2z_2^2})$ as $z_1,z_2 \rightarrow \infty$, and one can check from the definition of the two parts of the mean that it vanishes. Furthermore, $\boldsymbol{\omega}_{k_l}[f_l]$ is the vector with $k_l$-th entry $1$ and other entries zero. Hence, the $l$-th entry of the second term in the right-hand side of \eqref{Lin1km} is $\textnormal{\textsf{Gau\ss{}}}_{\amsmathbb{Z},k_l} - \N \hat{n}^{\boldsymbol{\mu}}_{k_l}$.

It is instructive to compare Theorem~\ref{Theorem_linear_statistics_fluctuation_ff} (as well as its extension, Theorem~\ref{Theorem_linear_statistics_fluctuation_sat} below) with Corollary~\ref{Corollary_CLT_relaxed}. The latter had deterministically fixed filling fractions and the answer there already contains the same Gaussian part --- although we describe the mean differently in the two statements. There is however a new part given by the discrete Gaussian.

\subsection{Proof of Theorem~\ref{Theorem_CLT_for_filling_fractions}}
\label{Secrigub}
The concentration inequality for filling fractions in Corollary~\ref{Corollary_a_priori_0} implies the existence of a constant $C > 0$ depending only on the constants in the assumptions, such that
\begin{equation}
\label{eq_x244}
\forall k \in [K] \qquad \amsmathbb{P}\big[ |N^{\circ}_k - \N \hat{n}^{\boldsymbol{\mu}}_{k}| \geq C\N^{\frac{1}{2}} \log \N\big] \leq \exp\bigg(- \frac{\N (\log \N)^2}{C}\bigg).\end{equation}
Therefore, to prove the statement \eqref{eq_discrete_filling_fractions_limit}, we can, and we will, restrict to studying the distribution of the vector $\boldsymbol{N}^{\circ}$ in the region where
\begin{equation}
\label{eq_x242}
\forall k \in [K] \qquad \big| N^{\circ}_k - \N \hat{n}^{\boldsymbol{\mu}}_{k}\big| > C \N^{\frac{1}{2}}\log \N,
\end{equation}
This implies $|\!|\hat{\boldsymbol{n}}^{\circ} - \hat{\boldsymbol{n}}^{\boldsymbol{\mu}}|\!|_{\infty} =O(\N^{-\frac{1}{4}})$.

\medskip

\noindent \textsc{Step 1: Choosing localization segments.} We would like to choose pairwise disjoint segments $[\hat a^{\circ}_k, \hat{b}^{\circ}_k]$ indexed by $k \in [K]$, which contain $[\alpha_k,\beta_k]$ in their interior, and then localize the ensemble to these segments so as to apply the results of Chapter~\ref{Chapter_partition_functions}. Although we may desire $[\hat a^{\circ}_k, \hat{b}^{\circ}_k]$ to be as close as possible to $(\alpha_k-\eps,\beta_k+\eps)$, there is a difficulty: because the particles in our discrete ensembles belong to complicated lattices that depend on the location of the particles to their left and satisfy the integrality conditions of Section~\ref{Section_configuration_space}, the segments $[\hat a^{\circ}_k, \hat{b}^{\circ}_k]$ should typically depend on the fluctuating filling fractions $\boldsymbol{N}^{\circ}$.

The state space of Section~\ref{Section_configuration_space} is slightly different depending on whether $a_1$ and $b_1$ are finite or infinite. We first exclude the situation where $a_1=-\infty$ and $b_1$ is finite, and later comment how the procedure should be adjusted to handle it as well. We introduce the set
\begin{equation}
\label{Lbbhatdef}
\hat{\amsmathbb{L}} = \begin{cases} \hat{a}_1 + \N^{-1}\amsmathbb{Z}_{\geq 0}, & \textnormal{if } a_1 \textnormal{ is finite},\\[2pt] \N^{-1}\amsmathbb{L},& \textnormal{if }a_1=-\infty\textnormal{ and }b_1=+\infty.
\end{cases}
\end{equation}
where $\amsmathbb{L}$ is the lattice from Section~\ref{Section_configuration_space}. We note that the rescaled leftmost particle $\frac{\ell_1}{\N}$ belongs to $\hat{\amsmathbb{L}}$ and choose the leftmost endpoint:
\begin{equation}
\label{ahatcircdef}
 \hat a^{\circ}_1 = \min\big\{x\in \hat{\amsmathbb{L}} \quad \big| \quad\, x>\alpha_1-\eps\big\}.
\end{equation}
In addition, we choose $2K-1$ positive real numbers $\mathfrak{t}^{\textnormal{b}}_1,\ldots,\mathfrak{t}^{\textnormal{b}}_K,\mathfrak{t}^{\textnormal{a}}_2,\ldots,\mathfrak{t}^{\textnormal{a}}_K$, which do not depend on $\boldsymbol{N}^{\circ}$ but may depend on $\N$. Then, we define by induction on $k$ all the remaining endpoints of $[\hat a^{\circ}_k, \hat{b}^{\circ}_k]$ as affine functions of these real numbers and the segment filling fractions:
\begin{equation}
\label{eq_x241}
 \hat{b}^{\circ}_k:= \hat a^{\circ}_k+\theta_{h^k,h^k} \frac{N^{\circ}_k-1}{\N} + \mathfrak{t}^{\textnormal{b}}_k, \qquad \hat{a}^{\circ}_{k+1}:= \hat{b}^{\circ}_k+ \mathfrak{t}^{\textnormal{a}}_{k+1}.
\end{equation}
In order for the choice \eqref{eq_x241} to be meaningful and helpful for us, we impose the following restrictions for any $k \in [K]$ and $l \in [K - 1]$.\begin{enumerate}
 \item $\N \mathfrak{t}^{\textnormal{b}}_k\in \amsmathbb Z_{\geq 0}$.
 \item If $h^{l+1}=h^l$, then $\N \mathfrak{t}^{\textnormal{a}}_{l+1}\in \amsmathbb{Z}_{\geq 0}$.
 \item If $h^{l+1}=h^l+1$, then $\N \mathfrak{t}^{\textnormal{a}}_{l+1} \in \N(\hat a_{h^{l+1}}-\hat{b}^{\circ}_l)+\amsmathbb{Z}_{\geq 0}$.
 \item For $\N$ large enough and any $\boldsymbol{N}^{\circ}$ satisfying \eqref{eq_x242}, the choices \eqref{eq_x241} are such that
 \begin{equation}
 \label{eq_x243}
 \hat a^{\circ}_k\in \bigg(\alpha_k-\eps,\alpha_k-\frac{\eps}{2}\bigg),\quad \textnormal{ and } \quad \hat b^{\circ}_k\in \bigg(\beta_k+\frac{\eps}{2},\beta_k+\eps\bigg).
 \end{equation}
\end{enumerate}
The first three restrictions are integrality constraints making sure that the localization of the discrete ensemble on the segments $[\hat a^{\circ}_k, \hat{b}^{\circ}_k]$ makes sense, \textit{cf.} the situation without saturations in Theorem~\ref{proposition_fluct_conditioning}. In more detail, the first restriction is a translation of the relation between the number of particles in $[\hat a^\circ_k, \hat b^{\circ}_k]$ and possible values for $\hat a^\circ_k$, $\hat b^\circ_k$; the second restriction is a translation of the fact that all the available sites between $\hat b^\circ_l$ and $\hat a^\circ_{l+1}$ should be filled with holes (\textit{cf.} Section~\ref{Sec : large_dev_saturated}) and the distance between two adjacent holes is $\frac{1}{\N}$; the third condition is a translation of the fact that all the available sites between $\hat a_{h_{l+1}}$ and $\hat a^\circ_{l+1}$ should be filled with holes at distance $\frac{1}{\N}$ from each other. The fourth restriction is the desired closeness of $[\hat a^{\circ}_k, \hat{b}^{\circ}_k]$ to $(\alpha_k-\eps,\beta_k+\eps)$. It is straightforward to check that for large $\N$ one can always find $\mathfrak{t}^{\textnormal{b}}_1,\ldots,\mathfrak{t}^{\textnormal{b}}_K,\mathfrak{t}^{\textnormal{a}}_2,\ldots,\mathfrak{t}^{\textnormal{a}}_K$ satisfying the above four conditions.

If $a_1=-\infty$ and $b_1$ is finite, then we amend the procedure for the definition of the segments $[\hat{a}_k^{\circ},\hat{b}_k^{\circ}]$ for $k \leq k^+(1)$, \textit{i.e.} containing in $[\hat{a}_1,\hat{b}_1]$. We start by setting
\[
\hat{b}^{\circ}_{k^+(1)} =\max\big\{x\in \hat b_1-\N^{-1}\amsmathbb{Z}_{\geq 0}\,\,\big| \,\,x<\beta_{k^+(1)}+\eps\big\},
\]
and then use the identities \eqref{eq_x241} backwards to define successively $\hat{a}_{k^+(1)}^{\circ},\hat{b}_{k^+(1) - 1}^{\circ},\ldots,\hat{b}_1^{\circ},\hat{a}_1^{\circ}$. The rest of the procedure for $k > k^+(1)$ remains unchanged.

\smallskip

\noindent \textsc{Step 2: Distribution of filling fractions in the localized ensemble.} By definition, the probability distribution of $\boldsymbol{N}^{\circ}$ can be computed as:
\begin{equation}
\label{eq_x245}
 \P(\boldsymbol{N}^{\circ}=\boldsymbol M)= \frac{\Z^{\boldsymbol M}_\N}{\Z_\N},
\end{equation}
where $\Z_\N^{\boldsymbol M}$ is the partition function of the ensemble restricted to the configurations satisfying $\boldsymbol{N}^{\circ}=\boldsymbol M$, \textit{i.e.} it is the part of the sum \eqref{eq_partition_function_definition} corresponding to such configurations. By Theorem~\ref{Theorem_ldpsup}, with overwhelming probability there are no particles outside ${\bigcup_{k=1}^K [\N(\alpha_k-\eps),\N(\beta_k+\eps)]}$; combining this with \eqref{eq_x244}, we conclude that up to multiplication by a negligible error term $\big(1+ O(e^{-\frac{\N}{C} (\log \N)^2})\big)$, the partition function $\Z_\N$ is the sum of $\Z_\N^{\boldsymbol M}$ over all possible choices of $\boldsymbol M$ satisfying \eqref{eq_x242} and the affine constraints \eqref{eq_equations_eqs} explained after \eqref{eq_x181}.

We compute $\Z^{\boldsymbol M}_\N$ in \eqref{eq_x245} for $\boldsymbol M$ satisfying \eqref{eq_x181} and \eqref{eq_x242} by localizing the ensemble to the segments $[\hat a^{\circ}_k, \hat{b}^{\circ}_k]$ for $k\in [K]$, and fixing simultaneously $\boldsymbol{N}^{\circ}=\boldsymbol M$, as in Theorem~\ref{proposition_fluct_conditioning}. Equation \eqref{eq_x181} implies that there are no particles outside ${\bigcup_{k=1}^K [\N(\alpha_k-\eps),\N(\beta_k+\eps)]}$ and therefore the partition function $\Z^{\boldsymbol M}_\N$ is unchanged in the localization procedure. Therefore, we can apply Theorem~\ref{Theorem_partition_multicut} to conclude that as $\N\rightarrow\infty$
\begin{equation}
\label{eq_x246}
 \Z^{\boldsymbol M}_\N=\exp\Bigg( \I^{\textnormal{loc}}[\boldsymbol{\mu}^{\hat{\boldsymbol m},\textnormal{loc}} ]\, \N^2  + \bigg(\sum_{k=1}^K \theta_{h^k,h^k} \hat m_k \bigg)\N \log \N + \mathbbm{Rest}_1(\hat{\boldsymbol{m}},\hat{\boldsymbol{a}}^{\circ},\hat{\boldsymbol{b}}^{\circ})\N +o(1)\Bigg),
\end{equation}
where we insisted on the dependence of $\mathbbm{Rest}_1$ on $\hat{\boldsymbol{m}} = \frac{\boldsymbol{M}}{\N}$ and the endpoints $\hat{\boldsymbol{a}}^{\circ},\hat{\boldsymbol{b}}^{\circ}$. Besides $-\I^{\textnormal{loc}}$ is the energy functional for the localized ensemble and $\boldsymbol{\mu}^{\hat{\boldsymbol m},\textnormal{loc}}$ is the corresponding equilibrium measure constrained to give mass $\hat{m}_k$ to the segment $[\hat{a}_k^{\circ,\prime},\hat{b}_k^{\circ,\prime}]$ for any $k \in [K]$. As usual, the prime notation include the shifts we always include when dealing with a variational datum and the corresponding equilibrium measure
\begin{equation}
\label{shiftedcirque}
\forall k \in [K] \qquad \hat{a}_k^{\circ\,\prime} = \hat{a}_k^{\circ} - \frac{\theta_{h^k,h^k} - \frac{1}{2}}{\N},\qquad \hat{b}_k^{\circ\,\prime} = \hat{b}_k^{\circ} + \frac{\theta_{h^k,h^k} - \frac{1}{2}}{\N}.
\end{equation}
In the present situation, since there are no saturations, the bands are surrounded by voids and equivalently $\hat{m}_k = \mu^{\hat{\boldsymbol{m}},\textnormal{loc}}([\hat{a}_k^{\circ},\hat{b}_k^{\circ}])$. For the same reason, the dependence of $\boldsymbol{\mu}^{\hat{\boldsymbol{m}},\textnormal{loc}}$ and $\mathbbm{Rest}_1$ in the endpoints $(\hat{\boldsymbol{a}}^{\circ},\hat{\boldsymbol{b}}^{\circ})$ is virtual, rather than meaningful: by Theorem~\ref{Theorem_ldpsup}, with overwhelming probability there are no particles near $\hat a^{\circ}_k$, $\hat{b}^{\circ}_k$, and therefore slightly changing $\hat a^{\circ}_k$, $\hat{b}^{\circ}_k$ --- \textit{e.g.} by modifying the choices of the numbers $\mathfrak{t}^{\textnormal{b}}_1,\ldots,\mathfrak{t}^{\textnormal{b}}_K,\mathfrak{t}^{\textnormal{a}}_2,\ldots,\mathfrak{t}^{\textnormal{a}}_K$ involved in their definition --- leads to an unchanged $\I^{\textnormal{loc}}[\boldsymbol{\mu}^{\hat{\boldsymbol{m}},\textnormal{loc}} ]$, while the changes in $\mathbbm{Rest}_1$ are negligible in the sense that they can be included in the $o(1)$ remainder term in \eqref{eq_x246}.

It is inconvenient that $\I^{\textnormal{loc}}[\boldsymbol{\mu}^{\hat{\boldsymbol{m}},\textnormal{loc}}]$ depends on the localization procedure and we would like to replace $\I^{\textnormal{loc}}$ in \eqref{eq_x246} by the energy functional $-\I$ of the original ensemble. Looking into the localization in Theorem~\ref{proposition_FFF_conditioning}, we see that the difference between $\I$ and $\I^{\textnormal{loc}}$ comes from the subleading $O(\frac{1}{\N})$ change in potential due to our treatment of Gamma functions, \textit{c.f} \eqref{eq_weight_subleading_3}. Using Proposition~\ref{proposition_Energy_series} with $t= \frac{1}{\N}$, we have
\begin{equation}
\label{eq_x254}
 \I^{\textnormal{loc}}[\boldsymbol{\mu}^{\hat{\boldsymbol{m}},\textnormal{loc}}]\,\N^2=\I[\boldsymbol{\mu}^{\hat{\boldsymbol{m}}}]\,\N^2+ \I'[\boldsymbol{\mu}^{\hat{\boldsymbol{m}}}]\,\N + \I''[\boldsymbol{\mu}^{\hat{\boldsymbol{m}}}]
+ o(1),
\end{equation}
where $\boldsymbol{\mu}^{\hat{\boldsymbol{m}}}$ is the minimizer of $-\I$ subject to the filling fractions constraints $\boldsymbol{\mu}^{\hat{\boldsymbol m}}([\hat a^{'\circ}_k, \hat{b}^{'\circ}_k])=\hat{m}_k$ for any $k\in [K]$, and the terms $\I'[\boldsymbol{\mu}^{\hat{\boldsymbol{m}}}], \I''[\boldsymbol{\mu}^{\hat{\boldsymbol{m}}}]$ depend smoothly on the rescaled filling fractions $\hat{\boldsymbol{m}}$.

\medskip

\noindent \textsc{Step 3: Taylor expansion near equilibrium filling fractions.} With help of Proposition~\ref{Proposition_Hessian_free_energy}, we perform the second-order Taylor expansion of $\I[\boldsymbol{\mu}^{\hat{\boldsymbol{m}}}]$ in $\hat{\boldsymbol{m}}$ near the point $\hat{\boldsymbol n}^{\boldsymbol{\mu}}$ given by \eqref{eq_def_equilibrium_ff}. The linear terms necessarily vanish because the equilibrium measure $\boldsymbol{\mu}$ at the filling fractions $\hat{\boldsymbol n}^{\boldsymbol{\mu}}$ is the minimizer of $-\I$, and the matrix $\Qu$ of \eqref{eq_discrete_covariance} gives the opposite of the Hessian. We get
\begin{equation}
\label{eq_x247}
 \I[\boldsymbol{\mu}^{\hat{\boldsymbol{m}}} ]\, \N^2 = \I[\boldsymbol{\mu} ]\, \N^2 -\frac{1}{2}\Qu({\boldsymbol M} - \N\hat{\boldsymbol n}^{\boldsymbol{\mu}},{\boldsymbol M} - \N\hat{\boldsymbol n}^{\boldsymbol{\mu}})+o(1).
\end{equation}
The first term in the right-hand side of \eqref{eq_x247} will eventually cancel between the numerator and the denominator of \eqref{eq_x245} and therefore can be ignored. We emphasize that the second term involves the integral vector ${\boldsymbol M}$, rather than its rescaled version $\hat{\boldsymbol{m}}$. Simultaneously, we expand the terms $\I'$ and $\I''$ of \eqref{eq_x254} near $\hat{\boldsymbol n}^{\boldsymbol{\mu}}$, recalling that $\boldsymbol{\mu}^{\hat{\boldsymbol{n}}^{\boldsymbol{\mu}}} = \boldsymbol{\mu}$
\begin{equation}
\label{eq_x255}
 \I'[\boldsymbol{\mu}^{\hat{\boldsymbol{m}}}]\,\N + \I''[\boldsymbol{\mu}^{\hat{\boldsymbol{m}}}] = \I'[\boldsymbol{\mu}]\,\N + \I''[\boldsymbol{\mu}] + \sum_{k=1}^K (M_k - \N\hat{n}^{\boldsymbol{\mu}}_k) \cdot \big(\partial_{\hat{m}_k} \I'[\boldsymbol{\mu}^{\hat{\boldsymbol{m}}}]|_{\hat{\boldsymbol{m}} = \hat{\boldsymbol{n}}^{\boldsymbol{\mu}}}\big) +o(1),
\end{equation}
The first two terms in the right-hand side of \eqref{eq_x255} will eventually cancel between the numerator and the denominator of \eqref{eq_x245}. Next, we transform the second term in the exponent in \eqref{eq_x246}:
\begin{equation}
\label{eq_x248}
 \bigg(\sum_{k=1}^K \theta_{h^k,h^k} \hat{m}_k \bigg)\N \log \N = \bigg(\sum_{k=1}^K \theta_{h^k,h^k}  \hat{n}^{\boldsymbol{\mu}}_k \bigg)\, \N \log \N + \bigg(\sum_{k=1}^K \theta_{h^k,h^k} (M_k-\N \hat{n}^{\boldsymbol{\mu}}_k )\bigg)\log \N.
\end{equation}
The first term in the right-hand side of \eqref{eq_x248} cancels between the numerator and the denominator of \eqref{eq_x245} and can be ignored. Finally, we write
\begin{equation}
\label{eq_x249}
\begin{split}
 & \quad \mathbbm{Rest}_1(\boldsymbol{\hat{m}},\hat{\boldsymbol{a}}^{\circ},\hat{\boldsymbol{b}}^{\circ})\,\N \\
 & = \mathbbm{Rest}_1(\hat{\boldsymbol{n}}^{\boldsymbol{\mu}},\hat{\boldsymbol{a}}^{\boldsymbol{\mu}}, \hat{\boldsymbol{b}}^{\boldsymbol{\mu}})\,\N +\sum_{k=1}^K (M_k-\N \hat{n}^{\boldsymbol{\mu}}_k ) \big(\partial_{\hat{m}_k} \mathbbm{Rest}_{1}(\boldsymbol{\hat{m}},\hat{\boldsymbol{a}}^{\circ},\hat{\boldsymbol{b}}^{\circ})\big)\big|_{\hat{\boldsymbol{m}} = \hat{\boldsymbol{n}}^{\boldsymbol{\mu}}} + o(1),
 \end{split}
\end{equation}
where in the first term $(\hat{\boldsymbol{a}}^{\boldsymbol{\mu}},\hat{\boldsymbol{b}}^{\boldsymbol{\mu}})$ are obtained as $(\hat{\boldsymbol{a}}^{\circ},\hat{\boldsymbol{b}}^{\circ})$ through the formulae \eqref{eq_x241}, but with $\boldsymbol{N}^{\circ}$ replaced with $\N \hat{\boldsymbol{n}}^{\boldsymbol{\mu}}$. Note that the derivatives with respect to the endpoints of the segments do not appear, because, as we have already mentioned, $\mathbbm{Rest}_1(\boldsymbol{\hat{m}},\hat{\boldsymbol{a}}^{\circ},\hat{\boldsymbol{b}}^{\circ})$ does not really depend on them. And, as before, the first term in the right-hand side of \eqref{eq_x248} cancels between the numerator and the denominator of \eqref{eq_x245}. Plugging \eqref{eq_x254}--\eqref{eq_x249} into \eqref{eq_x245}--\eqref{eq_x246}, we get as $\N\rightarrow\infty$
\begin{equation}
\label{eq_x250}
\begin{split}
& \quad \P(\boldsymbol{N}^{\circ}=\boldsymbol M)\\
& = \frac{1}{\mathscr{Q}} \exp\bigg(-\frac{1}{2}\Qu({\boldsymbol M} - \N\hat{\boldsymbol n}^{\boldsymbol{\mu}} ,{\boldsymbol M} - \N\hat{\boldsymbol n}^{\boldsymbol{\mu}} )
 +\sum_{k=1}^K (\theta_{h^k,h^k} \log \N+\mathbbm{Shift}_k) (M_k-\N \hat{n}^{\boldsymbol{\mu}}_k )+o(1) \,
 \bigg),
\end{split}
\end{equation}
where $\mathscr{Q}$ is a normalization constant making the expression a probability distribution and
\[
\mathbbm{Shift}_k= \partial_{\hat{m}_k} \I'[\boldsymbol{\mu}^{\hat{\boldsymbol{m}}}]|_{\hat{\boldsymbol{m}} = \hat{\boldsymbol{n}}^{\boldsymbol{\mu}}} + \partial_{\hat{m}_k} \mathbbm{Rest}_{1}(\boldsymbol{\hat{m}},\hat{\boldsymbol{a}}^{\circ},\hat{\boldsymbol{b}}^{\circ})\big|_{\hat{\boldsymbol{m}} = \hat{\boldsymbol{n}}^{\boldsymbol{\mu}}},
\]
 while the error $o(1)$ is uniform in $\boldsymbol M$ satisfying \eqref{eq_x242}. Comparing with the statement of Theorem~\ref{Theorem_CLT_for_filling_fractions} and the definition of the shift $\boldsymbol u$ in \eqref{eq_x233}, \eqref{eq_x234} and noting that the closeness of probabilities implies the convergence in the sense of Definition~\ref{Definition_asymptotic_equivalence}, we are done.

\medskip

The reader might worry about the fact that we started by conditioning on the restriction \eqref{eq_x242}, and therefore \eqref{eq_x250} is also obtained only under this restriction. However, for the random variable before taking a $\N \rightarrow \infty$ limit, the inequality \eqref{eq_x242} implies that we condition on the event of a probability tending to $1$ as $\N \rightarrow \infty$. For the limiting discrete Gaussian of \eqref{eq_Discrete_Gaussian_main} the probability of the condition similarly tends to $1$ because the restriction $\Qu$ to $\L_0$ is positive definite.

\subsection{Proof of Theorem~\ref{Theorem_linear_statistics_fluctuation_ff}}

\label{sec:proofth105}

For notational simplicity we only consider the $L=1$ case; the general case is obtained by the same argument. Thus, we fix some $k \in [K]$ and a holomorphic function $f$ in a complex neighborhood of the $k$-th band, and study the asymptotics of the random variable
\[
\textsf{Lin}_k[f] = \sum_{i = 1}^{N} \mathbbm{1}_{[\N(\alpha_k - \varepsilon),\N(\beta_k + \varepsilon)]}(\ell_i)\,f\bigg(\frac{\ell_i}{\N}\bigg).
\]
The first step is to study the conditional distribution of $\textsf{Lin}_k[f]$ given the values of the filling fractions $\boldsymbol{N}^{\circ}$. We use the segment $[\hat a^{\circ}_k, \hat{b}^{\circ}_k]$ introduced in the course of the proof of Theorem~\ref{Theorem_CLT_for_filling_fractions}, and note that by Corollary~\ref{Corollary_a_priori_0} and Theorem~\ref{Theorem_ldpsup} it is sufficient to only deal with filling fractions satisfying \eqref{eq_x242}, that is
\[
|\!|\boldsymbol{N}^{\circ} - \N \hat{\boldsymbol{n}}^{\boldsymbol{\mu}}|\!|_{\infty} \leq C \N^{\frac{1}{2}}\log \N,
\]
 and to assume that there are no particles $\ell_i$ outside $\bigcup_{k=1}^K [\N \hat a^{\circ}_k, \N \hat{b}^{\circ}_k]$.

Let $\textsf{Lin}_k^{\hat{\boldsymbol{m}}}[f]$ denote the random variable $\textsf{Lin}_k[f]$ conditioned on the intersection of the event $\boldsymbol{N}^{\circ}=\boldsymbol M = \N \hat{\boldsymbol{m}}$, with $\boldsymbol M$ satisfying \eqref{eq_x242}, and of the event that all particles are inside $\bigcup_{k=1}^K [\N \hat a^{\circ}_k, \N \hat{b}^{\circ}_k]$. Localizing the ensemble to $\bigcup_{k = 1}^{K} [\hat a^{\circ}_k, \hat{b}^{\circ}_k]$ as in Theorem~\ref{proposition_FFF_conditioning}, we can use Corollary~\ref{Corollary_CLT_relaxed}. Hence,
\[
\textsf{Lin}_k^{\hat{\boldsymbol{m}}}[f] - \amsmathbb{E}\big[\textsf{Lin}_k^{\hat{\boldsymbol{m}}}[f]\big]
\]
is asymptotically Gaussian with variance as in \eqref{eq_x240}. Note that the dependence on $\hat{\boldsymbol{m}}$ in this part disappears as $\N\rightarrow\infty$, because the function $\mathcal{F}_{k,l}$ in \eqref{eq_x240} depends only on the endpoints of the bands and on the matrix $\boldsymbol{\Theta}$, and the dependence of the endpoints on $\hat{\boldsymbol{m}}$ is negligible due to \eqref{eq_x242} and Theorem~\ref{Theorem_differentiability_full}. This gives the random part of the first term $\textnormal{\textsf{\textbf{Gau\ss{}}}}[f]$ in \eqref{Lin1km} and also proves the independence of the two terms in \eqref{Lin1km}.

The deterministic shift (expectation) for the first term $\textnormal{\textsf{\textbf{Gau\ss{}}}}[f]$ in \eqref{Lin1km} and the full second term in \eqref{Lin1km} arises from the asymptotic behavior of the expectation value of $\textsf{Lin}_k^{\hat{\boldsymbol{m}}}[f]$. By \eqref{eq_x164}, we have
\begin{equation}
\label{eq_x252}
\E\big[\textsf{Lin}_{k}^{\hat{\boldsymbol{m}}}[f]\big] =\N
\int_{\hat a^{\circ\,\prime}_k}^{\hat{b}^{\circ\,\prime}_k} f(x)\,\mu_{h^k}^{\hat {\boldsymbol{m}},\textnormal{loc}}(x)\dd x + \mathbbm{Rest}[f] + O(\N^{-\frac{3}{2}+\eps}),
\end{equation}
Due to \eqref{eq_x242}, the dependence of $\mathbbm{Rest}[f]$ on $\boldsymbol{\hat{m}}$ becomes negligible as $\N\rightarrow\infty$. This term gives the first part of the expectation for the term $\textnormal{\textsf{\textbf{Gau\ss{}}}}[f]$ in \eqref{Lin1km}.

It is convenient to re-express the first term in \eqref{eq_x252} using the Stieltjes transform of the equilibrium measure:
\begin{equation}
 \N
\int_{\hat a^{\circ\,\prime}_k}^{\hat{b}^{\circ\,\prime}_k} f(x)\,\mu_{h^k}^{\hat {\boldsymbol{m}},\textnormal{loc}}(x)\dd x=\frac{\N}{2\pi \ii} \oint_{\gamma_{h^k}} f(z)
 \Gm_{\boldsymbol{\mu}^{\hat {\boldsymbol{m}},\textnormal{loc}}}(z)\dd z,
\end{equation}
Recall that $\boldsymbol{\mu}^{\hat{\boldsymbol{m}},\textnormal{loc}}$ is the equilibrium measure for the potential of the original ensemble modified by the addition of subleading $O(\frac{1}{\N})$ terms coming from our treatment of Gamma functions in the localization procedure, \textit{cf.} \eqref{eq_weight_subleading_3}. Using Theorem~\ref{Theorem_differentiability_full}, we re-express its Stieltjes transform in terms of the equilibrium measure $\boldsymbol{\mu}^{\hat {\boldsymbol{m}}}$ corresponding to the potential of the original ensemble and constrained to segment filling fractions $\mu^{\hat{\boldsymbol{m}}}([ \hat a^{\circ\,\prime}_k, \hat{b}^{\circ\,\prime}_k])=\hat{m}_k$ for any $k\in[K]$:
\begin{equation}
\label{eq_x256}
\N \oint_{\gamma_{h^k}} f(z)
 \Gm_{\mu^{\hat {\boldsymbol{m}},\textnormal{loc}}}(z) \frac{\dd z}{2\ii\pi} = \N \oint_{\gamma_{h^k}} f(z)
 \Gm_{\mu^{\hat {\boldsymbol{m}}}}(z) \frac{\dd z}{2\ii\pi} + \oint_{\gamma_{h^k}} f(z)
 \mathbbm{Err}^{\hat{\boldsymbol{m}}}(z)\frac{\dd z}{2\ii\pi} +o(1).
\end{equation}
We will not specify the term $\mathbbm{Err}^{\hat{\boldsymbol{m}}}(z)$, although some information on it can be extracted from the proof of Theorem~\ref{Theorem_differentiability_full}. The only relevant information for us is that it depends smoothly on the segment filling fractions $\hat {\boldsymbol{m}}$. Hence, using \eqref{eq_x242}, we rewrite the second term in the right-hand side of \eqref{eq_x256} as
\begin{equation}
\label{eq_x257}
 \mathbbm{Rest}'[f_m] :=  \oint f(z)
 \mathbbm{Err}^{\hat{\boldsymbol{n}}^{\boldsymbol{\mu}}}(z) \frac{\dd z}{2\ii\pi} +o(1).
\end{equation}
This is the term appearing in the second part of the expectation for $\textnormal{\textsf{\textbf{Gau\ss{}}}}[f]$ in \eqref{Lin1km}.

The first integral in the right-hand side of \eqref{eq_x256} is the origin of another term in \eqref{Lin1km}. Using Theorem~\ref{Theorem_differentiability_full} --- see also the proof of Proposition~\ref{Proposition_Hessian_free_energy} --- we perform a Taylor expansion of this integral in \eqref{eq_x252} around the point $\hat{\boldsymbol{m}} = \hat{\boldsymbol{n}}^{\boldsymbol{\mu}}$ which corresponds to the filling fractions of the unconstrained equilibrium measure. These filling fractions do not have to satisfy the integrality constraints, but this is not important for the Taylor expansion. In this way,
\begin{equation}
 \label{eq_x253}
\begin{split}
& \quad \N \oint_{\gamma_{h^k}} f(z)
 \Gm_{\mu^{\hat {\boldsymbol{m}}}}(z) \frac{\dd z}{2\ii\pi} \\
 &=\N \oint_{\gamma_{h^k}} f(z)
 \Gm_{\mu}(z) \frac{\dd z}{2\ii\pi}  + \sum_{l=1}^K (M_{l}-\N \hat n^{\boldsymbol{\mu}}_{l}) \oint_{\gamma_{h^k}} f(z) \big(\partial_{\hat{m}_l}\Gm_{\mu^{\hat{\boldsymbol{m}}}}(z)\big)\big|_{\hat{\boldsymbol{m}} = \hat{\boldsymbol{n}}^{\boldsymbol{\mu}}}\, \frac{\dd z}{2\ii\pi} +o(1)
\\
& = \N\int_{\alpha_k}^{\beta_k} f(x)\,\mu_{h^k}(x)\dd x  + \sum_{l=1}^K (M_{l}-\N \hat n^{\boldsymbol{\mu}}_{l})\int_{\alpha_k}^{\beta_k} f(x) \big(\partial_{\hat{m}_{l}} \mu_{h^k}^{\hat{\boldsymbol{m}}}(x)\big)|_{\hat{\boldsymbol{m}} = \hat{\boldsymbol{n}}^{\boldsymbol{\mu}}}\, \dd x +o(1).
\end{split}
\end{equation}
In the last equality we used that the measure $\mu_{h^k}$ is supported on $[\alpha_k,\beta_k]$. In the first integral in the last formula we recognize the subtracted term in the left-hand side of \eqref{Lin1km}. The sum over $l$ in the last formula --- which is a part of \eqref{eq_x252} --- is a conditional average, and for the final answer, we should replace $\boldsymbol{M}$ by the random variable $\boldsymbol{N}^{\circ}$ whose value it was representing. Applying Theorem~\ref{Theorem_CLT_for_filling_fractions}, this gives us an affine function of the random integer-valued variable $\textnormal{\textsf{\textbf{Gau\ss{}}}}_{\amsmathbb{Z}}$, matching the second term in \eqref{Lin1km}. This concludes the proof of Theorem~\ref{Theorem_linear_statistics_fluctuation_ff}

\section{In presence of saturations}

\label{Section_Fillingfractionsatcasesec}

\subsection{Informal presentation}

In this section we generalize the results of Section~\ref{Section_FFF_no_sat} to the case when the equilibrium measure may have saturated regions. Informally, we still want to understand the asymptotics of the total number of particles in each band, however, a new difficulty arises: if we follow the same definition \eqref{Nbullek} of the filling fractions, then $\boldsymbol N^\circ$ becomes significantly dependent on the choice of $\eps>0$. Indeed an increase in $\eps$ leads to an --- essentially deterministic --- increase in the filling fractions, because a larger part of the saturated region ends up being in $(\alpha_k-\eps,\beta_k+\eps)$. An additional delicacy is caused by the complicated lattice on which particles $\ell_i$ sit if some of the diagonal elements $\theta_{h,h}$ of the interaction matrix are not equal to $1$. We overcome these difficulties by modifying the definition of the discrete random variables of interest: rather than looking only at the filling fractions in fixed segments as in \eqref{Nbullek}, we simultaneously start treating some of the endpoints of the segments as being random and looking at random filling fractions in these random segments.

\medskip

Concretely, we replace the $K$-tuple of random integers $\boldsymbol N^{\circ}$ of Section~\ref{Section_FFF_no_sat} by an extended random $4K$-tuple
\begin{equation}
\label{exttupledef}\overline{\boldsymbol N}^{\circ}=(\hat{\boldsymbol{a}}^\circ, \hat{\boldsymbol{b}}^\circ, \boldsymbol{N}^\circ,\boldsymbol{Q}^\circ),
\end{equation}
Here $\boldsymbol{Q}^\circ$ is an auxiliary $K$-tuple of integers. The $4K$-tuple $\overline{\boldsymbol N}^{\circ}$ is a function of $\boldsymbol{\ell} \in \W_\N$ which satisfies the following properties on an event of overwhelming probability as $\N\rightarrow\infty$.
\begin{enumerate}
 \item The segments $[\hat a^\circ_k, \hat b^\circ_k]$ indexed by $k \in [K]$ are pairwise disjoint. They contain in their interior the closure of the band $(\alpha_k,\beta_k)$, and are contained in the segments of the original ensemble that contain those bands
 \[
 [\alpha_k,\beta_k] \subset (\hat{a}^{\circ}_k,\hat{b}_k^{\circ}) \subset (\hat{a}_{h^k},\hat{b}_{h^k}).
 \]
 \item For any $k \in [K]$, the integer $N^\circ_k$ is the total number of particles in the segment $[\N \hat a^\circ_k, \N\hat b^\circ_k]$. In addition, the following three numbers are nonnegative integers:
 \begin{itemize}
 \item $\N\hat b^\circ_k-\N\hat a^\circ_k-\theta_{h^k,h^k} N^\circ_k$;
 \item the distance between $\N \hat a^{\circ}_k$ and the leftmost particle in $[\N \hat a^\circ_k, \N\hat b^\circ_k]$;
 \item the distance between the rightmost particle in $[\N \hat a^\circ_k, \N\hat b^\circ_k]$ and $\N \hat b^{\circ}_k$.
 \end{itemize}
 \item The complement of $\bigcup_{k = 1}^{K} [\N \hat{a}^{\circ}_k,\N\hat{b}_k^{\circ}]$ in $\bigcup_{h = 1}^{H} [a_h,b_h]$ splits into finitely many segments, each of which either has no particles or has no holes --- as in Theorem~\ref{Theorem_ldsaturated}.
 \item The random $3K$-tuple $(\hat{\boldsymbol{a}}^\circ, \hat{\boldsymbol{b}}^\circ, \boldsymbol{N}^\circ)$ is a deterministic affine function of the random $K$-tuple $\boldsymbol{Q}^\circ$.
\end{enumerate}

Informally, $N^\circ_k$ are the same filling fractions as before, but they refer to how the now \emph{random} segment $[\hat{a}^{\circ}_k,\hat{b}^{\circ}_k]$ is filled. Note that $\hat a^\circ_k$ and $\hat b^\circ_k$ were already present in Section~\ref{Section_FFF_no_sat} as a necessary ingredient of the proofs; the difference is that now we need them already to state the asymptotic theorems. The $K$-tuple $\boldsymbol{Q}$ is a convenient tool for formulating and proving the statements about the remaining $3K$ components of $\overline {\boldsymbol N}^\circ$, which due to Condition 4. we can consider as driven by $\boldsymbol{Q}^{\circ}$.

\smallskip

We are going to produce exact analogues of Theorems~\ref{Theorem_CLT_for_filling_fractions} and \ref{Theorem_linear_statistics_fluctuation_ff}: as $\N\rightarrow\infty$, the random vector $\overline{\boldsymbol N}^\circ$ is asymptotically distributed like a discrete Gaussian random variable, and linear statistics are asymptotically described as a sum of an affine function of this discrete Gaussian random variable and of an independent Gaussian component previously computed in Corollary~\ref{Corollary_CLT_relaxed}. The exact form of these theorems inevitably depends on the exact choices we make in defining the $4K$-tuple $\overline{\boldsymbol N}^\circ$. In contrast to the simpler case described in Section~\ref{Section_FFF_no_sat}, these choices are not canonical in the situation when saturations are present. Our way to make the choices is described below.

We remark that if $\theta_{h,h}=1$ for any $h\in[H]$, then we could have made a choice for which $\boldsymbol{a}^{\circ},\hat{\boldsymbol{b}}^{\circ}$ are deterministic, because the available sites for a particle only depends on the leftmost point on the segment it belongs to but not on the particles to its left. However, for generic $\theta_{h,h}$ this turns out to be impossible: indeed, if $\theta_{h^k,h^k}$ is irrational, then given $\hat a^\circ_k$, $\hat b^\circ_k$ and the integrality condition \[
\N(\hat b^\circ_k-\hat a^\circ_k)-\theta_{h^k,h^k} N^\circ_k\in\amsmathbb{Z}_{\geq 0},
\]
one can uniquely reconstruct the integer $N^\circ_k$. Therefore, it cannot be random unless one of the endpoints $\hat{a}^{\circ}_k,\hat{b}^{\circ}_k$ is random.

\subsection{Detailed formulation}

Let us now carefully define $\overline{\boldsymbol N}^\circ$. We recall that $[\hat a_h,\hat b_h]$ indexed by $h\in[H]$ are the defining segments of the ensemble, that $(\alpha_k,\beta_k)$ indexed by $k\in[K]$ are the bands of the equilibrium measure $\boldsymbol{\mu}$, and that the $k$-th band is contained in the $h^k$-th segment. Like in Definition~\ref{def:eq_rescaled_parameters}, we use the convention for the primed endpoints, to be used when dealing with the equilibrium measure):
\[
\forall k \in [K] \qquad \hat a^{\circ\,\prime}_k=\hat a_k^\circ - \frac{\theta_{h^k,h^k}-\frac{1}{2}}{\N},\qquad \hat b^{\circ\,\prime}_k=\hat b^\circ_k + \frac{\theta_{h^k,h^k}-\frac{1}{2}}{\N}.
\]
As in Section~\ref{Section_FFF_no_sat}, we note that the state space $\W_\N$ of Section~\ref{Section_configuration_space} is slightly different depending on whether $a_1$ and $b_1$ are finite or infinite. Throughout our detailed exposition we exclude the situation $a_1=-\infty$ and $b_1$ finite; it can be reduced by $x \mapsto -x$ to the situation $a_1$ finite and $b_1=+\infty$ that we do treat.

We start by fixing a deterministic small $\eps>0$. We choose $\hat a^\circ_1$ in a deterministic way:
\begin{equation}
\label{firstQcirc}\hat{a}^{\circ}_1 = \left\{\begin{array}{lcl} \hat{a}_1 && \textnormal{if}\,\,a_1 \,\,\textnormal{is finite}, \\ \max\big\{ x\in \N^{-1}\amsmathbb{L}\quad \big|\,\, \quad x< \alpha_1-\eps\big\} && \textnormal{if}\,\,a_1 = -\infty. \end{array} \right.
\end{equation}
Next, we sequentially apply the following rules for $k = 1,2,\ldots,K$. As in the proof of Theorem~\ref{Theorem_CLT_for_filling_fractions}, they will involve the choice of $(2K-1)$ deterministic parameters $\mathfrak{t}_{1}^{\textnormal{b}},\ldots,\mathfrak{t}_{K}^{\textnormal{b}},\mathfrak{t}_{2}^{\textnormal{a}},\ldots,\mathfrak{t}_{K}^{\textnormal{a}}$, that can be always made such that Lemma~\ref{Lemma_fluct_params_choice} below holds. To make the logic apparent, some of the required properties will already be mentioned in the rules. \label{index:Qcirck}
\begin{itemize}
 \item[(R1)] If $(\alpha_k,\beta_k)$ is the rightmost band in $[\hat a'_{h^k},\hat b'_{h^k}]$, then we set:
 \begin{equation*}
 Q^\circ_{k} = N_k^{\circ} = \#\big\{i \in [N]\quad \big|\quad \ell_i \in [\N\hat{a}_{k}^{\circ},\N \hat{b}_{h^k}]\big\},\qquad \hat b^\circ_{k}=\hat b_{h^k}.
 \end{equation*}
 \item[(R2)] If $(\alpha_k,\beta_k)$ is \emph{not} the rightmost band in $[\hat a'_{h^k},\hat b'_{h^k}]$ while $[\beta_k,\alpha_{k+1}]$ is a void, then we set
\begin{equation*}
Q^\circ_{k}= N^{\circ}_k = \#\big\{i \in [N]\quad \big|\quad \ell_i \in [\N\hat{a}_{k}^{\circ},\N(\beta_{k}+\eps)]\big\}, \qquad
\hat{b}^{\circ}_{k}= \hat{a}^{\circ}_{k} + \N^{-1}\theta_{h^k,h^k}Q^\circ_{k} + \mathfrak{t}^{\textnormal{b}}_k,
\end{equation*}
The choice of $\mathfrak{t}_k^{\textnormal{b}}$ should be such that $\N \mathfrak{t}^{\textnormal{b}}_k\in\amsmathbb{Z}_{\geq 0}$ and $\hat{b}^{\circ}_{k}$ is larger than but close to $\beta_k+\eps$.
\item[(R3)] If $(\alpha_k,\beta_k)$ is \emph{not} the rightmost band in $[\hat a'_{h^k},\hat b'_{h^k}]$ while $[\beta_k,\alpha_{k+1}]$ is a saturation, then we set
\[
 N_k^\circ:=\begin{dcases} \big\lfloor \N \mu([\hat a'_{h^k},\beta_{k}+2 (\theta_{h^k,h^k})^{-1}\eps])\big\rfloor,& \textnormal{if } (\alpha_k,\beta_k)\textnormal{ is the leftmost band in }[\hat a'_{h^k},\hat b'_{h^k}],\\[2pt] \big\lfloor \N \mu([\alpha_k,\beta_{k}+2 (\theta_{h^k,h^k})^{-1}\eps])\big\rfloor,&\textnormal{otherwise.} \end{dcases}
\]
\label{index:Qkcirc}We further set $\hat{b}_k^{\circ}$ to be the location of the $N^{\circ}_{k}$-th particle (starting from the left) in $[\hat{a}_k^{\circ},+\infty)$, and $Q_k^{\circ}$ to be determined by the equation
\begin{equation}
\label{eq_rule3}
\hat{b}_k^{\circ} = \hat{a}^{\circ}_{k} + \frac{\theta_{h^k,h^k}(N_k^\circ-1)}{\N} + \frac{Q^\circ_{k}}{\N}.
\end{equation}
As a matter of fact, the construction of $\hat{a}_k^{\circ}$ in the previously applied rule forces $Q_k^{\circ}$ to be an integer. More precisely, it is the number of holes in $[ \hat{a}_{k}^{\circ},\hat{b}_{k}^{\circ}]$ --- if $\theta_{h,h}=1$ this is transparent, in general it can be seen from the definition of holes (Definition~\ref{Definition_hole}). As the band has nonempty interior, $Q_k^{\circ}$ must in fact be positive.
\item[(R4)] If $(\alpha_k,\beta_k)$ is the rightmost band in $[\hat a'_{h^k},\hat b'_{h^k}]$, which equivalently means that $(\alpha_{k+1},\beta_{k+1})$ is the leftmost band in $[\hat a'_{h^{k+1}},\hat b'_{h^{k+1}}]$, then we set deterministically
  $\hat a^{\circ}_{k+1}=\hat a_{h^{k+1}}$.
\item[(R5)] If $(\alpha_k,\beta_k)$ is \emph{not} the rightmost band in $[\hat a'_{h^k},\hat b'_{h^k}]$ while $(\beta_{k},\alpha_{k+1})$ is a void, then we set ${\hat{a}^{\circ}_{k+1}:= \hat{b}^{\circ}_k+ \mathfrak{t}^{\textnormal{a}}_{k+1}}$. The choice of $\mathfrak{t}^{\textnormal{a}}_{k+1}$ should be such that $\N \mathfrak{t}^{\textnormal{a}}_{k+1}\in\amsmathbb{Z}_{\geq 0}$ and $\hat{a}^{\circ}_{k+1}$ is smaller than but close to $\alpha_{k+1}-\eps$.
\item[(R6)] If $(\alpha_k,\beta_k)$ is \emph{not} the rightmost band in $[\hat a'_{h^k},\hat b'_{h^k}]$ while $[\beta_{k},\alpha_{k+1}]$ is a saturation, then we set ${\hat{a}^{\circ}_{k+1}:= \hat{b}^{\circ}_k+ \theta_{h^{k},h^k}\mathfrak{t}^{\textnormal{a}}_{k+1}}$. The choice of $\mathfrak{t}^{\textnormal{a}}_{k+1}$ is made such that $\N \mathfrak{t}^{\textnormal{a}}_{k+1}\in\amsmathbb{Z}_{\geq 0}$ and $\hat{a}^{\circ}_{k+1}$ is smaller than but close to $\alpha_{k+1}-\eps$.
\end{itemize}
The last three rules are applicable only if $k<K$.

\begin{lemma} \label{Lemma_fluct_params_choice}
 Under Assumptions~\ref{Assumptions_Theta}, \ref{Assumptions_basic}, \ref{Assumptions_offcrit} and \ref{Assumptions_analyticity}, there exist $\eps>0$, $(2K - 1)$ real deterministic constants $\mathfrak{t}^{\textnormal{b}}_{1},\ldots,\mathfrak{t}^{\textnormal{b}}_K,\mathfrak{t}_{2}^{\textnormal{a}},\ldots,\mathfrak{t}_K^{\textnormal{a}}$ in the above rules, a $\N$-independent constant $C > 0$ depending only on the constants in the assumptions, and an event $\mathcal{A}$ of probability at least $1 - Ce^{-\frac{\N}{C}(\log \N)^2}$, having the following properties in combination with the above rules.
 \begin{itemize}
 \item[(P1)] The segments $[\hat a^\circ_k,\hat b^\circ_k]$ indexed by $k\in[K]$ are disjoint and at distance at least $\eps$ from each other.
 \item[(P2)] $\forall k \in [K] \qquad [\alpha_k-\eps,\beta_k+\eps]\subseteq [\hat a^\circ_k,\hat b^\circ_k]\subseteq [\hat a_{h^k},\hat b_{h^k}]$.
 \item[(P3)] $\forall k \in [K] \qquad \big|N^\circ_k-\N \mu([\hat a^\circ_k,\hat b^\circ_k])\big|< C \N^{\frac{1}{2}} \log \N$.
 \item[(P4)] $\forall k \in [K] \quad N^\circ_k = \# \big\{i \in [N] \quad \big| \quad \ell_i \in [\N \hat{a}_k^{\circ},\N \hat{b}_k^{\circ}]\big\}$.
 \item[(P5)] The complement of $\,\bigcup_{k = 1}^{K} (\hat a^\circ_k,\hat b^\circ_k)$ in $\bigcup_{h = 1}^{H} [\hat a_h,\hat b_h]$ splits into finitely many segments, such that on each segment there are either no particles or no holes in the configuration $\boldsymbol{\ell}$, corresponding to whether the equilibrium measure $\boldsymbol{\mu}$ is void or saturated on these segments.
 \item[(P6)] The random $3K$-tuple $(\hat{\boldsymbol{a}}^\circ, \hat{\boldsymbol{b}}^\circ, \boldsymbol{N}^\circ)$ is a deterministic affine function of the random $K$-tuple $\boldsymbol{Q}^\circ$. Furthermore, for any $k \in[K]$, the difference $N^{\circ}_k-\N \mu([\hat a^{\circ\,\prime}_k, \hat b^{\circ\,\prime}_k])$ is a deterministic affine function of $(Q_1^\circ,\ldots,Q_k^\circ)$ and the linear coefficient of $Q_k^\circ$ takes values in $\big\{1,-\frac{1}{\theta_{1,1}},\ldots,-\frac{1}{\theta_{H,H}}\big\}$.
\end{itemize}
\end{lemma}

\begin{proof}
 Let $\eps>0$ be an arbitrary positive number such that the segments $[\alpha_k-10\eps,\beta_k+10\eps]$ indexed by $k \in [K]$ are pairwise disjoint and $[\alpha_k-10\eps,\beta_k+10\eps]\subset [\hat a_{h^k},\hat b_{h^k}]$. Let $C_1$ be the constant from Corollary~\ref{Corollary_a_priori_0}. For the event $\mathcal A$ we choose the set of all particle configurations $ \boldsymbol{\ell} \in \W_\N$ such that
 \begin{itemize}
 \item for any $k \in [K]$, the total number of particles in the segment $[\N(\alpha_k-\eps),\N(\beta_k+\eps)]$ differs from $\N\mu([\alpha_k-\eps,\beta_k+\eps])$ at most by $C_1\N^{\frac{1}{2}} \log \N$;
 \item the complement of $\bigcup_{k = 1}^{K} (\alpha_k - \eps,\beta_k + \eps)$ in $\bigcup_{h = 1}^{H} [\hat{a}_h,\hat{b}_h]$ splits into finitely many segments, such that on each segment there are either no particles or no holes in the configuration $\boldsymbol{\ell}$, corresponding to whether the equilibrium measure $\boldsymbol{\mu}$ is void or saturated on these segments.
 \end{itemize}
Corollary~\ref{Corollary_a_priori_0}, Theorem~\ref{Theorem_ldpsup} and Theorem~\ref{Theorem_ldsaturated} imply that $\mathcal{A}$ has overwhelming probability as desired.

The delicate point in choosing the constants $\mathfrak{t}^{\textnormal{b}}_{1},\ldots,\mathfrak{t}_{K}^{\textnormal{b}},\mathfrak{t}_2^{\textnormal{a}},\ldots,\mathfrak{t}_{K}^{\textnormal{a}}$ is that they should all be deterministic, \textit{i.e.} they are not allowed to depend on $\boldsymbol{\ell} \in \mathcal{A}$. So we select \emph{some} configuration $\boldsymbol{\ell}^{\textnormal{ref}} \in\mathcal A$ in an arbitrary way and construct $\mathfrak{t}^{\textnormal{b}}_{1},\ldots,\mathfrak{t}_{K}^{\textnormal{b}},\mathfrak{t}_2^{\textnormal{a}},\ldots,\mathfrak{t}_{K}^{\textnormal{a}}$ using this configuration. For that we repeat the general construction of $\overline{\boldsymbol{N}}^\circ$ basing ourselves on $\boldsymbol{\ell}^{\textnormal{ref}}$ instead of a random $\boldsymbol{\ell}$. In particular:
\begin{itemize}
 \item whenever we use the rule (R2), we pick $\mathfrak{t}^{\textnormal{b}}_k$ so that $\hat{b}^{\circ}_{k}\in (\beta_k+2\eps,\beta_k+3\eps)$;
 \item whenever we use the rules (R5) or (R6), we pick $\mathfrak{t}^{\textnormal{a}}_{k+1}$ so that $\hat{a}^{\circ}_{k+1}\in (\alpha_k-3\eps,\alpha_k-2\eps)$.
\end{itemize}
Now, for $\boldsymbol{\ell} \in \mathcal{A}$, we reuse the constants $\mathfrak{t}^{\textnormal{b}}_{1},\ldots,\mathfrak{t}_{K}^{\textnormal{b}},\mathfrak{t}_2^{\textnormal{a}},\ldots,\mathfrak{t}_{K}^{\textnormal{a}}$ specified using $\boldsymbol{\ell}^{\textnormal{ref}}$, and proceed with the definition of the extended $4K$-tuple $\overline{\boldsymbol{N}}^{\circ}$ as before. The definition of $\mathcal A$ implies that for large $\N$ this extended tuple does not vary too much. Hence, in this way we achieve the following properties.
\begin{itemize}
 \item Whenever we use the rule (R2), we have $\hat{b}^{\circ}_{k}\in (\beta_k+\eps,\beta_k+4\eps)$.
 \item Whenever we use the rule (R5) or (R6), we have $\hat{a}^{\circ}_{k+1}\in (\alpha_k-4\eps,\alpha_k-\eps)$.
\end{itemize}
At this stage the first five properties in Lemma~\ref{Lemma_fluct_params_choice} are implied by the rules of construction of $\overline{\boldsymbol{N}}^{\circ}$ and the choice of $\mathcal A$.

Let us demonstrate that (P6) also holds. The affine structure follows from definitions, and the triangular structure follows from the fact that we were sequentially applying the six rules. We only need to compute the coefficient of $Q_k^{\circ}$ in the difference $N^{\circ}_k-\N \mu( [\hat a^{\circ\,\prime}_k, \hat b^{\circ\,\prime}_k])$. For that we go over the six (sequentially applied) rules for the construction of $\overline{\boldsymbol N}^\circ$. In (R1) we set $N^\circ_k=Q^\circ_{k}$ while $\hat a^\circ_{k}$, $\hat b^\circ_{k}$ do not depend on $Q^\circ_{k}$. Hence, the coefficient of $Q^{\circ}_k$ is $1$. In (R2), $\hat b^\circ_{k}$ now also depends on $Q^\circ_k$, but since this endpoint is in a void, $\mu( [\hat a^{\circ\, \prime}_k, \hat b^{\circ\,\prime}_k])$ does not depend on this endpoint. Hence, the coefficient of $Q^{\circ}_k$ is again $1$. In (R3), neither $N^\circ_k$ nor $\hat a^{\circ}_k$ depend on $Q^\circ_k$, and $\hat b^{\circ}_k$ has an affine dependence on $Q^\circ_k$ as in \eqref{eq_rule3}. Combining with the observation that $\hat b^{\circ}_k$ is in a saturated region, where the density of $\mu$ is $(\theta_{h^k,h^k})^{-1}$, we conclude that the coefficient of $Q^\circ_k$ in $N^{\circ}_k-\N \mu([\hat a^{\circ\,\prime}_k, \hat b^{\circ\,\prime}_k])$ is $-(\theta_{h^k,h^k})^{-1}$. In (R4)-(R5)-(R6) the sequential definition of $N^\circ_k= Q^\circ_{k}$ does not affect the coefficient of $Q^\circ_k$ in $N^{\circ}_k-\N \mu([\hat a^{\circ\,\prime}_k, \hat b^{\circ\,\prime}_k])$.
\end{proof}

If $a_1=-\infty$ and $b_1$ is finite, then the procedure must be adapted. Instead of choosing $\hat{a}^\circ_1$ in a deterministic way, we choose deterministically $\hat b^{\circ}_{k^+(1)} = \hat b_1$, where we recall that $k^+(1)$ is the index of the last band in the first segment $[\hat{a}_1',\hat{b}_1']$. Then for $k> k^+(1)$, we apply by the same above rules, while for $k< k^+(1)$ and for $\hat{a}^\circ_{k^+(1)}$, we proceed with similar rules but in the opposite direction of decreasing $k$. Then, Lemma~\ref{Lemma_fluct_params_choice} holds with (P6) taking into account this different sequential ordering. The details are left to the interested reader.

In any case, the particular way we have defined the extended $4K$-tuple $\overline{\boldsymbol{N}}^{\circ}$ is not really important. The only properties that we are going to use are the six properties listed in Lemma~\ref{Lemma_fluct_params_choice}, and Theorems~\ref{Theorem_CLT_for_filling_fractions_saturation} and \ref{Theorem_linear_statistics_fluctuation_sat} continue to hold for any other choice satisfying these six properties. In particular, if $\theta_{h,h}=1$ for any $h\in[H]$, then instead of the previous procedure a deterministic choice of $(\hat{\boldsymbol{a}}^\circ,\hat{\boldsymbol{b}}^\circ)$ would be possible as we have already remarked.

 \subsection{Discrete Gaussian part}
\label{Sec:gausdssat}
Theorem~\ref{Theorem_CLT_for_filling_fractions_saturation} below shows that the extended $4K$-tuple $\overline{\boldsymbol{N}}^{\circ}$ is asymptotically equal to a discrete Gaussian distribution in the sense of Definition~\ref{Definition_asymptotic_equivalence}. As in Section~\ref{Sec:gausds}, we first spend some time describing the parameters of this distribution.

We recall that the $h$-th segment $[\hat a'_h,\hat b'_h]$ contains the bands $(\alpha_k,\beta_k)$ with $k\in \llbracket k^-(h),k^{+}(h)\rrbracket$. We also recall that $\amsmathbb{S}_{h}$ is the union of the saturations in $[\hat a'_h,\hat b'_h]$. On the event $\mathcal A$ of Lemma~\ref{Lemma_fluct_params_choice}, due to (P5) we have
\begin{equation}
\label{eq_x258}
\begin{split}
 \N \hat{n}_h & = \sum_{k = k^-(h)}^{k^+(h)} N^{\circ}_{k}+ \frac{\N}{\theta_{h,h}} \Bigg(\mathbbm{1}_{\amsmathbb{S}_h}(\hat{a}'_h)\big(\hat a^\circ_{k^-(h)}-\hat a_h\big) + \sum_{k = k^-(h)}^{k^+(h)-1}\mathbbm{1}_{\amsmathbb{S}_h}(\beta_k)\Big(\hat a^\circ_{k+1}-\hat b^\circ_{k}-\frac{1}{\N}\Big) \\
 & \qquad\qquad \qquad \qquad \qquad + \mathbbm{1}_{\amsmathbb{S}_h}(\hat{b}_h') \big(\hat b_h-\hat b^\circ_{k^+(h)}\big)\bigg).
\end{split}
\end{equation}
Each of the three segments $[\hat{a}_h',\alpha_{k^-(h)}]$, $[\beta_k,\alpha_{k + 1}]$ and $[\beta_{k^+(h)},\hat{b}_h']$ are saturations if and only if one of their endpoints is in a saturation. We used this fact to write the formula in a shorter way, putting indicator for one of the endpoints only. Note that the right-hand side of \eqref{eq_x258} is an affine function in $(\hat{\boldsymbol{a}}^\circ,\hat{\boldsymbol{b}}^\circ, \boldsymbol{N}^\circ)$, and therefore is also an affine function in $\boldsymbol{Q}^\circ_{k}$. Accordingly, on the event $\mathcal A$ the affine constraints \eqref{eq_equations_eqs} in Section~\ref{DataS} on $(N_h)_{h = 1}^{H}$ can be rewritten as affine constraints on the numbers $\boldsymbol{Q}^{\circ}$. They define an affine subspace $\L$ in $\amsmathbb{R}^{K}$ and we call $\L_0$ the parallel subspace passing through the origin.

\medskip

Next, we recall that in Section~\ref{Section_FFF_no_sat} the normalized filling fractions were varying in a small neighborhood of $\hat{\boldsymbol n}^{\boldsymbol{\mu}}$ --- the filling fractions of the equilibrium measure $\boldsymbol{\mu}$. In the present situation we use $\boldsymbol{Q}^{\circ}$ instead of the filling fractions and we have to explain near which point the normalized versions of $\boldsymbol{Q}^\circ$ are varying. Recall that $(\hat{\boldsymbol{a}}^\circ,\hat{\boldsymbol{b}}^\circ, \boldsymbol{N}^\circ)$ is an affine function of $\boldsymbol{Q}^{\circ}$. Therefore, so are the rescaled version of the filling fractions and the shifted endpoints $(\hat{\boldsymbol{a}}^{\circ\,\prime},\hat{\boldsymbol{b}}^{\circ\,\prime},\hat{\boldsymbol{n}}^{\circ})$. We denote
\begin{equation}
\label{eq_ff_map}
 \boldsymbol{Q} \mapsto \big(\hat{\boldsymbol{a}}(\boldsymbol{Q}),\hat{\boldsymbol{b}}(\boldsymbol{Q}),\boldsymbol{N}(\boldsymbol{Q})\big) \qquad \textnormal{and} \qquad \boldsymbol{Q} \mapsto \big(\hat{\boldsymbol{a}}'(\boldsymbol{Q}),\hat{\boldsymbol{b}}'(\boldsymbol{Q}),\hat{\boldsymbol{n}}(\boldsymbol{Q})\big)
 \end{equation}
these affine functions. Note that the values of interest of $\boldsymbol {Q}$ are of order $\N$, while the hatted variables in \eqref{eq_ff_map} have coordinates of order $1$. Similarly to \eqref{eq_x258}, but now using the equilibrium measure $\boldsymbol{\mu}$, we have as long as the segments $[\hat a_k(\boldsymbol {Q}),\hat b_k (\boldsymbol {Q})]$ indexed by $k \in [K]$ are pairwise disjoint subsets of the intervals $[\hat a'_{h^k},\hat b'_{h^k}]$:
\begin{equation}
\label{eq_x259}
\begin{split}
 \mu([\hat a'_h, \hat{b}'_h]) &= \sum_{k = k^-(h)}^{k^+(h)} \mu([\hat a'_k(\boldsymbol {Q}),\hat b'_k (\boldsymbol {Q})]) + \frac{1}{\theta_{h,h}} \bigg(\mathbbm{1}_{\amsmathbb{S}_h}(\hat{a}'_h)\big(\hat a'_{k^-(h)}(\boldsymbol {Q}) -\hat a'_h\big)\\ & \quad + \sum_{k = k^-(h)}^{k^+(h)-1}\mathbbm{1}_{\amsmathbb{S}_h}(\beta_k)\big(\hat a_{k+1}'(\boldsymbol {Q})-\hat b_k'(\boldsymbol {Q})\big)
+ \mathbbm{1}_{\amsmathbb{S}_h}(\hat{b}'_h) \big(\hat b'_h-\hat b_{k^+(h)}'(\boldsymbol {Q})\big)\bigg).
\end{split}
\end{equation}
The indicator functions here select the situations where $[\hat{a}_h',\alpha_{k^-(h)}]$, $[\beta_{k},\alpha_{k + 1}]$ and $[\alpha_{k^+(h)},\hat{b}'_h]$ are saturations, respectively. Therefore, using \eqref{eq_x259} the affine constraints \eqref{eq_equations_eqs} in Section~\ref{DataS} on $\big(\mu([\hat{a}_h',\hat{b}_h'])\big)_{h = 1}^{H}$ can be rewritten as affine constraints on the vectors $(\hat{\boldsymbol{a}}'(\boldsymbol {Q}),\hat{\boldsymbol{b}}'(\boldsymbol {Q}))$ and $\big(\mu([\hat a_k'(\boldsymbol {Q}),\hat b_k' (\boldsymbol {Q})])\big)_{k=1}^K$. They define an affine subspace $\L'$ in $\amsmathbb{R}^{K}$. The underlying linear subspace $\L_0$ is the same as the one underlying $\L$. However, in general $\L\neq \L'$: they differ by a $O(1)$ translation which is caused by the $-\frac{1}{\N}$ term in \eqref{eq_x258}, which is missing in \eqref{eq_x259}, and by the difference between $\hat{\boldsymbol{a}}^\circ,\hat{\boldsymbol{b}}^\circ$ used for \eqref{eq_x258} versus their shifted version $\hat{\boldsymbol{a}}^{\circ\,\prime},\hat{\boldsymbol{b}}^{\circ\,\prime}$ relevant for \eqref{eq_x259}. Note however that these $O(1)$ difference should be thought as being small here, because the values of interest for $\boldsymbol{Q}$ are of order $\N$.

\begin{lemma} \label{Lemma_center_of_expansion}
 There exists a unique $K$-dimensional vector $\boldsymbol {Q}^{\boldsymbol{\mu}}\in\L'$, such that
 \begin{itemize}
  \item the segments $[\hat a_k'(\boldsymbol{Q}^{\boldsymbol{\mu}}),\hat b_k' (\boldsymbol {Q}^{\boldsymbol{\mu}})]$ indexed by $k \in [K]$ are pairwise disjoint;
  \item $\forall k \in [K] \quad [\alpha_k,\beta_k]\subseteq (\hat a_k'(\boldsymbol{Q}^{\boldsymbol{\mu}}),\hat b_k' (\boldsymbol{Q}^{\boldsymbol{\mu}})) \subseteq (\hat a'_{h^k}, \hat b'_{h^k})$;
 \item $\forall k \in [K] \quad \hat n_k(\boldsymbol{Q}^{\boldsymbol{\mu}}) = \mu([\hat a'_k(\boldsymbol{Q}^{\boldsymbol{\mu}}),\hat b'_k (\boldsymbol{Q}^{\boldsymbol{\mu}})])$.
 \end{itemize}
 \end{lemma}
Unlike $\boldsymbol{Q}^{\circ}$, the coordinates of $\boldsymbol{Q}^{\boldsymbol{\mu}}$ do not have to be integers.

\begin{proof}[Proof of Lemma~\ref{Lemma_center_of_expansion}] (P3) in Lemma~\ref{Lemma_fluct_params_choice} tells us that for random integral $\boldsymbol{Q}^\circ$ the third point of Lemma~\ref{Lemma_center_of_expansion} is satisfied up to a small error, with overwhelming probability. Simultaneously, for this integral $\boldsymbol{Q}^{\circ}$, (P1)-(P2) in Lemma~\ref{Lemma_fluct_params_choice} translate into the two first properties in Lemma~\ref{Lemma_center_of_expansion}.

On the other hand, (P6) in Lemma~\ref{Lemma_fluct_params_choice} implies that the difference ${\hat n_k(\boldsymbol {Q})-\mu([\hat a_k'(\boldsymbol {Q}),\hat b_k' (\boldsymbol {Q})])}$
is an affine function of $Q_1,\ldots,Q_k$ with linear coefficient of $Q_k$ equal to $\frac{1}{\N}$, for any $k \in [K]$. Hence, deforming a typical random integral $\boldsymbol{Q}^{\circ}$, we find a unique $\boldsymbol {Q}^{\boldsymbol{\mu}} \in\amsmathbb{R}^{K}$, such that $\mu([\hat a_k'(\boldsymbol {Q}^{\boldsymbol{\mu}}),\hat b_k' (\boldsymbol {Q}^{\boldsymbol{\mu}})])= \hat n_k(\boldsymbol {Q}^{\boldsymbol{\mu}})$ for any $k\in [K]$. The definition of the affine subspace $\L'$ was designed precisely to guarantee $\boldsymbol{Q}^{\boldsymbol{\mu}} \in\L'$. As explained in the previous paragraph, the desired deformation is small, and, therefore, the endpoints of the various segments do not change much. Therefore, (P1)-(P2)-(P3) in Lemma~\ref{Lemma_fluct_params_choice} guarantee that the first two points in the claim hold for this $\boldsymbol{Q}^{\boldsymbol{\mu}}$.
\end{proof}

We now turn to the definition of the quadratic bilinear form $\Qu$ on $\amsmathbb{R}^K$ entering into \eqref{eq_Discrete_Gaussian_main}. Recall that the equilibrium measure was defined in Section~\ref{Section_Energy_functional} as the minimizer of the energy functional $-\I$ over the set $\mathscr{P}_{\star}$ of positive measures on $\bigcup_{h = 1}^{H} [\hat{a}_h',\hat{b}_h']$ whose segment filling fractions obey the affine constraints \eqref{eq_equations_eqs}. Let us take an auxiliary $K$-tuple $\hat{\boldsymbol p}$, allowed to vary in a small $\N$-independent neighborhood of $\frac{\boldsymbol{Q}^{\boldsymbol{\mu}}}{\N}$ of Lemma~\ref{Lemma_center_of_expansion}, and consider the minimizer $\boldsymbol{\mu}^{\hat{\boldsymbol{p}}}$ of $-\I$ over the $H$-tuple of measures $\boldsymbol{\nu} \in \mathscr{P}_{\star}$ satisfying as well the constraints
\[
\forall k \in [K] \qquad \nu([\hat{a}_k'(\N\hat{\boldsymbol{p}}),\hat{b}_k'(\N\hat{\boldsymbol{p}})]) = \hat{n}_k(\N \hat{\boldsymbol{p}})
\]
using the affine map \eqref{eq_ff_map}. Since the latter describe a convex subset of $\mathscr{P}_{\star}$, this minimizer exists and is unique as in Proposition~\ref{Lemma_maximizer}. We then define $\Qu$ as the Hessian of $\hat{\boldsymbol{p}} \mapsto -\I[\boldsymbol{\mu}^{\hat{\boldsymbol{p}}}]$ at the point $\boldsymbol{Q}^{}$:
\begin{equation}
\label{eq_discrete_covariance_sat}
\forall k,l \in [K] \qquad \Qu_{k,l} = - \partial_{\hat{p}_k}\partial_{\hat{p}_l} \I[\boldsymbol{\mu}^{\hat{\boldsymbol{p}}}] \big|_{\hat{\boldsymbol{p}} = \frac{\boldsymbol{Q}^{\boldsymbol{\mu}}}{\N}} .
\end{equation}

\begin{proposition} \label{Lemma_quadratic_non_degen_sat}
The restriction of the quadratic form $\Qu$ to $\L_0$ is positive definite and its values on the intersection of $\L_0$ with the unit sphere in $\amsmathbb{R}^K$ is uniformly bounded and uniformly bounded away from $0$ as $\N\rightarrow\infty$.
\end{proposition}
\begin{proof}
 Let us choose small $\eps_1>0$, which does not depend on $\hat{\boldsymbol p}$ and such that for any $k \in [K]$ and $\hat{\boldsymbol{p}}$ in a neighborhood of $\frac{\boldsymbol{Q}^{\boldsymbol{\mu}}}{\N}$, we have $(\alpha_k-\eps_1,\beta_k+\eps_1)\subset [\hat a_k'(\N\hat{\boldsymbol{p}}), \hat b_k'(\N\hat{\boldsymbol{p}})]$. We then localize $\mu^{\hat{\boldsymbol p}}$ to the segments $\bigcup_{k=1}^K [\alpha_k-\eps_1,\beta_k+\eps_1]$ by the procedure of Chapter~\ref{Chapter_conditioning}. The localized measure has support in the new segments and minimizes an energy functional, in which the interaction matrix is the composition of the matrix $\boldsymbol{\Theta}$ of the original ensemble with the map $k\mapsto h^k$, and the potential $V_k(x)$ is the sum of $V_{h^k}(x)$ and the additional terms
 \begin{equation}
 -2 \frac{\theta_{h^k,g}}{\theta_{h^k,h^k}} \int_{\mathfrak a}^{\mathfrak b} \log|x-y| \dd y,
\end{equation}
for each $g \in [H]$ and each saturation $[\mathfrak{a},\mathfrak{b}]$ of $\mu$ in the complement of $\bigcup_{k = 1}^{K} (\alpha_k - \eps_1,\beta_k + \eps_1)$ in $[\hat a'_{g},\hat b'_{g}]$.

Arguing as in the first claim of the proof of Lemma~\ref{Lemma_energyconditioner}, the Hessian we are interested in is the same as the Hessian of the energy functional for the localized measures. Because both $[\hat a_k'(\N \hat{\boldsymbol{p}}), \alpha_k-\eps_1]$ and $[\beta_k+\eps_1, \hat b_k'(\N\hat{\boldsymbol{p}})]$ are either void or saturated, the constraints $\mu^{\hat{\boldsymbol p}}([\hat a_k'(\N \hat{\boldsymbol{p}}),\hat{b}_k'(\N\hat{\boldsymbol{p}})])=\hat{n}_k(\N\hat{\boldsymbol{p}})$ can be rewritten as filling fraction constraints on the localization of $\boldsymbol{\mu}^{\hat{\boldsymbol p}}$. In other words, we have constraints on $\mu^{\hat{\boldsymbol p}}([\alpha_k-\eps_1,\beta_k+\eps_1])$, which linearly depend on $\hat{\boldsymbol p}$. Therefore, the matrix $\Qu$ is obtained from the Hessian of Proposition~\ref{Proposition_Hessian_free_energy} by a linear change of variables. This change of variables is non-degenerate by (P6) in Lemma~\ref{Lemma_fluct_params_choice}. Hence, Proposition~\ref{Lemma_quadratic_non_degen_sat} follows from Proposition~\ref{Proposition_Hessian_free_energy}.
\end{proof}

The last ingredient to define the Gaussian variable relevant for Theorem~\ref{Theorem_CLT_for_filling_fractions_saturation} is the centering vector $\boldsymbol u$. It is very similar to Section~\ref{Section_FFF_no_sat}. We start by defining a quadratic polynomial function $P(\boldsymbol{x})$ of $\boldsymbol{x} \in \L\subseteq \amsmathbb{R}^K$ through the formula
\begin{equation}
\label{eq_x260}
P(\boldsymbol{x})=-\frac{1}{2}\Qu(\boldsymbol{x}- \boldsymbol{Q}^{\boldsymbol{\mu}}, \boldsymbol{x}-\boldsymbol{Q}^{\boldsymbol{\mu}}) +  \sum_{h=1}^H \theta_{h,h}\big(N_h(\x)-N_h(\boldsymbol{Q}^{\boldsymbol{\mu}})\big)\log \N  + \sum_{k=1}^K x_k\, \mathbbm{Shift}_k.
\end{equation}
Here $N_h(\boldsymbol{x})$ is an affine function of $\x$, obtained by expressing the total number of particles in $[\hat a_h,\hat b_h]$ through the vector $\boldsymbol Q^\circ$ via the formula \eqref{eq_x258}, and replacing $\boldsymbol{Q}^{\boldsymbol{\mu}}$ with $\boldsymbol{x}$. The constants $\mathbbm{Shift}_k$ in \eqref{eq_x260} are defined as partial derivatives with respect to $Q^\circ_k$ --- evaluated at $\boldsymbol{Q}^{\boldsymbol{\mu}}$ --- of the sum of two terms: first, the remainder $\mathbbm{Rest}_1$ in \eqref{eq_partition_multicut} for the expansion of the partition function of the discrete ensemble localized to $\bigcup_{l = 1}^{K} [\hat a_l^\circ,\hat b_l^\circ]$ and with filling fractions fixed to be $\boldsymbol{N}(\boldsymbol{Q}^{\boldsymbol{\mu}})$; second, the prefactor of the order $\N$ term in the difference of the logarithms of the partition functions of the original ensemble and the localized one, \textit{cf.} \eqref{eq_x272} and \eqref{eq_x273} for more details. Again, both terms are bounded as $\N$ becomes large, with upper bound depending only on the constants in the assumptions.

By completing the square, $P(\boldsymbol{x})$ can be rewritten as
\begin{equation}
\label{eq_x261}
 P(\boldsymbol{x})=-\frac{1}{2}\Qu(\boldsymbol{x}-\boldsymbol u, \boldsymbol{x}-\boldsymbol u) + c,
\end{equation}
where $c$ does not depend on $\boldsymbol{x}$. The determination of $\boldsymbol u$ through \eqref{eq_x260}-\eqref{eq_x261} may not be unique, because in this identity $\boldsymbol{x}$ ranges only over $\L$ and not over $\amsmathbb{R}^K$. In particular, we can specify $\boldsymbol{u}$ using the orthogonal projection onto $\L_0$. Both $\Qu$ and $\boldsymbol u$ may depend on $\N$, and they clearly depend on the equilibrium measure. The center $\boldsymbol u$ of the discrete Gaussian distribution always contains a term of order $\N$ coming from $\boldsymbol{Q}^{\boldsymbol{\mu}}$ and, in general, it also contains a correction of order $\log \N$. However, the latter becomes absent if all $\theta_{h,h}$ are equal for $h \in [H]$, because $(1,\ldots,1)$ is always orthogonal to $\L_0$, as the total number of particles is always deterministically fixed in our framework.

\subsection{Results on the fluctuations of extended parameters and linear statistics}

Our first result describes the asymptotic distribution of the random auxiliary $K$-tuple $\boldsymbol{Q}^{\circ}$ in terms of the discrete Gaussian random variable $\textnormal{\textsf{\textbf{Gau\ss{}}}}_{\amsmathbb{Z}} := \textnormal{\textsf{\textbf{Gau\ss{}}}}_{\amsmathbb{Z}}[\Qu,\L,\boldsymbol u]$ with its parameters described in Section~\ref{Sec:gausdssat}.
\begin{theorem}
 \label{Theorem_CLT_for_filling_fractions_saturation}
 If the discrete ensemble satisfies Assumptions~\ref{Assumptions_Theta}, \ref{Assumptions_basic}, \ref{Assumptions_offcrit} and \ref{Assumptions_analyticity}, then we have as $\N \rightarrow \infty$:
\begin{equation}
\label{eq_discrete_filling_fractions_limit_sat} \boldsymbol{Q}^{\circ} \stackrel{\textnormal{d}}{\sim} \textnormal{\textsf{\textbf{Gau\ss{}}}}_{\amsmathbb{Z}},
\end{equation}
where the convergence in the sense of Definition~\ref{Definition_asymptotic_equivalence} is uniform for fixed constants in the assumptions.
\end{theorem}
By construction, the random $3K$-tuple $(\hat{\boldsymbol{a}}^{\circ},\hat{\boldsymbol{b}}^{\circ},\hat{\boldsymbol{N}}^{\circ})$ is a deterministic affine functions of $\boldsymbol{Q}^{\circ}$. Therefore, the distribution of $(\hat{\boldsymbol{a}}^{\circ},\hat{\boldsymbol{b}}^{\circ},\hat{\boldsymbol{N}}^{\circ})$ is also asymptotically equal to a random discrete Gaussian vector.

\medskip

Similarly to the way Theorem~\ref{Theorem_CLT_for_filling_fractions} implied the asymptotic behavior of linear statistics in Theorem~\ref{Theorem_linear_statistics_fluctuation_ff}, we can now completely analyze the linear statistics in presence of saturations. We recall that the $h$-th shifted segment $[\hat a'_h,\hat b'_h]$ contains the bands $(\alpha_k,\beta_k)$ with $k\in \llbracket k^-(h),k^+(h)\rrbracket$, and that $\boldsymbol{Q}^{\boldsymbol{\mu}}$ is the deterministic $K$-tuple of integers introduced in Lemma~\ref{Lemma_center_of_expansion}.

\begin{definition}
\label{DhfL} Let $h \in [H]$ and $k \in \llbracket k^-(h),k^+(h) \rrbracket$. We introduce
\begin{equation*}
\begin{split}
& \mathfrak{a}_k = \hat{a}_k(\textnormal{\textsf{\textbf{Gau\ss{}}}}_{\amsmathbb{Z}}),\qquad \mathfrak{a}_k' = \mathfrak{a}_k  + \frac{\theta_{h,h} - \frac{1}{2}}{\N} \\
& \mathfrak{b}_k = \hat{b}_k(\textnormal{\textsf{\textbf{Gau\ss{}}}}_{\amsmathbb{Z}}),\qquad \mathfrak{b}_k' = \mathfrak{b}_k - \frac{\theta_{h,h} - \frac{1}{2}}{\N}
\end{split}
\end{equation*}
If $f$ is a continuous function defined in a neighborhood of the $h$-th segment $[\hat{a}_h,\hat{b}_h]$, we set
\begin{equation}
\label{def_shiftsfh}
\begin{split}
\mathsf{Shift}_h[f] & = \sum_{l = 1}^{K} (\textnormal{\textsf{Gau\ss{}}}_{\amsmathbb{Z},l} - Q_l^{\boldsymbol{\mu}}) \int_{\hat{a}_h'}^{\hat{b}'_h} f(x) \big(\partial_{\hat{p}_l} \mu_h^{\hat{\boldsymbol{p}}}(x)\big)\big|_{\hat{\boldsymbol{p}} = \frac{\boldsymbol{Q}^{\boldsymbol{\mu}}}{\N}}\,\dd x \\
& \quad + \frac{\theta_{h,h} - 1}{2\theta_{h,h}}\bigg[\mathbbm{1}_{\amsmathbb{S}_h}(\hat{a}_h')\big(f(\mathfrak{a}_{k^-(h)}') - f(\hat{a}_{h}')\big)  + \sum_{k = k^-(h)}^{k^+(h) - 1} \mathbbm{1}_{\amsmathbb{S}_h}(\beta_k) \big(f(\mathfrak{b}_k') + f(\mathfrak{a}'_{k + 1})\big)  \\
& \quad\qquad \qquad\qquad + \mathbbm{1}_{\amsmathbb{S}_h}(\hat{b}'_h) \big(f(\mathfrak{b}_{k^+(h)}') - f(\hat{b}'_h)\big)\bigg]
\end{split}
\end{equation}
As in \eqref{eq_x258}, the indicator functions indicate respectively that $[\hat{a}_h',\alpha_{k^+(h)}]$, $[\beta_{k},\alpha_{k + 1}]$ and $[\beta_{k^+(h)},\hat{b}_h']$ are saturations.  The equilibrium measure $\boldsymbol{\mu}^{\hat{\boldsymbol{p}}}$ is the one already encountered in Section~\ref{Sec:gausdssat}.
\end{definition}

In this definition, the first term in the right-hand side is the analogue of the second term of \eqref{Lin1km} where the expansion point $\hat{\boldsymbol{n}}^{\boldsymbol{\mu}}$ is replaced with $\frac{\boldsymbol{Q}^{\boldsymbol{\mu}}}{\N}$. The terms in brackets are new: they will arise from the difference between the sum of evaluations of $f$ at the densely packed particles in saturations and the integral of $f$, using an Euler--Maclaurin approximation like in Lemma~\ref{lem:EulerMaclaurin}.

To formulate the next result, we choose arbitrarily for each $k \in [K]$ a deterministic and $\N$-independent segment $[\overline{\alpha}_k,\overline{\beta}_k]$ containing $[\hat{a}_k^{\circ},\hat{b}_k^{\circ}]$ in its interior with probability tending to $1$ as $\N \rightarrow \infty$. Such a choice is always possible because Theorem~\ref{Theorem_CLT_for_filling_fractions_saturation} guarantees that the fluctuations of $\hat{a}_k^{\circ},\hat{b}_k^{\circ}$ are of order $O(\frac{1}{\N})$. Then, the $k$-th band $(\alpha_k,\beta_k)$ is automatically contained in the interior of $(\overline{\alpha}_k,\overline{\beta}_k)$.

\begin{theorem}
\label{Theorem_linear_statistics_fluctuation_sat}

Consider a discrete ensemble satisfying Assumptions~\ref{Assumptions_Theta}, \ref{Assumptions_basic}, \ref{Assumptions_offcrit} and \ref{Assumptions_analyticity}. Let $L \in \amsmathbb{Z}_{> 0}$ and a $L$-tuple of integers $\boldsymbol{h} \in [H]^L$ both independent of $\N$, let a (possibly $\N$-dependent) $L$-tuple of functions $\boldsymbol{f}(z)$ such that $f_l(z)$ is twice continuously differentiable in $\amsmathbb{A}_{h_l}^{\mathfrak{m}}$ and is a holomorphic function of $z$ in a $\N$-independent complex neighborhood $\amsmathbb{K}_h$ of $\bigcup_{k = k^-(h_l)}^{k^+(h_l)} [\overline{\alpha}_k,\overline{\beta}_k]$. Assume there exists a constant $C > 0$ such that $\max_{l} \sup_{z \in \amsmathbb{K}_h} |f_l(z)| \leq C$.

Let $\textnormal{\textsf{\textbf{Gau\ss{}}}}[\boldsymbol{f}]$ be a $L$-dimensional random Gaussian vector independent of the $K$-dimensional discrete Gaussian vector $\textnormal{\textsf{\textbf{Gau\ss{}}}}_{\amsmathbb{Z}}$ appearing in Theorem~\ref{Theorem_CLT_for_filling_fractions_saturation}, having covariance
\begin{equation}
\label{eq_x264}
\textnormal{\textsf{Cov}}_{l_1,l_2}[\boldsymbol{f}] = \sum_{k_1=k^-(h_{l_1})}^{k^+(h_{l_1})} \, \sum_{k_2=k^-(h_{l_2})}^{k^+(h_{l_2})} \oint_{\gamma_{l_1}} \oint_{\gamma_{l_2}}
 f_{l_1}(z_1) f_{l_2}(z_2)
\,\mathcal{F}_{k_1,k_2}(z_1,z_2)\,\frac{\dd z_1 \dd z_2}{(2\ii\pi)^2}.
\end{equation}
and $l$-th component of the mean equal to the sum of $\mathbbm{Rest}[f_l]$ in \eqref{Elink1} and $\mathbbm{Rest}'[f_l]$ in \eqref{Elinkm2} below.
Then, we have as $\N \rightarrow \infty$
\begin{equation}
\label{eq_x274}
\left(\sum_{i = 1}^{N} \mathbbm{1}_{[a_{h_l},b_{h_l}]}(\ell_i) f_l\bigg(\frac{\ell_i}{\N}\bigg) - \N \int_{\hat{a}'_{h_l}}^{\hat{b}'_{h_l}} f_l(x)\,\mu_{h_l}(x)\,\dd x\right)_{l = 1}^{L} \,\,\mathop{\sim}^{\textnormal{d}} \,\,\textnormal{\textsf{\textbf{Gau\ss{}}}}[\boldsymbol{f}] + \big(\mathsf{Shift}_{h_l}[f_l]\big)_{l = 1}^{L}.
\end{equation}
The convergence in Definition~\ref{Definition_asymptotic_equivalence} is uniform for fixed constants in the assumptions.
\end{theorem}

As in Theorem~\ref{Theorem_linear_statistics_fluctuation_ff}, the asymptotic distribution in Theorem~\ref{Theorem_linear_statistics_fluctuation_sat} contains two implicit shifts: in the mean for $\textnormal{\textsf{\textbf{Gau\ss{}}}}$ via $\mathbbm{Shift}_h$ and in the center $\boldsymbol{u}$ used in the definition of $\textnormal{\textsf{\textbf{Gau\ss{}}}}_{\amsmathbb{Z}}$ and thus $\mathsf{Shift}_h[f]$ in Definition~\ref{DhfL}. It would be interesting to find more explicit formulae for these shifts. Note however the simplification in Definition~\ref{DhfL} in case $\theta_{g,g} = 1$ for every $g \in [H]$.

There is a noteworthy difference between our formulations of Theorems~\ref{Theorem_linear_statistics_fluctuation_ff} and \ref{Theorem_linear_statistics_fluctuation_sat}: in the former the functions $f_l$ were defined in a neighborhood of the bands, while in the latter they are defined on a neighborhood of the full segments $[\hat a_h,\hat b_h]$. This is because saturations between bands prevent us from singling out contributions of each band individually. In Theorem~\ref{Theorem_linear_statistics_fluctuation_sat} It is possible to allow $f_l$ any function (not necessarily twice continuously differentiable) on $\amsmathbb{A}^{\mathfrak{m}}_{h_l}$ which is holomorphic in a complex neighborhood $\amsmathbb{K}_h$ of the bands of the $h_l$-th segment, provided we replace \eqref{def_shiftsfh} in Definition~\ref{DhfL} with
\begin{equation*}
\begin{split}
\textnormal{\textsf{Shift}}_h[f] & = \sum_{l = 1}^{K} (\textnormal{\textsf{Gau\ss{}}}_{\amsmathbb{Z},l} - Q_l^{\boldsymbol{\mu}}) \int_{\hat{a}_h'}^{\hat{b}'_h} f(x) \big(\partial_{\hat{p}_l} \mu_h^{\hat{\boldsymbol{p}}}(x)\big)\big|_{\hat{\boldsymbol{p}} = \frac{\boldsymbol{Q}^{\boldsymbol{\mu}}}{\N}}\,\dd x \\
& \quad + \mathbbm{1}_{\amsmathbb{S}_h}(\hat{a}'_h)\Bigg[f(\hat{a}_h) + \cdots + f\bigg(\mathfrak{a}_{k^-(h)}- \frac{\theta_{h,h}}{\N}\bigg)\Bigg] \\
& \quad + \sum_{k = k^-(h)}^{k^+(h) - 1} \mathbbm{1}_{\amsmathbb{S}_h}(\beta_k)\Bigg[f\bigg(\mathfrak{b}_k + \frac{\theta_{h,h}}{\N}\bigg) + \cdots + f\bigg(\mathfrak{a}_{k + 1} - \frac{\theta_{h,h}}{\N}\bigg)\Bigg] \\
& \quad + \mathbbm{1}_{\amsmathbb{S}_h}(\hat{b}'_h) \Bigg[f\bigg(\mathfrak{b}_{k^+(h)} + \frac{\theta_{h,h}}{\N}\bigg) + \cdots + f(\hat{b}_h)\Bigg],
\end{split}
\end{equation*}
where the sums in the three last lines involve arguments of $f$ increasing by steps of $\frac{\theta_{h,h}}{\N}$. Indeed, in the proof this expression shows up in \eqref{Lin3km2} and we use the regularity of $f_l$ outside the bands only to approximate these Riemann sums via Lemma~\ref{lem:EulerMaclaurin}, leading us to the simpler Definition~\ref{DhfL}.

\subsection{Proof of Theorem~\ref{Theorem_CLT_for_filling_fractions_saturation}}

The proof will follow the strategy of Step 2 in Section~\ref{Secrigub}. By definition, the probability distribution of $\boldsymbol{Q}^{\circ}$ can be computed as:
\begin{equation}
\label{eq_x265}
 \P(\boldsymbol{Q}^{\circ}=\boldsymbol M)= \frac{\Z^{\boldsymbol M}_\N}{\Z_\N},
\end{equation}
where $\Z_\N^{\boldsymbol M}$ is the partition function of the ensemble restricted to the configurations satisfying $\boldsymbol{Q}^{\circ}=\boldsymbol M$, \textit{i.e.} it is the part of the sum \eqref{eq_partition_function_definition} corresponding to such configurations. It is sufficient to only study $\boldsymbol M\in\L$ which corresponds to the configurations from the set $\mathcal A$ of Lemma~\ref{Lemma_fluct_params_choice}, because others have negligible probabilities. We claim that such configurations satisfy for a constant $C>0$:
\begin{equation}
\label{eq_x266}
 |\!|\boldsymbol{M} - \boldsymbol{Q}^{\boldsymbol{\mu}}|\!|_{\infty} < C \N^{\frac{1}{2}} \log \N.
\end{equation}
Indeed, \eqref{eq_x266} is obtained by comparing (P3) in Lemma~\ref{Lemma_fluct_params_choice} with the third point in Lemma~\ref{Lemma_center_of_expansion} and using the triangularity stated in (P6) in Lemma~\ref{Lemma_fluct_params_choice}. Furthermore, using the bound on the probability of $\mathcal A$ of Theorem~\ref{proposition_fluct_conditioning}, we can replace $\Z_\N$ in denominator of \eqref{eq_x265} by the sum of $\Z^{\boldsymbol M}_\N$ over all $\boldsymbol M\in \L$ satisfying \eqref{eq_x266}. This replacement leads to an exponentially small asymptotic error.

We compute $\Z^{\boldsymbol M}_\N$ in \eqref{eq_x265} by localizing the ensemble to $\bigcup_{k = 1}^{K} [\hat a^{\circ}_k, \hat{b}^{\circ}_k]$ and fixing simultaneously $\overline{\boldsymbol{N}}^{\circ}$ by the same procedure as in Theorem~\ref{proposition_fluct_conditioning}, recalling that $(\hat{\boldsymbol{a}}^{\circ}, \hat{\boldsymbol{b}}^{\circ},\boldsymbol{N}^\circ)$ are deterministic functions of $\boldsymbol{Q}^{\circ} = \boldsymbol{M}$. Let $\Z^{\boldsymbol M,\textnormal{loc}}_\N$, $-\I^{\textnormal{loc}}$, and $\boldsymbol{\mu}^{\hat{\boldsymbol{m}},\textnormal{loc}}$ be the partition function, the energy functional, and the equilibrium measure of the localized ensemble, with $\hat{\boldsymbol{m}} = \frac{\boldsymbol{M}}{\N}$. Let us compare $\Z^{\boldsymbol M,\textnormal{loc}}_\N$ with $\Z^{\boldsymbol M}_\N$.
On the event $\mathcal A$, a configuration of the localized ensemble can be obtained from the corresponding configuration of the original ensemble by removing densely packed collections of particles in the saturated parts outside $\bigcup_{k=1}^K [\hat a^{\circ}_k, \hat{b}^{\circ}_k]$. This leads to a change in the partition function by a certain factor, as in \eqref{eq_changeofZ}. As we have already seen in \eqref{eq_x220}, \eqref{eq_x223}, \eqref{eq_x268} --- the argument is exactly the same --- the logarithm of this partition function matches the change in the energy functional $-\I$ in the leading order. In more details, we have:
\begin{equation}
\label{eq_x269}
 \log\left(\frac{\Z^{\boldsymbol{M}}_\N}{\Z^{\boldsymbol{M},\textnormal{loc}}_\N}\right) = \big(\I[\boldsymbol{\mu}^{\hat{\boldsymbol{m}}}]-\I^{\textnormal{loc}}[\boldsymbol{\mu}^{\hat{\boldsymbol{m}},\textnormal{loc}}]\big)\N^2  +\sum_{h=1}^H \theta_{h,h} (N_h-N^{\textnormal{loc}}_h) \log \N
+ \mathbbm{Rest}^{[1]}_1(\hat{\boldsymbol{m}})\,\N + o(1),
\end{equation}
Here, $N_h = \N \hat{n}_h$ is the total number of particles in $[a_h,b_h]$ for the original ensemble and $N^{\textnormal{loc}}_h = \sum_{k = k^-(h)}^{k^+(h)} N_k^{\circ}$ is the total number of particles in $\bigcup_{k = k^-(h)}^{k^+(h)} [\hat{a}_k^{\circ},\hat{b}_k^{\circ}]$ for the localized ensemble. Both numbers are deterministic after we conditioned by $\boldsymbol{Q}^{\circ} = \boldsymbol{M}$, and their difference is the number of particles in the saturated regions of the original ensemble that are outside the segments of the localized ensemble, \textit{cf.} \eqref{eq_x258}. We know that the remainder $\mathbbm{Rest}_1^{[1]}(\hat{\boldsymbol{m}})$ depends smoothly on $\hat{\boldsymbol{m}}$ and the remainder $o(1)$ is uniform over $\boldsymbol{M}$ satisfying \eqref{eq_x266}.

Next, we apply Theorem~\ref{Theorem_partition_multicut} to compute $\log \Z^{\boldsymbol{M},\textnormal{loc}}_\N$ as
\begin{equation}
\label{eq_x270}
 \log \Z^{\boldsymbol M,\textnormal{loc}}_\N=\I^{\textnormal{loc}}[\boldsymbol{\mu}^{\hat{\boldsymbol{m}},\textnormal{loc}} ]\, \N^2  + \left(\sum_{k=1}^K \theta_{h^k,h^k} N_k^\circ \right) \log \N + \mathbbm{Rest}_1^{[2]}(\hat{\boldsymbol{m}})\, \N +o(1).
\end{equation}

We proceed as in the proof of Theorem~\ref{Theorem_CLT_for_filling_fractions} and expand all $\hat{\boldsymbol{m}}$-dependent ingredients in \eqref{eq_x269} and \eqref{eq_x270} in Taylor series near the point $\hat{\boldsymbol{m}}= \frac{\boldsymbol{Q}^{\boldsymbol{\mu}}}{\N}$, up to $o(1)$ terms. For the expansion of $\I[\boldsymbol{\mu}^{\hat{\boldsymbol{m}}}]$ we notice that the linear terms necessarily vanish because the equilibrium measure $\boldsymbol{\mu}$ is the minimizer of $-\I$ precisely at $\hat{\boldsymbol{m}} = \frac{\boldsymbol{Q}^{\boldsymbol{\mu}}}{\N}$, due to the third point in Lemma~\ref{Lemma_center_of_expansion}. Comparing with the definition of the matrix $\Qu$ in \eqref{eq_discrete_covariance_sat}, we obtain
\begin{equation}
\label{eq_x271}
 \I[\boldsymbol{\mu}^{\hat{\boldsymbol{m}}}]\, \N^2 = \I[\boldsymbol{\mu}]\, \N^2 -\frac{1}{2}\Qu({\boldsymbol M} - \boldsymbol{Q}^{\boldsymbol{\mu}},{\boldsymbol M} - {\boldsymbol Q}^{\boldsymbol{\mu}})+o(1).
\end{equation}
Note that the first term in the right-hand side of \eqref{eq_x247} will eventually cancel between numerator and denominator of \eqref{eq_x265} and therefore can be ignored. Also note that the second term involves the non-normalized vector ${\boldsymbol M}$, rather than $\hat {\boldsymbol{m}}$.
Let
\[
\mathbbm{Rest}_1^{[3]}(\hat{\boldsymbol{m}}) = \mathbbm{Rest}_1^{[1]}(\hat{\boldsymbol{m}}) + \mathbbm{Rest}_1^{[2]}(\hat{\boldsymbol{m}})
\]
and expand near $\frac{\boldsymbol{Q}^{\boldsymbol{\mu}}}{\N}$:
\begin{equation}
\label{eq_x272}
\mathbbm{Rest}_1^{[3]}(\hat{\boldsymbol{m}}) \,\N = \mathbbm{Rest}_1^{[3]}\bigg(\frac{{\boldsymbol Q}^{\boldsymbol{\mu}}}{\N}\bigg)\,\N + \sum_{k=1}^K (M_k - Q^{\boldsymbol{\mu}}_k) \big( \partial_{\hat{m}_k}\mathbbm{Rest}_1^{[3]}(\hat{\boldsymbol{m}})\big)\big|_{\hat{\boldsymbol{m}} = \frac{\boldsymbol{Q}^{\boldsymbol{\mu}}}{\N}} +o(1).
\end{equation}
The first term in the right-hand side of \eqref{eq_x272} will eventually cancel between numerator and denominator of \eqref{eq_x265}. Combining \eqref{eq_x269}--\eqref{eq_x272} with \eqref{eq_x265} we get as $\N\rightarrow\infty$
\begin{equation}
\label{eq_x273}
\begin{split}
 \P({\boldsymbol Q}^\circ = \boldsymbol M) & = \frac{1}{\mathcal {Q}} \exp\Bigg(-\frac{1}{2}\Qu({\boldsymbol M} - {\boldsymbol Q}^{\boldsymbol{\mu}},{\boldsymbol M} - {\boldsymbol Q}^{\boldsymbol{\mu}}) + \bigg(\sum_{h=1}^H \theta_{h,h}\big(N_h(\boldsymbol{M})-N_h(\boldsymbol{Q}^{\boldsymbol{\mu}})\big)\bigg)\log \N \\
 & \quad \qquad +\sum_{k=1}^K (\theta_{h^k,h^k} \log \N+\mathbbm{Shift}_k) (M_k-Q^{\boldsymbol{\mu}}_k)+o(1) \Bigg).
\end{split}
\end{equation}
Here, $\mathscr{Q}$ is a normalization constant making the right-hand side a probability distribution, we set
\[
\mathbbm{Shift}_k := \partial_{\hat{m}_k}\mathbbm{Rest}_1^{[3]}(\hat{\boldsymbol{m}})\big|_{\hat{\boldsymbol{m}} = \frac{\boldsymbol{Q}^{\boldsymbol{\mu}}}{\N}},
\]
and $\boldsymbol{Q} \mapsto N_h(\boldsymbol{Q})$ is the affine function expressing the total number of particles in $[\hat{a}_h,\hat{b}_h]$, taken from \eqref{eq_x258}
\begin{equation}
\begin{split}
N_h(\boldsymbol{Q}) & = \sum_{k = k^-(h)}^{k^+(h)} N_k(\boldsymbol{Q}) + \frac{\N}{\theta_{h,h}}\bigg(\mathbbm{1}_{\amsmathbb{S}_h}(\hat{a}_h')\big(\hat{a}_{k^-(h)}(\boldsymbol{Q}) - \hat{a}_h\big) \\
& \quad + \sum_{k = k^-(h)}^{k^+(h) - 1} \mathbbm{1}_{\amsmathbb{S}_h}(\beta_k)\Big(\hat{a}_{k + 1}(\boldsymbol{Q}) - \hat{b}_k(\boldsymbol{Q}) - \frac{1}{\N}\Big) + \mathbbm{1}_{\amsmathbb{S}_h}(\hat{b}'_h)\big(\hat{b}_h - \hat{b}_{k^+(h)}(\boldsymbol{Q})\big)\bigg).
\end{split}
\end{equation}
The error $o(1)$ in \eqref{eq_x273} is uniform for $\boldsymbol{M}$ satisfying \eqref{eq_x266}. Comparing with the statement of Theorem~\ref{Theorem_CLT_for_filling_fractions_saturation} and the definition of the center $\boldsymbol u$ in \eqref{eq_x260}-\eqref{eq_x261} and noting that the closeness of probabilities implies the convergence in the sense of Definition~\ref{Definition_asymptotic_equivalence}, we are done.

Just like at the end of the proof of Theorem~\ref{Theorem_CLT_for_filling_fractions}, we remark that our analysis was restricted to $\boldsymbol M$ satisfying \eqref{eq_x266}, which is sufficient, because the probability of all other values of $\boldsymbol M$ is negligible both before and after taking the limit $\N \rightarrow \infty$.

\subsection{Proof of Theorem~\ref{Theorem_linear_statistics_fluctuation_sat}}

As the proof is similar to Section~\ref{sec:proofth105}, we omit some details. For notational simplicity we only consider the $M=1$ case; the general case is obtained by a similar argument. Thus, we fix $h\in [H]$ and a function $f$ defined in $\amsmathbb{A}_h^{\mathfrak{m}}$ and study the asymptotic distribution of the random variable
\begin{equation}
\label{eq_x276}
\underline{\textsf{Lin}}_{h}[f] = \sum_{i = 1}^{N_h} f\bigg(\frac{\ell_i^{h}}{\N}\bigg).
\end{equation}

The first step is to study the conditional distribution of $\underline{\textsf{Lin}}_h[f]$ given the values of the extended filling fractions $\overline{\boldsymbol N}^{\circ}$ or, equivalently, given the value of $\boldsymbol{Q}^{\circ}$. As in the proof of Theorem~\ref{Theorem_CLT_for_filling_fractions_saturation}, we use the segments $[\hat a^{\circ}_k, \hat{b}^{\circ}_k]$ indexed by $k\in [K]$ and by Corollary~\ref{Corollary_a_priori_0} and Theorem~\ref{Theorem_ldpsup} and \ref{Theorem_ldsaturated}, it is sufficient to work on the event where $\boldsymbol{Q}^{\circ}$ satisfies \eqref{eq_x266} and each segment of the complement $\bigcup_{k=k^-(h)}^{k^+(h)} [\N \hat a^{\circ}_k, \N \hat{b}^{\circ}_k]$ in $[a_h,b_h]$ either has no particles or no holes as in Definition~\ref{Definition_no_hole}. Let $\underline{\textsf{Lin}}_h^{\hat{\boldsymbol{m}}}[f]$ be the random variable conditioned on this event and on $\boldsymbol{Q}^{\circ}=\boldsymbol M = \N \hat{\boldsymbol{m}}$ for a fixed value of $\boldsymbol M$ satisfying \eqref{eq_x266}. We write separately the contribution of the densely packed particles outside $\bigcup_{k = k^-(h)}^{k^+(h)}[\N\hat a^{\circ}_k, \N\hat{b}^{\circ}_k]$ in $[a_h,b_h]$:
\begin{equation}
\label{Lin3km2}
\begin{split}
\underline{\textsf{Lin}}_h^{\hat{\boldsymbol{m}}}[f] & = \sum_{k = k^-(h)}^{k^+(h)} \textsf{Lin}_k^{\hat{\boldsymbol{m}}}[f] + \mathbbm{1}_{\amsmathbb{S}_h}(\hat{a}'_h)\Bigg[f(\hat{a}_h) + \cdots + f\bigg(\hat{a}_{k^-(h)}(\boldsymbol{M}) - \frac{\theta_{h,h}}{\N}\bigg)\Bigg] \\
& \quad + \sum_{k = k^-(h)}^{k^+(h) - 1} \mathbbm{1}_{\amsmathbb{S}_h}(\beta_k)\Bigg[f\bigg(\hat{b}_k(\boldsymbol{M}) + \frac{\theta_{h,h}}{\N}\bigg) + \cdots + f\bigg(\hat{a}_{k + 1}(\boldsymbol{M}) - \frac{\theta_{h,h}}{\N}\bigg)\Bigg] \\
& \quad + \mathbbm{1}_{\amsmathbb{S}_h}(\hat{b}'_h) \Bigg[f\bigg(\hat{b}_{k^+(h)}(\boldsymbol{M}) + \frac{\theta_{h,h}}{\N}\bigg) + \cdots + f(\hat{b}_h)\Bigg],
\end{split}
\end{equation}
where $\cdots$ indicate summations over arguments of $f$ with increments of $\frac{\theta_{h,h}}{\N}$ between the two given extreme values. The random variables
\begin{equation}
\label{eq_x3001} \forall k \in \llbracket k^-(h),k^+(h)\rrbracket \qquad \textsf{Lin}_k^{\hat{\boldsymbol{m}}}[f] := \sum_{i = 1}^{N_h} \mathbbm{1}_{[\N\hat{a}_k(\boldsymbol{M}),\N\hat{b}_k(\boldsymbol{M})]}(\ell_i^h)\,f\bigg(\frac{\ell_i^h}{\N}\bigg)
\end{equation}
make sense and can be studied in the ensemble localized to $\bigcup_{g \neq h} [\hat{a}_g,\hat{b}_g] \cup \bigcup_{k = k^-(h)}^{k^+(h)} [\hat{a}_k^{\circ}(\boldsymbol{M}),\hat{b}_k^{\circ}(\boldsymbol{M})]$. For them, we can use Corollary~\ref{Corollary_CLT_relaxed}, which requires $f$ to be holomorphic in a $\N$-independent complex neighborhood of the bands of the $h$-th segment in the localized ensemble and to be uniformly bounded in this neighborhood. This explains the role of $(\overline{\alpha}_k,\overline{\beta}_k)$ and its complex neighborhood $\amsmathbb{K}_h$ in the assumptions of Theorem~\ref{Theorem_linear_statistics_fluctuation_sat}. We conclude from Corollary~\ref{Corollary_CLT_relaxed} that
\begin{equation}
\label{eq_x3000}
\sum_{k = k^-(h)}^{k^+(h)} \Big(\textsf{Lin}_{k}^{\hat{\boldsymbol{m}}}[f] - \amsmathbb{E}\big[\textsf{Lin}_{k}^{\hat{\boldsymbol{m}}}[f]\big]\Big)
\end{equation}
is asymptotically Gaussian and its variance agrees with \eqref{eq_x264} --- recall that the expectation value in \eqref{eq_x3000} is conditional to $\boldsymbol{Q}^{\circ} = \boldsymbol{M}$. The dependence in $\hat{\boldsymbol{m}}$ in this Gaussian distribution disappears as $\N\rightarrow\infty$, because the functions $\mathcal{F}_{k_1,k_2}$ in \eqref{eq_x264} depend only on the endpoints of the bands and the interaction matrix $\boldsymbol{\Theta}$, while the dependence of the endpoints on $\boldsymbol M$ is negligible due to \eqref{eq_x266} and Theorem~\ref{Theorem_differentiability_full}. This gives the random part of the first term $\textnormal{\textsf{\textbf{Gau\ss{}}}}$ in \eqref{eq_x274} and shows independence of the two terms in \eqref{eq_x274}.

The $\boldsymbol{M}$-dependent asymptotics of the conditional expectation value in \eqref{eq_x3000} is described in the second part of Corollary~\ref{Corollary_CLT_relaxed} applied to the localized ensemble. Namely, \eqref{eq_x164} yields
\begin{equation}
\label{Elink1}
\amsmathbb{E}\big[\textsf{Lin}_k^{\hat{\boldsymbol{m}}}[f]\big] = \N \int_{\hat{a}'_k(\boldsymbol{M})}^{\hat{b}'_k(\boldsymbol{M})} f(x)\,\mu_{h}^{\hat{\boldsymbol{m}},\textnormal{loc}}(x) \dd x + \mathbbm{Rest}[f] + O(\N^{-\frac{3}{2} + \varepsilon}),
\end{equation}
where the remainder $\mathbbm{Rest}[f]$ depends smoothly on $\hat{\boldsymbol{m}}$ near $\frac{\boldsymbol{Q}^{\boldsymbol{\mu}}}{\N}$ and the $O(\cdots)$ remainder is uniform over $\boldsymbol{M}$ satisfying \eqref{eq_x266}. In fact, due to \eqref{eq_x266}, the dependence in $\hat{\boldsymbol{m}}$ of the term $\mathbbm{Rest}[f]$ becomes negligible as $\N \rightarrow \infty$. Besides, $\boldsymbol{\mu}^{\hat{\boldsymbol{m}},\textnormal{loc}}$ is the equilibrium measure in the conditioned and localized ensemble. We would like to rewrite the corresponding term in terms of the equilibrium measure $\boldsymbol{\mu}^{\hat{\boldsymbol{m}}}$ of the ensemble before localization but after conditioning to fixed extended $4K$-tuple $\overline{\boldsymbol{N}}(\boldsymbol{M})$ --- as introduced in \eqref{eq_discrete_covariance_sat}. As in the proofs of Theorems~\ref{Theorem_linear_statistics_fluctuation_ff} and \ref{Theorem_CLT_for_filling_fractions_saturation}, we have
\begin{equation}
\label{Elinkm2}
\amsmathbb{E}\big[\textsf{Lin}_k^{\hat{\boldsymbol{m}}}[f]\big] = \N \int_{\hat{a}_k'(\boldsymbol{M})}^{\hat{b}_k'(\boldsymbol{M})} f(x)\,\mu_{k}^{\hat{\boldsymbol{m}}}(x) \dd x + \mathbbm{Rest}[f] + \mathbbm{Rest}'[f] + o(1).
\end{equation}
We recall that the key difference between $\mu_{k}^{\hat{\boldsymbol{m}},\textnormal{loc}}$ and $\mu_{k}^{\hat{\boldsymbol{m}}}$ is in the subleading terms \eqref{eq_weight_subleading_1}, \eqref{eq_weight_subleading_2}, and \eqref{eq_weight_subleading_3} in the potential, arising from our treatment of the Gamma functions in the localization procedure. This is the origin of the term $\mathbbm{Rest}'[f]$ in \eqref{Elinkm2}, and its sum with the $\mathbbm{Rest}[f]$ term from \eqref{Elink1} gives the mean of $\textnormal{\textsf{\textbf{Gau\ss{}}}}$ in \eqref{eq_x274}.

It remains to identify the second term $\textsf{Shift}_h[f]$ in \eqref{eq_x274}. For this, we approximate the three Riemann sums in \eqref{Lin3km2} with the integral between shifted endpoints, using the three formulae of Lemma~\ref{lem:EulerMaclaurin} and the assumption that $f$ is twice continuously differentiable in $\amsmathbb{A}_h^{\mathfrak{m}}$.  The integral contribution combines with the integral  terms of \eqref{Elinkm2} to reconstruct
\begin{equation}
\label{theintegralterm}
\N \int_{\hat{a}_h'}^{\hat{b}'_h} f(x) \mu_{h}^{\hat{\boldsymbol{m}}}(x) \dd x
\end{equation}
and there remains a contribution
\begin{equation}
\label{threthreesum}
\begin{split}
&  \frac{\theta_{h,h} - 1}{2\theta_{h,h}} \bigg[\mathbbm{1}_{\amsmathbb{S}_h}(\hat{a}'_h) \Big(f\big(\hat{a}_{k^-(h)}'(\boldsymbol{M})\big) - f(\hat{a}'_h)\Big) + \sum_{k = k^-(h)}^{k^+(h) - 1} \mathbbm{1}_{\amsmathbb{S}_h}(\beta_k)\Big(f\big(\hat{b}_{k}'(\boldsymbol{M})\big) + f\big(\hat{a}_{k + 1}'(\boldsymbol{M})\big)\Big) \\
& \quad \qquad \quad + \mathbbm{1}_{\amsmathbb{S}_h}(\hat{b}'_h) \Big(f\big(\hat{b}_{k^+(h)}'(\boldsymbol{M})\big)  - f(\hat{b}'_h)\Big)\bigg]  + O\bigg(\frac{1}{\N}\bigg).
\end{split}
\end{equation}

We now proceed to remove the conditioning on $\boldsymbol{Q}^{\circ} = \boldsymbol{M}$. Then for each for $k \in \llbracket k^-(h),k^+(h) \rrbracket$ the quantity  $\underline{\mathsf{Lin}}_{k}^{\hat{\boldsymbol{m}}}[f]$ becomes a random variable. We use the asymptotic distribution of $\boldsymbol{Q}^{\circ}$ known to be $\textnormal{\textsf{\textbf{Gau\ss{}}}}_{\amsmathbb{Z}}$ from Theorem~\ref{Theorem_CLT_for_filling_fractions_saturation} and perform a Taylor expansion of the right-hand side of \eqref{Elink1} around the point $\boldsymbol{M} = \boldsymbol{Q}^{\boldsymbol{\mu}}$, recalling that at this point $\boldsymbol{\mu}^{\hat{\boldsymbol{m}}}$ coincides with the equilibrium measure $\boldsymbol{\mu}$ of the original ensemble. The integral term \eqref{theintegralterm} then give rises to
\[
\N \int_{\hat{a}_h'}^{\hat{b}'_h} f(x)\mu_h(x)\dd x +  \sum_{l = 1}^{K} \big(\textnormal{\textsf{Gau\ss{}}}_{\amsmathbb{Z},l} - Q^{\boldsymbol{\mu}}_l\big) \int_{\hat{a}_h'}^{\hat{b}'_h} f(x) \big(\partial_{\hat{p}_l} \mu_h^{\hat{\boldsymbol{p}}}(x)\big)\big|_{\hat{\boldsymbol{p}} = \frac{\boldsymbol{Q}^{\boldsymbol{\mu}}}{\N}} \,\dd x + o(1)
\]
as $\N \rightarrow \infty$. The first-order term in the Taylor expansion gives rise to the contributions of the first line of $\mathsf{Shift}_h[f]$ in Definition~\ref{DhfL}, and contributions beyond first order are negligible as $\N \rightarrow \infty$. Eventually, plugging the discrete Gaussian asymptotic distribution of $\boldsymbol{Q}^{\circ} = \boldsymbol{M}$ in the three sums in \eqref{threthreesum} gives rise to the remaining terms in $\mathsf{Shift}_h[f]$ of \eqref{def_shiftsfh}, and we arrive to the claimed asymptotic distribution for $\underline{\textsf{Lin}}_h[f]$.

s
\chapter{Application to tiling models}
\label{Chap11}

In this section we study tilings of polygonal domains drawn on a triangular lattice. The allowed tiles are lozenges {\scalebox{0.16}{\includegraphics{lozenge_hor.pdf}}}, {\scalebox{0.16}{\includegraphics{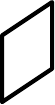}}} or
{\scalebox{0.16}{\includegraphics{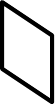}}} (\textit{cf.} Figure~\ref{Fig_trapezoid_tiling}), in other words they are rhombi formed by gluing two adjacent triangles. Our task is to show how the technology developed in the previous chapters leads to results describing macroscopic fluctuations of uniformly random lozenge tilings. We refer to \cite{Vadimlecture} for a general introduction to random tilings.

We start in Section~\ref{Section_simple_domains} from discussing lozenge tilings of the simplest domains: trapezoids, hexagons, hexagons with holes. In Section~\ref{Section_gluing_def} we introduce a more general class of domains which can be obtained by gluing several trapezoids together along a single vertical axis. In Section~\ref{sec:gen} we show how to apply our theorems to uniformly random lozenge tilings of these domains and deduce that the fluctuations of the filling fractions (which can be identified with relative heights of different components of the boundary of tiled domains) are asymptotically given by discrete Gaussian random variables. In Section~\ref{Section_KO_conjecture} we investigate the asymptotics of the two-dimensional field of fluctuations of the height functions of tilings. We state the Kenyon--Okounkov conjecture linking the fluctuations to the Gaussian free field in an appropriate complex structure, as well as its novel version for tilings on non-orientable surfaces. In Section~\ref{Section_KO_proofs} we prove these statements for the domains obtained by gluing trapezoids.

The connection between discrete ensembles with $H=1$ and $\theta=1$ and lozenge tilings was first noticed in \cite{CLP} who studied tilings of hexagons; later, discrete ensembles were also used in \cite{Johansson2,BKMM,Gor,BG,nordenstam2009interlaced,duits2018global,BuGo3}. A parallel story exists for domino tilings, see, \textit{e.g.}, \cite{Johansson2,fleming2011interlaced,Colomo_Pronko_2013_third_order,Colomo_Pronko_2015_thermodynamics,BuGo3,Colomo_Pronko_Sportiello_2019_arctic}. While we do not go in the domino direction in this book, we expect our methods to be relevant there as well.

\section{Lozenge tilings in simple domains}
\label{Section_simple_domains}
In this section we discuss the three simplest classes of domains to be tiled: trapezoids, which are basic building blocks for all our theorems; hexagons, which are the most well-known and well-studied domains; and hexagons with holes, which are the simplest domains where our methods provide new results that could not be reached by any of the methods previously existing in the literature.

\subsection{Lozenge tilings of trapezoids}
\label{Section_trapezoid}
\begin{figure}[t]
\includegraphics[width=0.35\linewidth]{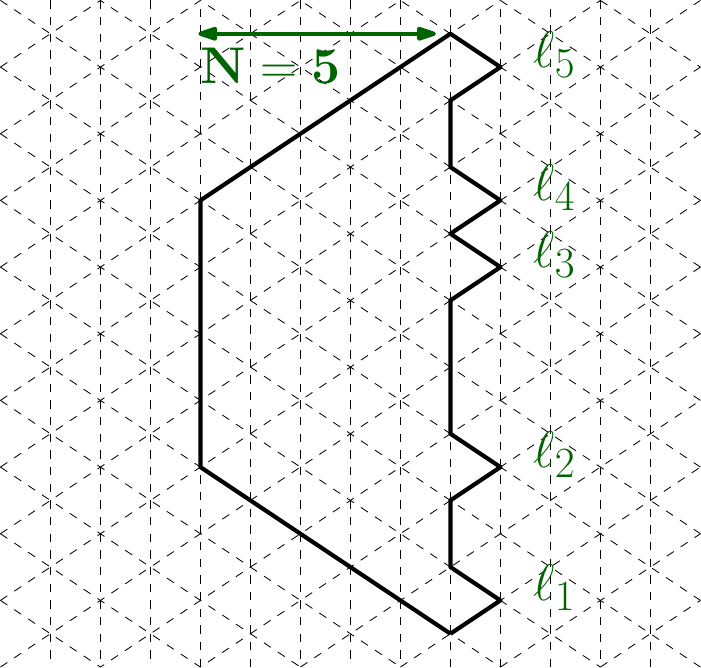} \hfill \includegraphics[width=0.17\linewidth]{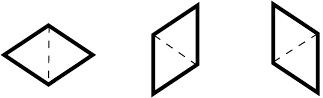} \hfill \includegraphics[width=0.35\linewidth]{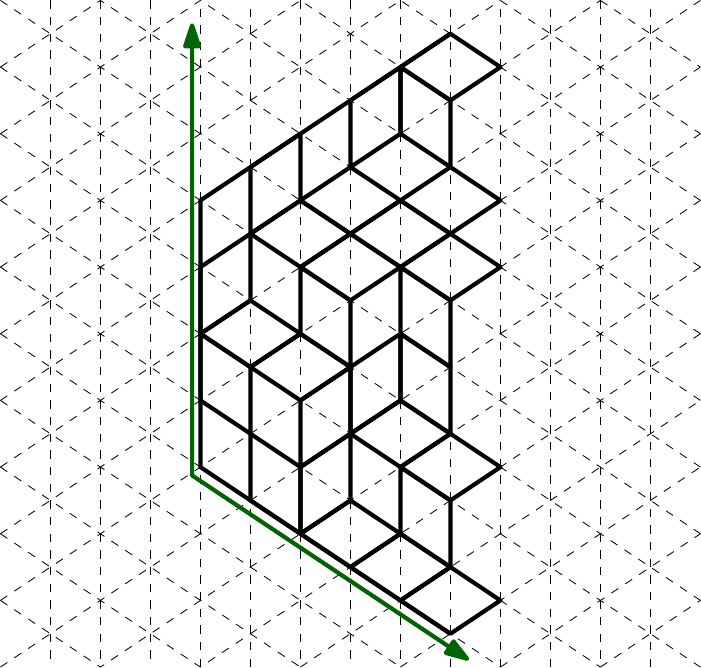} \caption{Trapezoid domain with dents at positions $1<3<6<7<9$, three types of lozenges, and one possible tiling. \label{Fig_trapezoid_tiling}}
\end{figure}

A trapezoid is a polygonal domain on the triangular lattice, parameterized by its width $N$ and integer coordinates $\ell_1<\ell_2<\dots<\ell_N$ of ``dents'' sticking out of its long base. The domain has three straight boundaries: the short base, the upper and lower sides, and the long base formed by dents, see Figure~\ref{Fig_trapezoid_tiling}. When we tile this domain, each of the $N$ dents is necessarily covered by a horizontal lozenge.

Our interest in lozenge tilings of trapezoids is two-fold. First, the total number of tilings turns out to be explicit and given by a formula reminiscent of the weight for discrete ensembles \eqref{eq_general_measure}.

 \begin{proposition}\label{Proposition_number_tilings_trapezoid}
The total number of lozenge tilings of the trapezoid of width $N$ and with dents at $\ell_1 < \ell_2 < \dots < \ell_N$ is
\begin{equation}
\label{eq_Weyl_dimension_main}
\prod_{1 \leq i < j \leq N}{\frac{\ell_j - \ell_i}{j - i}}.
\end{equation}
\end{proposition}
\begin{proof} This statement is due to \cite{CLP}. The expression \eqref{eq_Weyl_dimension_main} coincides with the Weyl formula for the dimension of an irreducible representation of the unitary group $\textnormal{U}(N)$ with a given highest weight: tilings themselves are in bijection with Gelfand--Tsetlin patterns which enumerate the bases in irreducible representations, \textit{cf.}\ \cite[Section 2]{BP_lectures}.
\end{proof}

Second, uniformly random tilings of trapezoids are very well-understood, due to a variety of methods which can be applied for their study: determinantal point processes \cite{Petrov_Airy,petrov2015asymptotics,gorin2017bulk}; Schur generating functions \cite{BuGo2,BuGo3}; dynamical loop equations \cite{huang2020height,gorin2024dynamical}. Our strategy for analyzing tilings of more complicated domains proceeds by gluing them out of trapezoids, and we will use results on trapezoids as building blocks in Sections \ref{Section_KO_conjecture}--\ref{Section_KO_proofs}.

\subsection{Lozenge tilings of the hexagon}
\label{Section_Hexagon}

 \begin{figure}[t]
\includegraphics[width=0.4\linewidth]{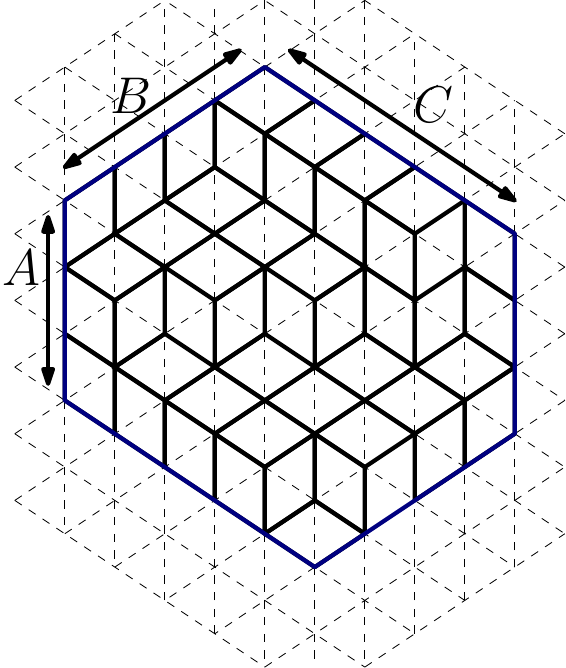} \hfill \includegraphics[width=0.4\linewidth]{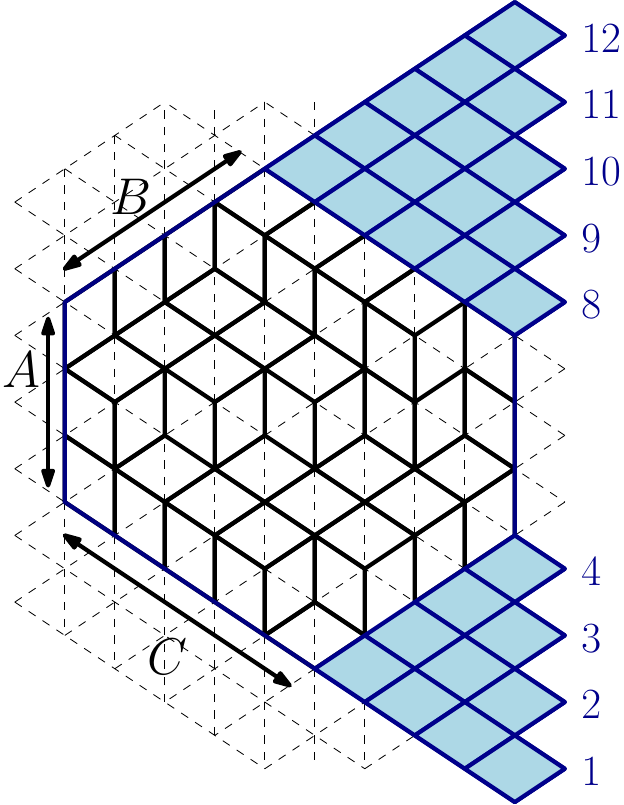} \caption{A lozenge tiling of $3\times 4 \times 5$ hexagon and of the corresponding trapezoid. \label{Fig_hexagon_trapezoid}}
\end{figure}

\begin{figure}[t]
\center \scalebox{0.7}{\includegraphics{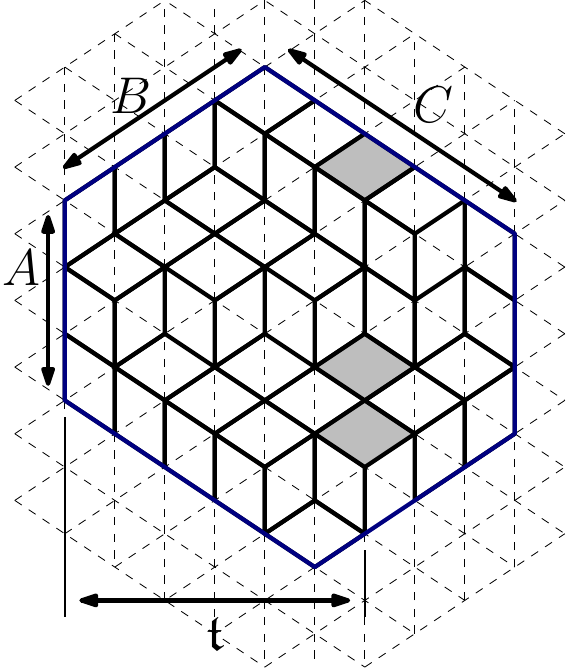}} \qquad
\scalebox{0.7}{\includegraphics{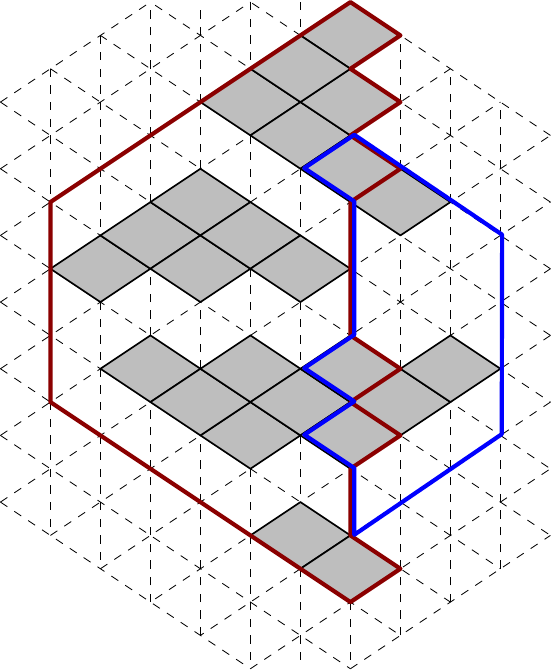}} \caption{Left panel: Lozenge tiling
of the $3\times 4\times 5$ hexagon and three horizontal lozenges on the sixth from the
left vertical line. Right panel: two tilings of trapezoids (inside blue and red
contours) corresponding to each tiling. \label{Fig_tiling_hex_section}}
\end{figure}

Our next domain is $A\times B\times C$ hexagon, as in Figure~\ref{Fig_hexagon_trapezoid}. Two points of view on lozenge tilings of the hexagon are important for us: the first one treats hexagon tilings as a particular case of trapezoid tilings; the second one considers the discrete ensemble \eqref{eq_general_measure} appearing as the marginal of a uniformly random tiling along a vertical section.

\begin{lemma} \label{Lemma_hexagon_as_trapezoid}
 Lozenge tilings of the $A\times B\times C$ hexagon are in bijection with lozenge tilings of the trapezoid of width $B+C$ and dents at $1<2<\dots<B<A+B+1<A+B+2<\dots<A+B+C$.
\end{lemma}
\begin{proof}
 The desired trapezoid of width $B+C$ is shown in the right panel of Figure~\ref{Fig_hexagon_trapezoid}. Observe that the dents at positions $1<2<\dots<B$ create in each tiling of this trapezoid a macroscopic triangle densely packed with $\frac{B(B+1)}{2}$ horizontal lozenges. Similarly, the dents at positions $A+B+1<A+B+2<\dots<A+B+C$ create a macroscopic triangle densely packed with $\frac{C(C+1)}{2}$ horizontal lozenges. If we remove all these $\frac{B(B+1)}{2}+\frac{C(C+1)}{2}$ horizontal lozenges, then each tiling of the trapezoid turns into a tiling of the $A\times B\times C$ hexagon.
\end{proof}

For the second point of view, we consider a uniformly random tiling and focus on positions of the horizontal lozenges on the vertical line at distance $\mathfrak{t}$ from the left border, as in Figure~\ref{Fig_tiling_hex_section}. The number of horizontal lozenges on the line depends on $\mathfrak{t}$ and equals
\begin{equation}
\label{eq_N_hexagon}
N=\left\{
 \begin{array}{lll}\mathfrak{t}&& \textnormal{if}\,\, 0 < \mathfrak{t}\leq \min(B,C),\\ \min(B,C)&& \textnormal{if}\,\, \min(B,C)\leq \mathfrak{t}\leq \max(B,C),\\B+C-\mathfrak{t}&& \textnormal{if}\,\,\max(B,C)\leq \mathfrak{t} < B+C.
 \end{array}\right.
\end{equation}
Let us introduce the coordinate system, so that the minimal possible coordinate of the horizontal lozenge on this vertical line is $1$ and the maximal possible coordinate is the length of the section of the hexagon by the vertical line explicitly given by
\begin{equation}
\label{eq_L_hexagon}
L=
\left\{
 \begin{array}{lll}A+\mathfrak{t}&& \textnormal{if}\,\, 0 < \mathfrak{t}\leq \min(B,C),\\ A+\min(B,C)&& \textnormal{if}\,\, \min(B,C)\leq \mathfrak{t}\leq \max(B,C),\\A+B+C-\mathfrak{t}&& \textnormal{if}\,\,\max(B,C)\leq \mathfrak{t}<B+C.
 \end{array}\right.
\end{equation}
It was first noticed in \cite{CLP} that the distribution of the $N$ horizontal lozenges along the vertical section is a particular case of the discrete ensemble \eqref{eq_general_measure}.

\begin{lemma}\label{lem:discretemodelhex}
 Consider a uniformly random tiling of $A\times B\times C$ hexagon. Using the notations \eqref{eq_N_hexagon} and \eqref{eq_L_hexagon}, the $N$ horizontal lozenges on the vertical section at distance $\mathfrak{t}$ from the left border with coordinates $1\leq \ell_1<\ell_2<\dots<\ell_N\leq L$ have probability distribution
\begin{equation}
\label{eq_Hahn_ensemble}
 \amsmathbb P_N(\ell_1,\ldots,\ell_N)=\frac{1}{\mathscr{Z}_N} \cdot \prod_{1 \leq i<j \leq N} (\ell_j-\ell_i)^2 \cdot \prod_{i=1}^N w(\ell_i),
\end{equation}
with the weight
\begin{equation}
\label{eq_Hahn_ensemble_2}
w(\ell)=
 \left\{\begin{array}{lll}
  (\ell)_{C-\mathfrak{t}}\cdot (A+\mathfrak{t}+1-\ell)_{B-\mathfrak{t}}&& \textnormal{if}\,\,0 < \mathfrak{t}\leq \min(B,C),\\
  (\ell)_{C-\mathfrak{t}}\cdot (A+B+1-\ell)_{\mathfrak{t}-B} && \textnormal{if}\,\,B\leq \mathfrak{t}\leq C,\\
  (\ell)_{\mathfrak{t}-C}\cdot (A+C+1-\ell)_{B-\mathfrak{t}}&& \textnormal{if}\,\,C\leq \mathfrak{t}\leq B,\\
  (\ell)_{t-C}\cdot (A+B+C+1-\mathfrak{t}-\ell)_{\mathfrak{t}-B} && \textnormal{if}\,\, \max(B,C)\leq \mathfrak{t} < B+C,\end{array}\right.
\end{equation}
where\label{index:Poch} $(a)_n$ is the Pochhammer symbol,
\[
(a)_n=a(a+1)\cdots (a+n-1)=\frac{\Gamma(a+n)}{\Gamma(a)},
\]
and $\mathscr{Z}_N$ is a normalizing
constant.
\end{lemma}
\begin{proof}
 $\amsmathbb P_N(\ell_1,\ldots,\ell_N)$ is proportional to the product of the total numbers of lozenge tilings to the left and to the right of the vertical line with fixed horizontal lozenges $\ell_1,\dots,\ell_N$. Repeating the argument of Lemma~\ref{Lemma_hexagon_as_trapezoid}, we identify tilings to the left with tilings of a trapezoid of width $\mathfrak{t}$ and tilings to the right with tilings of a trapezoid of width $B+C-\mathfrak{t}$, as in Figure~\ref{Fig_tiling_hex_section}. It remains to use Proposition~\ref{Proposition_number_tilings_trapezoid} twice and multiply the results.
\end{proof}

Upon closer inspection, one realizes that \eqref{eq_Hahn_ensemble} is an instance of the $zw$-measure of Chapter~\ref{Chapterzw}, with $\theta=1$. More specifically, after we multiply the weight $w(\ell)$ by a power of $\N$ not depending on $\ell$ (and including the same factor in the partition function $\mathscr{Z}_{N}$), we recognize the $zw$-discrete ensemble of the second kind (Section~\ref{S:sec_second_type}) with parameters $B_2 < A_1 < B_1 < A_2$ given in the table below.

\begin{equation}
\label{BCzwpar}
\begin{array}{|l||c|c|c|c|}
\hline
\textnormal{Case} & A_1 & A_2 & B_1 & B_2 \\
\hline\hline
 0 < \mathfrak{t} \leq \min(B,C) & 1 & A + B + 1 & A + \mathfrak{t}& \mathfrak{t} - C \\
 \hline
 B \leq \mathfrak{t} \leq C & 1 & A + \mathfrak{t} + 1 & A + B & \mathfrak{t} - C \\
 \hline
 C \leq \mathfrak{t} \leq B & 1 & A + B + C - \mathfrak{t} + 1 & A + C & C - \mathfrak{t} \\
 \hline
 \max(B,C) \leq \mathfrak{t} < B + C & 1 & A + C + 1 & A + B + C - \mathfrak{t} & C - \mathfrak{t} \\
 \hline
 \end{array}
 \end{equation}

\medskip

\begin{figure}[t]
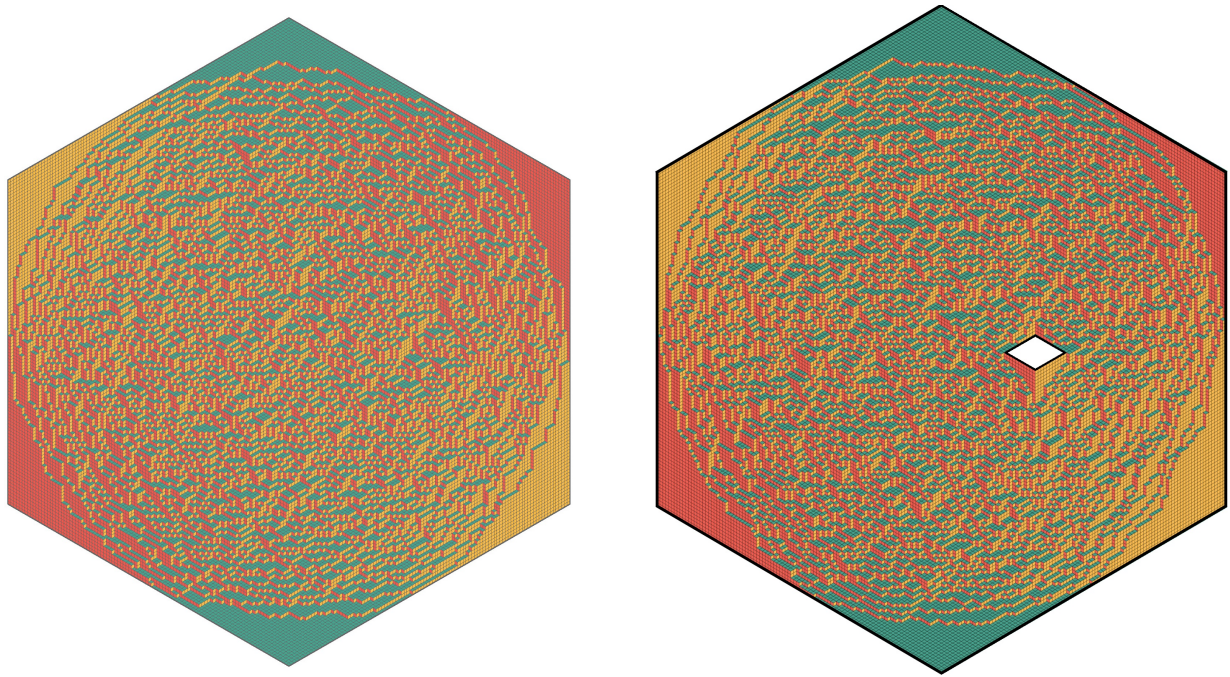

\includegraphics[width=0.47\linewidth]{hex100.jpg} \hfill \includegraphics[width=0.47\linewidth]{hex100_hole.jpg} \caption{Two simulations of uniformly random tilings with lozenges shown in three colors. Left panel: $100\times 100\times 100$ hexagon. Right panel: same hexagon with a small hole. \label{Fig_hex_simulation}}
\end{figure}

\begin{lemma}
\label{lem:checktil} Suppose that $A,B,C,\mathfrak{t}$ depend on $\N\rightarrow\infty$ and denote $\frac{A}{\N}=\hat A$, $\frac{B}{\N}=\hat B$, $\frac{C}{\N}=\hat C$, $\frac{\mathfrak{t}}{\N}=\hat{\mathfrak{t}}$. Then Assumptions \ref{Assumptions_Theta}, \ref{Assumptions_basic} and \ref{Assumptions_analyticity} hold for \eqref{eq_Hahn_ensemble} and the (simplified) equilibrium measure is given by the explicit formula of Theorem~\ref{Proposition_ZW_equilibrium_measure}. Assumption~\ref{Assumptions_offcrit} also hold for $\N$ large enough provided the distances between $\hat{\mathfrak{t}} \in (0,\hat{B} + \hat{C})$ and each of the following values
\begin{equation}
\label{6value}
0,\quad \frac{\hat{A}\hat{C}}{\hat{A} + \hat{B}},\quad \frac{\hat{A}\hat{B}}{\hat{A} + \hat{C}},\quad \hat{B} + \frac{\hat{B}\hat{C}}{\hat{A} + \hat{B}},\quad \hat{C} + \frac{\hat{B}\hat{C}}{\hat{A} + \hat{C}},\quad \hat{B} + \hat{C}
\end{equation}
remain bounded from below by a positive constant independent of $\N$.
\end{lemma}
\begin{proof} Except for Assumption~\ref{Assumptions_offcrit}, this was already verified for the $zw$-ensemble in Lemma~\ref{Checkzw} (the parameters parameters \eqref{BCzwpar} satisfy the assumptions of that lemma). The off-criticality Assumption~\ref{Assumptions_offcrit} was verified in Corollary~\ref{Corollary_ZW_off-critical_1} under the extra condition that
\begin{equation}
\label{awayfromoff}\big|(\hat{B}_1 - \hat{A}_1 - \hat{n})(\hat{B}_2 - \hat{A}_1 - \hat{n}) + \hat{n}(\hat{A}_2 - \hat{A}_1)\big| \cdot \big|(\hat{B}_1 - \hat{A}_1 - \hat{n})(\hat{B}_1 - \hat{A}_2 - \hat{n}) - \hat{n}(\hat{B}_2 - \hat{B}_1)\big| > \varepsilon
\end{equation}
for some constant $\varepsilon > 0$ independent of $\N$. Inserting \eqref{BCzwpar} for each of the four cases and neglecting the $O(1/\N)$ terms, we see that \eqref{awayfromoff} holds when $\hat{\mathfrak{t}}$ remains away from the values given in \eqref{6value}.
\end{proof}
Therefore, all our main results hold for the uniformly random lozenge tilings of hexagons, for generic sections (with the degenerate situations of Definitions~\ref{Definition_ZW_set} and \ref{Definition_off_critical_set_ZW} excluded). In particular, for the expansion of correlators we have Theorem~\ref{Theorem_correlators_expansion_relaxed}, in fact in the stronger form of Theorem~\ref{Theorem_correlators_expansion_relaxed_theta1} since $\theta =1$, and for the expansion of the partition function we have Theorem~\ref{Theorem_partition_one_band}. The left panel of Figure~\ref{Fig_hex_simulation} shows a simulation. It is known (see \cite{CLP} or \cite{Vadimlecture}) that the endpoints of the bands form ( the inscribed ellipse in the hexagon as we vary the vertical section. The critical vertical sections are those passing through the tangency points of this ellipse and in these situations (corresponding to the values \eqref{6value} for $\mathfrak{t}$) our theorems do not apply. For all other vertical sections the assumptions and theorems hold.

 We remark that in the setting of lozenge tilings of hexagons, multiple other approaches can be used to obtain the same results as ours, such as those in \cite{petrov2015asymptotics,breuer2017central,BuGo2,huang2020height,gorin2024dynamical}.

\subsection{Hexagon with a hole}

\label{Section_hex_hole}

\begin{figure}[t]
\center \scalebox{0.7}{\includegraphics{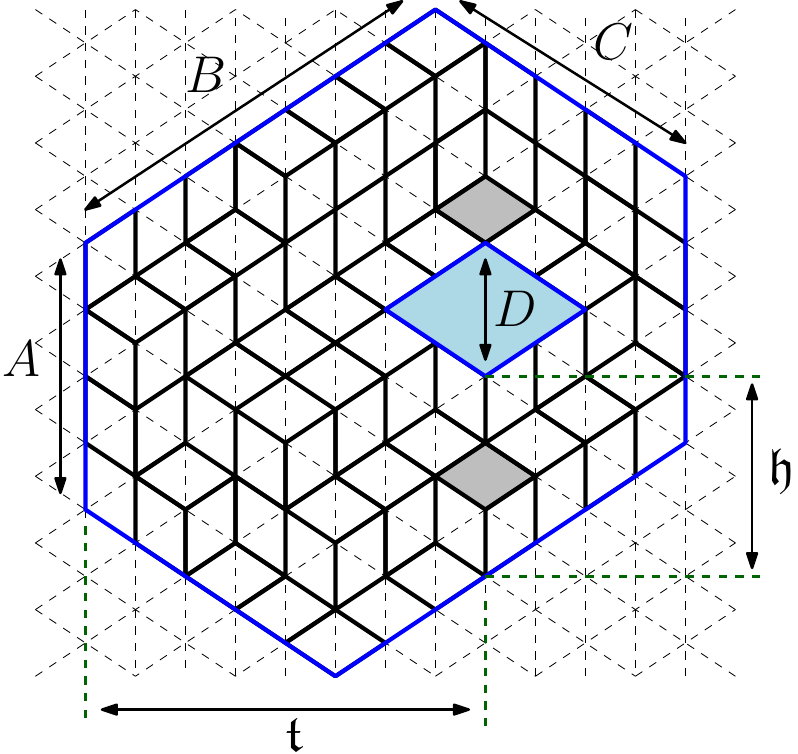}} \caption{Lozenge tiling of
the $4\times 7\times 5$ hexagon with a rhombic $2\times 2$ hole (shown in blue). The
center of the hole is at distance $\mathfrak{t}$ from the left side of the hexagon.
The remaining horizontal lozenges on $\mathfrak{t}$-th vertical line are shown in
gray. \label{Fig_tiling_hex_2}}
\end{figure}

Our next stop is the domain obtained by cutting a rhombic $D\times D$ hole in the $A\times B\times C$ hexagon, as shown in the right panel of Figure~\ref{Fig_hex_simulation} and in Figure~\ref{Fig_tiling_hex_2}. We assume that the bottom point of the hole is at fixed distance $\mathfrak{t}$ from the left
side of the hexagon and at fixed distance $\mathfrak{h}$ from the bottom side of the
hexagon along the $\mathfrak{t}$-th vertical line. We are interested in the asymptotic behavior of uniformly random lozenge tiling as
 $A,B,C,D,\mathfrak{t},\mathfrak{h} \rightarrow \infty$.

Let $\amsmathbb P_N$ be the probability distribution of the horizontal lozenges --- outside the
hole --- on the $\mathfrak{t}$-th vertical line induced by the uniform measure on all tilings of the
hexagon with the hole. Deciding that the lowest possible position for a horizontal lozenge is $1$, then the highest possible position is $L$ as in \eqref{eq_L_hexagon}, and the positions of the $N$ horizontal lozenges are specified by
\[
1 \leq \ell_1 < \cdots < \ell_N \leq L,
\]
where $N$ is given by the difference of \eqref{eq_N_hexagon} and $D$, \textit{i.e.}
\begin{equation}
\label{eq_N_hexagon_hole}
N=
 \left\{\begin{array}{lll}\mathfrak{t}-D&& \textnormal{if}\,\,0<\mathfrak{t}\leq \min(B,C),\\ \min(B,C)-D&& \textnormal{if}\,\,\min(B,C)\leq \mathfrak{t}\leq \max(B,C),\\B+C-D - \mathfrak{t} && \textnormal{if}\,\,\max(B,C)\leq \mathfrak{t}<B+C.\end{array}\right.
\end{equation}
It was noticed in \cite[Section 9.2]{BGG} that $\amsmathbb P_N$ is a particular case of the discrete ensemble of Section~\ref{Section_general_model} with $H=2$.

\begin{lemma}
 \label{lem:discretemodelhexhole}Consider a uniformly random tilings of $A\times B\times C$ hexagon with $D\times D$ hole (the bottom point of the hole is at fixed distance $\mathfrak{t}$ from the left side of the hexagon and at fixed distance $\mathfrak{h}$ from the bottom side).
 Using the notations \eqref{eq_N_hexagon_hole} and \eqref{eq_L_hexagon}, the $N$ horizontal lozenges on the vertical section at distance $\mathfrak{t}$ from the left border with coordinates $1\leq \ell_1<\ell_2<\dots<\ell_N\leq L$ have probability distribution
\begin{equation}
\label{eq_Hahn_hole_ensemble}
 \amsmathbb P_N(\ell_1,\ldots,\ell_N)=\frac{1}{\mathscr{Z}_N}\cdot \prod_{1 \leq i<j \leq N} (\ell_j-\ell_i)^2 \cdot \prod_{i=1}^N w(\ell_i),
\end{equation}
\begin{equation}
\label{eq_Hahn_hole_ensemble_2}
w(\ell)=
\left\{\begin{array}{lll}
  (\ell)_{C-\mathfrak{t}}\cdot(A+\mathfrak{t}+1-\ell)_{B-\mathfrak{t}}\cdot(\mathfrak{h}+1-\ell)_D^2&& \textnormal{if}\,\,0<\mathfrak{t}\leq \min(B,C),\\
  (\ell)_{C-\mathfrak{t}}\cdot(A+B+1-\ell)_{\mathfrak{t}-B}\cdot(\mathfrak{h}+1-\ell)_D^2 && \textnormal{if}\,\,B\leq \mathfrak{t}\leq C,\\
  (\ell)_{\mathfrak{t}-C}\cdot(A+C+1-\ell)_{B-\mathfrak{t}}\cdot(\mathfrak{h}+1-\ell)_D^2&& \textnormal{if}\,\,C\leq \mathfrak{t}\leq B,\\
  (\ell)_{\mathfrak{t}-C}\cdot(A+B+C+1-\mathfrak{t}-\ell)_{\mathfrak{t}-B}\cdot(\mathfrak{h}+1-\ell)_D^2 && \textnormal{if}\,\,\max(B,C)\leq \mathfrak{t}<B+C,
 \end{array}\right.
\end{equation}
where $(a)_n$ is the Pochhammer symbol and $\mathscr{Z}_N$ is a normalizing
constant.
\end{lemma}
\begin{proof}
 Tilings of the hexagon with a hole can be identified with tilings of the same hexagon without a hole, in which we freeze $D$ horizontal lozenges along the vertical line at distance $\mathfrak{t}$ from the left side at positions $\mathfrak{h}+1, \mathfrak{h}+2,\ldots,\mathfrak{h}+D$. Fixing these lozenges forces the $D\times D$ rhombus to be densely packed with $D^2$ horizontal lozenges. Removing them from the tiling, we get the desired hole.

 Hence, we can obtain \eqref{eq_Hahn_hole_ensemble} by fixing $D$ coordinates $\ell_i$ in \eqref{eq_Hahn_ensemble} to be equal to $\mathfrak{h}+1, \mathfrak{h}+2,\ldots,\mathfrak{h}+D$. The contribution of these fixed coordinates to the product over pairs is moved to the weight, resulting in the extra factor $(\mathfrak{h} + 1 - \ell)_D^2$ to go from \eqref{eq_Hahn_ensemble_2} to \eqref{eq_Hahn_hole_ensemble_2}. Note that the partition function also changes and the constants $\mathscr{Z}_N$ are different between \eqref{eq_Hahn_ensemble} and \eqref{eq_Hahn_hole_ensemble}.
\end{proof}

The particles in \eqref{eq_Hahn_ensemble} are naturally split into two groups --- we have $H=2$ in the notations of Section~\ref{Section_general_model} --- depending on whether they are below or above the hole. The first case is equivalent to $\ell_i\in [a_1,b_1]= [1,\mathfrak{h}]$ and the second to $\ell_i\in [a_2,b_2]=[ \mathfrak{h}+D+1,A+B+C-\mathfrak{t}]$. The filling fractions $(N_1,N_2)$ give by definition the number of horizontal lozenges below and above the hole, and $N_1 + N_2 = N$. For an alternative interpretation, we note that by viewing three types of lozenges as projections of three sides of unit cube, each lozenge tiling is realized as the projection of a stepped surface in the three-dimensional space. Then, the hole in the hexagon is viewed as a flat $D\times D$ square in the three-dimensional space and $N_1$ encodes the position of this square, or, equivalently, the height of the hole.

For this reason, the model with fixed filling fractions $(N_1,N_2)$ has an interpretation as tilings with a hole of fixed height; as such, this model was studied in \cite[Section 9.2]{BGG} and in \cite{BuGo3}, where global fluctuations in this model were computed and identified with the Gaussian free field. Here we can pursue another option and leave $N_1$ and $N_2$ free, instead: only their sum $N = N_1 + N_2$ is
deterministic while $N_1,N_2$ are random variables. Then, we can ask a new question: what is the asymptotic
distribution of $N_1$ as the size of the hexagon and its hole become very large? Theorem~\ref{Theorem_CLT_for_filling_fractionsi} shows that it converges towards a discrete Gaussian provided we can verify our assumptions.

Let us assume that $A,B,C,D,\mathfrak{h},\mathfrak{t}$ go to infinity at the same speed, which means that for a large master parameter $\mathcal N$, all other parameters grow linearly in it. In other words
\[
\left(\frac{A}{\N} , \frac{B}{\mathcal N},\frac{C}{\mathcal N}, \frac{D}{\N}, \frac{\mathfrak{h}}{\N}, \frac{\mathfrak{t}}{\mathcal N}\right)= (\hat A,\hat B,\hat C,\hat D, \hat{\mathfrak{h}},\hat{\mathfrak{t}}).
\]
remain $O(1)$ as $\N \rightarrow \infty$. The same is true for $\hat{n} = \frac{N}{\N}$ with $N$ given by \eqref{eq_N_hexagon_hole} in terms of the other parameters. There are certain extreme choices of $(\hat A,\hat B,\hat C,\hat D, \hat{\mathfrak{h}},\hat{\mathfrak{t}})$ for which the ensemble does not satisfy Assumptions \ref{Assumptions_Theta}, \ref{Assumptions_basic}, \ref{Assumptions_offcrit} or \ref{Assumptions_analyticity} of Section~\ref{Section_list_of_assumptions} (a simple example occurs when one of the parameters is close to $0$), yet we can identify situations when it does.
\begin{definition} \label{Definition_hex_hole_non_degenerate}
 We say that $(\hat A,\hat B,\hat C,\hat D, \hat{\mathfrak{h}},\hat{\mathfrak{t}})$ is \emph{non-degenerate} if there exists $\eps>0$ such that
 \[
 \hat{A},\hat{B},\hat{C} \in \big(\eps,\tfrac{1}{\eps}\big),\qquad \hat{D},\hat{\mathfrak{h}},\hat{n} > \eps,\qquad \hat{D} + \eps < \hat{\mathfrak{t}} < \hat{B} + \hat{C} - \hat{D} - \eps,
 \]
 and
 \[
 \hat{\mathfrak{h}} + \hat{D} + \eps < \left\{\begin{array}{lll}
 \hat A+\hat{\mathfrak{t}}&& \textnormal{if}\,\,0<\hat{\mathfrak{t}}\leq \min(\hat B,\hat C),\\
 \hat A+\min(\hat B,\hat C)&& \textnormal{if}\,\,\min(\hat B,\hat C)\leq \hat{\mathfrak{t}}\leq \max(\hat B,\hat C),\\
 \hat A+\hat B+\hat C-\hat{\mathfrak{t}}&& \textnormal{if}\,\,\max(\hat B,\hat C)\leq \hat{\mathfrak{t}}<\hat B+\hat C, \end{array} \right.
\]
 and the equilibrium measure of the discrete ensemble \eqref{eq_Hahn_hole_ensemble} has a band in each of the segments $[\hat a'_1+\eps,\hat b'_1-\eps]$ and $[\hat a'_2+\eps,\hat b'_2-\eps]$.
\end{definition}
 In words, the non-degeneracy condition means that all possible relevant sizes involved in the problem (such as the proportions of the hexagon, the size of the hole, the distance from the hole to the borders of the hexagon, the number $N$ of horizontal lozenges along the axis of the hole) stay bounded and bounded away from zero. In practice, the most tricky part to check is the very last condition of Definition~\ref{Definition_hex_hole_non_degenerate}, because it involves the equilibrium measure, which one needs to compute. Unfortunately, it is impossible to get rid of this condition: one can design a counterexample where a band starts at one of the endpoints $\hat a'_1$, $\hat b'_1$, $\hat a'_2$, $\hat b'_2$, while all other conditions hold; also one can design a counterexample where one of the segments $[\hat a'_1, \hat b'_1]$ or $[\hat a'_2,\hat b'_2]$ has no bands at all (\textit{i.e.}\ it is fully saturated, or fully void). At the end of this section we demonstrate some examples of the situations where we can check that the last condition holds.

\begin{proposition} \label{Proposition_hex_hole_offcrit}
Suppose that a $\N$-dependent sequence of hexagons with a hole has its parameters non-degenerate for all large enough $\N$ and some constant $\eps > 0$ not depending on $\N$. Then the corresponding discrete ensembles \eqref{eq_Hahn_hole_ensemble} satisfy Assumptions \ref{Assumptions_Theta}, \ref{Assumptions_basic}, \ref{Assumptions_offcrit} and \ref{Assumptions_analyticity} of Section~\ref{Section_list_of_assumptions} with $H=2$,
\[
\boldsymbol{\Theta}= \left(\begin{array}{cc} 1&1 \\ 1&1\end{array}\right),
\]
and the constants in these assumptions can be chosen to depend only on $\eps$.
\end{proposition}
\begin{remark} \label{rem:degeneratcomment}The condition which may fail for degenerate parameters is the first part of Assumption~\ref{Assumptions_offcrit}: one of the two segments might have no band (\textit{i.e.} might be fully void or fully saturated) and an endpoint of a band might coincide with an endpoint of one of the segments.
\end{remark}

\begin{corollary}
\label{Corollary_discrete_Gauss_in_hex}
 Under the conditions of Proposition~\ref{Proposition_hex_hole_offcrit}, the random filling fraction $N_1$ is asymptotically equal to a (one-dimensional) discrete Gaussian random variable, as in Theorems~\ref{Theorem_CLT_for_filling_fractions} and \ref{Theorem_CLT_for_filling_fractions_saturation}.
\end{corollary}
This corollary is a direct application of Theorem~\ref{Theorem_CLT_for_filling_fractions_saturation} (or Theorem~\ref{Theorem_CLT_for_filling_fractions} in the absence of saturations). We refer to \cite[Chapter 24]{Vadimlecture} for an interpretation of the parameters of the limiting discrete Gaussian random variable in terms of the Dirichlet energy of specific harmonic functions.

\begin{proof}[Proof of Proposition~\ref{Proposition_hex_hole_offcrit}] Since the weight \eqref{eq_Hahn_hole_ensemble_2} is given by four different formulae depending on the range of parameters, the details of the proof depend on the choice of the case. We only deal with the case $\hat{\mathfrak{t}} -\max(\hat{B},\hat{C})\geq \varepsilon$, the other cases can be handled similarly. In this situation the weight function in \eqref{eq_Hahn_ensemble} is given by
\begin{equation}
 w(\ell)=(\ell)_{\N\hat {\mathfrak{t}}-\N\hat C} \cdot (\N\hat A+\N\hat B+\N\hat C+1-\N\hat{\mathfrak{t}}-\ell)_{\N\hat{\mathfrak{t}}-\N\hat{B}} \cdot (\N\hat{\mathfrak{h}}+1-\ell)_{\N \hat D}^2.
\end{equation}
This weight can be recast through $\Gamma$-functions in two different forms. For $\ell \in \llbracket 1,\mathfrak{h} \rrbracket$ we can write
\begin{equation}
 \label{eq_x222}
w(\ell) = \frac{\Gamma\big(\N(\hat {\mathfrak{t}}-\hat C)+\ell\big)\cdot \Gamma\big(\N(\hat A+\hat C)+1-\ell\big)\cdot \Gamma\big(\N(\hat{\mathfrak{h}}+\hat D)+1-\ell\big)^2}{\Gamma(\ell) \cdot \Gamma\big(\N(\hat A+\hat B+\hat C-\hat{\mathfrak{t}})+1-\ell\big) \cdot \Gamma\big(\N \hat{\mathfrak{h}}+1-\ell\big)^2},
\end{equation}
while for $\ell \in \llbracket \mathfrak{h}+D+1,A+B+C-\mathfrak{t}\rrbracket$ we rather write
\begin{equation}
\label{eq_x222snd}
w(\ell) = \frac{\Gamma\big(\N(\hat {\mathfrak{t}}-\hat C)+\ell\big) \cdot \Gamma\big(\N(\hat A+\hat C)+1-\ell\big) \cdot \Gamma\big(-\N \hat{\mathfrak{h}}+\ell\big)^2}{\Gamma(\ell)\cdot \Gamma\big(\N(\hat A+\hat B+\hat C-\hat{\mathfrak{t}})+1-\ell\big)\cdot \Gamma\big(-\N(\hat{\mathfrak{h}}+\hat D)+\ell\big)^2}.
\end{equation}
We are in the setting of Section~\ref{Section_general_model} with two segments:
\[
H=2, \qquad [a_1,b_1]=[1,\mathfrak h], \qquad [a_2,b_2]={[\mathfrak h+D+1,
A+B+C-\mathfrak{t}]}.
\]
The first segment corresponds to the horizontal lozenges below the hole and the second to those above the hole. The interaction matrix is the $2 \times 2$ matrix with unit entries. The segment endpoints \eqref{eq_shifted_parameters} for the equilibrium measure are
\begin{equation}
\label{eq_x251}
[\hat a'_1,\hat b'_1]=\left[\frac{1}{2\N},\, \hat {\mathfrak{h}}+\frac{1}{2\N} \right], \qquad [\hat a'_2,\hat b'_2]=\left[\hat{\mathfrak h}+\hat D+\frac{1}{2\N},\, \hat A+\hat B+\hat C-\hat{\mathfrak{t}}+\frac{1}{2\N} \right].
\end{equation}
The single constraint \eqref{eq_equations_eqs} on the filling fractions $(N_1,N_2)$ is $N_1+N_2=N$, as in \eqref{eq_number_of_particles_contraint0}.

The validity of Assumption~\ref{Assumptions_Theta} is straightforward from the definitions. For the fourth point in the assumption, we see that for every $\boldsymbol{X}=(X_{1},X_{2})\in \amsmathbb R^{2}$ we have
 \[
 \boldsymbol{X}^{T}\cdot\boldsymbol{\Theta}\cdot \boldsymbol{X}=(X_{1}+X_{2})^{2}=\mathfrak r(\boldsymbol{X})^{2}
 \]
 with $\mathfrak r(\boldsymbol{X})=X_{1}+X_{2}$. We can therefore take $\mathfrak e=1$, the linear form $\mathfrak r_{1}=\mathfrak{r}$, the right-hand side $r_1 = \hat{n}$ and $\boldsymbol{\Theta}'$ the unit matrix of size $1$ in \eqref{eq_Theta_through_Theta_prime}. For the last point in Assumption~\ref{Assumptions_Theta}, we notice that the total length of the segments \eqref{eq_x251} is $\hat{A} + \hat{B} + \hat{C} + \hat{D} - \hat{\mathfrak{t}} = \hat{A} + \hat{n}$ (using the last case in \eqref{eq_N_hexagon_hole}), and the non-degeneracy assumption $\hat{A} > \eps$ also guarantees this last point.

Clearly, after a simple renormalization, \eqref{eq_x222} matches into the form \eqref{eq_ansatzw} of Assumption~\ref{Assumptions_analyticity} with
\[
\iota_1^- = \iota_2^+ = 1,\qquad \iota_1^+ = \iota_2^- = 2,\qquad \rho_{h,j}^{\pm} = 1.
\]
Hence, using the non-degeneracy of the parameters, both Assumptions~\ref{Assumptions_basic} and \ref{Assumptions_analyticity} follow from the definitions. The potentials $V_1$ and $V_2$ are explicitly given by:
\begin{align*}
 V_1(x)=& \mathrm{Llog}(x-\hat{a}'_1) +2\,\mathrm{Llog}(\hat{b}'_1-x) + U_1(x),
 \\
 U_1(x)=& - \frac{1}{\N} \log\left(\frac{\Gamma\big(\N(\hat {\mathfrak{t}}-\hat C+x)\big)}{\N^{\N(\hat {\mathfrak{t}}-\hat C+x)-\frac{1}{2}}} \cdot \frac{\Gamma\big(\N(\hat A+\hat C-x)+1\big)}{\N^{\N(\hat A+\hat C-x)+\frac{1}{2}}}\cdot \frac{\Gamma\big(\N (\hat{\mathfrak{h}}+\hat D-x)+1\big)^2}{\N^{2\N (\hat{\mathfrak{h}}+\hat D-x)+1}}\right)
 \\ &+ \frac{1}{\N} \log\left(\frac{\Gamma\big(\N(\hat A+\hat B+\hat C-\hat{\mathfrak{t}}-x)+1\big)}{\N^{\N(\hat A+\hat B+\hat C-\hat{\mathfrak{t}}-x)+\frac{1}{2}}}\right),
 \\
 V_2(x)=& 2\,\mathrm{Llog}(x-\hat{a}'_2) + \mathrm{Llog}(\hat{b}'_2-x) + U_2(x),\\
 U_2(x)=& - \frac{1}{\N}\log\left(\frac{\Gamma\big(\N(\hat {\mathfrak{t}}-\hat C+x)\big)}{\N^{\N(\hat {\mathfrak{t}}-\hat C+x)-\frac{1}{2}}} \cdot \frac{\Gamma\big(\N(\hat A+\hat C-x)+1\big)}{\N^{\N(\hat A+\hat C-x)+\frac{1}{2}}}\cdot \frac{\Gamma\big(\N (x-\hat{\mathfrak{h}})\big)^2}{\N^{2\N (x-\hat{\mathfrak{h}})-1}} \right)
 \\ &+ \frac{1}{\N} \log\left(\frac{\Gamma(\N x)}{\N^{\N x-\frac{1}{2}}}\right).
\end{align*}
It is somewhat inconvenient that $\Gamma$ factors in \eqref{eq_x222} are treated in two different ways as it is prescribed by the form of Assumption~\ref{Assumptions_analyticity} when we pass from the weight $w(\ell)$ to the potentials $V_1(x)$, $V_2(x)$: some of them become a part of $U_1$, $U_2$, while the others are replaced with a Stirling approximation to become a part of $V_1(x)-U_1(x)$ and $V_2(x)-U_2(x)$. Hence, we also introduce simplified weights $\tilde V_1(x)$, $\tilde V_2(x)$, obtained by replacing all $\Gamma$ factors by a Stirling approximation. While formally $x\in[\hat a'_1,\hat b'_1]$ for $\tilde V_1(x)$ and $x\in [\hat a'_2,\hat b'_2]$ for $\tilde V_2(x)$, both are given by the same formula
\begin{equation}
\label{eq_x284}
\begin{split}
& \quad \tilde{V}_1(x)=\tilde{V}_2(x)=\tilde{V}(x) \\
& = - \mathrm{Llog}\big(x + \hat {\mathfrak{t}}-\hat C - \tfrac{1}{2\N}\big)- \mathrm{Llog}\big(\hat A+\hat C + \tfrac{1}{2\N} - x\big) \\
& \quad + \mathrm{Llog}(x-\hat{a}'_1) + 2\,\mathrm{Llog}(x-\hat{a}'_2) +2\,\mathrm{Llog}(\hat{b}'_1-x) + \mathrm{Llog}(\hat{b}'_2-x).
\end{split}
\end{equation}
We introduced the shifts $\pm \frac{1}{2\N}$ in the two first terms for convenience in Lemma~\ref{Lemma_symmetric_hex}. We let $\boldsymbol{\mu}=(\mu_1,\mu_2)$ denote the equilibrium measure for the potential $(V_1(x),V_2(x))$ and $\tilde{\boldsymbol{\mu}}=(\tilde{\mu}_1,\tilde{\mu}_2)$ for the simplified potential $(\tilde{V}_1(x),\tilde{V}_2(x))$. In both cases we use the same matrix $\boldsymbol{\Theta}$ with unit entries and the segments \eqref{eq_x251}. The filling fractions $\hat n_1$, $\hat n_2$ are subjected to the constraint $\hat n_1+\hat n_2=\hat{B}+\hat{C}-\hat{\mathfrak{t}}-\hat{D}$, which is the normalized version of \eqref{eq_N_hexagon_hole} in the case we study.

Our task is to establish the off-criticality Assumption~\ref{Assumptions_offcrit} for $\boldsymbol{\mu}$. However, noticing that $V_h(x)-\tilde{V}_h(x)$ differ by a $O(\frac{1}{\N})$ term for $h \in \{1,2\}$, we can rely on Corollary~\ref{co:Lemregwithout} and it suffices to establish the off-criticality Assumption~\ref{Assumptions_offcrit} for the simplified equilibrium measure $\tilde{\boldsymbol{\mu}}$ under the conditions of Proposition~\ref{Proposition_hex_hole_offcrit}. This occupies us for the rest of the proof.

\medskip

Similarly to the constructions of Chapter~\ref{Chapter_smoothness} we introduce the two functions
\begin{equation}
\label{eq_x283}
\begin{split}
 \phi^+(x) & =\big(\hat{\mathfrak{t}}-\hat{C}+x - \tfrac{1}{2\N}\big) \cdot (\hat{b}'_1-x)^2 \cdot (\hat{b}'_2-x), \\
 \phi^-(x) & =\big(\hat{A}+\hat{C} + \tfrac{1}{2\N}-x\big) \cdot (x-\hat{a}'_1) \cdot (x-\hat{a}'_2)^2.
 \end{split}
\end{equation}
They are such that $\frac{\phi^+(x)}{\phi^-(x)}=e^{-\tilde{V}'(x)}$. We further let $\Gm(z)$ be the Stieltjes transform of the measure $\tilde \mu=\tilde\mu_1+\tilde \mu_2$ and set
\[
 q^{\pm}(z)=\phi^+(z) \cdot e^{\Gm(z)} \pm \phi^-(z) \cdot e^{-\Gm(z)}.
\]

\begin{lemma} \label{Lemma_qpm_hex_hole}
 The function $q^+(z)$ is a polynomial in $z$ of degree $4$ with leading coefficient $-2$. The function $(q^-(z))^2$ is a polynomial in $z$ of degree $6$.
\end{lemma}
\begin{proof}
Comparing with Definitions \ref{Definition_phi_functions_2} and \ref{GQdef}, we see that
\begin{equation}
\label{eq_x280}
\phi_1^{\pm}(z)=\frac{\phi^+(z)}{(z + \hat {\mathfrak{t}}-\hat C - \frac{1}{2\N}) \cdot (\hat{b}'_2-z)},\qquad q_1^\pm(z)=\frac{q^\pm(z)}{(z + \hat {\mathfrak{t}}-\hat C - \frac{1}{2\N}) \cdot (\hat{b}'_2-z)}.
\end{equation}
Hence, by Theorem~\ref{Theorem_regularity_density}-(iii) for $h=1$, $q^+(z)$ has no singularities on $[\hat a'_1,\hat b'_1]$. The same argument with $h=2$ shows that $q^+(z)$ has no singularities on $[\hat a'_2,\hat b'_2]$. Directly from its definition, $q^+(z)$ is holomorphic over $z\in \amsmathbb{C}\setminus \big([\hat a'_1,\hat b'_1]\cup [\hat a'_2,\hat b'_2]\big)$. Hence, $q^+(z)$ is an entire function. Directly from the definition and $\Gm(z)=O(\frac{1}{z})$, we get $q^+(z)\sim -2 z^4$ as $z\rightarrow \infty$. By Liouville theorem, $q^+(z)$ is a polynomial of degree $4$ with leading coefficient $-2$.

Switching to the properties of $q^-(z)$, we note $(q^+(z))^2-(q^-(z))^2=4\phi^+(z)\phi^-(z)$. Therefore, $(q^-(z))^2$ is a polynomial. Directly from its definition, $q^-(z)=O(z^3)$ as $z\rightarrow \infty$. Thus, $(q^-(z))^2$ is a polynomial of degree $6$.
\end{proof}

We now recall Proposition~\ref{Proposition_density}. Given the relation \eqref{eq_x280} and similar identities for $q_2^\pm$, \eqref{eq_density_tan} leads to the following formula for the density $\tilde \mu(x)$ of the simplified equilibrium measure:
\begin{equation}\label{eq_density_tan_hex_hole}
 \forall x \in (\hat{a}'_1,\hat{b}'_1)\cup (\hat{a}'_2,\hat{b}'_2) \qquad\tan\big(\pi \tilde \mu(x)\big)=\frac{ q^{-}(x^-) - q^{-}(x^+)}{2\ii \, q^+(x)}.
\end{equation}
We use this formula to check the conditions of Assumption~\ref{Assumptions_offcrit} through an auxiliary lemma.

\begin{lemma} \label{Lemma_hex_hole_bands}
 Under the conditions of Proposition~\ref{Proposition_hex_hole_offcrit}, $\tilde{\boldsymbol{\mu}}$ has exactly two bands and $q^-(z)$ factorizes as:
 \begin{equation}
 \label{eq_x282}
 q^-(z)=p(z) \sqrt{(z-\alpha_1)(z-\beta_1)(z-\alpha_2)(z-\beta_2)} ,
 \end{equation}
 where $(\alpha_1,\beta_1)$ and $(\alpha_2,\beta_2)$ are two bands, and $p(z)$ is a polynomial of degree $1$ with a zero inside $[\hat b'_1,\hat a'_2]$.
\end{lemma}
\begin{proof}
Suppose that $\hat b'_1$ is an endpoint of a saturation. In this situation $\Gm(\hat b'_1+y)$ behaves as $-\log(y)+O(1)$ for small positive $y$ and therefore $\exp(\Gm(\hat b'_1+y))$ has a simple pole at $y=0$ while $\exp(-\Gm(\hat b'_1+y))$ has a simple zero. Since $\phi^+(z)$ has a double zero at $z=\hat b'_1$, it follows that $q^{\pm}(z)$ are $O(z - \hat{b}_1')$ as $z \rightarrow \hat{b}_1'$. Similarly, if $\hat{a}'_2$ is an endpoint of a saturation, then $q^{\pm}(z)$ have a simple zero at $\hat{a}'_2$.

 If instead $\hat b'_1$ and $\hat a'_2$ are both endpoints of voids, then directly from the definition of $q^-(z)$, we have $q^-(\hat{b}'_1)<0<q^-(\hat{a}'_2)$. As $q^-(z)$ is real-valued and continuous over $z\in [\hat b'_1,\hat a'_2]$, it must have a zero at some point $z_* \in (\hat b'_1,\hat a'_2)$. As it is analytic near $z_*$, we must have $q^-(z) = O(z - z_*)$ as $z \rightarrow z_*$. Therefore, the polynomial $(q^-(z))^2$ of total degree $6$ factorizes as
 \[
 c\cdot (z-\alpha_1) \cdot (z-\beta_1) \cdot (z-\alpha_2) \cdot (z-\beta_2) \cdot (z-z_*)^2,
 \]
 where $c$ is a non-zero constant, $z_* \in [\hat b'_1,\hat a'_2]$ and $\alpha_1,\beta_1,\alpha_2,\beta_2$ are four complex numbers. They are identified with endpoints of the bands after we observe that \eqref{eq_density_tan_hex_hole} implies the necessity of a branchcut for $q^-(z)$ along each band. For the same reason, if there were more bands, their endpoints would have to be other zeros of $q^-(z)$. So, there are only two bands. The non-degeneracy of parameters included in the assumption that each segment contained at least a band. So, up to relabeling these are $(\alpha_1,\beta_1)\subset(\hat a'_1,\hat b'_1)$ and $(\alpha_2,\beta_2)\subset(\hat a'_2,\hat b'_2)$. We have established \eqref{eq_x282}. The fact that $\tilde \mu$ has exactly two bands then follows from \eqref{eq_density_tan_hex_hole}.
\end{proof}

We are now ready to check the conditions of Assumption~\ref{Assumptions_offcrit}. The first condition follows from the assumed non-degeneracy of parameters. The fourth, fifth and sixth conditions follow from \eqref{eq_density_tan_hex_hole}-\eqref{eq_x282} and the observation that zeros of $q^+(z)$ cannot be close to zeros of $q^{-}(z)$ at $\alpha_1,\beta_1,\alpha_2,\beta_2$ because of the relation $(q^+(z))^2-(q^{-}(z))^2=4\phi^+(z)\phi^-(z)$. It remains to check the second and the third conditions. Note that in a void we have directly from \eqref{eq_x283} and the Definition~\ref{def_eff_pot} of the effective potential
\begin{equation}\label{eq_effective_qm}
q^{-}(x)=\phi^{+}(x) e^{\Gm(x)}(1-e^{\partial_{x }V^{\textnormal{eff}}(x)}).
\end{equation}
Eq.\ \eqref{eq_x282} implies that $q^{-}(x)$ has no zeros in the voids inside $(\hat a'_1,\hat b'_1)$ and $(\hat a'_2,\hat b'_2)$. Hence, by \eqref{eq_effective_qm}, $\partial_{x }V^{\textnormal{eff}}(x)$ has fixed sign, that is, $V^{\textnormal{eff}}(x)$ is strictly monotone, which readily implies the bound on the effective potential in voids of the second condition in Assumption~\ref{Assumptions_offcrit}.

When $x$ is in a saturation, there is ambiguity for the meaning of $\Gm(x)$. A proper way to differentiate $V^{\textnormal{eff}}(x)$ in this situation is by using the principal value integral, which means:
\begin{equation}
 \partial_{x }V^{\textnormal{eff}}(x)= V'(x)- \Gm(x^+)-\Gm(x^-).
\end{equation}
However, because the density of the equilibrium measure in the saturated region is equal to $1$, $\Gm(x^+)$ differs from $\Gm(x^-)$ by $2\ii\pi$, which implies $\exp(\Gm(x^+))=\exp(\Gm(x^-))$. Hence, the choice of $x^+$ versus $x^-$ becomes irrelevant in the exponential, \eqref{eq_effective_qm} continues to hold, and $V^{\textnormal{eff}}(x)$ is again a monotone function, implying the third condition in Assumption~\ref{Assumptions_offcrit}. This finishes the verification of Assumptions~\ref{Assumptions_Theta}, \ref{Assumptions_basic}, \ref{Assumptions_offcrit} and \ref{Assumptions_analyticity} and completes the proof of Proposition~\ref{Proposition_hex_hole_offcrit}.
\end{proof}

\medskip

We conclude our discussion of the hexagon with a hole by giving a class of parameters for which the assumptions of Definition~\ref{Definition_hex_hole_non_degenerate} hold.

\begin{definition} \label{Definition_hex_hole_symmetric}
We say that the parameters $(\hat{A},\hat{B},\hat{C},\hat{D},\hat{\mathfrak{h}},\hat{\mathfrak{t}})$ are \emph{symmetric} if $\hat{B}=\hat{C}$ and
\[
 2 \hat {\mathfrak h}+\hat D = \left\{\begin{array}{ll} \hat A+\hat{\mathfrak{t}} &\textnormal{if } \,0<\hat{\mathfrak t}\leq \hat B,\\
                       \hat A+ 2\hat B -\hat{\mathfrak{t}} &\textnormal{if }\, \hat B\leq \hat{\mathfrak t}\leq 2\hat B. \end{array}\right.
 \]
\end{definition}
In words, symmetric parameters are those for which the hexagon and the hole are invariant under the reflection with respect to the horizontal line passing through the middle in the hole.

\begin{lemma} \label{Lemma_symmetric_hex}
 Fix $\hat{A} > 0$, $\hat{B}=\hat{C}>0$, and suppose that the vertical line at distance $\hat{\mathfrak t} \in (0,\hat{B} + \hat{C})$ from the left border of $\hat{A}\times \hat{B}\times \hat{C}$ hexagon does not pass through any of the tangency point of the inscribed ellipse. There exists $\tilde D_0>0$ such that, for any $\hat{D} \in [0,\hat{D}_0]$ and with the choice of $\hat{\mathfrak{h}}$ making the parameters symmetric, the simplified equilibrium measure $\tilde{\boldsymbol{\mu}}$ for the potential \eqref{eq_x284} has two bands $(\alpha_1,\beta_1)$ and $(\alpha_2,\beta_2)$ satisfying strict inequalities
 \begin{equation}
\label{eq_hex_hole_strict}
 \hat a'_1<\alpha_1<\beta_1<\hat b'_1<\hat a'_2<\alpha_2<\beta_2<\hat b'_2.
\end{equation}
 In addition, the segments $[\beta_1,\hat b_1']$ and $[\hat a'_2,\alpha_2]$ are void, while the segments $[\hat a'_1,\alpha_1]$ and $[\beta_2,\hat b'_2]$ are either both void or both saturated.
\end{lemma}
\begin{proof}
 The symmetry of the parameters implies that the variational datum for $\tilde{\boldsymbol{\mu}}$ is invariant under the reflection with respect to the point $\hat{\mathfrak{h}}+\frac{\hat D}{2}+\frac{1}{2\N}$. By the uniqueness in Proposition~\ref{Lemma_maximizer}, the simplified equilibrium measure $\tilde{\boldsymbol{\mu}}$ must be reflection-invariant as well. Its total mass remains strictly between $0$ and $2\mathfrak{h}$ because the vertical line remains away from the left and right borders of the hexagon. Therefore, $\tilde{\mu}$ cannot be identically $0$ neither can be fully saturated. As the continuity of its density established in Theorem~\ref{Theorem_regularity_density} forbids the direct adjacency of a void with a saturation, there must be at least a band. Since the symmetry point is between the two segments, the number of bands is even. We deduce there are at least two bands.

From the definition of $q^-(z)$ and the symmetry of the parameters, we have
 \[
 q^-\big(2\hat{\mathfrak{h}} + \hat{D} + \tfrac{1}{\N} - z\big) = -q^-(z).
 \]
 Hence, $q^-(\hat{\mathfrak{h}}+\frac{\hat{D}}{2}+\frac{1}{2\N})=0$. Moreover, combining with Lemma~\ref{Lemma_qpm_hex_hole} and the analyticity of $q^-(z)$ in a complex neighborhood of $(\hat b'_1,\hat a'_2)$, we conclude that
 \begin{equation}
 \label{eq_x281}
 \left( \frac{q^-(z)}{z- \hat{\mathfrak{h}} - \frac{\hat{D}}{2} - \frac{1}{2\N}}\right)^2
 \end{equation}
 is a polynomial of degree $4$, that we can factorize as $c\cdot(z-\alpha_1)(z-\alpha_2)(z-\beta_1)(z-\beta_2)$. Equation~\eqref{eq_density_tan_hex_hole} implies that $q^-(z)$ should have branch cuts along the bands and we know there are at least two bands. We deduce that there are exactly two bands, the four zeros of $(q^-(z))^2$, that we denote $\alpha_1,\beta_1,\alpha_2,\beta_2$, must be simple and give the endpoints of these bands. We can label them such that
\[
(\alpha_1,\beta_1)\subset [\hat a_1',\hat b_1'] \qquad \textnormal{and} \qquad (\alpha_2,\beta_2)\subset [\hat a_2',\hat b_2'].
\]

 Next, let us justify the strict inequalities $\beta_1<\hat{b}'_1$ and $\hat{a}'_2<\alpha$. By contradiction, suppose that $\beta_1=\hat{b}'_1$ (which implies also $\alpha_2=\hat{a}'_2$ by symmetry). Then $q^-(\hat{b}'_1)=0$ together with $(q^+(z)q^-(z))^2=4\phi^+(z)\phi^-(z)$ implies $q^+(\hat{b}'_1)=0$. Therefore, by polynomiality of $q^+(z)$, the point $\hat{b}'_1$ is a double zero of $(q^+(z))^2$; but it is also a double zero of $4\phi^+(z)\phi^-(z)$. Hence, applying $(q^+(z)q^-(z))^2=4\phi^+(z)\phi^-(z)$ again, we conclude that $\hat{b}'_1$ is a double zero of the polynomial $(q^-(z))^2$, which we already excluded.

 Further, let us show that $[\beta_1,\hat b'_1]$ and $[\hat a'_2,\alpha_2]$ are void. By contradiction, suppose that they are saturated instead. As we explained at the beginning of the proof of Lemma~\ref{Lemma_hex_hole_bands}, this would imply $q^-(\hat{b}'_1)=q^-(\hat{a}'_2)=0$, which is contradiction with $\alpha_1,\beta_1,\alpha_2,\beta_2$ being the only zeros of $(q^-(z))^2$.

 It remains to justify the strict inequalities $ \hat{a}'_1<\alpha_1$ and $\beta_2<\hat{b}'_2$. For this purpose we analyze the limit $\hat{D} \rightarrow 0$. If $\hat{D}=0$, then there is no hole and we deal with uniformly random tilings of the $\tilde A\times \tilde B\times \tilde C$ hexagon from Section~\ref{Section_Hexagon}. In this situation $ \hat{a}'_1<\alpha_1$ and $\beta_2<\hat{b}'_2$ follow from the fact that the endpoints of the bands form the inscribed ellipse (\textit{cf.} \cite{CLP} or \cite[Section 10.2 or Section 18.2 or Section 21.2]{Vadimlecture}). Due to Lemma~\ref{Lemma_continuity_varyinf_ff}, the simplified equilibrium measure $\tilde{\boldsymbol{\mu}}$ continuously depends on $\hat{D}$. Hence, so does the function $q^{-}(z)$, the ratio \eqref{eq_x281}, and its zeros. Hence, for small $\hat{D}$, the values of $\alpha_1$ and $\beta_2$ are close to their $\hat{D}=0$ counterparts and the inequalities $ \hat{a}'_1<\alpha_1$ and $\beta_2<\hat{b}'_2$ continue to hold.
\end{proof}

Lemma~\ref{Lemma_symmetric_hex} delivers a range of non-degenerate parameters, for which Proposition~\ref{Proposition_hex_hole_offcrit} and Corollary~\ref{Corollary_discrete_Gauss_in_hex} hold. Further, by the continuous dependence of the equilibrium measure on the parameters, as in Lemma~\ref{Lemma_continuity_varyinf_ff} and Theorem~\ref{Theorem_off_critical_neighborhood}, if we slightly deform the parameters from the symmetric case, the conditions remain. For instance, this is the case for a close-to-symmetric simulation on the right panel of Figure~\ref{Fig_hex_simulation}.

\section{Domains obtained as gluings of trapezoids}
\label{Section_gluing_def}

We present a class of domains drawn on a two-dimensional surface, which are more general than the domains of Section~\ref{Section_simple_domains} and whose random lozenge tilings can be analyzed by the methods of this book. The study of the resulting tiling models serve as a motivation for the general definitions given in Section~\ref{Section_general_model}.

\subsection{An example}
\label{sec:extile}

\begin{figure}[h]
\begin{center}
 {\scalebox{0.62}{\includegraphics{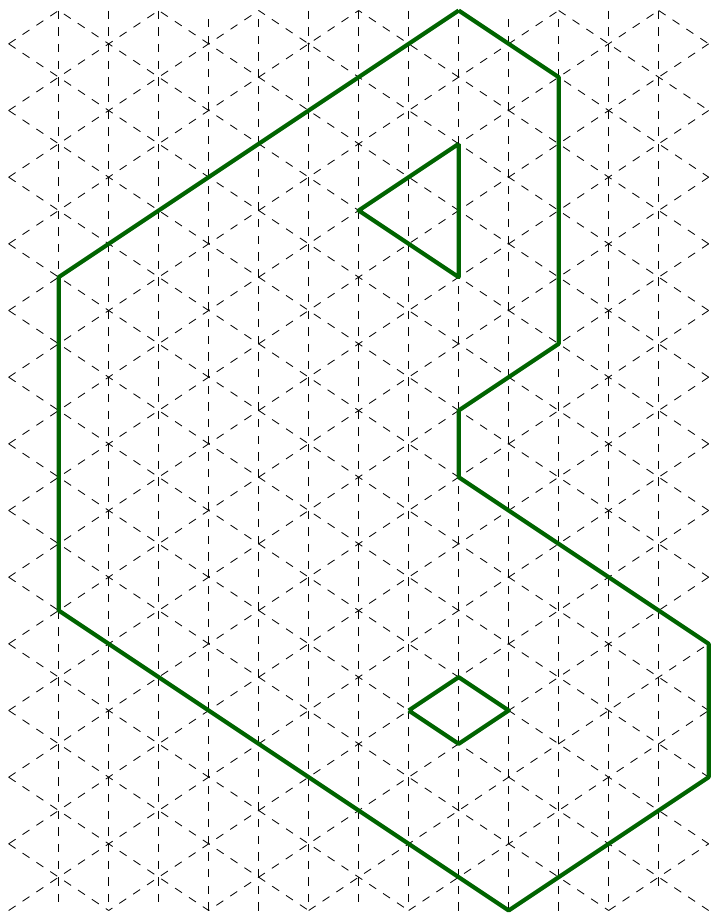}}} \hfill
 {\scalebox{0.62}{\includegraphics{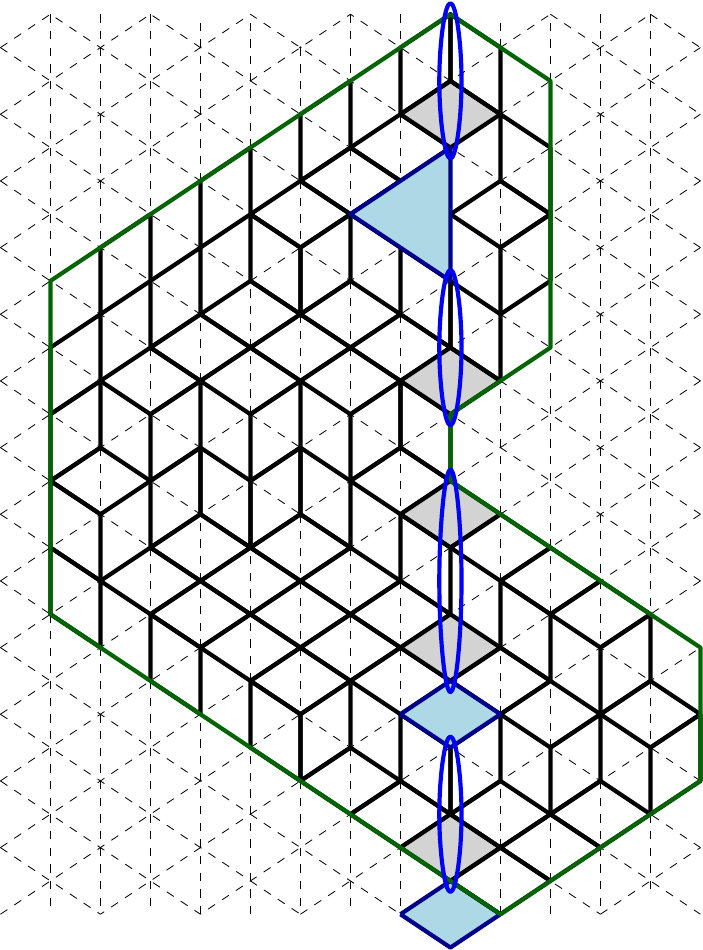}}}
\end{center}

\caption{Left panel: A polygon resulting in 4 segment filling fractions with 2 relations.
 Right panel: an example of tiling and its vertical
section with populations $N_1 = 1$, $N_2 = 2$, $N_3 = 1$ and $N_4 = 1$. Blue ovals indicate 4 groups of horizontal lozenges along the section.
\label{Fig_crazy_polygon}}
\end{figure}

\begin{figure}[h]
\begin{center}
 {\scalebox{0.6}{\includegraphics{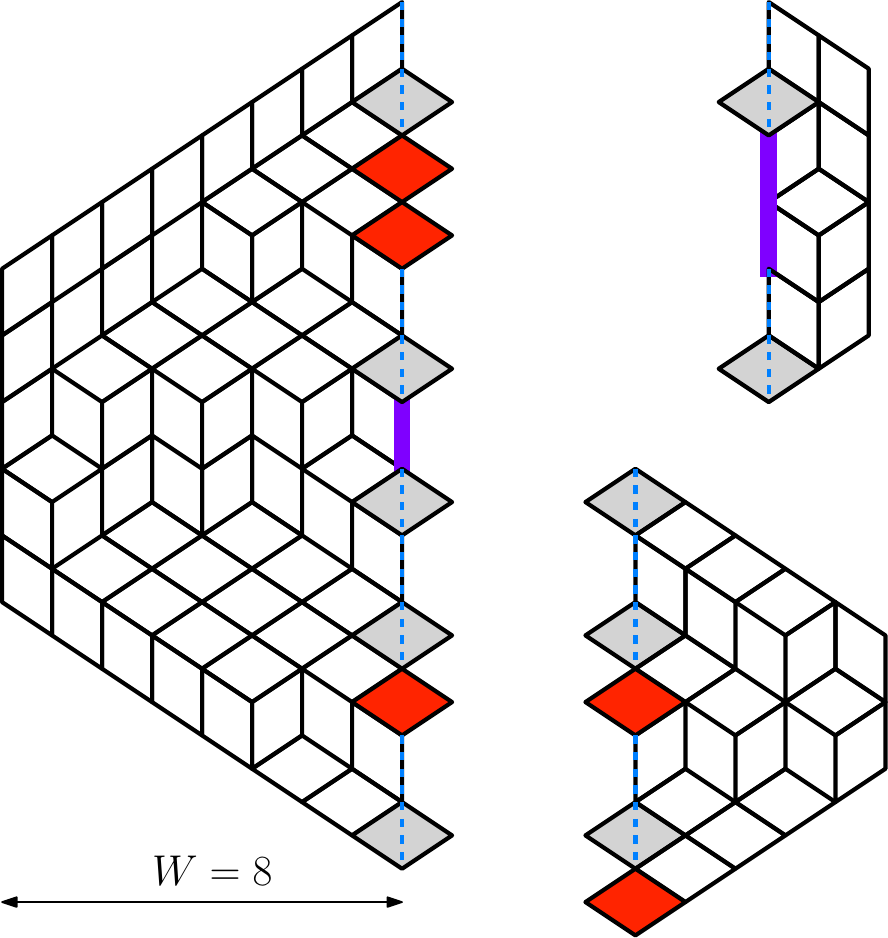}}}
\end{center}

\caption{The polygon splits into three trapezoids of respective widths $8$, $5$, and $2$. The dashed blue line is our vertical section. Horizontal lozenges are prohibited to occupy the thick violet part of the section. In scarlet we indicate the frozen horizontal lozenges (position $-1$, $2$, $9$, $10$, with the convention that the bottom of the left trapezoid is located at position $0$). The remaining segments of the vertical sections define the four groups.\label{Fig_polygon_split}}
\end{figure}

 Let us start by presenting a relatively complicated and representative \emph{planar} case. In Sections \ref{Section_Hexagon} and \ref{Section_hex_hole} the domains were represented as gluings of two trapezoids, while in this section we glue three of them. Consider a uniformly random lozenge tiling of the polygon of Figure
 \ref{Fig_crazy_polygon} and its section along the vertical line going through the
 axis of the two holes. We identify each tiling of this domain with three tilings of three trapezoids agreeing by their long bases: one left trapezoid and two right trapezoids, as in the right panel of Figure~\ref{Fig_crazy_polygon} and in Figure~\ref{Fig_polygon_split}. Along this vertical line, each tiling has four horizontal lozenges at prescribed (frozen) positions, and five horizontal lozenges at free positions split in four groups, \textit{cf.} Figure~\ref{Fig_crazy_polygon}. The four groups correspond to possible positions of free horizontal lozenges in the four integer intervals
 \[
 \llbracket A_1,B_1\rrbracket = \{0,1\},\qquad \llbracket A_2,B_2\rrbracket =\{3,4,5\},\qquad \llbracket A_3,B_3\rrbracket = \{7,8\},\qquad \llbracket A_4,B_4\rrbracket =\{11,12\}.
 \]
inside the vertical base of the trapezoid on the left panel of Figure~\ref{Fig_polygon_split}. Let $(N_h)_{h = 1}^{4}$ denote the numbers of
free horizontal lozenges in the four groups. The geometry dictates the two constraints $N_1 + N_2 = 3$ (the bottom right trapezoid has width $5$ but two horizontal lozenges are frozen), $N_3 + N_4 = 2$ (the top right trapezoid has width $2$ and no frozen horizontal lozenges). The third constraint $N_1+ N_2 + N_3 + N_4 = 5$ (the left trapezoid has width $8$ but three frozen horizontal lozenges) is a consequence of the two others.

 We would like to compute the probability distribution on the five horizontal lozenges induced by the uniform measure on all lozenge tilings of the domain. Once the positions of the five free lozenges are fixed, the domain splits into three simpler domains, trapezoids of Figure~\ref{Fig_polygon_split}. Hence, the distribution of the five free lozenges is proportional to the product of the number of tilings of three trapezoids with fixed boundaries. To enumerate tilings of trapezoids we use the Proposition~\ref{Proposition_number_tilings_trapezoid}. Multiplying \eqref{eq_Weyl_dimension_main} over the three trapezoids, we arrive at the following formula for the distribution of the five free lozenges
 \begin{equation}
 \label{eq_crazy_tiling_section}
 \begin{split}
 \amsmathbb{P}(\ell_1,\ldots,\ell_5) & = \frac{1}{\Z} \cdot \prod_{\substack{i \in I_1 \\ j\in I_2}}
 (\ell_j-\ell_i)^2 \cdot
 \prod_{\substack{i \in I_1 \\ j \in I_3}}
 (\ell_j-\ell_i) \cdot
 \prod_{\substack{i \in I_1 \\ j \in I_4}}
 (\ell_j-\ell_i) \cdot
 \\
& \quad \times \prod_{\substack{i \in I_2 \\ j \in I_3}}
 (\ell_j-\ell_i) \cdot
 \prod_{\substack{i \in I_2 \\ j \in I_4}}
 (\ell_j-\ell_i) \cdot
\prod_{\substack{i \in I_3 \\ j \in I_4}}
 (\ell_j-\ell_i)^2 \cdot
  \prod_{h=1}^4 \prod_{\substack{ i<j\, : \, \\ i,j  \in I_h}}
 (\ell_j-\ell_i)^2 \\
 & \quad \times
 \prod_{i \in I_1} (10 - \ell_i)(9 - \ell_i)(2 - \ell_i)^2(\ell_i+1) \cdot
 \prod_{i \in I_2} (10 - \ell_i)(9 - \ell_i)(\ell_i-2)^2(\ell_i+1)
 \\ & \quad \times \prod_{i \in I_3} (10 - \ell_i)(9 - \ell_i)(\ell_i-2) \cdot
 \prod_{i \in I_4} (\ell_i-10)(\ell_i-9)(\ell_i-2).
\end{split}
\end{equation}
with $\amsmathbb{I}_h =\big\{i\,\,|\,\,\ell_i \in [A_h,B_h]\big\}$. Further, if we take a positive integer parameter $\N$, multiply all the
proportions of the polygon (and holes) in Figure~\ref{Fig_crazy_polygon} by $\N$, and replace each frozen particle at position $k$ by frozen particles at positions $(k-1)\N+1,(k-1)\N+2 + 1,\ldots,k\N$, then we obtain a similar
probability distribution on $5\N$ particles, thought as encoding the position of the free horizontal lozenges on the vertical section. The only differences are that $[A_h,B_h]$ is replaced with $[a_h,b_h] = [\N\hat{a}_h,\N\hat{b}_h] := [\N (A_h-1)+1,\N(B_h+1) - 1]$, the number of particles in the $h$-th group (segment filling fractions) is $N_h = \N\hat{n}_h$, but we keep the relations
\begin{equation}
\label{starfortileex} (\star) \qquad \hat{n}_1 + \hat{n}_2 = 3\quad \textnormal{and}\quad \hat{n}_3 + \hat{n}_4 = 2.
\end{equation}
Besides, each individual factor $(\ell_i - k)$ appearing in \eqref{eq_crazy_tiling_section} for some $k \in \amsmathbb{Z}$ is replaced with a Pochhammer symbol to take into account the presence of $\N$ times many frozen horizontal lozenges. More precisely,
we should make in \eqref{eq_crazy_tiling_section} the substitutions
\begin{equation}
\label{substituN}
(k - \ell_i) \longrightarrow \big(((k-1)\N+1 - \ell_i\big)_{\N},\qquad (\ell_i - k) \longrightarrow (\ell_i - k\N)_{\N}.
\end{equation}
It is not necessary to distinguish the cases $(k - \ell_i)$ and $(\ell_i - k)$ as these two expressions are valid independently of the relative position $k$ and $\ell_i$, but we often prefer to do so in order to keep arguments of the Pochhammer symbols positive. In particular, all factors in the right-hand side of \eqref{eq_crazy_tiling_section} are positive. Note that the factors corresponding to consecutive frozen horizontal lozenges can be collected into a single Pochhammer symbol, for instance
\[
(10 - \ell_i)(9 - \ell_i) \longrightarrow (9\N+1 - \ell_i)_{\N} \cdot (8\N+1 - \ell_i)_{\N} = (8\N+1 - \ell_i)_{2\N}.
\]
We can rewrite the resulting formula for the probability distribution in the form
\begin{equation}
\label{eq_general_measure_int}
 \P(\boldsymbol{\ell})= \frac{1}{\Z_\N} \cdot \prod_{1\leq i<j \leq N}\bigg[
 \frac{1}{\N^{2\theta_{h(i),h(j)}}}\cdot \frac{\Gamma\big(\ell_j-\ell_i+1\big)\cdot
 \Gamma\big(\ell_j-\ell_i+\theta_{h(i),h(j)}\big)}{\Gamma\big(\ell_j-\ell_i\big) \cdot
 \Gamma\big(\ell_j-\ell_i+1-\theta_{h(i),h(j)}\big)}\bigg] \cdot
  \prod_{i=1}^N w_{h(i)}(\ell_i),
\end{equation}
where $h(i)=g$ whenever $\ell_i\in [a_g,b_g]$. The interaction matrix $\boldsymbol{\Theta} = (\theta_{g,h})_{g,h \in [H]}$ reads
\[
 \boldsymbol{\Theta} = \left(\begin{array}{cccc} 1 & 1 & \frac{1}{2} & \frac{1}{2}\\[0.6ex] 1 & 1 & \frac{1}{2}& \frac{1}{2}\\[0.6ex] \frac{1}{2} & \frac{1}{2}& 1
 &1\\[0.6ex] \frac{1}{2} & \frac{1}{2} &1 & 1 \end{array} \right)
 \]
 and the weights are
 \begin{equation}
\label{weightchoiceex} \begin{split}
w_1(\ell) & = \N^{-4\N}\cdot (8\N+1 - \ell_i)_{2\N}\cdot (\N+1 - \ell_i)_{\N}^2 \cdot (\ell_i + \N)_{\N}, \\
w_2(\ell) & = \N^{-4\N}\cdot (8\N+1 - \ell_i)_{2\N}\cdot (\ell_i - 2\N)_{\N}^2\cdot (\ell_i + \N)_{\N}, \\
w_3(\ell) & = \N^{-3\N} \cdot (8\N+1 - \ell_i)_{2\N} \cdot (\ell_i - 2\N)_{\N}, \\
w_4(\ell) & = \N^{-3\N}\cdot (\ell_i - 10\N)_{2\N} \cdot (\ell_i - 2\N)_{\N}.
 \end{split}
 \end{equation}
We included factors of $\N$ in \eqref{eq_general_measure_int} and in the definition \eqref{weightchoiceex} of the weights by convenience in anticipation of the large $\N$ analysis: they only affect the value of the partition function $\Z_{\N}$, which is always chosen such that \eqref{eq_general_measure_int} is a probability measure. Note that for special values of $\theta$, the ratio of four Gamma functions in \eqref{eq_general_measure_int} simplifies:
 \[
\forall \theta \in \big\{0,\tfrac{1}{2},1\big\}\qquad \frac{\Gamma\big(\ell_j-\ell_i+1\big) \cdot
 \Gamma\big(\ell_j-\ell_i+\theta\big)}{\Gamma\big(\ell_j-\ell_i\big) \cdot
 \Gamma\big(\ell_j-\ell_i+1-\theta\big)}=(\ell_j-\ell_i)^{2\theta}.
 \]
 Nevertheless, we use the form with Gamma functions, as this gives a proper generalization for general choices of the matrix $\boldsymbol{\Theta}$ making the results of Chapter~\ref{ChapterNekra} available.

For the resulting stochastic system of $5\N$ particles, we can ask the asymptotic
$\N \rightarrow\infty$ questions of the same flavor as before: expansion of the partition
function $\Z_\N$, law of large numbers for the position of particles and central limit theorem for global fluctuations,
asymptotic fluctuations for the segment filling fractions, \textit{i.e.} the number of horizontal tiles in each group (if we do not fix them). All these
questions can be addressed by our methods. Moreover, we can and will use the asymptotic behavior of these $5\N$ particles as an input for analyzing the full tiling.

\subsection{The general gluing procedure}
\label{sec:genglu}
There is a whole class of domains generalizing the hexagon with a hole and polygon of Figure~\ref{Fig_crazy_polygon},
 and on which tiling models still fit in the discrete ensembles studied in this book. The domains of this class are obtained by gluing several trapezoids along a single vertical axis, and do not have to be planar or orientable. Let us describe the procedure in more details.

\begin{figure}[h!]
\begin{center}
 {\scalebox{0.6}{\includegraphics{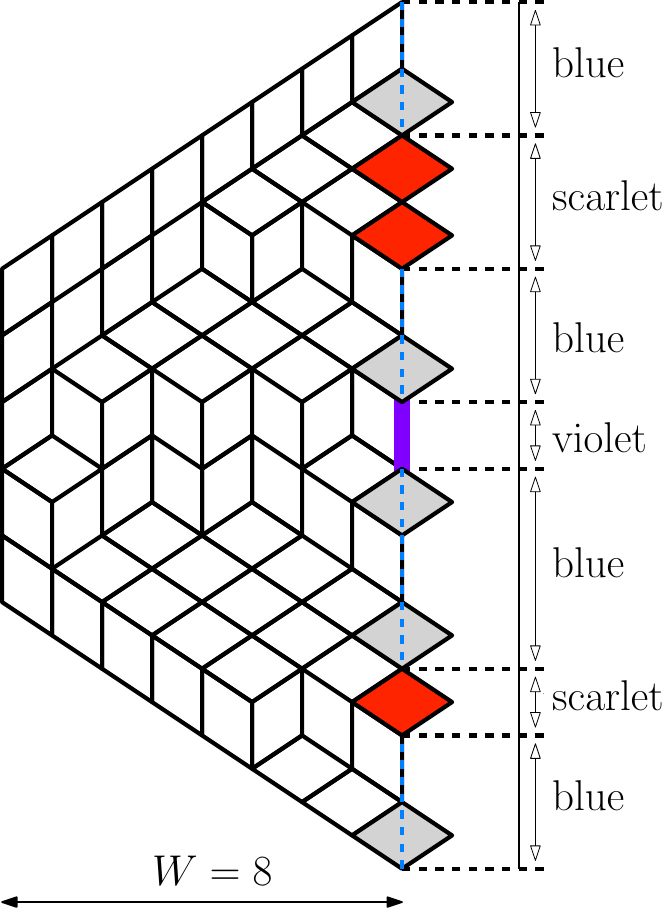}}}
\end{center}
\caption{Parameterization of the left trapezoid from Figures \ref{Fig_crazy_polygon}--\ref{Fig_polygon_split}. The segment $\llbracket D,U\rrbracket =\llbracket 0,12\rrbracket$ is split into seven segments: four blue, two scarlet and one violet.
\label{Fig_trapezoid_parameterization}}
\end{figure}

We start by parameterizing a single trapezoid. For that we take a segment of integers and a map
\begin{equation}
\label{eq:typemap}
\textnormal{type}\,:\,\, \llbracket D,U\rrbracket \rightarrow \{\textnormal{s},\textnormal{v},\textnormal{b}\},
\end{equation}
associating to each site in the segment one the values $\textnormal{s}$ = scarlet, $\textnormal{v}$ = violet, $\textnormal{b}$ = blue. Here $D/U$ stands for ``down''/''up''. In addition we fix a positive integer $W$. An admissible configuration of dents on a trapezoid of width $W$ with long base $\llbracket D,U\rrbracket$ is encoded by a configuration of $W$ sites in $\llbracket D,U\rrbracket$ such that
\begin{itemize}
\item a scarlet site $x$ must belong to the configuration: it corresponds to a frozen horizontal lozenge containing $[x,x + 1]$ as a diagonal;
\item a violet site $x$ is not allowed in the configuration: it corresponds to frozen segment $[x,x+1]$ on the long base.
\end{itemize}
An example of configuration is shown in Figure~\ref{Fig_trapezoid_parameterization}. There are no constraints for the blue sites: they correspond to possible positions for free dents (covered by horizontal lozenges), and will eventually become internal points of the tiled domain obtained after gluing several trapezoids together. The scarlet sites lead to the formation of frozen tiled subdomains common to all tilings compatible with the configuration of sites. These frozen subdomains can be removed from the two-dimensional domain, as we did in Figure~\ref{Fig_crazy_polygon}. The violet sites indicate parts of the long base of the trapezoid which are prohibited from crossing by lozenges. They will eventually become part of the boundary of the domain after gluing several trapezoids. Maximal sequences of consecutive blue sites form what we call \emph{blue segments}.

Let $S$ denote the total number of scarlet sites in $\llbracket D,U\rrbracket$. The geometry of the trapezoid dictates that there must be precisely $W-S$ horizontal lozenges on blue sites. Summing up, $\llbracket D,U\rrbracket$ comes with a finite collection of maximal blue segments of integers $\llbracket A_h,B_h\rrbracket$ (we call them simply blue segments), and the $h$-th segment can have $N_h$ horizontal lozenges (this number is not necessarily deterministic) subject to the condition $\sum_h N_h=W-S$.

Geometrically, shifting $\llbracket D,U\rrbracket $ and s/v/b sites by an integer leads to exactly the same trapezoid. However, a particular choice of shifts becomes important in the following definition, which explains our notion of gluing.

\begin{definition}[Gluing of trapezoids] \label{Definition_gluing}
 Take trapezoids $\mathcal{T}_1,\ldots,\mathcal{T}_m$ with long bases $\llbracket D_i,U_i\rrbracket$, $i\in[m]$, equipped with type maps \eqref{eq:typemap} and containing a total of $2H$ (maximal) blue segments. We say that the trapezoids are glued into a domain $\mathcal D$, if for each trapezoid $\mathcal T_i$ and each of its blue segment $\llbracket A_h,B_h\rrbracket\subset \llbracket D_i,U_i\rrbracket$ there exists exactly one other trapezoid $\mathcal T_{i'}$, $i'\neq i$ and a blue segment in it $\llbracket A_{h'},B_{h'}\rrbracket\subset \llbracket D_{i'},U_{i'}\rrbracket$, such that $\llbracket A_{h'},B_{h'}\rrbracket$ coincides with $\llbracket A_{h},B_{h}\rrbracket$. In this situation, the trapezoids $\mathcal T_i$ and $\mathcal{T}_{i'}$ are \emph{adjacent}, and we identify indices $h$ and $h'$. We require the gluing to be connected: any two trapezoids can be linked by a chain of adjacent ones.
\end{definition}

\begin{remark} \label{Example_gluing}
  Sometimes the domain $\mathcal{D}$ is planar, \textit{i.e.} it comes with an embedding into the triangular lattice drawn on $\amsmathbb R^2$. For that, let us designate each trapezoid as ``left'' or ``right'' and for each collection of adjacent $S$ scarlet sites remove the triangles whose base is formed by these sites from the corresponding trapezoid (this removes $\frac{S(S+1)}{2}$ horizontal lozenges). In addition, draw a vertical line $\mathsf L$ on the lattice and equip it with a coordinate system. Next, glue the trapezoids to $\mathsf L$ along their long bases at coordinates of these bases, on the left or on the right depending on the designation. Assuming that all trapezoids (with triangles corresponding to scarlet sites removed) do not overlap as subsets of the plane after such a gluing, the union of trapezoids is the desired domain $\mathcal D$. The example of Figure~\ref{Fig_crazy_polygon} is of this type.
\end{remark}

 The most general situation of Definition~\ref{Definition_gluing} allows two generalizations compared to the planar example. We do not require $\mathcal D$ to admit an embedding into the plane (there always exists an immersion but self-intersections are not excluded). Instead it might be only possible to embed it into a more complicated two-dimensional manifold. Second, the gluing does not have to be bipartite, \textit{i.e.} there is no designation left/right that would make blue segments in left trapezoids to be always glued to blue segments in right trapezoids. In this case, $\mathcal{D}$ can only be embedded in a non-orientable two-dimensional manifold. For instance, the non-orientable example in Figure~\ref{Fig:Domains23} correspond to a M\"obius band.

Given a gluing of trapezoids, we can speak about tilings of the resulting domain $\mathcal{D}$. By definition, these are tilings of each of the $\mathcal{T}_i$ satisfying the following conditions.
\begin{itemize}
\item For each $\mathcal{T}_i$, there must be horizontal lozenges at all the scarlet sites and horizontal lozenges are prohibited at the violet sites.
 \item For each $h \in [H]$, the configuration of lozenges in identified blue segments $\llbracket A_h,B_h \rrbracket$ must be the same in the two adjacent trapezoids to which it belongs. In particular, we have a well-defined number $N_h$ counting the number of horizontal lozenges in the $h$-th segment.
\item In each $\mathcal{T}_i$, if $W_i$ denotes the width and $S_i$ the number of scarlet sites, we have \begin{equation}
   \label{eq_trapezoid_restriction}
   \sum_{h \in \mathcal{H}_i} N_{h} = W_i - S_i,
  \end{equation}
  where $\mathcal{H}_i$ is the set of labels of blue segments in $\mathcal T_i$.
 \end{itemize}

There are finitely many tilings of a glued domain and we can speak about uniform random lozenge tilings of $\mathcal{D}$. The numbers $(N_h)_{h = 1}^H$ become \emph{random} filling fractions of the segments in $\mathcal{D}$, subjected to one affine constraint \eqref{eq_trapezoid_restriction} per trapezoid. Note that we should be careful in adjusting various parameters in the gluings, as, in principle, the system of equations \eqref{eq_trapezoid_restriction} might fail to have any positive integer solution, meaning that the glued domain has no lozenge tilings. Of course we restrict ourselves to tileable domains. For planar domains the question of tileability is well-understood, \textit{cf.} \cite[Lecture 1]{Vadimlecture} for a review. The resulting probability measure on the bottom positions of free horizontal lozenges in the segments takes the form
\begin{equation}
\label{PPPfinite}
\amsmathbb{P}(\boldsymbol{\ell}) = \frac{1}{\mathscr{Z}_{\N}} \cdot \prod_{(g,i) < (h,j)} \bigg(\frac{\ell_j^{h} - \ell_{i}^{g}}{\N}\bigg)^{2\theta_{g,h}} \cdot \prod_{(h,j),x} |x - \ell_j^h|.
\end{equation}
The second product ranges over scarlet sites $x$ (\textit{i.e.} the position of frozen lozenges containing $[x,x+1]$ as a diagonal) and occupied blue sites $\ell_j^h$ (\textit{i.e.} the position of fluctuating lozenges) in the same trapezoid. The interaction matrix $\boldsymbol{\Theta}$ has size $H \times H$ and its entry $\theta_{g,h} \in \big\{0,\frac{1}{2},1\big\}$ is equal to half the number of trapezoids that contain both the $g$-th and the $h$-th segment.
It can also be written
\begin{equation}
\label{ThetaRTR}
\boldsymbol{\Theta} = \frac{1}{2} \mathbf{R}^{T}\mathbf{R},
\end{equation}
where $\mathbf{R}$ is the $m \times H$ matrix expressing the adjacency between trapezoids and segments. In other words, this is the matrix whose rows are the linear forms indexed by $i \in [m]$
\[
\forall \boldsymbol{X} \in \amsmathbb{R}^H\qquad \mathfrak{r}_i(\boldsymbol{X}) = \sum_{h \in \mathcal{H}_i} X_{h}.
\]
In particular, $\boldsymbol{\Theta}$ is positive semi-definite. The $H$ filling fractions are subjected to $\mathfrak{e} = \textnormal{rank}(\mathbf{R})$ affine constraints. If we want (it is required for the setting of Chapter~\ref{Chapter_Setup_and_Examples} for the last point of Assumptions~\ref{Assumptions_Theta}), we can select $\mathfrak{e}$ of them which are linearly independent.

Once a tileable domain has been described, we can multiply all sides of the trapezoids by a positive integer $\N$. The uniform measure on all tilings of the rescaled domain takes the same form \eqref{PPPfinite}, but now particles live in the larger segments $\llbracket a_h,b_h\rrbracket = \llbracket \N \hat a_{h},\N\hat{b}_{h}\rrbracket := \llbracket \N (A_{h}-1)+1,\N (B_{h}+1) - 1\rrbracket$. This measure can be rewritten in the form
\begin{equation}
\label{PPPN}
\amsmathbb{P}(\boldsymbol{\ell}) = \frac{1}{\mathscr{Z}_{\N}} \cdot \prod_{(g,i) < (h,j)} \bigg(\frac{\ell_j^{h} - \ell_{i}^{g}}{\N}\bigg)^{2\theta_{g,h}} \cdot \prod_{(h,j)} w_{h}(\ell_j^h).
\end{equation}
The weights $w_{h}$ are obtained by the procedure we have illustrated in Section~\ref{sec:extile} and are products of Pochhammer symbols. They can be deduced from the weights of the probability measure for the non-rescaled domain using the substitution \eqref{substituN}. This measure fits in the scheme of discrete ensembles presented in Section~\ref{Section_general_model} and studied in this book.

\begin{figure}[!t]
\begin{center}
 {\scalebox{0.47}{\includegraphics{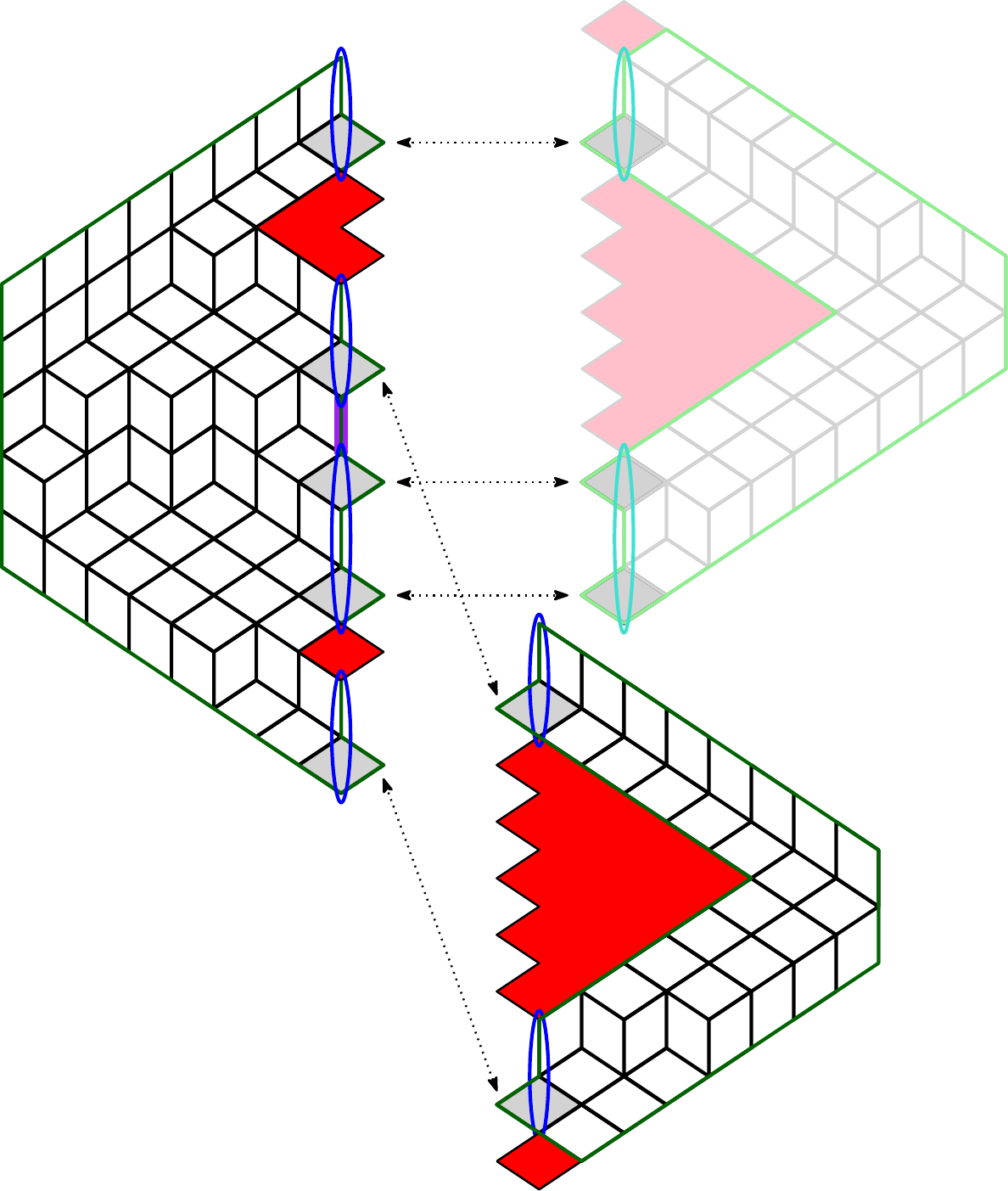}}}

 \medskip

 {\scalebox{0.47}{\includegraphics{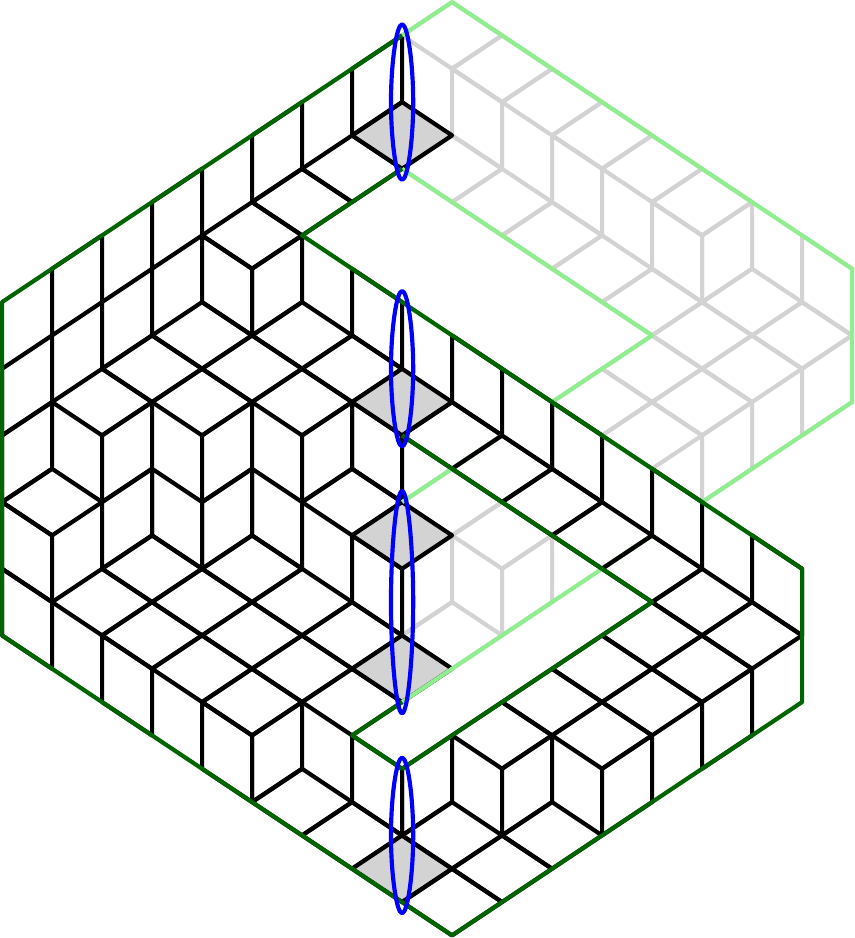}}}
 \end{center}
 \caption{A gluing of $m = 3$ tiled trapezoids (one left, two right) with $H = 4$ and resulting non-planar domain. Blue segments glued together are indicated by ovals. The scarlet frozen parts were removed for readability in the bottom panel. \label{Fig_nonplanar}}
 \end{figure}

\begin{figure}[!t]
\begin{center}
\includegraphics[width=0.9\textwidth]{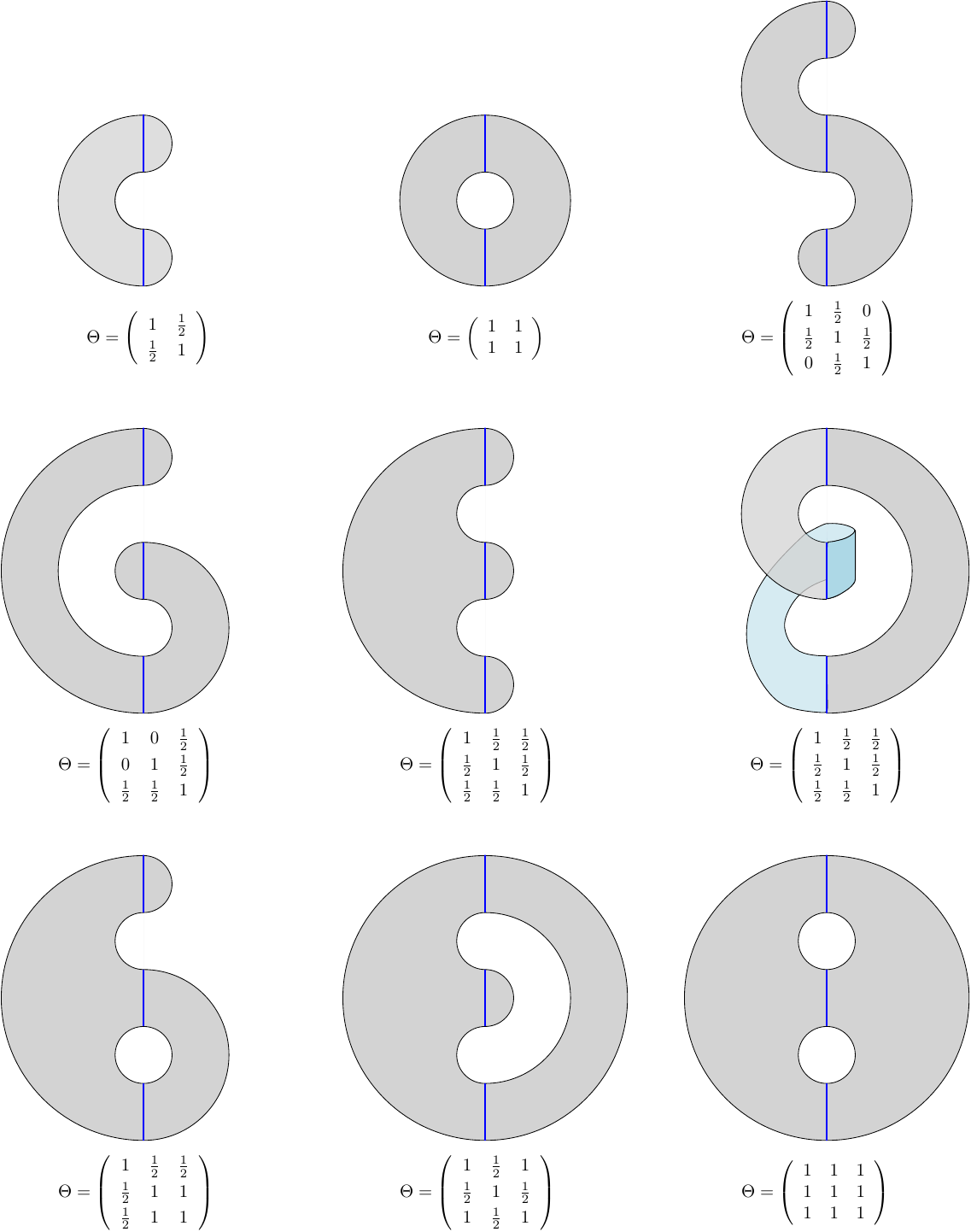}
\end{center}
\caption{\label{Fig:Domains23}Topology of domains corresponding to the gluings with $H = 2$ and $H = 3$. The E-shaped domain and the non-bipartite domain have the same interaction matrix $\boldsymbol{\Theta}$. The S-shaped and G-shaped domains differ from each other only by a reordering of blue segments, likewise for the two first domains of the last row.}
\end{figure}

\medskip

\begin{figure}[!t]
\begin{center}
\includegraphics[width=0.6\textwidth]{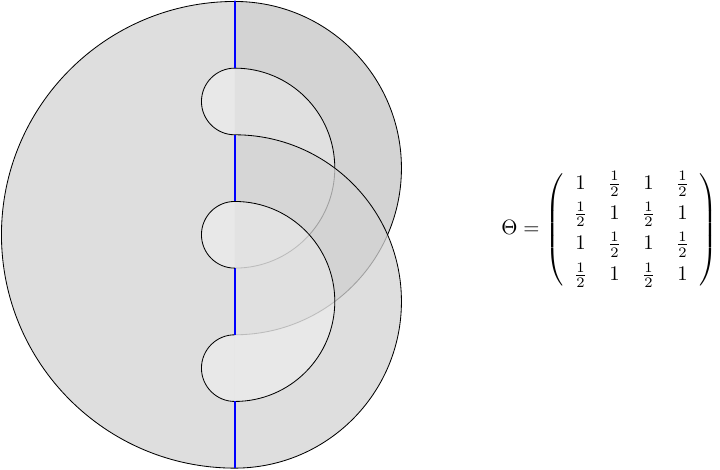}
\end{center}
\caption{A non-planar orientable domain with $H = 4$ segments and $m = 3$ trapezoids, describing the topology of Figure~\ref{Fig_nonplanar}.}
\label{Fig:Domains4}
\end{figure}

Many planar domains can be identified with glued domains $\mathcal{D}$ in Definition~\ref{Definition_gluing}. For instance, lozenge tilings of the $A\times B\times C$ hexagon are obtained in Figure~\ref{Fig_tiling_hex_section} by gluing two trapezoids. Here $H=1$ and we glue the two trapezoids along a single blue segment $\llbracket A_1,B_1\rrbracket$. The hexagon with one hole of Figure~\ref{Fig_tiling_hex_2} is obtained by gluing two trapezoids. This time we have $H=2$ groups and we glue along two blue segments $\llbracket A_1,B_1\rrbracket$ and $\llbracket A_2,B_2\rrbracket$. In the two-dimensional domain they are located along the vertical on each side of the hole. In general the domain may fail to be planar or orientable (equivalently, bipartite) but this does not prevent tileability and our study of uniformly random tilings. In Figure~\ref{Fig:Domains23} we describe all possible topologies for domains made from gluing trapezoids along $H = 2$ or $H = 3$ segments. This list contains a single non-orientable example. An orientable but non-planar example obtained by gluing $m = 3$ trapezoids along $H = 4$ blue segments is shown in Figure~\ref{Fig_nonplanar}, see also Figure~\ref{Fig:Domains4} where the topology appears more clearly. In Section~\ref{sec:C13tiling} we will study in greater details the topology and combinatorics of glued domains.

\medskip

It is possible to generalize Definition~\ref{Definition_gluing} and consider gluings in which only blue segments are being identified between trapezoids, but all the long bases together do not have to be embedded into a common oriented vertical line in a metric way. For instance, one could imagine allowing that the distance between two blue segments inside one trapezoid differs from the distance between the same segment inside another trapezoid. One could also imagine a situation where all segments are metrically embedded in the same common line but one allows changes of orientation (long bases from bottom to top may be mapped upside down). In these generalizations, the resulting probability measure for uniformly random tilings would not take the form \eqref{PPPN} as it would contain pairwise interactions of the forms $(c - \ell_i - \ell_j)$ or $(c + \ell_i - \ell_j)$. Such factors can appear as well when considering tilings with symmetries, or in random matrix theory in the $O(-2)$-matrix model \cite{GaudinKostov,BGK}. It is plausible that most arguments developed in this book can be adapted for such models, but the analogue of the Nekrasov equations of Chapter~\ref{ChapterNekra} in presence of such interactions needs to be found. We have not investigated this question further and stick to the setting of Definition~\ref{Definition_gluing}.

\section{Fluctuations of the filling fractions in tiling models}\label{sec:gen}

In this section we state and prove Lemma~\ref{Lemma_glued_tilings_assumptions} and Corollary~\ref{Corollary_discrete_Gaussian_tilings} describing asymptotics of fluctuating filling fractions for lozenge tilings of a class of domains obtained by gluing together several trapezoids.

\subsection{Discrete ensembles corresponding to tilings}

\label{Section_tiling_as_discrete_model}

We proceed to studying the tiling model defined in the general framework described in Section~\ref{sec:genglu}. We are studying uniformly random tilings of trapezoids $\mathcal{T}_1,\ldots,\mathcal{T}_m$ glued together to form a connected domain. The long base of each trapezoid is decomposed in blue, violet and scarlet segments. Blue segments of pairs of trapezoids are glued to a common vertical and the lozenges on both side should match. Along violet segments we impose our tilings to have lozenges prohibited from crossing them. We fill the scarlet segments with horizontal lozenges crossing them. We consider horizontal lozenges crossing the blue segments as particles, and the other integral sites as holes. The configuration of particles is described by a discrete ensemble whose weight functions are obtained by multiplying the result of Proposition~\ref{Proposition_number_tilings_trapezoid} over all trapezoids. Let us specify the parameters of the discrete ensemble in details.

\medskip

\noindent \textsc{Parameters and weights.} There is a total of $H$ blue segments along the common vertical. The position along the common vertical of the $h$-th blue segment is $[a_h,b_h] = [\N \hat{a}_h,\N \hat{b}_h]$ for $h \in [H]$ with integral $a_h,b_h$. The segment filling fractions indicate the number of particles they contain. They are subjected to the constraint that the number of particles in blue segments of a given trapezoid is equal to its width minus the number of frozen particles it contains. The intensity of interactions $\theta_{g,h} \in \{0,\frac{1}{2},1\}$ between particles in the $g$-th and the $h$-th segment is half the number of trapezoids containing both segments (in particular $\theta_{h,h} = 1$). The weight in the $h$-th segment $w_h(\ell)$ is proportional to the product, over all frozen particles in trapezoids that contain this segment, of their distance to the particle at position $\ell$. The proportionality factor is $\N$ to the power of total number of frozen particles. In other words, a scarlet segment $[p,q]$ in the same trapezoid as the $h$-th segment contributes to $w_h(\ell)$ by a factor
\begin{equation}
\begin{array}{ll}
\dfrac{\Gamma(\ell - p + 1)}{\N^{p - q + 1}\Gamma(\ell - q)} & \quad \textnormal{if}\,\,q < a_h, \\[2ex]
\dfrac{\Gamma(q - \ell + 1)}{\N^{p - q + 1}\Gamma(p - \ell)} & \quad \textnormal{if}\,\,p > b_h.
\end{array}
\end{equation}

\medskip

The equations \eqref{eq_equations_eqs} are produced from the following combinatorial constraint: In each trapezoid $\mathcal{T}_i$, if $W_i$ is the width and $S_i$ the number of scarlet sites, $N_h$ is the horizontal lozenges in the $h$-th segment, then we should have
\begin{equation}
   \label{eq_trapezoid_restriction_2}
   \sum_{h \in \mathcal{H}_i} N_h = W_i - S_i,
  \end{equation}
where we recall that $\mathcal{H}_i$ is the set of $h$ such that $[a_h,b_h]$ belongs to the trapezoid $\mathcal T_i$. The equations \eqref{eq_trapezoid_restriction_2} are silently assumed to be non-contradicting each other, \textit{i.e.} there should exist non-negative integers $N_h$ solving them. On the other hand, the equations might be excessive (for instance, in the hexagon case of Section~\ref{Section_Hexagon}, $H=1$, but we have two equations corresponding to two trapezoids), in which case we should omit some equations so as to get a subset of linearly independent equations with the same set of solutions as \eqref{eq_trapezoid_restriction_2}.

\medskip

\noindent \textsc{Dependence on a scale $\N$.} We consider all parameters to be large and depend on a master parameter $\N$, with some exceptions. Namely, we assume that the number of trapezoids $m$, the topological way they are glued (specified by $\boldsymbol{\Theta}$), the number and succession of blue, violet and scarlet segments along the long base of each trapezoid are independent of $\N$. We also assume that the size of each of these segments is bounded from below by $c\N$ and from above by $C\N$ for some constants $c,C > 0$, and that $|a_1| \leq c\N$.

For instance, one can start with a fixed domain and rescale all its proportions by an integer $\N$. To fully specify the $\N$-dependent discrete ensemble one should give the corresponding rule specifying the $\N$-dependent blue and scarlet segments along each long base of a trapezoid. A possible rule is to rescale the endpoints of the $\N$-dependent blue segments as $ [\N (A_h-1)+1,\N(B_h+1) - 1]$ and to replace each frozen particle at position $p$ with $\N$ frozen particles at positions $\N(p-1)+1,\N(p-1)+2,\ldots,\N p$. This is what we did in the examples of Sections~\ref{Section_Hexagon}, \ref{Section_hex_hole}, and \ref{sec:extile}.

However, our assumption also covers domains that are not obtain by a simple integer rescaling, and in general all $\hat{a}_h$ and $\hat{b}_h$ may depend on $\N$ in a complicated way. This flexibility allows for instance to get irrational ratios of proportions for the domain as $\N \rightarrow \infty$.

\medskip

\noindent \textsc{Potentials.} The potential in the $h$-th segment $V_h(x)$ receives an additive contribution from each scarlet segment $[p,q] = [\N\hat{p},\N\hat{q}]$ in the same trapezoid as the $h$-th segment. To get to the form \eqref{eq_potential_weight_match}, we treat differently the Gamma factors coming from $q$ or $p$ adjacent to $a_h$ or $b_h$. The formula for the additive contribution is
\begin{equation}
\label{eq_poteqp}
\begin{array}{ll}
\big(\hat{q} - \hat{p} + \tfrac{1}{\N}\big)\log \N - \frac{1}{\N}\log\Gamma\big(\N(x - \hat{p}) + 1\big) + \frac{1}{\N}\log \Gamma\big(\N(x - \hat{q})\big) & \quad \textnormal{if}\,\,q < a_h - 1, \\[2ex]
\big(x - \hat{p} + \frac{1}{2\N}\big) - \tfrac{1}{\N}\log \Gamma\big(\N(x - \hat{p}) + 1\big) + \mathrm{Llog}(x - \hat{a}_h') & \quad \textnormal{if}\,\,q = a_h - 1, \\[2ex]
\big(\hat{q} - x + \tfrac{1}{2\N}\big)\log\N -\tfrac{1}{\N}\log\Gamma\big(\N(\hat{q} - x) + 1\big) + \mathrm{Llog}(\hat{b}_h' - x) & \quad \textnormal{if}\,\,p = b_h + 1, \\[2ex]
 \big(\hat{q} - \hat{p} + \tfrac{1}{\N}\big)\log \N - \tfrac{1}{\N}\log\Gamma\big(\N(\hat{q} - x) + 1\big) + \tfrac{1}{\N}\log\Gamma\big(\N(\hat{p} - x)\big) & \quad \textnormal{if}\,\,p > b_h + 1,
\end{array}
\end{equation}
where $\hat{a}_h' = \hat{a}_h - \frac{1}{2\N}$ and $\hat{b}_h' = \hat{b}_h + \frac{1}{2\N}$. Consequently, the integer $\iota_h^{-}$ (resp. $\iota_h^+$) is the number in $\{0,1,2\}$ of scarlet segments adjacent to $a_h$ (resp. $b_h$) while the parameters $\rho_{h,j}^{\pm}$ are all equal to $1$. The constant $\log \N$ terms could also be moved to the error term $\mathbbm{e}_h(x)$.

\medskip

\noindent \textsc{Factorization of weights.} We can now specify the functions appearing in Definition~\ref{Definition_phi_functions}. The functions $\Phi_h^-(z)$ and $\phi_h^-(z)$ are both equal to $(z - \hat{a}_h')^{\iota_h^-}$, because $\rho_{h,j}^- = 1$. The function $\Phi_h^+(z)$ is obtained by multiplying the factors
\[
\frac{z - \hat{p} + \frac{1}{2\N}}{z - \hat{q} - \frac{1}{2\N}}
\]
for each scarlet segment $[p,q] = [\N\hat{p},\N\hat{q}]$ in a trapezoid containing the $h$-th blue segment, except that we should omit the denominator when $q = a_h - 1$ because it was already included in $\Phi_h^-(z)$. By construction, $\Phi_h^+(z)$ either does not vanish at $\hat{b}_h'$ (if $\iota_h^+ = 0$) or has a zero of order $\iota_h^+ \in \{1,2\}$ at $\hat{b}_h'$. The function $\phi_h^+(z)$ has a more complicated expression, namely it is exponential minus the sum of the $x$-derivative of all contributions to the potential \eqref{eq_poteqp}, after setting $x = z$ and analytically continuing to a complex domain.

\medskip

\noindent \textsc{The simplified variational datum.} It is inconvenient that different scarlet segments lead to different contributions to the potential in \eqref{eq_poteqp}. However, the distinction becomes negligible as $\N\rightarrow\infty$. Hence, along the same lines as in \eqref{eq_x284}, we also introduce simplified potentials $\tilde V_h$ and the corresponding equilibrium measure. The equilibrium measure $\boldsymbol{\mu}$ of the original ensemble is then approximated by the simplified equilibrium measure $\tilde{\boldsymbol{\mu}}$ as $\N\rightarrow\infty$. The latter is easier to work with. It involves the same segments $[\hat{a}_h,\hat{b}_h]$ for $h \in [H]$, and the potentials obtained from \eqref{eq_poteqp} by replacing all Gamma factors by a Stirling approximation. Concretely, for each scarlet segment $[p,q]$ in the same trapezoid as the $h$-th blue segment, $\tilde V_h(x)$ receives an additive contribution
\begin{equation}
\label{potpartinfini}
\begin{array}{ll}
-\mathrm{Llog}(x - \hat{p}') + \mathrm{Llog}(x - \hat{q}') &\quad \textnormal{if}\,\,q < a_h, \\[2pt]
-\mathrm{Llog}(\hat{q}' - x) + \mathrm{Llog}(\hat{p}' - x) & \quad \textnormal{if}\,\,p > b_h,
\end{array}
\end{equation}
where $\hat{p}' = \frac{p}{\N} - \frac{1}{2\N}$ and $\hat{q}' = \frac{q}{\N} + \frac{1}{2\N}$. The functions $\tilde{\phi}_h^{\pm}(z)$ provide the functions of Definition~\ref{Definition_phi_functions_2} for this variational datum. Namely, we have $\tilde{\phi}_h^{-}(z) = (z - \hat{a}'_h)^{\iota_h^-}$, and $\tilde{\phi}_h^{+}(z)$ is obtained by multiplying the factors
 \[
\frac{z - \hat{p}'}{z - \hat{q}'}
\]
for each scarlet segment $[p,q]$ in a trapezoid containing the $h$-th blue segment, except that we should omit the denominator when $q = a_h - 1$ because it was already included in $\tilde{\phi}_h^{-}(z)$.

\begin{lemma}
\label{Lemma_glued_tilings_assumptions} Consider a tiling model as above and choose $C>0$. Assumed that all the involved rescaled by $\N$ parameters of the model are smaller than $C$. Further assume that equations \eqref{eq_trapezoid_restriction_2} are non-contradictory, that each segment $[\hat a'_h,\hat b'_h]$, $h\in[H]$ contains at least one band of the simplified equilibrium measure $\tilde{\boldsymbol{\mu}}$ of length at least $C^{-1}$, and that the bands of $\tilde{\boldsymbol{\mu}}$ are at distance at least $C^{-1}$ from the endpoints of the segment. Then for $\N$ large enough (depending on the constant $C$), both the equilibrium measure $\boldsymbol{\mu}$ and the simplified equilibrium measure $\tilde{\boldsymbol{\mu}}$ have one band per segment. Furthermore, Assumptions~\ref{Assumptions_Theta}, \ref{Assumptions_basic}, \ref{Assumptions_offcrit}, and \ref{Assumptions_analyticity} hold for the associated discrete ensemble. And, for any $\delta > 0$ there exists a constant $C_2 > 0$ depending only on $\delta > 0$ and the constants in the assumptions such that, for any $h \in [H]$ and $z \in \amsmathbb{C}$ at distance at least $\delta$ from $[\hat{a}_h',\hat{b}_h']$ we have
\[
\big|\Gm_{\mu_h}(z) - \Gm_{\tilde{\mu}_h}(z)\big| \leq \frac{C_2}{\N^2}.
\]
\end{lemma}

The condition that equations \eqref{eq_trapezoid_restriction_2} are non-contradictory is a restatement of the condition for the existence of lozenge tilings of the domain. As we discussed in Section~\ref{Section_hex_hole}, the conditions on existence of bands in each $[\hat a'_h,\hat b'_h]$ and on the distance from the endpoints of the bands to endpoints of the segments are necessary: one can design a tiling model where they fail and the corresponding discrete ensemble does not satisfy Assumption~\ref{Assumptions_offcrit}. As we demonstrated in Section~\ref{Section_hex_hole}, with additional effort one can check these conditions in particular situations of interest. A numerical plot or Monte-Carlo simulations of the equilibrium measure can also help in this regard. In some situations, one can use the algorithm of \cite{kenyon2007limit} to compute explicitly the limit shape for the tilings and accordingly the equilibrium measure, and to check the conditions of Lemma~\ref{Lemma_glued_tilings_assumptions}. Note that Lemma~\ref{Lemma_glued_tilings_assumptions} applies both for fixed or for fluctuating filling fractions.

\begin{corollary}
\label{CLTfixedtiling}
Let $(f_h(z))_{h = 1}^{H}$ be an  tuple of functions independent of $\N$ and such that $f_h(z)$ is holomorphic in a complex neighborhood of $[\hat{a}_h,\hat{b}_h]$ for any $h \in [H]$. Under the conditions of Lemma~\ref{Lemma_glued_tilings_assumptions} and assuming that filling fractions are deterministically fixed, the random vector
\[
\left(\sum_{i = 1}^{N_h} f\bigg(\frac{\ell_i^{h}}{\N}\bigg) - \N \int_{\alpha_h}^{\beta_h} \tilde{\mu}_h(x)\dd x\right)_{h=1}^{H}
\]
is approximated in the sense of moments (Definition~\ref{Definition_convergence_moments}) by a centered Gaussian random vector with covariance
\[
\textnormal{\textsf{Cov}}[f_{h_1},f_{h_2}] = \oint_{\gamma_{h_1}}\oint_{\gamma_{h_2}} \frac{\dd z_1\dd z_2}{(2\ii\pi)^2}\,\mathcal{F}_{h_1,h_2}(z_1,z_2) f_{h_1}(z_1) f_{h_2}(z_2),
\]
where $\boldsymbol{\mathcal{F}}(z_1,z_2)$ is given by \eqref{eq_covariancepre} in terms of the bands of the simplified equilibrium measure.
\end{corollary}
\begin{corollary} \label{Corollary_discrete_Gaussian_tilings}
Under the conditions of Lemma~\ref{Lemma_glued_tilings_assumptions} and not assuming that filling fractions are deterministically fixed, the random filling fraction $(N_h)_{h=1}^H$ are asymptotically equal to a discrete Gaussian random variable, as in Theorems~\ref{Theorem_CLT_for_filling_fractions} and \ref{Theorem_CLT_for_filling_fractions_saturation}.
\end{corollary}
Corollary~\ref{CLTfixedtiling} is a direct application of Theorem~\ref{Theorem_correlators_expansion_relaxed_theta1}  and \ref{Corollary_CLT_relaxed}. A geometric understanding of the covariance will be achieved in Sections~\ref{Section_KO_conjecture}-\ref{Section_KO_proofs}, based on the results of Chapter~\ref{Chapter_AG} (see in particular Corollaries~\ref{cor:Greennormal} and \ref{cor:Greenfermion}). Corollary~\ref{Corollary_discrete_Gaussian_tilings} is a direct application of Theorem~\ref{Theorem_CLT_for_filling_fractions_saturation}. We refer to \cite[Chapter 24]{Vadimlecture} for a discussion of the parameters of the limiting (multivariate) discrete Gaussian random variable.

\subsection{Proof of Lemma~\ref{Lemma_glued_tilings_assumptions}}

The validity of Assumptions~\ref{Assumptions_Theta}, \ref{Assumptions_basic}, and \ref{Assumptions_analyticity} directly follows from the definition of the discrete ensemble, and we only prove that there is precisely one band per segment and check Assumption~\ref{Assumptions_offcrit}. Combining Lemma~\ref{Lemma_continuity_varyinf_ff} and Theorem~\ref{Theorem_off_critical_neighborhood} makes it sufficient to do so for the simplified potential and the simplified equilibrium measure $\tilde{\boldsymbol{\mu}}$. The rest of the proof focuses on this.

\smallskip

\noindent \textsc{Step 1: Convexity and connectedness.} Let $\tilde V^{\textnormal{eff}}_h(x)$ be the effective potential for the simplified equilibrium measure $\tilde{\boldsymbol{\mu}}$. We claim that for each $h\in[H]$, the function $x\mapsto \tilde V^{\textnormal{eff}}_h(x)$ is convex for $x$ in each connected component of $[\hat a_h,\hat b_h]\setminus \textnormal{supp}(\mu_h)$ . Indeed, the second derivative of \eqref{potpartinfini} is
\begin{equation}
\label{eq_tiling_convexity}
\begin{array}{ll}
-\dfrac{1}{x - \hat{p}'} + \dfrac{1}{x - \hat{q}'} &\quad \textnormal{if}\,\,q < a_h, \\[8pt]
-\dfrac{1}{\hat{q}' - x}+ \dfrac{1}{\hat{p}' - x} & \quad \textnormal{if}\,\,p > b_h,
\end{array}
\end{equation}
which is positive for $x\in[\hat a_h,\hat b_h]$. Because of the convexity of the functions $x\mapsto -\log(x-a)$ for $x>a$, and $x\mapsto -\log(a-x)$ for $x<a$, the remaining terms in the definition of the effective potential \eqref{eq_V_eff} are also convex outside the support of $\tilde{\mu}_h$.

Since $\tilde{V}^{\textnormal{eff}}_h(x)$ is a continuous function of $x$, convexity together with the characterization of Theorem~\ref{Theorem_equi_charact_repeat_2} implies that the support of the measure $\tilde{\mu}_h$ is connected (we already used this idea in the proof of Proposition~\ref{prop:convexmueq}). Moreover, the second derivative of $ \tilde{V}^{\textnormal{eff}}_h(x)$ is bounded away from $0$: if there is at least one scarlet segment, then this follows from positivity of \eqref{eq_tiling_convexity}; otherwise this follows from the fact that the band contributing to the effective potential has length at least $C^{-1}$. Hence, the second condition of Assumption~\ref{Assumptions_offcrit} necessarily holds.

\medskip

\noindent \textsc{Step 2: Hole perspective and dual convexity.} The previous step showed that the equilibrium measure $\tilde \mu_h$ cannot have bands separated by void regions and next we show that bands cannot be separated by saturated regions. To see this, we change the point of view by using the particle-hole involution. The discrete ensemble governing the hole configurations (\textit{i.e.} complements to the particles in segments $[a_h,b_h]$) is similar to the one for the particles, except for the weights and potentials. To derive it from the discrete ensemble governing the particle configurations, we use the following identity.

\begin{lemma}
\label{complementlem}Let $\llbracket c,d \rrbracket$ be a segment of integers, decomposed into two pairwise disjoint subsets $L$ and $\overline{L}$:
$\llbracket c,d \rrbracket=L\sqcup\overline{L}$. Then we have
\[
\prod_{\substack{i < j \\ i,j \in L}} (j - i) = \frac{\prod_{k = 1}^{d - c} k!}{\prod_{i \in \overline{L}} (d - i)!(i - c)!} \cdot \prod_{\substack{i < j \\ i,j \in \overline{L}}} (j - i).
\]
\end{lemma}
\begin{proof}
Another expression for the right-hand side is
\begin{equation}
\label{cijdji}
\prod_{c \leq i < j \leq d} (j - i) \cdot \prod_{\substack{i < j \\ i,j \in \overline{L}}} (j - i) \cdot \prod_{\substack{i \in \overline{L} \\ j \in \llbracket c,d \rrbracket - \{i\}}} |j - i|^{-1}.
\end{equation}
Let us count the number of factors (taking into account cancellations between numerator and denominator) involving the various types of pairs in this expression. We find one factor involving each given pair in $L$, $1 + 1 - 2 = 0$ factors involving each given pair in $\overline{L}$ (the $2$ comes from the fact that the pair of indices appears ordered in the denominator), $1 - 1 = 0$ factors involving each given pair made of one element of $L$ and one element of $\overline{L}$. Therefore, \eqref{cijdji} is also equal to the product of $(j - i)$ over all pairs $\{i,j\}\subset L$.
\end{proof}

In the tiling model, the probability $\amsmathbb{P}_\N(\boldsymbol{\ell})$ of a configuration of (fluctuating) particles $\boldsymbol{\ell}$ is obtained by multiplying over all trapezoids the products of distances between pairs of (fluctuating or frozen) particles belonging to the long base of the same trapezoid, as in Proposition~\ref{Proposition_number_tilings_trapezoid}. We use the identity in Lemma~\ref{complementlem} taking $L$ to be the sites occupied by (fluctuating or frozen) particles and $\overline{L}$ the sites occupied by (fluctuating or frozen) holes in a given trapezoid. Taking the product over all trapezoids, this rewrites $\amsmathbb{P}_\N(\boldsymbol{\ell})$, up to a proportionality constant, as the product of distances between pairs of (fluctuating or frozen) holes, times a product over all fluctuating holes $\overline{\ell}$ of
\begin{equation}
\label{dcdenom}
\frac{1}{(d - \overline{\ell})!(\overline{\ell} - c)!(d' - \overline{\ell})!(\overline{\ell} - c')!},
\end{equation}
where $[c,d]$ and $[c',d']$ are the long bases of the two trapezoids to which $\overline{\ell}$ belongs. As a result, the probability distribution of a configuration of (fluctuating) holes $\overline{\ell}$ in blue segments is also governed by a discrete ensemble. It is described in a similar way as the discrete ensemble for the particles $\ell$ with the following modifications:
\begin{itemize}
\item The roles of scarlet and violet segments are swapped. In particular, the additive contributions to the potential \eqref{eq_poteqp} are now based on violet segments.
\item Instead of \eqref{eq_trapezoid_restriction_2}, the segment filling fractions of the holes are subject to the constraints that the total number of holes in blue segments of each trapezoid is the length of the long base minus the width of the trapezoid, minus the length of violet segments.
\item The weight $\overline{w}_h(\overline{\ell})$ for a hole $\overline{\ell}$ in the $h$-th blue segment receives an extra factor \eqref{dcdenom};
\item the potential $\overline{V}_h(\overline{\ell})$ for a hole $\overline{\ell}$ in the $h$-th blue segment receives two extra terms, one for each trapezoid with long base $[c,d] = [\N\hat{c},\N\hat{d}]$ to which this segment belongs:
\begin{equation}
\label{eq_x291}
\begin{array}{ll}
\tfrac{1}{\N} \log \Gamma\big(\N(\hat{d} - x) + 1\big) + \tfrac{1}{\N}\log\Gamma\big(\N(x - \hat{c}) + 1\big) & \quad \text{if }c < a_h < b_h < d, \\[0.2ex]
\tfrac{1}{\N} \log\Gamma\big(\N(\hat{d} - x) + 1\big) + \mathrm{Llog}(x - \hat{a}_h') & \quad \text{if }c = a_h < b_h < d, \\[0.2ex]
\mathrm{Llog}(\hat{b}_h' - x) + \tfrac{1}{\N}\log\Gamma\big(\N(x - \hat{c}) + 1\big) & \quad \text{if }c < a_h < b_h = d, \\[0.2ex]
\mathrm{Llog}(\hat{b}_h' - x) + \mathrm{Llog}(x - \hat{a}_h') & \quad \text{if }c = a_h < b_h =d.
\end{array}
\end{equation}
\end{itemize}

For the simplified equilibrium measure of the holes, expressions \eqref{eq_x291} become in all four cases
\begin{equation}
\label{potholeinfiniplus}
 \mathrm{Llog}(x - \hat{c}') + \mathrm{Llog}(\hat{d}' - x),
 \end{equation}
where $\hat{c}' = \frac{c}{\N} + \frac{1}{2\N}$ and $\hat{d}' = \frac{d}{\N} - \frac{1}{2\N}$.

By the same argument as in Step 1, the potential $\overline{V}_h(\overline{\ell})$ is convex, and the simplified equilibrium measure of the holes has connected support. Since the density of (simplified) equilibrium measure of the holes is one minus the density for the particles, it follows that for the particles bands cannot be separated by saturated regions. Combined with the results of Step 1, this implies that each segment $[\hat a'_h, \hat b'_h]$ has exactly one band of the simplified equilibrium measure. In addition, the third condition of Assumption~\ref{Assumptions_offcrit} for the measure of the particles is the same as the second condition for the measure of the holes. Therefore, the latter (established as in Step 1) implies the former.

\medskip

\noindent \textsc{Step 3: Square-root behavior.} The previous results and the conditions of the theorem readily imply that the first condition of Assumption~\ref{Assumptions_offcrit} also holds, and it remains to check the conditions 4., 5., and 6. For that we recall the Tricomi formula \eqref{Tricomun} for the density of the simplified equilibrium measure: inside the band $(\alpha_h,\beta_h)$ of $\tilde \mu_h$ for $h\in[H]$, we have
\begin{equation}
\label{Tricomun_2}
\tilde \mu_h(x) = \frac{\sqrt{(\beta_h - x)(x - \alpha_h)}}{2\pi^2}\int_{\alpha_h}^{\beta_h} \frac{\tilde{V}'(x) - \tilde{V}'(y)}{x - y} \frac{\dd y}{\sqrt{(\beta - y)(y - \alpha)}}, \qquad x\in(\alpha_h,\beta_h),
\end{equation}
where $\tilde{V}$ is the modified potential defined in \eqref{eq_x292}. From that formula it follows that $\tilde{V}$ is uniformly strictly convex, and therefore the integrand in \eqref{Tricomun_2} is bounded away from $0$. Hence, \eqref{Tricomun_2} implies conditions 4.\ and 5.\ in Assumption~\ref{Assumptions_offcrit}. Condition 6.\ follows by the same argument applied to the equilibrium measure of the holes from Step 2.

\medskip

\noindent \textsc{Step 4: Mismatch of equilibrium measures after simplification.} The difference between $\boldsymbol{\mu}$ and $\tilde{\boldsymbol{\mu}}$ come from the difference between their potentials \eqref{eq_poteqp} and \eqref{potpartinfini}. The large $\N$ asymptotics of such differences were already computed around \eqref{Stirlingsimpli} which one should apply with $\rho = 1$ and $\theta = 1$. It shows that the difference in potentials is $O(\frac{1}{\N^2})$ uniformly in the segments. Theorem~\ref{Theorem_differentiability_full} then implies that for any $h \in [H]$, the difference of Stieltjes transforms $\Gm_{\mu_h}(z)- \Gm_{\tilde{\mu}_h}(z)$ is uniformly $O(\frac{1}{\N^2})$ for $z$ in any compact of $\widehat{\amsmathbb{C}} \setminus [\hat{a}_h',\hat{b}_h']$.

\section{Tilings and Gaussian free field: the Kenyon--Okounkov conjecture}

\label{Section_KO_conjecture}

In this section we continue our investigation of uniformly random lozenge tilings of domains obtained by gluing trapezoids, under the additional assumption that filling fractions are deterministically fixed. We study macroscopic fluctuations of tilings, encoded through height functions. Our goal is to show that the asymptotics is given by the Gaussian free field for a suitable complex structure, thus verifying the Kenyon--Okounkov conjecture. Our discussion is parallel to \cite[Section 4]{BuGo3}, where the case of the hexagon with a hole was investigated we also rely on the machinery of Schur generating functions developed in that paper.

\subsection{Complex structures, Gaussian free field, and complex Burgers equation}
\label{sec:complexstr}

Before formulating the Kenyon--Okounkov conjecture, it is useful to recall basic facts about complex structures, associated harmonic functions and Green functions, and how this interacts with the Gaussian free field and with solutions of the complex Burgers equation.

A \emph{complex structure} on a two-dimensional real surface\footnote{What we define here is an almost complex structure, but in two dimensions almost complex structure are automatically complex structures, \textit{i.e.} they imply the existence of holomorphic coordinates, see, \textit{e.g.} \cite[Section VII.11]{demailly1997complex}.} with local real coordinates $(x,y)$ amounts to the data of an endomorphism $\mathbf{J}$ of the (real) tangent space that depends smoothly on $(x,y)$ and such that $\mathbf{J}^2 = - \textbf{Id}$. A complex-valued function $f$ is then $\mathbf{J}$-\emph{holomorphic} when
\begin{equation}
\label{Jvder}
\partial_{\mathsf{J}(\boldsymbol{v})} f = \ii \partial_{\boldsymbol{v}} f
\end{equation}
for any tangent vector $\boldsymbol{v}$. We can represent $\mathbf{J}$ by a $2 \times 2$ matrix in the basis of unit vectors in the $x$- and $y$-direction. Note that the $2 \times 2$ real matrices squaring to $-\textbf{Id}$ have the general form
\begin{equation}
\label{eq_J_matrix}
\mathbf{J} = \left(\begin{array}{cc} \mathsf{a} & \mathsf{b} \\ \mathsf{c} & -\mathsf{a} \end{array}\right) \quad \textnormal{where}\,\, \mathsf{a}^2 + \mathsf{bc} = -1.
\end{equation}
Since $\mathbf{J}^2 = -\textbf{Id}$, a function $f$ is $\mathbf{J}$-holomorphic if and only if \eqref{Jvder} is satisfied in a single direction $\boldsymbol{v}$. For instance taking the $x$-direction, this gives the holomorphicity equation
\begin{equation}
\label{genCR}
\mathsf{a} \cdot \partial_x f + \mathsf{b} \cdot \partial_y f = \ii \partial_x f.
\end{equation}
In general, we can repackage \eqref{genCR} as a (possibly $(x,y)$-dependent) proportionality rule between the partial derivatives of $f$, and get the following equivalent definition of holomorphic functions.
\begin{definition}
\label{def:Jholo}A function $f$ is $\mathbf{J}$-holomorphic (we will simply say ``holomorphic'' when there is no confusion about the complex structure being used) if it satisfies for $\tau=\tau(x,y)\in\amsmathbb C$
\[
\partial_y f = \tau \cdot \partial_x f\qquad \text{where}\,\,\tau = \frac{\ii - \mathsf{a}}{\mathsf{b}}.
\]
\end{definition}
Conversely, the data of $\tau = \tau(x,y)$ with non-zero imaginary part uniquely determines the complex structure through \eqref{eq_J_matrix} with
\[
\mathsf{a} = - \frac{\textnormal{Re}(\tau)}{\textnormal{Im}(\tau)},\qquad \mathsf{b} = \frac{1}{\textnormal{Im}(\tau)},\qquad \mathsf{c} = -\frac{|\tau|^2}{\textnormal{Im}(\tau)}.
\]
Generally speaking, the entries of $\mathbf{J}$ are arbitrary and $\tau$ does not have to satisfy any particular equation apart from $\textnormal{Im}(\tau) \neq 0$ and smooth dependence on the surface point. The standard complex structure in the $(x,y)$-plane corresponds to $\mathbf{J}$ being the rotation by angle $\frac{\pi}{2}$, that is $\mathsf{a} = 0$, $\mathsf{b} = -\mathsf{c} = -1$ and $\tau = -\ii$. In this case, decomposing Definition~\ref{def:Jholo} into real and imaginary parts retrieves the well-known Cauchy--Riemann equations.

Associated to a complex structure, there is a notion of harmonic function. Given a holomorphic coordinate $z$, we can define locally the Laplace operator as $4\partial_z\partial_{z}^*$. If we take another holomorphic coordinate $w$, we get a different Laplace operator $|\partial_z w|^{-2} \cdot 4\partial_{w}\partial_{w}^*$. A function is \emph{harmonic} if it is annihilated by a local Laplace operator; given the previous observation, this notion does not depend on the choice of local holomorphic coordinate. It is also insensitive to replacing the complex structure $\mathbf{J}$ with $-\mathbf{J}$ (corresponding to replacing $\ii$ by its complex conjugate $-\ii$ in \eqref{Jvder}). More abstractly, the above transformation law allows the definition of the Laplace operator $\Delta$ so that its action on a function $f$ is a differential $2$-form, denoted $\Delta f$.

A \emph{Green function} $\textnormal{Gr}(p_1,p_2)$ is a real-valued function depending smoothly on two points $p_1,p_2$ in our surface and satisfying
\begin{equation}
\label{Laplaceeqn}
-\Delta_{p_1} \text{Gr}(p_1,p_2) = \delta(p_1 - p_2)
\end{equation}
where $\Delta_{p_1}$ is the Laplace operator with respect to the first point $p_1$. This equation does not depends on choice of local holomorphic coordinates because the Laplace operator transforms by the same Jacobian square factor as the Dirac distribution on the right-hand side. A local analysis shows that in any local holomorphic coordinate $z$, solutions of \eqref{Laplaceeqn} have a logarithmic singularity
\[
\textnormal{Gr}(p_1,p_2) = - \frac{1}{2\pi} \log|z_1 - z_2| + O(1) \qquad z_1 \rightarrow z_2
\]
If we have a compact surface with smooth non-empty boundary, there exists a unique Green function with Dirichlet boundary conditions (\textit{i.e.} vanishing when $p_1$ approaches the boundary) and it is symmetric $\textnormal{Gr}(p_1,p_2) = \textnormal{Gr}(p_2,p_1)$. This allows us to define the \emph{Gaussian free field}. We present a brief discussion and refer \textit{e.g.} to \cite{sheffield2007gaussian,werner2020lecture} for further technical details.
\begin{definition}\label{GFFdefinition}
 The Gaussian free field on a surface with smooth one-dimensional non-empty boundary is a (generalized) centered Gaussian field $\textnormal{GFF}(p)$ whose covariance is $\mathrm{Gr}(p_1,p_2)$.
\end{definition}

Although the Gaussian free field can be viewed as a random element of a certain functional space, the notation $\textnormal{GFF}(p)$ is slightly abusive as there is no such thing as its value at a point $p$; this is related to the singularity of $\mathrm{Gr}(p_1,p_2)$ at $p_1=p_2$. Despite that, the Gaussian free field inherits an important property of conventional functions: it can be integrated with respect to sufficiently regular measures. Avoiding an optimal description of the required regularity, we record this property in two special cases of relevance, presented without proofs.

\begin{lemma}
\label{Lemma_GFF_by_surface_measure}
 Let $\mu$ be an absolutely continuous compactly-supported complex measure on a surface, such that its density with respect to the Lebesgue measure (relative to coordinate charts) is smooth. If we have a holomorphic function $z$ on the surface which is a local coordinate near every point, we can also write the measure element $\dd\mu(p) = \mu(p)|\dd z(p)|^2$, where $\mu(p)$ is a smooth function of a point $p$ on the surface.

The Gaussian free field on this surface is such that $\int \textnormal{GFF}(p)\dd \mu(p) = \int \textnormal{GFF}(p)\mu(p)|\dd z(p)|^2$ is a well-defined centered Gaussian random variable. Moreover, if we have two such measures $\mu_1$ and $\mu_2$, then
\begin{equation*}
\begin{split}
 \amsmathbb E \left[ \bigg(\int \textnormal{GFF}(p_1)\dd\mu_1(p_1)\bigg) \cdot \bigg(\int\textnormal{GFF}(p_2)\dd\mu_2(p_2)\bigg)\right] & =\iint \textnormal{Gr}(p_1,p_2) \dd \mu_1(p_1) \dd \mu_2(p_2) \\
 & = - \int \mu_1(p) \big(\Delta^{-1} \mu_2(p)\big) |\dd z(p)|^2,
\end{split}
\end{equation*}
where $\Delta^{-1}$ is the inverse of the Laplace operator $4\partial_z\partial_{z^*}$ with Dirichlet boundary conditions on the surface.
\end{lemma}

\begin{lemma}
\label{Lemma_GFF_by_curve_measure}
Let $\nu$ be a complex measure on a surface whose support is a smooth compact curve $\gamma$, which has smooth density with respect to the Lebesgue measure relative to coordinate charts for $\gamma$, and which satisfies
\[
\int_{\gamma \times \gamma} \big|\textnormal{Gr}(p_1,p_2) \dd\nu(p_1)\dd\nu(p_2)\big| < +\infty.
\]
If we have a holomorphic function $z$ on the surface which is a local coordinate near every point of $\gamma$, we can also write $\dd\nu(p) = \mathbbm{1}_{\gamma}(p)\nu(p)\dd z(p)$ for some ``density'' $\nu(p)$ which is a smooth function of $p \in \gamma$.

The Gaussian free field on the surface is such that $\int \textnormal{GFF}(p) \dd \nu(p) = \int_{\gamma} \textnormal{GFF}(p) \nu(p) \dd z(p)$ is a well-defined Gaussian centered random variable. Moreover, if we have two such measures $\nu_1$ and $\nu_2$ supported respectively on curves $\gamma_1$ and $\gamma_2$, then
\[
\amsmathbb{E}\left[\bigg(\int \textnormal{GFF}(p_1)\dd\nu_1(p_1)\bigg) \cdot \bigg(\int \textnormal{GFF}(p_2) \dd\nu_2(p_2)\bigg)\right] = \iint_{\gamma_1 \times \gamma_2} \textnormal{Gr}(p_1,p_2) \nu_1(p_1)\nu_2(p_2) \dd z(p_1) \dd z(p_2).
\]
\end{lemma}
In principle, the above two lemmata can be taken as an alternative
definition of the Gaussian free field as a random functional on (regular enough) measures. The assumptions and the result do not depend on the choice of local holomorphic coordinates --- this is a manifestation of the fact that the Gaussian free field is essentially a function. We only mentioned local coordinates to offer a more concrete way to read formulae. As we will use these lemmata for surfaces obtained by gluing half-planes, there will indeed be a natural choice of holomorphic function $z$ which is a local coordinate near every point of the surface.

If the surface is the upper half-plane $\amsmathbb H=\{z\in\amsmathbb{C}\mid \mathrm{Im}(z)>0\}$, the Green function is fully explicit:
\begin{equation}
\label{eq_GFF_halfplane}
\textnormal{Gr}(z_1,z_2)=-\frac{1}{2\pi}\log\left|\frac{z_1-z_2}{z_1-z_2^*}\right|.
\end{equation}
The complex structure on $\amsmathbb H$ is preserved under conformal bijections, which are real M\"obius transformations. It is straightforward to check that \eqref{eq_GFF_halfplane} is also unchanged under M\"obius transformations --- a reflection of the \emph{conformal invariance} of the Gaussian free field.

\medskip

Another insightful way to think about the Gaussian free field is through uniformization. The Koebe
general uniformization theorem (see \textit{e.g.} \cite[Chapter III, Section 4]{ahlfors1960riemann}) implies that there exists a bijective
conformal (in the sense of \eqref{genCR}) map $T$
from our surface to a domain $\mathcal{U} \subset \amsmathbb{C}$ of the same topology.
The domain $\mathcal{U}$ is equipped with the standard complex structure from $\amsmathbb{C}$, and we can consider the Gaussian free field on
$\mathcal{U}$ which is defined in terms of the standard Laplace operator $\Delta = \partial_x^2 + \partial_y^2$. Then the Gaussian free field on our surface is simply the
pullback via $T$ of the Gaussian free field on $\mathcal{U}$. If our surface is simply-connected, one can choose $\mathcal{U}$ to be the upper half-plane $\amsmathbb H$. Combined with \eqref{eq_GFF_halfplane}, this gives a somewhat explicit formula for the covariance of the Gaussian free field on the surface:
\[
 \textnormal{Gr}(p_1,p_2)=-\frac{1}{2\pi}\log\left|\frac{T(p_1)-T(p_2)}{T(p_1)-(T(p_2))^*}\right|.
\]
For multiply-connected surfaces such simple formulae are not expected. In cases of direct relevance for tilings, the tools of Part~\ref{Part_Master_equation} also allow obtaining some expressions for the Green function.

\medskip

In the situations related to random tilings, the complex structure relevant for the definition of Gaussian free field arises from a distinguished function $\xi$ on the surface. This function called \emph{the complex slope} is defined in more detail in the next section. In the context of complex structures, its key property it that is solves the \emph{complex Burgers equation}:
\[
\partial_y \xi = \xi \cdot \partial_x \xi.
\]
Given any solution of this equation such that $\textnormal{Im}(\xi) \neq 0$, we can construct a complex structure $\mathbf{J}$ together with a distinguished $\mathbf{J}$-holomorphic function (namely $\xi$ itself), and thus consider the associated Gaussian free field. This is achieved simply by imposing
\[
\xi = \tau = \frac{\ii - \mathsf{a}}{\mathsf{b}},
\]
which makes the $\mathbf{J}$-holomorphicity equation of Definition~\ref{def:Jholo} for $f = \xi$ equivalent to the complex Burgers equation. In other words, $\xi = \mathsf{p} + \ii \mathsf{q}$ appears as a holomorphic function for the complex structure
\begin{equation}
\label{Jxipq}
\mathbf{J} = \left(\begin{array}{cc} -\frac{\mathsf{p}}{\mathsf{q}} & \frac{1}{\mathsf{q}} \\[2pt] -\mathsf{q} - \frac{\mathsf{p}^2}{\mathsf{q}} & \frac{\mathsf{p}}{\mathsf{q}} \end{array}\right).
\end{equation}
Note that this construction depends on the choice of local coordinates $(x,y)$, on the form of the complex Burgers equation and its solution $\xi$. Among all $\mathbf{J}$-holomorphic functions $\xi$ is special. For instance, for any choice of holomorphic function $F$ in the standard sense, $ F(\xi)$ is still $\mathbf{J}$-holomorphic (it satisfies \eqref{genCR}) but $\mathbf{J}$-holomorphic functions do not satisfy $\partial_y f = F(\xi) \cdot \partial_x f$ unless $F(\xi) = \xi$.

\subsection{Limit shapes, height functions, and complex slopes}
\label{sec:complexslopel}

Consider a polygonal domain $\Omega$ with boundary drawn on the triangular grid on the plane. We assume that all proportions of $\Omega$ scale linearly with an integer parameter $\N$, playing the same role as in the discrete ensembles. We denote $\hat{\Omega} = \N^{-1}\Omega$ the rescaled domain. We assume $\Omega$ to be tileable and encode its lozenge tilings through their height functions.

\begin{definition}
\label{Def_height_function}
For a given tiling of $\Omega$ and given point $(x,y) \in \hat{\Omega}$ represented in the coordinate system of the right panel of Figure~\ref{Figure_triangle}, the \emph{global} height function $\mathsf{Ht}_\N(x,y)$ counts the number of horizontal lozenges {\scalebox{0.16}{\includegraphics{lozenge_hor.pdf}}} directly above the closest point to $(\N x,\N y)$ in the triangular lattice inside $\Omega$.
\end{definition}

The definition admits some variants. Some authors instead count the number of lozenges of types
{\scalebox{0.16}{\includegraphics{lozenge_v_up.pdf}}} and
{\scalebox{0.16}{\includegraphics{lozenge_v_down.pdf}}} below the point $(\N x,\N y)$; others treat a tiling as the orthogonal projection of a stepped surface in $\amsmathbb{R}^3$ onto the plane $x+y+z=0$, the height being then the orthogonal projection of the stepped surface along the normal vector $(1,1,1)$. The last version can be restated as a local definition of the height function (\textit{cf.} \cite[Section 1.4]{Vadimlecture}).

\begin{figure}[t]
\centering
\includegraphics[width=0.5\textwidth]{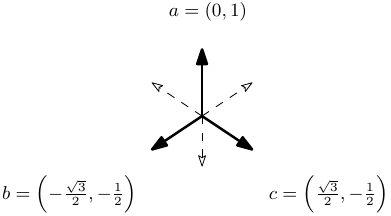}
\caption{Out of the six triangular lattice directions, three are chosen to be positive (in bold).}
\label{Figure_positive_dir}
\end{figure}

\begin{definition} \label{Def_local_height} Given a tiling of $\Omega$ and a vertex $v$ of the triangular grid inside $\Omega$, the \emph{local} height function $\mathsf{Ht}^{\textnormal{loc}}_\N(v)$ is defined by a local rule. If $u \rightarrow v$
is a positive direction, as in Figure~\ref{Figure_positive_dir}, then
\begin{equation}\label{eq_height_local}
  \mathsf{Ht}^\textnormal{loc}_\N(v)-\mathsf{Ht}^\textnormal{loc}_\N(u)=
  \left\{ \begin{array}{lll}
  1 & &\textnormal{if}\,\ u \rightarrow v\,\,\textnormal{is the edge of a lozenge,} \\
  -2 && \textnormal{if}\,\,u \rightarrow v\,\,\textnormal{crosses a lozenge diagonally.}
  \end{array}\right.
\end{equation}
\end{definition}

Note that $\mathsf{Ht}^\textnormal{loc}_\N(v)$ is determined up to a constant shift. We may assume without loss of generality that our
favorite vertex $v_0$ has $\mathsf{Ht}^\textnormal{loc}_\N(v_0)=0$. When the domain $\Omega$ is simply-connected, it may be easily checked that the local height function is defined consistently, \textit{i.e.} is single-valued. This is because the local rules are consistent for a single lozenge, and the definition extends consistently across unions. Along the boundary of a simply-connected domain the values of the local height function do not depend on the tiling, and the existence of a uniquely defined boundary-height function is equivalent to tileability of the domain, as explained in \cite{thurston1990conway} and \cite[Section 1.4]{Vadimlecture}.

If the domain is multiply-connected but holes are tileable (\textit{e.g.} as in Figure~\ref{Fig_tiling_hex_2}), then the local height function is still consistently defined, yet the boundary values of the height function are no longer uniquely determined. For general holes, the local height function can be multivalued, \textit{i.e.} make a jump when looping around the hole. This already happens for the triangular hole at the top of Figure~\ref{Fig_crazy_polygon}, which is indeed not tileable by lozenges. Note, however, that such a jump is deterministic: it depends on the shape of the hole but not on the choice of a particular tiling. Hence, if we subtract from the local height function of a random tiling its expectation, the difference is again consistently defined.

Out of all tilings of $\Omega$, we would like to study a subset singled out by the following assumption.

\begin{assumptionT}
\label{Assumption_T}
 For each connected component of the boundary of the tiled domain $\Omega$, we choose deterministically the values of the global height function $\mathsf{Ht}_\N(x,y)$ along it\footnote{Choosing instead the values of the local height function is an equivalent setting.}. We consider uniformly random lozenge tilings of $\Omega$ whose height functions agree with those prescribed boundary values.
\end{assumptionT}

When the domain $\Omega$ is simply-connected and tileable (\textit{e.g.} the $A\times B \times C$ hexagon of Section~\ref{Section_Hexagon}) and when we choose the unique boundary-height function compatible with tilings as just explained, Assumption~\ref{Assumption_T} does not bring any restriction. When $\Omega$ is multiply-connected (\textit{e.g.} the hexagon with a hole of Section~\ref{Section_hex_hole}), the assumption deterministically fixes the heights of all the holes --- although for non-tileable holes the notion of their heights is less intuitive. In the setting of Section~\ref{sec:gen} this is the same as choosing and deterministically fixing the filling fractions.

Under Assumption~\ref{Assumption_T}, the local and global height functions differ by multiplication by $3$ and a deterministic shift. This can be seen by inductively applying the rules \eqref{eq_height_local} to move from $(x,y)$ in the up direction until reaching the boundary of the domain. Hence, it does not matter which height function to study. We stick to the global height function $\mathsf{Ht}_\N(x,y)$ in our theorems, but they can be readily recast in terms of $\mathsf{Ht}^\text{loc}_\N(x,y)$.

Under Assumption~\ref{Assumption_T}, we will treat $\mathsf{Ht}_\N(x,y)$ as a random height function of a uniformly random tiling. As the linear size $\N$ of the domain grows, the height function concentrate around its mean. More precisely, there exist three non-negative continuous functions $\rho^{{\scalebox{0.16}{\includegraphics{lozenge_hor.pdf}}}}(x,y)$,
$\rho^{{\scalebox{0.16}{\includegraphics{lozenge_v_up.pdf}}}}(x,y)$,
$\rho^{{\scalebox{0.16}{\includegraphics{lozenge_v_down.pdf}}}}(x,y)$ defined over $(x,y) \in \hat{\Omega}$ and such that
\[
\rho^{{\scalebox{0.16}{\includegraphics{lozenge_hor.pdf}}}}(x,y)+\rho^{{\scalebox{0.16}{\includegraphics{lozenge_v_up.pdf}}}}(x,y)+
\rho^{{\scalebox{0.16}{\includegraphics{lozenge_v_down.pdf}}}}(x,y)=1.
\]
These functions encode the asymptotic local proportions of three types of lozenges in a large tiling. In the same way as $\mathsf{Ht}_\N(x,y)$ is defined by counting ${\scalebox{0.16}{\includegraphics{lozenge_hor.pdf}}}$ lozenges above $(\N x,\N y)$, its asymptotic version $\mathfrak{Ht}(x,y)$ is defined as
\begin{equation}
\label{eq_x293}
 \mathfrak{Ht}(x,y)=\int_y^{+\infty} \rho^{{\scalebox{0.16}{\includegraphics{lozenge_hor.pdf}}}}(x,\tilde{y}) \dd\tilde{y}.
\end{equation}
The density $\rho^{{\scalebox{0.16}{\includegraphics{lozenge_hor.pdf}}}}(x,y)$ vanishes outside $\hat{\Omega}$ and the integration segment in the above integral is finite.

\begin{theorem} \label{Theorem_tilings_LLN}
 There exists a deterministic function $\rho^{{\scalebox{0.16}{\includegraphics{lozenge_hor.pdf}}}}$ defined on $\Omega$ such that
\[
\lim_{\N \rightarrow \infty} \sup_{(x,y) \in \hat{\Omega}} \left|\frac{\mathsf{Ht}_\N(x,y)}{\N}-\mathfrak{Ht}(x,y)\right|=0\quad \textnormal{ in probability},
\]
using \eqref{eq_x293} to define the \emph{limit shape} $\mathfrak{Ht}(x,y)$.
\end{theorem}
\begin{remark} \label{Remark_L_independence}
Although we omitted $\N$ from most notations, in general the densities of lozenges $\rho$ and the limit shape $\mathfrak{Ht}$ depend on $\Omega$ and $\N$, which is the reason why we did not state a convergence in probability of $\N^{-1}\mathsf{Ht}_{\N}$ to $\mathfrak{Ht}$. This is in the spirit of all asymptotic results of this book. Yet, if the domains $\Omega$ grow regularly with $\N$ (\textit{e.g.} if the rescaled domain $\hat{\Omega}$ approaches a limit independent of $\N$), the functions $\rho$ and $\mathfrak{Ht}$ have $\N \rightarrow \infty$ limits which can be used in Theorem~\ref{Theorem_tilings_LLN} and all results mentioned later. In what follows, we will silently adapt the same convention, assuming that the domains grow regularly with $\N$ so that asymptotic proportions are $\N$-independent. Note that the results for $\N$-dependent asymptotic quantities still hold, if we keep all dependencies on $\N$ and replace evaluations of the limits by the statement that the difference between the quantities and their $\N$-dependent asymptotic values tends to $0$ as $\N \rightarrow \infty$.
\end{remark}

For simply-connected domains Theorem~\ref{Theorem_tilings_LLN} was proven in \cite{cohn2001variational} by establishing a variational principle: $\mathfrak{Ht}(x,y)$ is the minimizer of a surface tension functional, see also \cite[Lectures 5-10]{Vadimlecture}. For multiply-connected domains it is similar, yet the jumps in the height function complicate the proofs, \textit{cf.} \cite{kuchumov2021variational}. For gluings of trapezoids as introduced in Section~\ref{Section_gluing_def}, we provide an independent proof of Theorem~\ref{Theorem_tilings_LLN} in Section~\ref{Section_limit_shape_gluings}.

All three densities $\rho^{{\scalebox{0.16}{\includegraphics{lozenge_hor.pdf}}}}(x,y)$, $\rho^{{\scalebox{0.16}{\includegraphics{lozenge_v_up.pdf}}}}(x,y)$, $\rho^{{\scalebox{0.16}{\includegraphics{lozenge_v_down.pdf}}}}(x,y)$ can be read off $\mathfrak{Ht}(x,y)$: from \eqref{eq_height_local}, one needs to pass to the limit of the local height function, and then compute its derivatives in three directions indicated in Figure~\ref{Figure_triangle}. One can also restate Theorem~\ref{Theorem_tilings_LLN} without directly referring to the height functions: it is equivalent to the statement that normalized counts for lozenges of each of the three types in (regular enough) subdomains of $\hat{\Omega}$ converge towards integrals of densities over the subdomains.

The limit shape $\mathfrak{Ht}(x,y)$ has two types of local behavior, clearly seen in
Figure~\ref{Fig_hex_simulation}. In \emph{frozen regions} $\mathfrak{Ht}(x,y)$ is linear and its slope takes one of the three extreme values corresponding to three
different types of lozenges (this is parallel to voids and saturated regions for discrete ensembles). In the \emph{liquid region} of $\hat{\Omega}$, the limit shape $\mathfrak{Ht}(x,y)$ is curved and all three types of lozenges are present with varying density (some authors call it the \emph{rough phase}).

\begin{figure}[t]
\begin{center}
 {\scalebox{1.0}{\includegraphics{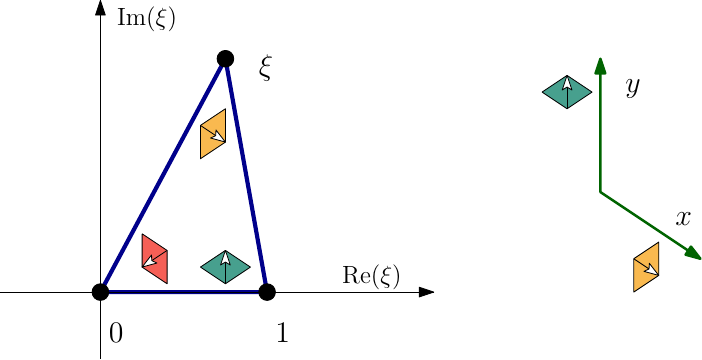}}}
 \caption{The complex slope $\xi$ encodes the three local proportions of lozenges through a geometric construction. Proportions can be read off partial derivatives of the limit shape in three directions. The coordinate directions for the complex Burgers equation can be also encoded by lozenges.
 \label{Figure_triangle}}
\end{center}
\end{figure}

Consider a triangle on the complex plane with angles $\pi\rho^{{\scalebox{0.16}{\includegraphics{lozenge_v_down.pdf}}}}(x,y)$, $\pi\rho^{{\scalebox{0.16}{\includegraphics{lozenge_hor.pdf}}}}(x,y)$, $\pi\rho^{{\scalebox{0.16}{\includegraphics{lozenge_v_up.pdf}}}}(x,y)$ at the respective vertices $0,1,\xi$, imposing that $\xi$ lies in the upper half-plane (\textit{cf.} Figure~\ref{Figure_triangle}). This $\xi$ is called the \emph{complex slope} of the limit shape $\mathfrak{Ht}(x,y)$. Kenyon and Okounkov \cite{kenyon2007limit} showed that the complex slope satisfies the complex Burgers equation everywhere inside the liquid region $\pounds \subset \hat{\Omega}$
\begin{equation}
\label{eq_Burgers}
 \partial_y \xi = \xi \cdot \partial_x \xi.
\end{equation}
It can be used to define a complex structure on $\pounds$ making $\xi$ holomorphic in the sense of Definition~\ref{def:Jholo}, as explained in Section~\ref{sec:complexstr}.

Some choices have been made to define the complex slope, and we may wonder how they affect the complex structure and the Gaussian free field. The short answer is that we always get either the same complex structure $\mathbf{J}$, or $-\mathbf{J}$ corresponding to complex conjugation. First, the assignment between a type of lozenge and the angles of the triangle of Figure~\ref{Figure_triangle} was rather arbitrary. If we reshuffle them, we get a modified complex slope
\begin{equation}
\label{eq_alternative_form}
\tilde{\xi} \in \left\{1 - \frac{1}{\xi}\,,\,\frac{1}{1 - \xi}\,,\,\frac{1}{\xi^*}\,,\,1 - \xi^*\,,\,\frac{\xi^*}{\xi^* - 1}\right\}.
\end{equation}
Although these modified slopes do not satisfy the same complex Burgers equation, they can be used as well to define complex structures in which they are holomorphic functions (although not the distinguished one of Section~\ref{sec:complexstr}). For instance, \cite{Petrov_Airy,petrov2015asymptotics} works with another convention where the coordinates are denoted $(\chi,\eta)$ and the complex slope $\Omega$ taking value in the upper-half plane satisfies the modified equation
\[
\partial_{\chi}\Omega = \bigg(1 - \frac{1}{\Omega}\bigg)\partial_{\eta}\Omega,
\]
so we cannot use it exactly as before to define a complex structure. Instead, we notice that $\xi = 1 - \frac{1}{\Omega}$ satisfies the usual Burgers equation \eqref{eq_Burgers} with $(x,y) = (\eta,\chi)$. The procedure of Section~\ref{sec:complexstr} gives a complex structure making $\xi$ holomorphic, but then $\Omega = \frac{1}{1 - \xi}$ is also holomorphic. We are back with the same complex structure $\mathbf{J}$. It is in the list \eqref{eq_alternative_form}, so this different convention corresponds indeed to reshuffling of angles. Our choice of the coordinate system in Figure~\ref{Figure_triangle} was also quite arbitrary and we could have chosen two other triangular lattice directions as $x$ and $y$. If we had swapped coordinates by setting $(x,y) = (\chi,\eta)$, we would rather get the complex structure $-\mathbf{J}$. More generally, changing the coordinate system $(x,y) \mapsto (\tilde{x},\tilde{y})$ by a symmetry of the triangular grid corresponds to change of bases and simultaneously transforms the matrix for $\mathbf{J}$ as well as the Burgers equation, in ways that compensate each other up to a sign.

\subsection{Kenyon--Okounkov conjecture and main results}

The Kenyon--Okounkov conjecture is a description of the next-order stochastic term\footnote{In contrast, there is no good exact prediction in the literature on the deterministic shift $\N^{-1}\E
[\mathsf{Ht}_\N(x,y)]-\mathfrak{Ht}(x,y)$.} in Theorem~\ref{Theorem_tilings_LLN}, governing the fluctuations of $\mathsf{Ht}_\N(x,y)$ around the limit shape. In parallel to Theorems~\ref{Theorem_ldpsup} and \ref{Theorem_ldsaturated}, as $\N \rightarrow \infty$ frozen regions exhibit only one type of lozenges among the three, and therefore there are no fluctuations. On the other hand, the prediction says
that in the liquid region $\pounds$ the recentered height function $\mathsf{Ht}_\N(x,y)-\E[
\mathsf{Ht}_\N(x,y)]$ converges to the \emph{Gaussian free field} for the complex structure defined by the complex slope, itself specified by the limit shape. Although the definition of the complex slope depends on some choices, we have explained that different choices either lead to the same complex structure or its complex conjugate, so that the notion of Gaussian free field and the final statement is independent of these choices.

\begin{conjecture}
\label{Conjecture_GFF}
 The random field $\sqrt{\pi}\big(\mathsf{Ht}_\N-\E[\mathsf{Ht}_\N]\big)$ converges as $\N \rightarrow \infty$ in
 the liquid region $\pounds$ to the Gaussian free field with respect to the complex slope
 $\xi$ and with Dirichlet boundary
 conditions.
\end{conjecture}

Conjecture~\ref{Conjecture_GFF} can be traced back to \cite[Section 2.3]{kenyon2007limit}. A detailed exposition and heuristics of the Kenyon--Okounkov predictions can be found in \cite[Lectures 11-12]{Vadimlecture}. Our main result is a partial proof of Conjecture~\ref{Conjecture_GFF}.

\begin{theorem}
\label{Theorem_GFF}
For \emph{planar} domains represented as gluings of trapezoids (see \textit{e.g.} Remark~\ref{Example_gluing}), Conjecture~\ref{Conjecture_GFF} is true under conditions of Lemma~\ref{Lemma_glued_tilings_assumptions} and Assumption~\ref{Assumption_T}.
\end{theorem}
Decoding the conditions of the theorem, we impose Assumption~\ref{Assumption_T} (saying that filling fractions in the discrete ensemble are deterministically fixed) and further require that all the involved rescaled by $\N$ parameters are bounded and bounded away from $0$, that each segment $[\hat a'_h,\hat b'_h]$, $h\in[H]$ contains at least one band of the simplified equilibrium measure $\tilde{\boldsymbol{\mu}}$ and that the bands of $\tilde{\boldsymbol{\mu}}$ are not touching the endpoints of the segments.

\medskip

In Section~\ref{Section_gluing_def} we have introduced lozenge tiling models on domains obtained by gluings of trapezoids which do not have to be planar. Our results also apply to this generalization. The height function can still be naturally defined: for a point $p$ in a given trapezoid $\mathcal T_j$, we count horizontal lozenges above $p$ inside $\mathcal T_j$, ignoring all other trapezoids. This height function has a discontinuity when we pass from one trapezoid to another and might differ from \eqref{Def_height_function} by deterministic shifts. However, both discontinuities and shifts cancel out when we compute the difference $\mathsf{Ht}_\N -\E[\mathsf{Ht}_\N]$.

For the definition of the Gaussian free field we need a complex structure. As we do not assume planarity, our choice of coordinate $(x,y)$ in a trapezoid may not extend to global coordinates in the domain. To define the complex structure we just need to extend smoothly the definition of partial derivatives $\partial_x$ and $\partial_y$ from a starting trapezoid to the whole domain, so that they match at junctions of two trapezoids. A property of our gluings is that the unit vector pointing up in the grid is globally defined, so $\partial_y$ is globally defined. If the domain is orientable (\textit{i.e.} bipartite), we can choose a grid direction defining $\partial_x$ in a starting trapezoid (like in Figure~\ref{Figure_triangle}) and propagate it smoothly to the whole domain. We then get a complex slope $\xi$ defining a complex structure in the liquid region $\pounds$ and can consider the associated Gaussian free field. If the domain is not orientable, this is impossible as it would force a discontinuity at least of the segments joining two trapezoids. Let us first state the result in the orientable case.

\begin{theorem}
\label{Theorem_GFF_general}
 Under conditions of Lemma~\ref{Lemma_glued_tilings_assumptions}, assume that all the filling fractions of the corresponding discrete ensemble are deterministically fixed, and that the glued domain is bipartite (trapezoids $\mathcal T_1,\ldots,\mathcal T_m$ can be partitioned into ``left'' and ``right'', so that left trapezoids can only be adjacent to right trapezoids). Then the random field $\sqrt{\pi}\big(\mathsf{Ht}_\N-\E[\mathsf{Ht}_\N]\big)$ converges as $\N \rightarrow \infty$ in
 the liquid region $\pounds$ to the Gaussian free field with respect to the complex slope $\xi$ and with Dirichlet boundary
 conditions.
\end{theorem}
We establish the convergence to the Gaussian free field in Theorems~\ref{Theorem_GFF} and \ref{Theorem_GFF_general} for finite-dimensional
distributions of pairings with a specific class of test-measures, which are supported on $x=\mathrm{const}$ lines and whose density depends on the $y$-coordinate in a polynomial way (\textit{cf.} Lemma~\ref{Lemma_GFF_by_curve_measure}). The proof of Theorem~\ref{Theorem_GFF_general} (from which Theorem~\ref{Theorem_GFF} is a particular case) is given in Section~\ref{Section_GFF_proof}, after a detailed discussion of the limit shapes and complex structure for gluings of trapezoids in Sections~\ref{Section_limit_shape_gluings} and \ref{Section_complex_structure}.

\begin{figure}[t]
 \begin{center}
\includegraphics[width=0.27\linewidth]{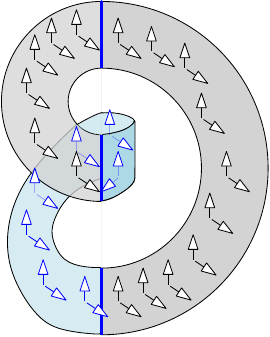}
\end{center}
\caption{Example of a non-orientable gluing: Propagating $(\partial_x,\partial_y)$ continuously from a trapezoid leads to a discontinuity. In the picture, the basis of vectors turns like the blue ribbon. \label{fig:nonorient3-4}}
\end{figure}

For non-orientable gluings (see \textit{e.g.} Figure~\ref{fig:nonorient3-4}) an extra ingredient is needed: we have to look at the orientation covering of the domain. This is a double covering $\widetilde{\mathcal{D}} \rightarrow \mathcal{D}$ such that $\widetilde{\mathcal{D}}$ is connected and orientable (such a covering exists uniquely for any non-orientable manifold). Each trapezoid $\mathcal{T}_i$ is lifted to two twin trapezoids in $\widetilde{\mathcal{D}}$: a 'left' one $\mathcal{T}_i^{\textnormal{L}}$ and a 'right' one $\mathcal{T}_i^{\textnormal{R}}$. When $\mathcal{T}_i$ and $\mathcal{T}_j$ share a blue segment, two copies of this segment appear in $\widetilde{\mathcal{D}}$: one of them is shared by $\mathcal{T}_i^{\textnormal{L}}$ and $\mathcal{T}_j^{\textnormal{R}}$, and the other by $\mathcal{T}_i^{\textnormal{R}}$ and $\mathcal{T}_j^{\textnormal{L}}$. Instead of uniformly random lozenge tilings in $\mathcal{D}$, we can equivalently work with uniformly random lozenge tilings in $\widetilde{\mathcal{D}}$ conditioned to be identical in each pair of twin trapezoids. As before, we can define in $\widetilde{\mathcal{D}}$ the height field $\widetilde{\mathsf{Ht}}_{\N}$, the complex slope $\xi$ and the associated complex structure on the lift $\widetilde{\pounds} \subset \widetilde{\mathcal{D}}$ of the liquid region $\pounds \subset \mathcal{D}$. To match the symmetry of the height field, we should compare to a symmetric conditioning of the Gaussian free field on $\widetilde{\pounds}$.

\begin{theorem}
\label{Theorem_GFF_general_nonor}
Under conditions of Lemma~\ref{Lemma_glued_tilings_assumptions}, let us assume that all the filling fractions of the corresponding discrete ensemble are deterministically fixed, and that the glued domain is not bipartite. Then the random field $\sqrt{\pi}\big(\widetilde{\mathsf{Ht}}_\N-\E[\widetilde{\mathsf{Ht}}_\N]\big)$ converges as $\N \rightarrow \infty$ in
 the orientation covering $\widetilde{\pounds}$ of the liquid region to $\sqrt{2}$ times the symmetric conditioning of the Gaussian free field with respect to the complex slope $\xi$ and with Dirichlet boundary conditions.
 \end{theorem}

As we explain in Section~\ref{sec:nonorGFFproof}. the proof of Theorem~\ref{Theorem_GFF_general_nonor} only requires minor adaptations compared to Theorem~\ref{Theorem_GFF_general}. This section will also describe in greater detail this symmetric conditioning of the Gaussian free field, in Definition~\ref{even:green} and Lemma~\ref{lem:evenGFFGreen}.

\begin{remark}\label{rem:sq2}
 The $\sqrt{2}$ prefactor in Theorem~\ref{Theorem_GFF_general_nonor} may appear surprising at first glance. $\widetilde{\mathsf{Ht}}_\N$ can be viewed as a height function for tilings of the orientable cover, conditioned to be symmetric, and one may have expected the limit to be exactly the symmetric conditioning of the Gaussian free field. However, this is misleading. Instead, a convenient way to describe the normalization is through the local behavior of the covariance, which should be $ - \frac{1}{2\pi} \log|z_1 - z_2| + O(1)$ as $z_1 \rightarrow z_2$ in all lozenge tilings models. This matches the $\sqrt{2}$ prefactor, see Definition~\ref{even:green} and Lemma~\ref{lem:evenGFFGreen} for some details.

A similar phenomenon can be found in a much simpler situation of lozenge tilings of the $A\times A\times A$ hexagon which are symmetric with respect to the vertical axis. These were studied \textit{e.g.} in \cite{panova2015lozenge}, which proved local convergence towards the GUE-corners process, and the latter can be used to predict (\textit{e.g.} through \cite{borodin2014clt}) the local behavior $-\frac{1}{2\pi} \log|z_1 - z_2| + O(1)$ for the covariance of the height function. Using Proposition~\ref{Proposition_number_tilings_trapezoid}, restriction of such symmetric tilings on the vertical axis can be described by a probability measure on the integer positions $1 \leq \ell_1 < \cdots < \ell_A \leq 2A$ of horizontal lozenges of the form
\begin{equation} \label{eq_x306}
\amsmathbb{P}(\ell_1,\ldots,\ell_A) = \frac{1}{\mathscr{Z}} \prod_{1 \leq i < j \leq A} (\ell_j - \ell_i).
\end{equation}
 Compared to \eqref{eq_Hahn_ensemble}, note the exponent $2\theta = 1$ (instead of $2$) for pairwise interaction, because the tiling on the right trapezoid is fixed by the tiling on the left. This looks like a $H=1$ discrete ensemble \eqref{eq_lattice_restriction},\eqref{eq_general_measure_one_cut}, except that we have $(\ell_{i + 1} - \ell_i) \in \{1,2,3,\ldots\}$ instead of $(\ell_{i + 1} - \ell_i) \in \{\frac{1}{2},\frac{3}{2},\frac{5}{2},\ldots\}$. Some adaptations of our proofs would be necessary to handle this change of lattice. Yet, we believe that it would not affect the leading covariance: it should still be given as in Theorem~\ref{Theorem_correlators_expansion} by the fundamental solution presented in Definition~\ref{def:Berg}. This is a reasonable expectation because the ensemble allowing $\ell_i \in \amsmathbb{R}$ to vary continuously has same leading covariance as the discrete ensemble, see Remark~\ref{cdiscrem}. From Proposition~\ref{thmBfund} (see as well \cite{Johansson}, or \cite{BGG} for this simple $H = 1$ case) it is clear that multiplying all entries of $\boldsymbol{\Theta}$ by the same factor $t > 0$ corresponds to multiplying the field of fluctuations by $t^{-1/2}$. Therefore, from \eqref{eq_x306} and the result for tilings without symmetries, for symmetric tilings of the hexagon we expect the height field along the vertical axis to converge to $(1/2)^{-1/2} = \sqrt{2}$ times the restriction of the Gaussian free field. Since the vertical is the symmetry axis, this coincides with $\sqrt{2}$ times the symmetric conditioning of the Gaussian free field. The full result on convergence of the two-dimensional field of fluctuations of heights for symmetric tilings of $A\times A\times A$ hexagon can probably be proven by combining the results of \cite{BuGo2} and \cite{panova2015lozenge}.

 An important difference between symmetric tilings of the hexagon and symmetric tilings of an orientation covering in Theorem~\ref{Theorem_GFF_general_nonor} is that the involution in the first case has a reflection axis, while in the second case it has no fixed points. Yet, in both cases the field of fluctuations of the height function is $\sqrt{2}$ times the symmetric conditioning of the Gaussian free field.

 A similar occurrence of $\sqrt{2}$ times an even conditioning of a Gaussian free field has been discussed in \cite[Section 1.4]{Besqrt2}.
\end{remark}

\medskip

Let us outline the existing literature related to Conjecture~\ref{Conjecture_GFF}. The seminal contributions \cite{kenyon2000conformal,kenyon2001dominos,kenyon2008height} developed
a method based on determinantal structure of correlation functions and exploiting the convergence of discrete
harmonic functions to their continuous counterparts. Later, these ideas were further developed in \cite{li2013conformal}, \cite{berestycki2016note,berestycki2020dimers,berestycki2024dimers1,berestycki2024dimers}, \cite{russkikh2018dimers,russkikh2020dominos}, \cite{basok2023dimers,nicoletti2025temperleyan}. The method works well for a particular class of domains with \emph{no} frozen regions, yet polygonal domains of our interest remain out of reach with it.

A new crucial set of ideas was added to the above methodology more recently through the introduction in \cite{chelkak2023dimer,chelkak2021bipartite} of t-embeddings and origami maps, which are expected to serve as discrete versions of the complex structure in the Kenyon--Okounkov conjecture. While there is a hope that this approach will eventually prove the convergence to Gaussian free field for very general domains, verifications of technical details in general case remain elusive, although they were checked in special situations in \cite{berggren2024perfect,berggren2024perfect2}.

A different approach towards the convergence to the Gaussian free field is pursued in \cite{borodin2014anisotropic,duits2013gaussian,petrov2015asymptotics,duits2018global,berggren2025gaussian}. These articles exploit exact formulae for the correlation functions in terms of contour integrals or orthogonal polynomials and perform direct asymptotic analysis of those formulae. Two more methods were fruitful: \cite{huang2020height,gorin2024dynamical,dimitrov2024global} use dynamical loop equations, which are relatives of the Nekrasov equations of Chapter~\ref{ChapterNekra}; \cite{BuGo2, bufetov2018asymptotics, BuGo3, boutillier2021limit, ahn2020global, ahn2022lozenge} rely on actions of differential/difference operators of Schur generating function on the probability measures of our interest.

Despite numerous previous approaches and results, our Theorem~\ref{Theorem_GFF} is significantly different. The main novelty is our ability to handle the situations in which the tiled domain both has a non-trivial topology and leads to rich (not flat) limit shapes for the height functions. We also have not seen detailed studies of the asymptotic heights for random tilings of non-orientable domains.

If we want to generalize the Kenyon--Okounkov conjecture and its non-orientable version to tilings on surfaces which are not necessarily gluings of trapezoids, it is important to stick to the situations with a well-defined height function. For tilings on a general surface a proper definition of the height function can be tricky: Definition~\ref{Def_local_height} requires a global consistent choice of the positive directions on the lattice and also does not guarantee that the height function is unchanged when we follow a non-contractible loop in the domain. Our gluings of trapezoids are somewhat special: despite the absence of global orientation, we can still define a global height function, whose jumps on the gluing line are deterministic due to Assumption~\ref{Assumption_T}. This happens, because both the ``up'' direction is preserved in gluings (in fact, all three positive directions of Figure~\ref{Figure_positive_dir} are preserved), and we fixed the filling fractions in blue segments. For general domains that are not necessarily gluings of trapezoids, the analogue of filling fractions are the jumps of the height function along loops in the domain. So, we conjecture that for all domains with a well-defined height function that is conditioned to have deterministic jumps along any loop, the fluctuations of the height field in the limit where the domain becomes large is given by the Gaussian free field (if the surface is orientable) or $\sqrt{2}$ times the symmetric conditioning of the Gaussian free field after lifting to the orientation covering (if the surface is not orientable).

\subsection{Gaussian free field in presence of fluctuating filling fractions}

In all the theorems of Section~\ref{Section_KO_conjecture}, the filling fractions were fixed through Assumption~\ref{Assumption_T}. One could ask about extensions of Theorems~\ref{Theorem_GFF_general} and \ref{Theorem_GFF_general_nonor} to the situation when filling fractions (equivalently, the relative heights of different connected components of the boundary of the domain) are allowed to fluctuate. In principle, the answer is delivered by the combination of Corollary~\ref{Corollary_discrete_Gaussian_tilings} (giving the asymptotics of the filling fractions) with Theorems~\ref{Theorem_GFF_general} and \ref{Theorem_GFF_general_nonor} (giving the asymptotics of the height function conditionally to the values of the discrete components). Eventually, like in the theorems of Chapter~\ref{Chapter_filling_fractions}, the fluctuations of heights would be asymptotically equivalent to a sum of two random independent terms: the Gaussian free field from Theorems~\ref{Theorem_GFF_general}-\ref{Theorem_GFF_general_nonor} and a discrete Gaussian from Corollary~\ref{Corollary_discrete_Gaussian_tilings} extended to a two-dimensional field.

Note that in Theorems~\ref{Theorem_GFF_general}-\ref{Theorem_GFF_general_nonor} the field was centered with respect to its expectation, and when we combine with Corollary~\ref{Corollary_discrete_Gaussian_tilings}, the expectation turns into the conditional expectation with respect to the values of the filling fraction. The discrete Gaussian from Corollary~\ref{Corollary_discrete_Gaussian_tilings} is being extended to the full domain as this conditional expectation. We believe that the difference between the conditional and unconditional expectation of the height in the liquid region should be asymptotically given by the unique harmonic function (with respect to the complex structure of Section~\ref{sec:complexslopel}) whose boundary values are prescribed by the height increments of the connected components of the boundary of the domain (these increments asymptotically governed by discrete Gaussian laws), see the related discussion in \cite[Section 24.4]{Vadimlecture}. For our general domains obtained by gluings of trapezoids, the link to harmonic extensions has not been rigorously justified yet, but it is very plausible that this can be done using the tools of \cite{BuGo,BuGo2,BuGo3}.

A more delicate part of the asymptotic fluctuations of the height function is the shift in discrete Gaussians of Corollary~\ref{Corollary_discrete_Gaussian_tilings}, which was $\boldsymbol u$ in $\textnormal{\textsf{\textbf{Gau\ss{}}}}_{\amsmathbb{Z}}[\Qu,\L,\boldsymbol u]$ of Theorem~\ref{Theorem_CLT_for_filling_fractions_saturation}. We do not have any explicit descriptions of $\boldsymbol u$ and the random tilings literature also does not provide any guidance on its identification. It would extremely interesting to have a better understanding of this part of the answer.

In the setting of random tilings on the torus, the asymptotic height field governed by the sum of the Gaussian free field and a discrete component was called the \emph{compactified Gaussian free field} in \cite{Dubedator}. In some examples where the limit shape is trivial (this is the case for the full torus, \textit{cf.} \cite[Lecture 3--6]{Vadimlecture}) the situations leading to the compactified Gaussian free field were analyzed for instance in \cite{basok2023dimers,nicoletti2025temperleyan}. For periodically weighted domino tilings of the Aztec diamond the same limiting object appeared very recently in \cite{berggren2025gaussian}.

\section{Proofs for the Kenyon--Okounkov conjecture}

\label{Section_KO_proofs}

\subsection{Limit shape in individual trapezoids}

\label{Section_limit_shape_gluings}

The goal of this section is to prove Theorem~\ref{Theorem_tilings_LLN} for gluings of trapezoids by combining Theorem~\ref{Theorem_main_LLN} and Corollary~\ref{Corollary_CLT_relaxed} with the results of \cite{BuGo3}. We also establish notations which are helpful for the next sections. The reader might find it helpful to review the gluing procedure of Section~\ref{sec:genglu}.

We start by looking at the particles $\ell_1<\ell_2<\cdots<\ell_N$ which encode the varying positions of horizontal lozenges in blue segments along the vertical line where the trapezoids are glued together. In Section~\ref{Section_tiling_as_discrete_model} we explained that their distribution is governed by the discrete ensemble of Section~\ref{Section_general_model}.

We identify $\boldsymbol{\ell} = (\ell_i)_{i=1}^N$ with $m$ signatures --- each signature corresponding to one of the trapezoids $\mathcal T_q$ with $q \in [m]$. The $q$-th signature is a $W_q$-tuple of integers
\[
\boldsymbol{\lambda}^{(q)} = (\lambda_1^{(q)}\geq \lambda_2^{(q)}\geq \cdots\geq \lambda_{W_q}^{(q)}),
\]
where $W_q$ is the width of the $q$-th trapezoid and $\lambda_j^{(q)}+W_q-j$ for $j \in [W_q]$ represent the positions of the $W_q$ horizontal lozenges along the long base of $\mathcal{T}_q$, as prescribed by the configuration $\boldsymbol{\ell}$ and the scarlet lozenges.

For example, let us revisit Figure~\ref{Fig_polygon_split}. The tiling in the right panel has five particles, corresponding to the five gray lozenges. Counted from the bottom of the vertical section of the domain, their coordinates are $1<4<6<8<12$. Additionally, there are three scarlet lozenges at positions $3 < 10 < 11$. The left trapezoid has width $W_1=8$ and its long base contains horizontal lozenges at the eight coordinates $1<3<4<6<8<10<11<12$, thus encoded in the signature $\boldsymbol{\lambda}^{(1)}=(5\geq 5 \geq 5 \geq 4\geq 3\geq 2\geq 2\geq 1)$. The top-right trapezoid has width $W_2=2$ and its long base contains the horizontal lozenges at two positions $8<12$, encoded in the signature $\boldsymbol{\lambda}^{(2)}=(11\geq 8)$. The bottom-right trapezoid has width $W_3=5$ and contains the horizontal lozenges at positions $(0<1<3<4<6)$, encoded in the signature $\boldsymbol{\lambda}^{(3)}=(2\geq 1\geq 1 \geq 0 \geq 0)$.

\medskip

Following \cite{BuGo3}, we encode each $\boldsymbol{\lambda}^{(q)}$ through its \emph{Schur generating function}. We recall that for $\boldsymbol{\lambda}=(\lambda_1\geq \lambda_2\geq\ldots\geq\lambda_W)$, the Schur function $s_{\boldsymbol{\lambda}}(\zeta_1,\ldots,\zeta_n)$ is a symmetric function defined by
\[
s_{\boldsymbol{\lambda}}(\zeta_1,\ldots,\zeta_W)=\frac{\det\big([ \zeta_i^{\lambda_j+W-j}\bigr]_{i,j=1}^n}{\prod_{1\leq i<j\leq W}(\zeta_i-\zeta_j)}.
\]
For each $q \in [m]$, we introduce the Schur generating function
\[
S^{(q)}(\zeta_1,\ldots,\zeta_{W_q})=\E\left[\frac{s_{\boldsymbol{\lambda}^{(q)}}(\zeta_1,\ldots,\zeta_{W_q})}{s_{\boldsymbol{\lambda}^{(q)}}(1,\ldots,1)}\right],
\]
where the expectation value is taken in the random lozenge tiling model. Under the conditions of Theorem~\ref{Theorem_GFF}, the particles $\boldsymbol{\ell}$ satisfy the law of large numbers and the central limit theorem as $\N \rightarrow \infty$ by Theorem~\ref{Theorem_main_LLN} and Corollary~\ref{Corollary_CLT_relaxed}. Hence, the same is true for the signature $\boldsymbol{\lambda}^{(q)}$ for each $q\in[m]$. This probabilistic information is recast in \cite[Theorem 2.7]{BuGo3} in terms of Schur generating functions, showing the existence of the following limits for any $k,k_1,k_2 \geq 1$:
\begin{equation}
\label{eq_CLT_appropriate_1}\begin{split}
c_k & := \lim_{\N \rightarrow \infty} W_q^{-1} \partial_{\zeta_i}^{k} \log S^{(q)} \big|_{\zeta_1 = \cdots = \zeta_{W_q} = 1}, \\
d_{k_1,k_2} & := \lim_{\N \rightarrow \infty} \partial_{\zeta_{i_1}}^{k_1}\partial_{\zeta_{i_2}}^{k_2} \log S^{(q)} \big|_{\zeta_1 = \cdots = \zeta_{W_q} = 1}\qquad \textnormal{for}\,\,i_1 \neq i_2.
\end{split}
\end{equation}
Besides, for any $i_1,\ldots,i_n \in [W_q]$ with at least three distinct elements
\begin{equation}
\label{eq_CLT_appropriate_3}
0 = \lim_{\N \rightarrow \infty} \partial_{\zeta_{i_1}}\cdots \partial_{\zeta_{i_n}} \log S^{(q)} \big|_{\zeta_1 = \cdots = \zeta_{W_q} = 1}.
\end{equation}
Random signatures satisfying these three properties are called \emph{CLT-appropriate}. In turn \cite[Theorem 3.13]{BuGo3} (see also \cite[Section 3.2]{BuGo}, \cite[Theorem 2.9 and Section 3.5]{BuGo2}) transforms \eqref{eq_CLT_appropriate_1}-\eqref{eq_CLT_appropriate_3} into the law of large numbers and central limit theorem for the (random) horizontal lozenges in the whole trapezoid $\mathcal{T}_q$, in particular reaching Theorem~\ref{Theorem_tilings_LLN}. For this section we only need the law of large numbers, so we will recast only this result.

 \begin{figure}[t]
 \begin{center}
\includegraphics[width=0.8\linewidth]{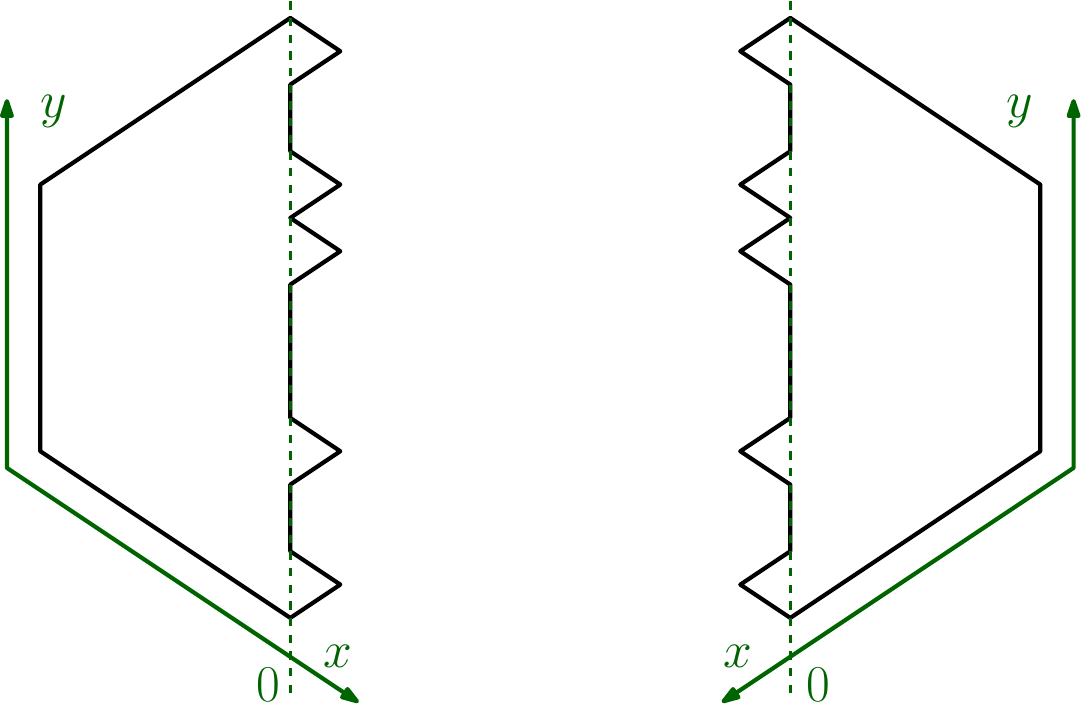}
\end{center}
\caption{Coordinate systems for left and right trapezoids.\label{Fig_left_right_trapezoids}.}
\end{figure}

In the rest of this section we focus on a single trapezoid $\mathcal{T}:= \mathcal{T}_q$ and momentarily remove $q$ from the notations. It has width $W$ and (random) dents covered by horizontal lozenges at positions
\[
l_i =\lambda_i+W-i,\qquad i \in [W],
\]
along the long base. We use the coordinate system of the left panel of Figure~\ref{Fig_left_right_trapezoids}, rescaled by $\N$ and with vertical axis $x=0$ along the long base. We also denote $\hat{W} =\N^{-1}W$. Hence, the horizontal coordinate varies from $-\hat{W}$ to $0$ inside the trapezoid, while the vertical coordinate is constant along the bottom border of the trapezoid. Let us define (we continue using the simplifying convention of Remark~\ref{Remark_L_independence}):
\begin{equation}
\label{eq_x296}
 G(z)=\lim_{\N \rightarrow \infty} \frac{1}{W} \sum_{i=1}^{W}\frac{1}{z-\frac{l_i}{W}}.
\end{equation}
As a corollary of Theorem~\ref{Theorem_main_LLN}, the function $G(z)$ is deterministic. Further, let
\[
\boldsymbol{l}_j = (l_{j,1} >l_{j,2} > \ldots>l_{j,j})
\]
denote the coordinates of the $j$ (random) horizontal lozenges on the vertical line at distance $j$ from the left boundary of the trapezoid (for $j=W$ this is simply $\boldsymbol{l}$). Then \cite[Theorem 3.13]{BuGo3} implies that the empirical measures
\[
\frac{1}{W}\sum_{i=1}^j \delta_{l_{j,i}/W}
\]
converge to a deterministic limit, which smoothly depends on the asymptotic parameter $\frac{j}{W}$. In turn, this implies Theorem~\ref{Theorem_tilings_LLN} for the values of the height function inside the trapezoid $\mathcal{T}$. The limit shape and liquid region inside $\mathcal{T}$ can be described in terms of $G(z)$. For horizontal coordinate $-\hat{W} \leq x \leq 0$ and vertical coordinate $y$, as in Figure~\ref{Figure_triangle}, consider the equation
\begin{equation}
\label{eq_critical_equation_rescaled}
 y= z + \frac{-x}{e^{-G(z/\hat{W})}-1}
 \end{equation}
for the unknown $z \in \amsmathbb{C}$.
Let $\pounds \subset \mathcal{T}$ denote the set of points $(x,y)$ in
$(-\hat{W},0)\times \amsmathbb R$, such that the equation has a solution $z$ in the upper
half-plane. It is proven in\footnote{Formally, these articles only analyze the case of deterministic $\boldsymbol{l}$. However, \cite[Theorem 3.13]{BuGo3} shows that for random $\boldsymbol{l}$ the limit shape is the same.}
\cite{petrov2015asymptotics,duse2015asymptotic,BuGo2,gorin2017bulk} that $\pounds$ is precisely the liquid region of the trapezoid $\mathcal{T}$ (compared too those article we have rescaled
$x$, $y$ and $z$ by $\hat{W}$). Moreover, for each
$(x,y)\in \pounds$ the equation \eqref{eq_critical_equation_rescaled} in fact has a
unique solution $z(x,y)$ in the upper-half plane $\amsmathbb{H}$ and the map $z :
\pounds \rightarrow \amsmathbb{H}$ is a smooth bijection, see \cite{duse2015asymptotic}.

Differentiating \eqref{eq_critical_equation_rescaled} in $x$ and $y$, we get
\begin{equation}
\label{eq_critical_equation_rescaled_xy}
\begin{split} 0 & = \left(1+\frac{G'(z/\hat{W})}{\hat{W}} e^{-G(z/\hat{W})} \cdot \frac{-x}{(e^{-G(z/\hat{W})} - 1)^2}\right) \partial_x z - \frac{1}{e^{-G(z/\hat{W})} - 1}, \\
1 & = \left(1+\frac{G'(z/\hat{W})}{\hat{W}} e^{-G(z/\hat{W})} \cdot \frac{-x}{(e^{-G(z/\hat{W})} - 1)^2}\right) \partial_y z,
\end{split}
\end{equation}
 which implies the relation
\begin{equation}
\label{eq_burgers_for_map}
\partial_yz = \big(e^{-G(z/\hat{W})} - 1\big) \partial_x z.
\end{equation}
In addition, the same articles \cite{petrov2015asymptotics,duse2015asymptotic,BuGo2,gorin2017bulk} explain that for left trapezoids $z(x,y)$ is related to the complex slope $\xi=\xi(x,y)$ of Figure~\ref{Figure_triangle} through
\begin{equation}
\label{eq_complex_structure_def}
\forall x \in [-\hat{W},0]\qquad \xi(x,y)= \frac{-x}{y-z(x,y)}= e^{-G(z/\hat{W})} - 1.
\end{equation}
\begin{remark}
\label{Remark_complex_slope_amb}
The articles
\cite{petrov2015asymptotics,duse2015asymptotic,BuGo2,gorin2017bulk} where the relations between $z$ and $\xi$
are established all use slightly different coordinate systems and definitions of
$\xi$. This is the ambiguity we discussed in Section~\ref{sec:complexslopel}. For instance, in
\cite{petrov2015asymptotics,petrov2015asymptotics}, the coordinate system $(\eta,\chi)$ is used and we have explained that the complex slope
$\Omega$ used there is related to our $\xi$ by $\Omega= \frac{1}{1 - \xi}$. Along the same lines, if $\mathcal{T}$ is a right, rather than left trapezoid, and we use the coordinate system of the right panel of Figure~\ref{Fig_left_right_trapezoids}, then the roles of the angles at $0$ and at $\xi$ in the definition of the complex slope by Figure~\ref{Figure_triangle} should be swapped, corresponding to the last transformation in \eqref{eq_alternative_form}, $\tilde{\xi} =\frac{\xi^*}{\xi^*-1}$. This is because the transformation between the left and right panels of Figure~\ref{Fig_left_right_trapezoids} is the symmetry with respect to $x=0$ axis, which swaps the roles of two lozenges {\scalebox{0.16}{\includegraphics{lozenge_v_up.pdf}}} and
{\scalebox{0.16}{\includegraphics{lozenge_v_down.pdf}}}.
\end{remark}

Combining the formulae \eqref{eq_critical_equation_rescaled} and \eqref{eq_complex_structure_def} gives the desired description of the limit shape inside the liquid region $\pounds$. The rest of the trapezoid $\mathcal{T}$ is a frozen region and one can figure out which of the three types of lozenges is observed in each part of it (\textit{cf.} Figure~\ref{Fig_hex_simulation}) by continuity of the densities extended from the liquid region. The combination of \eqref{eq_burgers_for_map} and \eqref{eq_complex_structure_def} is the complex Burgers equation for the complex slope $\xi$.

\subsection{Complex structure for gluings of trapezoids}

\label{Section_complex_structure}

The complex structure relevant for the Gaussian free field in Conjecture~\ref{Conjecture_GFF} is given in an abstract way, yet it becomes more explicit for gluings of trapezoids, as we explain in this section. In the setting of Theorem~\ref{Theorem_GFF_general}, we can still rely on Lemma~\ref{Lemma_glued_tilings_assumptions}: when we fix the filling fractions, the proofs in Lemma~\ref{Lemma_glued_tilings_assumptions} remain valid without any changes. Hence, the equilibrium measure has precisely one band $(\alpha_h,\beta_h)$ in the $h$-th blue segment for each $h \in [H]$. We use these bands to construct a bordered Riemann surface $\Sigma^{\textnormal{half}}$ which is\footnote{``Half'' in $\Sigma^{\textnormal{half}}$ indicates that is obtained by gluing half-planes.} a gluing of $m$ half-planes $\Sigma^{(q)}$ with a standard complex structure.

\begin{itemize}
 \item If $\mathcal T_q$ is a left trapezoid, then $\Sigma^{(q)}$ is the union of the open \emph{upper} half-plane $\{z\in\amsmathbb{C}\mid \mathrm{Im}(z)>0\}$ with all the open segments $(\alpha_h,\beta_h)$ over all bands inside the blue segments of $\mathcal{T}_q$.
 \item If $\mathcal{T}_q$ is a right trapezoid, then $\Sigma^{(q)}$ is the union of the open \emph{lower} half-plane $\{z\in\amsmathbb{C}\mid \mathrm{Im}(z)<0\}$ with open segments $(\alpha_h,\beta_h)$ over all bands inside the blue segments of $\mathcal{T}_q$.
\end{itemize}

\begin{figure}[t]
 \begin{center}
\includegraphics[width=0.8\linewidth]{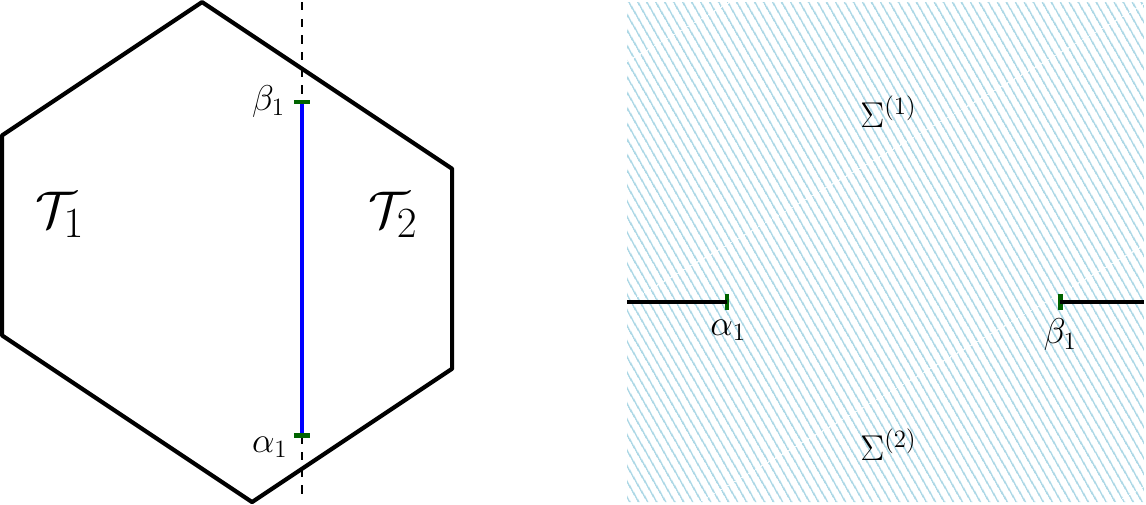}

\medskip

\includegraphics[width=0.8\linewidth]{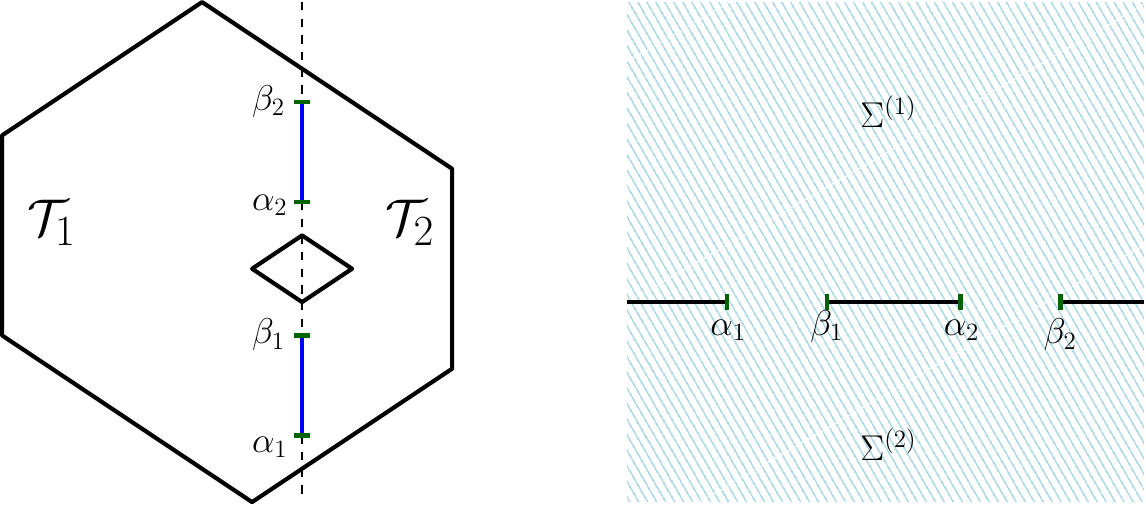}
\end{center}
\caption{Domains $\Sigma^{\textnormal{half}} =\Sigma^{(1)} \cup \Sigma^{(2)}$ for hexagon and hexagon with a hole. \label{Fig_uniformization_domains}}
\end{figure}

The gluings of $\Sigma^{(q)}$ into $\Sigma^{\textnormal{half}}$ happen for adjacent trapezoids: if two distinct trapezoids $\mathcal{T}_q$ and $\mathcal{T}_{q'}$ share a common blue segment $\llbracket A_{h},B_{h}\rrbracket$, then $(\alpha_h,\beta_h)$ in $\Sigma^{(q)}$ is identified with $(\alpha_h,\beta_h)$ in $\Sigma^{(q')}$. Figure~\ref{Fig_uniformization_domains} shows which subsets of $\amsmathbb{C}$ we get in the examples with the gluings of two trapezoids featured in Sections~\ref{Section_Hexagon} and \ref{Section_hex_hole}. When $m>2$, the surface $\Sigma^{\textnormal{half}}$ constructed as a gluing of $\Sigma^{(q)}$s would not be embedded into a complex plane directly. Yet, if the tiled domain $\mathcal{D}$ of Definition~\ref{Definition_gluing} is simply-connected, there would still exist a conformal bijection between $\Sigma^{\textnormal{half}}$ and a subset of $\amsmathbb{C}$ by the uniformization theorem.

In the general bipartite case, let us construct a uniformization map identifying the liquid region $\pounds$ for tilings of $\mathcal{D}$ with the bordered Riemann surface $\Sigma^{\textnormal{half}}$. For this we use Section~\ref{Section_limit_shape_gluings} with indices $q$ of trapezoids restored.

\begin{definition} \label{def:zqzq} For each left trapezoid $\mathcal{T}_q$, we define the map $z^{(q)} : \pounds_q \rightarrow \Sigma^{(q)}$ as the unique solution in the \emph{upper} half-plane of \eqref{eq_critical_equation_rescaled}, that is
\begin{equation}
\label{eq_critical_equation_rescaled2} \frac{-x}{y - z^{(q)}} = e^{-G^{(q)}(z^{(q)})} - 1,
\end{equation}
where $G^{(q)}$ is the function \eqref{eq_x296} for the $q$-th trapezoid. We denote this solution $z^{(q)}$ and treat it as a point in $\Sigma^{(q)}$. For each right trapezoid $\mathcal{T}_q$, we map its liquid region $\pounds_q$ into a unique solution of \eqref{eq_critical_equation_rescaled2} in the \emph{lower} half-plane. We again denote this solution $z^{(q)}$ and treat it as a point in $\Sigma^{(q)}$.
\end{definition}
\begin{proposition} \label{Proposition_uniformization_map}
 The maps $z^{(q)} :\pounds_q \rightarrow \Sigma^{(q)}$ extend by continuity to a bijection $z : \pounds \rightarrow \Sigma^{\textnormal{half}}$. This bijection is conformal with respect to the complex structure associated to the complex slope on $\pounds$, and the standard complex structure inherited from $\amsmathbb{C}$ on $\Sigma^{\textnormal{half}}$.
\end{proposition}
\begin{proof}
 The bijectivity of the map $z^{(q)}$ for trapezoids was proven in \cite{duse2015asymptotic}. For the left trapezoids, the formulae \eqref{eq_burgers_for_map}-\eqref{eq_complex_structure_def} imply that the map is conformal for the complex structure coming from the complex slope (\textit{cf.} Section~\ref{sec:complexstr} and Section~\ref{sec:complexslopel}). For the right trapezoids, as in Remark~\ref{Remark_complex_slope_amb} we need to replace the complex slope by $\frac{\xi^*}{\xi^*-1}$ and also conjugate $z^{(q)}$ since we now map to the lower rather than upper half-plane. Hence, the map is again conformal.

 It remains to show that the same properties continue to hold when we glue the maps $z^{(q)}$ together into a single map $z : \pounds \rightarrow \Sigma^{\textnormal{half}}$. For that we need to analyze the behavior of the map as a point in $\pounds_q\subset \mathcal{T}_q$ approaches the long base of the trapezoid. As we come from the liquid region, it is only possible to approach the parts of the long base which belong to the bands of the equilibrium measure (of the discrete ensemble describing the particles $\boldsymbol{\ell}$ on the vertical line of gluing), the rest of the long base being part $\partial\Sigma^{\textnormal{half}}$. In the equation \eqref{eq_critical_equation_rescaled}, the point approaching the long base corresponds to the $x$-coordinate approaching $0$. When the $y$-coordinate corresponds to the point inside the band, $G^{(q)}(y/\hat{W}_q)$ is non-real, which implies non-vanishing of the denominator in \eqref{eq_critical_equation_rescaled} and the solution of this equation at $x=0$ becomes real: it is $z=y$.

 We conclude that for any two adjacent trapezoids $\mathcal{T}_q$ and $\mathcal T_{q'}$, the corresponding maps $z^{(q)}$ and $z^{(q')}$ have the same real boundary values on the bands of the gluing vertical line, and therefore extend continuously across the line. Since any holomorphic function extends in a holomorphic over a cut where it is continuous (\textit{e.g.} by Morera theorem), this implies the holomorphicity of the glued map $z : \pounds \rightarrow \Sigma^{\textnormal{half}}$.
 \end{proof}

\subsection{Gaussian free field: proof of Theorem~\ref{Theorem_GFF_general}}
\label{Section_GFF_proof}

We first show that the field of fluctuations of the height function is asymptotically Gaussian. Then we identify the covariance with the one of the Gaussian free field.

\medskip

\noindent \textsc{Step 1.} As in Section~\ref{Section_limit_shape_gluings}, we would like to combine Theorem~\ref{Theorem_main_LLN} and Corollary~\ref{Corollary_CLT_relaxed} with the results of \cite{BuGo3}. While in Section~\ref{Section_limit_shape_gluings}, for each trapezoid the argument was completely independent and relied on its own Schur generating function, this time we also need to study correlations between height functions of different trapezoids. Hence, we need a generalization of the Schur generating function which captures information about all $m$ trapezoids simultaneously. We use the same notation as in Section~\ref{Section_limit_shape_gluings} and package the information about the random $m$ signatures $\boldsymbol{\lambda}^{(1)}, \ldots, \boldsymbol{\lambda}^{(m)}$ through their $m$-dimensional Schur generating function:
\[
S\big(\zeta_1^{(1)},\ldots,\zeta_{W_1}^{(1)};\zeta_1^{(2)},\ldots,\zeta_{W_2}^{(2)};\ldots; \zeta_1^{(m)},\ldots,\zeta_{W_m}^{(m)}\big)=\E \left[\prod_{q=1}^m\frac{s_{\boldsymbol{\lambda}^{(q)}}(\zeta^{(q)}_1,\ldots, \zeta^{(q)}_{W_q})}{s_{\boldsymbol{\lambda}^{(q)}}(1,\ldots,1)}\right].
\]
Under the assumptions of Theorem~\ref{Theorem_GFF}, the particle configuration $\boldsymbol{\ell}$ satisfies the law of large numbers and the central limit theorem as $\N \rightarrow \infty$ by Theorem~\ref{Theorem_main_LLN} and Corollary~\ref{Corollary_CLT_relaxed}. Hence, the same is true for joint distribution of all the signatures $(\boldsymbol{\lambda}^{(q)})_{q = 1}^{m}$. This probabilistic information is then recast in \cite[Theorem 2.7]{BuGo3} in terms of the $m$-dimensional Schur generating functions. The conclusion is that the following limits exist for any $q,q_1,q_2 \in [m]$ and $k,k_1,k_2 \geq 1$
\begin{equation}
\label{eq_CLT_ap_multi_1}
\begin{split}
c_{k}^{(q)} & = \lim_{\N \rightarrow \infty} W_q^{-1} \partial_{\zeta_i^{(q)}}^k \log S \big|_{\zeta_1^{(1)} = \cdots = \zeta_{W_m}^{(m)} = 1}, \\
d_{k_1,k_2}^{(q_1,q_2)} & = \lim_{\N \rightarrow \infty} \partial_{\zeta_{i_1}^{(q_1)}}^{k_1} \partial_{\zeta_{i_2}^{(q_2)}}^{k_2} \log S \big|_{\zeta_1^{(1)} = \cdots = \zeta_{W_m}^{(m)} = 1} \qquad \textnormal{for}\,\,(i_1,q_1) \neq (i_2,q_2).
\end{split}
\end{equation}
Besides, for any tuples $(i_1,q_1),\ldots,(i_n,q_n)$ with at least three distinct elements
\begin{equation}
\label{eq_CLT_ap_multi_3}
0 = \lim_{\N \rightarrow \infty} \partial_{\zeta_{i_1}^{(q_1)}} \cdots \partial_{\zeta_{i_n}^{(q_n)}} \log S \big|_{\zeta_1^{(1)} = \cdots = \zeta_{W_m}^{(m)} = 1}.
\end{equation}
We can extract from the $m$-dimensional Schur generating function $S$ information about random variables $p_{k}^{(q)}(\nu)$ which give the rescaled $k$-th power sum of positions of horizontal lozenges along the vertical line at distance $(1-\nu)W_q$ from the long base of the trapezoid $\mathcal{T}_q$. Formally, we denote
\[
\boldsymbol{l}_{j}^{(q)} = (l_{j,1}^{(q)} > l_{j,2}^{(q)} > \cdots > l_{j,j}^{(q)})
\]
the coordinates of (random) horizontal lozenges on the vertical line at distance $j$ from the short base of $\mathcal{T}_q$. For $j=W_q$ we retrieve $l_{W_q,i}^{(q)} = \lambda_i^{(q)}+W_q-i$ corresponding to the signature $\boldsymbol{\lambda}^{(q)}$. We define
\[
p_{k}^{(q)}(\nu) = \sum_{i=1}^{\lfloor \nu W_q\rfloor} \Bigg(\frac{l_{\lfloor \nu W_q \rfloor,i}^{(q)}}{\lfloor \nu W_q\rfloor}\Bigg)^k.
\]
This random variable can also be expressed as an integral of the height function $\mathsf{Ht}^{(q)}_\N(x,y)$, which counts the horizontal lozenges above $(\N x,\N y)$ inside the trapezoid $\mathcal{T}_q$, using the coordinate system of Figure~\ref{Fig_left_right_trapezoids}. Choosing an arbitrary large enough $M$ so that the entire trapezoid (rescaled by $\N$) is inside $[-\hat{W}_q,0]\times[-M,M]$ and integrating by parts, we have:
\begin{equation}
\label{eq_x295}
\begin{split}
 p_{k}^{(q)}(\nu) \cdot \frac{\lfloor \nu W_q\rfloor^{k}}{\N^k} & = \sum_{i=1}^{\lfloor \nu W_q\rfloor} \Bigg(\frac{l_{\lfloor \nu W_q\rfloor,i}^{(q)}}{\N}\Bigg)^k \\
 & = - k \int\limits_{-M}^M \left( \sum_{i=1}^{\lfloor \nu W_q\rfloor} \mathbbm{1}_{\{l_{\lfloor \nu W_q \rfloor,i}^{(q)} \geq \N y\}} \right) y^{k-1} \dd y+ \lfloor \nu W_q\rfloor (-M)^k \\
 & =  - k \int\limits_{-M}^M \mathsf{Ht}_{\N}^{(q)}\bigg(\frac{\lfloor \nu W_q\rfloor-W_q}{\N}\,,\,y\bigg) y^{k-1} \dd y+ \lfloor \nu W_q\rfloor (-M)^k.
\end{split}
\end{equation}
Hence, up to a constant shift, $p_{k}^{(q)}(\nu)$ is a pairing of the height function $\mathsf{Ht}^{(q)}_\N(x,y)$ with a specific function, which is polynomial in the $y$-direction and a $\delta$-function in the $x$-direction. Hence, tested against such functions, the asymptotic Gaussianity of the field $\mathsf{Ht}_\N^{(q)}(x,y)-\E[\mathsf{Ht}_\N^{(q)}(x,y)]$ is equivalent to the asymptotic (joint) Gaussianity of $p_{k}^{(q)}(\nu)$.

This asymptotic Gaussianity follows from conditions \eqref{eq_CLT_ap_multi_1}-\eqref{eq_CLT_ap_multi_3} by \cite[Theorem 3.13]{BuGo3}, in which we take all functions $g_{h,m}$ equal to $1$. The asymptotic covariance of $p_{k}^{(q)}(\nu)$ is also delivered by that theorem:
\begin{equation}
\label{eq_x294}
\begin{split}
& \quad \lim\limits_{\N \rightarrow \infty} \E^{(\textnormal{c})}\big[p_{k_1}^{(q_1)}(\nu_1),p_{k_2}^{(q_2)}(\nu_2)\big] \\
& = [\zeta_1^{-1} \zeta_2^{-1}] \prod_{i = 1}^{2}
 \left( (1+\zeta_i)
\bigg(\frac{1}{\zeta_i} + \frac{1}{\nu_i} \sum_{n_i=1}^\infty \frac{c_{n_i}^{(q_i)}\zeta_i^{n_i-1}
}{(n_i-1)!}\bigg)\right)^{k_i} \\
& \quad \phantom{ [\zeta_1^{-1} \zeta_2^{-1}]} \times \left( \delta_{q_1,q_2}\left(\sum_{n=0}^{\infty}
\frac{\zeta_1^n}{\zeta_2^{n+1}}\right)^2 +\sum_{n_1,n_2 =1}^{\infty}
\frac{d_{n_1,n_2}^{(q_1,q_2)}\zeta_1^{n_1-1} \zeta_2^{n_2-1} }{(n_1-1)! (n_2-1)!} \right),
\end{split}
\end{equation}
where $[\zeta_1^{-1} \zeta_2^{-1}]$ stands for the coefficient
$\zeta_1^{-1}\zeta_2^{-1}$ in the series to its right. The factors in the second line of \eqref{eq_x294} are related to the leading term deterministic asymptotics of $p_{k}^{(q)}(\nu)$ by the same \cite[Theorem 3.13]{BuGo3}, that is
\begin{equation}
\label{eq_x297}
 \lim_{\N \rightarrow \infty} p_{k}^{(q)}(\nu) = [\zeta^{-1}] \frac{1}{(k+1)(\zeta+1)}
\left( (1+\zeta) \left(\frac{1}{\zeta} + \frac{1}{\nu}\sum_{n=1}^\infty \frac{c_{n}^{(q)}\zeta^{n-1}
}{(n-1)!}\right) \right)^{k+1},
\end{equation}
where $[\zeta^{-1}]$ is the coefficient of $\zeta^{-1}$ in the series to its right. In the rest of the proof we will identify \eqref{eq_x294} with the covariance of pairings of the Gaussian free field against appropriate test measures.

\medskip

\noindent \textsc{Step 2.} Our next task is to write down formulae for the covariance of the Gaussian free field which is the limiting object in Theorem~\ref{Theorem_GFF_general}, solely in terms of the asymptotic covariance along the gluing line for our domain. Recall that the equilibrium measure of the particles $\ell_1<\ell_2<\cdots<\ell_N$ on the gluing line has $H$ bands according to the assumptions of Lemma~\ref{Lemma_glued_tilings_assumptions}, that is $(\alpha_h,\beta_h) \subseteq [\hat a_h,\hat b_h]$ for each $h \in [H]$. Recall the notation $\mathcal{H}_q \subseteq [H]$ for the set of blue segments in the trapezoid $\mathcal{T}_q$. The asymptotic covariance of the Stieltjes transforms for the $h_1$- and $h_2$-th segments is obtained as $\mathcal{F}_{h_1,h_2}(z_1,z_2)$ in Corollary~\ref{Corollary_CLT_relaxed} and \eqref{eq_covariancepre}. Taking into account all bands appearing in two given trapezoids, the asymptotic covariance of the Stieltjes transform for the horizontal lozenges in those trapezoids is given by the following function.
\begin{definition} \label{Definition_Covariance_tilings_St}
 For $q_1,q_2 \in [m]$ and $z_1,z_2 \in \amsmathbb{C}$ outside the segments, we define
\[
\mathcal{F}^{(q_1,q_2)}(z_1,z_2)= \sum_{\substack{h_1 \in \mathcal{H}_{q_1} \\ h_2 \in \mathcal{H}_{q_2}}} \mathcal{F}_{h_1,h_2}(z_1,z_2).
\]
\end{definition}
Clearly, $\mathcal{F}^{(q_1,q_2)}(z_1,z_2)$ is a holomorphic function of $z_1$ outside the bands $(\alpha_{h_1},\beta_{h_1})$ of $\mathcal T_q$ and a holomorphic function of $z_2$ outside the bands $(\alpha_{h_2},\beta_{h_2})$ of $\mathcal{T}_{q_2}$. We find the following relation to the Green function.

\begin{proposition} \label{Proposition_Green_function_explicit}
 The Green function $\textnormal{Gr}$ on $\Sigma^{\textnormal{half}}$ with Dirichlet boundary conditions has the following expression\footnote{One should compare this with (4.21), (4.22) in \cite{BuGo3} and note a sign error in that paper --- one needs to multiply $\frac{1}{4\pi}$ by $(-1)$ in both formulae. In the proofs of \cite{BuGo3} this sign error compensates with wrong signs in (4.25), (4.26) and in the computation of mixed partial derivative above (4.58).} for non-real $z_i\in \Sigma^{(q_i)}$ ($i = 1,2$):
 \begin{equation} \label{eq_Green_function_explicit}
 \textnormal{Gr}(z_1,z_2)= \left\{\begin{array}{lll} -\dfrac{1}{4\pi} \displaystyle\int_{z_1^*}^{z_1}\displaystyle\int_{z_2^*}^{z_2} \mathcal{F}^{(q,q)}(w_1,w_2) \dd w_1 \dd w_2- \dfrac{1}{2\pi} \log\left|\frac{z_1-z_2}{z_1-z_2^*}\right|& & \textnormal{ if }\, q_1 = q_2 = q,\\[12pt] - (-1)^{\textnormal{par}(q_1,q_2)}\dfrac{1}{4\pi} \displaystyle\int_{z_1^*}^{z_1}\displaystyle\int_{z_2^*}^{z_2} \mathcal{F}^{(q_1,q_2)}(w_1,w_2) \dd w_1 \dd w_2 && \textnormal{ if }\,\, q_1\neq q_2,\end{array}\right.
 \end{equation}
 where the parity $(-1)^{\textnormal{par}(q_1,q_2)}$ is $(+1)$ if $\mathcal T_{q_1},\mathcal T_{q_2}$ are both left, or both right trapezoids, and $(-1)$ otherwise. In particular, \eqref{eq_Green_function_explicit} extends to a continuous function of $\Sigma^{\textnormal{half}} \times \Sigma^{\textnormal{half}}$ minus the diagonal.
\end{proposition}
The integration contours in \eqref{eq_Green_function_explicit} can be arbitrary as long as they do not cross the bands where $\mathcal{F}^{(q_1,q_2)}$ is not well-defined. Since integrals of $\mathcal{F}_{h_1,h_2}(z_1,z_2)$ around the bands vanish by its definition in \eqref{eq_covariancepre}, the choice of the contours does not affect the value of the integral.

\begin{proof}[Proof of Proposition~\ref{Proposition_Green_function_explicit}] The Green function in $\Sigma^{\textnormal{half}}$ is characterized by three properties:
\begin{enumerate}
\item It should vanish when $z_i$ approaches the boundary of $\Sigma^{\textnormal{half}}$, which includes $\infty$ and the real axis minus the bands in each $\Sigma^{(q_i)}$.
\item It should be a real harmonic function outside the diagonal, \textit{i.e.} the locus $z_1 = z_2$ in the same $\Sigma^{(q)}$ for $q \in [m]$.
\item Near the diagonal, it can be locally represented as the sum of $-\frac{1}{2\pi} \log|z_1-z_2|$ and a bounded harmonic function.
\end{enumerate}

Let us call $\textnormal{gr}(z_1,z_2)$ the right-hand side of \eqref{eq_Green_function_explicit}. The first property is clear from the definition: when $z_1 \in\Sigma^{(q_1)}$ approaches the real axis outside the bands of $\mathcal{T}_{q_1}$, all terms in \eqref{eq_Green_function_explicit} vanish. The same is true as $z_1 \rightarrow \infty$ because $\mathcal{F}^{(q_1,q_2)}(z_1,z_2)$ decays as $O(\frac{1}{z_1^2})$, coming from the fact that the filling fractions are deterministically fixed (in \eqref{eq_covariancepre} it corresponds to $\boldsymbol{\kappa}=\boldsymbol{0}$). A similar vanishing holds as $z_2$ approaches the boundary of $\Sigma^{\textnormal{half}}$. The expression in \eqref{eq_Green_function_explicit} is real: for the integral part the preservation under conjugations can be seen by choosing the integration contours to be symmetric with respect to the real axis and using $(\mathcal{F}^{(q_1,q_2)}(z_1,z_2))^* = \mathcal{F}^{(q_1,q_2)}( z_1^*,z_2^*)$ from its definition. For non-real $z_1,z_2$, the double integrals in \eqref{eq_Green_function_explicit} define harmonic functions, as linear combinations of holomorphic and anti-holomorphic function. The desired singularity at $z_1=z_2$ is also clearly visible in $\textnormal{gr}(z_1,z_2)$ in the second term when $q_1=q_2$. Hence, it satisfies the second and third defining properties of the Green function for non-real $z_1,z_2$.

It remains to study the behavior of $\mathrm{gr}$ as one of the points approaches the bands connecting adjacent $\Sigma^{(q)}$s: we need to show that $\mathrm{gr}(z_1,z_2)$ extends continuously to a sum of a harmonic function and a singular term $-\frac{1}{2\pi} \log|z_1-z_2|$. For the simplest cases, such as those of Figure~\ref{Fig_uniformization_domains}, the check can be done directly using the explicit formulae for $\mathcal{F}^{(q_1,q_2)}(z_1,z_2)$, see \cite[Propositions 4.13, 4.14]{BuGo3}. For general glued domains, this check is done by identifying the jumps of $\mathcal{F}^{(q_1,q_2)}(z_1,z_2)$ on the real axis, \textit{i.e.} by representing it as a solution of a Riemann--Hilbert problem. This analysis is done in greater generality (for any discrete ensemble, not necessarily related to tilings) and in details in Section~\ref{Section_Master_by_covariance}, see in particular Proposition~\ref{thmBfund}. Double integrals of fundamental solutions like in \eqref{eq_Green_function_explicit} are introduced and studied systematically in Section~\ref{sec:Green13}, see in particular Proposition~\ref{prop:Greengen} (to compare with Definition~\ref{def:Hvv}, note that taking the real part of \eqref{eq_Green_function_explicit} does not change the result). Corollary~\ref{Greenwon} gives a criterion for relating it to the Green function. For tilings in the orientable glued domains that we are focusing on, the relevant results are Theorem~\ref{thm:Omegav} and Corollary~\ref{cor:Greennormal}, where $\textnormal{gr}$ corresponds to $\textnormal{gr}_{\boldsymbol{v}}$: they show continuity of the extension to $\Sigma^{\textnormal{half}}$ and finish the proof that $\textnormal{gr}$ coincides with the Green function $\textnormal{Gr}$.
\end{proof}

Note that the $- \frac{1}{2\pi} \log\big|\frac{z_1-z_2}{z_1-z_2^*}\big|$ term in the formula for $\mathrm{Gr}(z,w)$ is the Green function of the upper (or lower) half-plane. Hence, the decomposition of \eqref{eq_Green_function_explicit} can be viewed as a decomposition of the Gaussian free field in $\Sigma^{\textnormal{half}}$ into a sum of $m+1$ independent components: $m$ independent Gaussian free fields in each half-plane $\Sigma^{(q)}$ for $q\in[m]$, and the one-dimensional random field on $\bigcup_{h=1}^H (\alpha_h,\beta_h)$
extended as a harmonic function to the interior of each $\Sigma^{(q)}$. This is a shadow of the well-known domain Markov property of the Gaussian free field, see \textit{e.g.} \cite[Section 2.6]{sheffield2007gaussian} or, equivalently, the restriction property of the Green function encoding its covariance. A similar decomposition is already visible in the last line of \eqref{eq_x294} and we will make it more explicit in the next steps of the proof.

\medskip

\noindent \textsc{Step 3.} We massage the formulae of the previous two steps to relate them to each other. Introduce the notation
\[
 F^{(q)}(\zeta):=(1+\zeta) \left(\frac{1}{\zeta} + \sum_{n=1}^\infty \frac{c_n^{(q)} \zeta^{n-1}}{(n-1)!}\right),
\]
and for $x\in [-\hat{W}_q,0]$,
\[
F^{(q)}(\zeta;x):=(1+\zeta) \left(\frac{1}{\zeta} + \frac{\hat{W}_q}{\hat{W}_q+x}\sum_{n=1}^\infty \frac{c_{n}^{(q)}\zeta^{n-1}}{(n - 1)!}\right),
\]
so that $F^{(q)}(\zeta;0)=F^{(q)}(\zeta)$. As a corollary of \eqref{eq_x297} with $\nu=\frac{\hat{W}_q+ x}{\hat{W}_q}$ and recalling the function $G^{(q)}$ appearing in Definition~\ref{def:zqzq}, one can derive the identity
\begin{equation}
\label{eq_x302}
 F^{(q)}(\zeta;x)= \frac{1}{\hat{W}_q+x} \left(Z^{(q)}(\zeta) +
\frac{-x}{e^{-G^{(q)}(Z^{(q)}(\zeta)/\hat{W}_q)}-1} \right),\qquad \textnormal{where}\,\, \zeta = e^{G^{(q)}(Z^{(q)}(\zeta)/\hat{W}_q)} - 1.
\end{equation}
More precisely, we should take for $Z^{(q)}(\zeta)$ in this formula the unique analytic solution of the implicit equation defined near $\zeta = 0$ and such that $Z^{(q)}(\zeta) \sim \frac{1}{\zeta}$ as $\zeta \rightarrow 0$. We refer to \cite[proof of Lemma 9.2]{BuGo2} for more details of this computation. For $x=0$, \eqref{eq_x302} becomes $\hat{W}_q F^{(q)}(\zeta) = Z^{(q)}(\zeta)$.

\begin{lemma} \label{Lemma_Q_case_1}
 Using Definition~\ref{Definition_Covariance_tilings_St}, for $q_1\neq q_2$ the last line of \eqref{eq_x294} can be computed as
\begin{equation}
 \label{eq_x299}
 \begin{split}
 D^{(q_1,q_2)}(\zeta_1,\zeta_2) & :=\sum_{n_1,n_2 = 1}^{\infty} \frac{d_{n_1,n_2}^{(q_1,q_2)} \zeta_1^{n_1 - 1}\zeta_2^{n_2 - 1}}{(n_1-1)! (n_2-1)!} \\
& = \hat{W}_{q_1}\hat{W}_{q_2}\mathcal{F}^{(q_1,q_2)}\big(\hat{W}_{q_1} F^{(q_1)}(\zeta_1)\,,\,\hat{W}_{q_2}F^{(q_2)}(\zeta_2)\big) \cdot \partial_{\zeta_1}F^{(q_1)}(\zeta_1) \partial_{\zeta_2} F^{(q_2)}(\zeta_2).
\end{split}
\end{equation}
\end{lemma}
\begin{proof}
 We take \eqref{eq_x294} for $\nu_1=\nu_2=1$, multiply by $u_1^{-k_1-1}u_2^{-k_2-1}$ for larger complex numbers $u_1,u_2$, sum over all $k_1,k_2 \geq 0$ and express $[\zeta_1^{-1} \zeta_2^{-1}]$ as a contour integral over small contours enclosing $0,\frac{1}{u_1},\frac{1}{u_2}$ to get
\begin{equation}
\label{eq_x298}
\begin{split}
 & \quad \lim_{\N\rightarrow\infty} \E^{(\textnormal{c})}\left[\sum_{i_1=1}^{W_{q_1}} \frac{1}{u_1-l_{i_1}^{(q_1)}/W_{q_1}}\,,\,\sum_{i_2=1}^{W_{q_2}} \frac{1}{u_2-l_{i_2}^{(q_2)}/W_{q_2}} \right] \\
 & = \oint\oint \frac{D^{(q_1,q_2)}(\zeta_1,\zeta_2)}{(u_1-F^{(q_1)}(\zeta_1))(u_2-F^{(q_2)}(\zeta_2))}\,\frac{\dd \zeta_1\dd \zeta_2}{(2\ii\pi)^2},
\end{split}
\end{equation}
where $l_{i}^{(q)} := l_{W_q,i}^{(q)} = \lambda_i^{(q)} + W_q - i$. Comparing with Definition~\ref{Definition_Covariance_tilings_St} and Corollary~\ref{Corollary_CLT_relaxed}, the left-hand side of \eqref{eq_x298}
is $\hat{W}_{q_1}\hat{W}_{q_2}\mathcal{F}^{(q_1,q_2)}(\hat{W}_{q_1}u_1,\hat{W}_{q_2}u_2)$. The integrals in the right-hand side are computed as the residues at the (unique) close to $0$ solutions of the equations $u_1=F^{(q_1)}(\zeta_1)$ and $u_2=F^{(q_2)}(\zeta_2)$. Hence, \eqref{eq_x299} follows from \eqref{eq_x298}.
\end{proof}

If $q_1 = q_2 = q$ the analogue of Lemma~\ref{Lemma_Q_case_1} is a bit trickier. We first introduce another function:
\[
\tilde{D}^{(q)}(\zeta_1,\zeta_2)= \sum_{n_1,n_2 = 1}^{\infty} \frac{\tilde{d}_{n_1,n_2}^{(q)}\zeta_1^{n_1-1}\zeta_2^{n_2-1}}{(n_1-1)! (n_2-1)!}.
\]
In this expression the numbers $\tilde{d}_{n_1,n_2}^{(q_1,q_2)}$ are chosen so that \eqref{eq_x294} with $q_1 = q_2 = q$, $\nu_1=\nu_2=1$ and $d$ replaced with $\tilde{d}$ gives zero; in other words, these are the numbers corresponding to tilings of a trapezoid of width $W_q$ and whose horizontal lozenges along the long base have the same limit shape as the ones for $\mathcal{T}_q$, but no fluctuations. The existence of such numbers $\tilde{d}_{n_1,n_2}^{(q_1,q_2)}$ is guaranteed by \cite[Theorem 2.7]{BuGo3} and $\tilde{D}^{(q)}(\zeta_1,\zeta_2)$ matches $Q_\rho$ from \cite[Section 9.1]{BuGo2}. Further, set
\[
 \Delta D^{(q)}(\zeta_1,\zeta_2)= D^{(q,q)}(\zeta_1,\zeta_2)- \tilde{D}^{(q)}(\zeta_1,\zeta_2)=\sum_{n_1,n_2 = 1}^{\infty}
\frac{d_{n_1,n_2}^{(q,q)} - \tilde{d}_{n_1,n_2}^{(q)}}{(n_1-1)! (n_2-1)!} \zeta_1^{n_1-1} \zeta_2^{n_2-1}.
\]
\begin{lemma} \label{Lemma_Q_case_2}
 Using Definition~\ref{Definition_Covariance_tilings_St}, we have
 \begin{equation}
 \label{eq_x300}
 \Delta D^{(q)}(\zeta_1,\zeta_2)= \hat{W}_{q_1}\hat{W}_{q_2} \mathcal{F}^{(q,q)}\big(\hat{W}_{q}F^{(q)}(\zeta_1)\,,\,\hat{W}_{q}F^{(q)}(\zeta_2)\big) \cdot \partial_{\zeta_1}F^{(q)}(\zeta_1) \partial_{\zeta_2}F^{(q)}(\zeta_2).
 \end{equation}
\end{lemma}
\begin{proof}
 Subtracting \eqref{eq_x294} with $q_1=q_2 = q$ for $d_{n_1,n_2}^{(q,q)}$ and for $\tilde{d}_{n_1,n_2}^{(q)}$, the terms $\left(\sum_{n=0}^{\infty}
\frac{\zeta_1^n}{\zeta_2^{n+1}}\right)^2$ cancel out and afterwards the argument becomes the same as for $q_1\neq q_2$ in Lemma~\ref{Lemma_Q_case_1}.
\end{proof}

\medskip

\noindent \textsc{Step 4.} Using Proposition~\ref{Proposition_Green_function_explicit} and the height functions $\mathsf{Ht}^{(q)}_\N(x,y)$ inside the trapezoids $\mathcal{T}_q$, we can now restate Theorem~\ref{Theorem_GFF_general} in a more precise form.

\begin{proposition}
 \label{HtGre} For each $q\in[m]$, $k \geq 1$ and $x\in [-\hat{W}_q,0]$, the random variables
 \begin{equation}
 \label{eq_x301}
 \int_{\amsmathbb{R}} \left(\mathsf{Ht}^{(q)}_\N(\N^{-1}\lfloor x \N \rfloor, y)-\E\big[\mathsf{Ht}^{(q)}_\N(\N^{-1}\lfloor x \N \rfloor, y)\big]\right) y^k\dd y
 \end{equation}
 are asymptotically Gaussian (in the sense that the joint moments of these random variables converge towards expressions satisfying the Wick formula) with asymptotic covariance for the triplets $(k_1,q_1,x_1)$ and $(k_2,q_2,x_2)$ given by
 \begin{equation}
 \label{eq_GFF_sections_covariance}
 \frac{1}{\pi} \int_{\{y_1 \mid (x_1,y_1)\in\pounds_{q_1}\}}\int_{\{y_2 \mid (x_2,y_2)\in\pounds_{q_2}\}} \mathrm{Gr}\big(z^{(q_1)}(x_1,y_1),z^{(q_2)}(x_2,y_2)\big)\,y_1^{k_1}y_2^{k_2}\dd y_1\dd y_2,
 \end{equation}
 where the maps $z^{(q)}$ are as in Section~\ref{Section_complex_structure}.
\end{proposition}
\begin{proof}
 Note that the integrand in \eqref{eq_x301} vanishes for large enough $y$ and therefore the integration range can be replaced with a segment $[-M,M]$ for $M$ large enough. Afterwards, we can use \eqref{eq_x295}, which transforms \eqref{eq_x301} (after omitting the integer parts, which will eventually not affect the $\N \rightarrow \infty$ limit) into
 \begin{equation}
 -\frac{(\hat{W}_{q} + x)^{k + 1}}{k+1} \left(p_{k+1}^{(q)}\bigg(1 +\frac{x}{\hat{W}_q}\bigg) -\E\bigg[p_{k+1}^{(q)}\bigg(1 + \frac{x}{\hat{W}_q}\bigg)\bigg]\right).
 \end{equation}
Hence, by Step 1 these random variables are asymptotically Gaussian and the covariance between $(k_1,q_1,x_1)$ and $(k_2,q_2,x_2)$ triplets is given by
 \begin{equation*}
 \begin{split}
 & [\zeta_1^{-1} \zeta_2^{-1}] \prod_{i = 1}^{2} \frac{(\hat{W}_{q_i} + x_i)^{k_i + 1}}{k_i+1} \left( (1+\zeta_i)
\left(\frac{1}{\zeta_i} + \frac{\hat{W}_{q_i}}{\hat{W}_{q_i}+ x_i} \sum_{n_i=1}^\infty \frac{c_{n_i}^{(q_i)}\zeta_i^{n_i-1}}{(n_i-1)!}\right)\right)^{k_i+1} \\
& \phantom{[\zeta_1^{-1} \zeta_2^{-1}] } \times \left( \delta_{q_1,q_2}\left(\sum_{n=0}^{\infty}
\frac{\zeta_1^{n}}{\zeta_2^{n+1}}\right)^2 +\sum_{n_1,n_2=1}^{\infty}
\frac{d_{n_1,n_2}^{(q_1,q_2)}\zeta_1^{n_1-1}\zeta_2^{n_2-1}}{(n_1-1)! (n_2-1)!} \right).
\end{split}
\end{equation*}
Writing the extraction of coefficients as a double integral and using transformations of Step 3, the covariance becomes for the case $q_1\neq q_2$
 \begin{equation}
 \label{eq_x303}
 \begin{split}
& \quad \oint\oint\prod_{i = 1}^{2} \left[\frac{\dd \zeta_i}{2\ii\pi}\cdot\frac{\big((\hat{W}_{q_i} + x_i)F^{(q_i)}(\zeta_i;x_i)\big)^{k_i + 1}}{k_i+1}\cdot \hat{W}_{q_i} \partial_{\zeta_i}F^{(q_i)}(\zeta_i)\right] \\
&\qquad \qquad \qquad \mathcal{F}^{(q_1,q_2)}\big(\hat{W}_{q_1}F^{(q_1)}(\zeta_1),\hat{W}_{q_2}F^{(q_2)}(\zeta_2)\big),
\end{split}
\end{equation}
where the integration goes over two small circles around $0$. Let us make a change of variables in the last integral:
\[
z_i := Z^{(q_i)}(\zeta_i) = \hat{W}_{q_i} F^{(q_i)}(\zeta_i)\qquad i = 1,2.
\]
Since $F^{(q)}$ maps $0$ to $\infty$, the integration now goes over large contours. Comparing with \eqref{eq_x302} we conclude that the covariance becomes:
\[
\oint\oint \prod_{i = 1}^{2} \left[\frac{\dd z_i}{2\ii\pi}\cdot \frac{1}{k_i+1}\left(z_i + \frac{-x_i}{e^{-G^{(q_i)}(z_i/\hat{W}_{q_i})} - 1}\right)^{k_i + 1}\right] \mathcal{F}^{(q_1,q_2)}(z_1,z_2).
\]
For $i = 1,2$, we deform the $z_i$-contour in the last expression, so that its upper/lower half-plane part (depending on the way $\Sigma^{(q_i)}$ was specified in Section~\ref{Section_complex_structure}) becomes the image of $\{(x_i,y_i)\in\pounds_{q_i}\}$, where $x_i$ is fixed, under the uniformization map $(x_i,y_i)\mapsto z^{(q_i)}(x_i,y_i)$ of Section~\ref{Section_complex_structure}. Comparing with \eqref{eq_critical_equation_rescaled}, this implies the power factor between brackets to be $y_i^{k_i + 1}$ on the contour. Otherwise, the contours are kept symmetric with respect to the real axis.

After the contour deformation, we split the contour into four parts, according to the positive/negative imaginary parts of the arguments and integrate by parts in $z_1$ and in $z_2$. Note that for an arbitrary point $o$ in the complex plane, we can decompose
\[
 \int_{z_1^*}^{z_1} \int_{z_2^*}^{z_2} = \int_{o}^{z_1} \int_{o}^{z_1}-\int_o^{z_1^*} \int_o^{z_2} - \int_{o}^{z_1} \int_o^{z_2^*}+\int_o^{z_1^*}\int_o^{z_2^*}.
\]
Hence, the mixed partial derivative of $\mathrm{Gr}$ matches the term $\frac{1}{4\pi}\mathcal{F}^{(q_1,q_2)}$ when $q_1 \neq q_2$ in \eqref{eq_Green_function_explicit}. The sign is aligned with the image of $z^{(q)}$ being in the lower or upper half-plane depending on whether $\mathcal{T}_q$ is left or right trapezoid, and there is an additional minus sign coming from $\ii^2$. Recombining the four contours into a single one after integration by parts, we obtain the final formula for the covariance:
 \begin{equation*}
\frac{1}{\pi} \iint \left[ \prod_{i = 1}^{2} \left( z_i+
\frac{-x}{e^{-G^{(q_i)}( z_i/\hat{W}_{q_i})}-1} \right)^{k_i} \partial_{z_i} \left( z_i+
\frac{-x}{e^{-G^{(q_i)}( z_i/\hat{W}_{q_i})}-1} \right)\right] \textnormal{Gr}(z_1,z_2) \dd z_1 \dd z_2,
\end{equation*}
where the $i$-th integration range is the image of $z^{(q)}(x_i,\cdot)$ in $\Sigma^{(q_i)}$, for $i = 1,2$. Changing the variables to $y_i$ from \eqref{eq_critical_equation_rescaled}, we arrive at \eqref{eq_GFF_sections_covariance} with $q_1\neq q_2$.

We proceed to the case $q_1=q_2 = q$. Using Lemma~\ref{Lemma_Q_case_2} instead of Lemma~\ref{Lemma_Q_case_1}, the expression for the covariance \eqref{eq_x303} is replaced with
 \begin{equation}
 \label{eq_x304}
 \begin{split}
& \quad \oint\oint \left[\prod_{i = 1}^{2} \frac{\dd\zeta_i}{2\ii\pi}\cdot \frac{\big((\hat{W}_{q} + x)F^{(q)}(\zeta_i;x)\big)^{k_i+1}}{k_i +1}\cdot \hat{W}_{q}\partial_{z_i}F^{(q)}(z_i)\right] \mathcal{F}^{(q,q)}\big(\hat{W}_{q}F^{(q)}(\zeta_1),\hat{W}_{q}F^{(q)}(\zeta_2)\big) \\
 & \quad +
\oint\oint \left[\prod_{i = 1}^{2} \frac{\dd \zeta_i}{2\ii\pi}\cdot \frac{\big((\hat{W}_{q} + x)F^{(q)}(\zeta_i;x)\big)^{k_i + 1}}{k_i + 1}\right] \Delta D^{(q,q)}(\zeta_1,\zeta_2).
 \end{split}
\end{equation}
The first term in \eqref{eq_x304} gives rise to the part of the integral \eqref{eq_GFF_sections_covariance} corresponding to the first term in \eqref{eq_Green_function_explicit}. The second term in \eqref{eq_x304} corresponds to the field of fluctuations for the trapezoid with deterministic boundary. The latter fluctuations are given by the Gaussian free field in the upper half-plane by \cite[Theorem 3.14]{BuGo2} and the covariance for this Gaussian free field is precisely $-\frac{1}{2\pi}\log\big|\frac{z_1-z_2}{z_1-z_2^*}\big|$. The computation transforming the second term of \eqref{eq_x304} to the part of the integral \eqref{eq_GFF_sections_covariance} corresponding to the second term in \eqref{eq_Green_function_explicit} is \cite[Section 9.1]{BuGo2}.
\end{proof}

\subsection{Non-orientable case: proof of Theorem~\ref{Theorem_GFF_general_nonor}}
\label{sec:nonorGFFproof}
Most steps for the proofs of Theorem~\ref{Theorem_GFF} for bipartite glued domains are asymptotic statements about individual trapezoids that do not rely on gluing and therefore hold without modification for non-orientable (\textit{i.e.} non-bipartite) glued domains. This is the case with Section~\ref{Section_limit_shape_gluings} establishing the limit shape and the complex structure on the liquid region, and in Section~\ref{Section_GFF_proof} with Step 1 (asymptotic Gaussianity of height functions and formula in terms of coefficients $c$ and $d$ coming from Schur generating functions) and the computations of Step 4.

What requires adaptation in the non-orientable case is the construction of the bordered surface $\Sigma^{\textnormal{half}}$ in Section~\ref{Section_complex_structure}: in absence of left and right distinction of trapezoids, we do not have a consistent way to choose upper and lower half-planes that would glue together to form a Riemann surface. Instead, we apply the construction of Section~\ref{Section_complex_structure} to the orientation double covering $\widetilde{\mathcal{D}}$, where each trapezoid $\mathcal{T}_q$ rather appears in twin pairs consisting of a left trapezoid $\mathcal{T}_{q}^{\textnormal{L}}$ and a right trapezoid $\mathcal{T}_q^{\textnormal{R}}$. This defines $\Sigma^{\textnormal{half}}$ as gluing of upper-half planes $\Sigma^{(q,\textnormal{L})}$ and lower-half planes $\Sigma^{(q,\textnormal{R})}$. More precisely, if $\mathcal{T}_q$ is adjacent to $\mathcal{T}_{q'}$ via the $h$-th segment, we glue $\Sigma^{(q,\textnormal{L})}$ to $\Sigma^{(q',\textnormal{R})}$ and $\Sigma^{(q',\textnormal{L})}$ to $\Sigma^{(q,\textnormal{R})}$ along $(\alpha_h,\beta_h)$. For each trapezoid in $\widetilde{\mathcal{D}}$, we construct by the method of Section~\ref{Section_limit_shape_gluings} uniformizing maps $z^{(q,\textnormal{L})}$ (resp. $z^{(q,\textnormal{R})}$) sending the liquid region $\pounds_{q,\textnormal{L}} \subset \mathcal{T}^{\textnormal{L}}_q$ (resp. $\pounds_{q,\textnormal{R}} \subset \mathcal{T}^{\textnormal{R}}_q$) bijectively and conformally to the upper half-plane (resp. the lower half-plane). Checking that these glue to a uniformizing map from the liquid region $\widetilde{\pounds} \subset \widetilde{\mathcal{D}}$ to $\Sigma^{\textnormal{half}}$ is a statement only involving a left and a right trapezoid at bands joining them, so the proof of Proposition~\ref{Proposition_uniformization_map} remains valid.

\medskip

This being secured, we can continue with Step 2 of Section~\ref{Section_GFF_proof} and relate the asymptotic covariance of Stieltjes transforms of horizontal lozenges in two given trapezoids in $\widetilde{\mathcal{D}}$ to a new kind of Green function.
\begin{definition}
\label{even:green}
The \emph{even} Green function is defined for any two points $p_1,p_2 \in \Sigma^{\textnormal{half}}$
\begin{equation}
\label{GrplusGr}
\textnormal{Gr}^+(p_1,p_2) = \textnormal{Gr}(p_1,p_2) + \textnormal{Gr}(\overline{\varsigma}(p_1),p_2).
\end{equation}
where $\textnormal{Gr}$ is the usual Green function for $\Sigma^{\textnormal{half}}$ and $\overline{\varsigma}$ is the anti-holomorphic involution on $\Sigma^{\textnormal{half}}$ that sends a point $p \in \Sigma^{(q,\textnormal{L})}$ with coordinate $z$ to the point $\overline{\varsigma}(p)$ with coordinate $z^* \in \Sigma^{(q,\textnormal{R})}$. We define the \emph{even Gaussian free field} to be the centered (generalized) Gaussian field with covariance $\textnormal{Gr}^+$.
\end{definition}
Since the involution $z \mapsto z^*$ commutes with the local Laplace operator and by uniqueness of the Green function, the usual Green function on $\Sigma^{\textnormal{half}}$ has the symmetry $\textnormal{Gr}(p_1,p_2) = \textnormal{Gr}(\overline{\varsigma}(p_1),\overline{\varsigma}(p_2))$. Therefore, \eqref{GrplusGr} is a symmetric function of its two variables. Besides, the even Green function has a logarithmic singularity with coefficient $-\frac{1}{2\pi}$ at $p_1 = p_2$ and $p_1 = \overline{\varsigma}(p_2)$. In particular, it solves
\[
-\Delta_{p_1} \textnormal{Gr}^+(p_1,p_2) = \delta(p_1 - p_2) + \delta(p_1 - \overline{\varsigma}(p_2)).
\]
The following result clarifies the nature of the even Gaussian free field.
\begin{lemma}
\label{lem:evenGFFGreen}
The even Gaussian free field is $\sqrt{2}$ times the Gaussian free field conditioned to be invariant under the involution $\overline{\varsigma}$.
\end{lemma}

\begin{proof}
 Let $\textnormal{GFF}(p)$ denote the usual Gaussian free field, and decompose
 \begin{equation}
 \label{eq_x305}
 \textnormal{GFF}(p)= \frac{\textnormal{GFF}(p)+ \textnormal{GFF}(\overline{\varsigma}(p))}{2}+ \frac{\textnormal{GFF}(p)-\textnormal{GFF}(\overline{\varsigma}(p))}{2}.
 \end{equation}
 Using the symmetry $\textnormal{Gr}(p_1,p_2) = \textnormal{Gr}(\overline{\varsigma}(p_1),\overline{\varsigma}(p_2))$, one computes that the two fields in the right-hand side are uncorrelated and, hence, by Gaussianity also independent. Conditioning $\textnormal{GFF}(p)$ to be invariant under the involution $\overline{\varsigma}$ is the same as conditioning on the vanishing of the second term in \eqref{eq_x305}. Due to independence, the conditioning returns the first term $\frac{1}{2}\big(\textnormal{GFF}(p)+\textnormal{GFF}(\overline{\varsigma}(p)\big)$, whose covariance function is readily computed to be
\[
\frac{\textnormal{Gr}(p_1,p_2) + \textnormal{Gr}(\overline{\varsigma}(p_1),p_2)}{2}
\]
Comparing with \eqref{GrplusGr} concludes the proof.
\end{proof}

We are ready to give the replacement of Proposition~\ref{Proposition_Green_function_explicit}, which relates the asymptotic covariance of Stieltjes transform of horizontal lozenges (\textit{i.e.} $\mathcal{F}^{(q_1,q_2)}$ of Definition~\ref{Definition_Covariance_tilings_St}) and the even Green function of Definition~\ref{even:green}.
\begin{proposition}
\label{Proposition_Green_function_explicit_nonor}
The even Green function $\textnormal{Gr}^+$ is given for non-real $z_i \in \Sigma^{(q_i,\textnormal{X}_i)}$ with $q_i \in [m]$ and $\textnormal{X}_i \in \{\textnormal{L},\textnormal{R}\}$ for $i = 1,2$ by
\begin{equation}
\label{gr1222}
\left\{\begin{array}{lll} (-1)^{\textnormal{par}(\textnormal{X}_1,\textnormal{X}_2)}\bigg( -\dfrac{1}{4\pi} \displaystyle\int_{z_1^*}^{z_1}\int_{z_2^*}^{z_2} \mathcal{F}^{(q,q)}(w_1,w_2)\dd w_1 \dd w_2 - \dfrac{1}{2\pi}\log\left|\frac{z_1 - z_2}{z_1 - z_2^*}\right|\bigg) && \textnormal{if}\,\,q_1 = q_2 = q, \\[12pt]
- \dfrac{(-1)^{\textnormal{par}(\textnormal{X}_1,\textnormal{X}_2)}}{4\pi} \displaystyle\int_{z_1^*}^{z_1}\int_{z_2^*}^{z_2} \mathcal{F}^{(q_1,q_2)}(w_1,w_2)\dd w_1\dd w_2 && \textnormal{if}\,\,q_1 \neq q_2, \end{array}\right.
\end{equation}
where the parity $(-1)^{\textnormal{par}(\textnormal{X}_1,\textnormal{X}_2)}$ is $(+1)$ if $\textnormal{X}_1 = \textnormal{X}_2$ and $-1$ otherwise.
In particular, the right-hand side extends to a harmonic function on $\Sigma^{\textnormal{half}}\times \Sigma^{\textnormal{half}}$ minus the diagonal
\end{proposition}
 Note that if $q_i$ and $\textnormal{X}_i$ were coupled like in the bipartite case, \eqref{gr1222} would retrieve \eqref{eq_Green_function_explicit}. Proposition~\ref{Proposition_Green_function_explicit_nonor} is proved (similarly to Proposition~\ref{Proposition_Green_function_explicit}) by relying on Theorem~\ref{thm:spcurvenonbip} and Corollary~\ref{cor:Greenfermion} instead of Theorem~\ref{thm:Omegav} and Corollary~\ref{cor:Greennormal} to check the harmonic extension of the right-hand side of \eqref{gr1222} and Lemma~\ref{Lem:greenferm} for the link with the even Green function.

\medskip

Then, Step 3 in Section~\ref{Section_GFF_proof} is a computation relating the asymptotic covariance of height fluctuations in the trapezoids $\mathcal{T}_{q_1}$ and $\mathcal{T}_{q_2}$ to the function $\mathcal{F}^{(q_1,q_2)}$ of Definition~\ref{Definition_Covariance_tilings_St}, and inserting Proposition~\ref{Proposition_Green_function_explicit} from Step 2 to repackage the outcome in terms of the Green function (Proposition~\ref{HtGre}). This computation only involves individual trapezoids and not their gluing. Therefore, it holds without modification for the fluctuations of the height field in $\widetilde{\mathcal{D}}$, and using \eqref{gr1222} instead of Proposition~\ref{Proposition_Green_function_explicit} we conclude that Proposition~\ref{HtGre} is replaced with the following one exhibiting the even Gaussian free field on $\Sigma^{\textnormal{half}}$.

\begin{proposition}
\label{GrHttnon}
For each $q \in [m]$, $\textnormal{X} \in \{\textnormal{L},\textnormal{R}\}$, $k \geq 1$ and $x \in [-\hat{W}_q,0]$, the random variables
\[
\int_{\amsmathbb{R}} \Big(\widetilde{\mathsf{Ht}}_\N^{(q,\textnormal{X})}(\N^{-1}\lfloor \N x\rfloor,y) - \amsmathbb{E}\big[\widetilde{\mathsf{Ht}}^{(q,\textnormal{X})}_\N(\N^{-1}\lfloor \N x \rfloor,y)\big]\Big)y^k\dd y
\]
are asymptotically Gaussian and the asymptotic covariance between the quadruples $(k_1,q_1,\textnormal{X}_1,x_1)$ and $(k_2,q_2,\textnormal{X}_2,x_2)$ is given by
\[
\frac{1}{\pi} \int_{\{y_1 | (x_1,y_1) \in \pounds_{q_1,\textnormal{X}_1}\}}\int_{\{y_2|(x_2,y_2) \in \pounds_{q_2,\textnormal{X}_2}\}} \textnormal{Gr}^+\big(z^{(q_1,\textnormal{X}_1)}(x_1,y_1),z^{(q_2,\textnormal{X}_2)}(x_2,y_2)\big) y_1^{k_1}y_2^{k_2}\dd y_1\dd y_2.
\]
\end{proposition}
\noindent As already mentioned, the last Step 4 applies verbatim and this concludes the proof of Theorem~\ref{Theorem_GFF_general_nonor}.

\newpage

\part{COMPLEX ANALYSIS FOR THE MASTER EQUATION}

\label{Part_Master_equation}

In this part we study the master problem and a related Riemann--Hilbert problem that we need to solve asymptotically the Nekrasov equations of Chapter~\ref{ChapterNekra}. For that we resort to methods of complex and functional analysis and complex geometry, which are of a rather different nature than Part~\ref{Part1} and \ref{Part_Asymptotic}. The key theorems from this part have already been used by us in Chapter~\ref{Chapter_smoothness} (for derivatives of the equilibrium measure with respect to various parameters), in Chapter~\ref{Chapter_fff_expansions} (for properties of the solution operator $\boldsymbol{\Upsilon}$ for the master problem) and in Chapter~\ref{Chap11} (for identification of the covariance structure with the Gaussian free field).

In contrast to Parts~\ref{Part1} and \ref{Part_Asymptotic}, there is no need for a large parameter $\N$ in this Part~\ref{Part_Master_equation}. We keep notations from the previous parts but the important ones will be recalled. The application of the results of Part~\ref{Part_Master_equation} to Parts~\ref{Part1} and \ref{Part_Asymptotic} follow by specialization to possibly $\N$-dependent parameters, and in this setting the assumptions in Part~\ref{Part_Master_equation} are covered by the ones of Section~\ref{Section_list_of_assumptions}.

The results of this Part are of independent interest. They can also be used for the continuous $\sbeta$-ensembles with several groups of particles (mentioned in Chapter~\ref{SIntro}), because the same Riemann--Hilbert problem appears there. This fact will be exploited to give an alternative proof of the symmetry result of Theorem~\ref{thm:Bsym}, which in the context of random ensembles means symmetry of the expression we give for the leading covariance.

\chapter{Solving the master problem}
\label{Chapter_SolvingN}

Nekrasov equations in Chapter~\ref{ChapterNekra} can be seen as functional relations involving the correlators
\[
(\boldsymbol{W}_1,\boldsymbol{W}_2,\boldsymbol{W}_3,\ldots).
\]
They have an essential role in the bootstrap approach to the asymptotic expansion of the correlators at fixed filling fractions developed in Chapter~\ref{Chapter_fff_expansions}, allowing to solve the linearization of these equations near $(\N \Gm_{\boldsymbol{\mu}},0,0,\ldots)$. The analytic properties of the Stieltjes transform of the equilibrium measure $\Gm_{\boldsymbol{\mu}}$ --- near which we linearize the equations --- were established in Chapter~\ref{Chapter_smoothness}. Using the structural properties of the Nekrasov equations, it then becomes sufficient to solve a \emph{master problem} specified only by a minimal amount of data, the most important one being the endpoints of bands of the equilibrium measure and the interaction matrix $\boldsymbol{\Theta}$.

In Section~\ref{sec:Masterpb}, we describe the master problem, construct its solution operator and establish its continuity in the appropriate topology and its smooth parametric dependence: this is our main result Theorem~\ref{Theorem_Master_equation_12}, which was announced and used in Chapter~\ref{Chapter_fff_expansions}. The regularity properties of the solution operator were also used in the proofs of the parametric regularity for the equilibrium measure (Section~\ref{Section_Smoothnessparam}).

In Section~\ref{Section_Master_by_covariance}, we show that under stronger regularity assumptions, the master problem is equivalent to a vector-valued Riemann--Hilbert problem, and that solutions of the latter can be reconstructed from the knowledge of a fundamental solution $\boldsymbol{\mathcal{F}}(z_1,z_2)$. This fundamental solution appeared in the formulation of the central limit theorems of Chapter~\ref{Chapter_fff_expansions} and \ref{Chapter_filling_fractions} and we study its properties in detail. In particular, we derive its expression for simple matrices $\boldsymbol{\Theta}$, retrieving known formulae. Its systematic study that will turn effective for a larger class of matrices $\boldsymbol{\Theta}$ will be addressed in Chapter~\ref{Chapter_AG}.

\medskip

\section{Master problem and solution operator}
\label{sec:Masterpb}
\label{TheopUpsec}

The defining data of a master problem consist of
\begin{itemize}
\item[(i)] a positive integer $H$ and a real symmetric matrix $\boldsymbol{\Theta}$ of size $H$;
\item[(ii)] pairwise disjoint real compact intervals $\amsmathbb{A}^{\mathfrak{m}}_h := [\hat{a}_h^{\mathfrak{m}},\hat{b}_h^{\mathfrak{m}}]$ and subintervals $[\alpha_h,\beta_h] \subset (\hat{a}_h^{\mathfrak{m}},\hat{b}_h^{\mathfrak{m}})$ indexed by $h \in [H]$, appearing in increasing order with $h$, \textit{i.e.} $\alpha_1 < \beta_1 < \cdots < \alpha_H < \beta_H$;
\item[(iii)] pairwise disjoint open subsets $\amsmathbb{M}_h \subset \amsmathbb{C}$ indexed by $h \in [H]$ and containing $\amsmathbb{A}^{\mathfrak{m}}_h$ but not $\bigcup_{g \neq h} \amsmathbb{A}_{g}^{\mathfrak{m}}$;
\item[(iv)] an $H$-tuple of functions $\boldsymbol{E}(z)=(E_h(z))_{h=1}^H$, called sources, such that $E_h$ is holomorphic in
$\amsmathbb{M}_h\setminus \amsmathbb{A}_h^{\mathfrak{m}}$ for any $h \in [H]$;
\item[(v)] an $H$-tuple of complex numbers $\boldsymbol{\kappa} = (\kappa_h)_{h = 1}^H$.
\end{itemize}
Recall the notations for $\bth_h$ and scalar products in Section~\ref{Theopr}, in particular
\[
\big\langle \bth_h \cdot \boldsymbol{F}(z) \big\rangle = \sum_{g = 1}^{H} \theta_{h,g} F_g(z).
\]

We call $\gamma_h$ any contour surrounding $\amsmathbb{A}_h^{\mathfrak{m}}$ but not $\bigcup_{g \neq h} \amsmathbb{A}_g^{\mathfrak{m}}$ and having counterclockwise orientation. We are going to compute contour integrals of holomorphic functions on such contours $\gamma_h$ and these integrals will not depend on the choice of the particular $\gamma_h$, as long as the integrand is well-defined on the contour. In most cases one can think of $\gamma_h$ as a small loop tightly following the segment $\amsmathbb{A}_h^{\mathfrak{m}}$. We will explicitly comment whenever some particular properties of $\gamma_h$ start being important.

\begin{definition} \label{def:master} The master problem asks to find an $H$-tuple of functions $\boldsymbol{F}(z)=(F_h(z))_{h=1}^H$ with the following properties for any $h \in [H]$.
\begin{enumerate}
\item $F_h(z)$ is a holomorphic function of $z \in \amsmathbb C\setminus \amsmathbb{A}_h^{\mathfrak{m}}$.
\item There exists a holomorphic function $A_h(z)$ of $z \in \amsmathbb{M}_h$ such that
\begin{equation}
\label{eq_Master_equation}
 A_h(z) = \sqrt{(z-\alpha_h) (z-\beta_h)} \,\big\langle \bth_h \cdot \boldsymbol{F}(z)\big\rangle + E_h(z).
 \end{equation}
 \item $F_h(z) = \frac{\kappa_h}{z}+O\big(\frac{1}{z^2}\big)$ as $z \rightarrow \infty$. This can be equivalently restated as
\[
F_h(z)\,\, \mathop{=}_{z \rightarrow \infty} \,\, O\bigg(\frac{1}{z}\bigg) \qquad \textnormal{and} \qquad \oint_{\gamma_h} \frac{\dd z}{2\ii\pi}\,F_h(z) = \kappa_h.
\]
\end{enumerate}
\end{definition}

We recall that the master problem appeared previously as an essential step to solve asymptotically the system of equations satisfied by the correlators, \textit{cf.} Section~\ref{TheopUpsec2}. For the condition \eqref{eq_Master_equation}, let us emphasize that we have no information about $A_h(z)$ except that it is a holomorphic function in $\amsmathbb{M}_h$. This is a non-trivial property because the right-hand side is \textit{a priori} holomorphic in $\amsmathbb{M}_h \setminus \amsmathbb{A}_h^{\mathfrak{m}}$ only. It turns out that this information is sufficient to reconstruct uniquely the function $\boldsymbol{F}(z)$.

\begin{figure}[t]
\begin{center}
\includegraphics[width=0.42\textwidth]{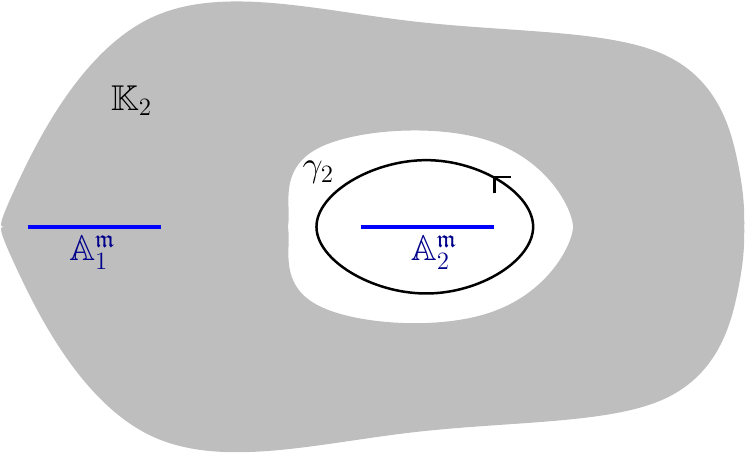} \hfill
\includegraphics[width=0.42\textwidth]{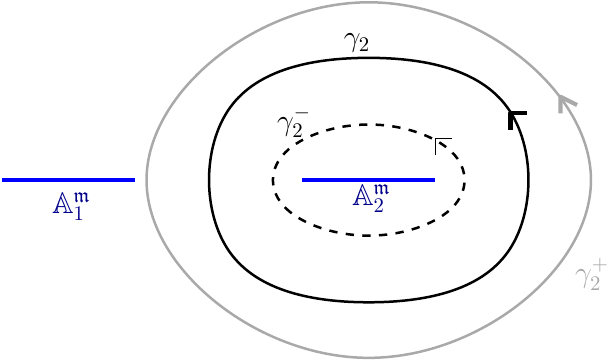}

\caption{\label{Fig:Kdomain} Left panel: Contour and domain in \eqref{ContinuityUp_eq} of Theorem~\ref{Theorem_Master_equation_12}. Right panel: Larger and smaller contours $\gamma_h^+$, $\gamma_h$, $\gamma_h^-$ in Theorem~\ref{Theorem_Masterspecial} and Section~\ref{Section_master_through_Fredholm}.}
\end{center}
\end{figure}

\begin{theorem}
\label{Theorem_Master_equation_12}
Let $C > 0$. Consider the defining data for a master problem such that
\begin{enumerate}
\item $\boldsymbol{\Theta}$ is positive semi-definite;
\item $\forall h \in [H] \quad \theta_{h,h} \geq \frac{1}{C}$;
\item $H \leq C$ and $|\!|\boldsymbol{\Theta}|\!|_{\infty} \leq C$;
\item the distance between any pair of points in $\bigcup_{h = 1}^{H} \big\{\hat{a}_h^{\mathfrak{m}},\alpha_h,\beta_h,\hat{b}^{\mathfrak{m}}_h\big\}$ is larger than $\frac{1}{C}$ and the absolute value of all these points is smaller than $C$.
\end{enumerate}
Then, the master problem of Definition~\ref{def:master} has a unique solution
 \begin{equation}
 \label{eq_solution_Master_equation}
 \forall h \in [H]\qquad \forall z \in \amsmathbb{C} \setminus \amsmathbb{A}_h^{\mathfrak{m}}  \qquad F_h(z) = \Op_h\big[\boldsymbol{E}\,;\,\boldsymbol{\kappa}\big](z).
 \end{equation}
For any $h \in [H]$, the solution operator $\Op_h$ is linear in the arguments $(\boldsymbol{E},\boldsymbol{\kappa}$), depends smoothly on the parameters $\boldsymbol{\Theta}$ and $\boldsymbol{\alpha},\boldsymbol{\beta}$ (\textit{i.e.}, for each $z$ the dependence on these parameters is infinitely differentiable, locally uniformly in $z$), depends on the segments $(\amsmathbb{A}_g^{\mathfrak{m}})_{g = 1}^{H}$ and the domains $(\amsmathbb{M}_g)_{g = 1}^{H}$ only through the domain of definition --- the precise meaning of this statement is explained below --- and is continuous in the following sense. Consider any $H$--tuple of contours $(\gamma_h)_{h=1}^H$, with each $\gamma_h\subset \amsmathbb{M}_h$ being a simple loop surrounding $\amsmathbb{A}_h^{\mathfrak{m}}$, and partitioning $\amsmathbb C$ into bounded interior of $\gamma_h$ and unbounded exterior. For any $H$-tuple of compact subsets $(\amsmathbb{K}_h)_{h=1}^H$, with each $\amsmathbb K_h$ being in the exterior of $\gamma_h$ (\textit{cf.} Figure~\ref{Fig:Kdomain})
 there exists a constant $C' > 0$ depending only on these contours, compacts, and on the constant $C$ --- in particular, independent of $\boldsymbol{E},\boldsymbol{\kappa},h$ --- such that for any $\boldsymbol{E},\boldsymbol{\kappa}$ we have
 \begin{equation}
\label{ContinuityUp_eq}\max_{h \in [H]} \sup_{z \in \amsmathbb{K}_h} \big|\Op_h\big[\boldsymbol{E}\,;\,\boldsymbol{\kappa}\big](z)\big| \leq C' \cdot \Big( \max_{h \in [H]} \sup_{z \in \gamma_h} |E_h(z)| + |\!|\boldsymbol{\kappa}|\!|_{\infty}\Big).
\end{equation}
Besides, if $E_h(z)$ is holomorphic in $\amsmathbb{M}_h$ for each $h \in [H]$, then
\begin{equation}
\label{eq_zero_solution_for_holomorphic}
\forall h \in [H]\qquad \Upsilon_h[\boldsymbol{E}\,;\,\boldsymbol{0}] = 0.
\end{equation}
\end{theorem}

An important feature of the solution is that, while for each $h \in [H]$ the functions $E_h$ can fail to be defined outside $\amsmathbb{M}_h$, the functions $F_h$ are holomorphic everywhere outside $\amsmathbb{A}_h^{\mathfrak{m}}$, and in particular in the complement of $\amsmathbb{M}_h$.

When we say that $\Upsilon_h$ depends on the segments $(\amsmathbb{A}_g^{\mathfrak{m}})_{g = 1}^{H}$ and the complex domains $(\amsmathbb{M}_g)_{g = 1}^{H}$ only through its domain of definition, we mean the following property. Take two master problems differing only by the segments $(\amsmathbb{A}_g^{\mathfrak{m},\textnormal{I}})_{g = 1}^{H}$ or $(\amsmathbb{A}_g^{\mathfrak{m},\textnormal{II}})_{g = 1}^{H}$ and by the complex domains $(\amsmathbb{M}_g^{\textnormal{I}})_{g = 1}^H$ or $(\amsmathbb{M}_g^{\textnormal{II}})_{g = 1}^{H}$. Set
$\amsmathbb{A}_{g}^{\mathfrak{m}} := \amsmathbb{A}_{g}^{\mathfrak{m},\textnormal{I}} \cup \amsmathbb{A}_{g}^{\mathfrak{m},\textnormal{II}}$ for any $g \in [H]$, and assume it is a segment included in $\amsmathbb{M}_{g} := \amsmathbb{M}_{g}^{\textnormal{I}} \cap \amsmathbb{M}_{g}^{\textnormal{II}}$. Then, for any $H$-tuples $\boldsymbol{E}^{\textnormal{I}}$ and $\boldsymbol{E}^{\textnormal{II}}$ such that for any $h \in [H]$ the function $E_h^{\textnormal{J}}$ is holomorphic in $\amsmathbb{M}_h^{\textnormal{J}} \setminus \amsmathbb{A}_{h}^{\mathfrak{m},\textnormal{J}}$ for $\textnormal{J} \in \{\textnormal{I},\textnormal{II}\}$ and the restrictions of $E_h^{\textnormal{I}}$ and $E_h^{\textnormal{II}}$ to $\amsmathbb{M}_h$ agree and equal to a function $E_h$, the following identity holds for any $\boldsymbol{\kappa} \in \amsmathbb{C}^H$:
\[
\forall h \in [H]\qquad \forall z \in \amsmathbb{C} \setminus \amsmathbb{A}_h^{\mathfrak{m}} \qquad \Upsilon^{\textnormal{I}}_h\big[\boldsymbol{E}^{\textnormal{I}}\,;\,\boldsymbol{\kappa}\big](z) = \Upsilon^{\textnormal{II}}_h\big[\boldsymbol{E}^{\textnormal{II}}\,;\,\boldsymbol{\kappa}\big](z) = \Upsilon_h\big[\boldsymbol{E}\,;\,\boldsymbol{\kappa}\big](z),
\]
where $\Upsilon_h$ is the solution operator for the master problem with segments $(\amsmathbb{A}_g^{\mathfrak{m}})_{g = 1}^{H}$ and complex domains $(\amsmathbb{M}_g)_{g = 1}^{H}$, and $\Upsilon^{\textnormal{J}}_h$ is the solution operator for the master problem with segments $(\amsmathbb{A}_g^{\mathfrak{m},\textnormal{J}})_{g = 1}^{H}$ and complex domains $(\amsmathbb{M}_{g}^{\textnormal{J}})_{g = 1}^{H}$ for $\textnormal{J} \in \{\textnormal{I},\textnormal{II}\}$ and $h \in [H]$. In particular, if $E_h$ happens to be holomorphic in $\amsmathbb{C} \setminus [\alpha_h,\beta_h]$, the solution $F_h$ is holomorphic in the same domain.

\medskip

We use Fredholm theory to construct the solution operator $\boldsymbol{\Upsilon}$ for general $\boldsymbol{\Theta}$ in Section~\ref{Section_master_through_Fredholm}. This gives a (complicated) formula for $\Upsilon_h[\boldsymbol{E}\,;\,\boldsymbol{\kappa}](z)$ whose virtue is to make the continuity and regularity with respect to all parameters manifest. An alternative construction of the solution through the Riemann-Hilbert theory is presented in Theorem~\ref{thm:genRHPsol} --- this construction might be more suitable for computing examples, but regularity properties are much harder to see through it.
In general, the solution operator can be very complicated. There are however two simpler cases allowing for a more explicit solution.

\begin{theorem}
\label{Theorem_Masterspecial} If $\boldsymbol{\Theta}$ is a diagonal matrix with positive diagonal elements, then the solution operator is
\begin{equation}
\label{Upsildecoupled}
\Upsilon_h[\boldsymbol{E}\,;\,\boldsymbol{\kappa}](z) = \frac{1}{\sqrt{(z - \alpha_h)(z - \beta_h)}}\left(\kappa_h + \oint_{\gamma_h} \frac{\dd\zeta}{2\ii\pi}\,\frac{E_h(\zeta)}{\theta_{h,h}(\zeta - z)} \right),
\end{equation}
where $\gamma_h\subset \amsmathbb{M}_h$ encloses $\amsmathbb{A}_h^{\mathfrak{m}}$, but not $z$.

If $\theta_{g,h} = \theta > 0$ is independent of $g,h \in [H]$, then the solution operator is uniquely determined by the identity
\[
\sum_{h = 1}^{H} \Op_h[\boldsymbol{E}\,;\,\boldsymbol{\kappa}](z) = \frac{1}{\sigma(z)}\Pi^{-1}\left[\boldsymbol{\kappa} - \bigg(\oint_{\gamma_g^+} \frac{\dd \zeta}{2\ii\pi}\,\mathcal{K}[\boldsymbol{E}](\zeta)\bigg)_{g = 1}^{H}\right](z) + \mathcal{K}[\boldsymbol{E}](z),
\]
where
\[
\mathcal{K}[\boldsymbol{E}](z) = \sum_{h = 1}^{H} \frac{1}{\sigma(z)} \oint_{\gamma_h^-} \frac{\dd\zeta}{2\ii\pi}\,\frac{\sigma(\zeta)}{\zeta - z}\,\frac{E_h(\zeta)}{\theta \sqrt{(\zeta - \alpha_h)(\zeta - \beta_h)}}
\]
and $\Pi^{-1}$ term is a polynomial in $z$ of degree at most $H-1$, defined through the hyperelliptic period map $\Pi$, whose more detailed definition is in \eqref{periodmap}. Contours $\gamma_h^-\subset\amsmathbb{M}_h$ enclose $\amsmathbb{A}_h^{\mathfrak{m}}$, but not $z$. Contours $\gamma_h^+\subset \amsmathbb{M}_h$ enclose the respective $\gamma_h^-$, as Figure~\ref{Fig:Kdomain}.\end{theorem}

\section{Construction of the solution operator (Proof of Theorems~\ref{Theorem_Master_equation_12} and \ref{Theorem_Masterspecial})} \label{Section_master_through_Fredholm}

The proof is decomposed in six parts. In the first part, we prove the theorems for positive diagonal $\boldsymbol{\Theta}$. In the second part for general $\boldsymbol{\Theta}$, we transform the master equation \eqref{eq_Master_equation} to get a functional relation between $\boldsymbol{E}(z)$ and $\boldsymbol{F}(z)$ which does not involve the unknown $\boldsymbol{A}(z)$. In the third part, we solve this equation explicitly by complex analysis techniques when all the entries of $\boldsymbol{\Theta}$ are equal: this completes the proof of Theorem~\ref{Theorem_Masterspecial}. In the fourth part, for general $\boldsymbol{\Theta}$ but $\boldsymbol{\kappa} = 0$ we show invertibility of the master problem using complex analysis and the convexity of the energy functional established in Lemma~\ref{Lemma_I_quadratic_Fourier}. In the fifth part we reinterpret the master problem as the question of inverting an operator in $\mathscr L^2$ space on the integration contour and use this to prove continuity of the inverse via functional analysis and Fredholm theory. In the sixth part, we explain how continuity on the contour implies continuity in the space of holomorphic functions. Finally, in the last part we turn on arbitrary $\boldsymbol{\kappa}$ to reach the solution in full generality. The arguments of this section follow some ideas from \cite{BG11,BEO,BGK,BGG,BG_multicut}.

\medskip

We start by introducing some notations. In the spirit of Definition~\ref{GQdef2}, we set
\[
\sigma(z) = \prod_{h = 1}^{H} \sqrt{(z - \alpha_h)(z - \beta_h)},
\]
which changes the sign when crossing $[\alpha_h,\beta_h]$ segments, and is otherwise a holomorphic function of $z\in\amsmathbb C\setminus \left(\bigcup_{h=1}^H [\alpha_h,\beta_h]\right)$.
If $\boldsymbol{F}$ is an $H$-tuple of functions, we denote
\begin{equation}
\label{FFhsum}
F(z) = \sum_{h = 1}^{H} F_h(z).
\end{equation}
This type of notation was already used in Definition~\ref{Definition_equilibrium_measure} to go from an $H$-tuple of measures to a single measure. The function $F(z)$ is defined in $\amsmathbb{C} \setminus \amsmathbb{A}^{\mathfrak{m}}$, where
\[
\amsmathbb{A}^{\mathfrak{m}} = \bigcup_{h = 1}^{H} \amsmathbb{A}_h^{\mathfrak{m}}.
\]
We recall that $\gamma_h\subset \amsmathbb M_h$ is a positively oriented contour enclosing $\amsmathbb{A}_h^{\mathfrak{m}}$. We will sometimes also need a slightly larger contour $\gamma_h^+\subset \amsmathbb M_h$ which encloses $\gamma_h$ or slightly smaller contour $\gamma_h^-$, which is inside $\gamma_h$, but still encloses $\amsmathbb{A}_h^{\mathfrak{m}}$, as in Figure~\ref{Fig:Kdomain}.

If in the expansion $F=\sum_{h=1}^H F_h$, each term $F_h$ is holomorphic in $\amsmathbb C\setminus \amsmathbb{A}_h^{\mathfrak{m}} $ and we know $F(z)=O(\frac{1}{z})$, $z\rightarrow\infty$, then we can reconstruct $F_h(z)$ from $F(z)$ by using the following projection operator:
\begin{equation}
\label{Fh_projection} \textnormal{Proj}_h[F](z):= \oint_{\gamma_h}
\frac{\dd\zeta}{2\ii\pi}\,\frac{F(\zeta)}{z - \zeta} = F_h(z),
\end{equation}
where $z$ should be outside the contour $\gamma_h$. Thus, the knowledge of $F(z)= \sum_{h = 1}^{H} F_h(z)$ is equivalent to the knowledge of $\boldsymbol{F}(z) = (F_h(z))_{h = 1}^H$.

\medskip

\noindent \textsc{Part 1 of the proof.} The case $\theta_{g,h} = 0$ for $g \neq h$ is simple and is handled in the same way as $H=1$ situation of \cite{BGG}. Let us explain how it works. For any $h \in [H]$, we start by rewriting
\begin{equation*}
\begin{split}
F_h(z) & = \Res_{\zeta = z} \left[\frac{\sqrt{(\zeta - \alpha_h)(\zeta - \beta_h)}}{\sqrt{(z - \alpha_h)(z - \beta_h)}}\,\frac{F_h(\zeta)}{\zeta - z}\right] \\
& = \frac{1}{\sqrt{(z - \alpha_h)(z - \beta_h)}}\left( \Res_{\zeta = \infty} + \oint_{\gamma_h} \frac{1}{2\ii\pi}\right) \frac{\sqrt{(\zeta - \alpha_h)(\zeta - \beta_h)}\,F_h(\zeta)}{z - \zeta}\dd\zeta,
\end{split}
\end{equation*}
where the second equality is by deforming the contour in the definition of the residue at $\zeta=z$ to surround $\infty$ and $\bigcup_{h = 1}^{H} \amsmathbb{A}_h^{\mathfrak{m}}$. Since $F_h(\zeta) \sim \frac{\kappa_h}{\zeta}$ as $\zeta \rightarrow \infty$, the residue at $\infty$ yields $\kappa_h$. The master equation \eqref{eq_Master_equation} says that
\begin{equation}
\label{premaster000}
\theta_{h,h}\sqrt{(\zeta - \alpha_h)(\zeta - \beta_h)}\,F_h(\zeta) = A_h(\zeta) - E_h(\zeta),
\end{equation}
where $A_h$ is holomorphic in $\amsmathbb{M}_h$. Inserting \eqref{premaster000} in the contour integral over $\gamma_h$, the term $A_h$ does not contribute and we obtain
\[
F_h(z) = \frac{1}{\sqrt{(z - \alpha_h)(z - \beta_h)}}\left(\kappa_h + \oint_{\gamma_h} \frac{\dd\zeta}{2\ii\pi}\,\frac{E_h(\zeta)}{\theta_{h,h}(\zeta - z)} \right).
\]
All continuity and parametric smoothness statements are obvious by inspection of this explicit integral formula. In particular, \eqref{ContinuityUp_eq} follows.  This proves the desired Theorems~\ref{Theorem_Master_equation_12} and \ref{Theorem_Masterspecial} for diagonal $\boldsymbol{\Theta}$.

\bigskip

\noindent \textsc{Part 2.} For general $\boldsymbol{\Theta}$, we rather start with
\begin{equation}\label{labF}
F(z) = \Res_{\zeta = z}\left[ \frac{F(\zeta)}{\zeta - z}\,\frac{\sigma(\zeta)}{\sigma(z)} \right].
\end{equation}
Since $\sigma(z)$ is holomorphic for $z \in \amsmathbb{C} \setminus \bigcup_{h = 1}^{H} [\alpha_h,\beta_h]$, we can again deform the contour in the definition of the residue \eqref{labF} to surround $\infty$ and $\bigcup_{h = 1}^{H} \amsmathbb{A}_h^{\mathfrak{m}}$, and get
\begin{equation}
\label{FCauchy} F(z) = \frac{\mathcal{P}[F](z)}{\sigma(z)} + \sum_{h = 1}^H \oint_{\gamma_h}
\frac{\dd\zeta}{2\ii\pi}\,\frac{\sigma(\zeta)}{\sigma(z)}\,\frac{F(\zeta)}{z
- \zeta},
\end{equation}
where $z$ is outside all the contours $\gamma_h$ and
\begin{equation}
\label{Polop}
\mathcal{P}[F](z) = \Res_{\zeta = \infty} \left[\frac{\sigma(\zeta)}{z - \zeta}\,F(\zeta)\right].
\end{equation}
Since $\sigma(\zeta) = \zeta^{H} + O(\zeta^{H - 1})$ and $F(\zeta) = O(\zeta^{-1})$ as $\zeta \rightarrow \infty$ we see (by expanding $\frac{1}{z-\zeta}=\frac{1}{\zeta}+\frac{z}{\zeta^2}+\frac{z^2}{\zeta^3}+\cdots$) that $\mathcal{P}[F](z)$ is a polynomial of degree at most $H - 1$. This is a new feature compared to the case treated in Part 1. In case $\sum_{h = 1}^{H} \kappa_h = 0$, we have $F(\zeta) = O(\zeta^{-2})$ as $\zeta \rightarrow\infty$, so there is a small improvement: $\mathcal{P}[F](z)$ has degree at most $H - 2$. In any case, for each $h \in [H]$ we can write
\[
F(z) = \frac{1}{\theta_{h,h}}\bigg(\big\langle \bth_h\cdot \boldsymbol{F}(z) \big\rangle + \sum_{\substack{g = 1 \\ g \neq h}}^{H}
(\theta_{h,h} - \theta_{h,g})F_{g}(z)\bigg).
\]
We then insert the master equation \eqref{eq_Master_equation} into \eqref{FCauchy}, to find that
\begin{equation}
\label{eq_x128}
F(z) = \frac{\mathcal{P}[F](z)}{\sigma(z)} + \frac{1}{\sigma(z)} \sum_{h = 1}^{H}
\oint_{\gamma_h} \frac{\dd\zeta}{2\ii\pi}\,\frac{\sigma(\zeta)}{\theta_{h,h} ({z-\zeta}
)}\bigg(\frac{ A_{h}(\zeta)-E_h(\zeta) }{\sqrt{(\zeta - \alpha_h)(\zeta - \beta_h)}} + \sum_{\substack{g = 1 \\ g \neq h}}^H
(\theta_{h,h} - \theta_{h,g})F_{g}(\zeta)\bigg).
\end{equation}
As $A_h(\zeta)$ and $\sigma(\zeta)/\sqrt{(\zeta - \alpha_h)(\zeta - \beta_h)}$ are
 holomorphic functions of $\zeta \in \amsmathbb{M}_{h}$, the term involving $A_h(\zeta)$ does not contribute to the contour integral. We now have to solve for
$F(z)$ in terms of $\boldsymbol{E}(\zeta)$, and the unknown quantity $\boldsymbol{A}(\zeta)$ disappeared. We can rewrite the equation \eqref{eq_x128} as
\begin{equation}
\label{rewritingEF}
\begin{split}
& \quad F(z)  - \frac{\mathcal{P}[F](z)}{\sigma(z)} + \sum_{1 \leq g \neq h \leq H} \frac{1}{\sigma(z)} \oint_{\gamma_h} \frac{\dd\zeta}{2\ii\pi} \frac{\sigma(\zeta)}{\zeta - z} \frac{(\theta_{h,h} - \theta_{h,g})F_g(\zeta)}{\theta_{h,h}} \\
& = \sum_{h = 1}^{H} \frac{1}{\sigma(z)}\oint_{\gamma_h} \frac{\dd\zeta}{2\ii\pi}\,\frac{\sigma(\zeta)}{\zeta - z} \,\frac{E_h(\zeta)}{\theta_{h,h}\sqrt{(\zeta - \alpha_h)(\zeta - \beta_h)}},
\end{split}
\end{equation}
where $z$ is outside the integration contours.

\bigskip

\noindent \textsc{Part 3.} We next complete the solution for $\theta_{g,h} = \theta > 0$ for all $g,h \in [H]$. In this case, we have $\theta_{h,h} - \theta_{h,g} = 0$ for any $g,h \in [H]$, hence, \eqref{rewritingEF} becomes
\begin{equation}
\label{FFFFF} F(z) =
\frac{\mathcal{P}[F](z)}{\sigma(z)} + \frac{1}{\sigma(z)} \sum_{h = 1}^H \oint_{\gamma_h}
\frac{\dd\zeta}{2\ii\pi} \frac{\sigma(\zeta)}{\zeta - z} \,\frac{E_h(\zeta)}{\theta
\sqrt{(\zeta-\alpha_h)(\zeta-\beta_h)}}.
\end{equation}
We claim that the coefficients of the polynomial $\mathcal{P}[F](z)$ can be determined by the $H$ conditions
\begin{equation}
\label{FFFFF1}
\forall h \in [H] \qquad \oint_{\gamma_{h}} \frac{\dd z}{2\ii\pi}\,F(z) = \kappa_h.
\end{equation}
Indeed, let $\Pi: \,\amsmathbb{C}_{H - 1}[\zeta] \rightarrow \amsmathbb{C}^{H}$ be the period map taking polynomials to the values of the contour integrals:
\begin{equation}
\label{periodmap} \Pi[P] = \bigg(\oint_{\gamma_h}
\frac{\dd\zeta}{2\ii\pi}\,\frac{P(\zeta)}{\sigma(\zeta)}\bigg)_{h = 1}^{H}.
\end{equation}

\begin{lemma}
\label{Pinvdet} The absolute value of the determinant of $\Pi$ in the basis $z^0,z^1,\ldots,z^{H - 1}$ is bounded from below by a positive constant that can be chosen to depend only on $H$ and $C$ in Theorem~\ref{Theorem_Master_equation_12}. Hence, $\Pi$ is invertible and the norm of its inverse is uniformly controlled by $H$ and $C$.
\end{lemma}
\begin{proof}
Let $h \in [H]$. Let us squeeze the contour $\gamma_h$ to the segment $[\alpha_h,\beta_h]$. There exists $\tau_h \in \{\pm 1\}$ such that for any $x \in [\alpha_h,\beta_h]$ we have $\sigma(x^-) = - \sigma(x^+) = \ii \tau_h |\sigma(x)|$. Thus
\[
\Pi_h[P] = \tau_h \int_{\alpha_{h}}^{\beta_h} \frac{P(x)\dd x}{\pi |\sigma(x)|}.
\]
We can then compute the determinant
\begin{equation}
\label{deting}\Big|\det_{1 \leq g,h \leq H} \Pi_h[z^{g - 1}]\Big|= \biggl|\prod_{h=1}^H \tau_h\biggr| \cdot \biggl|\int_{\alpha_1}^{\beta_1}\cdots\int_{\alpha_{H}}^{\beta_{H}} \Delta(x_1,\ldots,x_{H}) \prod_{h = 1}^{H} \frac{\dd x_h}{\pi |\sigma(x_h)|}\biggr|
\end{equation}
in terms of the Vandermonde determinant
\[
\Delta(x_1,\ldots,x_{H}) = \det_{1 \leq g,h \leq H} x_h^{g - 1} = \prod_{1 \leq g < h \leq H} (x_{h} - x_{g}).
\]
This quantity is positive in the whole range of integration and bounded from below by the product $\prod_{1 \leq g < h \leq H} (\alpha_h - \beta_g)$, which is itself bounded from below by some $c > 0$ depending only on the constant $C$ in the assumptions. Therefore
\begin{equation}
\label{detlower} \Big|\det_{1 \leq g,h \leq H} \Pi_h[z^{g - 1}]\Big| \geq c \prod_{h = 1}^{H} \int_{\alpha_h}^{\beta_h} \frac{\dd x_h}{\pi |\sigma(x_h)|}.
\end{equation}
By assumption, the distance between $\alpha_h$s and $\beta_g$s is bounded from below by $\frac{1}{C}$ while $\alpha_h$ and $\beta_g$ remain in the compact $[-C,C]$. Therefore \eqref{detlower} is bounded from below by a constant which only depends on $H$ and $C$.
\end{proof}

Since $\Pi$ is invertible, the system formed by \eqref{FFFFF}-\eqref{FFFFF1} has a unique solution given by
\begin{equation}
\label{Ups_special}
\forall h \in [H]\qquad F(z) = \frac{1}{\sigma(z)} \Pi^{-1}\left[\boldsymbol{\kappa} - \bigg(\oint_{\gamma_g^+} \frac{\dd\zeta}{2\ii\pi}\, \mathcal{K}[\boldsymbol{E}](\zeta)\bigg)_{g = 1}^{H}\right](z) + \mathcal{K}[\boldsymbol{E}](z),
\end{equation}
where
\begin{equation}
\label{Khdef}
\mathcal{K}[\boldsymbol{E}](z)= \frac{1}{\sigma(z)} \sum_{h = 1}^{H} \oint_{\gamma_h}
\frac{\dd\zeta}{2\ii\pi}\,\frac{\sigma(\zeta)}{\zeta - z}\,\frac{E_h(\zeta)}{\theta
\sqrt{(\zeta - \alpha_h)(\zeta - \beta_h)}}
\end{equation}
for $z$ outside of the contour of integration. If needed this contour can be moved closer to $\amsmathbb{A}_h^{\mathfrak{m}}$. This justifies the second part of Theorem~\ref{Theorem_Masterspecial}. We can also decompose back $F(z) = \sum_{h = 1}^{H} F_h(z)$ with the projection operator from \eqref{Fh_projection}. From this explicit solution and Lemma~\ref{Pinvdet} we can easily bound $\mathcal{K}$, and then $\Op_h$, uniformly in the constant $C$ of the assumptions to obtain the operator continuity and parametric smoothness statements. This proves Theorems~\ref{Theorem_Master_equation_12} and \ref{Theorem_Masterspecial} when all the entries of $\boldsymbol{\Theta}$ are equal.

\bigskip

\noindent \textsc{Part 4.} We proceed to the general $\boldsymbol{\Theta}$ case. In this situation, we are no longer able to give explicit formula for a solution. Instead, we treat the left-hand side of \eqref{rewritingEF} as a linear operator applied to function $F(z)$ and show that this operator is invertible. Note that for any polynomial $p(z)$, the transformation $F(z)\mapsto F(z)+\frac{p(z)}{\sigma(z)}$ does not change the left-hand side of \eqref{rewritingEF}. Hence, the linear operator of interest has a kernel and we need to get rid of it before studying invertibility. We start the analysis by assuming $\boldsymbol{\kappa} = 0$, which would imply that on the linear subspace of functions $F$ satisfying period conditions \eqref{FFFFF1} the invertibility is restored.

In order to make this precise, we introduce certain function spaces and operators on them.

\begin{definition}
\label{Definition_complex_dom}
If $\amsmathbb{O} \subset \amsmathbb{C}$ is an open set, we denote $\mathscr{O}(\amsmathbb{O})$ the linear space of holomorphic functions on $\amsmathbb{O}$. If $m \in \amsmathbb{Z}_{> 0}$ we denote
\begin{equation*}
\begin{split}
\mathscr{H}^{[m]} & = \bigg\{f \in \mathscr{O}(\amsmathbb{C} \setminus \amsmathbb{A}^{\mathfrak{m}}) \quad \bigg| \quad f(z) \mathop{=}_{z \rightarrow \infty} O\bigg(\frac{1}{z^m}\bigg) \bigg\}, \\
\mathscr{H}^{[m],0} & = \bigg\{f \in \mathscr{H}^{[m]} \quad \bigg| \quad \forall h \in [H] \quad \oint_{\gamma_h} \frac{\dd z}{2\ii\pi}\, f(z) = 0 \bigg\}.
\end{split}
\end{equation*}
\end{definition}

Notice that $\mathscr{H}^{[1],0} = \mathscr{H}^{[2],0}$ because the vanishing of the contour integrals prevents any residue at $\infty$.

\begin{definition}
\label{defofop} We introduce four linear operators
\begin{equation*}
\begin{split}
\mathcal{K}_h \,\,:\,\, \mathscr{O}(\amsmathbb{M}_h \setminus \amsmathbb{A}_h^{\mathfrak{m}}) \longrightarrow \mathscr{H}^{[H + 1]}, & \qquad \mathcal{K}\,\, :\,\, \prod_{h = 1}^{H} \mathscr{O}(\amsmathbb{M}_h \setminus \amsmathbb{A}_h^{\mathfrak{m}}) \longrightarrow \mathscr{H}^{[H + 1]}, \\
\mathcal{U}_{h}\,\,:\,\, \mathscr{H}^{[1]} \longrightarrow \mathscr{O}(\amsmathbb{M}_h \setminus \amsmathbb{A}_h^{\mathfrak{m}}), & \qquad \mathcal{P}\,\, :\,\, \mathscr{H}^{[1]} \longrightarrow \amsmathbb{C}_{H - 1}[z],
\end{split}
\end{equation*}
defined by the formulae ($z$ is outside the integration contours)
\begin{equation*}
\label{UhPl}
\begin{split}
\mathcal{K}_h[f](z) & = \frac{1}{\sigma(z)} \oint_{\gamma_h} \frac{\dd \zeta}{2\ii\pi} \,\frac{\sigma(\zeta)}{\zeta - z}\,\frac{f(\zeta)}{\theta_{h,h} \sqrt{(\zeta - \alpha_h)(\zeta - \beta_h)}}, \\
\mathcal{K}[\boldsymbol{f}](z) & = \sum_{h = 1}^H \mathcal{K}_h[f_h], \\
\mathcal{U}_h[f](z) & = \bigg(\sum_{\substack{g = 1 \\ g \neq h}}^H (\theta_{h,h} - \theta_{h,g})f_g(z)\bigg) \sqrt{(z - \alpha_h)(z - \beta_h)}, \\
\mathcal{P}[f](z) & = \Res_{\zeta = \infty}\left[ \frac{\sigma(\zeta)}{z - \zeta}\,f(\zeta)\right].
\end{split}
\end{equation*}
\end{definition}
The operator $\mathcal{K}$ already appeared in \eqref{Khdef} when $\theta_{h,h} = \theta$, and $\mathcal{K}_h$ are simply its summands. The fact that $\mathcal{K}$ and $\mathcal{K}_h$ land in $\mathscr{H}^{[H + 1]}$ comes from the asymptotic behavior $\sigma(z) \sim z^H$ as $z \rightarrow \infty$. The 'polynomial part' operator $\mathcal{P}$ was already discussed in \eqref{Polop}. In the formula for $\mathcal{U}_h[f]$, we used the fact that $f \in \mathscr{H}^{[1]}$ determines uniquely an $H$-tuple of functions $\boldsymbol{f} = (f_h)_{h = 1}^{H}$ such that $f(z) = \sum_{h = 1}^{H} f_h(z)$ by the projection formula $f_h(z) = \textnormal{Proj}_h[f](z)$ of \eqref{Fh_projection}.
The equation \eqref{rewritingEF} we want to solve then takes the form
\begin{equation}
\label{Idequ}\bigg(\textnormal{Id} - \sigma^{-1}\cdot \mathcal{P} + \sum_{h = 1}^H \mathcal{K}_{h} \circ \mathcal{U}_{h} \bigg)[F](z) = \mathcal{K}[\boldsymbol{E}](z),
\end{equation}
where $E_h \in \mathscr{O}(\amsmathbb{M}_{h}^{o})$ are given functions indexed by $h \in [H]$. As we impose $\boldsymbol{\kappa} = 0$, we look for unknown $F$ belonging to the space $\mathscr{H}^{[2],0}$.

\begin{lemma} \label{Lemma_eq_zero_solution}
 Choose a polynomial $p(z)$ of degree at most $H-1$ and suppose that $f \in \mathscr{H}^{[2],0}$ satisfies
\begin{equation}
\label{Gder} f(z) =
\frac{p(z)}{\sigma(z)} + \sum_{1 \leq g \neq h \leq H} \frac{1}{\sigma(z)}
\oint_{\gamma_h} \frac{\dd\zeta}{2\ii\pi}\,\frac{\sigma(\zeta)}{z - \zeta}\,\frac{(\theta_{h,h} -
\theta_{h,g})f_{g}(\zeta)}{\theta_{h,h}},
\end{equation}
where $z$ is outside all the contours $\gamma_h$. Then $f$ is identically equal to zero.
\end{lemma}

\begin{corollary}
\label{Corollary_Inj} The linear operator $\textnormal{Id} - \sigma^{-1}\cdot \mathcal{P} + \sum_{h = 1}^H \mathcal{K}_{h}\circ \mathcal{U}_h$ with domain $\mathscr{H}^{[2],0}$ is injective.
\end{corollary}
\begin{proof}[Proof of Corollary~\ref{Corollary_Inj}]
 If we choose $p(z)=\mathcal{P}[f](z)$, then \eqref{Gder} is equivalent to $f$ being in the kernel of $\textnormal{Id} - \sigma^{-1}\cdot \mathcal{P} + \sum_{h = 1}^H \mathcal{K}_{h}\circ \mathcal{U}_h$.
\end{proof}

\begin{proof}[Proof of Lemma~\ref{Lemma_eq_zero_solution}]
Let $f \in \mathscr{H}^{[2],0}$ satisfy \eqref{Gder}. For real $x$ we define
\begin{equation} \label{eq_x129}
 \nu(x)=\frac{1}{2\pi} \left( f(x^-)-f(x^+) \right).
\end{equation}
We claim that $\nu(x)$, $x\in\amsmathbb R$, satisfies two properties:
\begin{enumerate}
 \item $\nu_h(x)\dd x$ is integrable; $\nu(x)\sigma(x)$ is smooth inside $\bigcup_{h=1}^H [\alpha_h,\beta_h]$ and vanishes outside.
 \item Denoting $\nu_h(x)$ the restriction of $\nu(x)$ onto $[\alpha_h,\beta_h]$, the function $f_h(z)$ is reconstructed as the Stieltjes transform of the complex-valued finite measure $\nu_h(x)\dd x$:
   \begin{equation}
   \label{eq_x143}
   f_h(z)=\int_{\alpha_h}^{\beta_h} \frac{\nu_h(x)\dd x}{z-x}.
   \end{equation}
\end{enumerate}
In order to see the first property, we compute $\nu(x)$ by combining \eqref{eq_x129} with the right-hand side of \eqref{Gder}. The contribution of the first term in \eqref{Gder} satisfies the desired property, because $p(z)$ is a polynomial in $z$. In the second term, the integrand, and hence, the integral depend on $z$ in a holomorphic way, and therefore, when combined with $\frac{1}{\sigma(z)}$ prefactor satisfies the desired property. Note that when taking $z$ to be $x^-$ or $x^+$ in the integral in \eqref{Gder}, we might violate the condition that $z$ is outside $\gamma_h$. Thus, we need to additionally take into account the residue of the integral at $\zeta=z$ when dealing with $x$ inside $\gamma_h$. This residue is $\sigma(z)$ times a holomorphic function of $x$ inside $\gamma_h$ (here it is important that $\theta_{h,h}-\theta_{h,g}=0$ for $h=g$), and, hence, combined with $\frac{1}{\sigma(z)}$ prefactor, it has no jump on the real axis and does not contribute to $\nu(x)$.

In order to see the second property of $\nu(x)$, we consider the difference
\[
 g_h(z)= f_h(z)-\int_{\alpha_h}^{\beta_h} \frac{\nu_h(x)\dd x}{z-x}.
\]
Projection formula \eqref{Fh_projection} and the definition $f(z)=\sum_{h=1}^H f_h(z)$ imply that the jump of $f_h(z)$ when crossing the real axis coincides with that of $f(z)$ on the interval $[\alpha_h,\beta_h]$; there are no jumps outside the interval. By explicitly computing the jump of the subtraction term in the definition of $g_h(z)$, we get the same answer as for $f_h(z)$. We conclude that $g_h(z)$ has no jump on $[\alpha_h,\beta_h]$ and the singularity along this segment is removable (\textit{e.g.} by an application of the Morera theorem inside $(\alpha_h,\beta_h)$ and Riemann theorem on removable singularities at $\alpha_h$ and $\beta_h$). Hence, $g_h(z)$ is an entire function. In addition, by \eqref{Fh_projection} (for its applicability the $O(1/z)$ decay of $f(z)$, which follows from \eqref{Gder}, is used), $g_h(z)$ decays as $O(1/z)$ as $z \rightarrow \infty$. By Liouville theorem then $g_h(z)=0$.

\medskip

In the rest of the proof we show that in fact the complex measure $\nu(x)\dd x$ is necessary equal to identical zero, and therefore so is $f(z)$. The idea is to use Corollary~\ref{Corollary_I_positive}, for which we need to connect the equation \eqref{Gder} to the quadratic logarithmic energy functional of Definition~\ref{DEFI2}. This is achieved through the following auxiliary identity: for each $h\in[H]$ there exists a constant $v_h$, such that
\begin{equation}
\label{eq_x192} \forall x \in [\alpha_{h},\beta_{h}] \qquad \sum_{g=1}^H 2\theta_{h,g} \int_{\alpha_g}^{\beta_g} \log|x -
y|\,\nu_{g}(y)\dd y = v_h.
\end{equation}
Let us prove this identity. Since \eqref{eq_x192} is continuous in $x$, it is sufficient to take two values $\alpha_h<x_1<x_2<\beta_h$ and show that the difference of \eqref{eq_x192} at $x=x_1$ and at $x=x_2$ vanishes. It is helpful to smoothen the singularity of the logarithm by adding a small positive term to its argument, and rewrite the desired identity as
\[
\lim_{\eps\rightarrow 0} \Biggl[\sum_{g=1}^H \theta_{h,g} \int_{\alpha_g}^{\beta_g} \log\big(( (x_2 - y)^2+\eps^2\big)\,\nu_{g}(y)\dd y
- \sum_{g = 1}^H \theta_{h,g} \int_{\alpha_g}^{\beta_g} \log\big(( (x_1 - y)^2+\eps^2\big)\,\nu_{g}(y)\dd y \Biggr]\stackrel{?}{=} 0.
\]
Expressing the difference as the integral of derivative, we equivalently want to show
\[
\lim_{\eps \rightarrow 0} \int_{x_1}^{x_2} \Biggl[\sum_{g=1}^H \theta_{h,g} \int_{\alpha_g}^{\beta_g} \frac{2(x-y)}{(x - y)^2+\eps^2}\,\nu_{g}(y)\dd y \Biggr] \dd x \stackrel{?}{=} 0.
\]
Using \eqref{eq_x143}, this is equivalent to
\begin{equation}
\label{eq_x184}
\lim_{\eps \rightarrow 0} \int_{x_1}^{x_2} \Biggl[\theta_{h,h}
\big((f_h(x+\ii \eps)+f_h(x-\ii \eps)\big) + \sum_{\substack{g = 1 \\ g \neq h}}^{H} \theta_{h,g} \big(( f_g(x+\ii\eps)+f_g(x-\ii \eps)\big) \Biggr] \dd x \stackrel{?}{=} 0.
\end{equation}
We claim that the last integrand vanishes as $\eps \rightarrow 0$. Indeed, arguing as for the first property of $\nu(x)$, we add the upper- and lower-boundary values of \eqref{Gder} to get
\begin{equation}
\label{eq_x133}
 \forall x \in (\alpha_{h},\beta_{h})\qquad f(x^+) + f(x^-)= 2 \sum_{g \neq h} \frac{\theta_{h,h}-\theta_{h,g}}{\theta_{h,h}} f_{g}(x),
\end{equation}
where the right-hand side comes the contribution of the residue of the integral in \eqref{Gder} at $\zeta=z$. Multiplying \eqref{eq_x133} by $\theta_{h,h}$ and using $f(x)=\sum_{h=1}^H f_h(x)$, we get
\[
\forall h \in [H]\qquad \forall x \in (\alpha_{h},\beta_{h})\qquad \theta_{h,h}\big(f_h(x^+) + f_h(x^-)\big) + \sum_{g \neq h} 2\theta_{h,g}f_{g}(x) = 0,
\]
which implies \eqref{eq_x184} and finishes the proof of \eqref{eq_x192}.

\smallskip

We then integrate \eqref{eq_x192} over $x$ with respect to conjugated measure $\nu^*_{h}(x)\dd x$. Since $f \in \mathscr{H}^{[2],0}$, using \eqref{eq_x143} we have
\[
\forall h \in [H] \qquad 0 = \oint_{\gamma_h} \frac{\dd z}{2\ii\pi}\,f(z) = \oint_{\gamma_h} \frac{\dd z}{2\ii\pi}\,f_{h}(z) = \int_{\alpha_h}^{\beta_h} \nu_h(x)\dd x.
\]
Therefore, the constant in the right-hand side of \eqref{eq_x192} disappears after this integration. Summing over all $h \in [H]$, we conclude
\[
\sum_{g,h = 1}^H \theta_{h,g} \iint_{\amsmathbb{R}^2} \log|x - y|\,\nu_{h}(x)\nu_{g}^*(y)\dd x \dd y = 0.
\]
Using Lemma~\ref{Lemma_I_quadratic_Fourier}, this quadratic functional can be rewritten in Fourier space\footnote{In Lemma~\ref{Lemma_I_quadratic_Fourier}, $\nu$ was a real signed measure, while here $\nu$ might be complex. However, by bilinearity the complex identity is implied by the real version as soon as we add conjugations. The condition \eqref{eq_nu_mass} in Lemma~\ref{Lemma_I_quadratic_Fourier} follows from the fact that $\nu_h$ has zero mass on $[\alpha_h,\beta_h]$ for each $h \in [H]$.} as
\begin{equation}
\label{therightshad} - \sum_{g,h = 1}^H \theta_{h,g} \int_{\amsmathbb{R}} \frac{\dd s}{|s|}\,\hat{\nu}_{g}(s) \big(\hat{\nu}_{h}(s)\big)^* = 0.
\end{equation}
As $\boldsymbol{\Theta}$ satisfies the relevant parts of Assumption~\ref{Assumptions_Theta}, we can apply Corollary~\ref{Corollary_I_positive}, implying that \eqref{therightshad} vanishes if and only if $\nu_h = 0$ for all $h \in [H]$. We conclude that $f(z) = 0$.
\end{proof}

\bigskip

\noindent \textsc{Part 5.} We have learned in Corollary~\ref{Corollary_Inj} that $\textnormal{Id} - \sigma^{-1}\cdot\mathcal{P} + \sum_{h
= 1}^H \mathcal{K}_{h}\circ \mathcal{U}_{h} $ with domain $\mathscr{H}^{[2],0}$ is invertible on its image, and the next step is to show that the inverse is a bounded operator in an appropriate norm. We would like to use the Fredholm alternative, and for this we need to adjust the space on which the operator acts. Two features are important: first, rather than space of holomorphic functions $\mathscr{H}^{[m]}$, we need to switch to $\mathscr{L}^2$ space on a contour in order to identify $- \sigma^{-1}\cdot\mathcal{P} + \sum_{h
= 1}^H \mathcal{K}_{h}\circ \mathcal{U}_{h} $ with a compact integral operator; second, the functions $\frac{p(z)}{\sigma(z)}$ would be in such a $\mathscr{L}^2$ space, but also in the kernel of $\textnormal{Id} - \sigma^{-1}\cdot\mathcal{P} + \sum_{h
= 1}^H \mathcal{K}_{h}\circ \mathcal{U}_{h}$, and we need to adjust the definition of the operator on these functions. We proceed to the details.

We fix two sets of the contours: slightly larger $\gamma^+ = \bigcup_{h = 1}^{H} \gamma_h^+$ and slightly smaller $\gamma^- = \bigcup_{h = 1}^{H} \gamma_h^-$, as in Figure~\ref{Fig:Kdomain}. We would like to equip the union of the larger contours $\gamma^+$ with the curvilinear Lebesgue measure $|\dd z|$ and treat the operator as an endomorphism of $\mathscr{L}^2(\gamma^+)$ (the space of complex-valued functions whose modulus squared is integrable).

\begin{lemma}
\label{Lem:ConversionL2op}
 Consider the integral operator
\begin{equation}
\label{Lintegral}
 \begin{array}{llcll}\mathcal{S} & : & \mathscr{L}^2(\gamma^+) & \longrightarrow & \mathscr{L}^2(\gamma^+) \\[2pt] && f& \longmapsto & \displaystyle \oint_{\gamma^+} f(w)\, \dfrac{ S(z,w)}{2\ii\pi}\, \dd w \end{array}
\end{equation}
with integral kernel for $(z,w) \in \gamma^+ \times \gamma^+_{g}$ and $g \in [H]$ given by
\begin{equation}
\label{Szwkernel}
S(z,w) = \frac{1}{\sigma(z)}\bigg(
\Res_{\zeta = \infty} \left[\frac{\sigma(\zeta)}{(z - \zeta)(w - \zeta)}\right] + \sum_{\substack{h = 1 \\ h \neq g}}^{H} \frac{\theta_{h,h} -
\theta_{h,g}}{\theta_{h,h}}\oint_{\gamma^-_{h}}
\frac{\dd \zeta}{2\ii\pi}\,\frac{\sigma(\zeta)}{(z - \zeta)(w - \zeta)}\bigg).
\end{equation}
Then for each $f \in \mathscr{H}^{[2],0}$, which restricts onto $\gamma^+$ as $\hat f$, we have an identification with Definition~\ref{defofop}:
\begin{equation}
\label{eq_x193}
\mathcal{S}[\hat f](z) = - \frac{\mathcal{P}[f](z)}{\sigma(z)} +
\sum_{h = 1}^H \mathcal{K}_{h}\circ \mathcal{U}_{h}[f](z), \qquad z\in \gamma^+.
\end{equation}
\end{lemma}
\begin{proof} Using the projection formula \eqref{Fh_projection} we rewrite for $f \in \mathscr{H}^{[2],0}$ and $h \in [H]$
\begin{equation}
\begin{split}
\mathcal{K}_{h}\circ \mathcal{U}_{h}[f](z) & = \sum_{1\leq g \neq h \leq H} \frac{\theta_{h,h} - \theta_{h,g}}{\theta_{h,h} \sigma(z)} \oint_{\gamma_{h}} \frac{\dd\zeta}{2\ii\pi}\, \frac{\sigma(\zeta)}{\zeta - z} \oint_{\gamma_g} \frac{\dd w}{2\ii \pi} \frac{f(w)}{\zeta - w} \\
& = \sum_{1\leq g \neq h \leq H} \oint_{\gamma_g} \frac{\dd w}{2\ii\pi}\, \frac{f(w)}{\sigma(z)} \oint_{\gamma_{h}} \frac{\dd \zeta}{2\ii\pi}\, \frac{(\theta_{h,h} - \theta_{h,g})\sigma(\zeta)}{\theta_{h,h} (\zeta - z)(\zeta - w)},
\end{split}
\end{equation}
where both $z$ and $w$ should be outside the integration contours. We also rewrite
\[
\frac{\mathcal{P}[f](z)}{\sigma(z)} = \Res_{\zeta = \infty} \left[\frac{f(\zeta)\sigma(\zeta)}{(z - \zeta)\sigma(z)}\right] = \sum_{g = 1}^{H} \oint_{\gamma_g} \frac{\dd w}{2\ii\pi}\, \frac{f(w)}{\sigma(z)} \Res_{\zeta = \infty}\left[ \frac{\sigma(\zeta)}{(z - \zeta)(\zeta - w)}\right],
\]
where $z$ is outside the integration contour. The last equality comes writing the $w$-integral as the sum of the residues outside the contour: vanishing residue at $\infty$ and residue at $w=\zeta$. Therefore, deforming the $\gamma_g$ contours to $\gamma_g^+$ and deforming the $\gamma_h$ contours to $\gamma_h^-$, we get for $z\in\gamma^+$:
\begin{equation*}
\begin{split}
 & \quad - \frac{\mathcal{P}[f](z)}{\sigma(z)} + \sum_{h = 1}^H \mathcal{K}_{h}\circ \mathcal{U}_{h}[f](z) \\
 & = \sum_{g = 1}^{H} \oint_{\gamma_g^+} f(w)\cdot \frac{\dd w}{2\ii\pi}\, \frac{1}{\sigma(z)} \Res_{\zeta = \infty}\left[ \frac{\sigma(\zeta)}{(z - \zeta)(w-\zeta)}\right] \\
 & + \sum_{g \neq h} \oint_{\gamma_g^+} f(w) \, \frac{\dd w}{2\ii\pi}\, \frac{1}{\sigma(z)} \oint_{\gamma_{h}^-} \frac{\dd \zeta}{2\ii\pi}\, \frac{(\theta_{h,h} - \theta_{h,g})\sigma(\zeta)}{\theta_{h,h} (\zeta - z)(\zeta - w)},
\end{split}
\end{equation*}
which is the same as $\displaystyle \oint_{\gamma^+} \hat f(w)\, \dfrac{ S(z,w)}{2\ii\pi}\, \dd w$.
\end{proof}

\begin{lemma}
\label{Lem:ConversionL2op2}
 Let $p_1(z),\ldots,p_H(z)$ be polynomials of degree (at most) $H-1$, such that
 \begin{equation}
 \label{eq_x236}
  \oint_{\gamma_h} \frac{p_g(\zeta)}{\sigma(\zeta)} \frac{\dd \zeta}{2\ii\pi}= \delta_{g,h}.
 \end{equation}
 Consider the integral operator $\mathcal{P}^H:\mathscr{L}^2(\gamma^+) \longrightarrow \mathscr{L}^2(\gamma^+)$ given by $\displaystyle f \longmapsto \oint_{\gamma^+} f(w)\, \dfrac{ P^H(z,w)}{2\ii\pi}\dd w$, where
 \begin{equation}
 \label{PolopH}
 P^H(z,w)= \frac{p_h(z)}{\sigma(z)},\quad \textnormal{ for } w\in \gamma_h^+, \quad h\in [H].
 \end{equation}
 Then for any polynomial $f(z)$ of degree at most $H-1$, we have
 \begin{equation}
 \label{eq_x237}
 \mathcal P^H \left[\frac{f}{\sigma}\right](z)=\frac{f(z)}{\sigma(z)}.
 \end{equation}
 Further, for any $f\in \mathscr{L}^2(\gamma^+)$ with vanishing integrals around each $\gamma_h$, $h\in[H]$, we have $\mathcal P^H f=0$.
\end{lemma}
\begin{proof}
 The existence of the polynomials satisfying \eqref{eq_x236} follows from Lemma~\ref{Pinvdet}. The property \eqref{eq_x237} follows directly from \eqref{eq_x236} for $f= p_h(z)$, $h\in[H]$, and is extended to all polynomials of degree at most $H-1$ by linearity of $\mathcal P^H$. The identity $\mathcal P^H[f]=0$ for $f$ with vanishing integrals around $\gamma_h$ directly follows from the definition of $P^H(z,w)$.
\end{proof}
\begin{remark}
\label{Remark_about_op}
 The operator $\mathcal P$ in Definition~\ref{defofop} satisfies a similar property $\left(\sigma^{-1}\cdot \mathcal P\right) \frac{f(z)}{\sigma(z)}=\frac{f(z)}{\sigma(z)}$ for polynomials $f(z)$. However, the difference is that $\mathcal P^H$ preserves only polynomials of degree smaller than $H$ and its image has dimension $H$.
\end{remark}

The relevance of operators $\mathcal S$ and $\mathcal P^H$ to solution of \eqref{Idequ} is based on the following statement.

\begin{lemma} \label{Lemma_invertibility_Fredholm}
 $\mathcal S-\mathcal S \circ \mathcal P^H$ is a compact integral operator in $\mathscr{L}^2(\gamma^+)$. The operator $\textnormal{Id} + \mathcal S-\mathcal S \circ \mathcal P^H$ is invertible and its inverse can be represented as $\textnormal{Id}- \mathcal R$, where $\mathcal R$ is an integral operator whose kernel is smooth and smoothly depends on the parameters $\boldsymbol{\Theta}$ and $\boldsymbol{\alpha},\boldsymbol{\beta}$, with each partial derivative uniformly bounded in terms of $H$ and constant $C$ from Theorem~\ref{Theorem_Master_equation_12}. Finally, for each $f \in \mathscr{H}^{[2],0}$, which restricts onto $\gamma^+$ as $\hat f$, we have an identification with Definition~\ref{defofop}:
\begin{equation}
\label{eq_x238}
(\textnormal{Id} + \mathcal S- \mathcal S \circ \mathcal P^H)[\hat f](z) = f(z) - \frac{\mathcal{P}[f](z)}{\sigma(z)} +
\sum_{h = 1}^H \mathcal{K}_{h}\circ \mathcal{U}_{h}[f](z), \qquad z\in \gamma^+.
\end{equation}
\end{lemma}
\begin{proof} Let us first prove \eqref{eq_x238}. For $f\in \mathscr{H}^{[2],0}$ all contour integrals of $\hat f$ around $\gamma_h$ vanish and therefore $(\mathcal S \circ \mathcal P^H)[\hat f]=0$. Then \eqref{eq_x193} implies the desired statement.

Next, the kernels for both operators $\mathcal S$ and $\mathcal P^H$ are smooth, and so is the kernel for their product and for $\mathcal S-\mathcal S \circ \mathcal P^H$. Hence, the latter operator is compact. By the Fredholm alternative, $\textnormal{Id} + \mathcal S-\mathcal S \circ \mathcal P^H$ is invertible if and only if it is injective. Let us check that the last operator has no kernel. Indeed, suppose that for $f\in \mathscr{L}^2(\gamma^+)$, we have
 \begin{equation}
 \label{eq_x262}
 f+ \mathcal{S}[f] - (\mathcal S \circ \mathcal{P}^H)[f] =0.
 \end{equation}
 We rewrite the last identity as
 \begin{equation}
 \label{eq_x263}
 (\textnormal{Id}+\mathcal{S})\big[f-\mathcal P^H[f]\big]+\mathcal P^H[f] =0.
 \end{equation}
 By definition, $\mathcal P^H[f]$ is a function on $\gamma^+$ whose analytic continuation has the form $\frac{p(z)}{\sigma(z)}$ for a polynomial $p(z)$ of degree at most $H-1$. Let $\Delta f = f-\mathcal P^H[f]$. Since the integral kernel $S(z,w)$ has an analytic continuation to a holomorphic function of $z$ in $\amsmathbb{C} \setminus \amsmathbb{A}^{\mathfrak{m}}$ (one might need to make the contour $\gamma_h^-$ smaller to see that), by \eqref{eq_x262}, so does $f$, and therefore so does $\Delta f$. We denote these analytic continuations by the same letters $f(z)$ and $\Delta f(z)$, $z\in \amsmathbb{C} \setminus \amsmathbb{A}^{\mathfrak{m}}$. The kernel $S(z,w)$ decays as $O(\frac{1}{z})$ as $z \rightarrow \infty$. Hence, so do both $f(z)$ and $\Delta f(z)$. On the other hand, $\Delta f =f-\mathcal P^H[f]$ implies that all integrals of $\Delta f$ around $\gamma_h$ vanish. We conclude that $\Delta f$ belongs to $\mathscr{H}^{[2],0}$. Therefore, we can apply \eqref{eq_x238} to rewrite \eqref{eq_x263} as
 \begin{equation}
 \Delta f(z) - \frac{\mathcal{P}[\Delta f](z)}{\sigma(z)} +
\sum_{h = 1}^H (\mathcal{K}_{h}\circ \mathcal{U}_{h})[\Delta f](z) +\mathcal P^H[f] =0.
 \end{equation}
 The last identity together with Lemma~\ref{Lemma_eq_zero_solution} implies $\Delta f=0$. Then \eqref{eq_x263} implies $\mathcal{P}^H[f]=0$ and consequently $f=\Delta f+\mathcal{P}^H[\Delta f]=0$. This finishes the proof of invertibility of the operator $\textnormal{Id} + \mathcal{S} \circ (\textnormal{Id} - \mathcal{P}^H)$.

 \medskip

 Let $K(z,w)$ denote the kernel of the integral operator $\mathcal{S} \circ (\textnormal{Id} - \mathcal{P}^H)$. This kernel is a smooth function of $z$, $w$, interaction matrix $\boldsymbol{\Theta}$, and vectors of endpoints $\boldsymbol{\alpha},\boldsymbol{\beta}$. Writing\label{index:resolvent} $(\textnormal{Id} +\mathcal{S}-\mathcal{S}\circ \mathcal{P}^H)^{-1}=\textnormal{Id} - \mathcal{R}$, the kernel $R(z,w)$ of the integral operator $\mathcal{R}$ can be expressed as (see, \textit{e.g.}, \cite[Theorem 9.3]{gohberg2012traces})
\begin{equation}
\label{RSSSS}
\begin{split}
& R(z,w) \\
& = \frac{1}{\textnormal{Det}(\textnormal{Id} + \mathcal{S}-\mathcal{S} \circ \mathcal{P}^H)} \sum_{n \geq 0} \frac{1}{n!}\,\int_{(\gamma^+)^n} \det \left[\begin{array}{cccc} K(z,w) & K(z,\zeta_1) & \cdots & K(z,\zeta_n) \\ K(\zeta_1,w) & K(\zeta_{1},\zeta_{1}) & \cdots & K(\zeta_1,\zeta_n) \\ \vdots & \vdots & & \vdots \\ K(\zeta_n,w) & K(\zeta_n,\zeta_1) & \cdots & K(\zeta_n,\zeta_n) \end{array}\right] \prod_{i = 1}^{n} \frac{\dd \zeta_i}{2\ii\pi},
\end{split}
\end{equation}
where the Fredholm determinant is
\begin{equation}
\label{RSSSSF} \textnormal{Det}(\textnormal{Id} + \mathcal{S}-\mathcal{S} \circ \mathcal{P}^H) = \sum_{n \geq 0} \frac{1}{n!} \oint_{(\gamma^+)^{n}} \det_{1 \leq i,j \leq n} \big([K(\zeta_i,\zeta_j)\bigr] \prod_{i = 1}^n \frac{\dd \zeta_i}{2\ii\pi}.
\end{equation}
The series \eqref{RSSSS}-\eqref{RSSSSF} are absolutely convergent by application of the Hadamard inequality, first as
\begin{equation*}
\bigg|\oint_{\gamma^n} \det_{1 \leq i,j \leq n} \big([K(\zeta_i,\zeta_j)\bigr] \prod_{i = 1}^{n} \frac{\dd \zeta_i}{2\ii\pi} \bigg| \leq n^{\frac{n}{2}}\bigg(\frac{\textnormal{l}(\gamma^+)}{2\pi}\,|\!| K |\!|_{\infty}\bigg)^{n},
\end{equation*}
where $\textnormal{l}(\gamma^+)$ is the length of the contour, and then as
\begin{equation*}
\left|\oint_{\gamma^n} \det\left[\begin{array}{cccc} K(z,w) & K(z,\zeta_1) & \cdots & K(z,\zeta_n) \\ K(\xi_1,w) & K(\zeta_{1},\zeta_{1}) & \cdots & K(\zeta_1,\zeta_n) \\ \vdots & \vdots & & \vdots \\ K(\zeta_n,w) & K(\zeta_n,\zeta_1) & \cdots & K(\zeta_n,\zeta_n) \end{array}\right]\,\prod_{i = 1}^n \frac{\dd \zeta_i}{2\ii\pi}\right| \leq (n + 1)^{\frac{n + 1}{2}}\bigg( \frac{\textnormal{l}(\gamma^+)}{2\pi}\bigg)^{n}\big( |\!|K|\!|_{\infty}\big)^{n + 1}.
\end{equation*}

Since the kernel $R(z,w)$ is a smooth function of $z$, $w$, $\boldsymbol{\Theta}$,$\boldsymbol{\alpha},\boldsymbol{\beta}$, such that each partial derivative of finite order is bounded on $(\gamma^+)^2$, by the dominated convergence theorem the same smoothness is inherited by $R(z,w)$. Thus the operator norm of $\textnormal{Id} - \mathcal{R}$ is uniformly bounded, with operator norm depending only on the constants of the assumptions. \end{proof}

\bigskip

\noindent \textsc{Part 6.} We now come back to our problem in the form \eqref{Idequ} and finish the proof of Theorem~\ref{Theorem_Master_equation_12} for the case $\boldsymbol{\kappa} = 0$. We fix the domains $(\amsmathbb{K}_h)_{h=1}^H$ and contours $(\gamma_h)_{h=1}^H$ for \eqref{ContinuityUp_eq} and apply Lemma~\ref{Lemma_invertibility_Fredholm} on the contours $\gamma^+_h$ set to be equal to $\gamma_h$ from Theorem~\ref{Theorem_Master_equation_12}. This yields the existence of the unique solution $F_h(z)$ on the contour $\gamma_h$, as well as its smoothness and continuous dependence on $E_h$ in the supremum norm (because of the kernel of $\mathcal R$ being uniformly bounded
and smooth).

In order to extend the smoothness and bound \eqref{ContinuityUp_eq} from $\gamma_h$ to the set $\amsmathbb{K}_h$, we use the Cauchy integral formula:
\[
F_h(z) = \frac{1}{2\ii\pi} \oint_{\gamma_h} \frac{F_h(\zeta)}{z-\zeta} \dd\zeta, \qquad z\in \amsmathbb{K}_h.
\]
Smoothness and bounds on the integrand imply the same for the integral, which completes the proof of Theorem~\ref{Theorem_Master_equation_12}.

\bigskip

\noindent \textsc{Part 7.} It remains to turn on arbitrary $\boldsymbol{\kappa} \in \amsmathbb{C}^H$. In other words we want to solve \eqref{Idequ} for $F \in \mathscr{H}^{[1]}$ under the conditions
\begin{equation}
\label{periodkap}
\forall h \in [H]\qquad\oint_{\gamma_h} \frac{\dd z}{2\ii\pi}\,F(z) = \kappa_h.
\end{equation}
The latter implies that $\widetilde{F} = \big(F - \sigma^{-1} \cdot \Pi^{-1}[\boldsymbol{\kappa}]\big)$ belongs to the space $\mathscr{H}^{[2],0}$. We can then rewrite \eqref{Idequ} solely in terms of $\widetilde{F}$. Observing that
\[
\mathcal{P}\big[\sigma^{-1}\cdot \Pi^{-1}[\boldsymbol{\kappa}]\big] = \Pi^{-1}[\boldsymbol{\kappa}],
\]
we deduce that
\[
(\textnormal{Id} - \sigma^{-1} \cdot \mathcal{P})\big[\sigma^{-1}\cdot \Pi^{-1}[\boldsymbol{\kappa}]\big] = 0,
\]
and thus we obtain the equation
\[
\bigg(\textnormal{Id} - \sigma^{-1}\cdot \mathcal{P} + \sum_{h = 1}^H \mathcal{K}_{h} \circ \mathcal{U}_{h} \bigg)[\widetilde F](z) = \mathcal{K}[\boldsymbol{E}](z) - \sum_{h = 1}^{H} \mathcal{K}_h \circ \mathcal{U}_h\big[\sigma^{-1}\cdot\Pi^{-1}[\boldsymbol{\kappa}]\big](z).
\]
This is the same \eqref{Idequ}, but with a modified right-hand side. Hence, Part 6 implies existence of the unique $\widetilde{F}$ and its smooth dependence on all the involved parameters. Writing $F=\widetilde F + \sigma^{-1} \cdot \Pi^{-1}[\boldsymbol{\kappa}]$ we get the same statements for $F$.

\bigskip

We have established the existence of $ \Op_h\big[\boldsymbol{E}\,;\,\boldsymbol{\kappa}\big]$ operator, its smoothness and the bound \eqref{ContinuityUp_eq}. For the remaining statements of Theorem~\ref{Theorem_Master_equation_12}, note that the contours $\gamma_h$ can be moved, and therefore the weak dependence of $\boldsymbol{\Upsilon}$ in the segments $(\amsmathbb{A}_h^{\mathfrak{m}})_{h = 1}^{H}$ and the complex domains $(\amsmathbb{M}_h)_{h = 1}^{H}$ stated after Theorem~\ref{Theorem_Master_equation_12} is clear.

Finally, if $E_h$ is holomorphic in $\amsmathbb{M}_h$ (instead of $\amsmathbb{M}_h^o = \amsmathbb{M}_h \setminus \amsmathbb{A}_h^{\mathfrak{m}}$), looking at the formula for $\mathcal{K}_h$ in Definition~\ref{defofop}, we see that $\mathcal{K}_{h}[E_h] = 0$ because $\sigma(\zeta)/\sqrt{(\zeta - \alpha_h)(\zeta - \beta_h)}$ is holomorphic in $\amsmathbb{M}_h$. In this case, the right-hand side in \eqref{Idequ} is zero, and therefore by Corollary~\ref{Corollary_Inj}, the function $F$ is also zero.
This justifies the last property \eqref{eq_zero_solution_for_holomorphic} stated in Theorem~\ref{Theorem_Master_equation_12} and concludes the proof.

\section{From the fundamental solution to tame solutions}

\label{Section_Master_by_covariance}

In this section we develop alternative points of view on $\Op_h$ which solve the master problem of Definition~\ref{def:master}. First, we recast $\Op_h$ as a solution operator for a Riemann--Hilbert problem. Second, we express it in terms of the fundamental solution, which is similar in spirit to writing solutions to differential equations in terms of the corresponding Green function. The main result of this section is Theorem~\ref{thm:genRHPsol}.

\subsection{The master problem as a Riemann--Hilbert problem}

We are often interested in solutions of the master problem for sources $\boldsymbol{E}$ enjoying stronger properties.

\begin{definition}
\label{def:tame}
We say that $\boldsymbol{F}(z) = (F_h(z))_{h = 1}^H$ is \emph{tame} if for each $h \in [H]$ there exists a complex neighborhood $\amsmathbb{M}_h$ of $[\alpha_h,\beta_h]$ such that
\begin{itemize}
\item $F_h(z)$ is a holomorphic function of $z \in \amsmathbb{M}_h \setminus [\alpha_h,\beta_h]$;
\item $F_h(x^{\pm}) = \lim_{\epsilon \rightarrow 0^+} F_h(x \pm \ii \epsilon)$ exists for any $x \in (\alpha_h,\beta_h)$;
\item there exists an integer $m_{\boldsymbol{F}}$, called tameness exponent, and a real number $M_{\boldsymbol{F}}$ such that
\[
\sup_{z \in \amsmathbb{M}_h \setminus [\alpha_h,\beta_h]} \big| \big((z - \alpha_h)(z - \beta_h)\big)^{\frac{1}{2}m_{\boldsymbol{F}}} F_h(z)\big| \leq M_{\boldsymbol{F}};
\]
\item the function $x \mapsto F_h(x^+) + F_h(x^-)$ extends to a meromorphic function in $\amsmathbb{M}_h$.
\end{itemize}
\end{definition}
Let $\sigma_h$ be the holomorphic function in $\amsmathbb{C} \setminus [\alpha_h,\beta_h]$ such that
\begin{equation}
\label{sigmh}
(\sigma_h(z))^2 = (z - \alpha_h)(z - \beta_h),\qquad \sigma_h(z) \mathop{\sim}_{z \rightarrow \infty} z.
\end{equation}
Tameness exponents are not unique, but the optimal choice of exponent is irrelevant for our purposes. We are interested in the information they give about the worst possible divergence, not assuming that the divergence occurs exactly with this power.

For later use, we remark that since $\sigma_h(x^+) = - \sigma_h(x^-)$ for any $x \in (\alpha_h,\beta_h)$, tameness of $\boldsymbol{F}(z)/\sigma(z)$ is equivalent to the first three properties for $\boldsymbol{F}$ and existence of a meromorphic extension of $(F_h(x^+) - F_h(x^-))/\sigma_h(x^+)$ in a complex neighborhood of $[\alpha_h,\beta_h]$.

\begin{proposition}
\label{prop:tameRHP} Adopt the assumptions of Theorem~\ref{Theorem_Master_equation_12}, assume that $\boldsymbol{E}(z)/\sigma(z)$ is tame, and let $\boldsymbol{\kappa}$ vary in a fixed compact of $\amsmathbb{C}^H$.

Then the solution of the master problem $\boldsymbol{F}(z) = (\Upsilon_h[\boldsymbol{\boldsymbol{E}}\,;\,\boldsymbol{\kappa}](z))_{h = 1}^H$ is tame, with tameness exponent $m_{\boldsymbol{F}} = \max(1,m_{\boldsymbol{E}/\sigma})$ and constant $M_{\boldsymbol{F}}$ depending only on $H$ and $C$ in the assumptions. Besides, we have for any $h \in [H]$
\begin{equation}
\label{RHPeq} \forall x \in (\alpha_h,\beta_h) \qquad \theta_{h,h}\big(F_h(x^+) + F_h(x^-)\big) + \sum_{g \neq h} 2\theta_{h,g} F_g(x) = -\frac{E_h(x^+) - E_h(x^-)}{\sigma_h(x^+)}.
\end{equation}

Conversely, assume that $\boldsymbol{F}(z)$ is a tame solution of \eqref{RHPeq}, that both $\boldsymbol{F}$ and $\boldsymbol{E}/\sigma$ have tameness exponent $1$ and that for each $h \in [H]$ the function $F_h(z)$ is holomorphic in $\amsmathbb{C}\setminus [\alpha_h,\beta_h]$ and satisfies $F_h(z) \sim \frac{\kappa_h}{z}$ as $z \rightarrow \infty$. Then $\boldsymbol{F}$ solves the master problem, \textit{i.e.} $F_h(z) = \Upsilon_h[\boldsymbol{E}\,;\,\boldsymbol{\kappa}](z)$ for any $h \in [H]$.
\end{proposition}
\begin{proof}
We rewrite \eqref{rewritingEF} as

\begin{equation}
\begin{split}
\label{eq_x277}
 F(z) & =  \frac{\mathcal{P}[F](z)}{\sigma(z)}
 - \sum_{1 \leq g \neq h \leq H} \frac{1}{\sigma(z)} \oint_{\gamma_h} \frac{\dd\zeta}{2\ii\pi} \frac{\sigma(\zeta)}{\zeta - z} \frac{(\theta_{h,h} - \theta_{h,g})F_g(\zeta)}{\theta_{h,h}} \\ & \quad + \sum_{h = 1}^{H} \frac{1}{\sigma(z)}\oint_{\gamma_h} \frac{\dd\zeta}{2\ii\pi}\,\frac{\sigma(\zeta)}{\zeta - z} \,\frac{E_h(\zeta)}{\theta_{h,h}\sqrt{(\zeta - \alpha_h)(\zeta - \beta_h)}},
\end{split}
\end{equation}
where $z$ is outside the integration contour, and analyze the right-hand side. We are interested in the boundary values of $F(z)$ on the real axis, and for that need to move $z$ through the integration contour both in the double sum and in the single sum in \eqref{eq_x277}, picking up the residue at $\zeta=z$ in the process. Since $\boldsymbol{E}(z)/\sigma(z)$ is tame, this residue is also tame with the same exponent. The integrals themselves in \eqref{eq_x277} are tame with exponent $1$ because of the $\frac{1}{\sigma(z)}$ prefactors. Tameness of $F(z)$ implies tameness for $F_h(z) = \Upsilon_h[\boldsymbol{E}\,;\,\boldsymbol{\kappa}](z)$ of Theorem~\ref{Theorem_Master_equation_12} through \eqref{FFhsum} and \eqref{Fh_projection}.

To get more information on the boundary values of $F_h(z)$, we return to the master equation
\[
\sigma_h(z) \left( \sum_{g = 1}^H \theta_{h,g} F_g(z)\right) + E_h(z) = A_h(z),\qquad \sigma_h(z) = \sqrt{(z - \alpha_h)(z - \beta_h)}.
\]
We divide by $\sigma_h(z)$ and evaluate the equation at $z = x^{\pm}$ for $x \in (\alpha_h,\beta_h)$. Taking the difference of the two evaluations and using the fact that $\sigma_h(x^{+}) = - \sigma_h(x^-)$ while $A_h$ and $F_g$ for $g \neq h$ are holomorphic in a complex neighborhood of $[\alpha_h,\beta_h]$, we get \eqref{RHPeq}.

\medskip

For the converse part, take a solution of \eqref{RHPeq} with the announced properties. By multiplying \eqref{RHPeq} with $\sigma_h(x^+)$, we get for any $x \in (\alpha_h,\beta_h)$
\[
(\sigma_h \cdot F_h)(x^+) - (\sigma_h \cdot F_h)(x^-) + \sum_{g \neq h} \theta_{h,g} F_g(x)\big(\sigma_h(x^+) - \sigma_h(x^-)\big) = -\big(E_h(x^+) - E_h(x^-)\big).
\]
This is equivalent to
\begin{equation}
\label{Aeqn}
A_h(x^+) = A_h(x^-),\qquad \textnormal{where}\qquad A_h(z) = \sigma_h(z)\left(\sum_{g = 1}^{H} \theta_{h,g} F_g(z)\right) + E_h(z).
\end{equation}
The first property in the tameness assumptions for $\boldsymbol{E}/\sigma$ and $\boldsymbol{F}$ implies that $A_h$ is defined and holomorphic in $\amsmathbb{M}_h \setminus [\alpha_h,\beta_h]$. In conjunction with \eqref{Aeqn} this shows (say, by Morera theorem) that $A_h$ is holomorphic in $\amsmathbb{M}_h \setminus \{\alpha_h,\beta_h\}$. Having tameness exponents $1$ for $\boldsymbol{E}/\sigma$ and $\boldsymbol{F}$ implies that $A_h$ is bounded near $\alpha_h$ and $\beta_h$, therefore $A_h$ is holomorphic in full set $\amsmathbb{M}_h$ and we are in the setting of the master problem of Definition~\ref{def:master}.
\end{proof}

If we rather assume that the (optimal) tameness constants are greater than $1$ in the last part of Proposition~\ref{prop:tameRHP}, then $\boldsymbol{F}$ would be a solution of a variant of the master problem, namely one where the unknown function $A_h$ is meromorphic with possible poles at $\alpha_h$ and $\beta_h$, whose order is controlled by the tameness exponents. This meromorphic master problem in general does not have a unique solution. In fact, our uniqueness proof (essentially coming from Lemma~\ref{Lemma_eq_zero_solution}) relied on the construction of integrable measures $\nu_h$ on $[\alpha_h,\beta_h]$ whose Stieltjes transform is $F_h$. This is possible when $F_h$ diverges as inverse square roots at $\alpha_h,\beta_h$, but not if $F_h$ diverges with higher powers of the inverse square root. Solutions of this meromorphic variant of the master problem depending on the singular behavior near $\alpha_h,\beta_h$ will be discussed later in Theorem~\ref{thm:genRHPsol}.

Variants of Proposition~\ref{prop:tameRHP} with weaker regularity assumptions than 'tameness' could be formulated and proved along the same line. We restrict ourselves to an observation that the equilibrium measure of discrete ensembles can be reconstructed from the solution of the master Riemann--Hilbert problem with a right-hand side that may have logarithmic singularities.

\begin{lemma}
\label{Lem:app_RHP_mu}
Consider a discrete ensemble as described in Section~\ref{Section_general_model}, satisfying Assumptions~\ref{Assumptions_Theta}, \ref{Assumptions_basic} and \ref{Assumptions_analyticity}. This means in particular that the derivative of the potential in the $h$-th segment is
\[
V_h'(x) = U_h'(x) + \iota_h^- \log|x - \hat{a}_h'| + \iota_h^+ \log |x - \hat{b}_h'|,
\]
where $U_h$ has a holomorphic extension in a complex neighborhood of $[\hat{a}_h',\hat{b}_h']$. Then, the Stieltjes transform of the equilibrium measure provided by Theorem~\ref{Theorem_equi_charact} satisfies for any $h \in [H]$ and $x$ in a band of the $h$-th segment
\[
\theta_{h,h}\big(\mathcal{G}_{\mu_h}(x^+) + \mathcal{G}_{\mu_h}(x^-)\big) + \sum_{g \neq h} 2\theta_{h,g} \mathcal{G}_{\mu_g}(x) = V'_h(x).
\]
\end{lemma}
\begin{proof}[Sketch of the proof]
Fix $h \in [H]$. By Theorem~\ref{Theorem_equi_charact} the effective potential in the $h$-th segment
\[
V_h^{\textnormal{eff}}(x) = V_h(x) - 2 \sum_{g = 1}^H \theta_{h,g} \int_{\hat{a}_g'}^{\hat{b}_g'} \log|x - y|\mu_g(y)\dd y
\]
is constant equal to $v_h$ for $x$ in a band. Differentiating with respect to $x$ (we omit the justification) we get
\[
0 = V_h'(x) - 2 \sum_{g = 1}^H \theta_{h,g}\,\textnormal{p.v.} \int_{\hat{a}_g'}^{\hat{b}_g'} \frac{\mu_g(y) \dd y}{x - y}.
\]
The Cauchy principal value integral in the $g$-th term is equal to $\frac{1}{2}\big(\mathcal{G}_{\mu_g}(x^+) + \mathcal{G}_{\mu_g}(x^-)\big)$. For $g \neq h$ the measure $\mu_g$ has support disjoint from the $h$-th segment so the integral is also equal to $\mathcal{G}_{\mu_g}(x)$.
\end{proof}

\subsection{Master Riemann--Hilbert problem and fundamental solution}
\label{sec:RHPfun}
Proposition~\ref{prop:tameRHP} motivates a study of the Riemann--Hilbert problem \eqref{RHPeq}. The setting requires the data of
\begin{itemize}
\item[(i)] a positive integer $H$ and a real symmetric matrix $\boldsymbol{\Theta}$ of size $H$ which is positive semi-definite and has positive diagonal entries;
\item[(ii)] a $2H$-tuple of real numbers $\alpha_1 < \beta_1 < \alpha_2 < \cdots < \alpha_H < \beta_H$;
\item[(iii)] an $H$-tuple of functions $\boldsymbol{D}(z) = (D_h(z))_{h = 1}^{H}$, called source, such that $D_h$ is meromorphic in a complex neighborhood of $[\alpha_h,\beta_h]$ for any $h \in [H]$.
\end{itemize}

\begin{definition}
\label{def:masterRHP} The master Riemann--Hilbert problem asks to find an $H$-tuple of functions $\boldsymbol{F}(z) = (F_h(z))_{h = 1}^H$ with the following properties for any $h \in [H]$.
\begin{enumerate}
\item $F_h(z)$ is a meromorphic function of $z \in \amsmathbb{C} \setminus [\alpha_h,\beta_h]$.
\item $F_h(x^{\pm}) = \lim_{\epsilon \rightarrow 0^+} F_h(x \pm \ii \epsilon)$ exist for any $x \in (\alpha_h,\beta_h)$.
\item $F_h(z) = O(\frac{1}{z})$ as $z \rightarrow \infty$.
\item $\forall x \in (\alpha_h,\beta_h) \qquad \theta_{h,h}\big(F_h(x^+) + F_h(x^-)\big) + \sum_{g \neq h} 2\theta_{h,g}F_g(x) = D_h(x)$.
\end{enumerate}
\end{definition}

We are going to construct all solutions of the master Riemann--Hilbert problem from a fundamental solution. This fundamental solution already appeared in the context of the master problem as it was giving the leading covariance in the discrete ensembles (Theorem~\ref{Theorem_correlators_expansion}).

\begin{definition}
\label{def:Berg}The fundamental solution of the master problem is the $H^2$-tuple of functions $\boldsymbol{\mathcal{F}}(z_1,z_2) = \big(\mathcal{F}_{h_1,h_2}(z_1,z_2)\big)_{h_1,h_2 = 1}^{H}$ given in terms of the operator \eqref{eq_solution_Master_equation} by
\begin{equation}
\label{eq_covariance_11}
\forall h_1,h_2 \in [H] \qquad \mathcal{F}_{h_1,h_2}(z_1,z_2) =
\Op_{h_1}\left[\frac{\sqrt{(* - \alpha_{h_2})(* - \beta_{h_2})}}{2(*- z_2)^2}\boldsymbol{e}^{(h_2)}\,;\,\boldsymbol{\kappa} = \boldsymbol{0}\right](z_1),
\end{equation}
where $*$ is a placeholder for the variable in which this tuple is seen as a function on which the operator $\Upsilon_{h_1}$ is applied, and $\boldsymbol{e}^{(1)},\ldots,\boldsymbol{e}^{(H)}$ is the canonical basis of $\amsmathbb{R}^H$.
\end{definition}

A non-trivial property of the fundamental solution is the symmetry in its two variables.

\begin{theorem}
\label{thm:Bsym}
 The fundamental solution $\boldsymbol{\mathcal{F}}(z_1,z_2)$ depends only on $H$, $\boldsymbol{\Theta}$ and $\alpha_1 < \beta_1 < \cdots < \alpha_H < \beta_H$. It does so in a smooth way, and satisfies
 \begin{equation}
 \label{eq_Fundamental_symmetry}
 \mathcal{F}_{h_1,h_2}(z_1,z_2) = \mathcal{F}_{h_2,h_1}(z_2,z_1).
 \end{equation}
 for any $h_i \in [H]$ and $z_i \in \widehat{\amsmathbb{C}} \setminus [\alpha_{h_i},\beta_{h_i}]$ with $i = 1,2$. In particular, $\mathcal{F}_{h_1,h_2}(z_1,z_2)$ is a holomorphic function of $(z_1,z_2) \in \big({\amsmathbb{C}} \setminus [\alpha_{h_1},\beta_{h_1}]\big) \times \big({\amsmathbb{C}} \setminus [\alpha_{h_2},\beta_{h_2}]\big)$.
\end{theorem}
This property is important and will be used many times in the forthcoming arguments. The symmetry should have a direct analytic proof, perhaps similar to the proof of symmetry of Green function, but we have not been able to find one. Instead, we offer below an argument based on approximation by the (by definition symmetric) covariance in a statistical mechanics ensemble. There is some freedom in the choice of the approximating ensemble, and we in fact give two different proofs. The first one is based on approximation by discrete ensembles and uses Theorem~\ref{Theorem_correlators_expansion_relaxed} saying that the symmetric covariance ${W}_{2; h_1,h_2}(z_1,z_2)$ approaches the fundamental solution $\mathcal{F}_{h_1,h_2}(z_1,z_2)$ as $\N \rightarrow \infty$. Note that while we relied on the properties of the operator $\Op_h$ of Theorem~\ref{Theorem_Master_equation_12} in Chapter~\ref{Chapter_fff_expansions}, we never used the symmetry of $\mathcal{F}_{h_1,h_2}(z_1,z_2)$. Hence, the arguments in this first proof are not cyclic. The second proof instead uses an approximation by continuous $\sbeta$-ensembles and results about their covariance at leading order existing in the literature.

\begin{proof}[Proof of Theorem~\ref{thm:Bsym}] The analytic properties of the source and the properties stated in Theorem~\ref{Theorem_Master_equation_12} imply that $\boldsymbol{\mathcal{F}}(z_1,z_2)$ only depends on $\boldsymbol{\alpha},\boldsymbol{\beta},\boldsymbol{\Theta}$ and it does so in a smooth way. After we have justified symmetry, the last statement of the theorem comes from the following argument. For fixed $z_2$, $\mathcal{F}_{h_1,h_2}(z_1,z_2)$ is holomorphic with respect to $z_1 \in \widehat{\amsmathbb{C}} \setminus [\alpha_{h_1},\beta_{h_1}]$ because it is the solution of the master problem. For fixed $z_1$, it is also holomorphic with respect to $z_2 \in \widehat{\amsmathbb{C}} \setminus [\alpha_{h_2},\beta_{h_2}]$. By Hartogs' theorem on separate analyticity, this implies holomorphicity in $(z_1,z_2)$. In the rest of the proof we show \eqref{eq_Fundamental_symmetry}.

\medskip

\noindent \textsc{First proof.} We apply the approximation strategy, by constructing a discrete ensemble with deterministically fixed filling fractions, satisfying Assumptions~\ref{Assumptions_Theta}, \ref{Assumptions_basic}, \ref{Assumptions_offcrit}, and \ref{Assumptions_analyticity} and whose equilibrium measure has the desired $\N$-independent bands $(\alpha_h,\beta_h)$. Then, we will apply \eqref{eq_x81_2} from Theorem~\ref{Theorem_correlators_expansion_relaxed} together with the obvious symmetry \begin{equation}
\label{eq_x278}
 W_{2;h_1,h_2}(z_1,z_2)=W_{2;h_2,h_1}(z_2,z_1),
\end{equation}
coming from the definition \eqref{eq_correlators_def} of $W_{2;h_1,h_2}(z_1,z_2)$ as the covariance of $G_{h_1}(z_1)$ and $G_{h_2}(z_2)$. Note that the filling fractions $\hat{n}_h$ are not important, as $\boldsymbol{\mathcal{F}}(z_1,z_2)$ does not depend on them but only on the prescribed position of bands, and therefore we can choose them at our convenience.

When $H=1$, the desired discrete ensemble is delivered by the $zw$-measures of Chapter~\ref{Chapterzw}.
Proposition~\ref{Proposition_zw_band_and_type_1} guarantees that $zw$-discrete ensemble can have an arbitrary band $(\alpha_1,\beta_1)$; note that the required inequalities in the statement of the proposition do not restrict us, because we are free to choose any positive $\hat n$.

For general $H$ and $\boldsymbol{\Theta}$ we start from a collection of $H$ independent $zw$-discrete ensembles, so that their bands match $(\alpha_1,\beta_1)$, $(\alpha_2,\beta_2),\ldots,(\alpha_H,\beta_H)$ and with $\theta=\theta_{h,h}$ for the $h$-th ensemble. Since by Proposition~\ref{Proposition_zw_band_and_type_1}, the segments $[\hat A_1',\hat B_1']$ for the ensembles can be chosen arbitrary close to the bands, these $zw$-discrete ensembles can be chosen to live on pairwise disjoint segments denoted $[\hat a_h,\hat b_h]$ for $h\in [H]$. Let $w_h(x)$ denote the weight for the $h$-th $zw$-discrete ensemble, $V_h(x)$ the potential, $\mu_h$ the equilibrium measure, $\hat n_h$ the filling fraction, \textit{i.e.} the total mass of $\mu_h$. If $\boldsymbol{\Theta}$ is diagonal, then we are done, as the $H$ independent $zw$-discrete ensembles combine into a single discrete ensemble on $H$ segments with the equilibrium measures $\boldsymbol{\mu}=(\mu_h)_{h=1}^H$ and the desired bands $\{(\alpha_h,\beta_h)\}_{h=1}^H$. Let us refer to the last discrete ensemble as ensemble $X$.

If $\boldsymbol{\Theta}$ is not diagonal, then we need to modify the potentials $V_h(x)$ to get the same equilibrium measures for an ensemble with interactions given by $\boldsymbol{\Theta}$. Such modifications were previously used in the proof of Theorem~\ref{Theorem_partition_multicut}). Namely, we take disjoint complex neighborhoods $\amsmathbb{M}_h$ of $[\hat a_h,\hat b_h]$ for $h\in[H]$ and define the holomorphic functions
\begin{equation}
\forall h \in [H] \quad \forall z \in \amsmathbb{M}_h\qquad \Delta V_h(z)=\sum_{g \neq h} 2\theta_{h,g} \int_{\hat a'_{g}}^{\hat b'_{g}}
 \log\big(((x-y)\textnormal{sgn}(h - g)\big)\mu_{g}(y) \dd y.
\end{equation}
We define a discrete ensemble $Y$ on the same segments $[\hat a_h,\hat b_h]$, $h\in [H]$, and with the same filling fractions $\hat n_h$ as the ensemble $X$ by modifying the potential and its regular part:
\begin{equation*}
\begin{split}
 w_h^{Y}(x) & = w_h(x) \cdot \exp\left(-\N \Delta V_h\bigg(\frac{x}{\N}\bigg)\right), \\
 V_h^{Y}(x) & = V_h(x) + \Delta V_h(x), \\
 U_h^{Y}(x) & = U_h(x) + \Delta V_h(x).
\end{split}
\end{equation*}
The auxiliary functions for Assumption~\ref{Assumptions_analyticity} follow from these expressions by Definition~\ref{Definition_phi_functions}.
We claim that the equilibrium measure of $Y$
coincides with $\boldsymbol{\mu}$. This can be seen by applying Theorem~\ref{Theorem_equi_charact_repeat_2} and
noting that the effective potentials of Definition~\ref{def_eff_pot} computed with the measure $\boldsymbol{\mu}$ have the same characterizing properties for ensembles $X$ and $Y$. Indeed, the change in $\boldsymbol{\Theta}$ is exactly compensated by the change in $V_h$. Hence, the ensemble $Y$ is off-critical and satisfies Assumptions~\ref{Assumptions_Theta}, \ref{Assumptions_basic}, \ref{Assumptions_offcrit} and \ref{Assumptions_analyticity}, because $X$ satisfies them by Lemma~\ref{Checkzw} and Corollary~\ref{Corollary_ZW_off-critical_2}.

Now we can write using the covariance in the ensemble $Y$:
\begin{equation}
\label{eq_x279}
\begin{split}
 \left|\mathcal{F}_{h_1,h_2}(z_1,z_2) - \mathcal{F}_{h_2,h_1}(z_2,z_1)\right| & \leq \left|W_{2;h_1,h_2}(z_1,z_2)-W_{2;h_2,h_1}(z_2,z_1)\right| \\ &\quad +
 \left|\mathcal{F}_{h_1,h_2}(z_1,z_2) -W_{2;h_1,h_2}(z_1,z_2)\right| \\
 & \quad + \left|\mathcal{F}_{h_2,h_1}(z_2,z_1)-W_{2;h_2,h_1}(z_2,z_1)\right|.
 \end{split}
\end{equation}
The first term in the right-hand side of \eqref{eq_x279} is zero by \eqref{eq_x278}, while two other terms tend to $0$ as $\N \rightarrow \infty$ by Theorem~\ref{Theorem_correlators_expansion_relaxed}. Since the left-hand side of \eqref{eq_x279} does not depend on $\N$, we conclude that it must be equal to zero.

\medskip

\noindent \textsc{Second proof.} For each $h \in [H]$ let
\[
\rho_h(x) = \frac{2\sqrt{(x - \alpha_h)(\beta_h - x)}}{H\pi(\beta_h - \alpha_h)} \mathbbm{1}_{[\alpha_h,\beta_h]}(x)
\]
be the semi-circular density of mass $\frac{1}{H}$ with support $[\alpha_h,\beta_h]$. It is well-known (as the semi-circle law minimizes the energy functional $-\mathcal{I}$ defined in \eqref{eq_functional_general_repeat} when $H=1$ and $V_{h}$ is quadratic, under the sole condition that $\rho$ is a probability measure so that Theorem~\ref{Theorem_equi_charact} applies with no saturated regions)
 that there exists a constant $c_h$ such that
\begin{equation}
\label{8HHH}
\forall x \in \amsmathbb{R},\qquad \frac{8}{H(\beta_h - \alpha_h)^2} \left(x - \frac{\alpha_h + \beta_h}{2}\right)^2 - 2\int_{\alpha_h}^{\beta_h} \log|x - y|\rho_h(y)\dd y - c_h \geq 0,
\end{equation}
with equality if and only if $x \in [\alpha_h,\beta_h]$, and that the left-hand side of \eqref{8HHH} behaves likes $d_{h,-}\sqrt{\alpha_h - x}$ as $x \rightarrow \alpha_h^-$ and $d_{h,+}\sqrt{x - \beta_h}$ as $x \rightarrow \beta_h^-$ for some constants $d_{h,\pm} > 0$. We define an analytic function of $x$ in a real neighborhood of $[\alpha_h,\beta_h]$ by:
\begin{equation}
\label{Vhgauss}
V_h(x) = \frac{8\theta_{h,h}}{H(\beta_h - \alpha_h)^2} \left(x - \frac{\alpha_h + \beta_h}{2}\right)^2 + \sum_{g \neq h} 2\theta_{h,g}\int_{\alpha_g}^{\beta_g} \log|x - y|\rho_g(y) \dd y.
\end{equation}
Then $V_h^{\textnormal{eff}}(x) := V_h(x) - \sum_{g = 1}^{H} 2\theta_{h,g} \int_{\alpha_g}^{\beta_g} \log|x - y| \rho_g(y)\dd y -\theta_{h,h}c_{h}$ is nonnegative everywhere, equal to zero on $[\alpha_h,\beta_h]$ and is positive and behaves like a squareroot to the left of $\alpha_h$ and to the right of $\beta_h$. This shows that $(\rho_h(x)\dd x)_{h =1}^{H}$ is the unrestricted equilibrium measure associated to a variational datum (\textit{cf.} Remark~\ref{remark:continuous}), taking as segments small pairwise disjoint neighborhoods $[\hat{a}_h',\hat{b}_h']$ of $[\alpha_h,\beta_h]$. The main point of this construction is that this unrestricted equilibrium measure satisfies the off-criticality Assumption~\ref{Assumption_C}.

We introduce the statistical mechanics ensemble of $H$ groups of $M$ particles each, with positions $\boldsymbol{\ell}$ that can vary continuously and are drawn according to the probability measure on $\amsmathbb{R}^{MH}$
\begin{equation}
\label{continuousPPP}
\dd\amsmathbb{P}_M(\boldsymbol{\ell}) = \frac{1}{Z} \prod_{\substack{1 \leq h \leq H \\ 1 \leq j \leq M}} \mathbbm{1}_{[\hat{a}_h',\hat{b}_h']}(\ell_j^h) e^{-MV_{h}(\ell_i^h)} \dd \ell_j^h \prod_{g,h}
\prod_{i \neq j } |\ell_i^g - \ell_j^h|^{\theta_{g,h}}.
\end{equation}
Call $\amsmathbb{E}_M$ the expectation value for the probability measure $\amsmathbb{P}_M$, and define the $n$-point correlators as in Definition~\ref{def_correlators}
\[
W_{n;h_1,\ldots,h_n}^{\amsmathbb{R}}(z_1,\ldots,z_n) = \amsmathbb{E}^{(\textnormal{c})}_M\left[\sum_{i_1 = 1}^{M} \frac{1}{z_1 - \ell_{i_1}^{h_1}},\ldots,\sum_{i_n = 1}^{H} \frac{1}{z_n - \ell_{i_n}^{h_n}}\right] .
\]

The equilibrium measure for this ensemble is by construction $(\rho_h(x)\dd x)_{h =1}^{H}$. The strict convexity of the functional $\mathcal{I}$ proved in Corollary~\ref{Corollary_I_positive}, the existence of an analytic continuation of the potential \eqref{Vhgauss} in a complex neighborhood of $[\hat{a}_h',\hat{b}'_h]$, and the aforementioned off-criticality allow applying the results of \cite{BGK}. This guarantees the existence of an asymptotic expansion as $M \rightarrow \infty$
\begin{equation}
\label{expaWr}
\begin{split}
W_{1;h}(z) & = M \mathcal{G}_{\rho_h}(z) + O(1) \\
W_{2;h_1,h_2}(z_1,z_2) & = \mathcal{F}^{\amsmathbb{R}}_{h_1,h_2}(z_1,z_2) + o(1), \\
W_{n;h_1,\ldots,h_n}(z_1,\ldots,z_n) & = o(1),
\end{split}
\end{equation}
which is uniform for $h_i \in [H]$ and $z_i$ in the complex plane away from $[\hat{a}_{h_i}',\hat{b}_{h_i}']$. The leading covariance $\mathcal{F}^{\amsmathbb{R}}_{h_1,h_2}(z_1,z_2)$ is independent of $M$ and \cite{BEO} shows that it satisfies the same master Riemann--Hilbert problem with respect to $z_1$ than the fundamental solution $\mathcal{F}_{h_1,h_2}(z_1,z_2)$ of Definition~\ref{def:Berg}. By uniqueness, we must have
\[
\mathcal{F}_{h_1,h_2}(z_1,z_2) = \mathcal{F}_{h_1,h_2}^{\amsmathbb{R}}(z_1,z_2) = \mathcal{F}_{h_2,h_1}^{\amsmathbb{R}}(z_2,z_1) = \mathcal{F}_{h_1,h_2}(z_1,z_2).
\]

For completeness, we sketch the argument of \cite{BEO} leading to the master Riemann--Hilbert problem. First, we write the Dyson--Schwinger equation of the ensemble $\amsmathbb{P}_M$.
\begin{equation}
\begin{split}
0 & = \frac{1}{M} \amsmathbb{E}_M\left[\sum_{g \neq h} \sum_{i,j = 1}^M \frac{\theta_{g,h_1}}{(\ell_i^h - \ell_j^g)(z - \ell_i^h)} +\sum_{i,j = 1}^M \frac{\theta_{h,h}}{(z - \ell_i^h)(z - \ell_j^h)} + \sum_{i = 1}^M \frac{1 - \theta_{h,h}}{(z - \ell_i^h)^2}\right] \\
& \quad - \amsmathbb{E}_M\left[\sum_{i = 1}^{M} \frac{V_{h}'(\ell_i^h)}{z - \ell_i^h}\right].
\end{split}
\end{equation}
Varying the potential (we already used this method in Section~\ref{higher} to derive higher-order Nekrasov equation from the first Nekrasov equation), we arrive to the two-variable Dyson--Schwinger equation
\begin{equation}
\begin{split}
0 & = \frac{1}{M} \amsmathbb{E}_M^{(\textnormal{c})}\left[\left(\sum_{g \neq h} \frac{\theta_{g,h}}{(\ell_i^h - \ell_j^g)(z - \ell_i^h)} + \sum_{i,j = 1}^M \frac{\theta_{h,h}}{(z - \ell_i^h)(z - \ell_j^h)} + \sum_{i = 1}^M \frac{1 - \theta_{h,h}}{(z - \ell_i^h)^2}\right)\,,\,\sum_{i_2 = 1}^{M} \frac{1}{z_2 - \ell_{i_2}^{h_2}}\right] \\
& \quad - \amsmathbb{E}_M^{(\textnormal{c})}\left[\sum_{i = 1}^{M} \frac{V_{h}'(\ell_i^h)}{z - \ell_i^h}\,,\,\sum_{i_2 = 1}^{M} \frac{1}{z_2 - \ell_{i_2}^{h_2}}\right] + \frac{1}{M} \amsmathbb{E}_M\left[\sum_{i = 1}^M \frac{\delta_{h,h_2}}{(z - \ell_i^h)(z_2 - \ell_{i}^{h})^2}\right].
\end{split}
\end{equation}
In the limit $M \rightarrow \infty$ and due to \eqref{expaWr}, the $(1 - \theta_{h,h})$-term disappears and other terms are expressed only in terms of the equilibrium measure and the leading covariance:
\begin{equation}
\label{2ndSDeqnR}\begin{split}
0 & = \oint_{\gamma_h}  \frac{\dd \zeta}{2\ii\pi (z - \zeta)} \left(\sum_{g \neq h} \theta_{g,h} \oint_{\gamma_g} \frac{\dd \tilde{\zeta}}{2\ii\pi} \frac{G_{\rho_h}(\zeta) \mathcal{F}^{\amsmathbb{R}}_{g,h_2}(\tilde{\zeta},z_2) + G_{\rho_{g}}(\tilde{\zeta})\mathcal{F}^{\amsmathbb{R}}_{h,h_2}(\zeta,z_2)}{\zeta - \tilde{\zeta}}\right) \\
& \quad + 2\theta_{h,h} G_{\rho_h}(z) \mathcal{F}^{\amsmathbb{R}}_{h,h_2}(z,z_2) - \oint_{\gamma_h} \frac{\dd\zeta}{2\ii\pi}\,\frac{V'_h(\zeta)\mathcal{F}^{\amsmathbb{R}}_{h,h_2}(\zeta,z_2)}{z - \zeta} + \delta_{h,h_2} \partial_{z_2}\left(\frac{G_{\rho_h}(z) - G_{\rho_h}(z_2)}{z - z_2}\right),
\end{split}
\end{equation}
where $z,z_2$ are outside the integration contours. It is not hard to see from this equation that $(\mathcal{F}^{\amsmathbb{R}}_{h_1,h_2}(z_1,z_2))_{h_1 =1}^H$ is tame with respect to $z_1$ for each $h_2,z_2$. To obtain the master Riemann--Hilbert problem, for $x \in (\alpha_h,\beta_h)$ we send $z \rightarrow x^{+}$ or $z \rightarrow x^-$ in \eqref{2ndSDeqnR} and compute the difference. For this purpose we observe that
\begin{equation}
\begin{split}
& \quad 2\big(F_1(x^+)F_2(x^+) - 2F_1(x^-)F_2(x^-)\big) \\
& = \big(F_1(x^+) - F_1(x^-)\big)\big(F_2(x^+) + F_2(x^-)\big) - \big(F_1(x^+) + F_1(x^-)\big)\big(F_2(x^+) - F_2(x^-)\big),
\end{split}
\end{equation}
and that $V_h'(z)$, $G_{\rho_g}(z)$ and $\mathcal{F}^{\amsmathbb{R}}_{g,h_2}(z,z_2)$ are holomorphic functions of $z$ in a neighborhood of $[\alpha_h,\beta_h]$ as soon as $g \neq h$. This leads us to the equation
{\footnotesize \begin{equation}
\begin{split}
\label{disc2ndSDeqn}0 & = \big(G_{\rho_h}(x^+) - G_{\rho_h}(x^-)\big)\left(\sum_{g \neq h} 2\theta_{g,h} \oint_{\gamma_g} \frac{\dd\tilde{\zeta}}{2\ii\pi}\,\frac{\mathcal{F}_{g,h_2}^{\amsmathbb{R}}(\tilde{\zeta},z_2)}{x - \tilde{\zeta}} + \theta_{h,h}(\mathcal{F}^{\amsmathbb{R}}_{h,h_2}(x^+,z_2) + \mathcal{F}^{\amsmathbb{R}}_{h,h_2}(x^-,z_2)\big) + \partial_{z_2}\bigg(\frac{\delta_{h,h_2}}{x - z_2}\bigg)\right) \\
& \quad + \big(\mathcal{F}_{h,h_2}^{\amsmathbb{R}}(x^+,z_2) - \mathcal{F}_{h,h_2}^{\amsmathbb{R}}(x^-,z_2)\big) \left(\sum_{g \neq h} 2\theta_{g,h} \oint_{\gamma_g} \frac{\dd\tilde{\zeta}}{2\ii\pi}\,\frac{G_{\rho_g}(\tilde{\zeta})}{x - \tilde{\zeta}} + \theta_{h,h}\big(G_{\rho_h}(x^+) + G_{\rho_h}(x^-)\big) - V_h'(x)\right).
\end{split}
\end{equation}}
\! In both lines, the contour $\gamma_g$ can be moved to infinity, picking up the only pole of the integrands at $\tilde{\zeta} = x$. The quantity in bracket in the second line becomes
\begin{equation}
\begin{split}
& \quad \theta_{h,h}\big(G_{\rho_h}(x^+) + G_{\rho_h}(x^-)\big) + \sum_{g \neq h} 2\theta_{g,h} G_{\rho_g}(x) - V'_h(x) \\
& = 2\sum_{g = 1}^{H} \partial_x\left( \int_{\alpha_g}^{\beta_g} \rho_g(y)\log|x - y|\right) - V_h'(x) \\
& = -\partial_x V_h^{\textnormal{eff}}(x),
\end{split}
\end{equation}
which vanishes as $x$ belongs to the interior of a band. Then, $G_{\rho_h}(x^+) \neq G_{\rho_h}(x^-)$ by off-criticality, so the quantity in bracket in the first line of \eqref{disc2ndSDeqn} vanishes, that is after the aforementioned contour manipulation
\[
\theta_{h,h}\big(\mathcal{F}^{\amsmathbb{R}}_{h,h_2}(x^+,z_2) + \mathcal{F}^{\amsmathbb{R}}_{h,h_2}(x^-,z_2)\big) + \sum_{g \neq h} 2\theta_{g,h} \mathcal{F}^{\amsmathbb{R}}_{g,h_2}(x,z_2) + \frac{\delta_{h,h_2}}{(x - z_2)^2} = 0.
\]
This is exactly the Riemann--Hilbert problem of Definition~\ref{def:masterRHP} satisfied by the fundamental solution of Definition~\ref{def:Berg}. As we see, it is a consequence of the Dyson--Schwinger equations of the continuous ensemble in the same way that the master Riemann--Hilbert problem for the leading covariance of the discrete ensembles was a consequence of Nekrasov equations of discrete ensembles.
\end{proof}

Here are other properties of $\boldsymbol{\mathcal{F}}(z_1,z_2)$, which we prove directly.

\begin{proposition}
\label{thmBfund} The fundamental solution $\boldsymbol{\mathcal{F}}(z_1,z_2)$ has the following properties.
\begin{enumerate}
\item[(F1)] Writing $\{i,j\} = \{1,2\}$, it is tame with exponent $1$ with respect to $z_i$, uniformly so for $z_j$ in any compact subset of ${\amsmathbb{C}} \setminus [\alpha_{h_j},\beta_{h_j}]$.
\item[(F2)] It is a solution of the master Riemann--Hilbert problem with respect to both its variables. Namely, for any $h_1,h_2 \in [H]$, $z_i \in \amsmathbb{C} \setminus [\alpha_{h_i},\beta_{h_i}]$ and $x_i \in (\alpha_{h_i},\beta_{h_i})$ we have
\begin{equation*}
\begin{split}
\theta_{h_1,h_2}\big(\mathcal{F}_{h_1,h_2}(x_1^+,z_2) + \mathcal{F}_{h_1,h_2}(x_1^-,z_2)\big) + \sum_{g_1 \neq h_1} 2\theta_{h_1,g_1} \mathcal{F}_{g_1,h_2}(x_1,z_2) & = - \frac{\delta_{h_1,h_2}}{(x_1 - z_2)^2}, \\
\theta_{h_1,h_2}\big(\mathcal{F}_{h_1,h_2}(z_1,x_2^+) + \mathcal{F}_{h_1,h_2}(z_1,x_2^-)\big) + \sum_{g_2 \neq h_2} 2\theta_{h_2,g_2} \mathcal{F}_{h_1,g_2}(z_1,x_2) & = - \frac{1}{(z_1 - x_2)^2} .
\end{split}
\end{equation*}
Writing $\{i,j\} = \{1,2\}$, the limit indicated by $x_i^{\pm}$ is uniform with respect to $z_j$ in any compact of ${\amsmathbb{C}} \setminus [\alpha_{h_j},\beta_{h_j}]$, and a continuous function of $x_i$. Note that the source function in both cases is $D_{h_1,h_2}(z_1,z_2) = -\frac{\delta_{h_1,h_2}}{(z_1 - z_2)^2}$. \item[(F3)] It has a joint tameness property, namely there exists a constant $M > 0$ such that, for any $z_1 \in {\amsmathbb{C}} \setminus [\alpha_{h_1},\beta_{h_1}]$ and $z_2 \in {\amsmathbb{C}} \setminus [\alpha_{h_2},\beta_{h_2}]$
\[
\big|(z_1 - z_2)^{2\delta_{h_1,h_2}}\sigma_{h_1}(z_1)\sigma_{h_2}(z_2)\mathcal{F}_{h_1,h_2}(z_1,z_2)\big| \leq M
\]
\item[(F4)] For each $\varepsilon_i \in \{\pm 1\}$ and $x_i \in (\alpha_{h_i},\beta_{h_i})$ with $i = 1,2$, the joint boundary values
\[
\lim_{(t_1,t_2) \rightarrow (0^+,0^+)} \bigg(\mathcal{F}_{h_1,h_2}( x_1 + \ii \varepsilon_1 t_1,\,x_2 + \ii \varepsilon_2 t_2) + \frac{\delta_{h_1,h_2}\delta_{\varepsilon_1 + \varepsilon_2,0}}{( x_1 + \ii \varepsilon_1 t_1 - x_2 - \ii \varepsilon_2 t_2)^2}\bigg)
\]
exist uniformly for $x_i$ in any compact subset of $(\alpha_{h_i},\beta_{h_i})$, and are continuous functions of $x_1,x_2$.
\item[(F5)] for any $h, h_i \in [H]$ and $z_i \in \amsmathbb{C} \setminus [\alpha_{h_i},\beta_{h_i}]$ with $i = 1,2$, we have
\begin{equation}
\label{eq_F_integrals_vanish}
\oint_{\gamma_{h}} \mathcal{F}_{h_1,h_2}({z}_1,z_2)\dd z_1 = \oint_{\gamma_{h}} \mathcal{F}_{h_1,h_2}(z_1,{z}_2)\dd z_2 = 0.
\end{equation}
\item[(F6)] $\mathcal{F}_{h_1,h_2}(z_1^*,z_2^*) = \big(\mathcal{F}_{h_1,h_2}(z_1,z_2)\big)^*$.
\end{enumerate}
\end{proposition}

\begin{proof}
We start by fixing $h_2 \in [H]$ and $z_2 \in {\amsmathbb{C}} \setminus [\alpha_{h_2},\beta_{h_2}]$. The source
\begin{equation}
\label{Ehh2zz2}
E_{h_1;h_2}(z;z_2) = \frac{\delta_{h_1,h_2}\sigma_{h_2}(z)}{2(z_1 - z_2)^2}
\end{equation}
is such that $\boldsymbol{E}(z_1)/\sigma(z_1)$ is tame (using a neighborhood of the segments that does not contain $z_2$) with tameness exponent $0$, uniformly with respect to $z_2$ in any compact of ${\amsmathbb{C}} \setminus [\alpha_{h_2},\beta_{h_2}]$. The analytic properties of the source and the properties stated in Theorem~\ref{Theorem_Master_equation_12} imply (F1) with respect to $z_1$. The tameness property of $\boldsymbol{E}/\sigma$ and Proposition~\ref{prop:tameRHP} imply that $(\mathcal{F}_{h_1,h_2}(z_1,z_2))_{h_1 = 1}^{H}$ as a function of $z_1$ is tame and satisfies the master Riemann--Hilbert problem with source functions $D_{h_1}(z_1) = - \frac{\delta_{h_1,h_2}}{(z_1 - z_2)^2}$. This shows (F2) with respect to $z_1$. (F1) and (F2) with respect to $z_2$ follow from the symmetry provided by Theorem~\ref{thm:Bsym}.

To establish (F3) and (F4), we only need to examine the behavior of $\mathcal{F}_{h_1,h_2}(z_1,z_2)$ as $z_1$ and $z_2$ approach $[\alpha_{h_1},\beta_{h_1}]$ and $[\alpha_{h_2},\beta_{h_2}]$ because the aforementioned uniform tameness in each variable guarantees the bound (F3) elsewhere. We return to \eqref{rewritingEF} in the proof of Theorem~\ref{Theorem_Master_equation_12} that we specialize to the fundamental solution as given in Definition~\ref{def:Berg}. With the specific source \eqref{Ehh2zz2} this yields for any $h_2 \in [H]$ and $z_2 \in {\amsmathbb{C}} \setminus [\alpha_{h_2},\beta_{h_2}]$ the relation
\begin{equation}
\label{soHs7un}
\begin{split}
 \sum_{h = 1}^{H} \mathcal{F}_{h,h_2}(z_1,z_2) & = \frac{1}{\sigma(z_1)} \sum_{h = 1}^{H} \Res_{\zeta = \infty} \left[ \frac{\sigma(\zeta)\mathcal{F}_{h,h_2}(\zeta,z_2)}{z_1 - \zeta}\right] \\
& \quad + \sum_{1 \leq g \neq h \leq H} \frac{1}{\sigma(z_1)} \oint_{\gamma_{h}} \frac{\dd \zeta}{2\ii\pi} \frac{\sigma(\zeta)}{\zeta - z_1} \,\frac{\theta_{h,h} - \theta_{h,g}}{\theta_{h,h}} \mathcal{F}_{g,h_2}(\zeta,z_2) \\
& \quad + \frac{1}{\sigma(z_1)} \oint_{\gamma_{h_2}} \frac{\dd \zeta}{2\ii\pi} \frac{\sigma(\zeta)}{\zeta - z_1} \frac{1}{2\theta_{h_2,h_2}(\zeta - z_2)^2},
\end{split}
\end{equation}
where $|\zeta| > |z_1|$ in the integration contour realizing the residue near $\infty$ in the first-line of the right-hand side, and $z_1,z_2$ are outside the contours $\gamma_h$ and $\gamma_{h_2}$ in the last two lines. The contribution of the first line does not cause trouble for (F3) and (F4), because we know from (F1) the uniform tameness of $\boldsymbol{\mathcal{F}}$ with respect to its second variable. For the contribution of the second line, we should check for each $h_1 \in [H]$ what happens as $z_1$ approaches $[\alpha_{h_1},\beta_{h_1}]$. For this purpose, we can move the contour in the $h_1$-th term so that $z_1$ is now inside the contour, creating an additional term
\begin{equation}
\label{additionalr}
\sum_{g \neq h_1} \frac{\theta_{h_1,h_1} - \theta_{h_1,g_1}}{\theta_{h_1,h_1}} \mathcal{F}_{g,h_2}(z_1,z_2).
\end{equation}
The contour integral is now a holomorphic function of $z_1$ in a neighborhood of $[\alpha_{h_1},\beta_{h_1}]$, and so is the additional term \eqref{additionalr} because $g \neq h_1$. Combining with the tameness of $\boldsymbol{\mathcal{F}}(z_1,z_2)$ in the variable $z_2$ uniformly in $z_1$ away from $[\alpha_{g},\beta_{g}]$ that is known by (F1), we see that the two terms do not cause problems for (F3) and (F4). Now consider the last line of \eqref{soHs7un}. If $h_1 \neq h_2$, we move the contour so that $z_2$ is inside, getting a holomorphic function of $(z_1,z_2)$ in a complex neighborhood of $[\alpha_{h_1},\beta_{h_1}] \times [\alpha_{h_2},\beta_{h_2}]$. The price to pay is an additional term
\[
-\frac{1}{2\theta_{h_2,h_2}\sigma(z_1)}\partial_{z_2}\left(\frac{\sigma(z_2)}{z_2 - z_1}\right).
\]
If $h_1 = h_2$, we move the contour so that $z_1$ and $z_2$ are inside, a holomorphic function of $(z_1,z_2)$ in a neighborhood of $[\alpha_{h_1},\beta_{h_1}] \times [\alpha_{h_2},\beta_{h_2}]$. The price to pay is an additional term
\[
-\frac{1}{2\theta_{h_2,h_2}\sigma(z_1)} \partial_{z_2}\left(\frac{\sigma(z_1) - \sigma(z_2)}{z_1 - z_2}\right).
\]
In both cases, one can see that (F3) and (F4) hold. In particular, the fact that $\sigma(x^+) = - \sigma(x^-)$ for any $x \in [\alpha_{h},\beta_h]$ and $h \in [H]$ is responsible for the double pole of $\mathcal{F}_{h,h}(z_1,z_2)$ as $z_1$ and $z_2$ approach the same point on $(\alpha_h,\beta_h)$ with imaginary parts of opposite signs tending to $0$.

The vanishing (F5) comes from $\boldsymbol{\kappa} = \boldsymbol{0}$ in the definition \eqref{eq_covariance_11}. The reality property (F6) comes from the fact that $\alpha_h,\beta_h$ are real, so $\big(\boldsymbol{\mathcal{F}}(z_1^*,z_2^*)\big)^*$ satisfies the same master problem and enjoys the same analytic properties as $\boldsymbol{\mathcal{F}}(z_1,z_2)$. By uniqueness of the solution of the master problem, they must coincide.
\end{proof}

\subsection{General solution of the master Riemann--Hilbert problem}
 \label{sec:1233}

The meromorphic forms on a compact Riemann surface are classified according to three types: holomorphic (first kind), meromorphic but residueless (second kind), meromorphic with simple poles (third kind). Any meromorphic form can be decomposed into the sum of the forms of these types. We will introduce particular solutions of the master Riemann--Hilbert problem of Definition~\ref{def:masterRHP} which roughly mimic this classification (this will become clearer in Section~\ref{sec:3kinds}). Their properties are summarized in Lemma~\ref{lem:genRHPsol}. They will be the building blocks to describe general solutions of the master Riemann--Hilbert problem in Theorem~\ref{thm:genRHPsol}.

Let $\overline{\amsmathbb{B}}_h = [\alpha_h,\beta_h]$ and $\overline{\amsmathbb{B}} = \bigcup_{h = 1}^H [\alpha_h,\beta_h]$. This matches the notation of Section~\ref{sec:auxiliary_assum} when $\amsmathbb{B}_h = (\alpha_h,\beta_h)$ was the $h$-th band of an equilibrium measure, but here ($\overline{\amsmathbb{B}}_h)_{h = 1}^H$ are arbitrary pairwise disjoint segments. We recall that $\gamma_h$ denotes a contour surrounding $\overline{\amsmathbb{B}}_h$ but not $\bigcup_{g \neq h} \overline{\amsmathbb{B}}_{g}$, that can be taken arbitrarily close to $\overline{\amsmathbb{B}}_h$.

\medskip

\begin{definition}[Functions of the first kind] \label{def:1stkind} Let $o_g$ be a point in $(\alpha_g,\beta_g)$ and let $c_{o_g^+}$ (respectively $c_{o_g^-}$) be a path from $\infty$ to $o_g$ in ${\amsmathbb{C}} \setminus \overline{\amsmathbb{B}}$ approaching $o_g$ from the upper (respectively, lower) half-plane. We define for any $g,h \in [H]$
\begin{equation}
\label{ufrakgh}
\mathfrak{c}_{h;g}^{\textnormal{1st}}(z) = \frac{\delta_{h,g}}{z - o_g} + \int_{c_{o_g^+}} \langle \bth_g \cdot \boldsymbol{\mathcal{F}}_{h,\bullet}(z,\zeta) \rangle \dd \zeta+ \int_{ c_{o_g^-}} \langle \bth_g \cdot \boldsymbol{\mathcal{F}}_{h,\bullet}(z,\zeta) \rangle \dd \zeta.
\end{equation}
\end{definition}
\begin{lemma}
\label{Lem:1stkinddef}
 $\mathfrak{c}_{h;g}^{\textnormal{1st}}(z)$ is a holomorphic function of $z \in \amsmathbb{C} \setminus \overline{\amsmathbb{B}}_h$, which does not does not depend on the choice of $o_g$ and paths $c_{o_g}^{\pm}$. For each $g \in [H]$ the tuple $\boldsymbol{\mathfrak{c}}_{\bullet;g}^{\textnormal{1st}}(z) = (\mathfrak{c}^{\textnormal{1st}}_{h;g}(z))_{h = 1}^H$ is tame with exponent $1$.
\end{lemma}
\begin{proof}
 The integrals in \eqref{ufrakgh} exist because $\boldsymbol{\mathcal{F}}(z,\zeta)$ decays as $O(\frac{1}{\zeta^2})$, $\zeta \rightarrow \infty$, by definition \eqref{eq_covariance_11}.

 For fixed $o_g$, if we choose other paths $c_{o_g^{\pm}}$, the right-hand side of \eqref{ufrakgh} changes by an integral of $\sum_{g' = 1}^H \theta_{g,g'} \mathcal{F}_{h,g'}(z,\zeta)$ over $\zeta$ in some closed contour, which might touch, but does not intersect $\overline{\amsmathbb{B}}$. All such integrals are zero by (F5) in Proposition~\ref{thmBfund}. Besides, for fixed $z \in \amsmathbb{C} \setminus \overline{\amsmathbb{B}}_h$, the derivative of \eqref{ufrakgh} with respect to $o_g \in (\alpha_g,\beta_g)$ is
 \begin{equation}
 \label{RHPeqB}
 \frac{\delta_{h,g}}{(z - o_g)^2} + \theta_{g,g}\big(\mathcal{F}_{h,g}(z,o_g^+) + \mathcal{F}_{h,g}(z,o_g^-)\big) + \sum_{g' \neq g} 2\theta_{g,g'} \mathcal{F}_{h,g'}(z,o_g),
 \end{equation}
 which is zero since by (F2) in Proposition~\ref{thmBfund}.

 $\mathfrak{c}_{h;g}^{\textnormal{1st}}(z)$ is holomorphic because so is $\boldsymbol{\mathcal{F}}(z,\zeta)$. Tameness follows from (F3) in Proposition~\ref{thmBfund}.
\end{proof}

We will define three types of functions of the second kind, which differ by their singularities and the singularities of the source for the master Riemann--Hilbert problem they solve. Functions of the first type have a source with singularities away from the endpoints; functions of the second type have singularities at the edge but zero source; functions of the third type and the corresponding source have singularities at the edges.

\begin{definition}[Functions of the second kind, type I] \label{def:2ndkindI} For each $h \in [H]$ and $H$-tuple $\boldsymbol{f}(z)$ such that $f_g$ is a holomorphic function in a neighborhood of $\gamma_g$, we define
\[
\mathfrak{F}_h^{\textnormal{I}}[\boldsymbol{f}](z) = \sum_{g = 1}^H \oint_{\gamma_g} \frac{\dd\zeta}{2\ii\pi} \mathcal{F}_{h,g}(z,\zeta) f_g(\zeta).
\]
\end{definition}
This defines a holomorphic function of $z \in \amsmathbb{C} \setminus \overline{\amsmathbb{B}}_h$, and the corresponding $H$-tuple of functions $\boldsymbol{\mathfrak{\mathfrak{F}}}^{\textnormal{I}}[\boldsymbol{f}](z) = (\mathfrak{F}_h^{\textnormal{I}}[\boldsymbol{f}](z))_{h = 1}^H$ is tame with exponent $1$ by (F3) in Proposition~\ref{thmBfund}.

\begin{definition}[Functions of the second kind, type II and III] \label{def:2ndkindII}
For each $g \in [H]$, let $\gamma(\alpha_g)$ (respectively $\gamma(\beta_g)$) be a small positively oriented circular loop around $\alpha_g$ (respectively $\beta_g$). If $f_g(z)$ is a meromorphic function in a complex neighborhood $\amsmathbb{M}_g$ of $[\alpha_g,\beta_g]$ with possible poles only at $\alpha_g,\beta_g$, then we define for any $h \in [H]$
\begin{equation}
\label{BIIBIII}
\begin{split}
\mathfrak{F}_h^{\textnormal{II}}[\boldsymbol{f}](z) & = \sum_{g = 1}^H \left[\oint_{\gamma(\alpha_g)} + \oint_{\gamma(\beta_g)}\right] \frac{\dd \zeta}{2\ii\pi}\left(\frac{\delta_{g,h}}{2(z - \zeta)^2} + \big \langle \bth_g \cdot \mathcal{F}_{h,\bullet}(z,\zeta)\big \rangle\right) \sigma_g(\zeta) f_g(\zeta), \\
\mathfrak{F}_h^{\textnormal{III}}[\boldsymbol{f}](z) & = \sum_{g = 1}^H \left[\oint_{\gamma(\alpha_g)} + \oint_{\gamma(\beta_g)}\right] \frac{\dd\zeta}{2\ii\pi} \left(\frac{1}{2(z - \zeta)^2} + \sum_{g' \neq g} \theta_{g,g'} \mathcal{F}_{h,g'}(z,\zeta) \right)\frac{f_g(\zeta)}{\theta_{g,g}},
\end{split}
\end{equation}
where $z$ is outside the contours of integration.
\end{definition}
\begin{lemma}
\label{Lem:2ndkinddef}
 $\mathfrak{F}_h^{\textnormal{II}/\textnormal{III}}[\boldsymbol{f}](z)= (\mathfrak{F}_h^{\textnormal{II}/\textnormal{III}}[\boldsymbol{f}](z))_{h = 1}^H$ defines a holomorphic function of $z \in \amsmathbb{C} \setminus \overline{\amsmathbb{B}}_h$, which does not depend on the choice of $\gamma(\alpha_g)$ and $\gamma(\beta_g)$. If $\boldsymbol{f}(z)$ is tame, then so is  $\boldsymbol{\mathfrak{F}}^{\textnormal{II}/\textnormal{III}}[\boldsymbol{f}](z)$.
\end{lemma}
\begin{proof}
Analyticity of the integrands of \eqref{BIIBIII} in $\amsmathbb{M}_g \setminus \overline{\amsmathbb{B}}_g$ implies that the integrals can only depend on the intersection point $o_g$ of the loop with $\amsmathbb{B}_g$: the loop is then realized as an open path from $o_g^+$ to $o_g^-$ in $\amsmathbb{M}_g$. If we deform the loop so that $o_g$ is replaced with $\tilde{o}_g \in \amsmathbb{B}_g$, the integral in $\mathfrak{F}_h^{\textnormal{II}/\textnormal{III}}$ gets an extra contribution equal to the integral from $o_g$ and $\tilde{o}_g$ of the jump of the original integrand when crossing the interval $\amsmathbb{B}_g$ from the upper half-plane to the lower half-plane. In type III this jump is manifestly $0$, hence $\mathfrak{F}_h^{\textnormal{III}}[\boldsymbol{f}]$ does not depend on the choice of small loops around $\alpha_g$ and $\beta_g$. The same is true in type II, because (F2) in Proposition~\ref{thmBfund} leads to the equality valid for any $x \in \amsmathbb{B}_g$
\[
\big\langle \bth_{g} \cdot \mathcal{F}_{h,\bullet}(z,x^+)\big\rangle + \big \langle \bth_g \cdot \mathcal{F}_{h,\bullet}(z,x^-)\big \rangle + \frac{\delta_{h,g}}{(x - z)^2} = 0,
\]
and we have $\sigma_g(x^+)f_g(x^+) = - \sigma_g(x^-)f_g(x^-)$.

Tameness of $\boldsymbol{\mathfrak{F}}^{\textnormal{II}/\textnormal{III}}[\boldsymbol{f}](z)$ can be seen by moving $z$ inside the contour of integration and collecting the residue as $\zeta=z$, which is responsible for the growth near $z=\alpha_g$ or $z=\beta_g$.
\end{proof}

\begin{definition}[Functions of the third kind] \label{def:3rdkind}
Let $h,g \in [H]$, $z \in \amsmathbb{C}\setminus \overline{\amsmathbb{B}}_h$, $\zeta \in \amsmathbb{C} \setminus \amsmathbb{B}$ and $c_{\zeta}$ be a path from $\infty$ to $\zeta$ in $\amsmathbb{C} \setminus \overline{\amsmathbb{B}}$, we define
\[
\mathfrak{c}_{h;g}^{\textnormal{3rd}}(z;\zeta) = \int_{c_{\zeta}} \mathcal{F}_{h,g}(z,w)\dd w.
\]
\end{definition}
The integral is well-defined due to the decay of $\mathcal F_{h,g}$ at $\infty$ and there is no dependence on the choice of the path $c_{\zeta}$, because contour integrals of $ \mathcal{F}_{h,g}(z,w)$ vanish due to (F5) in Proposition~\ref{thmBfund}. The function $\mathfrak{c}^{\textnormal{3rd}}_{h;g}(z;\zeta)$ is holomorphic in $z \in \amsmathbb{C} \setminus \overline{\amsmathbb{B}}_h$ and $\zeta \in \amsmathbb{C} \setminus \overline{\amsmathbb{B}}_g$. Note that $\mathfrak{c}_{h;g}^{\textnormal{3rd}}(z;\zeta)$ is also well-defined at $\zeta = \alpha_g$ and $\zeta=\beta_g$, because $\zeta \mapsto \mathcal{F}_{h,g}(z,\zeta)$ diverges there at most like an inverse square root.

\medskip

We now explain in three statements the relevance of function $\mathfrak{c}$ to an alternative construction of the solution operator $\Upsilon_h$ recast as a solution to the Riemann-Hilbert problem of Definition~\ref{def:masterRHP} in Proposition~\ref{prop:tameRHP}. The proofs are in the next subsection.

\begin{lemma} Functions $\mathfrak{c}$ have the following properties (see also Table~\ref{table:sumf}):
\label{lem:genRHPsol}
\begin{itemize}
\item For each $g \in [H]$, the $H$-tuple of functions $\boldsymbol{\mathfrak{c}}^{\textnormal{1st}}_{\bullet;g}(z)$ satisfies the master Riemann--Hilbert problem with zero source. Besides, we have $\mathfrak{c}^{\textnormal{1st}}_{h;g}(z) = \frac{\delta_{h,g}}{z} + O(\frac{1}{z^2})$ as $z \rightarrow \infty$ for any $g,h \in [H]$; equivalently: $\oint_{\gamma_h} \frac{\dd z}{2\ii\pi} \mathfrak{c}^{\textnormal{1st}}_{h;g}(z) = \delta_{h,g}$. And, it is tame with exponent $1$.
\item For each $\boldsymbol{f}$, the $H$-tuple of functions $\boldsymbol{\mathfrak{F}}^{\textnormal{I}}[\boldsymbol{f}](z)$ satisfies the master Riemann--Hilbert problem with source
\[
\forall h \in [H]\qquad D_h^{\textnormal{I}}[\boldsymbol{f}](z) = - \oint_{\gamma_h} \frac{\dd\zeta}{2\ii\pi}\, \frac{f_h(\zeta)}{(z - \zeta)^2}.
\]
where $z$ is \emph{inside} the contour of integration. Besides, it behaves as $O(\frac{1}{z^2})$ as $z \rightarrow \infty$. And, it is tame with exponent $1$.
\item For each $\boldsymbol{f}$, the $H$-tuple of functions $\boldsymbol{\mathfrak{F}}^{\textnormal{II}}[\boldsymbol{f}](z)$ satisfies the master Riemann--Hilbert problem with zero source. Besides, it is $O(\frac{1}{z^2})$ as $z \rightarrow \infty$, and it becomes tame with exponent $1$ after subtraction of $(\frac{1}{2} \partial_z(\sigma_h(z)f_h(z)))_{h = 1}^H$.
\item For each $\boldsymbol{f}$, the $H$-tuple of functions $\boldsymbol{\mathfrak{F}}^{\textnormal{III}}[\boldsymbol{f}](z)$ satisfies the master Riemann--Hilbert problem with source
\[
\forall h \in [H] \qquad D_h^{\textnormal{III}}[\boldsymbol{f}](z) = \oint_{\gamma(\alpha_h)} \frac{\dd\zeta}{2\ii\pi}\,\frac{f_h(\zeta)}{(z - \zeta)^2}+\oint_{\gamma(\beta_h)} \frac{\dd\zeta}{2\ii\pi}\,\frac{f_h(\zeta)}{(z - \zeta)^2}.
\]
Besides, it is $O(\frac{1}{z^2})$ as $z \rightarrow \infty$, and it becomes tame with exponent $1$ after subtraction of $(\frac{1}{2\theta_{h,h}} \partial_z f_h(z))_{h = 1}^H$.
\item For each $g \in [H]$ and $\zeta \in \amsmathbb{C} \setminus \amsmathbb{B}$, the $H$-tuple of functions $\boldsymbol{\mathfrak{c}}^{\textnormal{3rd}}_{\bullet;g}(z;\zeta)$ satisfies the master Riemann--Hilbert problem with respect to $z$ and with source
\[
\forall h \in [H] \qquad D_{h;g}(z;\zeta) = -\frac{\delta_{h,g}}{z - \zeta}.
\]
Besides, it is $O(\frac{1}{z^2})$ as $z \rightarrow \infty$. If $\zeta$ is away from the endpoints $\alpha_g,\beta_g$ it is tame with exponent $1$. If $\zeta$ is located at $\alpha_g$ or $\beta_g$, it becomes tame (with respect to the variable $z$) with exponent $1$ after we subtract $\frac{1}{2\theta_{g,g}}\boldsymbol{D}_{\bullet;g}(z;\zeta)$.
\end{itemize}
\end{lemma}

\begin{table}
\begin{center}
\begin{tabular}{cccc}
\vspace{0.1cm}
$\boldsymbol{F}(z)$ & Tame exponent & Source & $\displaystyle\oint_{\gamma_g} F_h(z)\frac{\dd z}{2\ii\pi}$ \\
\hline\hline
\vspace{0.07cm}
1st kind & $1$ & $0$ & $\delta_{h,g}$ \\
\hline
\vspace{0.07cm}
2nd kind, I & $1$ & holomorphic & $0$ \\
\hline
\vspace{0.07cm}
2nd kind, II & odd $> 1$ & $0$ & $0$ \\
\hline
\vspace{0.07cm}
2nd kind, III & even $ > 1$ & $\begin{smallmatrix} \textnormal{meromorphic} \\ \textnormal{no residues}\end{smallmatrix}$ & $0$ \\
\hline
\vspace{0.07cm}
3rd kind & $1\,\,\textnormal{or}\,\,2$ & simple poles & $0$ \\
\hline
\hline
\end{tabular}
\end{center}
\caption{\label{table:sumf} Summary of the properties of the functions of Lemma~\ref{lem:genRHPsol}.}
\end{table}

\begin{theorem}
\label{thm:genRHPsol}
Consider the master Riemann--Hilbert problem with sources $(D_h(z))_{h = 1}^H$. In addition, choose $\boldsymbol{\kappa} \in \amsmathbb{C}$ and rational functions $(L_h(z))_{h=1}^H$ that can have poles at $z = \alpha_h,\beta_h$ and such that $L_h(z) = O(\frac{1}{z})$ as $z \rightarrow \infty$ for each $h\in[H]$.

Then there exists a unique solution $(F_h(z))_{h = 1}^H$ to the master Riemann--Hilbert problem such that:
\begin{enumerate}
\item $\displaystyle \left(F_h(z) - \sigma_h^{-1}(z) L_h(z) - \frac{1}{2\theta_{h,h}}D_h(z)\right)_{h = 1}^H$ is tame with exponent $1$; and
\item $\displaystyle\oint_{\gamma_h} \frac{\dd z}{2\ii\pi}\,F_h(z) = \kappa_h$ for all $h\in[H]$.
\end{enumerate}
The solution is given by the formula:
\[
F_h(z) = \sum_{g = 1}^H \big(\kappa_g \mathfrak{c}^{\textnormal{1st}}_{h;g}(z) - r(\alpha_g)\mathfrak{c}_{h;g}^{\textnormal{3rd}}(z;\alpha_g) - r(\beta_g)\mathfrak{c}^{\textnormal{3rd}}_{h;g}(\z;\beta_g)\big) - \sum_{\textnormal{J} \in \{\textnormal{I},\textnormal{II},\textnormal{III}\}} \mathfrak{F}_h^{\textnormal{J}}[\boldsymbol{f}^{\textnormal{J}}](z).
\]
The $H$-tuple of functions $\boldsymbol{f}^{\textnormal{J}}$ and the constants $(r(\alpha_g),r(\beta_g))_{g = 1}^{H}$ are determined by
\[
\forall g \in [H]\qquad f_g^{\textnormal{II}}(z) = \frac{2}{\sigma_g(z)} \int_{z}^{\infty} \frac{L_g(\zeta)}{\sigma_g(\zeta)}\dd \zeta,
\]
and the (unique up to addition of a constant to $f^{\textnormal{I}}_g$) decomposition of the source
\begin{equation}
\label{sourcedec}
\forall g \in [H]\qquad D_g(z) = \partial_z f^{\textnormal{I}}_g(z) + \frac{r(\alpha_g)}{z - \alpha_g} + \frac{r(\beta_g)}{z -\beta_g} + \partial_z f^{\textnormal{III}}_g(z),
\end{equation}
where the function $f_g^{\textnormal{I}}(z)$ is holomorphic in a complex neighborhood of $\overline{\amsmathbb{B}}_g$ and $f^{\textnormal{III}}_g(z)$ is a rational function with possible poles at $\alpha_g,\beta_g$ and such that $f^{\textnormal{III}}_g(z) = O(\frac{1}{z})$ as $z \rightarrow \infty$.
\end{theorem}
\begin{corollary}
 \label{cor:1stkindapp} We have for any $h \in [H]$ and $z \in \amsmathbb{C} \setminus \overline{\amsmathbb{B}}_h$
 \[
 \Upsilon_h[\boldsymbol{0}\,;\,\boldsymbol{\kappa}](z) = \sum_{g = 1}^H \kappa_g \mathfrak{c}^{\textnormal{1st}}_{h;g}(z).
 \]
\end{corollary}

From the point of the master problem of Definition~\ref{def:master}, Theorem~\ref{thm:genRHPsol} (combined with Proposition~\ref{prop:tameRHP}) means that tame solutions can be expressed in terms of the leading covariance identified in Theorem~\ref{Theorem_correlators_expansion}. This does not give an independent construction of the solution operator of the master problem in Theorem~\ref{Theorem_Master_equation_12}, because the fundamental solution itself was defined via the solution operator, but also due to the restriction to tame solutions in Theorem~\ref{thm:genRHPsol}: the latter theorem requires an additional data in functions $L_h(z)$ controlling tameness, while the former does not. This is an important distinction for our applications of Theorem~\ref{Theorem_Master_equation_12} in Chapter~\ref{Chapter_fff_expansions}, as we need to apply the solution operator $\boldsymbol{\Upsilon}$ to functions that contain small remainders that are not a priori tame (\textit{e.g.} only defined away from $[\hat{a}_h^{\mathfrak{m}},\hat{b}_h^{\mathfrak{m}}]$ instead of away from the bands) and exploit the corresponding operator continuity. However, Theorem~\ref{thm:genRHPsol} could be useful to obtain alternative expressions for the coefficients of asymptotic expansion of the correlators mentioned in Theorem~\ref{Theorem_correlators_expansion}. Besides, Corollary~\ref{cor:1stkindapp} was used in Proposition~\ref{Proposition_differentiability_filling_fraction} to give a more aesthetically pleasing formula\footnote{However, it was not necessary: the relevant proofs of Chapter~\ref{Chapter_smoothness} still work if we only use Theorem~\ref{Theorem_Master_equation_12}.} for the first-order partial derivatives with respect to filling fractions of the equilibrium measure.

\begin{remark} \label{cdiscrem} The same master Riemann--Hilbert problem of Definition~\ref{def:masterRHP} appears in the study of the large $\N$ behavior of continuous $\sbeta$-ensemble with several groups of particles (\textit{cf.} \cite{CE06,BEO} or Chapter~\ref{SIntro}) having several groups of particles with different intensities of interaction. The results of the present chapter have an interest in this context. In particular, Theorem~\ref{thm:Bsym} resolves the question of symmetry of the fundamental solution that was left open in \cite{BEO}.
\end{remark}

\subsection{Proofs of Lemma~\ref{lem:genRHPsol}, Theorem~\ref{thm:genRHPsol}, and Corollary~\ref{cor:1stkindapp}}

\begin{proof}[Proof of Lemma~\ref{lem:genRHPsol}]

The behavior of all functions as $z \rightarrow \infty$ is clear from the definitions and the fact that $\boldsymbol{\mathcal{F}}(z,\zeta) = O(z^{-2})$. Tameness was already briefly discussed after definitions, and we add some details below. To establish the master Riemann--Hilbert problems satisfied by our functions, the starting point is the master Riemann--Hilbert problem satisfied by $z \mapsto (\mathcal{F}_{h;g}(z,\zeta))_{h = 1}^H$ (\textit{cf.} Proposition~\ref{thmBfund}), namely for any $g,h \in [H]$, $x \in \amsmathbb{B}_h$ and $\zeta \in \amsmathbb{C} \setminus \overline{\amsmathbb{B}}_g$ we have the equality
\begin{equation}
\label{BthRHOP}
\theta_{h,h}\big(\mathcal{F}_{h,g}(x^+,\zeta) + \mathcal{F}_{h,g}(x^-,\zeta)\big) + \sum_{h' \neq h} 2\theta_{h,h'} \mathcal{F}_{h',g}(x,\zeta) = - \frac{\delta_{h,g}}{(x - \zeta)^2}.
\end{equation}

\medskip

\noindent \textsc{First kind.} In \eqref{BthRHOP} we replace $g$ with $g'$, then multiply by $\theta_{g,g'}$, sum over $g' \in [H]$ and integrate the variable $\zeta$ on the path $c_{o_g}^+ + c_{o_g}^-$ (we can always choose $o_g$ distinct from $x$). The right-hand side is equal to $-\frac{2\theta_{g,h}}{x - o_g}$, while the left-hand side is equal to
\[
\theta_{h,h}\left(\mathfrak{c}^{\textnormal{1st}}_{h;g}(x^+) + \mathfrak{c}^{\textnormal{1st}}_{h;g}(x^-) - \frac{2\delta_{h,g}}{x - o_g}\right) + \sum_{h' \neq h} 2\theta_{h,h'}\left(\mathfrak{c}^{\textnormal{1st}}_{h';g}(x) - \frac{\delta_{h',g}}{x - o_g}\right).
\]
Therefore, for any $g,h \in [H]$ and $x \in \amsmathbb{B}_h$ we have
\[
\theta_{h,h}\big(\mathfrak{c}^{\textnormal{1st}}_{h;g}(x^+) + \mathfrak{c}^{\textnormal{1st}}_{h;g}(x^-)\big) + \sum_{h' \neq h} 2\theta_{h,h'} \mathfrak{c}^{\textnormal{1st}}_{h';g}(x) = 0.
\]

\medskip

\noindent \textsc{Second kind.} For type I, we multiply \eqref{BthRHOP} by $f_{g}(\zeta) \frac{\dd\zeta}{2\ii\pi}$, integrate $\zeta$ over $\gamma_{g}$, and sum over $g \in [H]$. We get for any $h \in [H]$ and $x \in \amsmathbb{B}_h$
\[
\theta_{h,h}\big(\mathfrak{F}_h^{\textnormal{I}}[\boldsymbol{f}](x^+) + \mathfrak{F}_h^{\textnormal{I}}[\boldsymbol{f}](x^-)\big) + \sum_{h' \neq h} 2\theta_{h,h'} \mathfrak{F}_{h'}^{\textnormal{I}}[\boldsymbol{f}](x) = - \oint_{\gamma_h} \frac{\dd \zeta}{2\ii\pi}\frac{f_h(\zeta)}{(x - \zeta)^2},
\]
For type II, in \eqref{BthRHOP} we replace $g$ with $g'$, multiply by $\theta_{g,g'}$ and sum over $g' \in [H]$. We get for any $g,h \in [H]$ and $x \in (\alpha_h,\beta_h)$
\[
\theta_{h,h}\left(\big \langle \bth_{g} \cdot \mathcal{F}_{h,\bullet}(x^+,\zeta)\big\rangle + \big \langle \bth_g \cdot \mathcal{F}_{h,\bullet}(x^-,\zeta)\big\rangle\right) + \sum_{h' \neq h} 2\theta_{h,h'} \big\langle \bth_{g} \cdot \mathcal{F}_{h',\bullet}(x,\zeta)\big\rangle = - \frac{\theta_{h,g}}{(x - \zeta)^2}.
\]
We can include the right-hand side into the left-hand side in the following way:
\begin{equation}
\label{BthRHOP2}
\begin{split}
& \quad \theta_{h,h}\left(\big \langle \bth_g \cdot \mathcal{F}_{h,\bullet}(x^+,\zeta)\big\rangle + \big \langle \bth_g \cdot \mathcal{F}_{h,\bullet}(x^-,\zeta)\big\rangle + \frac{\delta_{h,g}}{(x - \zeta)^2}\right) \\
& \quad + \sum_{h' \neq h} 2\theta_{h,h'}\left(\big\langle \bth_g \cdot \mathcal{F}_{h',\bullet}(x,\zeta)\big\rangle + \frac{\delta_{h',g}}{2(x - \zeta)^2}\right) = 0.
\end{split}
\end{equation}
Then, we multiply by $\sigma_g(\zeta)f_g(\zeta)\frac{\dd\zeta}{2\ii\pi}$, integrate $\zeta$ over loops $\gamma(\alpha_g),\gamma(\beta_g)$ that leave $x$ outside, and sum over $g \in [H]$. This leads to
\[
\theta_{h,h}\big(\mathfrak{F}_h^{\textnormal{II}}[\boldsymbol{f}](x^+) + \mathfrak{F}_{h}^{\textnormal{II}}[\boldsymbol{f}](x^-)\big) + \sum_{h' \neq h} 2\theta_{h,h'} \mathfrak{F}_{h'}^{\textnormal{II}}[\boldsymbol{f}](x) = 0.
\]
For type III, we define the vector $\widetilde{\bth}_g = ((1 - \delta_{g,g'})\theta_{g,g'})_{g = 1}^H$. Then
\[
\big\langle\widetilde{\bth}_g \cdot \mathcal{F}_{h,\bullet}(z,\zeta)\big\rangle = \sum_{g' \neq g} \theta_{g,g'} \mathcal{F}_{h,g'}(z,\zeta).
\]
Taking the difference between \eqref{BthRHOP2} on one hand, and \eqref{BthRHOP} multiplied by $\theta_{g,g}$ on the other hand, yields for any $g,h \in [H]$ and $x \in \amsmathbb{B}_h$
\begin{equation*}
\begin{split}
& \quad \theta_{h,h}\left(\big\langle \widetilde{\bth}_g \cdot \mathcal{F}_{h,\bullet}(x^+,\zeta)\big\rangle + \big\langle \widetilde{\bth}_g \cdot \mathcal{F}_{h,\bullet}(x^+,\zeta)\big\rangle + \frac{\delta_{h,g}}{(x - \zeta)^2}\right) \\
& \quad + \sum_{h' \neq h} 2\theta_{h,h'} \left(\big\langle \widetilde{\bth}_g \cdot \mathcal{F}_{h',\bullet}(x,\zeta)\big\rangle + \frac{\delta_{h',g}}{2(x - \zeta)^2}\right) = \frac{\delta_{h,g}\theta_{g,g}}{(x - \zeta)^2}.
\end{split}
\end{equation*}
We multiply by $\frac{f_g(\zeta)}{\theta_{g,g}} \frac{\dd\zeta}{2\ii\pi}$, integrate $\zeta$ over loops $\gamma(\alpha_g),\gamma(\beta_g)$ that leave $x$ outside, and sum over $g \in [H]$. This yields
\[
\theta_{h,h}\big(\mathfrak{F}_h^{\textnormal{III}}[\boldsymbol{f}](x^+) + \mathfrak{F}_h^{\textnormal{III}}[\boldsymbol{f}](x^-)\big) + \sum_{h' \neq h} 2\theta_{h,h'} \mathfrak{F}_{h'}^{\textnormal{III}}[\boldsymbol{f}](x) = \left[ \oint_{\gamma(\alpha_h)}+\oint_{ \gamma(\beta_h)}\right] \frac{\dd\zeta}{2\ii\pi}\,\frac{f_h(\zeta)}{(x - \zeta)^2}.
\]

For the tameness statements in type II (respectively III), we move $z$ inside the contour of integrations in Definition~\ref{def:3rdkind} at the price of adding compensating terms
\begin{equation*}
\begin{split}
\partial_z\left(\frac{\sigma_h(z)f_h(z)}{2}\right),\qquad \textnormal{respectively} \qquad \partial_z\left(\frac{f_h(z)}{2\theta_{h,h}}\right).
\end{split}
\end{equation*}
The claimed tameness property of $\boldsymbol{\mathfrak{F}}^{\textnormal{II}/\textnormal{III}}(z)$ then come from the observation that the right-hand sides of \eqref{BIIBIII} with $z$ inside the contour of integration diverge at most like an inverse square root as $z$ approaches $\alpha_h,\beta_h$ because $\boldsymbol{\mathcal{F}}(z,\zeta)$ has this property.

\medskip

\noindent \textsc{Third kind.} In \eqref{BthRHOP}, we replace $\zeta$ with $w$, and then integrate $w$ on the path $c_{\zeta}$. This results in
\[
\theta_{h,h}\big(\mathfrak{c}_{h;g}^{\textnormal{3rd}}(x^+;\zeta) + \mathfrak{c}_{h;g}^{\textnormal{3rd}}(x^-;\zeta)\big) + \sum_{h' \neq h} 2\theta_{h,h'} \mathfrak{c}^{\textnormal{3rd}}_{h';g}(x;\zeta) = -\frac{\delta_{h,g}}{x - \zeta}.
\]
If $\zeta \notin \{\alpha_g,\beta_g\}$, the $H$-tuple of functions $z \mapsto \boldsymbol{\mathfrak{c}}^{\textnormal{3rd}}_{\bullet;g}(z;\zeta)$ is tame with exponent $1$ because $\mathcal{F}_{h,g}(z,\zeta)$ diverges at most like an inverse square root as $z \rightarrow \alpha_h,\beta_h$ in a uniform way for $\zeta$ away from $\overline{\amsmathbb{B}}_g$. In the case $\zeta = \alpha_{g}$, we need to analyze the behavior of $\mathfrak{c}^{\textnormal{3rd}}_{h;g}(z;\alpha_g)$ as $z$ approaches $\alpha_h,\beta_h$. If $h \neq g$ there is an inverse square root divergence as before, but the case $h = g$ needs more care. We come back to the master problem solved by $\boldsymbol{\mathcal{F}}(z,w)$. The rewriting \eqref{rewritingEF} gives in this case, for any $g \in [H]$
\begin{equation*}
\begin{split}
& \quad \sum_{h = 1}^H \mathcal{F}_{h,g}(z,\zeta) - \frac{\mathcal{P}[\sum_{h = 1}^H \mathcal{F}_{h,g}(*,\zeta)](z)}{\sigma(z)} + \sum_{1 \leq h \neq h' \leq H} \frac{1}{\sigma(z)} \oint\frac{\dd w}{2\ii\pi} \frac{\sigma(w)}{w - z} \frac{\theta_{h,h} - \theta_{h,h'}}{\theta_{h,h}} \mathcal{F}_{h',g}(w,\zeta) \\
& = \sum_{h = 1}^H \frac{1}{2\theta_{h,h}\sigma(z)} \oint_{\gamma_h} \frac{\dd w}{2\ii\pi} \frac{\sigma(w)}{w - z} \frac{\delta_{h,g}}{(w - \zeta)^2} \\
& = \frac{1}{2\theta_{g,g}\sigma(z)} \partial_{\zeta} \bigg(- \frac{\sigma(z) - \sigma(\zeta)}{z - \zeta} + \sum_{\substack{h = 1 \\ h \neq g}}^H \oint_{\gamma_h} \frac{\dd w}{2\ii\pi} \frac{\sigma(w)}{(w - z)(w - \zeta)}\bigg),
\end{split}
\end{equation*}
where $z$ and $\zeta$ are outside the integration contours, we recall $\sigma^2(z) = \prod_{h = 1}^H (z - \alpha_h)(z - \beta_h)$ and the operator $\mathcal{P}$ returns a polynomial in $z$. Let us choose a point $o \in \amsmathbb{C} \setminus \overline{\amsmathbb{B}}$. Integrating $\zeta$ on a path from $o$ to $\alpha_g$ in this equation yields
\begin{equation*}
\begin{split}
& \sum_{h = 1}^H \big(\mathfrak{c}^{\textnormal{3rd}}_{h;g}(z;\alpha_g) -\mathfrak{c}^{\textnormal{3rd}}_{h;g}(z;o)\big) \!\! + \sum_{1 \leq h \neq h' \leq H}\!\! \frac{1}{\sigma(z)}\oint_{\gamma_h} \frac{\dd w}{2\ii\pi} \frac{\sigma(w)}{w - z} \frac{\theta_{h,h} - \theta_{h,h'}}{\theta_{h,h}} \big(\mathfrak{c}^{\textnormal{3rd}}_{h';g}(w;\alpha_g) - \mathfrak{c}^{\textnormal{3rd}}_{h';g}(w;o)\big) \\
& = \frac{\textnormal{Pol}(z)}{\sigma(z)} + \frac{1}{2\theta_{g,g}\sigma(z)}\bigg( - \frac{\sigma(z)}{z - \alpha_g} + \frac{\sigma(z) - \sigma(o)}{z - o} + \sum_{\substack{h = 1 \\ h \neq g}}^{H} \oint_{\gamma_h} \frac{\dd w}{2\ii\pi} \frac{\sigma(w)}{w - z}\Big(\frac{1}{w - \alpha_g} - \frac{1}{w - o}\Big)\bigg).
\end{split}
\end{equation*}
Since we know that $\boldsymbol{\mathfrak{c}}^{\textnormal{3rd}}_{\bullet;g}(z;o)$ is tame with exponent $1$, we deduce that for any $h,g \in [H]$ we have
\[
\mathfrak{c}^{\textnormal{3rd}}_{h;g}(z;\alpha_g) = -\frac{\delta_{h,g}}{2\theta_{g,g}(z - \alpha_g)} + O\left(\frac{1}{\sigma(z)}\right)
\]
as $z \rightarrow \alpha_h,\beta_h$, which is the tameness statement we needed to justify.
\end{proof}

\begin{proof}[Proof of Theorem~\ref{thm:genRHPsol}] We first justify uniqueness. Assume that $\boldsymbol{F},\widetilde{\boldsymbol{F}}$ are two solutions of this master Riemann--Hilbert problem. Then $\boldsymbol{F} - \widetilde{\boldsymbol{F}}$ is a tame solution of the master Riemann--Hilbert problem with vanishing source, tameness exponent $1$ and such that $\boldsymbol{F}(z) - \widetilde{\boldsymbol{F}}(z) = O(z^{-2})$ as $z \rightarrow \infty$. By the converse part in Proposition~\ref{prop:tameRHP}, it also solves the master problem of Section~\ref{sec:Masterpb} with $\boldsymbol{E} = 0$ and $\boldsymbol{\kappa} = 0$, hence it must be the zero function.

For existence, let $h \in [H]$ and consider the $h$-th component of the source $D_h(z)$. We fix a simply-connected complex neighborhood $\amsmathbb{M}_h$ of $\overline{\amsmathbb{B}}_h$ bounded by a closed curve with counterclockwise orientation $\gamma_h^+$, such that $D_h(z)$ is meromorphic for $z \in \amsmathbb{M}_h$ and extends continuously to $\gamma_h^+$. By assumptions the possible poles of $D_h(z)$ are located at $\alpha_h,\beta_h$. Let us denote
\[
r(\alpha_h) := \Res_{\zeta = \alpha_h} D_h(\zeta)\dd\zeta, \qquad r(\beta_h) := \Res_{\zeta = \beta_h} D_h(\zeta)\dd\zeta,\qquad \tilde{D}_h(z) := D_h(z) - \frac{r(\alpha_h)}{z - \alpha_h} - \frac{r(\beta_h)}{z - \beta_h}.
\]
By Cauchy residue formula, we have
\begin{equation}
\label{fornunqueanh}\forall z \in \amsmathbb{M}_h \setminus \{\alpha_h,\beta_h\}\qquad \tilde{D}_h(z) = -\oint_{\gamma_h^+} \frac{\dd\zeta}{2\ii\pi} \frac{\tilde{D}_h(\zeta)}{z - \zeta} + \Res_{\zeta = \alpha_h,\beta_h}\left[ \frac{\tilde{D}_h(\zeta) }{z - \zeta}\right].
\end{equation}
Since $\tilde{D}_h(\zeta)$ has no residues and $\amsmathbb{M}_h$ is simply-connected, it admits a primitive $\tilde{f}^{\textnormal{I}}_h(\zeta)$ which is a meromorphic function in $\amsmathbb{M}_h$. Then, integration by parts allows rewriting the first term of \eqref{fornunqueanh}
\begin{equation}
\label{holoterm}
\oint_{\gamma_h^+} \frac{\dd\zeta}{2\ii\pi}\frac{\tilde{f}_h^{\textnormal{I}}(\zeta)}{(z - \zeta)^2}.
 \end{equation}
This is a holomorphic function of $z \in \amsmathbb{M}_h$ since $z$ is inside the contour of integration. We may also rewrite it as $\partial_zf_h^{\textnormal{I}}(z)$ with
\[
f^{\textnormal{I}}_h(z) = \oint_{\gamma_h^+} \frac{\dd\zeta}{2\ii\pi} \frac{\tilde{f}_h^{\textnormal{I}}(\zeta)}{z - \zeta}.
\]
The advantage of the last description is that the function $f^{\textnormal{I}}(z)$ is holomorphic in $\amsmathbb{M}_h$ (this is not necessary to use Definition~\ref{def:2ndkindI} but it brings us closer to the decomposition of $D_h$ given in \eqref{sourcedec}) and we also have
\[
\partial_z f_h^{\textnormal{I}}(z) = \oint_{\gamma_h^+} \frac{\dd\zeta}{2\ii\pi} \frac{f_h^{\textnormal{I}}(\zeta)}{(z - \zeta)^2}.
\]

The last term in \eqref{fornunqueanh} is a rational function of $z$ with possible poles at $\alpha_h,\beta_h$ and no residues. Therefore, it admits a primitive that we call $f_h^{\textnormal{III}}(z)$ and which is again a rational function. Substituting \eqref{fornunqueanh} in the residue computation, the holomorphic term \eqref{holoterm} does not contribute and we are left with
\[
\partial_z f_h^{\textnormal{III}}(z) = \Res_{\zeta = \alpha_h,\beta_h} \left[\frac{\partial_{\zeta} f_h^{\textnormal{III}}(\zeta)}{z - \zeta} \right] = -\left[\oint_{\gamma(\alpha_h)}+\oint_{\gamma(\beta_h)}\right] \frac{\dd\zeta}{2\ii\pi} \frac{f_h^{\textnormal{III}}(\zeta)}{(z - \zeta)^2},
\]
where we used integration by parts and wrote the residue as integral over small loops (this is valid when $z$ remains outside the loops). This yields the desired decomposition \eqref{sourcedec} for the source:
\[
D_h(z) = \frac{r(\alpha_h)}{z - \alpha_h} + \frac{r(\beta_h)}{z - \beta_h} + \underbrace{\oint_{\gamma_h^+} \frac{\dd\zeta}{2\ii\pi} \frac{f_h^{\textnormal{I}}(\zeta)}{(z - \zeta)^2}}_{\partial_z f_h^{\textnormal{I}}(z)} \,\,\, \underbrace{ - \left[\oint_{\gamma(\alpha_h) }+\oint_{\gamma(\beta_h)} \right] \frac{\dd\zeta}{2\ii\pi} \frac{f^{\textnormal{III}}_h(\zeta)}{(z - \zeta)^2}}_{\partial_z f_h^{\textnormal{III}}(z)}.
\]

Let us now introduce
\[
F_h(z) := \sum_{g = 1}^H \big(\kappa_g \mathfrak{c}^{\textnormal{1st}}_{h;g}(z) - r(\alpha_g) \mathfrak{c}_{h;g}^{\textnormal{3rd}}(z;\alpha_g) - r(\beta_g) \mathfrak{c}^{\textnormal{3rd}}_{h;g}(z;\beta_g)\big) - \sum_{\textnormal{J} \in \{\textnormal{I},\textnormal{II},\textnormal{III}\}} \mathfrak{F}_h^{\textnormal{J}}[\boldsymbol{f}^{\textnormal{J}}](z),
\]
for a choice of $H$-tuple of functions $\boldsymbol{f}^{\textnormal{II}}(z)$ yet to be specified. The properties of each term are described in Lemma~\ref{lem:genRHPsol}. We find that $\boldsymbol{F}(z) = (F_h(z))_{h = 1}^H$ satisfies the master Riemann--Hilbert problem with the desired source, that is for any $h \in [H]$ and $x \in \amsmathbb{B}_h$
\[
\theta_{h,h}\big(F_h(x^+) + F_h(x^-) \big) + \sum_{h' \neq h} 2\theta_{h,h'} F_{h'}(x) = D_h(x).
\]
Note that the functions of first kind and second kind type II do not contribute to the source. For the behavior at $z \rightarrow \infty$, only the first-kind functions contribute and we get $F_h(z) = \frac{\kappa_h}{z} + O(\frac{1}{z^2})$. For the behavior near $z \rightarrow \alpha_h,\beta_h$, the functions of first kind and second kind type I are tame with exponent $1$, while the other terms create divergences of higher-order. Namely, $(F_h(z) - F_h^{\textnormal{sing}}(z))_{h = 1}^H$ is tame with exponent $1$, where for any $h \in [H]$
\begin{equation*}
\begin{split}
F_h^{\textnormal{sing}}(z) & = \frac{1}{2\theta_{h,h}} \left(\frac{r(\alpha_h)}{z - \alpha_h} + \frac{r(\beta_h)}{z - \beta_h} - \left[\oint_{\gamma(\alpha_h)}+\oint_{ \gamma(\beta_h)}\right] \frac{\dd\zeta}{2\ii\pi}\frac{f_h^{\textnormal{III}}(\zeta)}{(z - \zeta)^2}\right) - \partial_z\left(\frac{\sigma_h(z)f_h^{\textnormal{II}}(z)}{2}\right) \\
& = \frac{D_{h}(z)}{2\theta_{h,h}} + O(1) - \partial_z\left(\frac{\sigma_h(z)f_h^{\textnormal{II}}(z)}{2}\right),
\end{split}
\end{equation*}
and the $O(1)$ is holomorphic in $\amsmathbb{M}_h$. The last term can be made equal to the prescribed $\frac{L_h(z)}{\sigma_h(z)}$ by choosing
\[
f^{\textnormal{II}}_h(z) = \frac{2}{\sigma(z)} \int_{z}^{\infty} \frac{L_h(\zeta)}{\sigma_h(\zeta)} \dd\zeta. \qedhere
\]
\end{proof}

\begin{proof}[Proof of Corollary~\ref{cor:1stkindapp}] According to Proposition~\ref{prop:tameRHP}, the $H$-tuple of functions $(\Upsilon_h[\boldsymbol{0}\,;\,\boldsymbol{\kappa}](z))_{h = 1}^H$ solves the master Riemann--Hilbert problem with zero source and is tame with exponent $1$. By definition in Theorem~\ref{Theorem_Master_equation_12}, we have
\[
\oint_{\gamma_h} \frac{\dd z}{2\ii\pi} \Upsilon_h[\boldsymbol{0}\,;\,\boldsymbol{\kappa}](z) = \kappa_g.
\]
The $H$-tuple of functions $(\sum_{g = 1}^H \kappa_g \mathfrak{c}^{\textnormal{1st}}_{h;g}(z))_{h = 1}^H$ has the same properties. We conclude it must be equal to $(\Upsilon_h[\boldsymbol{0}\,;\,\boldsymbol{\kappa}](z))_{h = 1}^H$ by the uniqueness part in Proposition~\ref{prop:tameRHP}.
\end{proof}

\section{Fundamental solution for special \texorpdfstring{$\boldsymbol{\Theta}$}{Theta}}

\label{Section_Explicit}

\subsection{Explicit formulae}
\label{Sec1241}
In Theorem~\ref{Theorem_Masterspecial} we gave an explicit formula for the solution operator $\boldsymbol{\Upsilon}$ for two special cases: $\boldsymbol{\Theta}$ diagonal, or filled with the same positive number. In this section, we use it to obtain the corresponding formulae for the fundamental solution of Definition~\ref{def:Berg}. These formulae are familiar in the study of matrix models \cite{ACM92,Ake96}. These results will also be indirectly useful in Chapter~\ref{Chapter_AG}.

\begin{proposition}
\label{W20diagonal} If $\boldsymbol{\Theta}$ is diagonal, we have
\begin{equation}
\label{Bdiagth}
\mathcal{F}_{h,g}(z,w) = \frac{\delta_{h,g}}{2\theta_{h,h}(z - w)^2}\left(-1 + \frac{zw - \frac{\alpha_h + \beta_h}{2}(z + w) + \alpha_h\beta_h}{\sigma_h(z)\sigma_h(w)}\right).
\end{equation}
\end{proposition}

Like in Theorem~\ref{Theorem_Masterspecial}, the case of $\theta_{g,h} = \theta$ for all $g,h \in [H]$ can be obtained explicitly up to inversion of the restricted period map $\hat{\Pi} : \amsmathbb{C}_{H - 2}[z] \rightarrow \amsmathbb{C}^{H - 1}$, defined by
\[
\hat{\Pi}[P] = \left(\oint_{\gamma_h} \frac{\dd\zeta}{2\ii\pi}\,\frac{P(\zeta)}{\sigma(\zeta)}\right)_{h = 1}^{H - 1}
\]
This map is invertible by the arguments employed for Lemma~\ref{Pinvdet}. We use the notation
\[
\mathcal{F}(z,w) = \sum_{h,g = 1}^H \mathcal{F}_{h,g}(z,w),\qquad S(z) := \sigma^2(z) = \prod_{h = 1}^H (z - \alpha_h)(z - \beta_h).
\]
\begin{proposition}
\label{W20second} If $\theta_{g,h} = \theta > 0$ for all $g,h \in [H]$, we have
\begin{equation}
\label{W20equal}
\begin{split}
\mathcal{F}(z,w) & = \frac{(\sigma(z) - \sigma(w))^2 - \frac{z - w}{2}(S'(z) - S'(w))}{4\theta(z - w)^2\sigma(z)\sigma(w)} \\
& \quad + \frac{P(z,w) + \sum_{j,i = 0}^{H - 2} c_{j,i}z^{j}w^{i}}{\theta\sigma(z)\sigma(w)}.
\end{split}
\end{equation}
Here, $P(z,w)$ is the unique symmetric polynomial of degree $2H - 2$ in both its variable such that each of its monomial contains at least $z^{H - 1}$ or $w^{H - 1}$, and
\begin{equation}
\label{equationPP}
\frac{(\sigma(z) - \sigma(w))^2 - \frac{z - w}{2}(S'(z) - S'(w))}{4(z - w)^2} + P(z,w) = O(z^{H - 2}w^{H - 2})
\end{equation}
as $z,w \rightarrow \infty$. The constants $c_{j,i}$ are given in terms of the inverse of the restricted period map
\begin{equation}
\label{equationc}\forall i \in \llbracket 0,H-2\rrbracket\qquad \sum_{j = 0}^{H - 2} c_{j,i}z^{j} = \hat{\Pi}^{-1}\Bigg[\bigg(\oint_{\gamma_h} \frac{\dd\zeta}{2\ii\pi} \frac{Q_{i}(\zeta)}{\sigma(\zeta)}\bigg)_{h = 1}^{H - 1}\Bigg],
\end{equation}
where $Q_{i}(\zeta)$ are polynomials of degree at most $2H - 2$, defined by the decomposition as $w \rightarrow \infty$
\begin{equation}
\label{equationc2}
\frac{S(z) + S(w) - \frac{z - w}{2}(S'(z) - S'(w))}{4(z - w)^2} = - \sum_{i = 0}^{2H - 2} Q_{i}(z)w^{i} + O\bigg(\frac{1}{w}\bigg).
\end{equation}
\end{proposition}

We remark that $P(z,w)$ and the term coming from $\frac{S'(z) - S'(w)}{z - w}$ are polynomials and in practice could be combined. Before coming to the proof, we illustrate the result for low values of $H$. For $H = 1$ it does coincide with \eqref{Bdiagth}. Let us introduce the elementary symmetric functions \label{index:elemsym}
\[
\mathsf{e}_k = \sum_{1 \leq i_1 < \cdots < i_k \leq 2H} \lambda_{i_1} \cdots \lambda_{i_k},\qquad \textnormal{where}\quad \{\lambda_i\,\,|\,\,i \in [2H]\} = \{\alpha_1,\beta_1,\ldots,\alpha_H,\beta_H\}.
\]
\begin{corollary}
\label{CoH2fund} If $H = 2$ and $\theta_{g,h} = \theta > 0$ for all $g,h \in \{1,2\}$, we have
\begin{equation}
\label{Hequal2}
\mathcal{F}(z_1,z_2) = \frac{1}{4\theta \sigma(z_1)\sigma(z_2)}\left(\bigg(\frac{\sigma(z_1) - \sigma(z_2)}{z_1 - z_2}\bigg)^2 - (z_1 + z_2)(z_1 + z_2 - \mathsf{e}_1) + 4c_{0,0} - \mathsf{e}_2\right).
\end{equation}
where
\begin{equation}
\label{claimedc}
c_{0,0} = \frac{1}{4}\left(\alpha_1\beta_2 + \alpha_2\beta_1 + (\beta_2 - \beta_1)(\alpha_2 - \alpha_1)\frac{E(\mathsf{k})}{K(\mathsf{k})}\right)
\end{equation}
in terms of the complete elliptic integrals
\begin{equation*}
\begin{split}
E(\mathsf{k}) = \int_{0}^{1} \sqrt{\frac{1 - \mathsf{k}^2t^2}{1 - t^2}},\qquad K(\mathsf{k}) = \int_{0}^{1} \frac{\dd t}{\sqrt{(1 - t^2)(1 - \mathsf{k}^2t^2)}}.
\end{split}
\end{equation*}
specialized to
\[
\mathsf{k} = \sqrt{\frac{(\beta_2 - \alpha_1)(\alpha_2 - \beta_1)}{(\beta_2 - \beta_1)(\alpha_2 - \alpha_1)}} < 1,
\]
\end{corollary}
\begin{proof}
From \eqref{equationPP} we find
\[
P(z,w) = \frac{z^2 + w^2}{4} + \frac{\mathsf{e}_1(z + w)}{8}.
\]
On the other hand
\[
-\frac{S'(z) - S'(w)}{8(z - w)} = - \frac{z^2 + zw + w^2}{2} + \frac{3\mathsf{e}_1(z + w)}{8} - \frac{\mathsf{e}_2}{4}.
\]
Adding these two polynomials gives the polynomial term in \eqref{Hequal2}, up to the constant $c_{0,0}$ which is defined by
\[
c_{0,0} = \frac{\oint_{\gamma_1} \frac{Q_0(z)\dd z}{\sigma(z)}}{\oint_{\gamma_1} \frac{\dd z}{\sigma(z)}} = \frac{\int_{\alpha_1}^{\beta_1} \frac{Q_0(x)\dd x}{|\sigma(x)|}}{\int_{\alpha_1}^{\beta_1} \frac{\dd x}{|\sigma(x)|}}.
\]
The method of reduction of these integrals to elliptic integrals can be found in \cite[p 550-551]{Speci}. The result can be found in \cite{Ake96}\footnote{Note that \cite[Equation B.6]{Ake96} can be heavily simplified: it is equal to $(x_1 - x_3)(x_2 - x_4)$, and in our notations $x_1 = \beta_2$, $x_2 = \alpha_2$, $x_3 = \beta_1$, $x_4 = \alpha_1$.} and leads to \eqref{claimedc}.
\end{proof}

We bother writing the case $H = 3$ without detailing the proof, as it shows up indirectly for the example of Section~\ref{H38shap}.

\begin{corollary}
\label{corHequal3} If $H = 3$ and $\theta_{g,h} = \theta > 0$ for all $g,h \in \{1,2,3\}$, we have
\[
\mathcal{F}(z,w) = \frac{1}{4\theta\sigma(z)\sigma(w)}\left(\frac{\sigma(z) - \sigma(w)}{z - w}\right)^2 + \frac{\hat{P}(z,w) + c_{0,0} + c_{0,1}(z + w) + c_{1,1}zw}{\theta\sigma(z)\sigma(w)}.
\]
with the polynomials
\begin{equation}
\label{Phatequation} \begin{split}
\hat{P}(z,w) & = \frac{1}{4}\left(-(z^2 + wz + w^2)^2 + \mathsf{e}_1(z + w)(z^2 + zw + w^2) - \mathsf{e}_2(z + w)^2 + \frac{3\mathsf{e}_3}{2}(z + w) - \mathsf{e}_4\right), \\
Q_0(z) & = \frac{z^4}{2} - \frac{3\mathsf{e}_1z^3}{8} + \frac{\mathsf{e}_2z^2}{4} - \frac{\mathsf{e}_3z}{8}, \qquad Q_1(z) = \frac{z^3}{4} - \frac{\mathsf{e}_1 z^2}{8} + \frac{\mathsf{e}_3}{8}.
\end{split}
\end{equation}
The remaining constants are solution of the system of equation for $h \in \{1,2\}$:
\begin{equation}
\oint_{\gamma_h} \frac{(c_{0,0} + c_{0,1}z)\dd z}{\sigma(z)} = \oint_{\gamma_h} \frac{Q_0(z) \dd z }{\sigma(z)} \qquad \textnormal{and}\qquad \oint_{\gamma_h} \frac{(c_{0,1} + c_{1,1}z) \dd z}{\sigma(z)} = \oint_{\gamma_h} \frac{Q_1(z)\dd z}{\sigma(z)}.
\end{equation}
\end{corollary}

\subsection{Proof of Propositions~\ref{W20diagonal} and \ref{W20second}}

In both cases, the starting point is Theorem~\ref{Theorem_Masterspecial} with $E_h(z) = \frac{\delta_{h,g}\sigma_g(\zeta)}{2(z - w)^2}$ and $\kappa_h = 0$, which gives access to $\mathcal{F}_{h,g}(z,w)$.

\begin{proof}[Proof of Proposition~\ref{W20diagonal}] Assume $\theta_{g,h} = 0$ for $g \neq h$. We then find:
\[
\mathcal{F}_{h,g}(z,w) = \frac{1}{\sigma_h(z)} \oint_{\gamma_h} \frac{\dd\zeta}{2\ii\pi}\frac{\delta_{h,g}\sigma_g(\zeta)}{2\theta_{h,h}(\zeta - z)(\zeta - w)^2},
\]
where $z,w$ are outside the contour of integration. We move the contour to $\infty$. There is no residue at $\infty$ but we pick the poles at $z$ and $w$, resulting in
\[
\mathcal{F}_{h,g}(z,w) = \frac{\delta_{g,h}}{2\theta_{h,h}}\partial_w\bigg(\frac{-1 + \frac{\sigma_h(w)}{\sigma_h(z)}}{z - w}\bigg) = \frac{\delta_{g,h}}{2\theta_{h,h}(z - w)^2}\left(-1 + \frac{\sigma_h(w)^2 + \sigma_h'(w)\sigma_h(w)(z - w)}{\sigma_h(z)\sigma_h(w)}\right).
\]
Inserting $\sigma_h(z)^2 = (z - \alpha_h)(z - \beta_h)$, it simplifies to the announced formula.
\end{proof}

\begin{proof}[Proof of Proposition~\ref{W20second}] Assume $\theta_{h,g} = \theta$ for any $h,g \in [H]$. We then find:
\[
\mathcal{F}(z,w) = \frac{1}{2\theta \sigma(z)} \sum_{h = 1}^H \oint_{\gamma_h}
\frac{\dd\zeta}{2\ii\pi}\,\frac{\sigma(\zeta)}{(\zeta - z)(\zeta - w)^2} +
\frac{P_0(z;w)}{\sigma(z)},
\]
where $z,w$ are outside the integration contours and $P_0(z;w)$ is a polynomial function of $z$ of degree at most $(H - 1)$. Actually we know that $P_0(z;w)$ has degree at most $(H - 2)$, because $\sigma(z) \sim z^H$ and $\mathcal{F}(z,w) = O(z^{-2})$ as $z \rightarrow \infty$. We then move the contour to pick up the residues at $\zeta = z,w$ and $\infty$. This yields
\[
\mathcal{F}(z,w) = -\frac{1}{2\theta \sigma(z)} \partial_{w}\bigg(\frac{\sigma_{-}(z) - \sigma_{-}(w)}{z - w}\bigg) + \frac{P_0(z;w)}{\sigma(z)},
\]
where we have decomposed
\begin{equation}
\label{sigmaplusmoins}\sigma(z) := \sigma_{-}(z) + z\sigma_{+}(z)
\end{equation}
such that $\sigma_{-}(z) = O(1)$ as $z \rightarrow \infty$, and $\sigma_+(z)$ is a polynomial of degree $H - 1$. The residue at $\zeta = \infty$ is the difference coming from the replacement of $\sigma(z)$ with $\sigma_-(z)$. The fact that $\mathcal{F}(z,w) = O(w^{-2})$ as $w \rightarrow \infty$ imposes that $P_0(z;w) = O(w^{-2})$ as $w \rightarrow \infty$. We can also rewrite
\begin{equation}
\label{W20exp} \mathcal{F}(z,w) = \frac{S(w) - \sigma(z)\sigma(w) + \frac{z - w}{2}\,S'(w)}{2\theta(z - w)^2\sigma(z)\sigma(w)} + \frac{P_1(z;w)}{\sigma(z)},
\end{equation}
where $S(w) = \sigma(w)^2$ and
\begin{equation}
\label{P0P1} P_1(z;w) = \partial_{w}\bigg(\frac{z\sigma_+(z) - w\sigma_+(w)}{2\theta(z - w)}\bigg) + P_0(z;w)
\end{equation}
is still a polynomial function of $z$, of degree at most $(H - 2)$, and such that $P_1(z;w) = O(w^{H - 2})$ as $w \rightarrow \infty$ due to the first term. Now, we use the \textit{a priori} symmetry $\mathcal{F}(z,w) = \mathcal{F}(w,z)$ to obtain (\textit{cf.} Theorem~\ref{thm:Bsym})
\begin{equation}
\label{preP2}\frac{1}{2\theta\sigma(z)\sigma(w)}\cdot \frac{S(w) - S(z) + \frac{z - w}{2}(S'(w) + S'(z))}{(z - w)^2} = \frac{P_1(w;z)}{\sigma(w)} - \frac{P_1(z;w)}{\sigma(z)}.
\end{equation}
One checks that the second ratio in the left-hand side is a polynomial in $z$ and $w$. Therefore, we must have
\begin{equation}
\label{P1P2} P_1(z;w) = \frac{P_2(z;w)}{\sigma(w)}
\end{equation}
for some $P_2(z;w)$ which is a polynomial function of $z$ of degree at most $H - 2$ and of $w$ of degree at most $2H - 2$. We decompose it
\[
P_{2}(z;w) = \frac{1}{\theta}\big(P_{2}^{\textnormal{s}}(z,w) + P_{2}^{\textnormal{a}}(z;w)\big),\qquad P_{2}^{\textnormal{s}}(z,w) = P_2^{\textnormal{s}}(w,z),\qquad P_2^{\textnormal{a}}(z;w) = -P_2^{\textnormal{a}}(w;z),
\]
with $P_2^{\textnormal{a}/\textnormal{s}}(z,w)$ polynomials of degree at most $2H - 2$ in both of its variables. Equation \eqref{preP2} then prescribes
\begin{equation}
\label{P2moins} P_{2}^{\textnormal{a}}(z;w) = \frac{1}{4}\,\frac{S(z) - S(w) - \frac{z - w}{2}(S'(z) + S'(w))}{(z - w)^2},
\end{equation}
and coming back to \eqref{W20exp} we find
\[
\mathcal{F}(z,w) = \frac{(\sigma(z) - \sigma(w))^2 -\frac{z - w}{2}(S'(z) - S'(w))}{4\theta \sigma(z)\sigma(w) (z - w)^2} + \frac{P_{2}^{\textnormal{s}}(z,w)}{\theta \sigma(z)\sigma(w)}.
\]
We can further determine part of $P_2^{\textnormal{s}}(z,w)$. Indeed, let us decompose
\[
P_2^{\textnormal{s}}(z,w) := \tilde{P}(z;w) + P(z;w)
\]
where $\tilde{P}(z;w)$ has degree at most $H - 2$ in both of its variables, and $P(z;w)$ has degree at most $2H - 2$ in both its variables but each of its monomials contain at least a power of $z^{H - 1}$ or of $w^{H - 1}$. The requirement that $\mathcal{F}(z,w) = O(z^{-2}w^{-2})$ as $z,w \rightarrow \infty$ imposes
\[
\frac{(\sigma(z) - \sigma(w))^2 - \frac{z - w}{2}(S'(z) - S'(w))}{4(z - w)^2} + P(z;w) = O(z^{H - 2}w^{H - 2}),
\]
which allows the computation of $P(z;w)$. In particular, $P(z;w) = P(w;z)$ and is simply denoted $P(z,w)$: this yields the desired \eqref{equationPP}. We still have to determine
\[
\tilde{P}(z;w) := \sum_{j,i = 0}^{H - 2} c_{j,i}z^{j}w^{i}.
\]
For this purpose, we impose the transcendental constraints $\oint_{\gamma_h} \mathcal{F}(z,w)\dd z = 0$ for any $h \in [H - 1]$ and $w$ outside the contour of integration. The remaining integral around $\gamma_H$ automatically vanishes because $\mathcal{F}(z,w) = O(z^{-2})$ has no residue as $z \rightarrow \infty$. We remark that the term $-\tfrac{1}{2\theta(z - w)^2}$ coming from the cross-term in $(\sigma(z) - \sigma(w))^2$ disappears in the contour integrals. It is actually enough to impose for any $h \in [H - 1]$ and $i \in \llbracket 0,H-2\rrbracket$
\begin{equation}
\label{forihH}
\forall i,h \in [H - 1]\qquad \oint_{\gamma_h} \dd z \left(\Res_{w = \infty} \frac{\dd w}{w^{i + 1}} \sigma(w)\mathcal{F}(z,w) \right) = 0.
\end{equation}
Since $i \leq H - 2$ the contribution of $P(z,w)$ also disappears from \eqref{forihH}. Taking the previous remark into account, it yields
\[
\oint_{\gamma_h} \bigg(\Res_{w = \infty} \frac{\dd w}{w^{i + 1}}\,\frac{S(z) + S(w) - \frac{z - w}{2}(S'(z) - S'(w))}{4(z - w)^2}\bigg)\,\frac{\dd z}{\sigma(z)} - \oint_{\gamma_h} \frac{\dd z}{\sigma(z)}\bigg(\sum_{j = 0}^{H - 2} c_{j,i}z^{j}\bigg) = 0,
\]
and entails \eqref{equationc}-\eqref{equationc2}.
\end{proof}

\chapter{Constructing the spectral curve}

\label{Chapter_AG}

We will construct \emph{spectral curves} associated with the master Riemann--Hilbert problem. This is a Riemann surface depending only on $\alpha_1,\beta_1,\ldots,\alpha_H,\beta_H$ and $\boldsymbol{\Theta}$, with the property that solutions of the master Riemann--Hilbert problem with meromorphic source can be analytically continued to the spectral curve. More precisely, we exploit the freedom to consider linear combinations $\langle \boldsymbol{v} \cdot \boldsymbol{F}(z)\rangle = \sum_{h = 1}^H v_h F_h(z)$, and different choices of $\boldsymbol{v} \in \amsmathbb{C}^H$ lead to different spectral curves. The master Riemann--Hilbert problem gives a direct access to the monodromy of such solutions. It is possible to handle all $\boldsymbol{v}$ at once by introducing a subgroup $\mathfrak{G} \subset \textnormal{GL}(\amsmathbb{R}^H)$ generated by (pseudo-)reflections depending only on $\boldsymbol{\Theta}$.

For special matrices $\boldsymbol{\Theta}$ it may happen that $\mathfrak{G}$ admits a finite orbit. When this is the case, we have vectors $\boldsymbol{v}$ such that the function t $\langle \boldsymbol{v} \cdot \boldsymbol{F}(z) \rangle$ is algebraic and the corresponding spectral curve is compact. This often leads to formulae for the fundamental solution $\boldsymbol{\mathcal{F}}(z,w)$ and accordingly for the solution operator $\boldsymbol{\Upsilon}$. In contrast, for general $\boldsymbol{\Theta}$, $\boldsymbol{\Upsilon}$ and $\boldsymbol{\mathcal{F}}(z,w)$ were constructed in Chapter~\ref{Chapter_SolvingN} only as series via Fredholm theory.

Among the specially nice cases are matrices $\boldsymbol{\Theta}$ which are diagonal, or filled with $1$s, or coming from tiling models. The latter are studied in greater detail in Section~\ref{sec:C13tiling}, unveiling a close relation between the topology of the spectral curve and the topology of the tiled domain on the one hand, and a close relation between the fundamental solution and the Green function of the Laplace operator on the spectral curve on the other hand. These relations play a key role in the proof of the Kenyon--Okounkov conjecture in Chapter~\ref{Chap11}.

The procedure to obtain (in good situations, compact) spectral curves has been used in a case-by-case approach in the theoretical physics literature about matrix models since the early 90s. An early occurrence of a master Riemann--Hilbert problem in the sense of Definition~\ref{def:masterRHP} is in the works of Kostov with $\boldsymbol{\Theta}$ being half the Cartan matrix of a root system (and segments are related by symmetries) \cite{KosADE0,KosADE}. The $O(n)$-matrix model \cite{GaudinKostov} features the Riemann--Hilbert problem with $H = 2$ with a symmetry between the two segments, and (non-)algebraicity and construction of general solutions was addressed in \cite{EynardZJ,EKOn}. An algebraic case without symmetry was discussed for the Cauchy two-matrix model \cite{BertolaCauchy}. A more systematic study of the master Riemann--Hilbert problem and algebraic cases appeared in \cite[Section 3]{BESeifert}, and this chapter is an adaptation and a further development of these techniques to study the master Riemann--Hilbert problem.

\section{Algebraic approach to the master Riemann--Hilbert problem}
\label{sec:Algapproach0}

\subsection{Spectral curves from monodromy}
\label{sec:Algapproach}

Assume that we are given $2H$ real numbers $\alpha_1 < \beta_1 < \cdots < \alpha_H < \beta_H$ and a real symmetric matrix $\boldsymbol{\Theta}$ with positive diagonal. As in Section~\ref{sec:1233} we denote $\amsmathbb{B}_h = (\alpha_h,\beta_h)$ and $\amsmathbb{B} = \bigcup_{h = 1}^H \amsmathbb{B}_h$, We recall that $\gamma_h$ denotes a contour in $\amsmathbb{C} \setminus \overline{\amsmathbb{B}}$ surrounding $\amsmathbb{B}_h$ but not $\amsmathbb{B}_g$ for $g \neq h$. For each $h \in [H]$ we recall the notation $\bth_{h} = (\theta_{h,g})_{g = 1}^H \in \amsmathbb{R}^H$ . If $\boldsymbol{w} \in \amsmathbb{R}^H$, we denote $\textnormal{supp}(\boldsymbol{w}) = \big\{h \in [H]\,\,|\,\,w_h \neq 0 \big\}$.

\begin{definition} \label{unimono} To these data, we associate the subgroup $\mathfrak{G} \subset \textnormal{GL}(\amsmathbb{R}^H)$ generated by $T^{(1)},\ldots,T^{(H)}$ defined as
\[
\forall h \in [H]\qquad \forall \boldsymbol{v} \in \amsmathbb{R}^H\qquad T^{(h)}(\boldsymbol{v}) = \boldsymbol{v} - \frac{2v_h}{\theta_{h,h}} \bth_h.
\]
\end{definition}

This is a pseudo-reflection group: the transformation $T^{(h)}$ is involutive and fixes the hyperplane $\{v_h = 0\}$. In particular, for any non-zero $\boldsymbol{v} \in \amsmathbb{R}^H$ and $h \in [H]$, we have $h \in \textnormal{supp}(\boldsymbol{v})$ if and only if $T^{(h)}(\boldsymbol{v}) \neq \boldsymbol{v}$. We will be interested in the orbit of $\boldsymbol{v} \in \amsmathbb{R}^H$ under the group $\mathfrak{G}$, which is denoted $\mathfrak{G}.\boldsymbol{v}$.

\begin{definition}
\label{def:spcurv1}
The \emph{spectral curve} associated to the above data and a non-zero $\boldsymbol{v} \in \amsmathbb{R}^H$ is the Riemann surface
\begin{equation}
\label{wcopydoubleplae}
\Sigma_{\boldsymbol{v}} := \left(\bigsqcup_{\boldsymbol{w} \in \mathfrak{G}.\boldsymbol{v}} \widehat{\amsmathbb{H}}^{+} \sqcup \widehat{\amsmathbb{H}}^-\right) \Big/ \sim.
\end{equation}
Here, $\widehat{\amsmathbb{H}}^{\pm} \subset \widehat{\amsmathbb{C}}$ is the upper (respectively, lower) hemisphere of the Riemann sphere in which we include $\widehat{\amsmathbb{R}}^{\pm} := \amsmathbb{R} \cup \{\infty\} = \partial \widehat{\amsmathbb{H}}^{\pm}$. We call $\amsmathbb{B}^\pm$ the segment $\amsmathbb{B}$ inside $\widehat{\amsmathbb{R}}^{\pm}$, and likewise for $\overline{\amsmathbb{B}}_{h}^{\pm}$. The equivalence relation consists in identifying
\begin{itemize}
\item $\widehat{\amsmathbb{R}}^+ \setminus \amsmathbb{B}^+$ with $\widehat{\amsmathbb{R}}^- \setminus \amsmathbb{B}^-$ in the same copies;
\item $\overline{\amsmathbb{B}}_{h}^{+}$ in the $\boldsymbol{w}$-th copy with $\overline{\amsmathbb{B}}_{h}^{-}$ in the $T^{(h)}(\boldsymbol{w})$-th copy, for each $\boldsymbol{w} \in \mathfrak{G}.\boldsymbol{v}$ and $h \in [H]$.
\end{itemize}
For $\boldsymbol{w} \in \mathfrak{G}.\boldsymbol{v}$ we denote
\begin{equation}
\label{Bamhbw}
\overline{\amsmathbb{B}}_{\boldsymbol{w}}^{\pm} = \bigcup_{g \in \textnormal{supp}(\boldsymbol{w})} \overline{\amsmathbb{B}}_{g}^{\pm} .
\end{equation}
We use the same notation for the image of \eqref{Bamhbw} in $\Sigma_{\boldsymbol{v}}$, and call
\begin{equation}
\label{wgluloc}
\overline{\amsmathbb{B}}_{\boldsymbol{w}} := \overline{\amsmathbb{B}}_{\boldsymbol{w}}^+ \cup \overline{\amsmathbb{B}}_{\boldsymbol{w}}^- \subset \Sigma_{\boldsymbol{v}}.
\end{equation}
the $\boldsymbol{w}$-th \emph{gluing locus}. The $\boldsymbol{w}$-th \emph{sheet}, denoted $\Sigma_{\boldsymbol{v},\boldsymbol{w}}$ is the image in $\Sigma_{\boldsymbol{v}}$ of the $\boldsymbol{w}$-th copy of $\widehat{\amsmathbb{H}}^+ \sqcup \widehat{\amsmathbb{H}}^-$ in \eqref{wcopydoubleplae}, from which we remove the $\boldsymbol{w}$-th gluing locus \eqref{wgluloc}. The $\boldsymbol{v}$-th sheet is called the \emph{principal sheet}.
\end{definition}
\begin{remark} Contrarily to Section~\ref{Section_complex_structure}, we glue here Riemann spheres and not half-planes. We will return to the half-construction in Definition~\ref{def:halfspcurvebip}.
\end{remark}

\begin{figure}[h!]
\begin{center}
\includegraphics[width=0.95\textwidth]{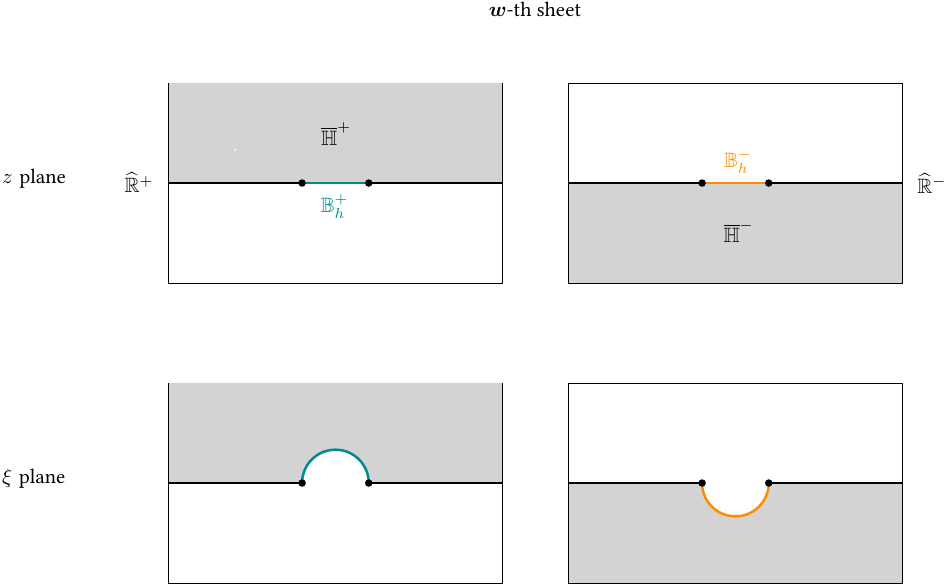}
\caption{\label{Fig:Half1} The two upper half-planes with the designated $h$-th segment and their image via the Zhukovsky map \eqref{Joukov}. The point at infinity is hard to depict here but it is included to get hemispheres $\widehat{\amsmathbb{H}}^{\pm}$.}
\end{center}
\end{figure}

Roughly speaking, the equivalence relation glues together copies of the Riemann sphere along some of the segments $\overline{\amsmathbb{B}}_h$. A more accurate description is that it opens up the slits $\amsmathbb{B}_h$ to form a topological circle, whose upper and lower parts correspond to $\amsmathbb{B}_h \pm \ii 0$, and then glues the upper part in the $\boldsymbol{w}$-th sheet to the lower part in the $T^{(h)}(\boldsymbol{w})$-th sheet. When $w_h = 0$, we have $T^{(h)}(\boldsymbol{w}) = \boldsymbol{w}$ so the result after gluing is the same as if the $h$-th slit was ignored; otherwise, we are gluing two distinct sheets (\textit{cf.} Figure~\ref{Fig:Half1}, and the example in Figure~\ref{Fig:Half4}).

\label{index:Zbv}The natural coordinate on the Riemann sphere induces a branched covering map $Z_{\boldsymbol{v}} : \Sigma_{\boldsymbol{v}} \rightarrow \widehat{\amsmathbb{C}}$ consisting in forgetting the sheet label. It admits simple ramification points at the points of the gluing locus that are endpoints of the identified segments. Its branchpoints are among $\alpha_1,\beta_1,\ldots,\alpha_H,\beta_H$. Note that the structure of the group $\mathfrak{G}$, its orbits and the topology of the spectral curve only depend on $\boldsymbol{\Theta}$, while the complex structure of the spectral curve depends as well on $\alpha_1,\beta_1,\ldots,\alpha_H,\beta_H$.

This construction naturally appears when studying the solutions of the master Riemann--Hilbert problem of Definition~\ref{def:masterRHP} with zero source.

\begin{theorem}
\label{spcurveRHP}Consider a solution $\boldsymbol{F}(z) = (F_h(z))_{h = 1}^H$ of the master Riemann--Hilbert problem with zero source, such that $\boldsymbol{F}(z) = O(z^d)$ as $z \rightarrow \infty$ for some $d \in \amsmathbb{Z}$. Then, for any non-zero $\boldsymbol{v} \in \amsmathbb{R}^H$, there exists a unique meromorphic function $F_{\boldsymbol{v}}$ on $\Sigma_{\boldsymbol{v}}$ such that, for any $\boldsymbol{w} \in \mathfrak{G}.\boldsymbol{v}$ and $p \in \Sigma_{\boldsymbol{v},\boldsymbol{w}}$ we have
\[
F_{\boldsymbol{v}}(p) = \big \langle \boldsymbol{w} \cdot \boldsymbol{F}(z) \big \rangle = \sum_{h = 1}^H w_h F_h(z) \qquad \textnormal{with}\quad z = Z_{\boldsymbol{v}}(p).
\]
The function $F_{\boldsymbol{v}}$ can have poles only at ramification points with order at most equal to the tameness exponent of $\boldsymbol{F}$, and at $Z^{-1}_{\boldsymbol{v}}(\{\infty\})$ with order at most $d$. If the $\mathfrak{G}$-orbit of $\boldsymbol{v}$ is finite, then $\Sigma_{\boldsymbol{v}}$ is a compact Riemann surface, and
\[
\big\{\langle \boldsymbol{w} \cdot \boldsymbol{F}(z)\rangle \,\,\big|\,\,\boldsymbol{w} \in \mathfrak{G}.\boldsymbol{v}\big\}
\]
are the various branches of an algebraic function.
\end{theorem}
\label{index:monodromy}The assumption on the behavior at $\infty$ can be dropped: in that case we only have a meromorphic function on $Z_{\boldsymbol{v}}^{-1}(\amsmathbb{C}) \subset \Sigma_{\boldsymbol{v}}$. The group $\mathfrak{G}$ is related but not always equal to the monodromy group $\mathfrak{M}_{\boldsymbol{v}}$ of the function $\langle \boldsymbol{v} \cdot \boldsymbol{F}(z) \rangle$. By definition, $\mathfrak{M}_{\boldsymbol{v}}$ is the image of the monodromy representation of the fundamental group of $\widehat{\amsmathbb{C}} \setminus \{\alpha_1,\beta_1,\ldots,\alpha_H,\beta_H\}$ in the permutation group of a generic fiber of the map $Z_{\boldsymbol{v}}$. Given the construction of our branched covering $Z_{\boldsymbol{v}}$, the group $\mathfrak{M}_{\boldsymbol{v}}$ is the image of $\mathfrak{G}$ in the permutation group of the set $\mathfrak{G}.\boldsymbol{v}$. Equivalently, $\mathfrak{M}_{\boldsymbol{v}}$ is the quotient of $\mathfrak{G}$ by the kernel of the action, \textit{i.e.} by the intersection of all stabilizers of elements of $\mathfrak{G}.\boldsymbol{v}$. Since $\mathfrak{G} \subset \textnormal{GL}(\amsmathbb{R}^H)$ acts linearly, $\mathfrak{M}_{\boldsymbol{v}}$ is the quotient of $\mathfrak{G}$ by the kernel of the restriction map $\mathfrak{G} \rightarrow \textnormal{GL}(\textnormal{span}_{\amsmathbb{R}}(\mathfrak{G}.\boldsymbol{v}))$. If the action of $\mathfrak{G}$ on $\amsmathbb{R}^H$ happens to be irreducible, the span of all vectors in any given non-zero $\mathfrak{G}$-orbit must be the full space, therefore $\mathfrak{M}_{\boldsymbol{v}} = \mathfrak{G}$ independently of $\boldsymbol{v}$.

\begin{definition} \label{def:goodor} An orbit $\mathfrak{G}.\boldsymbol{v}$ is good if for any $h \in [H]$, there exists $\boldsymbol{w} \in \mathfrak{G}.\boldsymbol{v}$ such that $w_h \neq 0$.
\end{definition}

If $\mathfrak{G}.\boldsymbol{v}$ is a good orbit, we can reconstruct the individual components of $\boldsymbol{F}(z)$ from the knowledge of $F_{\boldsymbol{v}}$. For each $h \in [H]$, this is done via the projection formula
\[
F_h(z) = \frac{1}{w_h} \oint_{\gamma_h} \frac{\dd \zeta}{2\ii\pi} \frac{\big\langle \boldsymbol{w} \cdot \boldsymbol{F}(\zeta)\big\rangle}{z - \zeta},
\]
using any $\boldsymbol{w} \in \mathfrak{G}.\boldsymbol{v}$ such that $w_h \neq 0$. In this formula, $z$ is outside the integration contour $\gamma_h$. For instance, if the action of $\mathfrak{G}$ on $\amsmathbb{R}^H$ is irreducible, non-zero orbits are always good.

\begin{proof}[Proof of Theorem~\ref{spcurveRHP}]
Let $\boldsymbol{F}(z)$ be a solution of the master Riemann--Hilbert problem with zero source. The analytic properties required for $\boldsymbol{F}(z)$ show that for any $\boldsymbol{w} \in \amsmathbb{R}^H$, the function $\big\langle\boldsymbol{w}\cdot\boldsymbol{F}(z)\big\rangle$ is holomorphic in $\amsmathbb{C} \setminus \bigcup_{h \in \textnormal{supp}(\boldsymbol{w})} \overline{\amsmathbb{B}}_h$. The master Riemann--Hilbert problem (Definition~\ref{def:masterRHP}) itself implies that for any $h \in \textnormal{supp}(\boldsymbol{w})$ and $x \in \amsmathbb{B}_h$ we have
\begin{equation}
\label{wcross}
\big\langle \boldsymbol{w} \cdot \boldsymbol{F}(x^+)\big\rangle = \big\langle T^{(h)}(\boldsymbol{w}) \cdot \boldsymbol{F}(x^-)\big\rangle.
\end{equation}
If $h \notin \textnormal{supp}(\boldsymbol{w})$ we have $T^{(h)}(\boldsymbol{w}) = \boldsymbol{w}$ so the equation \eqref{wcross} remains valid.

For a given $h \in \textnormal{supp}(\boldsymbol{w})$, consider the Zhukovsky map
\begin{equation}
\label{Joukov}
Z^{\textnormal{J}}(\xi) = \frac{\beta_h + \alpha_h}{2} + \frac{\beta_h - \alpha_h}{4}\left(\xi + \frac{1}{\xi}\right),
\end{equation}
and let $\amsmathbb{D} \subset \amsmathbb{C}$ be the unit disk centered at $0$. The map $Z$ sends the unit circle to $\overline{\amsmathbb{B}}_h = [\alpha_h,\beta_h]$, and $+1$ (respectively $-1$) to the point $\beta_h$ (respectively $\alpha_h$). Its restriction to $\amsmathbb{D}$, respectively to $\widehat{\amsmathbb{C}} \setminus \overline{\amsmathbb{D}}$, is a biholomorphic map to $\widehat{\amsmathbb{C}} \setminus \overline{\amsmathbb{B}}_h$ whose inverse is $\Xi_+(z)$, respectively $\Xi_-(z)$, with
\[
\Xi_\pm(z) = \frac{1}{2}\left(z - \frac{\alpha_h + \beta_h}{2} \pm \sigma_h(z)\right).
\]
Then, in the neighborhood of the unit circle in $\widehat{\amsmathbb{C}} \setminus \overline{\amsmathbb{D}}$ we have two holomorphic functions
\[
f_+(\xi) = \big\langle\boldsymbol{w} \cdot \boldsymbol{F}(Z^{\textnormal{J}}(\xi))\big\rangle,\qquad f_-(\xi) = \big\langle T^{(h)}(\boldsymbol{w}) \cdot \boldsymbol{F}(Z^{\textnormal{J}}(\xi))\big\rangle
\]
admitting a continuous extension to the unit circle except perhaps at $\{\pm 1\}$. Notice that $\Xi_{+}(z)\cdot \Xi_-(z) = 1$ for $z \in \widehat{\amsmathbb{C}}\setminus \overline{\amsmathbb{B}}_h$. Besides, for any $x \in \amsmathbb{B}_h$ we have $\Xi_{\pm}(x^+) = \frac{1}{\Xi_{\pm}(x^-)}$ because $\sigma_h$ takes a minus sign when $x$ crosses the segment from the upper to the lower half-plane. This construction then allows to rewrite \eqref{wcross} in the equivalent way
\[
\forall \xi \in \partial\amsmathbb{D} \qquad f_+(\xi) = f_-\left(\frac{1}{\xi}\right).
\]
Therefore, we can define a function $f_0(\xi)$ in the neighborhood of the unit circle in $\amsmathbb{C} \setminus \{\pm 1\}$ by declaring
\[
f_0(\xi) = \left\{\begin{array}{lll} f_+(\xi) && \textnormal{if}\,\, |\xi| \geq 1, \\ f_-(\xi) && \textnormal{if}\,\,|\xi| < 1. \end{array}\right.
\]
The tameness assumption implies that $\sigma_h(Z(\xi))^{m_{\boldsymbol{F}}} f_0(\xi) = (2\xi - Z(\xi) - \alpha_h - \beta_h)^{m_{\boldsymbol{F}}} f_0(\xi)$ is bounded near $\xi = \pm 1$, therefore $f_0$ is meromorphic with possible poles at $\xi = \pm 1$. The obstruction to extend $f_0$ as a holomorphic function away from the unit circle comes from the discontinuity of $F_g$ for $g \neq h$ on the segment $\overline{\amsmathbb{B}}_g$, that may be involved in $\langle\boldsymbol{w} \cdot \boldsymbol{F}\rangle$ and $\langle T^{(h)}(\boldsymbol{w}) \cdot \boldsymbol{F}\rangle$.

Let us bring back this construction to the $z$-plane. For this purpose we consider the topological space
\begin{equation}
\label{spsharp}
\Sigma^{\sharp,h} = \big(\widehat{\amsmathbb{H}}^+ \sqcup \widehat{\amsmathbb{H}}^-\big) \sqcup \big(\widehat{\amsmathbb{H}}^+ \sqcup \widehat{\amsmathbb{H}}^-\big)\big/ \sim
\end{equation}
where the equivalence relation identifies $\widehat{\amsmathbb{R}}^+ \setminus \amsmathbb{B}^+$ with $\widehat{\amsmathbb{R}}^+ \setminus \amsmathbb{B}^+$ in the same copies, and $\overline{\amsmathbb{B}}_{h}^+$ in one copy with $\overline{\amsmathbb{B}}_{h}^{-}$ in the other copy (the image of the latter in $\Sigma^{\sharp,h}$ is what we call the gluing locus). In each copy, the two hemispheres have their standard coordinate denoted $z$, which provides a local holomorphic coordinate in $\Sigma^{\sharp,h}$ away from the gluing locus. Defining $\Xi_{\sharp}(z) = \Xi_+(z)$ for $z$ in the first copy and $\Xi_{\sharp}(z) = \Xi_-(z)$ for $z$ in the second copy provides a local holomorphic coordinate $\xi = \Xi_{\sharp}(z)$ in the neighborhood of the gluing locus, and $\sigma_h(z)$ also becomes a local holomorphic coordinate. This turns $\Sigma^{\sharp,h}$ into a Riemann surface. Then, $f_{\sharp} \circ \Xi_{\sharp}$ is a holomorphic function in the neighborhood of the gluing locus in $\Sigma^{\sharp,h}$, which is equal to $\big\langle \boldsymbol{w} \cdot \boldsymbol{F}(z) \big\rangle$ for $z$ in the image in $\Sigma^{\sharp,h}$ of the first copy, and is equal to $\big\langle T^{(h)}(\boldsymbol{w}) \cdot \boldsymbol{F}(z) \big\rangle$ in the image of the second copy. Therefore, $\big\langle \boldsymbol{w} \cdot \boldsymbol{F}(z)\big\rangle$ initially defined and holomorphic on $\amsmathbb{C} \setminus \overline{\amsmathbb{B}}$ admits an analytic continuation to a meromorphic function on $\Sigma^{\sharp,h}$ away from the common image of $\overline{\amsmathbb{B}}_{g}^{\pm}$ for $g \neq h$, with possible poles at the image of $\infty$ in $\Sigma^{\sharp,h}$ (this gives two points, one for each copy, since within a copy the infinity points on each hemisphere were identified in $\Sigma^{\sharp,h}$) and at the image of $\alpha_h,\beta_h$ in $\Sigma^{\sharp,h}$.

\begin{figure}[h!]
\begin{center}
\includegraphics[width=0.95\textwidth]{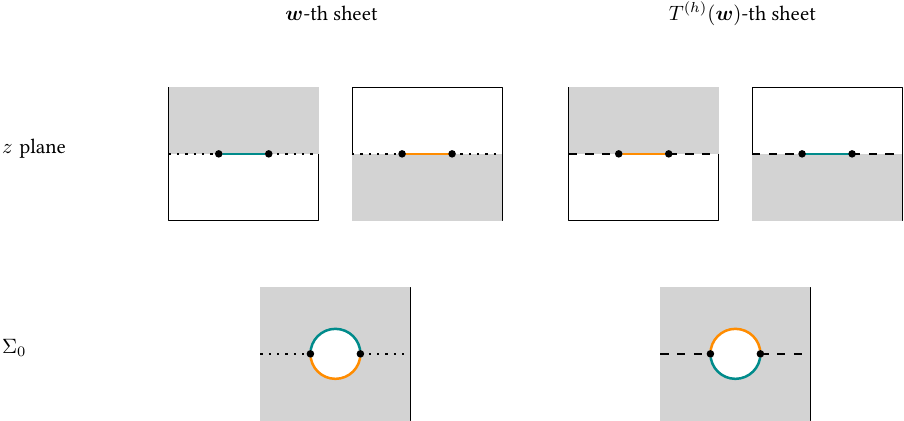}
\caption{\label{Fig:Half2} The Riemann surface $\Sigma^{\sharp,h}$. The orange curves (resp. the blue curves) are identified in a natural way in $\Sigma^{\sharp,h}$. Their endpoints are the ramification points of the branched covering map $Z$.}
\end{center}
\end{figure}

Starting from a non-zero vector $\boldsymbol{v} \in \amsmathbb{R}^H$, we can repeat this construction to get an analytic continuation for $\big\langle\boldsymbol{v} \cdot \boldsymbol{F}(z)\big\rangle$ by crossing all the possible segments. The outcome is a meromorphic continuation $\widetilde{F}_{\boldsymbol{v}}(p)$ of the function $\big\langle \boldsymbol{v} \cdot \boldsymbol{F}(z)\big\rangle$ to the Riemann surface $\widetilde{\Sigma}_{\boldsymbol{v}}$ obtained by gluing copies of $\widehat{\amsmathbb{H}}_+ \sqcup \widehat{\amsmathbb{H}}_-$ labeled by finite words formed with the alphabet $[H]$ modulo the following identifications (\textit{cf.} Figure~\ref{Fig:Half2})
\begin{itemize}
\item $\widehat{\amsmathbb{R}}^+ \setminus \amsmathbb{B}^+$ with $\widehat{\amsmathbb{R}}^- \setminus \amsmathbb{B}^-$ in the same copy;
\item $\overline{\amsmathbb{B}}_h^{\pm}$ in the copy labeled $g_m \cdots g_1$ with $\overline{\amsmathbb{B}}_h^{\mp}$ in the copy labeled $hg_m \cdots g_1$.
\end{itemize}
The empty word labels the principal sheet, \textit{i.e.} a copy where $F_{\boldsymbol{v}}(z)$ is equal to $\big\langle \boldsymbol{v} \cdot \boldsymbol{F}(z)\big\rangle$. The possible poles of this function are located at the points in $\widetilde{\Sigma}_{\boldsymbol{v}}$. with coordinate $z = \infty$. An immediate property of the construction is that $\widetilde{F}_{\boldsymbol{v}}$ takes the same values on two sheets $g_m \cdots g_1$ and $h_n \cdots h_1$ as soon as $T^{(g_m)} \circ \cdots \circ T^{(g_1)}(\boldsymbol{v}) = T^{(h_n)} \circ \cdots \circ T^{(h_1)}(\boldsymbol{v})$. As a consequence, it descends to a meromorphic function on the quotient of $\widetilde{\Sigma}_{\boldsymbol{v}}$ by this equivalence relation. The quotient Riemann surface is exactly $\Sigma_{\boldsymbol{v}}$ and its sheets are labeled by the orbit $\mathfrak{G}.\boldsymbol{v}$. An example of picture of $\Sigma_{\boldsymbol{v}}$ is given in Figure~\ref{Fig:Half4}.

If $\mathfrak{G}.\boldsymbol{v}$ is finite, the function
\[
P(y,z) = \prod_{\boldsymbol{w} \in \mathfrak{G}.\boldsymbol{v}} \big(y - \langle \boldsymbol{w} \cdot \boldsymbol{F}(z)\rangle\big)
\]
is a polynomial in $y$, and a holomorphic function of $z \in \amsmathbb{C} \setminus \overline{\amsmathbb{B}}$. Thanks to \eqref{wcross}, for any $h \in [H]$ and $x \in \amsmathbb{B}_h$ the values $P(y,x^+)$ and $P(y,x^-)$ are equal since they are related by a permutation of factors. Tameness and the assumption on the behavior of $\boldsymbol{F}(z)$ at $z=\infty$ implies that $P(y,z)$ is a rational function of $z$ whose poles are independent of $y$. Its numerator is then a polynomial $P_0(y,z)$ in the variables $y$ and $z$. Its roots with respect to $y$ are the branches of algebraic functions in $z$, and by construction these roots are $y = \langle \boldsymbol{w} \cdot \boldsymbol{F}(z)\rangle$ indexed by $\boldsymbol{w} \in \mathfrak{G}.\boldsymbol{v}$. We claim they are actually the various branches of a single algebraic function, \textit{i.e.} that $P_0(y,z)$ is an irreducible polynomial in $y$. Indeed, if $P_0(y,z)$ were not irreducible, there would be a subset $\mathfrak{W} \subset \mathfrak{G}.\boldsymbol{v}$ containing $\boldsymbol{v}$ such that
\[
\prod_{\boldsymbol{w} \in \mathfrak{W}} \big(y - \langle \boldsymbol{w} \cdot \boldsymbol{F}(z)\rangle\big)
\]
is a rational function of $z$, \textit{i.e.} has no discontinuity on $\amsmathbb{B}$. By \eqref{wcross} this would imply that $\mathfrak{W}$ is stable under $\mathfrak{G}$, which contradicts the fact that the orbit $\mathfrak{G}.\boldsymbol{v}$ is the minimal subset of $\amsmathbb{R}^H$ containing $\boldsymbol{v}$ and stable under the action of $\mathfrak{G}$.
\end{proof}

\begin{figure}[h!]
\begin{center}
\includegraphics[width=0.95\textwidth]{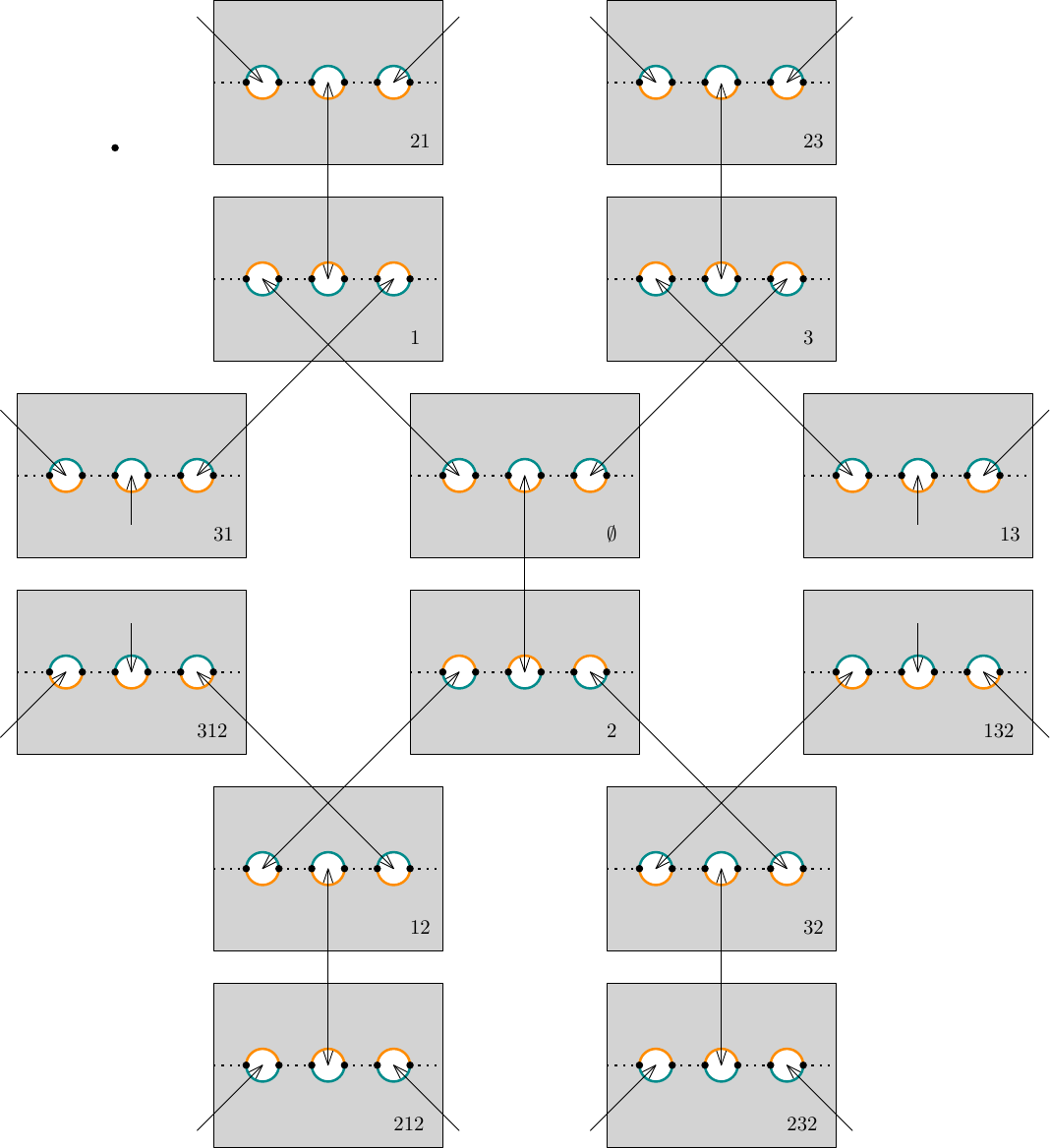}
\caption{\label{Fig:Half3} A part of the Riemann surface $\tilde{\Sigma}_{\boldsymbol{v}}$ for $H = 3$. The sheets are labeled by words in the letters $1,2,3$ labeling the three segments. The orange part (resp. blue part) of two circles related by an arrow are identified. The diagram should be continued, with one sheet for each vertex of the infinite $3$-regular tree.}
\end{center}
\end{figure}

\begin{figure}[h!]
\begin{center}
\includegraphics[width=0.6\textwidth]{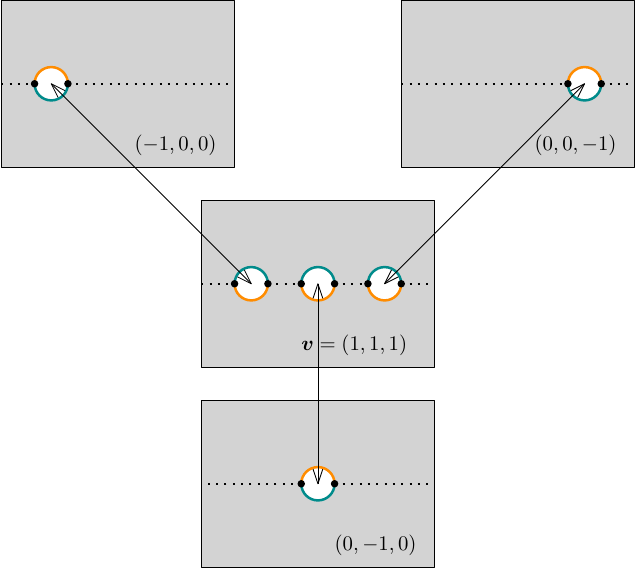}
\caption{\label{Fig:Half4} For $H = 3$ and $\theta_{h,h} = 1 $ and $\theta_{g,h} = \frac{1}{2}$ for any distinct $g,h \in [3]$, the vector $\boldsymbol{v} = (1,1,1)$ has an orbit of size $4$, containing as other vectors $(-1,0,0)$,$(0,-1,0)$ and $(0,0,-1)$. We depict the construction of the spectral curve $\Sigma_{\boldsymbol{v}}$ by gluing four copies of pairs of hemispheres. The union of blue and orange curves in the closure of the $\boldsymbol{w}$-sheet is the $\boldsymbol{w}$-th gluing locus $\overline{\amsmathbb{B}}_{\boldsymbol{w}}$.}
\end{center}
\end{figure}

\subsection{Structure of the group \texorpdfstring{$\mathfrak{G}$}{G}}
\label{sec:bigstr}
In this section, we investigate the structure of the group $\mathfrak{G}$ and its orbits. Recall the notation $\boldsymbol{e}^{(1)},\ldots,\boldsymbol{e}^{(H)}$ for the canonical basis of $\amsmathbb{R}^H$.

 \begin{proposition}
\label{lemrefl} The subspace $\mathscr{V} = \textnormal{Im}(\boldsymbol{\Theta}) = \textnormal{Ker}(\boldsymbol{\Theta})^{\bot}$ is stable under the action of $\mathfrak{G}$. The symmetric bilinear form $(\boldsymbol{v},\boldsymbol{v}') \mapsto \big\langle \boldsymbol{v} \cdot \boldsymbol{\Theta}(\boldsymbol{v}')\big\rangle$ is non-degenerate on $\mathscr{V}$ (call it the modified pairing). The restriction of $\boldsymbol{\Theta}$ to $\mathscr{V}$ gives $\hat{\boldsymbol{\Theta}} \in \textnormal{GL}(\mathscr{V})$ and for each $h \in [H]$ the map $\hat{T}^{(h)} := \hat{\boldsymbol{\Theta}}^{-1} \circ T^{(h)} \circ \boldsymbol{\Theta}$ is the orthogonal reflection (for the modified pairing in $\mathscr{V}$) with respect to the projection $\hat{\boldsymbol{e}}^{(h)} \in \mathscr{V}$ of $\boldsymbol{e}^{(h)}$ orthogonally to $\textnormal{Ker}(\boldsymbol{\Theta})$.
\end{proposition}
\begin{proof}
Equip $\amsmathbb{R}^H$ with its standard scalar product $\langle *\cdot *\rangle$. Since $\boldsymbol{\Theta}$ is symmetric we have the decomposition $\amsmathbb{R}^{H} = \mathscr{V} \oplus \mathscr{V}^{\bot}$ for
\begin{equation}
\label{VImKer}\mathscr{V} = \textnormal{Im}(\boldsymbol{\Theta}) = \textnormal{Ker}(\boldsymbol{\Theta})^{\bot} \qquad \textnormal{and}\qquad \mathscr{V}^{\bot} = \textnormal{Ker}(\boldsymbol{\Theta}) = \textnormal{Im}(\boldsymbol{\Theta})^{\bot}.
\end{equation}
 If $\boldsymbol{v},\boldsymbol{v}' \in \amsmathbb{R}^H$ and $h \in [H]$, we compute
\[
\langle \boldsymbol{v}' \cdot T^{(h)}(\boldsymbol{v}) \rangle = \langle \boldsymbol{v}' \cdot \boldsymbol{v}\rangle - \frac{2\lambda_{h}(\boldsymbol{v})}{\theta_{h,h}} \langle \boldsymbol{v}' \cdot \boldsymbol{\Theta}(\boldsymbol{e}^{(h)}) \rangle .
\]
where $\lambda_h : \amsmathbb{R}^H \rightarrow \amsmathbb{R}$ extracts the $h$-th coefficient. In particular, if $\boldsymbol{v} \in \mathscr{V}$ and $\boldsymbol{v}' \in \mathscr{V}^{\bot}$ the two terms in the right-hand side vanish due to \eqref{VImKer}, showing that $T^{(h)}$ leaves $\mathscr{V}$ stable. For $\boldsymbol{v} \in \amsmathbb{R}^H$ we compute
\[
T^{(h)}(\boldsymbol{\Theta}(\boldsymbol{v})) = \boldsymbol{\Theta}\left[\boldsymbol{v} - \frac{2 \langle \boldsymbol{e}^{(h)} \cdot \boldsymbol{\Theta}(\boldsymbol{v})\rangle}{\theta_{h,h}} \boldsymbol{e}^{(h)}\right].
\]
We can rewrite the right-hand side solely in terms of the orthogonal projection $\hat{\boldsymbol{e}}^{(h)}$ of the vector $\boldsymbol{e}^{(h)}$ onto $\mathscr{V}$:
\[
T^{(h)}(\boldsymbol{\Theta}(\boldsymbol{v})) = \boldsymbol{\Theta}\left[\boldsymbol{v} - \frac{2 \langle \hat{\boldsymbol{e}}^{(h)} \cdot \boldsymbol{\Theta}(\boldsymbol{v})\rangle}{\langle \hat{\boldsymbol{e}}^{(h)} \cdot \boldsymbol{\Theta}(\boldsymbol{e}^{(h)}) \rangle} \hat{\boldsymbol{e}}^{(h)}\right].
\]
From this formula the remaining properties are obvious.
\end{proof}

\begin{definition}
\label{def:smallmono}
We denote $\hat{\mathfrak{G}} \subset \textnormal{GL}(\mathscr{V})$ the group generated by $\hat{T}^{(1)},\ldots,\hat{T}^{(H)}$.
\end{definition}

If $\boldsymbol{\Theta}$ is positive semi-definite (as in Assumption~\ref{Assumptions_Theta} of Chapter~\ref{Chapter_Setup_and_Examples}), the modified pairing has Euclidean signature and the reflections $\hat{T}^{(h)}$ are elements of $\textnormal{O}(\mathscr{V})$. In absence of positivity assumptions, the signature may not be Euclidean but we can still talk about orthogonal reflection with respect to the modified pairing. Although $\hat{\mathfrak{G}}$ is a reflection group, the subgroup of $\textnormal{GL}(\mathscr{V})$ generated by the restrictions of $T^{(1)},\ldots,T^{(H)}$ to $\mathscr{V}$ may not be a reflection group because $\hat{\boldsymbol{\Theta}}$ need not be an orthogonal transformation of $\mathscr{V}$ equipped with the modified pairing.

\begin{proposition}\label{lem:extG}
$\mathfrak{G}$ is an extension of the reflection group $\hat{\mathfrak{G}}$ by a torsion-free abelian group $\mathfrak{K}$.
\end{proposition}
\begin{proof}
For any $h \in [H]$, we can decompose according to $\amsmathbb{R}^H = \mathscr{V} \oplus \mathscr{V}^\bot$:
\begin{equation}
\label{generatorsplit}
T^{(h)} = \left(\begin{array}{cc} \hat{\boldsymbol{\Theta}}\hat{T}^{(h)} \hat{\boldsymbol{\Theta}}^{-1} & -\dfrac{2\lambda_h}{\theta_{h,h}} \bth_h \\[2ex] 0 & \textnormal{id} \end{array}\right),
\end{equation}
where the linear form $\lambda_h : \amsmathbb{R}^H \rightarrow \amsmathbb{R}$ extracts the $h$-th component of a vector. Therefore, any element $G \in \mathfrak{G}$ takes the form
\[
G = \left(\begin{array}{cc} \hat{\boldsymbol{\Theta}}\hat{G} \hat{\boldsymbol{\Theta}}^{-1} & \iota_G \\0 & \textnormal{id}\end{array}\right),\qquad \iota_G \in \textnormal{Hom}(\mathscr{V}^{\bot},\mathscr{V}),
\]
and for any $G_1,G_2 \in\mathfrak{G}$ we have the composition rules
\begin{equation}
\label{comporule}
\widehat{G_1G_2} = \hat{G}_1 \circ \hat{G}_2,\qquad \iota_{G_1G_2} = \hat{\boldsymbol{\Theta}}\hat{G}_1\hat{\boldsymbol{\Theta}}^{-1} \circ \iota_{G_2} + \iota_{G_1}.
\end{equation}
Obviously, we have $\hat{\textnormal{id}} = \textnormal{id}$ and $\iota_{\textnormal{id}} = 0$. This shows that the map $G \mapsto \hat{G}$ defines a surjective morphism of groups $p : \mathfrak{G} \rightarrow \hat{\mathfrak{G}}$. Specializing \eqref{comporule} to elements $G_1,G_2 \in \textnormal{Ker}(p)$ yields $\iota_{G_1 G_2} = \iota_{G_1} + \iota_{G_2}$. Therefore, $\iota : \mathfrak{K}:= \textnormal{Ker}(p) \rightarrow \textnormal{Hom}(\mathscr{V}^{\bot},\mathscr{V})$ is a injective morphism of abelian groups (using the additive structure of the target). Since the target is torsion-free, the source must be as well. To sum up, we have a short exact sequence
\[
1 \longrightarrow \mathfrak{K} \mathop{\longrightarrow}^{\iota} \mathfrak{G} \mathop{\longrightarrow}^{p} \hat{\mathfrak{G}} \longrightarrow 1.
\]
\end{proof}

We stress that $\hat{\boldsymbol{\Theta}}$ maps bijectively the $\hat{\mathfrak{G}}$-orbit of a vector $\hat{\boldsymbol{v}} \in \mathscr{V}$ to the $\mathfrak{G}$-orbit of $\boldsymbol{v} = \hat{\boldsymbol{\Theta}}(\hat{\boldsymbol{v}})$. We are particularly interested in looking for finite orbits, possibly of minimal size. For this we can use the following observation.

\begin{lemma}
\label{lem:comparor} $\mathfrak{G}$ admits a non-zero finite orbit if and only if $\hat{\mathfrak{G}}$ does. In this case, the minimal sizes of non-zero finite orbits are the same for $\hat{\mathfrak{G}}$ and $\mathfrak{G}$.
\end{lemma}
\begin{proof}
Let $d$ (respectively $\hat{d}$) be the minimal size of a non-zero $\mathfrak{G}$-orbit (respectively $\hat{\mathfrak{G}}$-orbit) --- if it does not exist we set $d = \infty$ (respectively $\hat{d} = \infty$). If $\hat{\boldsymbol{v}} \in \mathscr{V}$ is non-zero and has finite $\hat{\mathfrak{G}}$-orbit of size $\hat{d}$, then $\hat{\boldsymbol{\Theta}}(\hat{\boldsymbol{v}})$ is non-zero and has finite $\mathfrak{G}$-orbit of same size. Hence $d \leq \hat{d}$. Conversely, assume we have a non-zero $\boldsymbol{v} \in \amsmathbb{R}^H$ with finite $\mathfrak{G}$-orbit of size $d$, and decompose it $\boldsymbol{v} = \hat{\boldsymbol{v}} + \boldsymbol{v}_0$ with $\hat{\boldsymbol{v}} \in \mathscr{V}$ while $\boldsymbol{v}_0 \in \mathscr{V}^{\bot}$. Since $\mathfrak{G}$ leaves $\mathscr{V}$ stable, $\hat{\boldsymbol{v}}$ has finite $\mathfrak{G}$-orbit of size $\leq d$, \textit{i.e.} $\hat{\boldsymbol{\Theta}}^{-1}(\hat{\boldsymbol{v}})$ has finite $\hat{\mathfrak{G}}$-orbit of size $\leq d$. If $\hat{\boldsymbol{v}} = 0$ we are not satisfied as this is the zero $\hat{\mathfrak{G}}$-orbit, but the form \eqref{generatorsplit} of the generators and the fact that $\boldsymbol{v}$ has at least a non-zero component implies that we can find $h \in [H]$ such that $\boldsymbol{w} = T^{(h)}(\boldsymbol{v})$ has projection $\hat{\boldsymbol{w}} \neq 0$ to $\mathscr{V}$ along $\mathscr{V}^{\bot}$. Since $\mathfrak{G}.\boldsymbol{v} = \mathfrak{G}.\boldsymbol{w}$, the $\hat{\mathfrak{G}}$-orbit of $\hat{\boldsymbol{\Theta}}^{-1}(\hat{\boldsymbol{w}})$ is non-zero and finite of size $\leq d$. In any case, we get $\hat{d} \leq d$.
\end{proof}

In the case of positive semi-definite $\boldsymbol{\Theta}$ (as in Assumption~\ref{Assumptions_Theta}), we can be more precise.

 \begin{lemma}
 \label{lem:classfinorb}
If $\boldsymbol{\Theta}$ is positive semi-definite, there exists an orthogonal decomposition $\mathscr{V} \simeq \bigoplus_{i = 1}^n \mathscr{V}_i$ such that $\hat{\mathfrak{G}} \simeq \prod_{i = 1}^n \hat{\mathfrak{G}}_i$, where for each $i \in [n]$ the reflection group $\hat{\mathfrak{G}}_i \subset \textnormal{O}(\mathscr{V}_i)$ is irreducible (\textit{i.e.} the only non-zero subspace of $\mathscr{V}_i$ stable under the action of $\hat{\mathfrak{G}}_i$ is $\mathscr{V}_i$ itself). Moreover, this decomposition is unique up to permutation of the factors, and for each $i \in [n]$ the following properties are equivalent.
\begin{itemize}
\item[(i)] $\hat{\mathfrak{G}}_i$ is discrete.
\item[(ii)] $\hat{\mathfrak{G}}_i$ is finite.
\item[(iii)] $\hat{\mathfrak{G}}_i$ admits a finite non-zero orbit.
\end{itemize}
\end{lemma}
\begin{proof}
For the decomposition of reflection groups in irreducible factors, see \textit{e.g.} \cite[Theorem 1.4.1]{Popov}. Let $\hat{\mathfrak{G}}_i \subset \textnormal{O}(\mathscr{V}_i)$ be one of the factors. As the orthogonal group $\textnormal{O}(\mathscr{V}_i)$ is compact, discrete subgroups are finite, justifying (i) $\Rightarrow$ (ii). The implication (ii) $\Rightarrow$ (iii) is clear, and it remains to show (iii) $\Rightarrow$ (i). Assume that $\hat{\mathfrak{G}}_i$ is not discrete and that there exists a non-zero $\boldsymbol{v} \in\mathscr{V}_i$ with finite orbit. Call $d_i$ the size of this orbit. Then, we must have $\dim(\mathscr{V}_i) \geq 2$ and one can approximate $\textnormal{id}_{\mathscr{V}_i}$ by elements of $\hat{\mathfrak{G}}_i$. In particular, using the normal form of orthogonal transformations, one can find an element of $\hat{\mathfrak{G}}_i$ whose restriction on a plane $\mathscr{V}_{i,0} \subseteq \mathscr{V}_i$ is a rotation of angle $\chi_i > 0$ such that $\frac{\pi}{\chi_i} > d_i$. Since this rotation has (possibly infinite) order larger than the size of the orbit, all vectors in this orbit must belong to the orthogonal subspace of $\mathscr{V}_{i,0}$ inside $\mathscr{V}_i$. Then, $\textnormal{span}_{\amsmathbb{R}}(\hat{\mathfrak{G}}_i.\boldsymbol{v})$ would be a strict non-zero subspace of $\mathscr{V}_i$ which is stable under the action of $\hat{\mathfrak{G}}_i$, contradicting its irreducibility. This proves (iii) $\Rightarrow$ (i).
\end{proof}

Finite irreducible Euclidean reflection groups are exactly the finite irreducible Coxeter groups, and the latter are described by their Coxeter graph, \textit{cf.} \cite{Humphreys}. Explicitly, they are isomorphic to: either the Weyl group of a simple Lie algebra $A_n$ for $n \geq 1$, $B_n$ for $n \geq 2$, $D_n$ for $n \geq 4$, $E_6$, $E_7$, $E_8$, $G_2$, $F_4$; or the isometry group $H_3$ of the icosahedron with $20$ triangular faces; or the isometry group $H_4$ of the regular polyhedron in $\amsmathbb{R}^4$ with $600$ tetrahedral faces; or a dihedral group $I_2(n)$ for $n \geq 5$ (the only redundancy in this list is $I_2(6) \simeq G_2$). Besides, their orbits are classified up to isomorphism by the standard parabolic subgroups, which can easily found by considering subgraphs of the Coxeter graphs obtained by deleting some vertices.

Note that for each irreducible factor, either all their orbits are finite, or all orbits are infinite. In practice, this gives a rather efficient way to find out whether $\hat{\mathfrak{G}}$ is infinite or finite, and in the latter case, find its isomorphism class, its orbits and construct an orbit of minimal size.
Orbits can easily be generated on a computer. If it takes too much time, one can compute the size of some orbits, or a lower bound for them. In practice, comparing these sizes with orders of Euclidean reflection groups or the (known) sizes for their orbits quickly leads to identify the isomorphism class of the group (if finite) or indicate that it is infinite. This method will be used for instance in Theorems~\ref{thm:Omegav}-\ref{thm:spcurvenonbip} to short cut computations. To come back to finite orbits for $\mathfrak{G}$ we use Lemma~\ref{lem:comparor} and the comment preceding it.

\begin{figure}[h!]
\begin{center}
\begin{tabular}{|c|c||c|c|}
\hline
$\mathfrak{G}$ & $\#\mathfrak{G}$ & $\mathfrak{P}$ & $\#\mathfrak{G}.\boldsymbol{v}$  \\
\hline
$A_n$ ($n \geq 1$) & $(n + 1)!$ & $A_{n - 1}$ & $n + 1$ \\
\hline
$B_n$ ($n \geq 2$) & $2^n n!$ & $B_{n - 1}$ & $2n$  \\
\hline
$D_n$ ($n \geq 4$) & $2^{n - 1}n!$ & $D_{n - 1}$ & $2n$ \\
\hline
$E_6$ & $51840$ & $D_5$ & $27$ \\
\hline
$E_7$ & $2903040$ & $E_6$ & $56$ \\
\hline
$E_8$ & $696729600$ & $E_7$ & $240$ \\
\hline
$G_2$ & $12$ & $A_1$ & $6$ \\
\hline
$F_4$ & $1152$ & $B_3$ & $24$ \\
\hline
$H_3$ & $120$ & $I_2(5)$ & $12$ \\
\hline
$H_4$ & $14400$ & $H_3$ & $120$ \\
\hline
$I_2(n)$ ($n \geq 5$) & $2n$ & $A_1$ & $n$ \\
\hline
\end{tabular}
\end{center}
\caption{\label{Fig:minor} Each line correspond to a finite irreducible Euclidean reflection group: we give its order, the type of a maximal parabolic subgroup $\mathfrak{P}$ (with conventions $A_0 = \{1\}$, $B_1 = A_1$, $D_3 = A_3$), the order of an orbit $\mathfrak{G}.\boldsymbol{v}$ of minimal size (\textit{i.e.} with $\boldsymbol{v}$ having stabilizer $\mathfrak{P}$).}
\end{figure}

\subsection{Geometry of compact spectral curves}
\label{sec:finiteor}

The situation where $\mathfrak{G}$ admits finite orbits is interesting for computational purposes, and deserves a discussion of its own. If $\mathfrak{G}.\boldsymbol{v}$ is finite, we can describe more precisely the geometry of Riemann surface $\Sigma_{\boldsymbol{v}}$, which is compact because it is obtained by gluing finitely many closed hemispheres. Its genus is determined by the Riemann--Hurwitz formula
\begin{equation}
\label{RHformula}\mathsf{g}_{\boldsymbol{v}} = 1 + \sum_{\boldsymbol{w} \in\mathfrak{G}.\boldsymbol{v}} \Big(-1 + \frac{1}{2} \# \textnormal{supp}(\boldsymbol{w})\Big).
\end{equation}
Above each branch point $\alpha_h,\beta_h$, there are exactly $\frac{1}{2}\#\{\boldsymbol{w} \in \mathfrak{G}.\boldsymbol{v} \,\,|\,\,w_h \neq 0\}$ ramification points in $\Sigma_{\boldsymbol{v}}$. The homology $H_1(\Sigma_{\boldsymbol{v}},\amsmathbb{Z})$ equipped with the intersection pairing is a symplectic vector space. By construction of $\Sigma_{\boldsymbol{v}}$, it admits a canonical Lagrangian sublattice $\mathcal{L}$ spanned by the cycles $\gamma_h \subset \Sigma_{\boldsymbol{v},\boldsymbol{w}}$ for some $h \in \textnormal{supp}(\boldsymbol{w})$ and $\boldsymbol{w} \in \mathfrak{G}.\boldsymbol{v}$. These cycles are not independent, as for fixed $\boldsymbol{w} \in\mathfrak{G}.\boldsymbol{v}$ their sum over $h \in \textnormal{supp}(\boldsymbol{w})$ is zero-homologous. We can define a basis $(\mathsf{a}_i)_{i = 1}^{\mathsf{g}_{\boldsymbol{v}}}$ for $\mathcal{L}$ as follows. The choice of order $T^{(1)},\ldots,T^{(H)}$ on the generators induces a total order on $\mathfrak{G}.\boldsymbol{v}$, and on the cycles $\gamma_h \subset \Sigma_{\boldsymbol{v},\boldsymbol{w}}$. We list them by increasing order the elements $(\boldsymbol{w}_i,h_i)_{i = 1}^{\mathsf{g}_{\boldsymbol{v}}}$ such that $\# \textnormal{supp}(\boldsymbol{w}_i) \geq 2$ and $h_i$ is a non-maximal element of $\textnormal{supp}(\boldsymbol{w}_i)$. We set $\mathsf{a}_i = \gamma_{h_i} \subset \Sigma_{\boldsymbol{v},\boldsymbol{w}_i}$.

All in all, this provides a point in the smooth part of the Hurwitz space $\widetilde{\mathcal{H}}_{\mathsf{g}_{\boldsymbol{v}},\# \mathfrak{G}.\boldsymbol{v}}$. By definition, this is the moduli space of branched coverings of the Riemann sphere by a Riemann surface of fixed genus and degree equipped with a Lagrangian sublattice of the first homology (and an integral basis). The smooth part of this moduli space parameterizes coverings with simple ramification points. The geometry of Hurwitz spaces is very rich. For instance, they carry the structure of a Frobenius manifold \cite{Dubrovinfrob}. There are various natural coordinates on Hurwitz spaces, which can be used to describe $\Sigma_{\boldsymbol{v}}$ and thus get formulae for solutions of master Riemann--Hilbert problems that are ``as explicit as possible''. We will see examples in Chapter~\ref{Chap14}.

 \subsection{Extension to inhomogeneous master Riemann--Hilbert problems}
\label{inHomo}

To generalize the construction of spectral curves associated to master Riemann--Hilbert problems with non-zero source (Theorem~\ref{spcurveRHP}), we have to extend the group $\mathfrak{G}$ of Definition~\ref{unimono}. We focus on the case of sources which are rational functions. Let $\boldsymbol{F}(z) = (F_h(z))_{h = 1}^H$ be a solution of the master Riemann--Hilbert problem such that $\boldsymbol{F}(z) = O(z^{d})$ as $z \rightarrow\infty$ for some $d \in \amsmathbb{Z}$. For $\boldsymbol{v},\boldsymbol{v}' \in \amsmathbb{R}^H$, we denote
\begin{equation}
\label{doubleba}
\big\langle\!\!\big\langle (\boldsymbol{v},\boldsymbol{v}') \cdot (\boldsymbol{F},\boldsymbol{D})(z)\big\rangle\!\!\big\rangle = \sum_{h = 1}^H \big(v_h F_h(z) + v_h' D_h(z)\big).
\end{equation}
We denote $\widetilde{\mathfrak{G}}$ the subgroup of $\textnormal{GL}(\amsmathbb{R}^{H} \oplus \amsmathbb{R}^H)$ generated by $\widetilde{T}^{(1)},\ldots,\widetilde{T}^{(H)}$ defined by
\[
\forall h \in [H] \qquad \forall \boldsymbol{v},\boldsymbol{v}' \in \amsmathbb{R}^H \qquad \widetilde{T}^{(h)}(\boldsymbol{v},\boldsymbol{v}') = \left(T^{(h)}(\boldsymbol{v}),\boldsymbol{v}' + \frac{v_h}{\theta_{h,h}} \boldsymbol{e}^{(h)}\right),
\]
This definition is tailored so that the functional relation in the master Riemann--Hilbert problem is equivalent to the identity for any $\boldsymbol{v},\boldsymbol{v}' \in \amsmathbb{R}^H$, any $h \in [H]$ and $x \in \amsmathbb{B}_h$
\[
\big\langle\!\!\big\langle (\boldsymbol{v},\boldsymbol{v}') \cdot (\boldsymbol{F},\boldsymbol{D})(x^+)\big\rangle\!\!\big\rangle = \big\langle\!\!\big\langle \widetilde{T}^{(h)}(\boldsymbol{v},\boldsymbol{v}') \cdot (\boldsymbol{F},\boldsymbol{D})(x^-) \big\rangle\!\!\big\rangle.
\]
Of course, we have $\boldsymbol{D}(x) = \boldsymbol{D}(x^{+}) = \boldsymbol{D}(x^-)$. With this equation replacing \eqref{wcross}, the remaining of the construction of Theorem~\ref{spcurveRHP} is valid. For any non-zero $\widetilde{\boldsymbol{v}} = (\boldsymbol{v},\boldsymbol{v}') \in \amsmathbb{R}^{2H}$, it gives a Riemann surface $\Sigma_{\widetilde{\boldsymbol{v}}}$ and a branched covering map $Z_{\widetilde{\boldsymbol{v}}} : \Sigma_{\widetilde{\boldsymbol{v}}} \rightarrow \widehat{\amsmathbb{C}}$ whose sheets correspond to the orbit $\widetilde{\mathfrak{G}}.\widetilde{\boldsymbol{v}}$, and such that the function $\langle\!\langle\widetilde{\boldsymbol{v}} \cdot (\boldsymbol{F},\boldsymbol{D})(z)\rangle\!\rangle$ admits an analytic continuation as a meromorphic function on $\Sigma_{\widetilde{\boldsymbol{v}}}$. The group $\widetilde{\mathfrak{G}}$ fits in the exact sequence
\[
1 \longrightarrow \textnormal{Ker}(\pi) \longrightarrow \widetilde{\mathfrak{G}} \mathop{\longrightarrow}^{\pi} \mathfrak{G} \longrightarrow 1 ,
\]
where $\pi$ is the restriction of endomorphisms of $\amsmathbb{R}^H \oplus \amsmathbb{R}^H$ to the first summand. As in Proposition~\ref{lem:extG}, $\textnormal{Ker}(\pi)$ can be identified with a subgroup of $\big(\textnormal{End}(\amsmathbb{R}^H),+\big)$, hence it is a torsion-free abelian group. If $\textnormal{Ker}(\pi)$ is trivial, then $\widetilde{\mathfrak{G}}.\widetilde{\boldsymbol{v}} \simeq \mathfrak{G}.\boldsymbol{v}$ and $\Sigma_{\widetilde{\boldsymbol{v}}}$ is the Riemann surface $\Sigma_{\boldsymbol{v}}$ already considered in Theorem~\ref{spcurveRHP}.

\begin{remark}
\label{Remark_about_rational_source_dependent}
If the rational functions $D_1(z),\ldots,D_H(z)$ are not linearly independent, we may replace the second copy of $\amsmathbb{R}^H$ by the smaller vector space
\[
\mathscr{V}_{\boldsymbol{D}} = \amsmathbb{R}^{H} \big/ \big\{\boldsymbol{v} \in \amsmathbb{R}^H\quad |\quad \forall z \in \widehat{\amsmathbb{C}} \quad \langle \boldsymbol{v} \cdot \boldsymbol{D}(z)\rangle = 0\big\}.
\]
The notation \eqref{doubleba} can then used for any $(\boldsymbol{v},\boldsymbol{v}') \in \amsmathbb{R}^H \oplus \mathscr{V}_{\boldsymbol{D}}$, and the group $\widetilde{\mathfrak{G}}$ is rather realized as a subgroup of $\textnormal{GL}(\amsmathbb{R}^H \oplus \mathscr{V}_{\boldsymbol{D}})$.
\end{remark}

\section{Geometry of the fundamental solution}

We apply the discussion of Section~\ref{sec:Algapproach0} to the study the fundamental solution $\boldsymbol{\mathcal{F}}(z,w)$ of the master Riemann--Hilbert problem, \textit{cf.} Definition~\ref{def:Berg}. Throughout this section we assume that $\boldsymbol{\Theta}$ is positive semi-definite, as the existence of the fundamental solution was established in Theorem~\ref{Theorem_Master_equation_12} under this assumption.

\subsection{Analytic continuation on spectral surfaces}
\label{sec:anaBext}

We recall that $\boldsymbol{\mathcal{F}}(z_1,z_2) = (\mathcal{F}_{h_1,h_2}(z_1,z_2))_{h_1,h_2 = 1}^H$ solves the master Riemann--Hilbert problem with respect to $z_1$ \emph{and} with respect to $z_2$ with source
\[
D_{h_1,h_2}(z_1,z_2) = - \frac{\delta_{h_1,h_2}}{(z_1 - z_2)^2}.
\]
The adaptation of Section~\ref{sec:Algapproach} to this situation suggests using an extended double of the group $\mathfrak{G}$.

\begin{definition}
\label{extgroups}
Let $\mathfrak{G}^{[2]}$ be the subgroup of $\textnormal{GL}\big((\amsmathbb{R}^{H} \otimes \amsmathbb{R}^H) \oplus \amsmathbb{R}\big)$ generated by $T^{(h)}_i$ indexed by $i = 1,2$ and $h \in [H]$ and defined by
\begin{equation}
\label{TtildeB}
\begin{split}
T^{(h)}_1\big(\boldsymbol{w}_1\otimes \boldsymbol{w}_2\,,\,t\big) & = \left(T^{(h)}(\boldsymbol{w}_1) \otimes \boldsymbol{w}_2\,,\, t - \frac{\lambda_h(\boldsymbol{w}_1)\lambda_h(\boldsymbol{w}_2)}{\theta_{h,h}}\right) , \\
T^{(h)}_2\big(\boldsymbol{w}_1\otimes \boldsymbol{w}_2\,,\,t\big) & = \left(\boldsymbol{w}_1 \otimes T^{(h)}(\boldsymbol{w}_2)\,,\,t - \frac{\lambda_h(\boldsymbol{w}_1)\lambda_h(\boldsymbol{w}_2)}{\theta_{h,h}}\right) ,
\end{split}
\end{equation}
where the linear form $\lambda_h : \amsmathbb{R}^H \rightarrow \amsmathbb{R}$ extracts the $h$-th component of a vector.
\end{definition}
Adapting Definition~\ref{def:spcurv1} to a situation with two variables, we associate \emph{spectral surfaces} to $\mathfrak{G}^{[2]}$-orbits in the following way.
\begin{definition}
\label{def:spsurf}
If $\boldsymbol{v} \in \amsmathbb{R}^H$ is non-zero, we set $\widetilde{\boldsymbol{v}} = (\boldsymbol{v} \otimes \boldsymbol{v},0)\in (\amsmathbb{R}^{H} \otimes \amsmathbb{R}^{H}) \oplus \amsmathbb{R}$ and construct the complex surface
\[
\Sigma^{[2]}_{\widetilde{\boldsymbol{v}}} = \left( \bigsqcup_{\widetilde{\boldsymbol{w}} \in \mathfrak{G}^{[2]}.\widetilde{\boldsymbol{v}}} (\widehat{\amsmathbb{H}}^{+} \sqcup \widehat{\amsmathbb{H}}^-) \times (\widehat{\amsmathbb{H}}^+ \sqcup \widehat{\amsmathbb{H}}^-) \right)\Big/\sim.
\]
The equivalence relation identifies for each $h \in [H]$, $\epsilon_1,\epsilon_2,\epsilon_1',\epsilon_2' \in \{\pm 1\}$ and $\widetilde{\boldsymbol{w}} \in \mathfrak{G}^{[2]}.\boldsymbol{v}$:
\begin{itemize}
\item $\prod_{i = 1}^2 (\widehat{\amsmathbb{R}}^{\epsilon_i} \setminus \amsmathbb{B}^{\epsilon_i})$ with $\prod_{i = 1}^{2} (\widehat{\amsmathbb{R}}^{\epsilon_i'} \setminus \amsmathbb{B}^{\epsilon_i'})$ both in the $\widetilde{\boldsymbol{w}}$-th copy;
\item $\overline{\amsmathbb{B}}_{h}^{+} \times \widehat{\amsmathbb{H}}^{\epsilon}$ in the $\widetilde{\boldsymbol{w}}$-th copy with $\overline{\amsmathbb{B}}_{h}^- \times \widehat{\amsmathbb{H}}^{\epsilon}$ in the $T_1^{(h)}(\widetilde{\boldsymbol{w}})$-th copy;
\item $\widehat{\amsmathbb{H}}^{\epsilon} \times \overline{\amsmathbb{B}}_h^+$ in the $\widetilde{\boldsymbol{w}}$-th copy with $\widehat{\amsmathbb{H}}^{\epsilon} \times \overline{\amsmathbb{B}}_h^-$ in the $T_2^{(h)}(\widetilde{\boldsymbol{w}})$-th copy.
\end{itemize}

Let us denote for $\textnormal{supp}_i(\widetilde{\boldsymbol{w}}) = \big\{h \in [H] \quad \big|\quad T_i^{(h)}(\widetilde{\boldsymbol{w}}) \neq \widetilde{\boldsymbol{w}}\big\}$ for $i = 1,2$. The image in $\Sigma^{[2]}_{\widetilde{\boldsymbol{v}}}$ of the $\widetilde{\boldsymbol{w}}$-th copy of
\[
\bigcup_{ \epsilon,\epsilon' \in \{\pm 1\}} \bigg(\widehat{\amsmathbb{H}}^{\epsilon} \setminus \bigcup_{g \in \textnormal{supp}_1(\widetilde{\boldsymbol{w}})} \overline{\amsmathbb{B}}_g^{\epsilon}\bigg) \times \bigg(\widehat{\amsmathbb{H}}^{\epsilon'} \setminus \bigcup_{g \in \textnormal{supp}_2(\widetilde{\boldsymbol{w}})} \overline{\amsmathbb{B}}_g^{\epsilon'}\bigg)
\]
is denoted $\Sigma^{[2]}_{\widetilde{\boldsymbol{v}},\widetilde{\boldsymbol{w}}}$ and called the $\widetilde{\boldsymbol{w}}$-th sheet. The $\widetilde{\boldsymbol{v}}$-th sheet is called the principal sheet.
\end{definition}
\label{index:Zvtilde}We have a branched covering map $Z_{\widetilde{\boldsymbol{v}}} : \Sigma^{[2]}_{\widetilde{\boldsymbol{v}}} \rightarrow \widehat{\amsmathbb{C}} \times \widehat{\amsmathbb{C}}$ forgetting the sheet label. We also have an involution $\textnormal{P}$ exchanging the two tensor copies in the space $(\amsmathbb{R}^{H} \otimes \amsmathbb{R}^H) \oplus \amsmathbb{R}$. It is easy to see that, if $\widetilde{\boldsymbol{v}}$ is invariant under this involution, its $\mathfrak{G}^{[2]}$-orbit is also stable under the involution. This induces a holomorphic involution on $\Sigma_{\widetilde{v}}^{[2]}$, by sending $(z_1,z_2)$ in the $\widetilde{\boldsymbol{w}}$-th sheet to $(z_1,z_2)$ in the $\textnormal{P}(\widetilde{\boldsymbol{w}})$-th sheet. This allows talking about symmetric functions or bidifferentials on $\Sigma^{[2]}_{\widetilde{\boldsymbol{v}}}$, namely the ones invariant under this involution. Spectral surfaces are the natural home for the analytic continuation of the fundamental solution.

\begin{theorem}
\label{thm:Bcont} Assume $\boldsymbol{\Theta}$ is positive semi-definite, take $\boldsymbol{v} \in \amsmathbb{R}^H$ non-zero and set $\widetilde{\boldsymbol{v}} = (\boldsymbol{v} \otimes \boldsymbol{v},0)$. Then, there exists a unique meromorphic bidifferential $\mathcal{B}_{\widetilde{\boldsymbol{v}}}$ on $\Sigma^{[2]}_{\widetilde{\boldsymbol{v}}}$ such that for any $\widetilde{\boldsymbol{w}} = (\boldsymbol{w}_1 \otimes \boldsymbol{w}_2,t) \in \mathfrak{G}^{[2]}.\widetilde{\boldsymbol{v}}$ and $p \in \Sigma^{[2]}_{\widetilde{\boldsymbol{v}},\widetilde{\boldsymbol{w}}}$ with $Z_{\widetilde{\boldsymbol{v}}}(p) = (z_1,z_2)$ we have
\begin{equation}
\label{Omegatildevnun}
\mathcal{B}_{\widetilde{\boldsymbol{v}}}(p) = \left(\big\langle (\boldsymbol{w}_1 \otimes \boldsymbol{w}_2) \cdot \boldsymbol{\mathcal{F}}(z_1,z_2)\big\rangle + \frac{t}{(z_1 - z_2)^2}\right)\dd z_1\dd z_2.
\end{equation}
Besides, $\mathcal{B}_{\widetilde{\boldsymbol{v}}}$ is symmetric, it has no other poles than the double poles of the second term in \eqref{Omegatildevnun}, and integrating one of its variables around a contour surrounding the $h$-th segment in a fixed sheet gives zero.
\end{theorem}

\begin{remark}
\label{Remark_bidiff_use}
A meromorphic bidifferential on $\Sigma_{\widetilde{\boldsymbol{v}}}^{[2]}$ in a local holomorphic coordinate chart $(\xi_1,\xi_2)$ is an object of the form $f(\xi_1,\xi_2) \dd \xi_1\dd \xi_2$. The function $f$ representing it depends on the choice of the chart. In another overlapping holomorphic coordinate chart $(\tilde{\xi}_1,\tilde{\xi}_2) = (\psi_1(\xi_1),\psi_2(\xi_2))$, the function becomes
\[
\tilde{f}(\tilde{\xi}_1,\tilde{\xi}_2) = f(\xi_1,\xi_2)\,\frac{\dd \xi_1}{\dd \tilde{\xi}_1}\frac{\dd \xi_2}{\dd\tilde{\xi}_2},
\]
so that it represents the same bidifferential $\tilde{f}(\tilde{\xi}_1,\tilde{\xi}_2)\dd \tilde{\xi}_1 \dd \tilde{\xi}_2 = f(\xi_1,\xi_2) \dd \xi_1 \dd \xi_2$.

On the spectral curve $\Sigma_{\boldsymbol{v}}$ a local holomorphic coordinate near the ramification points above $\alpha_h$ is the square root $\xi = \sqrt{(z - \alpha_h)}$ with $z = Z_{\boldsymbol{v}}$. As a result, a meromorphic function $F$ on $\Sigma_{\boldsymbol{v}}$ taking the form $f(z)$ in some sheet containing such a ramification point, and diverging like an inverse square-root as $z \rightarrow \alpha_h$, actually has a simple pole at this ramification point in $\Sigma_{\boldsymbol{v}}$. However, the $1$-form $\dd z = 2\xi\dd\xi$ defined on $\Sigma_{\boldsymbol{v}}$ has a simple zero at those ramification points, so the $1$-form $F \dd Z_{\boldsymbol{v}}$, represented by $f(z)\dd z$ in the given sheet, is regular at ramification points above $\alpha_h$ (a similar argument would hold above $\beta_h$).

As the fundamental solution $\boldsymbol{\mathcal{F}}(z_1,z_2)$ has at most inverse square-root divergence as $z_i \rightarrow \alpha_h$ or $\beta_h$ (tameness exponent $1$), we are in the same situation with respect to two variables. The analytic continuation $\mathcal{F}_{\widetilde{\boldsymbol{v}}}$ of the function $\langle (\boldsymbol{v} \otimes \boldsymbol{v}) \cdot \boldsymbol{\mathcal{F}}(z_1,z_2) \rangle$ to the spectral surface $\Sigma_{\widetilde{\boldsymbol{v}}}^{[2]}$ is meromorphic with simple poles on the ramification locus, but the bidifferential $\mathcal{B}_{\widetilde{\boldsymbol{v}}} = \mathcal{F}_{\widetilde{\boldsymbol{v}}}\,\dd Z_{\widetilde{\boldsymbol{v}},1} \dd Z_{\widetilde{\boldsymbol{v}},2}$ is regular on the ramification locus, as we explain in the proof. This is the reason why it is convenient to use bidifferentials.
\end{remark}

\begin{proof}
Let $\boldsymbol{w}_1,\boldsymbol{w}_2 \in \amsmathbb{R}^H$ and $z_1,z_2 \in \widehat{\amsmathbb{C}} \setminus \overline{\amsmathbb{B}}$. The functional relations of the master Riemann--Hilbert problem for $\boldsymbol{\mathcal{F}}(z_1,z_2)$ with respect to $z_1$ and with respect to $z_2$ yield for any $h \in [H]$ and $x \in \amsmathbb{B}_h$
\begin{equation}
\label{2RHPsf}
\begin{split}
\big\langle (\boldsymbol{w}_1 \otimes \boldsymbol{w}_2) \cdot \boldsymbol{\mathcal{F}}(x^+,z_2) \big\rangle & = \big\langle (T^{(h)}(\boldsymbol{w}_1) \otimes \boldsymbol{w}_2) \cdot \boldsymbol{\mathcal{F}}(x^-,z_2) \big\rangle - \frac{\lambda_{h}(\boldsymbol{w}_1)\lambda_{h}(\boldsymbol{w}_2)}{\theta_{h,h}(x - z_2)^2}, \\
\big\langle (\boldsymbol{w}_1 \otimes \boldsymbol{w}_2) \cdot \boldsymbol{\mathcal{F}}(z_1,x^+) \big\rangle & = \big \langle (\boldsymbol{w}_1 \otimes T^{(h)}(\boldsymbol{w}_2))\cdot \boldsymbol{\mathcal{F}}(z_1,x^-) \big\rangle - \frac{\lambda_h(\boldsymbol{w}_1)\lambda_h(\boldsymbol{w}_2)}{\theta_{h,h}(z_1 - x)^2}.
\end{split}
\end{equation}
Introducing an extended tuple $\overline{\boldsymbol{\mathcal{F}}}(z_1,z_2) = \big(\boldsymbol{\mathcal{F}}(z_1,z_2)\,,\,\frac{1}{(z_1 - z_2)^2}\big)$, Definition~\ref{extgroups} was tailored so that these relations are equivalent to
\begin{equation*}
\begin{split}
\big\langle\!\!\big\langle (\boldsymbol{w}_1 \otimes \boldsymbol{w}_2\,,\,t) \cdot \overline{\boldsymbol{\mathcal{F}}}(x^+,z_2)\big\rangle\!\!\big\rangle = \big\langle\!\!\big\langle T_1^{(h)}(\boldsymbol{w}_1 \otimes \boldsymbol{w}_2\,,\,t) \cdot \overline{\boldsymbol{\mathcal{F}}}(x^-,z_2)\big\rangle\!\!\big\rangle, \\
\big\langle\!\!\big\langle (\boldsymbol{w}_1 \otimes \boldsymbol{w}_2\,,\,t) \cdot \overline{\boldsymbol{\mathcal{F}}}(z_1,x^+)\big\rangle\!\!\big\rangle = \big\langle\!\!\big\langle T_2^{(h)}(\boldsymbol{w}_1 \otimes \boldsymbol{w}_2\,,\,t) \cdot \overline{\boldsymbol{\mathcal{F}}}(z_1,x^-)\big\rangle\!\!\big\rangle,
\end{split}
\end{equation*}
where the standard scalar product in $(\amsmathbb{R}^H \otimes \amsmathbb{R}^H) \oplus \amsmathbb{R}$ was denoted $\langle\!\langle \cdot \rangle\!\rangle$. Following the construction in the proof of Theorem~\ref{spcurveRHP} but now in the two variables $z_1,z_2$, we obtain a meromorphic function $\mathcal{F}_{\widetilde{\boldsymbol{v}}}$ on $\Sigma^{[2]}_{\widetilde{\boldsymbol{v}}}$, such that for any $(\boldsymbol{w}_1 \otimes \boldsymbol{w}_2\,,\,t) \in \mathfrak{G}^{[2]}.\widetilde{\boldsymbol{v}}$ and $p \in \Sigma^{[2]}_{\widetilde{\boldsymbol{v}},\widetilde{\boldsymbol{w}}}$ we have
 \begin{equation}
 \label{Bwv}
\mathcal{F}_{\widetilde{\boldsymbol{v}}}(p) = \big\langle\!\!\big\langle \widetilde{\boldsymbol{w}} \cdot \overline{\boldsymbol{\mathcal{F}}}(z_1,z_2) \big\rangle\!\!\big\rangle = \big\langle (\boldsymbol{w}_1 \otimes \boldsymbol{w}_2) \cdot \boldsymbol{\mathcal{F}}(z_1,z_2)\big\rangle + \frac{t}{(z_1 - z_2)^2} \quad \textnormal{with} \quad Z_{\widetilde{\boldsymbol{v}}}(p) = (z_1,z_2)
\end{equation}
Multiplying by $\dd z_1\dd z_2$ yields a bidifferential $\mathcal{B}_{\widetilde{\boldsymbol{v}}}$ on $\Sigma^{[2]}_{\widetilde{\boldsymbol{v}}}$ with all the required properties. Uniqueness is clear.

We need to check for $i = 1,2$ that $\mathcal{B}_{\widetilde{\boldsymbol{v}}}$ has no poles as $z_i$ approaches $\alpha_1,\beta_1,\ldots,\alpha_H,\beta_H$ or $\infty$. At $\infty$ it comes from the fact (known from $\boldsymbol{\kappa} = 0$ in Definition~\ref{def:Berg}) that $\mathcal{F}_{h_1,h_2}(z_1,z_2) = O(z_i^{-2})$ as $z_i \rightarrow \infty$ for any $h_1,h_2 \in [H]$, and the observations that $z_i^{-2}\dd z_i = - \dd(z_i^{-1})$ and $z_i^{-1}$ is a local coordinate near $\infty$. Let $h \in [H]$, denote $L_h \subset \Sigma^{[2]}_{\widetilde{\boldsymbol{v}}}$ the locus $\{z_1 = \alpha_h\} \cup \{z_1 = \beta_h\}$ and for $i\in \{1,2\}$ write $Z_{\widetilde{\boldsymbol{v}},i}$ for the $i$-th component of the map $Z_{\widetilde{\boldsymbol{v}}}$. By construction of the spectral surface, for each $\widetilde{\boldsymbol{w}} = (\boldsymbol{w}_1 \otimes \boldsymbol{w}_2 \,,\,t) \in \mathfrak{G}^{[2]}.\widetilde{\boldsymbol{v}}$, the locus $L_{h,\widetilde{\boldsymbol{w}}}$ obtained by intersecting $L_h$ with the closure of $\Sigma_{\widetilde{\boldsymbol{v}},\widetilde{\boldsymbol{w}}}^{[2]}$ in $\Sigma_{\widetilde{\boldsymbol{v}}}$ either does not intersect or is a connected component of the zero locus of $\dd Z_{\widetilde{\boldsymbol{v}},1}$. In the former case, $h$ is not in $\textnormal{supp}_1(\widetilde{\boldsymbol{w}}) = \textnormal{supp}(\boldsymbol{w}_1)$, so $\langle\!\langle(\boldsymbol{w}_1 \otimes \boldsymbol{w}_2,t) \cdot \overline{\boldsymbol{\mathcal{F}}}(z_1,z_2)\rangle\!\rangle$ does not involve $\mathcal{F}_{h,h'}(z_1,z_2)$ for any $h' \in [H]$, in particular it is holomorphic on $L_{h,\widetilde{\boldsymbol{w}}}$. In the latter case, $\dd Z_{\widetilde{\boldsymbol{v}},1}$ has a simple zero on $L_{h,\widetilde{\boldsymbol{w}}}$. But $\sigma_h \circ Z_{\widetilde{\boldsymbol{v}},1}$ locally defines a holomorphic function near $L_{h,\widetilde{\boldsymbol{w}}}$ which has a simple zero on $L_{h,\widetilde{\boldsymbol{w}}}$, and we know from Proposition~\ref{thmBfund} that $z_1 \mapsto \boldsymbol{\mathcal{F}}(z_1,z_2)$ is tame with exponent $1$, meaning that it diverges at most like $(\sigma_h \circ Z_{\widetilde{\boldsymbol{v}},1})^{-1}$ on $L_{h,\widetilde{\boldsymbol{w}}}$ away from isolated points, namely the zeros of $\dd Z_{\widetilde{\boldsymbol{v}},2}$ that lie on $L_{h,\widetilde{\boldsymbol{w}}}$. Accordingly, the function $\mathcal{F}_{\widetilde{\boldsymbol{v}}}$ of \eqref{Bwv} has a simple pole on $L_{h,\boldsymbol{w}}$ which is compensated in the bidifferential by the simple zero in $\dd z_1 = \dd Z_{\widetilde{\boldsymbol{v}},1}$. Thus $\mathcal{B}_{\widetilde{\boldsymbol{v}}}$ is holomorphic, except perhaps at isolated points which are at the intersection of $L_{h}$ and zeros of $\dd Z_{\widetilde{\boldsymbol{v}},2}$. By Hartogs' principle (see, \textit{e.g.}, \cite{Ebeling}) isolated points are removable singularities, so $\mathcal{B}_{\widetilde{\boldsymbol{v}}}$ is holomorphic on $L_{h}$. The same argument works for $z_2$ approaching $\alpha_h,\beta_h$.
\end{proof}

\subsection{Relation between spectral surfaces and spectral curves}

\label{index:Extended2} Products of two Riemann surfaces are more convenient to work with than arbitrary complex surfaces, so we would like to understand when $\Sigma_{\widetilde{\boldsymbol{v}}}^{[2]}$ can be replaced with $\Sigma_{\boldsymbol{v}} \times \Sigma_{\boldsymbol{v}}$. Call $\mathfrak{G}^{[2],\circ}$ the subgroup of $\textnormal{GL}(\amsmathbb{R}^{H} \otimes \amsmathbb{R}^{H})$ generated by the restriction of the generators of $\mathfrak{G}^{[2]}$ to the subspace $\amsmathbb{R}^{H} \otimes \amsmathbb{R}^H$. This restriction induces a projection morphism
\begin{equation}
\label{picircdefdef}\pi^{\circ} : \mathfrak{G}^{[2]} \longrightarrow \mathfrak{G}^{[2],\circ}.
\end{equation}
 We also have a second projection morphism
 \begin{equation}
\label{pibulletdefdef}\pi^{\bullet} : \mathfrak{G} \times \mathfrak{G} \longrightarrow \mathfrak{G}^{[2],\circ}
\end{equation}
sending $(G_1,G_2) \in \mathfrak{G} \times \mathfrak{G}$ to the transformation $\boldsymbol{w}_1 \otimes \boldsymbol{w}_2 \mapsto G_1(\boldsymbol{w}_1) \otimes G_2(\boldsymbol{w}_2)$.

As in Section~\ref{sec:anaBext} we restrict our interest to orbits of vectors of the form $\widetilde{\boldsymbol{v}} = (\boldsymbol{v} \otimes \boldsymbol{v},0)$ where $\boldsymbol{v} \in \amsmathbb{R}^H$ is non-zero. We then know that all elements in such orbits are of the form $(\boldsymbol{w}_1 \otimes \boldsymbol{w}_2,t)$ for $\boldsymbol{w}_1,\boldsymbol{w}_2 \in \amsmathbb{R}^H$ and some $t \in \amsmathbb{R}$, and $(\boldsymbol{w}_1 \otimes \boldsymbol{w}_2,t)$ belongs to this orbit if and only if $(\boldsymbol{w}_2 \otimes \boldsymbol{w}_1\,,\,t)$ does. The map $\pi^{\circ}$ induces a projection map at the level of orbits
\[
\pi^{\circ}_{\boldsymbol{v}} : \mathfrak{G}^{[2]}.\widetilde{\boldsymbol{v}} \longrightarrow \mathfrak{G}^{[2],\circ}.(\boldsymbol{v} \otimes \boldsymbol{v})
\]
intertwining the $\mathfrak{G}^{[2]}$-action on the source and the $\mathfrak{G}^{[2],\circ}$-action on the target (by left multiplication). Likewise $\pi^{\bullet}$ induces a projection
\[
\pi^{\bullet}_{\boldsymbol{v}} : \mathfrak{G}.\boldsymbol{v} \times \mathfrak{G}.\boldsymbol{v} \longrightarrow \mathfrak{G}^{[2],\circ}.(\boldsymbol{v} \otimes \boldsymbol{v}).
\]
These two projection maps at the level of orbits yield two unramified covering maps
\begin{equation}
\label{Sv2eq}\Sigma_{\widetilde{\boldsymbol{v}}}^{[2]} \,\,\mathop{\longrightarrow}^{\mathfrak{p}^{\circ}_{\boldsymbol{v}}} \,\,\Sigma_{\boldsymbol{v} \otimes \boldsymbol{v}}^{[2],\circ} \,\,\mathop{\longleftarrow}^{\mathfrak{p}_{\boldsymbol{v}}^{\bullet}}\,\, \Sigma_{\boldsymbol{v}} \times \Sigma_{\boldsymbol{v}},
\end{equation}
where $\Sigma^{[2],\circ}_{\boldsymbol{v} \otimes \boldsymbol{v}}$ is the spectral surface constructed from the $\mathfrak{G}^{[2],\circ}$-orbit of $\boldsymbol{v} \otimes \boldsymbol{v}$ following a procedure similar to Definition~\ref{def:spsurf}. If $\pi^{\circ}_{\boldsymbol{v}}$ is a bijection, then we simply have $\Sigma_{\widetilde{\boldsymbol{v}}}^{[2]} = \Sigma_{\boldsymbol{v} \otimes \boldsymbol{v}}^{[2],\circ}$; if $\pi^{\bullet}$ is an isomorphism, we also have $\Sigma_{\boldsymbol{v} \otimes \boldsymbol{v}}^{[2],\circ} = \Sigma_{\boldsymbol{v}} \times \Sigma_{\boldsymbol{v}}$. The following facts allow comparing the three types of orbits and thus clarify the comparison between the three complex surfaces in \eqref{Sv2eq}.

\begin{proposition}
\label{G2GG} Let $\boldsymbol{v} \in \amsmathbb{R}^H$ be non-zero and $\widetilde{\boldsymbol{v}} = (\boldsymbol{v} \otimes \boldsymbol{v}\,,\,0)$.
\begin{enumerate}
\item For any $\boldsymbol{w} \in \mathfrak{G}.\boldsymbol{v}$ we have $(\boldsymbol{w} \otimes \boldsymbol{w}\,,\,0) \in \mathfrak{G}^{[2]}.\widetilde{\boldsymbol{v}}$.
\item $\textnormal{Ker}(\pi^{\circ})$ is a torsion-free abelian group.
\item $\pi^{\circ}_{\boldsymbol{v}}$ is bijective if and only if $(\pi^{\circ}_{\boldsymbol{v}})^{-1}(\{\boldsymbol{v} \otimes \boldsymbol{v}\}) = \{\widetilde{\boldsymbol{v}}\}$.
\item If $\widetilde{\boldsymbol{v}}$ has finite $\mathfrak{G}^{[2]}$-orbit, then $\pi^{\circ}_{\boldsymbol{v}}$ is bijective.
\item If $-\textnormal{\textbf{Id}} \in \mathfrak{G}$, then $\textnormal{Ker}(\pi^{\bullet}) = \{(\tau \textnormal{\textbf{Id}},\tau \textnormal{\textbf{Id}})\,\,|\,\,\tau \in \{\pm 1\}\}$ and the maps $\pi_{\boldsymbol{v}}^{\bullet}$ and $\mathfrak{p}^{\bullet}_{\boldsymbol{v}}$ have degree $2$ (examples where this is realized will appear in Theorem~\ref{thm:spcurvenonbip}).
\item If $-\textnormal{\textbf{Id}} \notin \mathfrak{G}$, then $\pi^{\bullet}$, $\pi^{\bullet}_{\boldsymbol{v}}$ and $\mathfrak{p}_{\boldsymbol{v}}^{\bullet}$ are isomorphisms.
\end{enumerate}
\end{proposition}
\begin{proof}
 Let $\boldsymbol{w} \in \mathfrak{G}.\boldsymbol{v}$ and assume that $(\boldsymbol{w} \otimes \boldsymbol{w}\,,\,0) \in \mathfrak{G}^{[2]}.\widetilde{\boldsymbol{v}}$. Let $h \in [H]$ and set $\boldsymbol{w}' = T^{(h)}(\boldsymbol{w})$. We compute
\begin{equation*}
\begin{split}
T_1^{(h)}(\boldsymbol{w} \otimes \boldsymbol{w}\,,\,0) & = \left(\boldsymbol{w}' \otimes \boldsymbol{w}\,,\, - \frac{w_h^2}{\theta_{h,h}}\right), \\
T_2^{(h)}\left(\boldsymbol{w}' \otimes \boldsymbol{w}\,,\,-\frac{w_h^2}{\theta_{h,h}}\right) & = \left(\boldsymbol{w}' \otimes \boldsymbol{w}'\,,\, - \frac{w_h^2}{\theta_{h,h}} - \frac{w_h'w_h}{\theta_{h,h}}\right).
\end{split}
\end{equation*}
We see in Definition~\ref{unimono} that $T^{(h)}$ flips the sign of the $h$-th coordinate, namely $w_h = -w'_h$. Therefore, $(T_2^{(h)} \circ T_1^{(h)})(\boldsymbol{w} \otimes \boldsymbol{w}\,,\,0) = (\boldsymbol{w}' \otimes \boldsymbol{w}'\,,\,0)$ belongs to the $\mathfrak{G}^{[2]}$-orbit of $\widetilde{\boldsymbol{v}}$. This justifies 1. by induction.

The justification of 2. goes like in the proof of Proposition~\ref{lem:extG}: $\textnormal{Ker}(\pi^{\circ})$ is a torsion-free abelian group since it is identified with a subgroup of $(\amsmathbb{R}^{H} \otimes \amsmathbb{R}^H,+)$.

In 3. the direct implication is clear. Let $\boldsymbol{w}_1,\boldsymbol{w}_2 \in \mathfrak{G}.\boldsymbol{v}$ and $t,t' \in \amsmathbb{R}$ such that $\widetilde{\boldsymbol{w}} := (\boldsymbol{w}_1 \otimes \boldsymbol{w}_2\,,\,t)$ and $\widetilde{\boldsymbol{w}}' := (\boldsymbol{w}_1 \otimes \boldsymbol{w}_2\,,\,t')$ are in the $\mathfrak{G}^{[2]}$-orbit of $\widetilde{\boldsymbol{v}}$. We can find $G,G' \in \mathfrak{G}^{[2]}$ such that $G(\widetilde{\boldsymbol{v}}) = \widetilde{\boldsymbol{w}}$ and $G'(\widetilde{\boldsymbol{v}}) = \widetilde{\boldsymbol{w}}'$. The form of the generators in Definition~\ref{extgroups} implies that the scalar component is additive under composition in $\mathfrak{G}^{[2]}$. Therefore:
\[
(G^{-1} \circ G')(\boldsymbol{v} \otimes \boldsymbol{v}\,,\,0) = (\boldsymbol{v} \otimes \boldsymbol{v}\,,\, t' - t).
\]
If we assume that $(\pi^{\circ}_{\boldsymbol{v}})^{-1}(\{\boldsymbol{v} \otimes \boldsymbol{v}\}) = \{(\boldsymbol{v} \otimes \boldsymbol{v}\,,\,0)\}$, we must have $t = t'$. This shows that $\pi^{\circ}_{\boldsymbol{v}}$ is injective and proves the converse implication in 2. If we rather assume that $\widetilde{\boldsymbol{v}}$ has finite $\mathfrak{G}^{[2]}$-orbit, we infer from
\[
\forall n \in \amsmathbb{Z} \qquad (G^{-1} \circ G')^{n}(\boldsymbol{v} \otimes \boldsymbol{v}\,,\,0) = (\boldsymbol{v} \otimes \boldsymbol{v}\,,\,n(t - t'))
\]
that $t$ must be equal to $t'$. With help of 3. this proves 4.

To prove the dichotomy 5. or 6., we observe that $\boldsymbol{w}_1 \otimes \boldsymbol{w}_2 = \boldsymbol{w}'_1 \otimes \boldsymbol{w}'_2$ if and only if $\boldsymbol{w}_1 = c \boldsymbol{w}_1'$ and $\boldsymbol{w}_2 = c^{-1} \boldsymbol{w}_2'$ for some $c \in \amsmathbb{R}^*$. This shows that
\[
\textnormal{Ker}(\pi^{\bullet}) = \big\{(c \textnormal{\textbf{Id}},c^{-1}\textnormal{\textbf{Id}})\,\,|\,\,c \in \amsmathbb{R}^*\big\} \cap (\mathfrak{G} \times \mathfrak{G}).
\]
Since elements of $\mathfrak{G}$ have determinant $\pm 1$, the only scalar elements that $\mathfrak{G}$ can contain are $\pm \textnormal{\textbf{Id}}$, giving the desired description of the kernel. The statement about the orbits is then obvious.
\end{proof}

\subsection{Fundamental bidifferentials on compact Riemann surfaces}
\label{sec:fundbidiff}
\label{sec:fund02}

Bidifferentials with double poles on the diagonal on compact Riemann surfaces are important and well-studied objects. Reviewing their construction and properties will turn useful in Section~\ref{sec:computefundfund} and \ref{sec:C13tiling} to compute the bidifferential $\mathcal{B}_{\tilde{\boldsymbol{v}}}$ of Theorem~\ref{thm:Bcont} for a certain class of matrices $\boldsymbol{\Theta}$.

\begin{theorem}\cite{MumTata}
\label{MumUM} Let $\Sigma$ be a compact Riemann surface of genus $\mathsf{g}$ equipped with a lagrangian sublattice $\mathcal{L} \subset H_1(\Sigma,\amsmathbb{Z})$. There exists a unique meromorphic bidifferential $\mathcal{B}^{\Sigma,\mathcal{L}}$ with the following properties.
\begin{itemize}
\item It is symmetric in its two variables.
\item It has double poles on the diagonal with biresidue $1$, \textit{i.e.} for any $p \in \Sigma$ and local coordinate $\zeta_p$ near $p$, we have
\[
\mathcal{B}^{\Sigma,\mathcal{L}}(z_1,z_2)\,\, \mathop{=}_{z_1,z_2 \rightarrow p}\,\, \frac{\dd\zeta_{p}(z_1)\dd\zeta_{p}(z_2)}{(\zeta_p(z_1) - \zeta_p(z_2))^2} + O(1).
\]
\item It is normalized on $\mathcal{L}$, \textit{i.e.} for any $\mathsf{a} \in \mathcal{L}$ we have $\oint_{z \in \mathsf{a}} \mathcal{B}^{\Sigma,\mathcal{L}}(z_0,z) = 0$.
\end{itemize}
\end{theorem}
This $\mathcal{B}^{\Sigma,\mathcal{L}}$ is called \emph{fundamental bidifferential of the second kind}, or sometimes Bergman kernel \cite{EORev}. The latter name is unfortunate, as kernels studied by Bergman \cite{Bergman} are related but do not coincide with $\mathcal{B}^{\Sigma,\mathcal{L}}$. There are many ways to present Riemann surfaces, and therefore various formulae for $\mathcal{B}^{\Sigma,\mathcal{L}}$. The proof in \cite{MumTata} gives a formula in terms of the Riemann theta function. If we have an algebraic equation for $\Sigma$, \cite{Eynb} describes an algorithm to get an algebraic formula for $\mathcal{B}^{\Sigma,\mathcal{L}}$. We review some formulae for low genus, which are well-known and easy to check. We also indicate a basis of the space of holomorphic $1$-forms, that will be useful in view of Proposition~\ref{Lem:Holo1form}.

\medskip

\noindent \textsc{Genus $0$.} $\Sigma$ is isomorphic to the Riemann sphere $\widehat{\amsmathbb{C}}$ and the first homology vanishes. If $z$ is a global coordinate, one can check that
\begin{equation}
\label{Bergg0}
\mathcal{B}^{\widehat{\amsmathbb{C}}}(z_1,z_2) = \frac{\dd z_1\dd z_2}{(z_1 - z_2)^2}.
\end{equation}
Any other global coordinate on $\widehat{\amsmathbb{C}}$ is obtained from $z$ by a M\"obius transformation, \textit{i.e.} $z \mapsto \frac{az + b}{cz + d}$ for $a,b,c,d \in \amsmathbb{C}$ such that $ad - bc = 1$, and such transformations applied simultaneously to $z_1$ and $z_2$ leave \eqref{Bergg0} invariant. On the Riemann sphere there are no non-zero holomorphic $1$-forms.

\vspace{0.2cm}

\noindent \textsc{Genus $1$.} $\Sigma$ is isomorphic to the torus\label{index:torus} $\amsmathbb{T}_{\tau} = \amsmathbb{C}/(\amsmathbb{Z} \oplus \tau \amsmathbb{Z})$ for some $\tau \in \amsmathbb{C}$ such that $\textnormal{Im}(\tau) > 0$. Let $\zeta$ be the natural coordinate $\amsmathbb{T}_{\tau}$. For the Lagrangian $\mathcal{L}$ chosen to be generated by the homology class of $[0,1]$, we have
\begin{equation}
\label{Omegaelliptic}
\mathcal{B}^{\amsmathbb{T}_{\tau},[0,1]}(\zeta_1,\zeta_2) = \bigg(\frac{\pi^2E_2(\tau)}{3} + \sum_{m \in \amsmathbb{Z}} \frac{\pi^2}{\sin^2\pi(\zeta_1 - \zeta_2 + m\tau)}\bigg)\dd \zeta_1\dd \zeta_2,
\end{equation}
where $E_2(\tau)$ is the second Eisenstein series\label{index:E2}. The sum in the right-hand side coincides with Weierstra\ss{}\label{index:wp} elliptic function $\wp(\zeta_1 - \zeta_2|\tau)$. The addition theorem gives for it the alternative expression
\begin{equation}
\label{additionwp}
\wp(\zeta_1 - \zeta_2|\tau) = \frac{1}{4}\left(\frac{\wp'(\zeta_1|\tau) + \wp'(\zeta_2|\tau)}{\wp(\zeta_1|\tau) - \wp(\zeta_2|\tau)}\right)^2 - \wp(\zeta_1|\tau) - \wp(\zeta_2|\tau).
\end{equation}
The Riemann bilinear relation\label{index:eta} $\eta_{\mathsf{b}} - \tau\eta_{\mathsf{a}} = 2\ii\pi$ relates the periods of the Weierstra\ss{} functions, defined as
\begin{equation}
\label{etaperiod}
\eta_{\mathsf{a}} := \int_{0}^{1} \wp(\zeta_1 - \zeta_2)\dd \zeta_1,\qquad \eta_{\mathsf{b}} = \int_{0}^{\tau} \wp(\zeta_1 - \zeta_2) \dd \zeta_1,
\end{equation}
and independent of $\zeta_2$. Normalization of \eqref{Omegaelliptic} on $[0,1]$ comes from the identity $\eta_{\mathsf{a}} = - \frac{\pi^2 E_2(\tau)}{3}$, while integrating on $[0,\tau]$ yields
\begin{equation}
\label{0tauellitpic}
\int_{0}^{\tau} \mathcal{B}^{\amsmathbb{T}_{\tau},[0,1]}(\cdot, \zeta_2) = 2\ii\pi \dd \zeta_2
\end{equation}
thanks to the Riemann bilinear relation. On $\amsmathbb{T}_{\tau}$ there is a unique holomorphic $1$-form up to scale, and the one with unit period along $[0,1]$ is $\dd\zeta$. Equation~\ref{0tauellitpic} gives an alternative expression for it. Other expressions in terms of Riemann theta functions will be discussed in Chapter~\ref{Chap14}.

If we rather present $\Sigma$ as the (compactification of the) algebraic curve $Y^2 = \prod_{a = 1}^4 (X - \lambda_a)$, we have the formula
\begin{equation}
\label{OmeX1X2X1} \mathcal{B}^{\Sigma,\mathcal{L}}(p_1,p_2) = \frac{\dd X_1\,\dd X_2}{4Y_1Y_2}\left(4c_0 - (X_1+X_2)\Big(X_1 + X_2 - \sum_{i = 1}^4\lambda_i\Big) + \frac{(Y_1 + Y_2)^2}{(X_1 - X_2)^2}\right),
\end{equation}
Here, $(X_i,Y_i)$ the coordinates of a point $p_i \in \Sigma$ and the scalar $c_0$ is determined by the normalization on the chosen Lagrangian. If the $\lambda_i$ are real and $\lambda_1 < \lambda_2 \lambda_3 < \lambda_4$, and if $\mathcal{L}$ is taken to be a homology cycle around $[\lambda_1,\lambda_2]$ or $[\lambda_3,\lambda_4]$, this formula coincides with $\big(\theta\mathcal{F}(X_1,X_2) + \frac{1}{(X_1 - X_2)^2}\big)\dd X_1\dd X_2$ where $\mathcal{F}$ is taken from \eqref{Hequal2} in which one should substitute
\[
(\alpha_{1},\beta_1,\alpha_2,\beta_2) \leftarrow (\lambda_1,\ldots,\lambda_4) \qquad \textnormal{and} \qquad c_{0,0} \leftarrow c_0 + \frac{1}{4}\sum_{1 \leq i < j \leq 4} \lambda_i\lambda_j,
\]
that is
\begin{equation}
\label{4c0}
4c_0 = - (\alpha_1 + \beta_2)(\alpha_2 + \beta_1) + (\beta_2 - \beta_1)(\alpha_2 - \alpha_1)\frac{E(\mathsf{k})}{K(\mathsf{k})}.
\end{equation}
The $1$-form $\frac{\dd X}{Y}$ is holomorphic. Therefore, it is proportional to $\dd \zeta$.

\vspace{0.2cm}

\noindent \textsc{Genus $2$.} We can always realize $\Sigma$ as the (compactification of the) an algebraic curve of the form $Y^2 = \prod_{a = 1}^6 (X - \lambda_a)$. Denoting $(X_i,Y_i)$ the coordinates of a point $p_i \in \Sigma$, we then have
\begin{equation}
\label{Bgenus2}
\mathcal{B}^{\Sigma,\mathcal{L}}(p_1,p_2) = \frac{\dd X_1 \dd X_2}{Y_1Y_2}\bigg(c_{0,0} + c_{0,1}(X_1 + X_2) + c_{1,1}X_1X_2 + \hat{P}(X_1,X_2) + \frac{(Y_1 + Y_2)^2}{4(X_1 - X_2)^2}\bigg),
\end{equation}
where $\hat{P}$ is a bivariate polynomial given by \eqref{Phatequation}. The coefficients $c_{j,i} = c_{i,j}$ depend on the choice of $\mathcal{L}$, and are fixed with the substitution
\[
(\alpha_{1},\beta_{1},\alpha_{2},\beta_{2},\alpha_{3},\beta_{3}) \leftarrow (\lambda_1,\ldots,\lambda_6).
\]
A basis for the two-dimensional space of holomorphic $1$-forms is
\[
\frac{\dd X}{Y},\qquad \frac{X\dd X}{Y}.
\]

An alternative presentation by equations is
\[
\tilde{Y}^2 = \prod_{a = 1}^{5} (\tilde{X} - \tilde{\lambda}_a).
\]
This form is related to the previous sextic by the birational transformation
\[
\tilde{X} = \frac{1}{X - \lambda_6},\qquad \tilde{Y} = \frac{\ii Y \prod_{a = 1}^{5} \tilde{\lambda}_i}{(X - \lambda_6)^3}\qquad \tilde{\lambda}_a = \frac{1}{\lambda_a - \lambda_6}.
\]
This results in the alternative form for the fundamental bidifferential
\[
\mathcal{B}^{\Sigma,\mathcal{L}}(p_1,p_2) = \frac{\dd \tilde{X}_1 \dd \tilde{X}_2}{\tilde{Y}_1\tilde{Y}_2} \left( \frac{(\tilde{Y}_1 + \tilde{Y}_2)^2}{4(\tilde{X}_1 - \tilde{X}_2)^2} + \tilde{P}(\tilde{X}_1,\tilde{X}_2) + \tilde{c}_{0,0} + \tilde{c}_{0,1}(\tilde{X}_1 + \tilde{X}_2) + \tilde{c}_{1,1}\tilde{X}_1\tilde{X}_2\right)
\]
for some constants $\tilde{c}_{0,0},\tilde{c}_{1,0},\tilde{c}_{1,1}$ depending on the Lagrangian $\mathcal{L}$.

\subsection{Fundamental solution and fundamental bidifferentials}
\label{sec:computefundfund}

Apart from being defined on the spectral surface, the key properties of the bidifferential $\mathcal{B}_{\widetilde{\boldsymbol{v}}}$ of Theorem~\ref{thm:Bcont} are the coefficients of its double poles. To continue the discussion, we place ourselves in the case where $\pi_{\boldsymbol{v}}^{\circ}$ is a bijection, so that the spectral surface $\Sigma_{\widetilde{\boldsymbol{v}}}^{[2]}$ is either equal to $\Sigma_{\boldsymbol{v}} \times \Sigma_{\boldsymbol{v}}$ or admits a degree $2$ unramified covering $p^{\bullet}_{\boldsymbol{v}} : \Sigma_{\boldsymbol{v}} \times \Sigma_{\boldsymbol{v}} \rightarrow \Sigma_{\widetilde{\boldsymbol{v}}}^{[2]}$.

\begin{definition}
\label{def:Bires}
We write $\mathcal{B}_{\boldsymbol{v}}$ for the bidifferential on $\Sigma_{\boldsymbol{v}} \times \Sigma_{\boldsymbol{v}}$ equal to the bidifferential $\mathcal{B}_{\widetilde{\boldsymbol{v}}}$ from Theorem~\ref{thm:Bcont} (if $\mathfrak{p}_{\boldsymbol{v}}^{\bullet}$ is an isomorphism) or to the pullback of $\mathcal{B}_{\widetilde{\boldsymbol{v}}}$ via $\mathfrak{p}_{\boldsymbol{v}}^{\bullet}$ (if $\pi_{\boldsymbol{v}}^{\bullet}$ is a degree $2$ map).

The \emph{matrix of biresidues} \textbf{\foreignlanguage{russian}{B}} is the real symmetric (possibly infinite) matrix whose rows and columns are indexed by $\mathfrak{G}.\boldsymbol{v}$ and entries are the (unique by Proposition~\ref{G2GG}, 3.) scalars such that
\[
\forall \boldsymbol{w}_1,\boldsymbol{w}_2 \in \mathfrak{G}.\boldsymbol{v},\qquad \big(\boldsymbol{w}_1 \otimes \boldsymbol{w}_2\,,\,\text{\textit{\foreignlanguage{russian}{B}}}_{\boldsymbol{w}_1,\boldsymbol{w}_2}\big) \in \mathfrak{G}^{[2]}.\widetilde{\boldsymbol{v}}.
\]
In other words, they give the coefficient of the double pole of $\mathcal{B}_{\boldsymbol{v}}$ at the locus $z_1 = z_2$ in $\Sigma_{\boldsymbol{v},\boldsymbol{w}_1} \times \Sigma_{\boldsymbol{v},\boldsymbol{w}_2}$.
\end{definition}

By Proposition~\ref{G2GG}, 1. and 3., the matrix of biresidues has zero diagonal. In situations where it has a particularly simple form, we can relate $\mathcal{B}_{\boldsymbol{v}}$ to the fundamental bidifferentials of Section~\ref{sec:fundbidiff}. This will be applied for instance in tiling models, \textit{cf.} Section~\ref{sec:C13tiling} and Chapter~\ref{Chap14}.

\begin{proposition}
\label{funddif}
If $\boldsymbol{v}$ has a finite orbit and the matrix of biresidues $\text{\textbf{{\foreignlanguage{russian}{B}}}} = c\,\textnormal{\textbf{Id}} + c'\,\textnormal{\textbf{1}}$ for some scalars $c,c'$ and $\textnormal{\textbf{1}}$ being the matrix full of $1$s, we have
\[
\mathcal{B}_{\boldsymbol{v}} = c\,\mathcal{B}^{\Sigma_{\boldsymbol{v}}} + c' \frac{\dd z_1\dd z_2}{(z_1 - z_2)^2}.
\]
Here $z_i$ is an abuse of notation for the map $Z_{\boldsymbol{v}}$ on the $i$-th copy of $\Sigma_{\boldsymbol{v}}$ in $\Sigma_{\boldsymbol{v}} \times \Sigma_{\boldsymbol{v}}$, with $i = 1,2$, and the Lagrangian $\mathcal{L} \subset H_1(\Sigma_{\boldsymbol{v}},\amsmathbb{Z})$ is spanned by the contours surrounding the segments $\overline{\amsmathbb{B}}_h$ in $\Sigma_{\boldsymbol{v}}$ (\textit{cf.} Section~\ref{sec:finiteor}).
\end{proposition}
\begin{proof} By design the difference of the two sides is a holomorphic bidifferential on $\Sigma_{\boldsymbol{v}} \times \Sigma_{\boldsymbol{v}}$. Integrating its first variable on a cycle in $\mathcal{L}$ gives zero. Indeed, for any such cycle, take a contour $\gamma$ representing it and $p_2$ on the open dense subset $Z_{\boldsymbol{v}}^{-1}(Z_{\boldsymbol{v}}(\gamma)) \subset \Sigma_{\boldsymbol{v}}$. The integral of the left-hand side and the first-term of the right-hand side is zero by construction $\mathcal{B}^{\Sigma_{\boldsymbol{v}},\mathcal{L}}$, while the integral of the second term in the right-hand side is zero because it is represented as an integral along $Z_{\boldsymbol{v}}(\gamma) \subset \widehat{\amsmathbb{C}}$ which is automatically null-homologous. We conclude from the fact that holomorphic bidifferentials are bilinear combinations of holomorphic $1$-forms and a holomorphic $1$-form with vanishing integral along all cycles in a given Lagrangian is zero.
\end{proof}

\subsection{Relation to Green functions}

\label{sec:Green13}

In this section we study a function obtained by integrating the fundamental solution with respect to both variables taking into account the real structure of the spectral curves. The result has properties similar to a Green function, in particular we shall prove that it is symmetric and admit logarithmic singularities. In special cases it does coincide with the Green function with Dirichlet boundary condition on a certain bordered Riemann surface. In this discussion we assume for simplicity that $\boldsymbol{v}$ is such that $\pi_{\boldsymbol{v}}^{\circ}$ is a bijection, so that we have the bidifferential $\mathcal{B}_{\boldsymbol{v}}$ defined on $\Sigma_{\boldsymbol{v}} \times \Sigma_{\boldsymbol{v}}$ (Definition~\ref{def:Bires}). The interested reader should have no difficulty writing down the appropriate statements without this assumption.

First, we remark that the spectral curve $\Sigma_{\boldsymbol{v}}$ is equipped with an operation of complex conjugation: if $p \in \Sigma_{\boldsymbol{v}}$ is not in a gluing locus, we call $p^* \in \Sigma_{\boldsymbol{v}}$ the point in the same sheet as $p$ and such that $Z_{\boldsymbol{v}}(p^*) = (Z_{\boldsymbol{v}}(p))^*$. This definition is extended to the gluing loci by continuity. We say that $p$ is a real point of $\Sigma_{\boldsymbol{v}}$ if $p = p^*$. In other words, real points are either endpoints of a gluing locus or those points such that $Z_{\boldsymbol{v}}(p) = (Z_{\boldsymbol{v}}(p))^*$ which are not in a gluing locus. If $p_i \in \Sigma_{\boldsymbol{v}}$ we often write $z_i = Z_{\boldsymbol{v}}(p_i)$. We introduce the big diagonal and anti-diagonal:
\begin{equation}
\label{diagantidiag}\begin{split}
\Delta_{\boldsymbol{v}} & = \big\{(p_1,p_2) \in \Sigma_{\boldsymbol{v}} \times \Sigma_{\boldsymbol{v}}\,\,|\,\,Z_{\boldsymbol{v}}(p_1) = Z_{\boldsymbol{v}}(p_2)\big\}, \\
\Delta_{\boldsymbol{v}}^* & = \big\{(p_1,p_2) \in \Sigma_{\boldsymbol{v}} \times \Sigma_{\boldsymbol{v}}\,\,|\,\,Z_{\boldsymbol{v}}(p_1) = (Z_{\boldsymbol{v}}(p_2))^*\big\}
\end{split}
\end{equation}

\begin{figure}[h!]
\begin{center}
\includegraphics[width=0.3\textwidth]{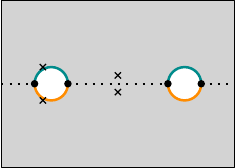}
\caption{\label{Fig:Half6} The closure of a sheet with $\# \textnormal{supp}(\boldsymbol{w}) = 2$. The real points are on the dotted lines. The crosses indicate pairs of points $\{p,p^*\}$. In particular, if $p \in \amsmathbb{B}_{\boldsymbol{w}}^+$, then $p^* \in \amsmathbb{B}_{\boldsymbol{w}}^-$ is a distinct point.}
\end{center}
\end{figure}

\begin{definition}
\label{def:Hvv} We define
\begin{equation}
\label{defHvp1p2}
\mathrm{gr}_{\boldsymbol{v}}(p_1,p_2) = \text{Re}\left( \int^{p_1}_{p_1^*} \int_{p_2^*}^{p_2} \mathcal{B}_{\boldsymbol{v}}\right),
\end{equation}
where the $i$-th integration ($i = 1,2$) runs over a path from $p_i$ to $p_i^*$ in the (closure of) sheet to which these points belong to, for $i = 1,2$.
\end{definition}
Since $\mathcal{B}_{\boldsymbol{v}}$ has double poles and $\Sigma_{\boldsymbol{v}}$ is in general not simply-connected, it is not obvious that this definition is well-posed on $\Sigma_{\boldsymbol{v}} \times \Sigma_{\boldsymbol{v}}$.

\begin{proposition}
\label{prop:Greengen}
Let $\boldsymbol{v} \in \amsmathbb{R}^H$ be non-zero and such that $\pi_{\boldsymbol{v}}^{\circ}$ is a bijection.
\begin{enumerate}
\item[(H1)] $\textnormal{gr}_{\boldsymbol{v}}(p_1,p_2)$ is a well-defined symmetric and pluriharmonic function of $(p_1,p_2) \in \Sigma_{\boldsymbol{v}} \times \Sigma_{\boldsymbol{v}} \setminus (\Delta_{\boldsymbol{v}} \cup \Delta_{\boldsymbol{v}}^*)$ which is independent on the choice of paths of integration;
\item[(H2)] $\textnormal{gr}_{\boldsymbol{v}}(p_1,p_2) + 2\text{\foreignlanguage{russian}{B}}_{\boldsymbol{w}_1,\boldsymbol{w}_2} \big( \log|z_1 - z_2^*| - \log|z_1 - z_2|\big)$ is a pluriharmonic function of $(p_1,p_2)$ near $\Delta_{\boldsymbol{v}} \cup \Delta_{\boldsymbol{v}}^*$ in the $(\boldsymbol{w}_1,\boldsymbol{w}_2)$-th sheet;
\item[(H3)] $\textnormal{gr}_{\boldsymbol{v}}(p_1,p_2) = - \textnormal{gr}_{\boldsymbol{v}}(p_1^*,p_2) = - \textnormal{gr}_{\boldsymbol{v}}(p_1,p_2^*)$;
\item[(H4)] if $p_1$ or $p_2$ is a real point, then $\textnormal{gr}_{\boldsymbol{v}}(p_1,p_2) = 0$.
\item[(H5)] Let $(f_{h}(z))_{h = 1}^{H}$ be an $H$-tuple of functions such that $f_{h}(z)$ is a holomorphic function of $z$ in a neighborhood of $[\alpha_h,\beta_h]$ and takes real values on the real line. Denoting $f_{\boldsymbol{w}_i}$ the function on $\amsmathbb{B}_{\boldsymbol{w}_i}^+ \subset \Sigma_{\boldsymbol{v}}$ and equal to $f_{h} \circ Z_{\boldsymbol{v}}$ on the $h$-th segment when it is part of $\amsmathbb{B}_{\boldsymbol{w}_i}^+$ (\textit{i.e.} for $h \in \textnormal{supp}(\boldsymbol{w}_i)$), we have
\begin{equation}
\label{BthecontourH}
\begin{split}
& \sum_{h_1,h_2 = 1}^{H} (w_{1})_{h_1} (w_{2})_{h_2}\oint_{\gamma_{h_1}} \oint_{\gamma_{h_2}} \mathcal{F}_{h_1,h_2}(z_1,z_2) f_{h_1}(z_1) f_{h_2}(z_2) \dd z_1 \dd z_2 \\
& = -\frac{1}{4\pi^2} \iint_{\amsmathbb{B}_{\boldsymbol{w}_1}^{+} \times \amsmathbb{B}_{\boldsymbol{w}_2}^{+}} \textnormal{gr}_{\boldsymbol{v}}(p_1,p_2) \dd f_{\boldsymbol{w}_1}(p_1) \dd f_{\boldsymbol{w}_2}(p_2).
\end{split}
\end{equation}
\item[(H6)] The result of 5. is unchanged if we add to $\mathcal{B}_{\boldsymbol{v}}$ in the Definition~\ref{def:Hvv} of $\textnormal{gr}_{\boldsymbol{v}}$ a bidifferential proportional to $\frac{\dd z_1\dd z_2}{(z_1 - z_2)^2}$.
\end{enumerate}
\end{proposition}
\begin{proof}
We know that $\mathcal{B}_{\boldsymbol{v}}$ is a meromorphic bidifferential on $\Sigma_{\boldsymbol{v}} \times \Sigma_{\boldsymbol{v}}$ with double poles on $\Delta_{\boldsymbol{v}}$ only, having real biresidues. There are two possible sources of ambiguities in the definition of $\textnormal{gr}_{\boldsymbol{v}}$. First, the double integral of the double pole yields a logarithm, which is only defined up to multiple of $2\ii \pi \amsmathbb{Z}$. Since the biresidues (\textit{i.e.} the coefficients of these double poles) are real, the resulting ambiguity disappears when we extract the real part in \eqref{defHvp1p2}. Second, one can add to the integration path from $p_1^*$ to $p_1$ a non-contractible contour in the same (closure of) sheet. These are linear combination of contours going around the $h$-th segment in the sheet in question. By Proposition~\ref{thm:Bcont}, integrating one variable of $\mathcal{B}_{\boldsymbol{v}}$ along such a contour gives zero. So, there is no ambiguity in the definition of $\textnormal{gr}_{\boldsymbol{v}}$. Then, the symmetry of $\textnormal{gr}_{\boldsymbol{v}}$ comes from the symmetry of $\mathcal{B}_{\boldsymbol{v}}$ (Theorem~\ref{thm:Bcont}). Besides, $\textnormal{gr}_{\boldsymbol{v}}$ is manifestly continuous on $\Sigma_{\boldsymbol{v}} \times \Sigma_{\boldsymbol{v}} \setminus (\Delta_{\boldsymbol{v}} \cup \Delta_{\boldsymbol{v}}^*)$. For any $\boldsymbol{w}_1,\boldsymbol{w}_2 \in \mathfrak{G}.\boldsymbol{v}$, the double integral in \eqref{defHvp1p2} is a sum of a locally defined holomorphic function (both locally defined) in $(\Sigma_{\boldsymbol{v},\boldsymbol{w}_1} \times \Sigma_{\boldsymbol{v},\boldsymbol{w}_2}) \setminus (\Delta_{\boldsymbol{v}} \cup \Delta_{\boldsymbol{v}}^*)$. Consider the dependence in $p_i$, with the other variable fixed. Since the said (anti-)holomorphic function is continuous when $p_i$ belongs to the intersection of two sheet closures, they are actually locally defined (anti-)holomorphic function of $p_i \in \Sigma_{\boldsymbol{v}}$. Hartogs' theorem on separate holomorphicity \cite{Ebeling} says that these functions are then locally defined holomorphic in $(\Sigma_{\boldsymbol{v}} \times \Sigma_{\boldsymbol{v}}) \setminus (\Delta_{\boldsymbol{v}} \cup \Delta_{\boldsymbol{v}}^*)$. Taking the real part, we obtain $\textnormal{gr}_{\boldsymbol{v}}$ which is thus pluriharmonic in $(\Sigma_{\boldsymbol{v}} \times \Sigma_{\boldsymbol{v}}) \setminus ( \Delta_{\boldsymbol{v}} \cup \Delta_{\boldsymbol{v}}^*)$. This finishes the justification of 1. Property 2. follows by the same arguments in the neighborhood of $\Delta_{\boldsymbol{v}} \cup \Delta_{\boldsymbol{v}}^*$ in the $(\boldsymbol{w}_1,\boldsymbol{w}_2)$-th sheet for
\[
\text{Re}\left(\int_{p_1^*}^{p_1} \int_{p_2^*}^{p_2} \mathcal{B}_{\boldsymbol{v}} - \frac{\text{\textit{\foreignlanguage{russian}{B}}}{}_{\boldsymbol{w}_1,\boldsymbol{w}_2}\,\dd Z_{\boldsymbol{v},1}\dd Z_{\boldsymbol{v},2}}{(Z_{\boldsymbol{v},1} - Z_{\boldsymbol{v},2})^2}\right) = \textnormal{gr}_{\boldsymbol{v}}(p_1,p_2) + 2\text{\textit{\foreignlanguage{russian}{B}}}{}_{\boldsymbol{w}_1,\boldsymbol{w}_2} \log \left|\frac{z_1 - z_2^*}{z_1 - z_2}\right|.
\]
where $Z_{\boldsymbol{v},i}(p_1,p_2) = z_i$ for $i = 1,2$. The properties 3. and 4. are clear, using the already established fact that the definition of $\textnormal{gr}_{\boldsymbol{v}}$ does not depend on the choice of integration paths between $p_i^*$ and $p_i$.

To get property 5., we start by fixing $h_1,h_2 \in [H]$ and defining
\[
\tilde{\mathcal{F}}_{h_1,h_2}(z_1,z_2) = \int_{\infty}^{z_1}\int_{\infty}^{z_2} \mathcal{F}_{h_1,h_2}(\tilde{z}_1,\tilde{z}_2) \dd \tilde{z}_1\dd \tilde{z}_2
\]
for $z_i \in \widehat{\amsmathbb{C}} \setminus \overline{\amsmathbb{B}}_{h_i}$ with $i = 1,2$. Integrability along path ending at $\infty$ is guaranteed because $\mathcal{F}_{h_1,h_2}(z_1,z_2) = O(\frac{1}{z_i^2})$ as $z_i \rightarrow \infty$. As the domain of definition of $\mathcal{F}_{h_1,h_2}$ is simply-connected, the definition does not depend on the choice of contour of integration from $\infty$ to $z_i$, and we can take for instance vertical paths. The fact that $\big(\mathcal{F}_{h_1,h_2}(z_1,z_2)\big)^* = \mathcal{F}_{h_1,h_2}(z_1^*,z_2^*)$ then implies the same property for $\tilde{\mathcal{F}}_{h_1,h_2}$.

Integration by parts yields
\begin{equation}
\label{ointFh1h2}
\oint_{\gamma_{h_1}} \oint_{\gamma_{h_2}} \mathcal{F}_{h_1,h_2}(z_1,z_2) f_{h_1}(z_1) f_{h_2}(z_2) \dd z_1 \dd z_2 = \oint_{\gamma_{h_1}} \oint_{\gamma_{h_2}} \tilde{\mathcal{F}}_{h_1,h_2}(z_1,z_2) \dd f_{h_1}(z_2) \dd f_{h_2}(z_2).
\end{equation}
Since $\mathcal{F}_{h_1,h_2}(z_1,z_2)$ diverges at most like an inverse square-root as $z_i$ approaches $\alpha_{h_i}$ or $\beta_{h_i}$, $\tilde{\mathcal{F}}_{h_1,h_2}$ has a well-defined finite value as $z_i = \alpha_{h_i}$ or $\beta_{h_i}$. From Proposition~\ref{thmBfund} (F4), we know that $\mathcal{F}_{h_1,h_2}$ admits simultaneous upper- and lower-boundary values when its variables approach the segments $(\alpha_{h_1},\beta_{h_1})$ and $(\alpha_{h_2},\beta_{h_2})$, apart from a double pole on the diagonal. In consequence, $\tilde{\mathcal{F}}$ also admits upper- and lower-boundary values when its variables approach the same segments, apart from a logarithmic singularity on the diagonal. Then, we can can squeeze the contour to the segment, and more precisely we get
\begin{equation}
\label{ah1ah2bh1bh2}
\int_{\alpha_{h_1}}^{\beta_{h_1}} \!\!\!\int_{\alpha_{h_2}}^{\beta_{h_2}} \big(\tilde{\mathcal{F}}_{h_1,h_2}(x_1^-,x_2^-) - \tilde{\mathcal{F}}_{h_1,h_2}(x_1^-,x_2^+) - \tilde{\mathcal{F}}_{h_1,h_2}(x_1^+,x_2^-) + \tilde{\mathcal{F}}_{h_1,h_2}(x_1^+,x_2^+)\big) \dd f_{h_1}(x_1) \dd f_{h_2}(x_2).
\end{equation}
The logarithmic singularity does not pose problem because it is integrable. The assumption that $f_{h_i}$ take real values on the real line allows rewriting \eqref{ah1ah2bh1bh2} as
\begin{equation}
\label{ah1ah2bh1bh2rr}
\int_{\alpha_{h_1}}^{\beta_{h_1}} \!\!\!\int_{\alpha_{h_2}}^{\beta_{h_2}} \text{Re}\big(\tilde{\mathcal{F}}_{h_1,h_2}(x_1^+,x_2^+) - \tilde{\mathcal{F}}_{h_1,h_2}(x_1^-,x_2^+) - \tilde{\mathcal{F}}_{h_1,h_2}(x_1^+,x_2^-) + \tilde{\mathcal{F}}_{h_1,h_2}(x_1^-,x_2^-)\big) \dd f_{h_1}(x_1) \dd f_{h_2}(x_2).
\end{equation}
The quantity in bracket in the integrand is the limit, as $z_i$ approaches a point on $(\alpha_{h_i},\beta_{h_i})$ from the upper half-plane (simultaneously for $i = 1,2$) of
\begin{equation}
\label{ReFtildeHtilde}
\begin{split}
& \text{Re}\big(\tilde{\mathcal{F}}_{h_1,h_2}(z_1,z_2) - \tilde{\mathcal{F}}_{h_1,h_2}(z_1^*,z_2) - \tilde{\mathcal{F}}_{h_1,h_2}(z_1,z_2^*) + \tilde{\mathcal{F}}_{h_1,h_2}(z_1^*,z_2^*)\big) \\
& = \text{Re}\left(\int_{z_1^*}^{z_1} \int_{z_2^*}^{z_2} \mathcal{F}_{h_1,h_2}(\tilde{z}_1,\tilde{z}_2) \dd \tilde{z}_1 \dd \tilde{z}_2\right),
\end{split}
\end{equation}
where the $i$-th path of integration are taken in $\widehat{\amsmathbb{C}} \setminus \overline{\amsmathbb{B}}_{h_i}$. Now multiply \eqref{ah1ah2bh1bh2rr} with $(w_1)_{h_1}(w_2)_{h_2}$ and sum over $h_1,h_2 \in [H]$. Applying the same operation on \eqref{ReFtildeHtilde} yields the value of
\begin{equation}
\label{Hvp1p2ru}
\textnormal{gr}_{\boldsymbol{v}}(p_1,p_2) - 2\text{\textit{\foreignlanguage{russian}{B}}}{}_{\boldsymbol{w}_1,\boldsymbol{w}_2} \log \left|\frac{z_1 - z_2^*}{z_1 - z_2}\right|.
\end{equation}
when $p_i$ is in the $\boldsymbol{w}_i$-th sheet and $z_i = Z_{\boldsymbol{v}}(p_i)$. This allows rewriting \eqref{ah1ah2bh1bh2rr} in terms of a double integral of \eqref{Hvp1p2ru} with range $\amsmathbb{B}_{\boldsymbol{w}_1}^+ \times \amsmathbb{B}_{\boldsymbol{w}_2}^+ \subset \Sigma_{\boldsymbol{v}} \times \Sigma_{\boldsymbol{v}}$. However, $z_i$ is real on $\amsmathbb{B}_{\boldsymbol{w}_i}^-$ so the second term in \eqref{Hvp1p2ru} disappears in the integral and we get the desired formula \eqref{BthecontourH}. The property 6. comes from the fact that adding a real multiple of $\frac{\dd z_1 \dd z_2}{(z_1 - z_2)^2}$ to $\mathcal{B}_{\boldsymbol{v}}$ in the right-hand side of \eqref{defHvp1p2} results in adding to $\textnormal{gr}_{\boldsymbol{v}}$ the same multiple of $2\log \big|\frac{z_1 - z_2}{z_1 - z_2^*}\big|$, and we have seen that it does not change \eqref{BthecontourH}.
\end{proof}

The logarithmic singularities on the anti-diagonal $\Delta_{\boldsymbol{v}}^*$ become irrelevant if we keep only a hemisphere in each sheet of the spectral curve. There is a however a consistency condition to do so.

\begin{definition}
\label{Definition_bipartite}
Let $\boldsymbol{v} \in \amsmathbb{R}^H$ be non-zero. We say that the $\mathfrak{G}$-orbit of $\boldsymbol{v}$ is bipartite if there exists a function $\epsilon : \mathfrak{G}.\boldsymbol{v} \rightarrow \{\pm 1\}$ such that $\epsilon(T^{(h)}(\boldsymbol{w})) = - \epsilon(\boldsymbol{w})$ for any $\boldsymbol{w} \in \mathfrak{G}.\boldsymbol{v}$ and $h \in [H]$.
\end{definition}

\begin{definition}[Half-spectral curve]
\label{def:halfspcurvebip}
Given $\boldsymbol{v} \in \amsmathbb{R}^H$ non-zero with bipartite $\mathfrak{G}$-orbit, we define $\Sigma_{\boldsymbol{v}}^{\textnormal{half}}$ to be the image in the quotient $\Sigma_{\boldsymbol{v}}$ (\textit{cf.} Definition~\ref{def:spcurv1}) of
\[
\bigcup_{\boldsymbol{w} \in \mathfrak{G}.\boldsymbol{v}} \widehat{\amsmathbb{H}}^{\epsilon(\boldsymbol{w})}.
\]
where $\epsilon$ is the unique bipartition of $\mathfrak{G}.\boldsymbol{v}$ such that $\epsilon(\boldsymbol{v}) = 1$. This is a bordered Riemann surface, whose boundary is the closure of the set of real points, that is the image in $\Sigma_{\boldsymbol{v}}$ of
\[
\bigcup_{\boldsymbol{w} \in \mathfrak{G}.\boldsymbol{v}} \big(\widehat{\amsmathbb{R}}^{\epsilon(\boldsymbol{w})} \setminus \amsmathbb{B}_{\boldsymbol{w}}^{\epsilon(\boldsymbol{w})}\big).
\]
See Figure~\ref{Fig:Half5} and \ref{Fig:Half7} for an example.
\end{definition}

Then, for special matrix of biresidues, we can relate $\textnormal{gr}_{\boldsymbol{v}}$ (and thus the fundamental solution of the master problem) to the usual notion of Green function.

\begin{corollary}
\label{Greenwon}
Assume that $\pi_{\boldsymbol{v}}^{\circ}$ is a bijection, that $\boldsymbol{v}$ has bipartite $\mathfrak{G}$-orbit and $\textbf{\text{\foreignlanguage{russian}{B}}} = \textnormal{\textbf{Id}} - \boldsymbol{1}$. Then, the Green function on $\Sigma_{\boldsymbol{v}}^{\textnormal{half}}$ with Dirichlet boundary conditions on $\partial\Sigma_{\boldsymbol{v}}^{\textnormal{half}}$ is
\begin{equation}
\label{Greenock}
\textnormal{Gr}(p_1,p_2) = -\frac{1}{4\pi}\left(\textnormal{gr}_{\boldsymbol{v}}(p_1,p_2) + 2\log\bigg|\frac{z_1 - z_2}{z_1 - z_2^*}\bigg|\right).
\end{equation}
\end{corollary}
\begin{proof}
The assumption on the matrix of biresidues implies together with Proposition~\ref{prop:Greengen} that $\textnormal{gr}_{\boldsymbol{v}}(p_1,p_2) + 2\log\big|\frac{z_1 - z_2}{z_1 - z_2^*}\big|$ on $\Sigma_{\boldsymbol{v}} \times \Sigma_{\boldsymbol{v}}$ is harmonic in both variables in $\Sigma_{\boldsymbol{v}} \times \Sigma_{\boldsymbol{v}}$ except for logarithmic singularities located on the diagonal $p_1 = p_2$ with coefficient $2$ and the anti-diagonal $p_1 = p_2^*$ with coefficient $-2$, irrespectively of the sheet closures in which $p_1 = p_2$ is located. The half spectral curve was constructed so that there are no pairs of points $p_1,p_2$ in the same sheet having $p_1 = p_2^*$, so that restricting to $\Sigma_{\boldsymbol{v}}^{\textnormal{half}} \times \Sigma_{\boldsymbol{v}}^{\textnormal{half}}$ we only see the logarithmic singularities on the diagonal. Multiplication by $-\frac{1}{4\pi}$ in the definition of \eqref{Greenock} sets the coefficient of this logarithmic singularity to $-\frac{1}{2\pi}$. Besides, Proposition~\ref{prop:Greengen} 4. indicates that the right-hand side of \eqref{Greenock} vanishes if $p_1$ or $p_2$ is a real point of a sheet of $\Sigma_{\boldsymbol{v}}$, that is of the boundary of $\Sigma_{\boldsymbol{v}}^{\textnormal{half}}$. This vanishing and the harmonicity with the said logarithmic singularities imply that \eqref{Greenock} is indeed a Green function.
\end{proof}

\subsection{Functions of the three kinds}
 \label{sec:3kinds}
In Section~\ref{sec:1233} we have described solutions of the master Riemann--Hilbert problem of first kind, second kind of type I,II, III, and third kind, and their properties were summarized in Lemma~\ref{lem:genRHPsol}. We revisit this classification by discussing how they define meromorphic quantities on (parts of) the spectral curve with specific properties. This explains a posteriori some peculiarities of the construction of Section~\ref{sec:1233}.

\medskip

Let $\boldsymbol{v} \in \amsmathbb{R}^H$ be a non-zero vector and $\Sigma_{\boldsymbol{v}}$ the spectral curve coming from Section~\ref{sec:Algapproach}. We recall that the global coordinate $z$ on the Riemann sphere glues to a branched covering map $Z : \Sigma_{\boldsymbol{v}} \rightarrow \widehat{\amsmathbb{C}}$. This means that $Z$ can be used as a local coordinate around any point of $\Sigma_{\boldsymbol{v}}$ which is a not a ramification point, \textit{i.e} a zero of $\dd Z$. By construction in Section~\ref{sec:Algapproach} the zeros of $\dd Z$ are simple and their $Z$-image (the branch points) are among $\alpha_1,\beta_1,\ldots,\alpha_H,\beta_H$. Near any ramification point above $\alpha_h$ or $\beta_h$, $Z$ is not a local coordinate on $\Sigma_{\boldsymbol{v}}$ but $\sqrt{(Z - \alpha_h)(Z - \beta_h)}$ is. This local coordinate vanishes at such ramification points. In particular, the $1$-form
\[
\frac{\dd Z}{\sqrt{(Z - \alpha_h)(Z - \beta_h)}}
\]
is holomorphic near such ramification points. In the following, we will also use the simpler Riemann surface $\Sigma^{\sharp, g}$ (met in the proof of Theorem~\ref{spcurveRHP}) obtained by gluings two pairs of hemispheres along $\overline{\amsmathbb{B}}_h$. It still admits the projection map $Z : \Sigma^{\sharp,h} \rightarrow \widehat{\amsmathbb{C}}$ and the same discussion about local coordinates on $\Sigma^{\sharp,h}$ applies.

\medskip

\noindent \textsc{First kind.} First kind differentials on compact Riemann surfaces are holomorphic differentials that are normalized on $\mathsf{a}$-cycles, and they can be obtained by $\mathsf{b}$-cycle integrals of the fundamental bidifferential of the second kind on the Riemann surface. The first kind functions of Definition~\ref{def:1stkind} have a similar flavor. First, from their definition we have
\[
\mathfrak{c}^{\textnormal{1st}}_{h;g}(z) = \int_{c_{o_g^+}} \big\langle \bth_{g} \cdot \boldsymbol{\mathcal{F}}_{h,\bullet}(z,\zeta)\big\rangle\dd \zeta + \int_{c_{o_g^-}} \left(\big\langle \bth_{g} \cdot \boldsymbol{\mathcal{F}}_{h,\bullet}(z,\zeta)\big\rangle + \frac{\delta_{g,h}}{(z - \zeta)^2}\right)\dd \zeta
\]
for $o_g \in \amsmathbb{B}_g$ and $c_{o_g^{\pm}}$ a contour from $\infty$ to $o_g$ approaching $o_g$ from the upper (respectively, lower) half-plane. On the other hand, for any $x \in \amsmathbb{B}_g$ due to the Riemann--Hilbert problem solved by $\boldsymbol{\mathcal{F}}$, \textit{cf.} Proposition~\ref{thmBfund}:
\begin{equation}
\label{lvirrgwrg}
\big\langle \bth_{g} \cdot \boldsymbol{\mathcal{F}}_{h,\bullet}(z,x^+)\big\rangle = - \left( \big\langle \bth_g \cdot \boldsymbol{\mathcal{F}}_{h,\bullet}(z,x^-)\big\rangle + \frac{\delta_{g,h}}{(z - x)^2}\right).
\end{equation}
Considering $z$ fixed and $\zeta$ as the variable, \eqref{lvirrgwrg} and the arguments at work in Section~\ref{sec:Algapproach} show the existence of a meromorphic $1$-form without residues on $\Sigma^{\sharp,g}$ which is equal to
\begin{equation*}
\begin{split}
\big\langle \bth_{g} \cdot \boldsymbol{\mathcal{F}}_{h,\bullet}(z,\zeta)\big\rangle \dd\zeta \qquad & \textnormal{for}\,\, \zeta \in \Sigma_+^{\sharp,g}, \\
-\left(\big\langle\bth_g \cdot \boldsymbol{\mathcal{F}}_{h,\bullet}(z,\zeta)\big\rangle + \frac{\delta_{g,h}}{(z - \zeta)^2}\right)\dd\zeta,\qquad & \textnormal{for}\,\,\zeta \in \Sigma_-^{\sharp,g}.
\end{split}
\end{equation*}
As a consequence, we may rewrite $\mathfrak{c}^{\textnormal{1st}}_{h;g}(z)$ as the integral of this meromorphic $1$-form on a path in $\Sigma^{\sharp,g}$ from $\infty \in \Sigma_+^{\sharp,g}$ to $\infty \in \Sigma_-^{\sharp,g}$. In this form the independence of the definition on the point $o_g$ at which the path intersects $\amsmathbb{B}_g$ is clear: moving $o_g$ just amounts to a deformation of the path in the same homotopy class relative to the endpoints in $\Sigma^{\sharp,g}$. The meromorphic $1$-form we constructed on $\Sigma^{\sharp,g}$ has a double pole at $\zeta = z$ for $\zeta \in \Sigma^{\sharp,g}_-$, sharing some similarity with fundamental bidifferentials of the second kind. The contours in $\Sigma_{\pm}^{\sharp,g}$ going counterclockwise around $\overline{\amsmathbb{B}}_h$ have intersection $\delta_{g,h}$ with the path between the two $\infty$. The latter plays the role of $\mathsf{b}_g$-cycle, while the former have the role of the $\mathsf{a}_g$-cycles.

The $H$-tuples $(\mathfrak{c}^{\textnormal{1st}}_{h;g}(z))_{h = 1}^{H}$ depending on $g \in [H]$ are solutions of the master Riemann--Hilbert problem with zero source, diverging at most like an inverse square root at $\alpha_h,\beta_h$ (tame exponent $1$) and behaving like $\frac{\delta_{h,g}}{z}$ as $z \rightarrow \infty$. Taking non-zero $\boldsymbol{v},\boldsymbol{v}' \in \amsmathbb{R}^H$, we get a $1$-form
\[
\big\langle \boldsymbol{v} \otimes \boldsymbol{v}' \cdot \mathfrak{c}^{\textnormal{1st}}_{\bullet;\bullet}(z)\big\rangle \dd z = \sum_{h,g = 1}^{H} v_h v'_g \mathfrak{c}^{\textnormal{1st}}_{h;g}(z) \dd z
\]
which continues analytically to a \emph{meromorphic} $1$-form on $\Sigma_{\boldsymbol{v}}$ which we denote\label{index:ufrak} $\mathfrak{u}_{\boldsymbol{v};\boldsymbol{v}'}$. It can have simple poles at the $Z$-preimages of $\infty$. Its value in the $\boldsymbol{w}$-th sheet is $\big\langle \boldsymbol{w} \otimes \boldsymbol{v}' \cdot \mathfrak{c}^{\textnormal{1st}}_{\bullet;\bullet}(z) \big\rangle \dd z$. By Lemma~\ref{lem:genRHPsol}, its integral on the contour $\gamma_{l} \subset \Sigma_{\boldsymbol{v},\boldsymbol{w}}$ is
\begin{equation}
\label{periodubar}\frac{1}{2\ii\pi} \oint_{\gamma_{l} \subset \Sigma_{\boldsymbol{v},\boldsymbol{w}}} \mathfrak{u}_{\boldsymbol{v};\boldsymbol{v}'} = w_{l}v'_l.
\end{equation}
This identity plays the role of the normalization on $\mathsf{a}$-cycles. The residue at $\infty$ in the $\boldsymbol{w}$-th sheet is
\begin{equation}
\label{residueatinf}
\sum_{l = 1}^{H} \frac{1}{2\ii\pi} \oint_{\gamma_{\boldsymbol{w},l}} \mathfrak{u}_{\boldsymbol{v};\boldsymbol{v}'} = \langle \boldsymbol{w} \cdot \boldsymbol{v}' \rangle.
\end{equation}
Since $\mathscr{V} = \textnormal{Ker}(\boldsymbol{\Theta})^{\bot}$ is preserved by the action of the group $\mathfrak{G}$, we deduce the following result, which is relevant in discrete ensembles with fixed filling fractions to compute the variation of the equilibrium measure with respect to filling fractions, \textnormal{cf.} Proposition~\ref{Proposition_differentiability_filling_fraction}.

\begin{proposition}
\label{Lem:Holo1form}
Let $\boldsymbol{v} \in \mathscr{V} = \textnormal{Ker}(\boldsymbol{\Theta})^{\bot}$ and $\boldsymbol{v}' \in \mathscr{V}^{\bot} = \textnormal{Ker}(\boldsymbol{\Theta})$. Then, $\mathfrak{u}_{\boldsymbol{v};\boldsymbol{v}'}$ is the unique \emph{holomorphic} $1$-form on $\Sigma_{\boldsymbol{v}}$, satisfying the period conditions \eqref{periodubar}.
\end{proposition}

\medskip

\noindent \textsc{Second kind.} Second-kind differentials on compact Riemann surfaces are meromorphic differentials having poles without residues, and having zero $\mathsf{a}$-periods. The functions of second kind in Definitions~\ref{def:2ndkindI}-\ref{def:2ndkindII} have similar properties but the fact the Riemann surface that we construct by gluing and on which these functions admit analytic continuation has local coordinate $z$ away from ramification points and $\sqrt{(z - \alpha_h)(z - \beta_h)}$ near ramification points forces us to distinguish between three types, depending on the position and nature of the poles.

Let us look at the second-kind function of type I of Definition~\ref{def:2ndkindI} for an $H$-tuple of \emph{rational functions} $\boldsymbol{f}(z) = (f_g(z))_{g = 1}^{H}$ having poles in $\widehat{\amsmathbb{C}} \setminus \overline{\amsmathbb{B}}$. By Lemma~\ref{lem:genRHPsol}, this is an $H$-tuple $\mathfrak{F}^{\textnormal{I}}[\boldsymbol{f}](z)$ solving the master Riemann--Hilbert problem with source
\begin{equation}
\label{Dhegrgugbh}
\forall h \in [H] \qquad D_h^{\textnormal{I}}[\boldsymbol{f}](z) = -\partial_z f_h^{\textnormal{I}}(z)
\end{equation}
and having tame exponent $1$ and decaying as $O(\frac{1}{z^2})$ as $z \rightarrow \infty$. The first consequence is that for each $h \in [H]$ the $1$-form $\mathfrak{F}_{h}^{\textnormal{I}}[\boldsymbol{f}](z)\dd z$ analytically continues to a meromorphic $1$-form $\Psi_h^{\textnormal{I}}[\boldsymbol{f}]$ on a Riemann surface $\Sigma_{(\boldsymbol{e}^{(h)},\boldsymbol{0})}$ constructed as in Section~\ref{inHomo} (recall that in general $\Sigma_{(\boldsymbol{e}^{(h)},\boldsymbol{0})}$ covers $\Sigma_{\boldsymbol{e}^{(h)}}$, although these two Riemann surfaces may coincide in special cases). Tameness with exponent $1$ implies that $\Psi_h^{\textnormal{I}}[\boldsymbol{f}]$ is regular at the ramification points in $\Sigma_{(\boldsymbol{e}^{(h)},\boldsymbol{0})}$. The fact that the source \eqref{Dhegrgugbh} has poles without residues leads to poles without residues for $\Psi_h^{\textnormal{I}}[\boldsymbol{f}]$ at the points of whose $Z$-projection are poles of $\boldsymbol{f}$. And, if $\gamma_l$ is a contour in the principal sheet of $\Sigma_{(\boldsymbol{e}^{(h)},\boldsymbol{0})}$ going counterclockwise around $[\alpha_l,\beta_l]$, we have
\[
\oint_{\gamma_l} \Psi_{h}^{\textnormal{I}}[\boldsymbol{f}] = 0.
\]
This identity comes from the analytic properties of $\mathfrak{F}_{h}^{\textnormal{I}}[\boldsymbol{f}](z)$ for $l \neq h$, and from its decay property at infinity for $l = h$. It is the analog of the normalization on $\mathsf{a}$-cycles.

Unlike in type I, the functions of second kind of type II and III have their singularities at $\alpha_h,\beta_h$. Let $\boldsymbol{f} = (f_g(z))_{g = 1}^{H}$ be an $H$-tuple of rational functions having poles only at $\alpha_g$ and $\beta_g$.
Consider the definition of the second-kind functions of type II in Definition~\ref{def:2ndkindII}, namely
\begin{equation}
\label{bIIII}
\mathfrak{F}_h^{\textnormal{II}}[\boldsymbol{f}](z) = \sum_{g = 1}^H \oint_{\gamma(\alpha_g) + \gamma(\beta_g)} \frac{\dd \zeta}{2\ii\pi}\left(\frac{\delta_{g,h}}{2(z - \zeta)^2} + \big \langle \bth_g \cdot \mathcal{F}_{h,\bullet}(z,\zeta)\big \rangle\right) \sigma_g(\zeta) f_g(\zeta),
\end{equation}
where $\gamma(\alpha_g)$ and $\gamma(\beta_g)$ are small counterclockwise loops around $\alpha_g$ and $\beta_g$ in the $\zeta$-plane. As a function of $\zeta$ and when $\zeta$ crosses $\amsmathbb{B}_g$ from the upper- to the lower half-plane, the $g$-th integrand is continuous because $f_g(\zeta)$ is continuous, $\sigma_g(\zeta)$ changes its sign and the quantity in bracket (due to the Riemann--Hilbert problem solved by $\boldsymbol{\mathcal{F}}$, \textit{cf.} Proposition~\ref{thmBfund}) also changes sign. Therefore, the $g$-th integrand is meromorphic in a neighborhood of $[\alpha_g,\beta_g]$, and the contour integrals are actually residues at the poles $\alpha_g$ and $\beta_g$. This explains why the definition \eqref{bIIII} does not depend on the choice of small loops around $\alpha_g$ and $\beta_g$. The properties of $\boldsymbol{\mathfrak{F}}^{\textnormal{II}}[\boldsymbol{f}]$ described in Lemma~\ref{lem:genRHPsol} imply that for each $h \in [H]$, $\mathfrak{F}_{h}^{\textnormal{II}}[\boldsymbol{f}](z)\dd z$ analytically continues to a meromorphic $1$-form on the Riemann surface $\Sigma_{\boldsymbol{e}^{(h)}}$, with poles only at the ramification points. We denote $\Psi_h^{\textnormal{II}}[\boldsymbol{f}]$ this meromorphic $1$-form. Near the ramification points $\alpha_h$ or $\beta_h$ in the principal sheet of $\Sigma_{\boldsymbol{e}^{(h)}}$ it is such that
\[
\Psi_h^{\textnormal{II}}[\boldsymbol{f}] = \frac{\dd(\sigma_h \cdot f_h)}{2} + \textnormal{regular}.
\]
In the local coordinate $\sigma_h$ the principal part (the first term in the right-hand side) is odd, because $f_h(z)$ is a rational function of $z$ and $z$ is locally an even function of $\sigma_h(z)$. Since this principal part is a differential of a (locally defined) meromorphic function, $\Psi_h^{\textnormal{II}}[\boldsymbol{f}]$ has no residue at these poles on $\Sigma_{\boldsymbol{e}^{(h)}}$. At ramification points in other sheets of $\Sigma_{\boldsymbol{e}^{(h)}}$, $\Psi_h^{\textnormal{II}}[\boldsymbol{f}]$ is a linear combination of $\mathfrak{F}_{g}^{\textnormal{II}}[\boldsymbol{f}](z)\dd z$ for $g \in [H]$ hence also has no residues and odd principal part.

One can make a similar discussion for the second-kind functions of type III. In Definition~\ref{def:2ndkindII} for $\boldsymbol{f} = (f_g(z))_{g = 1}^{H}$ an $H$-tuple of rational functions of $z$ with poles at $\alpha_g$ and $\beta_g$, these are given by
\[
\mathfrak{F}_h^{\textnormal{III}}[\boldsymbol{f}](z) = \sum_{g = 1}^H \oint_{\gamma(\alpha_g) + \gamma(\beta_g)} \frac{\dd\zeta}{2\ii\pi} \left(\frac{1}{2(z - \zeta)^2} + \sum_{g' \neq g} \theta_{g,g'} \mathcal{F}_{h,g'}(z,\zeta) \right)\frac{f_g(\zeta)}{\theta_{g,g}}.
\]
and the $g$-th integrand is manifestly a meromorphic $1$-form of $\zeta$ in a neighborhood of $\overline{\amsmathbb{B}}_g$. Using the properties established in Lemma~\ref{lem:genRHPsol}, $\mathfrak{F}_{h}^{\textnormal{III}}[\boldsymbol{f}](z)\dd z$ analytically continues to a meromorphic $1$-form on a Riemann surface covering $\Sigma_{\boldsymbol{e}^{(h)}}$ (because the master Riemann--Hilbert problem it satisfies has non-zero source). We denote $\Psi_h^{\textnormal{III}}[\boldsymbol{f}]$ this meromorphic $1$-form. At the ramification points $\alpha_h$ and $\beta_h$ in the principal sheet it satisfies
\[
\Psi_h^{\textnormal{III}}[\boldsymbol{f}] = \frac{\dd f_h}{2\theta_{h,h}} + \textnormal{regular}.
\]
Since $f_h(z)$ is a rational function of $z$, it is an even function of the local coordinate $\sigma_h(z)$. So, the poles of $\Psi_h^{\textnormal{III}}[\boldsymbol{f}]$ at these points have no residue and even principal part. At ramification points in other sheets, $\Psi_h^{\textnormal{III}}[\boldsymbol{f}]$ is a linear combination of $\mathfrak{F}_{g}^{\textnormal{III}}[\boldsymbol{f}](z)\dd z$ and of the sources
\[
D_g^{\textnormal{III}}[\boldsymbol{f}](z)\dd z = \left(\oint_{\gamma(\alpha_g) + \gamma(\beta_g)} \frac{\dd \zeta}{2\ii\pi}\,\frac{f_g(\zeta)}{(z - \zeta)^2}\right)\dd z = - \dd f_g(z),
\]
for $g \in [H]$ and this implies that these poles also have no residues and even principal part.

\medskip

\noindent \textsc{Third kind.} Differentials of the third kind on a compact Riemann surface are meromorphic with simple poles, and can be obtained by integrating one variable of the fundamental bidifferential of the second kind on path between various endpoints (that will become the poles). Definition~\ref{def:3rdkind} for our functions of the third kind is similar. Following Lemma~\ref{lem:genRHPsol}, the $1$-form $\mathfrak{c}_{h;g}^{\textnormal{3rd}}(z;\zeta) \dd z$ continues analytically to a meromorphic $1$-form on a Riemann surface covering $\Sigma_{\boldsymbol{e}^{(h)}}$ has simple poles at points whose $Z$-projection is $\zeta$. Due to the expression for the source in the master Riemann--Hilbert problem it satisfies, it can also have simple poles at points at $\infty$ in other sheets than the principal one.

\section{Matrices \texorpdfstring{$\boldsymbol{\Theta}$}{Theta} from gluing graphs}
\label{sec:green}

The tiling models defined in Section~\ref{Section_gluing_def} are special cases of the discrete ensembles of Section~\ref{Section_general_model}. Their interaction matrix $\boldsymbol{\Theta}$ describes the combinatorics of gluings of trapezoids from which the tiled domain is formed (Definition~\ref{Definition_gluing}). The topology of the gluing can be described as follows. We take $m$ subsets $\mathcal{H}_1,\ldots,\mathcal{H}_m$ of $[H]$ such that each $h \in [H]$ belongs to exactly two sets $\mathcal{H}_{i_h^+},\mathcal{H}_{i_h^-}$. Then, we take $m$ trapezoids $\mathcal{T}_1,\ldots,\mathcal{T}_m$ so that for each $i \in [H]$ the $i$-th trapezoid has segments labeled by $\mathcal{H}_i$ along its long base, and we glue for each $h \in [H]$ the trapezoids $\mathcal{T}_{i_h^-}$ and $\mathcal{T}_{i_h^+}$ along the $h$-th segment. Equivalently, this is encoded in the \emph{gluing graph}, denoted $\mathsf{G}$, whose vertices are trapezoids, and each segment shared by two trapezoids is indicated by an edge between the two corresponding vertices. We label the edges from $1$ to $H$ (we care about the order the edges appear), but the $m$ vertices are unlabeled (we do not care about the labeling of trapezoids).

There may be obstructions to realize this gluing metrically with trapezoids, and obstructions for the resulting domain to have a single vertical axis and to be tileable --- these are the conditions under which we get a discrete ensemble fitting in the setting of Section~\ref{Section_general_model}. However, the interaction matrix $\boldsymbol{\Theta}$ is defined for any undirected graph $\mathsf{G}$ with $H$ labeled edges and $m$ vertices:
\begin{equation}
\label{ThetaG}
\boldsymbol{\Theta} = \frac{1}{2} \sum_{i \in [m]} \textnormal{\textbf{1}}_{\mathcal{H}_i},
\end{equation}
where $\textnormal{\textbf{1}}_{\mathcal{H}_i}$ is the matrix whose $(g,h)$ entry is $1$ if $g,h \in \mathcal{H}_i$, and $0$ otherwise, and the subsets $\mathcal{H}_1,\ldots,\mathcal{H}_m$ record the edge-vertex adjacency. In this section, we want to study the master Riemann--Hilbert problems associated to such matrices $\boldsymbol{\Theta}$, and when we talk about domains $\mathcal{D}$ we only mean the ones obtained by a gluing of topological trapezoids (= closed disks, or equivalently, hemispheres) with the procedure that we have just described. We will make this discussion independent of whether there exists an actual tiling model (in view of the aforementioned obstructions) realizing this topological gluing or not. It is sufficient to restrict ourselves to connected gluing graphs $\mathsf{G}$, \textit{i.e.} connected domains $\mathcal{D}$.

In Section~\ref{sec:gluingcombi} we start by studying in greater detail the properties of matrices $\boldsymbol{\Theta}$ coming from gluing graphs. In Section~\ref{sec:C13tiling} we construct the spectral curve and the fundamental solutions of their associated master Riemann--Hilbert problem. Remarkably, we find that the spectral curves are compact and their topology is closely related to the topology of the glued domain.

\subsection{Topology and combinatorics of gluings}
\label{sec:gluingcombi}

Let us take a gluing graph $\mathsf{G}$. It is faithfully encoded in the vertex-edge adjacency matrix $\mathbf{R}$. This is a $m \times H$ matrix with entries in $\{0,1\}$, with non-zero rows and such that the sum in each column is $2$. Its lines are the linear forms
\begin{equation}
\label{imlinform}
\forall i \in [m] \qquad \forall \boldsymbol{X} \in \amsmathbb{R}^H\qquad \mathfrak{r}_i(\boldsymbol{X}) = \sum_{h \in \mathcal{H}_i} X_h
\end{equation}
When there is an underlying tiling model in such a domain, these are in fact the affine constraints \eqref{eq_trapezoid_restriction} for the segment filling fractions. The interaction matrix $\boldsymbol{\Theta}$ then takes the form
\begin{equation}
\label{ThetaGeq}
\boldsymbol{\Theta} = \frac{1}{2} \sum_{i = 1}^{m} \textnormal{\textbf{1}}_{\mathcal{H}_i \times \mathcal{H}_i} = \frac{1}{2}\, \mathbf{R}^{T}\mathbf{R}.
\end{equation}
In particular, it is positive semi-definite. We also introduce for each $i \in [m]$ the $H$-dimensional vector $\textnormal{\textbf{1}}_{\mathcal{H}_i} = \sum_{h \in \mathcal{H}_i} \boldsymbol{e}^{(h)}$. Then $\mathfrak{r}_i(\boldsymbol{X}) = \langle \textnormal{\textbf{1}}_{\mathcal{H}_i} \cdot \boldsymbol{X} \rangle$.

We noticed in Section~\ref{sec:genglu}, \textit{cf.} Figure~\ref{Fig:Domains23}, that many (but not all) domains could be drawn such that each segment has a trapezoid to its left and a trapezoid to its right. This is the case exactly when the gluing graph is bipartite, meaning that there is an assignment of signs $\pm 1$ to vertices so that adjacent vertices have opposite signs; we then say that the domain $\mathcal{D}$ is bipartite. In this situation we get a decomposition
\begin{equation}
\label{Thetapmdec}
\boldsymbol{\Theta} =\frac{1}{2}(\boldsymbol{\Theta}^{+} + \boldsymbol{\Theta}^-), \qquad \boldsymbol{\Theta}^\pm = \sum_{i \in M_{\pm}} \textnormal{\textbf{1}}_{\mathcal{H}_i \times \mathcal{H}_i},
\end{equation}
where $M_{+} \sqcup M_- = [m]$ is the bipartition and $\{\mathcal{H}_i\,\,|\,\,i \in M_{+}\}$ and $\{\mathcal{H}_i \,\,|\,\,i \in M_-\}$ give two partitions of $[H]$. In that case, each $h \in [H]$ belongs to exactly one $\mathcal{H}_{i_h^+}$ and one $\mathcal{H}_{i_h^-}$ for $i_h^{\pm} \in M_{\pm}$.

\begin{proposition}
\label{lem:H1G}We have $H_1(\mathcal{D},\amsmathbb{Z}) \simeq H_1(\mathsf{G},\amsmathbb{Z})$ and $\textnormal{Im}(\boldsymbol{\Theta}) = \textnormal{span}_{\amsmathbb{R}}(\mathbf{1}_{\mathcal{H}_1},\ldots,\mathbf{1}_{\mathcal{H}_m})$.
\begin{itemize}
\item If $\mathsf{G}$ is bipartite, then we have $m = \textnormal{rank}(\boldsymbol{\Theta}) + 1$ and $H_1(\mathsf{G},\amsmathbb{R}) \cong \textnormal{Ker}(\boldsymbol{\Theta})$.
\item If $\mathsf{G}$ is non-bipartite, then we have $m = \textnormal{rank}(\boldsymbol{\Theta})$ and $\dim H_1(\mathsf{G},\amsmathbb{R}) = \textnormal{corank}(\boldsymbol{\Theta}) + 1$.
\end{itemize}
\end{proposition}
\begin{proof}
By construction, the domain $\mathcal{D}$ strongly retracts onto $\mathsf{G}$, giving the first isomorphism. As $\mathsf{G}$ is connected the Euler relation gives $1 - \mathsf{b}_1(\mathsf{G}) = m - H$, where $\mathsf{b}_1(\mathsf{G}) = \dim H_1(\mathsf{G},\amsmathbb{Z})$ is the number of independent cycles in $\mathsf{G}$.

The computation of $\textnormal{Im}(\boldsymbol{\Theta})$ is easy. The image of any linear map and the kernel of its transpose have zero intersection, so we have $\textnormal{Im}(\boldsymbol{\Theta}) = \textnormal{Im}(\mathbf{R}^{T} \mathbf{R}) = \textnormal{Im}(\mathbf{R}^{T})$, and this is spanned by the columns of $\mathbf{R}^{T}$ which are the vectors $\mathbf{1}_{\mathcal{H}_1},\ldots,\mathbf{1}_{\mathcal{H}_m}$. We also have $\textnormal{Ker}(\boldsymbol{\Theta}) = \textnormal{Ker}(\mathbf{R}^T\mathbf{R}) = \textnormal{Ker}(\mathbf{R})$.

This does not give us yet the dimension of $\textnormal{Ker}(\boldsymbol{\Theta})$ and $\textnormal{Im}(\boldsymbol{\Theta})$, because there could be linear relations between the linear forms \eqref{imlinform}, or equivalently between the vectors $\mathbf{1}_{\mathcal{H}_1},\ldots,\mathbf{1}_{\mathcal{H}_m}$. Assume we have a non-trivial relation $\sum_{i = 1}^{m} v_i \mathbf{1}_{\mathcal{H}_i} = 0$ for $v_1,\ldots,v_m \in \amsmathbb{R}$. Then, for any $h \in [H]$ we have
\begin{equation}
\label{eq:balance} v_{i_h} + v_{i'_h} = 0,
\end{equation}
where the $h$-th segment is shared between the trapezoids $\mathcal{T}_{i_h}$ and $\mathcal{T}_{i'_h}$. If there existed a $j \in [H]$ such that $v_{j} = 0$, from the connectedness condition (c) we would deduce inductively that $v_i = 0$ for all $i \in [m]$. As we are looking at non-trivial relations, this shows that for any $i \in [m]$ we have $v_i \neq 0$. By the same connectedness argument we infer from \eqref{eq:balance} that $v_i = \epsilon_i v_1$ for some $\epsilon_i \in \{\pm 1\}$, and if $\mathcal{T}_i$ and $\mathcal{T}_{i'}$ share a segment we have $\epsilon_i + \epsilon_{i'} = 0$. In other words, $\mathsf{G}$ must be bipartite and there is (unique up to scale) a single relation between $(\mathbf{1}_{\mathcal{H}_1},\ldots,\mathbf{1}_{\mathcal{H}_m})$. Reading the argument from the end, we arrive to the following alternative.
\begin{itemize}
\item If $\mathsf{G}$ is bipartite, then $\textnormal{rank}(\mathbf{1}_{\mathcal{H}_1},\ldots,\mathbf{1}_{\mathcal{H}_m}) = m - 1$. Hence $\textnormal{rank}(\boldsymbol{\Theta}) = m -1$ and $\mathsf{b}_1(\mathsf{G}) = \textnormal{corank}(\boldsymbol{\Theta})$.
\item If $\mathsf{G}$ is not bipartite, then $\textnormal{rank}(\mathbf{1}_{\mathcal{H}_1},\ldots,\mathbf{1}_{\mathcal{H}_m}) = m$. Therefore $\textnormal{rank}(\boldsymbol{\Theta}) = m$ and $\mathsf{b}_1(\mathsf{G}) = \textnormal{corank}(\boldsymbol{\Theta}) + 1$.
\end{itemize}

Elements in $H_1(\mathsf{G},\amsmathbb{Z})$ are formal sums of oriented cycles with integral coefficients. Let us choose an arbitrary orientation for each edge. Then $H_1(\mathsf{G},\amsmathbb{Z})$ coincides with the set of maps $X : [H] \rightarrow \amsmathbb{Z}$ such that at each vertex the sum of $X_h$ for incoming edges equal to the sum of $X_h$ for outcoming edges ($X_h$ is the number of times the $h$-th edge is followed by a cycle, signed according to orientation). If $\mathsf{G}$ is bipartite, we can choose a bipartition $\epsilon : [m] \rightarrow \{\pm 1\}$ and then orient edges from negative to positive vertices. The advantage of the orientation is that the constraints at the $i$-th vertex takes the form
\[
\forall i \in [m] \qquad \epsilon_i \sum_{h \in \mathcal{H}_i} X_{h} = \epsilon_i \mathfrak{r}_i(\boldsymbol{X}) = 0.
\]
Therefore, $H_1(\mathsf{G},\amsmathbb{R}) \simeq \textnormal{Ker}(\mathbf{R})$, but $2\boldsymbol{\Theta} = \mathbf{R}^{T}\mathbf{R}$ also implies that $\textnormal{Ker}(\boldsymbol{\Theta}) = \textnormal{Ker}(\mathbf{R})$.
\end{proof}

 \begin{remark}\label{rem:numberfillingfrac} In terms of the corresponding tiling models, we have $H$ filling fractions subjected to $m$ affine constraints specified by $\mathbf{R}$. Therefore, the filling fractions vary in an affine space of dimension $b_1(\mathsf{G}) = \textnormal{corank}(\boldsymbol{\Theta})$ in the bipartite case, or $b_1(\mathsf{G}) = \textnormal{corank}(\boldsymbol{\Theta}) + 1$ in the non-bipartite case. In the bipartite case, to be in the setting of Chapter~\ref{Chapter_Setup_and_Examples} we need to select a linearly independent subset of $\mathfrak{e} = m - 1$ linear forms among $(\mathfrak{r}_1,\ldots,\mathfrak{r}_m)$.
 \end{remark}

Compared to the gluing graph $\mathsf{G}$, the matrix $\boldsymbol{\Theta}$ remembers the number $2\theta_{g,h}$ of vertices shared by the edges $g$ and $h$ in $\mathsf{G}$, \textit{i.e.} the number of trapezoids containing both the $g$-th and the $h$-th segment, with the convention $\theta_{g,h} = 1$ for $g = h$. We may ask if the graph can be reconstructed if we are only given the matrix $\boldsymbol{\Theta}$.

\begin{definition}
\label{def:connectoid}A \emph{connectoid} is a symmetric $H \times H$ positive semi-definite matrix $\boldsymbol{\Theta}$ with entries in $\{0,\frac{1}{2},1\}$ and having the following connectedness property: for any $h,h' \in [H]$ there exists a finite sequence $h_0,\ldots,h_s$ with $h_0 = h$, $h_d = h'$ and for any $c \in [d]$ we have $\theta_{h_{c - 1},h_c} \neq 0$. A connectoid is \emph{realizable} if it comes from a graph $\mathsf{G}$ as in \eqref{ThetaGeq}; we can also ask for bipartite or non-bipartite realizability.
\end{definition}

As the proof of Proposition~\ref{lem:H1G} only depends on the graph $\mathsf{G}$, to realize a connectoid in a bipartite (respectively, non-bipartite) way one should look for a graph with $m = \textnormal{rank}(\boldsymbol{\Theta}) + 1$ (respectively $m = \textnormal{rank}(\boldsymbol{\Theta})$) vertices. For $H \leq 3$ all connectoids are realizable, \textit{cf.} Figure~\ref{Fig:Domains23}. For $H = 4$, there are 70 connectoids: 60 of them are realizable, the 10 non-realizable ones have full rank and up to simultaneous row- and column-permutations they take two possible forms
\begin{equation}
\label{exnonrel}
\left(\begin{array}{cccc} 1 & 0 & 0 & \frac{1}{2} \\[0.6ex] 0 & 1 & 0 & \frac{1}{2} \\[0.6ex] 0 & 0 & 1 & \frac{1}{2} \\[0.6ex] \frac{1}{2} & \frac{1}{2} & \frac{1}{2} & 1 \end{array}\right)\qquad \textnormal{or}\qquad \left(\begin{array}{cccc} 1 & 0 & \frac{1}{2} & \frac{1}{2} \\[0.6ex] 0 & 1 & \frac{1}{2} & \frac{1}{2} \\[0.6ex] \frac{1}{2} & \frac{1}{2} & 1 & \frac{1}{2} \\[0.6ex] \frac{1}{2} & \frac{1}{2} & \frac{1}{2} & 1 \end{array}\right).
\end{equation}
The following realizability criterion was suggested by Tibor Szab\'o. We denote
\[
C_h = \{g \in [H] \,\,|\,\,\theta_{g,h} \neq 0\}.
\]

\begin{proposition}
\label{Prop:realiza}A connectoid $\boldsymbol{\Theta}$ is realizable if and only it is the matrix full of $1$s, or the conditions (1) and (2) below are simultaneously satisfied.
\begin{itemize}
\item[(1)] For any $g,h,k \in [H]$, having $\theta_{g,h} = 1$ implies $\theta_{g,k} = \theta_{h,k}$. In particular, the relation $\theta_{g,h} = 1$ defines an equivalence relation on $[H]$, and we denote $O_h$ the equivalence class of $h$.
\item[(2)] For each $h \in [H]$, the set of subsets $\mathfrak{C}_h := \big\{C_h \cap C_g \,\,|\,\,g \in C_h \setminus O_h \big\}$ has cardinality $1$ or $2$. If it has two distinct elements $C_{h}',C_h''$, then $C_h' \cap C_h''$ is the equivalence class containing $h$.
\end{itemize}
In this case, denoting $\mathfrak{C}_h^+ = \mathfrak{C}_h$ if the latter set has cardinality $2$, and $\mathfrak{C}_h^+ = \mathfrak{C}_h \sqcup \{O_h\}$ otherwise, and $\mathfrak{C}^+ = \bigcup_{h = 1}^H \mathfrak{C}_h^+$, we have
\begin{equation}
\label{thetareali}\boldsymbol{\Theta} = \frac{1}{2} \sum_{C \in \mathfrak{C}^+} \textnormal{\textbf{1}}_{C \times C}.
\end{equation}
\end{proposition}
\begin{proof}
Assume that $\boldsymbol{\Theta}$ is realizable, and look at a particular realization. For any $g,h \in [H]$, having $\theta_{g,h} = 1$ means that the $g$-th and $h$-th segment belongs to the same two trapezoids, so the $k$-th segment has the same number of trapezoids in common with $h$ and with $g$, proving the implication in (1). If $\theta_{g,h} < 1$, we can find a trapezoid containing the $h$-th segment. If there is a trapezoid containing the segments in $O_h$ and no other segment, then the assumption of connectedness and $\boldsymbol{\Theta} \neq \textnormal{\textbf{1}}_{[H] \times [H]}$ imply that $\mathfrak{C}_h$ has a unique element, namely the group of segments in the other trapezoid containing the segments in $O_h$ and at least another segment. Otherwise, $\mathfrak{C}_h$ has exactly two elements, corresponding to the group of segments in the two trapezoids containing the $h$-th segment (automatically: all segments in $O_h$). The intersection of these two groups is exactly $O_h$, proving (2).

Conversely, assume (1) and (2) are satisfied, and let $g,h \in [H]$. The $(g,h)$-th entry of the right-hand side of \eqref{thetareali} is half the number $c_{g,h}$ of elements of $\mathfrak{C}^+$ to which $g$ and $h$ belong jointly. We discuss separately all possible values of $\theta_{g,h}$.

We first consider the case $\theta_{g,h} = 0$. Then $g$ and $h$ cannot belong jointly to a set $C_{g} \cap C_{k}$ or $C_h \cap C_{k}$ for some $k \in [H]$. Neither can they belong to the same equivalence class defined by the equivalence relation in (1). If they belonged jointly to $C_k \cap C_{k'}$ for some $k,k' \in [H] \setminus \{g,h\}$ with $k' \notin O_k$, then $\mathfrak{C}_k$ would contain the three distinct sets $C_k \cap C_{k'}$ (contains $g$ and $h$), $C_k \cap C_g$ (contains $g$ but not $h$), $C_k \cap C_h$ (contains $h$ but not $g$), which contradicts (2). Hence $c_{g,h} = 0 = 2\theta_{g,h}$.

Next, we examine the case $\theta_{g,h} = \frac{1}{2}$. Then $g$ and $h$ are not in the same equivalence class, and $g$ and $h$ belong jointly to $C_g \cap C_h$, and the latter appears as an element of $\mathfrak{C}_h^+$ and of $\mathfrak{C}_g^+$. Since the intersection of the two elements of $\mathfrak{C}_h^+$ is $O_h$, the second element of $\mathfrak{C}_h^+$ contains $h$ but not $g$. For the same reason, the second element of $\mathfrak{C}_g^+$ contains $g$ but not $h$. If $g$ and $h$ belonged jointly to a set $C_k \cap C_{k'}$ with $k,k' \in [H] \setminus \{g,h\}$ and $k' \notin O_k$, then using the previous observation we see that both $C_h \cap C_k$ and $C_h \cap C_{k'}$ coincide with the element of $\mathfrak{C}_h^+$ containing $g$, and likewise exchanging the role of $g$ and $h$. Hence
\[
C_g \cap C_k = C_g \cap C_{k'} = C_h \cap C_{k} = C_h \cap C_{k'} = C_h \cap C_g.
\]
If $C_k \cap C_{k'}$ was distinct from this set, we would find $l \in (C_k \cap C_{k'}) \setminus (C_g \cap C_h)$, but then $\mathfrak{C}_k$ would contain the three distinct sets $C_{k} \cap C_{k'}$ (contains $g,h,k,k',l$), $C_k \cap C_g$ (contains $g,h,k,k'$ but not $l$), $C_k \cap C_l$ (contains $k,l$ but not $g$ nor $h$) and this is ruled out by (2). This shows that $C_g \cap C_h = C_k \cap C_{k'}$ and we conclude that the unique element of $\mathfrak{C}^+$ to which $g$ and $h$ jointly belongs is $C_g \cap C_h$. In particular $c_{g,h} = 1 = 2\theta_{g,h}$.

Last, we look at the case $\theta_{g,h} = 1$. Then $g$ and $h$ belongs to the same equivalence class, so $C_g = C_h$, and $g$ and $h$ belong jointly to the two distinct elements of $\mathfrak{C}_h^+$, and to no other distinct element of $\mathfrak{C}^+$. Hence $c_{g,h} = 2 = 2\theta_{g,h}$.
\end{proof}

\begin{figure}[h!]
\begin{center}
\includegraphics[width=0.4\textwidth]{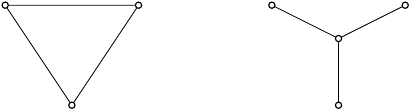}
\caption{\label{Fig:Startri}The triangle and the star have same interaction matrix $\boldsymbol{\Theta}$.}
\end{center}
\end{figure}

If a connectoid is realizable, Proposition~\ref{Prop:realiza} gives a canonical realization but it does not exclude the existence of other realizations. For instance, the connectoid
\[
\boldsymbol{\Theta} = \left(\begin{array}{ccc} 1 & \frac{1}{2} & \frac{1}{2} \\[0.6ex] \frac{1}{2} & 1 & \frac{1}{2} \\[0.6ex] \frac{1}{2} & \frac{1}{2} & 1 \end{array}\right)
\]
has one bipartite realization (the 3-branched star) and one non-bipartite realization (the triangle), shown in Figure~\ref{Fig:Startri}. We also see that some connectoids only have a non-bipartite realization (\textit{e.g.} the connectoid of the pentagon), and some only have a bipartite realization (\textit{e.g.} the $H = 2$ connectoids). In fact, we have the following uniqueness results, which we establish by elementary methods.

\begin{proposition}
\label{prop:uniquerelconn}If a connectoid is realized by a bipartite graph, it is so in a unique way. If a connectoid is realized by a non-bipartite graph, it is so in a unique way.
\end{proposition}

In Section~\ref{sec:extraconn} we will obtain more complete information: except for $m = 3$ due to the star-triangle equivalence, connectoids have at most one realization (given by Proposition~\ref{Prop:realiza}).

\begin{proof}
\textsc{Bipartite case.} Let $\boldsymbol{\Theta}$ be the connectoid coming in at least one way from a bipartite connected graph. Without lack of generality we can assume that the graph has at least two vertices. We construct the (undirected) connected graph $\mathsf{G}^{\vee}$ with vertex set $[H]$, having $2\theta_{g,h}$ edges between $g$ and $h$ for any $g \neq h$. This graph encodes faithfully the matrix $\boldsymbol{\Theta}$. We would like to assign signs $\pm 1$ to edges in $\mathsf{G}^{\vee}$, in such a way that the subgraph $\mathsf{G}_{+}^{\vee}$ (respectively $\mathsf{G}_-^{\vee}$) with vertex set $[H]$ and keeping only positive (respectively, negative) edges is a cluster graph, \textit{i.e.} a disjoint union of complete graphs with simple edges. The assumption that $\boldsymbol{\Theta}$ has a bipartite realization means there exists at least one way to do so, and we are going to examine the necessities of $\pm 1$ assignments for this. If $\epsilon \in \{\pm 1\}$, we write $\overline{\epsilon}$ for the flipped assignment. First, in each double edge one of them must be assigned $+1$ and the other $-1$. Then, if there are no simple edges, we are finished. Otherwise, let $e = e_0$ be the simple edge which is minimal for the lexicographic order on pairs of vertices. We have the freedom to assign $\epsilon_{e_{\min}} = +1$ or $\epsilon_{e_{\min}} = -1$ to this edge. With the data of $e = \{h_0,h_1\}$ and its assignment $\epsilon_{\{h_0,h_1\}}$, we enter the following loop in the algorithm. Look at any other simple edge of the form $\{h_0,h_2\}$ and not already visited. It is not possible that $\{h_1,h_2\}$ is a double edge; if $\{h_1,h_2\}$ is a simple edge, then we must declare $\epsilon_{\{h_0,h_2\}} := \epsilon_e$ and $\epsilon_{\{h_1,h_2\}} := \epsilon_e$. If we already visited $\{h_1,h_2\}$, this latter value should agree with the previously assigned one, because we assumed that the algorithm is successful in at least one way. If $\{h_1,h_2\}$ is not an edge, then we must declare $\epsilon_{\{h_0,h_2\}} := \overline{\epsilon_e}$. Then, for any simple edge $\{h_1,h_3\}$ that we have not already visited, $\{h_0,h_3\}$ is not an edge, so we must declare $\epsilon_{\{h_1,h_3\}} = \overline{\epsilon_e}$. If there are no simple edges adjacent to those already visited, we exit the loop and terminate the algorithm. Otherwise, we replace $e$ with the minimal (for the lexicographic order) edge that has not been visited yet but is adjacent to the edges we already visited, and re-enter the loop.

Since the graph is connected, the algorithm always terminates after all simple edges have been visited and yield the desired assignment where the only freedom was the choice of $\epsilon_{e_{\min}}$. The matrix $\boldsymbol{\Theta}^{+}$ (respectively $\boldsymbol{\Theta}^-$) is then the adjacency matrix of the graph obtained by keeping only positive (respectively, negative) edges, and the roles of $\boldsymbol{\Theta}^{+}$ and $\boldsymbol{\Theta}^-$ are exchanged if we flip the value of $\epsilon_{e_{\min}}$.

\medskip

\noindent \textsc{Non-bipartite case.} We show by induction on the size $H$ that connectoids admit at most one non-bipartite realization. For $H = 1$ there is nothing to prove as the only connectoid $\boldsymbol{\Theta} = (1)$ is only realized in a bipartite way. Take $H \geq 2$ and assume that all connectoids of size $< H$ admitting a non-bipartite realization admit a unique such realization. Let $\boldsymbol{\Theta}$ be a connectoid of size $H$ realizable by a non-bipartite graph. We examine the matrix $\boldsymbol{\Theta}^{[1]} = (\theta_{g,h})_{g,h = 2}^{H}$ of size $H - 1$ and the vector $(\theta_{1,h})_{h = 2}^{H}$; together they determine the matrix $\boldsymbol{\Theta}$. We write
\[
\mathcal{H}' = \big\{h \in \{2,\ldots,H\}\,\,|\,\,\theta_{h,1} \neq 0\big\}.
\]

Suppose that $\boldsymbol{\Theta}^{[1]}$ has the connectedness property (\textit{cf.} Definition~\ref{def:connectoid}). Then, it is a connectoid and $\mathcal{H}' \neq \emptyset$. Removing the edge labeled $1$ (and the vertex attached to it in case it was univalent) in a realization of $\boldsymbol{\Theta}$ gives a connected graph realizing $\boldsymbol{\Theta}^{[1]}$. We may not obtain all possible realizations of $\boldsymbol{\Theta}^{[1]}$ in this way, but any non-bipartite realization of $\boldsymbol{\Theta}$ is obtained from \emph{some} realization $\mathsf{G}^{[1]}$ of $\boldsymbol{\Theta}^{[1]}$ by adding an edge $e_1$ labeled $1$. There are three possibilities for the latter, which possibly overlap depending on the choice of realization.
\begin{itemize}
\item[(a1)] $\mathsf{G}^{[1]}$ is non-bipartite and $e_1$ connects an existing vertex $i$ ot $\mathsf{G}^{[1]}$ to a new vertex $0$.
\item[(a2)] $\mathsf{G}^{[1]}$ is non-bipartite and $e_1$ connects two vertices $i_1$ and $i_2$ of $\mathsf{G}^{[1]}$.
\item[(b)] $\mathsf{G}^{[1]}$ is bipartite and $e_1$ connects a negative vertex $i_1$ to a positive vertex $i_2$ of $\mathsf{G}^{[1]}$.
\end{itemize}
If case (a1) occurs, we have $\theta_{h,1} = \frac{1}{2} \delta_{h,\mathcal{H}'}$ for any $h \in \{2,\ldots,H\}$ and $\theta_{g,h} \geq \frac{1}{2}$ for any $g,h \in \mathcal{H}'$. By Proposition~\ref{lem:H1G} we must have $\textnormal{rank}(\boldsymbol{\Theta}) = \textnormal{rank}(\boldsymbol{\Theta}^{[1]}) + 1$ in cases (a1) and (b), and $\textnormal{rank}(\boldsymbol{\Theta}) = \textnormal{rank}(\boldsymbol{\Theta}^{[1]})$ in case (a2). In case (a2) or (b), we have $\mathcal{H}' = \mathcal{H}^{[1]}_{i_1} \cup \mathcal{H}^{[1]}_{i_2}$ where $\mathcal{H}^{[1]}_i$ labels the edges adjacent to $i$ in $\mathsf{G}^{[1]}$. Note that the two subsets are non-empty and the union may not be disjoint.

If $\textnormal{rank}(\boldsymbol{\Theta}) = \textnormal{rank}(\boldsymbol{\Theta}^{[1]})$, (a2) must occur and (a1) and (b) do not occur for any realization. Then, by the hypothesis assumption $\boldsymbol{\Theta}^{[1]}$ is realized by a non-bipartite graph $\mathsf{G}^{[1]}$ which is uniquely determined by $\boldsymbol{\Theta}^{[1]}$, hence by $\boldsymbol{\Theta}$, and obviously the set $\mathcal{H}'$ is determined by $\boldsymbol{\Theta}$. We claim that the vertices $i_1$ and $i_2$ are also determined uniquely by $\boldsymbol{\Theta}$. Indeed, if there is $h_0 \in \mathcal{H}'$ such that $\theta_{1,h_0} = 1$, $i_1$ and $i_2$ must be the vertices in $\mathsf{G}^{[1]}$ sharing the $h_0$-th segment. If $\theta_{1,h} = \frac{1}{2}$ for any $h \in \mathcal{H}'$, there exist $h_1,h_2 \in \mathcal{H}$ such that the $h_1$-th edge and the $h_2$-th edge are adjacent to distinct vertices in $\mathsf{G}^{[1]}$. Then these two vertices must be (up to order) $i_1$ and $i_2$, because $\mathcal{H}' = \mathcal{H}^{[1]}_{i_1} \cup \mathcal{H}^{[1]}_{i_2}$. This justifies the claim and shows that $\boldsymbol{\Theta}$ has a unique non-bipartite realization.

If $\textnormal{rank}(\boldsymbol{\Theta}) = \textnormal{rank}(\boldsymbol{\Theta}^{[1]}) + 1$, (a2) does not occur but (a1) or (b) could occur, without exclusivity. If $\boldsymbol{\Theta}^{[1]}$ has no non-bipartite realization, then (b) must occur for the unique (by the first part of the proof) bipartite realization of $\boldsymbol{\Theta}^{[1]}$. Then, the argument used to treat (a2) applies and shows that the vertices $i_1,i_2$ are determined by $\boldsymbol{\Theta}$, so that $\boldsymbol{\Theta}$ admits a unique non-bipartite realization. If $\boldsymbol{\Theta}^{[1]}$ has no bipartite realization, then (a1) must occur for the unique (by the induction hypothesis) non-bipartite realization of $\boldsymbol{\Theta}^{[1]}$, and the vertex $i$ must be the vertex in $\mathsf{G}^{[1]}$ whose set of adjacent edges in $\mathcal{H}'$, which is determined by $\boldsymbol{\Theta}$. If $\boldsymbol{\Theta}^{[1]}$ admits both a (unique) non-bipartite realization $\mathsf{G}^{[1]}_{\textnormal{a1}}$ and a (unique) bipartite realization $\textsf{G}^{[1]}_{\textnormal{b}}$, and (a1) is realized for the first while (b) is realized for the second. The previous arguments show that the vertex $i$ in (a1), or the vertices $i_1,i_2$ in (b) are determined by $\boldsymbol{\Theta}$. But for (a1) we must have $\theta_{1,h} = \frac{1}{2}$ and $\theta_{g,h} \geq \frac{1}{2}$ for any $g,h \in \mathcal{H}'$. The first equality implies in (b) that $i_1$ and $i_2$ are not related by an edge in $\mathsf{G}^{[1]}_{\textnormal{b}}$. Since $\mathcal{H}^{[1]}_{i_1}$ and $\mathcal{H}^{[1]}_{i_2}$ cannot be empty, this means that we can find $h_1,h_2 \in \mathcal{H}'$ such that $\theta_{h_1,h_2} = 0$, contradicting the second equality. So, (a1) and (b) cannot occur simultaneously, and we exhausted all possibilities. This proves that if $\boldsymbol{\Theta}^{[1]}$ satisfies the connectedness assumption, the non-bipartite realization of $\boldsymbol{\Theta}$ is unique.

Now suppose that $\boldsymbol{\Theta}^{[1]}$ does not have the connectedness property. Then, for any non-bipartite realization $\mathsf{G}$ of $\boldsymbol{\Theta}$, if we remove the edge labeled $1$ (and the adjacent vertex if it is univalent) we obtain a graph which is a disjoint union of two non-bipartite connected graphs $\mathsf{G}_1,\mathsf{G}_2$, whose edge label sets $C_1$ and $C_2$ form a partition of $\{2,\ldots,H\}$. Call $\boldsymbol{\Theta}_j$ the connectoid of $\mathsf{G}_j$, and define $\boldsymbol{\Theta}_{j}^{[1]}$ as the matrix whose rows and columns are indexed by $\{2,\ldots,H\}$, whose restriction to the rows and columns in $C_j$ coincides with $\boldsymbol{\Theta}_j$, and whose other entries are $0$. Then
\[
\boldsymbol{\Theta}^{[1]} = \boldsymbol{\Theta}_1^{[1]} + \boldsymbol{\Theta}_2^{[1]}.
\]
As the two summands have disjoint support (simultaneously in row and columns), the partition $C_1 \sqcup C_2 = \{2,\ldots,H\}$ and thus the connectoids $\boldsymbol{\Theta}_j$ are determined by $\boldsymbol{\Theta}^{[1]}$ and do not depend on the choice of realizations. By the induction hypothesis, $\boldsymbol{\Theta}_1$ and $\boldsymbol{\Theta}_2$ admit a unique non-bipartite realization, so $\mathsf{G}_1,\mathsf{G}_2$ (up to exchange) are determined by $\boldsymbol{\Theta}$. Besides, the edge labeled $1$ in $\mathsf{G}$ must connect the vertex of $\mathsf{G}_1$ adjacent to the edges labeled by $\{h \in C_1\,\,|\,\,\theta_{1,h} = \frac{1}{2}\}$ to the vertex of $\mathsf{G}_2$ adjacent to the edges labeled by $\{h \in C_2 \,\,|\,\,\theta_{1,h} = \frac{1}{2}\}$, showing that these vertices are uniquely determined by $\boldsymbol{\Theta}$ and therefore establishing the uniqueness of the non-bipartite realization of $\boldsymbol{\Theta}$.

We now have shown the induction step is valid in all situations, completing the proof.
\end{proof}

\subsection{Main result: spectral curves, fundamental solutions and Green functions}

\label{sec:C13tiling}

In this section we give the main results about spectral curves for the master Riemann--Hilbert problem of matrices $\boldsymbol{\Theta}$ which are realizable connectoids. The group $\mathfrak{G}$ admits in this case a finite and geometrically meaningful orbit. In the bipartite case, the spectral curve is obtained from the glued domain by replacing each trapezoid by a sheet. In the non-bipartite case, the spectral curve is the one of the bipartite covering of the glued domain, equipped with the flip involution of the covering. In particular it applies to the interaction matrices of tiling models. We also identify the fundamental solutions which relate to the leading covariance by Corollary~\ref{Corollary_CLT}, show that they relate to suitable Green functions, and identify the combinations of first-kind functions appearing (by Proposition~\ref{Proposition_differentiability_filling_fraction}) in the first variation of the Stieltjes transform of the equilibrium measure with respect to filling fractions.

If $i \in [m]$ and $h \in \mathcal{H}_i$, we denote $\mathcal{T}_{i_h'}$ the other trapezoid to which the $h$-th segment belong, \textit{i.e.} $i \neq i_h'$ and $h \in \mathcal{H}_i \cap \mathcal{H}_{i_h'}$.

\begin{theorem}
\label{thm:Omegav} Assume $\boldsymbol{\Theta}$ is realized by a bipartite gluing graph $\mathsf{G}$, with bipartition $\epsilon : [m] \rightarrow \{\pm 1\}$. Then the vector $\boldsymbol{v} = \epsilon_1\mathbf{1}_{\mathcal{H}_1}$ has $\mathfrak{G}$-orbit $\{\epsilon_i\mathbf{1}_{\mathcal{H}_i}\,\,|\,\,i\in [m]\}$, this orbit is good (\textit{cf.} Definition~\ref{def:goodor}) and has minimal size $m = \textnormal{rank}(\boldsymbol{\Theta}) + 1$. For any $i \in [m]$ and $h \in [H]$ we have
\[
T^{(h)}(\epsilon_i \mathbf{1}_{\mathcal{H}_i}) = \left\{\begin{array}{ccc}  \epsilon_i \mathbf{1}_{\mathcal{H}_i} && \textnormal{if}\,\,h \notin \mathcal{H}_i, \\[0.6ex] \epsilon_{i'_h}\mathbf{1}_{\mathcal{H}_{i'_h}} && \textnormal{if}\,\,h \in \mathcal{H}_i. \end{array}\right.
\]
The genus of the spectral curve $\Sigma_{\boldsymbol{v}}$ is $\mathsf{g}_{\boldsymbol{v}} = \textnormal{corank}(\boldsymbol{\Theta})$ and the group $\hat{\mathfrak{G}}$ is the Weyl group of type $A_{m - 1}$ (\textit{i.e.} the symmetric group in $m$ elements).

The bidifferential $\mathcal{B}_{\widetilde{\boldsymbol{v}}}$ of Theorem~\ref{thm:Bcont} is defined on $\Sigma^{[2],\circ}_{\boldsymbol{v}} = \Sigma_{\boldsymbol{v}} \times \Sigma_{\boldsymbol{v}}$ and can simply be denoted $\mathcal{B}_{\boldsymbol{v}}$. Its matrix of biresidues is \textbf{\foreignlanguage{russian}{B}} = $\textnormal{\textbf{Id}} - \boldsymbol{1}$, and the bidifferential itself is equal to
\begin{equation}
\label{Omev111}\mathcal{B}_{\boldsymbol{v}} = \mathcal{B}^{\Sigma_{\boldsymbol{v}},\mathcal{L}} - \frac{\dd z_1\dd z_2}{(z_1 - z_2)^2},
\end{equation}
where $(z_1,z_2) = Z_{\widetilde{\boldsymbol{v}}}$, the Lagrangian $\mathcal{L} \subset H_1(\Sigma_{\boldsymbol{v}},\amsmathbb{Z})$ is spanned by cycles around the components of the gluing locus (\textit{cf.} Definition~\ref{def:spcurv1}), and $\mathcal{B}^{\Sigma_{\boldsymbol{v}},\mathcal{L}}$ is the fundamental bidifferential of Section~\ref{sec:fund02}.

From the analytic continuation of first-kind functions in Proposition~\ref{Lem:Holo1form}, we get a linear isomorphism
\[
\begin{array}{ccc} \textnormal{Ker}(\boldsymbol{\Theta}) & \longrightarrow & H^1(\Sigma_{\boldsymbol{v}},\amsmathbb{C}) \\
\boldsymbol{v}' & \longmapsto & \mathfrak{u}_{\boldsymbol{v};\boldsymbol{v}'} \end{array}
\]
\end{theorem}

\begin{figure}[h!]
\begin{center}
\includegraphics[width=0.57\textwidth]{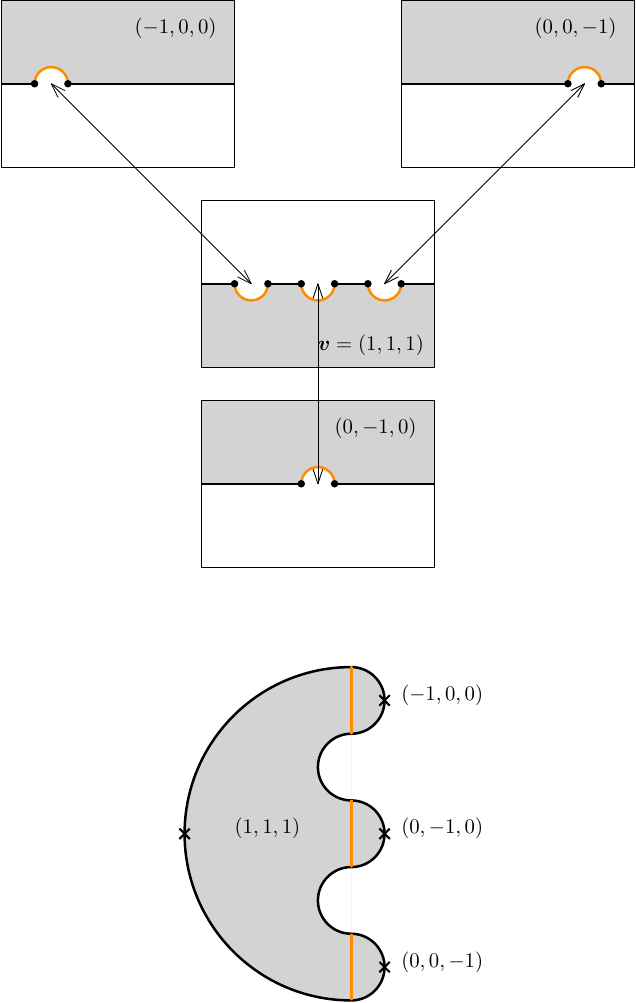}
\caption{\label{Fig:Half5} The half-spectral curve corresponding to the example of Figure~\ref{Fig:Half4}. The boundary is indicated by a fat black line. The orange curves correspond to the gluing locus: these are curves in the interior of the half-spectral curve joining two boundary points. Top: construction by gluing half-planes. Bottom: the corresponding glued domain, rotated for convenience (compare with Figure~\ref{ExEshaped}). The crosses indicate the points at infinity.}
\end{center}
\end{figure}

\begin{figure}[h!]
\begin{center}
\includegraphics[width=0.65\textwidth]{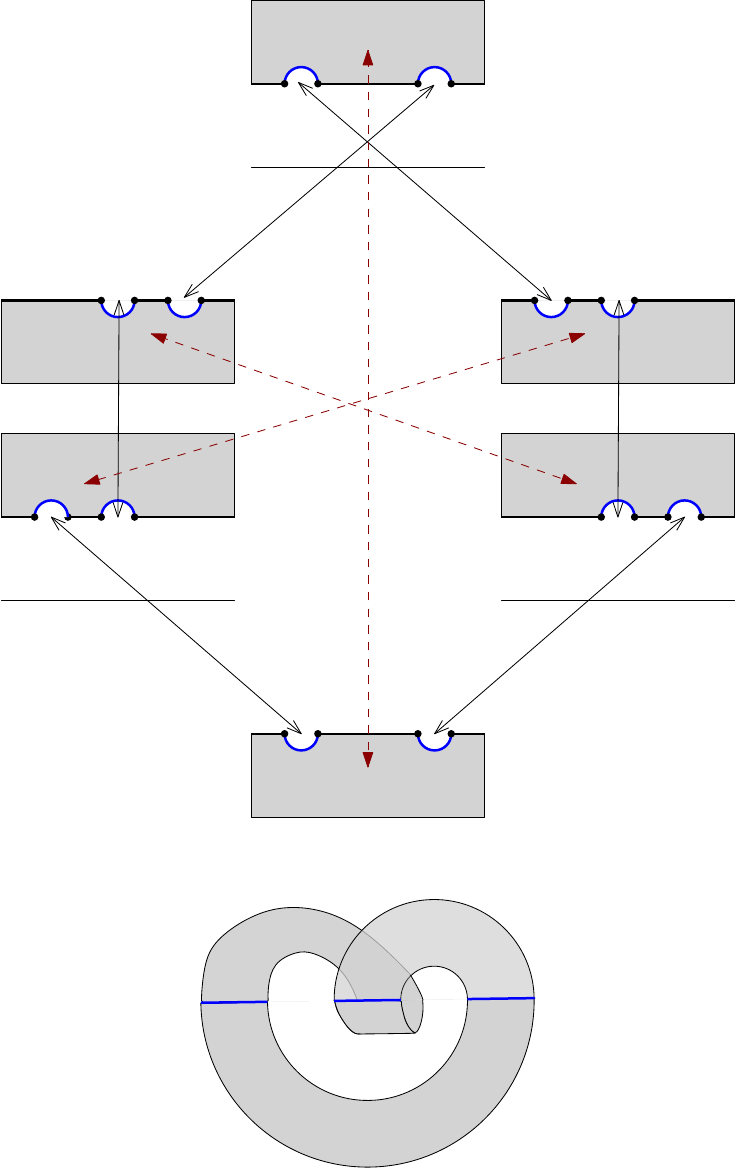}
\caption{\label{Fig:Half7} Bottom: a non-bipartite tiled domain (rotated for convenience) corresponding to $H = 3$, $\theta_{h,h} = 1$ and $\theta_{h,g} = \frac{1}{2}$ for any distinct $g,h \in [3]$. This $\boldsymbol{\Theta}$ is like in Figure~\ref{Fig:Half4} and \ref{Fig:Half5}, but we use a different seed vector $\boldsymbol{v} = (1,1,0)$, which has the (bipartite) orbit of size $6$ consisting of $\pm (1,1,0),\pm (1,0,1),\pm (0,1,1)$. Top: the half-spectral curve $\Sigma_{\boldsymbol{v}}^{\textnormal{half}}$ with boundary indicated by a fat black line and involution $\varsigma^*$ by dashed arrows.}
\end{center}
\end{figure}

We see that the glued domain is homeomorphic to the half-spectral curve $\Sigma_{\boldsymbol{v}}^{\textnormal{half}}$ of Definition~\ref{def:halfspcurvebip}: the $i$-th trapezoid is replaced with a hemisphere (upper or lower depending on the sign $\epsilon_i$) and these are glued in the same way specified by $\boldsymbol{\Theta}$ (\textit{cf.} Figure~\ref{Fig:Half5}). The structure of $\mathcal{B}_{\boldsymbol{v}}$ allows applying Corollary~\ref{Greenwon} and relating it to the Green function of $\Sigma_{\boldsymbol{v}}^{\textnormal{half}}$.

\begin{corollary}
\label{cor:Greennormal} Consider the situation of Theorem~\ref{thm:Omegav}. Identify for each $h \in [H]$ the interval $(\alpha_h,\beta_h)$ in the glued domain with its image in the interior of $\Sigma_{\boldsymbol{v}}^{\textnormal{half}}$. Let $(f_h(z))_{h = 1}^{H}$ be an $H$-tuple of functions such that $f_h(z)$ is a holomorphic function of $z$ in a complex neighborhood of $[\alpha_h,\beta_h]$, taking real values on the real axis, and call $\tilde{f}_{h}(p) = f_h(Z_{\boldsymbol{v}}(p))$ the function in the corresponding domain of $\Sigma_{\boldsymbol{v}}^{\textnormal{half}}$. Then, for any $h_1,h_2 \in [H]$
\[
\oint_{\gamma_{h_1}} \oint_{\gamma_{h_2}} \mathcal{F}_{h_1,h_2}(z_1,z_2) f_{h_1}(z_1) f_{h_2}(z_2) \frac{\dd z_1 \dd z_2}{(2\ii\pi)^2} = \frac{1}{\pi} \int_{\alpha_{h_1}}^{\beta_{h_1}}  \int_{\alpha_{h_2}}^{\beta_{h_2}} \textnormal{Gr}(p_1,p_2) \dd \tilde{f}_{h_1}(p_1) \dd \tilde{f}_{h_2}(p_2)
\]
involving the Green function on $\Sigma_{\boldsymbol{v}}^{\textnormal{half}}$.
\end{corollary}
\begin{proof} This is a direct consequence of Proposition~\ref{prop:Greengen} (H5) and Corollary~\ref{Greenwon}.
\end{proof}

\begin{remark} \label{rem:match1113} With the particular $\boldsymbol{v}$ chosen in Theorem~\ref{thm:Omegav}, $\Sigma^{\textnormal{half}}_{\boldsymbol{v}}$ matches the construction of $\Sigma^{\textnormal{half}}$ in Section~\ref{Section_complex_structure} by declaring that left trapezoids have $\epsilon = 1$ and right trapezoids have $\epsilon = -1$. The function of Definition~\ref{Definition_Covariance_tilings_St} is related to the bidifferential constructed here in the following way
\begin{equation}
\label{Fqqq12}\mathcal{F}^{(i_1,i_2)}(z_1,z_2) = \big\langle \mathbf{1}_{\mathcal{H}_{i_1}} \otimes\mathbf{1}_{\mathcal{H}_{i_2}} \cdot \boldsymbol{\mathcal{F}}(z_1,z_2) \big\rangle.
\end{equation}
Multiplying by $(-1)^{\textnormal{parity}(i_1,i_2)}$ and $\dd z_1\dd z_2$ and adding the biresidue terms is exactly the way $\mathcal{B}_{\boldsymbol{v}}$ is constructed in Theorem~\ref{thm:Omegav}. After double integration it led us to Proposition~\ref{Proposition_Green_function_explicit}, which was an essential ingredient in the proof of Theorem~\ref{Theorem_GFF_general}.
\end{remark}

\begin{theorem}
\label{thm:spcurvenonbip}
Assume $\boldsymbol{\Theta}$ is realized by a non-bipartite gluing graph $\mathsf{G}$. Then the vector $\boldsymbol{v} = \mathbf{1}_{\mathcal{H}_1}$ has $\mathfrak{G}$-orbit $\{\tau\,\mathbf{1}_{\mathcal{H}_i}\,\,|\,\,\tau \in \{\pm 1\}\,\,i \in [m]\}$, this orbit is good and has minimal size $2m = 2\,\textnormal{rank}(\boldsymbol{\Theta})$. For any $i \in [m]$ and $h \in \mathcal{H}_i$, we have
\[
T^{(h)}(\tau \mathbf{1}_{\mathcal{H}_i}) = \left\{\begin{array}{ccc} \tau \mathbf{1}_{\mathcal{H}_i} && \textnormal{if}\,\,h \notin \mathcal{H}_i \\[0.6ex] - \tau \mathbf{1}_{\mathcal{H}_{i'_h}} && \textnormal{if}\,\, h \in \mathcal{H}_i. \end{array}\right.
\]
The spectral curve $\Sigma_{\boldsymbol{v}}$ has genus $\mathsf{g}_{\boldsymbol{v}} = 2\,\textnormal{corank}(\boldsymbol{\Theta}) + 1$, and is equipped with a holomorphic involution \label{index:invs}$\varsigma : \Sigma_{\boldsymbol{v}} \rightarrow \Sigma_{\boldsymbol{v}}$ sending $p$ in the $\boldsymbol{w}$-th sheet to $p'$ in the $(-\boldsymbol{w})$-th sheet such that $Z_{\boldsymbol{v}}(p) = Z_{\boldsymbol{v}}(p')$. The group $\hat{\mathfrak{G}}$ is the Weyl group of type $D_{m}$.

The bidifferential $\mathcal{B}_{\widetilde{\boldsymbol{v}}}$ of Theorem~\ref{thm:Bcont} is defined on the spectral surface $\Sigma_{\boldsymbol{v} \otimes\boldsymbol{v}}^{[2],\circ} = \Sigma_{\boldsymbol{v}} \times \Sigma_{\boldsymbol{v}} / (\varsigma,\varsigma)$. Its pullback $\mathcal{B}_{\boldsymbol{v}}$ to $\Sigma_{\boldsymbol{v}} \times \Sigma_{\boldsymbol{v}}$ has biresidues \foreignlanguage{russian}{B}${}_{\boldsymbol{w}_1,\boldsymbol{w}_2} = \delta_{\boldsymbol{w}_1,\boldsymbol{w}_2} - 1 - \delta_{\boldsymbol{w}_1 + \boldsymbol{w}_2, \boldsymbol{0}}$, and is equal to
\[
\mathcal{B}_{\boldsymbol{v}} = \mathcal{B}^{\Sigma_{\boldsymbol{v}},\mathcal{L}} - (\textnormal{id},\varsigma)^*\mathcal{B}^{\Sigma,\mathcal{L}} - \frac{\dd z_1 \dd z_2}{(z_1 - z_2)^2},
\]
and we have $(\textnormal{id},\varsigma)^* \mathcal{B}^{\Sigma,\mathcal{L}} = (\varsigma, \textnormal{id})^* \mathcal{B}^{\Sigma,\mathcal{L}}$. Here $\mathcal{L}$ is the same Lagrangian as in Theorem~\ref{thm:Omegav} and $\mathcal{B}^{\Sigma,\mathcal{L}}$ is the fundamental bidifferential of Section~\ref{sec:fund02}.

From the analytic continuation of first-kind functions in Proposition~\ref{Lem:Holo1form}, we get an linear isomorphism
\[
\begin{array}{ccc} \textnormal{Ker}(\boldsymbol{\Theta}) & \longrightarrow & H^1(\Sigma_{\boldsymbol{v}},\amsmathbb{C})^- \\
\boldsymbol{v}' & \longmapsto & \mathfrak{u}_{\boldsymbol{v};\boldsymbol{v}'} \end{array}
\]
where $H^1(\Sigma_{\boldsymbol{v}},\amsmathbb{C})^-$ is the space of holomorphic $1$-forms that are anti-invariant under the involution.
\end{theorem}

In this non-bipartite case, the $\mathfrak{G}$-orbit of $\boldsymbol{v} = \mathbf{1}_{\mathcal{H}_1}$ is nevertheless bipartite (the signs $\tau$ give the bipartition) and we can introduce the half-spectral curve $\Sigma_{\boldsymbol{v}}^{\textnormal{half}} \subset \Sigma_{\boldsymbol{v}}$ as in Definition~\ref{def:halfspcurvebip}. It is built by gluing opposite hemispheres of adjacent pairs of sheets. The non-bipartite nature of $\mathsf{G}$ however gives us for each trapezoid $\mathcal{T}_i$ two hemispheres in the half-spectral curve: an upper half-plane in the sheet $\mathbf{1}_{\mathcal{H}_i}$, a lower half-plane in the sheet $-\mathbf{1}_{\mathcal{H}_i}$). Therefore, $\Sigma_{\boldsymbol{v}}^{\textnormal{half}}$ is homeomorphic to the orientation covering of the glued domain. The involution of this covering is realized in $\Sigma_{\boldsymbol{v}}^{\text{half}} \subset \Sigma_{\boldsymbol{v}}$ by the composition of $\varsigma$ with the complex conjugation. Exceptionally, we denote it $\overline{\varsigma}$, reserving the notation $\varsigma^{*}$ to the operation of pulling-back via $\varsigma$. The bidifferential $\mathcal{B}_{\boldsymbol{v}}$ has a structure different than the one of Corollary~\ref{Greenwon}, but we can relate it in a similar way to an even version of the Green function.

\begin{lemma}
\label{Lem:greenferm}
In the situation of Theorem~\ref{thm:spcurvenonbip}, define the even Green function as $\textnormal{Gr}^+(p_1,p_2) = \textnormal{Gr}(p_1,p_2) + \textnormal{Gr}(p_1,\overline{\varsigma}(p_2))$, where $\textnormal{Gr}$ is the usual Green function on $\Sigma_{\boldsymbol{v}}^{\textnormal{half}}$ with Dirichlet boundary conditions. We have
\begin{equation}
\label{eq:GRsigma}
\begin{split}
\textnormal{Gr}^{+}(p_1,p_2) & = - \frac{1}{4\pi}\left(\textnormal{gr}_{\boldsymbol{v}}(p_1,p_2) + 2\log\bigg|\frac{z_1 - z_2}{z_1 - z_2^*}\bigg| \right)
\end{split}
\end{equation}
\end{lemma}
\begin{proof} The proof strategy is the same than in Corollary~\ref{Greenwon}. We study the function
\begin{equation}
\label{grvp1p200}
\textnormal{gr}_{\boldsymbol{v}}(p_1,p_2) = \textnormal{Re}\left(\int_{z_1^*}^{z_1}\int_{z_2^*}^{z_2} \mathcal{B}_{\boldsymbol{v}}\right)
\end{equation}
of Definition~\ref{def:Hvv} with the bidifferential $\mathcal{B}_{\boldsymbol{v}}$ from Theorem~\ref{thm:spcurvenonbip}. By Proposition~\ref{prop:Greengen}, the function $\textnormal{gr}_{\boldsymbol{v}}$ is harmonic on $\Sigma_{\boldsymbol{v}} \times \Sigma_{\boldsymbol{v}}$ in both variables $p_1$ and $p_2$ except for logarithmic singularities as $z_1 = z_2$ and $z_1 = z_2^*$. The coefficient of the logarithmic singularity depends on the sheets $\boldsymbol{w}_1,\boldsymbol{w}_2$ to which $p_1,p_2$ belong, specified by the matrix of biresidue of the bidifferential given in Theorem~\ref{thm:spcurvenonbip}. Concretely
\begin{itemize}
\item If $\boldsymbol{w}_1 = \boldsymbol{w}_2$, there are no singularity;
\item If $\boldsymbol{w}_1 = - \boldsymbol{w}_2$, there is a singularity $-4 \log \big|\frac{z_1 - z_2}{z_1 - z_2^*}\big|$;
\item Otherwise, there is a singularity $-2\log\big|\frac{z_1 - z_2}{z_1 - z_2^*}\big|$.
\end{itemize}
Therefore, the function $-\frac{1}{4\pi}\big(\textnormal{gr}_{\boldsymbol{v}}(p_1,p_2) + 2\log\big|\frac{z_1 - z_2}{z_1 - z_2^*}\big|\big)$ has a singularity singularities $\mp \frac{1}{2\pi}\log\big|\frac{z_1 - z_2}{z_1 - z_2^*}\big|$ if $\boldsymbol{w}_1 = \pm \boldsymbol{w}_2$, and no singularity otherwise.

Let us now restrict $p_1,p_2$ to $\Sigma_{\boldsymbol{v}}^{\textnormal{half}}$. We find a logarithmic singularity as $p_1 = p_2$ (corresponding to $z_1 = z_2$ and $\boldsymbol{w}_1 = \boldsymbol{w}_2$) with coefficient $-\frac{1}{2\pi}$, and a logarithmic singularity as $p_1 = \overline{\varsigma}(p_2)$ (corresponding to $z_1 = z_2^*$ and $\boldsymbol{w}_1 = - \boldsymbol{w}_2$) with same coefficient $-\frac{1}{2\pi}$. In comparison, the usual Green function $\textnormal{Gr}(p_1,p_2)$ is also harmonic in both variables on $\Sigma_{\boldsymbol{v}}^{\textnormal{half}} \times \Sigma_{\boldsymbol{v}}^{\textnormal{half}}$ except for a logarithmic singularity with coefficient $-\frac{1}{2\pi}$ at $p_1 = p_2$. Therefore
\begin{equation}
\label{testevengreen}
-\frac{1}{4\pi}\bigg(\textnormal{gr}_{\boldsymbol{v}}(p_1,p_2) + 2\log\bigg|\frac{z_1 - z_2}{z_1 - z_2^*}\bigg|\bigg) - \big(\textnormal{Gr}(p_1,p_2) + \textnormal{Gr}(p_1,\overline{\varsigma}(p_2)\big)
\end{equation}
is a harmonic function in both variables in $\Sigma_{\boldsymbol{v}}^{\textnormal{half}} \times \Sigma_{\boldsymbol{v}}^{\textnormal{half}}$, without any singularity. By construction the Green function tends to zero when $p_1$ or $p_2$ approaches the boundary of $\Sigma_{\boldsymbol{v}}^{\textnormal{half}}$. As boundary points are real points, the same is true for the first bracket in \eqref{testevengreen} by direct inspection of the definition \eqref{grvp1p200} and the logarithmic term. As there is no non-zero harmonic function on $\Sigma_{\boldsymbol{v}}^{\textnormal{half}}$ having zero boundary values, we conclude that \eqref{testevengreen} is identically zero, justifying the desired formula.
\end{proof}

\begin{remark} \label{Remnonorlink} With the particular $\boldsymbol{v}$ chosen in Theorem~\ref{thm:spcurvenonbip}, $\Sigma_{\boldsymbol{v}}^{\textnormal{half}}$ matches the construction of $\Sigma^{\textnormal{half}}$ in Section~\ref{sec:nonorGFFproof}. The construction of the bidifferential $\mathcal{B}_{\boldsymbol{v}}$ (like in Remark~\ref{rem:match1113} but taking into account the different structure of biresidues in Theorem~\ref{thm:spcurvenonbip}) leads to the relation \eqref{gr1222} between the even Green function and the functions $\mathcal{F}^{(i_1,i_2)}$ given in \eqref{Fqqq12}. This was an essential ingredient in the proof of Theorem~\ref{Theorem_GFF_general_nonor}.
\end{remark}

\begin{corollary}
\label{cor:Greenfermion}
Consider the situation of Theorem~\ref{thm:spcurvenonbip}. Identify for each $h \in [H]$ the interval $(\alpha_h,\beta_h)$ in the glued domain with one (arbitrarily chosen) of its two images in $\Sigma_{\boldsymbol{v}}^{\textnormal{half}}$. Let $(f_h(z))_{h = 1}^{H}$ be an $H$-tuple of functions such that $f_h(z)$ is a holomorphic function of $z$ in a complex neighborhood of $[\alpha_h,\beta_h]$, taking real values on the real axis, and call $\tilde{f}_{h}(p) = f_h(Z_{\boldsymbol{v}}(p))$ the function in the corresponding domain of $\Sigma_{\boldsymbol{v}}^{\textnormal{half}}$. Then, for any $h_1,h_2 \in [H]$
\[
\oint_{\gamma_{h_1}} \oint_{\gamma_{h_2}} \mathcal{F}_{h_1,h_2}(z_1,z_2) f_{h_1}(z_1) f_{h_2}(z_2) \frac{\dd z_1 \dd z_2}{(2\ii\pi)^2} = \frac{1}{\pi} \int_{\alpha_{h_1}}^{\beta_{h_1}}  \int_{\alpha_{h_2}}^{\beta_{h_2}} \textnormal{Gr}^{+}(p_1,p_2) \dd \tilde{f}_{h_1}(p_1) \dd \tilde{f}_{h_2}(p_2).
\]
\end{corollary}
\begin{proof}
It is a direct consequence of Proposition~\ref{prop:Greengen} (H5) and Lemma~\ref{Lem:greenferm}.
\end{proof}

\subsection{Proof of Theorem~\ref{thm:Omegav}}

The definition of $\boldsymbol{\Theta}$ in \eqref{ThetaG} yields
\[
\bth_h = \boldsymbol{\Theta}(\boldsymbol{e}^{(h)}) = \frac{1}{2}\big(\mathbf{1}_{\mathcal{H}_{i_h^+}} + \mathbf{1}_{\mathcal{H}_{i_h^-}}\big),
\]
and we recall that $\theta_{h,h} = 1$. From Definition~\ref{unimono}, for any $i \in [m]$ and $h \in [H]$ we have
\begin{equation}
\label{Th111}T^{(h)}(\mathbf{1}_{\mathcal{H}_{i}}) = \left\{\begin{array}{lll} \mathbf{1}_{\mathcal{H}_i} && \textnormal{if}\,\,h \notin \mathcal{H}_i, \\[0.6ex] \mathbf{1}_{\mathcal{H}_i} - 2\bth_h = - \mathbf{1}_{\mathcal{H}_{i_h'}} && \textnormal{if}\,\,h \in \mathcal{H}_i. \end{array}\right.
\end{equation}
Since $\epsilon_i = - \epsilon_{i_h'}$ we deduce that the set $\{\epsilon_i \mathbf{1}_{\mathcal{H}_i}\,\,|\,\,i \in [m]\}$ is stable under the action of $\mathfrak{G}$. The action on this set is transitive because we assumed that the two-dimensional domain obtained from gluing trapezoids is connected, so it is an orbit. This orbit is good because each segment belongs at least to one trapezoid.

We would like to compute the group $\hat{\mathfrak{G}}$, which by definition is isomorphic to the group generated by the restriction of the generators of $\mathfrak{G}$ to the stable subspace $\mathscr{V} = \textnormal{Im}(\boldsymbol{\Theta})$. By Proposition~\ref{lem:H1G} we know that $\mathscr{V} = \textnormal{Im}(\boldsymbol{\Theta})$ is spanned by $\mathbf{1}_{\mathcal{H}_1},\ldots,\mathbf{1}_{\mathcal{H}_m}$. We exploit the bipartition and define vectors $\boldsymbol{f}^{(i)} = \epsilon_i \mathbf{1}_{\mathcal{H}_i}$ for $i \in [m]$. They are subjected to the unique relation $\sum_{i = 1}^m \boldsymbol{f}^{(i)} = 0$ because each segment belongs to exactly one positive and one negative trapezoid. So, the vector space $\mathscr{V}$ has the same presentation as the root space of type $A_{m - 1}$, and the corresponding Weyl group is the group of permutations of the vectors $\boldsymbol{f}^{(1)},\ldots,\boldsymbol{f}^{(m)}$. In \eqref{Th111} we saw that for each $h \in [H]$, the generator $T^{(h)}$ acts on $\mathscr{V}$ as the transposition of the vectors $\boldsymbol{f}^{(i_h^-)}$ and $\boldsymbol{f}^{(i_h^+)}$. To show that $\hat{\mathfrak{G}}$ is the Weyl group of type $A_{m - 1}$ it suffices to show that all transpositions of two $\boldsymbol{f}$-vectors can be obtained by composition of our generators. And indeed, let $i,j$ be distinct elements of $[m]$ and take a path from the $i$-th vertex to the $j$-th vertex in the connected gluing graph $\mathsf{G}$. Calling $h_0,\ldots,h_s$ be the sequence of edges along this path, the composition
\[
T^{(h_s)} \circ T^{(h_{s - 1})} \circ \cdots \circ T^{(h_1)} \circ T^{(h_0)} \circ T^{(h_1)} \circ \cdots \circ T^{(h_{s - 1})} \circ T^{(h_s)}
\]
acts on $\mathscr{V}$ as the transposition of $\boldsymbol{f}^{(i)}$ and $\boldsymbol{f}^{(j)}$.

We offer a second proof of $\hat{\mathfrak{G}} = A_{m - 1}$, which is less elegant as it relies on classification but illustrates the method of looking at orbit sizes to find the group. Proposition~\ref{lemrefl} tells us that the group $\hat{\mathfrak{G}}$ is a finite Euclidean reflection group, generated by orthogonal reflections with respect to vectors $\hat{\boldsymbol{e}}^{(h)}$ indexed by $h \in [H]$. The latter have scalar products
\begin{equation}
\label{scprodgh}
\big\langle \hat{\boldsymbol{e}}^{(h)} \cdot \boldsymbol{\Theta}(\hat{\boldsymbol{e}}^{(g)})\big\rangle = \theta_{g,h}
\end{equation}
and since $\boldsymbol{\Theta}$ has the connectedness property, we cannot find a partition $[H] = P_1 \sqcup P_2$ with non-empty $P_1,P_2$ such that \eqref{scprodgh} is zero for any $g \in P_1$ and $h \in P_2$. Therefore, $\hat{\mathfrak{G}}$ is irreducible. By Proposition~\ref{lem:H1G} it has rank equal to $\textnormal{rank}(\boldsymbol{\Theta}) = m - 1$. Since $\mathfrak{G}.\boldsymbol{v}$ has a orbit of size $m$, the general principles from Lemma~\ref{lem:comparor} tell us that $\hat{\mathfrak{G}}$ has an orbit of size $\leq m$. Among irreducible finite Euclidean reflections group of rank $m - 1$, only the Weyl group of type $A_{m - 1}$ has an orbit of this size. We also observe that this orbit is minimal (\textit{cf.} Figure~\ref{Fig:minor}).

We now turn to the fundamental solution of the master Riemann--Hilbert problem $\boldsymbol{\mathcal{F}}(z,w)$. Coming back to \eqref{2RHPsf}, it satisfies for any $z,w \in \amsmathbb{C} \setminus \overline{\amsmathbb{B}}$, any $i,j \in [m]$, $h \in [H]$ and $x \in (\alpha_h,\beta_h)$
\begin{equation}
\label{runHPS}
\begin{split}
\big\langle (\epsilon_{i}\mathbf{1}_{\mathcal{H}_{i}} \otimes \epsilon_{j}\mathbf{1}_{\mathcal{H}_{j}}) \cdot \boldsymbol{\mathcal{F}}(x^+,w) \big\rangle & = \big\langle \big(T^{(h)}(\epsilon_{i}\mathbf{1}_{\mathcal{H}_{i}}) \otimes \epsilon_{j}\mathbf{1}_{\mathcal{H}_{j}}\big) \cdot \boldsymbol{\mathcal{F}}(x^-,w)\big\rangle - \frac{\epsilon_{i}\epsilon_{j}\delta_{h,\mathcal{H}_{i} \cap \mathcal{H}_{j}}}{(x - w)^2}, \\
\big\langle(\epsilon_{i}\mathbf{1}_{\mathcal{H}_{i}} \otimes \epsilon_{j}\mathbf{1}_{\mathcal{H}_{j}}) \cdot \boldsymbol{\mathcal{F}}(z,x^+) \big\rangle & = \big\langle \big(\epsilon_{i}\mathbf{1}_{\mathcal{H}_{i}} \otimes T^{(h)}(\epsilon_{j}\mathbf{1}_{\mathcal{H}_{j}})\big) \cdot \boldsymbol{\mathcal{F}}(z,x^-)\big\rangle - \frac{\epsilon_{i}\epsilon_{j}\delta_{h,\mathcal{H}_{i} \cap \mathcal{H}_{j}}}{(z - x)^2} ,
\end{split}
\end{equation}
where $\delta_{h,\mathcal{H}}$ is equal to $1$ if $h \in \mathcal{H}$ and $0$ otherwise. Since each $h \in [H]$ belongs to exactly one negative and one positive trapezoid, there are only two situations where the shift in the right-hand side is zero:
\begin{itemize}
\item if $\mathcal{T}_{i} = \mathcal{T}_{j}$, \textit{i.e.} $i = j$, then we shift by $-\frac{1}{(z - w)^2}$;
\item if $\mathcal{T}_{i},\mathcal{T}_{j}$ are facing each other along the $h$-th segment, \textit{i.e.} $i_h' = j$, then we shift by $\frac{1}{(z - w)^2}$.
\end{itemize}
We stress that the shift does not depend whether $T^{(h)}$ is applied to the first or second variable. Consider a finite sequence
\[
(\mathcal{T}_{i_0},\mathcal{T}_{j_0})\,\, \mathop{\rightsquigarrow}^{h_1}\,\, (\mathcal{T}_{i_1},\mathcal{T}_{j_1}) \,\,\mathop{\rightsquigarrow}^{h_2} \,\,(\mathcal{T}_{i_2},\mathcal{T}_{j_2}) \,\, \rightsquigarrow \,\,\cdots\,\, \mathop{\rightsquigarrow}^{h_d} \,\, \mathcal{T}_{i_d} \times \mathcal{T}_{j_d}
\]
of ordered pairs of trapezoids starting from $(i_0,j_0) = (1,1)$ and such that at each step
\begin{itemize}
\item[(S1)] either $\mathcal{T}_{i_c}$ is the trapezoid sharing the $h_c$-th segment with $\mathcal{T}_{i_{c - 1}}$, and $\mathcal{T}_{j_c} = \mathcal{T}_{i_{c - 1}}$;
\item[(S2)] or $\mathcal{T}_{i_c} = \mathcal{T}_{i_{c - 1}}$ and $\mathcal{T}_{j_c}$ is the trapezoid sharing the $h_c$-th segment with $\mathcal{T}_{j_{c - 1}}$.
\end{itemize}
We want to compute the constant $t_c$ in the shift $\frac{t_c}{(z - w)^2}$ that results from analytic continuation of $\big\langle (\epsilon_1\mathbf{1}_{\mathcal{H}_1} \otimes \epsilon_1 \mathbf{1}_{\mathcal{H}_1}) \cdot \boldsymbol{\mathcal{F}}(z,w)\big\rangle$ by letting $z$ (in case S1) or $w$ (in case S2) crossing $\amsmathbb{B}_{h_c}$ from the upper half-plane to the lower half-plane using \eqref{runHPS}, successively for each $c \in [d]$. The only steps where $t_c$ is updated are those when $\mathcal{T}_{i_c}$ and $\mathcal{T}_{j_c}$ are either equal (we get a $-1$), or opposite to each other (we get a $+1$). These updates occur irrespectively of the segment across which we make the transition. We start with $t_0 = 0$. At the first step we replace one of the copy of $\mathcal{T}_1$ by a trapezoid facing it so we get $t_1 = -1$. The value of $t$ then remains constant until we meet a step where $\mathcal{T}_{i_{c - 1}}$ and $\mathcal{T}_{j_{c - 1}}$ are facing each other and we transition to $\mathcal{T}_{i_c} = \mathcal{T}_{j_c}$, so that $t_c = 0$. Repeating the argument starting with $\mathcal{T}_{i_c} \times \mathcal{T}_{j_c}$ instead of $\mathcal{T}_1 \times \mathcal{T}_1$ shows by induction that along any sequence as above, the value of the shift is $0$ if and only if $i_c = j_c$, and $-1$ in all other cases. Notice that for any $c \in [d]$ the value of $t_c$ only depends on $(i_c,j_c)$ and not on the sequence of steps it was reached starting from $(1,1)$. In the language of Section~\ref{sec:anaBext} and setting $\boldsymbol{v} = \epsilon_1 \mathbf{1}_{\mathcal{H}_1}$ and $\widetilde{\boldsymbol{v}} = (\boldsymbol{v} \otimes \boldsymbol{v}\,,\,0)$, we have just shown that
\begin{equation}
\label{G2orbitun92}
\mathfrak{G}^{[2]}.\widetilde{\boldsymbol{v}} = \big\{(\epsilon_i \mathbf{1}_{\mathcal{H}_i} \otimes \epsilon_j \mathbf{1}_{\mathcal{H}_j}\,,\,\delta_{i,j} - 1)\quad |\quad i,j \in [m]\big\}.
\end{equation}
In particular, the orbit-projection map $\pi^{[2]}_{\boldsymbol{v}}$ is a bijection and we get from Theorem~\ref{thm:Bcont} a bidifferential $\mathcal{B}_{\boldsymbol{v}}$ on $\Sigma^{[2]}_{\widetilde{\boldsymbol{v}}} = \Sigma_{\boldsymbol{v}} \times \Sigma_{\boldsymbol{v}}$ by analytic continuation of
\[
\big\langle (\epsilon_1 \mathbf{1}_{\mathcal{H}_1} \otimes \epsilon_1 \mathbf{1}_{\mathcal{H}_1}) \cdot \boldsymbol{\mathcal{F}}(z_1,z_2) \big\rangle \dd z_1 \dd z_2 = \sum_{g,h \in \mathcal{H}_1} \mathcal{F}_{g,h}(z_1,z_2)\dd z_1\dd z_2.
\]
The matrix of biresidues of \textbf{\foreignlanguage{russian}{B}} is read from the scalar component in \eqref{G2orbitun92}: it is $\textnormal{\textbf{Id}} - \boldsymbol{1}$. Proposition~\ref{funddif} then expresses $\mathcal{B}_{\boldsymbol{v}}$ in terms of fundamental bidifferentials.

\subsection{Proof of Theorem~\ref{thm:spcurvenonbip}}

We follow the same scheme of proof as in Theorem~\ref{thm:Omegav}. Equation~\eqref{Th111} is valid, but the non-existence of bipartitions now implies that the $\mathfrak{G}$-orbit of $\mathbf{1}_{\mathcal{H}_1}$ is the full $\{\tau\,\mathbf{1}_{\mathcal{H}_i}\,\,|\,\,\tau \in \{\pm 1\}\,\, i \in [m]\}$. The orbit is good for the same reasons. The same connectedness argument show that the group $\hat{\mathfrak{G}}$ is irreducible of rank $m$, and knowing that we have an orbit of size $2m$ only gives the two possibilities $\hat{\mathfrak{G}} \simeq B_m$ or $\hat{\mathfrak{G}} \simeq D_m$. By looking more carefully at the generators, we will see that only the last possibility is realized. Recall that the Weyl group of type $D_m$ is the group of signed permutations of $[m]$ with an even number of signs. This is the group of permutations $g$ of $\{-m,\ldots,-1,1,\ldots,m\}$ such that $g(-i) = -g(i)$ and $\#\{i \in [m]\,\,|\,\,g(i) < 0\}$ is even. It is a semi-direct product $(\amsmathbb{Z}_2)^{m - 1} \rtimes \mathfrak{S}_m$ and is generated by the signed transpositions
\[
\pm i \mapsto \mp j,\qquad j \mapsto \mp i,\qquad k \mapsto k \quad (k \neq \pm i,\pm j).
\]
Equivalently, the signed transposition is represented by the endomorphism of $\amsmathbb{R}^m$ transposing the $i$-th and the $j$-th basis vectors and flipping their sign.

Proposition~\ref{lem:H1G} in the non-bipartite case tells us that $\mathbf{1}_{\mathcal{H}_1},\ldots,\mathbf{1}_{\mathcal{H}_m}$ is a basis of $\mathscr{V} = \textnormal{Im}(\boldsymbol{\Theta})$. For each $h \in [H]$, calling $\mathcal{T}_{i_h},\mathcal{T}_{i_h'}$ the two trapezoids containing the $h$-th segment, \eqref{Th111} says that $T^{(h)}$ acts on the basis vectors as the linear map
 \[
\mathbf{1}_{\mathcal{H}_{i_h}} \mapsto -\mathbf{1}_{\mathcal{H}_{i_h'}},\qquad \mathbf{1}_{\mathcal{H}_{i_h'}} \mapsto - \mathbf{1}_{\mathcal{H}_{i_h}},\qquad \mathbf{1}_{\mathcal{H}_k} \mapsto \mathbf{1}_{\mathcal{H}_k} \quad (k \notin i_h,i_h')
\]
representing a signed transposition. Composing any number of such transformations certainly gives a signed permutation of $[m]$ with an even number of signs. To prove $\hat{\mathfrak{G}} \simeq D_m$ it suffices to check that any signed transposition of $[m]$ arises in this way. And indeed, if $i,j$ are two distinct elements of $[m]$, we can find an even-length path $h_0,\ldots,h_{2s + 1}$ from the $i$-th vertex and the $j$-th vertex in the gluing graph $\mathsf{G}$ because it is non-bipartite, and the resulting transformation
\[
 T^{(h_{2s + 1})} \circ T^{(h_{2s})} \circ \cdots \circ T^{(h_1)} \circ T^{(h_0)} \circ T^{(h_1)} \circ \cdots \circ T^{(h_{2s})} \circ T^{(h_{2s + 1})}
 \]
 acts on the basis of $\mathscr{V}$ as the signed transposition of the $i$-th and the $j$-th basis vector. In the Weyl group of type $D_m$ the minimal size of an orbit is $2m$, so the orbit we found is indeed minimal.

To compute the fundamental solution of the master Riemann--Hilbert problem, Equation~\ref{runHPS} is also valid and we can track the shifts of double pole as we did in the proof of Theorem~\ref{thm:Omegav}. Start with $\tau_0 = \eta_0 = 1$ and $\mathcal{T}_{i_0} = \mathcal{T}_{j_0} = \mathcal{T}_{1}$ and take a finite sequence of ordered pairs of elements of $\mathfrak{G}.\boldsymbol{v}$, here written as signed trapezoids
\[
(\tau_0,\mathcal{T}_{i_0}; \eta_0,\mathcal{T}_{j_0}) \,\, \mathop{\rightsquigarrow}^{h_1}\,\, (\tau_{1},\mathcal{T}_{i_1};\eta_1,\mathcal{T}_{j_1})  \,\,\mathop{\rightsquigarrow}^{h_2} \,\,(\tau_2,\mathcal{T}_{i_2};\eta_2 \mathcal{T}_{j_2}) \,\, \rightsquigarrow \,\,\cdots\,\, \mathop{\rightsquigarrow}^{h_d} \,\, (\tau_{d},\mathcal{T}_{i_d}; \eta_d,\mathcal{T}_{j_d})
\]
with $\tau_c,\eta_c \in \{\pm 1\}$ and $i_c,j_c \in [m]$ for any $c \in [d]$ and such that at each step
\begin{itemize}
\item[(S1)] either $\tau_{c - 1} = - \tau_c$ and $\mathcal{T}_{i_{c}}$ is the trapezoid sharing the $h_c$-th segment with $\mathcal{T}_{i_{c - 1}}$ while $(\eta_{c - 1},\mathcal{T}_{j_{c - 1}}) = (\eta_c,\mathcal{T}_{j_c})$;
\item[(S2)] or $(\tau_{c - 1},\mathcal{T}_{i_{c - 1}}) = (\tau_{c},\mathcal{T}_{i_c})$ while $\eta_c = - \eta_{c - 1}$ and $\mathcal{T}_{j_{c}}$ is the trapezoid sharing the $h_c$-th segment with $\mathcal{T}_{j_{c -1}}$.
\end{itemize}
As before, we let $t_c$ be the constant in the double pole after we have reached the $c$-th step. The value of $t_c$ is updated only when $\mathcal{T}_{i_{c - 1}}$ and $\mathcal{T}_{j_{c - 1}}$ are either equal or glued to each other: in both situations $t_{c} - t_{c - 1} = -\tau_{c - 1}\eta_{c - 1} = \tau_c\eta_c$. We start with $t_0 = 0$. At the first step we get $t_1 = -1$, and then $t_c$ remain constant until we reach the situation where $\mathcal{T}_{i_{c - 1}}$ and $\mathcal{T}_{j_{c - 1}}$ are glued to each other and we make the transition to $\mathcal{T}_{i_c} = \mathcal{T}_{j_c}$. Then we get $t_c = -1 + \tau_c \eta_c$. If $\tau_c = \eta_c$ we get back to $t_c = 0$; if $\tau_c = -\eta_c$ we rather get $t_c = -2$, and at the next step $t_{c + 1} = -1$. Continuing the process we see the value of $t_c$ only depends on $(\tau_{c},i_c;\eta_c,j_c)$ and not on the way it was reached from $(1,\mathcal{T}_1;1,\mathcal{T}_1)$, and
\[
\forall c \in [d] \qquad t_c = \left\{\begin{array}{lll} 0 && \textnormal{if}\,\,i_c = j_c \,\,\textnormal{and}\,\,\tau_c = \eta_c \\ -2 && \textnormal{if}\,\,i_c = j_c \,\,\textnormal{and}\,\,\tau_c = - \eta_c \\ -1 && \textnormal{if}\,\,i_c \neq j_c \end{array}\right.
\]
Even more, it only depends on $(i_c,j_c)$ and the product $\tau_c\eta_c \in \{\pm 1\}$. This proves
\begin{equation}
\label{G2333}\mathfrak{G}^{[2]}.\widetilde{\boldsymbol{v}} = \big\{(\epsilon \mathbf{1}_{\mathcal{H}_i} \otimes \mathbf{1}_{\mathcal{H}_j}\,,\,\epsilon \delta_{i,j} - 1)\quad | \quad \epsilon \in \{\pm 1\}\quad i,j \in [m]\big\}.
\end{equation}
in particular $\Sigma_{\widetilde{\boldsymbol{v}}}^{[2]} = \Sigma_{\boldsymbol{v} \otimes \boldsymbol{v}}^{[2],\circ}$ and the bidifferential $\mathcal{B}_{\widetilde{\boldsymbol{v}}}$ from Theorem~\ref{thm:Bcont} is defined on this spectral surface.

The orbit-projection map $\pi^{\bullet}_{\boldsymbol{v}}$ sends $(\tau\mathbf{1}_{\mathcal{H}_i}, \eta \mathbf{1}_{\mathcal{H}_j}\,,\,\tau\eta \delta_{i,j} - 1)$ to $(\tau\eta \mathbf{1}_{\mathcal{H}_i} \otimes \mathbf{1}_{\mathcal{H}_j}\,,\,\tau\eta\delta_{i,j} - 1)$. Comparing with Proposition~\ref{G2GG} we see that the group morphism $\pi^{\bullet} : \mathfrak{G} \times \mathfrak{G} \rightarrow \mathfrak{G}^{[2],\circ}$ is not a bijection and we must have $-\textnormal{\textbf{Id}} \in \mathfrak{G}$. In particular, we have
\[
\mathfrak{G}^{[2],\circ}.(\boldsymbol{v} \otimes \boldsymbol{v}) = \mathfrak{G}.\boldsymbol{v} \times \mathfrak{G}.\boldsymbol{v} \big/(-\textnormal{\textbf{Id}},-\textnormal{\textbf{Id}}).
\]
Denoting $\varsigma$ the holomorphic involution $-\textnormal{\textbf{Id}}$ induces on $\Sigma_{\boldsymbol{v}}$ by action on the sheet labels, we get
\[
\Sigma_{\boldsymbol{v} \otimes \boldsymbol{v}}^{[2],\circ} = \Sigma_{\boldsymbol{v}} \times \Sigma_{\boldsymbol{v}} \big/(\varsigma,\varsigma).
\]
The pullback of $\mathcal{B}_{\widetilde{\boldsymbol{v}}}$ to $\Sigma_{\boldsymbol{v}} \times \Sigma_{\boldsymbol{v}}$ then gives a meromorphic symmetric bidifferential with matrix of biresidues read from \eqref{G2333}, thus having entries $\text{\textit{\foreignlanguage{russian}{B}}}_{\boldsymbol{w}_1,\boldsymbol{w}_2} = \delta_{\boldsymbol{w}_1,\boldsymbol{w}_2} - \delta_{\boldsymbol{w}_1 + \boldsymbol{w}_2,0} - 1$ for any $\boldsymbol{w}_1,\boldsymbol{w}_2 \in \mathfrak{G}.\boldsymbol{v}$. The right-hand side of \eqref{Omev111} has exactly the same matrix of biresidues as $\mathcal{B}_{\boldsymbol{v}}$ and has zero integral with respect to first variable along cycles in $\mathcal{L}$, so must be equal to $\mathcal{B}_{\boldsymbol{v}}$ like in the proof of Proposition~\ref{funddif}. The equality $(\textnormal{id},\varsigma)^* \mathcal{B}^{\Sigma,\mathcal{L}} = (\varsigma,\textnormal{id})^*\mathcal{B}^{\Sigma,\mathcal{L}}$ comes as well from the fact that both sides have same matrix of biresidues and have zero integral with respect to the first variable along cycles in $\mathcal{L}$.

\subsection{Realizability and reflection groups for connectoids}
\label{sec:extraconn}

While the reflection group $\hat{\mathfrak{G}}$ is determined solely by $\boldsymbol{\Theta}$, we see in Theorems~\ref{thm:Omegav}-\ref{thm:spcurvenonbip} that the existence of bipartite or non-bipartite realizations imply that $\hat{\mathfrak{G}}$ is either or type A or D. Since these groups coincide only in rank $\leq 3$ and triangles are the only non-bipartite graphs with no more than three vertices, this has the following consequence.
\begin{corollary}
\label{cor:connsingle}
If a realizable connectoid admits both a bipartite and a non-bipartite realization, it must be the connectoid of a triangle ($m = 3$), possibly with multiple edges between its vertices.
\end{corollary}
Having reflection groups of type A or D is by no means a characterization of realizable connectoids. The non-realizable connectoids of size 4 (\textit{cf.} \eqref{exnonrel}) have $\hat{\mathfrak{G}} = D_4$, and the following non-realizable connectoid (found by an extensive search by Alain Goldberg) has $\hat{\mathfrak{G}} = E_6$:
\begin{equation*}
\begin{split}
\boldsymbol{\Theta} = \left(\begin{array}{cccccc} 1 & 0 & 0 & 0 & 0 & \frac{1}{2} \\[0.6ex] 0 & 1 & 0 & \frac{1}{2} & 0 & 0 \\[0.6ex] 0 & 0 & 1 & \frac{1}{2} & \frac{1}{2} & \frac{1}{2} \\[0.6ex] 0 & \frac{1}{2} & \frac{1}{2} & 1 & 0 & 0 \\[0.6ex] 0 & 0 & \frac{1}{2} & 0 & 1 & 0 \\[0.6ex] \frac{1}{2} & 0 & \frac{1}{2} & 0 & 0 & 1 \end{array}\right).
\end{split}
\end{equation*}

\chapter{Examples}

\label{Chap14}

This chapter illustrates the methods and results of Chapters~\ref{Chapter_SolvingN}-\ref{Chapter_AG}, and the logic to extract explicit (to the extent possible) formulae for discrete ensembles. The vocabulary and notions of these chapters will therefore be repeatedly used. Our exposition features mostly matrices of interaction $\boldsymbol{\Theta}$ relevant for the tiling models described in Section~\ref{Section_gluing_def} and aims at concrete descriptions of
\begin{itemize}
\item the spectral curve $\Sigma_{\boldsymbol{v}}$, on which all functions of interest are defined. Although it depends on a choice of seed vector $\boldsymbol{v} \in \amsmathbb{R}^H$, for tiling models there is a canonical choice dictated by the geometry of the domain, \textit{cf.} Theorems~\ref{thm:Omegav} and \ref{thm:spcurvenonbip};
\item the fundamental solution of the master Riemann--Hilbert problem, computing the leading covariance in the discrete ensemble with fixed filling fractions as found in Corollary~\ref{Corollary_CLT};
\item the first-kind functions, which can be extracted from holomorphic $1$-forms by Proposition~\ref{Lem:Holo1form} and computing the variation of the equilibrium measure with respect to filling fractions as found in Proposition~\ref{Proposition_differentiability_filling_fraction}.
\end{itemize}
We give a detailed treatment for all connected domains with three segments or less. We add the discussion of a few examples with more segments, that may or not be realized in tiling models, so as to stress the common features and the differences.

Many computations are reduced to the study of appropriate algebraic functions, and what is meant by ``explicit formulae'' in this context largely depends on the motivation and appreciation of the reader. There are indeed several layers of computational complexity to achieve this goal. Certainly the most explicit formulae can be obtained for simply-connected domains, for which the underlying spectral curve has genus $0$. In our case, it is easy to describe the general form of a uniformization map to the Riemann sphere depending on parameters, and the parameters are determined in terms of the endpoints of the bands in the discrete ensemble as solution of a system of algebraic equations. It would be part of an involved case-by-case analysis to get an ``explicit'' description of the solution of this system, in particular of the correct branch to be chosen among the many solutions of the algebraic system, and this is not attempted here. For domains with a hole, we need elliptic functions. The branches of parameters to be used in the elliptic functions is also subjected to more complicated systems of equations and, as we saw in Proposition~\ref{W20second}, some of these equations are of transcendental nature as they involve the period map. Domains with two holes require genus $2$ hyperelliptic functions, which are a level of complexity higher and about which much less is known. For domains with more than two holes more general algebraic functions can occur and although some formulae can be given --- \textit{e.g.} using fields of the form $\amsmathbb{C}(x,y)/P(x,y)$ for some polynomial $P$, or Riemann theta functions --- one must accept that their explicit character fades out.

\section{\texorpdfstring{$H = 2$}{H=2}: C-shaped domain}

\label{sec:Csha}
\begin{figure}[h!]
\begin{center}
\includegraphics[width=0.14\textwidth]{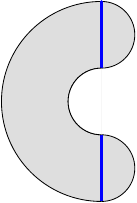}
\caption{\label{Fig:Cshape} Topology of the domain with matrix of interactions \eqref{Ex1}.}
\end{center}
\end{figure}

Consider the matrix $\boldsymbol{\Theta}$ for a C-shaped domain obtained by gluing three trapezoids (Figure~\ref{Fig:Cshape}).
\label{SecEx1}
\begin{equation}
\label{Ex1}\boldsymbol{\Theta} = \left(\begin{array}{cc} 1 & \tfrac{1}{2} \\ \tfrac{1}{2} & 1 \end{array}\right).
\end{equation}

Written in the canonical basis $(\boldsymbol{e}^{(1)},\boldsymbol{e}^{(2)})$ of $\amsmathbb{R}^2$, the two generators $T^{(1)}$ and $T^{(2)}$ of the group $\mathfrak{G}$ are
\[
T^{(1)} = \left(\begin{array}{cc} -1 & 0 \\ -1 & 1 \end{array}\right),\qquad T^{(2)} = \left(\begin{array}{cc} 1 & -1 \\ 0 & -1 \end{array}\right),
\]
and we have the relation $(T^{(1)} \circ T^{(2)})^3 = \textnormal{id}$. Therefore $\mathfrak{G}$ is the symmetric group in three elements, \textit{i.e.} the Weyl group of type $A_2$. Recalling Section~\ref{sec:bigstr}, it agrees with the fact that, $\boldsymbol{\Theta}$ being positive-definite, $\mathfrak{G}$ itself is conjugated to the reflection group $\hat{\mathfrak{G}}$ acting on $\mathscr{V} = \amsmathbb{R}^2$. The maximal parabolic subgroup is the symmetric group in two elements, hence the minimal size for an orbit is $3 = 3!/2!$. The $\mathfrak{G}$-orbit of $\boldsymbol{v} := \boldsymbol{e}^{(1)} = (1,0)$ has size $3$:
\begin{equation}
\label{eq:Cshapedvvv}
(1,0) \,\,\mathop{\longleftrightarrow}^{T^{(1)}} \,\, (-1,-1)\,\, \mathop{\longleftrightarrow}^{T^{(2)}} \,\,(0,1).
\end{equation}
This is the orbit exhibited by Theorem~\ref{thm:Omegav}.

\subsection{Rational parameterization of the spectral curve}

The spectral curve $\Sigma_{\boldsymbol{v}}$ has three sheets and four simple ramification points: by the Riemann--Hurwitz formula \eqref{RHformula} it has genus $0$. This agrees with the fact that the tiled domain in Figure~\ref{Fig:Cshape} is simply-connected. Therefore, there exists a biholomorphic map $\zeta : \Sigma_{\boldsymbol{v}} \rightarrow \widehat{\amsmathbb{C}}$ realizing the isomorphism to the Riemann sphere. The meromorphic function $Z := Z_{\boldsymbol{v}} : \Sigma_{\boldsymbol{v}} \rightarrow \widehat{\amsmathbb{C}}$ must then be a rational function of $\zeta$. To find its expression, we observe that $Z$ has three simple poles. Up adjusting $\zeta$ by application of a M\"obius transformation, we can choose the location of these poles to be $0,1,\infty$, respectively corresponding to the point with $Z = \infty$ in the sheets $(0,1)$, $(-1,-1)$, $(1,0)$. The partial fraction decomposition of $Z(\zeta)$ then takes the form
\begin{equation}
\label{Zvzeta}Z(\zeta) = \frac{u_0}{\zeta} + \frac{u_1}{\zeta - 1} + u' + u_{\infty}\zeta.
\end{equation}
for some coefficients $u_0,u_1,u_{\infty} \in \amsmathbb{C}^*$ and $u' \in \amsmathbb{C}$. The coefficients are determined by imposing that the location of the branchpoints of $Z$ are located at $\alpha_1 < \beta_1 < \alpha_2 < \beta_2$. Let us call $s_1,s_2,s_3,s_4$ the ramification points, \textit{i.e.} the zeros of $Z'(\zeta)$, that we label such that
\begin{equation}
\label{Zvalpha}
Z(s_1) = \alpha_1,\qquad Z(s_2) = \beta_1,\qquad Z(s_3) = \alpha_2,\qquad Z(s_4) = \beta_2.
\end{equation}
Since we know the zeros and poles of $Z'(\zeta)$, as well as its behavior near $\zeta = \infty$ for instance, we get
\begin{equation}
\label{Zvzetader}
Z'(\zeta) = u_{\infty} \frac{\prod_{j = 1}^{4} (\zeta - s_j)}{\zeta^2(\zeta - 1)^2}.
\end{equation}
Comparing with \eqref{Zvzeta} and its behavior at $\zeta = 0$ and $\zeta = 1$, we deduce
\[
u_0 = - u_{\infty}\prod_{j = 1}^{4} s_j,\qquad u_1 = - u_{\infty}\prod_{j = 1}^{4} (s_j - 1).
\]
There are additional constraints on the parameters appearing in \eqref{Zvzetader} coming from the fact that $Z'(\zeta)$ is the derivative of a rational function, thus cannot have residues. As the sum of all residues of any rational function vanishes, it is sufficient to examine the two constraints given by the vanishing of residues at $0$ and $\infty$, namely
\begin{equation}
\label{ssss}\sum_{j = 1}^4 s_j = \sum_{j = 1}^4 \frac{1}{s_j} = 2.
\end{equation}
Imposing \eqref{Zvalpha} in \eqref{Zvzeta}, we get parametric expressions for the four endpoints $\alpha_1,\beta_1,\alpha_2,\beta_2$ in terms of the six parameters $s_1,s_2,s_3,s_4,u',u_{\infty}$ subjected to the two constraints \eqref{ssss}
\begin{equation*}
\begin{split}
\alpha_1 - \beta_1 & = u_{\infty}(s_1 - s_2)(2 - s_3 - s_4 + s_3s_4) \\
\alpha_2 - \beta_2 & = u_{\infty}(s_3 - s_4)(2 - s_1 - s_2 + s_1s_2) \\
\alpha_1 - \alpha_2 & = u_{\infty}(s_1 - s_3)(2 - s_2 - s_4 + s_2s_4) \\
\alpha_1 + \beta_1 + \alpha_2 + \beta_2 & = u' + 2u_{\infty}\left(\sum_{i < j} s_is_j - \sum_{i < j <k} s_is_js_k\right).
\end{split}
\end{equation*}

\subsection{The fundamental solution and how to use it}

The unique fundamental bidifferential of the Riemann sphere is $\frac{\dd\zeta_1\dd \zeta_2}{(\zeta_1 - \zeta_2)^2}$. Theorem~\ref{thm:Omegav} then shows that the bidifferential $\mathcal{F}_{1,1}(z_1,z_2) \dd z_1 \dd z_2$ built from the fundamental solution of the master Riemann--Hilbert problem admits a meromorphic continuation to $\Sigma_{\boldsymbol{v}} \times \Sigma_{\boldsymbol{v}}$, given by
\begin{equation}
\label{1213eqnb}
\begin{split}
\mathcal{B}_{\boldsymbol{v}}(\zeta_1,\zeta_2) & = \left(\frac{1}{(\zeta_1 - \zeta_2)^2} - \frac{Z'(\zeta_1)Z'(\zeta_2)}{(Z(\zeta_1) - Z(\zeta_2))^2}\right)\dd \zeta_1\dd \zeta_2 \\
& = \dd_{\zeta_1}\dd_{\zeta_2}\log\left(\frac{\zeta_1 - \zeta_2}{Z(\zeta_1) - Z(\zeta_2)}\right) \\
& = -\dd_{\zeta_1}\dd_{\zeta_2}\log\left(1 - \frac{u_0}{u_{\infty}\zeta_1\zeta_2} - \frac{u_1}{u_{\infty}(\zeta_1 - 1)(\zeta_2 - 1)}\right)
\end{split}
\end{equation}
We summarize in Figure~\ref{fig:doublesheetBC} the way $\mathcal{B}_{\boldsymbol{v}}$ is defined from the fundamental solution when its variables are in various sheets. This is dictated by the structure of the orbit \eqref{eq:Cshapedvvv}, while the shift by double poles were derived in Theorem~\ref{thm:Omegav} from the master Riemann--Hilbert problem with non-trivial source satisfied by $\boldsymbol{\mathcal{F}}(z_1,z_2)$ in both variables, \textit{cf.} Proposition~\ref{thmBfund}.

\begin{figure}[h!]
\[
\begin{array}{|c|c|c|c|}
\hline
\textnormal{sheet} & (1,0) & (-1,-1) & (0,1) \\[4pt]
\hline
(1,0) & \mathcal{F}_{1,1} + \frac{1}{(Z_1 - Z_2)^2}& - \mathcal{F}_{1,1} - \mathcal{F}_{2,1} & \mathcal{F}_{2,1} \\[4pt]
\hline
(-1,-1) & -\mathcal{F}_{1,1} - \mathcal{F}_{1,2} & \sum_{h_1,h_2 = 1}^{2} \mathcal{F}_{h_1,h_2} + \frac{1}{(Z_1 - Z_2)^2} & - \mathcal{F}_{2,1} - \mathcal{F}_{2,2} \\[4pt]
\hline
(0,1) & \mathcal{F}_{1,2} & - \mathcal{F}_{1,2}- \mathcal{F}_{2,2} & \mathcal{F}_{2,2} + \frac{1}{(Z_1 - Z_2)^2} \\[4pt]
\hline
\end{array}
\]
\caption{\label{fig:doublesheetBC} The lines (respectively, columns) heads indicate the sheet to which the first (respectively, second) variable belongs to, and the entries are the expressions for the function $\frac{\mathcal{B}_{\boldsymbol{v}}}{\dd Z_1\dd Z_2} + \frac{1}{(Z_1 - Z_2)^2}$.}
\end{figure}

The last formula is adapted for expansion near $\zeta_i \rightarrow \infty$, corresponding to the pole of $Z$ in the sheet labeled $(1,0)$, which is the one where the bidifferential is equal to $\mathcal{F}_{1,1}(z_1,z_2)\dd z_1\dd z_2$. For expansion near $\zeta_i \rightarrow 1$ (pole of $Z$ in the sheet labeled $(-1,-1)$) or $\zeta_i \rightarrow 0$ (pole of $Z$ in the sheet labeled $(0,1)$) the following equivalent formulae are better suited:
\begin{equation*}
\begin{split}
\mathcal{B}_{\boldsymbol{v}}(\zeta_1,\zeta_2) & = -\dd_{\zeta_1}\dd_{\zeta_2}\log\left(1 + \frac{(\zeta_1 - 1)(\zeta_2 - 1)}{u_1}\bigg(\frac{u_0}{\zeta_1\zeta_2} - u_{\infty}\bigg)\right) \\
& = - \dd_{\zeta_1}\dd_{\zeta_2}\log\left(1 + \frac{\zeta_1\zeta_2}{u_{0}}\bigg( \frac{u_1}{(\zeta_1 - 1)(\zeta_2 - 1)} - u_{\infty}\bigg)\right).
\end{split}
\end{equation*}

\label{Csun02ex}
In the context of discrete ensembles with fixed filling fractions, the fundamental solution gives access to the leading covariance. Let us unfold how the above formulae can be used in practice, for the covariance involving polynomial test functions $f$. For fluctuations of position of particles in the first segment, we should use integrals around the segment $[\alpha_1,\beta_1]$ in the sheet labeled $(1,0)$ for both variables. Moving the contour to infinity in this sheet, this reduce to a residue computation near $\zeta = \infty$ in the spectral curve
\begin{equation*}
\begin{split}
& \quad \sum_{i,j = 1}^{N_1} \amsmathbb{E}^{(\textnormal{c})}\big[f(\ell^1_i)f(\ell^1_j)\big] + o(1) \\
& = \Res_{\zeta_1 = \infty} \Res_{\zeta_2 = \infty} f(Z(\zeta_1)) f(Z(\zeta_2))\,\mathcal{B}_{\boldsymbol{v}}(\zeta_1,\zeta_2) \\
& = -\Res_{\zeta_1 = \infty} \Res_{\zeta_2 = \infty} \dd f(Z(\zeta_1)) \dd f(Z(\zeta_2)) \log\left(1 - \frac{u_0}{u_{\infty}\zeta_1\zeta_2} - \frac{u_1}{u_{\infty}(\zeta_1 - 1)(\zeta_2 - 1)}\right).
\end{split}
\end{equation*}
In the last line we used integration by parts in both variables. Note that $\dd f(Z(\zeta)) = f'(Z(\zeta)) Z'(\zeta) \dd \zeta$ is a rational function with poles at $\zeta = 0,1,\infty$. The residues therefore extract a combination of coefficients in the series expansion of $\mathcal{B}_{\boldsymbol{v}}(\zeta_1,\zeta_2)$ near $\zeta_i \rightarrow \infty$. Likewise, we get the fluctuations of positions of particles in the second segment by a residue computation at $\zeta_i = 0$:
\begin{equation*}
\begin{split}
& \quad \sum_{i,j = 1}^{N_2} \amsmathbb{E}^{(\textnormal{c})}\big[f(\ell^2_i)f(\ell^2_j)\big] + o(1) \\
& = \Res_{\zeta_1 = 0} \Res_{\zeta_2 = 0} f(Z(\zeta_1)) f(Z(\zeta_2))\,\mathcal{B}_{\boldsymbol{v}}(\zeta_1,\zeta_2) \\
& = -\Res_{\zeta_1 = 0} \Res_{\zeta_2 = 0} \dd f(Z(\zeta_1)) \dd f(Z(\zeta_2))\log\left(1 + \frac{\zeta_1\zeta_2}{u_{0}}\bigg( \frac{u_1}{(\zeta_1 - 1)(\zeta_2 - 1)} - u_{\infty}\bigg)\right).
\end{split}
\end{equation*}
If we want the covariance between linear statistics of particles in the first segment and linear statistics of particles in the second segment, we rather compute, for polynomial test functions $f_1,f_2$
\begin{equation}
\begin{split}\nonumber
& \quad \sum_{i = 1}^{N_1} \sum_{j = 1}^{N_2} \amsmathbb{E}^{(\textnormal{c})}\big[f_1(\ell^1_i)f_2(\ell^2_j)\big] + o(1) \\
& = \Res_{\zeta_1 = \infty} \Res_{\zeta_2 = 0} f_1(Z(\zeta_1)) f_2(Z(\zeta_2))\,\mathcal{B}_{\boldsymbol{v}}(\zeta_1,\zeta_2)\nonumber \\
& = \Res_{\zeta_1 = \infty} \Res_{\zeta_2 = 0} \frac{f_1(Z(\zeta_1)) f_2(Z(\zeta_2))\dd \zeta_1 \dd \zeta_2}{(\zeta_1 - \zeta_2)^2}\nonumber \\
& = \sum_{m \geq 0} (m + 1) \Big(\Res_{\zeta_1 = \infty} f_1(Z(\zeta_1)) \zeta_1^{-(m + 2)}\dd\zeta_1\Big) \Big( \Res_{\zeta_2 = 0} f_2(Z(\zeta_2)) \zeta_2^{m}\dd\zeta_2\Big)\label{cov:Csha}
\end{split}
\end{equation}
To get the third line we removed the second term in \eqref{1213eqnb} from $\mathcal{B}_{\boldsymbol{v}}$ because it does not contribute to the residue. Unlike in the two previous cases, it leads here to a simpler formula because we are taking residues at two different points ($0$ and $\infty$) in the two variables and the first term of \eqref{1213eqnb} is already regular in this region. In the last line the sum over $m$ is finite because $f_1$ and $f_2$ are polynomials.

In comparison, the residues at $\zeta_i = 1$ record the (redundant) information on linear statistics of all particles
\begin{equation*}
\begin{split}
& \quad \sum_{i,j = 1}^{N} \amsmathbb{E}^{(\textnormal{c})}\big[f_1(\ell_i)f_1(\ell_j)\big] + o(1) \\
 & = \Res_{\zeta_1 = 1} \Res_{\zeta_2 = 1} f_1(Z(\zeta_1)) f_2(Z(\zeta_2))\,\mathcal{B}_{\boldsymbol{v}}(\zeta_1,\zeta_2) \\
 & = - \Res_{\zeta_1 = 1} \Res_{\zeta_2 = 1} \dd f_1(Z(\zeta_1)) \dd f_2(Z(\zeta_2)) \log\left(1 + \frac{(\zeta_1 - 1)(\zeta_2 - 1)}{u_1}\bigg(\frac{u_0}{\zeta_1\zeta_2} - u_{\infty}\bigg)\right).
\end{split}
\end{equation*}
The fact that the residues at $\zeta_i = 1$ encode redundant information can also be traced back from the fact that the sum of residues of a meromorphic form on the compact Riemann surface $\Sigma_{\boldsymbol{v}}$ vanishes.

The geometry of the domain fixes all filling fractions. Therefore, we do not have freedom to vary filling fractions in the setting of Proposition~\ref{Proposition_differentiability_filling_fraction}. This is also reflected in the fact that first-kind functions vanish. One way to see this is to remark that first-kind functions are encoded in holomorphic $1$-form by Proposition~\ref{Lem:Holo1form}, and on the Riemann sphere these must be zero.

\section{\texorpdfstring{$H = 2$}{H=2}: O-shaped domain (hexagon with a hole)}
\label{sec:Oshaped}
\begin{figure}[h!]
\begin{center}
\includegraphics[width=0.17\textwidth]{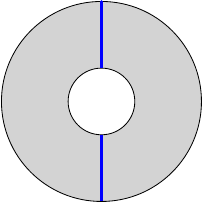}
\caption{\label{Fig:Oshaped} Topology of the domain with matrix of interactions \eqref{ExOshaped}.}
\end{center}
\end{figure}

The hexagon with a hole discussed in Section~\ref{Section_tiling_hole} (\textit{cf.} Figure~\ref{Fig:Oshaped}) corresponds to the gluing of two trapezoids with matrix of interactions
\begin{equation}
\label{ExOshaped}
\boldsymbol{\Theta} = \left(\begin{array}{cc} 1 & 1 \\ 1 & 1 \end{array}\right).
\end{equation}
The two generators of the group $\mathfrak{G}$ are
\[
T^{(1)} = \left(\begin{array}{cc} -1 & 0 \\ -2 & 1 \end{array}\right),\qquad T^{(2)} = \left(\begin{array}{cc} 1 & -2 \\ -1 & 0 \end{array}\right).
\]
Therefore, $\mathfrak{G}$ is the free group in two elements, while $\hat{\mathfrak{G}}$ is the Weyl group of rank $1$, \textit{i.e.} $\hat{\mathfrak{G}} = \amsmathbb{Z}_2$. The vector $\boldsymbol{v} := (1,1)$ has a very simple orbit
\[
(1,1)\,\, \mathop{\longleftrightarrow}^{T^{(1)},T^{(2)}}\,\, (-1,-1).
\]
In the notations of Proposition~\ref{lemrefl}, we have $\mathscr{V} = \textnormal{Im}(\boldsymbol{\Theta}) = \amsmathbb{R}.\boldsymbol{v}$.

The corresponding spectral curve $\Sigma_{\boldsymbol{v}}$ has two sheets and four ramification points, hence it has genus $1$. According to Theorem~\ref{thm:Omegav}, the bidifferential
\[
\mathcal{F}(z_1,z_2)\dd z_1\dd z_2 := \big(\mathcal{F}_{1,1}(z_1,z_2) + \mathcal{F}_{1,2}(z_1,z_2) + \mathcal{F}_{2,1}(z_1,z_2) + \mathcal{F}_{2,2}(z_1,z_2)\big)\dd z_1 \dd z_2
\]
admits an analytic continuation on $\Sigma_{\boldsymbol{v}} \times \Sigma_{\boldsymbol{v}}$ given by
\begin{equation}
\label{4gig4ugzbg}
\mathcal{B}_{\boldsymbol{v}} = \mathcal{B}^{\Sigma_{\boldsymbol{v}},\mathcal{L}} - \frac{\dd Z_1\dd Z_2}{(Z_1 - Z_2)^2}.
\end{equation}
The first term is the fundamental bidifferential on $\Sigma_{\boldsymbol{v}}$ normalized on the contour going counterclockwise around $[\alpha_2,\beta_2]$ in the $Z_{\boldsymbol{v}}$-plane, while $Z_1,Z_2$ refer to the function $Z := Z_{\boldsymbol{v}}$ on each of the two factors of $\Sigma_{\boldsymbol{v}} \times \Sigma_{\boldsymbol{v}}$. Since $\mathscr{V}^{\bot} = \textnormal{Ker}(\boldsymbol{\Theta}) = \amsmathbb{R}.(-1,1)$, Proposition~\ref{Lem:Holo1form} shows that the $1$-form
\[
\big(-\mathfrak{c}_{1;1}^{\textnormal{1st}}(z) + \mathfrak{c}_{1;2}^{\textnormal{1st}}(z) - \mathfrak{c}^{\textnormal{1st}}_{2;1}(z) + \mathfrak{c}^{\textnormal{1st}}_{2;2}(z)\big)\dd z
\]
analytically continues to a holomorphic $1$-form $\mathfrak{u}_{\boldsymbol{v};(-1,1)}$ such that
\begin{equation}
\label{normgubar} \frac{1}{2\ii\pi} \oint_{\gamma_2} \mathfrak{u}_{\boldsymbol{v};(1,-1)} = 1
\end{equation}
where $\gamma_2$ is a counterclockwise contour around $[\alpha_2,\beta_2]$ in the principal sheet of $\Sigma_{\boldsymbol{v}}$.

In the corresponding tiling model, only the sum of the filling fractions $\hat{n} := \hat{n}_1 + \hat{n}_2$ is fixed. Studying the corresponding discrete ensemble with both $\hat{n}_1,\hat{n}_2$ fixed gives an equilibrium measure $\mu^{\hat{n}_1,\hat{n}_2}$. According to Proposition~\ref{Proposition_differentiability_filling_fraction}, the first derivative of the Stieltjes transform $\mathcal{G}_{\boldsymbol{\mu}^{\hat{n}_1,\hat{n} - \hat{n}_1}}(z)$, with respect to $\hat{n}_1$ while keeping $\hat{n}$ fixed, is expressed in terms of the functions of the first kind
\[
\partial_{\hat{n}_1} \mathcal{G}_{\mu^{\hat{n}_1,\hat{n} - \hat{n}_1}}(z)\big|_{\hat{n}\,\,\textnormal{fixed}} = \mathfrak{c}^{\textnormal{1st}}_{1;1}(z) - \mathfrak{c}_{1;2}^{\textnormal{1st}}(z) + \mathfrak{c}^{\textnormal{1st}}_{2;1}(z) - \mathfrak{c}^{\textnormal{1st}}_{2;2}(z).
\]

There are several ways to present the elliptic curve $\Sigma_{\boldsymbol{v}}$, by equations or by identifying it with a complex torus $\amsmathbb{C}/(\amsmathbb{Z} \oplus \tau \amsmathbb{Z})$. By equations, we remark that we have two meromorphic functions $Y,Z$ on $\Sigma_{\boldsymbol{v}}$ satisfying
\[
Y^2 = (Z - \alpha_1)(Z - \beta_1)(Z - \alpha_2)(Z - \beta_2).
\]
Concretely, $Y$ in the physical sheet of $\Sigma_{\boldsymbol{v}}$ is the function $\sigma(z)$ already used many times in this book, namely the unique holomorphic function away from the two segments $[\alpha_1,\beta_1] \cup [\alpha_2,\beta_2]$ such that $\sigma^2(z) = (z - \alpha_1)(z - \beta_1)(z - \alpha_2)(z - \beta_2)$ and $\sigma(z) \sim z^2$ as $z \rightarrow \infty$. The fundamental bidifferential to be used in \eqref{4gig4ugzbg} was described in Section~\ref{sec:fund02}, from which we deduce
\[
\mathcal{B}_{\boldsymbol{v}} = \frac{\dd Z_1 \dd Z_2}{4Y_1Y_2}\left(c_{0} - (Z_1 + Z_2)(Z_1 + Z_2 - \alpha_1 - \beta_1 - \alpha_2 - \beta_2) + \frac{(Y_1 -Y_2)^2}{(Z_1 - Z_2)^2}\right).
\]
with the constant $c_0$ given in \eqref{4c0}. We have
\begin{equation}
\label{ubareqn}
\overline{u}_{\boldsymbol{v};(1,-1)} = \frac{2\pi}{\int_{\alpha_2}^{\beta_2} \frac{\dd z}{|\sigma(z)|}}\,\frac{\dd Z}{Y}
\end{equation}
because the right-hand side is the unique holomorphic $1$-form on $\Sigma_{\boldsymbol{v}}$ with the normalization property \eqref{normgubar}.

\medskip

\begin{figure}[h!]
\begin{center}
\includegraphics[width=0.4\textwidth]{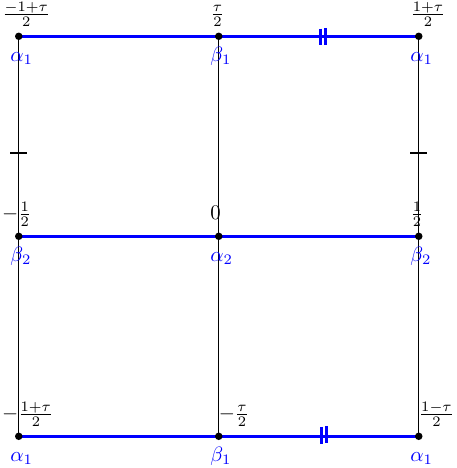}
\caption{\label{Fig:Complextoruscoord} The parameterization $\zeta$ (in black), with the corresponding values of $Z$ (in blue). The blue segments correspond to $Z^{-1}\big([\alpha_1,\beta_1] \cup [\alpha_2,\beta_2]\big)$, the remaining black segments correspond to the rest of $Z^{-1}(\amsmathbb{R})$. The bars indicate that the vertical (respectively, horizontal) sides of the rectangle that are identified in $\Sigma_{\boldsymbol{v}} \simeq \amsmathbb{C}/(\amsmathbb{Z} \oplus \tau\amsmathbb{Z})$.}
\end{center}
\end{figure}

The identification with a complex torus $\amsmathbb{T}_{\tau} = \amsmathbb{C}/(\amsmathbb{Z} \oplus \tau\amsmathbb{Z})$ can be done with the coordinate
\begin{equation}
\label{zetaZint2}
\zeta = \frac{\int_{\alpha_2}^{Z} \frac{\dd z}{\sigma(z)}}{-2\ii\int_{\alpha_2}^{\beta_2} \frac{\dd z}{|\sigma(z)|}},
\end{equation}
The possible choices of path of integration give $\zeta$ as a coordinate on the universal covering $\amsmathbb{C}$ of $\Sigma_{\boldsymbol{v}}$. Alternatively, the multivaluedness of $\zeta$ on $\Sigma_{\boldsymbol{v}}$ can be eliminated by computing it modulo the lattice $\amsmathbb{Z} \oplus \tau\amsmathbb{Z}$ with
\begin{equation}
\label{tauequation}
\tau = \ii \frac{\int_{\beta_1}^{\alpha_2} \frac{\dd z}{|\sigma(z)|}}{\int_{\alpha_2}^{\beta_2} \frac{\dd z}{|\sigma(z)|}},
\end{equation}
\textit{cf.} Figure~\ref{Fig:Complextoruscoord}. The normalization in \eqref{zetaZint2} was chosen so that the contour going counterclockwise around $[\alpha_2,\beta_2]$ is mapped to $[0,1]$ modulo lattice translation. In particular, we have
\[
\mathfrak{u}_{\boldsymbol{v};(1,-1)} = - 2\ii \pi \dd \zeta
\]
because both sides are holomorphic $1$-forms with the same period along $\gamma_2$.

With this presentation, $Z_{\boldsymbol{v}} := Z$ is an elliptic function with periods $1$ and $\tau$. It is even, has simple poles $\pm \zeta_{\infty} + \amsmathbb{Z} \oplus \tau \amsmathbb{Z}$ with residue $u_{\infty}$, where the parameters are given by
\begin{equation}
\label{Zpolareq}
\begin{split}
\zeta_{\infty} & := \frac{\int_{\alpha_2 + \ii \amsmathbb{R}_+} \frac{\dd z}{\sigma(z)}}{-2\ii \int_{\alpha_2}^{\beta_2} \frac{\dd z}{\sigma(z)}} \in \big(\tfrac{1}{2},\tfrac{1 + \tau}{2}\big), \\
u_{\infty} & := \Res_{\zeta = \zeta_{\infty}} Z(\zeta)\dd\zeta = - \Res_{\zeta = -\zeta_{\infty}} Z(\zeta)\dd \zeta = \frac{1}{2\ii \int_{\alpha_2}^{\beta_2} \frac{\dd z}{|\sigma(z)|}},
\end{split}
\end{equation}
and we have the special values
\begin{equation}
\label{Zspeci}
Z(0) = \alpha_2,\qquad Z\big(\pm \tfrac{1}{2}\big) = \beta_2,\qquad Z\big(\pm \tfrac{\tau}{2}\big) = \beta_1,\qquad Z\big(\pm \tfrac{1 + \tau}{2}\big) = \alpha_1.
\end{equation}
Both \eqref{Zpolareq} and \eqref{Zspeci} come from the normalization in \eqref{zetaZint2} together with the expression \eqref{tauequation} of the modulus $\tau$. For the first one, we have used
\[
\frac{\dd \zeta}{\dd Z} \,\,\mathop{\sim}_{\zeta \rightarrow \pm \zeta_{\infty}}\,\, \frac{1}{-2\ii Z^2\int_{\alpha_2}^{\beta_2} \frac{\dd z}{|\sigma(z)|}}.
\]
There are many equivalent ways to represent the function $Z$ characterized by these properties and we refer to \cite{Speci} for general principles and more information on elliptic functions. A first approach is via Weierstrass $\wp$-function. Even elliptic functions have two simple poles are of the form
\begin{equation}
\label{Zparac}
Z(\zeta) = \frac{c_{\infty}}{\wp(\zeta|\tau) - p_{\infty}} + c_{0}.
\end{equation}
and depend on four parameters $c_{0},p_{\infty} \in \amsmathbb{C}$, $c_{\infty} \in \amsmathbb{C}^*$ and $\tau$. The parameters $p_{\infty},c_{\infty}$ can be traded for $\zeta_{\infty},u_{\infty}$, since the comparison of \eqref{Zparac} with the announced behavior of $Z$ near its pole yields
\[
p_{\infty} = \wp(\zeta_{\infty}|\tau),\qquad c_{\infty} = u_{\infty}\wp'(\zeta_{\infty}|\tau).
\]
The four conditions \eqref{Zspeci} give parametric expressions for the endpoints $\alpha_1,\beta_1,\alpha_2,\beta_2$.  In the same style, the fundamental bidifferential to be used in \eqref{4gig4ugzbg} is the one normalized on the homology cycle $[0,1]$, and it can be expressed as
\[
\mathcal{B}^{\Sigma_{\boldsymbol{v}},[0,1]}(\zeta_1,\zeta_2) = \left(\wp(\zeta_1 - \zeta_2|\tau) + \frac{\pi^2E_2(\tau)}{3}\right)\dd\zeta_1\dd\zeta_2.
\]

A second approach is via ratios of the Jacobi theta function
\begin{equation}
\label{Jacobitheta11}
\vartheta_{11}(\zeta|\tau) = \sum_{m\in \amsmathbb{Z}} e^{\ii \pi \tau (m + \frac{1}{2})^2 + 2 \ii \pi (\zeta + \frac{1}{2})(m + \frac{1}{2})}.
\end{equation}
It is a entire function of $\zeta$ having a simple zero at $\zeta = 0$ and satisfying
\[
\vartheta_{11}(\zeta + 1 |\tau) = - \vartheta_{11}(\zeta|\tau),\qquad \vartheta_{11}(\zeta + \tau) = -e^{-2\ii\pi \zeta}\vartheta_{11}(\zeta|\tau).
\]
We then have
\begin{equation}
\label{Zzetatheta}
\begin{split}
Z(\zeta) - \beta_1 & = \frac{ u_{\infty}\,\vartheta_{11}'(0|\tau)\theta_{11}(2\zeta_{\infty}|\tau)}{\vartheta_{11}\big(\zeta_{\infty} - \frac{\tau}{2}\big|\tau\big)\vartheta_{11}\big(\zeta_{\infty} + \frac{\tau}{2}\big|\tau)}\, \frac{\vartheta_{11}\big(\zeta - \frac{\tau}{2}\big|\tau\big)\vartheta_{11}\big(\zeta + \frac{\tau}{2}\big|\tau\big)}{\vartheta_{11}(\zeta - \zeta_{\infty}|\tau)\theta_{11}(\zeta + \zeta_{\infty}|\tau)}, \\
\mathcal{B}^{\Sigma_{\boldsymbol{v}}}(\zeta_1,\zeta_2) & = \dd_{\zeta_1}\dd_{\zeta_2} \log \vartheta_{11}(\zeta_1 - \zeta_2).
\end{split}
\end{equation}
This first identity is proved by checking that the right-hand side is an elliptic function with the same zeros, and same simple poles with same residue at $\zeta = \zeta_{\infty}$ as the left-hand side. The second identity is also proved by checking that the right-hand side has all the properties (elliptic in both $\zeta_1$ and $\zeta_2$, double pole at $\zeta_1 = \zeta_2$ with biresidue $1$, zero period as $\zeta_1$ is integrated on $[0,1]$) characterizing the fundamental bidifferential.

A third approach is via Jacobi trigonometric functions. We do not develop it here but refer to \cite{Speci} for the basic techniques to do so.

\section{\texorpdfstring{$H = 3$}{H = 3}: S- and G-shaped domains}
\label{GSshapesec}
\begin{figure}[h!]
\begin{center}
\includegraphics[width=0.4\textwidth]{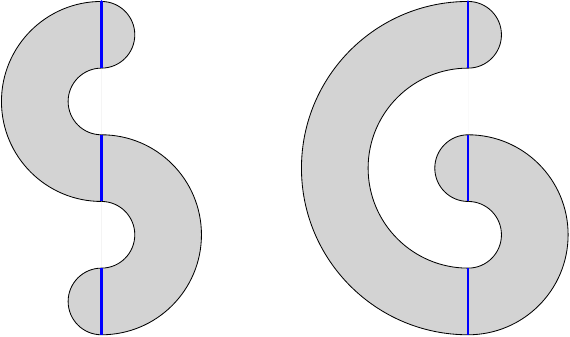}
\caption{\label{Fig:SGShaped} Topology of the domains with matrix of interactions \eqref{ExSGshape}.}
\end{center}
\end{figure}

The matrix of interactions for a S-shaped domain or an G-shaped domain obtained by gluing three trapezoids (\textit{cf.} Figure~\ref{Fig:SGShaped}) are respectively
\begin{equation}
\label{ExSGshape} \boldsymbol{\Theta} = \left(\begin{array}{ccc} 1 & \frac{1}{2} & 0 \\[0.6ex] \frac{1}{2} & 1 & \frac{1}{2} \\[0.6ex] 0 & \frac{1}{2} & 1 \end{array}\right),\qquad \boldsymbol{\Theta} = \left(\begin{array}{ccc} 1 & 0 & \frac{1}{2} \\[0.6ex] 0 & 1 & \frac{1}{2} \\[0.6ex] \frac{1}{2} & \frac{1}{2} & 1 \end{array}\right).
\end{equation}
They differ by exchange of the second and the third segment. Therefore the discussion for the two domains are obtained from each other by exchanging $(\alpha_2,\beta_2)$ with $(\alpha_3,\beta_3)$.

Let us treat the S-shaped domain. The generators are
\[
T^{(1)} = \left(\begin{array}{ccc} -1 & 0 & 0 \\ -1 & 1 & 0 \\ 0 & 0 & 1 \end{array}\right),\qquad T^{(2)} = \left(\begin{array}{ccc} 1 & -1 & 0 \\ 0 & -1 & 0 \\ 0 & -1 & 1 \end{array}\right),\qquad T^{(3)} = \left(\begin{array}{ccc} 1& 0 & 0 \\ 0 & 1 & -1 \\ 0 & 0 & -1 \end{array}\right).
\]
We can check that they satisfy the relations $(T^{(1)}T^{(2)})^3 = (T^{(2)}T^{(3)})^3 = (T^{(1)}T^{(3)})^2 = 1$, confirming the result of Theorem~\ref{thm:Omegav} (observe that $\boldsymbol{\Theta}$ has full rank) saying that $\mathfrak{G}$ is the Weyl group of type $A_3$, \textit{i.e} the symmetric group in four elements. We have a minimal orbit generated by $\boldsymbol{v} := \boldsymbol{e}^{(1)} = (1,0,0)$
\begin{equation}
\label{start1}
(1,0,0)\,\, \mathop{\longleftrightarrow}^{T^{(1)}}\,\, (-1,-1,0) \,\,\mathop{\longleftrightarrow}^{T^{(2)}}\,\,(0,1,1) \,\,\mathop{\longleftrightarrow}^{T^{(3)}}\,\, (0,0,-1).
\end{equation}

\subsection{Rational parameterization of the spectral curve}

The spectral curve $\Sigma_{\boldsymbol{v}}$ has genus zero and the function $Z := Z_{\boldsymbol{v}} : \Sigma_{\boldsymbol{v}} \rightarrow \widehat{\amsmathbb{C}}$ has four poles. There is a unique uniformizing coordinate $\zeta \in \widehat{\amsmathbb{C}} \simeq \Sigma_{\boldsymbol{v}}$ such that $Z$ has poles at $0,1,\infty$ and some $t \neq 0,1,\infty$. This gives a partial fraction expansion of the form
\begin{equation}
\label{ZzetaShape}
Z(\zeta) = \frac{u_0}{\zeta} + \frac{u_1}{\zeta - 1} + \frac{u_t}{\zeta - t} + u' + u_{\infty}\zeta.
\end{equation}
for some parameters $u_0,u_1,u_t,u_{\infty} \in \amsmathbb{C}^*$ and $u' \in \amsmathbb{C}$. We can arrange the poles so that $\zeta = \infty$ is in the sheet labeled $(1,0,0)$, $\zeta = t$ in the sheet labeled $(-1,-1,0)$, $\zeta = 1$ in the sheet labeled $(0,1,1)$, and $\zeta = 0$ in the sheet labeled $(0,0,-1)$. Parameterizing $Z'(\zeta)$ in terms of its six zeros $s_1,\ldots,s_6 \in \amsmathbb{C} \setminus \{0,1,t\}$ and $u_{\infty}$, we get
\begin{equation}
\label{ZderSh}
 Z'(\zeta) = \frac{u_{\infty}\prod_{j = 1}^6 (\zeta - s_j)}{\zeta^2(\zeta - 1)^2(\zeta - t)^2}.
\end{equation}
Comparing the behavior at the poles with \eqref{ZzetaShape}, we have
\[
u_0 = \frac{u_{\infty}}{t^{2}}\prod_{j = 1}^{6} s_j ,\qquad u_1 = \frac{u_{\infty}}{(1 - t)^{2}} \prod_{j = 1}^{6} (s_j - 1),\qquad u_t = \frac{u_{\infty}}{t^2(1 - t)^2} \prod_{j = 1}^{6} (t - s_j).
\]
which allow the elimination of the parameters $u_0,u_1,u_t$. The vanishing of residues of $Z'(\zeta)$ at $0,1,\infty$ yields three relations
\[
\frac{1}{2} \sum_{j = 1}^{6} \frac{1}{s_j} = 1 + \frac{1}{t},\qquad \frac{1}{2}\sum_{j = 1}^{6} \frac{1}{1 - s_j} = \frac{t- 2}{t - 1},\qquad \frac{1}{2}\sum_{j = 1}^{6} \frac{1}{1 - s_j/t} = \frac{2t - 1}{t - 1},
\]
for the nine parameters $u',u_{\infty},t,s_1,\ldots,s_6$, and writing
\[
\{Z(s_j)\,\,|\,\,j \in [6]\} = \{\alpha_1,\beta_1,\alpha_2,\beta_2,\alpha_3,\beta_3\}.
\]
gives a parametric expression for these six endpoints.

\subsection{The fundamental solution and how to use it}
\label{Csha02exp}
Like in Section~\ref{sec:Csha}, Theorem~\ref{thm:Omegav} and the uniqueness of the fundamental bidifferential on the Riemann sphere gives
\begin{equation}
\label{O02GSsh} \begin{split}
\mathcal{B}_{\boldsymbol{v}} (\zeta_1,\zeta_2) & = \dd_{\zeta_1}\dd_{\zeta_2}\log\left(\frac{\zeta_1 - \zeta_2}{Z(\zeta_1) - Z(\zeta_2)}\right) \\
& = -\dd_{\zeta_1} \dd_{\zeta_2} \log\left(1 - \frac{u_0}{u_{\infty}\zeta_1\zeta_2} - \frac{u_1}{u_{\infty}(\zeta_1 - 1)(\zeta_2 - 1)} - \frac{u_t}{u_{\infty}(\zeta_1 - t)(\zeta_2 - t)} \right)
\end{split}
\end{equation}
as analytic continuation of the fundamental solution $\mathcal{F}_{1,1}(z_1,z_2)\dd z_1 \dd z_2$ to $\Sigma_{\boldsymbol{v}} \times \Sigma_{\boldsymbol{v}}$. This formula differs from Section~\ref{sec:Csha} only by a more complicated expression for the uniformization map $Z(\zeta)$. The structure of the orbit \eqref{start1} dictates the expression of $\mathcal{B}_{\boldsymbol{v}}$ in terms of the fundamental solution $\boldsymbol{\mathcal{F}}(z_1,z_2)$ when its two variables are in various sheets. The reader can easily construct the analog of Figure~\ref{fig:doublesheetBC} for this case, and deduce formulae for the leading covariance of linear statistics of various groups of particles in discrete ensembles with matrix of interactions \eqref{ExSGshape}. The logic is the same as in Section~\ref{sec:Csha}. If $f$ is a polynomial test function, we have for instance as $\N \rightarrow \infty$
\begin{equation*}
\begin{split}
 \sum_{i,j = 1}^{N_1} \amsmathbb{E}^{(\textnormal{c})}\big[f(\ell^1_i)f(\ell^1_j)\big] + o(1) & = \Res_{\zeta_1 = \infty} \Res_{\zeta_2 = \infty} f(Z(\zeta_1)) f(Z(\zeta_2)) \mathcal{B}_{\boldsymbol{v}}(\zeta_1,\zeta_2) \\
 \sum_{i,j = 1}^{N_2} \amsmathbb{E}^{(\textnormal{c})}\big[f(\ell^2_i)f(\ell^2_j)\big] + o(1) & = \big(\Res_{\zeta_1 = 1} + \Res_{\zeta_1 = 0}\big)\big(\Res_{\zeta_2 = 1} + \Res_{\zeta_2 = 0}\big) f(Z(\zeta_1)) f(Z(\zeta_2)) \mathcal{B}_{\boldsymbol{v}}(\zeta_1,\zeta_2), \\
\sum_{i,j = 1}^{N_3} \amsmathbb{E}^{(\textnormal{c})}\big[f(\ell^3_i)f(\ell^3_j)\big] + o(1) & = \Res_{\zeta_1 = 0} \Res_{\zeta_2 = 0} f(Z(\zeta_1)) f(Z(\zeta_2)) \mathcal{B}_{\boldsymbol{v}}(\zeta_1,\zeta_2).
\end{split}
\end{equation*}
In the second identity we used the fact that $\boldsymbol{e}^{(2)} = (0,1,1) + (0,0,-1)$ to access covariance of linear statistics for particles in the second segment via a sum of residues at the poles of $Z$ in the sheet labeled $(0,1,1)$ and the sheet labeled $(0,0,-1)$.

\section{\texorpdfstring{$H = 3$}{H = 3}: E-shaped domain}
\label{H3Esec}
\begin{figure}[h!]
\begin{center}
\includegraphics[width=0.14\textwidth]{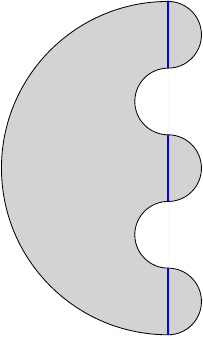}
\caption{\label{Fig:Eshaped} Topology of the domain with matrix of interactions \eqref{ExEshaped}.}
\end{center}
\end{figure}

The matrix of interactions for a E-shaped domain obtained by gluing four trapezoids (\textit{cf.} Figure~\ref{Fig:Eshaped}) is
\begin{equation}
\label{ExEshaped}
\boldsymbol{\Theta} = \left(\begin{array}{ccc} 1 & \frac{1}{2} & \frac{1}{2} \\[0.6ex] \frac{1}{2} & 1 & \frac{1}{2} \\[0.6ex] \frac{1}{2} & \frac{1}{2} & 1 \end{array}\right).
\end{equation}
It has full rank and by Theorem~\ref{thm:Omegav}, $\mathfrak{G}$ is the Weyl group of type $A_3$ (the symmetric group in four elements) and we have a minimal orbit
\begin{equation}
\label{start2}
\boldsymbol{v} := (1,1,1) \,\, \mathop{\longleftrightarrow}^{T^{(h)}}\,\, - \boldsymbol{e}^{(h)}\,\, \quad h = 1,2,3.
\end{equation}
The spectral curve $\Sigma_{\boldsymbol{v}}$ has again four sheets, six ramification points, and genus $0$. The function $Z_{\boldsymbol{v}}$ has four poles and can be parameterized exactly in the same way as in Section~\ref{GSshapesec}. The analytic continuation of
\[
\mathcal{F}(z_1,z_2) \dd z_1 \dd z_2 = \sum_{h_1,h_2 = 1}^{3} \mathcal{F}_{h_1,h_2}(z_1,z_2) \dd z_1 \dd z_2
\]
to $\Sigma_{\boldsymbol{v}} \times \Sigma_{\boldsymbol{v}}$ is given again by \eqref{O02GSsh}. Its values in the other sheets are indicated in Figure~\ref{fig:Eshapedsheet}.

\begin{figure}[h!]
\[
\begin{array}{|c|c|c|c|c|}
\hline
\textnormal{sheet} & (1,1,1) & (-1,0,0) & (0,-1,0) & (0,0,-1) \\[4pt]
\hline
(1,1,1) & \sum_{h_1,h_2 = 1}^{3} \mathcal{F}_{h_1,h_2} + \frac{1}{(Z_1 - Z_2)^2} & - \sum_{h_2 = 1}^{3} \mathcal{F}_{1,h_2} & - \sum_{h_2 = 1}^{3} \mathcal{F}_{2,h_2} & - \sum_{h_2 = 1}^{3} \mathcal{F}_{3,h_2} \\[4pt]
\hline
(-1,0,0) & -\sum_{h_1 = 1}^{3} \mathcal{F}_{h_1,1} & \mathcal{F}_{1,1}+ \frac{1}{(Z_1 - Z_2)^2} & \mathcal{F}_{2,1} & \mathcal{F}_{3,1} \\[4pt]
\hline
(0,-1,0) & - \sum_{h_1 = 1}^{3}\mathcal{F}_{h_1,2} & \mathcal{F}_{1,2} & \mathcal{F}_{2,2} + \frac{1}{(Z_1 - Z_2)^2}& \mathcal{F}_{3,2} \\[4pt]
\hline
(0,0,-1) & -\sum_{h_1 = 1}^{3} \mathcal{F}_{h_1,3} & \mathcal{F}_{1,3} & \mathcal{F}_{2,3} & \mathcal{F}_{3,3} + \frac{1}{(Z_1 - Z_2)^2}\\[4pt]
\hline
\end{array}
\]
\caption{\label{fig:Eshapedsheet} The lines (respectively, columns) heads indicate the sheet to which the first (respectively, second) variable belongs to, and the entries are the expressions for the function $\frac{\mathcal{B}_{\boldsymbol{v}}}{\dd Z_1\dd Z_2} + \frac{1}{(Z_1 - Z_2)^2}$.}
\end{figure}

The structure of the branched covering $Z_{\boldsymbol{v}}$ is however not the same as in Section~\ref{GSshapesec} (compare \eqref{start1} and \eqref{start2}). This is reflected in the branches of the algebraic functions to be chosen to extract the parameters of $Z_{\boldsymbol{v}}$ as functions of the six branch points $\alpha_1,\beta_1,\alpha_2,\beta_2,\alpha_3,\beta_3$. For the E-shaped domain one should choose the branch so that there is one sheet containing the six ramification points $s_1,\ldots,s_6$. In any case, as functions of $\alpha_1,\beta_1,\alpha_2,\beta_2,\alpha_3,\beta_3$, the expressions for the G, S and E-shaped domain related to each other by analytic continuation.

\section{\texorpdfstring{$H = 3$}{H = 3}: 6-shaped domain}
\label{sec:6yann}
The interaction matrix for the $6$-shaped domain obtained by gluing three trapezoids (\textit{cf.} Figure~\ref{Fig:6shaped}) is
\begin{equation}
\label{Ex6shaped}
\boldsymbol{\Theta} = \left(\begin{array}{ccc} 1 & \frac{1}{2} & \frac{1}{2} \\[0.6ex] \frac{1}{2} & 1 & 1 \\[0.6ex] \frac{1}{2} & 1 & 1 \end{array}\right).
\end{equation}
The generators are
\[
T^{(1)} = \left(\begin{array}{ccc} -1 & 0 & 0 \\ -1 & 1 & 0 \\ -1 & 0 & 1 \end{array}\right),\qquad T^{(2)} = \left(\begin{array}{ccc} 1 & -1 & 0 \\ 0 & -1 & 0 \\ 0 & -2 & 1 \end{array}\right),\qquad T^{(3)} = \left(\begin{array}{ccc} 1 & 0 & -1 \\ 0 & 1 & -2 \\ 0 & 0 & -1\end{array}\right).
\]
and they generate an infinite group. Indeed, $T^{(2)}$ and $T^{(3)}$ has an upper-triangular block form if we decompose in blocks of size $1$ and $2$, and their lower right block coincide with the generators in Section~\ref{sec:Oshaped} that were generating the free group in two elements.

As $\textnormal{rank}(\boldsymbol{\Theta}) = 2$, Theorem~\ref{thm:Omegav} says that the group $\hat{\mathfrak{G}}$ is the Weyl group of type $A_2$, \textit{i.e.} the symmetric group in three elements, and gives a minimal orbit of size $3$
\[
\boldsymbol{v} := (1,0,0) \,\,\mathop{\longleftrightarrow}^{T^{(1)}} \,\,(-1,-1,-1)\,\,\mathop{\longleftrightarrow}^{T^{(2)},T^{(3)}} \,\,(0,1,1).
\]

\subsection{Elliptic parameterization of the spectral curve}
\label{Ell6s}
The spectral curve $\Sigma_{\boldsymbol{v}}$ has three sheets and genus $1$. It is biholomorphic to the complex torus $\amsmathbb{T}_{\tau} := \amsmathbb{C}/(\amsmathbb{Z} \oplus \tau \amsmathbb{Z})$, and the biholomorphic map is unique up to postcomposition with a translation. We call $\zeta$ the natural coordinate in the complex torus. We look for an expression of the branched covering $Z_{\boldsymbol{v}} := Z : \Sigma_{\boldsymbol{v}} \rightarrow \widehat{\amsmathbb{C}}$. This is a meromorphic function on $\Sigma_{\boldsymbol{v}}$ having three simple poles $p_0,p_-,p_+$. We label them so that $p_0$ is in the sheet $(1,0,0)$, $p_-$ in the sheet $(-1,-1,-1)$ and $p_+$ in the sheet labeled $(0,1,1)$. Using the translation freedom, we can choose the above biholomorphism such that $\zeta = 0$ corresponds to the pole $p_0$, and we call $\zeta_{\pm} \in \amsmathbb{T}_{\tau}$ the points corresponding to $p_{\pm}$. In terms of the Weierstra\ss{} $\wp$-function, the general form of a meromorphic elliptic function with simple poles at $0,\zeta_{\pm}$ reads
\begin{equation}
\label{Zzetacccc}
Z(\zeta) = -\frac{c_+}{2} \frac{\wp'(\zeta|\tau) + \wp'(\zeta_+|\tau)}{\wp(\zeta|\tau) - \wp(\zeta_+|\tau)} - \frac{c_-}{2} \frac{\wp'(\zeta|\tau) + \wp'(\zeta_-|\tau)}{\wp(\zeta|\tau) - \wp(\zeta_-|\tau)} + c_0.
\end{equation}
for some constants $c_{\pm} \in \amsmathbb{C}^*$ and $c_0 \in \amsmathbb{C}$. The coefficients $c_{\pm}$ are the residues of $Z(\zeta)\dd \zeta$ at $\zeta_{\pm}$. Since $\wp(\zeta) \sim \frac{1}{\zeta^2}$ as $\zeta \rightarrow 0$, the residue at $0$ is $-(c_+ + c_-)$, in agreement with the fact that the sum of residues of a meromorphic function on $\amsmathbb{T}_{\tau}$ vanishes. The derivative $Z'(\zeta)$ is an elliptic function with three double poles with known leading coefficients, no residues, and zero periods along $[0,1]$ and $[0,\tau]$. There is an unique function with such properties and it can be exhibited as
\begin{equation}
\label{Zprimecccc}
\begin{split}
Z'(\zeta) & = (c_+ + c_-) \wp(\zeta|\tau) - c_+ \wp(\zeta - \zeta_+|\tau) - c_- \wp(\zeta - \zeta_-|\tau) \\
& = c_+ \left[\wp(\zeta_+|\tau) - \frac{1}{4} \bigg(\frac{\wp'(\zeta|\tau) + \wp'(\zeta_+|\tau)}{\wp(\zeta|\tau) - \wp(\zeta_+|\tau)}\bigg)^2\right] + c_- \left[\wp(\zeta_-|\tau) - \frac{1}{4} \bigg(\frac{\wp'(\zeta|\tau) + \wp'(\zeta_-|\tau)}{\wp(\zeta|\tau) - \wp(\zeta_-|\tau)}\bigg)^2\right].
\end{split}
\end{equation}
This expression is more convenient for computations than taking naively the derivative of \eqref{Zzetacccc}. The second equality in \eqref{Zprimecccc} comes from the addition formula, \textit{cf.} \eqref{additionwp}. An elliptic function admits as many poles as zeros (counted with multiplicity). Therefore, $Z'(\zeta)$ has $6$ zeros $s_1,\ldots,s_6 \in \amsmathbb{T}_{\tau}$. They correspond to the six (as expected) ramification points of $Z$. Writing
\begin{equation}
\label{alpha6j}
\{\alpha_1,\beta_1,\alpha_2,\beta_2,\alpha_3,\beta_3\} = \big\{Z(s_j)\,\,|\,\,j \in [6]\big\}
\end{equation}
gives six constraints for the branch points in terms of the six parameters $c_0,c_{\pm},\zeta_{\pm},\tau$. They however involve the auxiliary parameters $s_1,\ldots,s_6$ which are not easy to access.

\subsection{Algebraic presentation of the spectral curve}
\label{1452Sec}
There is an alternative way to obtain an algebraic equation for $Z$ without resorting to a parameterization of $\Sigma_{\boldsymbol{v}}$, which will give a simpler way to get formulae for the branch points. It relies on the fact (consequence of the Riemann--Roch theorem for elliptic curves) that the space $\mathscr{L}(d_0,d_+,d_-)$ of meromorphic functions on $\Sigma_{\boldsymbol{v}}$ that can have poles at $p_0,p_{\pm}$ with respective orders at most $d_0,d_{\pm}$, has exactly dimension $d_0 + d_+ + d_-$. Then, $\mathscr{L}(0,1,0)$ has dimension $1$ and contains only the constant functions. The space $\mathscr{L}(0,1,1)$ (respectively $\mathscr{L}(1,1,0)$) has dimension two, so contains a non-constant function $\xi$ (respectively $\eta$). Clearly, $1,\xi,\eta$ must be linearly independent. Then, $\mathscr{L}(1,1,1)$ contains $1,\xi,\eta,Z$ and has dimension $3$, so $Z$ must be a linear combination of $1,\xi,\eta$. By dimension matching arguments, we find bases
\begin{equation*}
\begin{split}
\mathscr{L}(1,2,1) & = \textnormal{span}_{\amsmathbb{C}}(1,\xi,\eta,\xi\eta), \\
\mathscr{L}(1,2,2) & = \textnormal{span}_{\amsmathbb{C}}(1,\xi,\eta,\xi\eta,\xi^2), \\
\mathscr{L}(2,2,2) & = \textnormal{span}_{\amsmathbb{C}}(1,\xi,\eta,\xi\eta,\xi^2,\eta^2) \\
\mathscr{L}(2,3,2) & = \textnormal{span}_{\amsmathbb{C}}(1,\xi,\eta,\xi\eta,\xi^2,\eta^2,\xi^2\eta)
\end{split}
\end{equation*}
Then, the space $\mathscr{L}(2,3,3)$ contains the seven basis elements of $\mathscr{L}(2,3,2)$ as well as $\xi \eta^2$ and $\xi^3$, but is only of dimension $8$, so we must have a linear relation between these nine elements. Eliminating $\eta$ to keep only $\xi$ and $Z$, we get an algebraic equation with same type of monomials. Exploiting the freedom of changing $\xi$ by an affine transformation and dividing the whole equation by a non-zero constant, we can bring it to the form
\begin{equation}
\label{Pxieta}
P(\xi,Z) := \xi Z^2 + sZ^2 + t_0Z + t_1\xi Z + t_2\xi^2 Z + u_0 + u_1 \xi + \xi^3 = 0
\end{equation}
This is an algebraic equation for $\Sigma_{\boldsymbol{v}}$ depending on six parameters $s,t_0,t_1,t_2,u_0,u_1$. For generic parameters, the ramification points of $Z$ correspond to the common solutions of $P(\xi,Z) = 0$, $\partial_{\xi} P(\xi,Z) = 0$ and $\partial_{Z} P(\xi,Z) \neq 0$. From $\partial_1 P(\xi,Z)$ we can extract $\xi$ as a rational function with numerator of degree $3$ and denominator of degree $2$ in $Z$. Inserting back in \eqref{Pxieta} yields an equation of degree $6$ for $Z$, namely
\begin{equation*}
\begin{split}
0 & = (4 - t_2^2) Z^6 + 2(6t_1 - 9st_2 - t_1t_2^2 + 2s t_2^3)Z^5 \\
& \quad + (27s^2 + 12t_1^2 - 18t_0t_2 - 18st_1t_2 - t_1^2t_2^2 + 4t_0t_2^3 + 12u_1 - 2t_2^2u_1)Z^4 \\
& \quad + (27st_0 + 2t_1^3 - 9t_0t_1t_2 - 9t_2u_0 + 2t_2^3u_0 + 12t_1u_1 - 9st_2u_1 - t_1t_2^2u_1)Z^3 \\
& \quad + (27t_0^2 + 54su_0 - 18t_1t_2u_0 + 12t_1^2u_1 - 18t_0t_2u_1 + 12u_1^2 - t_2^2u_1^2)Z^2 \\
& \quad + 2(27t_0u_0 - 9t_2u_0u_1 + 6t_1u_1^2)Z + 27u_0^2 + 4u_1^3.
\end{split}
\end{equation*}
Imposing that the six roots of this equation are $\alpha_1,\beta_1,\alpha_2,\beta_2,\alpha_3,\beta_3$ give a rational parameterization for them in terms of $s,t_0,t_1,t_2,u_0,u_1$. Equivalently, this gives an explicit algebraic system to solve in order to obtain the latter (appearing in the equation for $\Sigma_{\boldsymbol{v}}$) in terms of the former.

It is possible to relate directly this presentation with the first one we have obtained via the parameterization $\zeta$ and Weierstra\ss{} functions. Indeed, by construction $\xi$ is a meromorphic function on $\Sigma_{\boldsymbol{v}}$ with simple poles at $p_{\pm}$, and inspection of \eqref{Pxieta} reveals that it is normalized such that $\xi \sim \frac{\delta_{\pm}c_{\pm}}{\zeta - \zeta_{\pm}}$ near $p_{\pm}$ with
\[
\delta_{\pm} = \frac{-t_2 \pm \sqrt{t_2^2 - 4}}{2}.
\]
Since the sum of residues of $\xi \dd \zeta$ should vanish, we deduce that
\[
c_+ \delta_+ + c_- \delta_- = 0.
\]
Besides, at the remaining pole of $Z$ (namely, $p_0$), \eqref{Pxieta} shows that $\xi \sim -s$. By matching this special value and the behavior at poles we find
\[
\xi = \frac{\alpha_+ \delta_+}{2} \left(\frac{\wp'(\zeta|\tau) + \wp'(\zeta_+|\tau)}{\wp(\zeta|\tau) - \wp(\zeta_+|\tau)} - \frac{\wp'(\zeta|\tau) + \wp'(\zeta_-|\tau)}{\wp(\zeta|\tau) - \wp(\zeta_-|\tau)}\right) - s.
\]

\subsection{Elliptic parameterization of the fundamental solution}

We turn to the fundamental solution of the master Riemann--Hilbert problem. By Theorem~\ref{thm:Omegav}, the bidifferential
\begin{equation}
\label{Omegav6sha}
\mathcal{B}_{\boldsymbol{v}}(\zeta_1,\zeta_2) = \mathcal{B}^{\Sigma_{\boldsymbol{v}},\mathcal{L}}(\zeta_1,\zeta_2) - \frac{\dd Z(\zeta_1)\dd Z(\zeta_2)}{(Z(\zeta_1) - Z(\zeta_2))^2}
\end{equation}
is the analytic continuation of $\mathcal{F}_{1,1}(z_1,z_2)\dd z_1 \dd z_2$ to $\Sigma_{\boldsymbol{v}} \times \Sigma_{\boldsymbol{v}}$. Its values in the other sheets in terms of the fundamental solution $\boldsymbol{\mathcal{F}}$ are indicated in Figure~\ref{fig:doublesheetB6}. The first term in \eqref{Omegav6sha} is the fundamental bidifferential on $\Sigma_{\boldsymbol{v}}$ normalized on the cycle going counterclockwise around $[\alpha_2,\beta_2]$ in the sheet $(0,1,1)$. In the notation, the exponent $\mathcal{L}$ refers to the homology class of this cycle. Finding explicitly this homology class would require a more involved analysis of the structure of the branched covering and we do not enter this matter. It only affects the constant part in the expressions for the fundamental bidifferentials given in Section~\ref{sec:fund02} by a discrete ambiguity. In any case, in terms of the coordinate $\zeta$, we have from \eqref{OmeX1X2X1}
\begin{equation*}
\begin{split}
\mathcal{B}_{\boldsymbol{v}}(\zeta_1,\zeta_2) & = \big(\wp(\zeta_1 - \zeta_2|\tau) + c\big)\dd \zeta_1\dd \zeta_2 - \frac{\dd Z(\zeta_1)\dd Z(\zeta_2)}{(Z(\zeta_1) - Z(\zeta_2))^2} \\
& = \dd_{\zeta_1}\dd_{\zeta_2} \log\left(\frac{\vartheta_{11}(\zeta_1 - \zeta_2|\tau)}{Z(\zeta_1) - Z(\zeta_2)}\right) + c'\dd\zeta_1 \dd \zeta_2.
\end{split}
\end{equation*}
for some constants $c,c' \in \amsmathbb{C}$. Given the parameterization $Z(\zeta)$ of \eqref{Zzetacccc} and the description in Figure~\ref{fig:doublesheetB6} of the values taken by $\mathcal{B}_{\boldsymbol{v}}$ in all the sheets of $\Sigma_{\boldsymbol{v}}$, it is possible to carry out covariance computations in the style of Sections~\ref{Csun02ex}-\ref{Csha02exp}. The result involve residues at $\zeta = 0$ and $\zeta = p_{\pm}$ and can be expressed in terms of the (known) series expansion of Weierstra\ss{} $\wp$-functions or theta functions (instead of series expansion of rational functions).

\subsection{Algebraic presentation of the fundamental solution}
\label{thebingfsgun}
If we would present the elliptic curve $\Sigma_{\boldsymbol{v}}$ via an algebraic equation
\begin{equation}
\label{Weier4}
Y^2 = (X - \lambda_1)(X - \lambda_2)(X - \lambda_3)(X - \lambda_4)
\end{equation}
we get from \eqref{OmeX1X2X1} the formula
\begin{equation}
\label{fundsolfgrge}
\mathcal{B}_{\boldsymbol{v}} = \frac{\dd X_1\dd X_2}{4Y_1Y_2}\left(c'' - (X_1 + X_2)\Big(X_1 + X_2 - \sum_{i = 1}^{4} \lambda_i\Big) + \frac{(Y_1 + Y_2)^2}{(X_1 - X_2)^2}\right) - \frac{\dd Z_1\dd Z_2}{(Z_1 - Z_2)^2}
\end{equation}
for some constant $c'' \in \amsmathbb{C}$. The equation \eqref{Pxieta} for $Z$ can be brought in this form by setting
\[
Y = \frac{2(\xi + s)Z + t_0 + t_1 \xi + t_2\xi^2}{\sqrt{t_2^2 - 4}} = \frac{\partial_{Z} P(\xi,Z)}{\sqrt{t_2^2 - 4}},\qquad X = \xi
\]
We then have $\sum_{i = 1}^{4} \lambda_i = \frac{2s - t_1t_2}{4 - t_2^2}$. The relation $\partial_{\xi}P(\xi,Z) \dd \xi + \partial_{Z} P(\xi,Z) \dd Z = 0$ between meromorphic differentials on $\Sigma_{\boldsymbol{v}}$ allows rewriting
\[
\frac{\dd \xi}{Y} = - \frac{\sqrt{t_2^2 - 4}\,\dd Z}{Z^2 \Psi},\qquad \Psi := Z^{-2} \partial_{\xi}P(\xi,Z).
\]
and therefore
\begin{equation}
\label{Omegaunuggeggga}\mathcal{B}_{\boldsymbol{v}} = \frac{(t_2^2 - 4)\dd Z_1\dd Z_2}{4Z_1^2Z_2^2\Psi_1\Psi_2}\left(c'' - (X_1 + X_2)\Big(X_1 + X_2 - \sum_{i = 1}^{4} \lambda_i\Big) + \frac{(Y_1 + Y_2)^2}{(X_1 - X_2)^2}\right) - \frac{\dd Z_1\dd Z_2}{(Z_1 - Z_2)^2}.
\end{equation}

Examples of covariance computation would be
\begin{equation}
\label{Eccqun}\begin{split}
\amsmathbb{E}^{(\textnormal{c})}\big[f_1(\ell_i^1)f_2(\ell_j^1)\big] + o(1) & = \Res_{p_0} \Res_{p_0} f(Z_1)f(Z_2) \mathcal{B}_{\boldsymbol{v}} \\
\amsmathbb{E}^{(\textnormal{c})}\left[\Big(\sum_{i = 1}^{N_1} f_1(\ell_i^1)\Big)\Big(\sum_{j = 1}^{N_2} f_2(\ell_j^2) + \sum_{j = 1}^{N_3} f_2(\ell_j^3)\Big)\right] + o(1) & = \Res_{p_0} \Res_{p_+} f_1(Z_1)f_2(Z_2) \mathcal{B}_{\boldsymbol{v}}
\end{split}
\end{equation}
for polynomials $f_1,f_2$. If $f_l(Z)$ is a monomial $Z^{k_l}$ for $l = 1,2$, this is extracting the coefficient of $Z_1^{-(k_1 + 1)}Z_2^{-(k_2 + 1)}$ in the series expansion of $\frac{\mathcal{B}_{\boldsymbol{v}}}{\dd Z_1\dd Z_2}$ near $(p_0,p_0)$, respectively $(p_0,p_+)$. From the equation of the curve \eqref{Pxieta}, it is easy to compute order by order the Laurent series expansions of $\xi$ (and then of $Y$ and $\Psi$) near the poles $p_0,p_{\pm}$, in the variable $\frac{1}{Z}$. They take the form
\begin{equation}
\left\{\begin{array}{l} \textnormal{Near}\,\,p_0 : \\[4pt]
X = - s + \frac{t_0 + st_1 - s^2t_2}{Z}\big(1 + X^{[0]}\big) \\[4pt]
Y = \frac{t_0 + st_1 - s^2t_2}{\sqrt{t_2^2 - 4}}\big(1 + Y^{[0]}\big) \\[6pt]
\Psi = 1 + \Psi^{[0]}
\end{array}\right. \qquad \left\{\begin{array}{l} \textnormal{Near}\,\,p_{\pm} : \\ [4pt] X = \frac{\delta_{\pm}}{Z}\big(1 + X^{[\pm]}\big) \\[4pt] Y = -\delta_{\pm}^2Z^2\big(1 + Y^{[\pm]}\big) \\[4pt] \Psi = \delta_{\pm}\sqrt{t_2^2 - 4}\big(1 + \Psi^{[\pm]}\big) \end{array} \right.
\end{equation}
where
\[
\forall \varepsilon \in \{0,\pm\}\qquad X^{[\varepsilon]},Y^{[\varepsilon]},\Psi^{[\varepsilon]} = O(Z^{-1})
\]
are power series in the variable $\frac{1}{Z}$ that are explicitly computable order by order. One deduces without difficulty the series expansion of $\mathcal{B}_{\boldsymbol{v}}$ near $(p_{\varepsilon},p_{\varepsilon}')$ for $\varepsilon,\varepsilon' \in \{0,\pm\}$ in the variables $\frac{1}{Z_1}$ and $\frac{1}{Z_2}$ needed in \eqref{Eccqun}, in the case $\varepsilon \neq \varepsilon'$. In the case $\varepsilon = \varepsilon'$ the apparent double poles in the second line of \eqref{Omegaunuggeggga} seems to complicate the task. If we naively expanded each term, we would need to choose an ordering, say $|Z_1| < |Z_2|$, and we would obtain monomials involving large positive powers of $Z_1$, which is not desired. We however know that $\mathcal{B}_{\boldsymbol{v}}$ is regular at coinciding points in $\Sigma_{\boldsymbol{v}} \times \Sigma_{\boldsymbol{v}}$. Rewriting the formula to make this regularity apparent will facilitate the practical computation of expansions. A formula equivalent to \eqref{Omegaunuggeggga} and adapted for the expansion near $(p_0,p_0) \in \Sigma_{\boldsymbol{v}} \times \Sigma_{\boldsymbol{v}}$ is
\begin{equation}
\begin{split}
\mathcal{B}_{\boldsymbol{v}} & = \frac{\dd Z_1\dd Z_2}{4\Psi_1\Psi_2}\big((t_2^2 - 4)c'' + (t_1t_2 - 2s)(\xi_1 + \xi_2) + (4 - t_2^2)(\xi_1 + \xi_2)^2\big) \\
& \quad + \frac{\dd Z_1 \dd Z_2}{Z_1^2Z_2^2} \frac{1}{(Z_1^{-1} - Z_2^{-1})^2}\left[\frac{1}{\Psi_1\Psi_2} \left(\frac{1 + \frac{1}{2}(Y^{[0]}_1 + Y^{[0]}_2)}{1 + \frac{X_1^{[0]} - X_2^{[0]}}{Z_1^{-1} - Z_2^{-1}}}\right)^2 - 1 \right].
\end{split}
\end{equation}
It is indeed manifest that the double pole in the denominator of the second line cancel and this expression has a power series expansion near $(p_0,p_0)$ in the variables $\frac{1}{Z_1}$ and $\frac{1}{Z_2}$. A formula equivalent to \eqref{Omegaunuggeggga} and adapted for the expansion near $(p_{\varepsilon},p_{\varepsilon})$ is
\begin{equation}
\begin{split}
\mathcal{B}_{\boldsymbol{v}} & = \frac{\dd Z_1\dd Z_2}{4\Psi_1\Psi_2}\big((t_2^2 - 4)c'' + (t_1t_2 - 2s)(\xi_1 + \xi_2) + (4 - t_2^2)(\xi_1 + \xi_2)^2\big) \\
& \quad + \frac{\dd Z_1 \dd Z_2}{Z_1^2Z_2^2 (Z_1 - Z_2)^2}\left[\frac{1}{4(1 + \Psi_1^{[\varepsilon]})(1 + \Psi_2^{[\varepsilon]})} \left(\frac{Z_1^2(1 + Y_1^{[\varepsilon]}) + Z_2^2(1 + Y_2^{[\varepsilon]})}{1 - \frac{X_1^{[\varepsilon]} - X_2^{[\varepsilon]}}{Z_1^{-1}Z_2^{-1}(Z_1^{-1} - Z_2^{-1})}}\right)^2 - 1 \right].
\end{split}
\end{equation}
All symbols with exponent $[\varepsilon]$ being subleading, the expression in square bracket is manifestly regular near $(p_{\varepsilon},p_{\varepsilon})$ and has a power series expansion in the variables $\frac{1}{Z_1}$ and $\frac{1}{Z_2}$.

\begin{remark}
\label{rem:compare_spcurv}
An important difference with Section~\ref{sec:Oshaped}, where the spectral curve also had genus $1$, is that the equation of minimal degree \eqref{Pxieta} presenting $\Sigma_{\boldsymbol{v}}$ as an algebraic curve and $Z$ as an algebraic function, is more complicated: $Z$ has degree $3$ (= the number of trapezoids in the domain) and not $2$. This leads to more complicated expressions for $Z$ seen as function of the complex torus coordinate $\zeta$. Any two non-zero meromorphic functions on $\Sigma_{\boldsymbol{v}}$ are related by a polynomial equation. For us $\xi$ plays this role of this second meromorphic function, and we could express $Z$ in terms of $\xi$ by solving a quadratic equation.
\end{remark}

\subsection{Holomorphic \texorpdfstring{$1$}{1}-forms and variation with respect to filling fractions}

In the elliptic presentation, holomorphic $1$-forms are proportional to $\dd \zeta$. In the algebraic presentation \eqref{Pxieta}, it can be checked that the $1$-form $\frac{\dd Z}{\partial_\xi P(\xi,Z)}$ is holomorphic on $\Sigma_{\boldsymbol{v}}$. The proportionality constant relating it to $\dd \zeta$ can be determined by comparing the behavior near the point $p_0$. This yields
\begin{equation}
\label{dzetadxi}
\dd \zeta = \frac{(1-2t_2s)(c_+ + c_-) \dd Z}{\partial_{\xi} P(\xi,Z)} = \frac{(1 - 2t_2s)(c_+ + c_-)\dd Z}{Z^2 + t_1Z + 2t_2 \xi Z^2 + u_1 + 3 \xi^2}.
\end{equation}

The kernel of $\boldsymbol{\Theta}$ is spanned by $\boldsymbol{v}' = (0,1,-1)$. Proposition~\ref{Lem:Holo1form} shows that\
\[
\big\langle \boldsymbol{v} \otimes \boldsymbol{v}' \cdot \mathfrak{c}^{\textnormal{1st}}_{\bullet;\bullet}(z) \big\rangle \dd z = \big(\mathfrak{c}^{\textnormal{1st}}_{1;2}(z) - \mathfrak{c}^{\textnormal{1st}}_{1;3}(z)\big)\dd z
\]
continues to a holomorphic $1$-form $\mathfrak{u}_{\boldsymbol{v};\boldsymbol{v}'}$ such that
\[
\frac{1}{2\ii\pi} \int_{\mathcal{L}} \mathfrak{u}_{\boldsymbol{v};\boldsymbol{v}'} = 1.
\]
where the contour representing $\mathcal{L}$ goes counterclockwise around $[\alpha_2,\beta_2]$ in the $(0,1,1)$-th sheet. So, there must be a proportionality constant $M$ such that
\[
\mathfrak{u}_{\boldsymbol{v};\boldsymbol{v}'} = 2\ii\pi M\dd \zeta,\qquad M = \frac{1}{ \int_{\mathcal{L}} \dd \zeta}.
\]
Since $\mathcal{L}$ is a non-trivial primitive homology class on the torus $\amsmathbb{T}_{\tau} = \amsmathbb{C}/(\amsmathbb{Z} \oplus \tau \amsmathbb{Z})$, there exists coprime integers $(m,m')$ such that $M = \frac{1}{m + m'\tau}$.

Using the algebraic presentation, we can expand \eqref{dzetadxi} near $p_0$ or $p_{\pm}$ in the variable $Z$. Let us start the computation for instance near $p_0$. Recalling that $\xi \sim -s$ at this point, we obtain first from the equation \eqref{Pxieta} the expansion of $\xi$ near $p_0$
\[
\xi = -s - \frac{t_0 - st_1 + s^2t_2}{Z} + O\left(\frac{1}{Z^2}\right).
\]
This yields
\begin{equation}
\label{u12zc}
\quad \frac{\mathfrak{c}^{\textnormal{1st}}_{1;2}(z) - \mathfrak{c}^{\textnormal{1st}}_{1;3}(z)}{2\ii\pi} = \frac{1}{m + m'\tau}\left( \frac{c_- + c_+}{z^2} + \frac{(c_- + c_+)(2t_0t_2 - t_1 - 2st_1t_2 + 2t_2^2s^2)}{(1 - 2st_2)z^3} + O\left(\frac{1}{z^4}\right)\right)
\end{equation}
as $z \rightarrow \infty$. For the tiling model, Proposition~\ref{Proposition_differentiability_filling_fraction} tells us that the coefficient of $z^{-(k + 1)}$ in the right-hand side of \eqref{u12zc} is the first derivative with respect to $\hat{n}_2$ (at $\hat{n}_1$ and $\hat{n}_2 + \hat{n}_3$ fixed) of the $k$-th moment of the $h = 1$-th component of the equilibrium measure.

\begin{figure}[h!]
\[
\begin{array}{|c|c|c|c|}
\hline
\textnormal{sheet} & (1,0,0) & (-1,-1,-1) & (0,1,1) \\[4pt]
\hline
(1,0,0) & \mathcal{F}_{1,1} + \frac{1}{(Z_1 - Z_2)^2} & - \mathcal{F}_{1,1} - \mathcal{F}_{2,1} - \mathcal{F}_{3,1} & \mathcal{F}_{2,1} + \mathcal{F}_{3,1} \\[4pt]
\hline
(-1,-1,-1) & - \mathcal{F}_{1,1} - \mathcal{F}_{1,2} - \mathcal{F}_{1,3} & \sum_{h_1,h_2 = 1}^{3} \mathcal{F}_{h_1,h_2} + \frac{1}{(Z_1 - Z_2)^2} & - \sum_{h_1 = 2}^{3}\sum_{h_2 = 1}^{3} \mathcal{F}_{h_1,h_2}\\[4pt]
\hline
(0,1,1) & \mathcal{F}_{1,2} + \mathcal{F}_{1,3} & - \sum_{h_1 = 1}^{3} \sum_{h_2 = 2}^{3} \mathcal{F}_{h_1,h_2} & \sum_{h_1,h_2 = 2}^{3} \mathcal{F}_{h_1,h_2} + \frac{1}{(Z_1 - Z_2)^2} \\[4pt]
\hline
\end{array}
\]
\caption{\label{fig:doublesheetB6} The lines (respectively, columns) heads indicate the sheet to which the first (respectively, second) variable belongs to, and the entries are the expressions for the function $\frac{\mathcal{B}_{\boldsymbol{v}}}{\dd Z_1\dd Z_2} + \frac{1}{(Z_1 - Z_2)^2}$.}
\end{figure}

\begin{figure}[h!]
\begin{center}
\includegraphics[width=0.2\textwidth]{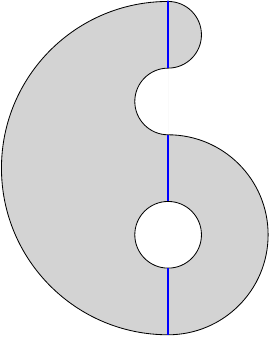}
\caption{\label{Fig:6shaped} Topology of a domain with matrix of interactions \eqref{Ex6shaped}.}
\end{center}
\end{figure}

\section{\texorpdfstring{$H = 3$}{H = 3}: triangle}
\label{sec:3yann}
\begin{figure}[h!]
\begin{center}
\includegraphics[width=0.4\textwidth]{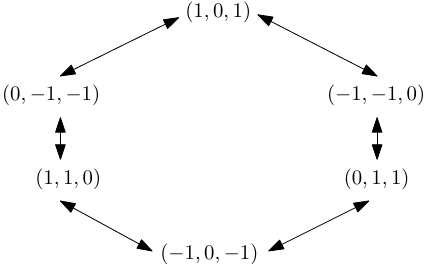}
\caption{\label{Fig:Cycle6} The $\mathfrak{G}$-orbit of $\boldsymbol{v} = (1,1,0)$.}
\end{center}
\end{figure}

\begin{figure}[h!]
\begin{center}
\includegraphics[width=0.2\textwidth]{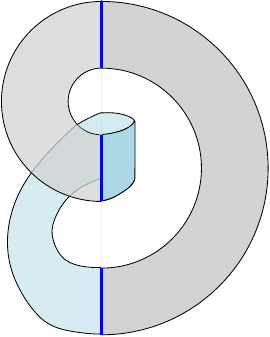}
\caption{\label{Fig:nonorient3} Topology of the non-bipartite domain with matrix of interactions \eqref{ExEshaped}}
\end{center}
\end{figure}

The same matrix of interaction that we used for the E-shaped domain
\begin{equation}
\label{Thetatrig}\boldsymbol{\Theta} = \left(\begin{array}{ccc} 1 & \frac{1}{2} & \frac{1}{2} \\[0.6ex] \frac{1}{2} & 1 & \frac{1}{2} \\[0.6ex] \frac{1}{2} & \frac{1}{2} & 1 \end{array}\right).
\end{equation}
has a realization as the matrix of interactions for a non-bipartite gluing of trapezoids, \textit{cf.} Figure~\ref{Fig:nonorient3}. Its gluing graph is a triangle. The only difference is that instead of the minimal orbit generated by $(1,1,1)$, we should use the non-minimal orbit of size $6$ generated by $\boldsymbol{v} = (1,1,0)$, which is adapted to the geometry of the domain. This is depicted in Figure~\ref{Fig:Cycle6}.

The spectral curve $\Sigma_{\boldsymbol{v}}$ has genus $1$ and is equipped with a non-trivial holomorphic involution $\varsigma$ having no fixed point and leaving $Z:= Z_{\boldsymbol{v}} : \Sigma_{\boldsymbol{v}} \rightarrow \widehat{\amsmathbb{C}}$ invariant. Taking the quotient gives again a genus $1$ Riemann surface together with a branched covering $\tilde{Z} : \Sigma_{\boldsymbol{v}}/\langle \varsigma \rangle \rightarrow \widehat{\amsmathbb{C}}$. It has three sheets, labeled by elements of the orbit of $\boldsymbol{v}$ modulo sign. In particular, the map $\tilde{Z}$ has six ramification points, one above each of the branch point $\alpha_1,\beta_1,\alpha_2,\beta_2,\alpha_3,\beta_3$, and it admits three simple poles, which we label by $p_0$ in the sheet corresponding to $\pm (1,1,0)$, $p_-$ in the sheet $\pm (1,0,1)$, $p_+$ in the sheet $\pm (1,1,1)$.

Using an identification $\Sigma_{\boldsymbol{v}}$ with a complex torus $\amsmathbb{T}_{\tau} = \amsmathbb{C}/(\amsmathbb{Z} \oplus \tau \amsmathbb{Z})$ with natural coordinate $\zeta$, we must have $\varsigma(\zeta) = \zeta + \varepsilon$ for some $\varepsilon \in \{\frac{1}{2},\frac{\tau}{2},\frac{1 + \tau}{2}\}$. Up to using an isomorphic complex torus, we can declare $\varepsilon = \frac{\tau}{2}$. Then $\Sigma_{\boldsymbol{v}} /\langle \varsigma \rangle$ is identified with $\amsmathbb{T}_{\frac{\tau}{2}}$, the maps $Z$ and $\tilde{Z}$ are represented by the same elliptic function of $\zeta$ with periods $1$ and $\frac{\tau}{2}$ having three simple poles a fundamental domain of $\amsmathbb{T}_{\frac{\tau}{2}}$. The only difference between $Z$ and $\tilde{Z}$ is that the former is seen as a function on $\amsmathbb{T}_{\tau}$ while the latter is seen as a function on $\amsmathbb{T}_{\frac{\tau}{2}}$. In Section~\ref{Ell6s} we already described the general shape of these functions depending on six parameters: we just have to replace $\tau$ by $\frac{\tau}{2}$ in \eqref{Zzetacccc} and the exact same conditions for the branch of the six parameters in terms of the branch points hold. This gives the parameters as algebraic functions of the branch points.

Given the properties of $\tilde{Z}$, the quotient curve $\Sigma_{\boldsymbol{v}}/\langle \varsigma\rangle$ is presented by an algebraic equation depending on six parameters in exactly the same way as in Section~\ref{1452Sec}: the equation is \eqref{Pxieta} with $\tilde{Z}$ replacing $Z$, and the parameterization of the branch points through these six parameters is identical. Reciprocally, this gives the parameters as algebraic functions of the branch points, but the branch to be chosen for these algebraic functions will be different than the one relevant in Section~\ref{1452Sec}, because the structure of the branched covering $\tilde{Z}$ (that can be visualized on Figure~\ref{Fig:nonorient3}) is different from the one the 6-shaped domain (visualized in Figure~\ref{Fig:6shaped}). The two branches are related by analytic continuation, as it was the case in the relation between Sections~\ref{GSshapesec} and \ref{H3Esec}.

From Theorem~\ref{thm:spcurvenonbip} and the expression of the fundamental bidifferential on genus $1$ curves in \eqref{Omegaelliptic}, the bidifferential
\[
\langle ( \boldsymbol{v} \otimes \boldsymbol{v}) \cdot \boldsymbol{\mathcal{F}}(z_1,z_2) \rangle \dd z_1 \dd z_2 = \big(\mathcal{F}_{1,1}(z_1,z_2) + \mathcal{F}_{1,2}(z_1,z_2) + \mathcal{F}_{2,1}(z_1,z_2) + \mathcal{F}_{2,2}(z_1,z_2)\big) \dd z_1\dd z_2
\]
 admits an analytic continuation $\mathcal{B}_{\boldsymbol{v}}$ as the meromorphic bidifferential on $\Sigma_{\boldsymbol{v}} \times \Sigma_{\boldsymbol{v}}$ with the following expression
\begin{equation}
\label{antiinvOme20}
\begin{split}
\mathcal{B}_{\boldsymbol{v}}(\zeta_1,\zeta_2) & = \big(\wp(\zeta_1 - \zeta_2|\tau) - \wp(\zeta_1 - \zeta_2 - \tfrac{\tau}{2}|\tau)\big)\dd \zeta_1 \dd \zeta_2 - \frac{\dd Z(\zeta_1) \dd Z(\zeta_2)}{(Z(\zeta_1) - Z(\zeta_2))^2} \\
& = \dd_{\zeta_1}\dd_{\zeta_2} \log\left(\frac{\vartheta_{11}(\zeta_1 - \zeta_2|\tau)}{\vartheta_{01}(\zeta_1 - \zeta_2|\tau) (Z(\zeta_1) - Z(\zeta_2))}\right),
\end{split}
\end{equation}
The theta function expression for the Weierstra\ss{} function was already mentioned in \eqref{Zzetatheta}, and we have used
\begin{equation}
\label{Jac01}
\vartheta_{01}(\zeta|\tau) := \sum_{n \in \amsmathbb{Z}} (-1)^n e^{\ii \pi n^2\tau + 2\ii\pi n \zeta} = -\ii e^{\ii\pi \frac{\tau}{4} - \ii\pi \zeta} \vartheta_{11}(\zeta - \tfrac{\tau}{2}|\tau).
\end{equation}
We recognize Jacobi sine function
\begin{equation}
\label{Jacobisnk}
\textnormal{sn}_{\mathsf{k}}(2K(\mathsf{k})\zeta) = - \frac{1}{\sqrt{\mathsf{k}}} \frac{\vartheta_{11}(\zeta|\tau)}{\vartheta_{01}(\zeta|\tau)},\qquad \tau = \frac{\ii K(\sqrt{1 - \mathsf{k}^2})}{K(\mathsf{k})}.
\end{equation}
This leads to the formula
\[
\mathcal{B}_{\boldsymbol{v}} = \dd_{\zeta_1}\dd_{\zeta_2} \log\left(\frac{\textnormal{sn}_{\mathsf{k}}[K(\mathsf{k})(\zeta_1 - \zeta_2)]}{Z(\zeta_1) - Z(\zeta_2)}\right).
\]

For tilings models in the domain of Figure~\ref{Fig:nonorient3} the filling fractions are deterministically fixed by the geometry of the domain. This is related to the fact that $\boldsymbol{\Theta}$ in \eqref{Thetatrig} is invertible. For the same reason, Proposition~\ref{Lem:Holo1form} does not construct any holomorphic $1$-forms on $\Sigma_{\boldsymbol{v}}$ from first-kind functions. In fact, the (unique up to scale) holomorphic $1$-form on the elliptic curve $\Sigma_{\boldsymbol{v}}$ does not play a role in the tiling model, because it is invariant under the involution while we care about anti-invariant ones (\textit{cf.} Theorem~\ref{thm:spcurvenonbip}).

\section{\texorpdfstring{$H = 3$}{H = 3}: 8-shaped domain}
\label{H38shap}
\begin{figure}[h!]
\begin{center}
\includegraphics[width=0.23\textwidth]{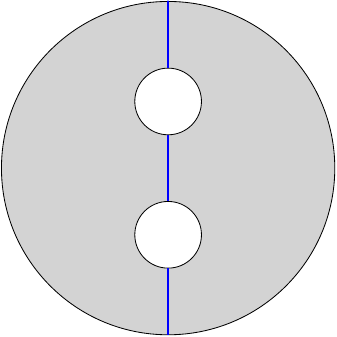}
\caption{\label{Fig:8shaped} The $8$-shaped domain}
\end{center}
\end{figure}

The $8$-shaped domain is obtained by gluing two trapezoids along three segments (\textit{cf.} Figure~\ref{Fig:8shaped}) with matrix of interactions full of $1$s. The vector $\boldsymbol{v} = (1,1,1)$ generates an orbit of size $2$, namely $\{\boldsymbol{v},-\boldsymbol{v}\}$. The corresponding spectral curve has genus $2$, it is hyperelliptic of equation
\[
\sigma^2 = (Z - \alpha_1)(Z - \beta_1)(Z - \alpha_2)(Z - \beta_2)(Z - \alpha_3)(Z - \beta_3).
\]
Since all $\theta_{g,h}$ are equal to $1$, the fundamental solution of the master Riemann--Hilbert problem has already been computed in Corollary~\ref{corHequal3}.

\section{\texorpdfstring{$H = 4$}{H=4}: B-shaped domain and its non-planar cousin}

\label{SecEx2}
\begin{figure}[h!]
\begin{center}
\includegraphics[width=0.5\textwidth]{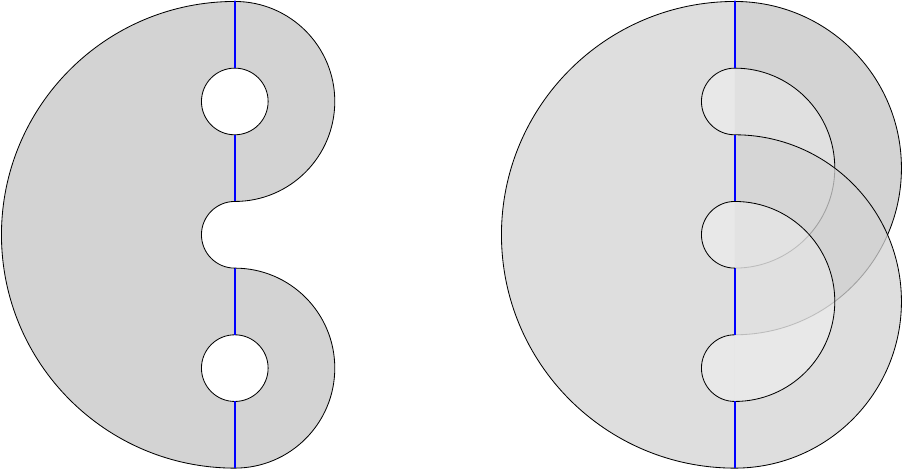}
\caption{\label{Fig:Bshaped} The $B$-shaped domain and its non-planar cousin.}
\end{center}
\end{figure}

Consider the matrix $\boldsymbol{\Theta}$ for a B-shaped domain obtained by gluing three trapezoids
\begin{equation}
\label{theta4eq}\boldsymbol{\Theta} = \left(\begin{array}{cccc} 1 & 1 & \frac{1}{2} & \tfrac{1}{2} \\[0.6ex] 1 & 1 & \frac{1}{2} & \frac{1}{2} \\[0.6ex] \frac{1}{2} & \frac{1}{2} & 1 & 1 \\[0.6ex] \frac{1}{2} & \frac{1}{2} & 1 & 1 \end{array}\right).
\end{equation}
If we exchange the role of the second and third segment (\textit{i.e.} the second and third rows and columns in $\boldsymbol{\Theta}$), we obtain a bipartite but non-planar domain, \textit{cf.} Figure~\ref{Fig:Bshaped}. The study we will make for the B-shaped domain is thus identical, except that at the end one should choose a branch of the algebraic functions giving the parameters of the spectral curve such that the order of $(\alpha_2,\beta_2)$ and $(\alpha_3,\beta_3)$ on the real line is exchanged.

The matrices for the four generators of the group $\mathfrak{G}$ in the canonical basis are
\begin{eqnarray*}
T^{(1)} = \left(\begin{array}{cccc} -1 & 0 & 0 & 0 \\ -2 & 1 & 0 & 0 \\ -1 & 0 & 1 & 0 \\ -1 & 0 & 0 & 1 \end{array}\right), & \qquad & T^{(2)} = \left(\begin{array}{cccc} 1 & -2 & 0 & 0 \\ 0 & -1 & 0 & 0 \\ 0 & -1 & 1 & 0 \\ 0 & -1 & 0 & 1\end{array}\right), \\
T^{(3)} = \left(\begin{array}{cccc} 1 & 0 & -1 & 0 \\ 0 & 1 & -1 & 0 \\ 0 & 0 & -1 & 0 \\ 0 & 0 & -2 & 1 \end{array}\right),& \qquad & T^{(4)} = \left(\begin{array}{cccc} 1 & 0 & 0 & -1 \\ 0 & 1 & 0 & -1 \\ 0 & 0 & 1 & -2 \\ 0 & 0 & 0 & -1 \end{array}\right).
\end{eqnarray*}
The group $\mathfrak{G}$ is infinite, as one can check
\[
(T^{(2)}\circ T^{(1)})^n(1,0,0,0) = (2n + 1,2n,n,n).
\]
The matrix $\boldsymbol{\Theta}$ having rank $2$, Theorem~\ref{thm:Omegav} tells us that $\hat{\mathfrak{G}}$ is the symmetric group in three elements. Let us check this directly from the definition. The space $\mathscr{V}$ from Proposition~\ref{lemrefl} is spanned by $(1,1,0,0)$ and $(0,0,1,1)$, and $\hat{\mathfrak{G}}$ is generated by the orthogonal reflections with respect to $\hat{\boldsymbol{e}}^{(1)} = \hat{\boldsymbol{e}}^{(2)} = (\frac{1}{2},\frac{1}{2},0,0)$ and $\hat{\boldsymbol{e}}^{(3)} = \hat{\boldsymbol{e}}^{(4)} = (0,0,\frac{1}{2},\frac{1}{2})$, using the scalar product
\begin{equation*}
\begin{split}
& \langle \hat{\boldsymbol{e}}^{(1)} \cdot \boldsymbol{\Theta}(\hat{\boldsymbol{e}}^{(1)}) \rangle = \langle \hat{\boldsymbol{e}}^{(3)} \cdot \boldsymbol{\Theta}(\hat{\boldsymbol{e}}^{(3)}) \rangle = 1,\\
& \langle \hat{\boldsymbol{e}}^{(1)} \cdot \boldsymbol{\Theta}(\hat{\boldsymbol{e}}^{(3)}) \rangle = \langle \hat{\boldsymbol{e}}^{(3)} \cdot \boldsymbol{\Theta}(\hat{\boldsymbol{e}}^{(1)}) \rangle = \frac{1}{2}.
\end{split}
\end{equation*}
We retrieve the example of Section~\ref{sec:Csha} with the substitutions $\boldsymbol{e}^{(1)} \rightarrow \hat{\boldsymbol{e}}^{(1)}$ and $\boldsymbol{e}^{(2)} \rightarrow \hat{\boldsymbol{e}}^{(3)}$, thus showing that $\hat{\mathfrak{G}}$ is the symmetric group in three elements. We are interested in the orbit of $\boldsymbol{v} = (1,1,0,0)$, that is
\begin{equation*}
(1,1,0,0) \,\,\mathop{\longleftrightarrow}^{T^{(1)},T^{(2)}} \,\, (-1,-1,-1,-1)\,\, \mathop{\longleftrightarrow}^{T^{(3)},T^{(4)}} \,\,(0,0,1,1).
\end{equation*}
The resulting spectral curve $\Sigma_{\boldsymbol{v}}$ has genus $2$ and the map $Z:= Z_{\boldsymbol{v}} : \Sigma_{\boldsymbol{v}} \rightarrow \widehat{\amsmathbb{C}}$ has three simple poles and eight ramification points. Let us call $p_0$ (respectively $p_+$ or $p_-$) the pole in the sheet labeled $(-1,-1,-1,-1)$ (respectively $(1,1,0,0)$ or $(0,0,1,1)$). the $(1,1,0,0)$-th sheet.

Any genus $2$ Riemann surface admits two meromorphic functions $\xi,\eta$ satisfying an equation of the form
\begin{equation}
\label{xietaF}
\eta^2 = F(\xi),\qquad F(\xi) = \xi(\xi - 1)(\xi - \lambda_1)(\xi - \lambda_2)(\xi - \lambda_3)
\end{equation}
for some $\lambda_1,\lambda_2,\lambda_3 \in \amsmathbb{C} \setminus \{0,1\}$ pairwise disjoint. Call $\textnormal{P}$ the set of Weierstrass points, \textit{i.e.} the five points at which $\eta = 0$ and the point $p_{\infty}:=(\infty,\infty)$. Let us choose three points on the curve $(\xi_0,\eta_0),(\xi_{\pm},\eta_{\pm})$ away from $\textnormal{P}$. The Riemann--Roch theorem says that the space of meromorphic functions with at most simple poles at those three points has dimension two. We can check this explicitly. The general form of a function with three simple poles is
\begin{equation}
\label{Zcccetac12}
Z = c \frac{\eta + \eta_+}{\xi - \xi_{+}} + c_+ \frac{\eta + \eta_+}{\xi - \xi_+} + c_- \frac{\eta + \eta_-}{\xi - \xi_-} + c',
\end{equation}
where the constants $c,c_{\pm}$ are subjected to the condition that there is no pole at $p_{\infty}$. Near $p_{\infty}$ we have $\eta^2 \sim \xi^5$, so we can find a centered local coordinate $\zeta$ such that $\xi = \zeta^{-2}$, and from \eqref{xietaF} it follows that $\eta = \zeta^{-5} + O(\zeta^{-1})$. The right-hand side of \eqref{Zcccetac12} has Laurent series expansion
\[
(c_0 + c_+ + c_-)\zeta^{-3} + (c_0\xi_0 + c_+\xi_+ + c_-\xi_-)\zeta^{-1} + O(1).
\]
The absence of pole at $p_{\infty}$ in \eqref{Zcccetac12} then yields
\begin{equation}
\label{Zxietac}
Z = c\left(\frac{\eta + \eta_0}{\xi - \xi_0} + \frac{\xi_- - \xi_0}{\xi_+ - \xi_-} \frac{\eta + \eta_+}{\xi - \xi_+} + \frac{\xi_0 - \xi_+}{\xi_+ - \xi_-} \frac{\eta + \eta_-}{\xi - \xi_-}\right) + c'.
\end{equation}
This expression depends on eight parameters $c,c',\xi_0,\xi_+,\xi_-,\lambda_1,\lambda_2,\lambda_3$, and the choice of signs of square roots to get $\eta_0,\eta_{\pm}$ from $\xi_0,\xi_{\pm}$. We can invert \eqref{Zxietac} to express the function $\eta$ in terms of $Z$ and $\xi$. The result takes the form $\eta = A(\xi)Z + B(\xi)$ where $A(\xi)$ and $B(\xi)$ are polynomials of degree $3$ and $2$, respectively. Inserting it in \eqref{xietaF} gives the equation
\begin{equation}
\label{igunegrg}
(A(\xi)Z + B(\xi))^2 = F(\xi) \quad \textnormal{with} \quad \left\{\begin{array}{lll} A(\xi) = \frac{(\xi - \xi_0)(\xi - \xi_-)(\xi - \xi_+)}{(\xi_0 - \xi_-)(\xi_0 - \xi_+)} \\[4pt] B(\xi_{\circ}) = \eta_{\circ} \,\,\textnormal{for}\,\,\circ \in \{0,\pm\} \end{array}\right.
\end{equation}
The eight zeros of $\dd Z$, which are to be found among the common solution of
\[
(A(\xi)Z + B(\xi))^2 = F(\xi),\qquad 2(A'(\xi)Z + B(\xi))(A(\xi)Z + B(\xi)) = F'(\xi)
\]
Eliminating $Z$ in this system to expressing it as a rational function of $\xi$ which are the roots of the polynomial $(-AF' + 2A'F)^2 = 4F(A'B - B'A)^2$. This polynomial has degree $14$, but it admits three double roots at $\xi_0,\xi_{\pm}$ due to the last property in \eqref{igunegrg}. This leaves us with eight roots, giving rise to eight values for $Z$, which are exactly $\alpha_1,\beta_1,\alpha_2,\beta_2,\alpha_3,\beta_3,\alpha_4,\beta_4$. This gives a parametric expression for symmetric functions of the latter in terms of the eight parameters $c,c',\xi_0,\xi_+,\xi_-,\lambda_1,\lambda_2,\lambda_3$. We do not attempt to be more explicit.

The birational transformation
\[
\xi = X^{-1} - 1,\qquad \eta = \qquad Y = \frac{\ii X^3 \eta}{\sqrt{2\prod_{a = 1}^{3} (1 + \lambda_a)}}
\]
transforms the algebraic equation \eqref{xietaF} for $\Sigma_{\boldsymbol{v}}$ into
\[
Y^2 = \tilde{F}(X),\qquad \tilde{F}(X) = X\Big(X - \frac{1}{2}\Big)(X - 1)\prod_{a =1}^{3} \Big(X - \frac{1}{1 + \lambda_a}\Big)
\]
With this parameterization, the fundamental bidifferential $\mathcal{B}^{\Sigma_{\boldsymbol{v}},\mathcal{L}}$ was described in \eqref{Bgenus2}. This allows in principle the computation of the analytic continuation of
\[
(\mathcal{F}_{1,1}(z_1,z_2) + \mathcal{F}_{1,2}(z_1,z_2) + \mathcal{F}_{2,1}(z_1,z_2) + \mathcal{F}_{2,2}(z_1,z_2))\dd z_1\dd z_2,
\]
since by Theorem~\ref{thm:Omegav} it is the bidifferential $\mathcal{B}_{\boldsymbol{v}}$ on $\Sigma_{\boldsymbol{v}} \times \Sigma_{\boldsymbol{v}}$ given by
\[
\mathcal{B}_{\boldsymbol{v}} = \mathcal{B}^{\Sigma_{\boldsymbol{v}},\mathcal{L}} - \frac{\dd Z_1\dd Z_2}{(Z_1 - Z_2)^2}.
\]
In order to access leading covariances in tiling models, we need the expansion of this bidifferential near the points $Z_1,Z_2 \rightarrow \infty$ in various sheets, that can be accessed with the same technique as in Section~\ref{thebingfsgun}.

\section{\texorpdfstring{$H = 4$}{H = 4}: two non-domain examples}

\begin{figure}[h!]
\begin{center}
\includegraphics[width=0.95\textwidth]{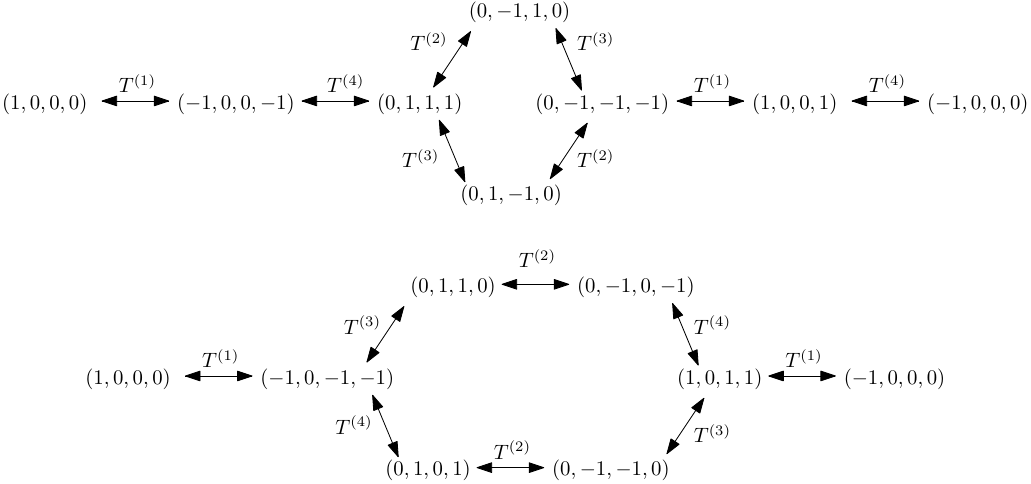}
\caption{\label{Fig:8orbit} The respective orbits of $\boldsymbol{v} = (1,0,0,0)$ for the matrices of interactions \eqref{EqThetanon}.}
\end{center}
\end{figure}

Let us consider the two representative non-realizable connectoids of size $4$, \textit{cf.} \eqref{exnonrel}.
\begin{equation}
\label{EqThetanon}
\boldsymbol{\Theta} = \left(\begin{array}{cccc} 1 & 0 & 0 & \frac{1}{2} \\[0.6ex] 0 & 1 & 0 & \frac{1}{2} \\[0.6ex] 0 & 0 & 1 & \frac{1}{2} \\[0.6ex] \frac{1}{2} & \frac{1}{2} & \frac{1}{2} & 1 \end{array}\right)\qquad \boldsymbol{\Theta} = \left(\begin{array}{cccc} 1 & 0 & \frac{1}{2} & \frac{1}{2} \\[0.6ex] 0 & 1 & \frac{1}{2} & \frac{1}{2} \\[0.6ex] \frac{1}{2} & \frac{1}{2} & 1 & \frac{1}{2} \\[0.6ex] \frac{1}{2} & \frac{1}{2} & \frac{1}{2} & 1 \end{array}\right).
\end{equation}
These matrices have full rank and we mentioned in Section~\ref{sec:extraconn} that the associated group $\mathfrak{G}$ is the Weyl group of type $D_4$. In both cases the orbit of $\boldsymbol{v} := \boldsymbol{e}^{(1)}$ has minimal size $8$, \textit{cf.} Figure~\ref{Fig:8orbit}, the spectral curve $\Sigma_{\boldsymbol{v}}$ has genus $1$ and is equipped with a holomorphic involution without fixed point and sending the $\boldsymbol{w}$-th sheet to the $(-\boldsymbol{w})$-th sheet for any $\boldsymbol{w} \in \mathfrak{G}.\boldsymbol{v}$. From the sheet structure, we see that the quotient curve $\Sigma_{\boldsymbol{v}}/\langle \varsigma\rangle$ also has genus $1$. The branched covering $Z:= Z_{\boldsymbol{v}} : \Sigma_{\boldsymbol{v}} \rightarrow \widehat{\amsmathbb{C}}$ has eight simple poles, and in the quotient it determines a branched covering $\tilde{Z} : \Sigma_{\boldsymbol{v}}/\langle \varsigma\rangle \rightarrow \widehat{\amsmathbb{C}}$ with four simple poles.

Letting $\widetilde{\boldsymbol{v}} := (\boldsymbol{v} \otimes \boldsymbol{v}\,,\,0)$, one can check that
\[
\mathfrak{G}^{[2]}.\widetilde{\boldsymbol{v}} = \big\{(\boldsymbol{w} \otimes \boldsymbol{w}'\,,\,-1 + \delta_{\boldsymbol{w},\boldsymbol{w}'} - \delta_{\boldsymbol{w},-\boldsymbol{w}'})\quad \big|\quad \boldsymbol{w},\boldsymbol{w}' \in \mathfrak{G}.\boldsymbol{v}\big\}
\]
Therefore, the bidifferential $\mathcal{F}_{1,1}(z_1,z_2)\dd z_1 \dd z_2$ admits an analytic continuation to the spectral surface $\Sigma_{\boldsymbol{v}} \times \Sigma_{\boldsymbol{v}}/\langle(\varsigma,\varsigma)\rangle$ and can be pulled back to a meromorphic bidifferential $\mathcal{B}_{\boldsymbol{v}}$ on $\Sigma_{\boldsymbol{v}} \times \Sigma_{\boldsymbol{v}}$ taking the form
\[
\mathcal{B}_{\boldsymbol{v}} = \mathcal{B}^{\Sigma_{\boldsymbol{v}}} - (\varsigma,\textnormal{Id})^*\mathcal{B}^{\Sigma_{\boldsymbol{v}}} - \frac{\dd Z_1\dd Z_2}{(Z_1 - Z_2)^2}.
\]
Here, $\mathcal{B}^{\Sigma_{\boldsymbol{v}}}$ is the fundamental bidifferential on $\Sigma_{\boldsymbol{e}^{(1)}}$ normalized to have zero periods around $[\alpha_h,\beta_h]$ for each $h \in [4]$. All of this is similar to what Theorem~\ref{thm:spcurvenonbip} would give for a connectoid of rank $4$ that is realizable in a non-bipartite way, except that the entries of the vectors in the orbit can have mixed signs. This difference does not affect the geometry of the spectral curve. We could continue like in Sections~\ref{sec:6yann} and \ref{sec:3yann} to get a parameterization of the branched points of $Z$ and some expressions of $\mathcal{B}_{\boldsymbol{v}}$.

\appendix

\part{Appendices}
\label{apppart}
\section*{List of definitions (D)}
{\footnotesize
\begin{tabular}{lll}
\hline\hline
D \ref{def:eq_rescaled_parameters} & Rescaled endpoints \dotfill & \pageref{def:eq_rescaled_parameters} \\
D \ref{def:eq_shifted_parameters} & Shifted rescaled endpoints \dotfill & \pageref{def:eq_shifted_parameters} \\
D \ref{def:eq_enlarged_parameters} & Enlarged rescaled segments \dotfill & \pageref{def:eq_enlarged_parameters} \\
D \ref{Definition_functional} & Energy functional \dotfill & \pageref{Definition_functional} \\
D \ref{Definition_equilibrium_measure} & Equilibrium measure \dotfill & \pageref{Definition_equilibrium_measure} \\
D \ref{Definition_void_saturated} & Bands, voids, saturations \dotfill & \pageref{Definition_void_saturated} \\
D \ref{def_eff_pot} & Effective potential \dotfill & \pageref{def_eff_pot} \\
D \ref{Definition_phi_functions} & Auxiliary functions for weight decomposition \dotfill & \pageref{Definition_phi_functions} \\
\hline\hline
\end{tabular}

\newpage

\begin{tabular}{lll}
\hline\hline
D \ref{unimono} & Group $\mathfrak{G}$ by generators \dotfill & \pageref{unimono} \\
D \ref{def:spcurv1} & Spectral curve of a seed vector \dotfill & \pageref{def:spcurv1} \\
D \ref{def:goodor} & Good orbit of $\mathfrak{G}$ \dotfill & \pageref{def:goodor} \\
D \ref{def:smallmono} & Reduced group $\hat{\mathfrak{G}}$ \dotfill & \pageref{def:smallmono} \\
D \ref{extgroups} & Extended two-variable group $\mathfrak{G}^{[2]}$ \dotfill & \pageref{extgroups} \\
D \ref{def:spsurf} & Spectral surface of seed vector \dotfill & \pageref{def:spsurf} \\
D \ref{def:Bires} & Matrix of biresidues \dotfill & \pageref{def:Bires} \\
D \ref{def:Hvv} & Double primitive of fundamental solution \dotfill & \pageref{def:Hvv} \\
D \ref{Definition_bipartite} & Bipartite orbits \dotfill & \pageref{Definition_bipartite} \\
D \ref{def:halfspcurvebip} & Half-spectral curve \dotfill & \pageref{def:halfspcurvebip} \\
D \ref{def:connectoid} & Connectoid \dotfill & \pageref{def:connectoid} \\
\hline\hline
\end{tabular}}

\newpage

\section*{List of results}

T = Theorem, P = Proposition, L = Lemma, C = Corollary

{\footnotesize
\begin{flushleft} \begin{tabular}{lll}
\hline\hline
 T \ref{Theorem_equi_charact} & \textit{cf.} T \ref{Theorem_equi_charact_repeat_2} \dotfill & \pageref{Theorem_equi_charact} \\
L \ref{Lemma_phi_properties} & Direct properties of auxiliary functions from Assumption~\ref{Assumptions_analyticity} & \pageref{Lemma_phi_properties} \\
\hline\hline
T \ref{Theorem_Nekrasov} & First-order Nekrasov equations \dotfill & \pageref{Theorem_Nekrasov} \\
C \ref{Corollary_higherNek} & Higher-order Nekrasov equations \dotfill & \pageref{Corollary_higherNek} \\
\hline\hline
L \ref{Lemma_I_welldef} & Energy functional $-\mathcal{I}$ is well-defined \dotfill & \pageref{Lemma_I_welldef} \\
L \ref{Lemma_I_quadratic_smoothing} & Quadratic part $\mathcal{I}^{(2)}$ is well-defined \dotfill & \pageref{Lemma_I_quadratic_smoothing} \\
L \ref{Lemma_I_quadratic_Fourier} \& C \ref{Corollary_I_positive} & Fourier representation and positivity of $-\mathcal{I}^{(2)}$ \dotfill & \pageref{Lemma_I_quadratic_Fourier} \& \pageref{Corollary_I_positive} \\
C \ref{Corollary_I_bound} \& L \ref{Lemma_J_alpha} & Quantitive bounds for $\mathcal{I}^{(2)}$ \dotfill & \pageref{Corollary_I_bound} \& \pageref{Lemma_J_alpha}\\
\hline\hline
\end{tabular}

\newpage

\begin{tabular}{lll}
\hline\hline
P \ref{Proposition_Olsh_ZW} & Partition function of the $zw$-ensemble (Olshanski) \dotfill & \pageref{Proposition_Olsh_ZW} \\
P \ref{Prop_62} & Partition function of reference ensemble \dotfill & \pageref{Prop_62} \\
L \ref{Checkzw} & Checking Assumptions for reference ensemble \dotfill & \pageref{Checkzw} \\
P \ref{Proposition_ZW_partition} & Free energy asymptotics up to $o(1)$ for reference ensemble \dotfill & \pageref{Proposition_ZW_partition} \\
L \ref{Lemma_densely_packed_expansion} & Asymptotics of interaction factor in saturations \dotfill & \pageref{Lemma_densely_packed_expansion} \\
T \ref{Proposition_ZW_equilibrium_measure} & Equilibrium measure of reference ensemble \dotfill & \pageref{Proposition_ZW_equilibrium_measure} \\
C \ref{Corollary_ZW_off-critical_1}-\ref{Corollary_ZW_off-critical_2} & Off-criticality for reference ensemble \dotfill & \pageref{Corollary_ZW_off-critical_1}-\pageref{Corollary_ZW_off-critical_2} \\
P \ref{Proposition_ZW_partition_leading} & Leading term in free energy is given by $\mathcal{I}$ \dotfill & \pageref{Proposition_ZW_partition_leading} \\
\hline\hline
\end{tabular}
\end{flushleft}}

\newpage
{\footnotesize
\begin{flushleft}
\begin{tabular}{lll}
\hline\hline
T \ref{Theorem_Master_equation_12} & Solution operation to master equation and its properties \dotfill & \pageref{Theorem_Master_equation_12} \\
T \ref{Theorem_Masterspecial} & Simple solution operators \dotfill & \pageref{Theorem_Masterspecial} \\
L \ref{Pinvdet} & Invertibility of period map \dotfill & \pageref{Pinvdet} \\
L \ref{Lemma_eq_zero_solution} \& C \ref{Corollary_Inj} & Injectivity of reduced operator \dotfill & \pageref{Lemma_eq_zero_solution} \\
L \ref{Lem:ConversionL2op}-\ref{Lem:ConversionL2op2} & Conversion to operators in $\mathscr{L}^2$ spaces \dotfill & \pageref{Lem:ConversionL2op}-\pageref{Lem:ConversionL2op2} \\
L \ref{Lemma_invertibility_Fredholm} & Invertibility of reduced operators in $\mathscr{L}^2$ spaces \dotfill & \pageref{Lemma_invertibility_Fredholm} \\
P \ref{prop:tameRHP} & Tameness for solutions of master RHP \dotfill & \pageref{prop:tameRHP} \\
L \ref{Lem:app_RHP_mu} & Application of the solution operator for equilibrium measures \dotfill & \pageref{Lem:app_RHP_mu} \\
\hline\hline
\end{tabular}
\end{flushleft}}

\newpage

\section*{List of remarks}

{\footnotesize
\begin{flushleft}
\begin{tabular}{lll}
\hline\hline
R \ref{rembounder} & Worst divergence of derivative of potential \dotfill & \pageref{rembounder} \\
R \ref{remtheta11} & Simplification if $\theta_{h,h} = 1$ \dotfill & \pageref{remtheta11} \\
R \ref{rem:vanishing_interpretation} & Vanishing of weights vs. decomposition factors \dotfill & \pageref{rem:vanishing_interpretation} \\
R \ref{rem:allway} & Fixing or not filing fractions in Nekrasov equations \dotfill & \pageref{rem:allway} \\
R \ref{remark:continuous} & Equilibrium measures without constraints on max. density \dotfill & \pageref{remark:continuous} \\
R \ref{remark:v_when_fully_void_or_sat} & Fully void or fully saturated segments \dotfill & \pageref{remark:v_when_fully_void_or_sat} \\
R \ref{rkcompactsupport} & Particles stay in a compact \dotfill & \pageref{rkcompactsupport} \\
R \ref{remextrasssh} & Simplification in definition of $\psi$ \dotfill & \pageref{remextrasssh} \\
\hline\hline
\end{tabular}
\end{flushleft}}

\section*{List of figures}

{\footnotesize
\begin{flushleft}
\begin{tabular}{lll}
\hline\hline
F \ref{Fig_tiling_hex_222} & Tiling of a hexagon \dotfill & \pageref{Fig_tiling_hex_222} \\
F \ref{Fig_hex_hole_small} & Tiling of a small hexagon with hole \dotfill & \pageref{Fig_hex_hole_small} \\
F \ref{Fig:Cshapeint} & C-shaped domains \dotfill & \pageref{Fig:Cshapeint} \\
F \ref{Fig:Cshape_simulation} & Tilings of C-shaped domains \dotfill & \pageref{Fig:Cshape_simulation} \\
F \ref{Fig:Cshape_simulation_fluct} & Simulation of height fluctuations in C-shaped domain \dotfill & \pageref{Fig:Cshape_simulation_fluct} \\
F \ref{Fig_Gaussian_intro} & Plot of equilibrium measure for Gaussian weight \dotfill & \pageref{Fig_Gaussian_intro} \\
F \ref{fig:Cgluefig} & Spectral curve for C-shaped domain \dotfill & \pageref{fig:Cgluefig} \\
\hline\hline
F \ref{Fig_Gaussian} & Plot of effective potential for Gaussian weights \dotfill & \pageref{Fig_Gaussian} \\
\hline\hline
\end{tabular}
\end{flushleft}}

\newpage

\section*{List of symbols}

Page numbers refer to the first occurrence or place of definition.

{\footnotesize \begin{flushleft}
\begin{tabular}{lll}
\multicolumn{3}{l}{\textsc{Special functions}} \\
\hline\hline
$\mathbbm{1}$ & Indicator function \dotfill & \pageref{index:indicfn} \\
\hline
$\mathsf{e}_k$ & Elementary symmetry polynomial \dotfill & \pageref{index:elemsym} \\
\hline
$\textnormal{B}_{2n}$ & Bernoulli numbers \dotfill & \pageref{index:Bernou} \\
\hline
$(x)_n$ & Pochhammer symbol \dotfill & \pageref{Pochdef} \\
$\mathrm{Llog}, \mathrm{LLlog}$ & Integrated logarithms \dotfill & \pageref{Llogdef} \\
$\Gamma_2, \chi_2$ & Double Gamma and related functions \dotfill & \pageref{index:Gamma2fn} \\
$\gimel_N(x), \tilde{\gimel}_N$ & Products and ratios of Gammas \dotfill & \pageref{phiexp} \\
$\mathsf{\Psi}$ & Digamma function \dotfill & \pageref{index:digam} \\
\hline
$E(\mathsf{k}), K(\mathsf{k})$ & Complete elliptic integrals \dotfill & \pageref{index:ellipticEK} \\
$\Pi_1[c;\mathsf{k}]$ & Incomplete elliptic integral \dotfill & \pageref{index:ellipticEK} \\
$\wp$ & Weierstra\ss{} function \dotfill & \pageref{index:wp} \\
$\eta_{\mathsf{a}},\eta_{\mathsf{b}}$ & Periods of $\wp$ \dotfill & \pageref{index:eta} \\
$E_2(\tau)$ & Second Eisenstein series \dotfill & \pageref{index:E2} \\
$\vartheta_{11}, \vartheta_{01}$ & Jacobi theta functions \dotfill & \pageref{Jacobitheta11}, \pageref{Jac01} \\
$\mathsf{sn}_{\mathsf{k}}$ & Jacobi sine function \dotfill & \pageref{Jacobisnk} \\
$Z^{\textnormal{J}}$ & Zhukovsky map \dotfill & \pageref{Joukov} \\
\hline\hline
\end{tabular}
\end{flushleft}

\vspace{0.2cm}

{\footnotesize \begin{flushleft}
\begin{tabular}{lll}
\multicolumn{3}{l}{\textsc{Sets}} \\
\hline\hline
$I^{\pm}(h)$ & Index of extreme particles \dotfill & \pageref{index:Iplusmoins} \\
$\amsmathbb{W}_{\N}$ & State space \dotfill & \pageref{Section_configuration_space} \\
$\amsmathbb{L}$ & Lattice in $\amsmathbb{R}$ \dotfill & \pageref{index:Iplusmoins} \\
$\Lambda_{\star},\hat{\Lambda}_{\star},\hat{\Lambda}_{\star}^{\amsmathbb{R}}$ & Filling fraction ranges \dotfill & \pageref{eq_Lambda_def} \\
\hline
$\mathcal{H}_i$ & Sets of segments in $i$-th trapezoid \dotfill & \pageref{Definition_gluing} \\
\hline
$\mathfrak{G}, \hat{\mathfrak{G}}$ & Groups of sheet transitions \dotfill & \pageref{unimono}, \pageref{def:smallmono} \\
$\mathfrak{M}_{\boldsymbol{v}}$ & Monodromy group \dotfill & \pageref{index:monodromy} \\
$\mathfrak{G}^{[2],\circ}, \mathfrak{G}^{[2]}$ & Two-variable versions of $\mathfrak{G}$ \dotfill & \pageref{extgroups}, \pageref{Sv2eq} \\
\hline
$\Delta_{\boldsymbol{v}}, \Delta^*_{\boldsymbol{v}}$ & Diagonal and anti-diagonal \dotfill & \pageref{diagantidiag} \\
\hline
$\mathsf{G}$ & Gluing graph \dotfill & \pageref{sec:green} \\
$\mathfrak{C}_h^{\pm}$ & Intersection of equivalence classes & \pageref{Prop:realiza} \\
\hline\hline
\end{tabular}
\end{flushleft}

\vspace{0.2cm}

{\footnotesize \begin{flushleft}
\begin{tabular}{lll}
\multicolumn{3}{l}{\textsc{Domains}} \\
\hline\hline
$\mathcal{D}$ & Tileable domain \dotfill & \pageref{Definition_gluing} \\
$\mathcal{T}_i$ & Trapezoid \dotfill & \pageref{Definition_gluing} \\
\hline
$\amsmathbb{A}$ & Segment for equilibrium measure \dotfill & \pageref{eq_measure_space} \\
$\amsmathbb{A}^{\mathfrak{m}}_h$ & Enlarged rescaled segments \dotfill & \pageref{def:eq_enlarged_parameters}, \pageref{sec:Masterpb} \\
$\amsmathbb{V},\amsmathbb{B},\amsmathbb{S}$ & Voids, bands, saturations \dotfill & \pageref{Definition_void_saturated} \\
\hline
$\mathcal{D}^{\pm}$ & Parameter ranges for reference ensemble \dotfill & \pageref{Definition_ZW_set} \\
$\mathcal{D}^{\textnormal{off}}$ & Off-critical locus \dotfill & \pageref{Definition_off_critical_set_ZW} \\
\hline \hline
\end{tabular}
\end{flushleft}}

\vspace{0.5cm}

{\footnotesize \begin{flushleft}
\begin{tabular}{lll}
\multicolumn{3}{l}{\textsc{Operators}} \\
\hline\hline
$\Upsilon_h$ & Solution operator \dotfill & \pageref{Theorem_Master_equation}, \pageref{Theorem_Master_equation_12} \\
$\textnormal{Proj}_h$ & Projection of analytic functions \dotfill & \pageref{Fh_projection} \\
$\mathcal{P}, \mathcal{P}^H$ & Polynomial part operators \dotfill & \pageref{Polop}, \pageref{PolopH} \\
$\Pi$ & Period map \dotfill & \pageref{periodmap} \\
$\mathcal{K},\mathcal{U}_h$ & Auxiliary operators \dotfill &\pageref{Khdef}, \pageref{UhPl} \\
$\mathcal{S}$ & Integral operator on $\mathscr{L}^2$-spaces \dotfill & \pageref{Lintegral} \\
$S(z,w)$ & Integral kernel of $\mathcal{S}$ \dotfill & \pageref{Szwkernel} \\
$\mathcal{R}$ & Resolvent operator \dotfill & \pageref{index:resolvent} \\
\hline
$T^{(h)}$ & Generators of $\mathfrak{G}$ \dotfill & \pageref{unimono} \\
$\hat{T}^{(h)}$ & Generators of $\hat{\mathfrak{G}}$ \dotfill & \pageref{lemrefl} \\
$T^{(h)}_{1}, T_2^{(h)}$ & Two-variable version of $T^{(h)}$ \dotfill & \pageref{extgroups} \\
\hline\hline
\end{tabular}
\end{flushleft}}

\vspace{0.5cm}

{\footnotesize \begin{flushleft}
\begin{tabular}{lll}
\multicolumn{3}{l}{\textsc{Probability}} \\
\hline\hline
$\amsmathbb{P}_{\N}$ & Probability measure \dotfill & \pageref{eq_general_measure_one_cut} \\
$\overline{\amsmathbb{P}}_{\N}$ & Unnormalized version of $\amsmathbb{P}_{\N}$ \dotfill & \pageref{index:notnorm} \\
$\amsmathbb{E}^{(\textnormal{c})}$ & Cumulant expectation value \dotfill & \pageref{eq_cumulant_def} \\
$\mathscr{Z}_\N$ & Partition function \dotfill &\pageref{eq_zw_1} \\
$\mathscr{Z}_{N}^{zw}$ & Partition function of $zw$-ensemble \dotfill & \pageref{zwintro}, \pageref{zwdiscretemd} \\
$\mathscr{Z}_{\N}^{\textnormal{ref}}$ & Partition function of reference ensemble \dotfill & \pageref{S:sec_first_type} \\
$\mathscr{F}^{[p],\textnormal{ref}}$ & Free energy expansion coefficients \dotfill & \pageref{Proposition_ZW_partition} \\
\hline\hline
\end{tabular}
\end{flushleft}}

\vspace{0.5cm}

{\footnotesize \begin{flushleft}
\begin{tabular}{lll}
\multicolumn{3}{l}{\textsc{Error terms}} \\
\hline\hline
$\mathbbm{e}_h$ & Error between potential and weight \dotfill & \pageref{eq_weight_form} \\
$\mathbbm{f}_h$ & Remainder in expansion of $\Phi_h^{\pm}$ & \pageref{Def_finer_Phi} \\
\hline
$\mathbbm{q}_h^{\pm}$ & Remainder in expansion of $Q_h^{\pm}$ \dotfill & \pageref{eq_Q_plus_determ_expansion} \\
$\mathbbm{Err}^{(l)}$ & Errors in Nekrasov equations \dotfill & \pageref{eq_First_Nekr_asymptotic} \\
\hline
$\mathbbm{Rest}_0^{\textnormal{ref}}$ & Subleading free energy of reference ensemble \dotfill & \pageref{Proposition_ZW_partition_leading} \\
$\mathbbm{Rest}_1^{\mathcal{I}}$ & Localization error in energy functional \dotfill & \pageref{RestI:index} \\
$\mathbbm{Rest}$ & Error in CLT \dotfill & \pageref{eq_x164} \\
$\mathbbm{Rest}_1$ & Error in expansion of $\mathscr{Z}_{\N}$ \dotfill & \pageref{eq_partition_one_band}, \pageref{eq_partition_multicut} \\
$\mathbbm{Shift}_k$ & Error in centering of discrete Gaussian \dotfill & \pageref{eq_x233} \\
\hline\hline
\end{tabular}
\end{flushleft}}

\vspace{0.5cm}

{\footnotesize \begin{flushleft}
\begin{tabular}{lll}
\multicolumn{3}{l}{\textsc{Vector spaces}} \\
\hline\hline
$\boldsymbol{e}^{(h)}$ & Canonical basis of $\amsmathbb{R}^H$ \dotfill & \pageref{eigenvalhhh}, \pageref{formulaehhh} \\
$\langle \cdot \rangle$ & Scalar product \dotfill & \pageref{eigenvalhhh} \\
$\mathscr{V}$ & Image of $\boldsymbol{\Theta}$ \dotfill & \pageref{lemrefl} \\
$\hat{\boldsymbol{e}}^{(h)}$ & Orthogonal projection of $\boldsymbol{e}^{(h)}$ onto $\mathscr{V}$ \dotfill & \pageref{lemrefl} \\
$\langle\!\langle \cdot \rangle\!\rangle$ & Extended scalar product \dotfill & \pageref{doubleba} \\
$\textnormal{\textbf{1}}_{\mathcal{H}_i}$ & Vectors with 0-1 entries & \pageref{ThetaG} \\
\hline
$|\!|\cdot|\!|_{\infty}$ & Sup norm \dotfill & \pageref{Assumptions_Theta} \\
$|\!|\cdot|\!|_{\textnormal{Lip}}$ & Lipschitz norm \dotfill & \pageref{deffdeminorm} \\
$|\!|\cdot|\!|_{\frac{1}{2}}$ & Sobolev norm \dotfill & \pageref{deffdeminorm} \\
$\mathscr{H}_{\text{Lip},\frac{1}{2}}$ & Lipschitz--Sobolev space \dotfill &\pageref{deffdeminorm} \\
$\text{Lip}_1$ & 1-bounded 1-Lipschitz functions \dotfill & \pageref{Definition_asymptotic_equivalence} \\
\hline
$\mathscr{O}$ & Space of holomorphic functions \dotfill & \pageref{Definition_complex_dom} \\
$\mathscr{H}^{[m]},\mathscr{H}^{[m],0}$ & Similar, with controlled growth \dotfill & \pageref{Definition_complex_dom} \\
\hline\hline
\end{tabular}
\end{flushleft}}

\vspace{0.5cm}

{\footnotesize \begin{flushleft}
\begin{tabular}{lll}
\multicolumn{3}{l}{\textsc{Functions}} \\
\hline\hline
$w_h$ & Weight \dotfill & \pageref{DataS} \\
$V_h, U_h$ & Potential and regular part \dotfill & \pageref{eq_weight_form} \\
\hline
$\mathcal{I}$ & Minus the energy functional \dotfill & \pageref{Definition_functional} \\
$\mathcal{I}^{(2)}$ & Quadratic part of $\mathcal{I}$ \dotfill & \pageref{eq_functional_quadratic} \\
$\mathcal{J}_{\varkappa}$ & Auxiliary functional on measures & \pageref{eq_J_alpha} \\
$\tilde{\mathcal{I}}$ & Minus the regularized energy functional \dotfill & \pageref{eq_x207} \\
$\mathfrak{D}^2$ & Squared pseudo-distance & \pageref{Definition_pseudodistance} \\
\hline\hline
\end{tabular}
\end{flushleft}}

\vspace{0.5cm}

{\footnotesize \begin{flushleft}
\begin{tabular}{lll}
\multicolumn{3}{l}{\textsc{Variables and parameters}} \\
\hline\hline
$\ell_i, \ell_i^h$ & Particles \dotfill & \pageref{Section_framework}, \pageref{eljah} \\
$\N$ & Master parameter \dotfill & \pageref{DataS} \\
\hline
$\theta_{g,h}, \boldsymbol{\Theta}$ & Interaction intensity \dotfill & \pageref{DataS} \\
$\boldsymbol{\Theta}'$ & Reduced interaction intensity \dotfill & \pageref{Assumptions_Theta} \\
$\boldsymbol{\Theta}^{\pm}$ & Bipartite decomposition of $\boldsymbol{\Theta}$ \dotfill & \pageref{Thetapmdec} \\
$\mathbf{R}$ & Adjacency segment-trapezoid \dotfill & \pageref{ThetaRTR}, \pageref{imlinform} \\
\hline
$a_h, b_h$ & Segment endpoints \dotfill & \pageref{DataS} \\
\hline\hline
\end{tabular}
\end{flushleft}}

\providecommand{\bysame}{\leavevmode\hbox to3em{\hrulefill}\thinspace}
\providecommand{\MR}{\relax\ifhmode\unskip\space\fi MR }
\providecommand{\MRhref}[2]{%
  \href{http://www.ams.org/mathscinet-getitem?mr=#1}{#2}
}
\providecommand{\href}[2]{#2}

\end{document}